\documentclass[12pt,a4paper]{ociamthesis}
\usepackage{graphicx}

\usepackage{amssymb,amsmath,stmaryrd,mathrsfs, mathpartir,color, cmll, tabulary, enumitem}
\usepackage[dvipsnames]{xcolor}

\usepackage{setspace}
\usepackage{changepage}

\usepackage{caption}
\usepackage{tikz}
\usetikzlibrary{positioning,calc,decorations.pathmorphing,shapes}
\usepackage[framemethod=TikZ]{mdframed}
\mdfdefinestyle{defcol}{%
  backgroundcolor=yellow!10,roundcorner=8pt,hidealllines=true, skipbelow = 5pt,
  innertopmargin=\the\abovedisplayskip,innerbottommargin=\the\belowdisplayskip}

\newenvironment{noverticalspace}
 {%
  \par %
  \offinterlineskip %
 }
{\par}

\DeclareMathAlphabet{\mathpzc}{OT1}{pzc}{m}{it}

\definecolor{brightgreen}{rgb}{0.4, 1.0, 0.0}
\definecolor{brightturquoise}{rgb}{0.03, 0.91, 0.87}
\definecolor{brightpink}{rgb}{1.0, 0.0, 0.5}
\definecolor{carrotorange}{rgb}{0.93, 0.57, 0.13}

\definecolor{CDcolor}{rgb}{0.0,0.5,0.75}	%
\definecolor{CLDcolor}{rgb}{0.8,0.0,0.2}		%

\def\definetac{\newif\iftac}    %
\ifx\tactrue\undefined
  \definetac
  \ifx\state\undefined\tacfalse\else\tactrue\fi
\fi

\def\definebeamer{\newif\ifbeamer}
\ifx\beamertrue\undefined
  \definebeamer
  \ifx\uncover\undefined\beamerfalse\else\beamertrue\fi
\fi

\def\definecref{\newif\ifcref}
\ifx\creftrue\undefined
  \definecref
  \creffalse
\fi

\iftac\else\usepackage{amsthm}\fi
\usepackage[all,2cell]{xy}
\ifbeamer\else
  \usepackage{enumitem}
  \usepackage{xcolor}
  \usepackage{hyperref}
  \hypersetup{colorlinks,breaklinks,
            urlcolor=[RGB]{7,7,96},
            linkcolor=[RGB]{7,7,96},
            citecolor=[RGB]{7,7,96}}
  \usepackage{bookmark}
\fi
\usepackage{mathtools}          %
\usepackage{graphics}           %
\usepackage{ifmtarg}            %
\usepackage{microtype}
\usepackage{braket}             %

\usepackage{url}                %
\usepackage{xspace}             %
\ifcref\usepackage{cleveref,aliascnt}\fi

\makeatletter
\let\ea\expandafter

\def\mdef#1#2{\ea\ea\ea\gdef\ea\ea\noexpand#1\ea{\ea\ensuremath\ea{#2}\xspace}}
\def\alwaysmath#1{\ea\ea\ea\global\ea\ea\ea\let\ea\ea\csname your@#1\endcsname\csname #1\endcsname
  \ea\def\csname #1\endcsname{\ensuremath{\csname your@#1\endcsname}\xspace}}

\DeclareRobustCommand\widecheck[1]{{\mathpalette\@widecheck{#1}}}
\def\@widecheck#1#2{%
    \setbox\z@\hbox{\m@th$#1#2$}%
    \setbox\tw@\hbox{\m@th$#1%
       \widehat{%
          \vrule\@width\z@\@height\ht\z@
          \vrule\@height\z@\@width\wd\z@}$}%
    \dp\tw@-\ht\z@
    \@tempdima\ht\z@ \advance\@tempdima2\ht\tw@ \divide\@tempdima\thr@@
    \setbox\tw@\hbox{%
       \raise\@tempdima\hbox{\scalebox{1}[-1]{\lower\@tempdima\box
\tw@}}}%
    {\ooalign{\box\tw@ \cr \box\z@}}}

\newcount\foreachcount

\def\foreachletter#1#2#3{\foreachcount=#1
  \ea\loop\ea\ea\ea#3\@alph\foreachcount
  \advance\foreachcount by 1
  \ifnum\foreachcount<#2\repeat}

\def\foreachLetter#1#2#3{\foreachcount=#1
  \ea\loop\ea\ea\ea#3\@Alph\foreachcount
  \advance\foreachcount by 1
  \ifnum\foreachcount<#2\repeat}
  
\let\oldit\it

\def\definescr#1{\ea\gdef\csname s#1\endcsname{\ensuremath{\mathscr{#1}}\xspace}}
\foreachLetter{1}{27}{\definescr}
\def\definecal#1{\ea\gdef\csname c#1\endcsname{\ensuremath{\mathcal{#1}}\xspace}}
\foreachLetter{1}{27}{\definecal}
\def\definebold#1{\ea\gdef\csname b#1\endcsname{\ensuremath{\mathbf{#1}}\xspace}}
\foreachLetter{1}{27}{\definebold}
\def\definebb#1{\ea\gdef\csname l#1\endcsname{\ensuremath{\mathbb{#1}}\xspace}}
\foreachLetter{1}{27}{\definebb}
\def\definefrak#1{\ea\gdef\csname k#1\endcsname{\ensuremath{\mathfrak{#1}}\xspace}}
\foreachletter{1}{27}{\definefrak}
\foreachLetter{1}{27}{\definefrak}
\def\definesf#1{\ea\gdef\csname i#1\endcsname{\ensuremath{\mathsf{#1}}\xspace}}
\foreachletter{1}{6}{\definesf}
\foreachletter{7}{14}{\definesf}
\foreachletter{15}{27}{\definesf}
\foreachLetter{1}{27}{\definesf}
\def\definebar#1{\ea\gdef\csname #1bar\endcsname{\ensuremath{\overline{#1}}\xspace}}
\foreachLetter{1}{27}{\definebar}
\foreachletter{1}{8}{\definebar} %
\foreachletter{9}{15}{\definebar} %
\foreachletter{16}{27}{\definebar}
\def\definetil#1{\ea\gdef\csname #1til\endcsname{\ensuremath{\widetilde{#1}}\xspace}}
\foreachLetter{1}{27}{\definetil}
\foreachletter{1}{27}{\definetil}
\def\definehat#1{\ea\gdef\csname #1hat\endcsname{\ensuremath{\widehat{#1}}\xspace}}
\foreachLetter{1}{27}{\definehat}
\foreachletter{1}{27}{\definehat}
\def\definechk#1{\ea\gdef\csname #1chk\endcsname{\ensuremath{\widecheck{#1}}\xspace}}
\foreachLetter{1}{27}{\definechk}
\foreachletter{1}{27}{\definechk}
\def\defineul#1{\ea\gdef\csname u#1\endcsname{\ensuremath{\underline{#1}}\xspace}}
\foreachLetter{1}{27}{\defineul}
\foreachletter{1}{27}{\defineul}

\let\it\oldit

\def\autofmt@b#1\autofmt@end{\mathbf{#1}}
\def\autofmt@l#1#2\autofmt@end{\mathbb{#1}\mathsf{#2}}
\def\autofmt@c#1#2\autofmt@end{\mathcal{#1}\mathit{#2}}
\def\autofmt@s#1#2\autofmt@end{\mathscr{#1}\mathit{#2}}
\def\autofmt@f#1\autofmt@end{\mathsf{#1}}
\def\autofmt@k#1\autofmt@end{\mathfrak{#1}}
\def\autofmt@u#1\autofmt@end{\underline{\smash{\mathsf{#1}}}}
\def\autofmt@U#1\autofmt@end{\underline{\underline{\smash{\mathsf{#1}}}}}
\def\autofmt@h#1\autofmt@end{\widehat{#1}}
\def\autofmt@r#1\autofmt@end{\overline{#1}}
\def\autofmt@t#1\autofmt@end{\widetilde{#1}}
\def\autofmt@k#1\autofmt@end{\check{#1}}

\def\auto@drop#1{}
\def\autodef#1{\ea\ea\ea\@autodef\ea\ea\ea#1\ea\auto@drop\string#1\autodef@end}
\def\@autodef#1#2#3\autodef@end{%
  \ea\def\ea#1\ea{\ea\ensuremath\ea{\csname autofmt@#2\endcsname#3\autofmt@end}\xspace}}
\def\autodefs@end{blarg!}
\def\autodefs#1{\@autodefs#1\autodefs@end}
\def\@autodefs#1{\ifx#1\autodefs@end%
  \def\autodefs@next{}%
  \else%
  \def\autodefs@next{\autodef#1\@autodefs}%
  \fi\autodefs@next}

\DeclareSymbolFont{bbold}{U}{bbold}{m}{n}
\DeclareSymbolFontAlphabet{\mathbbb}{bbold}

\renewcommand*{\sG}{\mathcal{G}}
\renewcommand*{\sJ}{\mathcal{J}}
\renewcommand*{\sR}{\mathcal{R}}
\renewcommand*{\sS}{\mathcal{S}}
\renewcommand*{\sT}{\mathcal{T}}
\renewcommand*{\sU}{\mathcal{U}}
\renewcommand*{\sV}{\mathcal{V}}
\renewcommand*{\cS}{\mathscr{S}}
\renewcommand*{\cF}{\mathscr{F}}
\renewcommand*{\cG}{\mathtt{GE}}

\definecolor{cblue}{RGB}{0,0,255}
\definecolor{cdarkblue}{RGB}{0, 38, 153}
\definecolor{clightblue}{RGB}{128, 170, 255}
\definecolor{cturquoise}{RGB}{0, 255, 255}
\definecolor{cgreen}{RGB}{0, 230, 0}
\definecolor{cdarkgreen}{RGB}{0, 128, 0}
\definecolor{corange}{RGB}{255, 153, 0}
\definecolor{cred}{RGB}{255, 51, 0}
\definecolor{cpurple}{RGB}{172, 0, 230} 
\definecolor{cyellow}{RGB}{230, 230, 0}
\definecolor{cpink}{RGB}{255, 51, 204} 
\newcommand*{\cpink}{\textcolor{cpink}{pink}}
\newcommand*{\cyellow}{\textcolor{cyellow}{yellow}}
\newcommand*{\cred}{\textcolor{cred}{red}}
\newcommand*{\cpurple}{\textcolor{cpurple}{purple}}
\newcommand*{\corange}{\textcolor{corange}{orange}}
\newcommand*{\cdarkgreen}{\textcolor{cdarkgreen}{green}}
\newcommand*{\cturquoise}{\textcolor{cturquoise}{turquoise}}
\newcommand*{\cgreen}{\textcolor{cgreen}{green}}
\newcommand*{\clightblue}{\textcolor{clightblue}{light blue}}
\newcommand*{\cdarkblue}{\textcolor{cdarkblue}{dark blue}}
\newcommand*{\cblue}{\textcolor{cblue}{blue}}
\newcommand*{\cgray}{\textcolor{gray}{gray}}
\newcommand*{\cblack}{\textbf{black}}

\newcommand{\bang}{\ensuremath{!}}

\newcommand{\vslash}{\mathbin{/\mkern-6mu/}}
\newcommand{\regop}{^{\mathrm{reg}}}
\newcommand{\singop}{^{\mathrm{sing}}}
\newcommand{\sing}{\ensuremath{\mathsf{sing}}}
\newcommand{\cone}{\mathrm{cone}}
\newcommand{\retr}{\mathrm{retr}}
\newcommand{\Fol}{^\mathrm{fol}}
\newcommand{\extsing}{\ensuremath{\widehat{\mathsf{sing}}}}
\newcommand{\reg}{\ensuremath{\mathsf{reg}}}
\newcommand{\SiR}{\ensuremath{\mathsf{R}}}
\newcommand{\edgeset}{\mathtt{E}}
\newcommand{\oedgeset}{_\mathtt{E}}
\renewcommand{\succ}[1]{\ensuremath{\mathtt{S}({#1})}}
\newcommand{\succfill}[1]{\ensuremath{\mathtt{F}({#1})}}
\newcommand{\succindex}[1]{\mathtt{I}({#1})}
\newcommand{\pred}[1]{\ensuremath{\mathtt{P}({#1})}}

\newcommand{\succn}[2]{\ensuremath{\mathtt{S}^{#2}({#1})}}
\newcommand{\predn}[2]{\ensuremath{\mathtt{P}^{#2}({#1})}}
\newcommand{\soe}{\mathtt{s}} %
\newcommand{\sop}{\mathtt{s}} %
\newcommand{\tae}{\mathtt{t}} %
\newcommand{\tap}{\mathtt{t}} %
\newcommand{\eda}{\mathtt{a}} %
\newcommand{\edb}{\mathtt{b}} %
\newcommand{\edc}{\mathtt{c}} %
\newcommand{\scA}{A} %
\newcommand{\scB}{B} %
\newcommand{\scC}{C} %
\newcommand{\scD}{D} %
\newcommand{\SIf}[1]{{\underline{#1}}}
\newcommand{\SIs}[1]{\ensuremath{\underrightarrow{#1}}}

\newcommand{\tgl}{\ensuremath{\mathsf{G}}}
\newcommand{\ppi}{\ensuremath{\mathsf{proj}}}
\newcommand{\const}{\ensuremath{\mathsf{const}}}
\newcommand{\msrc}{\ensuremath{\mathsf{s}}}
\newcommand{\mtgt}{\ensuremath{\mathsf{t}}}

\newcommand{\csrc}{\ensuremath{\mathsf{src}_\square}}
\newcommand{\ctgt}{\ensuremath{\mathsf{tgt}_\square}}
\newcommand{\gsrc}{\ensuremath{\mathsf{src}_\circ}}
\newcommand{\gtgt}{\ensuremath{\mathsf{tgt}_\circ}}
\newcommand{\stack}{\ensuremath{\triangleright}}
\newcommand{\glue}[1]{\ensuremath{\triangleright^{#1}}}
\newcommand{\whisker}[2]{\ensuremath{\triangleright^{#1}_{#2}}}

\newcommand{\igadd}[1]{+^\sim_{#1}}
\newcommand{\gadd}[1]{+_{#1}}
\newcommand{\icg}[1]{\ic^{#1}}
\newcommand{\dualdag}{^\dagger}
\newcommand{\ccirc}{\mathrm{c}^\circ}
\newcommand{\cell}{\mathrm{c}}
\newcommand{\regcont}{\mathsf{reg}}
\newcommand{\singcont}{\mathsf{sing}}
\newcommand{\quotg}{\sim}
\newcommand{\wwidehat}{\widehat}
\newcommand{\ttimes}{\times}
\newcommand{\truncleq}{\leq}
\renewcommand*{\sslash}{\slash}
\newcommand{\cpartial}{\mathrm{c}^\partial}
\newcommand{\secp}[1]{\ensuremath{\overline{#1}}}
\newcommand{\lGR}[1]{\ensuremath{\lG^{#1}}}
\newcommand{\topmon}{_\top}
\newcommand{\restsec}[1]{_{#1}}
\newcommand{\NF}[1]{[#1]}
\newcommand{\qmin}[1]{\iq^{#1}}
\newcommand{\intrel}[1]{\mathbf{I}[#1]}
\newcommand{\Poss}[1]{\Pos\slash_{\!{#1}}}
\newcommand{\kcoll}[1]{\coll^{#1}}
\newcommand{\coll}{\ensuremath{\epi}}
\newcommand{\mcoll}{\ensuremath{\epi}}
\newcommand{\kmcoll}[1]{\epi^{#1}}
\newcommand{\starcoll}{\coll^*}
\newcommand{\pbstar}{^*}
\newcommand{\postar}{_*}
\newcommand{\psstar}{_*}
\newcommand{\vvec}[1]{\vec{#1}}
\newcommand{\SIvertone}[1]{\ensuremath{\mathbf{SI}\vslash_{\! #1}}}
\newcommand{\SIvert}[2]{\ensuremath{\mathbf{SI}\vslash^{#1}_{\! #2}}}
\newcommand{\tsG}[1]{\sG^{#1}}
\newcommand{\upi}[2]{\ensuremath{\tsG {#2}({#1})}}
\newcommand{\tpi}[1]{\pi^{#1}}
\newcommand{\tsU}[1]{\ensuremath{\sU^{#1}}}
\newcommand{\tsR}[1]{\ensuremath{\sR^{#1}}}
\newcommand{\tusU}[1]{\ensuremath{\underline\sU^{#1}}}

\newcommand{\monosing}{\ensuremath{_{\mathsf{sing}}}} 
\newcommand{\restemb}{\ensuremath{\mathsf{rest}}} 
\newcommand{\nfc}{\ensuremath{\Lambda}} 
\newcommand{\vsS}[1]{\ensuremath{\vec\sS^{{#1}}}} 
\newcommand{\ssoe}{_\mathtt{s}} %
\newcommand{\ttae}{\mathtt{t}} %
\newcommand{\stfwd}{\textbf{straight-forward}} %
\newcommand{\stacklow}{\ensuremath{\mathsf{In}_1}}
\newcommand{\stackup}{\ensuremath{\mathsf{In}_2}}
\newcommand{\stackinc}[1]{\ensuremath{\mathsf{In}^{#1}}}

\newcommand{\Free}{Presented}
\newcommand{\free}{presented}
\newcommand{\rsDual}{Dual}

\newcommand{\rsdual}{dual}
\newcommand{\typed}{typed}

\newcommand{\obj}{\mathrm{obj}}
\newcommand{\PSk}{\mathrm{PSk}}
\newcommand{\elpsk}{\mathrm{SE}}
\newcommand{\discr}{\mathrm{dis}}
\newcommand{\Discr}{\mathrm{Dis}}
\newcommand{\converse}{\ensuremath{^\top}}

\newcommand{\grph}{\glssymbol{grph}}

\newcommand{\SI}{\ensuremath{\mathbf{SI}}}
\newcommand{\oSI}{_\mathbf{SI}}
\newcommand{\globset}{\ensuremath{\mathbf{\lG Set}}}
\newcommand{\PRel}{\ensuremath{\mathbf{PRel}}}
\newcommand{\xcoll}[2]{...}

\newcommand{\ata}{\ensuremath{\mathsf{f}}}
\newcommand{\bracket}{\ensuremath{\mathsf{brkt}}}
\newcommand{\bunit}{\ensuremath{\mathsf{unit}}}
\newcommand{\witness}[1]{\ensuremath{\mathsf{wit}(#1)}}
\newcommand{\tinv}{\ensuremath{\mathrm{inv}}}
\newcommand{\ttinv}{\ensuremath{\widetilde{\mathrm{inv}}}}

\newcommand{\Comp}{\ensuremath{\mathrm{Comp}}}
\newcommand{\GComp}{\ensuremath{\mathrm{GComp}}}
\newcommand{\GComps}{\ensuremath{\mathrm{GComp}^*}}
\newcommand{\extkM}{\ensuremath{\mathfrak{M}^\infty}}
\newcommand{\interchanger}{\ensuremath{\mathrm{IR}}}
\newcommand{\complaw}{associativity}

\newcommand{\cohlaw}{coherence}
\newcommand{\unitlaw}{unit}
\newcommand{\kiC}[1]{\ensuremath{\kR^{#1}}}

\newcommand{\infeq}[2]{\ensuremath{\mathsf{m}^{#1}_{#2}}}
\newcommand{\infinc}[1]{\ensuremath{\ii_{#1}}}
\newcommand{\ctyp}{\ensuremath{\mathrm{regtyp}}}
\newcommand{\ctypsum}{\ensuremath{\mathrm{codim}}}
\newcommand{\codim}{\ensuremath{\mathrm{codim}}}

\newcommand{\elcat}{\ensuremath{\mathrm{El}}}

\newcommand{\Spaces}{\ensuremath{\mathbf{Spaces}}}
\newcommand{\Bool}{\ensuremath{\mathbf{Bool}}}
\newcommand{\Rel}{\ensuremath{\mathbf{Rel}}}
\newcommand{\Pos}{\ensuremath{\mathbf{Pos}}}
\newcommand{\cMfld}{\ensuremath{\mathbf{cMfld}}}
\newcommand{\SetCat}{\ensuremath{\mathbf{Set}}}
\newcommand{\Cat}{\ensuremath{\mathbf{Cat}}}
\newcommand{\pCat}{\ensuremath{\mathbf{Cat}^{\mathrm{pres}}}}
\newcommand{\pwCat}{\ensuremath{\omega\mathbf{Cat}^{\mathrm{pres}}}}
\newcommand{\pgCat}{\ensuremath{\mathbf{Cat}^{\mathrm{pres,g}}}}

\newcommand{\BDelta}{\ensuremath{\boldsymbol{\Delta}}}
\newcommand{\stratatype}{region type}

\newcommand{\Bun}{\ensuremath{\mathbf{\lC ube}}}

\newcommand{\Cubeo}[2]{\ensuremath{\mathbf{Cube}^{\circ , #1}_{#2}}}
\newcommand{\Bunbc}{\ensuremath{\mathbf{\lC ube}}}
\newcommand{\Buno}[2]{\ensuremath{\mathbf{\lC ube}^{\circ , #1}_{#2}}}
\newcommand{\Buncoll}[1]{\ensuremath{\mathbf{\lC ube}^{\rightsquigarrow}_{#1}}}
\newcommand{\TI}{\ensuremath{\mathbf{TI}}}
\newcommand{\TTI}{\ensuremath{\redGamma{\mathbf{TI}}}}
\newcommand{\redGamma}[1]{\ensuremath{\widetilde{\mathsf{\Gamma}}(#1)}}
\newcommand{\GGamma}[2]{\ensuremath{\mathsf{\Gamma}^{#1}(#2)}}

\newcommand{\PProf}{\ensuremath{\mathbb{P}\mathbf{rof}}} 
\newcommand{\Prof}{\ensuremath{\mathbf{Prof}}} 
\newcommand{\bnum}[1]{\ensuremath{\mathbf{#1}}}

\newcommand{\und}[1]{{\underline{#1}}}
\newcommand{\singint}[1]{{\mathbf{I}_{#1}}}

\newcommand{\abs}[1]{\left|{#1}\right|} 
\newcommand{\abss}[1]{\left\llbracket{#1}\right\rrbracket} 
\newcommand{\norm}[1]{\ensuremath{\left\lVert{#1}\right\rVert}}

\newcommand{\eps}{\epsilon}

\newcommand{\avg}[1]{\left\langle{#1}\right\rangle} 
\newcommand{\rest}[2]{\ensuremath{\left.{#1}\right|_{#2}}}

\DeclareSymbolFont{extraup}{U}{zavm}{m}{n}
\DeclareMathSymbol{\varheart}{\mathalpha}{extraup}{86}
\DeclareMathSymbol{\vardiamond}{\mathalpha}{extraup}{87}

\mdef\delbar{\overline{\partial}}

\newcommand{\inv}{^{-1}}

\mdef\hf{\textstyle\frac12 }
\mdef\thrd{\textstyle\frac13 }
\mdef\qtr{\textstyle\frac14 }

\newcommand{\op}{^{\mathrm{op}}}

\SelectTips{cm}{}
\newdir{ >}{{}*!/-10pt/\dir{>}}    %

\newcommand{\pushoutfar}[1][dr]{\save*!/#1+1.7pc/#1:(1,-1)@^{|-}\restore}
\newcommand{\pullback}[1][dr]{\save*!/#1-1.2pc/#1:(-1,1)@^{|-}\restore}
\newcommand{\pullbacktest}[1][dr]{\save*!/va(-15) -2pc/va(135)@^{|-}\restore}
\newcommand{\pullbackfar}[1][dr]{\save*!/#1-1.7pc/#1:(-1,1)@^{|-}\restore}
\let\iso\cong

\mdef\Id{\mathrm{Id}}
\mdef\id{\mathrm{id}}
\mdef\uni{\mathrm{is\_uni}}
\mdef\sym{\mathrm{is\_sym}}
\mdef\rel{\mathrm{is\_rel}}
\mdef\Uni{\mathrm{Uni}}
\mdef\Sym{\mathrm{Sym}}
\mdef\SYM{\textsc{sym}}
\mdef\REL{\textsc{rel}}
\alwaysmath{ell}
\alwaysmath{infty}
\alwaysmath{odot}
\def\frc#1/#2.{\frac{#1}{#2}}   %
\mdef\ten{\mathrel{\otimes}}

\mdef\sqten{\mathrel{\boxtimes}}

\DeclareMathOperator\dom{dom}
\DeclareMathOperator\cod{cod}
\DeclareMathOperator\mor{mor}

\DeclareMathOperator\Fun{Fun}

\mdef\Im{\mathrm{Im}}
\mdef\im{\mathrm{im}}
\let\lim\relax
\DeclareMathOperator\lim{lim}

\DeclareMathOperator\colim{colim}

\DeclareMathOperator\coeq{coeq}

\DeclareMathOperator\Hom{Hom}
\DeclareMathOperator\Map{Func}

\newcommand{\ot}{\ensuremath{\leftarrow}}

\let\toot\rightleftarrows

\let\imp\Rightarrow
\let\limp\Leftarrow

\let\into\hookrightarrow
\let\xinto\xhookrightarrow
\mdef\we{\overset{\sim}{\longrightarrow}}
\mdef\leftwe{\overset{\sim}{\longleftarrow}}
\let\mono\rightarrowtail

\let\epi\twoheadrightarrow

\newcounter{sarrow}

\let\xto\xrightarrow

\def\rightarrowtailfill@{\arrowfill@{\Yright\joinrel\relbar}\relbar\rightarrow}
\newcommand\xrightarrowtail[2][]{\ext@arrow 0055{\rightarrowtailfill@}{#1}{#2}}

\def\twoheadrightarrowfill@{\arrowfill@{\relbar\joinrel\relbar}\relbar\twoheadrightarrow}
\newcommand\xtwoheadrightarrow[2][]{\ext@arrow 0055{\twoheadrightarrowfill@}{#1}{#2}}

\def\slashedarrowfill@#1#2#3#4#5{%
  $\m@th\thickmuskip0mu\medmuskip\thickmuskip\thinmuskip\thickmuskip
   \relax#5#1\mkern-7mu%
   \cleaders\hbox{$#5\mkern-2mu#2\mkern-2mu$}\hfill
   \mathclap{#3}\mathclap{#2}%
   \cleaders\hbox{$#5\mkern-2mu#2\mkern-2mu$}\hfill
   \mkern-7mu#4$%
}
\def\rightslashedarrowfill@{%
  \slashedarrowfill@\relbar\relbar\mapstochar\rightarrow}
\newcommand\xslashedrightarrow[2][]{%
  \ext@arrow 0055{\rightslashedarrowfill@}{#1}{#2}}
\mdef\hto{\xslashedrightarrow{}}
\mdef\htoo{\xslashedrightarrow{\quad}}

\def\xiso#1{\mathrel{\mathrlap{\smash{\xto[\smash{\raisebox{1.3mm}{$\scriptstyle\sim$}}]{#1}}}\hphantom{\xto{#1}}}}

\def\jd#1{\@jd#1\ej}
\def\@jd#1|-#2\ej{\@@jd#1,,\;\vdash\;\left(#2\right)}
\def\@@jd#1,{\@ifmtarg{#1}{\let\next=\relax}{\left(#1\right)\let\next=\@@@jd}\next}
\def\@@@jd#1,{\@ifmtarg{#1}{\let\next=\relax}{,\,\left(#1\right)\let\next=\@@@jd}\next}
\def\jdm#1{\@jdm#1\ej}
\def\@jdm#1|-#2\ej{\@@jd#1,,\\\vdash\;\left(#2\right)}

\long\def\my@drawfill#1#2;{%
\@skipfalse
\fill[#1,draw=none] #2;
\@skiptrue
\draw[#1,fill=none] #2;
}
\newif\if@skip
\newcommand{\skipit}[1]{\if@skip\else#1\fi}
\newcommand{\drawfill}[1][]{\my@drawfill{#1}}

\newif\ifhyperref
\@ifpackageloaded{hyperref}{\hyperreftrue}{\hyperreffalse}
  \ifhyperref
    \def\defthm#1#2#3{%
      \newtheorem{#1}{#2}[section]%
      \numberwithin{#1}{subsection}%
      \expandafter\def\csname #1autorefname\endcsname{#2}%
      \expandafter\let\csname c@#1\endcsname\c@thm}
\let\cref\autoref

\newtheorem{thm}{Theorem}[section]
\numberwithin{thm}{subsection}
\ifcref
  \crefname{thm}{Theorem}{Theorems}
\else
  
\fi
\defthm{cor}{Corollary}{Corollaries}
\defthm{prop}{Proposition}{Propositions}
\defthm{lem}{Lemma}{Lemmas}
\defthm{sch}{Scholium}{Scholia}
\defthm{assume}{Assumption}{Assumptions}
\defthm{claim}{Claim}{Claims}
\defthm{conj}{Conjecture}{Conjectures}
\defthm{hyp}{Hypothesis}{Hypotheses}
\theoremstyle{definition}
\defthm{indclaim}{Inductive Claim}{Inductive Hypotheses}
\defthm{defn}{Definition}{Definitions}
\defthm{conv}{Convention}{Conventions}
\defthm{constr}{Construction}{Constructions}
\defthm{eg}{Example}{Examples}
\defthm{egs}{Examples}{Examples}
\defthm{ex}{Exercise}{Exercises}
\defthm{ceg}{Counterexample}{Counterexamples}
\theoremstyle{remark}
\defthm{rmk}{Remark}{Remarks}
\defthm{notn}{Notation}{Notations}

\let\c@equation\c@thm
\numberwithin{equation}{subsection}

\iftac
  \let\your@endproof\endproof
  \def\my@endproof{\your@endproof}
  \def\endproof{\my@endproof\gdef\my@endproof{\your@endproof}}
  \def\qedhere{\tag*{\endproofbox}\gdef\my@endproof{\relax}}
\fi

\iftac
  \def\pr@@f[#1]{\subsubsection*{\sc #1.}}
\fi

\def\thmqedhere{\expandafter\csname\csname @currenvir\endcsname @qed\endcsname}

\ifbeamer\else

\fi

\ifbeamer\else
  \setitemize[1]{leftmargin=2em}
  \setenumerate[1]{leftmargin=*}
\fi

\ifcref\else
  \@ifpackageloaded{mathtools}{\mathtoolsset{showonlyrefs,showmanualtags}}{}
\fi

\providecommand\phantomsection{}

\makeatletter
\newenvironment{restoretext}%
    { \let\if@nobreak\iffalse
  \let\if@noskipsec\iffalse
  \let\par\@@par
  \let\-\@dischyph
  \let\'\@acci\let\`\@accii\let\=\@acciii
  \parindent\z@ \parskip\z@skip
  \everypar{}%
     \linewidth\hsize \@totalleftmargin\z@
     \leftskip\z@skip \rightskip\z@skip \@rightskip\z@skip
     \vspace{7pt}
     \begin{adjustwidth}{}{\leftmargin} \centering %
    }{\end{adjustwidth} \vspace{0pt}
     }
\makeatother

\alwaysmath{alpha}
\alwaysmath{beta}
\alwaysmath{gamma}
\alwaysmath{Gamma}
\alwaysmath{delta}
\alwaysmath{Delta}
\alwaysmath{epsilon}
\mdef\ep{\varepsilon}
\alwaysmath{zeta}
\alwaysmath{eta}
\alwaysmath{theta}
\alwaysmath{Theta}
\alwaysmath{iota}
\alwaysmath{kappa}
\alwaysmath{lambda}
\alwaysmath{Lambda}
\alwaysmath{mu}
\alwaysmath{nu}
\alwaysmath{xi}
\alwaysmath{pi}
\alwaysmath{rho}
\alwaysmath{sigma}
\alwaysmath{Sigma}
\alwaysmath{tau}
\alwaysmath{upsilon}
\alwaysmath{Upsilon}
\alwaysmath{phi}
\alwaysmath{Pi}
\alwaysmath{Phi}
\mdef\ph{\varphi}
\alwaysmath{chi}
\alwaysmath{psi}
\alwaysmath{Psi}
\alwaysmath{omega}
\alwaysmath{Omega}

\makeatother
\title{Associative $n$-categories}   %

\author{Christoph Dorn}             %
\college{St Catherine's College}  %

\degree{Doctor of Philosophy}     %
\degreedate{Michaelmas 2018}         %

\usepackage{lipsum}
\usepackage{booktabs}
\usepackage{calc}
\usepackage[nomain,toc,entrycounter]{glossaries}

\newglossary[tesy]{termsandsymbols}{tsi}{tsg}{Terms and symbols} 
\newglossary[var]{variables}{vai}{vag}{List of variables} 

\glsaddkey{vars}{\glsentrytext{\glslabel}}{\glsentryvars}{\GLsentryvars}{\glsvars}{\Glsvars}{\GLSvars}

\makeglossaries

\newlength\glstermwidth
\newlength\glssymbolwidth
\newlength\glsvariableswidth
\newlength\glstermdescwidth
\newlength\glsvardescwidth

\newglossarystyle{termdescsymb}{%
  \setlength{\glstermdescwidth}{\linewidth-\glstermwidth-\glssymbolwidth-6\tabcolsep}%
  \renewenvironment{theglossary}%
    {\begin{longtable}{p{\glstermwidth}p{\glstermdescwidth}p{\glssymbolwidth}}}%
    {\end{longtable}}%
  \renewcommand*{\glsgroupheading}[1]{}%
}

\newglossarystyle{vardesc}{%
  \setlength{\glsvardescwidth}{\linewidth-\glsvariableswidth-4\tabcolsep}%
    {\begin{longtable}{p{\glsvariableswidth}p{\glsvardescwidth}}}%
    {\end{longtable}}%
  \renewcommand*{\glsgroupheading}[1]{}%
}

\makeatletter
 \appto\@newglossaryentryposthook{%
   \settowidth{\dimen@}{\glsentrysymbol{\@glo@label}}%
   \ifdim\dimen@>\glssymbolwidth
     \setlength{\glssymbolwidth}{\dimen@}%
     \setlength{\glstermdescwidth}{\linewidth-\glstermwidth-\glssymbolwidth}
   \fi
   \settowidth{\dimen@}{\glsentryvars{\@glo@label}}%
 }%
\makeatother

\newglossaryentry{TI_gen}
{
name= {invertibility coherence},
description={A generator of the theory of $\infty$-invertibility $\TI$. See \autoref{sec:group_TI} and \autoref{ssec:com_cob}.},
symbol={\ensuremath{\ic_{S\equiv T}}},
type={termsandsymbols}
}

\newglossaryentry{TI2}
{
name= {theory of $\infty$-dualisability},
description={Usually called theory of $\infty$-invertibility. See \autoref{sec:group_TI} and \autoref{ssec:com_cob}.},
symbol={\ensuremath{\TI}},
type={termsandsymbols}
}

\newglossaryentry{TI}
{
name= {theory of $\infty$-invertibility},
description={This is a presented associative $n$-category $\TI$. Alternatively also called theory of $\infty$-dualisability (as invertibility and dualisability coincide if no $n$-truncation is performed). See \autoref{sec:group_TI} and \autoref{ssec:com_cob}.},
symbol={\ensuremath{\TI}},
type={termsandsymbols}
}

\newglossaryentry{adjoin_inv_gen_set}
{
name= {adjoining a set of invertible generators},
description={See \autoref{constr:adjoining_sets_of_inv_generators}. },
symbol={\ensuremath{ \igadd{x_\ig,y_\ig} \{\ig \in I\}}},
type={termsandsymbols}
}

\newglossaryentry{adjoin_inv_gen_coh_data2}
{
name= {coherence data of invertible generator},
description={Usually called ``singularities associated to invertible generator $\ig$". See \autoref{constr:adjoin_inv_gen}.},
symbol={\ensuremath{\icg\ig_{S \equiv T}}},
type={termsandsymbols}
}

\newglossaryentry{adjoin_inv_gen_coh_data}
{
name= {singularities associated to invertible generator},
description={Also called ``coherence data" of an invertible generator $\ig$. See \autoref{constr:adjoin_inv_gen}.},
symbol={\ensuremath{\icg\ig_{S \equiv T}}},
type={termsandsymbols}
}

\newglossaryentry{adjoin_inv_gen}
{
name= {adjoining an invertible generator},
description={This is similar to ``adjoining a generator $\ig$", but instead of just adding a single generator we also need to add ``singularities" associated to it, as for instance its inverse, witnesses of invertibility etc. These singularities $\icg\ig_{S \equiv T}$ correspond to generators of the theory of $\infty$-invertibility $\TI$. See \autoref{constr:adjoin_inv_gen}.},
symbol={\ensuremath{\igadd{x,y} ~\ig}},
type={termsandsymbols}
}

\newglossaryentry{pres_gpd}
{
name= {presented associative $n$-groupoid},
description={A presented associative $n$-category in which every generator is invertible. See \autoref{constr:pres_gpd}.},
symbol={\ensuremath{}},
type={termsandsymbols}
}

\newglossaryentry{ata_map}
{
name= {attaching map},
description={See \autoref{defn:CW_complex}},
symbol={\ensuremath{\ata}},
type={termsandsymbols}
}

\newglossaryentry{cell_int}
{
name= {cell interior},
description={See \autoref{defn:CW_complex}},
symbol={\ensuremath{\ccirc}},
type={termsandsymbols}
}

\newglossaryentry{cell_bound}
{
name= {cell boundary},
description={See \autoref{defn:CW_complex}},
symbol={\ensuremath{\cpartial}},
type={termsandsymbols}
}

\newglossaryentry{glob_fol}
{
name= {globular foliation},
description={See \autoref{ssec:sum_CW}.},
symbol={\ensuremath{}},
type={termsandsymbols}
}

\newglossaryentry{framed_strat}
{
name= {framed $n$-stratifications},
description={The set of stratifications of a manifold $M$, such that strata have framing (and are part of a namescope $\sN$), up to a notion of cobordism. See \autoref{defn:framed_strat}.},
symbol={\ensuremath{\Psi^{\mathrm{fr}}_{\sN}({\color{lightgray}M})}},
type={termsandsymbols}
}

\newglossaryentry{cone}
{
name= {cone of stratification},
description={Given a stratification $\bM$ of a space $X$, its cone $\cone(\cM)$ is the stratification of $\cone(X)$ in which all strata are extended up to (but excluding) the vertex point, which becomes its own $0$-dimensional stratum. See \autoref{constr:dual_complex}.},
symbol={\ensuremath{\cone({\color{lightgray}\bM})}},
type={termsandsymbols}
}

\newglossaryentry{gen_TP}
{
name= {generalised Thom-Pontryagin construction},
description={A correspondence between homotopy classes $f$ of maps from a manifold $M$ into a CW-complex $X$, to framed $X$-manifolds $\kP_M(f)$ in $M$. See \autoref{constr:dual_complex}.},
symbol={\ensuremath{\kP_M}},
type={termsandsymbols}
}

\newglossaryentry{dual_strat}
{
name= {dual stratification},
description={The dual stratification $\kD(X)$ of a CW-complex $X$ is a stratification of $X$, with strata corresponding to cells, and intersecting them ``transversally". See \autoref{constr:dual_complex}.},
symbol={\ensuremath{\kD({\color{lightgray}X})}},
type={termsandsymbols}
}

\newglossaryentry{dual_stratum}
{
name= {dual stratum},
description={A dual stratum $g^\dagger$ is a stratum in the dual stratification $\kD(X)$. See \autoref{constr:dual_complex}.},
symbol={\ensuremath{{\color{lightgray}g}^\dagger}},
type={termsandsymbols}
}

\newglossaryentry{ass_cat_multiplication}
{
name= {composition operation of associative $n$-category},
description={See \autoref{defn:ANC}.},
symbol={\ensuremath{\kM_\sC}},
type={termsandsymbols}
}

\newglossaryentry{assoc_n_cat}
{
name= {associative $n$-category},
description={A globular set with composition operation on its compositional shapes. See \autoref{defn:ANC}.},
symbol={\ensuremath{}},
type={termsandsymbols}
}

\newglossaryentry{namescope}
{
name= {namescope},
description={A dimension-ordered sequence of sets of ``names" (or ``labels"). See \autoref{ssec:sum_namescopes}},
symbol={\ensuremath{}},
type={termsandsymbols}
}

\newglossaryentry{realisation}
{
name= {realisation of syntax as singular cube},
description={A syntactic expression $t$ can have realisations as a singular cube $\abss{t}$. This notation is used for \hyperref[glsentry-types]{types} ($\abss{g}$) and \hyperref[glsentry-double_cone]{double cones} ($\abss{S \xto g T}$).},
symbol={\ensuremath{\abss{-}}},
type={termsandsymbols}
}

\newglossaryentry{types}
{
name= {type of generator},
description={Each generator $g$ of a higher structure has a type $\abss{g}$, describing its interaction with other generators. See for instance \autoref{defn:pres_ANC}.},
symbol={\ensuremath{\abss{\color{lightgray}g}}},
type={termsandsymbols}
}

\newglossaryentry{summit_type}
{
name= {vertex point of type},
description={The central $0$-dimensional region of a type $\abss{g}$ for a label $g \in \sC_k$. See \autoref{defn:pres_ANC}.},
symbol={\ensuremath{\ip_g}},
type={termsandsymbols}
}

\newglossaryentry{pres_ass_cat}
{
name= {presented associative $n$-category},
description={A notion of semi-strict higher categories given by generators and relations. See \autoref{defn:pres_ANC}.},
symbol={\ensuremath{}},
type={termsandsymbols}
}

\newglossaryentry{map_pres}
{
name= {map of presentations},
description={A mapping of generators between presented associative $n$-categories, which preserves types. See \autoref{defn:map_of_presentations}.},
symbol={\ensuremath{}},
type={termsandsymbols}
}

\newglossaryentry{colim_pres}
{
name= {colimit of presentations},
description={One can take colimits of inclusions of presentations for presented associative $n$-categories $\sC^i$. See \autoref{rmk:colimit_of_presentations}.},
symbol={\ensuremath{\colim_i(\sC^i)}},
type={termsandsymbols}
}

\newglossaryentry{mor_pres_cat}
{
name= {composites of presented associative $n$-category},
description={Given a presented associative $n$-categories, its ($k$-)composites are globular, normalised and well-typed $k$-cubes. See \autoref{defn:PANC_mor}.},
symbol={\ensuremath{\Comp(\sC)_k}},
type={termsandsymbols}
}

\newglossaryentry{quot_glob_set}
{
name= {quotient globular set},
description={See \autoref{constr:globular_quotients}.},
symbol={\ensuremath{S^{\sim}_{\leq n}}},
type={termsandsymbols}
}

\newglossaryentry{truncleq}
{
name= {truncated globular set},
description={See \autoref{notn:trunc_glob_sets}.},
symbol={\ensuremath{S_{\truncleq n}}},
type={termsandsymbols}
}

\newglossaryentry{adjoin_eq}
{
name= {adjoining an equality},
description={See \autoref{constr:adjoin_eq}.},
symbol={\ensuremath{\gadd{x,y}\ie_{x = y}}},
type={termsandsymbols}
}

\newglossaryentry{adjoin_gen_set}
{
name= {adjoining a set of generators},
description={See \autoref{constr:adjoining_sets_of_gen}. },
symbol={\ensuremath{ \gadd{x_\ig,y_\ig} \{\ig \in I\}}},
type={termsandsymbols}
}

\newglossaryentry{adjoin_gen}
{
name= {adjoining a generator},
description={Given a presented associative $n$-category $\sC$, and a fresh  name $\ig$ together with candidates $x,y$ for $(k-1)$-dimensional sources and targets of $\ig$, then we can form $\sC \gadd{x,y} \ig$ which is the presented associative $n$-category obtained from $\sC$ by adjoining a new generating $k$-morphism $\ig$ with source $x$ and target $y$. See \autoref{constr:adjoining_gen}.},
symbol={\ensuremath{\gadd{x,y} ~\ig}},
type={termsandsymbols}
}

\newglossaryentry{manifold_diag}
{
name= {manifold diagrams},
description={Manifold diagrams are a geometric model for morphisms in presented associative $n$-categories. A definition is given in \autoref{ssec:po_mfld_diag}.},
symbol={\ensuremath{}},
type={termsandsymbols}
}

\newglossaryentry{deform}
{
name= {equivalence of flag-foliation-compatible stratifications},
description={See \autoref{ssec:po_mfld_diag}.},
symbol={\ensuremath{\simeq}},
type={termsandsymbols}
}

\newglossaryentry{ordered_sum}
{
name= {ordered sum functor},
description={A two-variable functor on singular intervals. See \autoref{constr:SI_family_stacking}.},
symbol={\ensuremath{\uplus}},
type={termsandsymbols}
}

\newglossaryentry{cubical_src}
{
name= {cubical source},
description={See \autoref{constr:source_target}.},
symbol={\ensuremath{\csrc^k}},
type={termsandsymbols}
}

\newglossaryentry{cubical_tgt}
{
name= {cubical target},
description={See \autoref{constr:source_target}.},
symbol={\ensuremath{\ctgt^k}},
type={termsandsymbols}
}

\newglossaryentry{stack}
{
name= {stacking of interval families},
description={A composition operation for interval families which morally stacks them on top of one another. See \autoref{defn:stacking_labelled_families}},
symbol={\ensuremath{\stack}},
type={termsandsymbols}
}

\newglossaryentry{stack_emb}
{
name= {stacking embedding},
description={For $k$-level stacking of $\SIvert n \cC$-families see  \autoref{rmk:can_subbund_of_gluing}. For stacking of $\SI$-families see  \autoref{constr:stacking_gluing_iso}. Note that in both cases the subscript $i \in \Set{1,2}$ refers to the first or second component of a specified ($k$-level) stacking composition, and the reference to this composition is implicit in the notation.},
symbol={\ensuremath{\stackinc k_i}},
type={termsandsymbols}
} 

\newglossaryentry{klvl_stack}
{
name= {stacking of cube families},
description={A composition operation for $n$-cube families which morally stacks them on top of one another in direction $k$. See \autoref{constr:k_level_stacking}.},
symbol={\ensuremath{\glue k}},
type={termsandsymbols}
}

\newglossaryentry{whiskering}
{
name= {whiskering},
description={A composition operation for normalised globular cube modelling the usual whiskering operation in strict higher categories. Defined in \autoref{constr:glob_comp}.},
symbol={\ensuremath{\whisker k n}},
type={termsandsymbols}
}

\newglossaryentry{homotopy}
{
name= {homotopy},
description={In general this is a well-typed $n$-cube without $n$-dimensional singularities. A definition is given in \autoref{defn:homotopies}.},
symbol={\ensuremath{\abss{T_1 \xiso \scD T_2}}},
type={termsandsymbols}
}

\newglossaryentry{well_typedness}
{
name= {typability},
description={A manifold diagram is typable if its labelled region have globally constant local neighbourhoods, see \autoref{ssec:po_mfld_diag}. In the context of \textit{presented} associative $n$-categories see \autoref{defn:PANC_mor}.},
symbol={\ensuremath{\empty}},
type={termsandsymbols}
}

\newglossaryentry{proj_tow_iso}
{
name= {trivial product tower isomorphism},
description={Inductively defined isomorphisms between ``trivial product bundle and bundle product". See \autoref{constr:n_lvl_coll_lim_from_coll}.},
symbol={\ensuremath{\sW^{T,k}_Z}},
type={termsandsymbols}
}

\newglossaryentry{proj_tow}
{
name= {trivial product towers},
description={We can take a trivial product of a tower of $\SI$-families $T$ with a poset $Z$ and obtain a trivial product tower $T\times Z$. In this case $T$ is called the projected tower of the product tower $T\times Z$. See \autoref{constr:n_lvl_coll_lim_from_coll}.},
symbol={\ensuremath{T \times Z}},
type={termsandsymbols}
}

\newglossaryentry{seq_coll_lim}
{
name= {sequenced collapse limit from ordered collapse sequence},
description={A chain of morphisms of length $n$ in $\SIvert n \cC$ obtainable from an ordered collapse sequences of labelled $n$-cubes. See \autoref{constr:n_lvl_coll_lim_from_coll}.},
symbol={\ensuremath{\sL}},
type={termsandsymbols}
}

\newglossaryentry{loc_triv_fam}
{
name= {locally trivial family},
description={Family which normalises to the identity on every morphism of it basespace. See \autoref{defn:loc_triv_fam}.},
symbol={\ensuremath{\empty}},
type={termsandsymbols}
}

\newglossaryentry{sing_cont}
{
name= {singular content},
description={This refers to the singular heights in the total space of an $\SI$-bundle. See \autoref{defn:singular_content}. More specifically it can also refer to \hyperref[glsentry-sing]{singular heights} of a single interval.},
symbol={\ensuremath{\singcont(-)}},
type={termsandsymbols}
}

\newglossaryentry{reg_cont}
{
name= {regular content},
description={This refers to the regular segments in the total space of a $\SI$-bundle. See \autoref{defn:singular_content}. More specifically it can also refer to \hyperref[glsentry-reg]{regular segments} of a single interval.},
symbol={\ensuremath{\regcont(-)}},
type={termsandsymbols}
}

\newglossaryentry{globularity}
{
name= {globular $n$-cube families},
description={Globularity is a condition on an $n$-cube family guaranteeing ``constancy" on sides of the cube (and similarly for all sub-cubes). See \autoref{defn:globular_families}.},
symbol={\ensuremath{\empty}},
type={termsandsymbols}
}

\newglossaryentry{globe}
{
name= {$n$-globe},
description={Usually a reference to a globular $n$-cube. See \autoref{defn:globular_families}. Sometimes reference to classical $n$-globes. See \autoref{notn:globeset}.},
symbol={\ensuremath{\empty}},
type={termsandsymbols}
}

\newglossaryentry{glob_src}
{
name= {globular source},
description={The source of a globular $n$-cube, that is, the first value of its $1$-level labelling. See \autoref{defn:globular_src_tgt}. For the iterated version see \autoref{defn:klvl_source_target}.},
symbol={\ensuremath{\gsrc}},
type={termsandsymbols}
}

\newglossaryentry{glob_tgt}
{
name= {globular target},
description={The target of a globular $n$-cube, that is, the last value of its $1$-level labelling. See \autoref{defn:globular_src_tgt}. For the iterated version see \autoref{defn:klvl_source_target}},
symbol={\ensuremath{\gtgt}},
type={termsandsymbols}
}

\newglossaryentry{elcat_of_globe}
{
name= {category of elements of representable globular set},
description={See \autoref{constr:elcat}},
symbol={\ensuremath{\cG^n}},
type={termsandsymbols}
}

\newglossaryentry{cat_of_el}
{
name= {category of elements},
description={See \autoref{constr:elcat}},
symbol={\ensuremath{\elcat}},
type={termsandsymbols}
}

\newglossaryentry{globe_cat_rep}
{
name= {globe category representable},
description={See \autoref{notn:globeset}.},
symbol={\ensuremath{\lG^n}},
type={termsandsymbols}
}

\newglossaryentry{globe_cat}
{
name= {globe category},
description={See \autoref{notn:globeset}.},
symbol={\ensuremath{\lG}},
type={termsandsymbols}
}
\newglossaryentry{globeset_cat}
{
name= {globular set},
description={See \autoref{notn:globeset}.},
symbol={\ensuremath{\globset}},
type={termsandsymbols}
}

\newglossaryentry{term_globe}
{
name= {terminal $n$-globe},
description={A globular $n$-cube looking like a classical $n$-globe in a globular set. See \autoref{constr:terminal_n_globe}.},
symbol={\ensuremath{\tgl^n}},
type={termsandsymbols}
}

\newglossaryentry{top_monad}
{
name= {$\top$-monad},
description={The $\top$-monad, pronounced ``top monad", adjoins a terminal object to a category. This generalises to singular $n$-cubes. See \autoref{constr:top_monad}.},
symbol={\ensuremath{(-)_\top}},
type={termsandsymbols}
}

\newglossaryentry{double_cone}
{
name= {double cone of sources and targets},
description={Let $S, T$ be normalised globular $\SIvert {n-1} \cC$-cubes whose globular sources and targets coincide, then their double cone $\abss{S \to T}$ (or alternatively, their double cone  $\abss{S \xto g T}$ with vertex point $g$) is the $n$-cube with source and target $S$, $T$ and a single singularity (labelled by $g$) in its centre. See \autoref{constr:double_cones_of_src_and_tgt}.},
symbol={\ensuremath{\abss{S \xto g T}}},
type={termsandsymbols}
}

\newglossaryentry{embedding}
{
name= {embedding},
description={A (multi-level) embedding $\theta: \scA \mono \scB$ identifies $\scA$ as a subfamily of the cube family $\scB$. See \autoref{defn:embedding_fctr_of_int}.},
symbol={\ensuremath{\mono}},
type={termsandsymbols}
}

\newglossaryentry{multilevel_embedding}
{
name= {multi-level embedding},
description={See \gls{embedding}.},
symbol={\ensuremath{\mono}},
type={termsandsymbols}
}

\newglossaryentry{embedding_fctr}
{
name= {embedding functor},
description={An injective open functor of singular intervals. See \autoref{defn:subfamilies}.},
symbol={\ensuremath{\empty}},
type={termsandsymbols}
}

\newglossaryentry{rest_emb}
{
name= {restriction embedding},
description={An embedding obtained by a restriction of the basespace to a subposet. See \autoref{eg:subfamily_by_restriction}.},
symbol={\ensuremath{\restemb}},
type={termsandsymbols}
}

\newglossaryentry{emb_fact}
{
name= {factorisation of embeddings},
description={See \autoref{claim:factorisation_subfamilies}.},
symbol={\ensuremath{\theta\inv\chi}},
type={termsandsymbols}
}

\newglossaryentry{open_sec}
{
name= {open sections},
description={See \autoref{defn:sections}.},
symbol={\ensuremath{\empty}},
type={termsandsymbols}
}

\newglossaryentry{src_sec}
{
name= {source section},
description={Special open section defined in \autoref{constr:source_and_target_inclusion}.},
symbol={\ensuremath{\msrc_A}},
type={termsandsymbols}
}

\newglossaryentry{tgt_sec}
{
name= {target section},
description={Special open section defined in \autoref{constr:source_and_target_inclusion}.},
symbol={\ensuremath{\mtgt_A}},
type={termsandsymbols}
}

\newglossaryentry{endpoints}
{
name= {endpoints},
description={A pair of open sections. See \autoref{constr:endpoint_inclusions}. Several structures can be determined from other structure by specifying endpoints. For instance, we defined \hyperref[glsentry-subfam_from_end]{subfamilies from endpoints} (${{\color{lightgray}\scA}\restsec{[q_-,q_+]}}$), \hyperref[glsentry-emb_from_end]{family embedding functors from endpoints} (${\sJ^\scA\restsec{[q_-,q_+]}}$), \hyperref[glsentry-subfam_from_k_end]{subfamilies determined by $k$-level endpoints} (${\color{lightgray}\scA}^k\restsec{[q_-,q_+]}$) and \hyperref[glsentry-emb_from_k_end]{embeddings from $k$-level endpoints} ($\sJ^{\scA,k}\restsec{[q_-,q_+]}$).},
symbol={\ensuremath{[q_-,q_+]}},
type={termsandsymbols}
}

\newglossaryentry{subfam_from_end}
{
name= {subfamily from endpoints},
description={Generalised for cubes in \autoref{constr:subfamilies_from_endpoints}. Originally defined in \autoref{constr:endpoint_inclusions}.},
symbol={\ensuremath{{\color{lightgray}\scA}\restsec{[q_-,q_+]}}},
type={termsandsymbols}
}

\newglossaryentry{emb_from_end}
{
name= {family embedding functor},
description={Defined in \autoref{constr:subfamilies_from_endpoints}. Generalised in \autoref{constr:subfamilies_from_endpoints}. },
symbol={\ensuremath{\sJ^\scA\restsec{[q_-,q_+]}}},
type={termsandsymbols}
}

\newglossaryentry{subfam_from_k_end}
{
name= {subfamily determined by $k$-level endpoints},
description={Defined in \autoref{constr:subfamilies_from_endpoints}.},
symbol={\ensuremath{{\color{lightgray}\scA}^k\restsec{[q_-,q_+]}}},
type={termsandsymbols}
}

\newglossaryentry{emb_from_k_end}
{
name= {embedding from $k$-level endpoints},
description={Defined in \autoref{constr:subfamilies_from_endpoints}.},
symbol={\ensuremath{\sJ^{\scA,k}\restsec{[q_-,q_+]}}},
type={termsandsymbols}
}

\newglossaryentry{klvl_emb}
{
name= {$k$-level embedding},
description={See \gls{emb_from_k_end}.},
symbol={\ensuremath{\empty}},
type={termsandsymbols}
} 

\newglossaryentry{min_endpoints}
{
name= {minimal endpoints},
description={Open sections used in the construction of minimal embeddings. See \autoref{constr:minimal_endpoints}.},
symbol={\ensuremath{q^p_\pm}},
type={termsandsymbols}
}

\newglossaryentry{min_subfam}
{
name= {minimal subfamily},
description={Minimal ``open" neighbourhood of a given point in a singular $n$-cube. See \autoref{constr:minimal_subfamilies}.},
symbol={\ensuremath{ {\color{lightgray}\scA} \slash {\color{lightgray}p} }},
type={termsandsymbols}
}

\newglossaryentry{slash_notn}
{
name= {slash notation},
description={The slash notation either denotes over-categories (see \autoref{notn:overcats}), or it denotes \hyperref[glsentry-min_subfam]{minimal subfamilies}.},
symbol={\ensuremath{ {\color{lightgray}X} \slash {\color{lightgray}x} }},
type={termsandsymbols}
}

\newglossaryentry{min_emb}
{
name= {minimal embedding},
description={Given a cube family $\scA$, a minimal embedding is an embedding $\iota^p_\scA : (\scA \slash p) \scA$ whose subfamily is the minimal (non-trivial) subfamily containing a given point $p\in \sG^n(\scA)$. See \autoref{constr:minimal_subfamilies}.},
symbol={\ensuremath{\color{lightgray}{\color{black}\iota}^p_\scA}},
type={termsandsymbols}
}

\newglossaryentry{coll_emb}
{
name= {collapsed embedding},
description={Given a embedding $\theta : \scA \mono \scB$ and a collapse (sequence) $\vec \lambda : \scB \to \scC$, then there is a induced collapsed $\vec\lambda\postar\theta : \scD \mono \scC$ (for a uniquely determined $\scD$). See \autoref{constr:collapse_on_subfamilies} and (for sequences) \autoref{constr:restricting_collapse_seq}.},
symbol={\ensuremath{{\color{lightgray}\vec\lambda}\postar {\color{lightgray}\theta}}},
type={termsandsymbols}
}

\newglossaryentry{pullback}
{
name= {pullback notation},
description={The pullback notation $A\pbstar B$ can denote pullback of collapse along collapse ($\mu\pbstar\lambda$, cf. \autoref{lem:coll_pull}), pullback of collapse along embedding ($\theta\pbstar \vec\lambda$, cf. \hyperref[glsentry-restrict_coll]{restriction of collapse}), pullback of functors along inclusions ($H\pbstar F$, cf. \hyperref[glsentry-codrest]{restriction to codomain}), pullback of $n$-cubes along functors factoring their labelling ($F\pbstar \scA$, cf. \autoref{notn:restrictions_of_labels}), or the pullback of a dual stratum into a cell ($h^*g\dualdag$, cf. \autoref{constr:dual_complex}).},
symbol={\ensuremath{\empty^*}},
type={termsandsymbols}
}

\newglossaryentry{restrict_coll}
{
name= {restriction of collapse (sequences) to embedding},
description={Given a embedding $\theta : \scA \mono \scB$ and a collapse (sequence) $\vec \lambda : \scB \to \scC$, then there is a restriction of the collapse (sequence) to $\theta$, denoted by $\theta^*\vec\lambda : \scA \to \scD$ (for a uniquely determined $\scD$). See \autoref{constr:collapse_on_subfamilies} and (for sequences) \autoref{constr:restricting_collapse_seq}.},
symbol={\ensuremath{{\color{lightgray}\theta}^*{\color{lightgray}\vec\lambda}}},
type={termsandsymbols}
}

\newglossaryentry{k_lvl_collapse}
{
name= {$k$-level collapse},
description={A $k$-level collapse $\lambda : \scA \coll^k \scB$ of $n$-cube families $\scA, \scB : X \to \SIvert n \cC$ is a collapse of interval families $\lambda : \sU^{k-1}_\scA \coll \sU^{k-1}_\scB$. See \autoref{defn:singular_cube_families}.},
symbol={\ensuremath{\kcoll k}},
type={termsandsymbols}
}

\newglossaryentry{normal_form2}
{
name= {normalised cube family},
description={A cube family in normal form. See \autoref{defn:collapse_normal_form}.},
symbol={\ensuremath{\empty}},
type={termsandsymbols}
}

\newglossaryentry{normal_form_notn}
{
name= {normal form notation},
description={The notation $\NF{-}$ is used for \hyperref[glsentry-k_lvl_normal_form]{$k$-level normal forms of cube families}, \hyperref[glsentry-normal_form]{normal forms of cube families} and \hyperref[glsentry-normal_form_SI]{normal forms of of interval families}.},
symbol={\ensuremath{\NF{-}}},
type={termsandsymbols}
}

\newglossaryentry{k_lvl_normal_form}
{
name= {$k$-level normal form (up to level $n$)},
description={See \autoref{thm:normal_forms_unique}.},
symbol={\ensuremath{\NF{-}^n_k}},
type={termsandsymbols}
}

\newglossaryentry{normal_form}
{
name= {normal form (up to level $n$) of cube family},
description={The unique maximal collapse of an $n$-cube family. See \autoref{defn:collapse_normal_form}.},
symbol={\ensuremath{\NF{-}^n}},
type={termsandsymbols}
}

\newglossaryentry{normal_form_SI}
{
name= {normal form of interval family},
description={The unique maximal collapse of a interval families. See \autoref{defn:collapse_normal_form}.},
symbol={\ensuremath{\NF{-}}},
type={termsandsymbols}
}

\newglossaryentry{vector}
{
name= {vector notation},
description={We use the vector notation to either denote \hyperref[glsentry-ordered_collapse_seq]{ordered collapse sequences}, \hyperref[glsentry-ass_multi_lvl_collapse]{associated multi-level collapse} or \hyperref[glsentry-multi_lvl_collapse]{multi-level collapse} in general.},
symbol={\ensuremath{\vec{\empty}}},
type={termsandsymbols}
}

\newglossaryentry{ordered_collapse_seq}
{
name= {ordered collapse sequence},
description={A sequence of collapses ordered by level of the cube projection. See \autoref{defn:ordered_coll_seq}.},
symbol={\ensuremath{\vec {\color{lightgray}\lambda}}},
type={termsandsymbols}
}

\newglossaryentry{normal_form_collapse}
{
name= {normal form collapse},
description={The unique ordered collapse sequence $\vec \lambda_\scC$ leading to normal form of a $\SIvert n \cC$-family $\scC$. See \autoref{thm:normal_forms_unique}.},
symbol={\ensuremath{\vec {\nfc}_{\color{lightgray}\scC}}},
type={termsandsymbols}
}

\newglossaryentry{multi_lvl_collapse}
{
name= {multi-level collapse},
description={A collapse acting at all projection levels of the cube simultaneously. See \autoref{defn:multilevel_collapse}. Note that multi-level collapse is an instance of multi-level base change.},
symbol={\ensuremath{\vec {\color{lightgray}S} : A \coll B}},
type={termsandsymbols}
} 

\newglossaryentry{multi_lvl_bc}
{
name= {multi-level base change},
description={See \autoref{ssec:multi_bc}},
symbol={\ensuremath{\Bunbc^n_\cC}},
type={termsandsymbols}
} 

\newglossaryentry{klvl_bc}
{
name= {$k$-level base change},
description={Can refer to both a $k$-level base change functor $H$ (see \autoref{defn:klvl_base_change}) or the multi-level base change derived from it (see \autoref{ssec:multi_bc}).},
symbol={\ensuremath{\empty}},
type={termsandsymbols}
} 

\newglossaryentry{ass_multi_lvl_collapse}
{
name= {multi-level collapse associated to ordered collapse sequence},
description={Every ordered collapse sequence $\vec\lambda$ induces a multi-level collapse $\vec \sS ^{\vec\lambda}$. See \autoref{constr:multilevel_collapse}.},
symbol={\ensuremath{\vsS{\color{lightgray}\vec\lambda}}},
type={termsandsymbols}
}

\newglossaryentry{injections}
{
name= {injection of singular interval families},
description={An injection $\lambda : \scA \into \scB$ of families $\scA, \scB : X \to \SI$ is a monomorphic natural transformation. See \autoref{defn:natural_injections}},
symbol={\ensuremath{\into}},
type={termsandsymbols}
}

\newglossaryentry{ass_subset_section}
{
name= {singular subset sections associated to injections},
description={Given an injection $\lambda : \scA \coll \scB$ of families $\scA, \scB : X \to \SI$ its associated singular subset section $\cS^\lambda$ is the subset section defined by its image. See \autoref{defn:stability_vs_injections}.},
symbol={\ensuremath{\cS^{\color{lightgray}\lambda}}},
type={termsandsymbols}
}

\newglossaryentry{SI_ass_SSS}
{
name= {singular interval family associated to stable singular subset section},
description={See \autoref{defn:stability_vs_injections}.},
symbol={\ensuremath{\intrel{-}}},
type={termsandsymbols}
} 

\newglossaryentry{SI_inj_ass_SSS}
{
name= {injection associatied to subsets of singular heights},
description={For single intervals see \autoref{defn:eta_inclusion} and for families see \autoref{defn:stability_vs_injections}.},
symbol={\ensuremath{\eta}},
type={termsandsymbols}
} 

\newglossaryentry{stable_subset_section}
{
name= {stable},
description={Used in the context of \hyperref[glsentry-subset_section]{stable singular subset sections}.},
symbol={\ensuremath{\empty}},
type={termsandsymbols}
}

\newglossaryentry{subset_section}
{
name= {singular subset sections},
description={A singular subset section $\cF$ of a \SI-family is a section of subsets of singular heights of the family's bundle. The notion specialises to stable singular subset sections which are images of injections. See \autoref{defn:subset_sections}.},
symbol={\ensuremath{\empty}},
type={termsandsymbols}
}

\newglossaryentry{open_map}
{
name= {open functor},
description={An open functor of singular intervals, is a functor of singular intervals which preserves regular segments. See \autoref{eg:open_maps}.},
symbol={\ensuremath{\empty}},
type={termsandsymbols}
} 

\newglossaryentry{underlying_mono}
{
name= {underlying monomorphism},
description={See \autoref{constr:underlying_monomorphism_fctr}.},
symbol={\ensuremath{\monosing}},
type={termsandsymbols}
} 

\newglossaryentry{coll_map_SI}
{
name= {collapse bundle map},
description={A collapse bundle map $\sS^\lambda : \pi_\scB \coll \pi_\scA$ is a bundle map canonically constructed from an injection $\lambda : \scA \mono \scB$ of families $\scA, \scB : X \to \SIvert {} \cC$. See \autoref{defn:glambda}.},
symbol={\ensuremath{\sS^{\color{lightgray}\lambda}}},
type={termsandsymbols}
}

\newglossaryentry{coll_notn}
{
name= {collapse notation},
description={The notation ``$\coll$" is used both for \hyperref[glsentry-multi_lvl_collapse]{multi-level collapse} and \hyperref[glsentry-coll_SI]{collapse of interval families}.},
symbol={\ensuremath{\coll}},
type={termsandsymbols}
}

\newglossaryentry{coll_SI}
{
name= {collapse of interval families},
description={A collapse $\lambda : \scA \coll \scB$ for $\scA, \scB : X \to \SIvert {} \cC$ is a natural transformation $\und\scA \to \und\scB$, such that its associated collapse bundle map $\sS^\lambda$ factors the labelling functors of $\scA$ and $\scB$. See \autoref{defn:collapse_of_SI_families}.},
symbol={\ensuremath{\coll}},
type={termsandsymbols}
}

\newglossaryentry{bun_coll_cat}
{
name= {category of interval families and collapses},
description={A category of $\SIvert {} \cC$-bundles and their collapses. Defined in \autoref{constr:cat_of_bun_and_coll}. The category has pushouts as proven in \autoref{thm:family_pushouts}.},
symbol={\ensuremath{\Bun^\coll_\cC}},
type={termsandsymbols}
}

\newglossaryentry{pullback_coll}
{
name= {pullback of collapse},
description={Given a collapse $\lambda : \scB \to \widetilde \scA$ then this can be ``pre-whiskered" by functors $H$ that precompose with $\scB, \scA$. This is defined in \autoref{claim:collapse_dimension_interaction1}.},
symbol={\ensuremath{\lambda H}},
type={termsandsymbols}
}

\newglossaryentry{pushforward_coll}
{
name= {pushforward of collapse along functor},
description={Given a collapse $\lambda : \scA \to \widetilde \scA$ for $\scA : X \to \SIvert {} \cC$ and $H : X \to Y$ such that $\scB H = \scA$, then there is a pushforward of collapse $H_* : \scB \coll H_* \widetilde \scA$. This is defined in \autoref{claim:collapse_dimension_interaction2}.},
symbol={\ensuremath{{\color{lightgray}H}_* {\color{lightgray}\lambda}}},
type={termsandsymbols}
}

\newglossaryentry{pushforward_coll_coll}
{
name= {pushforward of collapse along collapse},
description={See \autoref{lem:coll_push}.},
symbol={\ensuremath{{\color{lightgray}\mu}_* {\color{lightgray}\lambda}}},
type={termsandsymbols}
}

\newglossaryentry{pushforward_notn}
{
name= {pushforward notation},
description={Pushforward notation can refer to \hyperref[glsentry-pushforward_coll]{pushforward of collapse along functors}, \hyperref[glsentry-pushforward_coll_coll]{pushforward of collapse along collapse} or pushforward of embedding along collapse ($(\vec\lambda)\pbstar \theta$ \hyperref[glsentry-coll_emb]{collapsed embedding}).},
symbol={\ensuremath{\empty_*}},
type={termsandsymbols}
}

\newglossaryentry{repack}
{
name= {repacking},
description={A construction applied to a certain tuple $(V,U)$ yielding a functor $\sR_{V,U}$ into $\SIvert {} \cC$. See \autoref{defn:unpacking_and_repacking}.},
symbol={\ensuremath{\sR}},
type={termsandsymbols}
} 

\newglossaryentry{unpack}
{
name= {unpacking},
description={A construction applied to $R : X \to \SIvert {} \cC$ producing a tuple $(\sV_R : X \to \SI, \sU_R : \sG(\sV_R) \to \cC)$. See \autoref{defn:unpacking_and_repacking}.},
symbol={\ensuremath{(\sV, \sU)}},
type={termsandsymbols}
}

\newglossaryentry{towers}
{
name= {towers of \SI-bundles},
description={See \autoref{defn:towers_of_bundles}.},
symbol={\ensuremath{\empty}},
type={termsandsymbols}
}

\newglossaryentry{SI_notn}
{
name= {singular cube notation},
description={The ``over double category" notation $\SIvert {} {}$ is used for \hyperref[glsentry-SIvert1cat]{$\cC$-labelled singular intervals}, \hyperref[glsentry-relabel]{relabelling of labelled singular intervals}, \hyperref[glsentry-SIvertcat]{$\cC$-labelled singular $n$-cubes} and \hyperref[glsentry-krelabel]{relabelling of labelled singular $n$-cubes}.},
symbol={\ensuremath{\SIvert {\bullet} {}}},
type={termsandsymbols}
}

\newglossaryentry{krelabel}
{
name= {relabelling of $\cC$-labelled singular $n$-cube},
description={Given a functor $F : \cC \to \cD$ the relabelling on labelled singular $n$-cube induced by $F$ is a functor $\SIvert n F : \SIvert n \cC \to \SIvert n \cD$ defined in \autoref{defn:grothendieck_cubes}.},
symbol={\ensuremath{\SIvert n {\color{lightgray}F}}},
type={termsandsymbols}
}

\newglossaryentry{SIvertcat}
{
name= {$\cC$-labelled singular $n$-cube category},
description={The category of $\cC$-labelled singular $n$-cubes is inductively defined in  \autoref{defn:grothendieck_cubes} by setting $\SIvert n \cC = \SIvert {} {\SIvert {n-1} \cC}$.},
symbol={\ensuremath{\SIvert n \cC}},
type={termsandsymbols}
}

\newglossaryentry{SIvert1cat}
{
name= {$\cC$-labelled singular interval category},
description={See \autoref{defn:SIvert}.},
symbol={\ensuremath{\SIvert {} \cC}},
type={termsandsymbols}
}

\newglossaryentry{SIvertfam2}
{
name= {$\SIvert n \cC$-family},
description={See \gls{SIvertfam}.},
symbol={\ensuremath{\empty}},
type={termsandsymbols}
}

\newglossaryentry{SIvertfam3}
{
name= {$n$-cube},
description={The topological $n$-cube $[0,1]^n$, or, an abbreviation for \gls{SIvertfam}.},
symbol={\ensuremath{\empty}},
type={termsandsymbols}
}

\newglossaryentry{SIvertfam}
{
name= {$\cC$-labelled singular $n$-cube family},
description={A functor $\scA : X \to \SIvert n \cC$ as defined in \autoref{defn:grothendieck_cubes}. If $X = \bnum 1$ we speak of a cube rather than a family.},
symbol={\ensuremath{\empty}},
type={termsandsymbols}
}

\newglossaryentry{krepack}
{
name= {$k$-repacking},
description={The mapping of a tower of $\SI$-bundles $T$ together with a  labelling $U$ of the tower's total space to a functor $\sR^k_{T, U}$. Defined in \autoref{defn:complete_unpacking_and_repacking}.},
symbol={\ensuremath{\tsR k_{\color{lightgray}T,U}}},
type={termsandsymbols}
}

\newglossaryentry{kunpack}
{
name= {$k$-unpacking},
description={The mapping from $\scA : X \to \SIvert n \cC$ to $(\sT^k_\scA, \sU^k_\scA)$ defined in   \autoref{defn:complete_unpacking_and_repacking}.},
symbol={\ensuremath{(\sT^k_{\color{lightgray}\scA}, \tsU k_{\color{lightgray}\scA})}},
type={termsandsymbols}
}

\newglossaryentry{ktower}
{
name= {height $k$ tower of cube family},
description={$\sT^k$ is applied to a family $\scA : X \to \SIvert n \cC$ ($n \geq k$) and denotes height $k$ tower of its first $k$ $\SI$-bundles. This is defined in \autoref{defn:complete_unpacking_and_repacking}.},
symbol={\ensuremath{\sT^k}},
type={termsandsymbols}
}

\newglossaryentry{kbund}
{
name= {$k$-level bundle},
description={The $k$-level bundle $\tpi k_\scA : \tsG k (\scA) \to \tsG {k-1} (\scA)$ of a family $\scA : X \to \SIvert n \cC$ ($n \geq k$). Defined in \autoref{defn:complete_unpacking_and_repacking}.},
symbol={\ensuremath{\tpi k}},
type={termsandsymbols}
}

\newglossaryentry{klabel}
{
name= {$k$-level labelling},
description={The $k$-level labelling $\tsU k _\scA : \tsG k (\scA) \to \SIvert {n-k} \cC$ of a family $\scA : X \to \SIvert n \cC$ ($n \geq k$). Defined in \autoref{defn:complete_unpacking_and_repacking}.},
symbol={\ensuremath{\tsU k}},
type={termsandsymbols}
} 

\newglossaryentry{multi_lvl_base_change}
{
name= {multi-level basechange},
description={See \autoref{defn:multilevel_base_change}.},
symbol={\ensuremath{\empty}},
type={termsandsymbols}
}  

\newglossaryentry{ktotalspace}
{
name= {$k$-level total space, total base change},
description={If $\tsG k$ is applied to a family $\scA : X \to \SIvert n \cC$ ($n \geq k$) then $\tsG k (\scA)$ refers to the $k$-level total space  $\tsG k(\scA)$ defined in \autoref{defn:complete_unpacking_and_repacking}. If it is applied to a functor of posets $H: X \to Y$, and $\scA = \scB H$ for $\scB : Y \to \SIvert n \cC$ then this refers to  the $k$-level total base change $\tsG k (H) : \tsG k(\scA) \to \tsG k(\scB)$ defined in \autoref{constr:unpacking_collapse}. Note that in this case, indices can become negative as explained in \autoref{notn:k_lvl_basechange}},
symbol={\ensuremath{\tsG k}},
type={termsandsymbols}
}

\newglossaryentry{identity_cube}
{
name= {identity cubes},
description={See \autoref{defn:identities}. Note that the notation $\Id$ is exclusively used for identity cubes. Categorical identities are denoted by $\id$.},
symbol={\ensuremath{\Id}},
type={termsandsymbols}
}

\newglossaryentry{region_type}
{
name= {region type},
description={See \autoref{defn:regions}.},
symbol={\ensuremath{\ctyp}},
type={termsandsymbols}
}

\newglossaryentry{region_dim}
{
name= {region dimension},
description={See \autoref{defn:regions}.},
symbol={\ensuremath{\ctypsum}},
type={termsandsymbols}
} 

\newglossaryentry{uncurrying}
{
name= {uncurrying},
description={See \autoref{notn:cartesian_closed}.},
symbol={\ensuremath{\cY}},
type={termsandsymbols}
}

\newglossaryentry{currying}
{
name= {currying},
description={See \autoref{notn:cartesian_closed}.},
symbol={\ensuremath{\cN}},
type={termsandsymbols}
}

\newglossaryentry{extended_sing_coll}
{
name= {extended singular collapse},
description={See \autoref{constr:uncurrying}.},
symbol={\ensuremath{\sE}},
type={termsandsymbols}
}

\newglossaryentry{relabel}
{
name={relabelling},
description={Given a functor $F : \cC \to \cD$ we obtain a functor $\SIvert {} F : \SIvert {} \cC \to \SIvert {} \cD$ which morally ``relabels labels by acting on them with $F$". See \autoref{defn:transfer_of_coloring}.},
symbol= {\ensuremath{\SIvert {} F}},
type={termsandsymbols}
}

\newglossaryentry{SIvertsec}
{
name={labelling of labelled singular interval (morphism)},
description={Objects and morphism in $\SIvert {} \cC$ are tuples. The labelling of these is the second component of the tuple. See \autoref{defn:SIvert}.},
symbol= {\ensuremath{\SIs{(-)}}},
type={termsandsymbols}
}

\newglossaryentry{SIvertforget}
{
name={forgetting labels},
description={Passing to underlying singular interval structure from labelled singular intervals $\SIvert {} \cC$ gives rise to a ``label-forgetting functor" $(\und{-}) : \SIvert {} \cC \to \SI$ mapping objects $I$ to $\SIf I$ and morphisms $f$ to $\SIf f$. See \autoref{defn:SIvert} and \autoref{defn:label_forgetting}. Note that by the inductive \autoref{defn:grothendieck_cubes} of $\SIvert n \cC$ we similarly obtain a functor $\und{-} : \SIvert n \cC \to \SI$.},
symbol= {\ensuremath{\und{(-)}}},
type={termsandsymbols}
}

\newglossaryentry{horcomp}
{
name={horizontal composition},
description={Composition of relations (see \autoref{defn:relations_and_profunctors}), of profunctorial relations (see \autoref{constr:PRel}), or of labels of morphisms in $\SIvert {} \cC$ (see \autoref{defn:SIvert}).},
symbol= {\ensuremath{\odot}},
type={termsandsymbols}
}

\newglossaryentry{funnat}
{
name={natural transformation associated to functor},
description={See \autoref{notn:functor_morphism_action}.},
symbol= {\ensuremath{\nu_F}},
type={termsandsymbols}
}

\newglossaryentry{SIbunpull}
{
name={singular interval bundle pullback},
description={See \autoref{rmk:grothendieck_span_construction_basechange}.},
symbol= {\ensuremath{\empty}},
type={termsandsymbols}
}

\newglossaryentry{SIfam2}
{
name={\SI-bundle},
description={See \gls{SIbun}.},
symbol= {\ensuremath{\empty}},
type={termsandsymbols}
}

\newglossaryentry{SIbun2}
{
name={\SI-family},
description={See \gls{SIfam}.},
symbol= {\ensuremath{\empty}},
type={termsandsymbols}
}

\newglossaryentry{SIfam}
{
name={singular interval family},
description={A functor from a poset $X$ to the category $\SI$. See \autoref{constr:SI_families}.},
symbol= {\ensuremath{\empty}},
type={termsandsymbols}
}

\newglossaryentry{SIbun}
{
name={singular interval bundle},
description={A bundle obtained from an \SI-family. See \autoref{constr:SI_families}.},
symbol= {\ensuremath{\empty}},
type={termsandsymbols}
}

\newglossaryentry{SIbunmap}
{
name={map of \SI-bundles},
description={A map of bundles which is fibrewise monotone. See \autoref{constr:SI_families}.},
symbol= {\ensuremath{\empty}},
type={termsandsymbols}
}

\newglossaryentry{bisurjective}
{
name={bisurjective},
description={A property of relations of between sets (``any element is related to some other element"), proven to hold for profunctorial realisations of singular height morphisms. See \autoref{cor:relation_fullness}.},
symbol= {\ensuremath{\empty}},
type={termsandsymbols}
}

\newglossaryentry{bimon}
{
name={bimonotone},
description={A property of relations of between sets of numbers (``monotonicity in both directions"), proven to hold for profunctorial realisations of singular height morphisms. See \autoref{claim:order_realisations_monotone}.},
symbol= {\ensuremath{\empty}},
type={termsandsymbols}
}

\newglossaryentry{mon}
{
name={monotone},
description={Preserving the order on integers.},
symbol= {\ensuremath{\empty}},
type={termsandsymbols}
}

\newglossaryentry{fillind}
{
name={filler's index},
description={See \autoref{defn:successor}.},
symbol= {\ensuremath{\succindex{-}}},
type={termsandsymbols}
}

\newglossaryentry{filled}
{
name={filler edge},
description={The edge set $\edgeset (f : I_\soe \to I_\tae)$ has a total order. If $\edb$ is the successor of $\eda$ in that order then the filler edge of $\eda$ is an arrow completing $\edb$ and $\eda$ to a ``triangle". See \autoref{defn:successor}.},
symbol= {\ensuremath{\succfill{-}}},
type={termsandsymbols}
}

\newglossaryentry{succed}
{
name={successor edge},
description={The edge set $\edgeset (f : I_\soe \to I_\tae)$ has a total order, the successor $\succ\eda$ of an edge $\eda$ is the successor with respect to that order. See \autoref{defn:successor}.},
symbol= {\ensuremath{\succ{-}}},
type={termsandsymbols}
}

\newglossaryentry{preded}
{
name={predecessor edge},
description={The edge set $\edgeset (f : I_\soe \to I_\tae)$ has a total order, the successor $\succ\eda$ of an edge $\eda$ is the successor with respect to that order. See \autoref{defn:predecessor}.},
symbol= {\ensuremath{\pred{-}}},
type={termsandsymbols}
}

\newglossaryentry{soed}
{
name={source of edge},
description={An edge $\eda \in \edgeset (f : I_\soe \to I_\tae)$ in a edge set is tuple of numbers, and its source is the first number of the tuple. See \autoref{defn:edge_sets}.},
symbol= {\ensuremath{{\color{lightgray}\eda}_\soe}},
type={termsandsymbols}
}

\newglossaryentry{taed}
{
name={target of edge},
description={An edge $\eda \in \edgeset (f : I_\soe \to I_\tae)$ in a edge set is tuple of numbers, and its source is the second number of the tuple. See \autoref{defn:edge_sets}.},
symbol= {\ensuremath{{\color{lightgray}\eda}_\tae}},
type={termsandsymbols}
}

\newglossaryentry{edgenorm}
{
name= {norm of an edge},
description={The norm $\avg{\eda}$ is sum of the numbers of the two endpoints of an edge $\eda \in \edgeset(f)$. See \autoref{defn:edge_sets}.},
symbol={\ensuremath{\avg{-}}},
type={termsandsymbols}
}

\newglossaryentry{edgeset}
{
name= {edge set},
description={The set $\edgeset(f)$ is the set of tuples in the profunctorial realisation of a singular height morphism $f$. See \autoref{defn:edge_sets}.},
symbol={\ensuremath{\edgeset(f)}},
type={termsandsymbols}
}

\newglossaryentry{prelreal}
{
name= {profunctorial realisation},
description={A canonical embedding $\SiR : \SI \into \PRel$, obtained by minimal completion of the joint graphs of singular interval morphisms and their \rsdual. See \autoref{defn:order_realisation}. See also \autoref{claim:order_realisation_symmetric_def}.},
symbol={\ensuremath{\SiR}},
type={termsandsymbols}
}

\newglossaryentry{ambidextcond}
{
name= {ambidexterity condition},
description={See \autoref{rmk:ambidext_cond}.},
symbol={\ensuremath{\empty}},
type={termsandsymbols}
}

\newglossaryentry{openint}
{
name= {open interval of integers},
description={For $a,b \in \lZ$, $]a,b[$ denotes the open interval with endpoints $a,b$. That is $]a,b[$ contained number $a + 1, a  + 2, ... , b -1$. See \autoref{eg:posets}.},
symbol={\ensuremath{]a,b[}},
type={termsandsymbols}
}

\newglossaryentry{closedint}
{
name= {closed interval},
description={For $a,b \in \lZ$, $[a,b]$ denotes the closed interval with endpoints $a,b$. See \autoref{eg:posets}.},
symbol={\ensuremath{[a,b]}},
type={termsandsymbols}
}

\newglossaryentry{singop}
{
name= {singular dual},
description={The singular \rsdual{} of a regular-segment morphism is a canonically associated singular-height morphism. See \autoref{defn:regop_singop}.},
symbol={\ensuremath{{\color{lightgray}g}\singop}},
type={termsandsymbols}
}

\newglossaryentry{regop}
{
name= {regular dual},
description={The regular \rsdual{} of a singular-height morphism is a canonically associated regular-segment morphism. See \autoref{defn:regop_singop}.},
symbol={\ensuremath{{\color{lightgray}f}\regop}},
type={termsandsymbols}
}

\newglossaryentry{extsingmor}
{
name= {extended singular height morphism},
description={The extension $\widehat{f}$ of a singular height morphism $f : \singint k \to\oSI \singint l$ is its extension to the extended singular heights $\extsing$, and this extension is defined to preserve the first and last extended singular height. See \autoref{notn:singular_morphism_boundary_cases}.},
symbol={\ensuremath{\widehat{\color{lightgray}f}}},
type={termsandsymbols}
}

\newglossaryentry{singmor}
{
name= {singular height morphism},
description={A morphism of singular interval is a mapping of its singular heights. See \autoref{defn:singular_intervals_morphism}.},
symbol={\ensuremath{\to\oSI}},
type={termsandsymbols}
}

\newglossaryentry{regmor}
{
name= {regular segment morphism},
description={A regular segment morphisms of singular intervals is a mapping of its regular segments. See \autoref{defn:singular_intervals_morphism}.},
symbol={\ensuremath{\to\regop\oSI }},
type={termsandsymbols}
}

\newglossaryentry{SI}
{
name= {category of singular intervals},
description={See \autoref{defn:singular_intervals}.},
symbol={\ensuremath{\SI}},
type={termsandsymbols}
}

\newglossaryentry{singint}
{
name= {singular interval},
description={A zig-zag poset of natural number $0$ to $2k$. See \autoref{defn:singular_intervals}.},
symbol={\ensuremath{\singint k}},
type={termsandsymbols}
}

\newglossaryentry{extsing}
{
name= {extended singular heights},
description={The extended singular heights of a singular interval $I$ are its odd numbers together with the ``extended singular heights" $(-1)$ and $2\iH_I +1$ (topologically these can be thought of as the endpoints of the interval).  See \autoref{defn:singular_intervals}.},
symbol={\ensuremath{\extsing}},
type={termsandsymbols}
}

\newglossaryentry{sing}
{
name= {singular heights},
description={The singular heights of a singular interval are its odd numbers. Topologically they correspond to points in the interior of the interval. See \autoref{defn:singular_intervals}.},
symbol={\ensuremath{\singcont(-)}},
type={termsandsymbols}
}

\newglossaryentry{reg}
{
name= {regular segments},
description={The regular segments of a singular interval are its even numbers. Topologically they correspond to open connected segments of the interval. See \autoref{defn:singular_intervals}.},
symbol={\ensuremath{\regcont(-)}},
type={termsandsymbols}
}

\newglossaryentry{totalorders}
{
name= {total finite orders},
description={The category of finite non-empty total orders (also known as the category of simplices), see \autoref{ssec:SIvert}.},
symbol={\ensuremath{\BDelta}},
type={termsandsymbols}
}

\newglossaryentry{face}
{
name= {face map},
description={The face maps in the simplicial category, see \autoref{ssec:coloring}.},
symbol={\ensuremath{\delta^i_k}},
type={termsandsymbols}
}

\newglossaryentry{degeneracy}
{
name= {degeneracy map},
description={The degeneracy maps in the simplicial category, see \autoref{ssec:coloring}.},
symbol={\ensuremath{\sigma^i_k}},
type={termsandsymbols}
}

\newglossaryentry{bnum}
{
name= {numeral poset},
description={The totally ordered poset with $n$ elements $\Set{0 < 1 < ... < n-1}$. See \autoref{eg:posets}.},
symbol={\ensuremath{\bnum n}},
type={termsandsymbols}
}

\newglossaryentry{bang}
{
name= {terminal functor},
description={The unique functor from a category to the terminal category, see \autoref{rmk:terminal_coloring}.},
symbol={\ensuremath{\bang}},
type={termsandsymbols}
}

\newglossaryentry{labelpushout}
{
name= {pushouts of labelled posets},
description={A pushout diagram in labelled posets $\Poss \cC$ whose legs have \gls{lifts} admits a pushout. See \autoref{constr:pushouts_in_labelled_posets}.},
symbol={\ensuremath{\cup}},
type={termsandsymbols}
}

\newglossaryentry{posslash}
{
name= {labelled posets},
description={Full subcategory of $\Cat \slash \cC$ consisting of functors into $\cC$ from posets. See \autoref{defn:posslash}.},
symbol={\ensuremath{\Poss \cC}},
type={termsandsymbols}
}

\newglossaryentry{downwardclosed}
{
name= {downward closed},
description={A condition on posets. See \autoref{defn:downward_closed_subposet}},
symbol={\empty},
type={termsandsymbols}
}

\newglossaryentry{lifts}
{
name= {lifts},
description={A ``map having lifts" is a condition on maps of posets given in \autoref{defn:having_lifts}},
symbol={\empty},
type={termsandsymbols}
}

\newglossaryentry{codrest}
{
name= {restriction to codomain},
description={Given a functor $H : \cC \to \cD$ injective on objects and morphism, and functor $F : \cE \to \cD$, then $H^*F : F\inv(\im(H)) \to \cD$ denotes a pullback of $F$ along $H$. See \autoref{notn:subsets_and_restrictions}.},
symbol={\ensuremath{{\color{lightgray}H}^*{\color{lightgray}F}}},
type={termsandsymbols}
}

\newglossaryentry{rest}
{
name= {restriction to subcategory},
description={If $F : \cC \to \cD$ is a functor, and $\cC_0 \subset \cC$ then $\rest F {\cC_0}$ denotes the restriction of a functor $F$ to $\cC_0$, as explained in \autoref{notn:subsets_and_restrictions}. Note that sometimes $\cC_0$ might take the form of a single object $x$ or morphism $x \to y$ as well.},
symbol={\ensuremath{\rest {\color{lightgray} F} {\cC_0}}},
type={termsandsymbols}
}

\newglossaryentry{rest2}
{
name= {restriction to fibre},
description={See \autoref{notn:fiber_restrictions}. The same notation is used for \hyperref[glsentry-rest]{restriction to subcategories}.},
symbol={\ensuremath{\rest {\color{lightgray} F} x}},
type={termsandsymbols}
}

\newglossaryentry{testfunc}
{
name= {test functors},
description={Given a category $\cC$, and $x \in \obj(\cC)$, $f \in \mor(\cC)$ then the test functors $\Delta_x : \bnum 1 \to \cC$ and $\Delta_f : \bnum 2 \to \cC$ have image $x$ respectively $f$. See \autoref{notn:subsets_and_restrictions}.},
symbol={\ensuremath{\Delta}},
type={termsandsymbols}
}

\newglossaryentry{secp}
{
name= {fiber value},
description={Given a point $p$ in a total space of a bundle $\pi_F$, $\secp p$ denotes its value in the fiber $\pi\inv_F(\pi_F(p))$.  See \autoref{notn:total_poset_projs}.},
symbol={\ensuremath{\secp p}},
type={termsandsymbols}
}

\newglossaryentry{bundleproj}
{
name= {bundle map},
description={Usually a bundle map $\pi_F : \sG(F) \to X$ associated to a \PRel-family $F : X \to \PRel$, projecting the total space $\sG(F)$ to the base space $X$. See \autoref{defn:grothendieck_construction}},
symbol={\ensuremath{\mathchar"119}},
type={termsandsymbols}
}

\newglossaryentry{prodproj}
{
name= {product projection},
description={The projection of a product $X \times Y$ of posets to one of its components. See \autoref{notn:product_categories}},
symbol={\ensuremath{\ppi}},
type={termsandsymbols}
}

\newglossaryentry{totalspace}
{
name= {total space, total base change},
description={If $\sG$ is applied to a family $F : X \to \PRel$ then $\sG(F)$ is the domain of a bundle map $\pi_F$, called the total space of $F$. This is explained \autoref{defn:grothendieck_construction}. If $\sG$ is applied to a base change $H$ of families $F,G$ then $\sG(H) : \sG(F) \to \sG(H)$ denotes the total base change as explained in \autoref{defn:grothendieck_base_change}.},
symbol={\ensuremath{\mathscr{G}}},
type={termsandsymbols}
}

\newglossaryentry{settruth}
{
name= {truth values as sets},
description={A functor $\iP : \Bool \into \SetCat$ mapping $\bot$ to the empty set $\emptyset$ and $\top$ to the singleton set $\Set{*}$. See \eqref{eq:Bool_to_Set}.},
symbol={\ensuremath{\mathsf{P}}},
type={termsandsymbols}
}

\newglossaryentry{converse}
{
name= {converse of (profunctorial) relation},
description={For $R : X \xslashedrightarrow{} Y$ a (profunctorial) relation, its converse $R^\top$ is obtained by precomposing with the contravariant symmetry $Y\op\times X \iso X\op \times Y$.},
symbol={\ensuremath{{\color{lightgray}R}^\top}},
type={termsandsymbols}
}

\newglossaryentry{grph}
{
name= {graph of poset maps},
description={For $F : X \to Y$, $\grph_F$ is $\Hom_Y(F-,-)$. See \autoref{defn:graph_of_functors}},
symbol={\ensuremath{\mathrm{grph}}},
type={termsandsymbols}
}

\newglossaryentry{PRel}
{
name= {profuncorial relation category},
description={$\PRel$ is the category of profunctorial relations. See \autoref{constr:PRel}},
symbol={\ensuremath{\mathbf{PRel}}},
type={termsandsymbols}
}

\newglossaryentry{Rel}
{
name= {relation set},
description={$\Rel(X,Y)$ is the set of relations between $X$ and $Y$. See \autoref{defn:relations_and_profunctors}},
symbol={\ensuremath{\mathrm{Rel}}},
type={termsandsymbols}
}

\newglossaryentry{obj}
{
name= {object set functor},
description={The functor $\obj: \Cat \to \SetCat$ maps categories to their object sets and functors to their maps of objects. It is often kept implicit. See \autoref{notn:basic_category_theory}},
symbol={\ensuremath{\obj}},
type={termsandsymbols}
}

\newglossaryentry{discr_cat}
{
name= {discrete category functor},
description={See \autoref{notn:basic_category_theory}},
symbol={\ensuremath{\discr}},
type={termsandsymbols}
} 

\newglossaryentry{discr_prof}
{
name= {discrete profunctorial relation functor},
description={See \autoref{rmk:set_valued_grothendieck_construction}},
symbol={\ensuremath{\Discr}},
type={termsandsymbols}
}

\newglossaryentry{mor_cat}
{
name= {morphism set},
description={If $\cC$ is an ordinary category then $\mor(\cC)$ denotes the set of all morphisms in $\cC$. A similar notation is used for \hyperref[glsentry-mor_pres_cat]{morphisms of presented associative $n$-categories}.},
symbol={\ensuremath{\mathrm{mor}}},
type={termsandsymbols}
}

\newglossaryentry{Setcat}
{
name= {category of sets},
description={See \autoref{notn:basic_category_theory}.},
symbol={\ensuremath{\SetCat}},
type={termsandsymbols}
}

\newglossaryentry{function_set}
{
name= {set of functions},
description={See \autoref{notn:basic_category_theory}.},
symbol={\ensuremath{\Map}},
type={termsandsymbols}
}

\newglossaryentry{Catcat}
{
name= {category of categories},
description={See \autoref{notn:basic_category_theory}.},
symbol={\ensuremath{\Cat}},
type={termsandsymbols}
}

\newglossaryentry{functor_category}
{
name= {category of functors},
description={See \autoref{notn:basic_category_theory}.},
symbol={\ensuremath{\Fun}},
type={termsandsymbols}
}

\newglossaryentry{Boolcat}
{
name= {category of booleans},
description={See \autoref{notn:basic_category_theory}.},
symbol={\ensuremath{\Bool}},
type={termsandsymbols}
} 

\newglossaryentry{Poscat}
{
name= {category of posets},
description={See \autoref{defn:posets}.},
symbol={\ensuremath{\Pos}},
type={termsandsymbols}
}

\newglossaryentry{constant_functor}
{
name= {constant functor},
description={See \autoref{notn:simplicial}.},
symbol={\ensuremath{\const}},
type={termsandsymbols}
}

\newglossaryentry{ps_stra}
{
name= {projection-stable stratification},
description={See \autoref{ssec:po_mfld_diag}.},
symbol={\ensuremath{}},
type={termsandsymbols}
}

\newglossaryentry{ffc_stra}
{
name= {flag-foliation-compatible stratification},
description={See \autoref{ssec:po_mfld_diag}.},
symbol={\ensuremath{}},
type={termsandsymbols}
}

\newglossaryentry{top_cub}
{
name= {topological $n$-cube},
description={See \autoref{ssec:po_mfld_diag}.},
symbol={\ensuremath{}},
type={termsandsymbols}
}

\newglossaryentry{dir_triang}
{
name= {triangulation, directed},
description={See \autoref{ssec:coloring}.},
symbol={\ensuremath{\abs{-}_\bullet}},
type={termsandsymbols}
}

\newglossaryentry{ggamma_tp}
{
name= {namescope realisation, total poset},
description={See \autoref{defn:namescope}},
symbol={\ensuremath{\GGamma{k}{-}}},
type={termsandsymbols}
}

\newglossaryentry{gen_comp}
{
name= {generic composites},
description={Manifold diagrams representing morphisms, with manifold in generic position. We give two characterisations: see \autoref{defn:top_down_gen_comp} or \autoref{constr:bottom_up_gcomp}.},
symbol={\ensuremath{\GComp}},
type={termsandsymbols}
}

\newglossaryentry{interchanger}
{
name= {interchanger},
description={The lowest-dimensional (generic) homotopy. See \autoref{notn:interchanger}.},
symbol={\ensuremath{\interchanger}},
type={termsandsymbols}
}

\newglossaryentry{poset_skel}
{
name= {posetal skeleton},
description={See \autoref{defn:posets}.},
symbol={\ensuremath{\PSk}},
type={termsandsymbols}
}

\newglossaryentry{poset_skel_label}
{
name= {posetal skeleton of category of elements},
description={See \autoref{constr:panc_from_glob}.},
symbol={\ensuremath{\elpsk}},
type={termsandsymbols}
}

\newglossaryentry{cubeo_cat}
{
name= {category of topological $\cC$-labelled cubes},
description={See \autoref{rmk:cube_cat}.},
symbol={\ensuremath{\Cubeo n \cC}},
type={termsandsymbols}
}

\newglossaryentry{bunbc_cat}
{
name= {category of cubes and multi-level base change},
description={See \autoref{sec:sum_Buno}.},
symbol={\ensuremath{\Bunbc^n_\cC}},
type={termsandsymbols}
}

\newglossaryentry{buno_cat}
{
name= {category of cubes and open multi-level base change},
description={See \autoref{sec:sum_Buno2}.},
symbol={\ensuremath{\Buno n \cC}},
type={termsandsymbols}
}

\newglossaryentry{cube_equiv}
{
name= {equivalence of topological (labelled) cubes},
description={See \autoref{ssec:po_mfld_diag}.},
symbol={\ensuremath{\simeq}},
type={termsandsymbols}
}

\newglossaryentry{geom_real}
{
name= {geometric realisation},
description={See \autoref{ssec:coloring}.},
symbol={\ensuremath{\norm{-}}},
type={termsandsymbols}
}

\newglossaryentry{min_label}
{
name= {minimal labelling category},
description={See \autoref{defn:sum_red_lab_cat}.},
symbol={\ensuremath{\redGamma \sC}},
type={termsandsymbols}
}

\newglossaryentry{pres_cat_n}
{
name= {category of $n$-presentations},
description={See \autoref{defn:cat_of_n_pres}.},
symbol={\ensuremath{\pCat_n}},
type={termsandsymbols}
}

\newglossaryentry{pres_wcat}
{
name= {category of weak $\infty$-presentations},
description={See \autoref{ssec:weakification}.},
symbol={\ensuremath{\pwCat}},
type={termsandsymbols}
}

\newglossaryentry{inf_eq_res}
{
name= {invertible elements resolving algebra equations},
description={See \autoref{constr:resolutions}.},
symbol={\ensuremath{\infeq {f} {k}}},
type={termsandsymbols}
}

\newglossaryentry{glob_to_pres}
{
name= {presented associative $n$-category from globular set},
description={See \autoref{ssec:pres_to_globe}.},
symbol={\ensuremath{\kC}},
type={termsandsymbols}
}

\newglossaryentry{alg_resol}
{
name= {resolutions of algebras},
description={See \autoref{constr:resolutions}.},
symbol={\ensuremath{\kiC{}(-)}},
type={termsandsymbols}
}

\newglossaryentry{V}
{
name={Edges},
description={These usually denote elements of edge sets  $\edgeset (f)$. See \autoref{defn:edge_sets}.},
vars= {\ensuremath{\eda,\edb,\edc,...}},
type={variables}
}

\newglossaryentry{fams}
{
name={Families},
description={These usually denote families (of singular intervals, labelled cubes, etc.)},
vars={$\scA, \scB, \scC, \scD$},
type={variables}
}

\newglossaryentry{injs}
{
name={Injections},
description={These usually denote injections, that is, the data for collapses.},
vars={$\lambda, \mu, \eps$},
type={variables}
}

\newglossaryentry{embdss}
{
name={Imbedding},
description={These usually denote embeddings.},
vars={$\theta, \psi, \phi$},
type={variables}
}

\newglossaryentry{pnts}
{
name={Points},
description={These usually denote points (or regions) in spaces (such as total spaces or base spaces of bundles). Sometimes (especially with additional decoration) they denote maps into these spaces.},
vars={$p,q,r$},
type={variables}
}

\newglossaryentry{pots}
{
name={Posets},
description={These usually denote posets and sometimes spaces.},
vars={$X,Y,Z$},
type={variables}
}

\newglossaryentry{acats}
{
name={cats higher},
description={These usually denote higher categories},
vars={$\sC, \sD, \sE$},
type={variables}
}

\newglossaryentry{cats}
{
name={cats higher},
description={These usually denote ordinary categories},
vars={$\cC,\cD,\cE$},
type={variables}
}

\newglossaryentry{nmbrs}
{
name={Numbers},
description={These usually denote natural numbers, such as objects of singular intervals. Sometimes, they are elements in posets.},
vars={$a,b,c,d$},
type={variables}
}

\let\Hic \ic
\renewcommand*{\ic}{\hyperref[glsentry-TI_gen]{\Hic}}

\let\HTI \TI
\renewcommand*{\TI}{\hyperref[glsentry-TI]{\HTI}}

\let\Higadd \igadd
\renewcommand*{\igadd}[1]{\hyperref[glsentry-adjoin_inv_gen_set]{\Higadd{#1}}}

\renewcommand*{\icg}[1]{\hyperref[glsentry-adjoin_inv_gen_coh_data]{\Hic}^{\hyperref[glsentry-adjoin_inv_gen_coh_data]{#1}}}

\let\Hata\ata
\renewcommand*{\ata}{\hyperref[glsentry-ata_map]{\Hata}}

\renewcommand*{\ccirc}{\hyperref[glsentry-cell_int]{\mathrm{c}}^{\hyperref[glsentry-cell_int]{\circ}}}

\renewcommand*{\cpartial}{\hyperref[glsentry-cell_bound]{\mathrm{c}}^{\hyperref[glsentry-cell_bound]{\partial}}}

\let\HPsi\Psi
\renewcommand*{\Psi}{\hyperref[glsentry-framed_strat]{\HPsi}}

\let\Hcone\cone
\renewcommand*{\cone}{\hyperref[glsentry-cone]{\Hcone}}

\let\HypkP\kP
\renewcommand*{\kP}{\hyperref[glsentry-gen_TP]{\HypkP}}

\let\HypkD\kD
\renewcommand*{\kD}{\hyperref[glsentry-dual_strat]{\HypkD}}

\renewcommand*{\dualdag}{^{\hyperref[glsentry-dual_stratum]{\dagger}}}

\let\HypkM\kM
\renewcommand*{\kM}{\hyperref[glsentry-ass_cat_multiplication]{\HypkM}}

\let\Hip\ip
\renewcommand*{\ip}{\hyperref[glsentry-summit_type]{\Hip}}

\let\Hcolim\colim
\renewcommand*{\colim}{\hyperref[glsentry-colim_pres]{\Hcolim}}

\let\HComp\Comp
\renewcommand*{\Comp}{\hyperref[glsentry-mor_pres_cat]{\HComp}}

\let\Hquotg\quotg
\renewcommand*{\quotg}{\hyperref[glsentry-quot_glob_set]{\Hquotg}}

\let\Htruncleq\truncleq
\renewcommand*{\truncleq}{\hyperref[glsentry-truncleq]{\Htruncleq}}

\let\Hypie\ie
\renewcommand*{\ie}{\hyperref[glsentry-adjoin_eq]{\Hypie}}

\let\Hypgadd\gadd
\renewcommand*{\gadd}[1]{\hyperref[glsentry-adjoin_gen_set]{\Hypgadd{#1}}}

\let\Hypuplus\uplus
\renewcommand*{\uplus}{\hyperref[glsentry-ordered_sum]{\Hypuplus}}

\renewcommand*{\csrc}{\hyperref[glsentry-cubical_src]{\mathsf{src}}_{\hyperref[glsentry-cubical_src]{\square}}}

\renewcommand*{\ctgt}{\hyperref[glsentry-cubical_tgt]{\mathsf{tgt}}_{\hyperref[glsentry-cubical_tgt]{\square}}}

\let\Hypstack\stack
\renewcommand*{\stack}{\hyperref[glsentry-stack]{\Hypstack}} 

\let\Hypstackup\stackup
\renewcommand*{\stackup}{\hyperref[glsentry-stack_emb]{\Hypstackup}} 

\let\Hypstacklow\stacklow
\renewcommand*{\stacklow}{\hyperref[glsentry-stack_emb]{\Hypstacklow}}

\renewcommand*{\stackinc}[1]{\hyperref[glsentry-stack_emb]{\mathsf{In}}^{\hyperref[glsentry-stack_emb]{#1}}}

\let\Hypglue\glue
\renewcommand*{\glue}[1]{\hyperref[glsentry-klvl_stack]{\Hypglue{#1}}}

\let\Hypwhisker\whisker
\renewcommand*{\whisker}[2]{\hyperref[glsentry-whiskering]{\Hypwhisker{#1}{#2}}}

\let\Hypxiso\xiso
\renewcommand*{\xiso}[1]{\hyperref[glsentry-homotopy]{\Hypxiso{#1}}}

\let\HypsW\sW
\renewcommand*{\sW}{\hyperref[glsentry-proj_tow_iso]{\HypsW}}

\let\HypcN\cN
\renewcommand*{\cN}{\hyperref[glsentry-currying]{\HypcN}}

\let\HypcY\cY
\renewcommand*{\cY}{\hyperref[glsentry-uncurrying]{\HypcY}}

\let\HypsE\sE
\renewcommand*{\sE}{\hyperref[glsentry-extended_sing_coll]{\HypsE}}

\let\Hypttimes\ttimes
\renewcommand*{\ttimes}{\hyperref[glsentry-proj_tow]{\Hypttimes}}

\let\HypsL\sL
\renewcommand*{\sL}{\hyperref[glsentry-seq_coll_lim]{\HypsL}}

\let\Hypsing\singcont
\renewcommand*{\singcont}{\hyperref[glsentry-sing_cont]{\Hypsing}}

\let\Hypreg\regcont
\renewcommand*{\regcont}{\hyperref[glsentry-reg_cont]{\Hypreg}}

\renewcommand*{\gsrc}{\hyperref[glsentry-glob_src]{\mathsf{src}}_{\hyperref[glsentry-glob_src]{\circ}}}

\renewcommand*{\gtgt}{\hyperref[glsentry-glob_tgt]{\mathsf{tgt}}_{\hyperref[glsentry-glob_tgt]{\circ}}}

\let\HypcG\cG
\renewcommand*{\cG}{\hyperref[glsentry-elcat_of_globe]{\HypcG}}

\let\Hypelcat\elcat
\renewcommand*{\elcat}{\hyperref[glsentry-cat_of_el]{\Hypelcat}}

\renewcommand*{\lGR}[1]{\hyperref[glsentry-globe_cat_rep]{\HyplG}^{\hyperref[glsentry-globe_cat_rep]{#1}}}

\let\HyplG\lG
\renewcommand*{\lG}{\hyperref[glsentry-globe_cat]{\HyplG}}

\let\Hypglobeset\globset
\renewcommand*{\globset}{\hyperref[glsentry-globeset_cat]{\Hypglobeset}}

\let\Hyptgl\tgl
\renewcommand*{\tgl}{\hyperref[glsentry-term_globe]{\Hyptgl}}

\renewcommand*{\topmon}{_{\hyperref[glsentry-top_monad]{\top}}}

\let\Hypabss\abss
\renewcommand*{\abss}[1]{\hyperref[glsentry-realisation]{\Hypabss{#1}}}

\let\Hypmono\mono
\renewcommand*{\mono}{\hyperref[glsentry-embedding]{\Hypmono}}

\let\Hypinto\into
\renewcommand*{\into}{\hyperref[glsentry-injections]{\Hypinto}}

\let\Hyprestemb\restemb
\renewcommand*{\restemb}{\hyperref[glsentry-rest_emb]{\Hyprestemb}}

\let\Hypmsrc\msrc
\renewcommand*{\msrc}{\hyperref[glsentry-src_sec]{\Hypmsrc}}

\let\Hypmtgt\mtgt
\renewcommand*{\mtgt}{\hyperref[glsentry-tgt_sec]{\Hypmtgt}}

\renewcommand*{\intrel}[1]{\hyperref[glsentry-SI_ass_SSS]{\mathbf{I}[}{#1} \hyperref[glsentry-SI_ass_SSS]{]}}

\let\Hypeta\eta
\renewcommand*{\eta}{\hyperref[glsentry-SI_inj_ass_SSS]{\Hypeta}}

\renewcommand*{\restsec}[1]{_{\hyperref[glsentry-endpoints]{#1}}}

\let\HypsJ\sJ
\renewcommand*{\sJ}{\hyperref[glsentry-emb_from_k_end]{\HypsJ}}

\renewcommand*{\qmin}[1]{\hyperref[glsentry-min_endpoints]{\iq}^{#1}}

\let\Hypsslash\sslash
\renewcommand*{\sslash}{\hyperref[glsentry-slash_notn]{\Hypsslash}}

\let\Hypiota\iota
\renewcommand*{\iota}{\hyperref[glsentry-min_emb]{\Hypiota}}

\renewcommand*{\pbstar}{^{\hyperref[glsentry-pullback]{*}}}

\renewcommand*{\psstar}{_{\hyperref[glsentry-pushforward_notn]{*}}}

\renewcommand*{\postar}{_{\hyperref[glsentry-pushforward_notn]{*}}}

\let\Hypvvec\vvec
\renewcommand*{\vvec}[1]{\hyperref[glsentry-vector]{\Hypvvec{#1}}}

\let\Hypkcoll\kcoll
\renewcommand*{\kcoll}[1]{\hyperref[glsentry-k_lvl_collapse]{\Hypkcoll{#1}}}

\let\Hypstarcoll\starcoll
\renewcommand*{\starcoll}{\hyperref[glsentry-k_lvl_collapse]{\Hypstarcoll}}

\renewcommand*{\NF}[1]{\hyperref[glsentry-normal_form_notn]{[}{#1}\hyperref[glsentry-normal_form_notn]{]}}

\let\Hypnfc\nfc
\renewcommand*{\nfc}{\hyperref[glsentry-normal_form_collapse]{\Hypnfc}}

\let\HypcS\cS
\renewcommand*{\cS}{\hyperref[glsentry-ass_subset_section]{\HypcS}}

\let\HypsS\sS
\renewcommand*{\sS}{\hyperref[glsentry-coll_map_SI]{\HypsS}}

\renewcommand*{\mcoll}{\hyperref[glsentry-coll_notn]{\epi}}
\renewcommand*{\kmcoll}[1]{\hyperref[glsentry-coll_notn]{\epi^{#1}}}

\renewcommand*{\SIvertone}[1]{\ensuremath{\hyperref[glsentry-SI_notn]{\mathbf{SI}}\vslash_{\! #1}}}

\renewcommand*{\SIvert}[2]{\ensuremath{\hyperref[glsentry-SI_notn]{\mathbf{SI}}\vslash^{\hyperref[glsentry-SI_notn]{#1}}_{\! #2}}}

\let\HypsR\sR
\renewcommand*{\sR}{\hyperref[glsentry-repack]{\HypsR}}

\let\HypsV\sV
\renewcommand*{\sV}{\hyperref[glsentry-unpack]{\HypsV}}

\let\HypsU\sU
\renewcommand*{\sU}{\hyperref[glsentry-unpack]{\HypsU}}

\let\HypsG\sG
\renewcommand*{\sG}{\hyperref[glsentry-totalspace]{\HypsG}}

\let\Hyppi\pi
\renewcommand*{\pi}{\hyperref[glsentry-bundleproj]{\Hyppi}}

\renewcommand*{\tsR}[1]{\hyperref[glsentry-krepack]{\HypsR}^{\hyperref[glsentry-krepack]{#1}}}

\let\HypsT\sT
\renewcommand*{\sT}{\hyperref[glsentry-ktower]{\HypsT}}

\renewcommand*{\tpi}[1]{\hyperref[glsentry-kbund]{\Hyppi}^{\hyperref[glsentry-kbund]{#1}}}

\renewcommand*{\tsU}[1]{\hyperref[glsentry-klabel]{\HypsU}^{\hyperref[glsentry-klabel]{#1}}}

\renewcommand*{\tusU}[1]{\und{\hyperref[glsentry-klabel]{\HypsU}}^{\hyperref[glsentry-klabel]{#1}}}

\renewcommand*{\tsG}[1]{\hyperref[glsentry-ktotalspace]{\HypsG}^{\hyperref[glsentry-ktotalspace]{#1}}}

\let\HypSIs\SIs
\renewcommand*{\SIs}[1]{\hyperref[glsentry-SIvertsec]{\HypSIs{#1}}}

\let\HypSIf\SIf
\renewcommand*{\SIf}[1]{\hyperref[glsentry-SIvertforget]{\HypSIf{#1}}}

\let\Hypund\und
\renewcommand*{\und}[1]{\hyperref[glsentry-SIvertforget]{\Hypund{#1}}}

\let\Hypodot\odot
\renewcommand*{\odot}{\hyperref[glsentry-horcomp]{\Hypodot}}

\let\Hypnu\nu
\renewcommand*{\nu}{\hyperref[glsentry-funnat]{\Hypnu}}

\renewcommand*{\succindex}[1]{\hyperref[glsentry-fillind]{\mathtt{I}}({#1})}

\renewcommand*{\succfill}[1]{\hyperref[glsentry-filled]{\mathtt{F}}({#1})}

\renewcommand*{\succ}[1]{\hyperref[glsentry-succed]{\mathtt{S}}({#1})}
\renewcommand*{\succn}[2]{\hyperref[glsentry-succed]{\mathtt{S}}^{#2}({#1})}

\renewcommand*{\pred}[1]{\hyperref[glsentry-preded]{\mathtt{P}}({#1})}
\renewcommand*{\predn}[2]{\hyperref[glsentry-preded]{\mathtt{P}}^{#2}({#1})}

\renewcommand*{\ssoe}{_{\hyperref[glsentry-soed]{\mathtt{s}}}}

\renewcommand*{\ttae}{_{\hyperref[glsentry-taed]{\mathtt{t}}}}

\let\Hypavg\avg
\renewcommand*{\avg}[1]{\hyperref[glsentry-edgenorm]{\Hypavg{#1}}}

\let\HypSiR\SiR
\renewcommand*{\SiR}{\hyperref[glsentry-prelreal]{\HypSiR}}

\renewcommand*{\singop}{^{\hyperref[glsentry-singop]{\mathsf{sing}}}}

\renewcommand*{\regop}{^{\hyperref[glsentry-regop]{\mathsf{reg}}}}

\let\Hypwwidehat\wwidehat
\renewcommand*{\wwidehat}[1]{\hyperref[glsentry-extsingmor]{\Hypwwidehat{#1}}}

\let\Hypsingint\singint
\renewcommand*{\singint}[1]{\hyperref[glsentry-singint]{\Hypsingint{#1}}}

\let\Hypextsing\extsing
\renewcommand*{\extsing}{\hyperref[glsentry-extsing]{\Hypextsing}}

\let\HypBDelta\BDelta
\renewcommand*{\BDelta}{\hyperref[glsentry-totalorders]{\HypBDelta}}

\let\Hypdelta\delta
\renewcommand*{\delta}{\hyperref[glsentry-face]{\Hypdelta}}

\let\Hypsigma\sigma
\renewcommand*{\sigma}{\hyperref[glsentry-degeneracy]{\Hypsigma}}

\let\Hypbnum\bnum
\renewcommand*{\bnum}[1]{\hyperref[glsentry-bnum]{\Hypbnum{#1}}}

\let\Hypbang\bang
\renewcommand*{\bang}{\hyperref[glsentry-bang]{\Hypbang}}

\let\HypPoss\Poss
\renewcommand*{\Poss}[1]{\hyperref[glsentry-posslash]{\HypPoss{#1}}}

\let\Hyprest\rest
\renewcommand*{\rest}[2]{\Hyprest{#1}{\hyperref[glsentry-rest2]{#2}}}

\let\HypDelta\Delta
\renewcommand*{\Delta}{\hyperref[glsentry-testfunc]{\HypDelta}}

\let\Hypsecp\secp
\renewcommand*{\secp}[1]{\hyperref[glsentry-secp]{\Hypsecp{#1}}}

\let\Hypppi\ppi
\renewcommand*{\ppi}{\hyperref[glsentry-prodproj]{\Hypppi}}

\let\HypiP\iP
\renewcommand*{\iP}{\hyperref[glsentry-settruth]{\HypiP}}

\renewcommand*{\converse}{^{\hyperref[glsentry-converse]{\top}}}

\let\HypPRel\PRel
\renewcommand*{\PRel}{\hyperref[glsentry-PRel]{\HypPRel}}

\let\HypRel\Rel
\renewcommand*{\Rel}{\hyperref[glsentry-Rel]{\HypRel}}

\let\Hypobj\obj
\renewcommand*{\obj}{\hyperref[glsentry-obj]{\Hypobj}}

\let\Hypdiscr\discr
\renewcommand*{\discr}{\hyperref[glsentry-discr_cat]{\Hypdiscr}}

\let\HypDiscr\Discr
\renewcommand*{\Discr}{\hyperref[glsentry-discr_prof]{\HypDiscr}} 

\let\HypBunbc\Bunbc
\renewcommand*{\Bunbc}{\hyperref[glsentry-bunbc_cat]{\HypBunbc}} 

\renewcommand*{\Buncoll}[1]{\ensuremath{\hyperref[glsentry-bun_coll_cat]{\mathbf{Bun}}^{\rightsquigarrow}_{#1}}}

\let\Hypctypsum\ctypsum
\renewcommand*{\ctypsum}{\hyperref[glsentry-region_dim]{\Hypctypsum}}

\let\Hypctyp\ctyp
\renewcommand*{\ctyp}{\hyperref[glsentry-region_type]{\Hypctyp}}

\let\HypId\Id
\renewcommand*{\Id}{\hyperref[glsentry-identity_cube]{\HypId}}

\let\Hypsing\sing
\renewcommand*{\sing}{\hyperref[glsentry-sing]{\Hypsing}}

\let\Hypreg\reg
\renewcommand*{\reg}{\hyperref[glsentry-reg]{\Hypreg}}

\let\Hypedgeset\edgeset
\renewcommand*{\edgeset}{\hyperref[glsentry-edgeset]{\Hypedgeset}}

\renewcommand*{\SI}{\ensuremath{\hyperref[glsentry-SI]{\mathbf{SI}}}}

\renewcommand*{\oSI}{_{\hyperref[glsentry-singmor]{\mathbf{SI}}}}

\renewcommand*{\monosing}{_{\hyperref[glsentry-underlying_mono]{\mathsf{sing}}}}

\let\HypPos\Pos
\renewcommand*{\Pos}{\hyperref[glsentry-Poscat]{\HypPos}}

\let\HypBool\Bool
\renewcommand*{\Bool}{\hyperref[glsentry-Boolcat]{\HypBool}}

\let\HypMap\Map
\renewcommand*{\Map}{\hyperref[glsentry-function_set]{\HypMap}}

\let\HypSetCat\SetCat
\renewcommand*{\SetCat}{\hyperref[glsentry-Setcat]{\HypSetCat}}

\let\HypFun\Fun
\renewcommand*{\Fun}{\hyperref[glsentry-functor_category]{\HypFun}}

\let\HypCat\Cat
\renewcommand*{\Cat}{\hyperref[glsentry-Catcat]{\HypCat}}

\let\Hypconst\const
\renewcommand*{\const}{\hyperref[glsentry-constant_functor]{\Hypconst}}

\let\Hypstfwd\stfwd
\renewcommand*{\stfwd}{\hyperref[sec:notes]{\Hypstfwd}}

\let\HPSk\PSk
\renewcommand*{\PSk}{\hyperref[glsentry-poset_skel]{\HPSk}}

\let\Helpsk\elpsk
\renewcommand*{\elpsk}{\hyperref[glsentry-poset_skel_label]{\Helpsk}}

\let\HCubeo\Cubeo
\renewcommand*{\Cubeo}[2]{\hyperref[glsentry-cubeo_cat]{\HCubeo{#1}{#2}}}

\let\HBuno\Buno
\renewcommand*{\Buno}[2]{\hyperref[glsentry-buno_cat]{\HBuno{#1}{#2}}}

\let\Hsimeq\simeq
\renewcommand*{\simeq}{\hyperref[glsentry-cube_equiv]{\Hsimeq}}

\let\Hnorm\norm
\renewcommand*{\norm}[1]{\hyperref[glsentry-geom_real]{\Hnorm{#1}}}

\let\HredGamma\redGamma
\renewcommand*{\redGamma}[1]{\hyperref[glsentry-min_label]{\HredGamma{#1}}}

\renewcommand*{\pCat}{\hyperref[glsentry-pres_cat_n]{\mathbf{Cat}}^{\hyperref[glsentry-pres_cat_n]{\mathrm{pres}}}}

\let\HpwCat\pwCat
\renewcommand*{\pwCat}{\hyperref[glsentry-pres_wcat]{\HpwCat}}

\let\Hinfeq\infeq
\renewcommand*{\infeq}[2]{\hyperref[glsentry-inf_eq_res]{\Hinfeq{#1}{#2}}}

\let\HkC\kC
\renewcommand*{\kC}{\hyperref[glsentry-glob_to_pres]{\HkC}}

\let\HkiC\kiC
\renewcommand*{\kiC}[1]{\hyperref[glsentry-alg_resol]{\HkiC{#1}}}

\let\Hinterchanger\interchanger
\renewcommand*{\interchanger}{\hyperref[glsentry-interchanger]{\Hinterchanger}}

\let\HGComp\GComp
\renewcommand*{\GComp}{\hyperref[glsentry-gen_comp]{\HGComp}}

\let\HGComps\GComps
\renewcommand*{\GComps}{\hyperref[glsentry-gen_comp]{\HGComps}}

\let\HGGamma\GGamma
\renewcommand*{\GGamma}[2]{\hyperref[glsentry-ggamma_tp]{\HGGamma{#1}{#2}}}

\begin{document}

\baselineskip=16pt plus2pt

\setcounter{secnumdepth}{3}
\setcounter{tocdepth}{3}

\maketitle                  %

\begin{dedication}
to my parents
\end{dedication}

\begin{abstractlong} \small{
We define novel fully combinatorial models of higher categories. Our definitions are based on a connection of higher categories to ``directed spaces". Directed spaces will be locally modelled on manifold diagrams, which are stratifications of the $n$-cube such that strata are transversal to the flag foliation of the $n$-cube. The first part of this thesis develops a combinatorial language for manifold diagrams called singular $n$-cubes. In the second part we apply this language to build our notions of higher categories.

Singular $n$-cubes can be thought of as ``flag-foliation-compatible" stratifications of the $n$-cube, such that strata are ``stable" under projections from the $(k+1)$- to the $k$-cube, together with a functorial assignment of data to strata. The definition of singular $n$-cubes is inductive, with $(n+1)$-cubes being defined as combinatorial bundles of $n$-cubes over the (stratified) interval. The combinatorial structure of singular $n$-cubes can be naturally organised into two categories: $\SIvert n \cC$, whose morphisms are bundles themselves, and $\Bunbc^n_\cC$, whose morphisms are inductively defined as base changes of bundles. The former category is used for the inductive construction of singular $n$-cubes. The latter category describes the following interactions of these cubes. There is a subcategory of ``open" base changes, which topologically correspond to open maps of bundles. We show this subcategory admits an (epi,mono) factorisation system. Monomorphism will be called embeddings and describe how cubes can be embedded in one another such that strata are preserved. Epimorphisms will be called collapses and describe how strata can be can be refined. Two cubes are equivalent if there is a cube that they both refine. We prove that each  ``equivalence class" (that is, the connected component of the subcategory generated by epimorphisms) has a terminal object, called the collapse normal form. Geometrically speaking the existence of collapse normal forms translates into saying that any combinatorially represented manifold diagram has a unique coarsest stratification, making the equality relation between manifold diagrams decidable and computer implementable. 

As the main application of the resulting combinatorial framework for manifold diagrams, we give algebraic definitions of various notions of higher categories. In particular, we define associative $n$-categories, \free{} associative $n$-categories and \free{} associative $n$-groupoids. The first depends on a theory of sets, while the latter two do not, making them a step towards a framework for working with general higher categories in a foundation-independent way. All three notions will have strict units and associators. The only ``weak" coherences which are present will be called homotopies. We propose that this is the right conceptual categorisation of coherence data: homotopies are essential coherences, while all other coherences can be uniformly derived from them. As evidence to this claim we define \free{} weak $n$-categories, and develop a mechanism for recovering the usual coherence data of weak $n$-categories, such as associators and pentagonators and their higher analogues. This motivates the conjecture that the theory of associative higher categories is equivalent to its fully weak counterpart. }
\end{abstractlong}
\newpage

\begin{acknowledgementslong}
I feel equally grateful, happy and surprised that this thesis came into being. When the direction of research first manifested itself at the end of my second year in Oxford, it was far from clear whether anything concrete would ever come out of it. The ambition behind this work, of finding ``right" definitions, can be dangerous to follow as it is often not decidable what ``right" means by itself, without the context of concrete problems and applications. However, the burden of decision is lifted if one is in the lucky situation of there being only one viable choice. In several ways, I felt I was in a lucky situation many times in the past years, and it was the many supportive and thoughtful friends and colleagues I was surrounded by who placed me in (and sometimes pushed me into) those situations. If it wasn't for them this work would evidently not exist.

A most central force behind the mathematically lucky situations during this project was without any doubt my supervisor, Christopher Douglas. There is (at least!) one parallel between being a supervisor and being a psychotherapist, namely, sometimes asking the right questions leads to the most progress in a patient's understanding. If anyone has perfected this skill, then I claim it is Chris. Both in mathematics, and in its human aspects, Chris has been a wonderful teacher. There were many other personal obstacles to overcome during my time in Oxford, some of which were not so easily resolved. In those situations Chris has been more supportive than I could have asked for, and more significantly, he has been a true friend to me. I was further in the lucky position of benefiting from the highly supportive voice of Samson Abramsky, under whose co-supervision I started my studies at Oxford. Samson's insights and (seemingly inexhaustible) wisdom had a lasting impact on my intellectual journey. I would like to thank both Chris and Samson for their invaluable guidance and for encouraging me to trust (more) in my own ideas. 

My growth as a mathematician and more generally as a human being is the result of the cumulative efforts of many wonderful people that I connected with at ETH, Cambridge and Oxford. Mathematically I would like to thank Stefano Gogioso, Amar Hadzihasanovic, Kohei Kishida, David Reutter, Matthjis Vakar and Jamie Vicary for many enlightening discussions. Andre Henriques and Mike Shulman had to carefully read this thesis---and they did a fantastic job. The union of their suggestions improved the presentation dramatically, saving me from bending the truth in several places. Michael Eichmair introduced me to the world of geometry during my undergraduate studies, and simultaneously conveyed an amount of enthusiasm for mathematics that inspired me until the end of my DPhil and beyond. I would also like to thank Julia Goedecke for introducing me to Category Theory at Cambridge and supervising my Part III essay, which would later on become one of the foundations for my thesis research. I further want to thank Martin Bendersky, Noson Yanofsky and Mahmoud Zeinalian for being very welcoming and gracious hosts during my time in New York, allowing me give many talks which benefited my own understanding of the topics immensely.

I gratefully acknowledge the funding I received from the Engineering and Physical Sciences Research Council, the Department of Computer Science and St Catherine's College, which I am proud to call my home at Oxford. I would like to thank the administrative staff at my college and my departments, particularly Keith Gillow, who provided the template for this thesis.
\thispagestyle{empty}
\end{acknowledgementslong}
\newpage

\begin{versionnotes}
\begin{center}
\vspace{-\baselineskip}
		\small\emph{added March 2023}
\vspace{2\baselineskip}
\end{center}
		This version of my thesis ``Associative $n$-categories'' is essentially identical to the version publicly available at the Oxford Research Archives (modulo minor changes in formatting which probably resulted from using a newer version of \texttt{pdflatex}). The latter document can be accessed via the following link:
	\begin{center}
			\url{https://ora.ox.ac.uk/objects/uuid:9f9d8668-a48a-4b30-9883-b11e4edbfc54} .
	\end{center}
It may further be helpful to point out that manifold diagrams continue to be an active area of research: by now, there have been substantial further developments in the area, in particular relating to the construction of stratified-geometric models of manifold diagrams. In many ways, this subsequent work improved, generalized and superseeded much of the material presented here.
\end{versionnotes}

\begin{romanpages}          %
\tableofcontents            %
\end{romanpages}            %

\cleardoublepage
\phantomsection
\addcontentsline{toc}{chapter}{Introduction}
\setcounter{section}{0}
\renewcommand{\thesection}{I.{\arabic{section}}}
\renewcommand{\theHsection}{I.{\arabic{section}}}

\chapter*{Introduction} \label{ch:introduction}

In \autoref{sec:overview} we give an overview of the contents of this thesis and briefly discuss the larger research programme and context that it is part of. In \autoref{sec:coherences} we give an informal introduction to higher category theory and coherences from two different perspectives that are relevant to this work. In the remainder of the chapter we discuss results, conjectures and related work. 

For the reader not interested in the philosophical context of this work, it will suffice to read \autoref{ssec:fast_facts} (listing the main results) and \autoref{sec:notes} (containing notes to the reader) before proceeding to the next chapter.

\section{Overview} \label{sec:overview}

This thesis is a first step towards a generic framework for working with general higher categories in a foundation-independent way. The work is based on building a novel bridge between algebra and geometry.

Higher categories can be regarded as a generalisation of topological spaces. While in spaces any path (or higher homotopy of such) can be travelled along in two directions, in higher categories this need not be the case. This is relevant for the description of irreversible processes, which are abundant in the mathematical and physical world. From this \textit{geometric} perspective, higher categories can thus be approached by formulating them as ``directed spaces". 

Higher categories can also be regarded as certain \textit{algebraic} objects generalising ordinary categories. The geometric and algebraic perspectives on higher categories should be equivalent in some way. In the special case of spaces, this equivalence is the content of the \textit{homotopy hypothesis}, which states that the theory of spaces is (in a certain sense) equivalent to the theory of so-called higher groupoids.

Neither ``directed spaces" nor ``algebraic higher categories" are straight-forward to define, and no single, canonical definition for them exists. For the former this is because the \textit{global} geometric models of spaces don not easily generalise to the directed case. For the latter this is because many algebraic relations in the form of so-called ``coherences" have to be accounted for in higher category theory.

This thesis develops novel algebraic approaches to higher categories by exploiting the combinatorial structure underlying a certain \textit{local} geometric model of directed spaces. The novelty and relevance of our geometric model are in turn based on the existence of such an elementary and canonical combinatorial description (amenable, for instance, to computer implementation) and the natural connection of algebra and geometry that this entails. The model also provides novel insights into the algebra of higher-categorical coherences, and motivates that coherences fall into one of two classes of either ``essential" coherences (those naturally present in the geometry) or ``derived" coherences (those which can be obtained by composition operations on the geometry). This will motivate our semi-strict (``associative") approach to higher categories.

The local geometric model mentioned above is given by so-called \textit{manifold diagrams}. They are geometrically ``dual" to the usual cellular perspective on directed spaces. Manifold diagrams also appear to have a natural connection to the classical theory of spaces via a generalised Thom-Pontryagin construction as developed e.g. by Buonchristiano, Rourke and Sanderson \cite{buonchristiano1976geometric}. This in turn points towards promising approaches to proving the aforementioned homotopy hypothesis for associative higher categories and will be further discussed in an appendix.

The combinatorial structures that we use to capture manifold diagrams will be called  \textit{singular $n$-cubes}. Most of this thesis will be dedicated to the study of their rich combinatorial theory. As an immediate application we will be able to define various algebraic notions of higher categories. The central notions will be that of \textit{associative $n$-categories}.

The models of higher categories that we introduce using manifold diagrams come in the following flavours: besides associative $n$-categories, we also define fully weak $n$-categories, presented associative $n$-categories, presented associative $n$-fold categories and presented fully weak $n$-categories. We briefly give an overview of the most important features of these flavours.

Firstly, the ``fully weak" flavour is indeed different from the ``associative" flavour. The former contains all coherence data, which means ``strict" equalities are being replaced by ``weak" isomorphisms as much as possible. The latter has strict \textit{units} and strict (higher) \textit{associators}. Its weakness only lies in coherences called \textit{homotopies}. We propose that this is the right conceptual categorisation of coherence data: homotopies are essential coherences, while all other coherences (including weak identities and associators) can be uniformly derived from them by introducing an explicit invertiable ``composition operation". Based on this we will argue that the associative and the fully weak flavour are equivalent. This equivalence can be regarded as a strengthening of a conjecture by Simpson \cite{simpson1998homotopy}. Notably, associative $n$-categories generalise Gray-categories (for $n = 3$) which were introduced by Gordon, Power and Street \cite{gordon1995coherence}, showing that their choice of definition was not arbitrary but follows general geometric principles. 

Secondly, the ``associative" and the ``presented associative" flavour are different. Associative $n$-categories are globular sets with composition operations satisfying equations. Presented associative $n$-categories, however, are given by \textit{presentations}---that is, lists of generating $k$-morphisms ($0\leq k \leq n+1$). Type-theoretically, the latter correspond to constructors of higher inductive types. Indeed, specialising to presented associative $\infty$-groupoids, one can translate $k$-cells of a CW-complex into generating $k$-morphisms of an $\infty$-groupoid, and this is analogous to translating $k$-cells to constructors of higher inductive types in Homotopy Type Theory \cite{univalent2013homotopy}. Subject to completing the proof of the homotopy hypothesis, this has many immediate applications. For instance, representatives of elements in the homotopy groups of spaces can be ``algorithmically listed" using the language developed in this thesis.

\subsection{Research programme}

The work in this thesis, which includes the fully algebraic definition of presented associative $n$-categories, is a step towards the definition of a generic framework for $\infty$-theories, or more precisely, $n$-theories for arbitrary $n$. 

An $n$-theory \cite{shulman2018ntheory} is a foundational language describing the rules governing a class of $(n-1)$-categorical objects (for instance, if $(n - 1) = 0$ then these objects are sets possibly with additional structure). However, classes of $n$-theories are described themselves by $(n+1)$-theories (For instance, the $1$-theory of monoids describes the class of monoids, but is itself e.g. part of the class of categories with finite products described by an according $2$-theory). This inductive raise of dimension leads to difficulties in letting an $n$-theory define or describe itself easily. 

Ideally however, we would want mathematical foundations to be self-descriptive.  Based on the above it is therefore natural to let $n \to \infty$ and ask whether there is a ``uniform" theory of $n$-theories for all $n$. In this work we define objects of such a theory in the form of presented associative $\infty$-categories, and a fuller description of the theory (including $k$-morphisms, for $k \geq 1$, and a small universe type) is work in progress. 

The following two subsections provide informal explanations of the ideas involved in $n$-theories and their relation to classical Set-theoretical thinking in mathematics \textit{without} the use of type-theoretical terminology. The remaining subsections will set the stage for the definition of presented associative $n$-categories. 

In \autoref{ssec:po_alg_vs_geo} we discuss (``foundation-dependent") Set-theoretical thinking in the context of more general (``foundation-independent") type-theoretical languages. In \autoref{ssec:po_alg_geo_higher_cat}, we transfer this discussion of Set Theory vs. Type Theory, or more concretely, Geometry vs. Algebra, to the context of higher category theory. In \autoref{sec:intro_mfld} we introduce our geometric model for (morphisms) in higher categories by means of examples, and sketch its connection to classical topology in \autoref{rmk:GTP}. In \autoref{ssec:po_alg_of_diag} we mention our approach to ``algebraisation" of manifold diagrams. Finally, in \autoref{ssec:po_anc_def} we can already sketch a definition of presented associative $\infty$-category.

\subsection{Foundation-dependent vs. foundation-independent models} \label{ssec:po_alg_vs_geo}

While spaces are ``less general" than higher categories, they have historically played the much more important role in mathematics. For instance, spaces appear as vector \textit{spaces}, phase \textit{spaces}, moduli \textit{spaces} etc. in many areas of research in mathematics and physics. One reason for this might be that spaces have an intuitive \textit{geometric model} as sets with (geometric) structure, but an analogous intuitive geometric definition has been difficult to give for higher categories. However, despite its great versatility, a drawback of the geometric approach to spaces is that sets underlying spaces such as the linear continuum can be ``large" and not directly amenable to finite axiomatisation or implementation on a computer with finite capacity. 

We briefly broaden our perspective, commenting on ideas in mathematical foundations relevant to this research. There are two approaches to phrasing a concept in mathematical language. The first of them is to encode a \textit{model} of the concept \textit{within} an existing foundational language (we assume this to be Set Theory for now). One can then try to derive new properties of the concept using that language, and in particular verify whether the concept's expected behaviour is modelled fully and faithfully. This approach could be called ``building \textit{set-theoretic models}", and if there is reference to the theory of spaces, we will also refer to it as ``building  \textit{geometric models}". Since the sets involved are often large, their manipulation can require some intuition and mathematical ingenuity. In contrast, the second approach proceeds as follows. One can alternatively try to find a finite collection of ``universal" rules that directly govern the concept's behaviour, possibly extending to other concepts' behaviour if there are relevant interactions. These rules can then form the grammar of an independent foundational language describing the concept. The advantage of this second approach is that it might capture the nature of the concept more efficiently, allowing for clearer, more systematic and possibly computer implementable proofs. A drawback is that it might be difficult to capture \textit{all} aspects of the set-theoretical model (cf. \cite{shulman2018brouwer}). We suggest to call this second approach a \textit{fully algebraic} model of the concept, since rules ``compose" to form proofs just as elements of an algebraic structure compose. Importantly, in actual mathematical work we often find a combination of both approaches, which leads to a spectrum of ``partially algebraic" models: that is, part of the concept's behaviour might be described in the context of other languages. For instance, this is the case for symbolic manipulations in Algebra (such as differentiation) which are governed by ``independent" (and computer implementable) rule sets while still being situated\footnote{The point here is \textit{not} to doubt that Set Theory can be used to formulate such algebraic computations, but just to point out that they can also be performed independently of such powerful foundational languages (for instance, in a minimal type theory capturing the necessary rules).} in the larger context of Set Theory.

For a long time, some topologists had the impression that the \textit{concept of spaces} might have its own fully algebraic formulation as well, independent of its geometric model. A first indication of this feeling might have been Whitehead's ``algebraic homotopy" programme presented at his 1950 ICM talk (quoted in \cite{baues1989algebraic}). Steps towards partial solutions were consequently taken by Kan \cite{kan1955abstract} and Quillen \cite{quillen2006homotopical} among others. In the 1980s Grothendieck took an important conceptual step towards fully realising the vision, hypothesising that the homotopy theory of spaces is in fact equivalent to a theory of algebraic structures called higher groupoids \cite{grothendieck1983pursuing}, based on translating a space into an $\infty$-groupoid by the so-called \textit{fundamental $\infty$-groupoid} construction. This was later on called the homotopy hypothesis by Baez \cite{baez2007homotopy}. A central hurdle towards a fully algebraic formulation, however, remained, namely, a sufficiently nice algebraic description of higher groupoids was missing. That is, until recently: the quest for a foundational (and in particular, computer-implementable) language of higher groupoids was given its first solution in 2013 in the form of homotopy type theory (HoTT), developed in a large collaborative effort, the Univalent Foundations Program at the IAS \cite{univalent2013homotopy}.

\subsection{Geometry and algebra of higher categories} \label{ssec:po_alg_geo_higher_cat}

A fully algebraic language for spaces, however, is only the beginning of a larger story. Namely, the translation of the geometric model of spaces into higher groupoids is a special case of the translation between geometric and algebraic models of higher categories. Hypothetically, the latter (that is, a fully algebraic formulation of higher categories) would be a powerful foundational language for all of mathematics. Indeed, being able to express higher directed paths allows higher category theory to formulate so-called universal properties of the most general kind, and in general, types of concepts with these properties cannot be easily formulated, for instance, in HoTT. Such a fully algebraic formulation of higher category theory has not yet been found. In fact, not only does the fully algebraic approach to higher categories appear difficult, but there also doesn't seem to be a canonical set-theoretical model, leading to a large variety of approaches. A first ``fully" geometric model was recently constructed by Ayala, Francis and Rozenblyum \cite{ayala2018factorization}. Many partially algebraic approaches have been given, some of which more on the geometric side (e.g. by Barwick  \cite{barwick2005infinity} and Rezk \cite{rezk2010cartesian}), and others more on the algebraic side (e.g. by Batanin \cite{batanin1998monoidal}, Leinster  \cite{leinster-operads} and Maltsiniotis \cite{maltsiniotis2010grothendieck}). The latter gave rise to a fully algebraic sub-theory of higher categories by Finster and Mimram \cite{finster2017type}.

None of the above models provides a translation between the geometric and fully algebraic world. We argue that this, however, is highly desirable because of the following: having a faithful geometric model of an algebraic language gives a \textit{correctness} test for the algebraic language. Conversely, having an algebraic description of a geometric model guarantees that this model is constrained with respect to possible pathologies (which are often uncontrollable in a Set-theoretic setting), which ultimately also supports its correctness. In summary, we can ask the following question:

\begin{center}
\textit{Is there a theory of higher categories that naturally \\* has both a geometric and an algebraic description?}
\end{center}

In this thesis we derive fully and partially algebraic models of higher categories based on geometric principles and argue that these models give a positive answer to the preceding question, with two caveats: firstly, we defer a rigorous proof of equivalence of the geometric and algebraic description to future work; secondly, we do not give an account of the full theory of higher categories (which includes functors, natural transformation, etc.) but only define the objects of that theory, that is, higher categories. 

The geometric description of (morphisms in) our algebraic higher categories will be modelled by so-called \textit{manifold diagrams}. 

\subsection{Examples of manifold diagrams} \label{sec:intro_mfld}

Manifold diagrams are a useful, but not widely used geometric tool for approaching higher category theory. The more common approach to higher categories is via cellular diagrams. Roughly speaking, cellular diagrams and manifold diagrams are Poincar\'e dual. Their relation will be explained in more detail in \autoref{sec:sum_CW} (and \autoref{sec:CWcomplexes}) in the case of higher groupoids, to which cellular diagrams apply the best.

For the reader not familiar with manifold diagrams, in this section we will briefly illustrate the translation between cellular diagrams and manifold diagrams by means of examples. We leave these examples without much explanation, but the reader should note that in each of these examples the following principles hold:
\begin{enumerate}
\item Cellular diagrams (representing compositions of morphisms in $n$-categories) look like $n$-skeleta of certain cell complexes. $k$-cells represent $k$-morphisms. Cells are given a direction from ``input" to ``output" parts of their boundary, which is indicated by an arrow.
\item Manifold diagrams (representing compositions of morphisms in $n$-categories) consist of certain manifolds embedded in the $n$-cube. $k$-morphisms are represented by $(n-k)$-dimensional manifolds. The distinction into ``input" and ``output" of adjacent (higher-dimensional) manifolds is determined by the $n$ directions of the $n$-cube, so no arrows are required.
\item If a cell $f$ of a cellular diagram contains another cell $A$ in its boundary (that is, either in its source or its target), then the submanifold $A$ in the corresponding manifold diagrams contains the submanifold $f$ in its boundary.
\end{enumerate}
We give examples up to dimension $n = 3$. In this case, $n$-manifold diagrams are also known as point diagrams ($n = 1$), string diagrams ($n = 2$) and surface diagrams ($n = 3$). 

\begin{enumerate}
\item In the case of point diagrams consider the following four examples
\begin{restoretext}
\begingroup\sbox0{\includegraphics{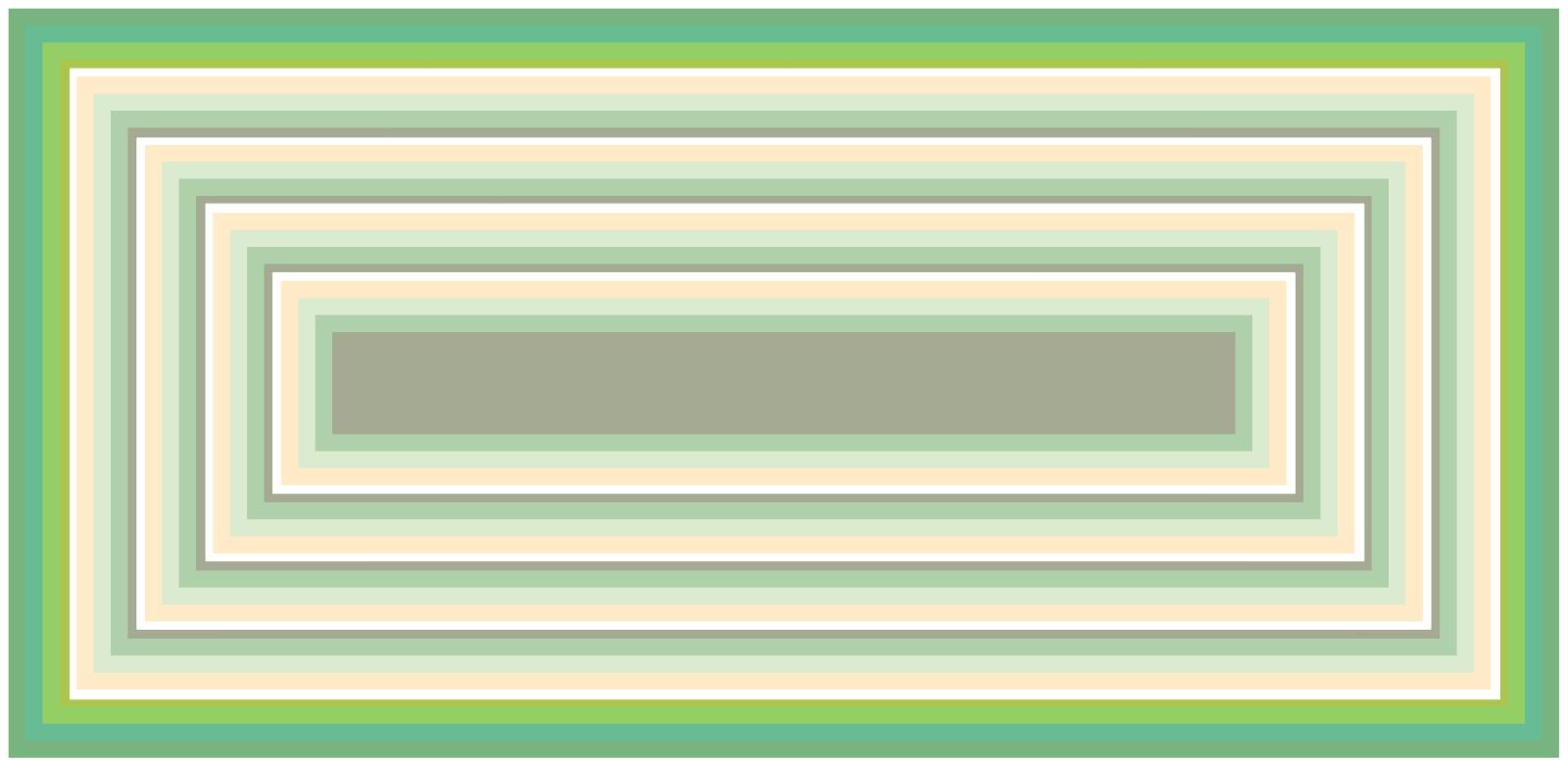}}\includegraphics[clip,trim={.2\ht0} {.15\ht0} {.2\ht0} {.1\ht0} ,width=\textwidth]{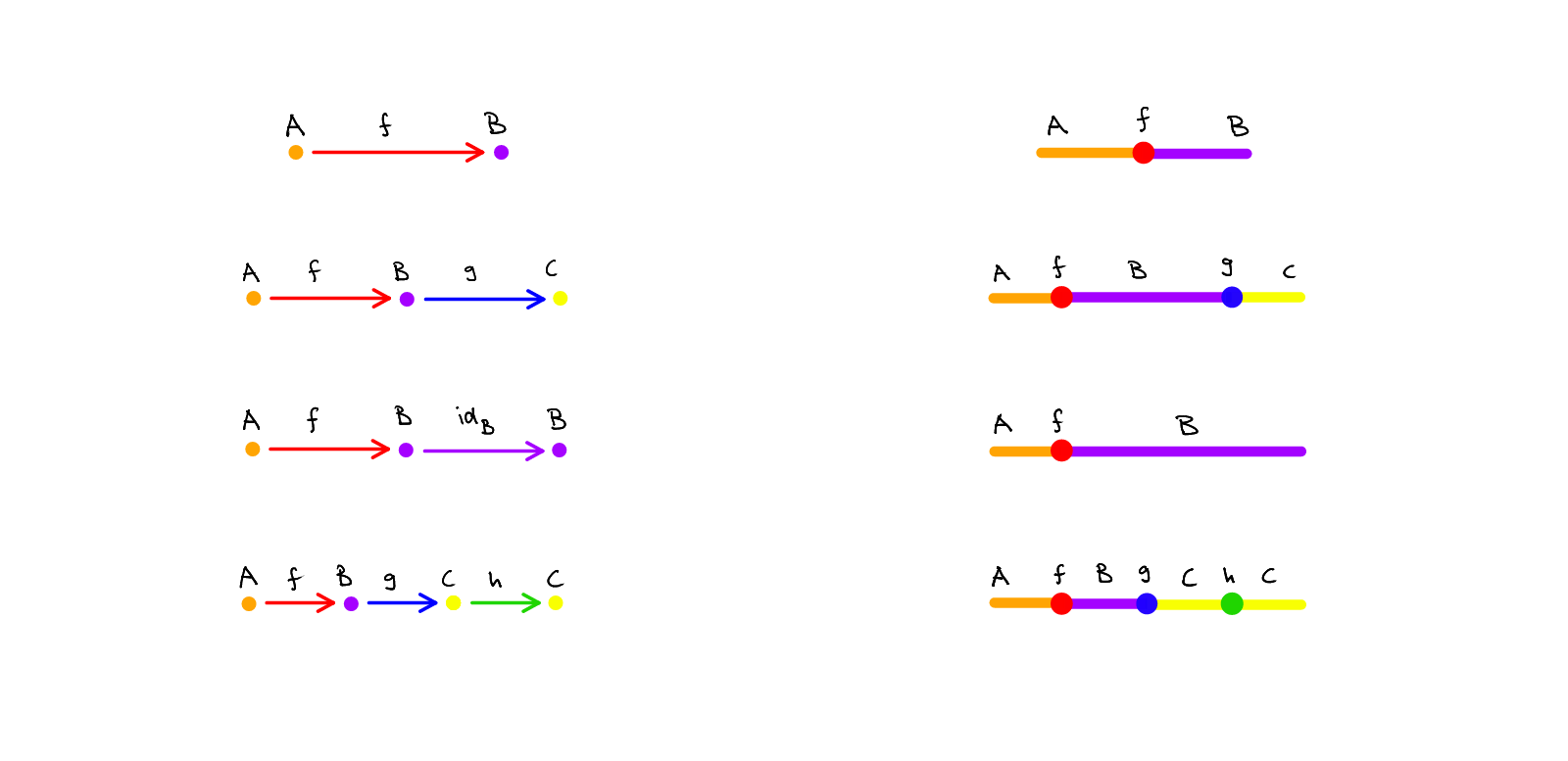}
\endgroup\end{restoretext}
The usual cellular notation for ordinary categories in the left columns is translated to point diagrams in the right column. Note that the third example shows how identity cells can be naturally interpreted by ``absence" of any $0$-manifolds.

\item For string diagrams we first mention the following two translations
\begin{restoretext}
\begingroup\sbox0{\includegraphics{test/page1.png}}\includegraphics[clip,trim={.2\ht0} {.05\ht0} {.2\ht0} {.05\ht0} ,width=\textwidth]{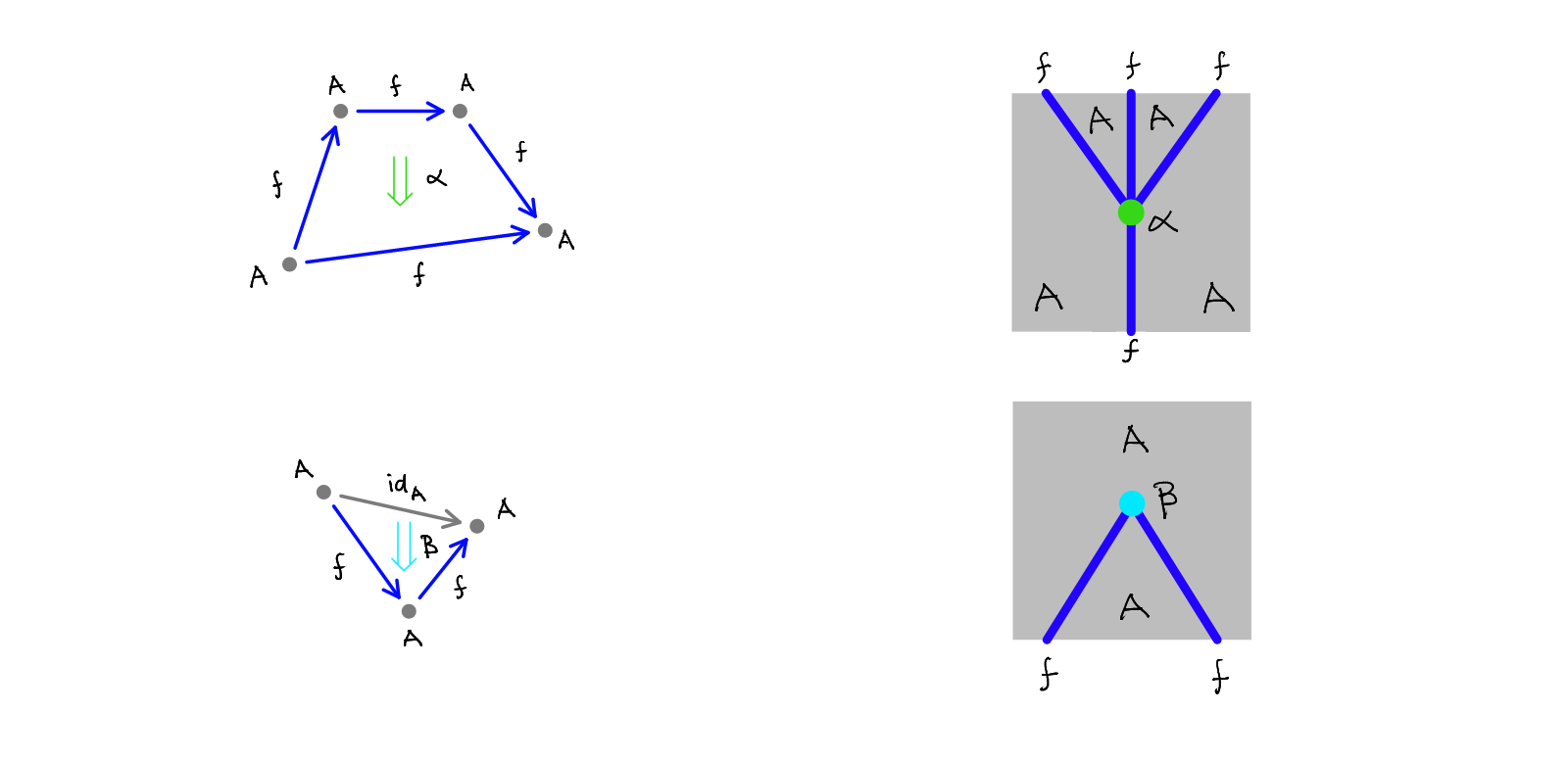}
\endgroup\end{restoretext}
In each case a single $2$-cell on the left is translated into a $0$-dimensional manifold on the right. Further, $1$-cells on the left become $1$-manifolds on the right, and $0$-cells turn into $2$-manifolds. A more complex example is then the following
\begin{restoretext}
\begingroup\sbox0{\includegraphics{test/page1.png}}\includegraphics[clip,trim={.2\ht0} {.2\ht0} {.2\ht0} {.2\ht0} ,width=\textwidth]{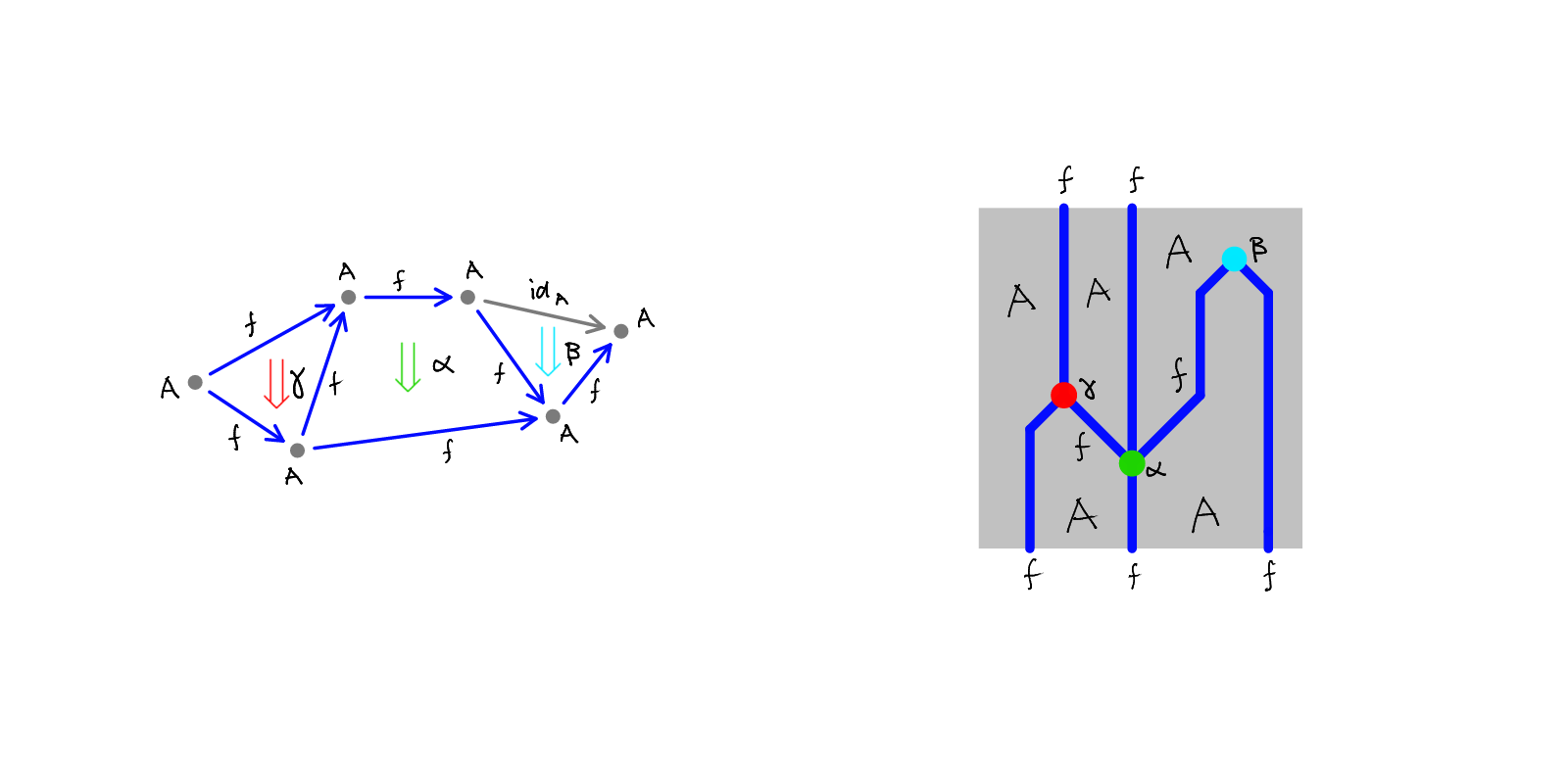}
\endgroup\end{restoretext}
Note the $\alpha$ and $\beta$ are the $2$-cells that we just defined individually.

Finally, to highlight our interpretation of identities again, consider the following translation from a composite of $2$-cells into a string diagram
\begin{restoretext}
\begingroup\sbox0{\includegraphics{test/page1.png}}\includegraphics[clip,trim={.2\ht0} {.3\ht0} {.2\ht0} {.2\ht0} ,width=\textwidth]{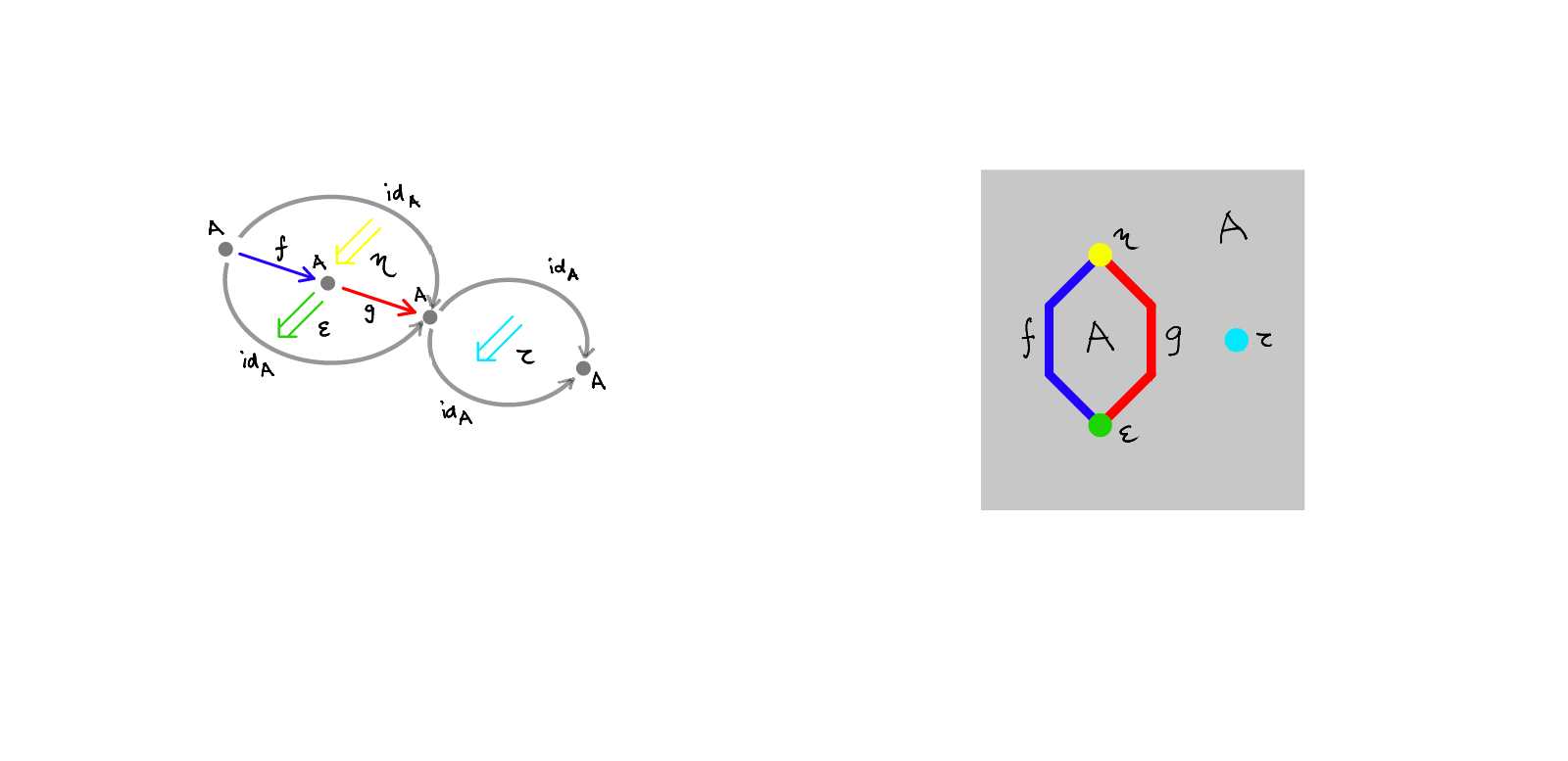}
\endgroup\end{restoretext}

\item For surface diagrams consider the following example
\begin{restoretext}
\begingroup\sbox0{\includegraphics{test/page1.png}}\includegraphics[clip,trim={.2\ht0} {.0\ht0} {.2\ht0} {.0\ht0} ,width=\textwidth]{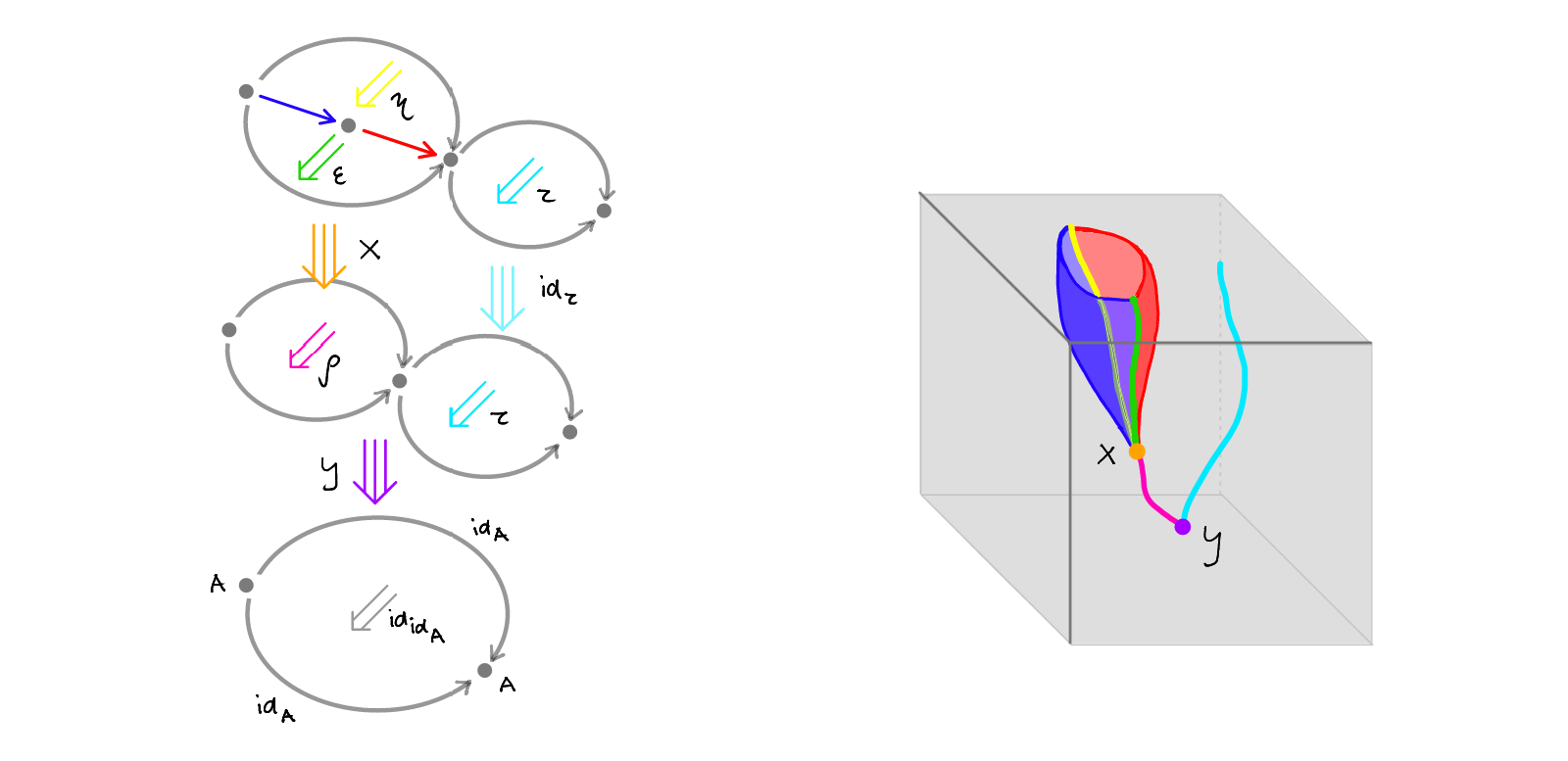}
\endgroup\end{restoretext}
Note how the top slice of the cube equals the last of our string diagram examples. Also note how the dualisation of dimension plays out in dimension $3$: $3$-cells turn into $0$-manifolds, $2$-cells into $1$-manifolds, $1$-cells into $2$-manifolds and $0$-cells into $3$-manifolds. 
\end{enumerate}

We end this section on manifold diagrams by recalling the following central point.

\begin{rmk}[Manifold diagrams as local models of higher categories] As mentioned before, higher categories don't have a canonical geometric model in the same way that higher groupoids have a model in spaces. This means theories of higher categories don't have a canonical benchmark for ``correctness". As proposed by various authors, geometric models for higher categories do exist (see for instance \cite{trimble1995} and \cite{ayala2017cobordism}). In agreement with \cite{trimble1995}, in this thesis we will treat manifold diagrams as a natural local geometric model for higher categories, analogous to pasting diagrams of cells being local models for groupoids.
\end{rmk}

\subsection{Recovering manifold diagrams} \label{rmk:GTP} 

In a special case of ``framed" manifold diagrams (meaning manifolds with a certain framing, which conjecturally is the case for manifold diagrams corresponding to morphisms of higher groupoids as discussed in \autoref{ch:geom}), manifold diagrams can be motivated from a yet different, more classical line of thinking. We will outline this here, as the idea might be of interest to the more topologically minded reader, with more details found in \autoref{ch:geom}.

We start from the classical Thom-Pontryagin construction \cite{pontrjagin2007smooth} which can be stated as 
\begin{equation}
[S^n,S^k] \iso \Omega^{\mathrm{fr}}_{n-k}(S^{n})
\end{equation}
where the left hand side denotes homotopy classes of continuous maps $S^n \to S^k$ from the $n$-sphere to the $k$-sphere, and the right hand side denotes closed framed $(n-k)$-manifolds in $S^n$. This was generalised in \cite{buonchristiano1976geometric} to give a correspondence of the form
\begin{center}
$[M,X]$ $\quad \cong \quad$ \{\text{framed $X$-stratification of $M$ up to ``cobordism"}\}
\end{center}
where $X$ is a CW-complex and $M$ is closed manifold. Conjecturally, this can be generalised even further to include compact manifolds $M$ with boundary, such as the $n$-disk $D^n$, to yield a correspondence
\begin{center}
$[M,X]^{\text{rel-}\partial}$ $\quad \cong \quad$ \{\text{framed $X$-stratification of $M$ up to ``rel-$\partial$ cobordism"}\}
\end{center}
(Here, the left-hand side denotes maps up to homotopies which preserve the boundary in a certain sense.) In the case that $M = D^n$ the left-hand side gives a good notion of ``$n$-morphisms in the groupoid $X$" and the right-hand side outputs a stratification of the $D^n$ disk in which $k$-strata correspond to $(n-k)$-cells of the CW-complex $X$. This is the dualisation of dimension that we have seen for string and more generally manifold diagrams in the previous section.

For instance, the Hopf map\footnote{This can be obtained from the usual Hopf map $S^3 \to S^2$ by removing a point from $S^3$. Note that $\mathrm{int}(Y)$ denotes the interior of a topological space $Y$} $\mathrm{int}(D^3) \to S^2$, leads to a (framed) $S^2$-stratification of $\mathrm{int}(D^3)$ of the form
\begin{restoretext}
\begingroup\sbox0{\includegraphics{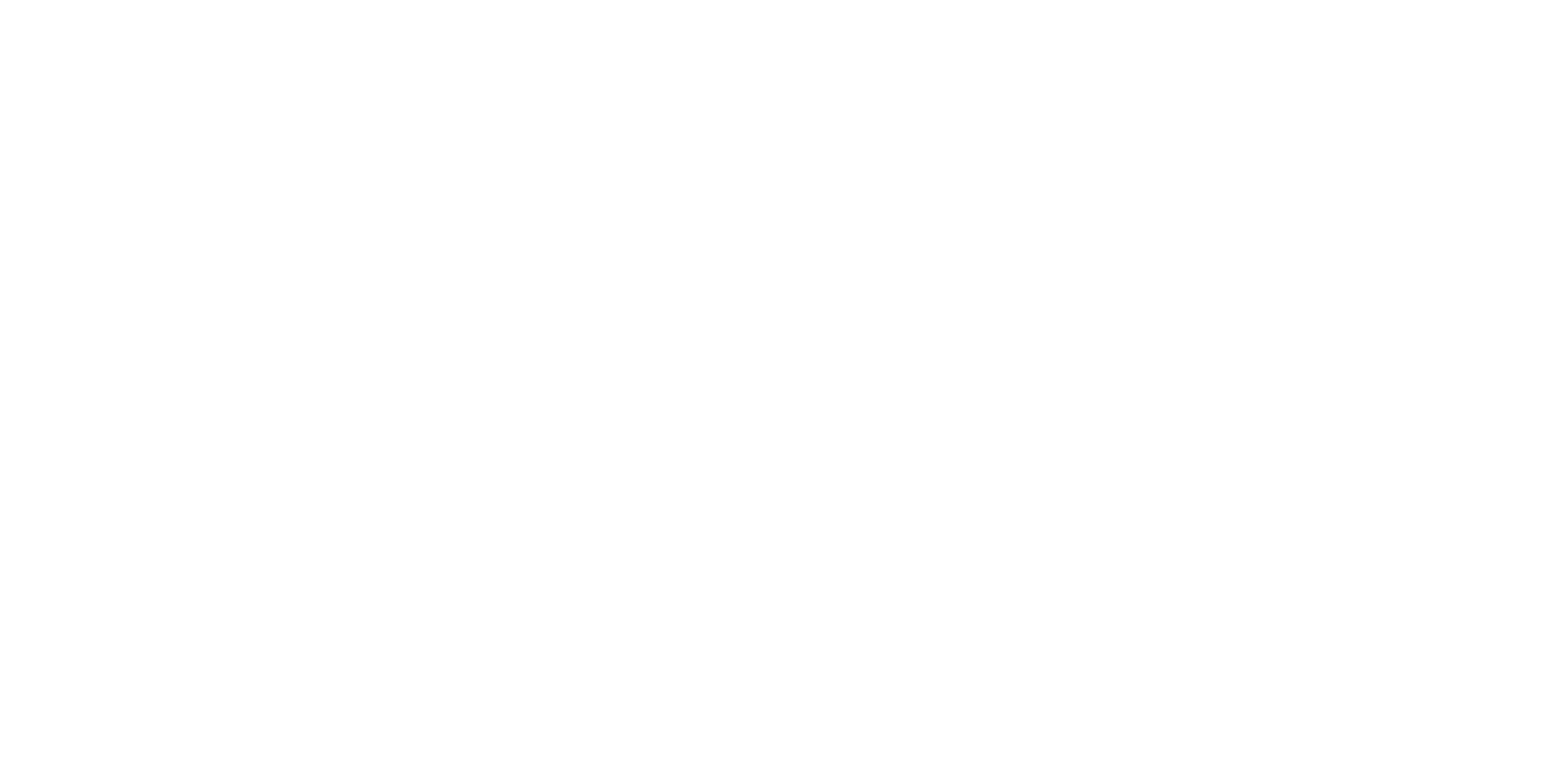}}\includegraphics[clip,trim={.0\ht0} {.25\ht0} {.0\ht0} {.25\ht0} ,width=\textwidth]{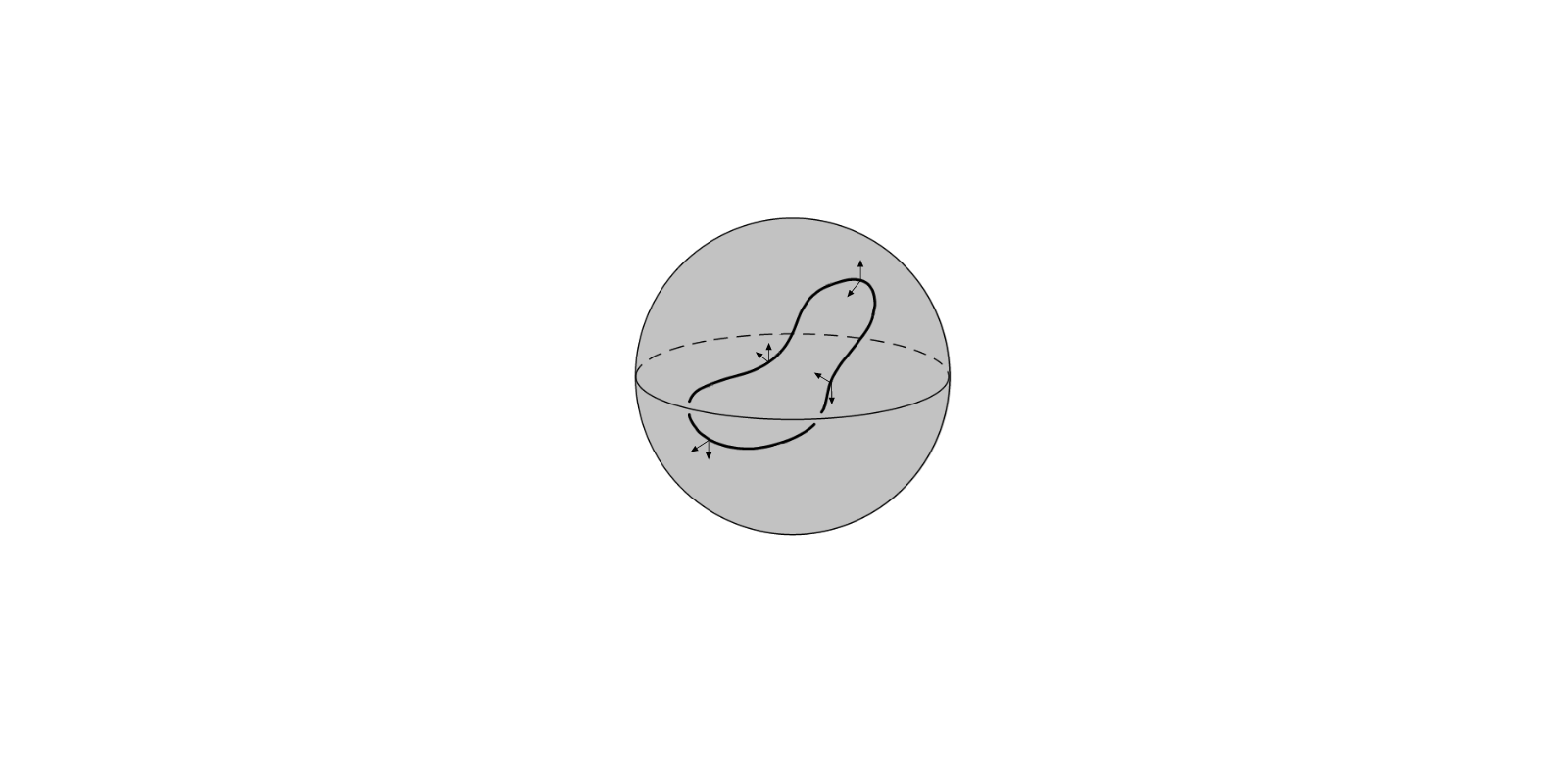}
\endgroup\end{restoretext}
Here, we are using the standard presentation of $S^2$ as a CW-complex which contains two cells (a $0$-cell and a $2$-cell) leading to two strata in the above picture (of dimension $3$ and $1$ respectively). A different presentation would  lead to a different stratification of $D^3$.

To fully recover our notion of $n$-manifold diagrams from a (framed) $X$-stratification of $D^n$, in \autoref{ch:geom}, we will pull back strata  along certain maps $F: (0,1)^n \to \mathrm{int}(D^n)$ later on called globular foliations, obtaining a stratification of the $n$-cube $(0,1)^n$. This in general is not (yet) a manifold diagram since the pullback strata might have certain ``singularities", which needs to be remedied by introducing new strata. For example, pulling back the $S^2$-stratification of $D^3$ above along a chosen globular foliation could yield any of the following $3$-manifold diagrams
\begin{restoretext}
\begingroup\sbox0{\includegraphics{ANCimg3/empty.png}}\includegraphics[clip,trim={.0\ht0} {.25\ht0} {.0\ht0} {.25\ht0} ,width=\textwidth]{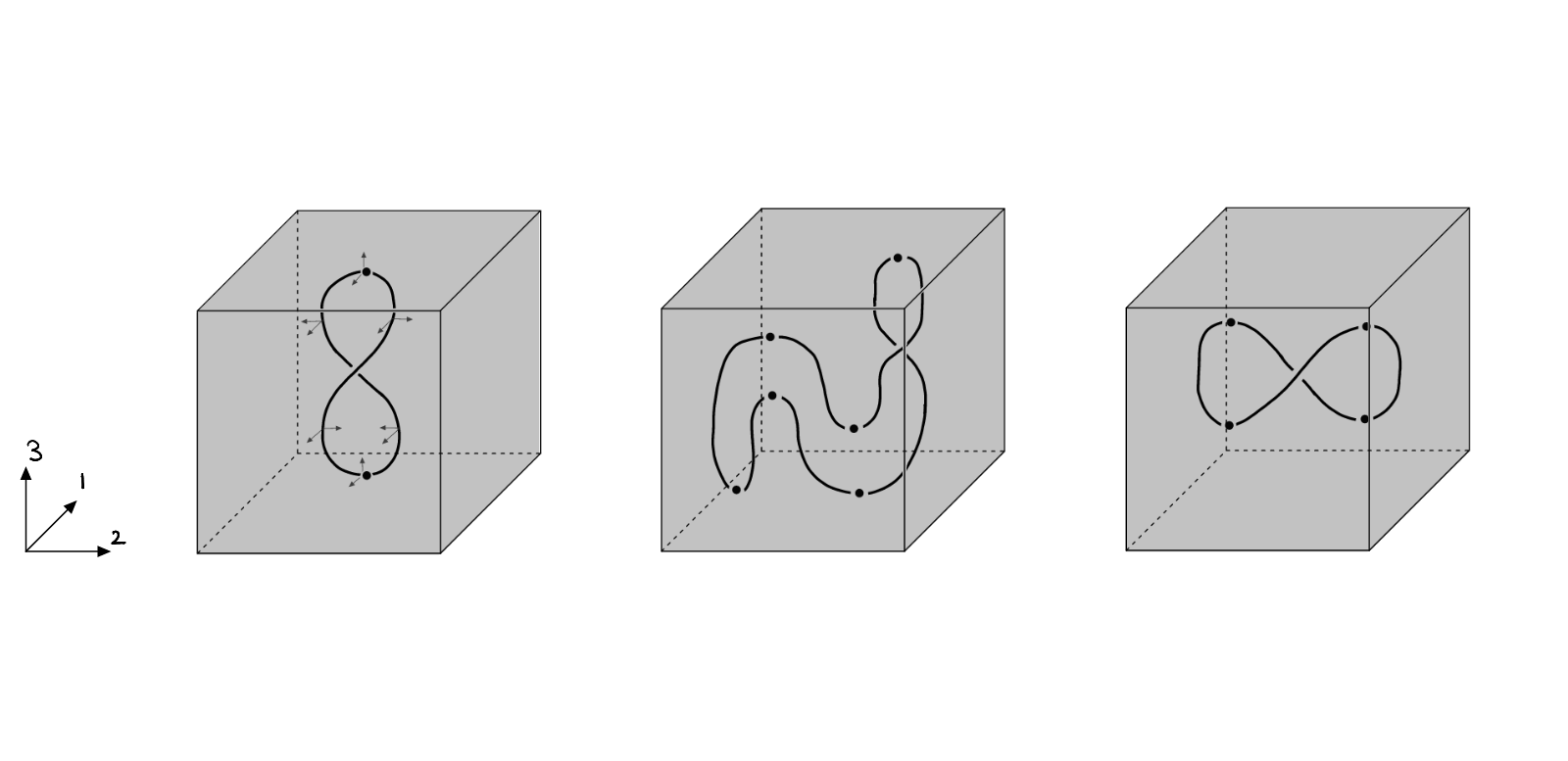}
\endgroup\end{restoretext}
In \autoref{ch:geom} we argue that all of these diagrams are in fact equivalent. There, we will also give a rough sketch of the conjecture that (groupoidal) manifold diagrams allow for a ``canonical choice" of framing (indicated for the cube on the left), which can represent the framing obtained in the generalised Thom-Pontryagin construction.

\subsection{Building an algebraic model} \label{ssec:po_alg_of_diag}

In this thesis we will work almost exclusively with a combinatorial analogue of manifold diagrams. The central observation leading towards this combinatorial structure is that restricting a $n$-manifold diagram to one of its ``$k$-sheets" ($k$-subcubes of the form $\Set{p} \times (0,1)^k$, $p \in (0,1)^{n-k}$, see \autoref{ssec:po_mfld_diag}) yields a $k$-manifold diagram. Consequently, we can recursively analyse a $n$-manifold diagram by looking at its sheets, and as it turns out for sufficiently nice such cubes, only finitely many sheets need to be recorded together with the transitions between them. This (recursively) leads to a finite combinatorial structure describing our original $n$-manifold diagram.

The resulting combinatorial structures will be called \textit{$\cC$-labelled singular $n$-cubes}, or $\SIvert n \cC$-cubes for short, where $\cC$ is a category whose objects are the ``labels" (also sometimes called ``colors" as they will often be represented as such). Their theory will be studied in \autoref{ch:intervals} to \autoref{ch:globes} of this thesis. There are various algebraic notions of higher categories that can be build from this, and we will do so starting in \autoref{ch:presented}.

\subsection{Definitions of higher categories} \label{ssec:po_anc_def}

We briefly sketch ideas involved in the various definitions of higher categories. To begin with, we first emphasize that $\cC$-labelled singular $n$-cubes are more expressive than a mere description of manifold diagrams and the additional generality is crucial, leading for instance to a definition of (presented associative) \textit{$n$-fold categories} in \autoref{ch:nfold} (also see  \autoref{rmk:po_nfold}).

The first higher categorical definition in this thesis will be a ``foundation-independent" (and in particular, computer implementable) definition of \textit{\free{} associative $n$-categories} in \autoref{ch:presented}. In fact, the definition is simple enough to be sketched as follows. We remark that the terminology ``\free{} associative $n$-categories"  is a red herring: it is not a classical $n$-category but a presentation thereof meaning it is given by sets of (higher) generators and relations\footnote{In higher category theory both generators and relations will be given in the form of higher morphisms. Thus the distinction between generators and relations becomes somewhat artificial.}.

\begin{defn}[Presented associative $n$-categories] A presented associative $n$-category $\sC$ is a list of sets $\sC_0, \sC_1, ... \sC_{n+1}$, where $\sC_k$ is called set of generating $k$-morphisms, together with an assignment of a ``conical" $k$-manifold diagrams $\abss{g}$ to each $g \in \sC_k$ (here, ``conical" means that the manifold diagram is obtainable as the cone of its boundary), such that the following is satisfied
\begin{enumerate}
\item $\abss{g}$ is \textit{globular} (which is a certain constancy condition of manifolds on the side of the $k$-cube)
\item $\abss{g}$ is \textit{well-typed}: If $p$ is a point in the cube $\abss{g}$ of color $f \in \sC_l$, then the minimal open neighbourhood of $p$ looks like $\abss{f}$ (or more precisely, the $(k-l)$-fold identity of $\abss{f}$)
\end{enumerate}
$\abss{g}$ is called the \textit{type} of $g$.
\end{defn}
Further to our remark above it is in fact not immediately clear whether a presented associative $n$-category is also an associative $n$-category, or even how to define the latter, which we will only do in later chapters. \textit{Morphisms} (also called \textit{composites}) of an presented associative $n$-category are manifold diagrams which are globular and well-typed (but, unlike types, they are not required to be conical).

After discussing theories of coherent invertibility we will specialise presented associative $n$-categories to a notion of \textit{presented associative $\infty$-groupoids} in \autoref{ch:groupoids}. The latter notion will be obtained from the former by requiring every generator to be invertible. Note that, with a view towards applications in topology, spaces (that is, CW-complexes) can also be regarded as presented structures: they are ``freely" obtained from gluings of their generating cells. In \autoref{ch:geom}, we sketch how CW-complexes can be translated into presented associative $\infty$-groupoids.

In \autoref{ch:composition} we will study the compositional behaviour of cubes, and how they can be glued together along their sides. This will lead us to formulate a notion of whiskering composition for morphisms in presented associative $n$-categories. We will then characterise ``perturbation-stable" \textit{generic composites}\footnote{In geometric terms, the idea should be taken literally: A cube is generic if a small directed perturbation of any subcube does not change the cube's equivalence class. Examples will be given in \autoref{ch:composition}.} in two ways: They can be described as the subset of all composites satisfying a genericity condition, or they can be characterised as being inductively built from generators and \textit{elementary homotopies} using the whiskering operation.

This inductive description of the set of generic composites will then finally allow us to give a definition of associative $n$-categories in \autoref{ch:associative}. An associative $n$-category is a globular set with an algebra structure, which in particular provides a weak coherence for each elementary homotopy. All other coherences (that is,  composition operations and ``higher-depth" homotopies involving such operations) will be chosen strictly. These strictness conditions will be mathematically packaged by a notion of algebra \textit{resolutions}. In fact, the idea of resolutions can be iterated, and the higher-depth homotopies arising in this process need not at all be chosen to act strictly in order to yield a valid definition. As a consequence we will actually discover a \textit{spectrum} of definitions, ranging from ``fully associative" to ``fully weak".

Finally, it is worthwile re-emphasizing that all ``associative" notions of higher categories above are conjecturally \textit{semi-strict}. This means that while part of the coherence data is strict the resulting notion of category can still capture all equivalence classes of weak $n$-categories (cf. \cite{simpson2011homotopy}) and in particular all homotopy types of spaces. Preliminary thoughts on this equivalence are offered after we give a definition \textit{\free{} weak $n$-categories} in \autoref{ch:weak}, but also the fact that there is a ``continuous" spectrum of definitions ranging from ``fully associative" to ``fully weak" might be a  promising lead.

\section{Coherences in higher category theory} \label{sec:coherences} 

This section is a short and informal introduction to certain basic ideas in higher categories. At the centre of our discussion will be the notion of ``coherence", which captures non-trivial relations between higher-dimensional objects in a general higher category.

Higher category theory has turned out to be a difficult subject. This is unsurprising, as higher categories should be at least as difficult as the study of spaces. Similar to higher dimensional homotopy groups being difficult to compute, ``homotopical phenomena" or \textit{coherences} complicate the development of fully algebraic approaches to higher category theory---unless one ignores all coherences which means working in the setting of strict higher categories. We substantiate this line of thought as follows. 

On the side of homotopy types, homotopy groups of spheres are trivial in dimension $n = 0$ and $n = 1$ (that is, $\pi_{k}(S^n) = 0$ for $k > n$). Only in dimension $n = 2$ we start to see non-trivial phenomena of homotopy groups. In lowest dimensions these non-trivial phenomena arise in $\pi_3(S^2) = \lZ$. This homotopy group is generated by the so-called Hopf map, which by the Thom-Pontryagin construction (and using certain conventions) can be visualised as
\begin{restoretext}
\begingroup\sbox0{\includegraphics{test/page1.png}}\includegraphics[clip,trim=0 {.35\ht0} 0 {.25\ht0} ,width=\textwidth]{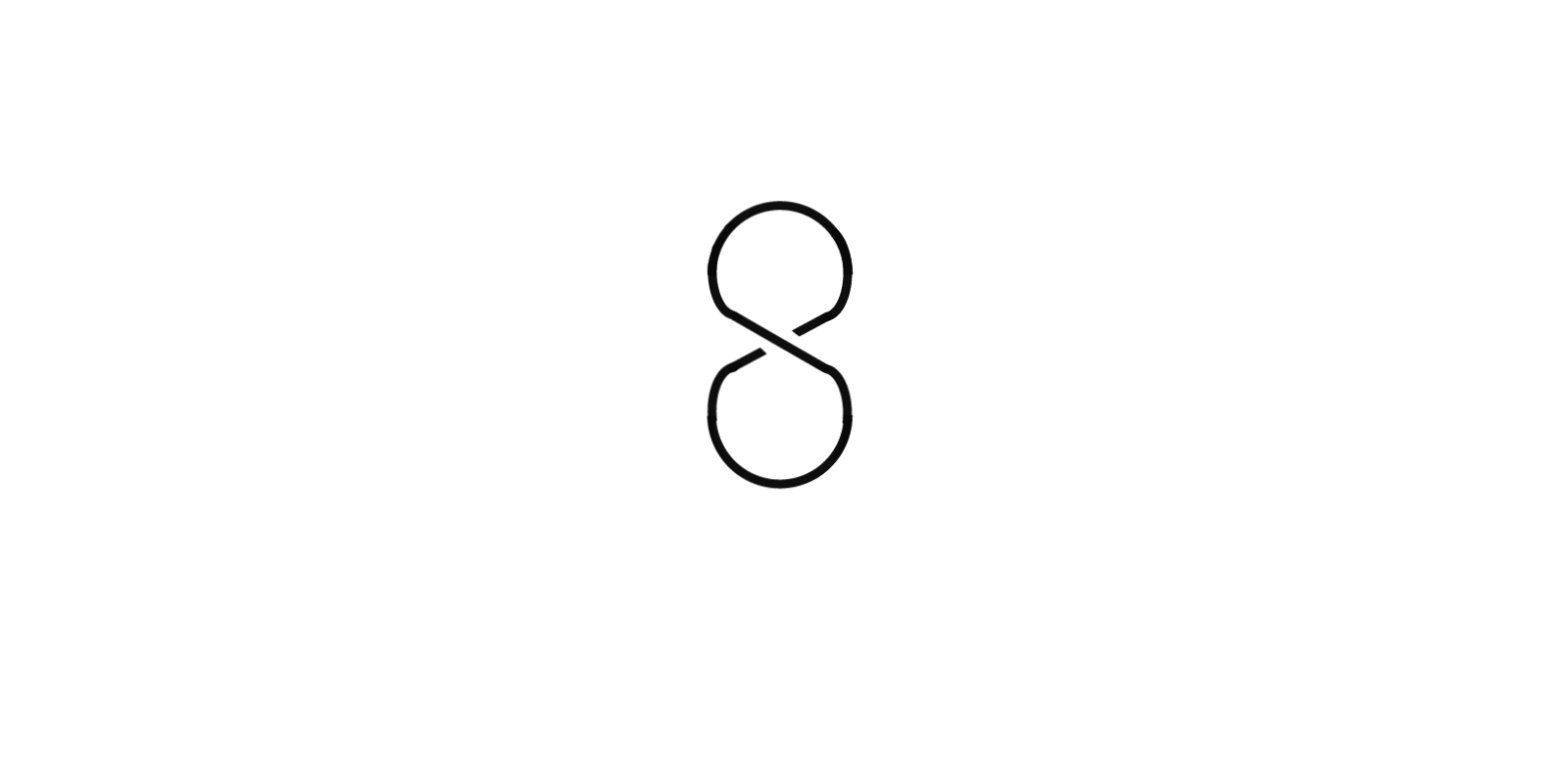}
\endgroup\end{restoretext}

On the side of categories, the study of $1$-categories and $2$-categories can be considered ``easy". For both it suffices to study their strict version. In dimension $3$ we find the first obstruction: not every $3$-category is equivalent to a strict $3$-category, and thus strict $3$-category theory is not fully general. In \cite{gordon1995coherence} it was shown that the ``strictest possible" fully general theory of $3$-categories is the theory of Gray categories. A Gray category is strict up to a cell called ``the interchanger" and commonly visualised as
\begin{restoretext}
\begingroup\sbox0{\includegraphics{test/page1.png}}\includegraphics[clip,trim=0 {.37\ht0} 0 {.37\ht0} ,width=\textwidth]{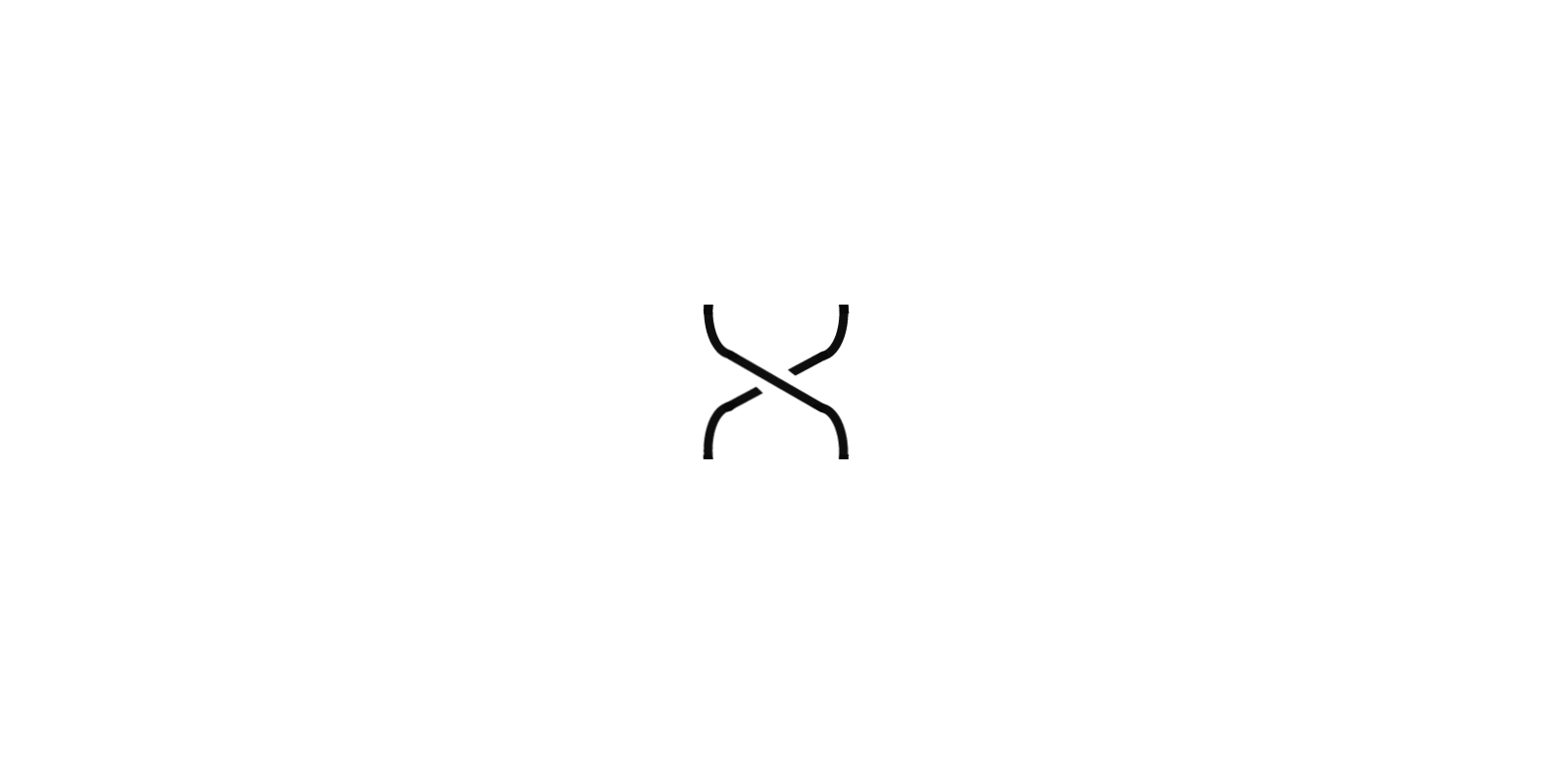}
\endgroup\end{restoretext}

In this work we will see that non-triviality of $\pi_3(S^2)$ is directly to the related non-triviality of the theory of $3$-categories. Moreover, we will give an algebraic and computer-implementable language to describe all homotopical phenomena (or as they will later on be called, \textit{homotopies}) such as the interchanger above.

The notion of associative $n$-category that we will describe, can be regarded as a generalisation of (unbiased\footnote{Here, the word ``unbiased" refers to axiomatising categories using general $n$-ary composition as opposed to the ``biased" approach, in which we only consider binary composites. A good discussion of the distinction can be found in \cite{leinster-operads}.}) Gray categories to arbitrary dimension $n$. By saying this, we in particular want to suggest that associative $n$-category theory is equivalent to weak $n$-category theory. This will turn out to be a substantial strengthening of Simpson's conjecture (cf. \cite{simpson1998homotopy}), reducing the requirement of ``weak identities" therein to merely ``weak homotopies". Before explaining our approach in more detail, we will briefly recall pathways to the study of higher categories and give a little bit of context for the notion of Gray categories and the aforementioned Simpson's conjecture.

\subsection{Algebraic perspective: Categorification} \label{ssec:categorification}

Equality $a = b$ of two objects $a$, $b$ plays a central role in mathematics. However, stating equality $a = a$ of two actually \textit{indistinguishable} objects is trivial, and commonly not of interest to mathematicians. Instead, mathematicians are usually concerned with equality $a = b$ of distinguishable objects, and the statement of their equality is reached by a computation, or proof, or some other process (all of ``computation", ``proof" and ``process" should really be understood synonymously here). As an example we consider equality of morphisms in a category: consider the statement $fg = h$. The terms $fg$ (``$f$ composed with $g$") and $h$ are distinguishable to the reader. A proof of the equality could for instance be derived from postulating $g = kl$ and $(fk)l = h$ together with the associativity law of categories.

\textit{Categorification} is the (informal) process of retrieving information that has been lost by omitting these computational steps, and making them explicit again. While equality $a = b$ usually allows to bidirectionally compute $a$ from $b$ and $b$ from $a$, categorification in principle also allows us to consider \textit{directed} processes $a \to b$. Roughly speaking, categorification replaces equalities by processes or \textit{(higher) morphisms}. In particular, proofs just as morphisms compose. For our example of equality proofs of morphism, we note the following modes of composition for such proofs
\begin{enumerate}
\item \textit{Vertical composition}: If $f = g$ and $g = l$ then $f = l$.
\item \textit{Horizontal composition}: If $f = g$ and $h = k$ then $fh = gk$ (if these composites exist)
\item \textit{Post-whiskering}: If $f = g$ then $fh = gh$ (if these composites exist)
\item \textit{Pre-whiskering} If $h = k$ then $fh = fk$ (if these composites exist)
\end{enumerate}
To introduce the usual categorical language for this context, we note that an equality proof $f = g$ can be regarded as a (bidirectional) morphism between morphism $f$ and $g$. This is called a $2$-morphism. We therefore found that $2$-morphisms have more than one mode of composition. These different modes of composition will later be fully explained by their underlying geometry. More generally, $(k+1)$-morphisms are morphisms between $k$-morphisms and will have yet other modes of composition.

Categories consisting of $1$-morphism and $2$-morphisms (and their respective compositional data) are called $2$-categories. Based on the above discussion, the theory of $2$-categories is a natural candidate for a categorification of the theory of $1$-categories. Similarly, $3$-categories additionally introduce $3$-morphisms between $2$-morphisms. Repeating the process leads to the concept of $n$-categories and as a ``colimit" of this process we obtain $\infty$-categories which consist of $k$-morphisms for any $k \in \lN$. However, the simplicity of this picture is deceiving and details have been omitted. In general, it not straight-forward to replace the ``structural equalities" which are part of the theory of categories by higher morphisms. 

An example of a structural equality is the associativity law in the theory of categories, which states the equality $(fg)h = f(gh)$ (more generally, we understand structural equalities to mean equalities that are given as part of the theory of $n$-categories). According to the principle of categorification this should be interpreted as a $2$-morphism $\alpha_{f,g,h} : (fg)h \to f(gh)$ after categorification. This $2$-morphism is usually called an \textit{associator}. The associator and other cells introduced by categorifying structural equalities will be referred to as \textit{coherence data}. Having introduced the associator as a new morphism we now need to answer the question what compositional data it has. This question is non-trivial, and for instance, it turns out that composites of associators naturally satisfy the so-called \textit{pentagon identity}. Categorifying this identity in turn leads to a cell called the pentagonator.

We will provide a general procedure for generating coherence data in \autoref{ch:weak}, where we for instance derive associators and pentagonators based on algebraic principles. For now, we will turn to geometric principles to explain the ``mysteries" of coherence data.

\subsection{Geometric perspective: Homotopy hypothesis} \label{ssec:intro_hom_hyp}

A second pathway to higher categories is of geometric nature and concerns a special case of higher categories: $\infty$-groupoids are $\infty$-categories in which all (higher) morphism are invertible. The homotopy hypothesis asserts that there is an equivalence of $\infty$-groupoids with (the homotopy theory of) spaces, with the following underlying idea: given an $\infty$-groupoid $\sX$ then we obtain a space $\abs{G}$ called \textit{geometric realisation} of $G$. This process involves translating all morphisms in $\sX$ into cells, and gluing them together according to their specified inputs and outputs. Conversely given a space $X$ we can form its \textit{fundamental $\infty$-groupoid} $\Pi_n X$, whose $k$-morphisms are given by ``$k$-cells in $X$", that is, by mappings of the $k$-disk\footnote{More precisely, we would have to chose a combinatorial model of the $k$-disk (e.g. the $k$-simplex $\Delta^k$, in which case the fundamental $\infty$-groupoid is called singular nerve of $X$).} $D^k$ into $X$. Composition of morphisms is given by gluing (and reparametrisation) of disks. For instance, consider two glueable $1$-cells $f , g : [0,1] \to X$ in a space $X$. 
Their composition gives a $1$-cell $(fg) : [0,1] \to X$ after identifying (or ``reparametrising") the gluing of two $1$-cells with a single $1$-cell. For example, we can set
\begin{equation}
(fg)(x) = \begin{cases} f(2x) & \text{~if~} x \leq \frac{1}{2} \\
g(2x - 1) & \text{~if~} x \geq \frac{1}{2} \end{cases}
\end{equation}
Choosing this reparametrisation implies that, given three glueable $1$-cells $f,g,h$ in $X$, $((fg)h)$ and $(f(gh))$ are not indistinguishable in general, but are different maps. However, they are equal up to a homotopy $\alpha : [0,1] \times [0,1] \to X$ such that $\alpha(-,0) = ((fg)h)$ and $\alpha(-,1) = (f(gh))$ while the endpoints $\alpha(0,-)$ and $\alpha(1,-)$ are held constant. Identifying $[0,1] \times [0,1]$ with the $2$-disk $D^2$, we observe that $\alpha$ is also a $2$-morphism in the groupoid corresponding to $X$. $\alpha$ can be seen to play the role of an associator. More generally, all coherence data (including for instance the pentagonator) can be derived from such homotopies. This is the ``geometric perspective" on coherence data.

\subsection{Weak, strict and semistrict notions} \label{ssec:intro_strict}

When categorifying a structure one is given a choice which equality to retain as equality and which to turn into morphisms. In the case of $n$-categories this leads to the following distinction.

The theory of \textit{weak} $n$-categories is obtained when turning ``all" structural equalities into morphisms. The theory of \textit{strict} $n$-categories is obtained by turning none of them into morphisms, and retaining all of them as equalities. The theory of strict $n$-categories is not equivalent to the theory of weak $n$-categories, and in particular the homotopy hypothesis does not hold when restricting to strict $n$-groupoids. This was shown by Simpson in \cite{simpson1998homotopy} refuting an earlier proof of the homotopy hypothesis by Kapranov and Vovoedsky in \cite{kapranov1991infty}. In the same work Simpson conjectured the following

\begin{conj}[Simpson's conjecture]  Every weak $\infty$-category is equivalent to an $\infty$-category in which composition and exchange laws are strict and only the unit laws are allowed to hold weakly.
\end{conj}

\noindent Subsequently to this conjecture, the idea of \textit{semistrict} higher categories gained traction: the theory of semi-strict $n$-categories is theory of higher categories with minimal coherence data (that is a minimal amount of morphisms that replace structural identities) such that each weak $n$-category is still equivalent to some semi-strict $n$-category. Semi-strict $n$-categories have only found a well-accepted formulation up to dimension $n = 3$. This formulation of semi-strict $3$-categories is called Gray categories and it is based on a single (weak) coherence called the \textit{interchanger} \cite{gordon1995coherence}. 

We will motivate that the idea underlying Gray-categories is not a random choice, but can in fact be motivated and generalised by geometric principles. Associative $n$-categories will be exactly this generalisation of (unbiased) Gray categories to dimension $n$.

\subsection{The role of $n$-fold categories} \label{sec:intro_nfold}

The framework we develop will first and foremost be a description of associative \textit{$n$-fold categories}---however, we will only give a definition of $n$-fold categories in the \free{} case. Roughly speaking, the $n$-morphisms in an $n$-fold category are (directed) $n$-cubes in place of $n$-globes. In this section we comment on the relation of $n$-fold categories to manifold diagrams and the cellular approach to $n$-categories.

To specialise from ``cube-like" morphisms in an $n$-fold categories to ``globe-like" morphisms $n$-categories one can introduce a certain constancy condition on the cube called \textit{globularity}. Topologically, this constancy condition can be enforced by quotienting with certain maps $[0,1]^n \to D^n$ (later called globular foliations). In dimension $n=2$, this can be illustrated as follows
\begin{restoretext}
\begingroup\sbox0{\includegraphics{test/page1.png}}\includegraphics[clip,trim=0 {.44\ht0} 0 {.36\ht0} ,width=\textwidth]{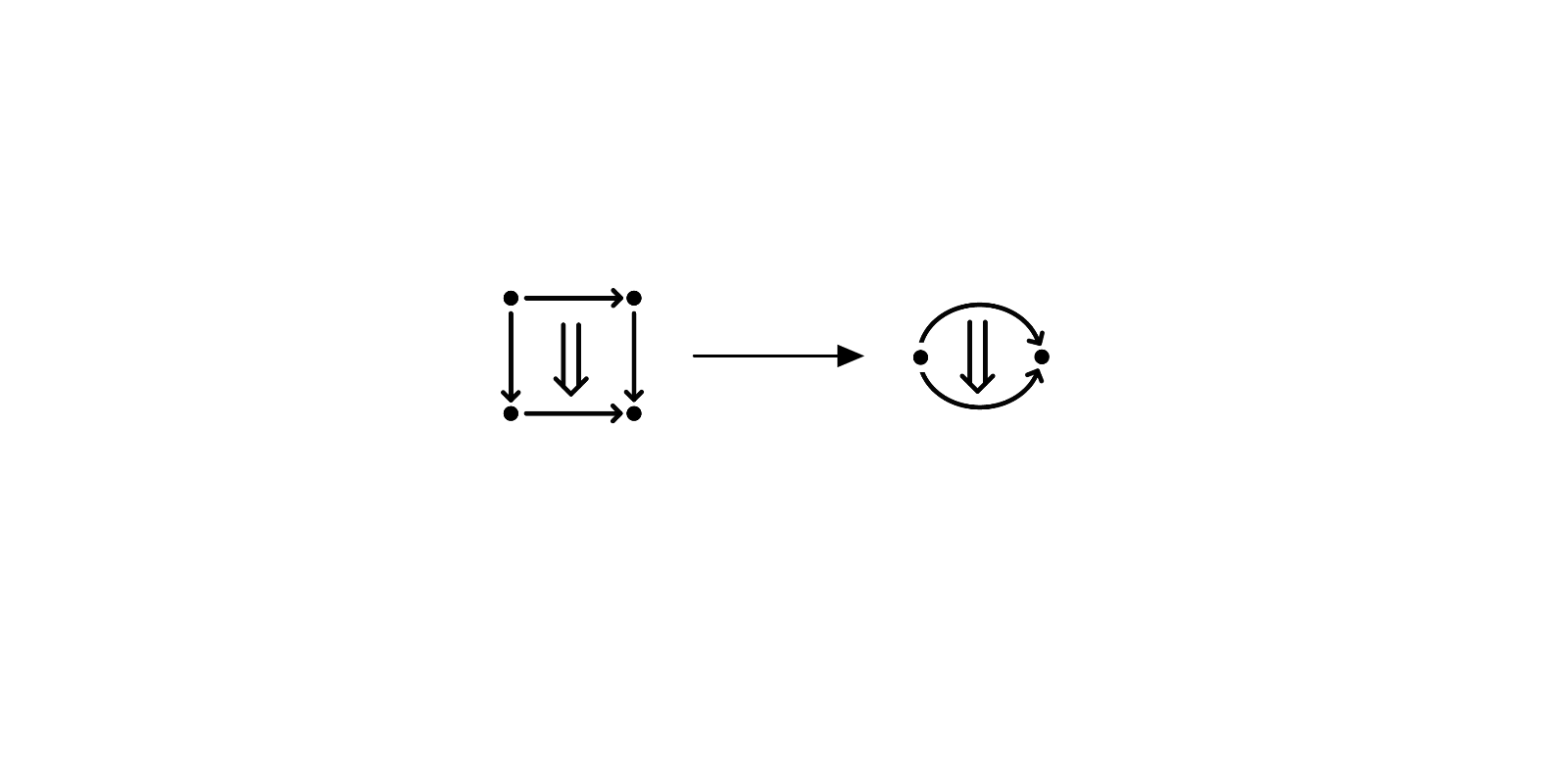}
\endgroup\end{restoretext}
Importantly, in this work we will not ``fully quotient" the sides of the cube and work with globes. Instead we will keep working with cubes and require the sides of the cube to be constant as forced by the globularity condition. This is illustrated as follows
\begin{restoretext}
\begingroup\sbox0{\includegraphics{test/page1.png}}\includegraphics[clip,trim=0 {.44\ht0} 0 {.36\ht0} ,width=\textwidth]{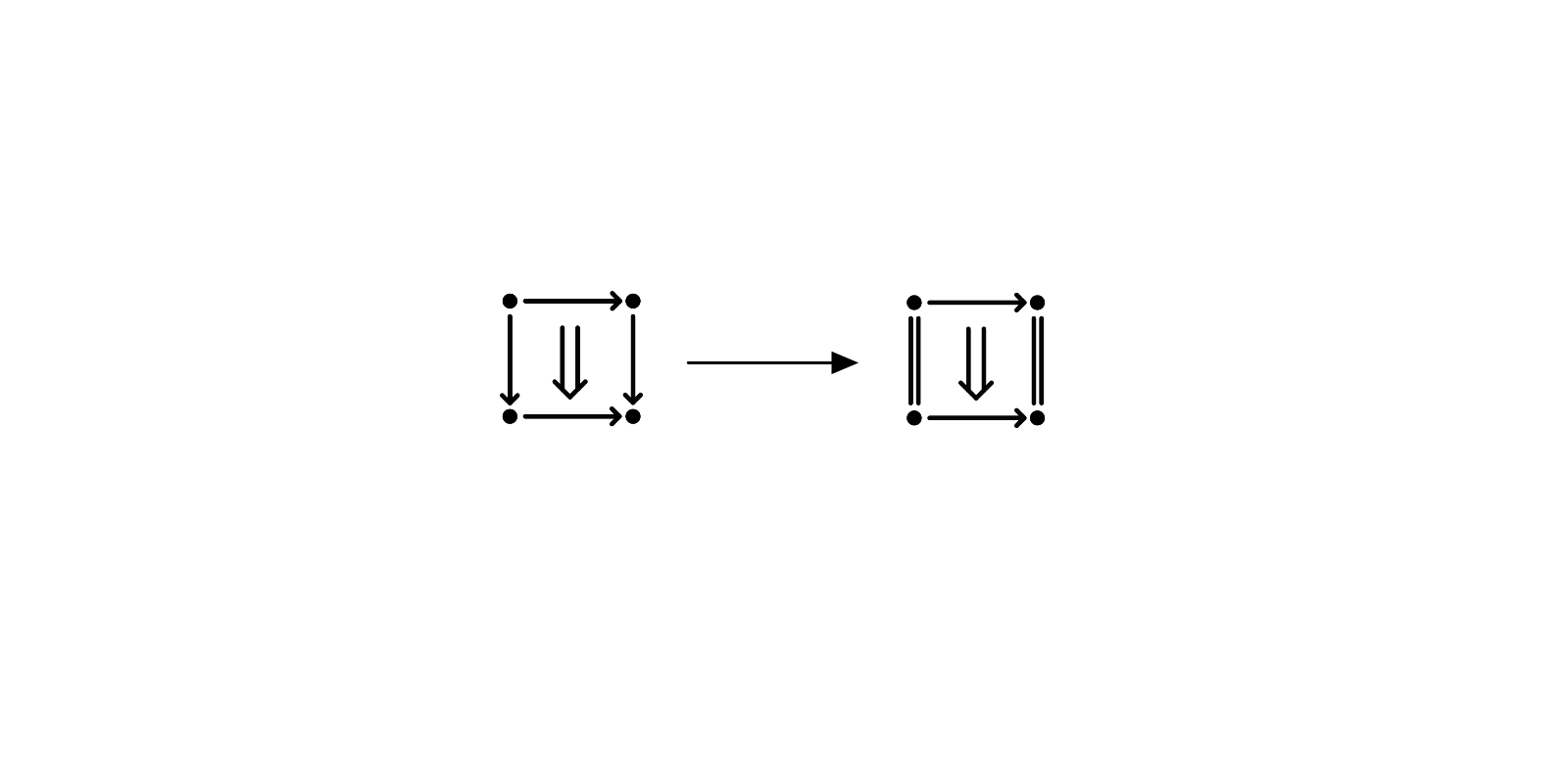}
\endgroup\end{restoretext}
Leaving the cubical structure present will play a crucial role in describing homotopies for coherence data. 

We substantiate this last statement with an example. We will consider two $2$-globes $\alpha : f \to h$ and $\beta :  g \to k$ (where $f,h : A \to B$, $g,k :  B \to C$) which can be represented as
\begin{restoretext}
\begingroup\sbox0{\includegraphics{test/page1.png}}\includegraphics[clip,trim=0 {.36\ht0} 0 {.35\ht0} ,width=\textwidth]{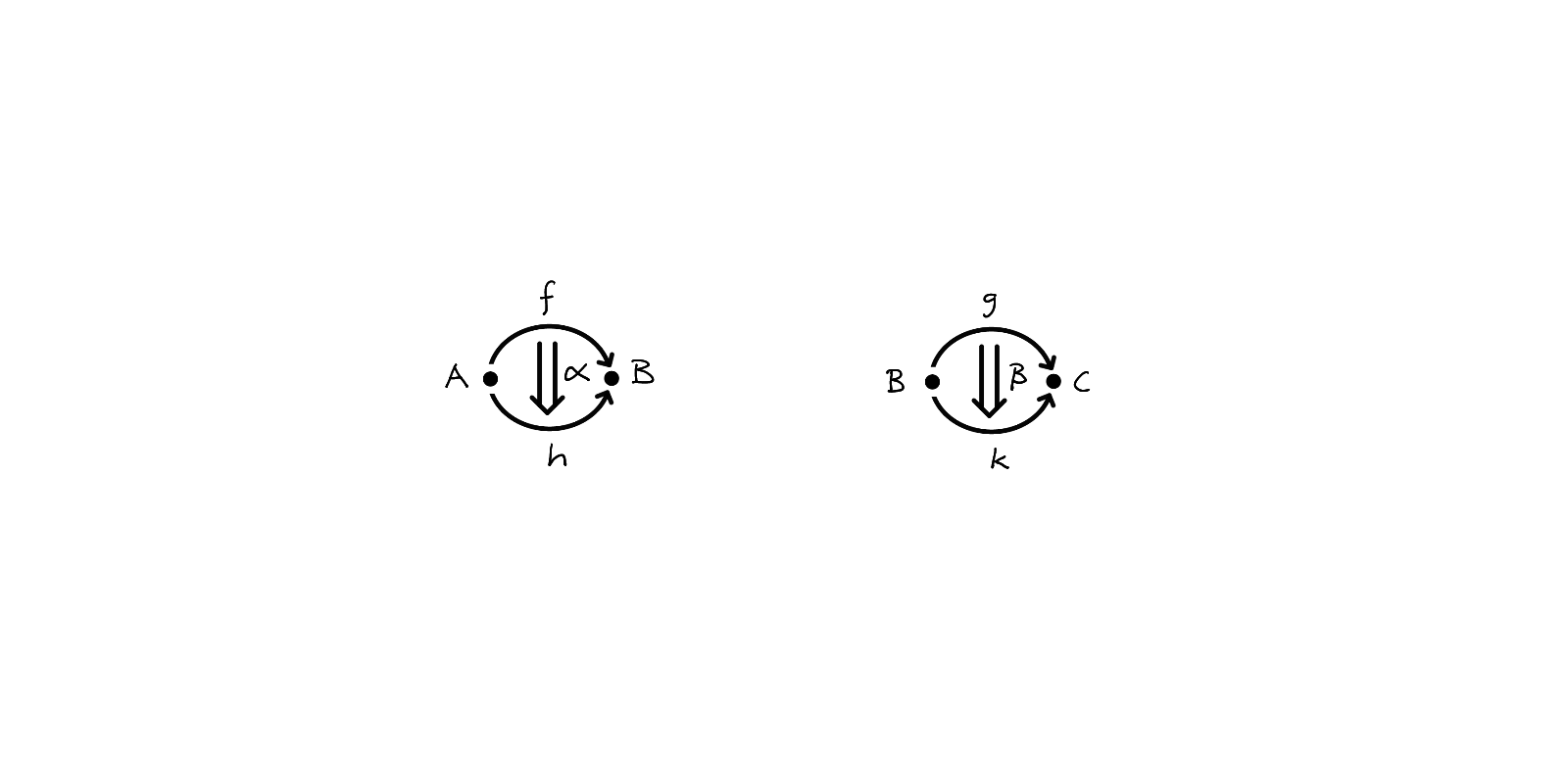}
\endgroup\end{restoretext}
The different modes of (categorical) composition of $2$-morphisms correspond to certain gluing operations of the underlying $2$-disks $D^2$ in this representation of $2$-globes. We now consider three sequences of such gluings (To clarify again: We are treating $2$-globes as topological disks when performing the gluing operations below).
\begin{enumerate}
\item Horizontal composition of $\alpha$ and $\beta$ corresponds to a gluing along their mutual $0$-boundary $B$ (indicated in \cred{}) yielding 
\begin{restoretext}
\begingroup\sbox0{\includegraphics{test/page1.png}}
\includegraphics[clip,trim=0 {.44\ht0} 0 {.4\ht0} ,width=\textwidth]{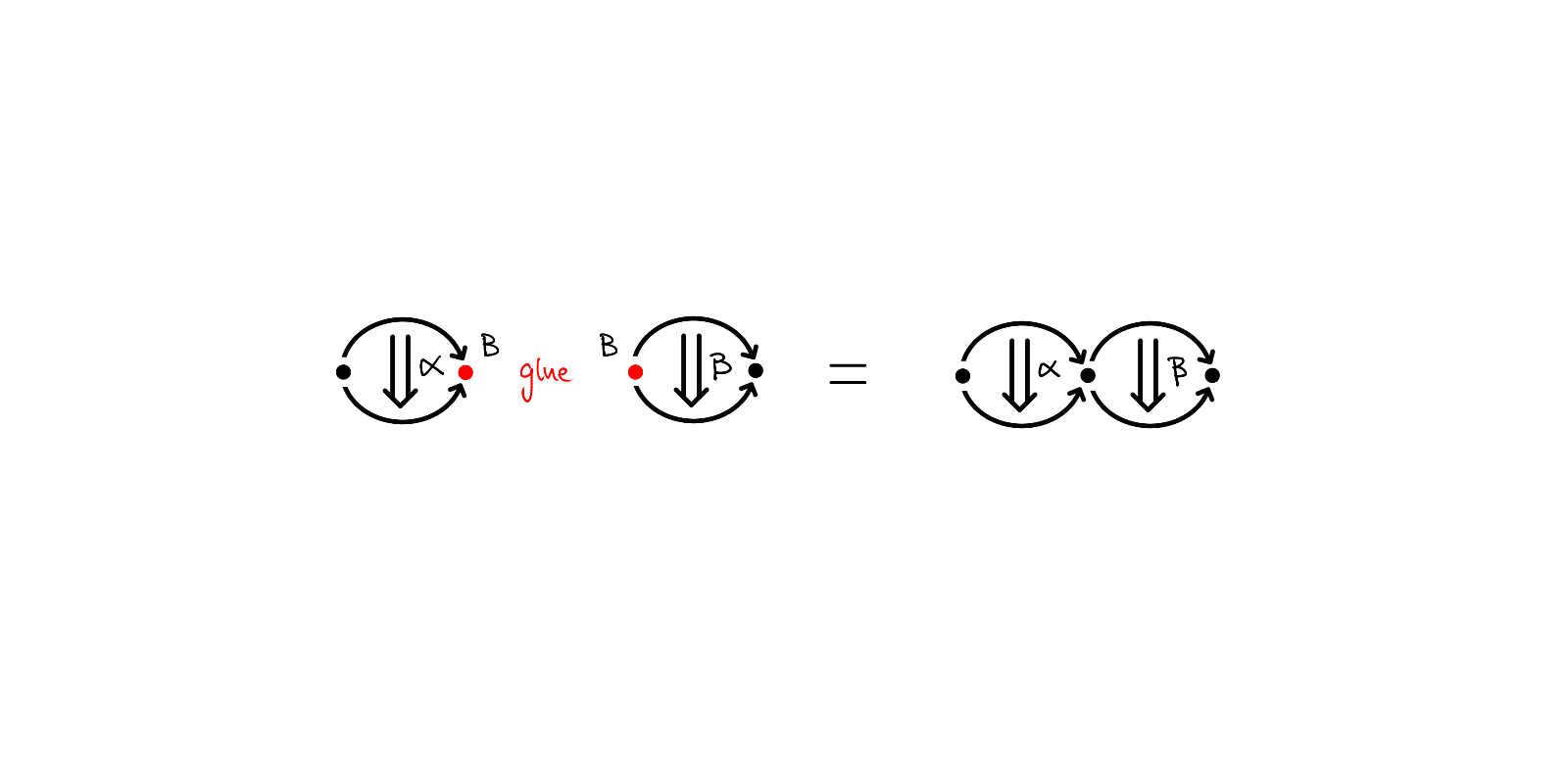}
\endgroup\end{restoretext}
\item Secondly, we can also glue $\alpha$ and $g$ along their mutual $0$-boundary $B$ (corresponding to post-whiskering of $\alpha$ with $f$) indicated in \cred{} on the left below, and we can glue $h$ and $\beta$ along their mutual $0$-boundary (corresponding to pre-whiskering $h$ and $\beta$) indicated in \cred{} on the right
\begin{restoretext}
\begingroup\sbox0{\includegraphics{test/page1.png}}\includegraphics[clip,trim=0 {.37\ht0} 0 {.12\ht0} ,width=\textwidth]{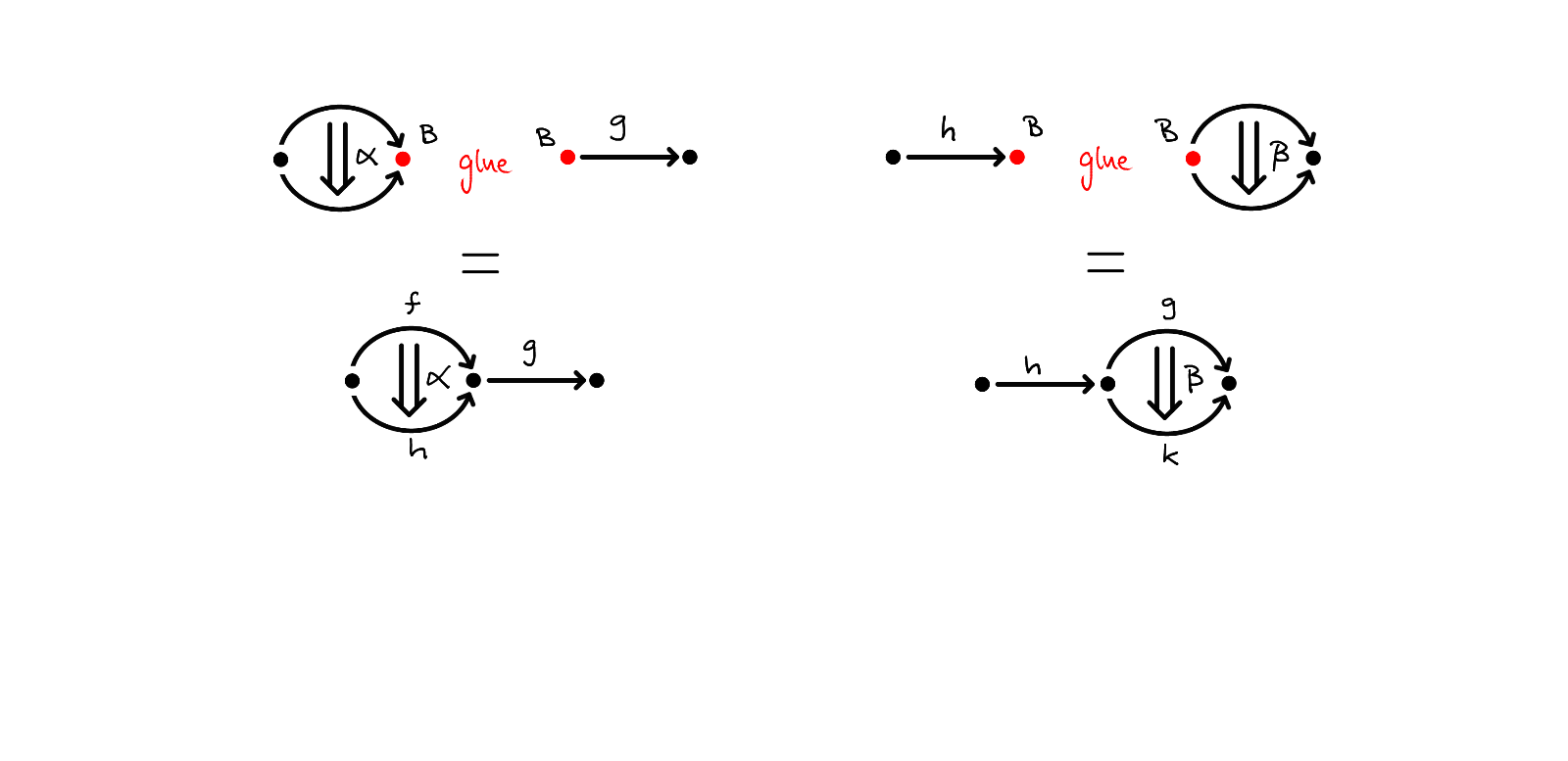}
\endgroup\end{restoretext}
The resulting compositions of cells share a mutual $1$-boundary $gh$. Gluing them along this boundary yields
\begin{restoretext}
\begingroup\sbox0{\includegraphics{test/page1.png}}\includegraphics[clip,trim=0 {.29\ht0} 0 {.3\ht0} ,width=\textwidth]{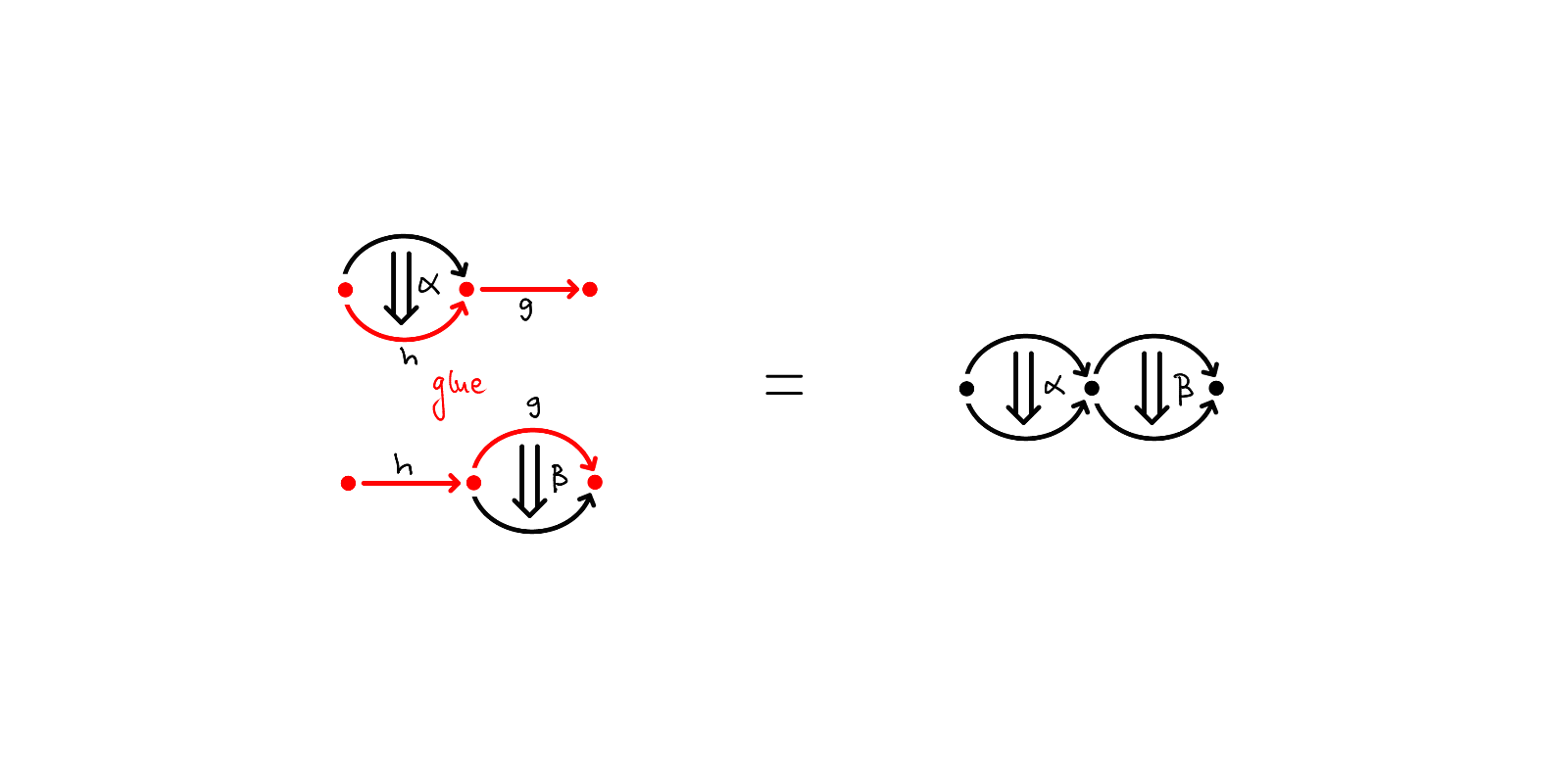}
\endgroup\end{restoretext}
\item Similarly, we can pre-whisker $f$ and $\beta$, post-whisker $\alpha$ and $k$ and glue the resulting composition along $kf$ (corresponding to vertical composition) to obtain
\begin{restoretext}
\begingroup\sbox0{\includegraphics{test/page1.png}}\includegraphics[clip,trim=0 {.25\ht0} 0 {.25\ht0} ,width=\textwidth]{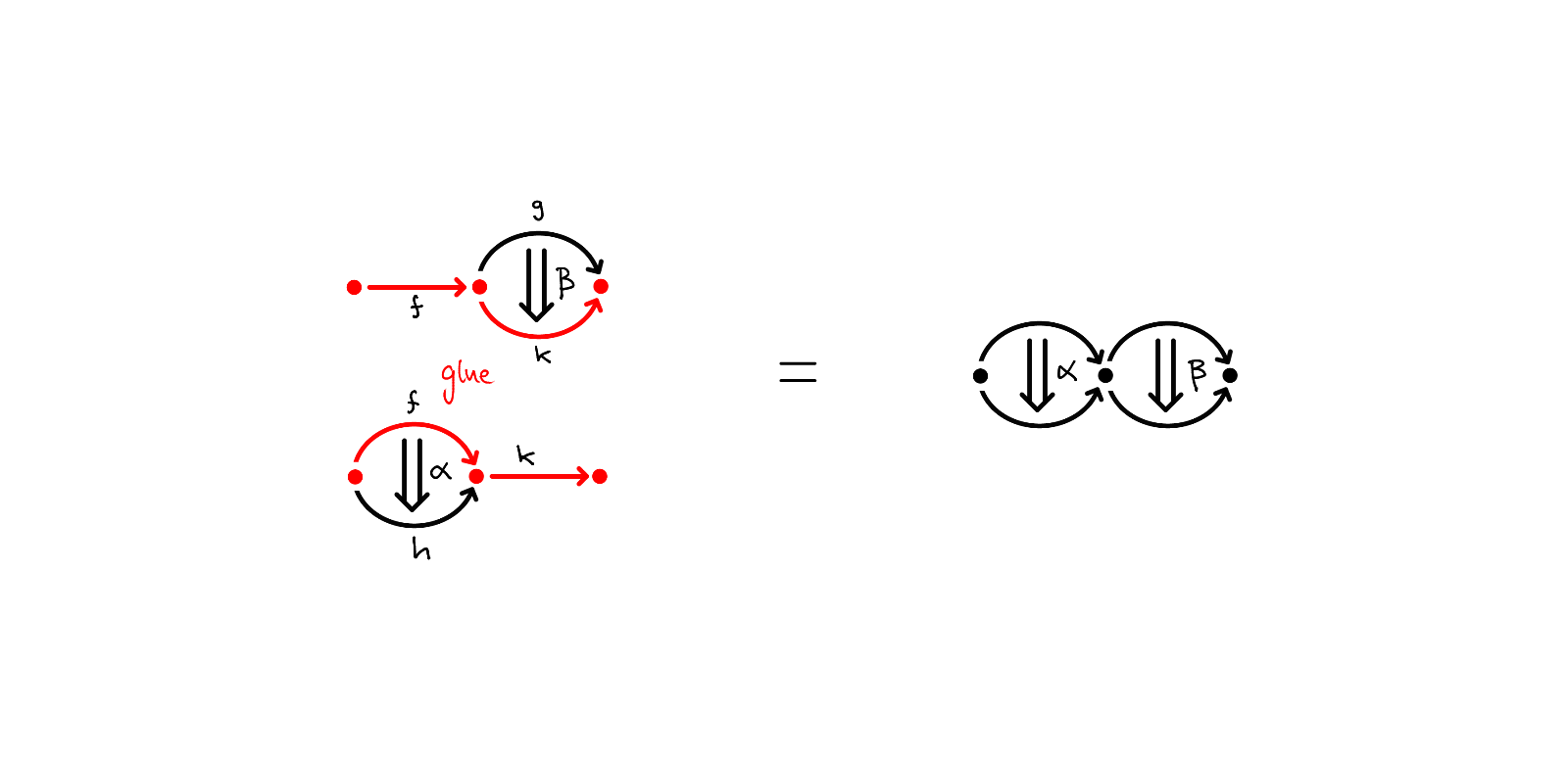}
\endgroup\end{restoretext}
\end{enumerate}
Importantly, we now observe that information about how these morphisms were glued is lost in the above, in the sense that \textit{the above cell complexes have the same presentation} in all three cases, which for our algebraic purposes will mean that they are \textit{indistinguishable}. As remarked in the beginning of \autoref{ssec:categorification} equalities are interesting if they are equalities between distinguishable objects, since in that case non-trivial proof is required. We consequently face the following dilemma: based on the principles of weak higher categories, we want there to be a non-trivial proof or witness that relates (certain) composites of morphisms instead of their strict equality. Since the above cell complexes are indistinguishable, it is not immediately clear what this proof could look like. This is a short-coming of the cellular approach. 

In the setting of ``globular" $2$-fold categories, $\alpha$ and $\beta$ take the form of $2$-cubes as follows
\begin{restoretext}
\begingroup\sbox0{\includegraphics{test/page1.png}}\includegraphics[clip,trim=0 {.35\ht0} 0 {.35\ht0} ,width=\textwidth]{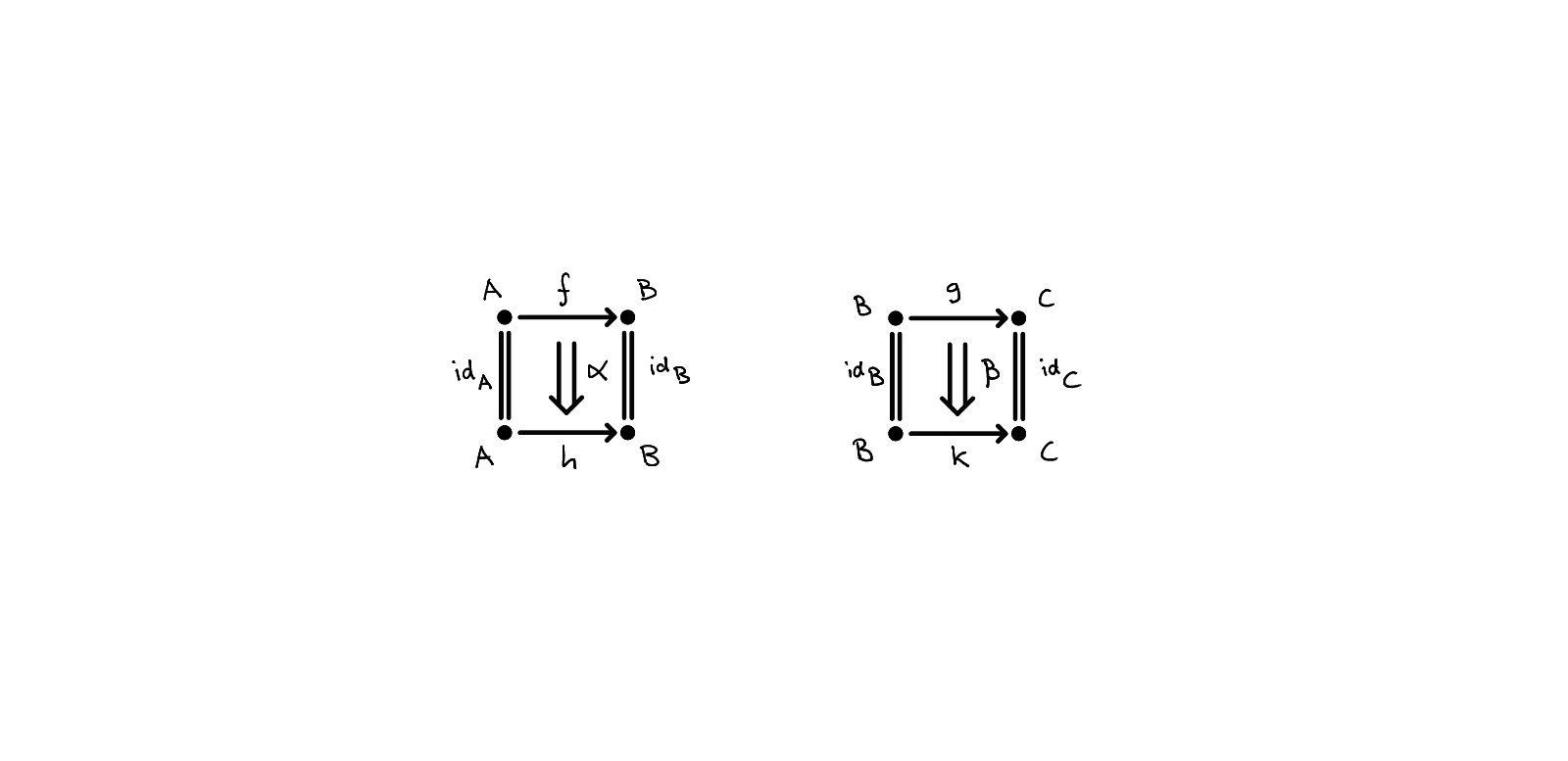}
\endgroup\end{restoretext}
Now, composing $\alpha$ and $\beta$ along their mutual $0$-boundary (or ``vertical" $1$-boundary) yields
\begin{restoretext}
\begingroup\sbox0{\includegraphics{test/page1.png}}\includegraphics[clip,trim=0 {.37\ht0} 0 {.35\ht0} ,width=\textwidth]{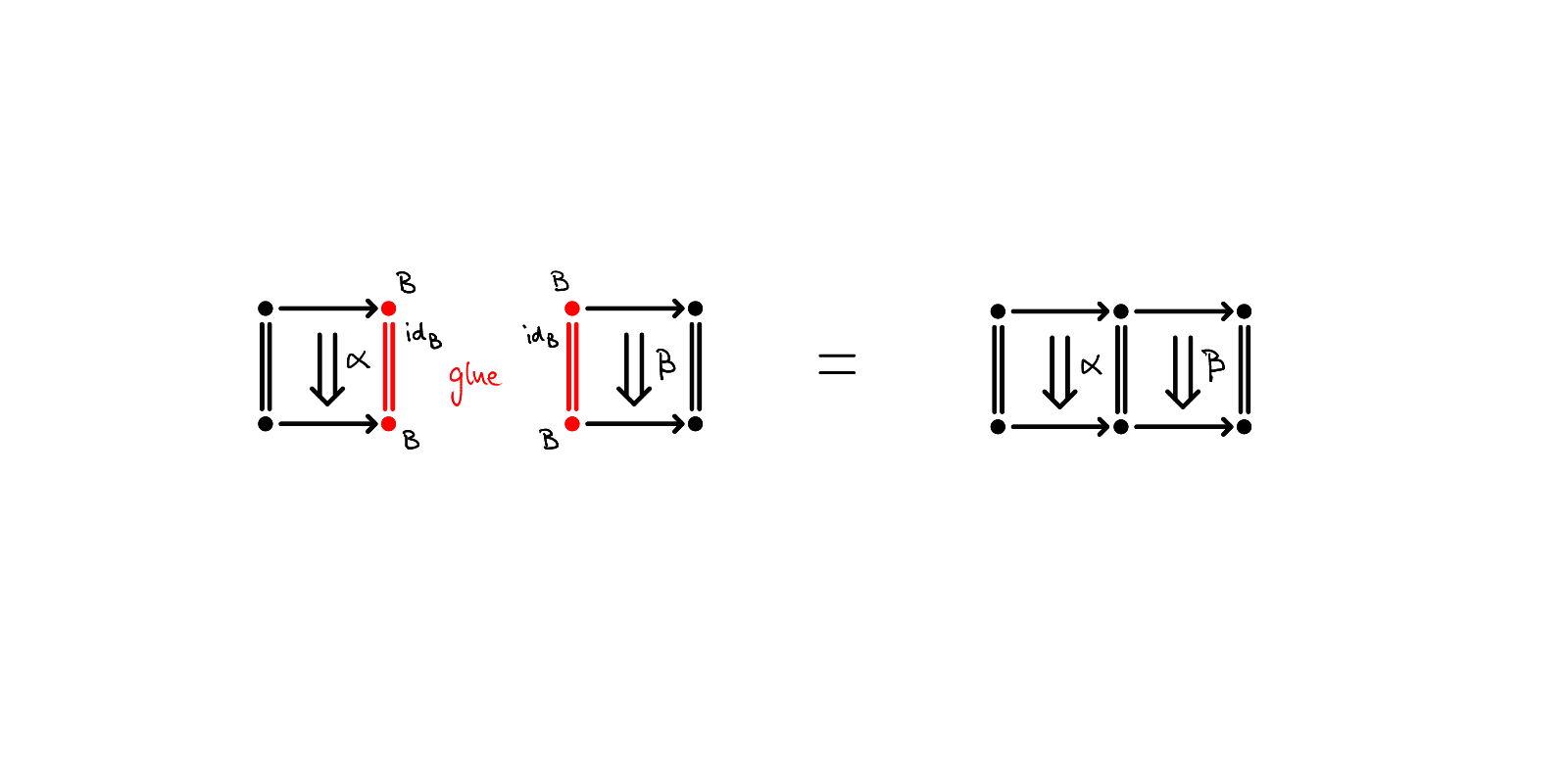}
\endgroup\end{restoretext}
On the other hand, post-whiskering $\alpha$ and $g$, and pre-whiskering $h$ and $\beta$ gives
\begin{restoretext}
\begingroup\sbox0{\includegraphics{test/page1.png}}\includegraphics[clip,trim=0 {.3\ht0} 0 {.1\ht0} ,width=\textwidth]{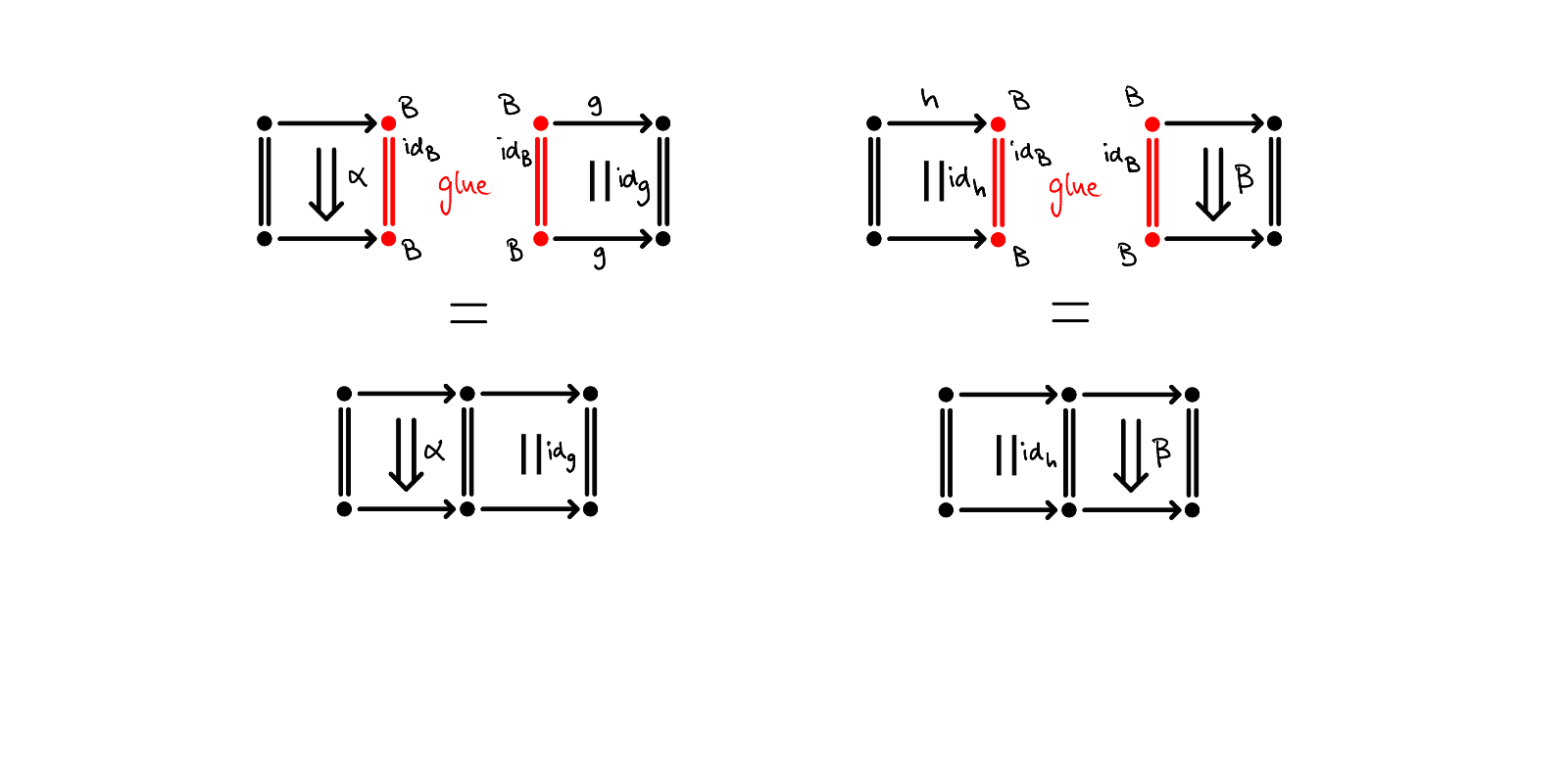}
\endgroup\end{restoretext}
which compose along their $1$-boundary $gh$ to give
\begin{restoretext}
\begingroup\sbox0{\includegraphics{test/page1.png}}\includegraphics[clip,trim=0 {.24\ht0} 0 {.24\ht0} ,width=\textwidth]{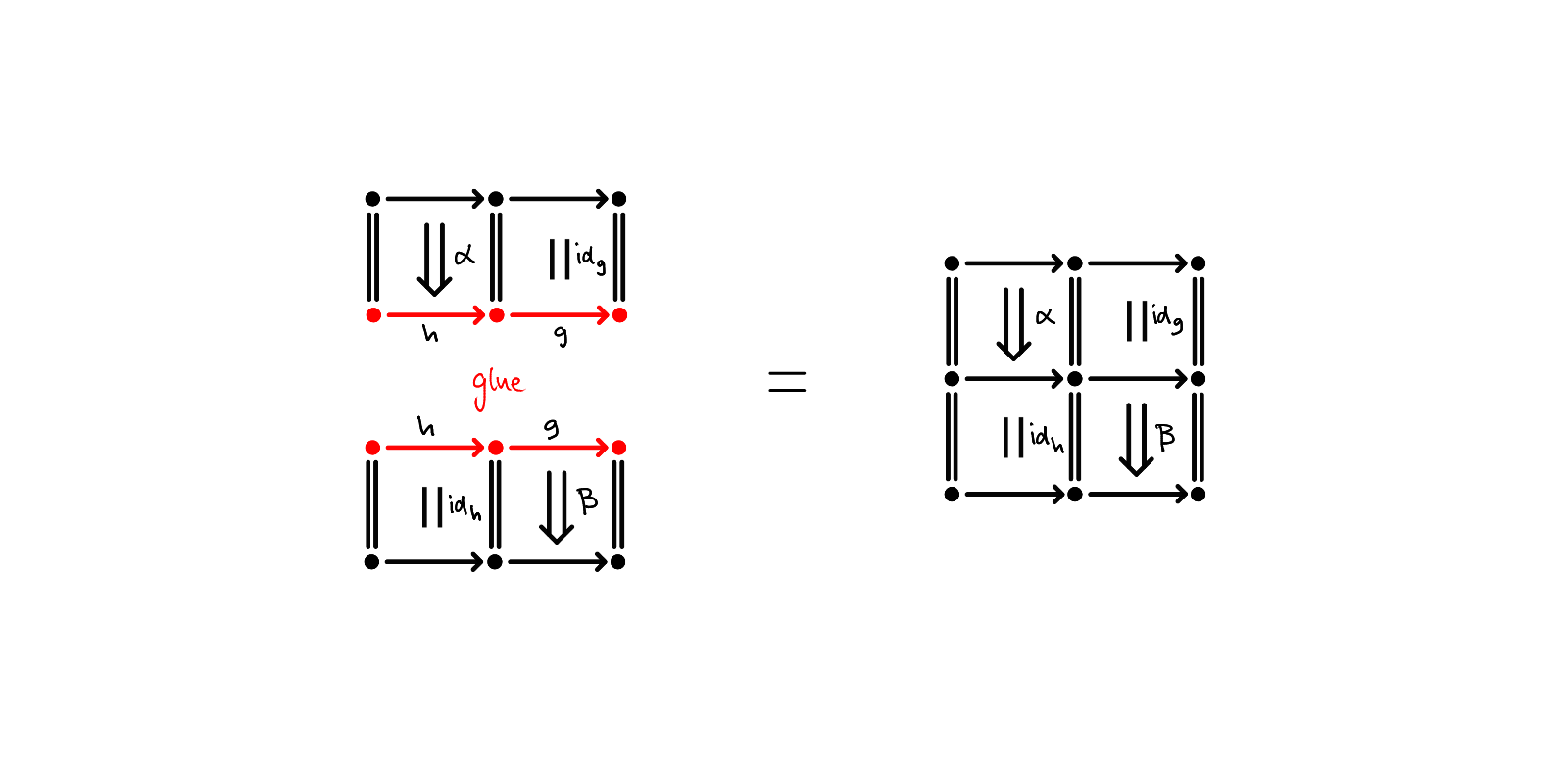}
\endgroup\end{restoretext}
Similarly, pre-whiskering $f$ with $\beta$, post-whiskering $\alpha$ with $f$ gives
\begin{restoretext}
\begingroup\sbox0{\includegraphics{test/page1.png}}\includegraphics[clip,trim=0 {.25\ht0} 0 {.1\ht0} ,width=\textwidth]{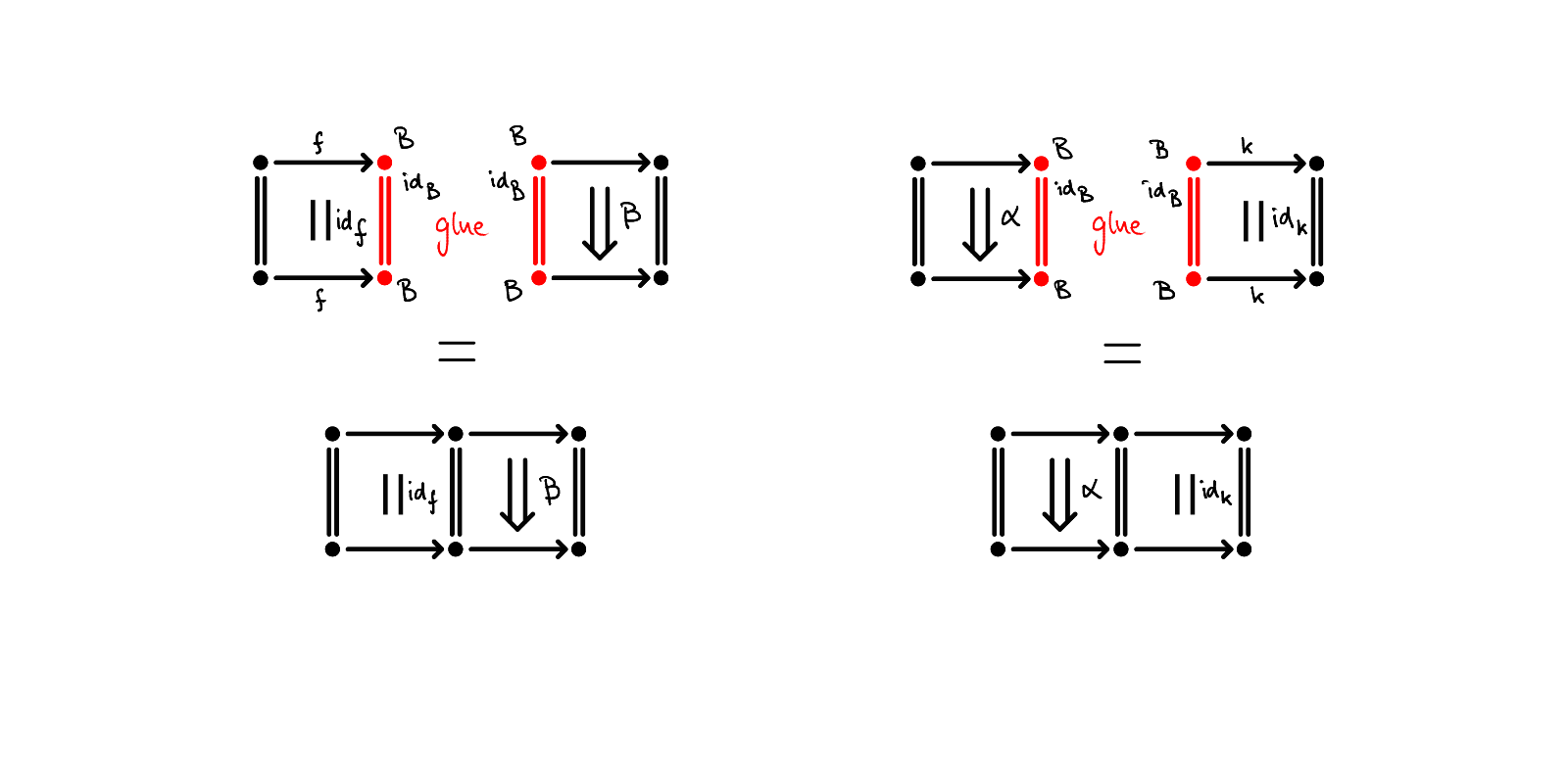}
\endgroup\end{restoretext}
And composing the results along their mutual $1$-boundary $kf$ we find
\begin{restoretext}
\begingroup\sbox0{\includegraphics{test/page1.png}}\includegraphics[clip,trim=0 {.25\ht0} 0 {.2\ht0} ,width=\textwidth]{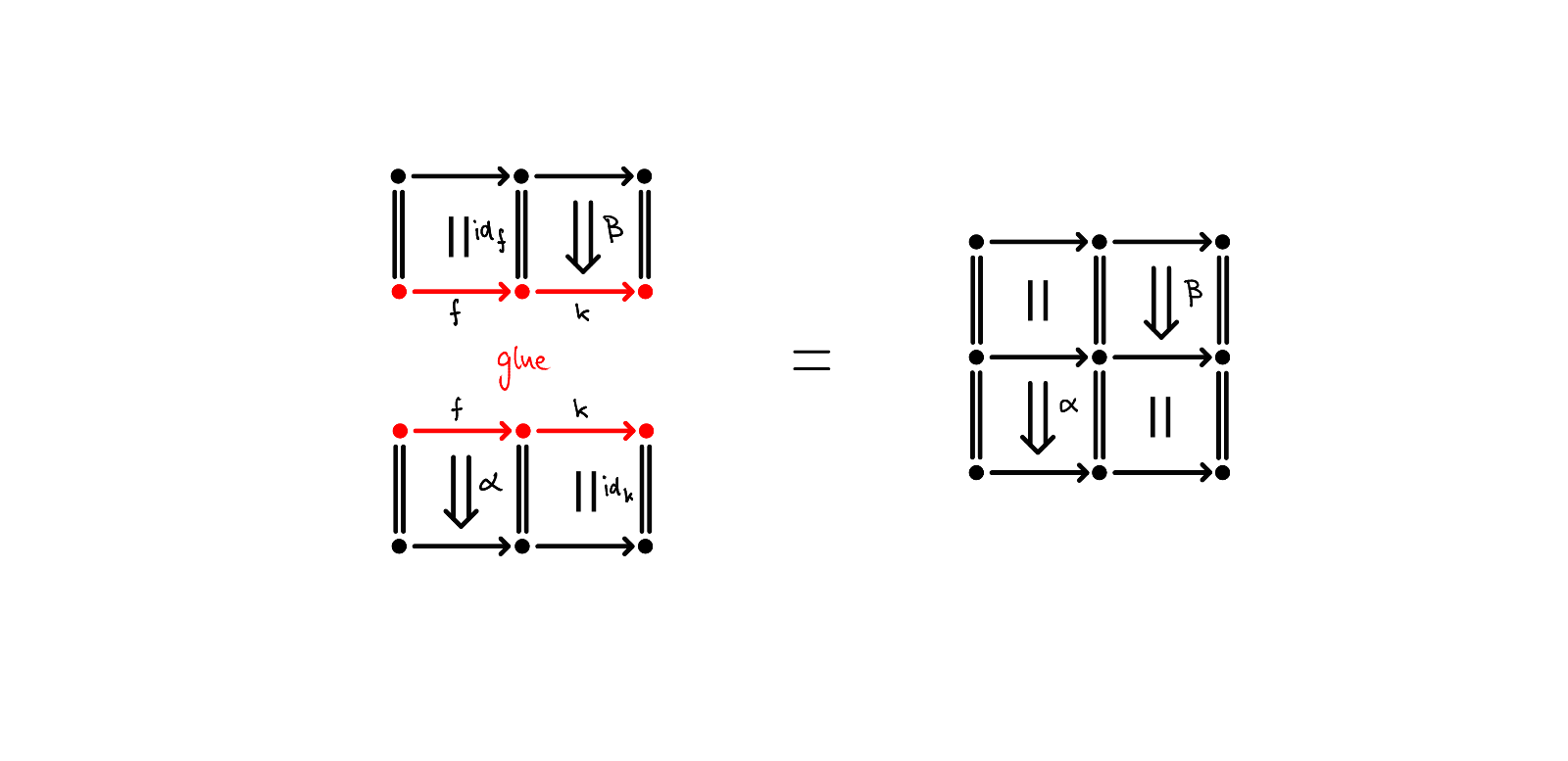}
\endgroup\end{restoretext}
The resulting composites thus have \textit{distinguishable} cell complexes. The coherence relating them is called the \textit{interchanger}. The ``underlying data" of the interchanger, which we proclaimed to be difficult to express in the cellular approach, will be given by proofs that the following transformations (indicated by arrows) are valid
\begin{restoretext}
\begingroup\sbox0{\includegraphics{test/page1.png}}\includegraphics[clip,trim=0 {.35\ht0} 0 {.3\ht0} ,width=\textwidth]{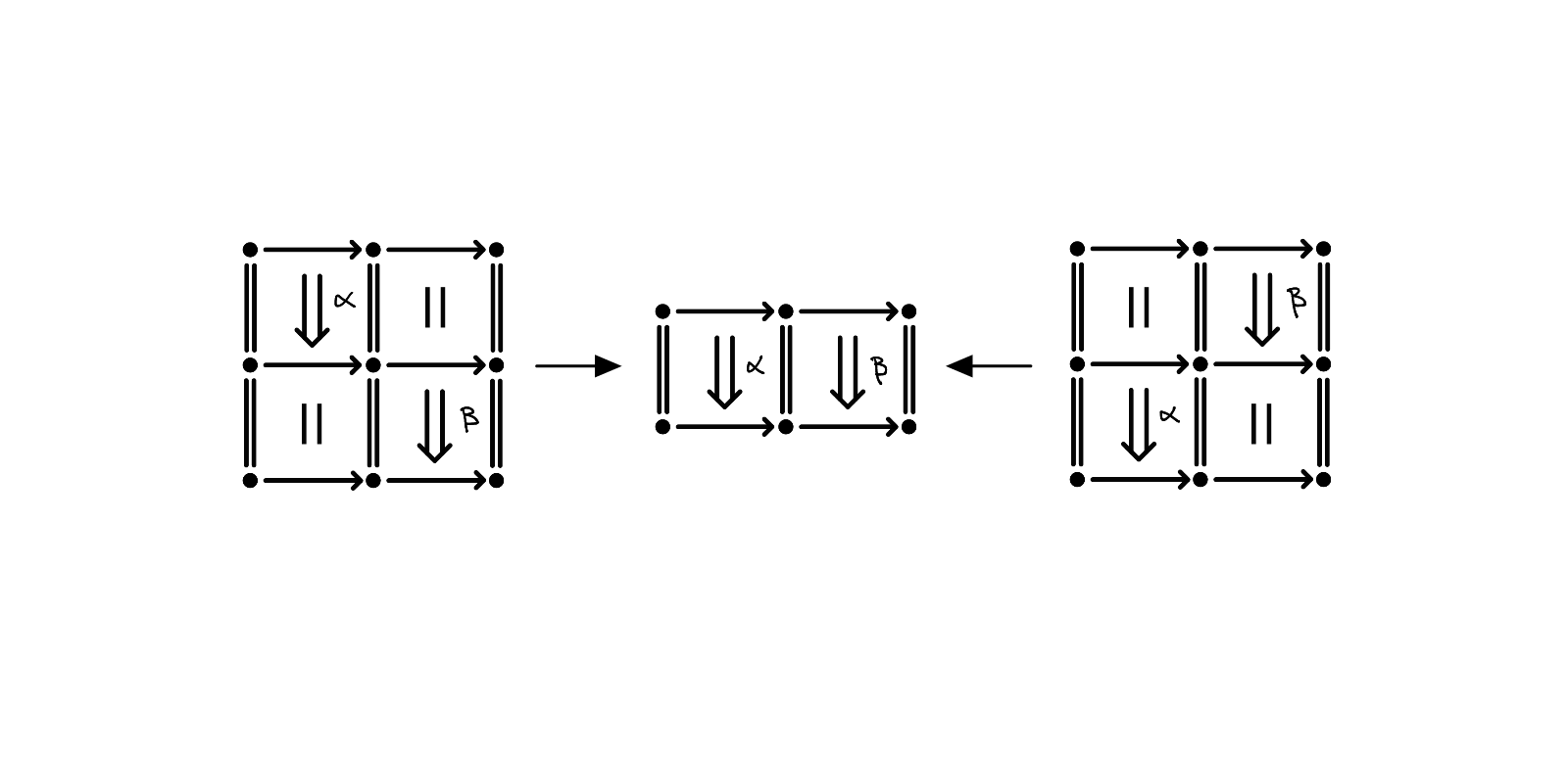}
\endgroup\end{restoretext}
We will formulate these proofs in a rigorous algebraic manner. However, they are best described in geometric terms. First observe, representing of the above composites in \textit{string diagrams} (by turning cells into strata of dual dimension, see e.g. \cite{selinger2010survey}) faithfully distinguishes the above three composites. Namely, using string diagrams, the three different composites above can be represented as
\begin{restoretext}
\begingroup\sbox0{\includegraphics{test/page1.png}}\includegraphics[clip,trim=0 {.3\ht0} 0 {.3\ht0} ,width=\textwidth]{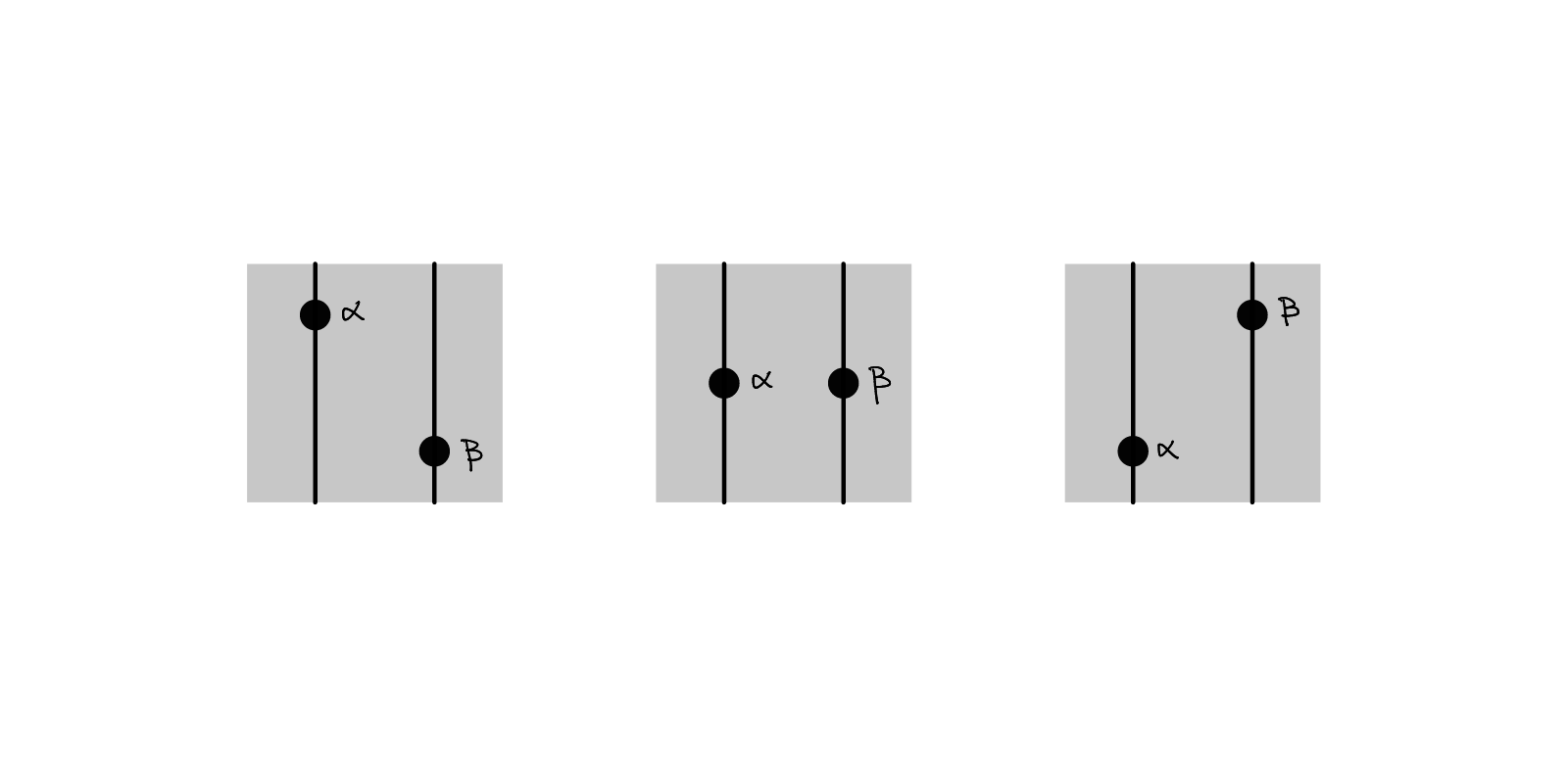}
\endgroup\end{restoretext}
While all three string diagrams are isotopic in the sense of the classical notion of string diagrams, differences in the coordinates describing the positions of $\alpha$ and $\beta$ allow us to distinguish the three cases that have previously been conflated in the cellular approach. 

Now the data of the proclaimed ``geometric proof" showing equivalence (but not strict equality) of the above three composites, takes the form
\begin{restoretext}
\begingroup\sbox0{\includegraphics{test/page1.png}}\includegraphics[clip,trim=0 {.3\ht0} 0 {.3\ht0} ,width=\textwidth]{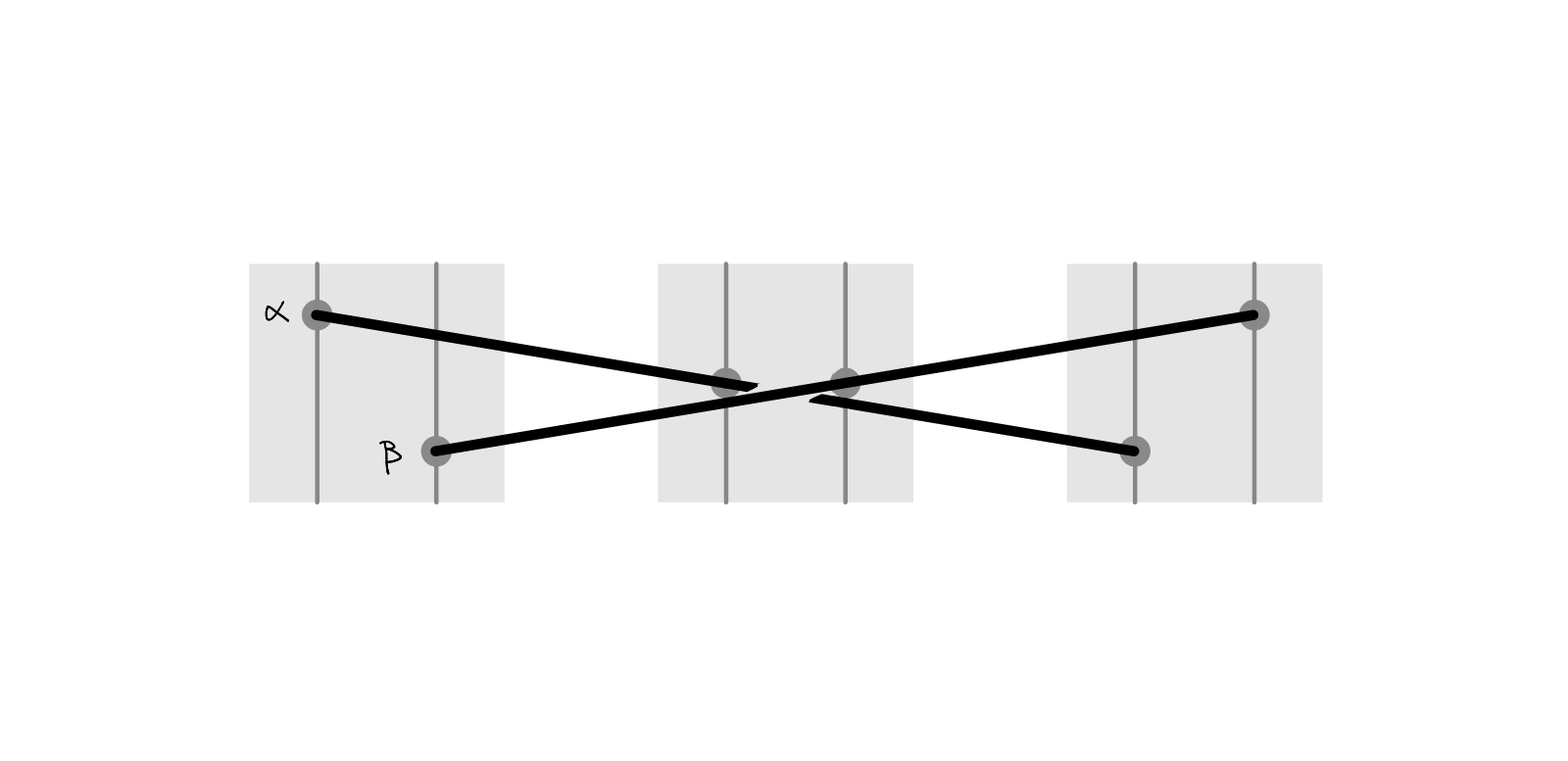}
\endgroup\end{restoretext}
This geometric object is an illustration of a \textit{surface diagram} (a manifold diagram in dimension $3$). This surface diagram is an instance of an ``interchange homotopy" which as mentioned plays a central role for Gray categories. 

In this work, we will generalise string diagrams to arbitrary dimensions $n$, which will lead to a notion of \textit{$n$-manifold diagram}. More importantly, we will develop a language to algebraically represent manifold diagrams, including coherences (or ``homotopies") such as the interchange. We re-emphasize that the ``cubical" structure of $n$-manifold diagrams, and the ``cube" terminology used throughout this work, ultimately is reflected in the fact that the language we develop is a language of $n$-fold categories.

\section{Results and conjectures}

The next two subsections respectively contain a less wordy and a more wordy summary of the results in this thesis. The third and final subsection lists the main conjectures we are making.

\subsection{Core results} \label{ssec:fast_facts}

This work contains the following

\begin{itemize}
\item A definition of $\Cubeo n \cC(\bnum 1)$ (given \autoref{rmk:cube_cat}), capturing the category of (topological) $\cC$-labelled $n$-cubes which are ``flag-foliation-compatible stratifications" of the $n$-cube together with functorial data on the strata in some fixed category $\cC$. The definition in particular specialises to a novel definition of \textit{$n$-manifold diagrams} (cf. \autoref{ssec:po_mfld_diag}).
\item A combinatorial analogue of $\Cubeo n \cC(\bnum 1)$, denoted by $\Buno n \cC (\bnum 1)$. Objects of $\Buno n \cC(X)$ are called $\cC$-labelled singular $n$-cube \textit{families} over (a poset) $X$, or $\SIvert n \cC$-bundles over $X$ for short. If we let $X$ vary we obtain a category $\Buno n \cC$ (cf. \autoref{sec:sum_Buno} and \autoref{sec:sum_Buno2}).
\item A definition of a geometric realisation $\norm{-} : \Buno n \cC (\bnum 1) \to \Cubeo n \cC(\bnum 1)$ at the level of \textit{objects} (cf. \autoref{ssec:coloring}). The description of a full functor should be equally possible (modulo slight modifications) but exceeds the scope of this thesis.
\item A proof that $\Buno n \cC$ has an (epi,mono) factorisation system. Monos are also called \textit{embeddings} (cf. \autoref{sec:sum_comm}), epis are also called \textit{collapse}. The former, denoted by $A \mono B$, describes when $A$ embeds as a subcube in $B$. The latter, denoted by $A \mcoll B$, describe when $A$ is a refinement (of the stratification) of $B$\footnote{Or conversely, strata of $A$ ``collapse" to yield strata of $B$, but the functorial data remains the same}. 
\item A combinatorial definition of ``geometric equivalence": Cubes $A$, $B$ are geometrically equivalent, written $A \simeq B$, if they refine the same cube, that is, there is a cospan $A \mcoll C \twoheadleftarrow B$. We introduce an analogous notion of equivalence $\simeq$ for topological $\cC$-labelled cubes and conjecturally $\norm{-}$ then preserves and reflects equivalence (cf. \autoref{ssec:coloring} and \autoref{sec:sum_NF}).
\item A proof of a ``normal form theorem": every equivalence class of cubes as a unique maximally collapsed representative called the normal form. Morally, this means that $\norm{A} \simeq \norm{B}$ if and only if $A$ and $B$ have the same normal form. This is of fundamental importance for our work being ``foundation-independent", since it makes geometric equivalence decidable (cf. \autoref{sec:sum_NF}).
\item A construction of minimal embeddings in $\Buno n \cC$ and proof of their minimality (cf. \autoref{sec:sum_min}).
\item An endofunctor $\SIvertone {-} : \Cat \to \Cat$ such that $\SIvert n \cC$ is the classifying category of $\cC$-labelled singular $n$-cube bundles (cf. \autoref{ssec:SIvert}).
\item A construction of an embedding $\sL$ of $(\Buno n \cC(\bnum 1),\mathrm{epi})$ (denoting the wide subcategory of all epis) into $\SIvert n \cC$ (this follows from results in \autoref{ch:globes}).
\item Several definitions of higher categories, both ``foundation-independent" (and computer-implementable) as well as ``foundation-dependent" (in a classical Set-theoretic setting) (cf. \autoref{sec:sum_pres} and \autoref{sec:sum_cat}).
\item An informal discussion of a connection of our notion of $\infty$-groupoids and CW-complexes (cf. \autoref{sec:sum_CW}).
\end{itemize}

We note that $\Buno n \cC$ itself will not be very prominent in the main body of the thesis, as for the most part we will be working with decompositions of its  morphisms (also called \textit{multi-level basechanges}) into their ``dimension-wise" components (also called \textit{$k$-level basechanges}). This will enable our (often very explicit and pedestrian but straight-forward) computations and proofs.

\subsection{Overview of work}

A more wordy, and more chronological, description of what is done in this thesis is the following: We define associative $n$-categories, \free{} associative $n$-categories and \free{} associative $n$-groupoids. In the appendices we further define notions of \free{} associative $n$-fold categories, and \free{} weak $n$-categories. On our way, we describe the following tools and observations:

\begin{enumerate}
\item We develop a theory of iterated Grothendieck fibrations over posets, called towers of $\SI$-bundles, whose individual fibers live in a category $\SI$. The latter is called the category of \textit{singular intervals} and is a subcategory of the category of $\Bool$-enriched profunctors (\autoref{ch:intervals}). Towers of $\SI$-bundles together with a ``labelling" functor (from the top bundle's total space to a category $\cC$)  give rise to the notion of \textit{$\cC$-labelled singular $n$-cube families} (\autoref{ch:cubes}), or $\SIvert n \cC$-cube families for short. We show that such families can be equivalently defined as functors from posets into a certain category $\SIvert n \cC$, called \textit{classifying category of $\SIvert n \cC$-cube bundles}. %

\item We note that $\SIvert n \cC$-cube families support a notion of \textit{multi-level base change} giving rise to a category $\Bunbc^n_\cC$ ($\Buno n \cC$ is a subcategory thereof, whose morphisms are required to be \textit{open}). As a first special instance of multi-level basechange we define a notion of  \textit{collapse} (whose individual base change maps will be required to be open and surjective). We show that any cube family has a unique ``maximal" collapse which leads to its \textit{normal form} (\autoref{ch:collapse_intervals} and \autoref{ch:collapse_cubes}). More abstractly, this means connected components of the subcategory of $\Bunbc^n_\cC$ generated by collapses have terminal objects. %

\item As a second special instance of multi-level basechange, we define a notion of \textit{embeddings} of $\SIvert n \cC$-cubes  (whose individual base change maps will be required to be open and injective). In particular, this will lead to a notion of \textit{minimal neighbourhoods} (\autoref{ch:emb}). %
 We prove that embedded cubes inherit collapses from their parent cube. More abstractly, this means the subcategory of $\Bunbc^n_\cC$ generated by collapses and embeddings admits an (epi,mono) factorisation system.

\item We introduce a condition of \textit{globularity} (\autoref{ch:globes}), which, in geometric terms, forces constancy of labelling data on the boundary of the $n$-cube so that the $n$-cube can be quotiented to an $n$-globe. In fact, our condition will be stable under taking subcubes, such that every subcube of a globular cube will be globular again. To facilitate proofs about the properties of globular cubes we will also construct an embedding of $(\Buno n \cC (\bnum 1),\mathrm{epi})$ into $\SIvert n \cC$. Finally, we show that the ``adjoin-a-terminal-object"-monad extends to $\SIvert n \cC$-cubes. We use it to define conical globular $n$-cubes, also called double cones, obtained by forming a ``cone" between two compatible globular $(n-1)$-cubes.

\item In \autoref{ch:presented} we define \free{} associative $n$-categories, which will be valid collections of generating morphisms (or, in type-theoretic lingo, constructors). Validity will be subject to a condition of \textit{well-typedness}. Geometrically, this condition guarantees regions labelled by the same generating morphisms have globally constant minimal local neighbourhoods, namely, they equal the so-called \textit{type} of the generating morphism.  

\item In \autoref{ch:groupoids} we then introduce candidates for a \textit{theory of coherent invertibility} $\TI$ as certain \free{} associative $\infty$-categories. This will allow us to speak about general invertible elements. \Free{} associative $n$-groupoids are \free{} associative $n$-categories with the condition that all generators are invertible. (In \autoref{ch:geom} we conjecture that morphisms of a specific candidate theory for $\TI$ yield combinatorial representations of framed $k$-tangles. We further outline a procedure for translating CW-complexes into \free{} associative $\infty$-groupoids using a generalised Thom-Pontryagin construction).

\item In \autoref{ch:composition} we describe ``perturbation-stable" morphisms, which will be called \textit{generic composites}. The notion intends to capture morphisms whose (collapse) equivalence class is stable under small perturbations of subcubes. We give two characterisations, the first being a simple ``genericity" condition on morphisms. For the second one we first show that cubes can be stacked (along their sides), and that globular cubes can be \textit{whiskered} (which corresponds to the usual whiskering compositions of globes). This then allows us to characterise generic composites as being inductively build from generators and ``elementary homotopies" via whiskering. We sketch a proof of equivalence of the two characterisations.

\item Finally, in \autoref{ch:associative} we suggest a definition for associative $n$-categories as unital algebras to the ``generic composites endofunctor" $\GComp$ on globular sets $\globset$, satisfying equations (which in particular, make the definition ``associative"). It might come as a disappointment to some that $\GComp$ does not have an  immediate monad structure, and thus the above mentioned equations are (concisely) enforced by so-called \textit{resolutions} of the algebra structure instead. A lucky feature of this technique is, however, that it allows to control the exact ``amount of weakness" that we want, opening up a spectrum of definitions from ``fully associative" to ``fully weak". No attempt is made to prove any form of correctness results for our chosen definition. However, by means of examples we motivate in which way associative $n$-categories can be thought of ``Gray $n$-categories".
\end{enumerate}

\subsection{Conjectures}

We also make the following conjectures

\begin{conj} Every \free{} associative $n$-category is naturally an associative $n$-category, up to some arbitrary (contractible) choices for morphisms witnessing homotopies.
\end{conj}

\noindent This is a conjecture because a priori \free{} associative $n$-category and associative $n$-category are very different structures. It will be further discussed in \autoref{ssec:pres_conj}.

\begin{conj} There is a one-to-one correspondence of presentations of CW-complexes (up to homotopy of the attaching maps) and \free{} associative $\infty$-groupoids (up to choices of ``generalised direction" for each generator).
\end{conj}

\noindent This conjecture will be discussed again in \autoref{ssec:homotopy_hyp} after a translation of CW-complexes into \free{} associative $\infty$-groupoids has been sketched.

\begin{conj} The theory of associative $n$-categories is equivalent to a theory weak $n$-categories.
\end{conj}

\noindent Preliminary evidence for this conjecture is discussed in \autoref{sec:weak_assoc_equiv}, where a theory of \free{} weak $n$-categories will be constructed.

\section{Related work}

\subsection{Algebra} In low dimensions, associative $n$-categories are related to existing concepts as follows
\begin{enumerate}
\item Associative $0$-categories are sets.
\item Associative $1$-categories are unbiased $1$-categories (cf. \cite{leinster-operads}). Here, the word ``unbiased" refers to the presence of composition operations of arbitrary arity (instead of just binary compositions).
\item Associative $2$-categories are unbiased strict $2$-categories (cf. \cite{leinster-operads}).
\item Associative $3$-categories are unbiased Gray-categories. Note that the ``unbiased" predicate here not only refers to composition but also extends to homotopies: instead of just the binary interchanger used in biased Gray-categories, we also have $n$-ary interchangers and other types of interchangers.
\end{enumerate}

\noindent Similarly, \free{} associative $n$-categories can be understood as follows in low dimension.

\begin{enumerate}
\item \Free{} associative $(-1)$-categories are booleans.
\item \Free{} associative $0$-categories are setoids (that is, a set with an equivalence relation).
\item \Free{} associative $1$-categories are presentations of categories given by objects, generating morphism with generating equality relations between their composites.
\item \Free{} associative $2$-categories are presentations of strict $2$-categories, given by objects, generating $1$- and $2$-morphisms and generating equality relations between their ($2$-dimensional) composites.
\item \Free{} associative $3$-categories are presentations of (unbiased) Gray-categories, given by objects, generating $1$-, $2$-  and $3$-morphisms and generating equality relations between their ($3$-dimensional) composites.
\end{enumerate}

\subsection{Geometry}

A notion of string diagram has been made precise by Joyal and Street in \cite{joyalGTC1} and \cite{joyal1991geometry}, and it was shown to correctly capture various notions of monoidal categories (that is one-object $2$-categories) and monoidal categories with duals. An attempt to formally capture the notion of surface diagrams using stratified Morse Theory was started by Trimble and McIntyre in 1996 (\cite{mcintyre1997surface}), but was left in an unfinished state in 1999 due to technical difficulties in their approach. The motivation for building such a theory of surface diagrams was the hope that they would provide a natural model for Gray categories, such that each composite of morphisms would have an associated surface diagram. Later on in work by Barrett, Meusburger and Schaumann \cite{barrett2012gray} as well as in work by Hummon \cite{hummon2012surface} notions of surface diagrams have been made precise. 

\subsection{Generalising the homotopy hypothesis}

As a generalisation of the homotopy hypothesis, we can ask: what notion of space do $n$-categories or $(\infty,n)$-categories correspond to? The answer should be some type of ``directed space" as discussed in the beginning of this thesis. In \cite{ayala2015stratified}, Ayala, Francis and Rozenblyum (see also \cite{ayala2017cobordism}) seek such a generalisation of the homotopy hypothesis for $(\infty,n)$-categories. Concretely, they establish a ``factorization homology" functor
\begin{equation}
\int : \Cat_{(\infty,n)} \to \Fun(\cMfld^{\mathsf{vfr}}_n, \Spaces)
\end{equation}
from the $(\infty,1)$-category $\Cat_{(\infty,n)}$ to the $(\infty,1)$-category of functors from $\cMfld^{\mathsf{vfr}}_n$ into spaces: here, $\cMfld^{\mathsf{vfr}}_n$ is the $(\infty,1)$-category of so-called \textit{vari-framed compact $n$-manifolds}, which roughly are certain stratified manifolds with compatible framing. They then prove that the factorization homology functor $\int$ is fully-faithful, and the image of $\int$ thus gives an answer to the generalised homotopy hypothesis. Comparing to our setting, the role of vari-framed manifolds is played by manifold diagrams: they are the building blocks from which directed spaces are built by gluing. A detailed comparison of vari-framed manifolds to manifold diagrams is left to future work.

\section{Notes to the reader} \label{sec:notes}

The reader might want to take note of the following:
\begin{itemize}
\item Reading the \hyperref[ch:summary]{Summary} should be sufficient for most readers. \autoref{ch:prerequisites} to \autoref{ch:globes} mainly develop technical tools to analyse the concepts introduced in the summary. Conversely, the \hyperref[ch:summary]{Summary} should not be skipped and provides the narrative for most of our developments. For space reasons a couple of interesting ideas and examples from the later chapters (\autoref{ch:presented} to \autoref{ch:associative}) and in particular the appendices have not been reproduced in the summary. Up to a few definitions, it is possible to look at these chapters, which are probably the most interesting to many readers, directly after having read the \hyperref[ch:summary]{Summary}.

\item The reader should expect to find completely elementary and largely self-contained mathematics in this work (maybe with the exception of \autoref{ch:geom}). Most proofs are straight-forward. Since the ``degree of straight-forwardness" is usually subjective, we decided to include all proofs, but marked proofs which can probably be safely skipped by most readers with the key-word ``\stfwd{}".

\item The theory of \free{} associative $n$-categories has a ``natural geometric model" by translating morphisms into manifold diagrams. We will not attempt to prove this statement, however we will make reference to manifold diagrams on multiple occasion, and some parts will be entirely devoted to explaining the connections of algebra and geometry. \textit{We emphasize that whenever we enter the realm of geometry in this work, arguments should not be expected to be fully formal, but rather regarded as providing a guiding intuition for the algebraic work}. Further, the word ``manifold" should be read as \textit{piecewise linear} manifold throughout the document with the exception of \autoref{ch:geom} where we will assume manifolds to be \textit{smooth} unless otherwise stated. In fact, it seems plausible to conjecture that in $\Cubeo n \cC(\bnum 1)$ ``smooth" and ``piecewise linear" can be used interchangeably (up to a notion of equivalence, see \autoref{ssec:po_mfld_diag}), and in illustrations we will often depict either notion depending on which is more convenient.

\item For the reader sensitive to foundational issues, we remark that throughout the document we have chosen a concrete representation (namely, by natural numbers) of the core objects (namely, singular intervals) at play. This allows us to write ``strict equalities" in many places where one would otherwise expect ``isomorphisms". The advantage of this concrete representation is that our definitions are directly amenable to computer implementation. The potential disadvantage is that some of our arguments could possibly be simplified by working more abstractly.

\item The length of this thesis should by no means deter the potential reader, as it is mainly a consequence of the presence of many pictures, examples and often very detailed (albeit elementary) proofs and computations.

\item To help readability, a glossary has been added to this thesis which contains the majority of terms and symbols used in it. The glossary can also be accessed by clicking on the respective symbols, which are hyper-linked to their corresponding glossary entries.

\end{itemize}

\cleardoublepage
\phantomsection
\addcontentsline{toc}{chapter}{Summary}
\setcounter{section}{0}
\renewcommand{\thesection}{S.{\arabic{section}}}
\renewcommand{\theHsection}{S.{\arabic{section}}}

\chapter*{Summary} \label{ch:summary}

The following is a concise overview of the theory of (presented) associative higher categories, and the tools needed for their development. An attempt is made to give a treatment which is mathematically as complete as possible, while focusing only on the most relevant and interesting aspects of the theory.

While some important details will have to be deferred to later chapters, most parts of the later chapters will be dedicated to constructing the sometimes complex language for otherwise completely straight-forward proofs and computations. The reader not interested in these details could skip directly to \autoref{ch:presented}, \autoref{ch:groupoids} and \autoref{ch:associative} after having read the summary. Conversely, some observations and definitions will not be repeated in later on chapters and are only stated here, which means this chapter should be read before proceeding to the main body of the thesis.

\section{Geometric model} \label{ssec:po_mfld_diag}

In this section we discuss the local geometric model that our definition of higher categories will be based on: manifold diagrams, and further, their generalisation to so-called (topological) $\cC$-labelled $n$-cubes. While the rest of this thesis (with the exception of \autoref{ch:geom}) will not explicitly depend on the geometric model, for expositional reasons we decided to include this brief discussion of the geometry as early as possible: it will provide a guiding intuition for the reader in the later fully algebraic arguments.

We will summarise most relevant ideas of the geometric model, but some further discussion of ideas mentioned here can be found in \autoref{ssec:coloring} (geometric realisation of combinatorial structures), \autoref{ssec:regions} (combinatorial \stratatype{}s) and \autoref{ch:geom} (generalised Thom-Pontryagin construction). 

\subsection{Flag-foliation-compatible stratifications}

We start with the basic observation that Euclidean space $\lR^n \iso (0,1)^n$ is \textit{flag foliated}: that is, it has a codimension 1 foliation as $(0,1)^n = (0,1)\times (0,1)^{n-1}$ and recursively each sheet $(0,1)^{n-1}$ is in turn flag foliated (for the base case, note that every $0$-manifold is by definition flag foliated). The \textit{$k$-sheets} $(0,1)^k \times \Set{p}$, $k \leq n$ of the flag foliation of $(0,1)^n$ are parametrised by $p \in (0,1)^{n-k}$. The order of components in $(0,1)^n = (0,1) \times ... \times (0,1)$ will be often indicated by coordinate axis in figures below.

\begin{notn}[Ordinal numbers] \label{notn:simplicies}
Let $\bnum n = \Set{0 < 1 < ... < n-1}$ be the $n$th ordinal category. 
\end{notn}
\noindent In particular, $\bnum 2 = (0 \to 1)$ is the interval category. Note for $n \in \lN$ that there is a functor $\codim : \bnum 2^n \to (\bnum {n+1})$ mapping $u = (u_1, ..., u_{n}) \in \bnum 2^n$ to the sum $\sum_i u_i$ of the components of $u$. 

\begin{defn}[Region types] An open submanifold $M$ of $(0,1)^n$ is said to be of \textit{\stratatype{}} $u \in \bnum 2^n$, $\codim(u) = n-m$, if 
\begin{enumerate}
\item For every $k$-sheet $S$, $S \cap M$ is a submanifold of $S$
\item $u_k = 0$ implies that, for any $k$-sheet $S$, $S \cap M$ intersects each $(k-1)$-sheet of $S$ transversally (in $S$)
\item $u_k = 1$ implies that there is $\eps > 0$ such that for all $x = (x_1, ..., x_{n}) \in M$ we have $\eps < x_k < 1-\eps$.
\end{enumerate}
Let $\overline{M}$ denote the closure of $M$, then the \textit{boundary} $\partial M$ is defined to be $\overline{M}\setminus M$ (and a similar definition can be given  replacing $(0,1)^n$ by a general manifold $X$).
\end{defn}
We give examples in $(0,1)^3$ as follows
\begin{restoretext}
\begingroup\sbox0{\includegraphics{ANCimg3/empty.png}}\includegraphics[clip,trim={.0\ht0} {.25\ht0} {.0\ht0} {.17\ht0} ,width=\textwidth]{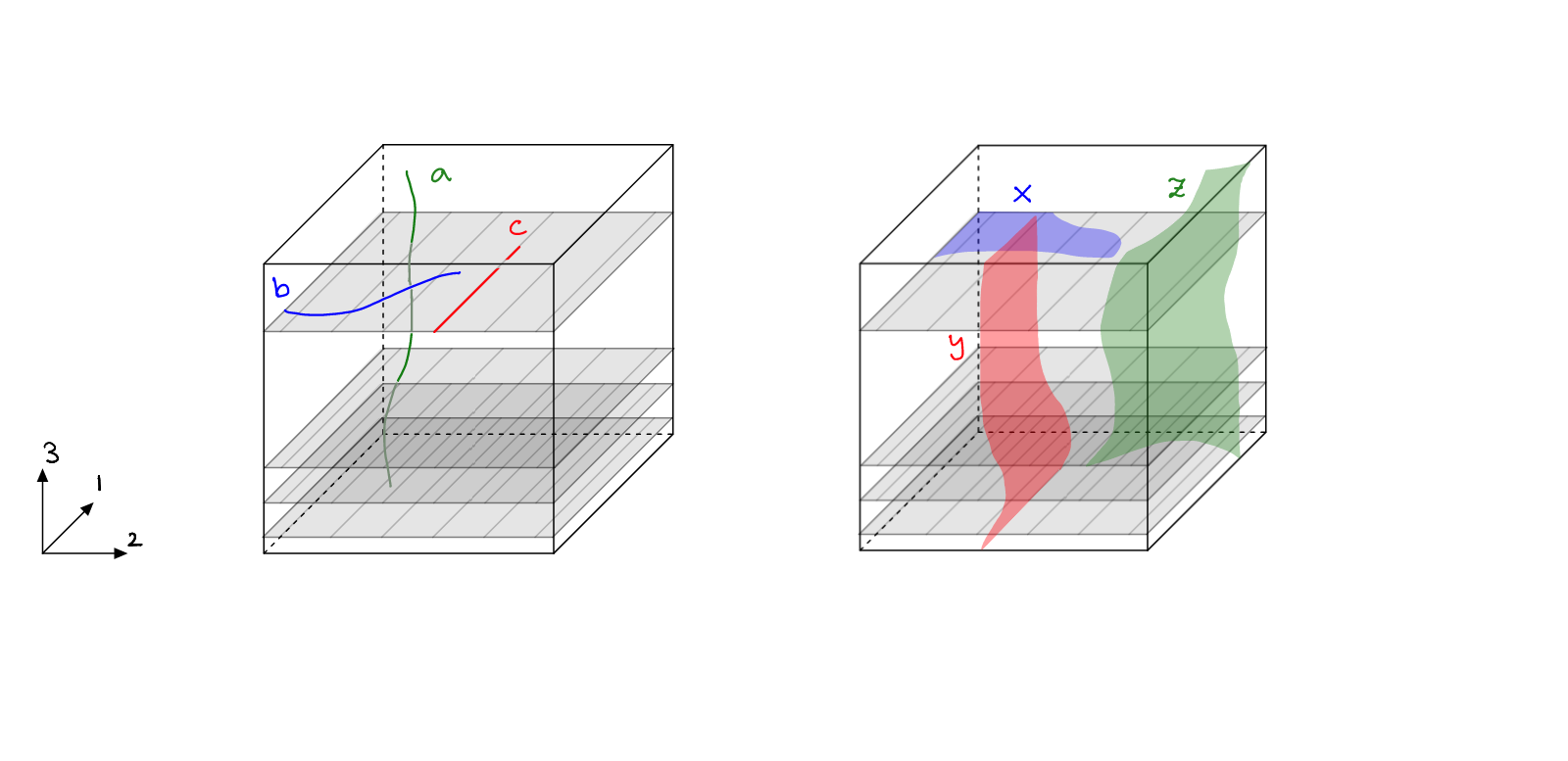}
\endgroup\end{restoretext}
On the left, the $1$-manifolds $a,b,c$ are of type $(1,1,0)$, $(1,0,1)$ and $(0,1,1)$ respectively. On the right, the $2$-manifolds $x,y,z$ are of type $(0,0,1)$, $(0,1,0)$, and $(1,0,0)$ respectively.

\begin{defn}[Flag-foliation-compatible stratifications] \label{defn:ffcs} A (finite) \textit{stratification} $A$ of a manifold $X$ is a decomposition of $X$ into finitely many connected, open submanifolds $\Set{p^A_i}_{i \in I}$, such that the boundary of each $p^A_i$ is given by a union of lower-dimensional $p^A_j$. The \textit{total poset} of the stratification $A$, denoted by $\sG(A)$, is defined to have object set $I$ and morphisms $i \to j$ are given if $p^A_j \subset \partial p^A_i$.

If $X = (0,1)^n$ then $A$ is said to be \textit{flag-foliation-compatible} stratification of the $n$-cube if each stratum $p^A_i$ is of some type $u(p^A_i) \in \bnum 2^n$.
\end{defn}
 \noindent As an example, consider the following flag-foliation compatible stratification of the $2$-cube on the left, which we will refer to as $A_0$, and its corresponding total poset $\sG(A_0)$ on the right,
\begin{restoretext}
\begingroup\sbox0{\includegraphics{ANCimg3/empty.png}}\includegraphics[clip,trim={.0\ht0} {.27\ht0} {.0\ht0} {.20\ht0} ,width=\textwidth]{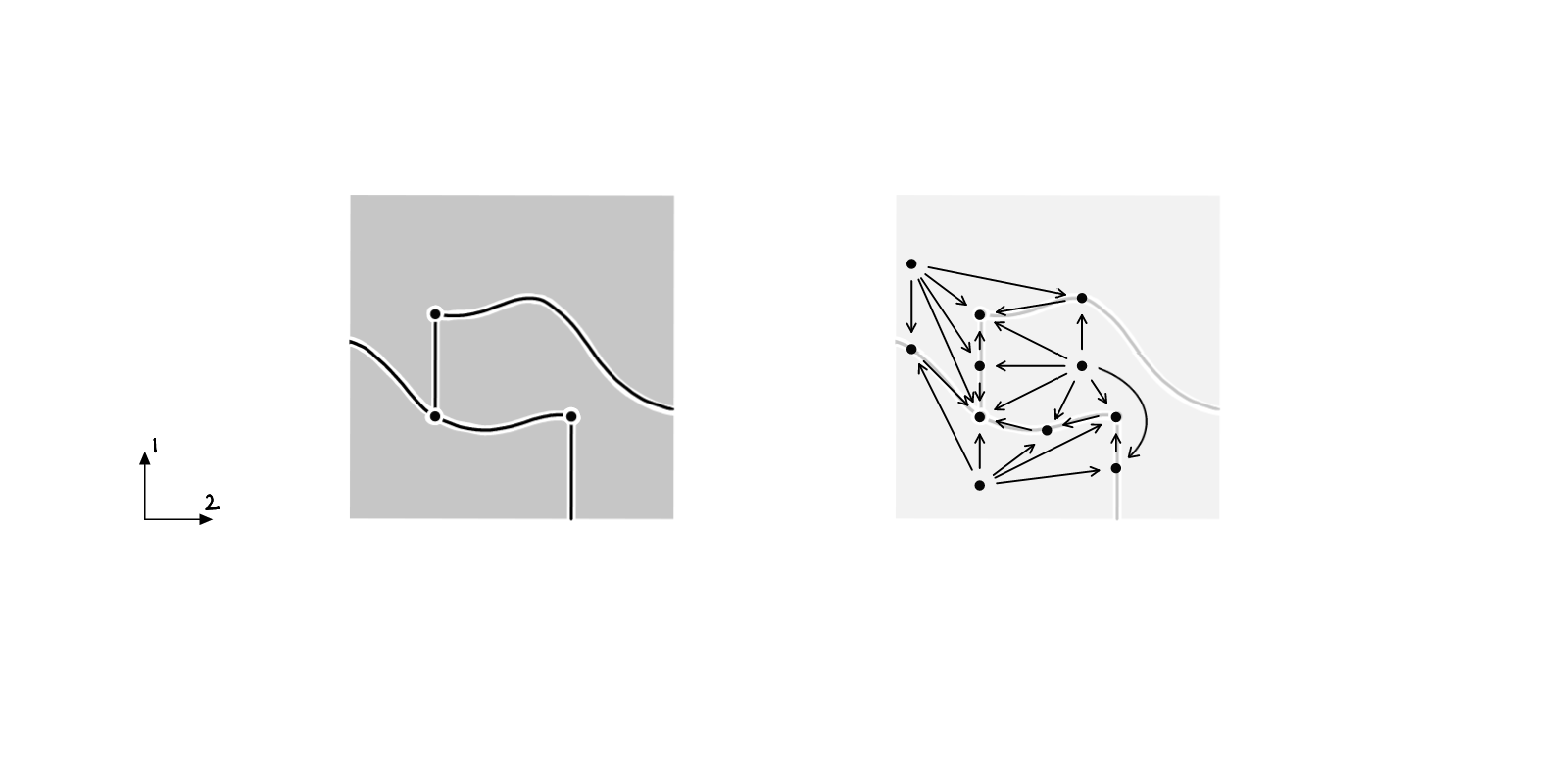}
\endgroup\end{restoretext}

\begin{defn}[Untyped $n$-manifold diagrams] An \textit{(untyped) $n$-manifold diagram} is a flag-foliation-compatible stratification of the $n$-cube whose manifolds have type $u$ only for $u$ which are monotonically decreasing sequences\footnote{In other words, the $k$-strata of a manifold diagram are required to intersect $(n-k)$-sheets transversally.}. 
\end{defn}
The above $A_0$ is \textit{not} an example of an untyped manifold diagram, but the following flag-foliation-compatible stratification is 
\begin{restoretext}
\begingroup\sbox0{\includegraphics{ANCimg3/empty.png}}\includegraphics[clip,trim={.0\ht0} {.3\ht0} {.0\ht0} {.18\ht0} ,width=\textwidth]{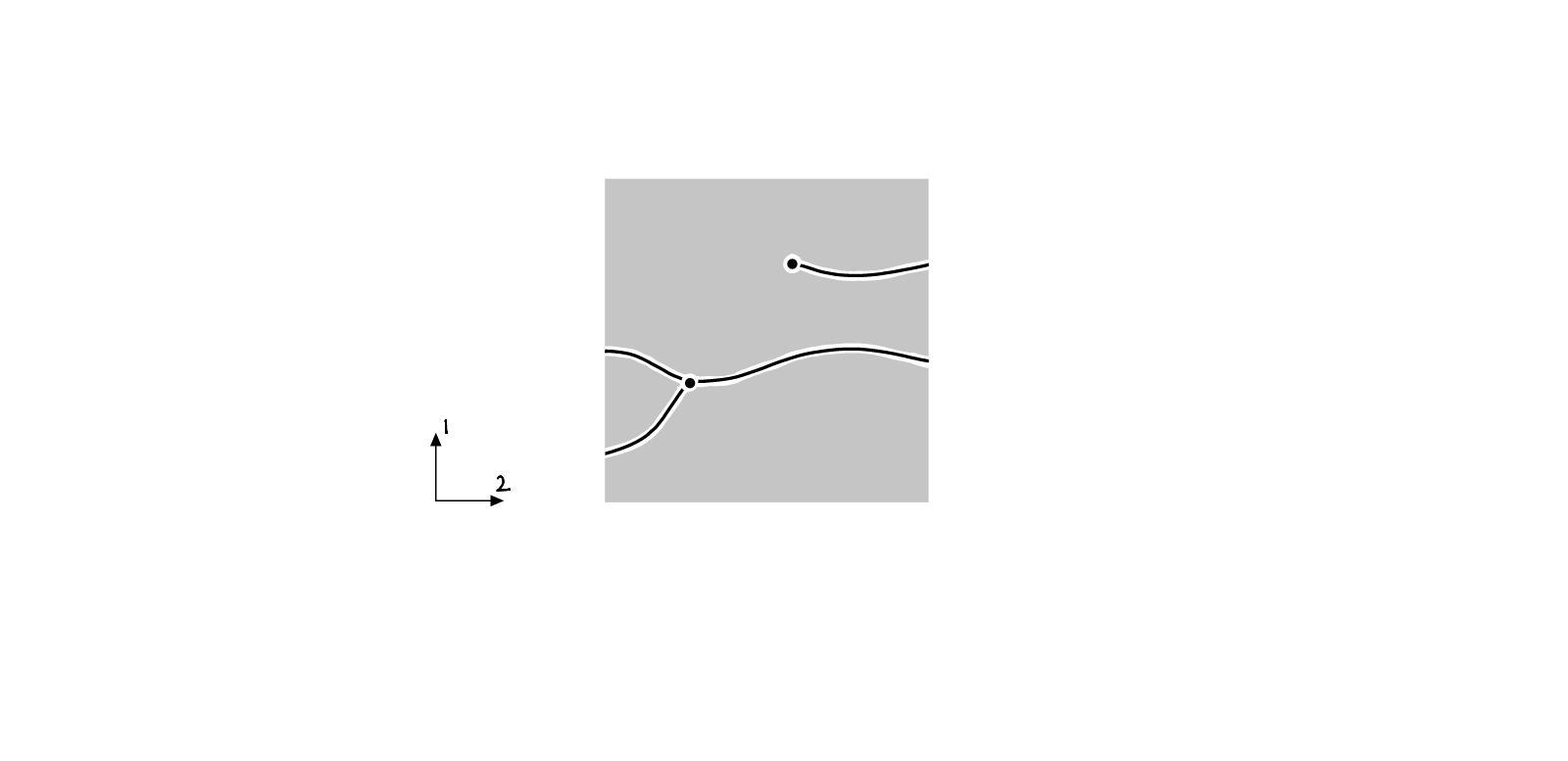}
\endgroup\end{restoretext}
\begin{defn}[Flag-foliation-preserving maps and equivalences]
A continuous map $f: (0,1)^n \to (0,1)^n$ is called \textit{flag-foliation-preserving} if it maps $k$-sheets into $k$-sheets and is monotone in each coordinate, meaning $f(x)_i$ depends monotonously on $x_i$.\footnote{In yet other words, if $f$ maps a $k$-sheet $\Set{p} \times (0,1)^k$ into the $k$-sheet $\Set{q}\times (0,1)^k$, then $(k-1)$-sheets $(p,a) \times (0,1)^{k-1}$ are mapped into $(k-1)$-sheets $(q,g(a))\times (0,1)^{k-1}$ such that $a \mapsto g(a)$ is monotonous.} 

Two flag-foliation compatible stratifications $A$, $B$ are \textit{equivalent}, written $A \simeq B$ if there is a flag-foliation-preserving homeomorphism mapping strata of one to strata of the other (thus in particular exhibiting a bijection of strata). 
\end{defn}
As an example, the following two manifold diagrams are equivalent
\begin{restoretext}
\begingroup\sbox0{\includegraphics{ANCimg3/empty.png}}\includegraphics[clip,trim={.0\ht0} {.25\ht0} {.0\ht0} {.22\ht0} ,width=\textwidth]{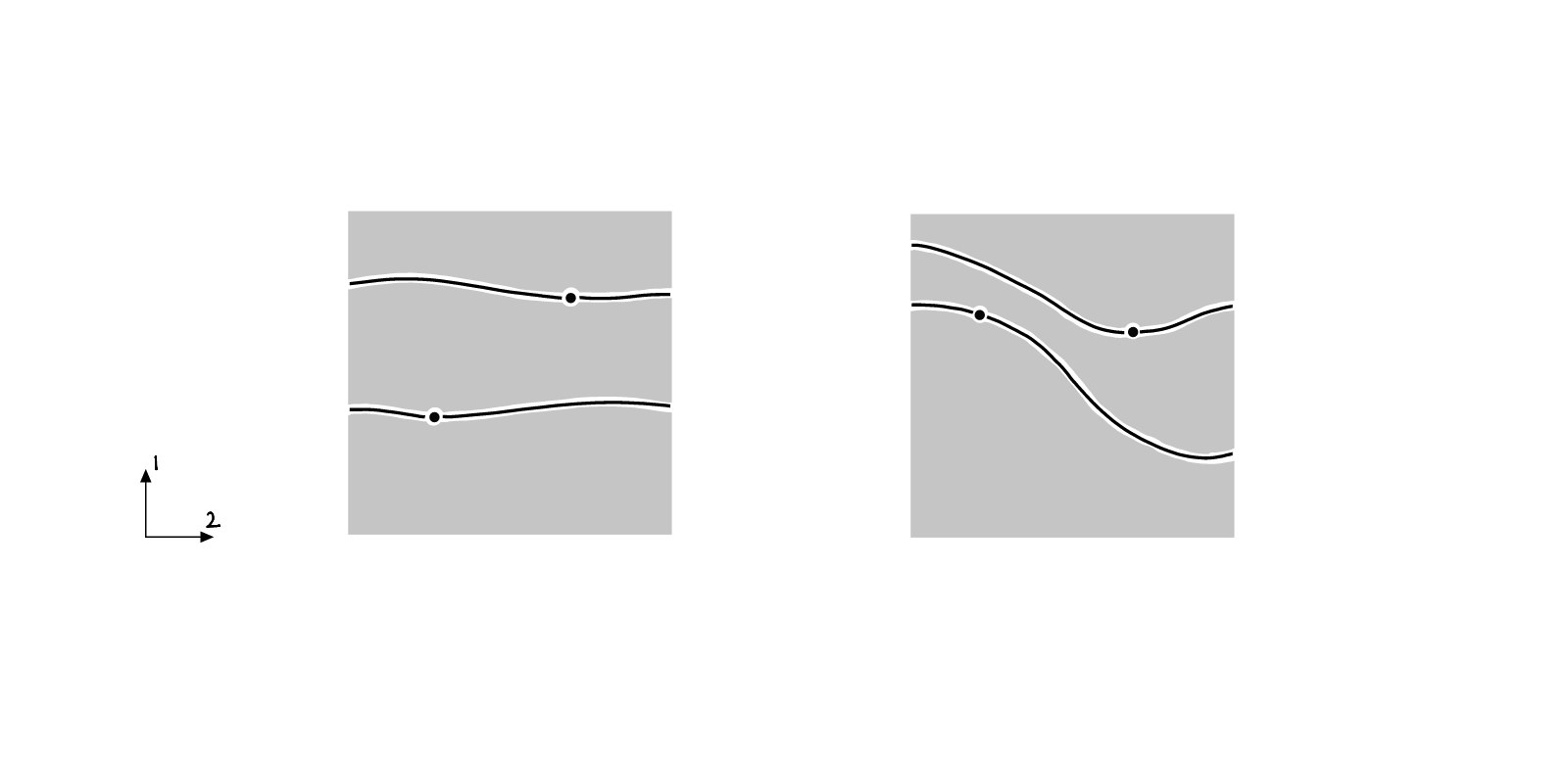}
\endgroup\end{restoretext}
but the next two diagrams are (crucially) not equivalent
\begin{restoretext}
\begingroup\sbox0{\includegraphics{ANCimg3/empty.png}}\includegraphics[clip,trim={.0\ht0} {.25\ht0} {.0\ht0} {.22\ht0} ,width=\textwidth]{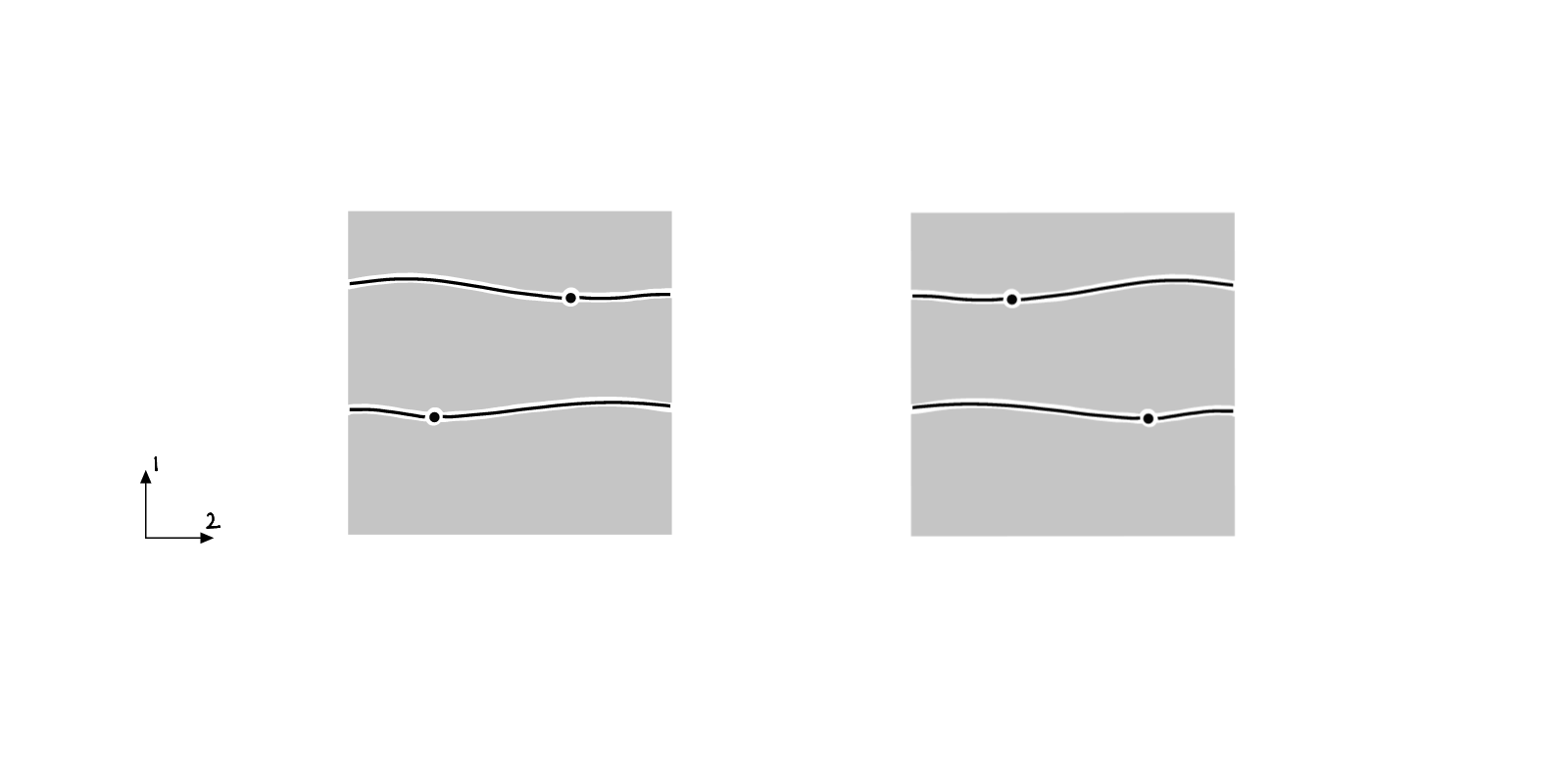}
\endgroup\end{restoretext}

\subsection{Topological labelled $n$-cubes}

Now, let $\cC$ be a category.
\begin{defn}[$\cC$-labelled cubes] A \textit{(topological) $\cC$-labelled $n$-cube} $A$ is (abusing notation) a flag-foliation-compatible stratification $A$ of the $n$-cube together with a functor $\sU_A : \sG (A) \to \cC$. We define 
\begin{equation}
c_A = \bigcup_{i \in (\sU _A)\inv(c)} p^A_i
\end{equation}
called the \textit{region of label} (or \textit{color}) $c$ in $A$.
\end{defn}
As an example of a $\bnum 2$-labelled 2-cube, which we will refer to as $A_1$, consider
\begin{restoretext}
\begingroup\sbox0{\includegraphics{ANCimg3/empty.png}}\includegraphics[clip,trim={.0\ht0} {.25\ht0} {.0\ht0} {.22\ht0} ,width=\textwidth]{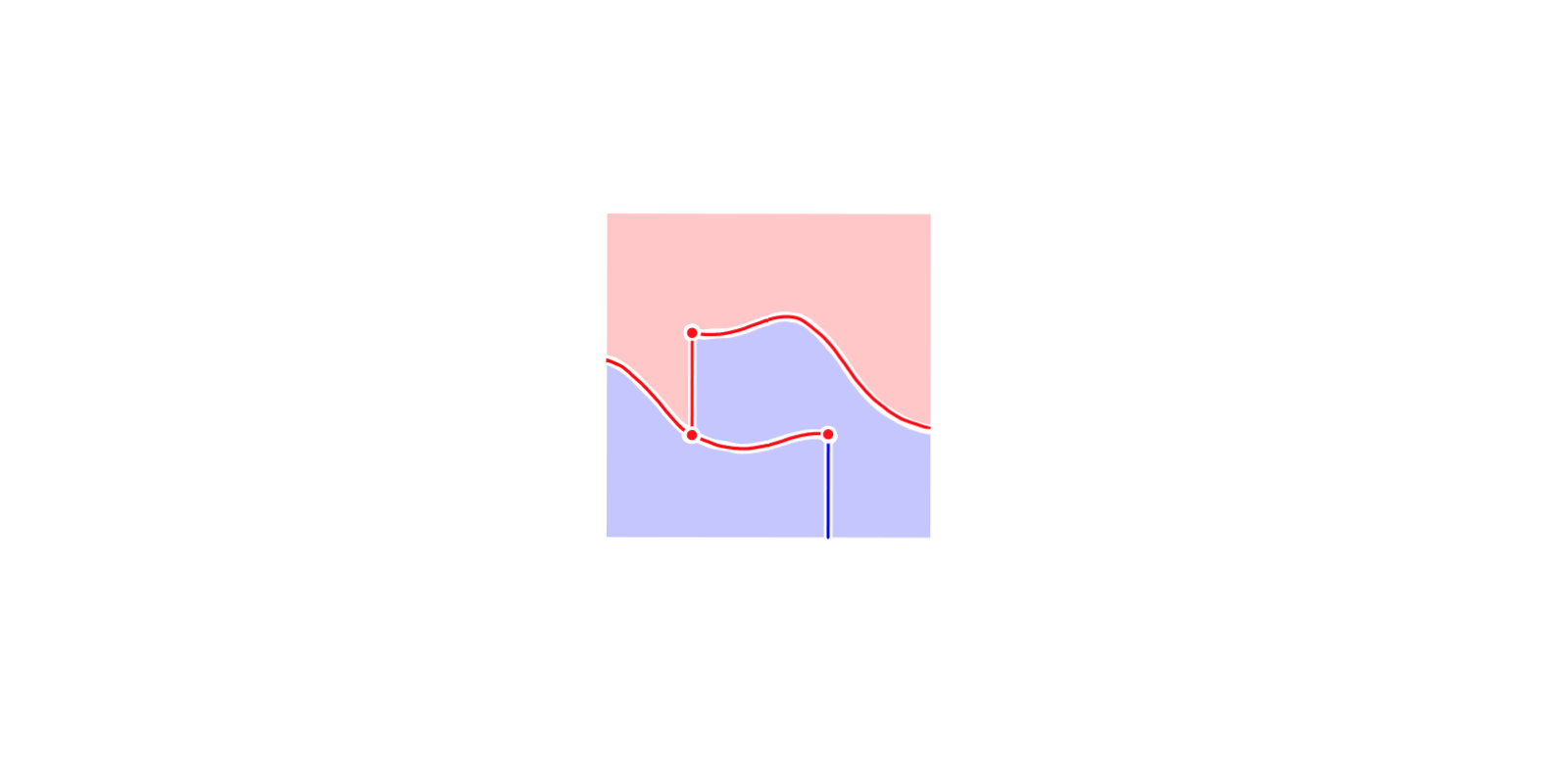}
\endgroup\end{restoretext}
The underlying flag-foliation-compatible stratification of $A_1$ equals our earlier example $A_0$. The regions of color $0$ and $1$ are indicated by \cblue{} and  \cred{} respectively. 

Later on we will be especially interested in the situation when a $\cC$-labelled cube $A$ is such that the collection $\Set{c_A}_{c \in \cC}$ forms an $n$-manifold diagram, in which case we speak of a \textit{$\cC$-labelled $n$-manifold diagram}. 

An example of a $\bnum 3$-labelled $2$-manifold diagram, which we will refer to as  $A_2$, is
\begin{restoretext}
\begingroup\sbox0{\includegraphics{ANCimg3/empty.png}}\includegraphics[clip,trim={.0\ht0} {.25\ht0} {.0\ht0} {.22\ht0} ,width=\textwidth]{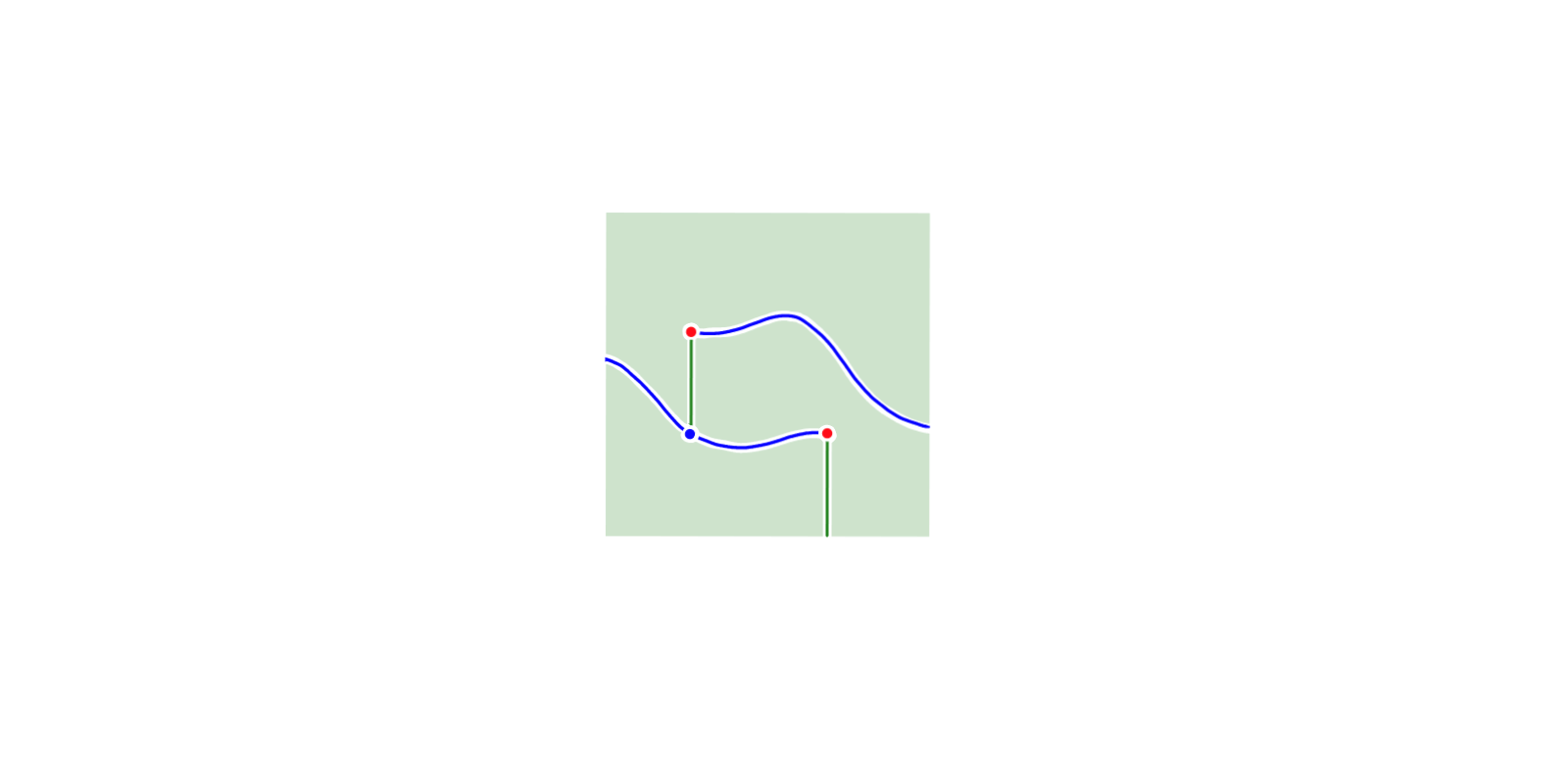}
\endgroup\end{restoretext}
The underlying flag-foliation-compatible stratification of $A_2$ is again $A_0$. The regions labelled by $0$, $1$ and $2$ are respectively indicated by colors \cdarkgreen{}, \cblue{} and \cred{} and arise from a functor $\sG(A_0) \to \bnum 3$ as required. 

Once we discuss the notion of refinements, it will become clear why one would want to keep the ``underlying stratification" and ``labelled regions" separate from one another.

\begin{rmk}[Singularity terminology] Note that a $\cC$-labelled manifold diagram $A$ has \textit{two} stratification: its underlying flag-foliation-compatible stratification $A$ with strata $\Set{p^A_i}_{i \in I}$, and its stratification as a manifold diagram with strata being the labelled regions $\Set{c_A}_{c \in \cC}$. To distinguish these to collections of strata, we will often refer to (the connected components of) $c_A$ as \textit{singularities} in $A$ instead of strata---in some situations this terminology nicely highlights the connection to classical singularity theory \cite{arnold1981singularity} that we will encounter.
\end{rmk}

There is also a notion of $\cC$-labelled $n$-cube \textit{families} indexed by a (stratified) space. We will only give a simplified account of this idea here.
\begin{defn}[Families of labelled cubes] \label{defn:cube_fam} Let $X$ be a stratification of (abusing notation) a $k$-manifold $X$. A $\cC$-labelled $n$-cube bundle $A$ over $X$ is (abusing notation) a stratification $A$ of $X \times (0,1)^n$ together with a ``global $\cC$-labelling" functor $\sU^n_A$ from the total poset $\sG(A)$ of $A$ into $\cC$, such that the following holds
\begin{enumerate}
\item Restricted to any point in $X$, $A$ yields a flag-foliation-compatible stratification of the $n$-cube
\item Restricted to any open neighbourhood in a stratum of $X$, $A$ is bundle homeomorphic to the trivial bundle\footnote{That is, there is a fibrewise flag-foliation-preserving bundle homeomorphism to the trivial product bundle.}
\end{enumerate}
\end{defn}
\noindent Note that $\sU_A$ in particular gives every fibre the structure of a $\cC$-labelled $n$-cube. The definition is simplified since we omitted conditions on $A$ for transitions between strata in $X$.

The right side of the depiction below shows an example of a $\bnum 3$-labelled $1$-cube bundle, which we will refer to as $A_3$, over a stratification $X_3$ of the open interval $(0,1)$. 
\begin{restoretext}
\begingroup\sbox0{\includegraphics{ANCimg3/empty.png}}\includegraphics[clip,trim={.0\ht0} {.07\ht0} {.0\ht0} {.15\ht0} ,width=\textwidth]{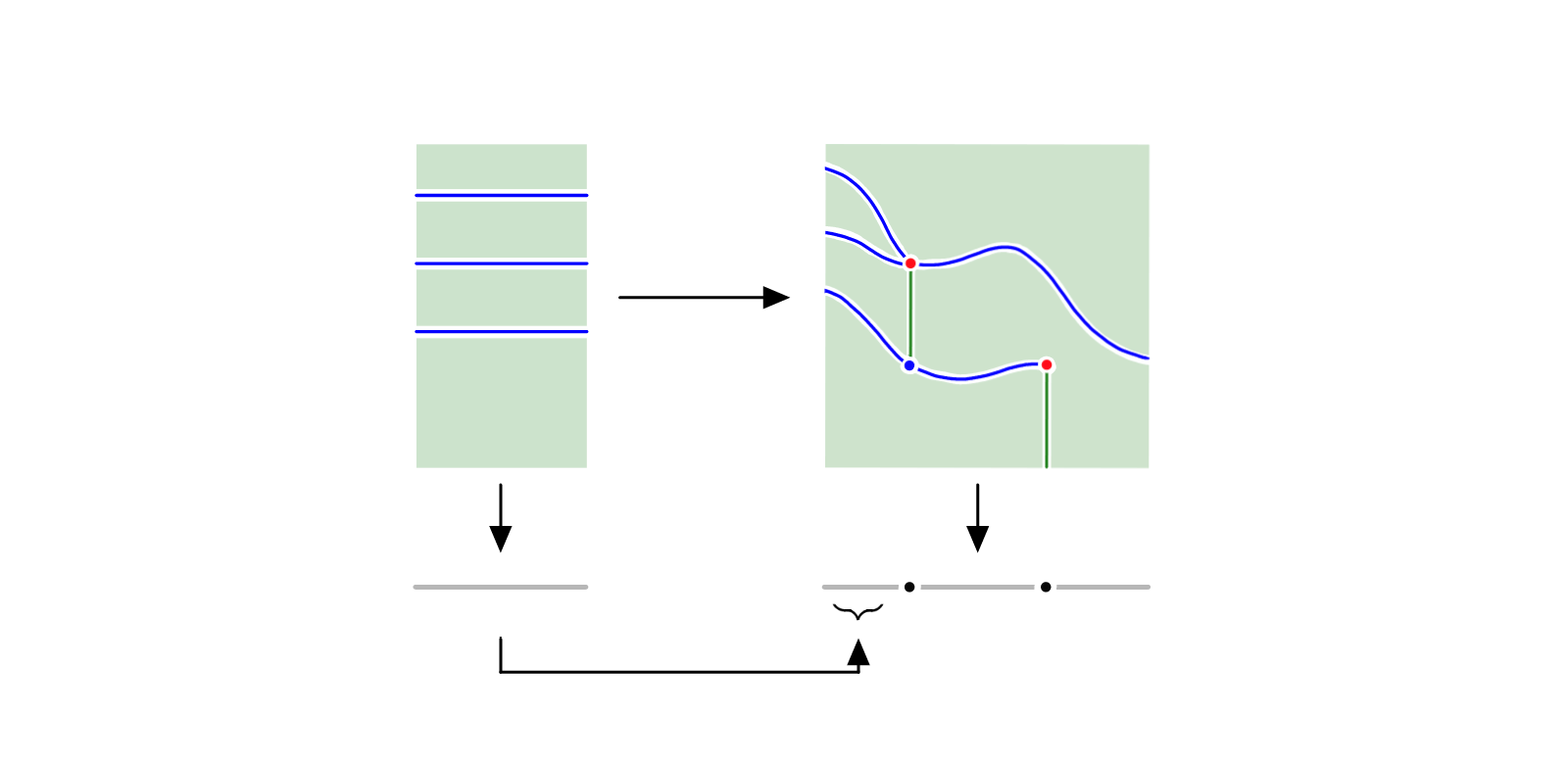}
\endgroup\end{restoretext}
The right vertical arrow is the projection map from $(0,1) \times (0,1)$ to the first component. $X_3$ stratifies $(0,1)$ into three $1$-dimensional strata and two $0$-dimensional strata as indicated. The stratification of $(0,1)\times (0,1)$ slighty modifies a previous example, and strata are indicated by connected components as before. The functor $\sU_{A_3} : \sG(A_3) \to \bnum 3$ takes \cdarkgreen{}, \cblue{} and \cred{} strata to $0$, $1$ and $2$ respectively. On the left we illustrated that restricting $A_3$ to an open neighbourhood in a stratum of $X_3$ yields (up to fibrewise flag-foliation-preserving bundle homeomorphism) the trivial bundle as required. Importantly, observe that $A_3$ is itself a flag-foliation-compatible stratification of the $(n+k)$-cube (here $n = k =1$).

Given labelled $n$-cubes we can always build trivial bundles over some manifold $X$, and the global labelling functor will be canonically inherited from the fibres. For instance in the case of the previous $\bnum 3$-labelled $2$-manifold diagram, we obtain a $\bnum 3$-labelled $2$-cube bundle as a trivial bundle over $X = (0,1)$ (or more precisely the stratification of $(0,1)$ with a single stratum)
\begin{restoretext}
\begingroup\sbox0{\includegraphics{ANCimg3/empty.png}}\includegraphics[clip,trim={.0\ht0} {.2\ht0} {.0\ht0} {.17\ht0} ,width=\textwidth]{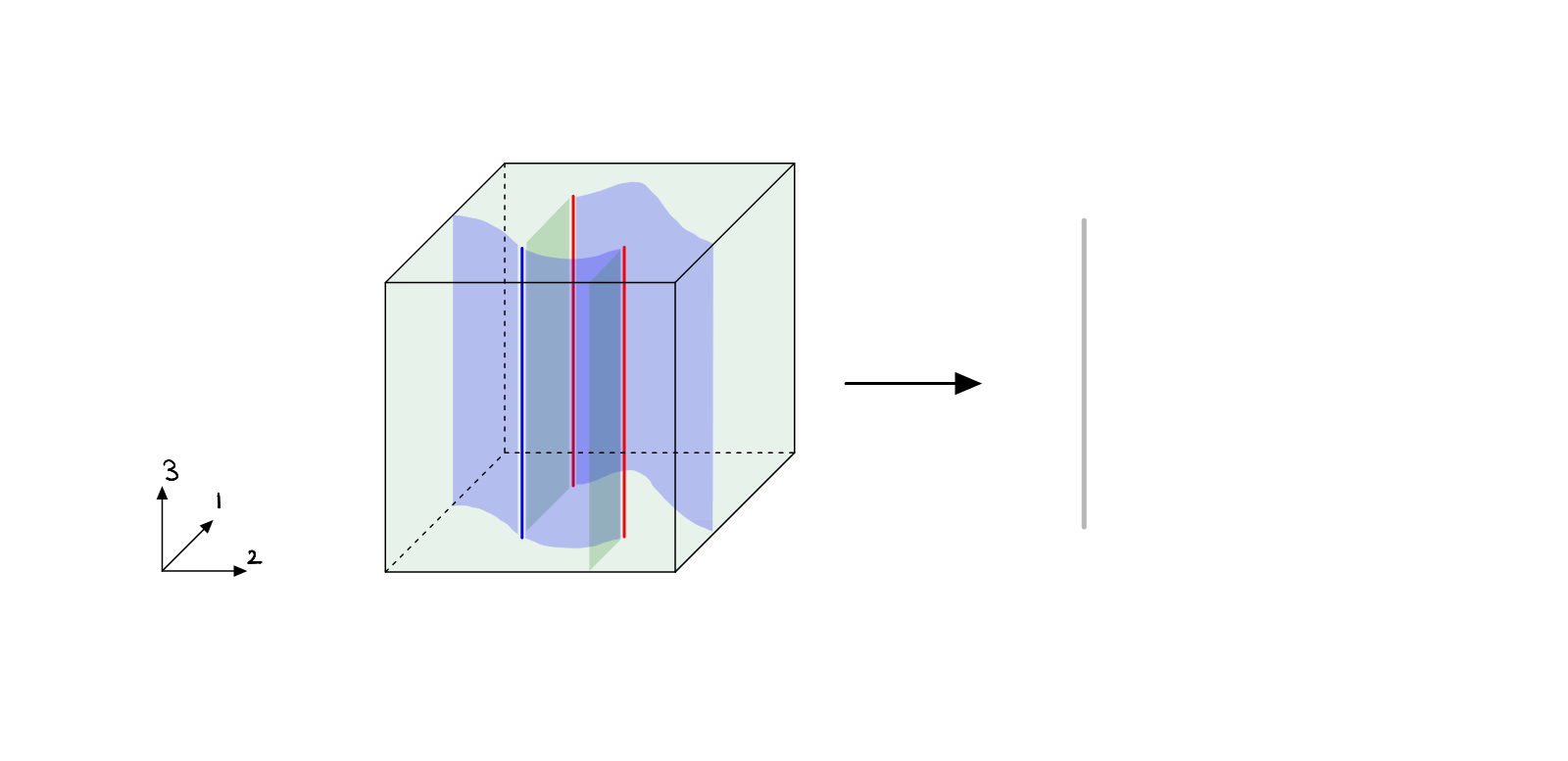}
\endgroup\end{restoretext}
This in particular is again a $\bnum 3$-labelled $3$-manifold diagram.

There is a natural categorical structure on $\cC$-labelled $n$-cubes (and thus on $n$-manifold diagrams) giving rise to a category $\Cubeo n \cC (\bnum 1)$, where $\bnum 1$ denotes the singleton space. The symbol $\circ$ indicates that we will restrict our discussion to open maps (but a definition without this restriction will only be given combinatorially). A similar definition could be given for families over $X$, giving rise to a category $\Cubeo n \cC (X)$ (however, this is subject to a more in-depth treatment of \autoref{defn:cube_fam}).

\begin{defn}[Category of (topological) $\cC$-labelled $n$-cubes and open cube maps] \label{rmk:cube_cat} An \textit{open $n$-cube map} $f : A \to B$ of $\cC$-labelled $n$-cubes $A, B$ is an open flag-foliation-preserving map such that $f$ maps strata of $A$ into strata of $B$, and the induced functor of total posets $f: \sG(A) \to \sG(B)$ factors $\sU_A$ through $\sU_B$. $\cC$-labelled cubes and open cube maps organise into a category denoted by $\Cubeo n \cC(\bnum 1)$.
\end{defn}

\begin{defn}[Refinements and equivalence] \label{defn:ref_and_eq} If $f : A \to B$ is a open $n$-cube map with its underlying map being a homeomorphism then we call $f$ a \textit{refinement} (or say ``it refines $B$ by $A$"). Given $\cC$-labelled $n$-cubes $A$, $B$ we say $A$ and $B$ are \textit{equivalent}, written $A \simeq B$, if there is a cospan $A \to C \ot B$ of refinements, in other words, if they have a ``mutual coarsening".
\end{defn}
\noindent Note that this notion of equivalence in a way extends the previous notion of equivalence of flag-foliation-compatible stratifications, now also allowing ``coarsening" and ``refinements" of the stratification.

Having briefly discussed bundles, a possible additional condition in their definition is to require the bundle map to descend to a map of total posets. In the specific context of cubes we can further ask for the following. Let $\pi : (0,1)^{k+1} \iso (0,1) \times (0,1)^{k} \to (0,1)^k$ denote the projection from the $(k+1)$-cube to the $k$-cube by omitting the first coordinate (thus mapping $l$-sheets to $(l-1)$-sheets).

\begin{defn}[Projection-stable cubes] \label{defn:projection_stable}  We say $A \in \Cubeo n \cC(\bnum 1)$ is \textit{projection-stable}, if there are flag-foliation-compatible stratifications $A^k$ of $(0,1)^k$, $k < n$, with total posets $\sG^k(A)$ such that
\begin{equation}
(0,1)^n \xto {\pi} (0,1)^{n-1} \xto {\pi} ... \xto{\pi} (0,1)^k
\end{equation}
induces a functor of posets $\sG^n(A) \to \sG^k(A)$ and exhibits $A$ to be a $\cC$-labelled $(n-k)$-cube bundle over $A^k$. 
\end{defn} 
\noindent The conditions in the previous definition imply that $A^k$ is unique if it exists (namely it means strata of $A^k$ are exactly projections of strata in $A$).

A guiding motivation of the algebraic work in this thesis is the observation that every sufficiently ``nice and finite" $\cC$-labelled $n$-cube has a projection-stable refinement. A procedure to obtain such a refinement will be sketched in the next section. For instance, the projection-stable refinements of our $\bnum 3$-labelled $2$-cube $A_3$ is given by the cube, which we will refer to as $A_4$, depicted at the top of the following illustration
\begin{restoretext}
\begingroup\sbox0{\includegraphics{ANCimg3/empty.png}}\includegraphics[clip,trim={.0\ht0} {.05\ht0} {.0\ht0} {.07\ht0} ,width=\textwidth]{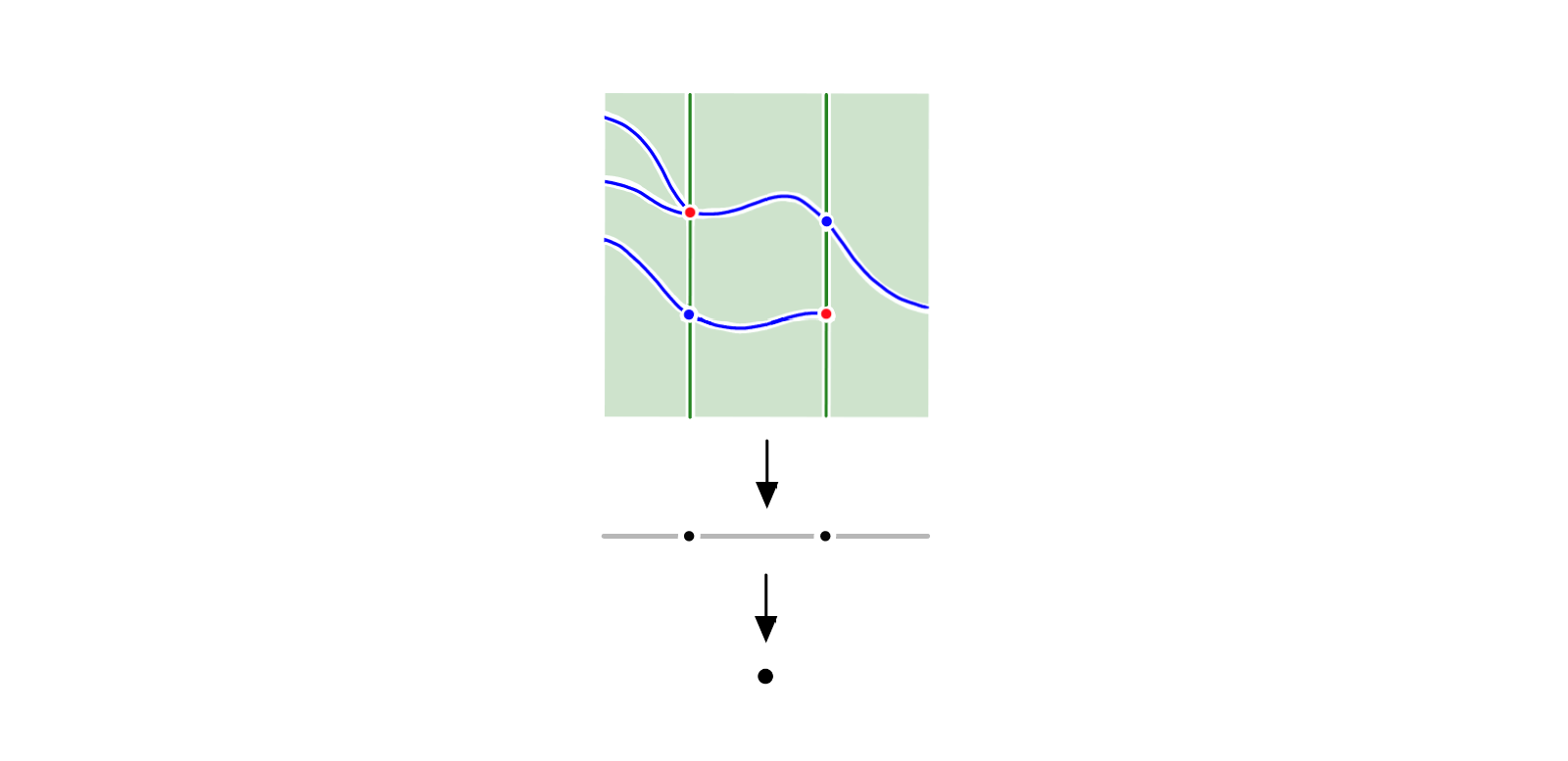}
\endgroup\end{restoretext}
Vertical arrows are projections $\pi : (0,1)^{k+1} \to (0,1)^{k}$ and as required these induce maps of total posets $\sG^2(A_4) \to \sG^1(A_4)$ and $\sG^1(A_4) \to \sG^0(A_4)$ respectively. One can further see that the condition for $A_4$ being a cube bundle over both $A^1_4$ and $A^0_4$ are satisfied.

The importance of projection-stable cubes derives from them having an elegant, fully faithful combinatorial description, forming the heart of this thesis. In particular, we will describe a combinatorial version of $\Cubeo n \cC(\bnum 1)$, denoted by $\Buno n \cC(\bnum 1)$. In \autoref{ssec:coloring} we will also provide a ``geometric realisation" function (at the level of \textit{object sets})
\begin{equation}
\norm{-} : \Buno n \cC(\bnum 1) \to \Cubeo n \cC(\bnum 1)
\end{equation}
Elements in the image will be exactly projection-stable ``finite" cubes (up to isomorphism), which then, conjecturally and up to equivalence, capture all ``finite" cubes. While it won't be relevant to us, the following remark suggests one way to intrinsically capture such a finiteness condition.
\begin{rmk}[Finiteness] \label{rmk:finiteness} A morphism in $\Cubeo n \cC(\bnum 1)$ is called an {embedding} if its induced functor of total posets is injective. Two morphisms $f,g : A \to B$ are called homotopic if they are related by a homotopy through open cube maps. We say $B$ is {finite} if $\sG^n(B)$ is finite and $B$ has finitely many embeddings up to homotopy.
\end{rmk}
\noindent From now on we will restrict our discussion to the class of sufficiently finite cubes (to which all examples drawn in this thesis belong), in other words, those living in the image of $\norm{-}$ up to equivalence.

We end this section with two further remarks: the first concerning ``types" of morphisms in $n$-categories represented by manifold diagrams, the second concerning the distinguishing feature of morphisms in $n$-fold and $n$-categories.

\begin{rmk}[Typability] \label{rmk:typability} So far we described \textit{untyped} (topological) $\cC$-labelled $n$-cubes $A$ (and resp. $n$-manifold diagrams): here, the word ``untyped" refers to the fact that regions $c_A$ need not be of the same type globally, that is, local neighbourhoods of points in the region $c_A$ can globally look very different. This is unwanted behaviour if we think of $\cC$ describing a set of higher morphisms with fixed types. The following suggests a notion of typability.

Fix $p \in (0,1)^n$ and let $A,B \in \Cubeo n \cC(\bnum 1)$. $A$ is called \textit{a minimal neighbourhood around $p$} if every embedding into $A$ containing $p$ in its image is an isomorphism up to homotopy. $B$ is called \textit{typable} if for all $c \in \cC$ we have a cube $\abss{c}$ (called the \textit{type} of $c$ in $B$) such that for all $q \in c_B$, $f : A \to B$, $f(p) = q$ and $A$ a minimal neighbourhood around $p$, we have $A \simeq \abss{c}$.

As an example, none of $A_1, A_2, A_3$ or $A_4$ were typable. But the following $\bnum 3$-labelled $2$-cube $A_5$ (depicted on the left) is
\begin{restoretext}
\begingroup\sbox0{\includegraphics{ANCimg3/empty.png}}\includegraphics[clip,trim={.0\ht0} {.25\ht0} {.0\ht0} {.25\ht0} ,width=\textwidth]{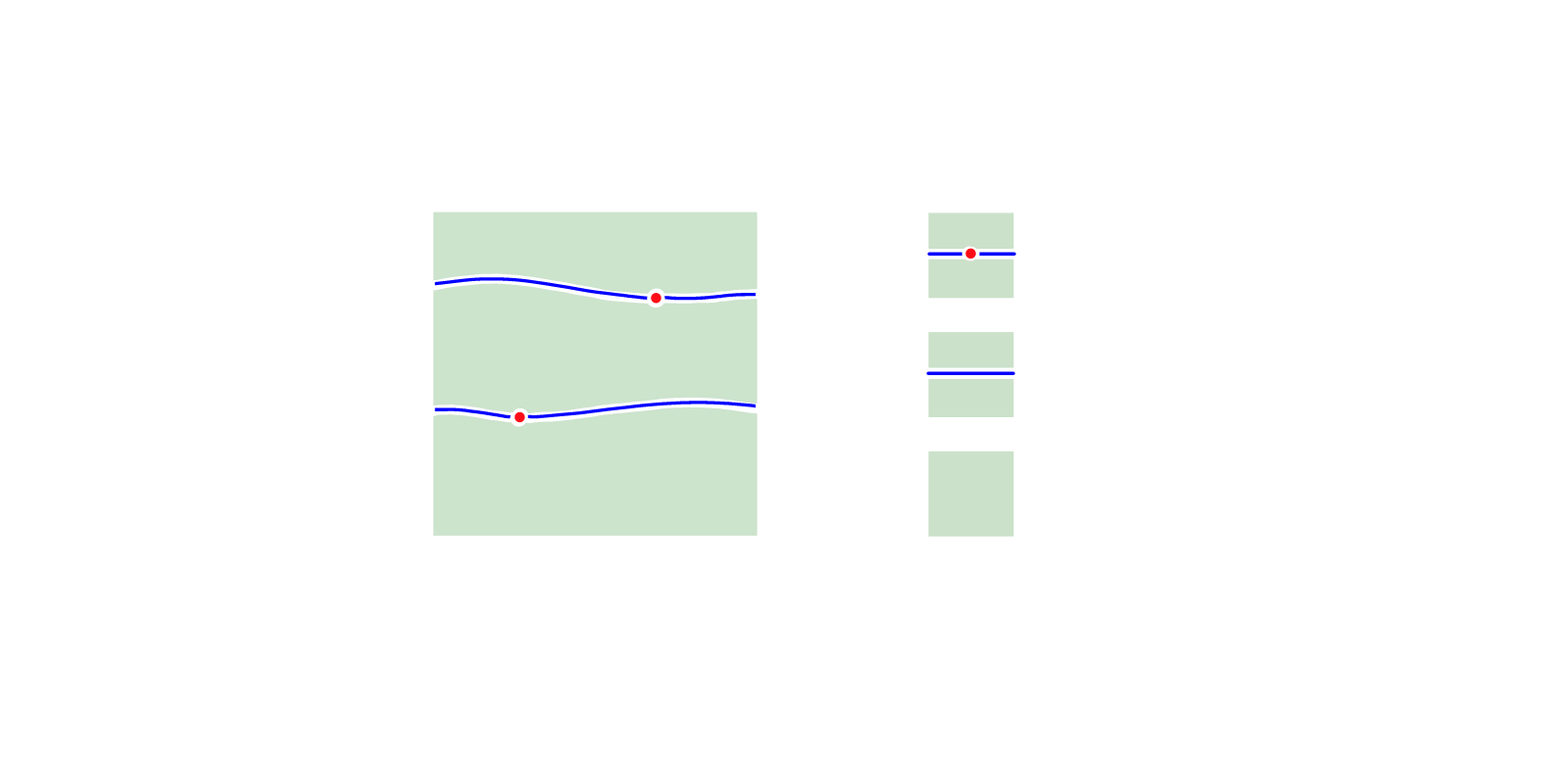}
\endgroup\end{restoretext}
\end{rmk}
\noindent The types (i.e. local neighbourhoods up to equivalence) of the \cdarkgreen{}, \cblue{} and \cred{} singularities can be seen to be the same everywhere and are shown on the right.

We further make a first remark on the importance of the generality of $\cC$-labelled $n$-cubes in contrast to just ``mere" manifold diagrams.

\begin{rmk}[$n$-fold categories] \label{rmk:po_nfold} The relationship of (typable) $n$-manifold diagrams and (typable) $\cC$-labelled $n$-cubes reflects the relationship of $n$-categories and $n$-fold categories. Concretely, there are $\bnum 2^n$ \stratatype{}s of morphisms in an $n$-fold category, but only $(n+1)$ of those yield ``monotonically increasing sequences" (as required in our definition of manifold diagrams above), thus there are only $(n+1)$ types of morphism in an $n$-category (corresponding to $0$-, $1$-, ... $(n-1)$ and $n$-morphisms). 
\end{rmk}

\subsection{Projection-stable refinements} \label{ssec:sum_rec_struct}

We will now sketch a procedure to obtain a projection-stable $n$-cube refinement for any given (sufficiently finite) $\cC$-labelled $n$-cube. This will motivate the algebraic definitions of the next sections.

We outline the procedure in the case of the following example of a $\bnum 3$-labelled\footnote{The color scheme for the objects in $\bnum 3$ in this and the following sections is different than that in previous sections.} singular $3$-cube $P$, which we already met in \autoref{rmk:GTP}, and which is depicted on the right below
\begin{restoretext}
\begingroup\sbox0{\includegraphics{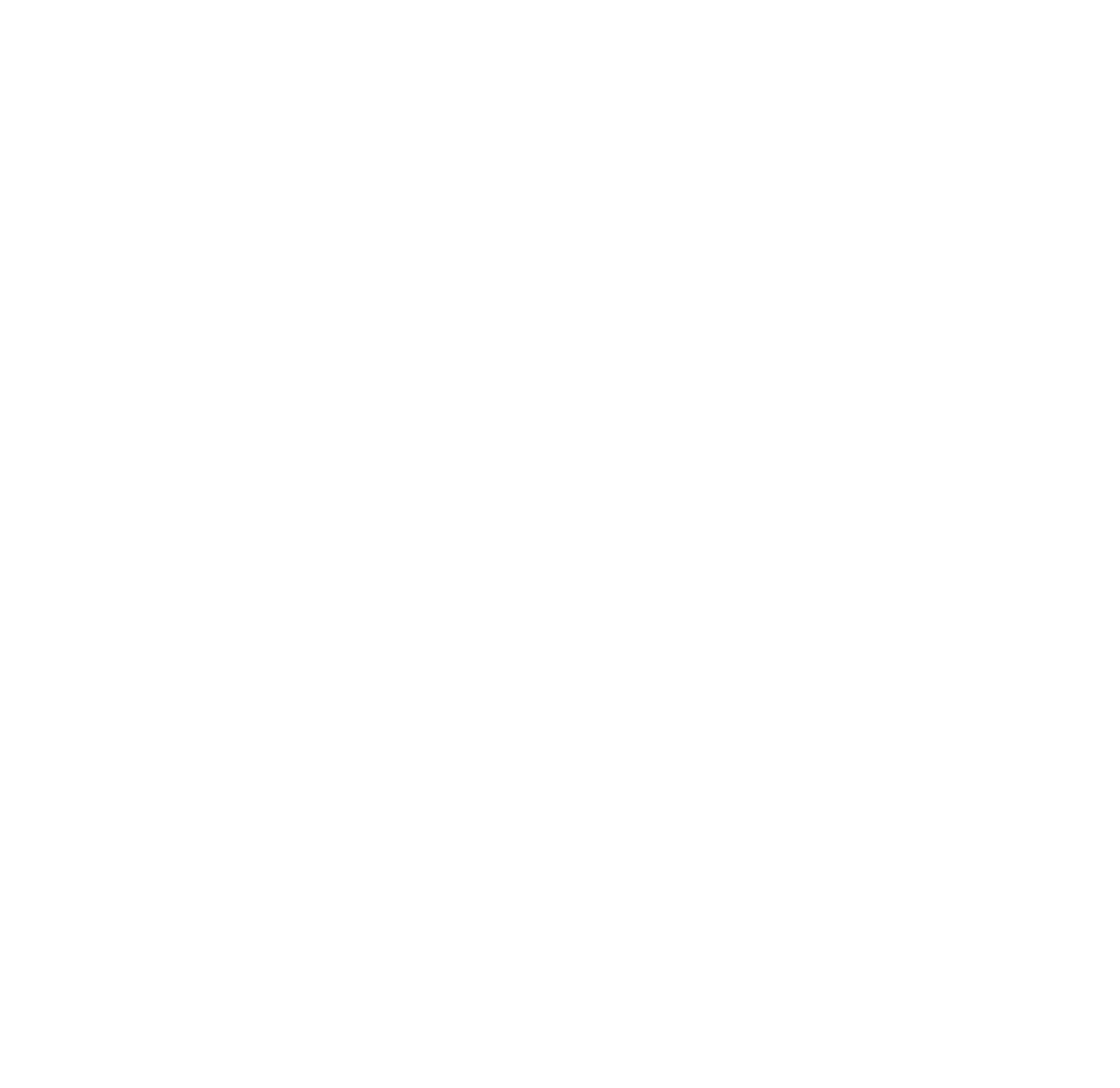}}\includegraphics[clip,trim=0 {.1\ht0} 0 {.1\ht0} ,width=\textwidth]{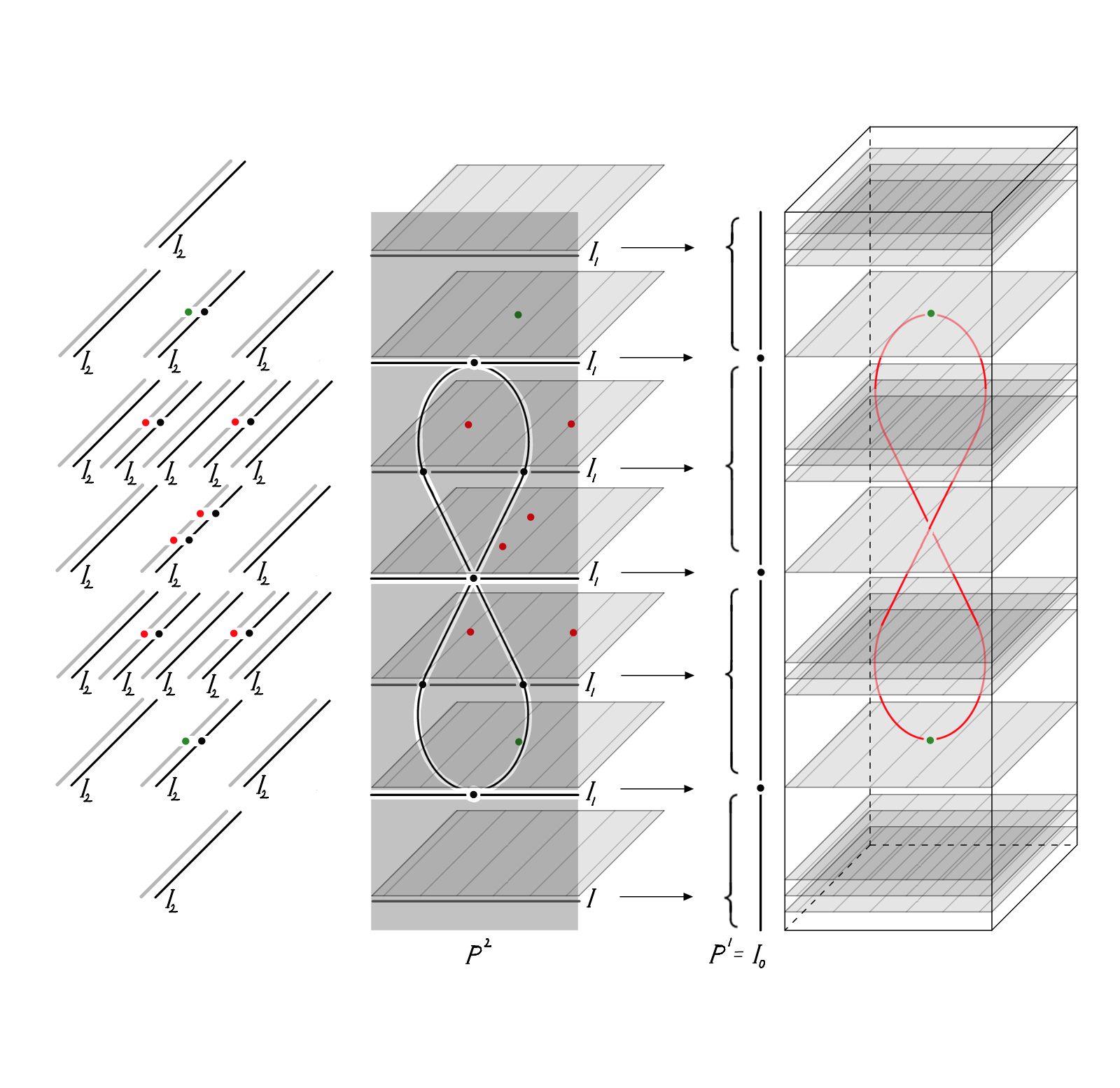}
\endgroup\end{restoretext}
Note that the cube consists of 5 strata in total: two $0$-dimensional strata (labelled by \cdarkgreen{}), two $1$-dimensional strata (labelled by \cred{}) and one $3$-dimensional stratum (labelled by \cgray{}, but not colored on the right). Note that $P$ is \textit{not} projection-stable. We want to show how to construct its projection-stable refinement. First, from the previous sections we recall that
\begin{enumerate}
\item by restricting to a $k$-sheet in a labelled $n$-cube, one can obtain a labelled $k$-cube. For instance, the figure above contains several examples of labelled $1$-cubes and $2$-cubes obtained by restriction to $1$-sheets and $2$-sheets in a labelled $3$-cube.
\item $(n-k)$-sheets are parametrised by $p \in (0,1)^{k}$. For instance, $2$-sheets of our $3$-cube above are parametrised by an interval $I_0$ (a stratified version of which is drawn on the left of the cube in the right column). 
\item $P$ being projection-stable is equivalent to the existence of flag-foliation-compatible stratifications  $P^k$ of $(0,1)^k$ for $k \leq n$ where strata of $P^k$ are projections of the strata of $P^{k+1}$ under $\pi : (0,1)\times (0,1)^k\to (0,1)^k$ and $P^n = P$, such that $P$ becomes a $(n-k)$-cube bundle over each $P^k$ under $\pi^{n-k} : (0,1)^n \to (0,1)^k$
\end{enumerate}
We sketch an inductive construction of $P^k$, $k \leq n$, such that the bundle condition is satisfied. We then obtain $P^n$ as the desired projection-stable refinement of $P$. We assume all cubes that we encounter are sufficiently ``finite" by requiring $P$ to satisfy an appropriate finiteness condition (such as \autoref{rmk:finiteness}).

$P^0$ is the trivial stratification of the point. Assume $P^k$ to be constructed with the bundle condition satisfied. The latter implies that restriction to $(n-k)$-sheets over each stratum in $q$ in $P^k$ are constant (up to equivalence). Let those restrictions equal cubes $P^k_{x}$ over a given point $x \in q$.  Since $P^k_x$ is finite, its $(n-k-1)$-sheets, parametrised by an open interval $I_k$, appear in an alternating pattern: namely, they are part of trivial bundles over open subintervals of $I_k$ with ``intermediate" $(n-k-1)$-sheets over isolated points in $I_k$ in between those subintervals. This makes $I_k$ into a (flag-foliation-preserving) stratification of $(0,1)$, with $0$-strata also called \textit{singular heights} and $1$-strata also called \textit{regular segments}. Note that since $(0,1)$ is linearly ordered we can speak of the $i$th singular height (resp. regular segment). Note also that the stratification $I_k$ tacitly depends on $x \in q$.

We observe that the stratification $I_k$ is not uniquely determined by the above conditions, because of the following fact: every trivial $(n-k-1)$-cube bundle can be sub-divided into two such trivial bundles with an intermediate $(n-k-1)$-cube. In other words, singular heights can be added arbitrarily. However, subject to a more in-depth treatment of our notion of bundles (cf. \autoref{defn:cube_fam}), we can guarantee that a choice and, in fact, a minimal choice exists with the following properties: we can choose $I_k$ for each $P^k_x$, $x \in q$, and each $q$, such that singular heights of $I_k$ depend continuously on $x \in q$ (in particular their number remains constant over all of $q$). Further, if $q'$ is a stratum in the boundary of stratum $q$ then, over any continuous ``entrance" path $r : (0,1] \to (0,1)^k$ with $r(a < 1) \in q$, $r(1) \in q'$, the $i$th singular height in $I_k$ over $r(a) \in q$ converges to the $f(i)$th singular height in $I_k$ over $r(1) \in q'$ for some fixed function $f$ which is \textit{independent} of the chosen $r$.

As a result, we can now obtain $P^{k+1}$ as follows: A stratum in $P^{k+1}$ is the union of all $i$th singular heights (resp. regular segments) lying over a stratum $q \in P^k$ for fixed $i$. We illustrate this with our concrete example above as follows.
\begin{itemize}
\item In the case of $P$, we first choose a stratification $P^1 = I_0$ of the interval parametrising $2$-sheets such as the one depicted in the right column above, satisfying that $2$-sheet bundles over regular segments (i.e. over $1$-strata in $P^1$) are trivial bundles.
\item $P^2$ is then obtained by inductively choosing stratified intervals $I_1$, indicated in the middle column above, lying over points in the the $0$- and $1$-strata in $P^1$ and satisfying the following: we choose $I_1$ over $1$-strata in $P^1$ such that singular heights vary continuously over points in those strata, and converge ``functionally" over entrance paths, that is, when passing from $1$- to $0$- strata in $P^1$, as required in our general construction above. By taking unions of $i$th singular heights (resp. $i$th regular segments) over each stratum $q \in P^0$ separately, we then build the strata of $P^2$, a stratification of $(0,1)^2$ as indicated in the middle column of the above picture.
\item  Similarly, for building $P^3$, over each stratum $q$ of $P^2$ and each point $x \in q$ in it, we choose interval stratifications $I_2$ (as indicated in the left column above) and require the same continuous and functional properties as before. Letting $x \in q$ vary and taking unions of the $i$th singular heights and regular segments, we then obtain the strata of $P^3$ (which are not drawn above, but a related drawing can be found in \autoref{ssec:coloring}).
\end{itemize}
$P^3$ is a stratification of $(0,1)^3$ refining our original stratification $P$, thus in particular we have a map $\sG(P^3) \to \sG(P)$ of total posets we which lets $P^3$ inherit a labelling from $P$. As a result, we constructed a $\bnum 3$-labelled $3$-cube $P^3$ refining $P$ such that $P^3$ is projection-stable.

Starting in the next section, we will study the \textit{total posets} $\sG^k(P) := \sG(P^k)$ (cf. \autoref{defn:ffcs}), which come with natural ``poset projections" $\sG^{k+1}(P) \to \sG^k(P)$. As it will turn out, this combinatorial structure has elegant features. For a start, these towers of poset projections can be described by an iterated Grothendieck bundle construction, whose classifying category can be derived from an endofunctor on $\Cat$, as we will see soon. There is also a combinatorial analogue of cube maps which leads to well-behaved categories of these combinatorial $n$-cubes. This opens up an interesting and novel combinatorial playground forming the basis of this thesis.

\begin{rmk}[Geometric realisation] Crucially, with view towards the next sections, (we claim that) all the information needed to build $P$ (or rather, its projection-stable refinement $P^3$) is, up to equivalence, already is encoded in these total posets\footnote{Technically, we also need a second ``direction" order recording the direction of the parametrising intervals} $\tsG k(P)$ and the projection functors between them. Thus, all of $P$ can be encoded combinatorially as a tower of posets projections together with a functor into $\cC$. A ``geometric realisation" procedure for recovering $P$ from its combinatorial data will be given in \autoref{ssec:coloring}.
\end{rmk}

\section{Algebraic model: singular cubes} 

Labelled singular $n$-cubes are our main object of study. They will specialise to the combinatorial structures describing manifold diagrams. They can be thought of as total posets of projection-stable labelled cubes (as discussed in the previous section).

Labelled singular $n$-cubes will be built inductively, in each step adding one dimension of the cube by means of a bundle whose fibers are ``1-dimensional" combinatorial objects called singular intervals. The bundle construction forms the core of the inductive step and will be discussed first. 

\subsection{The classifying map paradigm}

We preface the discussion of the bundle construction with a short reminder about the classifying map paradigm, which plays a role in many parts of mathematics. The paradigm is similar to constructing the graph of a function: given a \textit{family} of ``pieces" $F$ indexed by a \textit{base space} 
\begin{restoretext}
\begingroup\sbox0{\includegraphics{test/page1.png}}\includegraphics[clip,trim=0 {.1\ht0} 0 {.1\ht0} ,width=\textwidth]{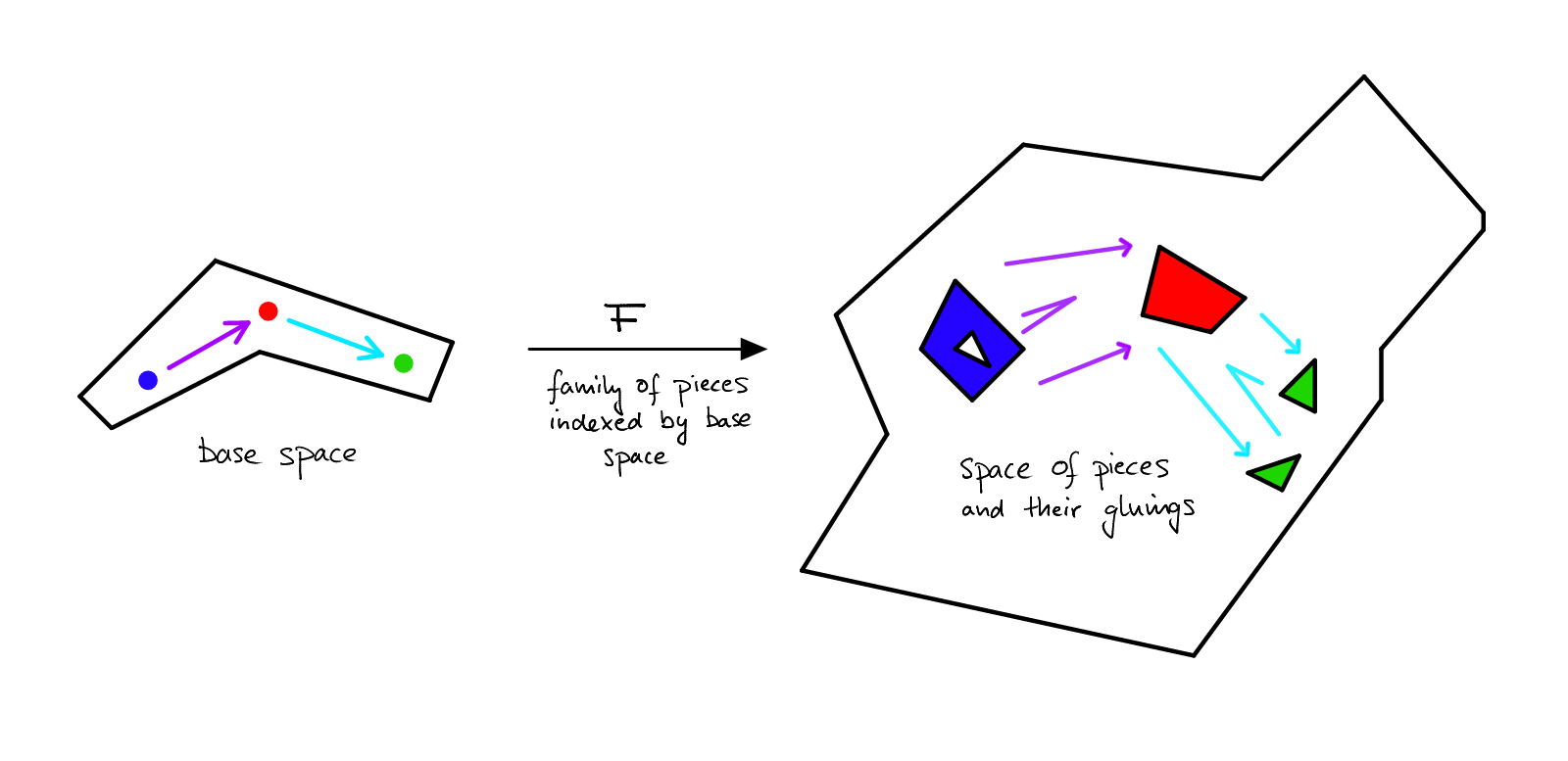}
\endgroup\end{restoretext}
we can bundle the pieces together, obtaining a \textit{bundle} $\pi_F$ over the base space
\begin{restoretext}
\begingroup\sbox0{\includegraphics{test/page1.png}}\includegraphics[clip,trim=0 {.1\ht0} 0 {.0\ht0} ,width=\textwidth]{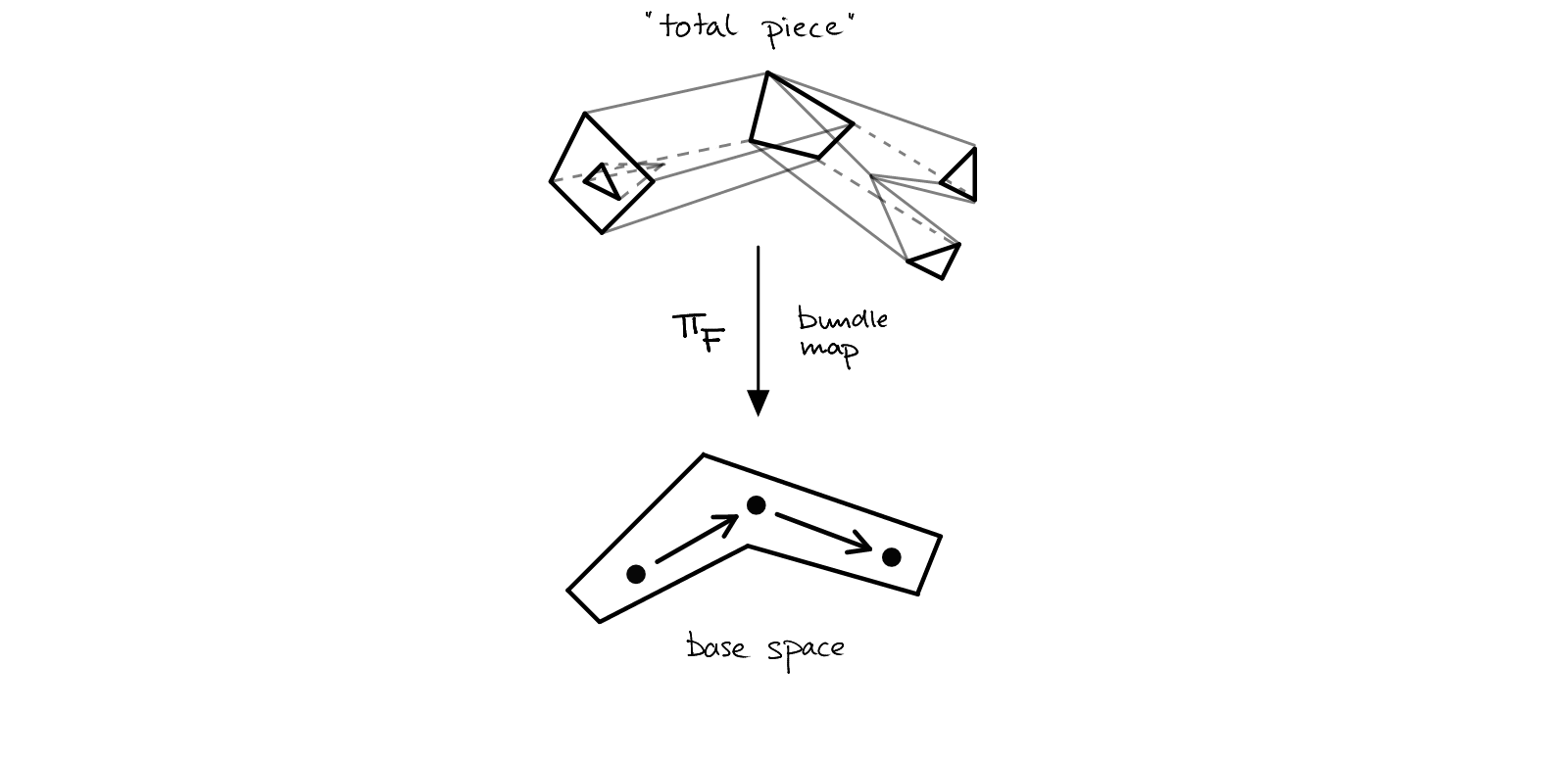}
\endgroup\end{restoretext}
Conversely, a total piece bundled over a base space can be broken apart into a family of its individual pieces indexed by the base space, which is often called the \textit{classifying map of the bundle}. However, we will usually only speak of families and their associated bundles. 

\subsection{Families of profunctorial relations and their bundles} \label{ssec:sum_fam_bun_prel}

We will now discuss an algebraic instance of a family-to-bundle construction. The goal of the section will be to define the following mapping
\begin{align} \label{eq:PRel_grothendieck}
\Fun(\cC, \PRel) &\to \Cat \sslash \cC \\
F &\mapsto (\pi_F : \sG(F) \to X)
\end{align}
Here, $\Cat$ is the category of categories and functors, $\cC$ is a category and $\left(\Cat \sslash \cC\right)$ denotes the overcategory over $\cC$. Denote by $\Bool$ the usual monoidal category of truth values $\bot \to \top$. $\Fun(\cC, \PRel)$ is the set of functors from $\cC$ into $\PRel$, the full subcategory of the category of $\Bool$-enriched profunctors whose objects are posets. Explicitly, objects in $\PRel$ are posets $Z, Y$ and morphism $R : Z \xslashedrightarrow{} Y$ are functors
\begin{align}
R : Z\op \times Y \to \Bool
\end{align}
Such $R : Z \xslashedrightarrow{} Y$ will be called a profunctorial relation to emphasise that morphisms in $\PRel$ are given by certain relations (namely, $R\inv(\top) \subset Z \times Y$ in the case of $R$) and compose like relations.

\begin{eg}[Profunctorial relations] Depicting elements $(x,y) \in R\inv(\top) \subset Z \times Y$ by edges from $x$ to $y$. The following is a profunctorial relation $R : \bnum 3 \to \bnum 2$
\begin{restoretext}
\begingroup\sbox0{\includegraphics{test/page1.png}}\includegraphics[clip,trim=0 {.3\ht0} 0 {.25\ht0} ,width=\textwidth]{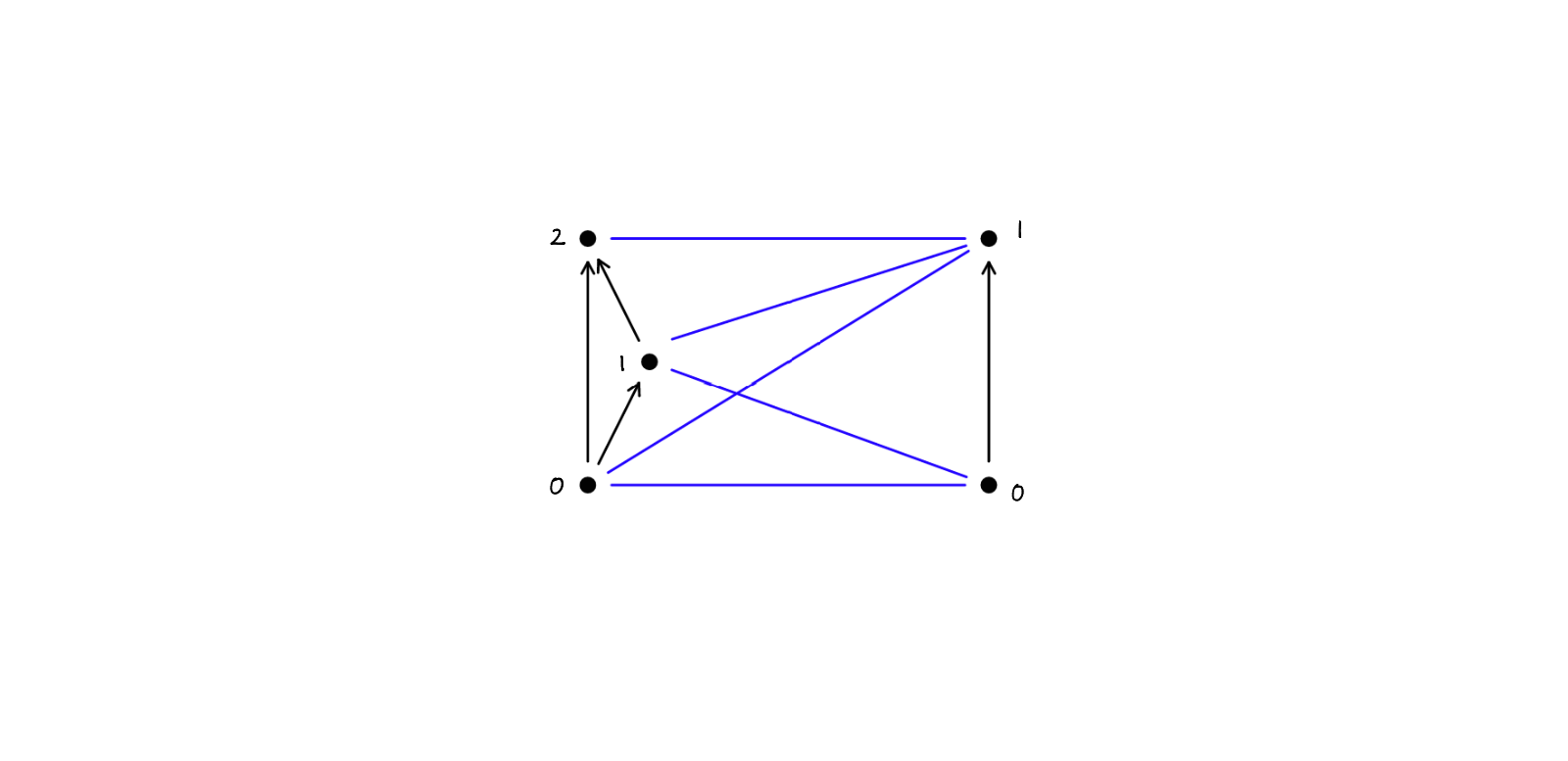}
\endgroup\end{restoretext}
The following is a relation $R : \bnum 3 \to \bnum 2$, but \textit{not} a profunctorial one
\begin{restoretext}
\begingroup\sbox0{\includegraphics{test/page1.png}}\includegraphics[clip,trim=0 {.3\ht0} 0 {.25\ht0} ,width=\textwidth]{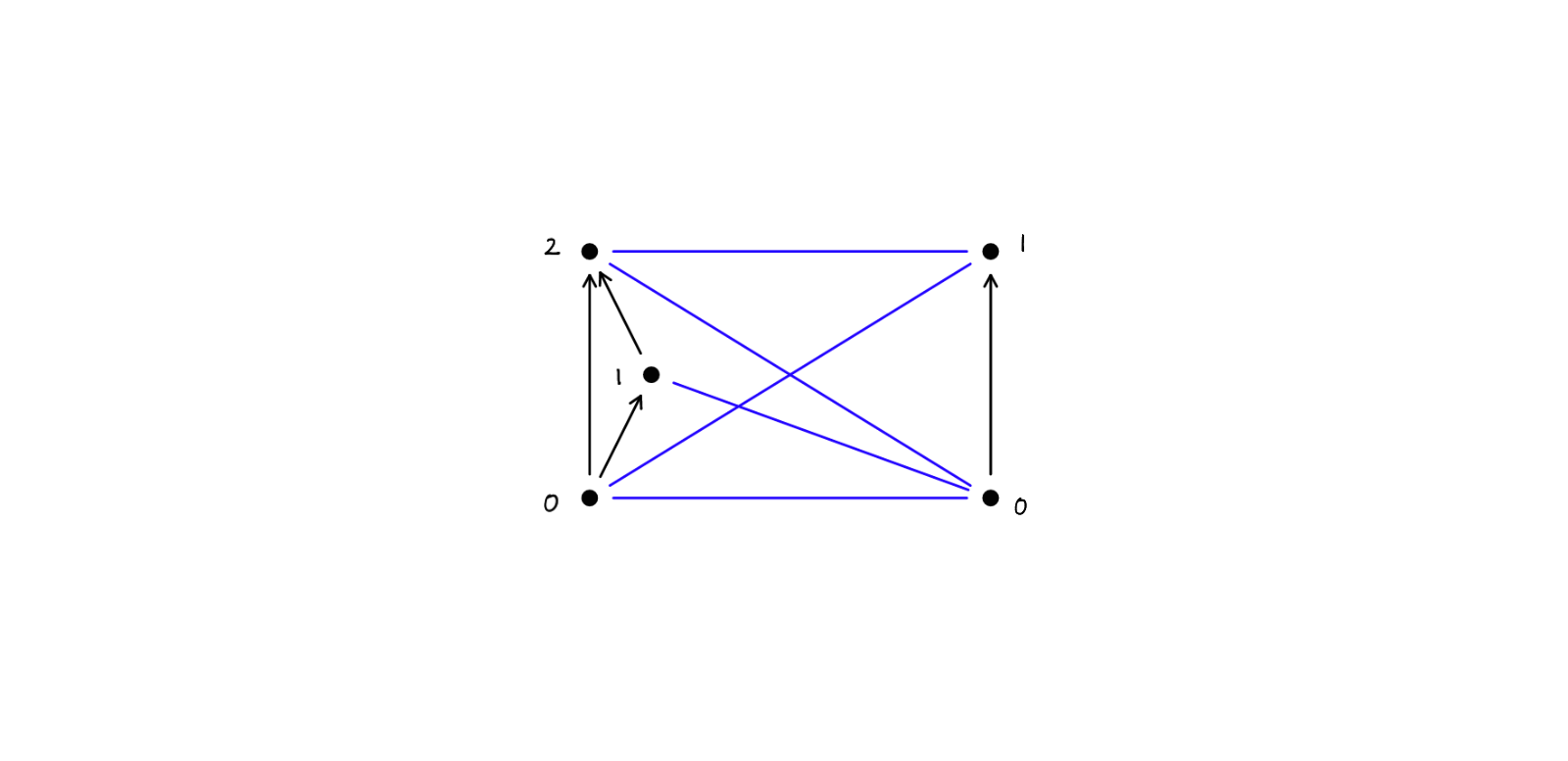}
\endgroup\end{restoretext}
For instance, functoriality would require an arrow $R(1,0 \to 1) : R(1,0) \to R(1,1)$ in \Bool{}, but $R(1,0) = \top$ and $R(1,1) = \bot$ makes this impossible.
\end{eg}

We now define the mapping \eqref{eq:PRel_grothendieck}. Given a family of profunctorial relations $F$ indexed by a category $\cC$ (that is, a functor $F : \cC \to \PRel$) we first construct $\sG(F) \in \Cat$ as follows: objects in $\sG(F)$ are given by
\begin{equation}
\obj(\sG(F)) = \Set{(c,a)~|~c \in \cC, a \in F(x)}
\end{equation}
morphisms in $\sG(F)$ are given by
\begin{align}
&(f,a,b) : (c,a) \to (d,b) \in \mor(\sG(F)) \\
&\iff \quad (f : c \to d) \in \mor(\cC) \text{~and~} F(f)(a,b)
\end{align}
Note that $F(f)$ is a profunctorial relation and can thus be evaluated on two values as done on the right hand-side.

If $g\circ f = h$ in $\cC$, then morphisms in $\sG(F)$ compose as
\begin{equation}
(g,b,c) \circ (f,a,b) = (h,a,c)
\end{equation}
We then construct a functor $\pi_F : \sG(F) \to \cC$ by defining
\begin{equation}
\pi_F (c,a) := c
\end{equation}
If $\cC$ is a poset then $\sG(F)$ is a poset which will be called \textit{total poset} of $F$, otherwise we speak of $\sG(F)$ as \textit{total category}\footnote{The letter $\sG$ is chosen to allude to both ``graph" and ``Grothendieck construction".}. $\pi_F : \sG(F) \to \cC$ is called the ($\PRel$-)\textit{bundle} associated to the family $F$. 

\begin{rmk}[Inverse of Grothendieck construction] The explicit construction of the mapping $F \mapsto \pi_F$ above is injective, and thus $F$ is fully determined by $\pi_F$.
\end{rmk}

As an example consider the following family $F$ of profunctorial relations
\begin{restoretext}
\begingroup\sbox0{\includegraphics{test/page1.png}}\includegraphics[clip,trim={.1\ht0} {.09\ht0} {.1\ht0} {.1\ht0} ,width=\textwidth]{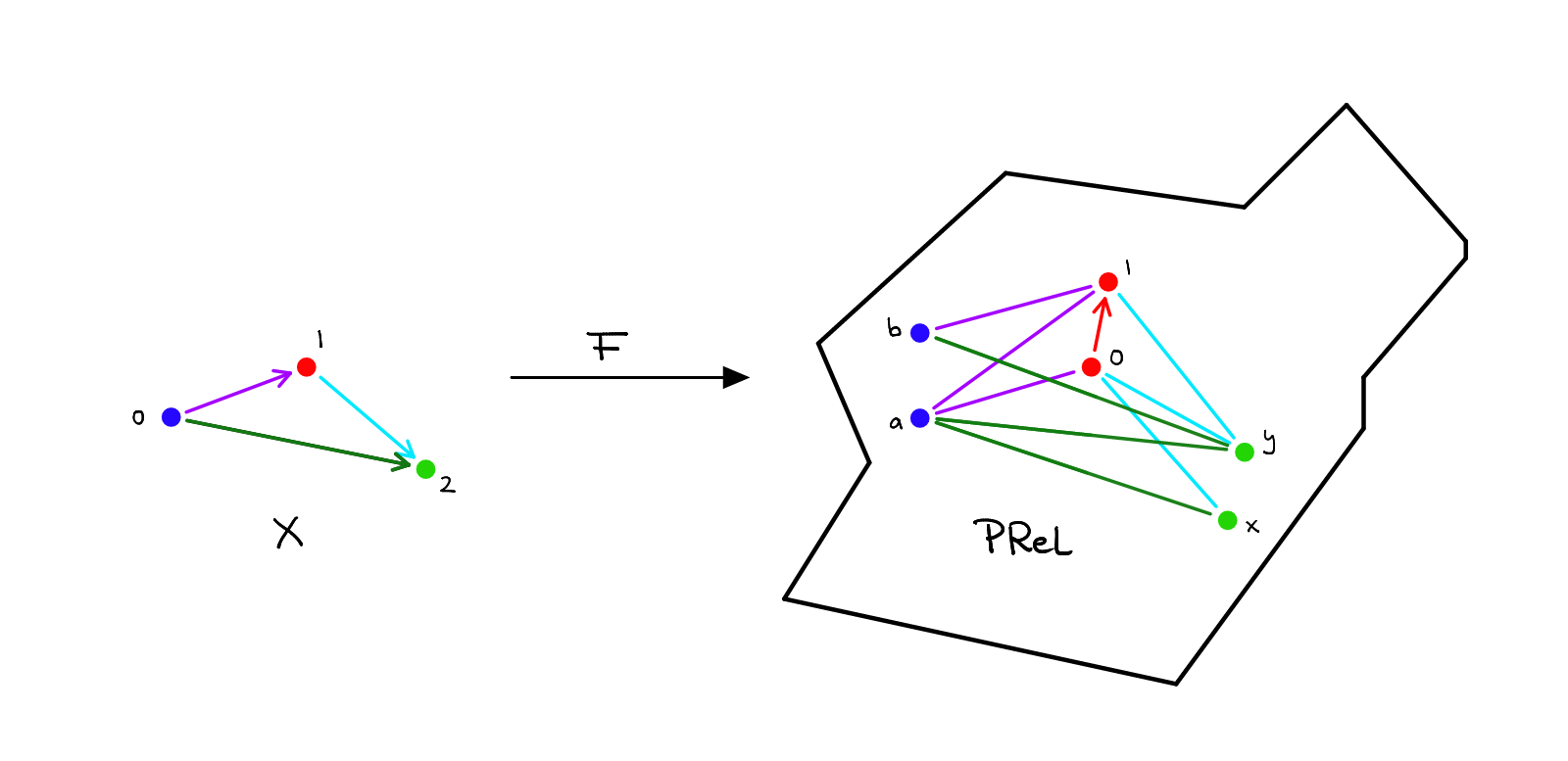}
\endgroup\end{restoretext}
Here, we colored posets and profunctorial relations in the image of $F$ correspondingly to their preimage (note that tuples in a profunctorial relation are indicated by edges, as opposed to arrows used for poset morphisms). Its associated \PRel-bundle is the map of posets given by
\begin{restoretext}
\begingroup\sbox0{\includegraphics{test/page1.png}}\includegraphics[clip,trim={.1\ht0} {.0\ht0} {.1\ht0} {.0\ht0} ,width=\textwidth]{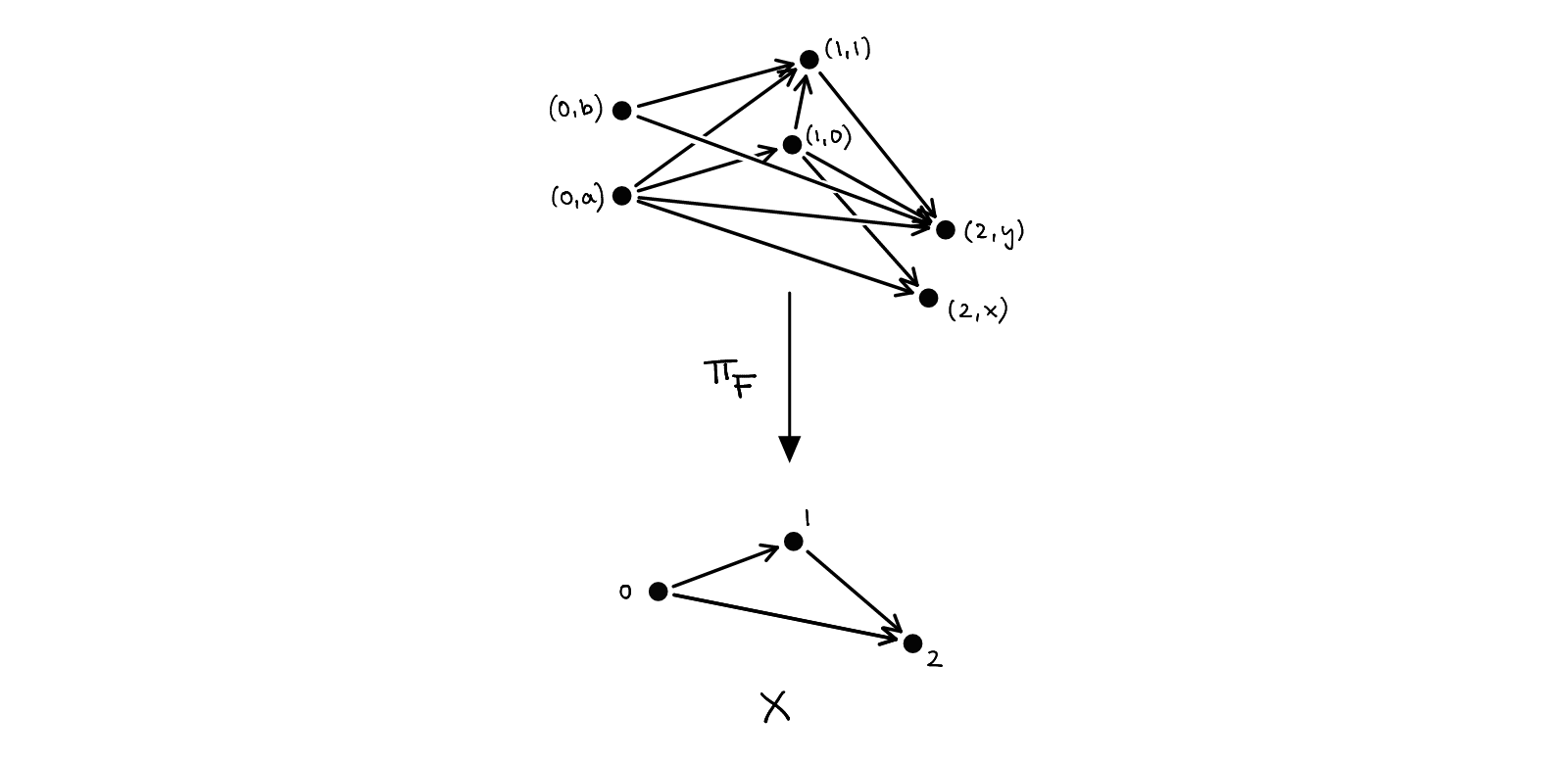}
\endgroup\end{restoretext}

\subsection{Families of singular intervals} \label{sec:sum_SI}

In this section we give definitions of combinatorial objects called singular intervals, and bundles thereof. A geometric reason for introducing them was given in \autoref{ssec:po_mfld_diag}. They play the role of ``1-dimensional fibers" in the inductive construction of $n$-dimensional space via bundles, which was mentioned in the beginning of the chapter.

\begin{defn} The category $\SI$ of singular intervals has as objects posets $\singint k$, $k \in \lN$. $\singint k$ has objects $\Set{0, 1, 2, \dots ,2k}$ and its order is generated by $a \to a + 1$ for even $a \in \singint k$, $a + 1 \to a$ for odd $a \in \singint k$. A morphisms $f : \singint k \to \singint l$ is a monotone map
\begin{equation} 
f : \Set{1, 3, \dots , 2k -1} \to \Set{1,3, \dots, 2l-1}
\end{equation}
Morphisms in $\SI$ compose as functions.
\end{defn}
\noindent Note that $\SI$ is isomorphic to $\BDelta$, the simplex category.

We next construct the natural embedding $\SiR: \SI \into \PRel$. Given a morphism $f : \singint k \to \singint l$ we first define its regular \rsdual{}
\begin{equation}
f\regop : \Set {0,2, \dots, 2l} \to \Set{0,2, \dots , 2k}
\end{equation}
by the ``ambidexterity condition": for $a \in \Set{-1, 1, 3, \dots , 2k -1, 2k+1}$ and $b  \in \Set {0,2, \dots, 2l}$, and setting $f(-1) := -1$ as well as $f(2k+1) := 2l+1$ then we define $f\regop$ by
\begin{align}
f(a) <  b &\iff a <  f\regop(b) \\
b <  f(a) &\iff f\regop(b) <  a
\end{align}
We note that actually one of those conditions is sufficient to define $f\regop$. Now, $\SiR(f)$ is the minimal profunctorial relation $\SiR(f) : \singint k \xslashedrightarrow{} \singint l$ satisfying that for all $a,c \in \Set{1, 3, \dots , 2k -1}$ and $b,d \in \Set{0,2, \dots, 2l}$
\begin{align}
f(a) = c &\imp \SiR(f)(a,c) \\
d = f\regop(b) &\imp \SiR(f)(d,b)
\end{align}
Here, minimality means that for any other $S : \singint k \xslashedrightarrow{} \singint l$ satisfying these condition we have $\SiR(a,b) \implies S(a,b)$, for all $a \in \singint k$, $b \in \singint I$. 

We give two examples of the construction
\begin{restoretext}
\begingroup\sbox0{\includegraphics{test/page1.png}}\includegraphics[clip,trim=0 {.15\ht0} 0 {.15\ht0} ,width=\textwidth]{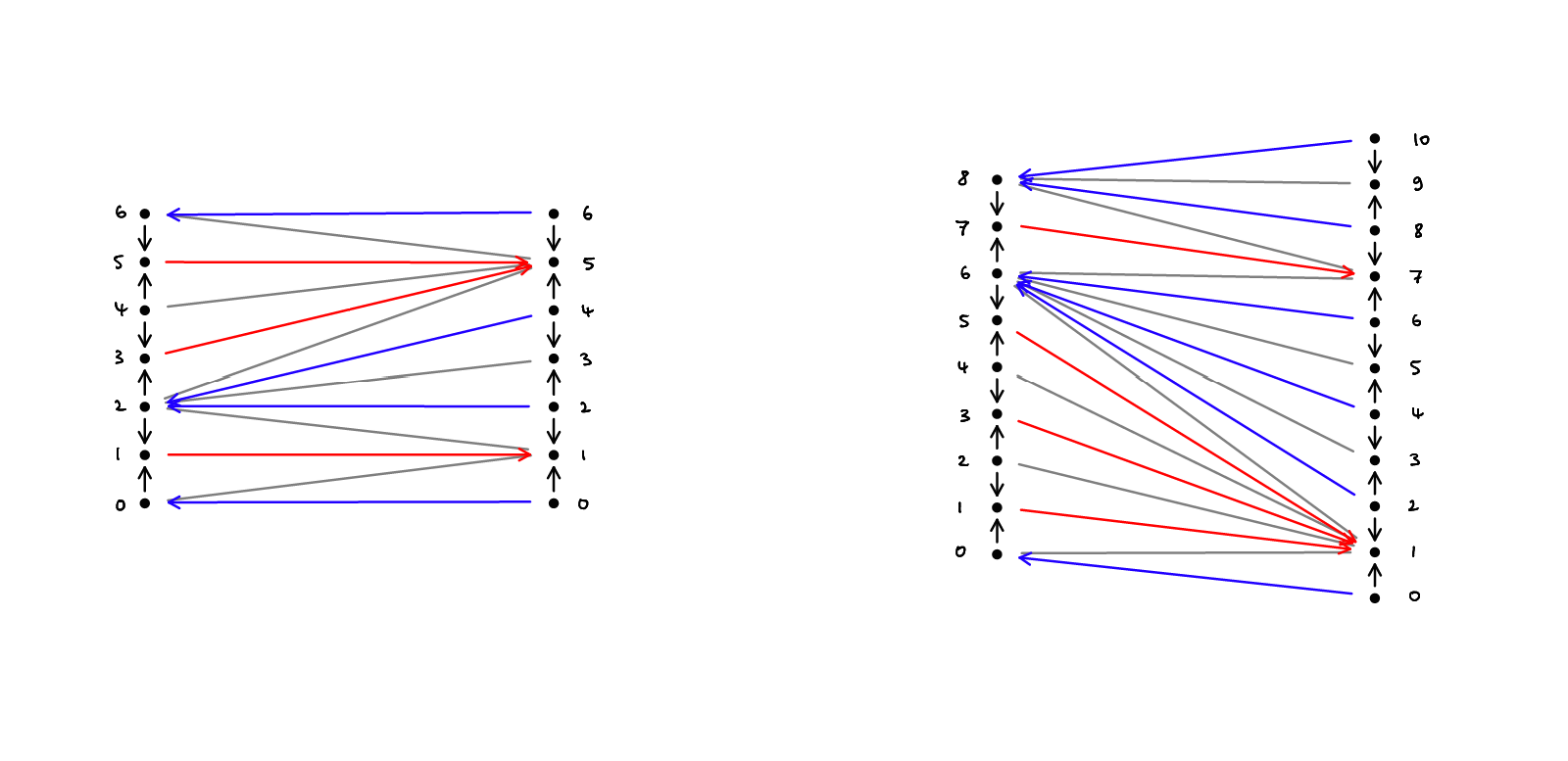}
\endgroup\end{restoretext}
Here, \cred{} arrows indicate morphisms of singular intervals. \textcolor{cblue}{Blue} arrows indicate their regular \rsdual{}s{}. Taking \cred{} arrows, \cblue{} arrows and \cgray{} edges together gives the edges of the resulting profunctorial relation in each case.

Given a poset $X$, a functor $F : X \to \SI$ is called a \textit{family of singular intervals}, or $\SI$-family for short. Such $F$ can be post-composed with $\SiR$ to yield a family $\SiR F : X \to \PRel$ of profunctorial relations.

\begin{notn}[\SI-bundles]\label{notn:SI_bundles_proj}  Recall the bundle construction from \autoref{ssec:sum_fam_bun_prel}. We will denote $\pi_{\SiR F} : \sG(\SiR F) \to X$ by $\pi_F : \sG(F) \to X$. $\pi_F$ is called a \textit{bundle of singular intervals}, or $\SI$-bundle for short.
\end{notn}

\subsection{Families of singular cubes}

We are now in the position to define labelled singular $n$-cube families, which consist of towers of singular interval bundles together with a ``labelling".

\begin{defn}[Singular $n$-cube families] Given a category $\cC$ and $n \in \lN$, a $\cC$-labelled singular $n$-cube family  $\scA$ indexed by a poset $X$ is a list of functors
\begin{align} \label{eq:tower_for_cube}
\tsU n_\scA : \tsG n(\scA) &\to \cC\\
\tusU {n-1}_\scA : \tsG {n-1}(\scA) &\to \SI \\
&\dots \\
\tusU {1}_\scA : \tsG {1}(\scA) &\to \SI \\
\tusU {0}_\scA : \tsG 0(\scA) &\to \SI 
\end{align}
satisfying $\tsG 0(\scA) = X$ and for $1 \leq k \leq n$ that
\begin{equation}
\tsG k(\scA) = \sG(\tusU {k-1}_\scA)
\end{equation}
Note that here, as stated in \autoref{notn:SI_bundles_proj}, $\sG(\tusU {k-1}_\scA)$ denotes the total space of the family $\SiR \tusU {k-1}_\scA : \tsG {k-1}(\scA) \to \PRel$ of profunctorial relations. $\tsU n_\scA$ is called the \textit{labelling functor} of $\scA$. 
\end{defn}

\begin{notn}[Bundle notation] Given a $\cC$-labelled singular $n$-cube family $\scA$ indexed by $X$ we denote for $0 \leq k < n$
\begin{equation}
\tpi {k+1}_\scA := \pi_{\tusU k_\scA}
\end{equation}
\end{notn}

\begin{rmk}[Defining $\cC$-labelled singular $n$-cube families in terms of bundles] \label{notn:SIn_bundles_proj}
Using the previous notation, we remark that $\scA$ can be defined by giving bundles $\tpi {k+1}_\scA$ (for $0 < k \leq n$) together with $\tsU n_\scA$, instead of giving $\tusU k_\scA$ (for $0 \leq k < n$) and $\tsU n_\scA$. This is because we can reconstruct $\tusU k_\scA$ from $\tpi {k+1}_\scA$ uniquely.
\end{rmk}

\begin{notn}[Shorthand for $n$-cube families] A $\cC$-labelled singular $n$-cube family $\scA$ indexed by $X$ is also called an $\SIvert n \cC$-family indexed by $X$. If further $X = \bnum{1}$ (where $\bnum{1}$ is the poset with a single object $0$), $\scA$ is usually just called a $\cC$-labelled singular \textit{$n$-cube}, or simply an $\SIvert n \cC$-cube.
\end{notn}

\begin{conv}[$(-1)$-cubes] \label{conv:minus_one_cubes} By convention there is a single $\SIvert {-1} \cC$-cube, namely the empty one, denoted by $\emptyset$. Note that this convention follows nicely from the above definition if we were to set $\tsG {-1} (\scA) = \emptyset$ (namely, the empty cube then takes the form of a map $\emptyset : \emptyset \to \cC$).
\end{conv}

As a first example recall the poset $\bnum{3}$ defined by
\begin{restoretext}
\begingroup\sbox0{\includegraphics{test/page1.png}}\includegraphics[clip,trim=0 {.35\ht0} 0 {.35\ht0} ,width=\textwidth]{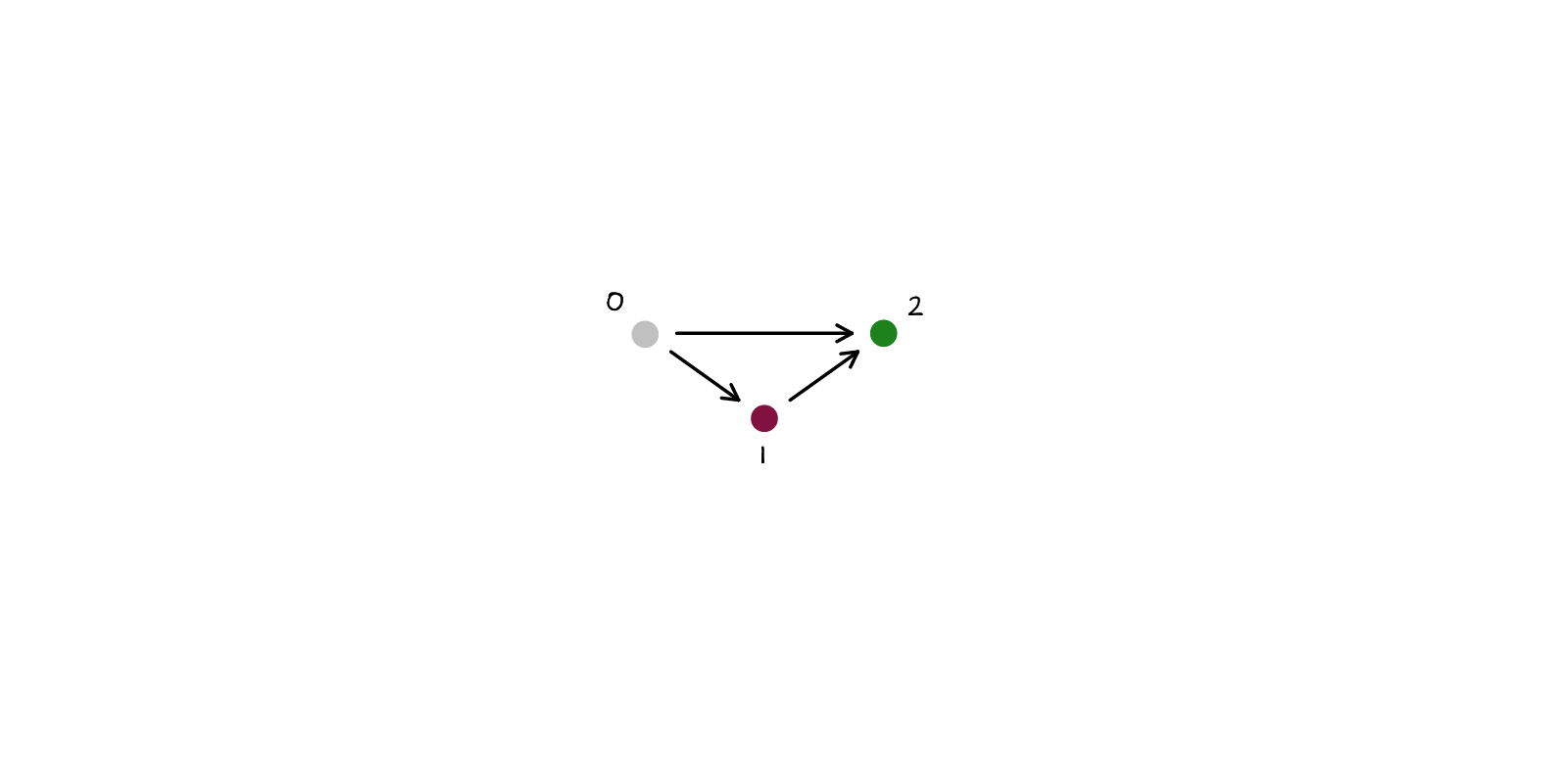}
\endgroup\end{restoretext}
Using \autoref{notn:SIn_bundles_proj}, this allows us to form the following $\bnum{3}$-labelled singular $2$-cube $Q$
\begin{restoretext}
\begingroup\sbox0{\includegraphics{test/page1.png}}\includegraphics[clip,trim=0 {.15\ht0} 0 {.05\ht0} ,width=\textwidth]{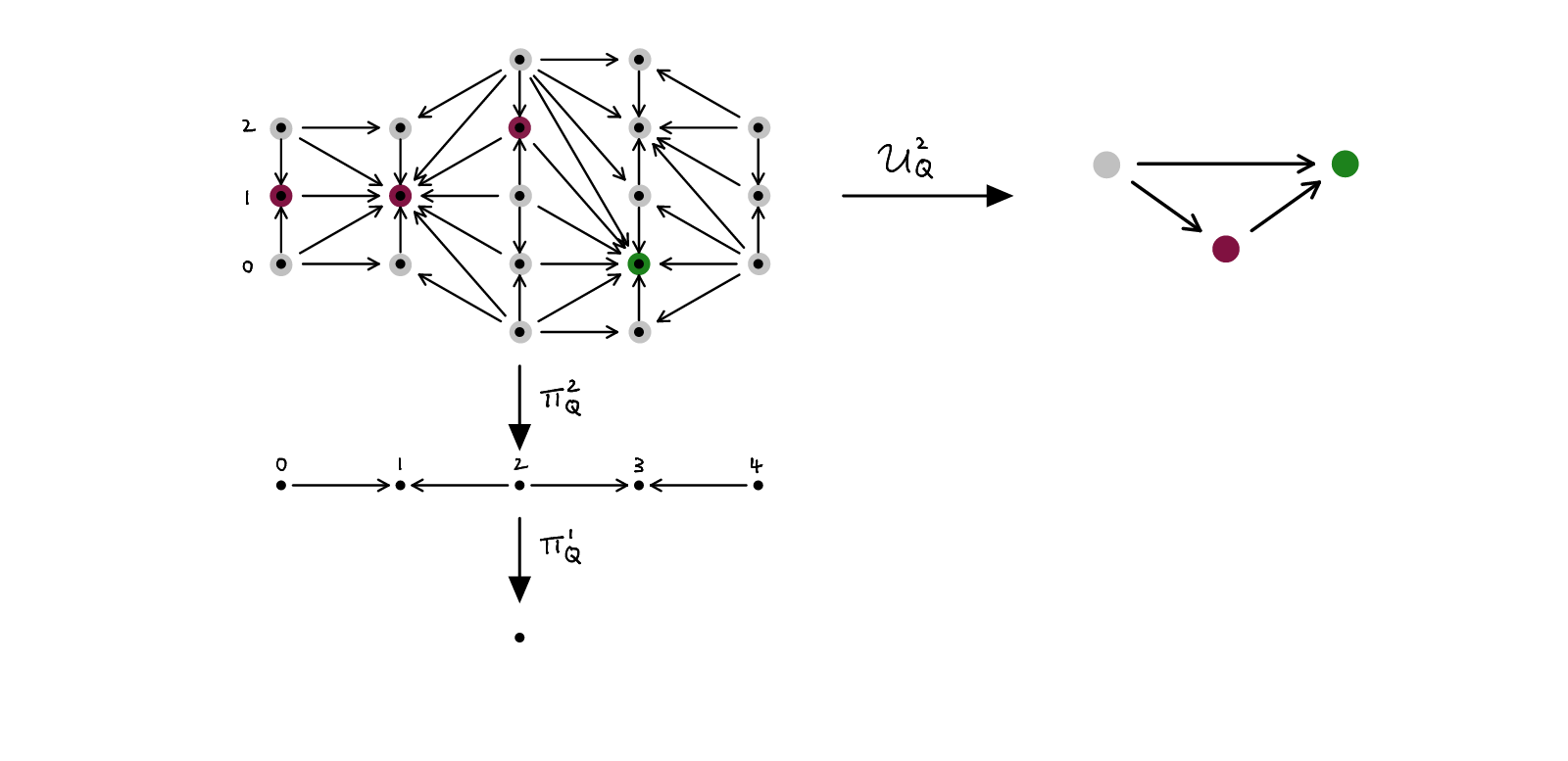}
\endgroup\end{restoretext}
Here and in the following, we often depict labelling functors by coloring preimages and their respective images by the same color.

As a second example example of a $\cC$-labelled singular $n$-cube we will algebraically reconstruct our example in \autoref{ssec:sum_rec_struct}. Using \autoref{notn:SIn_bundles_proj} and leaving certain poset arrows implicit for readability, we define a $\bnum{3}$-labelled singular $3$-cube $P$ as follows
\begin{restoretext}
\begingroup\sbox0{\includegraphics{test/page1.png}}\includegraphics[clip,trim=0 {.0\ht0} 0 {.0\ht0} ,width=\textwidth]{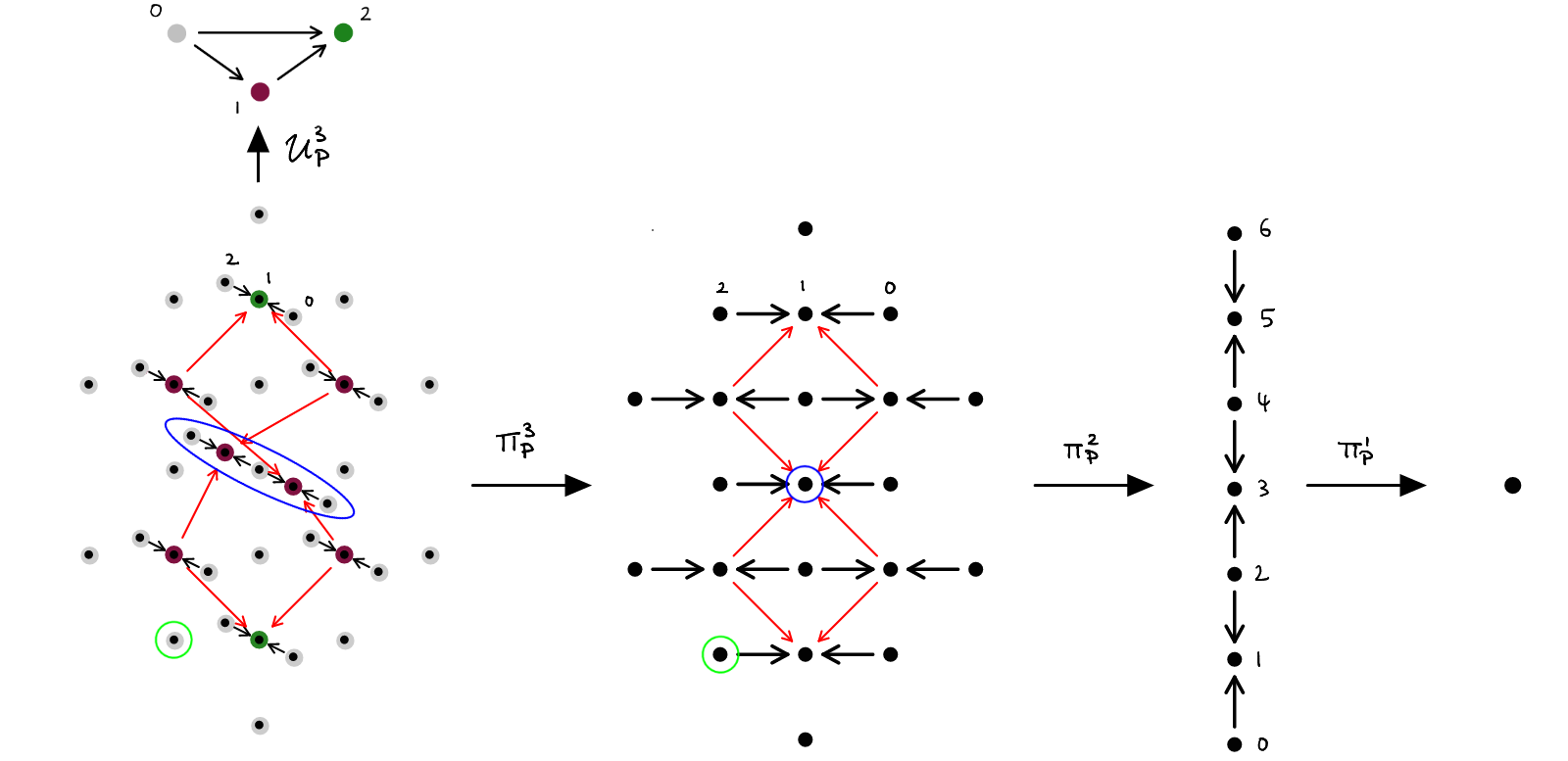}
\endgroup\end{restoretext}
Here, the maps $\tpi k_P$ should be understood from the spatial arrangement of their domain and codomain. For instance, the preimages of the \cblue{} and \cdarkgreen{}-circled points in the codomain of the bundle $\tpi 3_P$ are lying at the same position in the domain (and are indicated by \cblue{} and \cdarkgreen{} circles respectively). Also note that, instead of recording all arrows in the sets $\tsG k(P)$ we only recorded (in \cblack{}) arrows lying over objects (along $\tpi k_P$), or (in \cred{}) arrows lying over morphisms  which go between singular heights. Since profunctorial relations can be reconstructed from maps between singular heights (namely, by the functor $\SiR$) this notation scheme is sufficient to fully determine $\tsG k(P)$. This same notational scheme will be frequently applied later on.

\begin{rmk}[Non-typability] In fact, the above representation of the Hopf map will turn out to be a non-typable representation since local neighbourhoods of the same stratum are \textit{not} the same everywhere (cf. \autoref{rmk:typability} and \autoref{ssec:coloring}). A more correct and in-depth discussion of the Hopf map will be given in \autoref{ch:presented} and \autoref{ch:geom}.
\end{rmk}

Comparing to \autoref{ssec:sum_rec_struct}, the reader should observe that $\sG^k(P)$ is indeed the total poset of $P^k$ for $0 \leq k \leq 3$. Thus, the above is the combinatorial structure describing the projection-stable refinement of our original topological cube $P$. In the next section, we will discuss how to recover a topological labelled cube from this combinatorial structure by a process of ``geometric realisation". In particular, for our example we will find that the geometric realisation of the (combinatorial) singular cube $P$ recovers the topological cube $P$ up to equivalence.

\subsection{Geometric realisation} \label{ssec:coloring}

We now discuss a geometric realisation of (combinatorial) labelled singular $n$-cubes. The goal is to give a hands-on procedure to pass from labelled singular $n$-cubes to (topological) labelled $n$-cubes. We re-emphasize that mathematically this procedure does not carry any relevance for the development of the combinatorial ideas in this thesis. However, it provides a guiding intuition for many of the given constructions.

The procedure will consist of two steps. In the first step, we extract a from labelled singular $n$-cube $\scA$ the \textit{directed triangulation} $\abs{\scA}^n_\bullet$. Simplices of this triangulation will allow us to build a flag-foliation-compatible stratification of the $n$-cube (cf. \autoref{ssec:po_mfld_diag}) that is also stable under projection maps (cf. \autoref{ssec:sum_rec_struct}).  Namely, in the second step we obtain a $\cC$-labelled $n$-cube $\norm{\scA}$ by defining prestrata as certain unions of simplices of the directed triangulation, then ``pad" empty space in the cube linearly to turn prestrata into strata and add the labelling from $A$.

\begin{enumerate}
\item \textit{Step 1}: We first define the  \textit{directed triangulation} $\abs{\scA}^n_\bullet$ for a given $\cC$-labelled singular $n$-cube $\scA$. For a poset $X$ we say $f :  (\bnum {m+1}) \to X$ is a non-degenerate $m$-simplex of $X$ if $f$ is injective on objects. Inductively in $k=0,1,..., n$ we define $\abs{\scA}^k_\bullet$, the $k$-projected triangulation of $\scA$. This will be a sub-division\footnote{Here, ``sub-division" entails that simplices intersect only in a mutual face (or one is a face of the other). For details, we refer to the wealth of literature on triangulations and simplicial complexes.} of a subspace $\abs{\scA}^k$ of the $k$-cube $(0,1)^k$ into collections of (standard topological) $l$-simplices $\abs{\scA}^k_l$, $l \leq k$, which correspond exactly to non-degenerate $l$-simplices of $\sG^k(\scA)$, and satisfies that all (topological) simplices are convex. 

\begin{enumerate}
\item $\abs{\scA}^0_\bullet$ consists of the single $0$-simplex $0 \in \abs{\scA}^0_0$ sub-dividing the $0$-cube $\abs{\scA}^0 = (0,1)^0$. This corresponds to $\sG^0(\scA) = \bnum 1$ having a single $0$-simplex $0$.

\item Assume the sub-division $\abs{\scA}^k_\bullet$ of $\abs{\scA}^k \subset (0,1)^k$ has been defined. Let $\pi : (0,1) \times (0,1)^k \to (0,1)^k$ be the projection from $(k+1)$- to $k$-cube omitting the first coordinate. For each $0$-simplex $x \in \abs{\scA}^k_0$ (corresponding to $x \in \sG^k(\scA)$), distribute a $0$-simplex on $\pi\inv(x) \iso (0,1)$ for each $a \in \und\sU^k_\scA(x)$ in a strictly monotone (but otherwise arbitrary) way. Add these $0$-simplices in $\pi\inv(x)$ to $\abs{\scA}^{k+1}_0$. Note that this gives a correspondence of $0$-simplices in  $\abs{\scA}^{k+1}_0$ to $0$-simplices $(x,a)$ in $\sG^{k+1}(\scA)$ as claimed.

Now, for each non-degenerate $l$-simplex $f$ in $\sG^{k+1}(\scA)$ we add to $\abs{\scA}^k_l$ the (topological) $l$-simplex defined as the convex closure of the points $f(0), f(1), ... f(l-1)$. This indeed yields a (topological) $l$-simplex, since it can be shown that all $1$-simplices $f(i \to (i+1))$ are linearly independent\footnote{To prove this fact, note that odd numbers never have arrows to even numbers (cf. \autoref{sec:sum_SI}), that is, $f(0), f(1), .., f(l-1)$ must be a sequence of the form
\begin{equation}
(x_0,a_0) , (x_1,a_1), ... , (x_j,a_j) , (x_{j+1},b_{j+1}) , ... , (x_{l-2},b_{l-2}), (x_{l-1},b_{l-1})
\end{equation}
where $a_i$ are even numbers, $b_i$ are odd numbers and $-1 \leq j \leq l-1$. Non-degeneracy of $f$ implies $x_i \neq x_{i + 1}$ unless $i = j$. Linear independence thus follows inductively: $\pi(f(x_i \to x_{i+1}))$ (excluding $i = j$) are linearly independent inductively and $f(x_j \to x_{j+1})$ lies in the kernel of $\pi$. Note that the same argument inductively shows that $\sG^{k+1}(\scA)$ can have non-degenerate $l$-simplices only for $l \leq k + 1$.}. 

Finally, taking the union of all such topological simplices yields the space $\abs{\scA}^{k+1}$. Together with its subdivision into simplices we have thus defined $\abs{\scA}^{k+1}_\bullet$. 
\end{enumerate}
This completes the inductive construction of the directed triangulation $\abs{\scA}^n_\bullet$ of $\abs{\scA}^n \subset (0,1)^n$. The reader familiar with simplical sets will be able to verify that we have in fact constructed a special instance of the geometric realisation of the nerve of $\sG^n(\scA)$ (cf. \cite{riehl2014categorical}): this specific construction is necessitated by desired compatibility with the cube projection $\pi$. Namely, if we apply $\pi$ to $\abs{\scA}^{k+1}_\bullet$ we recover the triangulation $\abs{\scA}^k_\bullet$.

\item \textit{Step 2}: Note that $\abs{\scA}^n_\bullet$ is the disjoint union of the interiors $f^\circ$ of its non-degenerate simplices $f$. We define the set of \textit{prestrata} $\mathrm{Pre}(\scA)$ of $\scA$ by setting
\begin{equation}
\mathrm{Pre}(A) = \Set{\cup_{f(0) = x} f^\circ }_{x \in \abs{\scA}^n_0}
\end{equation}
where the union runs over non-degenerate $k$-simplices $f$ in $\abs{\scA}^n_\bullet$. We now define the set of \textit{strata} $\mathrm{Str}(A)$ by linearly (and inductively) extending prestrata along the directions $k = 1,2,.., n$ to fill the $n$-cube $(0,1)^n$. It is straight-forward to prove that this extension yields a flag-foliation-compatible stratification of the $n$-cube. One way to formally describe this extension construction is as follows.
\begin{enumerate}
\item We set $\mathrm{Str}_0(A) = \mathrm{Pre}(A)$
\item To define $\mathrm{Str}_{k+1}(A)$, let
\begin{equation}
I^{a,b}_{p,q} = \Set{p \in (0,1)^k} \times (a,b) \times \Set{q \in (0,1)^{n-k}}
\end{equation}
be an open $(a,b)$-interval pointing in direction $k$. Note $I^{a,b}_{p,q}  = I^{b,a}_{p,q}$. Define $A_{p,q}^{-}, A_{p,q}^{+}$ to be the lower respectively upper boundary of the intersection  
\begin{equation}
I^{0,1}_{p,q} \cap \abs{A}^n
\end{equation}
(if it is non-empty). For each $u \in \mathrm{Str}_k(A)$ define its ``linear extension in direction $k$" by setting
\begin{equation}
\hat u = u \cup \Set{I^{A_{p,q}^{\pm},\frac{1}{2}\pm\frac{1}{2}}_{p,q} ~|~ A_{p,q}^{\pm} \in u} 
\end{equation}
Then set
\begin{equation}
\mathrm{Str}_{k+1}(A) = \Set{ \hat u ~|~ u \in \mathrm{Str}_k(A)}
\end{equation}
\end{enumerate}
As a result we obtain $\mathrm{Str}(A) := \mathrm{Str}_n(A)$, which forms a flag-foliation-compatible statification of the $n$-cube $(0,1)^n$ and comes with a bijection $\mathrm{bp} : \mathrm{Str}(A) \iso \mathrm{Pre}(A) \iso \abs{A}^n_0$. 

Now we construct a (topological) $\cC$-labelled $n$-cube $\norm{A}$ as follows: using the flag-foliation-compatible stratification $\mathrm{Str}(A)$ that we already built, It remains to build a labelling functor from the underlying poset of the stratification to $\cC$. Using the isomorphism $\mathrm{bp}$ we see that the total poset (cf. \autoref{ssec:po_mfld_diag}) is in fact just $\sG^n(A)$, and consequently we can just use the existing functor $\sU^n_A : \sG^n(A) \to \cC$ as our labelling functor for $\norm{A}$. This completes the construction.
\end{enumerate}

We give examples of the preceding construction.

\begin{eg}[Directed triangulations and prestrata] We apply this procedure for obtaining prestrata to the two examples of labelled singular $n$-cubes $Q$ and $P$ that were previously discussed. Since prestrata are in correspondence with strata, the labelling functor associates to each prestratum a label (or ``color"). For $Q$ we obtain the following ($k$-projected) directed triangulations and (colored) prestrata
\begin{restoretext}
\begingroup\sbox0{\includegraphics{test/page1.png}}\includegraphics[clip,trim={.0\ht0} {.1\ht0} {.0\ht0} {.0\ht0} ,width=\textwidth]{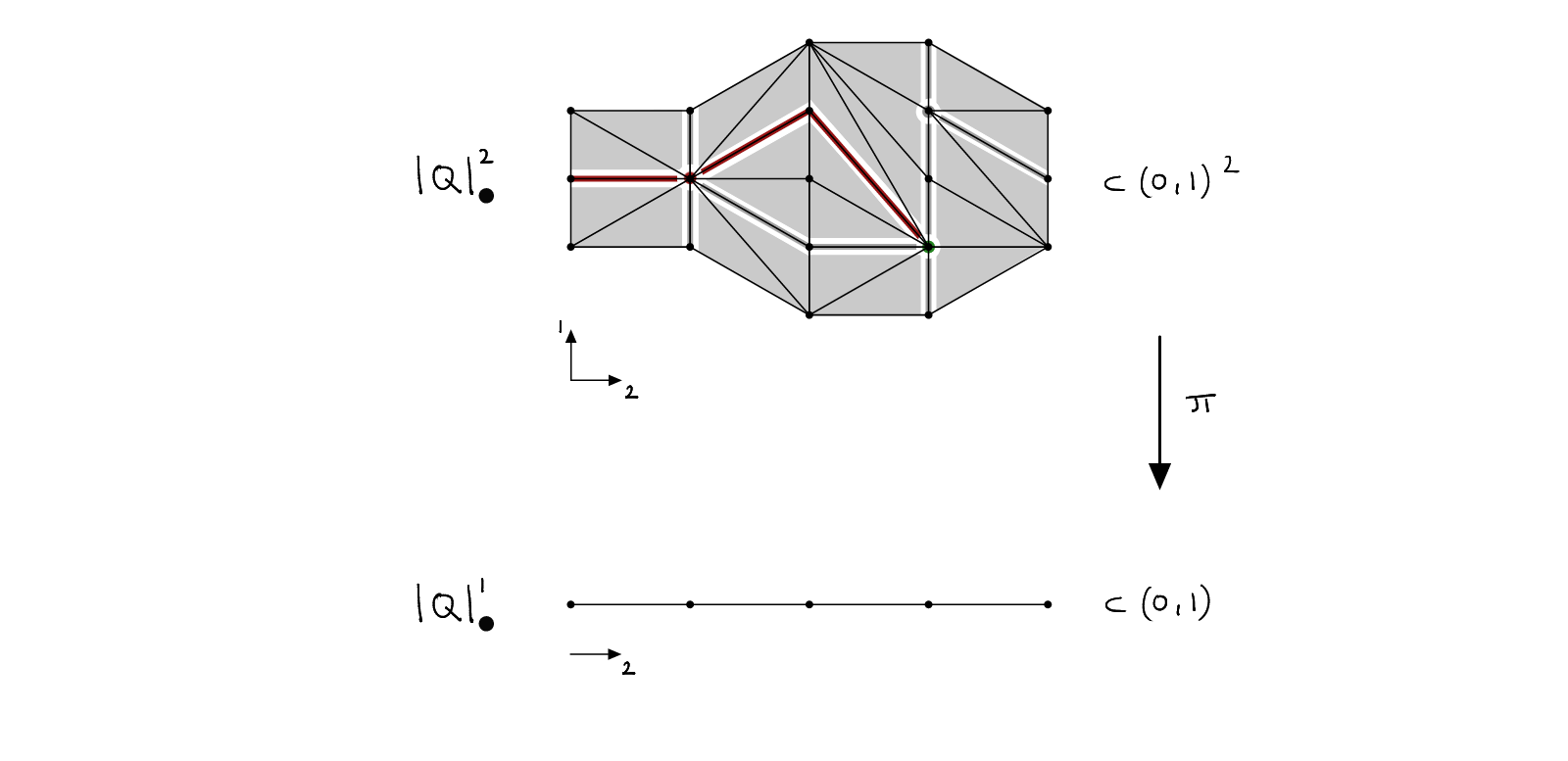}
\endgroup\end{restoretext}
Interiors of simplices in $\abs{Q}^2_\bullet$ were joined to obtain the set of prestrata $\mathrm{Pre}(Q)$ whose elements are indicated as ``connected colored regions" with coloring according to $\sU^2_Q$. The above should be compared to the previously defined posets $\tsG 2(Q)$, $\tsG 1(Q)$ and the labelling $\tsU 2_Q$. In particular, note that $\mathrm{Pre}(Q)$ is indeed in bijective correspondence with $\tsG 2(Q)$. The projection map $\pi : (0,1)^2 \to (0,1)$ has been indicated by an arrow. Note that we also added two coordinate axes 1 and 2 of the cube $(0,1)^2$ (of which $\abs{Q}^2$ is a subspace) and a single coordinate axis $2$ for the target $(0,1)$ of the projection $\pi$ (of which $\abs{Q}^1$ is a subspace). By virtue of the construction (namely the ``monotonicity" condition in the inductive step) these correspond to the directions of singular intervals at level 1 and 2 of the tower of bundles of $Q$, which were previously indicated by two sequences of numbers in the definition of $Q$.

The above subspace $\abs{Q}^2$ happens to be homeomorphic to the square $I^2$ but this does not hold in general as we will see in the case of our previous example $P$. For $P$, building (and coloring) prestrata according to the above procedure, we obtain $\abs{P}^2_\bullet$ and $\mathrm{Pre}(P)$ as follows 
\begin{restoretext}
\begingroup\sbox0{\includegraphics{test/page1.png}}\includegraphics[clip,trim={.2\ht0} {.05\ht0} {.2\ht0} {.15\ht0} ,width=\textwidth]{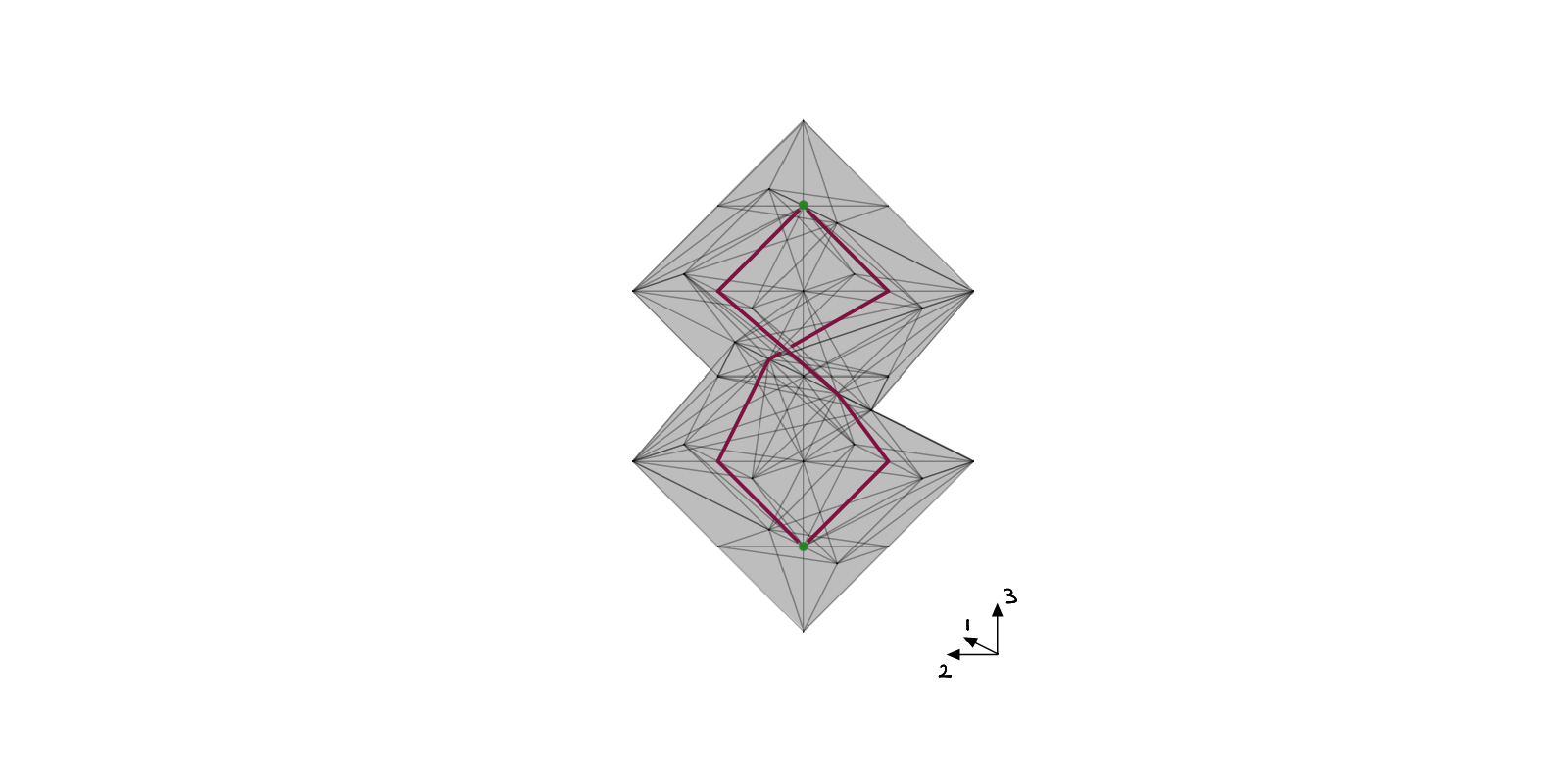}
\endgroup\end{restoretext}
This should be compared with the poset $\tsG 3(P)$ and the labelling map $\tsU 3_P$. Again, note that 3 coordinate axes of the $3$-cube $(0,1)^3$ have been added which agree with the directions of singular intervals (previously indicated by 3 sequences of numbers in the data of $P$). The picture should also be compared to the topological cube $P$ in \autoref{ssec:sum_rec_struct} (note that there, coordinate axis $2$ is shown in opposite direction).

Note that in this example, the space $\abs{P}^3$ is indeed not homeomorphic to the cube $(0,1)^3$. Such ``degeneracy" can occur when fibers of $\SI$-bundles equal the initial singular interval $\singint 0$ (which realises to a point and is thus not homeomorphic to the interval). Thus, in general, prestrata do not form a stratification of the cube (or of a space homeomorphic to it).
\end{eg}

This last remark underlines the importance of ``linear extending" prestrata to strata as described in the above procedure. This is visualised in the following example.

\begin{eg}[Linearly extending prestrata]
We illustrate the procedure in the case of $Q$. We first explicitly visualise the embedding $\abs{Q}^2 \subset (0,1)^2$ as follows
\begin{restoretext}
\begingroup\sbox0{\includegraphics{test/page1.png}}\includegraphics[clip,trim=0 {.0\ht0} 0 {.0\ht0} ,width=.8\textwidth]{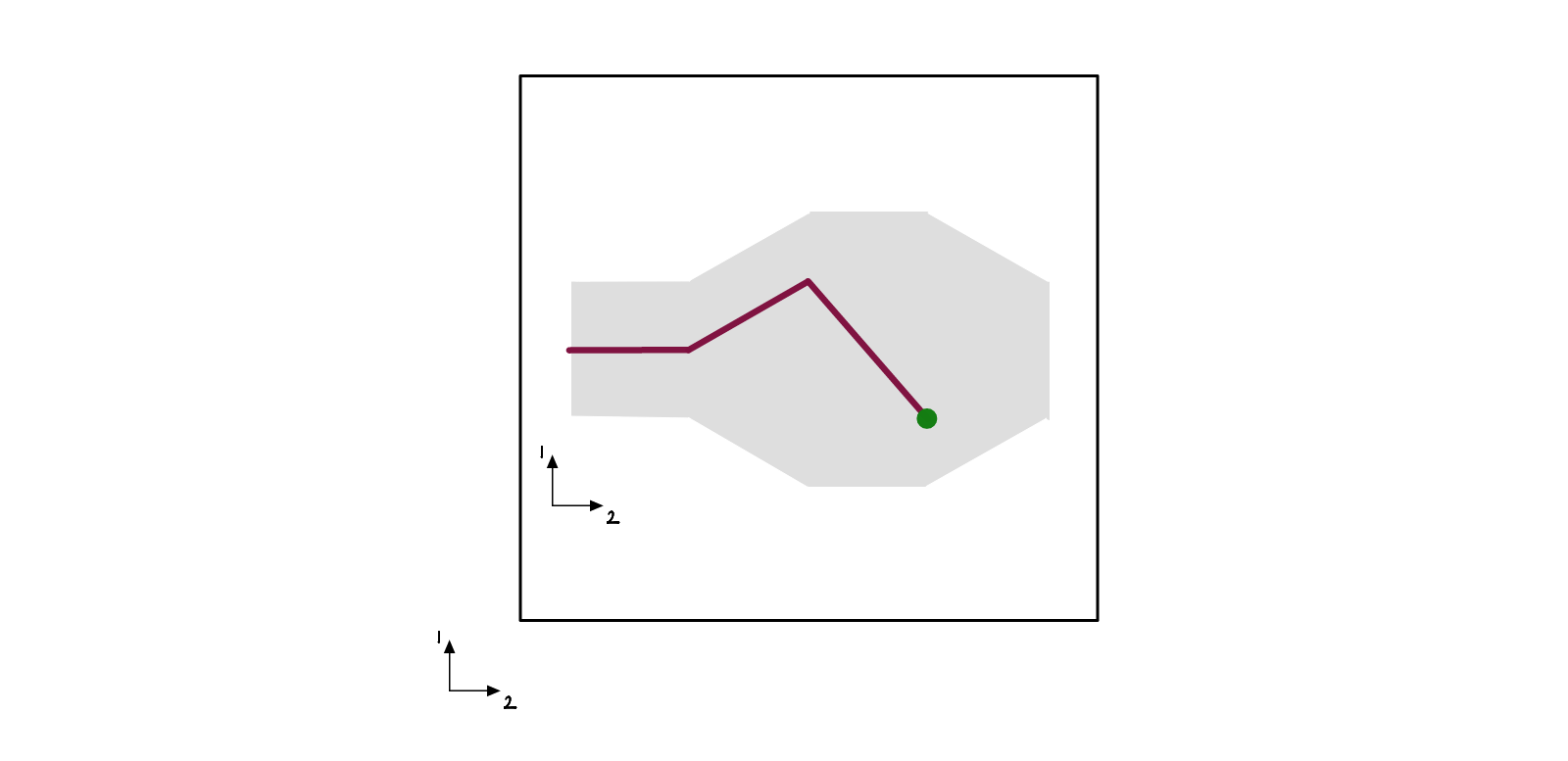}
\endgroup\end{restoretext}
Then we extend prestrata linearly in direction $1$, obtaining
\begin{restoretext}
\begingroup\sbox0{\includegraphics{test/page1.png}}\includegraphics[clip,trim=0 {.0\ht0} 0 {.0\ht0} ,width=.8\textwidth]{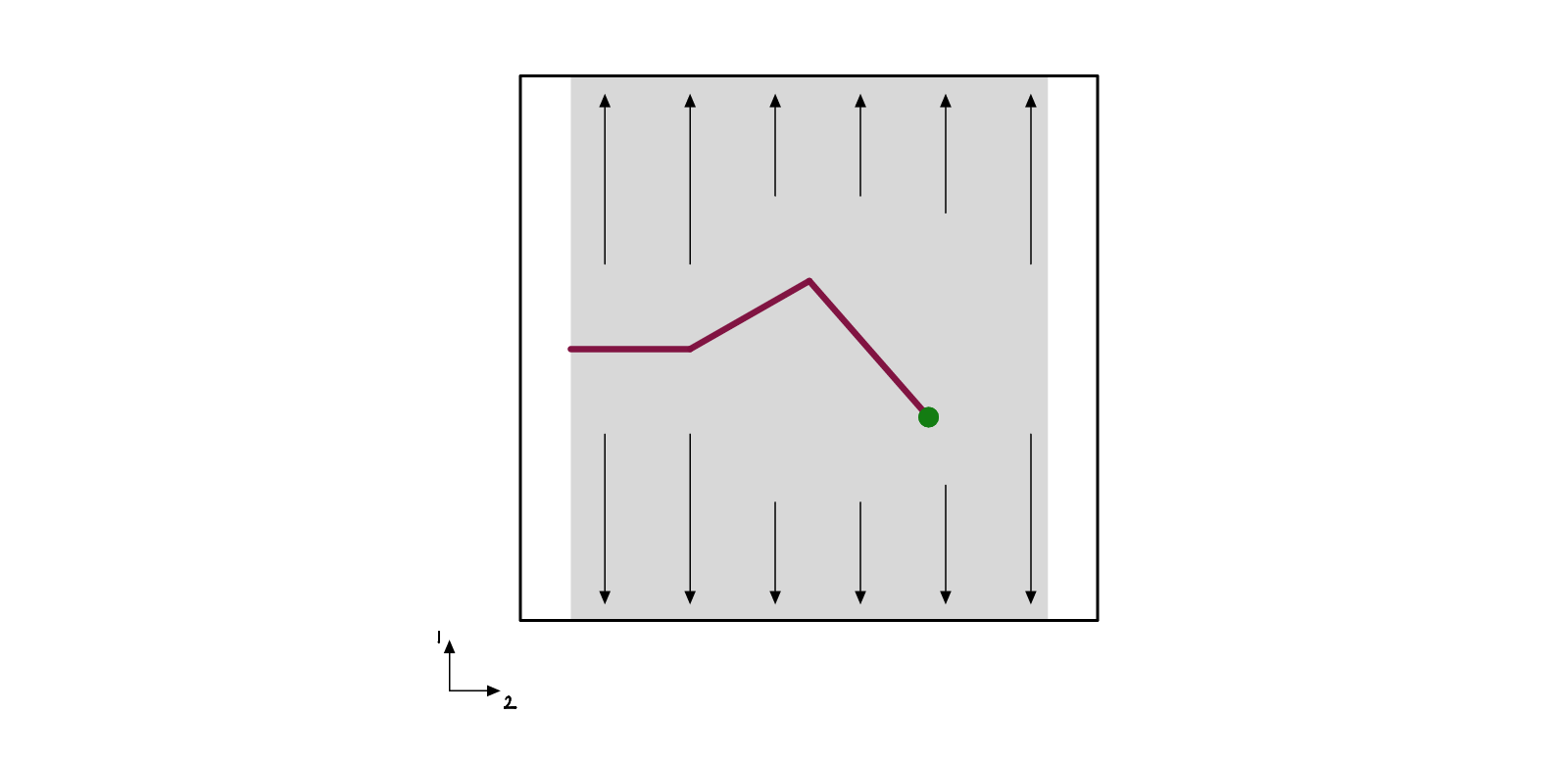}
\endgroup\end{restoretext}
Finally, we extend the resulting subspaces linearly in direction $2$, obtaining
\begin{restoretext}
\begingroup\sbox0{\includegraphics{test/page1.png}}\includegraphics[clip,trim=0 {.0\ht0} 0 {.0\ht0} ,width=.8\textwidth]{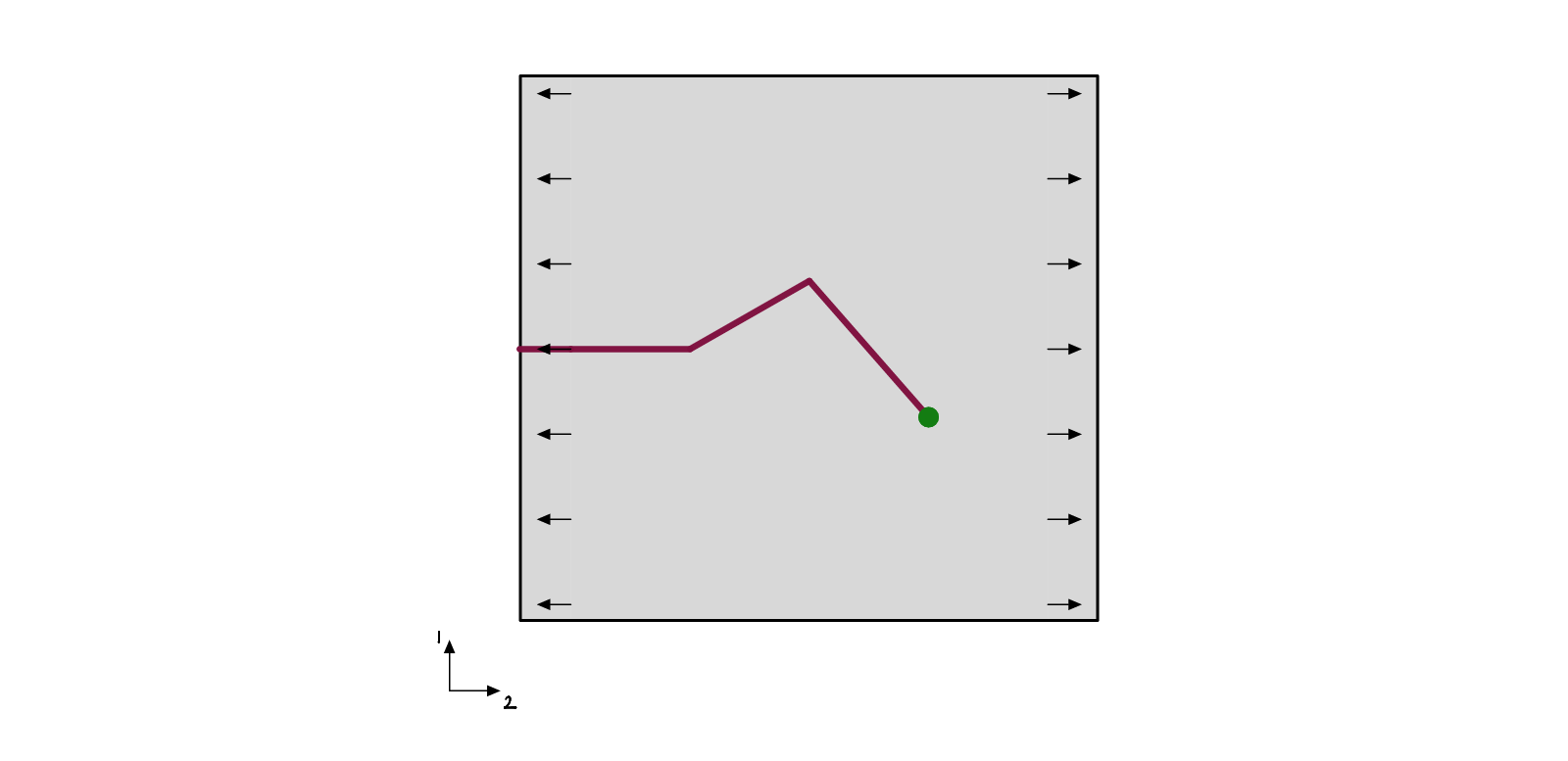}
\endgroup\end{restoretext}
This is the (topological) $\bnum 3$-labelled $2$-cube $\norm{Q}$ built from $Q$.

A similar extension of strata for $P$ recovers our labelled $3$-cube used in \autoref{ssec:sum_rec_struct}.
\end{eg}

\begin{rmk}[Different choices of $\norm{A}$ for given $A$ are isomorphic] Importantly, in the construction of the directed triangulation $\abs{A}^n_\bullet$ there were arbitrary choice to be made (when distributing singular heights over the topological interval). It is easy to see that all of these choices yield isomorphic cubes in $\Cubeo n \cC(\bnum 1)$. Thus $Q \mapsto \norm{Q}$ is well-defined up to isomorphism. 

For example, a different choice of triangulation for $Q$ could have resulted in the cube
\begin{restoretext}
\begingroup\sbox0{\includegraphics{test/page1.png}}\includegraphics[clip,trim=0 {.20\ht0} 0 {.25\ht0} ,width=\textwidth]{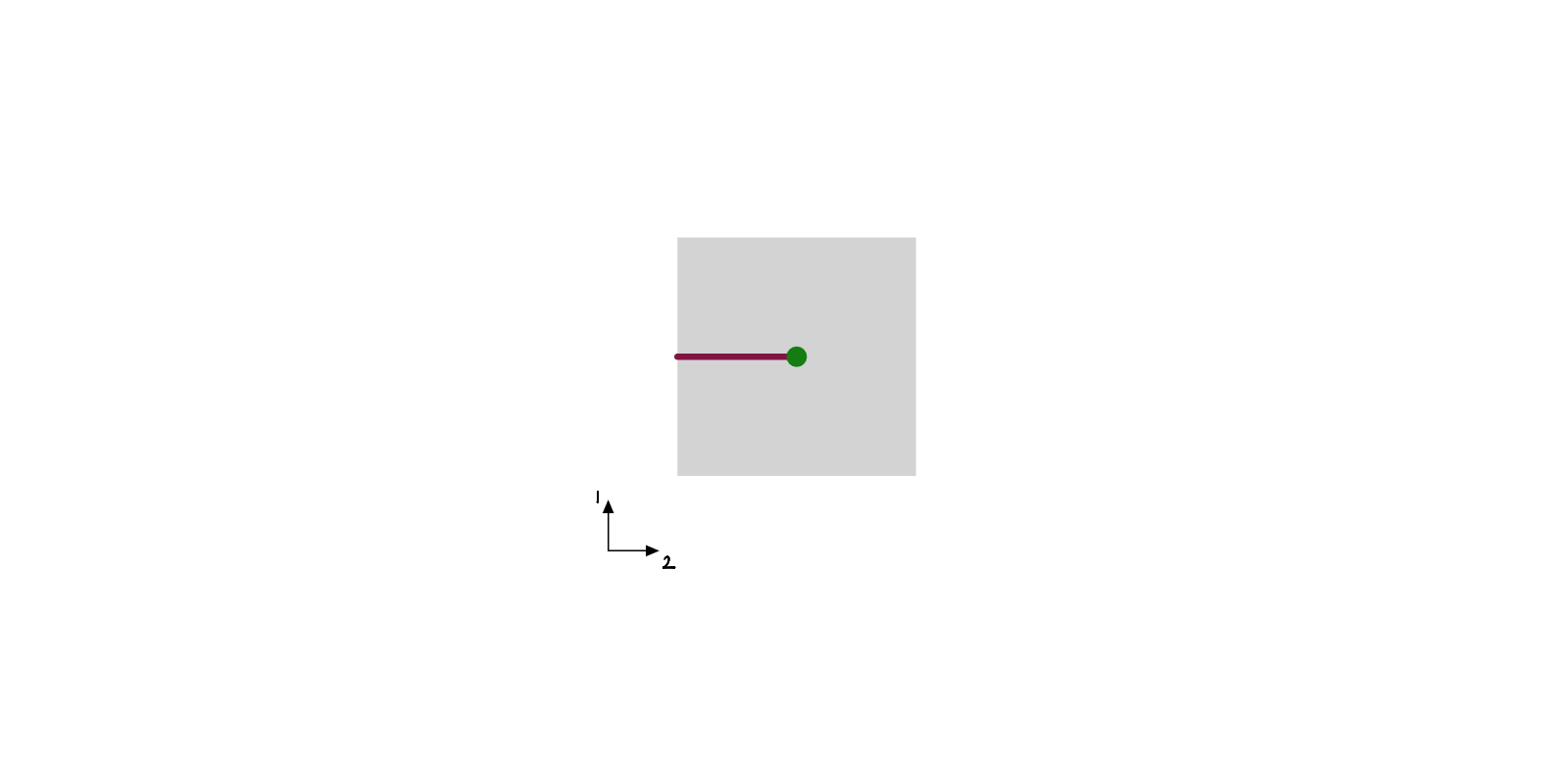}
\endgroup\end{restoretext}
\end{rmk}

\begin{rmk}[Combinatorial manifold diagrams] The reader will have noticed that the labelled $n$-cubes for $\norm{Q}$ and $\norm{P}$ look like manifold diagrams and not just mere labelled cubes. This will hold true in general only if we restrict to special $\SIvert n \cC$-cubes satisfying a combinatorial globularity condition which will be described in \autoref{sec:sum_glob}.
\end{rmk}

\begin{rmk}[Correspondence of topological cubes up to equivalence and combinatorial cubes up to normalisation] As remarked in \autoref{ssec:sum_rec_struct} already, there are (infinitely) many ways to ``refine" the structure of a topological $\cC$-labelled $n$-cube, obtaining ``projection-stable" underlying stratifications that can be captured combinatorially as a $\cC$-labelled singular $n$-cube. However, each such singular $n$-cube then geometrically realises to the same topological $n$-cube up to equivalence. What combinatorial relationship do all these singular $n$-cubes have? The answer is, that they all share the same \textit{normal form}, which can be thought of as the unique simplest (or ``coarsest") projection-stable stratification of a given topological cube. 

This notion can be captured in a completely combinatorial manner: two singular $n$-cubes $A$, $B$ have the same normal form if they are related by a cospan $A \to C \ot B$ of epimorphisms in a category $\Buno n \cC$, which will be defined shortly, and the normal form itself is the terminal object of the connected component of $A$ and $B$ (in the subcategory spanned by epimorphisms), see \autoref{sec:sum_NF}. The (constructive) existence of normal forms will guarantee that the question of ``geometric equivalence", meaning\footnote{up to a proof that our two notions cube equivalence (namely, the combinatorial and topological notion) are preserved and reflected by $\norm{-}$.} the question, whether $\norm{A} \simeq \norm{B}$, is decidable. This is the basic fact making our model of higher categories computer implementable.

As an example, we will see that $Q$ is \textit{not} normalised, and indeed it is \textit{not} an example of a coarsest directed triangulation (or equivalently, a coarsest projection-stable stratification) of its manifold diagram. We will meet its normal form $\NF{Q}$ shortly.
\end{rmk}

Further to the previous remark, we have the following.

\begin{rmk}[Is geometric realisation a functor?] The mapping $Q \mapsto \norm{Q}$ should give rise to a functor (up to slight modifications of the codomain)
\begin{equation}
\Buno n \cC (\bnum 1) \to \Cubeo n \cC(\bnum 1)
\end{equation}
where the domain $\Buno n \cC (\bnum 1)$ will be defined shortly. However, this construction lies beyond the goal of this thesis, which is to study the combinatorial side of the story.
\end{rmk}

\subsection{$k$-level labelling and relabelling}

In this short section we make two observations about the behaviour of labels of labelled singular $n$-cubes.

\begin{notn}[$k$-level labellings] \label{notn:sum_und_drop_notn} Note that by truncating the list \eqref{eq:tower_for_cube} after $k$ elements we obtain a $\cC$-labelled singular $k$-cube family indexed by $\tsG k(\scA)$
\begin{align}
\tsU n_\scA : \tsG n(\scA) &\to \cC\\
\tusU {n-1}_\scA : \tsG {n-1}(\scA) &\to \SI \\
&\dots \\
\tusU {k}_\scA : \tsG {k}(\scA) &\to \SI
\end{align}
which will be denoted by $\tsU k_\scA$, and called the \textit{$k$-level labelling of $\scA$}.
\end{notn}

For instance, using our above definition of $Q$, we find the family $\tsU 1_Q$ to equal
\begin{restoretext}
\begingroup\sbox0{\includegraphics{test/page1.png}}\includegraphics[clip,trim=0 {.4\ht0} 0 {.0\ht0} ,width=\textwidth]{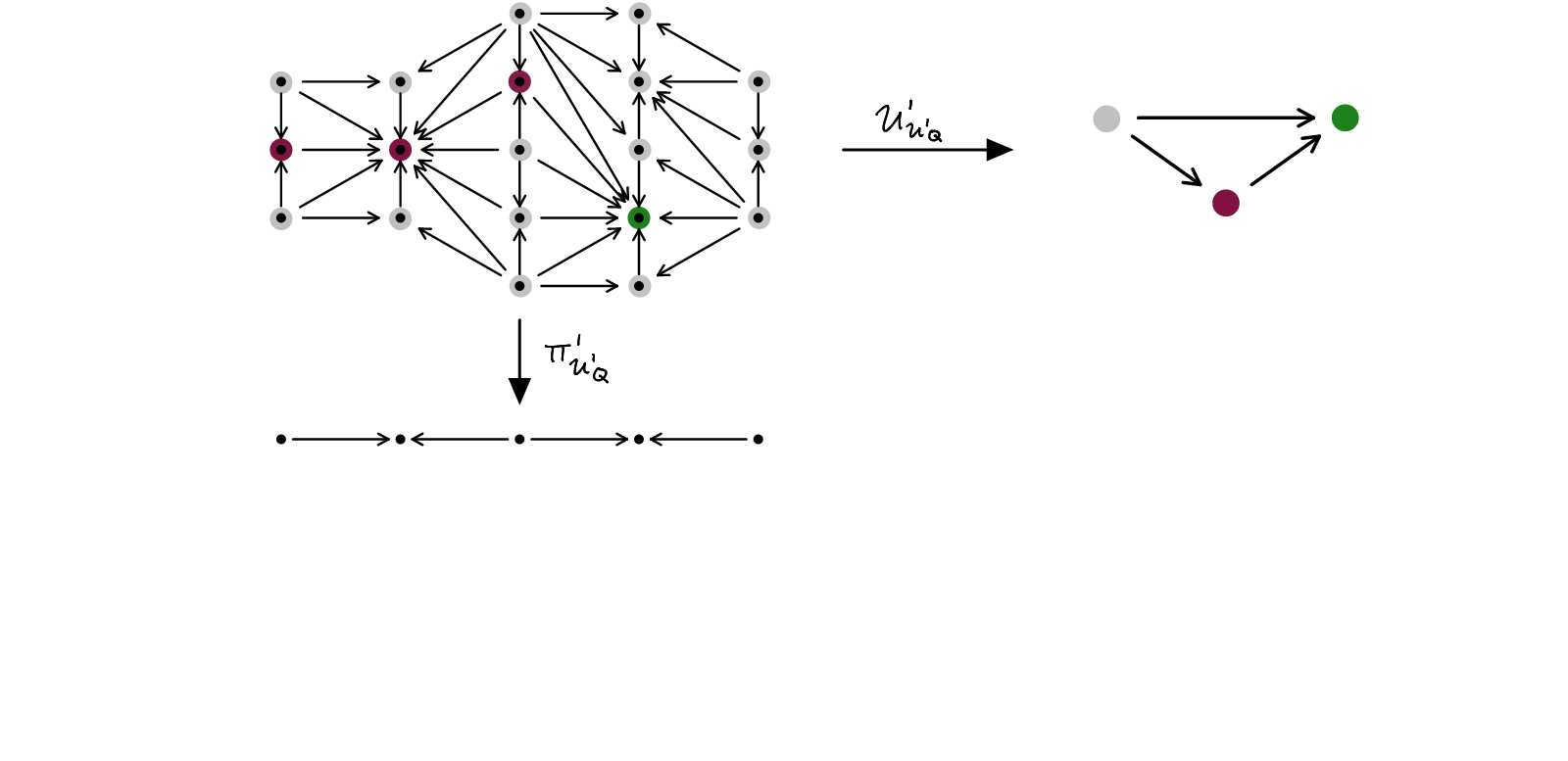}
\endgroup\end{restoretext}

\begin{defn}[Relabelling families] \label{def:sum:relabelling} Let $\scA$ be an $\SIvert n \cC$-family indexed by $X$. Let $F : \cC \to \cD$ be a functor of categories. Then $\SIvert n F \scA$ denotes the $\SIvert n \cD$-family $\scB$ indexed by $X$ with data $\tusU k_\scB = \tusU k_\scA$ (for $k < n$) and $\tsU n_\scB = F \tsU n_\scA$.
\end{defn}

Morally, relabelling by a functor $F$ leaves the geometric structure of a cube family in place but changes its labels according to $F$.

\subsection{Pullback families}

Often when working with bundles, given a map into the base space of a bundle, we can ``pullback" the bundle along this map. We now define this pullback in the context of towers of \SI-bundles.

\begin{defn}[Pullback of families] \label{def:sum:pullbacks_of_families} Given a $\cC$-labelled singular $n$-cube family $\scA$ indexed by a poset $X$ and a map $H : Y \to X$, then we define a $\cC$-labelled singular $n$-cube family $\scA H$ indexed by $Y$ as follows. Set $\tsG 0(H) = Y$, and inductively define $\tsG {k+1}(H)$ and $\tusU k_{\scA H}$ (for $0 \leq k < n$) by the pullbacks
\begin{equation}
\xymatrix{ \tsG k(\scA H) \pullbackfar \ar[r]^{\tsG k(H)} \ar[d]_{\pi_{\tusU {k-1}_{\scA H}}} & \tsG k(\scA) \ar[d]^{\pi_{\tusU {k-1}_\scA}} \\
\tsG {k-1}(\scA H) \ar[r]_{\tsG {k-1}(H)} & \tsG {k-1}(\scA) }
\end{equation}
Further set
\begin{equation}
\tsU n_{\scA H} = \tsU n_\scA \tsG n(H)
\end{equation}
This completes the definition of $\scA H$ (cf. \autoref{rmk:poset_data_objects}).
\end{defn}

As a first example, setting $V = \tsU 1_Q$ consider the map $H : \singint 1 \to \singint 2$ which maps $i \mapsto i + 2$. Then we obtain $VH$ as follows
\begin{restoretext}
\begingroup\sbox0{\includegraphics{test/page1.png}}\includegraphics[clip,trim=0 {.0\ht0} 0 {.0\ht0} ,width=\textwidth]{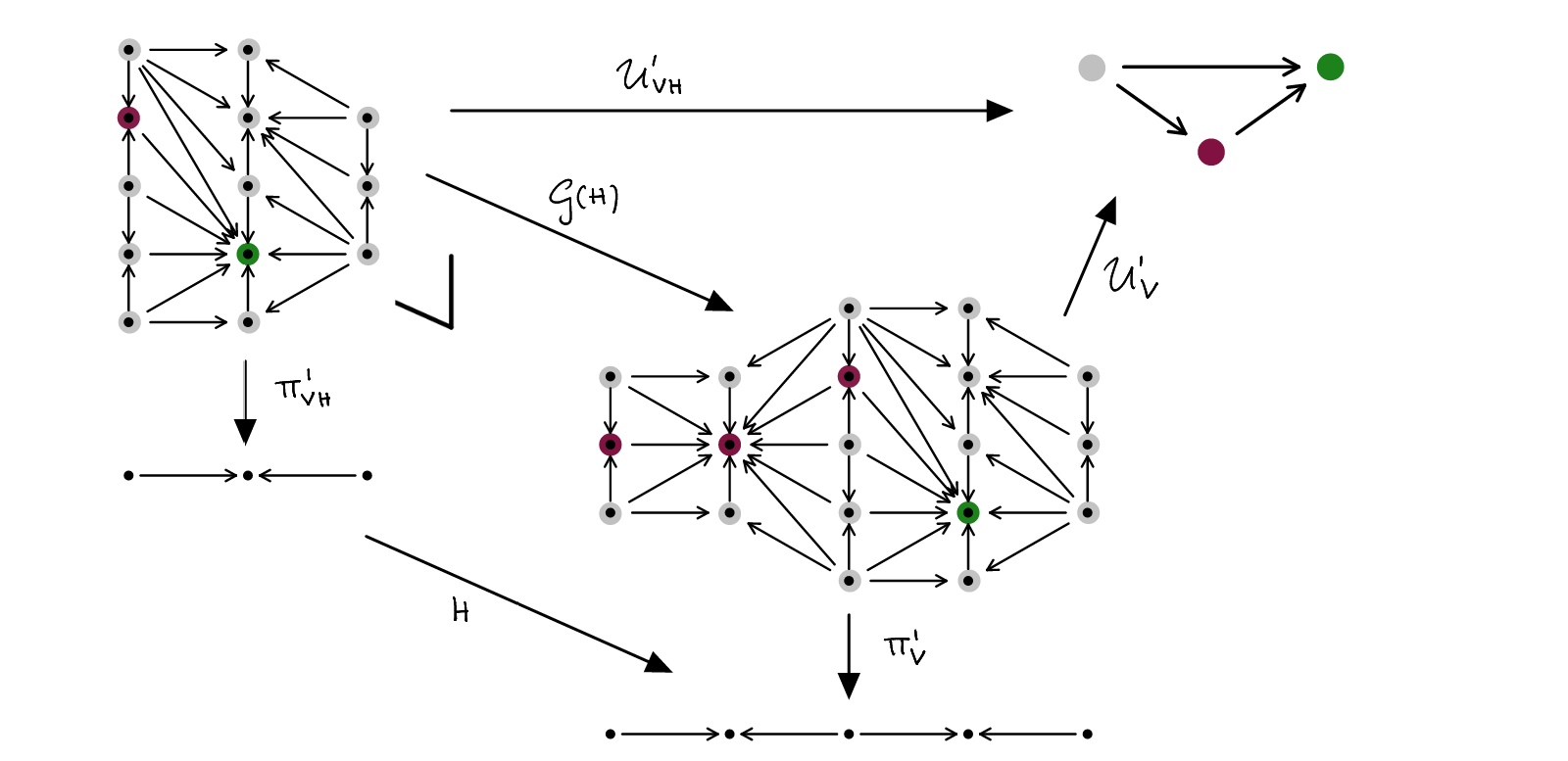}
\endgroup\end{restoretext}
Next consider the map $K : \singint 2 \to \singint 1$ mapping $i$ to $\max(i-2,0)$. Then $VHK$ is obtained as follows
\begin{restoretext}
\begingroup\sbox0{\includegraphics{test/page1.png}}\includegraphics[clip,trim=0 {.0\ht0} 0 {.0\ht0} ,width=\textwidth]{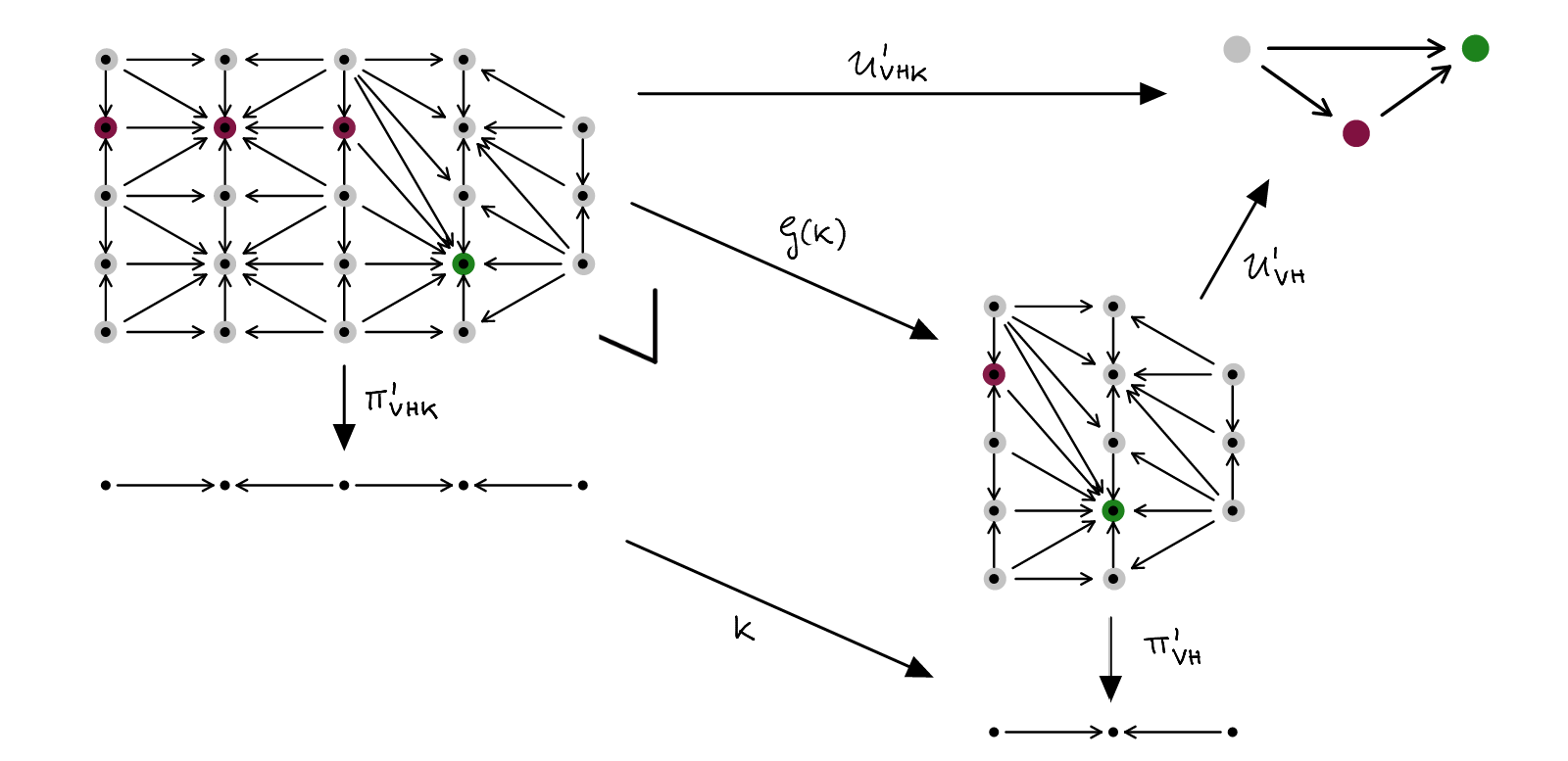}
\endgroup\end{restoretext}

\subsection{Identities and projections}

We introduce identities of singular $n$-cubes and projections of points (in the total poset). Geometrically, identities are labelled $(n+1)$-cubes obtained as trivial bundles of labelled $n$-cubes. Projections of a point in the $n$-cubes are the images of that point under (iterated) projections from the $(k+1)$-cube to the $k$-cube.

\begin{defn}[Identity bundles] Let $\scA$ be a $\cC$-labelled singular $n$-cube. We define the $\cC$-labelled singular $(n+1)$-cube $\Id_\scA$ by setting (cf. \autoref{notn:sum_und_drop_notn})
\begin{equation}
\tsU 1_{\Id_\scA} = \tsU 0_{\scA}
\end{equation}
and 
\begin{equation}
\tusU 0_{\Id_\scA} = \const_{\singint 0}
\end{equation}
where $\const_{\singint 0} : \bnum{1} \to \SI$ is the constant functor to $\singint 0 \in \SI$. We inductively set $\Id^0_\scA = \scA$ and $\Id^k_\scA = \Id_{\Id^{k-1}_\scA}$ for $k > 0$. 
\end{defn}

\begin{defn}[Projections] \label{defn:sum_proj} Given a $\cC$-labelled singular $n$-cube family $\scA$ indexed by a poset $X$, then for an object $p \in \tsG n(\scA)$ we define $p^n = p$ and further set for $0 \leq k < n$
\begin{equation}
p^k := \pi_{\tusU k_\scA} (p^{k+1}) \in \tsG k(\scA)
\end{equation}
$p^k$ is called the $k$-level projection of $p$.
\end{defn}

As an example, the identity bundle $\Id_Q$ is given by the data
\begin{restoretext}
\begingroup\sbox0{\includegraphics{test/page1.png}}\includegraphics[clip,trim=0 {.0\ht0} 0 {.0\ht0} ,width=\textwidth]{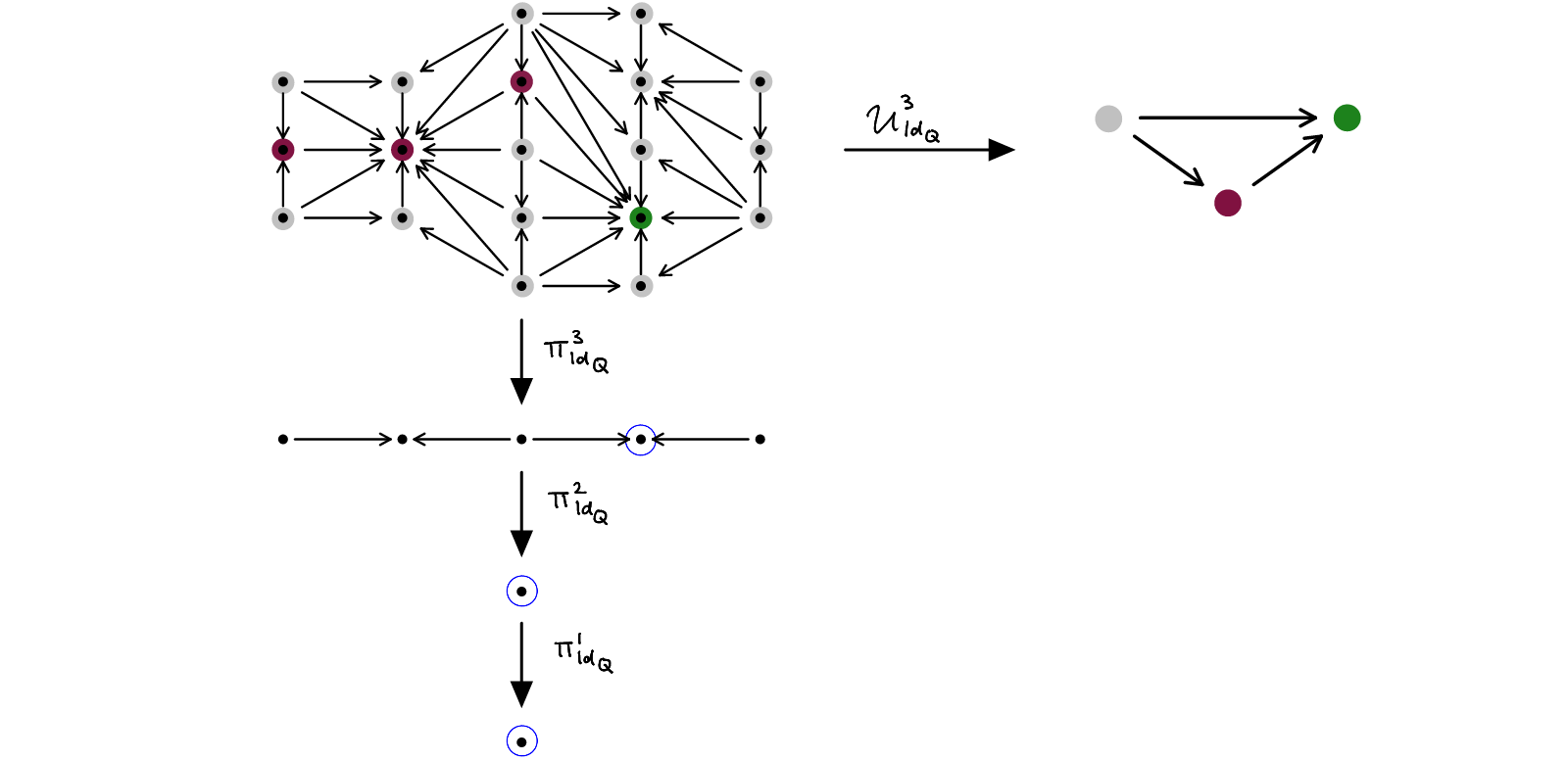}
\endgroup\end{restoretext}
Further, the projections of the \cdarkgreen{} point in $\tsG 3(\Id_Q)$ are marked by \cblue{} circles.

\section{Categorical structures of cubes} \label{sec:sum_norm}

\subsection{The classifying category of cube bundles: $\SIvert n \cC$} \label{ssec:SIvert}

An important and not immediately obvious feature of our $\cC$-labelled singular $n$-cube families over $X$ is that they have a classifying category $\SIvert n \cC$. That is, an $X$-indexed $\SIvert n \cC$-cube family corresponds to a functor $X \to \SIvert n \cC$. 

Let $i_{\BDelta} : \BDelta \into \Cat$ be the simplex category (resp. its inclusion into $\Cat$), the full subcategory of $\Cat$ of non-empty finite ordinals $\bnum 1, \bnum 2, ...$ . Recall that the simplicial nerve $N$ of a category is defined to map
\begin{align}
N : \Cat \to [\BDelta\op, \SetCat] \\
\cC \mapsto \Cat(i_{\BDelta} , \cC)
\end{align}
Note that $N$ is fully faithful and injective on objects. 

The classifying category of labelled $n$-cubes is then defined as follows. 
\begin{defn}[The classifying category of cube bundles] We define $\SIvert n \cC$ by defining $N (\SIvert n \cC)$ as follows: $N (\SIvert n \cC)(\bnum k)$ is the set $\bnum k$-indexed $\SIvert n \cC$-cube families. If $f : \bnum k \to \bnum m$ is a morphism in $\BDelta$ then $N (\SIvert n \cC)(\bnum k)(f)$ acts by pullback along that morphism (cf. \autoref{def:sum:pullbacks_of_families}).
\end{defn}

Importantly, $\SIvert n \cC$ has another construction as well: we will prove that the above coincides with an inductive ``over-double-category" construction, giving rise to an endofunctor
\begin{equation}
\SIvertone {-} : \Cat \to \Cat
\end{equation}
mapping $\cC$ to a category $\SIvertone \cC$ whose objects are functors $F : I \to \cC$ for some $I \in \SI$, and morphisms are natural transformations $\alpha : (\SiR(f) : I_1 \xslashedrightarrow{} I_2) \to \Hom_\cC(F_1 -, F_2-)$ for $f \in \mor(\SI)$. Composition is given by the canonical composition of squares
\begin{equation}
\xymatrix@C=1.5cm@R=0.8cm{  I_1 \ar[r]^-{\SiR(f_1)}|-*=0@{|} \ar[d]_{F_1} \ar@{}[dr]|{ \Downarrow \alpha_1} &  I_2 \ar[d]^{F_2} \ar[r]^-{\SiR(f_2)}|-*=0@{|}  \ar@{}[dr]|{ \Downarrow \alpha_2} & I_3 \ar[d]^{F_3} \\
\cC \ar[r]_{\Hom_{\cC}}|-*=0@{|} & \cC  \ar[r]_{\Hom_{\cC}}|-*=0@{|} & \cC }
\end{equation}
where the top two morphisms compose as relations, the bottom two as profunctors. We then define
\begin{defn}[The classifying category of cube bundles] $\SIvert n \cC$ is the $n$-fold application of the endofunctor $\SIvertone{-}$ to $\cC\in \Cat$.
\end{defn}
To verify that these two definitions coincide, we will show that there is a construction taking the data of a $\cC$-labelled singular $n$-cube family $\scA$ over $X$, and outputting a functor into $\SIvert n \cC$ (in the sense of the second definition) denoted by
\begin{equation}
\sR_{\sT_\scA, \sU^n_\scA} : X \to \SIvert n \cC
\end{equation}
where
\begin{equation}
\sT_\scA = \Set{\und\sU^0_\scA, \und\sU^1_\scA, ... ,\und\sU^{n-1}_\scA}
\end{equation}
This construction then gives a bijective correspondence between $\SIvert n \cC$-cube families over $X$ and functors $X \to \SIvert n \cC$. We note that the correspondence also behaves nicely for pullbacks as we will see.

\subsection{The category of cube bundles $\Bunbc^n_\cC$} \label{sec:sum_Buno}

We define the category $\Bunbc^n_\cC(\bnum 1)$ of ($\cC$-labelled $n$-)cubes and cube maps, which is a combinatorial model of the category $\Cubeo n \cC(\bnum 1)$ (minus the openness condition on maps) as defined in \autoref{ssec:po_mfld_diag}. The openness condition will be combinatorially formulated only in the next section. Note that further generality is added by considering \textit{families} of cubes indexed by some (stratified) base space $X$ yielding a category $\Bunbc^n_\cC (X)$---this also has a topological analogue which was mentioned in \autoref{ssec:po_mfld_diag}. Letting $X$ vary we will first define a category $\Bunbc^n_\cC$.

\begin{defn}[Category of cube bundles and multi-level base changes] The category $\Bunbc^n_\cC$ has as objects $\cC$-labelled singular $n$-cube families $\scA$, and a morphism $M : \scA \to \scB$ consists of a list of functors (for $0 \leq k \leq n$)
\begin{equation}
M^k : \tsG k(\scA) \to \tsG k(\scB)
\end{equation}
such that firstly
\begin{equation} \label{eq:sum_multi_lvl_bch_tri}
\xymatrix{ & \cC & \\
\tsG n(\scA) \ar[ur]^{\tsU n_\scA} \ar[rr]_{M^n} && \tsG n(\scB) \ar[ul]_{\tsU n_\scB} }
\end{equation}
commutes and (for $0 \leq k < n$) so does
\begin{equation} \label{eq:sum_multi_lvl_bch_sq}
\xymatrix{ \tsG {k+1}(\scA) \ar[r]^{M^{k+1}} \ar[d]_{\tpi {k+1}_\scA} & \tsG {k+1}(\scB) \ar[d]^{\tpi {k+1}_\scB} \\
\tsG k(\scA) \ar[r]^{M^k} & \tsG k(\scB)  }
\end{equation}
Secondly, we require that all $M^{k+1}$ are fibrewise monotone, meaning that for each $M^k(x) = y$ we have
\begin{equation}
M^{k+1} : (\tpi {k+1}_\scA)\inv(x) \to(\tpi {k+1}_\scB)\inv(y)
\end{equation}
is monotone. Morphisms are also called \textit{multi-level base changes}. The full subcategory of cube bundles over a fixed poset $\tsG 0 (A) = X$ is denoted by $\Bunbc^n_\cC (X)$.

Composition is given level-wise, that is, for $M : \scA \to \scB$ and $L : \scB \to \scC$ we have $LM : \scA \to \scC$ which is given by
\begin{equation}
(LM)^k = L^k M^k
\end{equation}
\end{defn}

\begin{constr}[Decomposition into $k$-level base changes] \label{constr:multilvl_decomp} A morphism $M_i : \scA \to \scB$ in $\Bunbc^n_\cC$ is called an $i$-level base change if 
\begin{enumerate}
\item for $k \geq i$, \eqref{eq:sum_multi_lvl_bch_sq} are pullbacks 
\item for $k < i$, $M^k$ are identities
\end{enumerate}
Any multi-level base change $M : \scA \to \scB$ has a unique decomposition as 
\begin{equation}
M = M_n M_{n-1}... M_1 M_0
\end{equation}
where $M_i$ is an $i$-level base change. This can be shown inductively: write $M = M' M_0$ (with $(M')^0 = \id$) using the universal property of pullback, and then recognise $M'$ as a morphism in $\Bunbc^{n-1}_\cC$ by forgetting level $0$, and use the inductive assumption.
\end{constr}

The importance of the preceding construction lies in the fact that most proofs about multi-level base changes can be reduced to proofs about $k$-level base changes: consequently, for the most part of the technical development in the next chapters we will actually be working with $k$-level base changes.

\subsection{The category of cube bundles and \textit{open} maps $\Buno n \cC$} \label{sec:sum_Buno2}

More precisely, we will be working with two special instances of multilevel base changes (and their $k$-level decompositions) called \textit{collapse} and \textit{embedding}. These are epis and monos in a special subcategory $\Buno n \cC$ of $\Bunbc^{n}_\cC$. The relationship of $\Bunbc^n_\cC$ and $\Buno n \cC$ can be understood (heuristically) in geometric terms as follows: $\Bunbc^n_\cC$ encodes maps of stratified cubes bundles (up to deformation) including possibly degenerate maps, where as $\Buno n \cC$ is the full sub-category of \textit{open} maps. We introduce the following terminology

\begin{defn}[Open, collapse and embedding functors] A functor of posets $F : \singint k \to \singint l$ is called
\begin{enumerate}
\item \textit{open} if it preserves regular segments (that is, on objects it maps even numbers to even numbers)
\item \textit{a collapse functor} if it is open, and also surjective on objects
\item \textit{an embedding functor} if it is open, and also injective on objects
\end{enumerate}
\end{defn}

We now define

\begin{defn}[Category of cube bundles and open multi-level base changes]
Let $\Buno n \cC$ be the wide subcategory of $\Bunbc^n_\cC$ whose morphisms $M : \scA \to \scB$ satisfy
\begin{enumerate}
\item For $k = 0$, $M^k$ is injective.
\item For $k > 0$, we have that $M^k$ is fibrewise open, meaning that for each $M^k(x) = y$ we require 
\begin{equation}
M^{k} : (\tpi {k}_\scA)\inv(x) \to(\tpi {k}_\scB)\inv(y)
\end{equation}
is open. 
\end{enumerate}
If $M^k$, $k > 0$, further is fibrewise a collapse functor, then $M$ is called a \textit{(multi-level) collapse}. Equally, if it further is fibrewise an embedding functor, then $M$ is called a \textit{(multi-level) embedding}.
\end{defn}

The injectivity condition can be understood as requiring ``non-degeneracy" in the base. Examples of collapses and embeddings are at the end of the section.

We will not need the following observation, but mention it for completeness.

\begin{thm}[Factorisation system in $\Buno n \cC$] $\Buno n \cC$ has a (epi,mono)-factorisation system and epis are exactly collapses, while monos are exactly embeddings.
\proof This will follow from \autoref{thm:restriction_of_collapse} (see also \autoref{rmk:epi_mono_pf}). \qed
\end{thm}

\begin{notn}[Symbols for later computations] In the later chapter, $k$-level collapses $\scA \to \scB$ will be characterised by certain natural injections $\lambda$. The corresponding $k$-level collapse has a map at level $k$ given by $\cS^\lambda : \sG^k(\scA) \mcoll \sG^k(\scB)$ (we also write $\lambda : \scA \kmcoll k \scB$). Multi-level collapses will be characterised by sequences $\vvec\lambda$ of natural injections. The corresponding multi-level collapse is written as $\vvec\cS^{\vvec\lambda} : \scA \mcoll \scB$ (we also write $\vvec\lambda : \scA \kmcoll * \scB$).

Embeddings will often be written as Greek letters $\theta : \scA \mono \scB$, and their $k$-level counterparts will be sometimes written as $\sJ^{\scA,k}\restsec{[q_-,q_+]} : \scA^k\restsec{[q_-,q_+]} \mono \scA$ as they are characterised by pairs of ``endpoint sections" $q_-, q_+$ (note that we also keep track of $\scA, k$ in this notation).
\end{notn}

Embeddings give us a notion of \textit{subfamilies}. That is, $S : \scA \mono \scB$ should be understood as $\scA$ being a (singular $n$-cube) subfamily of $\scB$. The following is an example of an embedding $\theta : M \mono Q$ which is given by the data
\begin{restoretext}
\begingroup\sbox0{\includegraphics{test/page1.png}}\includegraphics[clip,trim=0 {.0\ht0} 0 {.0\ht0} ,width=\textwidth]{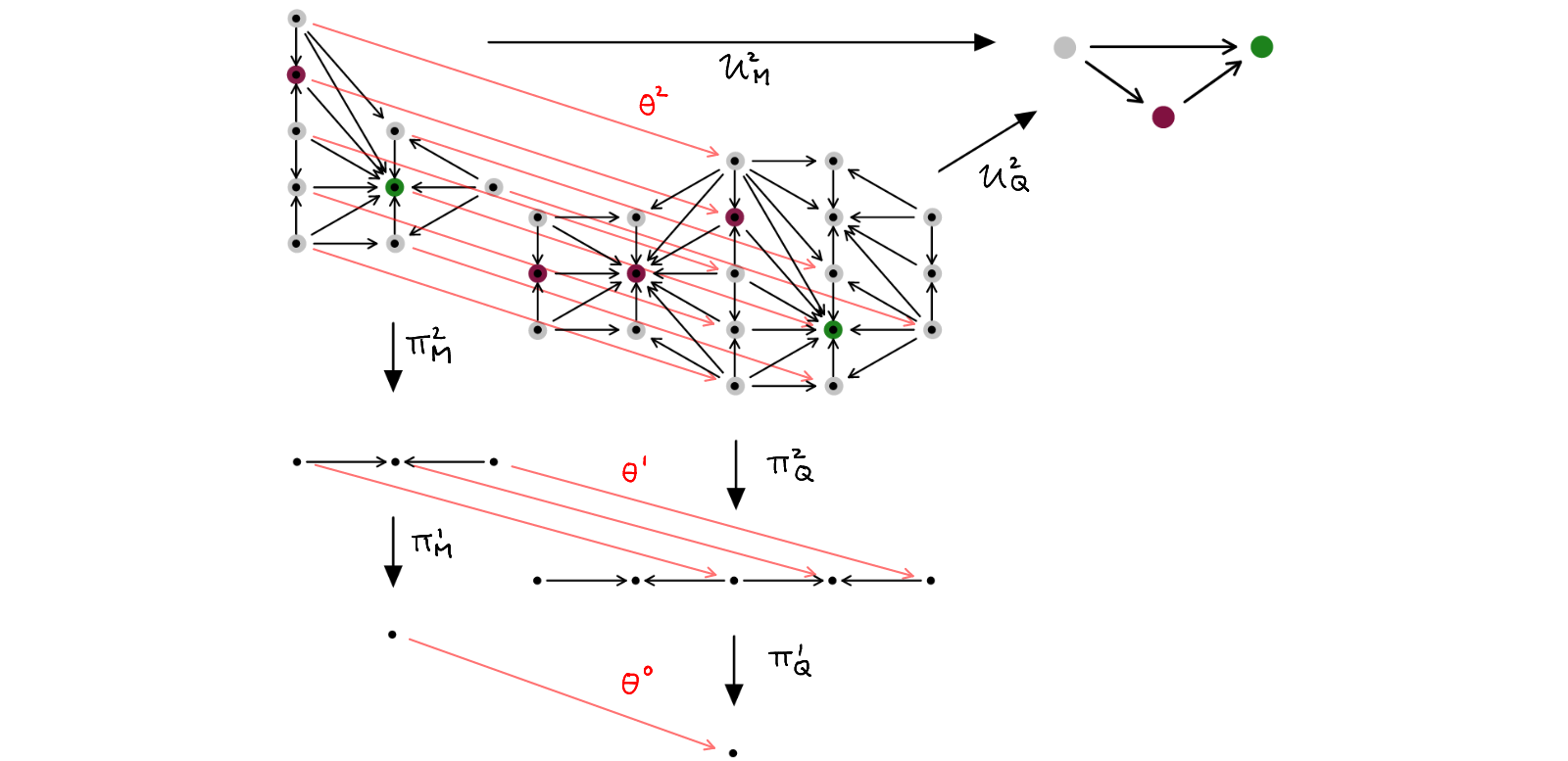}
\endgroup\end{restoretext}
Here, the level-wise maps $\theta^k$ are indicated by \cred{} arrows. Note that $M$ and $Q$ are cubes in this case and thus we call $M$ a subcube of $Q$.

On the other hand, collapses give us a notion of ``quotient families". That is, $S : \scA \mcoll \scB$ should be understood as $\scB$ being the cube $\scA$ but with singularities ``quotiented away". An example of a collapse $\lambda : Q \mcoll N$ is given by the following data
\begin{restoretext}
\begingroup\sbox0{\includegraphics{test/page1.png}}\includegraphics[clip,trim=0 {.0\ht0} 0 {.0\ht0} ,width=\textwidth]{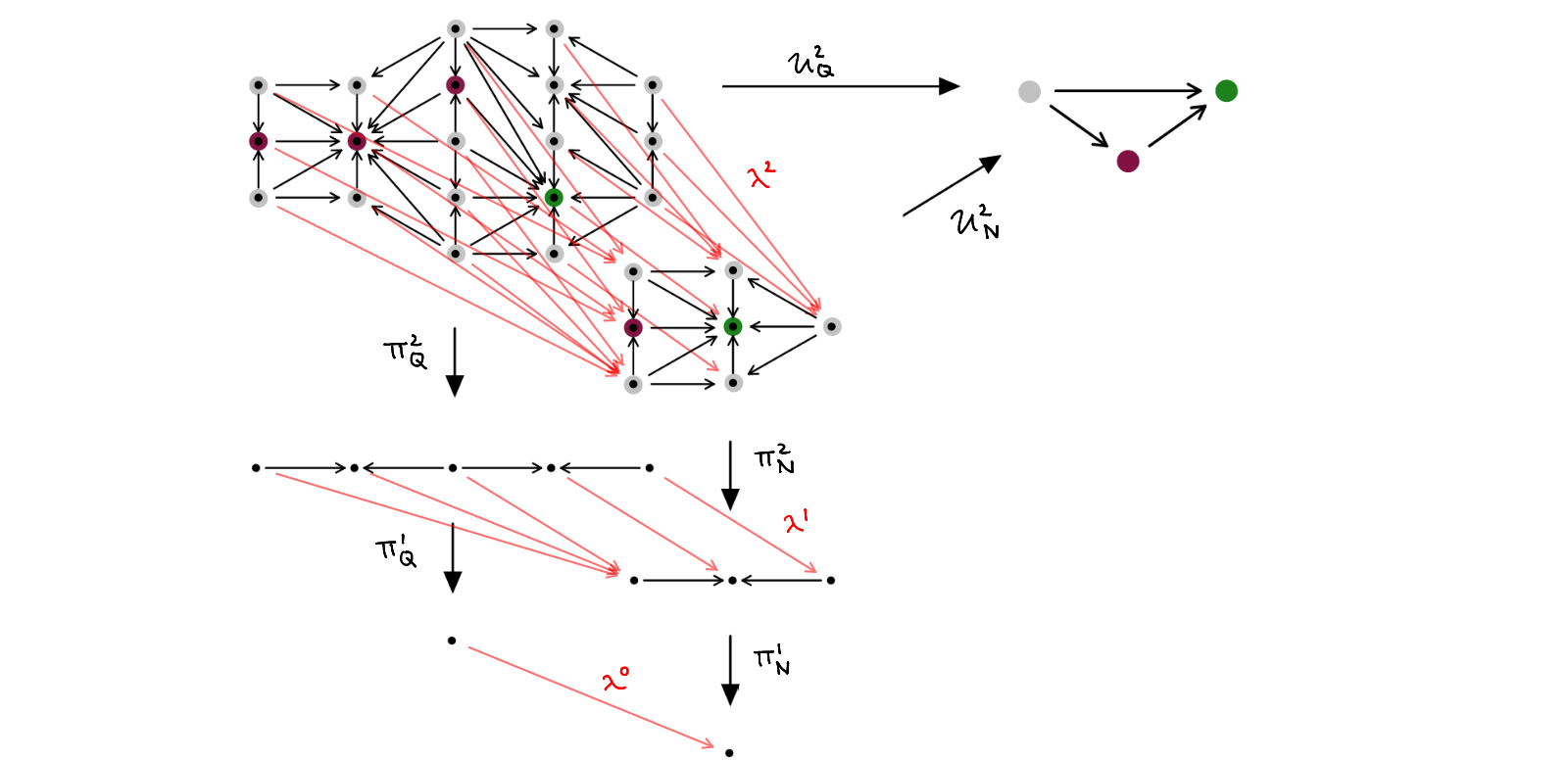}
\endgroup\end{restoretext}
The level-wise maps $\lambda^k$ are indicated by \cred{} arrows as before.

Observe that this quotienting process retains all the ``geometric information" of a cube, that is, the corresponding manifold diagram (up to equivalence) stays the same. In other words, collapse only quotients away ``redundant" singular heights. Collapsing a cube then eventually leads to its unique \textit{normal form} which is the cube with the same geometric information but only non-redundant singular heights (in other words, with the unique coarsest stratification). In this sense, collapse gives a confluent rewriting system on $n$-cubes. Given a collapse $S : \scA \mcoll \scB$ we call $\scA$ an \textit{expansion} of $\scB$, and $\scB$ a \textit{reduct} of $\scA$.

\subsection{Normal form theorem: Terminal objects in ($\Buno n \cC$, epi)} \label{sec:sum_NF}

We state the existence and uniqueness of normal forms. We will show the following.

\begin{thm}[Normal forms] Let $\scA$ be a $\cC$-labelled singular $n$-cube. Then there is a $\cC$-labelled singular $n$-cube $\NF{\scA}^n$ such that for any $S : \scA \mcoll \scB$ we have unique $S' : \scB \mcoll \NF{\scA}^n$. $\NF{\scA}^n$ is called the normal form of $\scA$ (up to level $n$)\footnote{The terminology ``up to level $n$" is not redundant even though $\scA$ is assumed to be an $n$-cube. As we will later on see every singular $n$-cube is in particular a singular $k$-cube, and thus can be normalised ``up to level $k$" as well.}.
\end{thm}

In particular, applying the previous theorem for $\id : A \mcoll A$ we find a collapse $\nfc^n_\scA : \scA \mcoll \NF{\scA}^n$. A cube in normal form is called \textit{normalised}.

\begin{rmk}[Terminal objects in ($\Buno n \cC$, epi)] Let ($\Buno n \cC$, epi) denote the wide subcategory of $\Buno n \cC$ generated by epimorphisms (that is, by collapses). The proof of the preceding theorem will then show that each connected component of ($\Buno n \cC$, epi) has a terminal object (namely, the normal form of any object in the component). 
\end{rmk}

Compare the following to \autoref{defn:ref_and_eq}.

\begin{defn}[Geometric equivalence] \label{defn:sum_geom_equiv} Given $\cC$-labelled singular $n$-cubes $A, B$ we say $A$ and $B$ are \textit{geometrically equivalent}, written $A \simeq B$, if there is a cospan $A \epi C \twoheadleftarrow B$, in other words, if they have a mutual collapse. The previous theorem guarantees that each equivalence class has a unique maximally collapsed representative, which makes equivalence decidable by computing normal forms.
\end{defn}

As an example the normal form of the singular $2$-cube $Q$ is given by the singular $2$-cube $N$ given above. As the reader can verify, no non-identity collapses apply to $N$. Using the discussion in \autoref{ssec:coloring} we find the (labelled) directed triangulation corresponding to $N$ to be
\begin{restoretext}
\begingroup\sbox0{\includegraphics{test/page1.png}}\includegraphics[clip,trim=0 {.3\ht0} 0 {.35\ht0} ,width=\textwidth]{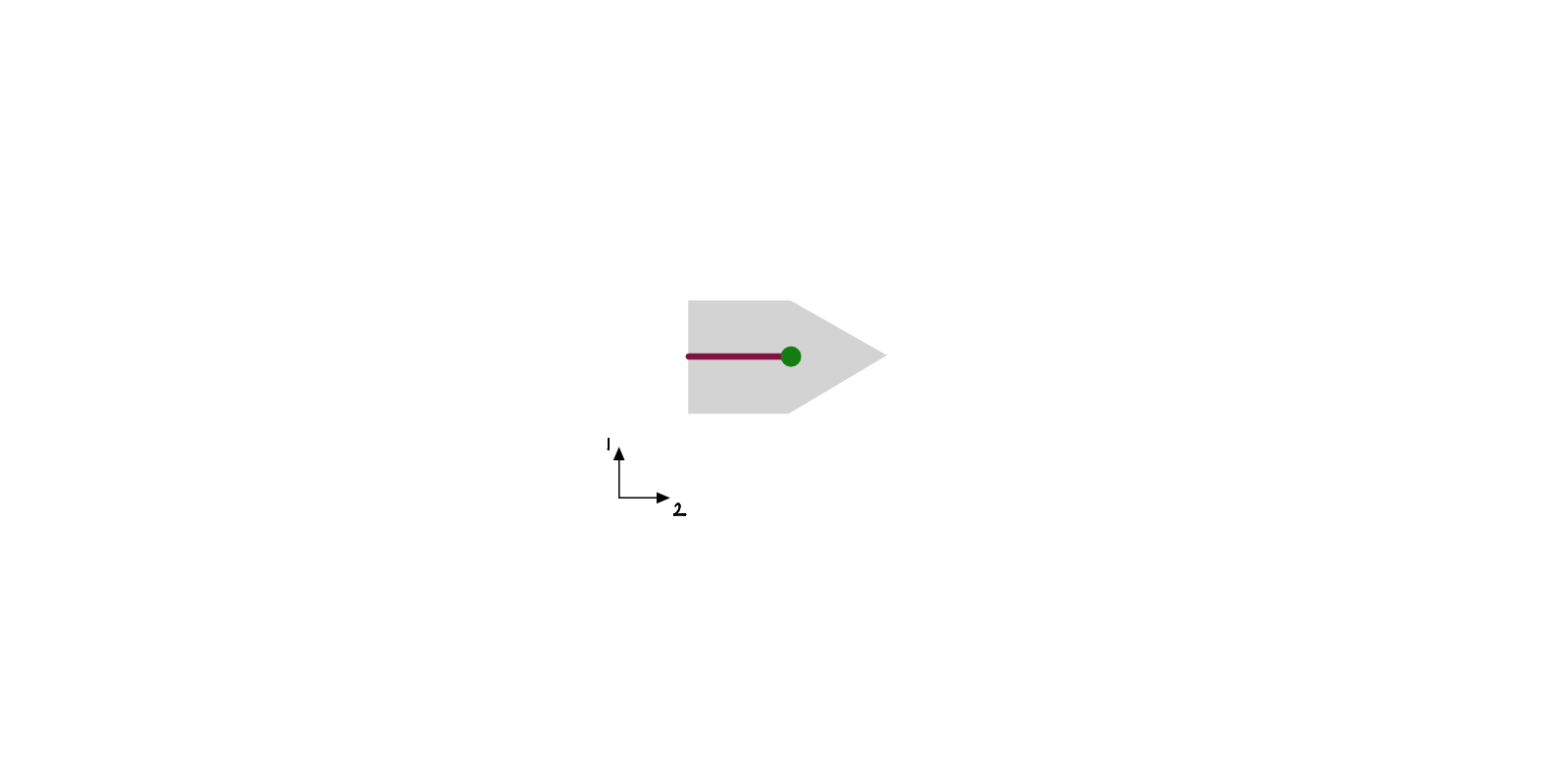}
\endgroup\end{restoretext}
This should be compared to the definition of $N$. Following the discussion in \autoref{ssec:coloring} we can see that the above also gives the manifold diagram corresponding to $Q$, but that $N$, unlike $Q$, provides the coarsest triangulation of the manifold diagram.

\subsection{Commutation theorem: (epi,mono)-factoisation system} \label{sec:sum_comm}

Given a collapse $S : \scA \mcoll \widetilde\scA$ of a cube $\scA$ which has a sub-cube $T : \scB \mono \scA$ then the former induces a collapse on its sub-cube. We will show the following.

\begin{thm}[Commutation of collapse and embedding] \label{thm:restriction_of_collapse} Given $S : \scA \mcoll \widetilde\scA$, $T : \scB \mono \scA$, then we have $T\pbstar  S : \scB \mcoll \widetilde\scB$ and $S\postar  T : \widetilde \scB \mono \widetilde \scA$, determined by the following being pullbacks for all $0 \leq k \leq n$
\begin{equation}
\xymatrix{\tsG k(\scB) \ar[r]^{T^k} \ar[d]_{(T\pbstar  S)^k} \pullbackfar & \tsG k(\scA) \ar[d]^{S^k} \\
\tsG k(\widetilde \scB) \ar[r]_{(S\postar  T)^k} & \tsG k(\widetilde \scA) }
\end{equation}
\end{thm}
\noindent Later on we will also talk about the \textit{restriction} of a collapse to a subcube. 

\begin{rmk}[(epi,mono)-factorisation system] \label{rmk:epi_mono_pf} As an immediate consequence of the above, since every $k$-level basechange can be written as an epi followed by a mono (or simply a mono for $k = 0$), we deduce using \autoref{constr:multilvl_decomp} that every morphism $\Buno n \cC$ can be written as an epi following by a mono. It is easy to see that this decomposition is necessarily unique and thus $\Buno n \cC$ has an (epi,mono)-factorisation system.
\end{rmk}

In previous examples we have defined $\theta : M \mono Q$ and $\lambda : Q \mcoll N$. In this case we have that $\lambda\postar \theta : N \mono N$ is in fact the identity sub-cube, and $\theta\pbstar \lambda : M \mcoll N$ is given by the collapse
\begin{restoretext}
\begingroup\sbox0{\includegraphics{test/page1.png}}\includegraphics[clip,trim=0 {.0\ht0} 0 {.0\ht0} ,width=\textwidth]{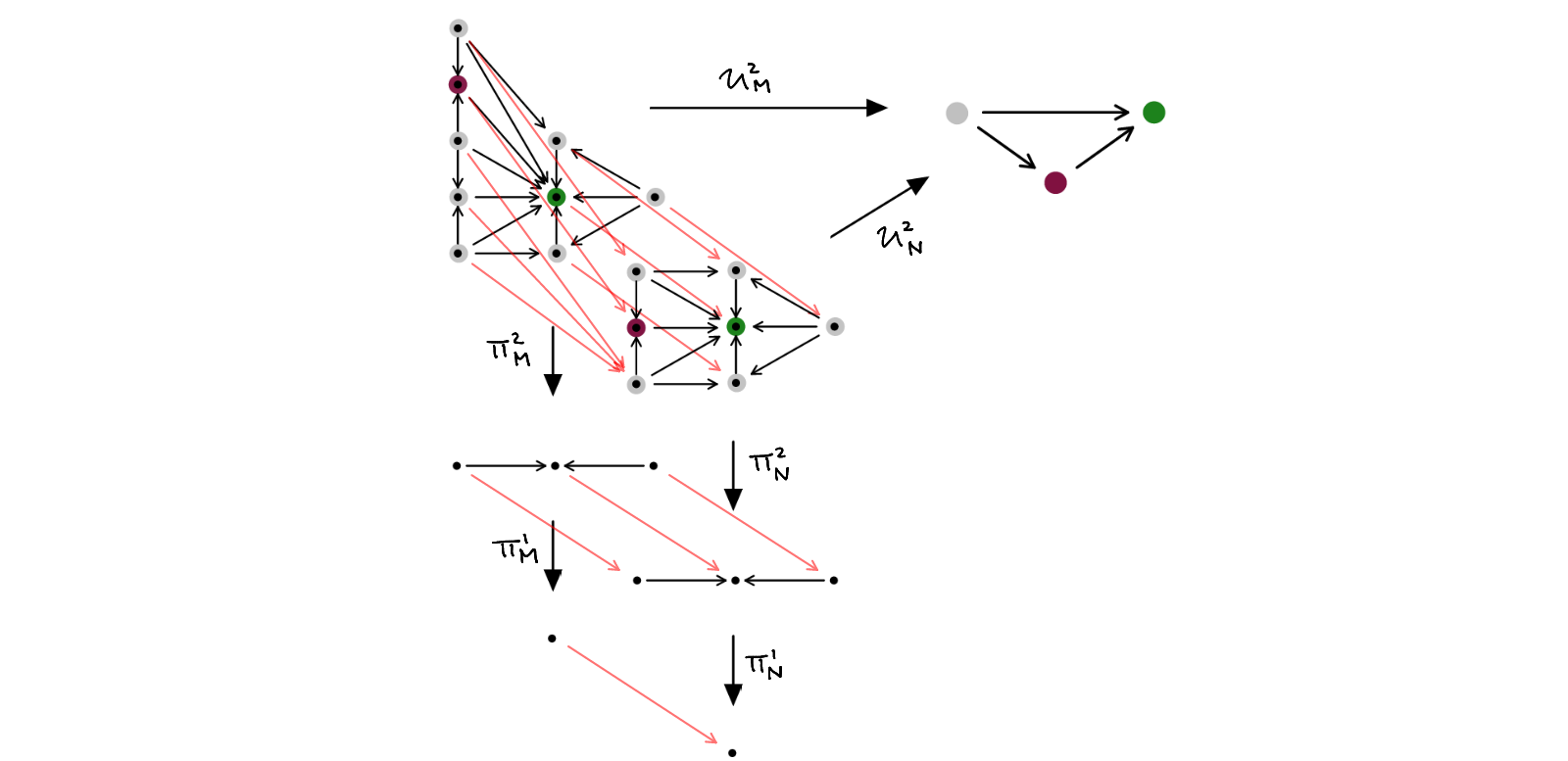}
\endgroup\end{restoretext}

\subsection{The minimal subcube construction} \label{sec:sum_min}

There is also a minimal (non-trivial) object among all subcubes which we record in this section. In terms of directed triangulations, the minimal subcube around a point of the triangulation, should be thought of as the subtriangulation bounded by the link of that point. We remark that the statement below uses the observation that given $x \in X$ for a poset $X$, then the overposet $X \sslash x$ is naturally a subposet of $X$.

\begin{thm} Let $\scA \in \SIvert n \cC$, and $p \in \tsG n(\scA)$. Then there is
\begin{enumerate}
\item $\scA \in \SIvert n \cC$, called minimal (open) neighbourhood of $p$
\item $(\iota^p_\scA) : (\scA \sslash p) \mono \scA$, called minimal neighbourhood embedding of $p$
\end{enumerate}
This satisfies
\begin{equation}
(\tsG n(\scA) \sslash p) \subset \im((\iota^p_\scA)^n)
\end{equation}
together with the following universal property: if $T : \scB \mono \scA$ is such that $\left(\tsG n(\scA) \sslash p\right) \subset \im(T^n)$ then there is a unique $F : (\scA \sslash p) \mono \scB$ such that $TF = (\iota^p_\scA)$.
\end{thm}

\noindent The reader can verify that our previous example  $\theta : M \mono Q$ is in fact also an example of a minimal sub-cube.

\section{Globes and cones} \label{sec:sum_glob}

So far we have allowed our labelled subsets of the cube to ``roam freely" in the $n$-cube. However, as discussed in the \autoref{sec:intro_nfold} we would like to impose constancy conditions on some of the sides of the cube, such that the $n$-cube factors through an $n$-globe. Technically, we will impose this condition not only on the side of a given $n$-cube, but also on all sides of possible embedded subcubes of it. This condition is then called ``globularity" and will be formulated in the following. Globularity will also allow for interesting modes of composition of cubes, which we call ``whiskering".

\begin{notn}[Test functors] \label{notn:sum_test_fctr} Let $X$ be a poset. Given $x \in X$ define $\Delta_x :  \bnum{1} \to X$ to be the functor mapping $0$ to $x$. Given $(x\to y) \in \mor(X)$, denote by $\Delta_{x \to y} : \bnum{2} \to X$ the functor mapping $0$ to $x$ and $1$ to $y$.
\end{notn}

\subsection{Globes}

For an $X$-labelled $\SI$-family $\scB$, denote by $\regcont(\scB) \subset \sG(\scB)$ the full subposet having objects $p = (x,a)$ where $a$ is a regular segment (that is, $a$ is even).

\begin{defn}[Globular cubes] A $\cC$-labelled singular $n$-cube family $\scA$ is globular if for every $1 \leq k \leq n$, $\tsU k_\scA$ normalises to a constant functor when restricted to any connected component of $\regcont(\tusU {k-1}_\scA)$.

A globular $\cC$-labelled singular $n$-cube is also called a $\cC$-labelled singular \textit{$n$-globe}.
\end{defn}

All previous examples (for instance the cubes $M, N, P, Q$) are in fact globes. A non-example can be produced for instance from $N$, by replacing the \cgray{} color of one of two (horizontally) centred points by the color \cdarkgreen{}. We would then find that $\tsU 2_N$ is non-constant on $\regcont(\tusU 1_N)$ and thus globularity is not satisfied (note that for $0$-cube families the constancy-after-normalisation condition simplifies to just requiring constancy).

Note further that all $(-1)$- and $0$-cubes are automatically globular.

\begin{defn}[Source and target of globes] Let $\scA$ be a $\cC$-labelled singular $n$-globe for $n \geq -1$. For $n \geq 1$. Define $k$ to satisfy $\tusU 0_\scA(0) = \singint k$. Then we define $\gsrc(\scA)$, called \textit{globular source} of $\scA$, to be the $(n-1)$-globe given by (cf. \autoref{notn:sum_test_fctr})
\begin{equation}
{\gsrc(\scA)} := \tsU 1_\scA \Delta_0
\end{equation}
Similarly, we define the \textit{globular target} $\gtgt(\scA)$ of $\scA$ by
\begin{equation}
{\gtgt(\scA)} := \tsU 1_\scA \Delta_{2k}
\end{equation}
For $n = 0$ define the globular source and target of $\scA$ to be $\emptyset$ (recall \autoref{conv:minus_one_cubes}). For $n = -1$, define globular source and target of $A = \emptyset$ to again equal $\emptyset$.

Note that the notation $\gsrc$ will later on distinguish \textit{globular} sources from \textit{cubical} sources ($\csrc$) and similarly for targets.
\end{defn}

\subsection{Double cones}

For the theory of \textit{\free{}} associative $n$-categories we need the following coning construction to concisely describe minimal neighbourhood of singularities (which will be of ``conical" shape). Geometrically speaking, this construction takes a stratification of the $(n-1)$-sphere $S^{n-1}$ and produces a stratification of the $n$-disk $D^n$, by taking the cone of $S^{n-1}$ and labelling the vertex point by a given label. Note that since we work with $n$-globes, $S^{n-1}$ will be cut into two disks: a source $S$ and a target $T$, both stratifications of the $(n-1)$-disk coinciding on their boundary.

\begin{defn}[Double cones] Let $n \geq 0$, and $S$ and $T$ be globular normalised $\SIvert {n-1} {\cC}$-cubes which coincide on their globular sources and targets. Let $g \in \cC$ be an object such that for each $s \in \im(\tsU {n-1}_S)$ and $t \in \im(\tsU {n-1}_T)$ there are unique morphisms $s \to g$ and $t \to g$. We define the \textit{double cone} $\abss{S \xto g T}$ to be the unique globular normalised $\SIvert n \cC$-cube determined by the properties
\begin{enumerate}
\item \textit{Source and target}: $\gsrc\abss{S \xto g T} = S$ and $\gtgt\abss{S \xto g T} = T$
\item \textit{Minimality} There is a unique $p_g \in \tsG n(\abss{S \xto g T})$ such that $\tsU n_{\abss{S \xto g T}}(p_g) = g$ and
\begin{equation}
\abss{S \xto g T} \sslash p_g = \abss{S \xto g T}
\end{equation}
\item \textit{Correctness of dimension}: $p^k_g = (p^{k-1}_g,1)$ for all $1 \leq k \leq n$ (for this condition recall notation for elements in total spaces as tuples from \autoref{ssec:sum_fam_bun_prel}. The condition essentially requires that the neighbourhood of $p_g$ is an open set in $n$-dimensional space and not a lower-dimensional projection of such.)
\end{enumerate}
Note that for $n = 0$ we have $\abss{\emptyset \xto g \emptyset} = \Delta_g$.
\end{defn}

Morally, double cones form an $n$-cube from two $(n-1)$-cubes (that agree on their boundary) by contracting their content into a single central point $p_g$ labelled by $g$.

\subsection{Terminal globes} \label{ssec:sum_term_glob}

For the theory of associative $n$-categories, which will be formulated as \textit{globular sets} with operations, we will need an $n$-cube that models elements in globular sets (including their iterated sources and targets). There is a natural choice for this, and it will be called a \textit{terminal globular $n$-cube}, or terminal $n$-globe for short.

To define this we first construct a ``generic" (terminal) $n$-globe. All other globes corresponding to elements of globular sets will be obtained from relabelling this generic globe.  

\begin{constr}[Terminal $n$-globes] Let $\lG$ denote the globe category. This has object set $\lN$ and generating morphisms $\sigma_{k-1}, \tau_{k-1} : (k-1) \to k$ subject to the relations (for $\rho_{k} \in \Set{\sigma_{k}, \tau_{k}}$)
\begin{equation}
\rho_k \sigma_{k-1} = \rho_k \tau_{k-1}
\end{equation}
Let $\globset = [\lG\op,\SetCat]$ be the category of globular sets (see e.g.  \cite{maltsiniotis2010grothendieck}). For $S \in \globset$ we sometimes write $S_k = S(k)$.

Given $S \in \globset$ we can construct its category of elements $\elcat(S)$. This has as objects all elements of $S$, and arrows go from elements to their sources and targets. Formally, using our bundle construction $\sG$ and the ``discrete" functor $\Discr :  \SetCat \to \PRel$, this can be defined by
\begin{equation}
\elcat(S) := \sG(\Discr S)
\end{equation}
Here, $\Discr :  \SetCat \to \PRel$ maps a set $X$ to the discrete poset $X$, and a function of sets $f : X \to Y$, to the profunctorial relation $\Hom_{\Discr Y} (f - , - ) : X \xslashedrightarrow{} Y$. 

As an example consider the globular set $S_{(3)}$ determined by having a single element in $S_{(3)}(0)$, $S_{(3)}(1)$ and $S_{(3)}(2)$ but no other elements. We call its elements $a$,$c$ and $h$ respectively. The opposite of its category of elements, $\elcat(S_{(3)})\op$, is then given by
\begin{restoretext}
\begingroup\sbox0{\includegraphics{test/page1.png}}\includegraphics[clip,trim=0 {.4\ht0} 0 {.4\ht0} ,width=\textwidth]{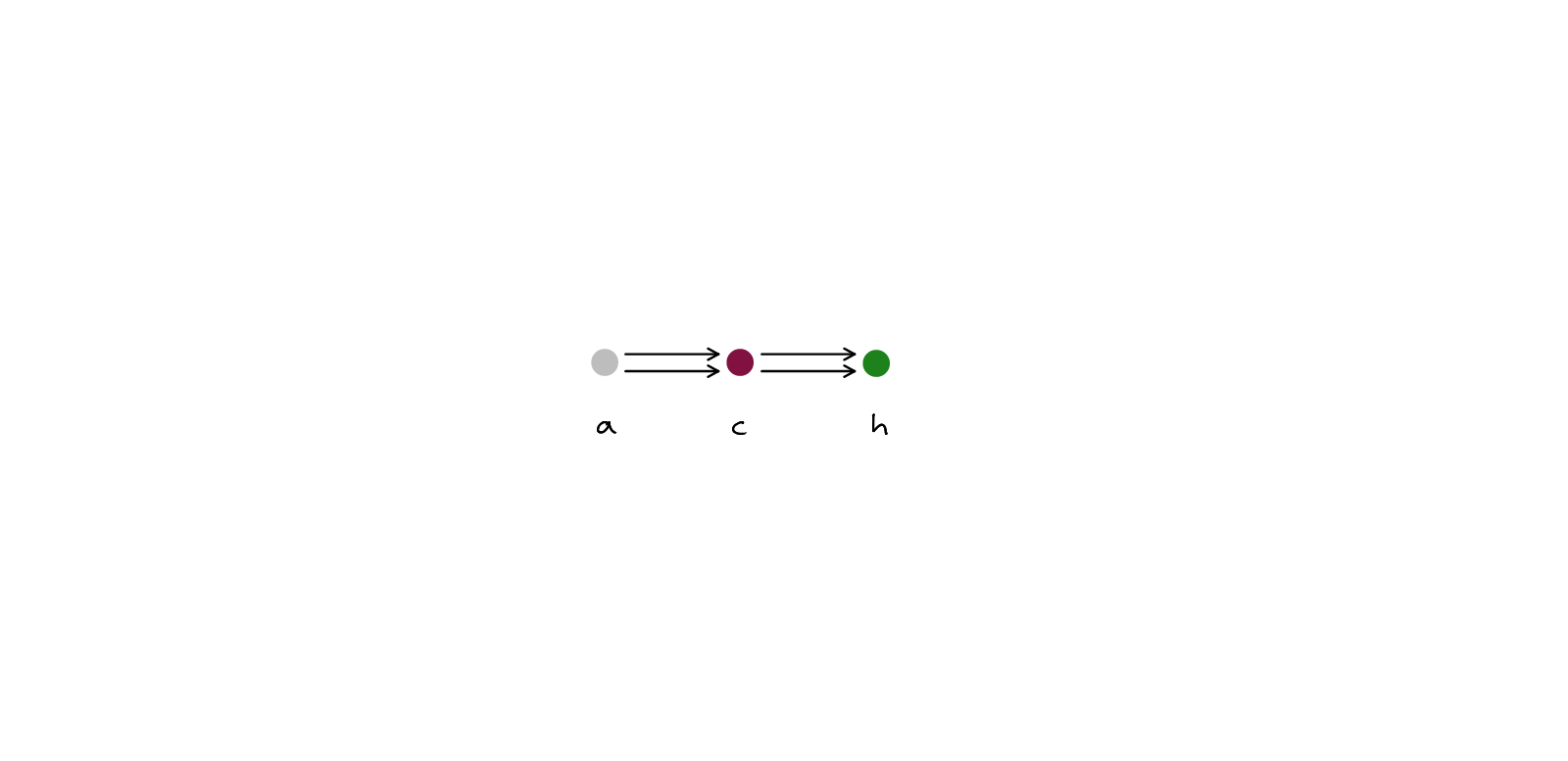}
\endgroup\end{restoretext}
Note that the category of elements construction is functorial in that a map of globular sets (that is, a natural transformation) $\alpha : S \to T$ induces a canonical functor
\begin{equation}
\elcat(\alpha) : \elcat(S) \to \elcat(T)
\end{equation}

For each $n \geq 0$, we define
\begin{equation}
\cG^n := \elcat(\lG(-,n))\op
\end{equation}
which has $2n + 1$ elements of the form (for $0 \leq k \leq n$)
\begin{align}
\sigma_{k,n} &= \sigma_{n-1} \sigma_{n-2} \dots \sigma_k  \\
\tau_{k,n} &= \tau_{n-1} \tau_{n-2} \dots \tau_k 
\end{align}
Note $\sigma_{n,n} = \tau_{n,n} = \id_n$.

We construct a $\cG^n$-labelled $n$-globe $\tgl^n$, called the terminal $n$-globe as follows. First define
\begin{align}
s &:= \elcat(\sigma_{n-1,n})\op : \cG^{n-1} \to \cG^n\\
t &:= \elcat(\tau_{n-1,n})\op : \cG^{n-1} \to \cG^n
\end{align}
where both $\sigma_{n-1,n}$ and $\tau_{n-1,n}$ are natural transformations by the Yoneda lemma
\begin{equation}
\lG(n-1,n) \iso \mathrm{Nat}(\lGR {n-1}, \lGR n)
\end{equation}
The definition of $\tgl^n$ is now inductive in $n$. $\tgl^0 : \bnum 1 \to \cG^0$ is uniquely determined since $\cG^0$ is the terminal category. We then define $\tgl^{n} : \bnum 1 \to \SIvert {n+1} {\cG^n}$ by setting
\begin{equation}
\abss{ \SIvert n s \tgl^{n-1} \xto {\id_n} \SIvert n t \tgl^{n-1} }
\end{equation}
where we used the double cone construction from the previous section.
\end{constr}

As an example the $\cG^2$-labelled singular $2$-cube $\tgl^2$ is given by the data
\begin{restoretext}
\begingroup\sbox0{\includegraphics{test/page1.png}}\includegraphics[clip,trim=0 {.25\ht0} 0 {.05\ht0} ,width=\textwidth]{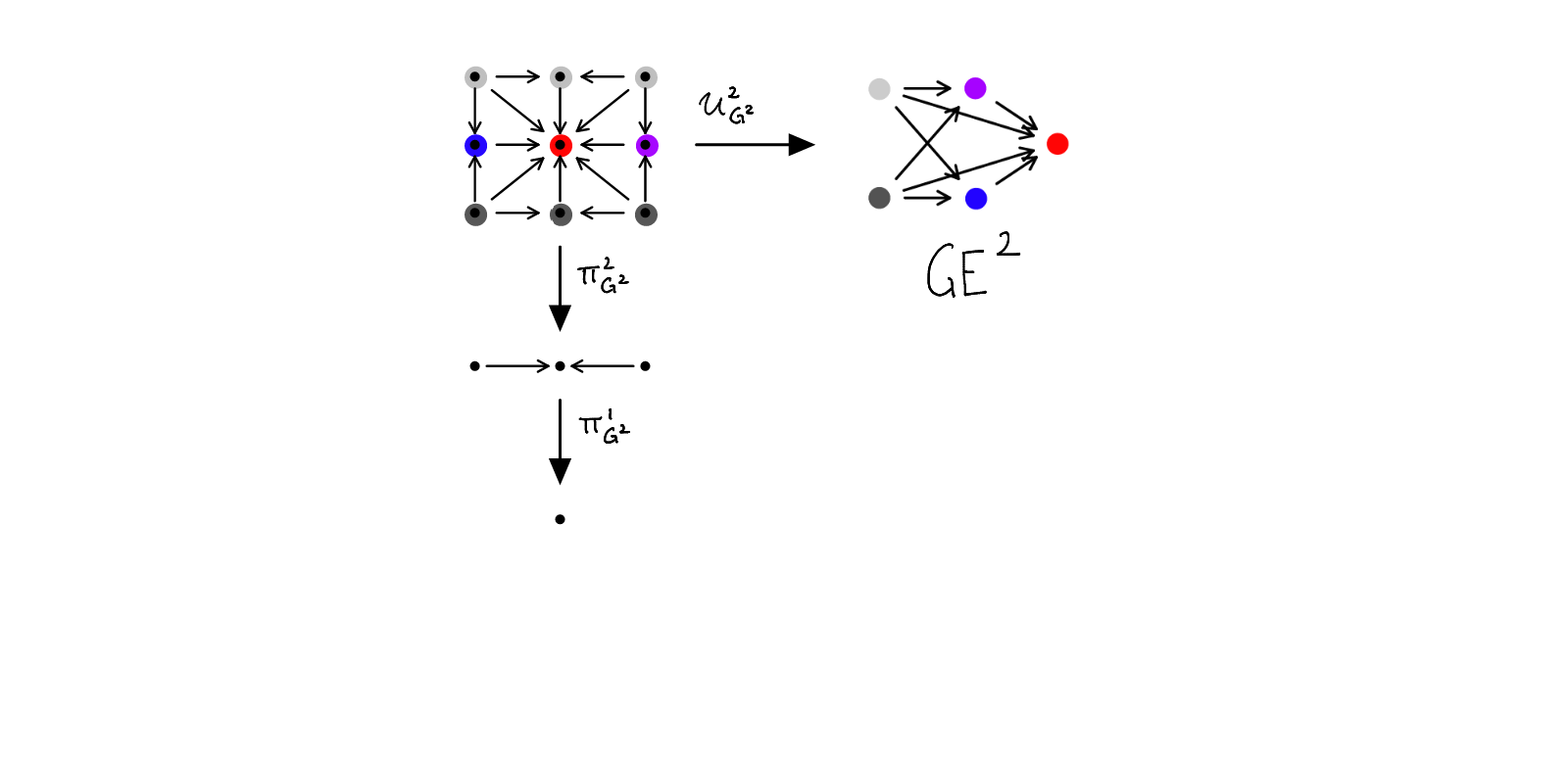}
\endgroup\end{restoretext}

\section{Presentations of higher categories} \label{sec:sum_pres}

We are now in the position to define our first model of higher categories. Recall from \ref{sec:overview} the heuristic distinction between ``fully algebraic", ``algebraic" and ``geometric" models: The former two terms referred to models whose ``algebraic" rules could be (fully or partially) interpreted independently of other foundational systems, and thus for instance implemented as a type theory on a computer. The latter referred to model which fully live within another foundational system, which in our case will be (Set-theory-based) classical Geometry.

The model in this section will be fully algebraic (however, it can also be interpreted within Set Theory). Namely, we will describe presented higher categories which are given by generators and relations. This comes in handy for computational purposes, such as in type theory when working with higher inductive types, or in topology when working with finitely generated complexes. In the latter context we will later rediscover the Hopf map as a generator of $\pi_3(S^2)$.

\subsection{Namescopes} \label{ssec:sum_namescopes}

We start with the notion of namescopes, which are lists of fixed names that will be useful on several occasions for the purpose of comparing structures.

\begin{notn} Given two posets $X$ and $Y$, there is a profunctorial relation $R^{\mathrm{full}}_{X,Y} : X \xslashedrightarrow{} Y$ defined by setting $R^{\mathrm{full}}_{X,Y}(x,y)$ to be true for all $x \in X$ and $y \in Y$.
\end{notn}

\begin{notn}[Namescopes] \label{notn:sum_namescopes} We refer to an $\lN$-ordered list of sets $\sN_0,\sN_1,\sN_2, ... $ as a \textit{namescope} (of dimension $\infty$). Given a namescope $\sN$, we define a functor
\begin{equation}
\Gamma_\sN : \lN \to \PRel
\end{equation}
to map $i \mapsto \Discr \sN_i$ and $(i \to j)$ to $R^{\mathrm{full}}_{\Discr \sN_i,\Discr \sN_j}$ (here $\lN$ is a poset with order $\leq$ and $\Discr$ is the discrete functor defined earlier). 

We denote $\GGamma{}\sN = \sG(\Gamma_\sN)$. We usually identify objects $(i,g) \in \GGamma{}\sN$ with objects $g \in \sN_i$.
\end{notn}

\subsection{\Free{} associative $n$-categories and their morphisms}

In this summary, we describe the definition of \free{} associative $n$-categories only in the case $n = \infty$  (the case of $n < \infty$ is similar and will be discussed in \autoref{ch:presented}).

\begin{defn}[\Free{} associative $\infty$-categories and their morphisms] \label{defn:sum_pres_anc} A \textit{\free{} associative $\infty$-category} $\sC$, is a namescope $\sC_0, \sC_1, \sC_2, ... $, where $\sC_k$ is called the set of \textit{generating $k$-morphisms}, together with data of a certain $k$-cube for each $g \in \sC_k$, called \textit{type} of $g$, and denoted by
\begin{equation}
\abss{g} : \bnum{1} \to \SIvert {k} {\GGamma{}{\sC}}
\end{equation}
such that
\begin{equation}
\abss{g} = \abss{S \xto g T}
\end{equation}
for two $(k-1)$-morphisms $S,T \in \Comp(\sC)_{k-1}$: (mutually) inductively, a \textit{$k$-morphism} (also called a \textit{$k$-composite}) $\scA$ of $\sC$, written $\scA \in \Comp(\sC)_k$, is defined to be a globular normalised $\SIvert k {\GGamma{}\sC}$-cube satisfying
\begin{itemize}
\item \textit{Well-typedness}: For each $p \in \tsG k(\scA)$, $\tsU k_{\scA}(p) = f \in \sC_l$, we require $l \leq k$ and
\begin{equation}
\NF{\scA \sslash p}^k = \Id^{k-l}_{\abss{f}}
\end{equation}
\end{itemize}
$\abss{g}$ is called the type of $g$, and the assignment $g \mapsto \abss{g}$ is called type data for $\sC$.
\end{defn}

We remark, that given a \free{} associative $n$-category $\sC$ there is a function 
\begin{equation}
\abss{-} : \sC_k \to \Comp(\sC)_k
\end{equation}
which maps $g \in \sC_k$ to $\abss{g} \in \Comp(\sC)_k$. We further introduce the following notation

\begin{notn}[Minimal labelling category] \label{defn:sum_red_lab_cat} Let $\sC$ be a presented associative $\infty$-category. Define its \textit{minimal labelling category of $\sC$}, denoted by $\redGamma \sC$, to be the wide subcategory (that is, containing all objects) of $\GGamma{}{\sC}$ with morphisms $c \to d$ (for $c\in \sC_l$, $d \in \sC_k$, $l < k$) whenever $c \in \im(\tsU k_{\abss{d}})$. 
\end{notn}
\noindent  The meaning of $\redGamma \sC$ is that it is the minimal labelling category through which $\sU^k_f$ of all morphisms $f \in \Comp(\sC)_k$ can be factored.

\begin{rmk}[$\Comp(\sC)$ is a globular set] \label{rmk:sum_mor_glob} We remark that there is a globular set $\Comp(\sC)$ with components
\begin{equation}
\Comp(\sC)(k) := \Comp(\sC)_k
\end{equation}
 and source and target maps being given by globular source and target. To see that this is a valid globular set, note that the well-typedness condition is inherited by globular sources and targets.
\end{rmk}

\begin{defn}[Category of presentations $\pCat_\infty$] Let $\sC$ and $\sD$ be \free{} associative $\infty$-categories. A map of presentations $\alpha : \sC \to \sD$, is a bundle functor $\alpha : \GGamma{}\sC \to \GGamma{}\sD$ (that is, a functor satisfying $\pi_{\Gamma_\sD}\alpha =  \pi_{\Gamma_\sC}$) such that for all $g \in \sC_k$, $k \leq n$, we have (cf. \autoref{def:sum:relabelling})
\begin{equation}
\SIvert k \alpha \abss{g} = \abss{\alpha(g)}
\end{equation}
If $\alpha$ is injective, then it is referred to as an inclusion of presentations $\sC \into \sD$.

The category $\pCat_\infty$ is defined to have presented associative $\infty$-categories as objects, and maps of presentations as morphisms.
\end{defn}

\begin{rmk}[Functors]
Maps of presentations are special instances of functors of higher categories, but not every general functor can be written as a map of presentation. A general account of functors, transformations and $k$-transformations will be given in future work.
\end{rmk}

\subsection{\Free{} associative $n$-categories from globular sets} \label{ssec:sum_panc_from_gset}

In this section we record the fact that every globular set canonically determines a presented associative $\infty$-category. 

\begin{defn}[\Free{} associative $n$-categories from globular sets] \label{defn:sum_panc_from_glob} Let $S$ be a globular set. We define a presented associative $\infty$-category $\kC(S)$ from $S$ by setting $\kC(S)_k = S_k$ for $k \geq 0$, and further for each $g \in S_k$ defining
\begin{equation}
\abss{g} = \SIvert n {\elcat(g)\op} \tgl^n
\end{equation}
where we employed the Yoneda lemma to find $g : G(-,n) \to S$ and then used functoriality of $\elcat(-)$ as well as a relabelling functor (cf. \autoref{def:sum:relabelling}). We further implicitly used relabelling by a functor $\elcat(S)\op \to \GGamma{}{\kC(S)}$ that acts as the identity on objects (which uniquely determines it). 

The above construction gives rise to a functor
\begin{equation}
\kC : [\lG\op,\SetCat] \to \pCat_\infty
\end{equation}
\end{defn}

In other words, $\kC(S)$ is the presented associative $\infty$-category whose generating $k$-morphisms are the elements of $S_k$ with types determined by the sources and targets in $S$. 

As an example, recall our globular set $S_{(3)}$ from \autoref{ssec:sum_term_glob}. The type $\abss{c}$ of $c$ in $\kC(S_{(3)})$ is 
\begin{restoretext}
\begingroup\sbox0{\includegraphics{test/page1.png}}\includegraphics[clip,trim=0 {.45\ht0} 0 {.1\ht0} ,width=\textwidth]{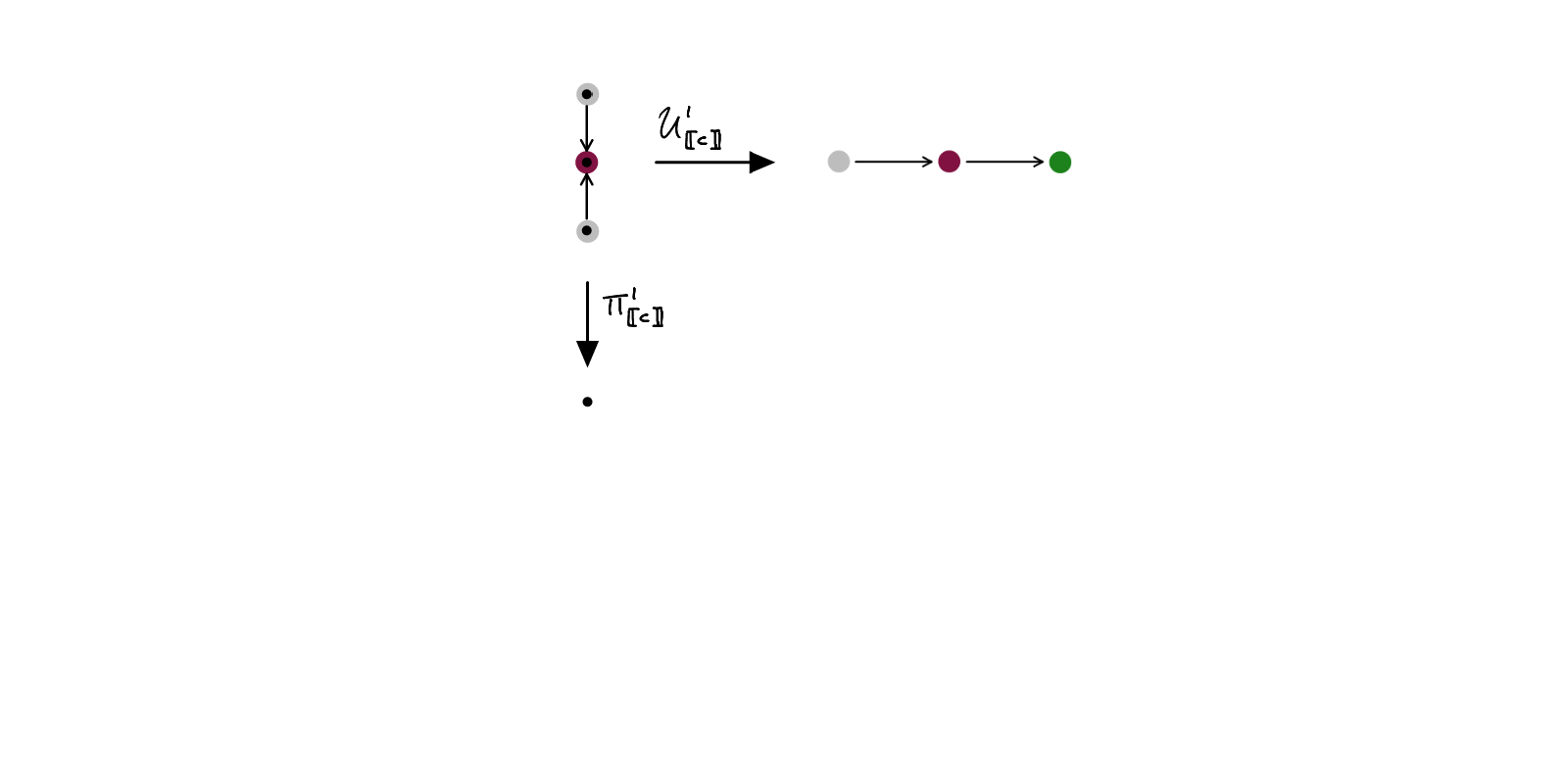}
\endgroup\end{restoretext}
Similarly the type $\abss{h}$ of $h$ in $\kC(S_{(3)})$ is
\begin{restoretext}
\begingroup\sbox0{\includegraphics{test/page1.png}}\includegraphics[clip,trim=0 {.15\ht0} 0 {.15\ht0} ,width=\textwidth]{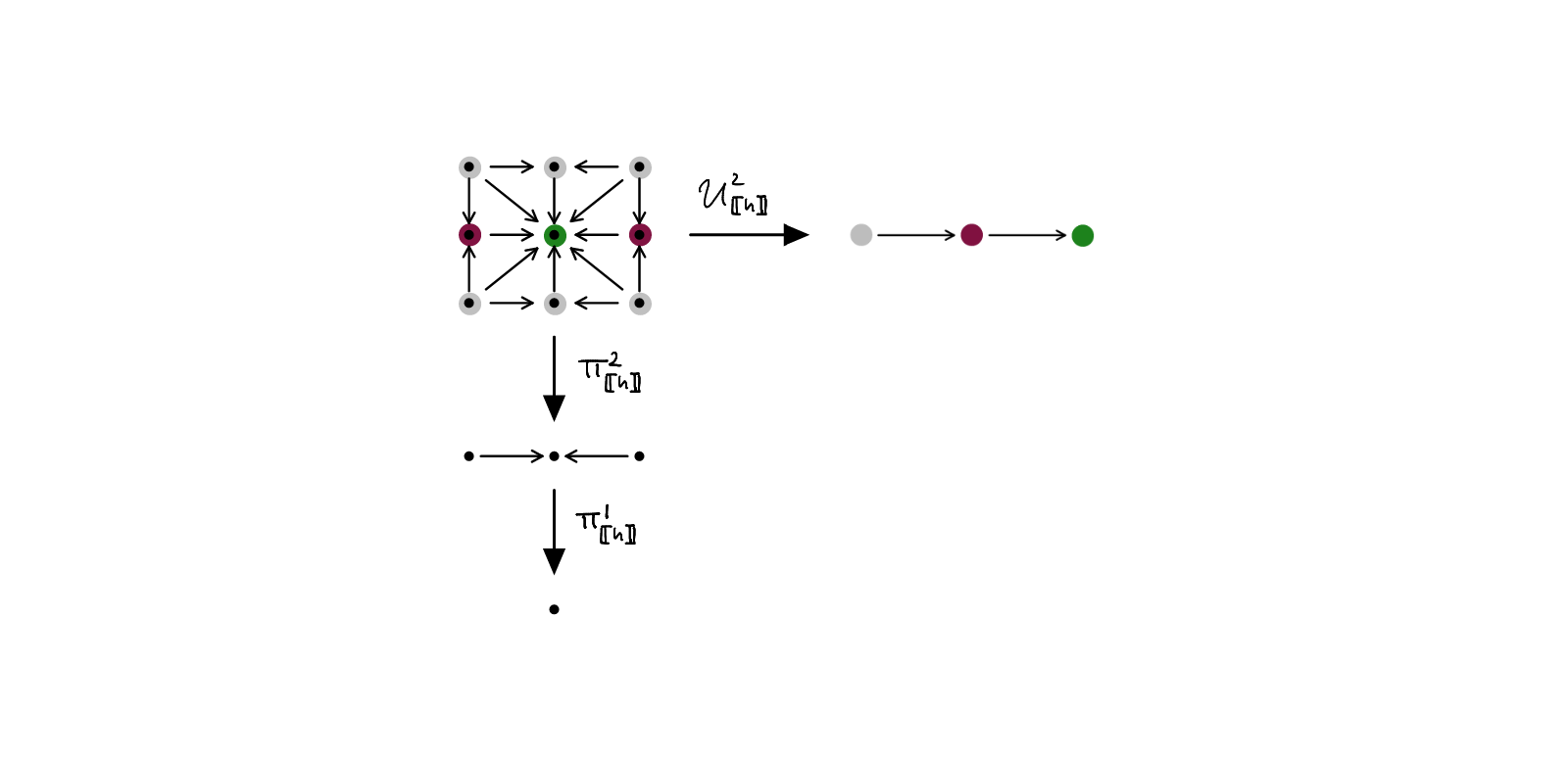}
\endgroup\end{restoretext}

\section{Coherent invertibility}

\Free{} associative $n$-groupoids are \free{} associative $n$-categories in which all generators are invertible. Usually when thinking about invertibility in higher category theory, we use a coinductive approach: ($k$-)morphisms $f$ and $g$ are invertible if there are invertible ($(k+1)$-)morphisms $\alpha : fg \to 1$ and $\beta : gf \to 1$. However, in reality these higher invertible morphisms ($\alpha, \beta...$) have many other interesting \textit{coherence} relations, just as associators are related by ``higher associator" coherences such as the pentagonator). In this section we define a candidate for the theory of (fully) coherent invertibility $\TI$ (in fact, we will see that there appears to be a spectrum of choices for such a theory, cf. \autoref{rmk:sum_TI_choices}). This will then lead us to the definition of presented associative (coherent) $\infty$-groupoids, but $\TI$ will also play a role in the definitions of associative $n$-categories.

\subsection{Unbiased, full presentation} \label{ssec:sum_groupoids}

The (full-coherent) \textit{theory of $\infty$-invertibility} $\TI$, freely generated by a single invertible $1$-morphism, is the \free{} associative $\infty$-category inductively defined as follows. Set
\begin{equation}
\TI_0 = \Set{-,+}
\end{equation}
Then, inductively, $\TI_k$ contains elements $\ic_{S \equiv T}$, called the coherence between $S$ and $T$ where $S, T \in \Comp(\TI)_{k-1}$, which satisfy that $(\tsU k_{\abss{S \xto c T}})\inv(x)$ is non-empty for $x \in \Set{-,+}$ (where we denote $c = \ic_{S \equiv T}$). We then set $\abss{c} = \abss{S \xto c T}$. This completes the inductive definition of $\TI$. 

\begin{rmk}[$r$-connected coherent invertibility] \label{rmk:sum_TI_choices} The condition of mere non-emptiness in the definition of $\TI$ can be strengthened to $r$-connectedness for $-1 \leq r \leq \infty$, or even to requiring a homeomorphism with the $(k-1)$-ball. As we will discuss later on, the latter has an elegant geometric interpretation since it guarantees that $k$-morphisms in $\TI$ look like $(k+1)$-framed $k$-tangles.
\end{rmk}

\begin{rmk}[Binary theory of invertibility] Further to the previous remark note that $\TI$ (even after strengthening the non-emptiness condition in its definition) contains infinitely many generators in dimension $n \geq 3$. Just as the extended $n$-cobordism categeory is believed to have a finite presentation, we expect $\TI$ to have a finite (binary) presentation. A candidate for the latter is the subject of work in progress.
\end{rmk}

\subsection{Freely adjoining coherently invertible generators}

Recall from \autoref{defn:sum_red_lab_cat} the \textit{minimal labelling category of $\TI$}, denoted by $\TTI$.

\begin{constr}[Freely adjoining a coherently invertible generator] \label{constr:sum_free_coh_inv_mor} Let $\sC$ be a \free{} associative $\infty$-category. Let $x,y \in \Comp(\sC)_m$ be $m$-morphisms with coinciding source and target. We define a \free{} associative $\infty$-category $\sC \igadd{x,y} \ig$ called the \textit{\free{} associative $\infty$-category obtained by adjoining an (coherently) invertible generator $\ig$ in between $x$ and $y$}. 

$\sC \igadd{x,y} \ig$ is fully determined by the following structure
\begin{enumerate}
\item For all $k$, there are injective functions
\begin{equation}
\tinv^\ig_k : \Comp(\TI)_k \into \Comp(\sC \igadd{x,y} \ig)_{k + m}
\end{equation}
such that
\begin{align}
\tinv^\ig_0(-) &= x \\
\tinv^\ig_0(+) &= y
\end{align}
and a functor
\begin{equation}
\ttinv^\ig : \TTI \to \SIvert m {\GGamma{}{\sC \igadd{x,y} \ig}}
\end{equation}
such that for any $\scD \in \Comp(\TI)_k$ we have
\begin{equation}
\SIvert k {\ttinv^\ig} \scD = \tinv^\ig_k (\scD)
\end{equation}

\item For all $k > 0$, the pullback
\begin{equation}
\xymatrix{ \TI_k \ar[r]^{i_k} \ar[d]_{\abss{-}} \pullback & (\sC \igadd{x,y} \ig)_{k + m} \ar[d]^{\abss{-}} \\
\Comp(\TI)_k \ar[r]_-{\tinv^\ig_k} & \Comp(\sC \igadd{x,y} \ig)_{k +m}}
\end{equation}
is such that $(\sC \igadd{x,y} \ig)_k$ is the disjoint union of $\sC_k$ and the image of $i_k$, and the inclusion $\sC_k \subset (\sC \igadd{x,y} \ig)_k$ is required to extend to an inclusion of presentations $\sC \into (\sC \igadd{x,y} \ig)$. 
\end{enumerate}

\noindent We note that later on $i_k(\ic_{S \equiv T})$ will be denoted by $\icg\ig_{S \equiv T}$, and $\sC \igadd{x,y} \ig$ will be constructed more explicitly.

Morally, $\sC\igadd{x,y} \ig$ is of course just the category $\sC$ with both a generator $\ig$ and all the generator's invertibility data adjoined to it.
\end{constr}

There is an analogous treatment for adjoining \textit{sets} of of invertible generators which uses the fact that given sequences of inclusions of presentations
\begin{equation}
\sC^0 \into \sC^1 \into \sC^2 \into \dots
\end{equation}
then there is a natural way to pass to the colimit $\colim_i (\sC^i)$ of this sequence. Adjoining a set $I$ of invertible generators $\ig \in I$ with source $x_\ig$ and target $y_\ig$ (both of them living in $\Comp(\sC)$) to $\sC$ will give an $\infty$-category denoted by $\sC \igadd{x_\ig,y_\ig} \Set{\ig \in I}$.

\subsection{\Free{} associative $n$-groupoids}

\begin{defn}[\Free{} associative $\infty$-groupoids] \label{def:sum_infty_group} An associative $\infty$-groupoid $\sX$ is a \free{} associative $\infty$-category obtained by the following inductive procedure. We start with a \free{} associative $\infty$-category $\sX^{(0)}$ consisting only of $0$-generators (called \textit{$0$-cells}).

In the $k$-th step ($1 \leq k \leq n$), we take a set $\sX^k$ called the set of \textit{$k$-cells}, together with $(k-1)$-morphisms $s_c, t_c \in \Comp(\sX^{(k-1)})_{k-1}$ for each $c \in \sX^k$, such that $s_c$, $t_c$ agree on their sources and targets. We then define 
\begin{equation}
\sX^{(k)} = \sX^{(k-1)} \igadd{s_c,t_c} \Set{c \in \sX^k}
\end{equation}
Finally, $\sX$ is defined as the colimit of
\begin{equation}
\sX^{(0)} \into \sX^{(1)} \into \sX^{(2)} \into \dots
\end{equation}
\end{defn}

Morally, the above procedure should be understood as an algebraic version of inductively attaching cells to a CW-complex. We will discuss this analogy more precisely in \autoref{sec:sum_CW} (and in \autoref{ch:geom}).

\section{Generic composites}

In geometric terms, the notion of \textit{generic} composites intends to capture composites (that is, manifold diagrams typed in a presented associative $n$-category) which are stable under perturbation: if we perturb a subcube of the $n$-cube in any direction then up to equivalence the manifold diagram stays constant. In yet other words, the manifolds of the manifold diagram are in ``generic position". The prototypical example and non-example are the following
\begin{restoretext}
\begingroup\sbox0{\includegraphics{test/page1.png}}\includegraphics[clip,trim=0 {.2\ht0} 0 {.25\ht0} ,width=\textwidth]{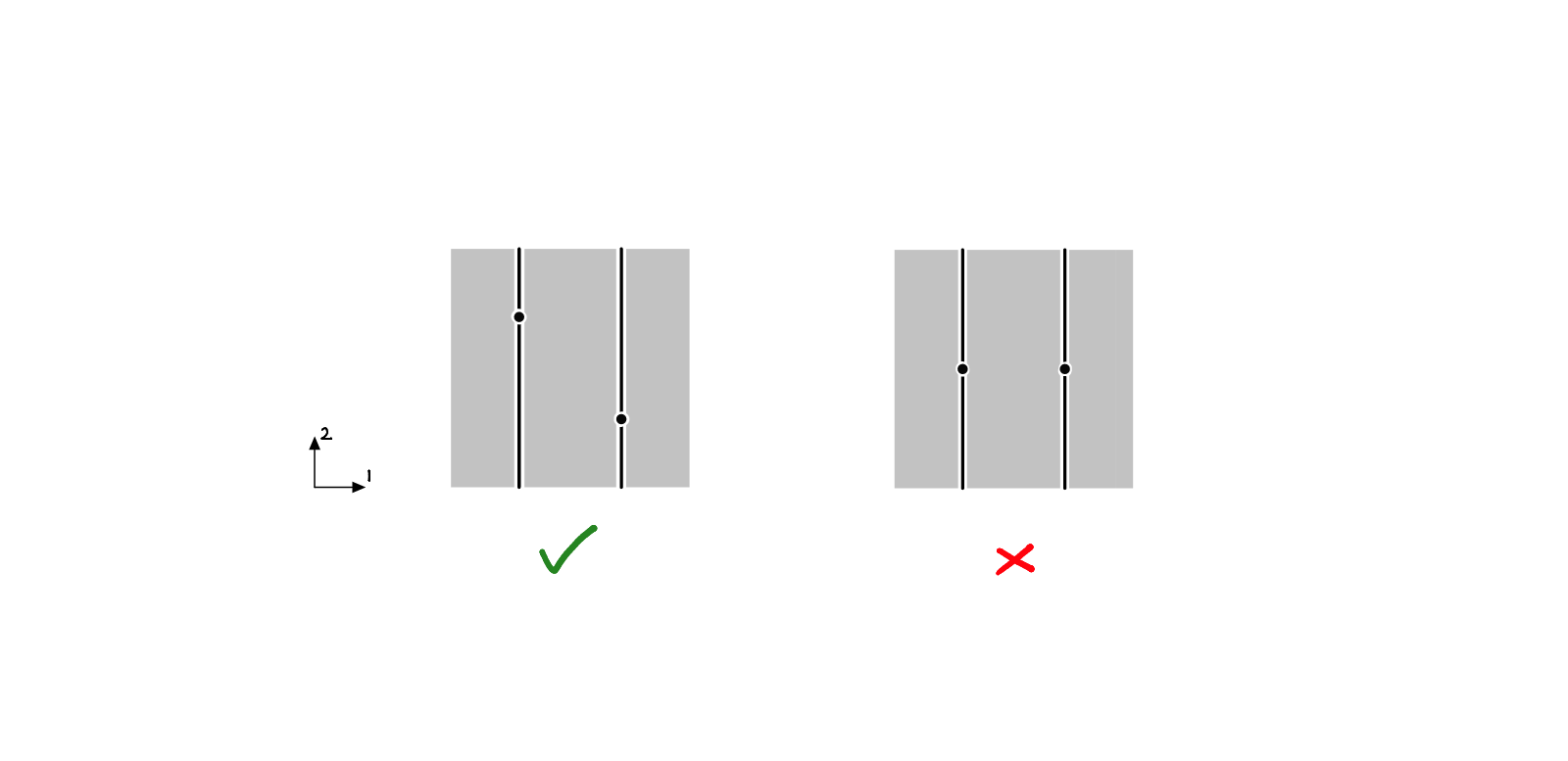}
\endgroup\end{restoretext}
with the left diagram being generic, and the right diagram not being generic (indeed perturbing one of the $0$-strata in it up or down changes the equivalence class of the diagram). More examples will be given in \autoref{ch:composition}.

We will give two (equivalent) characterisations of generic composites: The first characterisation will be ``top-down", and identify generic composites as composites satisfying a certain condition (suitably called the condition of ``being generic"). The second characterisation will be ``bottom-up", and build generic composites from composing generators and so-called elementary homotopies in a special way. This special way of composition will be called ``whiskering" (and coincides with the usual notion of whiskering of $n$-globes up to appropriately identifying $n$-globes and $n$-cubes).

\subsection{First characterisation} \label{ssec:sum_gcomps}

We start with the important and central notion of homotopy, which is a morphism in an associative $n$-category that is locally trivial. The notion can be formulated in even more general terms as we will see in \autoref{ch:composition}.

\begin{defn}[Homotopies] Let $\sC \in \pCat_\infty$. $H \in \Comp(\sC)_n$ is called a \textit{homotopy} if $\im(\tsU n_H) \cap \sC_n = \emptyset$.

$g \in \SIvert n {\GGamma{}{\sC}}$ with $\NF{g} \in \Comp(\sC)_n$ is called \textit{homotopically trivial} if there is a homotopy $H \in \Comp(\sC)_{n+1}$ and a morphism $h \in \Comp(\sC)_{n-1}$ such that $\gsrc(H) = \NF{g}$ and $\gtgt(H) = \Id(h)$. We say $\theta : g \mono f$ is homotopically non-trivial if $g$ is.
\end{defn}

\begin{defn}[Generic composites] \label{defn:sum_gcomps} Let $\sC \in \pCat_\infty$. Inductively in $k$, we say $f \in \Comp(\sC)_k$ is \textit{generic}, written $f \in \GComps(\sC)_k$ if either $k = 0$ or
\begin{enumerate}
\item for all $b \in \reg(\tsG 1 (f))$ we have $\tsU 1 _f \Delta_b \in \GComps(\sC)_{k-1}$
\item for all $a \in \sing(\tsG 1 (f))$ there is a unique homotopically non-trivial $\theta_a : f_a \mono f$ with $a \in \im(\theta^1_a)$ such that for any other homotopically non-trivial $\phi : h \mono f$ we have $\theta_a \mono \phi$.
\end{enumerate}
If only the second condition is satisfied, $g$ is called \textit{non-recursively} generic.
\end{defn}

Morally, the existence of a unique $\theta_a$ guarantees that there is at most one non-trivial ``event" happening at any given singular height.

\subsection{Composing cubes}

We now consider the question of how to compose cubes which will enable us to give a second characterisation of generic cubes. We start with single binary gluing operations. These operations will apply to general singular $n$-cubes, and no structure of a presented associative $\infty$-category will be required. After an initial definition which extends the ordered sum to $\SI$, we will introduce three notions of composition: stacking (denoted by $\stack$), $k$-level stacking (denoted by $\glue k$) and $k$-level whiskering (denoted by $\whisker k n$). Geometrically, $k$-level stacking takes two $n$-cubes and glues them together along a mutual side (facing into the $(n-k+1)$th direction) of the cube. On the other hand, whiskering (geometrically) takes an $n$-globe and $k$-globe and glues them together along a mutual $(k-1)$-boundary.

\begin{defn}[Ordered sum] The ordered sum functor $(-\uplus-) : \SI \times \SI \to \SI$ is defined to map
\begin{equation}
(\singint {k_1} \uplus \singint {k_2}) = \singint {k_1 + k_2}
\end{equation}
and for $f_i : \singint {k_i} \to \singint {l_i}$ we set
\begin{equation}
(f_1 \uplus f_2)(a) = \begin{cases} f_1(a) & \text{if~} a \leq 2k_1 \\ f_2(a - 2k_1) + 2l_1 & \text{if~} a \geq 2k_1 \end{cases}
\end{equation}
\end{defn}

\begin{constr}[Composite cubes and globes] Firstly, given two families $\scA, \scB : X \to \SI$ we define $\scA \stack \scB$ to be the family
\begin{equation}
X \xto {\Delta} X \times X \xto {\scA \times \scB} \SI \times \SI \xto {\uplus} \SI
\end{equation}
where $\Delta$ is the diagonal functor. $\scA \stack \scB$ is called the \textit{stacking of $\scA$ with $\scB$}. Note that there are canonical (fully-faithful) poset inclusions
\begin{align}
\stacklow &: \sG(\scA) \into \sG(\scA \stack \scB) \\
\stackup &: \sG(\scB) \into \sG(\scA \stack \scB)
\end{align}
defined by
\begin{align}
\stacklow(x,a) &= (x,a) \\
\stackup(x,b) &= (x,b + l^\scA_x)
\end{align}
where $\scA(x) = \singint {l^\scA_x}$. We remark that $\scA$, $\scB$ are implicit in the notation $\stacklow$, $\stackup$.

Secondly, for $1 \leq k \leq n$, a $\cC$-labelled singular $n$-cube $\scC$ is said to be the \textit{$k$-level stacking of $\cC$-labelled singular $n$-cubes $\scA$ and $\scB$}, written $\scC = \scA \glue k \scB$, if
\begin{itemize}
\item For $0 \leq l < k-1$
\begin{equation}
\tusU l_\scA = \tusU l_\scB = \tusU l_\scC
\end{equation}
\item We have 
\begin{equation}
\tusU {k-1}_\scC = \tusU {k-1}_\scA \stack \tusU {k-1}_\scB
\end{equation}
\item We have
\begin{align}
\tsU k_\scA &= \tsU k_{\scC} \stacklow \\
\tsU k_\scB &= \tsU k_{\scC} \stackup
\end{align}
\end{itemize}

Finally, for $1 \leq k \leq n$, a $\cC$-labelled singular $n$-globe $\scC$ is said to be the\textit{ $k$-level post-whiskering composition of a $\cC$-labelled singular $n$-globe $\scA$ and a $(n-k+1)$-globe $\scB$,} written $\scC = \scA \whisker k n \scB$, if $\scC = \scA \glue k \widetilde \scB$ where $\widetilde \scB$ is the $\cC$-labelled singular $n$-cube defined as follows: we set
\begin{equation}
\tsU {k-1}_{\widetilde \scB} := \tsU 0_\scB !
\end{equation}
where $! : \tsG {k-1}(\scA) \to \bnum{1}$ is the terminal functor, and further (for $0 \leq l < k-1$)
\begin{equation}
\tusU l_{\widetilde \scB} := \tusU l(\scA)
\end{equation}

Similarly, $\cC$ is called the $k$-level pre-whiskering composition of (an $(n-k+1)$-globe) $\scB$ and (an $n$-globe) $\scA$, written $\cC = \scB \whisker k n \scA$, if $\cC = \widetilde \scB \glue k \scA$ where $\widetilde \scB$ is defined as before.
\end{constr}
We remark that the notation $\scA \whisker k n \scB$ is chosen since the construction involves a stacking of cubes at level $k$. The fact that the dimensions of $\scA$ and $\scB$ might differ is recorded by the subscript $n$, which is the highest dimension among the dimensions of $\scA$ and $\scB$.

As an example, recall $S_{(3)}$ and $\kC(S_{(3)})$ from \autoref{ssec:sum_panc_from_gset} (as well as the types $\abss{c}$ and $\abss{h}$). Then $\abss{h} \glue 1 \abss{h}$ is given by the data
\begin{restoretext}
\begingroup\sbox0{\includegraphics{test/page1.png}}\includegraphics[clip,trim=0 {.15\ht0} 0 {.15\ht0} ,width=\textwidth]{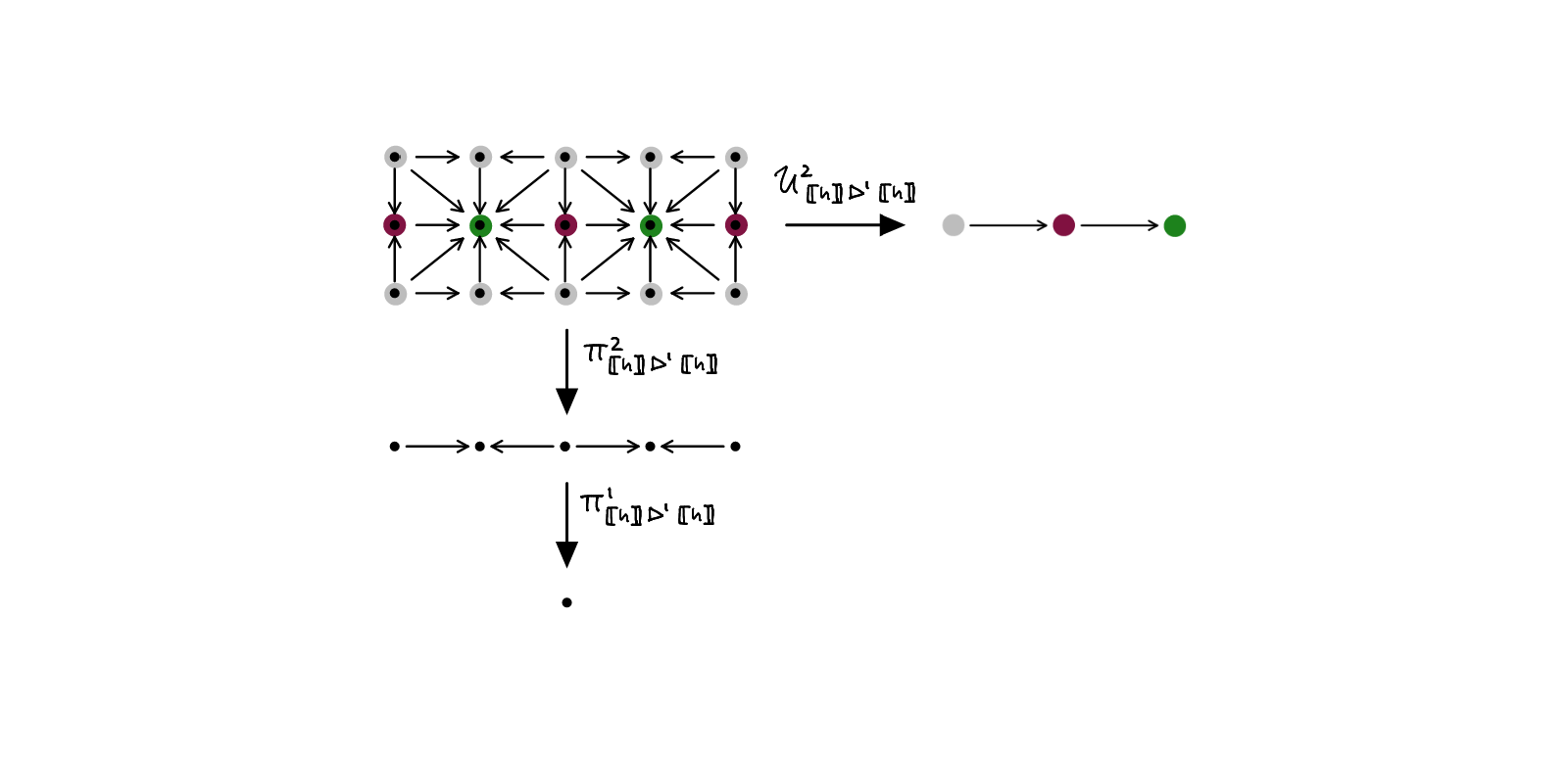}
\endgroup\end{restoretext}
Note that in this case the computation of $\abss{h} \whisker 1 2 \abss{h}$ would have yielded the same result. Similarly, we find $\abss{h} \whisker 2 2 \abss{c}$ is given by
\begin{restoretext}
\begingroup\sbox0{\includegraphics{test/page1.png}}\includegraphics[clip,trim=0 {.05\ht0} 0 {.1\ht0} ,width=\textwidth]{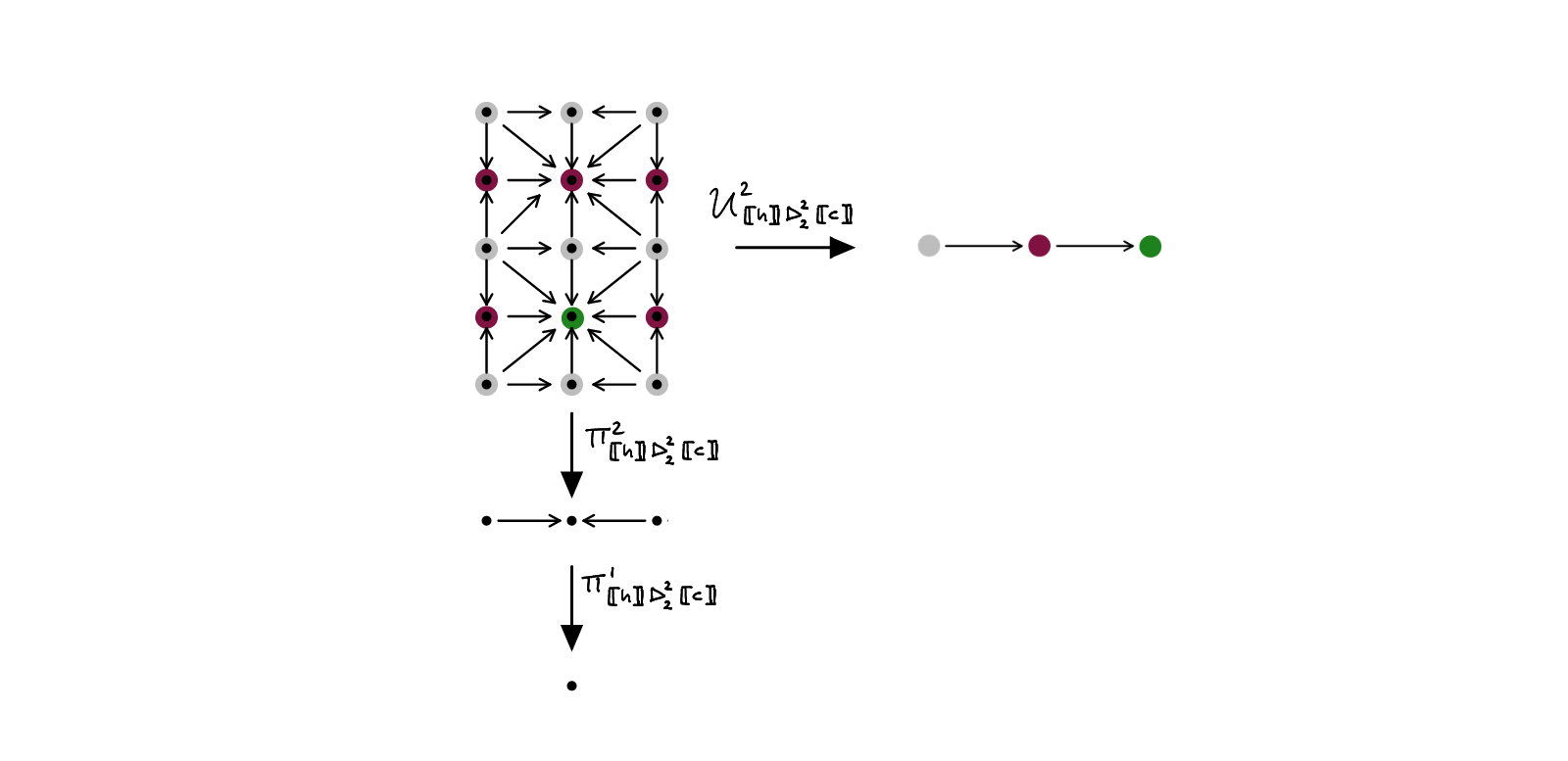}
\endgroup\end{restoretext}
Finally, we find $(\abss{h} \whisker 2 2 \abss{c}) \whisker 1 2 (\abss{1} \whisker 2 2 \abss{h})$ is given by
\begin{restoretext}
\begingroup\sbox0{\includegraphics{test/page1.png}}\includegraphics[clip,trim=0 {.05\ht0} 0 {.1\ht0} ,width=\textwidth]{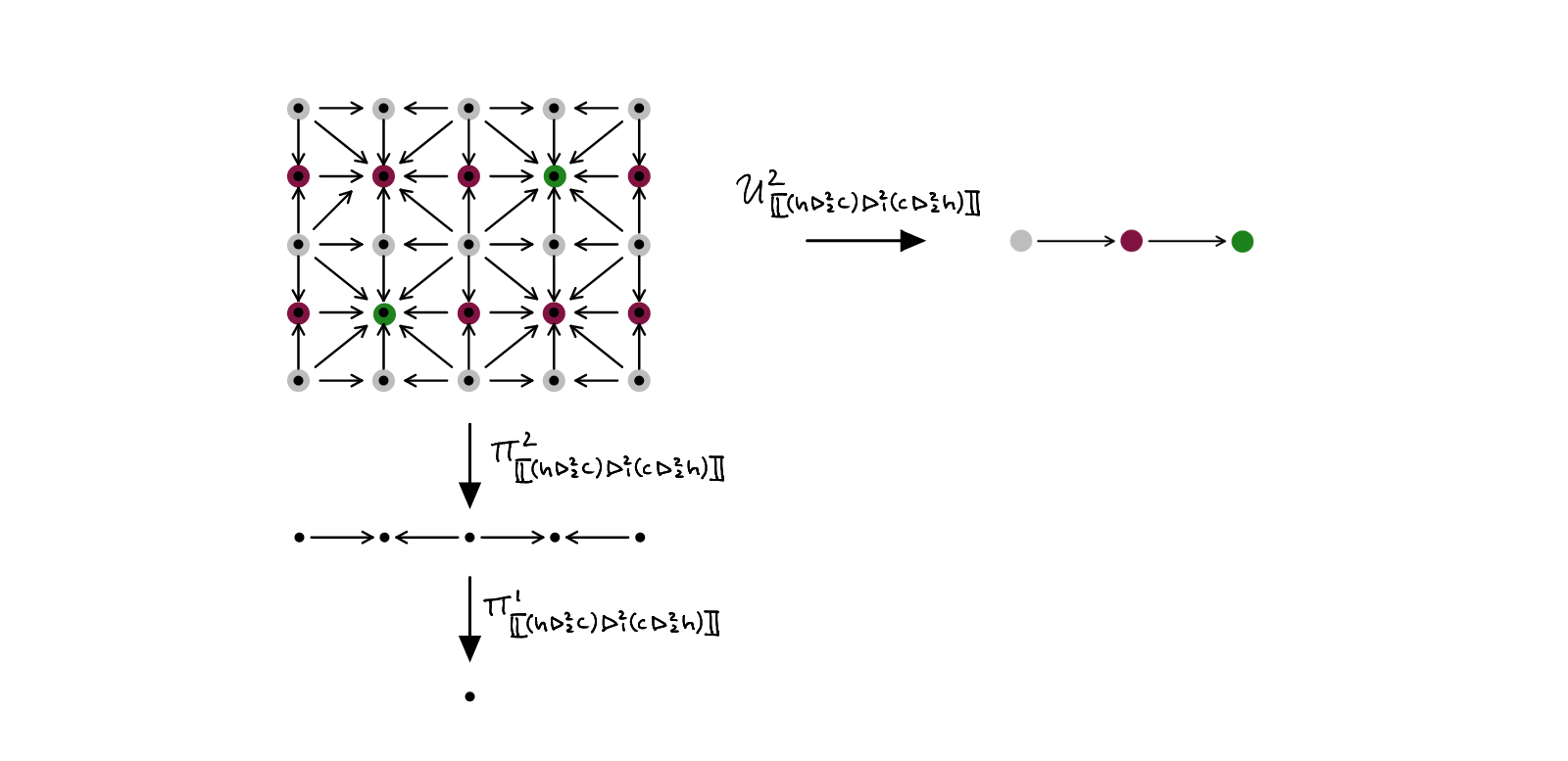}
\endgroup\end{restoretext}

\subsection{Second characterisation} \label{ssec:sum_gcomp}

We now define generic composites as being inductively built from generators and elementary homotopies and using only whiskering operations.

\begin{constr}[Generic composites functor] Let $\sC \in \pCat$. We define the globular subset $\GComp(\sC) \subset \Comp(\sC)$, whose elements are called \textit{generic composites}. Inductively, we first set $\GComp(\sC)_0 = \Comp(\sC)_0$ and then
\begin{enumerate}
\item For $g \in \sC_k$ we have $g \in \GComp(\sC)_k$ if $\gsrc(g),\gtgt(g) \in \GComp(\sC)_{k-1}$
\item For $H \in \Comp(\sC)_k$ an \textit{elementary} homotopy, meaning $H$ is a (non-recursively) generic homotopy (cf. \autoref{defn:sum_gcomps}) with $\gsrc(H),\gtgt(H) \in \GComp(\sC)$, and such that $H = h_1 \whisker j k h_2$ implies either $h_1$ or $h_2$ is an identity
\item If $f_1, f_2 \in \GComp(\sC)_k$ then (if it exists) $f_1 \whisker j k f_2 \in \GComp(\sC)_k$
\end{enumerate}
This completes the construction of $\GComp(\sC)$. $f \in \GComp(\sC)$ is called \textit{elementary} if $f = h_1 \whisker j k h_2$ implies either $h_1$ or $h_2$ is an identity.

We note that this gives rise to an functor
\begin{equation}
\GComp : \pCat_\infty \to \globset
\end{equation}
which is a subfunctor of $\Comp$. Using our functor $\kC : \globset \to \pCat_\infty$, abusing notation, we find an endofunctor
\begin{equation}
\GComp : \globset \to \globset
\end{equation}
which comes with a natural transformation $\abss{-} : \id \to \GComp$ with components $\abss{-}_S : g \mapsto \abss{g} \in \GComp(S)$.
\end{constr}

\begin{claim}[Equivalence of characterisations] Let $\sC \in \pCat_\infty$. Then 
\begin{equation}
\GComps(\sC) = \GComp(\sC) \subset \Comp(\sC)
\end{equation}
\end{claim}

\section{Higher categories} \label{sec:sum_cat}

We will now finally describe a candidate for the notion of associative $n$-categories that stands in line with other traditional Set-theory-based definitions of higher categories---in other words, unlike our previous ``foundation-independent" notion of presented associative $n$-categories, the following will crucially depend on the framework of Set Theory.

\subsection{Algebras of $\GComp$ and their resolution $\kiC{}$} \label{ssec:sum_alg_and_inf}

Recall the endofunctor $\GComp$ on globular sets. Associative $n$-categories will be special ``algebras" for this functor. 

\begin{defn}[Algebras of $\GComp$] An ($n$-truncated\footnote{An $n$-truncated globular set is a globular set with no elements in degree higher than $n$. Any globular set $S$ can be $n$-truncated to yield an $n$-truncated globular set with the same elements, sources and targets as $S$ in degree less or equal to $n$.}) algebra $(S,\iM)$ for $\GComp$ is an ($n$-truncated) globular set $S$, together with a map of ($n$-truncated) globular sets $\iM : \GComp(S) \to S$ which when precomposed with the unit $\abss{-}_S$ gives the identity.
\end{defn}

There is a special presented associative $\infty$-category associated with an algebra $S$, called the resolution of $S$, which replaces all algebraic equations $\iM(f) = g$ by coherently invertible generators.

\begin{constr}[Resolutions of algebras] Let $(S,\iM)$ be an $n$-truncated algebra of $\GComp$. Its \textit{resolution} $\kiC{}(S)$ is the presented associative $\infty$-category constructed as the colimit of a sequence $\kiC{k}(\sC)$ inductively constructed as follows: Define $\kiC{0}(S)  = \kC(S)$. Now assume we have constructed $\kiC{k}(S) $ together with generators $\infeq f k \in \kiC{k}(S)_k$ for each $f \in \GComp(S)_{k-1} \setminus S_{k-1}$. For $f \in S_{k-1}$ we define $\infeq f k := \Id(f)$). Let $f \in \GComp(S)_k \setminus S_k$, and set for $0 < i < k$
\begin{align}
s_i &= \infeq {\gsrc^i(f)} {k-i}\\
t_i &=  \infeq {\gtgt^i(f)} {k-i}
\end{align}
Define $x^f \in \GComp(\kiC{k}(S))_k \subset \Comp(\kiC{k}(S))_k$ to be the whiskering composite
\begin{equation}
x^f := (s_1 \whisker 1 n ... (s_{k-1} \whisker {k-1} k f \whisker {k-1} k t_{k-1}) ... \whisker 1 n t_1)
\end{equation}
(which we claim to exist). Next, set $y^f \in \GComp(\kiC{k}(S))_k$ to equal
\begin{equation}
y^f := \abss{\iM(f)}
\end{equation} 
Then define
\begin{equation}
\kiC{k+1}(S) = \kiC{k}(S) \igadd {x^f,y^f} \Set{\infeq f {k+1} ~|~ f \in \GComp(S)_k}
\end{equation}
which completes the inductive construction. Finally, $\kiC{}(S)$ is defined as the colimit of 
\begin{equation}
\kiC{0}(S)  \into \kiC{1}(S) \into \kiC{2}(S)  \into \dots
\end{equation}
\noindent A morphism $f \in \Comp(\kiC{}(S))$ is said to be of \textit{compositional depth 1} if it doesn't lie in the image of the canonical inclusion $\Comp(\kC(S)) \into \Comp(\kiC{}(S))$ of $n$-truncated globular sets. Otherwise it is of \textit{depth 0}.
\end{constr}

We re-emphasize that, morally, the \textit{resolution} $\kiC{}(S)$ of an algebra $(S,\iM)$ encodes every single equation $\iM(f) = g$ as a higher (coherently) invertible generator $\infeq f k$ which we call a composition witness. Morphisms of compositional depth 1 are exactly those involving composition witnesses.

\subsection{Associative $n$-categories} \label{ssec:sum_anc}

Our candidate notion for an associative $n$-category is now ``nothing but" an algebra of $\GComp$ compatible with the inductive structure of generic composites of its resolution. Here, ``compatibility with the inductive structure of $\GComp$" in particular entails strict equations for compatibility with whiskering composition, and strict equations for elementary homotopies of compositional depth $1$. The only weak coherence data lies in elementary homotopies of depth $0$. We note that the definition can be amended to be weak ``up to depth $k$", with $k = \infty$ yielding a ``fully weak" definition (see \autoref{ch:associative}).

We emphasize the importance of generic composites in this definition, which allow us to separate ``elementary" from ``composite" homotopies by fully classifying the generic composites of the resolution as being generated by generators and elementary homotopies using whiskering composition.

\begin{defn}[Associative $\infty$-categories] An \textit{associative $\infty$-category} $\sC$ is an $n$-truncated algebra $(\sC, \kM_\sC)$ of $\GComp$ where $\sC$ is an $n$-truncated globular such that $\kM_\sC$ factors through the following maps of $n$-truncated globular sets
\begin{equation} \label{eq:sum_complaw}
\xymatrix{ \GComp(\sC) \ar[d] \ar[r]^-{\kM_\sC} & \sC \\
 \GComp(\kiC{}(\sC)) \ar@{-->}[ur]_-{\extkM_\sC}& }
\end{equation}
The factorisation is unique subject to the condition that
\begin{equation}
\extkM_\sC(f_1 \whisker j k f_2) = \kM_\sC(\extkM(f_1) \whisker j k \extkM(f_2))
\end{equation}
and further for all elementary morphisms $f \in \GComp(\kiC{}(\sC))_{k \leq n+1}$ of compositional depth $1$ we have
\begin{equation}
\extkM_\sC(\gsrc(f)) = g = \extkM_\sC(\gtgt(f))
\end{equation}
and (if $ k \leq n$)
\begin{equation} \label{eq:sum_cohlaw}
\extkM_\sC(f) = \kM_\sC(\Id_{\abss{g}})
\end{equation}
\end{defn}
\noindent Morally, the first condition \eqref{eq:sum_complaw} makes our definition ``associative" whereas the second condition \eqref{eq:sum_cohlaw} asks coherences of depth 1 (i.e. elementary homotopies involving composition witnesses) to be strict. The only non-strict coherences are of depth 0, which are homotopies already present in $\GComp(\sC)$. Later on in \autoref{ch:associative} we will see that it is in fact possible to also iterate the idea of resolutions up to ``depth $k$", and amend the above definition to chose weak data for coherences up to depth $k$ instead of making them strict equations. As a result we will obtain a spectrum of definitions ranging from ``fully associative" ($k = 1$) to ``fully weak" ($k = \infty$). \\

We summarise the following for low dimensions.
\begin{itemize}
\item An associative $0$-category is a set
\item An associative $1$-category is an unbiased category
\item An associative $2$-category is an unbiased strict $2$-category
\item An associative $3$-category is an unbiased $\mathbf{Gray}$-category
\end{itemize}

\section{Appendix A: Connection of algebraic and geometric models}
 \label{sec:sum_CW}

We preface this section with the reminder that the following discussion, and similarly the discussion in \autoref{sec:CWcomplexes}, is not (yet) fully rigorous. It was included nonetheless as an outline of subtantial geometric ideas, which are ``precise enough" to allow for many interesting computations, and which also have relevance for potential proofs of the homotopy hypothesis \cite{baez2007homotopy} and the generalised tangle hypothesis \cite{baez1995higher} (in the setting of our fully algebraic definition of higher groupoids given in the previous section).

\subsection{Manifold-like theory of invertibility}

We introduce the following modification of $\TI$ by replacing non-emptiness with a  homeomorphism condition: let $\TI$ from now on denote the sub-presentation (of $\TI$ as defined before) given only by $k$-generators $g$ which satisfy that, for $x \in \Set{-,+}$, $(\sU^k_{\abss{g}})\inv(x)$ is homeomorphic to a $k$-ball. 

The motivation for introducing this strengthened condition can be understood after passing to the geometric realisation as described in \autoref{ssec:coloring}: it guarantees that geometric realisations of generators look like manifolds. Further discussion of this can be found in \autoref{ch:groupoids} and \autoref{ch:geom}.

\subsection{Disclaimer on geometry-related conjectures} \label{ssec:sum_mfdl_diag}

The discussion in \autoref{ch:geom} is based on several implicit conjectures. No honest attempt is made to prove these conjectures or even formulate them precisely. Here, we briefly discuss what type of conjectures could be expected as a step towards formalising the discussion of that appendix.

\begin{enumerate}
\item Firstly, \autoref{ch:geom} freely assumes a ``correspondence of geometry and combinatorics". In \autoref{ssec:coloring} we described a procedure to ``geometrically realise" a $\cC$-labelled singular $n$-cube $A$ as a (topological) $\cC$-labelled $n$-cube $\norm{A}$. Furthermore, in \autoref{ssec:sum_rec_struct} we essentially discussed steps towards constructing a ``refinement" function $\mathrm{Ref}$ in the opposite direction; at least on the subcategory of appropriately \textit{finite} topological labelled $n$-cubes, which we will assume to have chosen here without changing notation. With a little bit more thought one can see that this construction allows one to output a singular cube in normal form for any given (topological) $\cC$-labelled $n$-cube. We make the following conjecture
\begin{conj}[Correspondence of algebra and geometry] \label{conj:manifold_diag} Let $E_{alg}, E_{top}$ denote the class of epimorphisms in $\Buno n \cC (\bnum 1)$ and $\Cubeo n \cC(\bnum 1)$ respectively. Then 
\begin{equation}
\norm{-} :\Buno n \cC(\bnum 1) [E\inv_{alg}] \toot \Cubeo n \cC(\bnum 1) [E\inv_{top}] : \mathrm{Ref}
\end{equation}
induces bijection of objects.
\end{conj}

\item Secondly, \autoref{ch:geom} assumes that we can ``equate" the (framed) piecewise linear and (framed) smooth notion of manifolds in the setting of manifold diagrams for $\infty$-groupoids---the former notion naturally arises from the combinatorics, while the latter naturally arises in the setting of CW-complexes, or more precisely, of a generalised Thom-Pontryagin construction. We formulate 
\begin{conj}[Translating piecewise linearity into smoothness (with framing)] Let $A \in \Buno n \cC$. Then there is $B \simeq \norm{A}$ such that strata of the underlying (flag-foliation-compatible) stratification of $B$ are smooth manifolds.

Moreover and more specifically, if $\sC \in \pCat_\infty$ is a presented associative $\infty$-groupoid, and $f \in \Comp(\sC)_n$ is a morphism of $\sC$, then there is $B \simeq \norm{f}$ such that for all $g \in \sC_k$ the union of regions in $B$ with labels $\ic^g_{S\equiv T}$ is a ($n - k$)-manifold admitting a canonical normal framing.
\end{conj}

\item Finally, \autoref{ch:geom} assumes that for any framed stratification $A$ of $S^n$ obtained as an output of the generalised Thom-Pontryagin construction, there is a conical $(n+1)$-manifold diagram equivalent to the cone of $A$ up to identifying the $(n+1)$-cube and the $(n+1)$-disk.
\end{enumerate}

\subsection{Connection of $\infty$-groupoids and CW-complexes} \label{ssec:sum_CW}

In this section, we will (leisurely) sketch a translation from CW-complexes to presented associative $\infty$-groupoids. First, recall the homotopy hypothesis from \autoref{ssec:intro_hom_hyp} which states that the homotopy theory of groupoids is equivalent to the homotopy theory of topological spaces. We focus on spaces for a moment.
\begin{itemize}
\item How can we describe a given space $X$? One answer is to study how $k$-disks $D^k$ and $k$-spheres $S^k$ can be mapped into $X$, that is, the study of the homotopy classes (rel boundary) of maps $[D^k, X]$ respectively $[S^k,X]$ (cf. \cite{hatcher2001algebraic}). At least categorically, this approach can be naturally justified: elements of $[D^k,X]$ can be regarded as $k$-morphism of the corresponding higher groupoid $X$, whereas $[S^k,X]$ describes valid types, or Hom-spaces, of $(k+1)$-morphisms in $X$ (for instance the morphism $f : A \to B$ is of type $A \to B$).

\item How can we build a space $X$? One answer is to build spaces by inductively attaching $(k+1)$-disks of a given ``type", more commonly called attaching map, $f \in [S^k,X]$, thereby obtaining a new space $X \cup_f D^{k+1}$. Doing so inductively leads to the notion of CW-complexes (cf. \cite{hatcher2001algebraic}). Again, categorically this construction can be justified: attaching a disk corresponds to freely adding (a coherently invertible\footnote{Adding a coherently invertible morphism means that we not only add a morphism to a category, but also its inverse as well as ``higher coherence data" which witnesses how the morphism and its inverse interact. This coherence data is given by the higher category $\TI$ that was defined in \autoref{ssec:sum_groupoids}. The process of freely adding a coherently invertible morphism was described in \autoref{constr:sum_free_coh_inv_mor}. A more detailed discussion can be found in \autoref{ch:groupoids}.}) morphism to the Hom-space of type $f$. It is important to note that if $X$ is a CW-complex, then any $g \in [S^k,X]$ has a representative $g' : S^k \to X^{(k)}$ where $X^{(k)}$ is the $k$-skeleton of $X$. Categorically, this can simply be interpreted as saying that the source and target of a $(k+1)$-morphism are $k$-morphisms (and not $l$-morphisms for $l > k$),
\end{itemize}

Our translation of a CW-complex $X$ into a groupoid $\sX$ proceeds by the following three steps: first, translate an attaching map of a $(k+1)$-cell $f$ of $X$ into a stratification $\kP(f)$ of the $k$-sphere by taking inverse images of ``dual strata" (of the so-called dual stratification $\kD(X^{(k)})$ of $X^{(k)}$). Then form the ``cone" of $\kP(f)$ to build a stratification $\cone(\kP(f))$ of the $(k+1)$-disk (attaching this stratified disk back onto $X^{(k)}$ inductively constructs the dual stratification $\kD(f)$). Finally, chose a ``globular foliation" $f\Fol$ of $\cone(\kP(f))$ in order to pass from undirected to directed space. This outputs a manifold diagram which can then be translated into a generating $k$-morphism $f$ of $\sX$. In summary this translates cells of $X$ into generators of $\sX$.

More in detail, starting from a CW-complex $X$ we construct a groupoid $\sX$ as follows
\begin{enumerate}
\item \textit{Encoding attaching maps}: An attaching map $f \in [S^k,X^{(k)}]$ can be encoded in a ``dual form" using a generalised Thom-Pontryagin construction (a version can be found as Prop. VII.4.1 in \cite{buonchristiano1976geometric}) which for a given closed $k$-manifold $M$  provides us with a 1-to-1 correspondence
\begin{equation}
\kP_M : [M,X^{(k)}] \iso \Set{\text{framed $X^{(k)}$-stratifications of $M$ up to ``cobordism"}}
\end{equation}
Here, an $X^{(k)}$-stratification of $M$ is a decomposition of $M$ into (possibly disconnected) manifolds, also called \textit{strata}, which correspond to cells in $X^{(k)}$. We will mainly be interested in the case where $M$ equals $S^k$, we will usually denote $\kP_{S^k}$ just by $\kP$.

We give a brief sketch of the construction: one first inductively constructs a stratification $\kD(X^{(k)})$ of $X^{(k)}$ itself, which is such that each cell $g$ of $X^{(k)}$ corresponds to a stratum $g\dualdag$ of this stratification. This is called the \textit{dual stratification} $\kD(X^{(k)})$ of $X^{(k)}$ ($g$ and $g\dualdag$ coincide in the base case $k = 0$, and the inductive construction will be given in the next step). Up to certain assumptions on a representative $x$ of a same-named class $x \in [M,X^{(k)}]$, we can then define $\kP_M(x)$ to be the stratification of $M$ whose (framed) strata are inverse images $x\inv(g\dualdag)$ of strata  $g\dualdag$ of $\kD(X^{(k)})$. \textit{For simplicity we will omit discussion of framing here, until our more detailed discussion in \autoref{sec:CWcomplexes}.} 

Setting $M = S^k$, this completes the construction of the stratification $\kP(f)$ of $S^k$.

We remark that the classical Thom-Pontryagin construction \cite{pontrjagin2007smooth}
\begin{equation}
\kP : [S^k,S^{k-l}] \iso \Set{\text{framed $l$-manifolds in $S^k$ up to cobordism}} = \Omega^{\mathrm{fr}}_l(S^k)
\end{equation}
is in fact a special case of the generalised construction, when one chooses standard CW-structure for $S^{k-l}$. %

\item \textit{(Inductively) adding morphisms}: Note that
\begin{equation}
D^{k+1} \iso \frac{S^k \times [0,1]}{S^k \times \Set{1}} = \cone(S^k)
\end{equation}
We build a stratification of $D^{k+1}$. Using the above identification, we can regard $f\inv(g\dualdag) \times [0,1)$ as a (framed) submanifold of $D^{k+1}$. The vertex point $S^k \times \Set{1}$ of $\cone(S^k)$ is not contained in any of those strata and becomes its own $0$-dimensional stratum with name $f$. In this way we have constructed a stratification $\cone(\kP(f))$ of $D^{k+1}$. To complete the inductive step of defining the dual stratification, we define the dual stratification of $\kD(X^{(k)} \cup_f D^{k+1})$ to be obtained as a gluing of the (inductively assumed) stratification $\kD(X^{(k)})$ of $X^{(k)}$ and the (newly constructed) stratification $\cone(\kP(f))$ of $D^{k+1}$. Note that they coincide along the gluing of $S^k$ by definition of $\kP(f)$, which makes this gluing of strata possible. Also note that this process can be used to inductively extend the stratification to all $k$-cells $f_1, f_2, f_3, ...$, yielding the dual stratification of $X^{(k)} \cup_{f_1} D^k \cup_{f_2} D^k \cup_{f_3} ... ~ $ which we denote by $\kD(X^{({k+1})})$. This completes the inductive construction of the dual stratifications $\kD(X^{(k)})$, from which $\kD(X)$ is obtained as a colimit.

\begin{notn}[Strata names] \label{notn:strata_name} Abusing notation, we will usually denote both $f\inv(g\dualdag) \subset S^k$ and $f\inv(g\dualdag) \times [0,1) \subset D^{k+1}$ by $g\dualdag$. The respective meaning of $g\dualdag$ can be inferred from the space it is embedded in.
\end{notn}

We give an example of the construction of the dual stratification. Let $X$ be the torus, whose $k$-skeleta $X^{(k)}$ (for $k=0,1,2$) are given as follows . 
\begin{restoretext}
\begingroup\sbox0{\includegraphics{test/page1.png}}\includegraphics[clip,trim=0 {.15\ht0} 0 {.15\ht0} ,width=\textwidth]{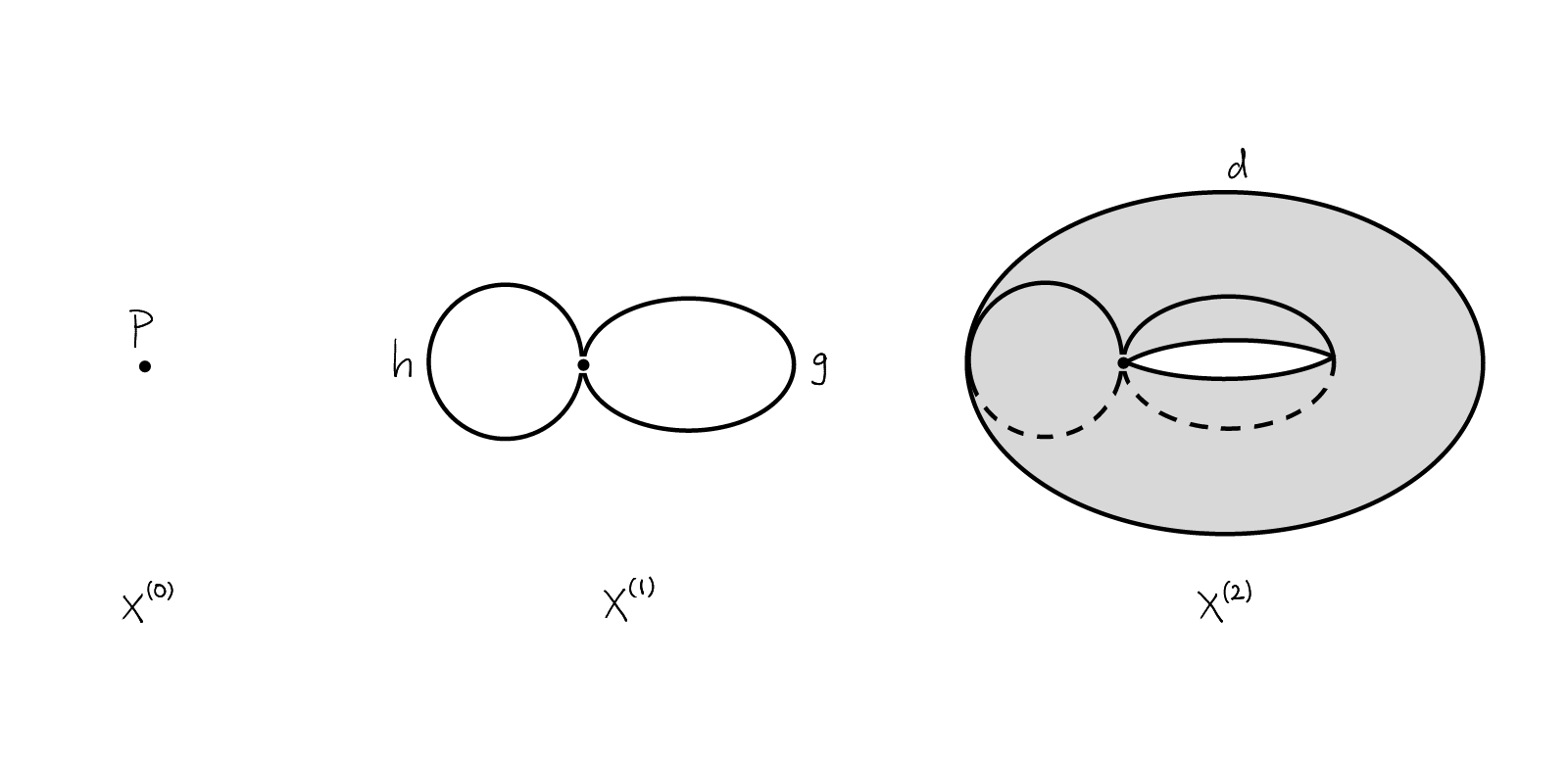}
\endgroup\end{restoretext}
We follow the inductive procedure above to find the dual complex $\kD(X)$. As stated, in $\kD(X^{(0)})$ $0$-cells $q$ equal their dual strata $q\dualdag$. The torus has a single $0$-cell $p$ and thus $\kD(X^{(0)})$ is the stratification with a single $0$-dimensional stratum $p\dualdag$
\begin{restoretext}
\begingroup\sbox0{\includegraphics{test/page1.png}}\includegraphics[clip,trim=0 {.3\ht0} 0 {.3\ht0} ,width=\textwidth]{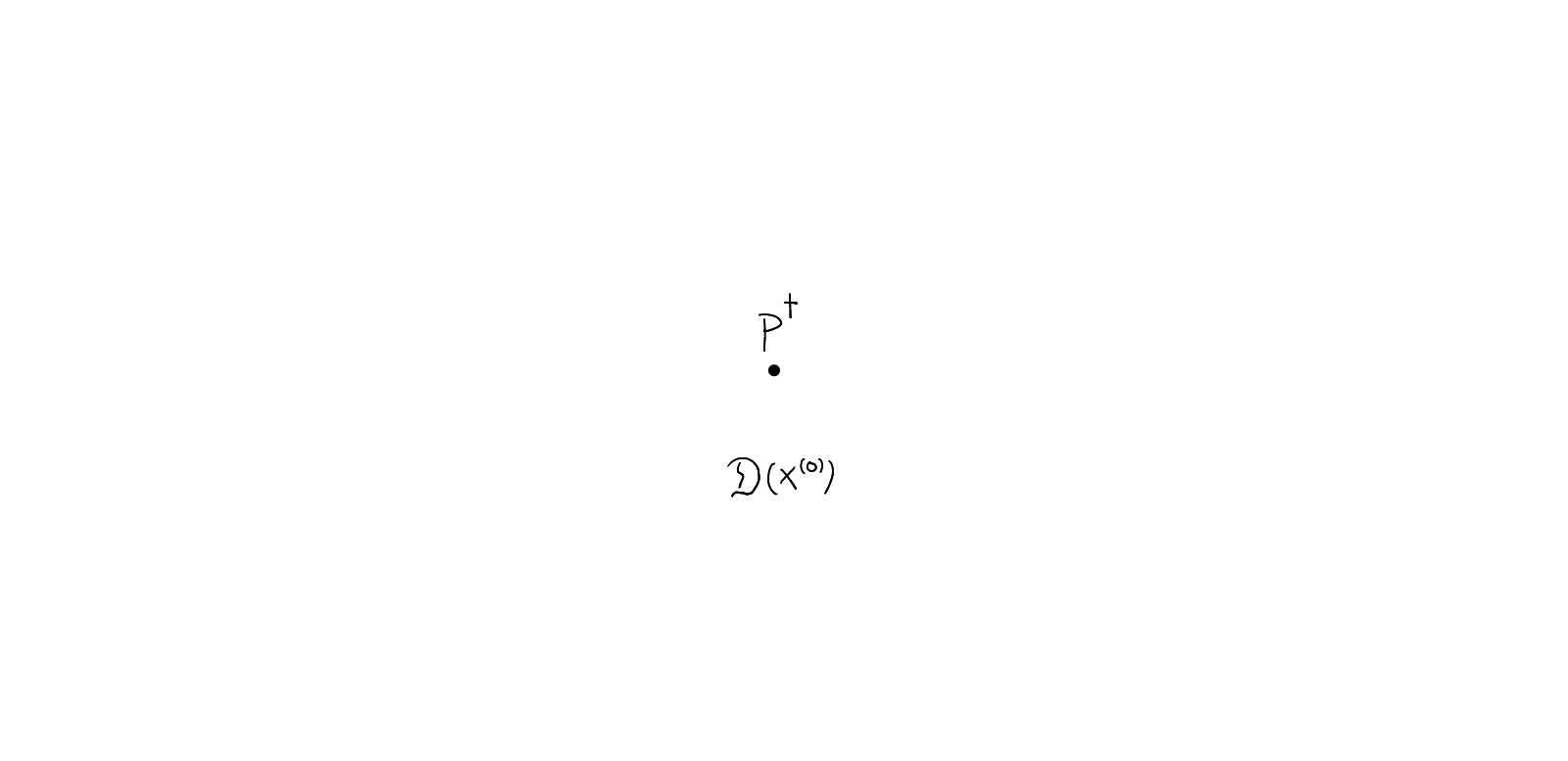}
\endgroup\end{restoretext}
Now, $h$ is attaching a $1$-cell $D^1$ by mapping its boundary $S^0$ to $X^{(0)}$ as shown on the left below
\begin{restoretext}
\begingroup\sbox0{\includegraphics{test/page1.png}}\includegraphics[clip,trim=0 {.25\ht0} 0 {.15\ht0} ,width=\textwidth]{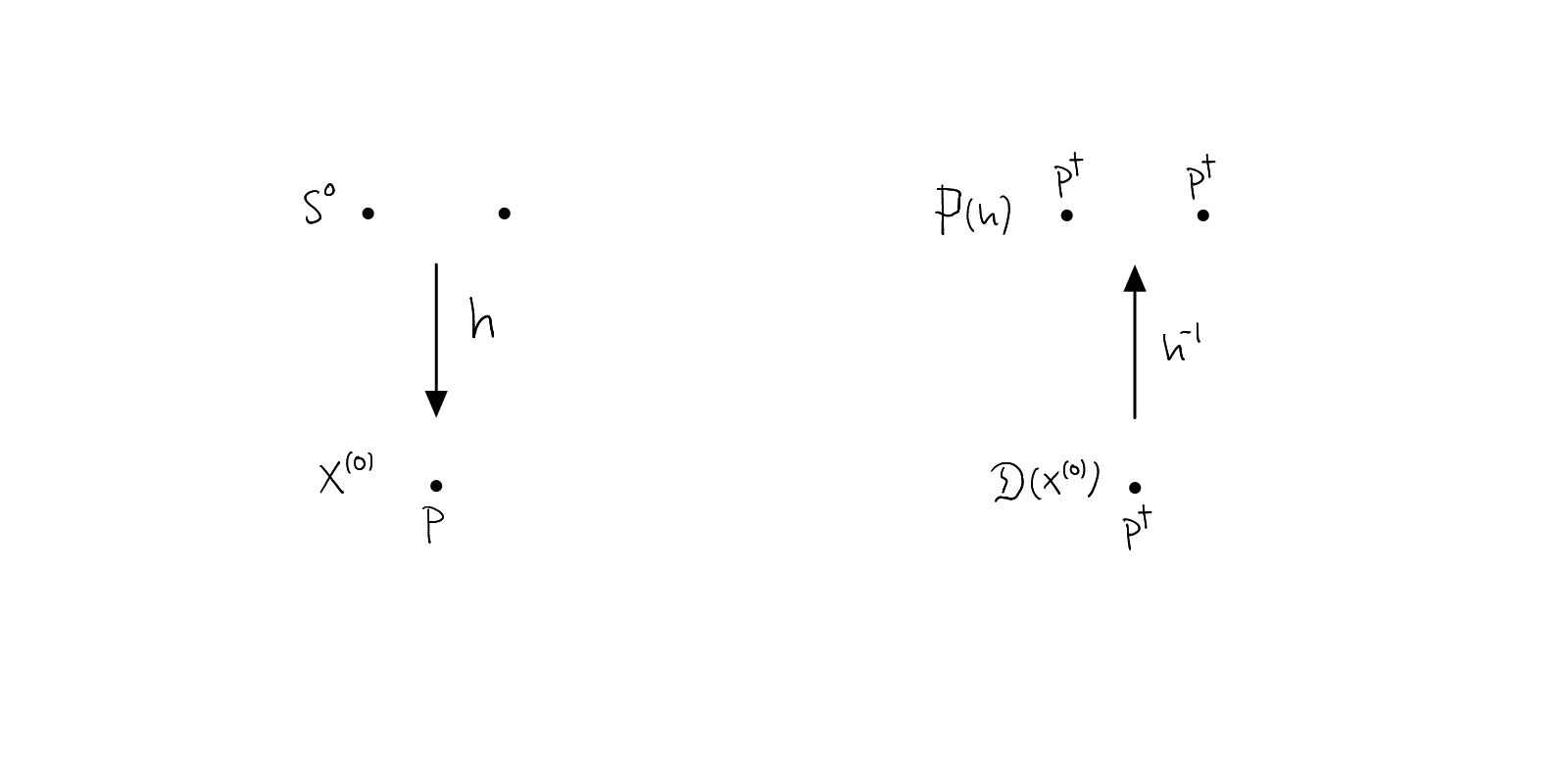}
\endgroup\end{restoretext}
On the right we depict the stratification $\kP(h)$ of $S^0$ obtained by taking inverse images of the strata of $\kD(X^{(0)})$. The resulting stratification $\kP(h)$ contains a single stratum $h\inv(p\dualdag) \equiv p\dualdag$ (cf. \autoref{notn:strata_name}) consisting of two disjoint points as shown above. Now, forming the cone of this gives a stratification $\cone(\kP(h))$ of the $1$-disk $D^1$ as follows
\begin{restoretext}
\begingroup\sbox0{\includegraphics{test/page1.png}}\includegraphics[clip,trim=0 {.3\ht0} 0 {.15\ht0} ,width=\textwidth]{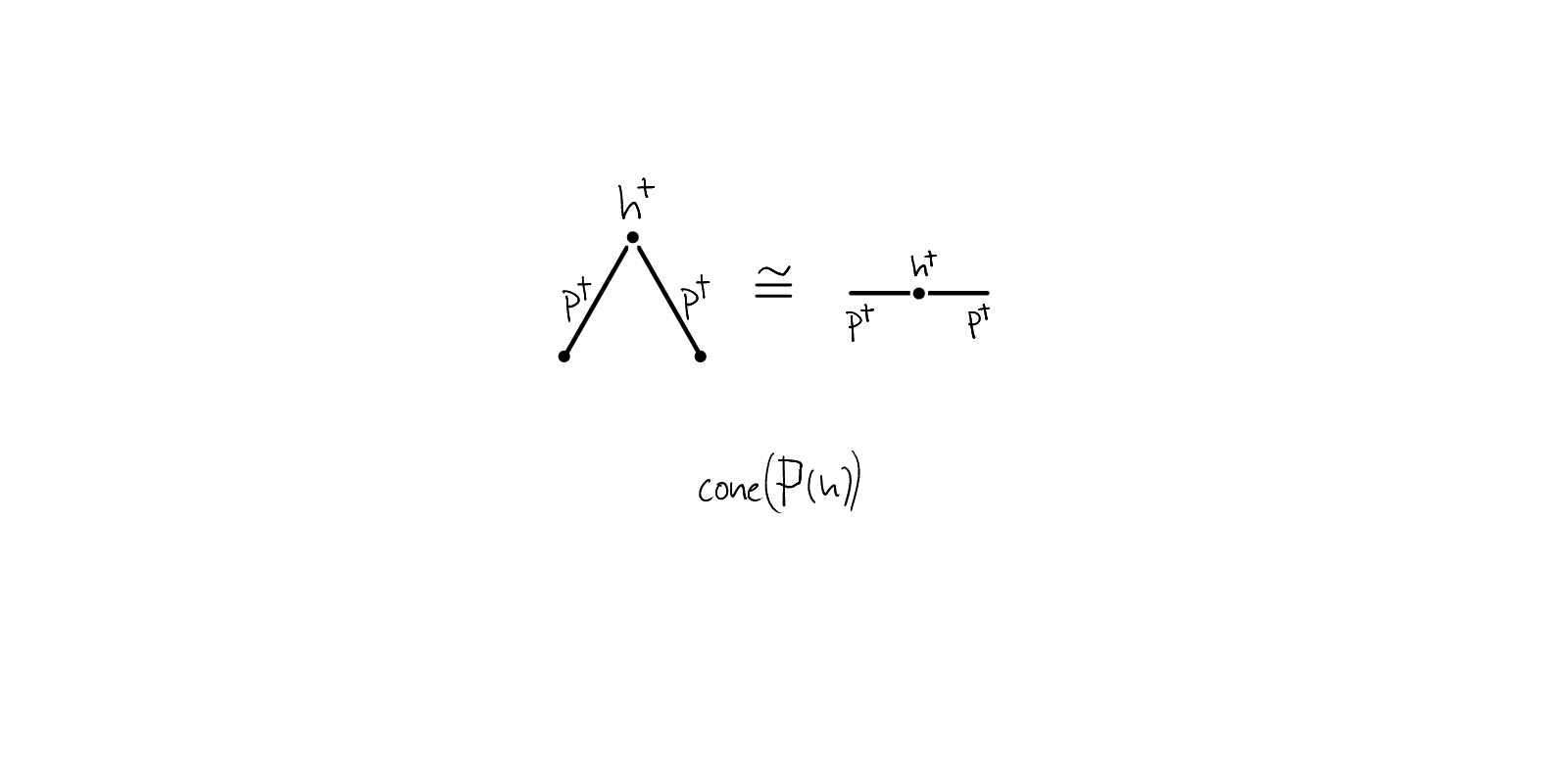}
\endgroup\end{restoretext}
Here, we gave two equivalent depictions of the stratification, illustrating the coning process on the left, and the the resulting stratified $1$-disk and the right. Note how a new $0$-dimensional stratum called $h\dualdag$ is added at the vertex point of the cone.

Now, the obtained stratified $1$-disk can be glued to $X^{(0)}$, extending the stratification $\kD(X^{(0)})$, which yields $\kD(X^{(0)} \cup_h D^1)$ as follows
\begin{restoretext}
\begingroup\sbox0{\includegraphics{test/page1.png}}\includegraphics[clip,trim=0 {.3\ht0} 0 {.3\ht0} ,width=\textwidth]{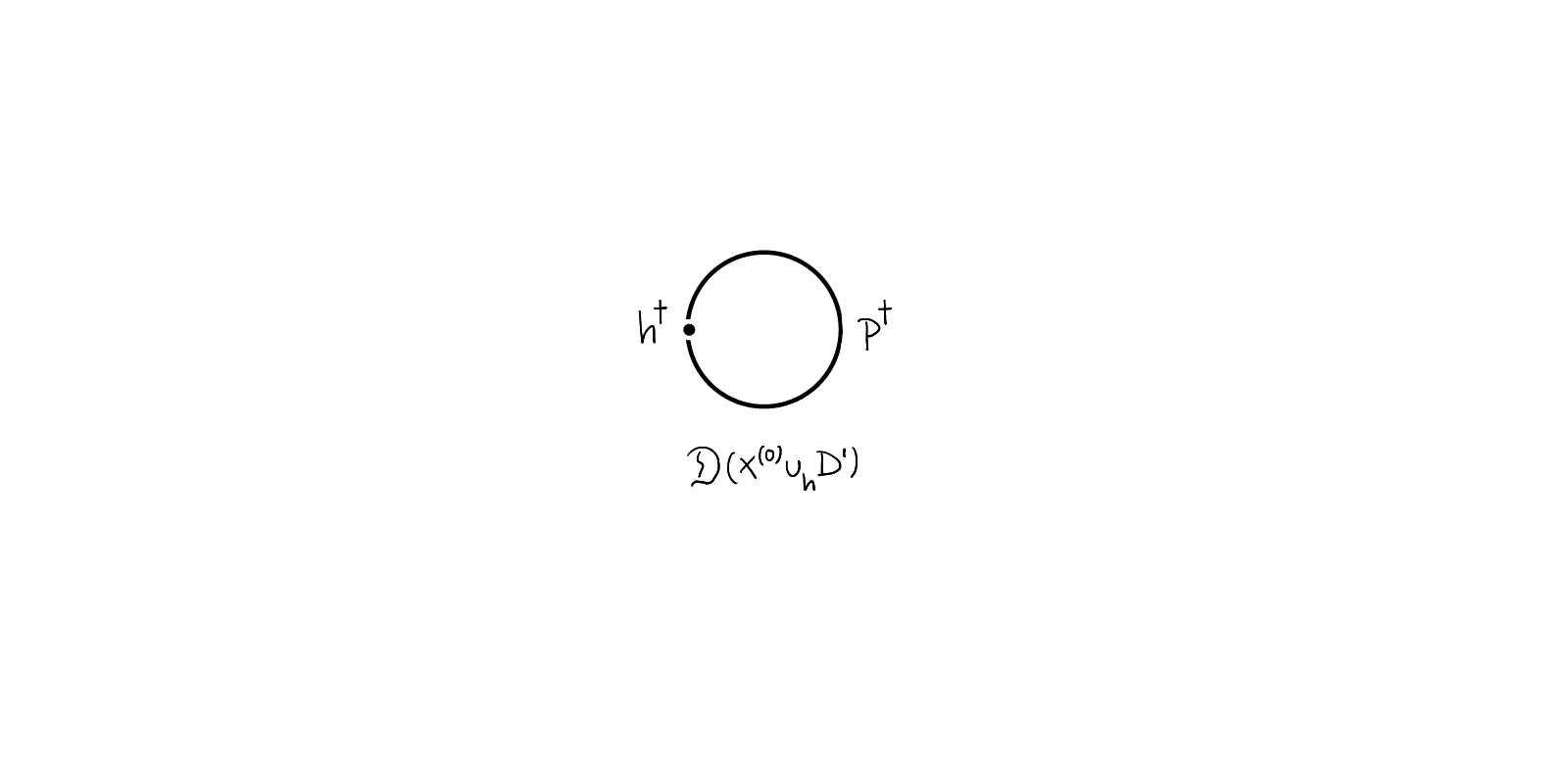}
\endgroup\end{restoretext}
Similarly, $\cone(\kP(g))$ is given by
\begin{restoretext}
\begingroup\sbox0{\includegraphics{test/page1.png}}\includegraphics[clip,trim=0 {.3\ht0} 0 {.15\ht0} ,width=\textwidth]{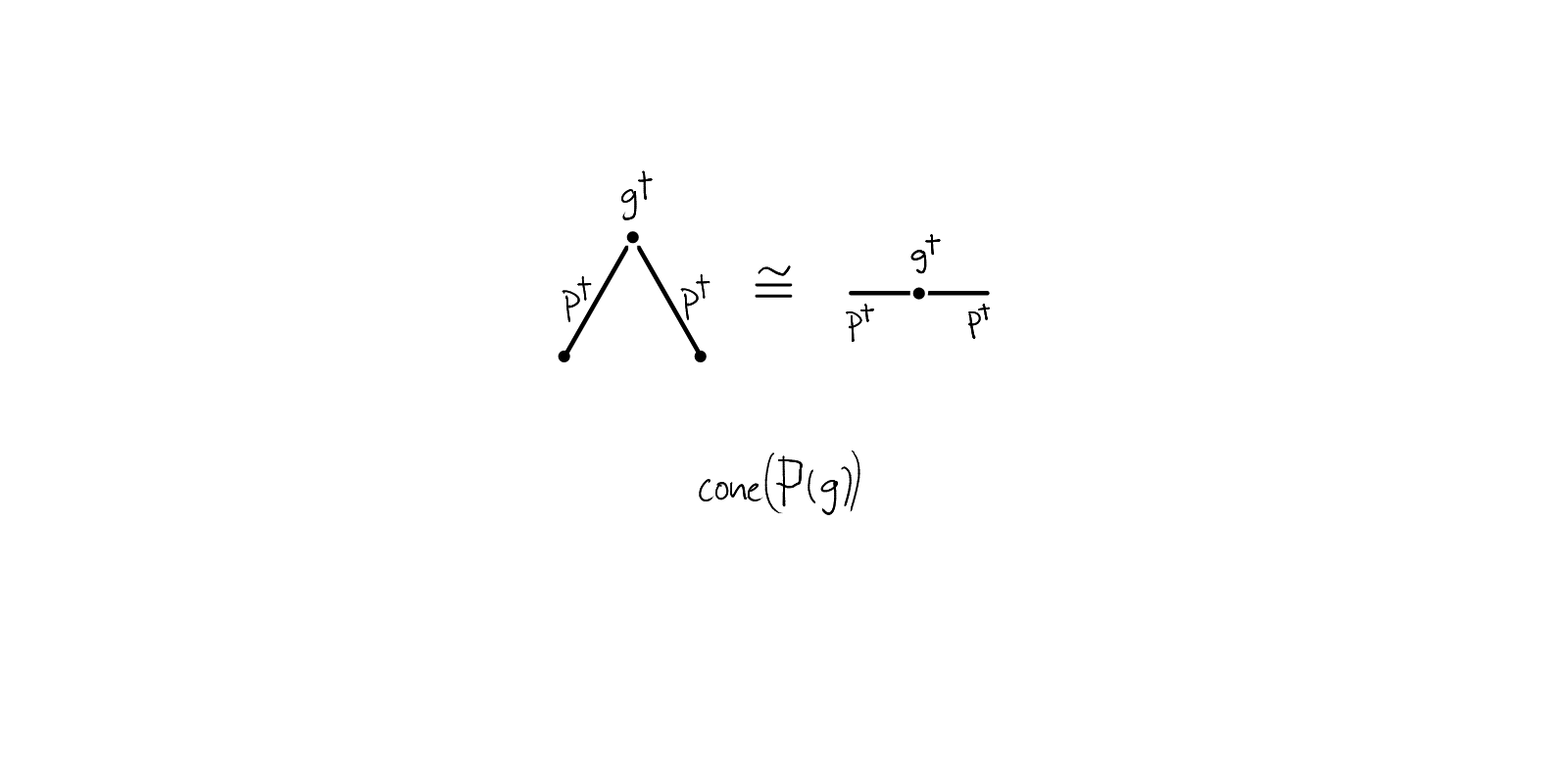}
\endgroup\end{restoretext}
and combining the strata by gluing we find $\kD(X^{(0)} \cup_h D^1 \cup_g D^1)$
\begin{restoretext}
\begingroup\sbox0{\includegraphics{test/page1.png}}\includegraphics[clip,trim=0 {.25\ht0} 0 {.3\ht0} ,width=\textwidth]{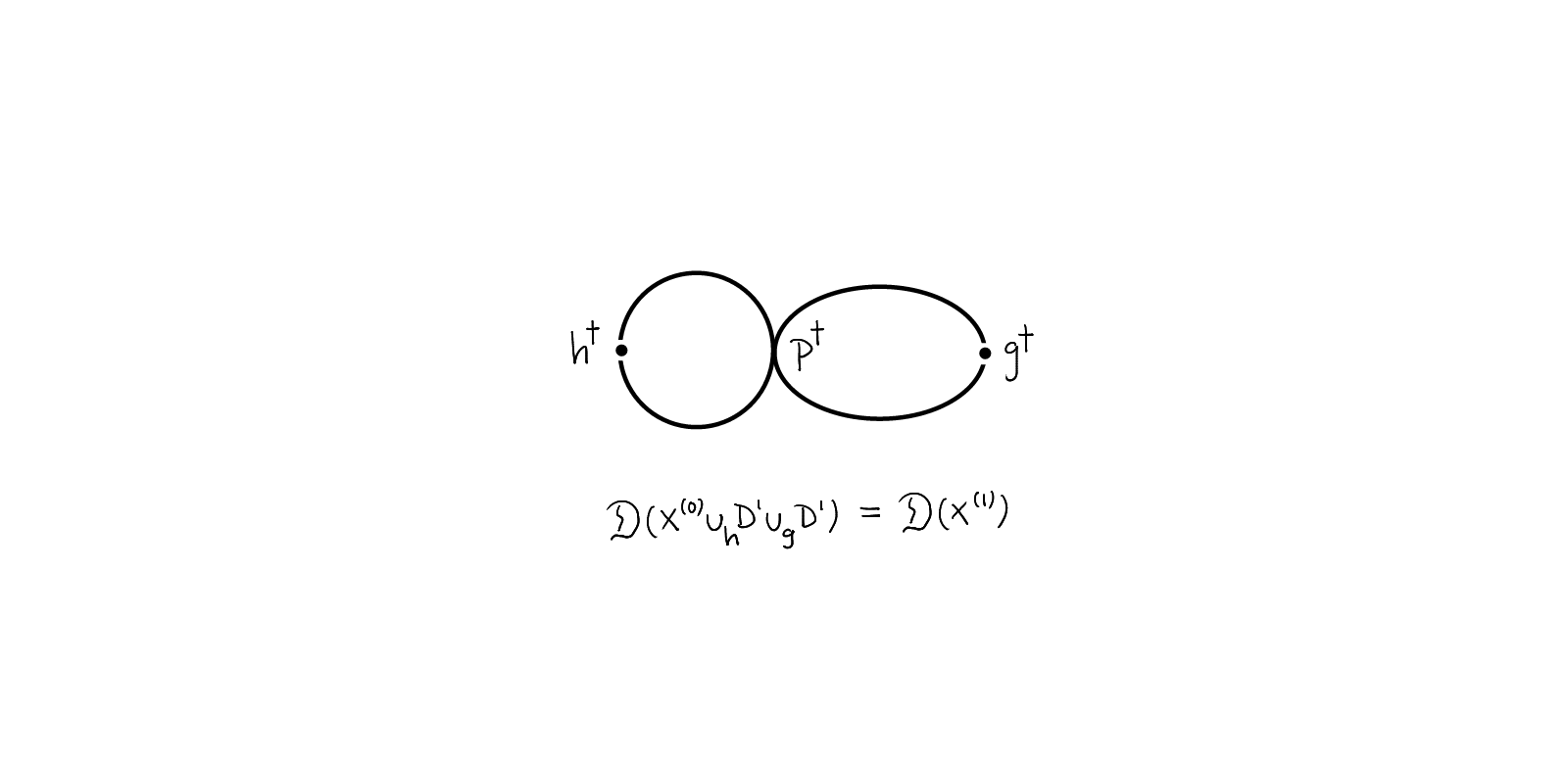}
\endgroup\end{restoretext}
which equals $\kD(X^{(1)})$. 

In the next inductive step we consider the attaching map $d$, which maps $S^1$ to $X^{(2)}$ as follows: starting at $p$, we first loop around $h$, then loop around $g$, then loop (in the other direction) around $h$ and finally loop (in the other direction) around $g$ ending again at $p$. $\kP(d)$ is obtained as the stratification of $S^1$ given by inverse images of strata in $\kD(X^{(1)})$. This is indicated on the right below
\begin{restoretext}
\begingroup\sbox0{\includegraphics{test/page1.png}}\includegraphics[clip,trim=0 {.0\ht0} 0 {.0\ht0} ,width=\textwidth]{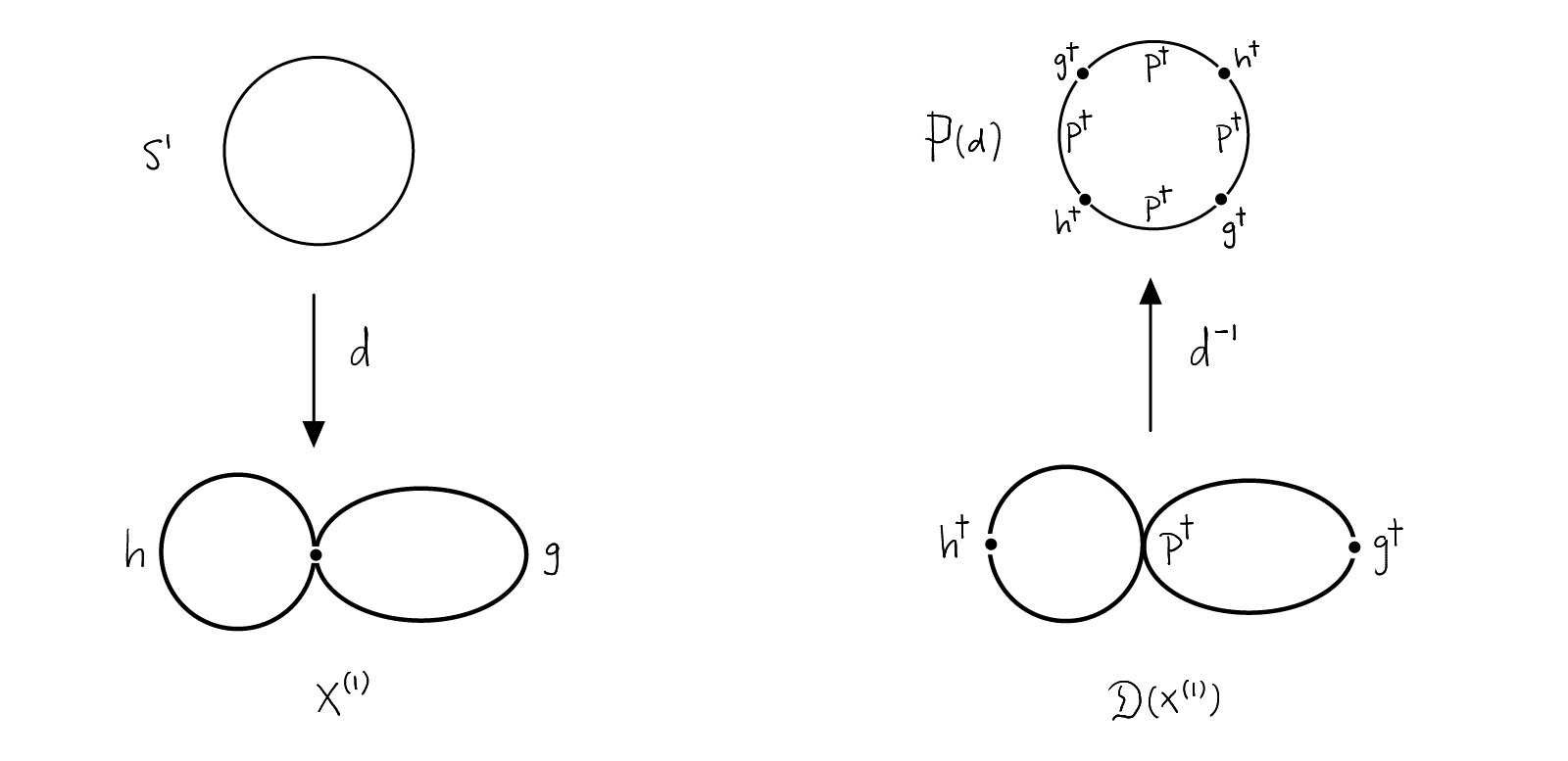}
\endgroup\end{restoretext}
Now $\cone(\kP(d))$ is obtained from this stratification $\kP(d)$ by extending strata all the way up to the vertex point of the cone (but excluding the vertex point), which yields
\begin{restoretext}
\begingroup\sbox0{\includegraphics{test/page1.png}}\includegraphics[clip,trim=0 {.25\ht0} 0 {.15\ht0} ,width=\textwidth]{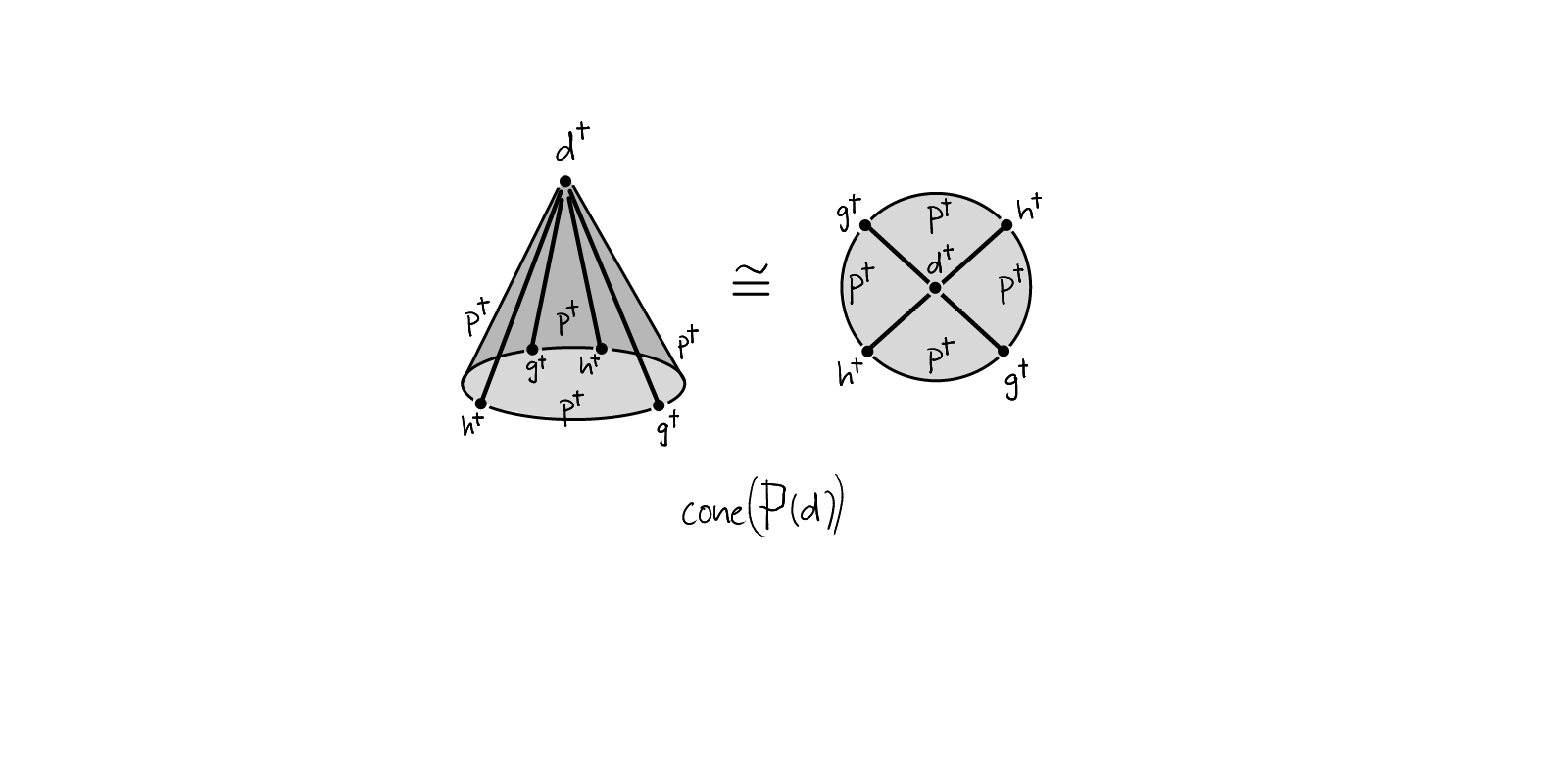}
\endgroup\end{restoretext}
Note how a new $0$-dimensional stratum called $d\dualdag$ is added at the vertex point of the cone. We can now glue the stratification $\cone(\kP(h))$ of the $2$-disk $D^2$ on $X^{(2)}$ combining their stratifications. This yields a stratification $\kD(X^{(1)} \cup_d D^2)$ of $X^{(1)} \cup_d D^2 = X^{(2)}$ as follows
\begin{restoretext}
\begingroup\sbox0{\includegraphics{test/page1.png}}\includegraphics[clip,trim=0 {.15\ht0} 0 {.25\ht0} ,width=\textwidth]{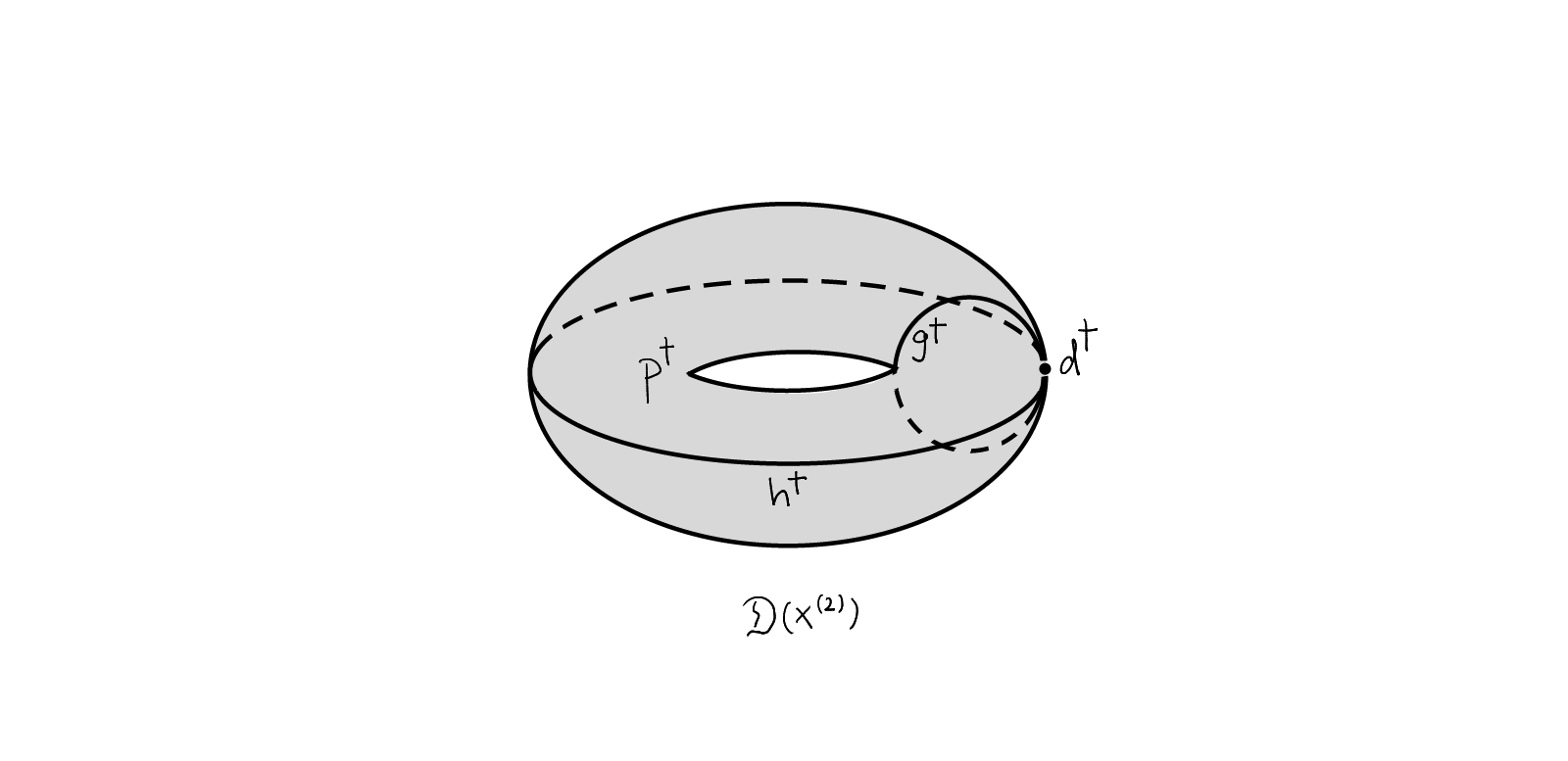}
\endgroup\end{restoretext}
Thus, $\kD(X^{(2)})$ contains four strata in total: a 2-dimensional stratum $p\dualdag$, two $1$-dimensional strata $g\dualdag$ and $h\dualdag$, and a $0$-dimensional stratum $d\dualdag$.

This completes our example, and we now turn to the last step of the translation.

\item \textit{Choosing direction}: We now want to add a notion of direction to the $k$-disk. The $n$-cube has a notion of direction by being flag-foliated (cf. \autoref{ssec:po_mfld_diag}). We want to ``push" this flag-foliation onto the disk now. 

A \textit{globular foliation} of the $k$-disk is a map from the (open) $k$-cube to the (open) $k$-disk, extending to a map from the closed cube to the closed disk, which maps the boundary of the $k$-cube to the boundary of the $k$-disk, quotienting certain sides of the cube in the usual sense of globularity (the sides of the $n$-cube facing in direction $k$ map to the $(n-k+1)$-fold source or target of the $n$-globe). In more precise terms, a globular foliation $F$ is a homeomorphism required to extend to a map $F_n : [0,1]^n \to D^n$ of the form (for $p \in \Set{0,1}$, $k \in \bnum n$)
\begin{equation}
(x,p,y) \in [0,1]^{k} \times [0,1] \times [0,1]^{n-k-1} \mapsto F_n(x,p,y) = F_k(x) \in D^n
\end{equation}
where $F_k$, $k \in (\bnum {n+1})$, restrict to homeomorphisms on $(0,1)^k$, and $F_k((0,1)^k)$ together disjointly cover $D^n$. In yet other words, we identify $D^n$ with an $n$-globe\footnote{Inductively, the $n$-globe is an $n$-disk $D^n$ whose boundary $S^{n-1}$ is the union of a ``source" and ``target" hemisphere, which are themselves $(n-1)$-globes, and the union is such that the source (resp. target) of the source coincides with the source (resp. target) of the target. (cf. \cite{leinster-operads})} and then $F_n$ takes the lower and upper side of the $n$-cube facing in the $k$th direction, and maps them to the $(n-k+1)$-fold source respectively target of the $n$-globe. 

For $k = 1$, a globular foliation (up to homotopy) means choosing an orientation of the $1$-disk, i.e., up to homotopy there are two globular foliations $\id_{\pm} : (0,1) \to D^1$. For $k = 2$, we give three examples $F_1$, $F_2$, $F_3$ of globular foliations of the $2$-disk in the following three pictures
\begin{restoretext}
\begingroup\sbox0{\includegraphics{test/page1.png}}\includegraphics[clip,trim=0 {.25\ht0} 0 {.25\ht0} ,width=\textwidth]{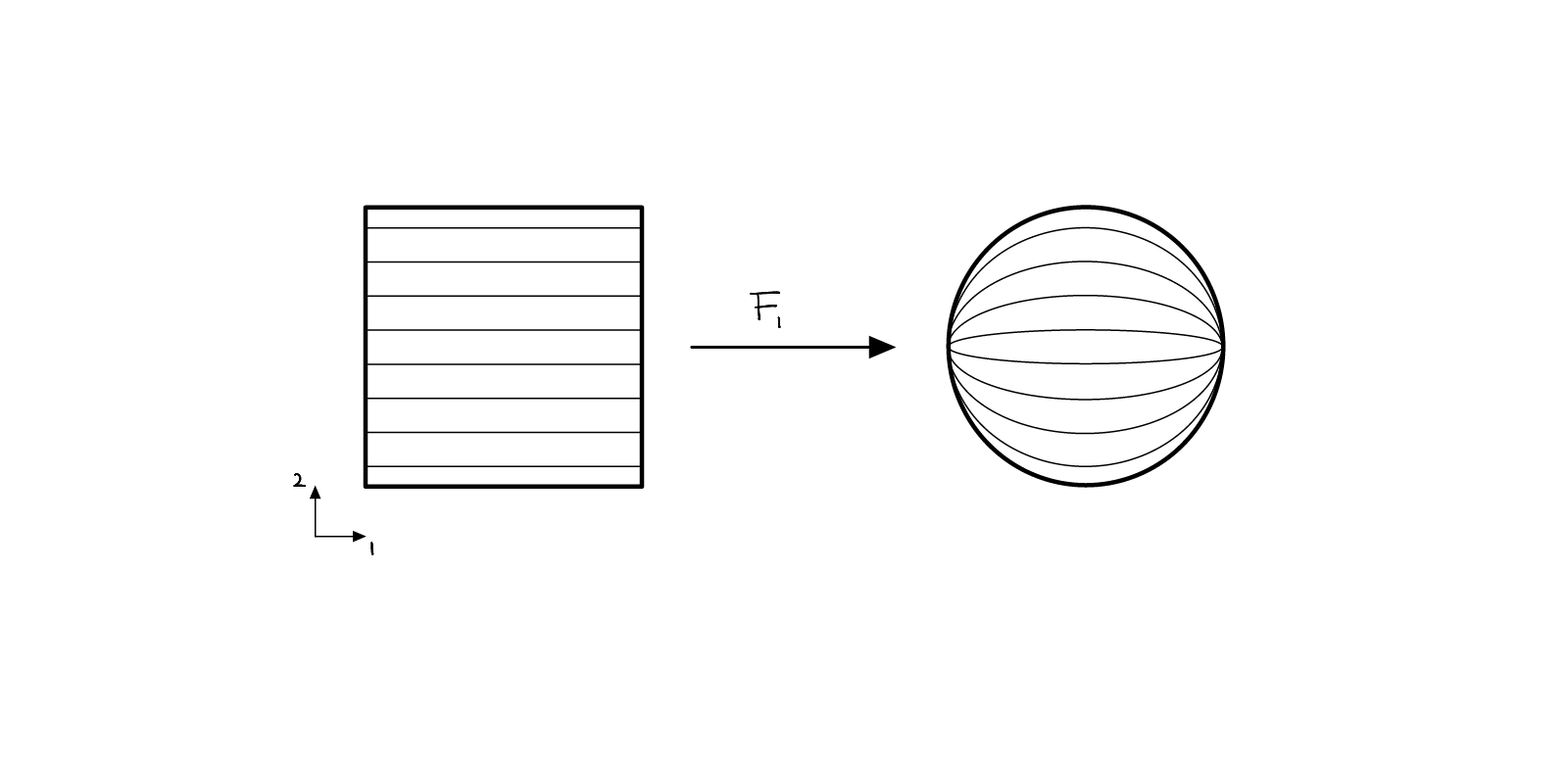}
\endgroup\end{restoretext}
\begin{restoretext}
\begingroup\sbox0{\includegraphics{test/page1.png}}\includegraphics[clip,trim=0 {.24\ht0} 0 {.26\ht0} ,width=\textwidth]{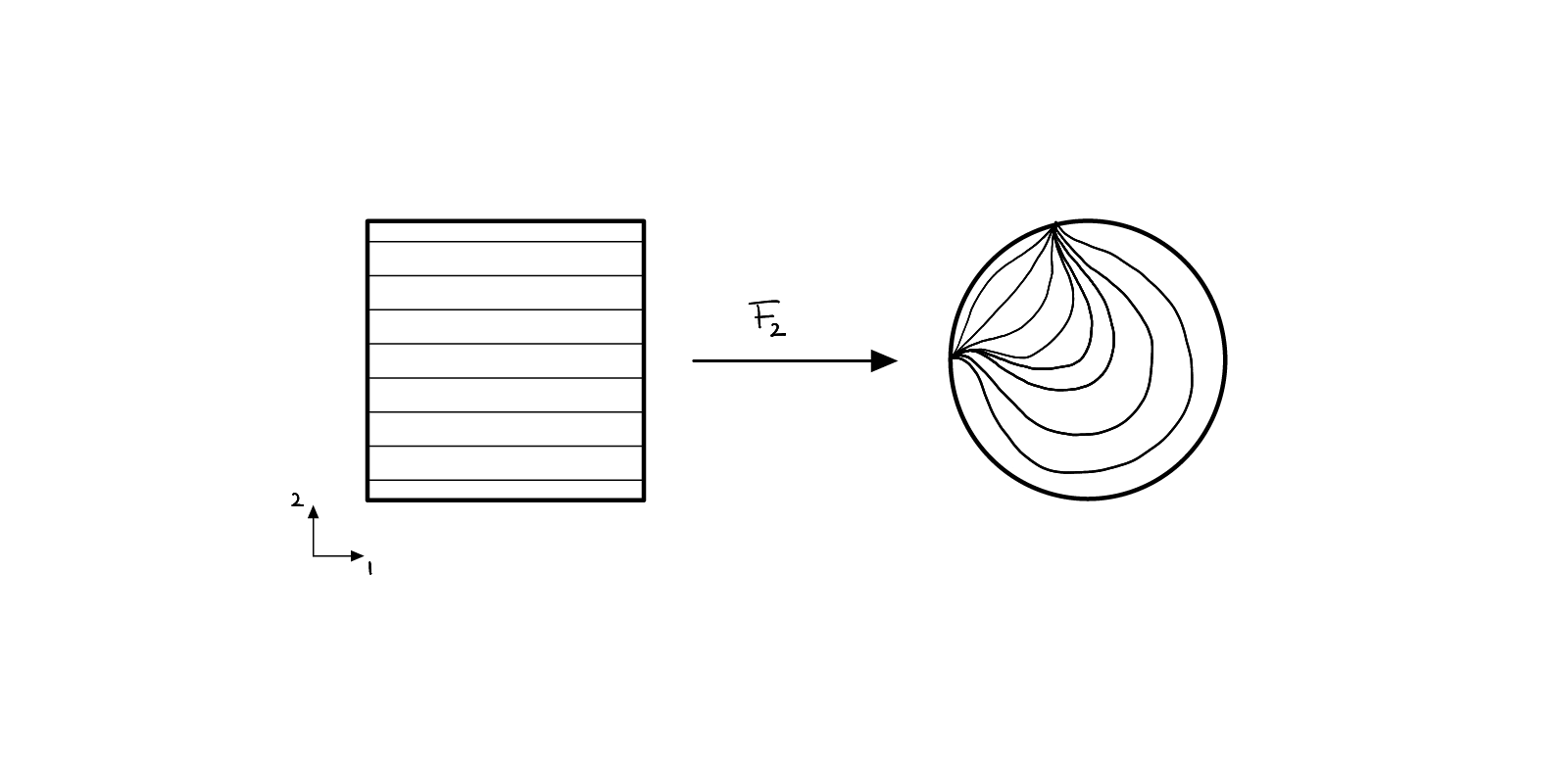}
\endgroup\end{restoretext}
\begin{restoretext}
\begingroup\sbox0{\includegraphics{test/page1.png}}\includegraphics[clip,trim=0 {.24\ht0} 0 {.26\ht0} ,width=\textwidth]{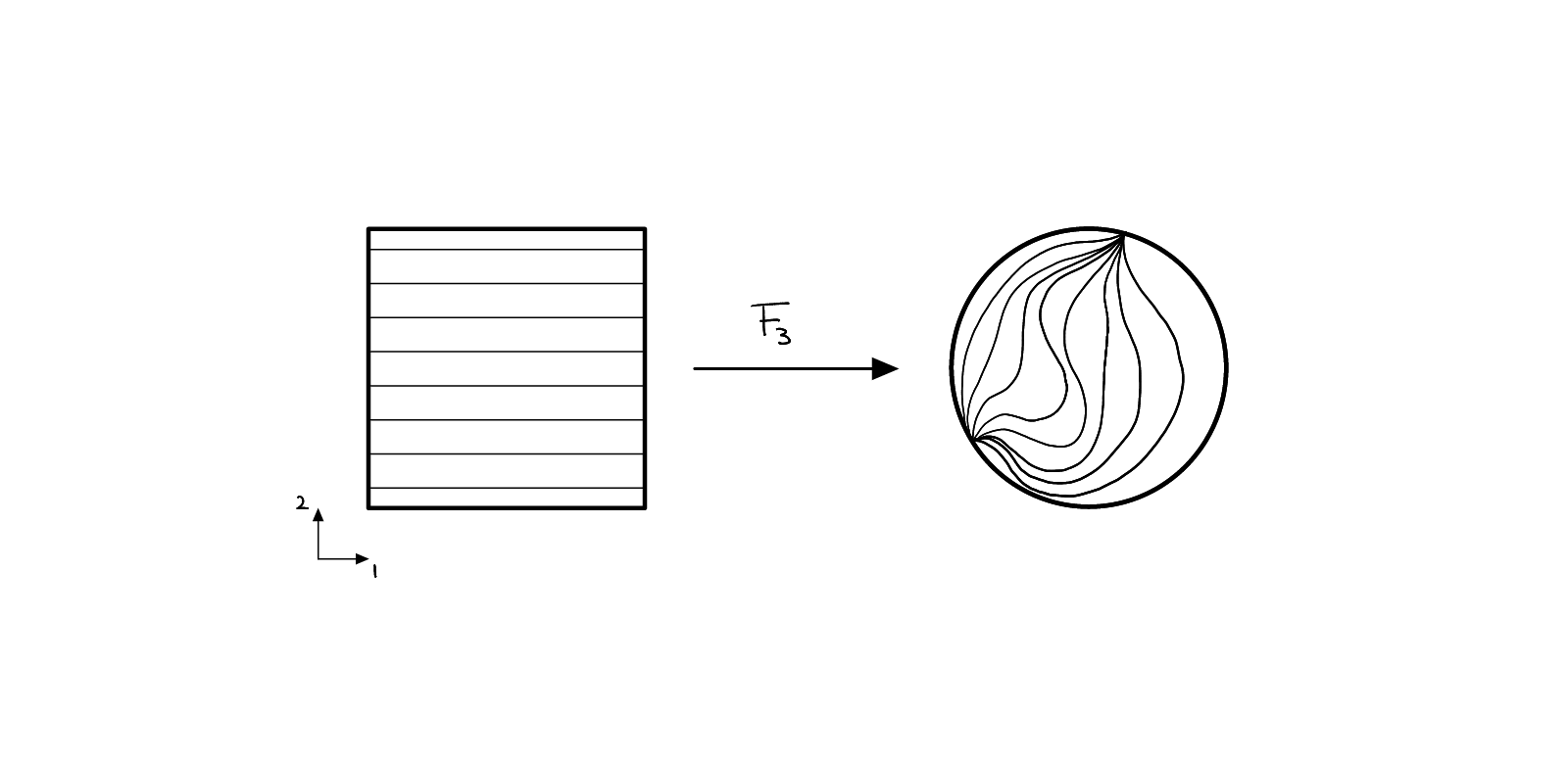}
\endgroup\end{restoretext}
Above, we depicted images of certain sheets (marked by lines) in the $2$-cube by lines in the $2$-disk. 

Now, let $\cone(\kP(f))$ be a (framed) stratification of the $k$-disk $D^k$. We say $F : (0,1)^k \to D^k$ is compatible with the cell attached via $f$ in the CW-complex $X$, if $F$ pulls back strata of $\cone(\kP(f))$ to strata in $(0,1)^k$ which, after possibly refining the latter strata using finitely many new strata, yields a manifold diagram (or more precisely, a $\GGamma{}{X}$-labelled manifold diagram, cf. \autoref{notn:sum_namescopes}, where $X_i$ is the set of $i$-cells in $X$).

We give examples. All of $F_1$, $F_2$, $F_3$ defined above are compatible with $d$ in the torus ($k = 2$). Indeed, pulling back we obtain the following stratification of the $2$-cube for $F_1$
\begin{restoretext}
\begingroup\sbox0{\includegraphics{test/page1.png}}\includegraphics[clip,trim=0 {.25\ht0} 0 {.25\ht0} ,width=\textwidth]{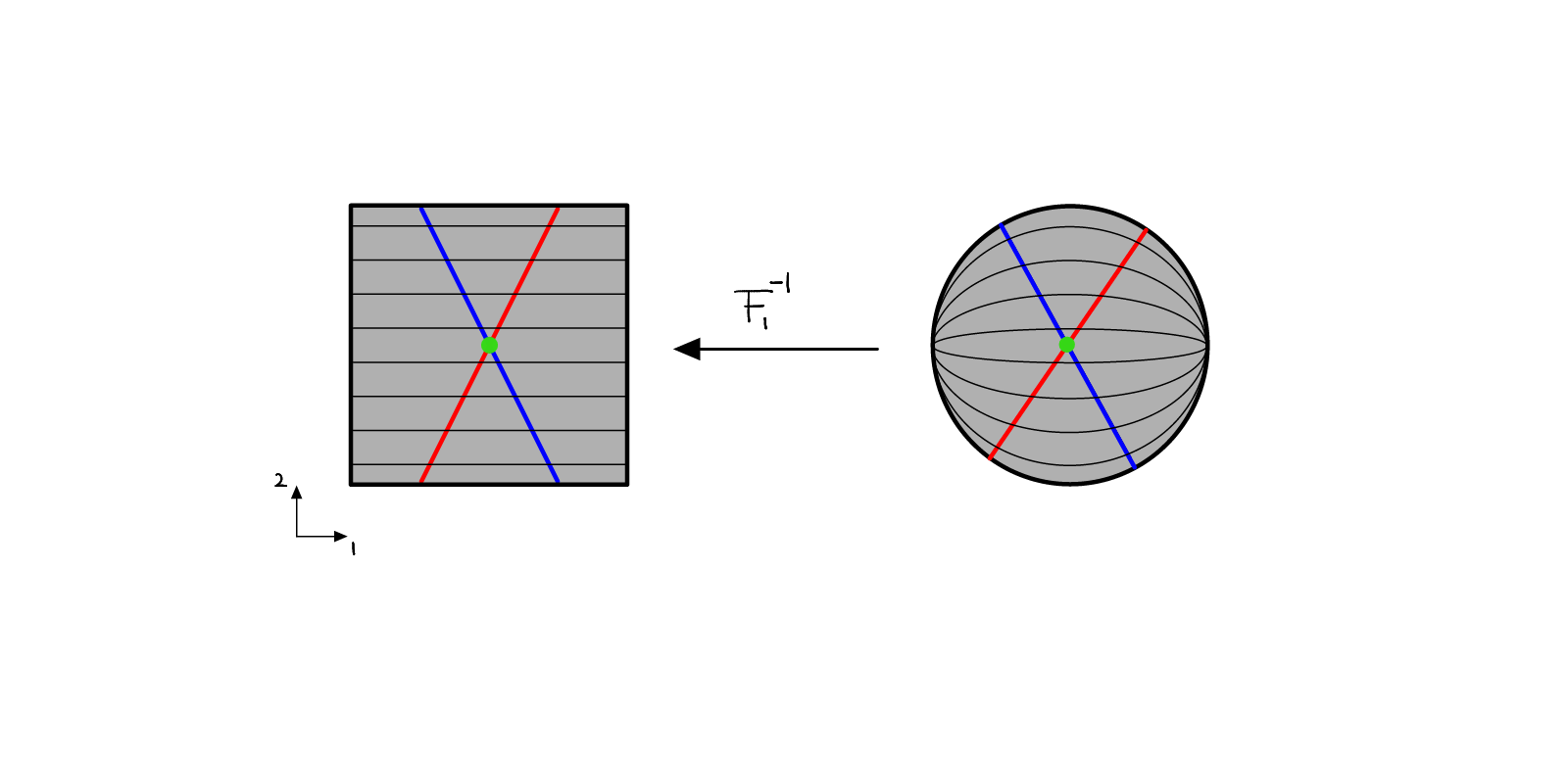}
\endgroup\end{restoretext}
For $F_2$ we find
\begin{restoretext}
\begingroup\sbox0{\includegraphics{test/page1.png}}\includegraphics[clip,trim=0 {.24\ht0} 0 {.26\ht0} ,width=\textwidth]{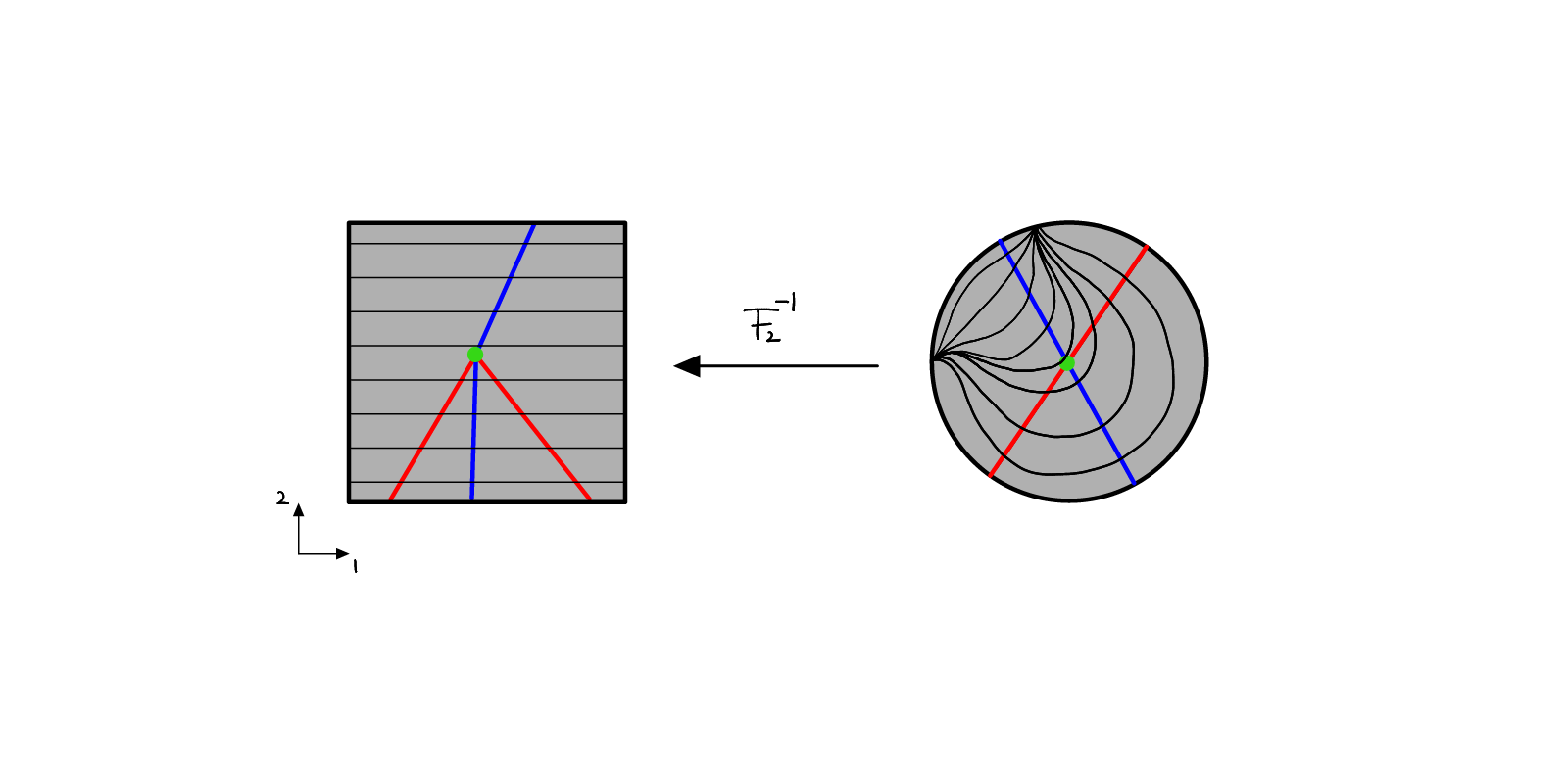}
\endgroup\end{restoretext}
and for $F_3$ we find
\begin{restoretext}
\begingroup\sbox0{\includegraphics{test/page1.png}}\includegraphics[clip,trim=0 {.25\ht0} 0 {.25\ht0} ,width=\textwidth]{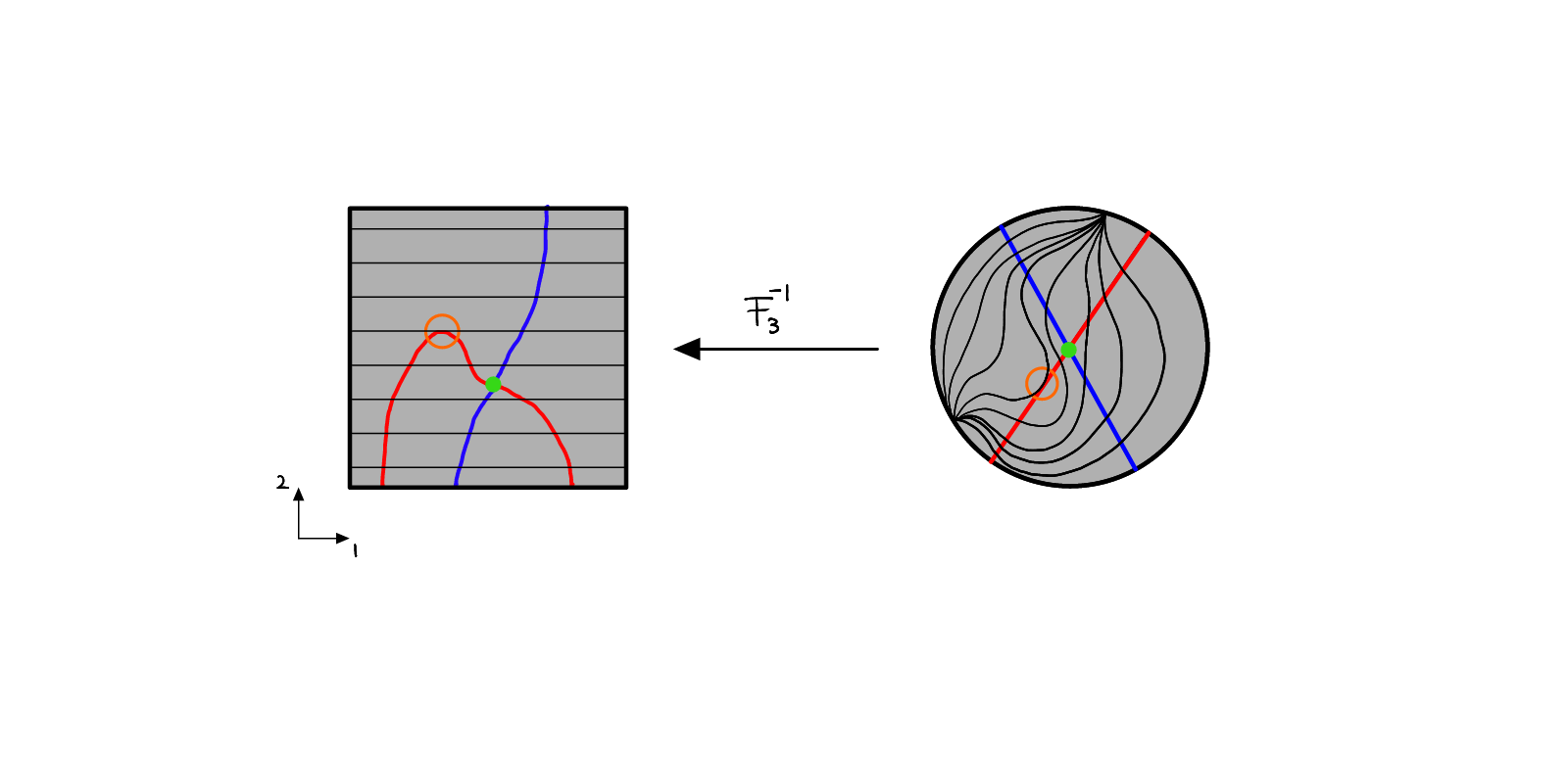}
\endgroup\end{restoretext}
This last example is special when compared to the others: the $l$-dimensional sheets are not everywhere  \textit{transversal} to $(k+1-l)$-dimensional strata of $\cone(\kP(d))$. This leads to a $0$-dimensional \textit{singularity} circled in \corange{}. Due to this singularity, the condition of \textit{transversality} in \autoref{ssec:po_mfld_diag} is not satisfied, unless we replace this singularity with its own $0$-dimensional stratum, thus refining the \cred{} stratum into $3$ disconnected components in total.

Importantly, in general such singularities cannot be avoided. However, a crucial but implicit conjecture is the following: in our coarse setting of ``\textit{directed topology}" these singularities can be enumerated, and more precisely, we claim they are enumerated by the theory of $\infty$-dualisability $\TI$ (for instance, the above singularity is the ``cap singularity" $\ic_{\abss{\ic_{-\equiv +}} \whisker 1 1 \abss{\ic_{+ \equiv -}} \equiv \Id_{-}} \in \TI_2$). This is in stark contrast with the finer \textit{classical setting of algebraic descriptions} of singularities, where such an enumeration is not possible. A formalisation and proof of this conjecture is left to future work.\\

To complete our translation from CW-complexes to $\infty$-groupoids, we now choose a globular foliation $f\Fol$ for each attaching map $f$ of a cell in $X$ such that $(f\Fol)\inv(\cone(\kP(f)))$ is itself of the form $\cone(S)$ for some stratification $S$ of the boundary of the $k$-cube. Such a stratification will be called conical (algebraically it translates to our previous definition of double cones). The choice of $f\Fol$ is non-unique but we claim it always exists. The resulting $(f\Fol)\inv(\cone(\kP(f)))$ is a manifold diagram with duals (that is, it might contain singularities classified by the theory of $\infty$-dualisability). Algebraically, this diagram now provides the type $\abss{f\dualdag}$ for a coherently invertible generating $(k+1)$-morphism $f\dualdag$ for the $\infty$-groupoid $\sX$ corresponding to $X$ (cf. \autoref{def:sum_infty_group}). Using this construction inductively, we are able to add invertible generating $(k+1)$-morphisms $f\dualdag$ to $\sX$ for each cell $f$ of $X$, thereby defining $\sX$ from $X$ as we set out to do.

We give an example. In the case of the torus $X$, choosing $p\Fol = \id$, $g\Fol = h\Fol = \id_+$, and $d\Fol = F_1$, we find a presented associative $\infty$-groupoid $\sX$ with the following invertible generators
\begin{restoretext}
\begingroup\sbox0{\includegraphics{test/page1.png}}\includegraphics[clip,trim=0 {.25\ht0} 0 {.25\ht0} ,width=\textwidth]{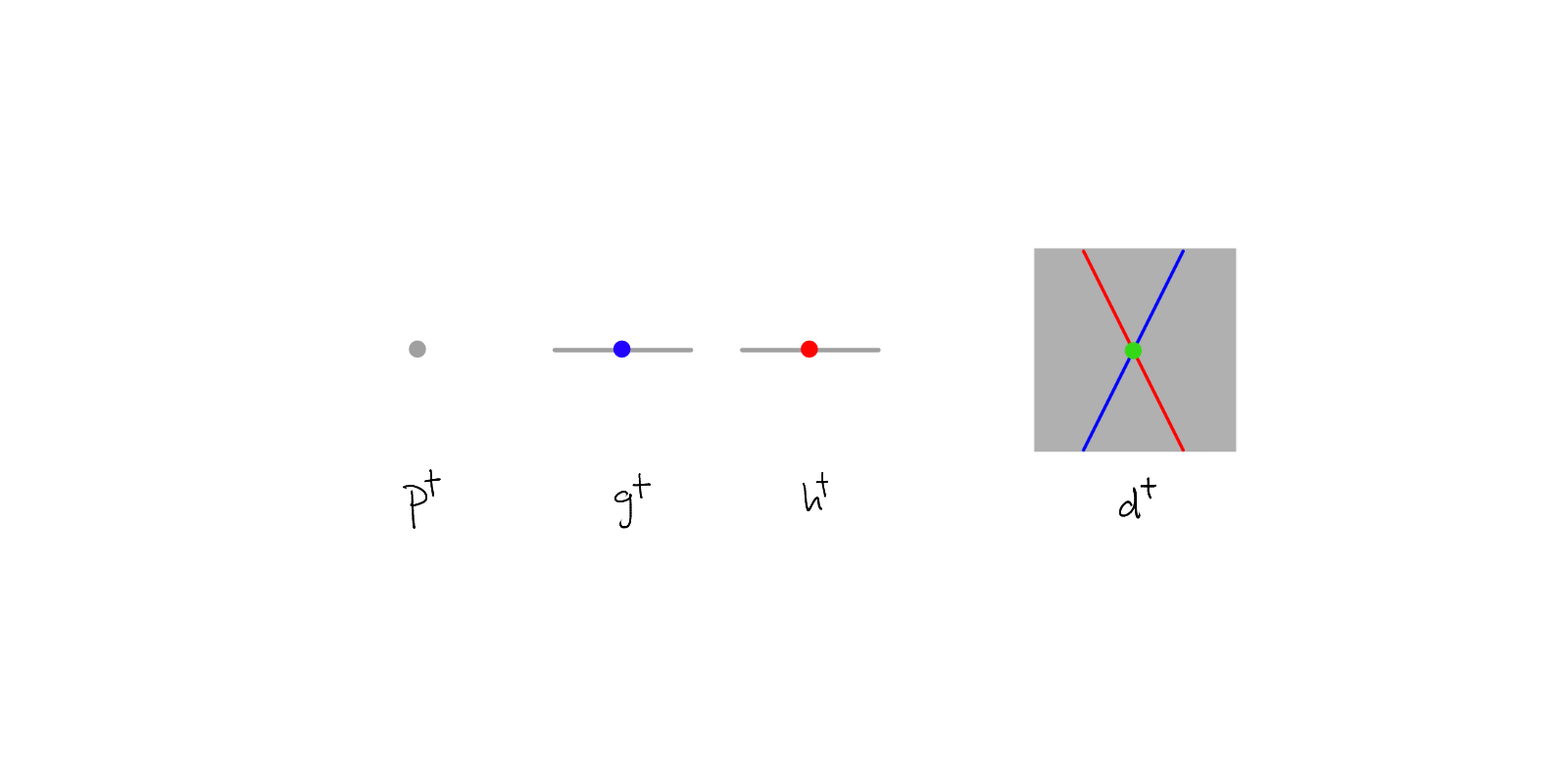}
\endgroup\end{restoretext}
Note that, if we had chosen $g\Fol = h\Fol = \id_-$ instead of $\id_+$, we would have obtained exactly the same groupoids but with $g$ and $h$ now switching roles with their inverses\footnote{In fact, since in this example there is only one object, $g,h$ and their respective inverses play exactly the same role. But if $g\dualdag : A \to B$ was a general morphism between objects $A$ and $B$, then, after changing the orientation of the foliation we would obtain a morphism $g\dualdag : B \to A$, now playing the role of its own former inverse.}, which are part of the coherent invertibility data for generators in a groupoid (cf \autoref{def:sum_infty_group}). 

More generally we make the following remark.

\begin{rmk}[Different choices of foliation yield equivalent groupoids] As mentioned above, types depend on foliation choice. Any two choices should yield equivalent groupoids, since making a different choice of foliation can be compensated for by the presence of coherent invertibility data from $\TI$. We will not attempt to make this precise but we illustrate the point in the following example case: if we had chosen $d\Fol = F_2$ instead of $F_1$ we would have obtained a different type $\abss{d\dualdag}$, and a different groupoid $\sX$, but the resulting higher groupoid would be equivalent in the sense that the new type appears as a morphism in the old type and vice versa, and thus again any morphism can be stated in either groupoid. For instance, the new groupoid's type for $d\dualdag$ (on the left) appears as the following morphism (on the right) in the old groupoid
\begin{restoretext}
\begingroup\sbox0{\includegraphics{test/page1.png}}\includegraphics[clip,trim=0 {.25\ht0} 0 {.23\ht0} ,width=\textwidth]{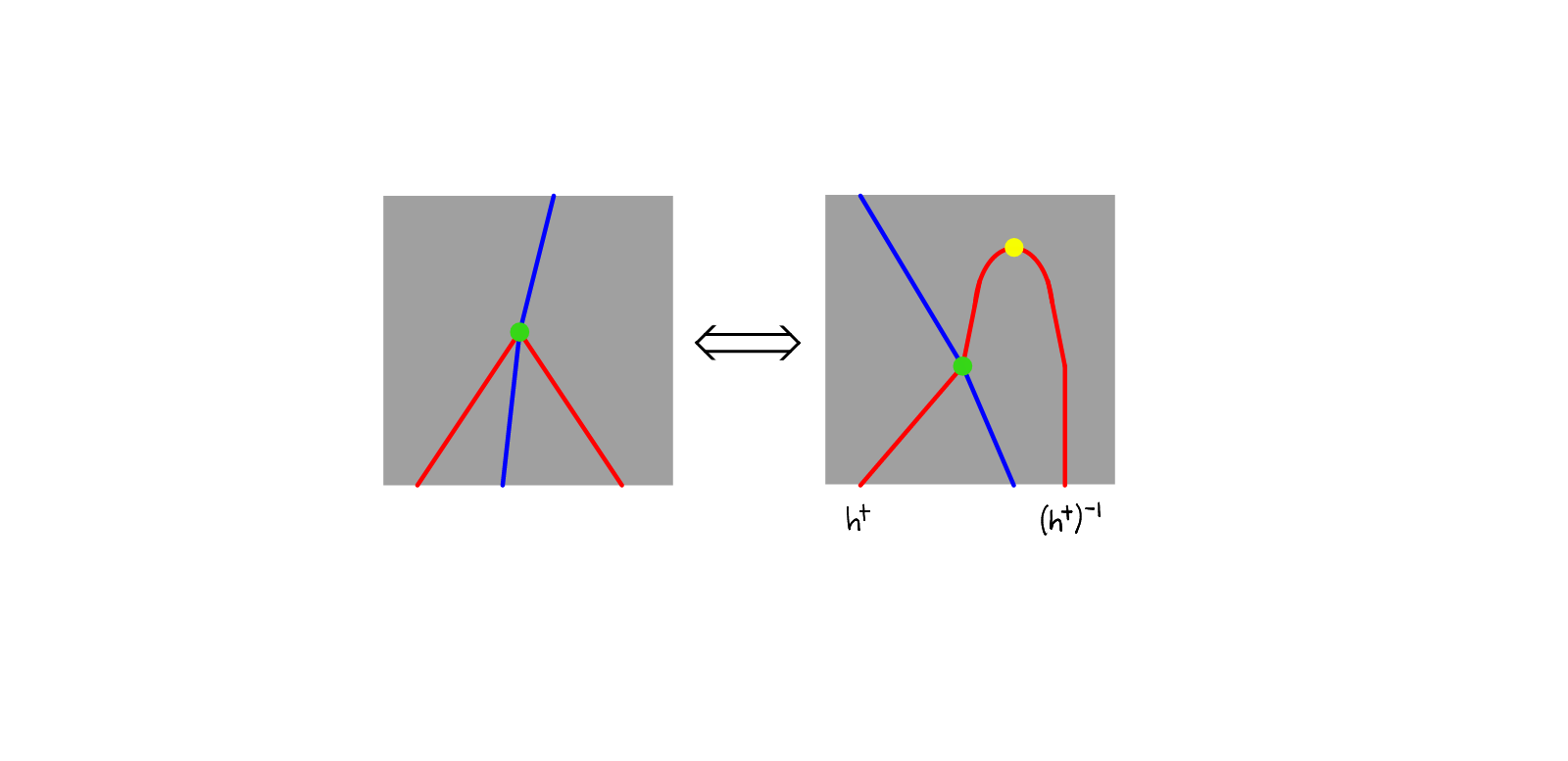}
\endgroup\end{restoretext}
Note that indeed, the right morphism contains a ``cap singularity" (marked in \cyellow{}) which is part of the coherent invertibility data for $h\dualdag$.  

As a concluding remark, with a pointer towards \autoref{sec:CWcomplexes}, note also that the ``cap singularity" witnesses $(h\dualdag)\inv \circ h\dualdag = \id$. The reader might wonder why inverses like $(h\dualdag)\inv$ did so far not show up in our discussion. This is because we decided to omit framing here, which will be remedied in \autoref{sec:CWcomplexes} .
\end{rmk}
\end{enumerate}
This completes the (non-unique) construction of $\sX$ from $X$. \\
 
To summarise the idea of this section: maps $[D^k, X]$ (up to homotopy with an appropriate ``rel-boundary" condition) are to spaces what manifold diagrams with duals are to presented associative $\infty$-groupoids---and using the generalised Thom-Pontryagin construction, roughly speaking, they differ only by a notion of ``direction" (in the sense of a choice of globular foliation for $\kP_{D^k}(f)$, $f \in [D^k,X]$, as discussed above). This in turn gives a clear idea on where future work on interesting questiona such as the homotopy hypotheses has to start.

\section{Appendix B: $n$-fold categories}

The framework presented here is (very crucially) intrinsically cubical. By working with cubes instead of globes (that is, omitting the globularity condition but leaving the well-typedness condition in place) we will also briefly define and discuss a notion of \free{} associative $n$-fold category in \autoref{ch:nfold}. This provides a nice graphical calculus for double, triple and $n$-fold categories.

\section{Appendix C: Weak categories}

In \autoref{ch:weak} we will present a straight-forward definition of \free{} weak $n$-categories which is obtained via a colimiting process from \free{} associative $n$-categories. Based on this definition we will recover the usual perspective on higher categories: compositions are defined ``up to equivalence" and coherences (associators, pentagonators etc.) are present as morphisms. Our presentation will also motivate the equivalence between the associative and the weak approach.

\renewcommand{\thesection}{\arabic{chapter}.{\arabic{section}}}
\renewcommand{\theHsection}{\arabic{chapter}.{\arabic{section}}}

\chapter{Notation and prerequisites} \label{ch:prerequisites}

In this chapter we recall basic definitions and fix notation needed later on. The reader is encouraged to skip this chapter and refer back to it if needed.

In \autoref{sec:basic} we recall elementary notions from category theory such as the cartesian closed structure of \Cat{}, \SetCat{} and \Bool{}, and fix notation for subcategories and restrictions of functors. We also recall posets, profunctorial relations (which are \Bool{}-enriched profunctors) and a Grothendieck construction for the latter. We also discuss our visual notation style for these notions. \autoref{sec:label_poset} further recalls properties of ``labelled posets" which are posets $P$ together with a labelling functor $F : P \to \cC$ to a category $\cC$.

\section{Very elementary concepts in category theory} \label{sec:basic}

\subsection{Cartesian closed structures}

In the following we discuss our choices of notation for (very) elementary concepts of category theory. In an attempt to make this work as self-contained as possible, we only assume minimal familiarity with categories, functors, natural transformations, the Yoneda lemma and adjoints. 

Given a category $\cC$, we denote the object set of $\cC$ by $\obj(\cC)$, and the morphism set by $\mor(\cC)$. 

\begin{notn}[Products, $\ppi_\cC$ and implicit isomorphisms] \label{notn:product_categories} Given two categories $\cC$, $\cD$ we denote by $\cC \times \cD$ their product category, whose objects are denoted by tuples $(c,d)$, where $c \in \cC$, $d \in \cD$ and whose morphisms are denoted by tuples $(f,g) : (a,b) \to (c,d)$, where $(f : a \to c) \in \mor(\cC)$, $(g :  b \to d) \in \mor(\cD)$. We denote the projection 
functors by $\ppi_\cC : \cC \times \cD \to \cC$ and $\ppi_\cD : \cC \times \cD \to \cD$, which are defined to map tuples $(x,y)$ (of objects or morphism) to their first component $x$ respectively to their second component $y$.

Given two functors $F : \cC_1 \to \cC_2$ and $G : \cD_1 \to \cD_2$ then their \textit{product} $F \times G : \cC_1 \times \cD_1 \to \cC_2\times \cD_2$ maps objects $(c,d)$ to $(Fc,Gd)$ and morphisms $(f,g)$ to $(Ff, Gg)$. This gives rise to a functor
\begin{equation}
(- \times -)  : \Cat \times \Cat \to \Cat
\end{equation}
We will usually implicitly use that $(- \times -)$ is associative and unital up to natural isomorphism, that is
\begin{equation}
(\cC_1 \times \cC_2) \times \cC_3 \iso \cC_1 \times (\cC_2 \times \cC_3)
\end{equation}
and (denoting by $\bnum{1}$ the terminal category)
\begin{equation}
\cC \times \bnum{1} \iso \cC
\end{equation}
\end{notn}

\begin{notn}[$\Bool$, $\SetCat$, $\Cat$ and their cartesian closed structure] \label{notn:basic_category_theory} \hfill
\begin{enumerate}
\item Let $\Bool$ be the category of truth values which has as objects the truth values $\bot$ (``false") and $\top$ (``true") together with a single non-identity morphisms $\bot \to \top$. Let $(-\vee-), (- \wedge -),(- \imp-)$ denote the usual conjunction, disjunction and implication operations on truth values (note that these operations are functorial).
\item  We denote by $\SetCat$ the category of sets and functions, and by $(-\cup-),(-\times-), \Map(-,-)$ the usual sum, product and function set operations on sets (note that these operations are functorial).
\item Given a set $C$, denote by $(-=_C-)$ the diagonal relation, i.e. the function $C \times C \to \Bool$ which maps $(c,c')$ to $\top$ if $c = c'$ and to $\bot$ otherwise.
Given a category $\cC$ denote by $\Hom_\cC : \cC\op \times \cC \to \SetCat$ the hom-functor on $\cC$.
\item Let $\Cat$ denote the category of categories and functors. Given categories $\cC,\cD$ and functors $F,G : \cC \to \cD$, denote by $\Fun(\cC,\cD)$ the set of functors from $\cC$ to $\cD$ and by $\mathrm{Nat}(F,G)$ the set of natural transformations from $F$ to $G$. The product functor $(- \times -)$ together with the \textit{internal hom}\footnote{Since ultimately all categorical structures in this work will be finite, we will be disregarding any size issues.}
\begin{equation}
[-,-] : \Cat\op \times \Cat \to \Cat
\end{equation}
equip $\Cat$ with \textit{cartesian closed} structure: $[\cC,\cD]$ is the category which has functors $F : \cC \to \cD$ as objects an natural transformations $\alpha : F \to G$ as morphisms. Given $H : \cC' \to \cC$ and $K : \cD \to \cD'$ then $[H,K] : [\cC,\cD] \to [\cC',\cD']$ acts on functors by pre-composition with $H$ and post-composition with $K$, that is,
\begin{equation}
(F : \cC \to \cD)\quad \mapsto\quad (KFH : \cC' \to \cD')
\end{equation}
and it acts on natural transformations by pre-whiskering with $H$ and post-whiskering with $K$, denoted by
\begin{equation}
(\lambda : F_1 \to F_2) \quad \mapsto \quad (K\lambda H : KF_1H \to KF_2H)
\end{equation}
Note that if either $H$ or $K$ are the identity functor $\id$ they will be omitted in this notation.
\item There is a \textit{object functor}
\begin{equation}
\obj : \Cat \to \SetCat
\end{equation}
mapping categories to their object sets, and functors to their underlying functions on object sets. Note that this functor has a left adjoint
\begin{equation}
\discr : \SetCat \to \Cat
\end{equation}
mapping sets to their corresponding discrete categories. Both $\obj$ and its left adjoint are sometimes kept implicit in our notation.
\end{enumerate} 
\end{notn}

\subsection{Special functors}

We introduce notation for a range of special functors used frequently later on. First, further to our notes to the reader in \autoref{sec:notes}, we make the following  remark.

\begin{rmk}[Foundation independence/computer implementability] \label{rmk:material_set_theory} Since we will ultimately be working with special (finite) structures that don't admit non-trivial isomorphism, it will be convenient take the standpoint of ``material set theory". That means, most elements in this work can be thought of as existing independently of the sets or categories that they live in (especially with respect to subsets and subcategories). 

Since we claim that at least parts of this thesis are ``foundation-independent" and/or computer-implementable, we add that we understand an ``element" to mean a concrete \textit{name} that we use to refer to it (for instance, a string of numbers or bits). Note the terms ``name", ``label" or ``color" are usually used synonymously.
\end{rmk}

\noindent The previous remark in particular influences how we think about ``subsets" and ``subcategories" in this thesis. Keeping this in mind, we introduce the following notation.

\begin{notn}[Subcategories, images and restrictions of functors] \label{notn:subsets_and_restrictions} \hfill 
\begin{enumerate}

\item Given a functor $F : \cC \to \cD$ denote by $\im(F)$ the image of $F$,  that is, the subgraph given by objects and morphism in the image of $F$. Note that $\im(F)$ is not necessarily a subcategory of $\cD$, unless $F$ is injective on objects.

\item A subcategory (or subset) $\cC_0 \subset \cC$ comes with a canonical inclusion functor $\cC_0 : \cC_0 \into \cC$ defined to map $d \mapsto d$ on objects and $f \mapsto f$ on morphisms. This entails the notation $\cC = \id_\cC : \cC \to \cC$.

\item Given a functor $F : \cC \to \cD$ and a subcategory $\cC_0 \subset \cC$ we define $\rest F {\cC_0} : \cC_0 \to \cD$ to be the composite
\begin{equation}
\cC_0 \xto {\cC_0} \cC \xto {F} \cD
\end{equation}
$\rest F {\cC_0}$ is called the \textit{domain restriction of $F$ to $\cC_0$}. Both of $\rest F {\cC_0}$ and $F\cC_0$ will be used depending on which notation is more convenient.

\item Let $\cC$ be a category with $c,d \in \cC$ and $(f : c \to d)$ a morphism in $\cC$. Then we denote by
\begin{equation}
\Delta_c : \bnum{1} \to \cC
\end{equation}
the functor mapping $0$ to $d$, and by
\begin{equation}
\Delta_f : \bnum{2} \to \cC
\end{equation}
the functor mapping $(0 \to 1)$ to $f$. $\Delta_c$ and $\Delta_f$ are called \textit{test functors}.%

\item Now assume a functor $F : \cC \to \cD$ and an injective and faithful functor $H : \cD_0 \to \cD$. Let $\cC_0 \subset \cC$ be the maximal subcategory such that $\im(F \cC_0) \subset \im(H)$. Define a functor, called the  \textit{codomain restriction of $F$ to $H$} by setting
\begin{align} \label{eq:codomain_restriction}
H\pbstar : \cC_0 ~& \to D_0 \\
x ~& \mapsto H\inv F(x)
\end{align}
where $x$ is either an object or a morphism. Note that this defines $H\pbstar F$ as a pullback of $F$ along $H$.

\end{enumerate}
\end{notn}

\begin{notn}[Constant functors] \label{notn:simplicial} If $F : \cC \to \cD$ is a functor, $d \in \cD$, then we write 
\begin{equation}
F = \const_d
\end{equation}
to mean that $F$ is the constant functor with image $d$. If $d$ is implicit, we write $F = \const$. %
\end{notn}

\subsection{Relations and profunctors}

We first recall the notion of relations. A slight generalisation of this notion, namely that of profunctorial relations, will play a central role in developing the combinatorics of associative $n$-categories. For completeness, we also briefly remind the reader of the definition of profunctors.

\begin{defn}[Relations and profunctors] \label{defn:relations_and_profunctors} \hfill
\begin{enumerate}
\item Given sets $C, D$, a \textit{relation} $R : C  \xslashedrightarrow{} D$ is a function 
\begin{equation}
R : C \times D \to \Bool
\end{equation}
Denote the set of relations $R : C  \xslashedrightarrow{} D$ by $\Rel(C,D)$. Given two relations $R : C \xslashedrightarrow{} D$, $S : D \xslashedrightarrow{} E$ their composite $R \odot S : C \xslashedrightarrow{} E$ is given by the mutual implication 
\begin{equation}
(R \odot S) (c,e) \iff \exists d \in D.~ R(c,d) \wedge S(d,e)
\end{equation}
An alternative description of composition is the following: introduce structure of a (posetal) category on $\Rel(C,D)$ by setting for $T,T' \in \Rel(C,D)$
\begin{equation}
T \to T' \quad \text{~if and only if~}\quad \forall	 c \in C, d \in D ~.~ T(c,d) \imp T'(c,d)
\end{equation}
Then we can define $(R \odot S)$ by the formula
\begin{equation}
(R \odot S) (c,e) = \int^{d : D}  R(c,d) \wedge S(d,e)
\end{equation}
where the functor $\int^{d : D} : \Rel(D,D) \to \Bool$, $T \mapsto \int^{d : D} T(d,d)$ is defined as the left adjoint
\begin{equation}
\int^{d : D} \vdash \big((- =_D -) \imp -\big) : \Bool \to \Rel(D,D)
\end{equation}

\item Given categories $\cC$, $\cD$, a \textit{profunctor} $R : \cC \xslashedrightarrow{} \cD$ is a functor 
\begin{equation}
R : \cC\op \times \cD \to \SetCat
\end{equation}
Denote the category of profunctors $R : \cC \xslashedrightarrow{} \cD$ and natural transformations by $\Prof(\cC,\cD)$. Given two profunctors $R : \cC \xslashedrightarrow{} \cD$, $S : \cD \xslashedrightarrow{} \cE$ their composite $R \odot S : \cC \xslashedrightarrow{} \cE$ (assuming enough colimits in $\cD$ exist) is given by the formula
\begin{equation}
(R \odot S) (-,-) \iso \int^{d : \cD} R(-,d) \times S(d,-)
\end{equation}
where the functor $\int^{\cD} : \Prof(\cD,\cD) \to \SetCat$ can be defined as the left adjoint
\begin{equation}
\int^{d : \cD} \dashv \Map(\Hom_{\cD}, -) : \SetCat \to \Prof(\cD,\cD)
\end{equation}
\end{enumerate}
\end{defn}

From the above definitions it should be clear that profunctors can be thought of as a categorification of relations. However, somewhat more naturally profunctors categorify the more general notion of $\Bool$-enriched profunctors, which (when restricted to posets) we will call profunctorial relations. This choice of terminology emphasizes that these structures compose as relations. To keep our discussion self-contained, we will however not refer to enriched category theory in the following definition. Instead we will use that there is a fully-faithful inclusion
\begin{equation} \label{eq:Bool_to_Set}
\iP : \Bool \into \SetCat
\end{equation}
mapping $\bot$ to the empty set $\emptyset$ and $\top$ to the singleton set $\Set{*}$. Usually this \Bool{}-to-\SetCat{} parsing (and its inverse) will be kept implicit.

\subsection{Posets} \label{sec:poset}

We recall the basic notion of posets and related concepts.

\begin{defn}[Preorders and posets] \label{defn:posets} \hfill
\begin{itemize}
\item A \textit{preorder} $X$ is a category whose hom sets are either the empty set $\emptyset$ or the singleton set $\Set{*}$. Let $x,y \in X$. If the hom-set $X(x,y)$ is non-empty, i.e. it equals $\Set{*}$, then $(x \to y)$ denotes  morphism $* : x \to y$. Otherwise we read $(x \to y)$ as false. %
\item A \textit{poset} $X$ is a preorder such that if $x, y \in {X}$ and both $(x \to y) \in \mor(X)$ and $(y \to x) \in \mor(X)$ then $x = y$. 
\item A \textit{functor of posets} $F : X \to Y$ is a functor $X\to Y$ between posets $X$ and $Y$. Poset and their functors form the category $\Pos$ which is a subcategory of $\Cat$. The inclusion $\Pos \into \Cat$ has a left adjoint $\PSk : \Cat \to \Pos$ associating to a category its \textit{posetal skeleton}.
\item A \textit{subposet} $Y$ of a poset $X$,  written $Y \subset X$, is a subcategory $Y$ of $X$. We say $Y$ is a \textit{full} subposet of $X$ if it is a full subcategory. An inclusion of posets $F : X \into Y$ is functor of posets which is injective on objects. 
\end{itemize}
\end{defn}

\begin{rmk}[Data for maps and full subposets determined on objects] \label{rmk:poset_data_objects} For maps and full subposets it suffices to give data for the above definitions on objects:
\begin{itemize}
\item We note that a functor of posets $F : X \to Y$ as a functor is already fully determined by its underlying function of objects sets $\obj(F) : \obj(X) \to \obj(Y)$.
\item We also note that a full subposets $Y \subset X$ is fully determined by its subset of objects $\obj(Y) \subset \obj(X)$.
\end{itemize}
\end{rmk}

We record relevant posets related to the integers and their notation.

\begin{constr}[$\lZ$, intervals and $\bnum{n}$] \label{eg:posets} \hfill
\begin{enumerate}
\item We denote by $\lZ$ the poset of integers. That is, $\obj(\lZ)$ is the set of integers and $(i \to j) \in \mor(\lZ)$ if $i \leq j$. We will use $\leq$ exclusively to refer to the poset of integers and its full subposets. A map of subposets of $\lZ$ is also called a \textit{monotone} map. 
\item Given $a, b \in \lZ$ with $a \leq b$, then $[a,b] \subset \lZ$ respectively $]a,b[ ~ \subset\lZ$ denote the closed interval respectively open interval with endpoints $a,b$.

\item For a natural number $n \in \lN$, we denote by $\bnum{n}$ the full subposet of $\lZ$ with object set $]-1,n[$ (thus, $\bnum{0}$ is the empty poset and $\bnum n$ is the $(n-1)$-simplex). We denote by $\delta^n_i : \bnum{n} \into (\bnum{n+1})$ and $\sigma^n_i : (\bnum{n+2}) \to (\bnum{n+1})$ the usual face and degeneracy maps for simplices (cf. \cite{riehl2011leisurely}), that is $\delta^n_i$ is injective but omits $i$ in its image, and $\sigma^n_i$ is surjective with $i$ having two preimages.
\end{enumerate}
\end{constr} 

\noindent We will give examples of these definitions in the next section.

\subsection{Depicting posets and functors on posets}

The following introduces a (sometimes) helpful visual notation for posets, poset maps and maps from posets into categories.

\begin{notn}[Depicting posets and functors on posets]  \label{notn:depicting_posets} Given a poset $X$, we depict an object $x$ of $X$ as a point labelled by $x$ (though often we will omit the labelling), and draw an arrow from point $x$ to point $y$ whenever $(x \to y) \in \mor(X)$ for $x \neq y$. For instance, the following
\begin{restoretext}
\begingroup\sbox0{\includegraphics{test/page1.png}}\includegraphics[clip,trim=0 {.3\ht0} 0 {.3\ht0} ,width=.8\textwidth]{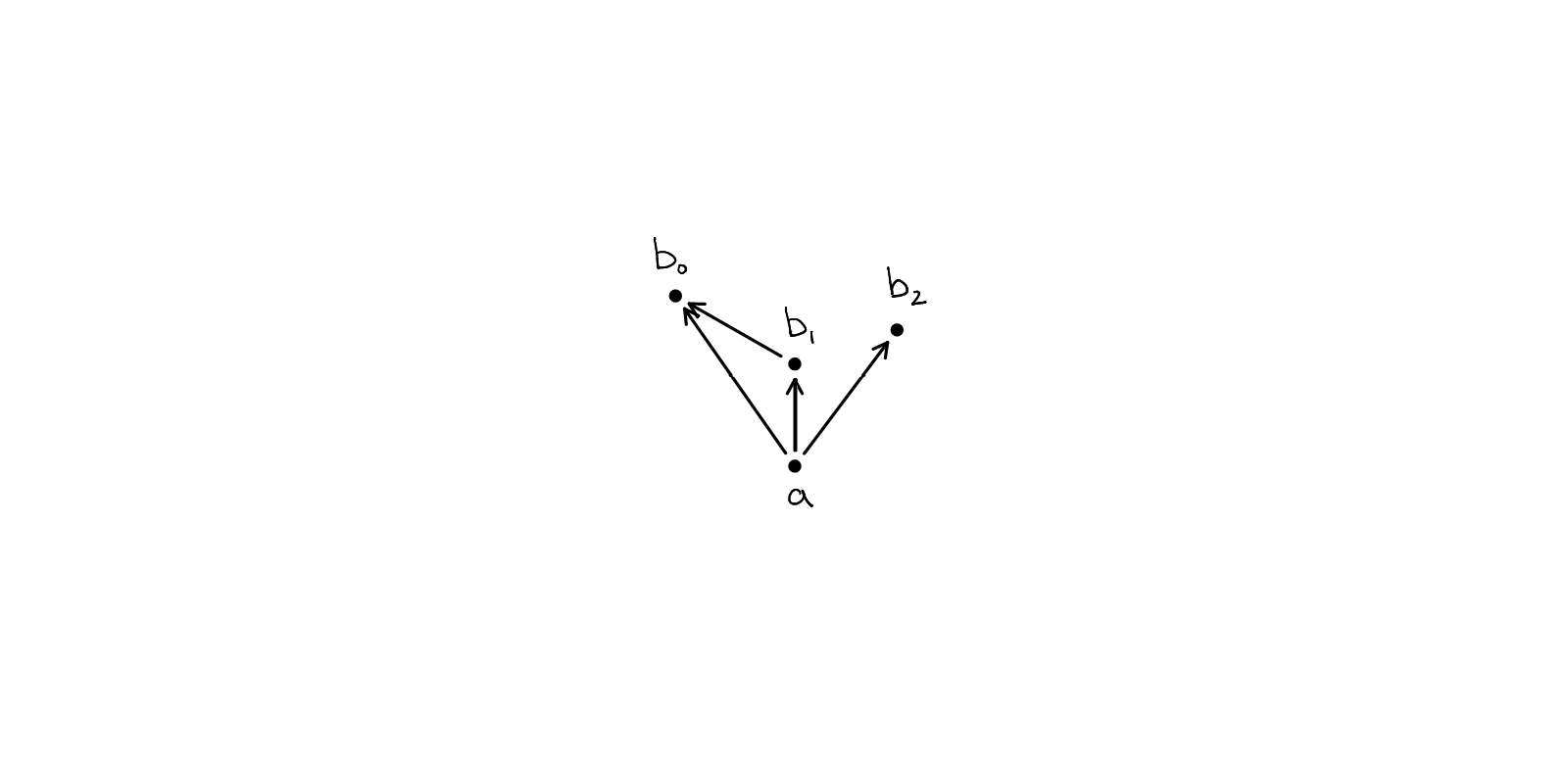}
\endgroup\end{restoretext}
is a depiction of a poset $X$ with four objects and four non-identity morphisms.

A functor of posets $F : X \to Y$ can be depicted by color coding its preimages $F\inv(y)$ for each $y \in Y$. For instance,
\begin{restoretext}
\begingroup\sbox0{\includegraphics{test/page1.png}}\includegraphics[clip,trim=0 {.2\ht0} 0 {.25\ht0} ,width=.8\textwidth]{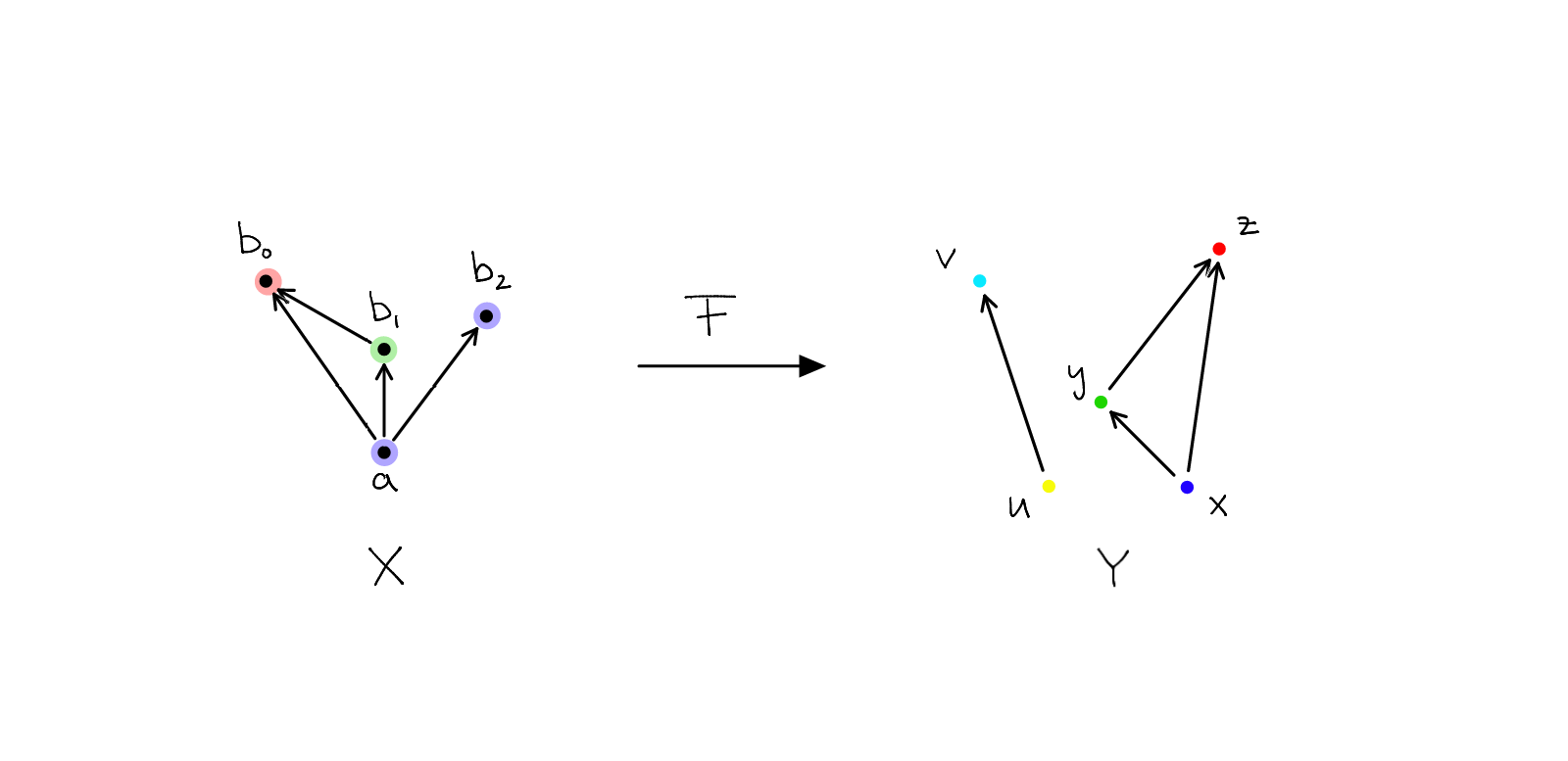}
\endgroup\end{restoretext}
depicts a map $F: X \to Y$, where $Y$ is a poset with five objects, and $F$ maps $a$ and $b_2$ to the same object in $Y$, while $a $, $b_1$ and $b_0$ are all mapped to different objects in $Y$. Since as remarked above poset functors are in fact fully determined by their mapping on objects, we can also omit the coloring of morphisms. In addition to (and sometimes in place of) the color coding, preimages will be indicated by labelling with symbols corresponding to symbols of their image points, or simply as ``spatially lying over" their imagine points relative to the direction of the map's arrow (labelled by $F$ in the above). The following
\begin{restoretext}
\begingroup\sbox0{\includegraphics{test/page1.png}}\includegraphics[clip,trim=0 {.1\ht0} 0 {.1\ht0} ,width=\textwidth]{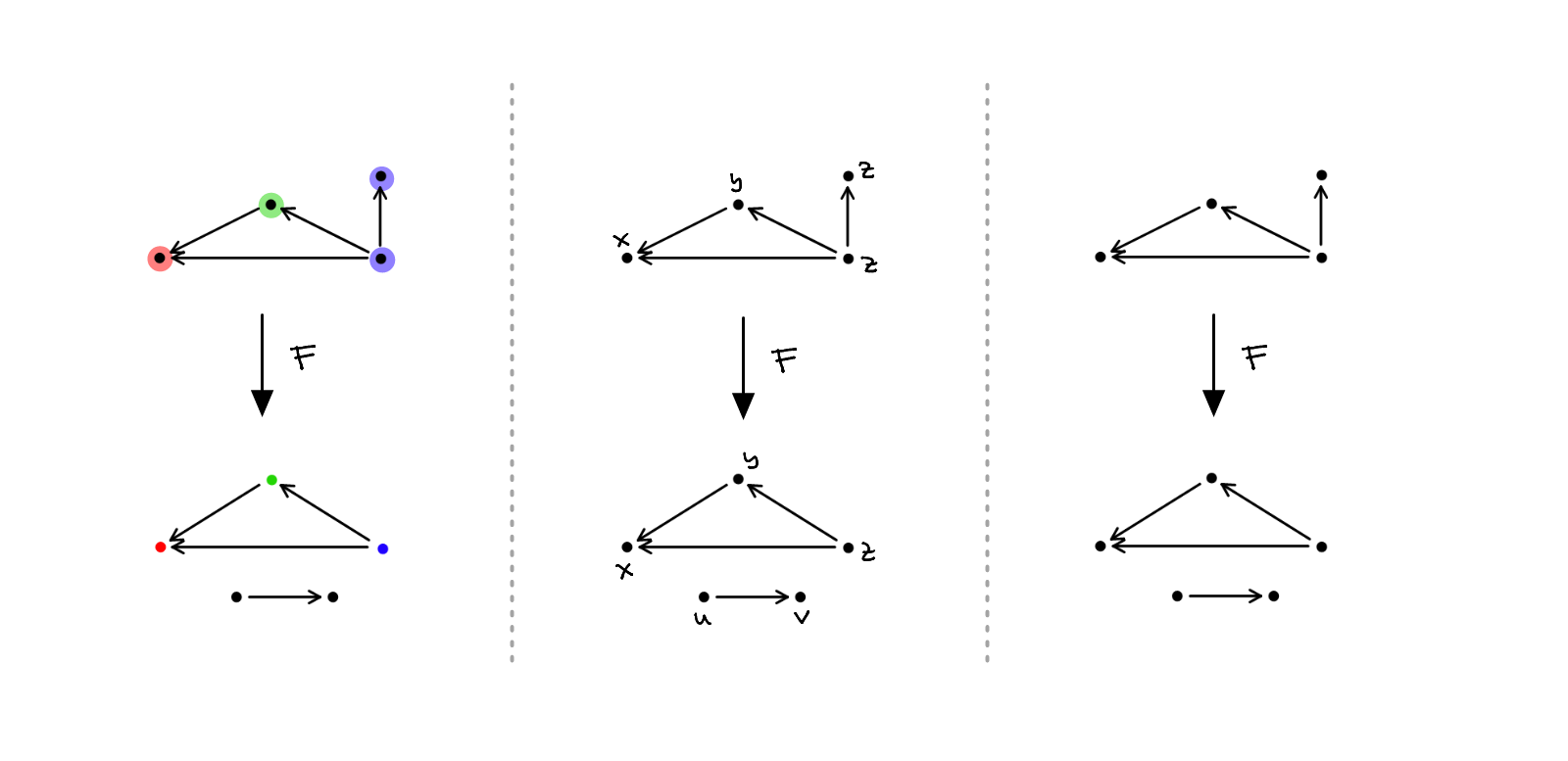}
\endgroup\end{restoretext}
shows three depictions (one coloured, one labelled and the last uncoloured and unlabelled but spatially oriented) of the same functor $F$ defined above.

A functor $F : X \to \cC$ from a poset $X$ to a category $\cC$ can be depicted by an $X$-shaped diagram in $\cC$. Since morphisms contain information about domain and codomain objects, it is also sufficient to only give labels for morphism. For instance, the following
\begin{restoretext}
\begingroup\sbox0{\includegraphics{test/page1.png}}\includegraphics[clip,trim=0 {.3\ht0} 0 {.25\ht0} ,width=.8\textwidth]{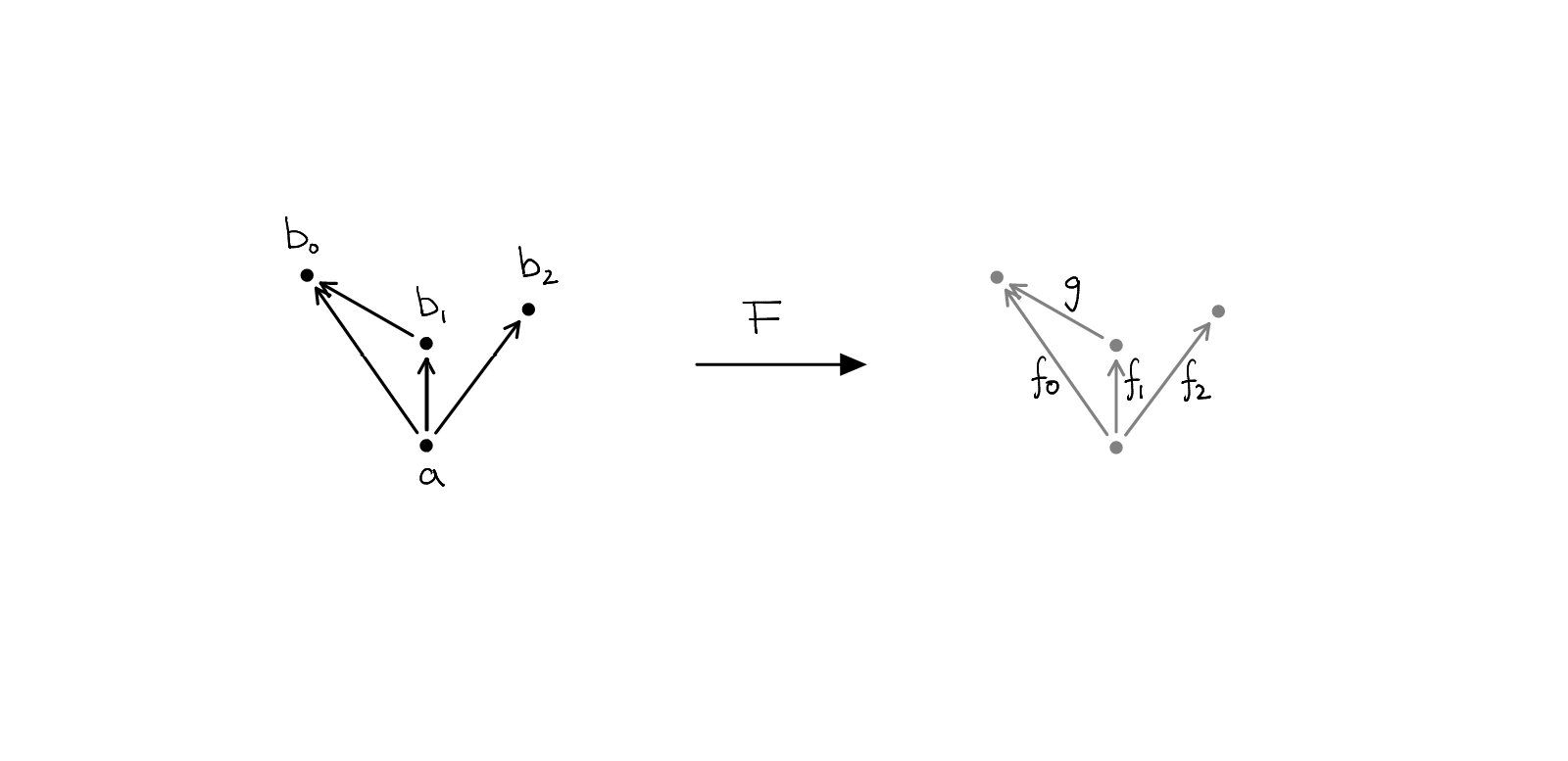}
\endgroup\end{restoretext}
depicts a functor $F : X \to \cC$, such that $F$ maps $a$, $b_i$ to objects $c, d_i$ in $\cC$ respectively, as well as morphism $(b_1 \to b_0)$, $(a \to b_i)$ to morphisms $(g : d_1 \to d_0), (f_i : c \to d_i)$ in $\cC$ (for $i \in \Set{1,2,3}$). 
\end{notn}

\begin{eg}[Finite totally ordered sets] Using our conventions, the posets $\bnum{3}$ and $\bnum{4}$ can be depicted as follows
\begin{restoretext}

\begingroup\sbox0{\includegraphics{test/page1.png}}\includegraphics[clip,trim=0 {.3\ht0} 0 {.3\ht0} ,width=.9\textwidth]{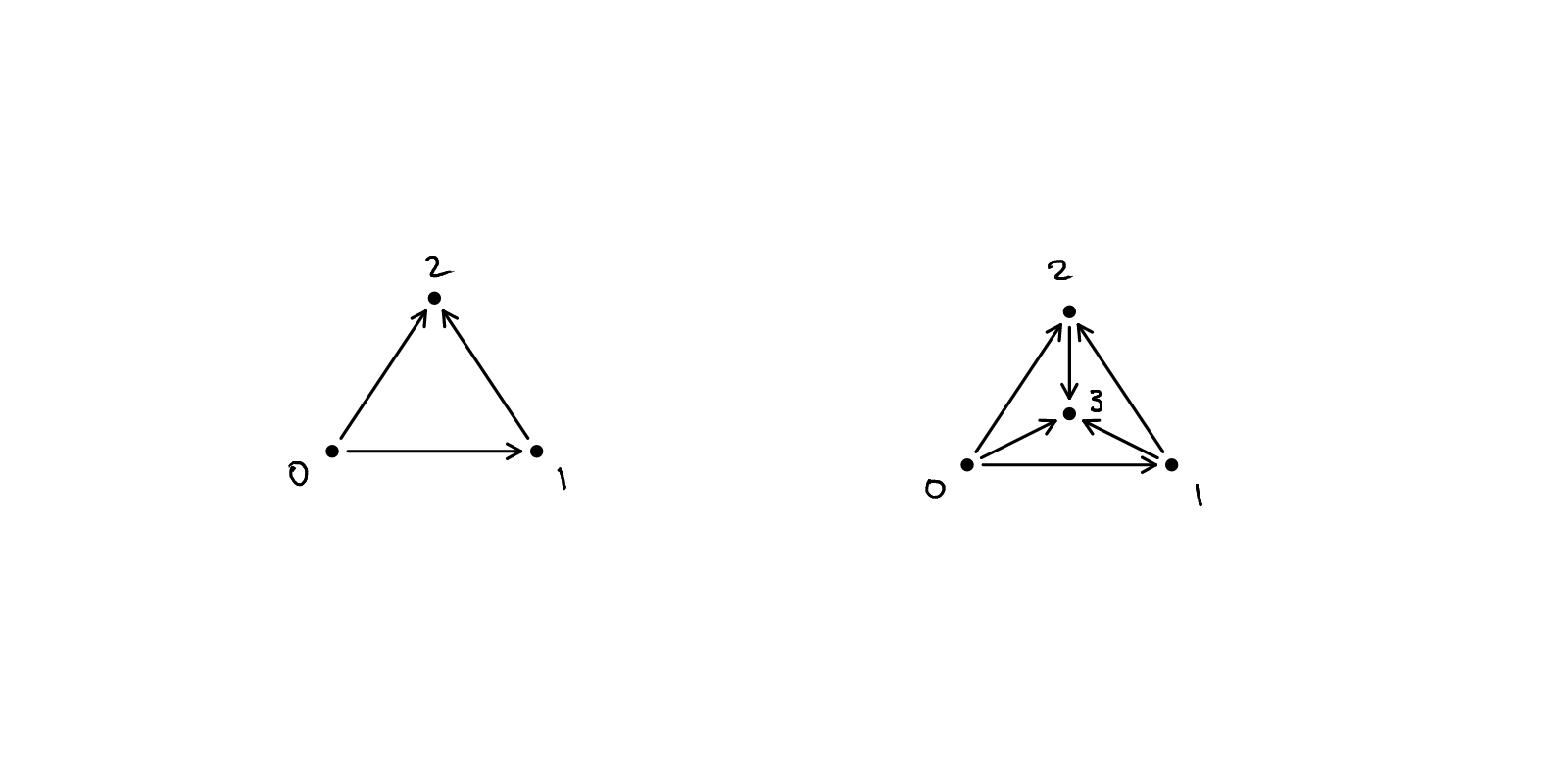}
\endgroup\end{restoretext}
Some face and degeneracy maps for $n = 3$ and $n = 1$ can be depicted as follows.
\begin{restoretext}
\begingroup\sbox0{\includegraphics{test/page1.png}}\includegraphics[clip,trim=0 {.15\ht0} 0 {.09\ht0} ,width=.9\textwidth]{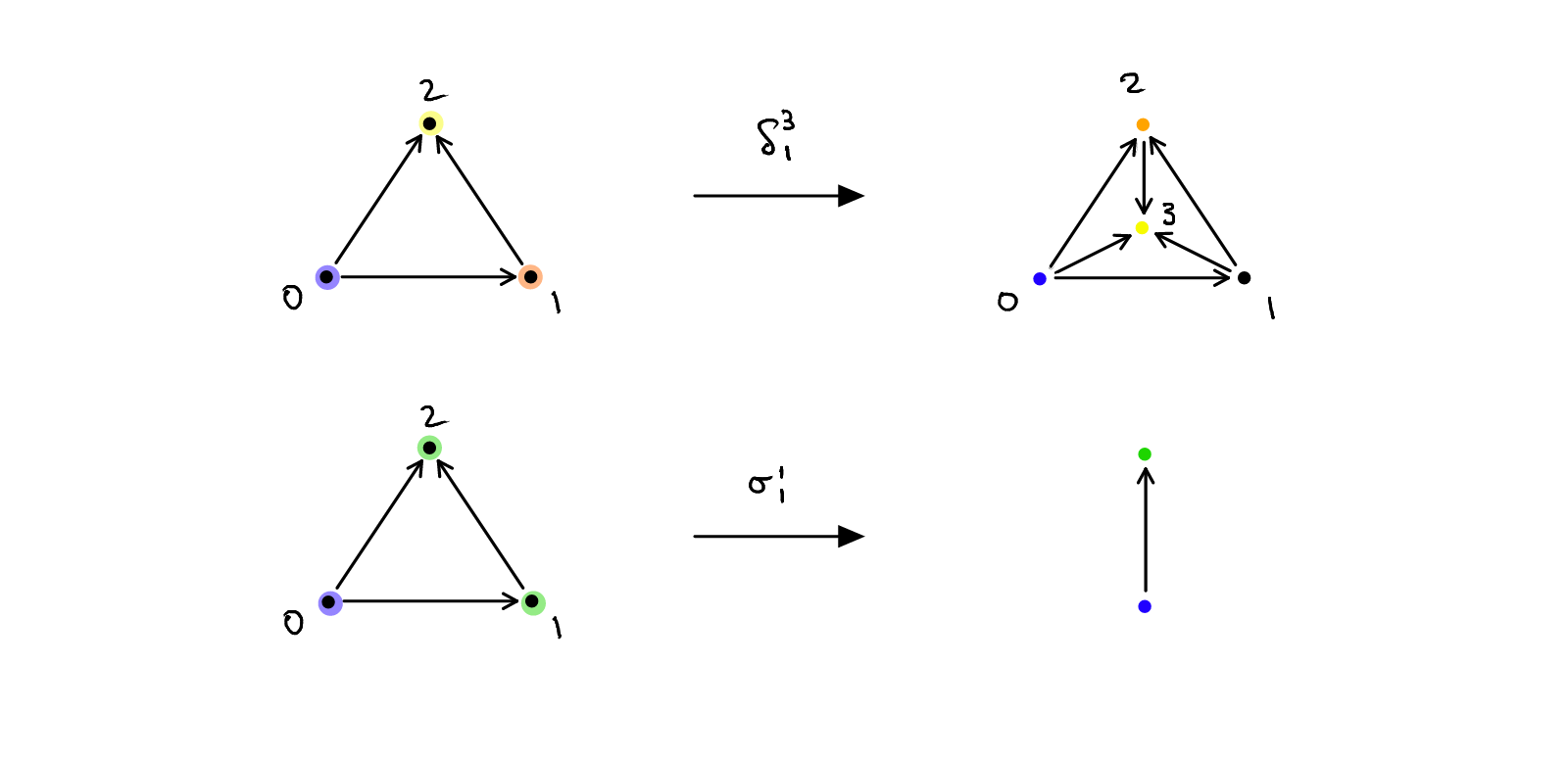}
\endgroup\end{restoretext}
\end{eg}

While for these example the notation seems quite cumbersome, it will turn out to be quite useful once we work with bigger posets.

\subsection{Profunctorial relations}

Having seen relations and posets, we now introduce $\Bool$-enriched profunctors as follows.

\begin{defn}[Profunctorial relations] \label{defn:prerequisites} Let $X$, $Y$, $Z$ be posets. A \textit{profunctorial relation} $R : X \xslashedrightarrow{} Y$ is a functor
\begin{equation} \label{eq:profunctorial_relation_defn}
R : X\op \times Y \to \Bool
\end{equation}
Denote the set of profunctorial relations $R : X \xslashedrightarrow{} Y$ by $\PRel(X,Y)$.
\end{defn}

\begin{notn}[Depicting (profunctorial) relations] \label{notn:depicting_prel} A (profunctorial) relation $R : X\op \times Y \to \Bool$ can by depicted as follows: we depict the (po)sets $X$ and $Y$ individually as points and arrows, and draw an edge between $x \in X$ and $y \in Y$ if and only if $R(x,y)$ is true. For instance
\begin{restoretext}
\begingroup\sbox0{\includegraphics{test/page1.png}}\includegraphics[clip,trim=0 {.25\ht0} 0 {.15\ht0} ,width=.8\textwidth]{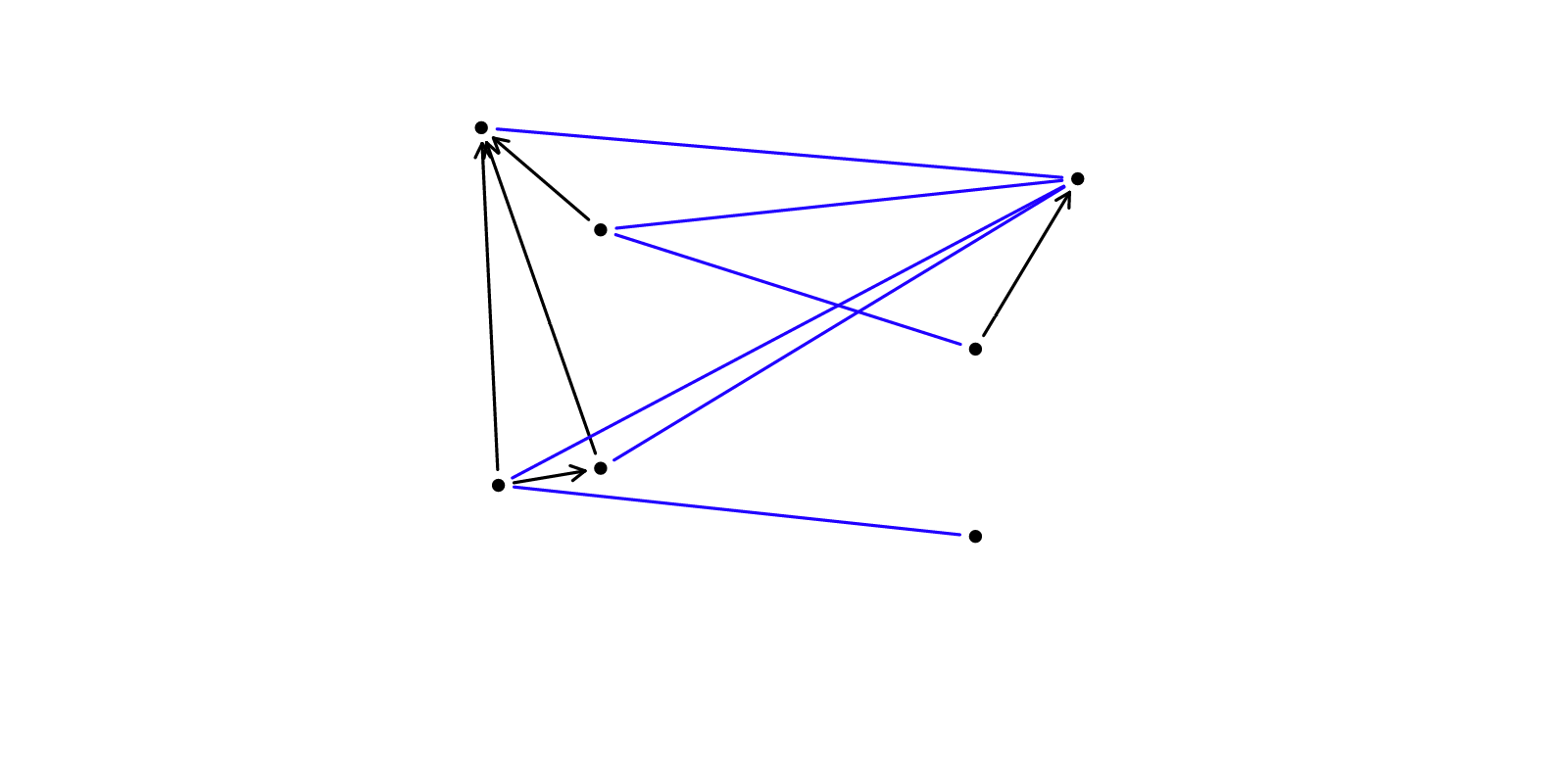}
\endgroup\end{restoretext}
depicts a profunctorial relation $R$ between a four element poset $X$ and a three elements poset $Y$. Note also that here we highlighted the edges in \cblue{}. However the colouring can be omitted since we depict edges only as lines (and not as arrows). That means, that the direction of the profunctorial relation itself needs to be inferred by other means: in the above, we use the convention that $Y$ being on the right of $X$ should indicate that $X$ is the domain and $Y$ the codomain of $R$.
\end{notn}

Composition of profunctorial relations can be derived from the general theory of enriched profunctors (which uses a definition analogous to the ones given in \autoref{defn:relations_and_profunctors}). It can also be characterised more directly, as done in the following construction.

\begin{constr}[Composition of profunctorial relations and the category $\PRel$]\label{constr:PRel} We construct the composite of two (composable) profunctorial relations and use this to define the category $\PRel$. Given $R \in \PRel(X,Y)$ first note that we have $\obj(R) \in \Rel(X,Y)$. Conversely, given a relation $R \in \Rel(X,Y)$, then there is a (necessarily unique) $R' \in \PRel(X,Y)$ with $\obj(R') = R$ if and only if $R$ satisfies the \textit{profunctoriality conditions}
\begin{gather}\label{eq:profunctorial_relation}
\text{if~} x \to x' \text{~and~} R(x',y) \text{~then~} R(x,y) \\
\text{if~} R(x,y') \text{~and~} y' \to y \text{~then~} R(x,y)
\end{gather}
Now, given two profunctorial relations $R : X \xslashedrightarrow{} Y$, $S : Y \xslashedrightarrow{} Z$, we observe that the composite relation $\obj(R)\odot \obj(S)$ again satisfies the profunctoriality conditions. We can thus define $R \odot  S$ by setting (cf. \autoref{rmk:poset_data_objects})
\begin{equation}
\obj(R \odot  S) = \obj(R)\odot \obj(S)
\end{equation}
This notion of composition gives rise to a category $\PRel$ with objects being posets and morphisms being profunctorial relations. Note that identities in this category are given by \textit{hom-relations} $\Hom_X : X \xslashedrightarrow{} X$ defined by $\Hom_X (x,y) \iff (x \to y) \in \mor(X)$. 
\end{constr} 

To illustrate the condition \eqref{eq:profunctorial_relation} we give examples of relations that satisfy and don't satisfy the condition.

\begin{eg}[Profunctorial relations] \label{eg:prel} The following are three examples of profunctorial relations $R_1 : \bnum{3} \xslashedrightarrow{} \bnum{4}$, $R_2 : \bnum{4} \xslashedrightarrow{} \bnum{3}$ and $R_3 : \bnum{3} \xslashedrightarrow{} \bnum{3}$
\begin{restoretext}
\begingroup\sbox0{\includegraphics{test/page1.png}}\includegraphics[clip,trim=0 {.15\ht0} 0 {.05\ht0} ,width=.8\textwidth]{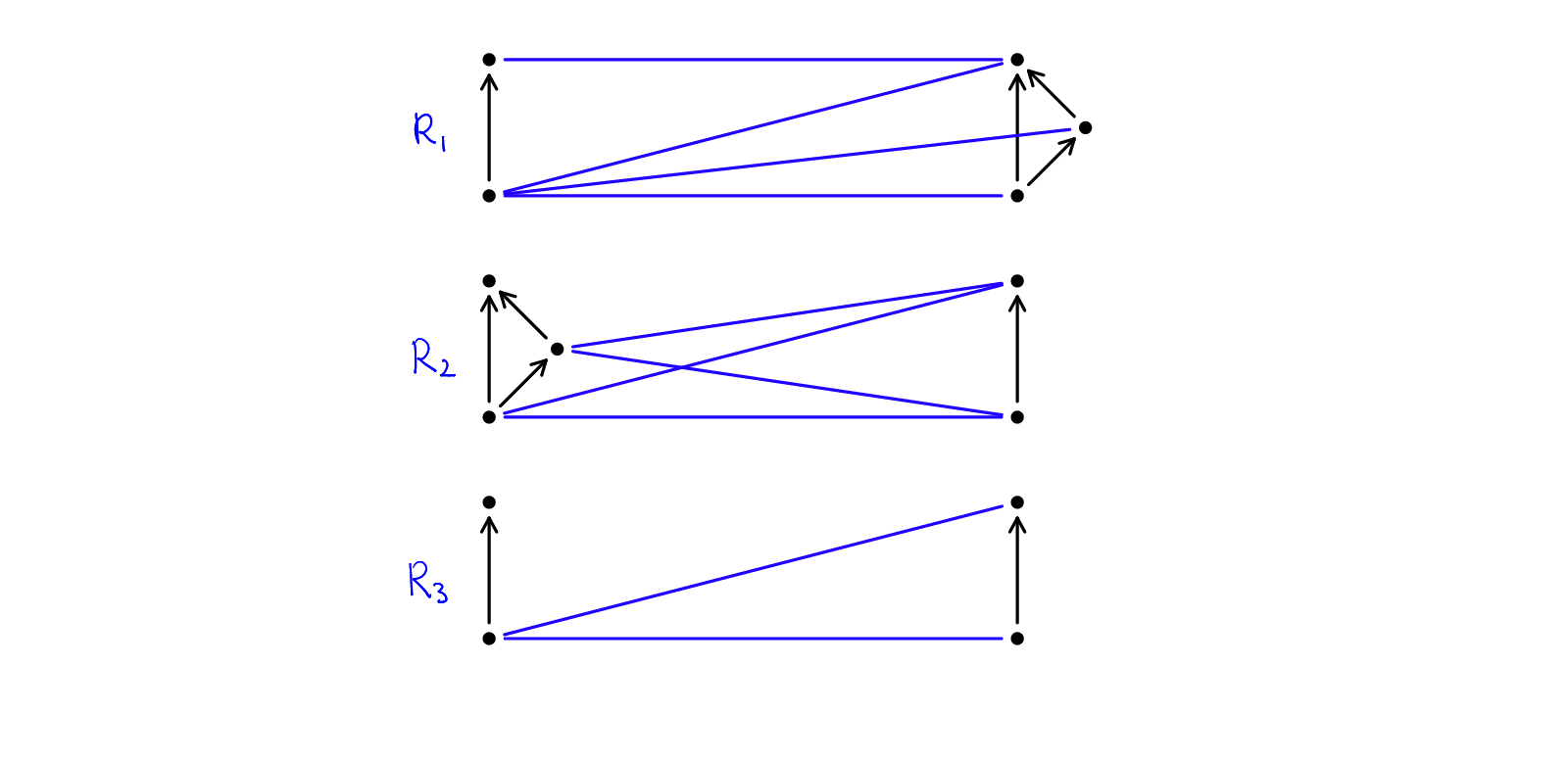}
\endgroup\end{restoretext}
Note that $R_3 = R_1 \odot R_2$. T

The following is a variation of the previous examples, defining $R'_2 : \bnum{4} \xslashedrightarrow{} \bnum{3}$ and $R'_3 : \bnum{3} \xslashedrightarrow{} \bnum{3}$.
\begin{restoretext}
\begingroup\sbox0{\includegraphics{test/page1.png}}\includegraphics[clip,trim=0 {.15\ht0} 0 {.05\ht0} ,width=.8\textwidth]{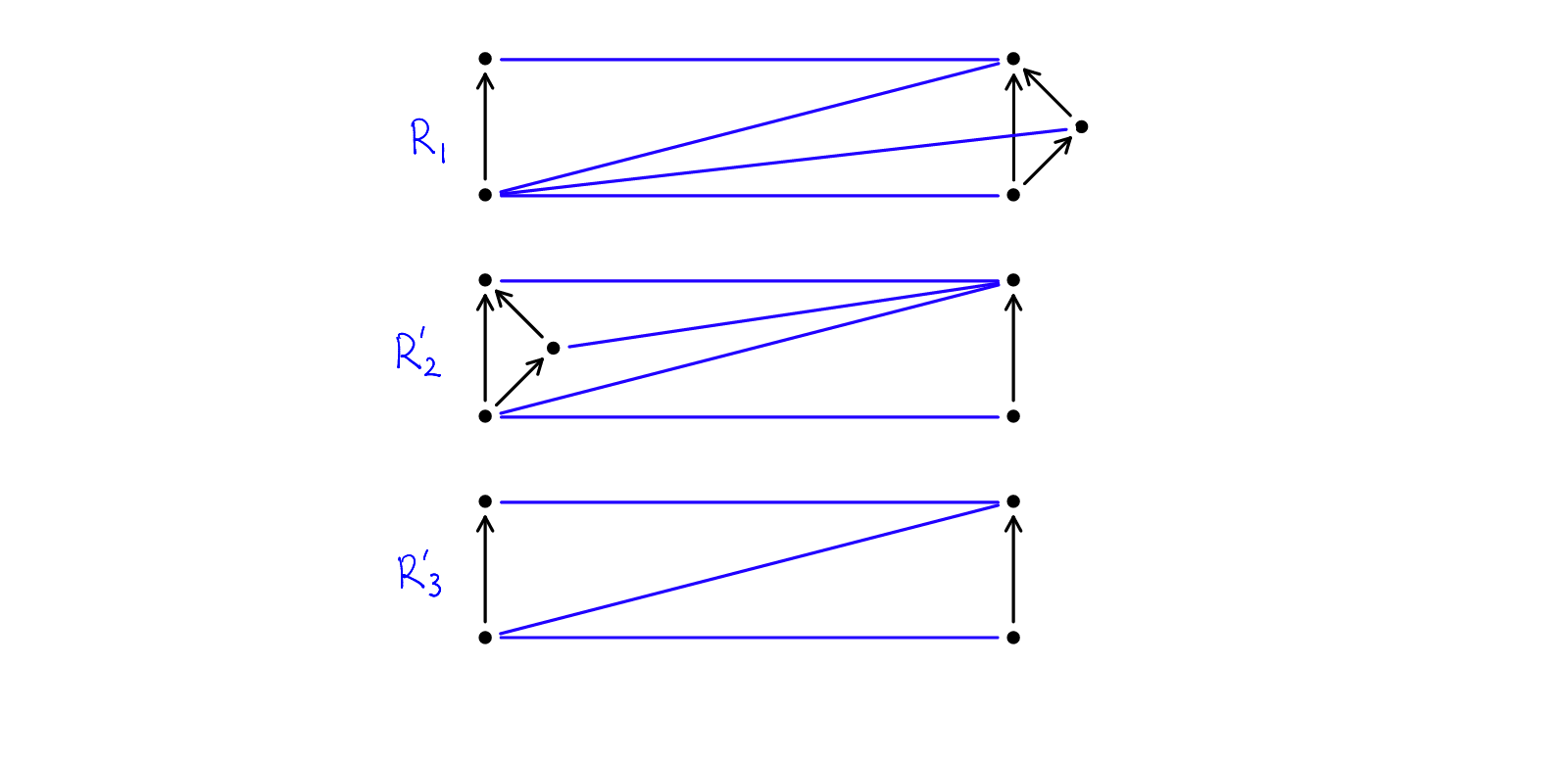}
\endgroup\end{restoretext}
Once more, note that $R'_3 = R_1 \odot R'_2$. 

Finally, we give two non-examples of profunctorial relations
\begin{restoretext}
\begingroup\sbox0{\includegraphics{test/page1.png}}\includegraphics[clip,trim=0 {.25\ht0} 0 {.25\ht0} ,width=.8\textwidth]{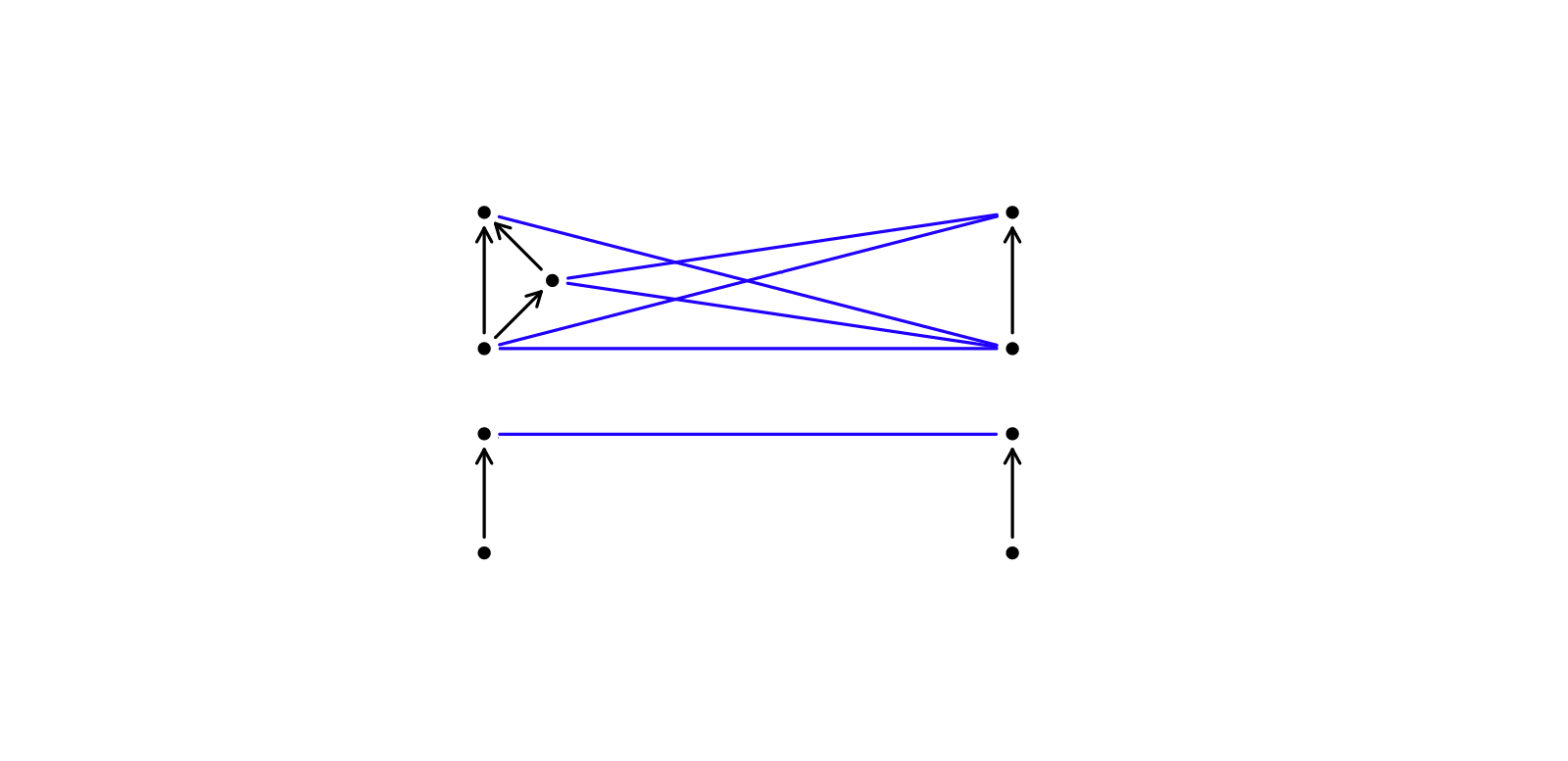}
\endgroup\end{restoretext}
However, in both cases the reader should check that there is one edge which, if added to the relation, would make the relation profunctorial.
\end{eg}
\noindent We conclude this section by remarking that, just as profunctors can be obtained as conjoints (or companions) from functors (cf. \cite{shulman2008framed}), profunctorial relations can be obtained from functors of posets. Explicitly this means the following.

\begin{defn}[Profunctorial relations from functors of posets] \label{defn:graph_of_functors} Given a functor of posets $F: X \to Y$, we define $\grph_F$ to be the profunctorial relation obtained as $\Hom_Y(F-,-)$, that is, the composite
\begin{equation}
X \times Y \xto {F \times \id} Y \times Y \xto {\Hom_Y} \SetCat \xto {\iP\inv} \Bool
\end{equation}
Here, we used that $\Hom_Y$ only takes values in $\im(\iP)$. $\grph_F$ is called the \textit{graph of $F$}. Similarly, for a function of sets $f : X \to Y$ we define $\grph_f$ to be the relation $\grph_{\discr(f)}$ which is the usual \textit{graph of a function} $f$.
\end{defn}

\subsection{Families of profunctorial relations} 

\noindent We recall (and extend terminology of) the construction in\autoref{ssec:sum_fam_bun_prel}. This is a Grothendieck construction for profunctorial relations which can be seen as a (de-categorified) analogue of the Grothendieck construction for profunctors. We can also use this definition to recover the usual discrete Grothendieck construction. 

\begin{constr}[Total poset of a family of profunctorial relations] \label{defn:grothendieck_construction}  Let $\cC$ be a category.
\begin{enumerate}
\item A functor $F : \cC \to \PRel$ is called an ($\cC$-indexed) \textit{family of profunctorial relations}. We define its \textit{total category} $\sG(F)$ as well as its \textit{(Grothendieck) bundle} $\pi_F$ \textit{over $\cC$} as follows: $\sG(F)$ is a poset with objects being pairs $(x \in \cC, a \in F(c))$. Morphisms of $\sG(F)$ are defined by
\begin{equation}
\Hom((x,a),(y,b)) = \Set{(f,a,b) ~|~ (f : x \to y) \in \mor(\cC) , F(f)(a,b)}
\end{equation}
and compose as
\begin{equation}
(g,b,c) \circ (f,a,b) = (gf,a,c)
\end{equation}
Recall that the bundle $\pi_F$ over $\cC$ is the functor $\pi_F : \sG(F) \to \cC$ defined to map $(x,a)$ to $x$.

Note that if $\cC$ is a poset then so is $\sG(F)$ (in which case it is called \textit{total poset}) which follows $F$ being a family of profunctorial relations, cf. \eqref{eq:profunctorial_relation}. 

\item Given two Grothendieck bundles $\pi_{F_1}, \pi_{F_2}$ over $\cC$, a \textit{map of bundles} $G : \pi_{F_1} \to \pi_{F_2}$ is a functor $G : \sG(F_1) \to \sG(F_2)$ such that $\pi_{F_2} G = \pi_{F_1}$.
\end{enumerate}
\end{constr}

\begin{rmk}[Poset-indexed families] In later chapters chapters (with an exception in \autoref{ch:globes}) all our families will be indexed by posets, for which we usually use the letters $X,Y,Z$ (in place of $\cC$, $\cD$).
\end{rmk}

\begin{rmk}[Fibers of bundles]\label{rmk:fiber_identification} Based on the above explicit construction, note that there is a canonical identification $F(x) \iso \pi\inv_F (x)$ by the functor of posets mapping $a \in F(x) \mapsto (x,a) \in \sG(F)$. 
\end{rmk}

\begin{notn}[Projecting to values in fiber] \label{notn:total_poset_projs} Let $F : \cC \to \PRel$ and $p \in \sG(F)$. By the preceding \autoref{defn:grothendieck_construction} we have $p = (x,a)$ for a unique tuple of $x \in \cC$ and $a \in F(x)$. This satisfies
\begin{equation}
\pi_F p = x
\end{equation}
We further denote
\begin{equation}
\secp p := a
\end{equation}
called the \textit{fiber value} of $p$. In other words, $p = (\pi_F p, \secp p)$. Given a functor $f : \cD \to \sG(F)$ for conciseness of notation we often write $\secp f (y)$ instead of $\secp {f(y)}$.
\end{notn}

\begin{eg}[Total posets and functors of bundles] We define a functor $F : \bnum{3} \to \PRel$ using the definition from \autoref{eg:prel} as follows
\begin{restoretext}
\begingroup\sbox0{\includegraphics{test/page1.png}}\includegraphics[clip,trim=0 {.35\ht0} 0 {.29\ht0} ,width=.8\textwidth]{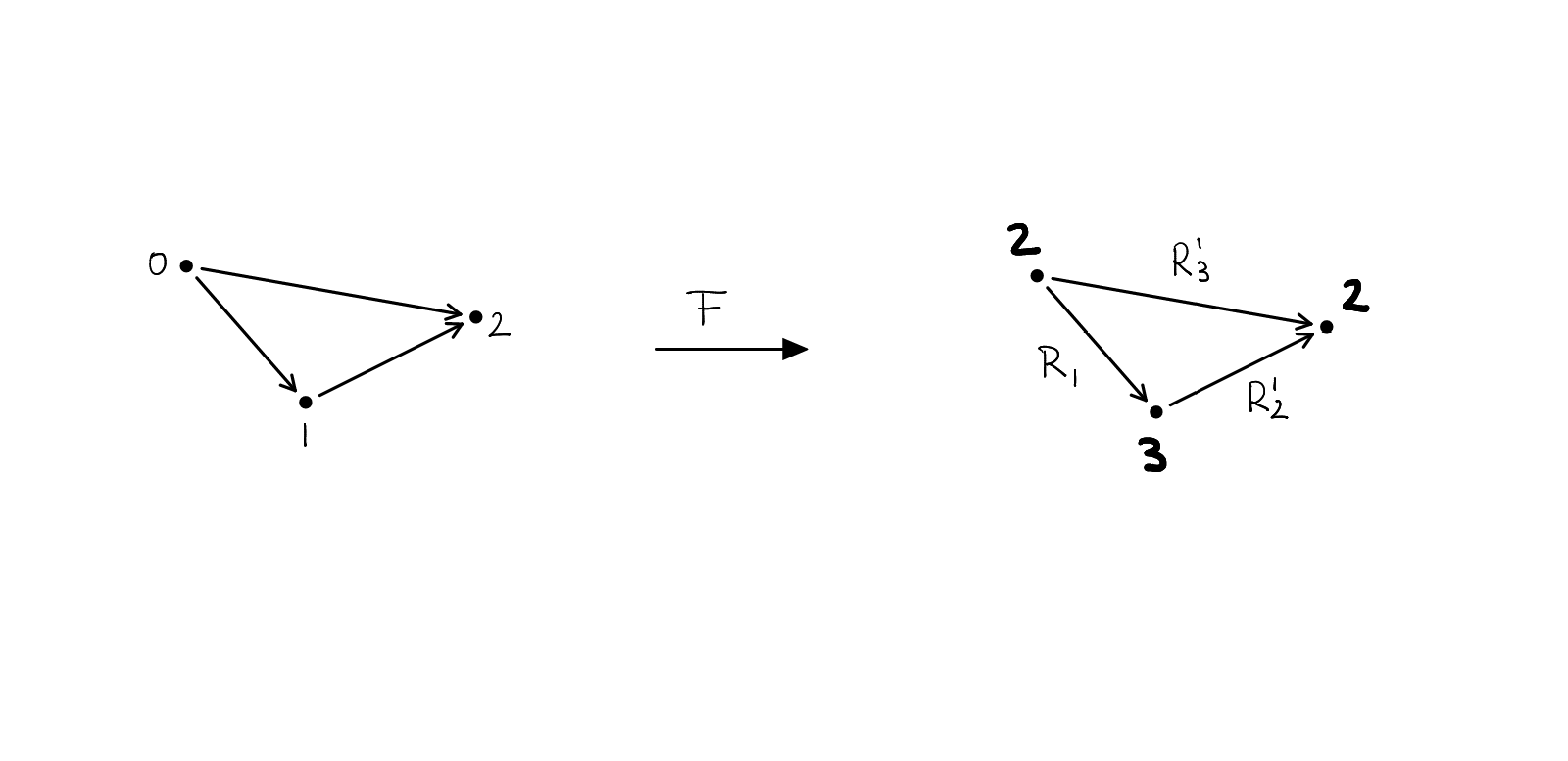}
\endgroup\end{restoretext}
The Grothendieck construction then yields a bundle $\pi_F : \sG(F) \to \bnum{3}$ defined by the data
\begin{restoretext}
\begingroup\sbox0{\includegraphics{test/page1.png}}\includegraphics[clip,trim=0 {.12\ht0} 0 {.04\ht0} ,width=.8\textwidth]{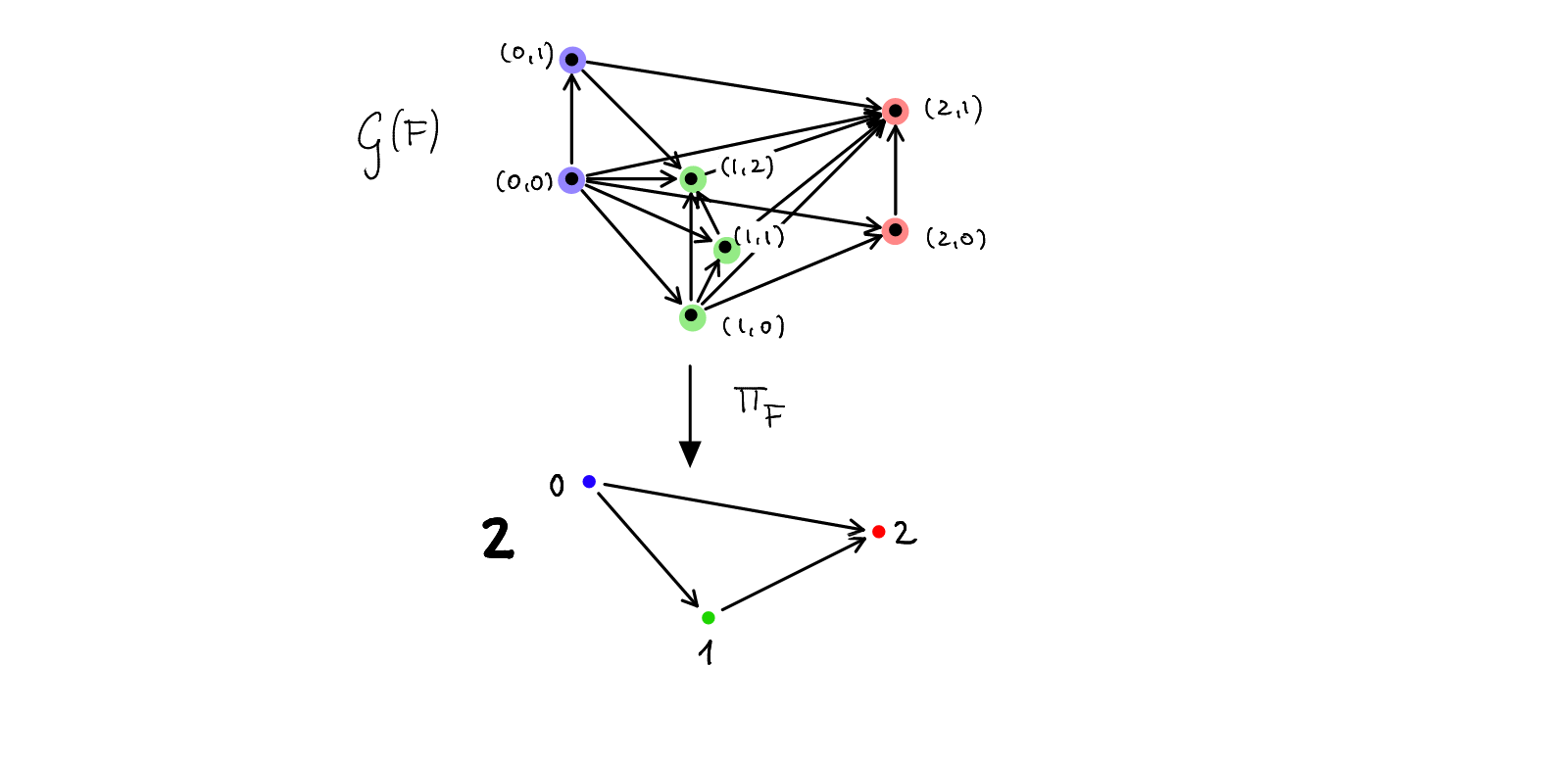}
\endgroup\end{restoretext}
As another example consider
\begin{restoretext}
\begingroup\sbox0{\includegraphics{test/page1.png}}\includegraphics[clip,trim=0 {.35\ht0} 0 {.3\ht0} ,width=.8\textwidth]{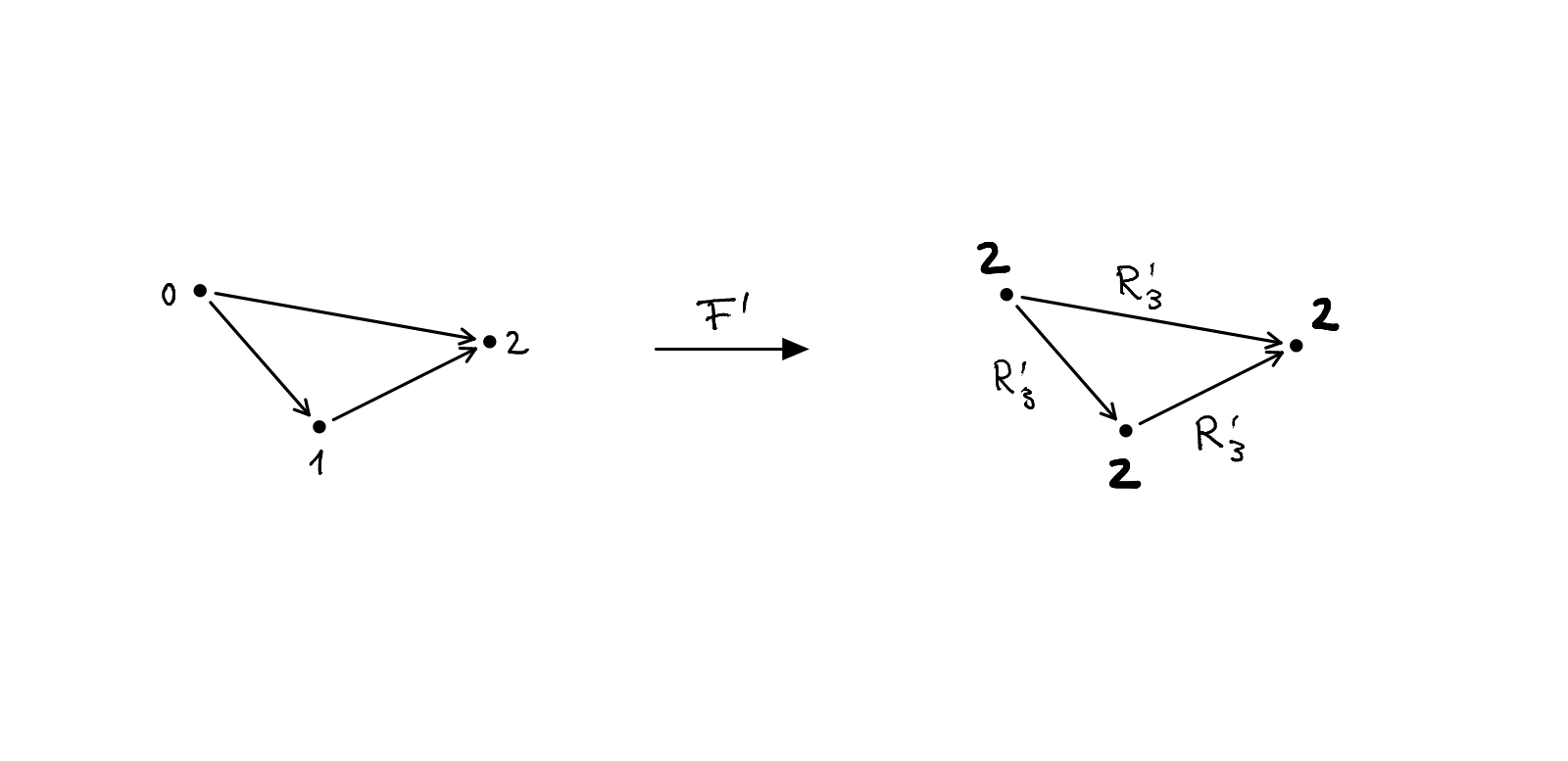}
\endgroup\end{restoretext}
defining $F' : \bnum{3} \to \PRel$, whose Grothendieck construction yields
\begin{restoretext}
\begingroup\sbox0{\includegraphics{test/page1.png}}\includegraphics[clip,trim=0 {.1\ht0} 0 {.1\ht0} ,width=.8\textwidth]{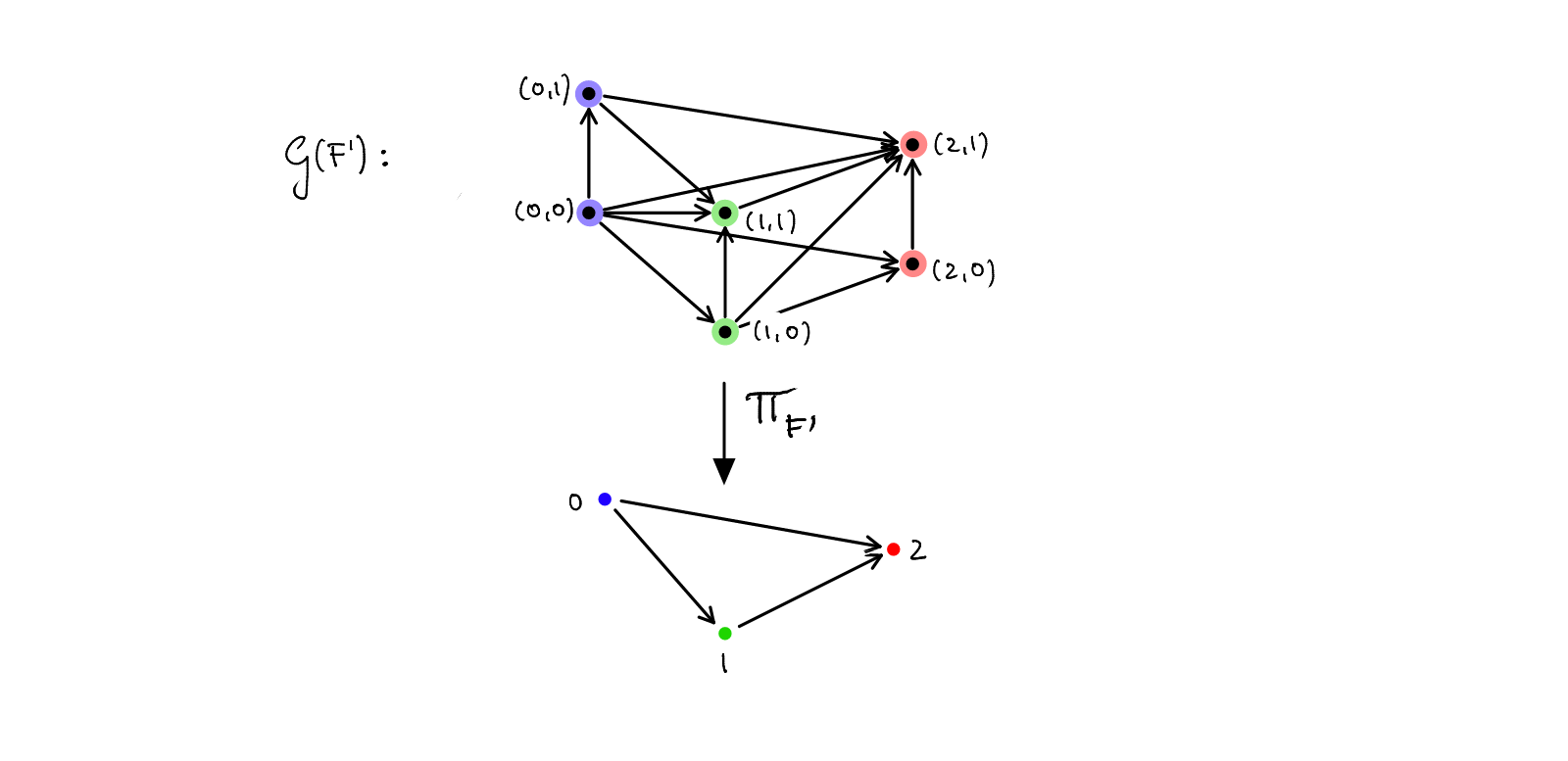}
\endgroup\end{restoretext}
Then there is a map of bundles $G : \pi_F \to \pi_{F'}$ defined as follows
\begin{restoretext}
\begingroup\sbox0{\includegraphics{test/page1.png}}\includegraphics[clip,trim=0 {.05\ht0} 0 {.05\ht0} ,width=.8\textwidth]{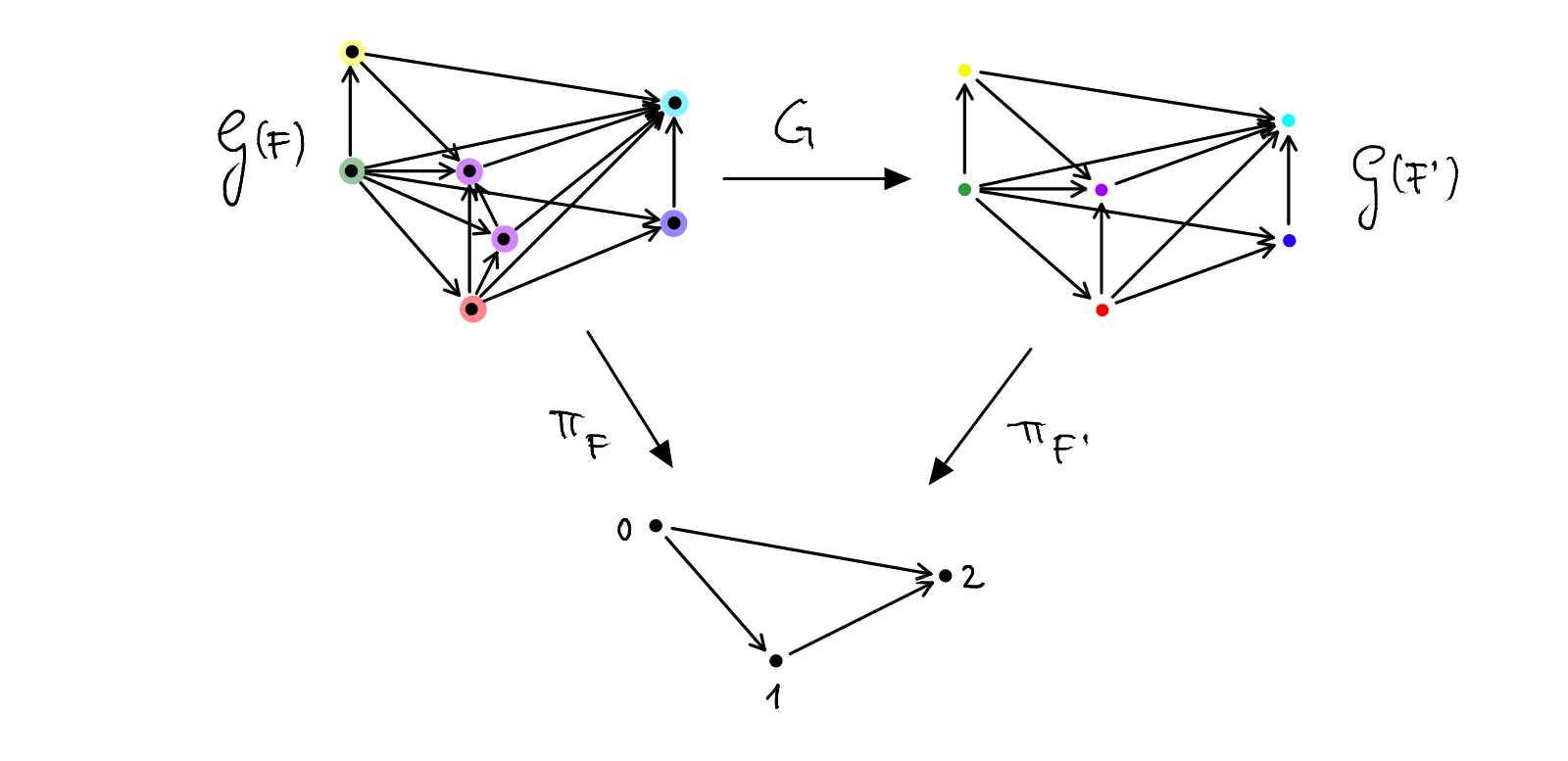}
\endgroup\end{restoretext}
Note that here we only defined $G$ on objects, which (cf. \autoref{rmk:poset_data_objects}) fully determines a functor of posets on morphisms.
\end{eg}

\noindent The following remark compares the usual discrete Grothendieck construction for functors $F : Z \to \SetCat$ (where $Z$ is a poset) to the Grothendieck construction that we defined here.

\begin{rmk}[Comparison to discrete Grothendieck construction] \label{rmk:set_valued_grothendieck_construction}
We first observe that there is a functor $\Discr : \SetCat \to \PRel$: Namely, $\Discr$ maps a set $X$ to the discrete poset $\discr(X)$, and a function $f : X \to Y$ of sets $X,Y$ to the profunctorial relation $\grph_f$. Given a poset $Z$, and a functor $F : Z  \to \SetCat$, we can thus apply the Grothendieck construction defined above to the functor $\Discr F : Z \to \PRel$ to obtain a bundle $\pi_{\Discr F} : \sG(\Discr F) \to Z$. This recovers the usual discrete Grothendieck construction for $\SetCat$-valued functors.
\end{rmk}

\subsection{Base change}

From now on we assume all our bundles to be indexed by posets (with a single exception in \autoref{ch:globes}). Given a classifying map $F$ of a bundle over $Y$, and a map $H$ into $Y$, then $F$ can be precomposed with $H$, and $FH$ is usually said to be obtained from $F$ by pullback along the base change $H$. The bundles build from families of profunctorial relations also have a notion of base change. Since we will work with explicit representations of posets (cf. \autoref{sec:notes}), we will give an explicit base change construction in this section. First, since we will frequently be working with bundles and fibers we introduce the following notation.

\begin{notn}[Restrictions and fiber restrictions] \label{notn:fiber_restrictions} Given families of profunctorial relations $F : X \to \PRel$, $G : Y \to \PRel$ and a commuting square
\begin{align} \label{eq:restr_notn}
\xymatrix@R=0.4cm{ \sG(F) \ar[r]^{K} \ar[d]_{\pi_F} & \sG(G) \ar[d]^{\pi_G} \\
X \ar[r]_{H} & Y}
\end{align}
then for $x \in X$ we denote by 
\begin{equation}
\rest K x : F(x) \to G(H(x))
\end{equation}
the canonical functor induced by $K$ on fibers, mapping $a$ to $\secp{K}(x,a)$ (cf. \autoref{notn:total_poset_projs}).
\end{notn}

\begin{constr}[Base-change of bundles] \label{defn:grothendieck_base_change}  Let $X,Y$ be posets, and assume families of profunctorial relations $F : X \to \PRel$, $G : Y \to \PRel$, $H : X \to Y$ such that $F = GH$. In this case, we say $F$ is the \textit{pullback family} of $G$ along $H$, and $H$ is called a \textit{base-change}. There is a functor $\sG(H) : \sG(F) \to \sG(G)$ defined by mapping $(x,a)$ to $(H(x),a)$. Note that, functoriality of this assignment follows since by \autoref{defn:grothendieck_construction} and functoriality of $H$ we have
\begin{align}
(x,a) \to (y,b) \in \mor(\sG(F)) &\iff x \to y \in \mor(X) \text{~and~} F(x \to y)(a,b) \\ 
&\quad \imp ~H(x) \to H(y) \in \mor(Y) \text{~and~}  G(H(x\to y))(a,b) \\
&\iff (H(x),a) \to (H(y),b) \in \mor(\sG(G))
\end{align}
More concisely, $\sG(H)$ is defined as a leg of the pullback (in the category $\Pos$) along $H$ as follows
\begin{align} \label{eq:grothendieck_base_change}
\xymatrix@R=0.4cm{ \sG(F) \ar[r]^{\sG(H)} \ar[d]_{\pi_{F}} \pullback & \sG(G) \ar[d]^{\pi_G} \\
X \ar[r]_{H} & Y }
\end{align}
We leave the verification of the pullback property to the reader. Note that $F,G$ are implicit in the notation $\sG(H)$.

Further, using the explicit definition of $\sG(H)$, one can see that for each $x \in X$ it restricts to the identity $\rest {\sG(H)} x : F(x) \to G(H(x))$ mapping $a$ to $a$ (cf. \autoref{notn:fiber_restrictions})
\end{constr}

\begin{eg}[Base-change of bundles] Further to the previous example, setting $H = \delta^2_0$ and redefining the functor $G : Y \to \PRel$ to equal $FH$, then we obtain $\sG(H)$ as the following map
\begin{restoretext}
\begingroup\sbox0{\includegraphics{test/page1.png}}\includegraphics[clip,trim=0 {.05\ht0} 0 {.1\ht0} ,width=.8\textwidth]{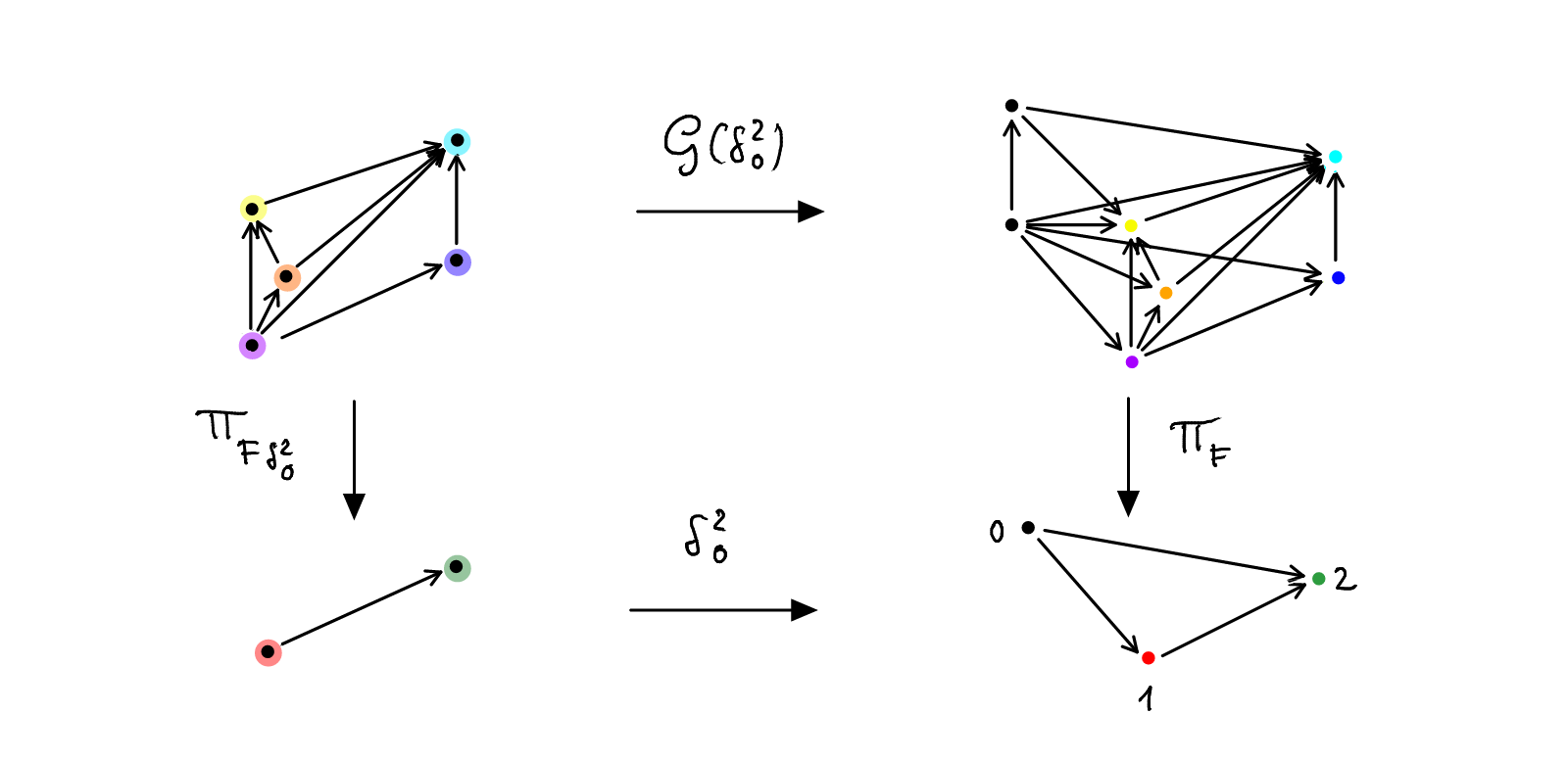}
\endgroup\end{restoretext}\end{eg}

\section{Labelled posets} \label{sec:label_poset}

We now extend our discussion of posets to labelled posets. These are tuples consisting of a poset together with a functor into some (fixed) category $\cC$, which is called the ``labelling functor". The idea of labelled poset structures will be ubiquitous in the rest of this work.

\begin{defn}[The category $\Poss \cC$] \label{defn:posslash} Let $\cC$ be a category. We define the category $\Poss \cC$, called the category of ($\cC$-)\textit{labelled posets} to be the full subcategory of the over-category $\Cat \sslash \cC$ whose image under the forgetful functor is $\Pos$. In other words: objects are tuples $(X,F)$ such that $F: X \to \cC$ is functor from the poset $X$ to $\cC$. $F$ is called the \textit{labelling} of $X$. Morphisms $h : (Y,G) \to (X,F)$ are maps of posets $h : Y \to X$ such that $G = Fh$. The forgetful functor $U : \Poss \cC \to \Pos$ is obtain as the restriction of the forgetful $U : \Cat\sslash{\cC} \to \Cat$. In other words: it maps $h : (Y,G) \to (X,F)$ to $h: Y \to X$.
\end{defn}

\begin{rmk} \label{rmk:terminal_coloring} There is an isomorphism of categories $\Pos$ and $\Pos\sslash_{\!\bnum{1}}$ where $\bnum{1}$ is the terminal category.  Denoting by $\bang : \cC \to \bnum{1}$ the terminal functor to $\bnum{1}$, this isomorphism identifies $X$ and $(X,\bang)$.
\end{rmk}

\subsection{Maps with lifts and downward closed subposets}

A central definition in the context of (labelled) posets is that of poset maps with \gls{lifts}. Under many different circumstances, this condition will allow us to construct certain labelling functors which otherwise would not exist in general.

\begin{defn}[Poset maps with lifts]\label{defn:having_lifts} A functor of posets $F : X \to Y$ \textit{has \gls{lifts}} if whenever $F(x) = y$ and $y' \to y \in \mor(Y)$, then there is $x' \in X$ such that $x' \to x \in \mor(X)$ and $F(x') = y'$.
\end{defn}

\begin{eg} $\delta^2_2 : \bnum{2} \to \bnum{3}$ is an example of a functor of posets with \gls{lifts}. As a non-example, neither $\delta^2_1 : \bnum{2} \to \bnum{3}$ nor $\delta^2_0 : \bnum{2} \to \bnum{3}$ has \gls{lifts}. \end{eg}

\begin{defn}[Downward closed subposets] \label{defn:downward_closed_subposet} Let $X$ be a poset and $Y \subset X$ be a subposet of $X$. We say $Y$ is \textit{\gls{downwardclosed}} if whenever $x \to y$, $x \in X$ and $y \in Y$ then $x \in Y$.
\end{defn}

\begin{rmk} \label{eg:image_of_map_with_lifts} Let $H : Z \to X$ be a functor of posets. Then $\im(H)$ is a \gls{downwardclosed} subposet of $X$ if and only if $H$ has \gls{lifts}.\end{rmk}

\begin{rmk}\label{rmk:downward_closed_implies_full} Every \gls{downwardclosed} subposet is full.
\end{rmk}

\begin{claim} \label{claim:pullback_preserves_having_lifts} Let $X,Y$ be posets, and assume functors $F : X \to \PRel$, $G : Y \to \PRel$, $H : X \to Y$ such that $F = GH$. If $H$ has \gls{lifts}, then $\sG(H)$ (as defined in \autoref{defn:grothendieck_base_change}) has \gls{lifts}.
\proof Assume $(y,a) \to (H(x),b) \in \mor(\sG(G))$ and note $(H(x),b) = \sG(H)(x,b)$. Then $(y \to H(x)) \in \mor(Y)$ and $G(y \to H(x))(a,b)$. Since $H$ has \gls{lifts} we find $x'$ such that $x' \to x \in \mor(X)$ and $H(x') = y$. Since $F = GH$ we deduce $F(x' \to x)(a,b)$. This implies $(x',a) \to (x,b) \in \mor(\sG(F))$ and $\sG(H)(x',a) = (y,a)$ as required. \qed
\end{claim}

\subsection{Gluing posets along downward closed subposets} 

We now explicitly construct pushouts in the category $\Poss \cC$ in a special case which is based on the conditions defined in the previous section.

\begin{constr}[Pushouts in $\Poss \cC$ of injective maps with lifts]  \label{constr:pushouts_in_labelled_posets} Assume a diagram
\begin{equation} \label{eq:labelled_poset_pushout_diag_explicit}
\xymatrix{ (Y,G) \ar[r]^{h_1} \ar[d]_{h_2} & (X_1,F_1) \\ (X_2, F_2) & }
\end{equation}
in $\Poss \cC$ such that $h_i : Y \into X_i$ are injective maps of posets with \gls{lifts}. Using \autoref{rmk:downward_closed_implies_full} and \autoref{eg:image_of_map_with_lifts} this implies that $Y \iso \im(h) \subset X_i$ is a \gls{downwardclosed} full subposet of $X_i$. We construct a pushout 
\begin{equation} \label{eq:labelled_poset_pushout_diag}
\xymatrix{ (Y,G) \ar[r]^{h_1} \ar[d]_{h_2} & (X_1,F_1) \ar[d]^{g_1} \\ (X_2, F_2) \ar[r]_-{g_2}& (X_1 \cup_Y X_2, F_1 \cup_Y F_2) \pushoutfar }
\end{equation}
which in fact will also be a pushout in $\Cat \sslash \cC$.

We first construct a pushout in $\Cat$
\begin{equation} \label{eq:labelled_poset_pushout_diag2}
\xymatrix{ Y \ar[r]^{h_1} \ar[d]_{h_2} & X_1 \ar[d]^{g_1} \\ X_2 \ar[r]_-{g_2}& X_1 \cup_Y X_2 \pushoutfar}
\end{equation}
On objects (that is, after applying $\obj : \Cat \to \SetCat$) this square is constructed as the pushout of $\obj(h_1)$ and $\obj(h_2)$ in $\SetCat$. Since pushouts in $\SetCat$ preserve monomorphisms, we will from now on regard $\obj(Y)$, $\obj(X_i)$ as subsets of $\obj(X_1 \cup_Y X_2)$
\begin{equation}
\obj(Y) = (\obj(X_1) \cap \obj(X_2)) \subset \obj(X_i) \subset \obj(X_1 \cup_Y X_2) 
\end{equation}
and thus identify $h_i(y) = y$, $g_i(x) = x$ in our notation.

Using this convention, the morphism set of $X_1 \cup_Y X_2$ is defined by
\begin{gather} \label{eq:morphism_set_quotient}
(x \to y) \in \mor(X_1 \cup_Y X_2) \iff \exists i \in \Set{1,2}.~ (x \to y) \in \mor(X_i)
\end{gather}
Note that this definition makes $g_i : X_i \subset X_1 \cup_Y X_2$ into a full subposet.

We now show the pushout property in $\Cat$. For this assume $G, F_1$ and $F_2$ as shown in \eqref{eq:labelled_poset_pushout_diag}. We then need to show there is a unique $F_1 \cup_X F_2$. This is constructed as follows.

We remark that since $\obj(Y) = (\obj(X_1) \cap \obj(X_2))$, and since $Y$ is \gls{downwardclosed} in $X_i$ we have that any composite morphism $(x \to y \to z)$ in $X_1 \cup_Y X_2$ is either fully contained in $X_1$ or in $X_2$ (since either $z \in X_1$ or $z \in X_2$). In other words,
\begin{equation} \label{eq:downward_closed}
(x \to y \to z) \in \mor(X_1 \cup_Y X_2) \imp \exists i \in \Set{1,2}.~ (x \to y \to z) \in \mor(X_i)
\end{equation}
As a consequence of the last observation, we can define a functor $F_1 \cup_Y F_2 : X_1 \cup_Y X_2 \to \cC$ by defining it to map morphisms $(x \to y) \in \mor(X_1 \cup_Y X_2)$ as follows
\begin{equation}
(F_1 \cup_Y F_2)(x \to y) = F_i(x \to y) \text{~if~} (x \to y) \in \mor(X_i)
\end{equation}
This is well-defined by our assumption \eqref{eq:labelled_poset_pushout_diag_explicit}, and functorial by observation \eqref{eq:downward_closed} and $F_i$ being functorial.  Conversely, it is the only possible choice of values for $F_1 \cup_Y F_2$ by commutativity of \autoref{eq:labelled_poset_pushout_diag}. 

The claim that this is a pushout in $\Cat \sslash \cC$ follows from the fact that $U : \Cat \sslash \cC \to \Cat$ creates limits. Since $\Pos \subset \Cat$ is a full subcategory it follows that our construction also gives a pushout in  $\Poss \cC$.
\end{constr}

\begin{eg} Consider the following maps $h_i : (Y,G) \to (X_i,F_i)$ ($i \in \Set{1,2}$) in $\Poss \cC$ (cf. \autoref{notn:depicting_posets})
\begin{restoretext}
\begingroup\sbox0{\includegraphics{test/page1.png}}\includegraphics[clip,trim=0 {.0\ht0} 0 {.0\ht0} ,width=\textwidth]{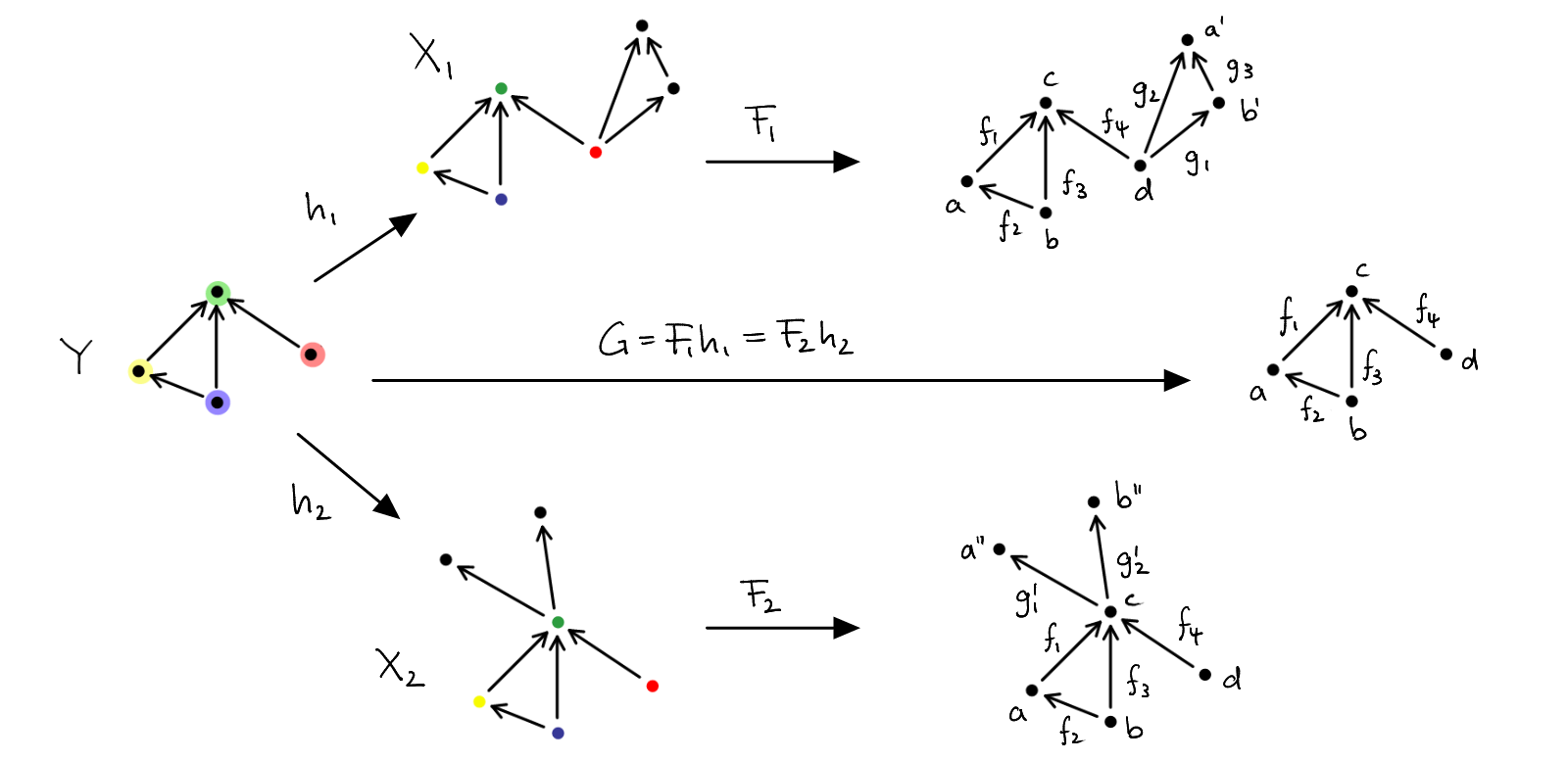}
\endgroup\end{restoretext}
Note that $h_i$ are indeed maps with \gls{lifts}. Then the pushout construction yields an object $(X_1 \cup_Y X_2, F_1 \cup_Y F_2)$ with the data
\begin{restoretext}
\begingroup\sbox0{\includegraphics{test/page1.png}}\includegraphics[clip,trim=0 {.2\ht0} 0 {.25\ht0} ,width=.8\textwidth]{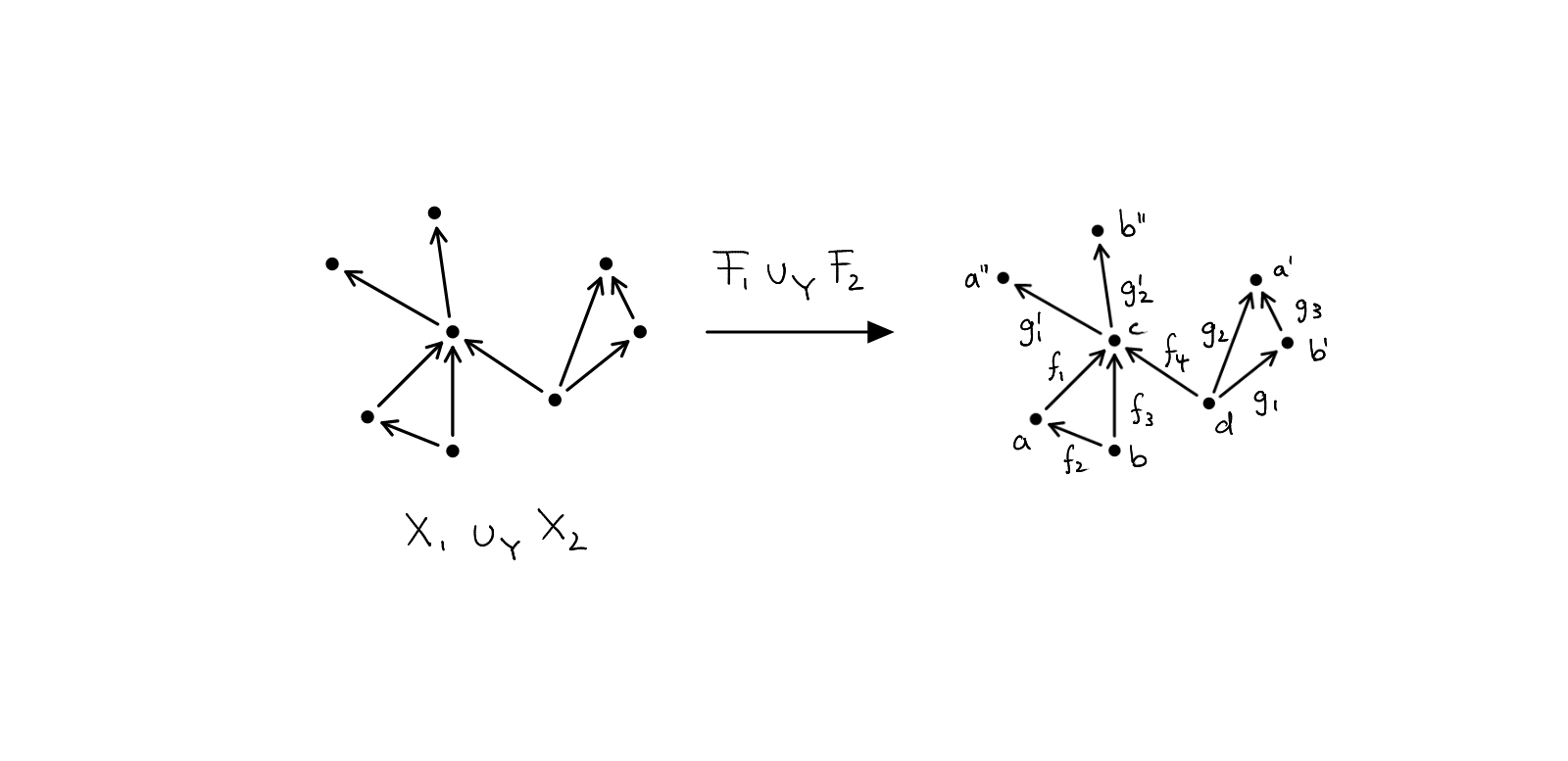}
\endgroup\end{restoretext}
\end{eg}

\subsection{Necessity of lift condition}

The reader might wonder whether the condition of having \gls{lifts} is needed for the previous construction to yield a pushout. The following example shows that it is indeed.

\begin{eg} Consider the following maps $h_i : (Y,G) \to (X_i,F_i)$ ($i \in \Set{1,2}$) in $\Poss \cC$
\begin{restoretext}
\begingroup\sbox0{\includegraphics{test/page1.png}}\includegraphics[clip,trim=0 {.1\ht0} 0 {.0\ht0} ,width=.8\textwidth]{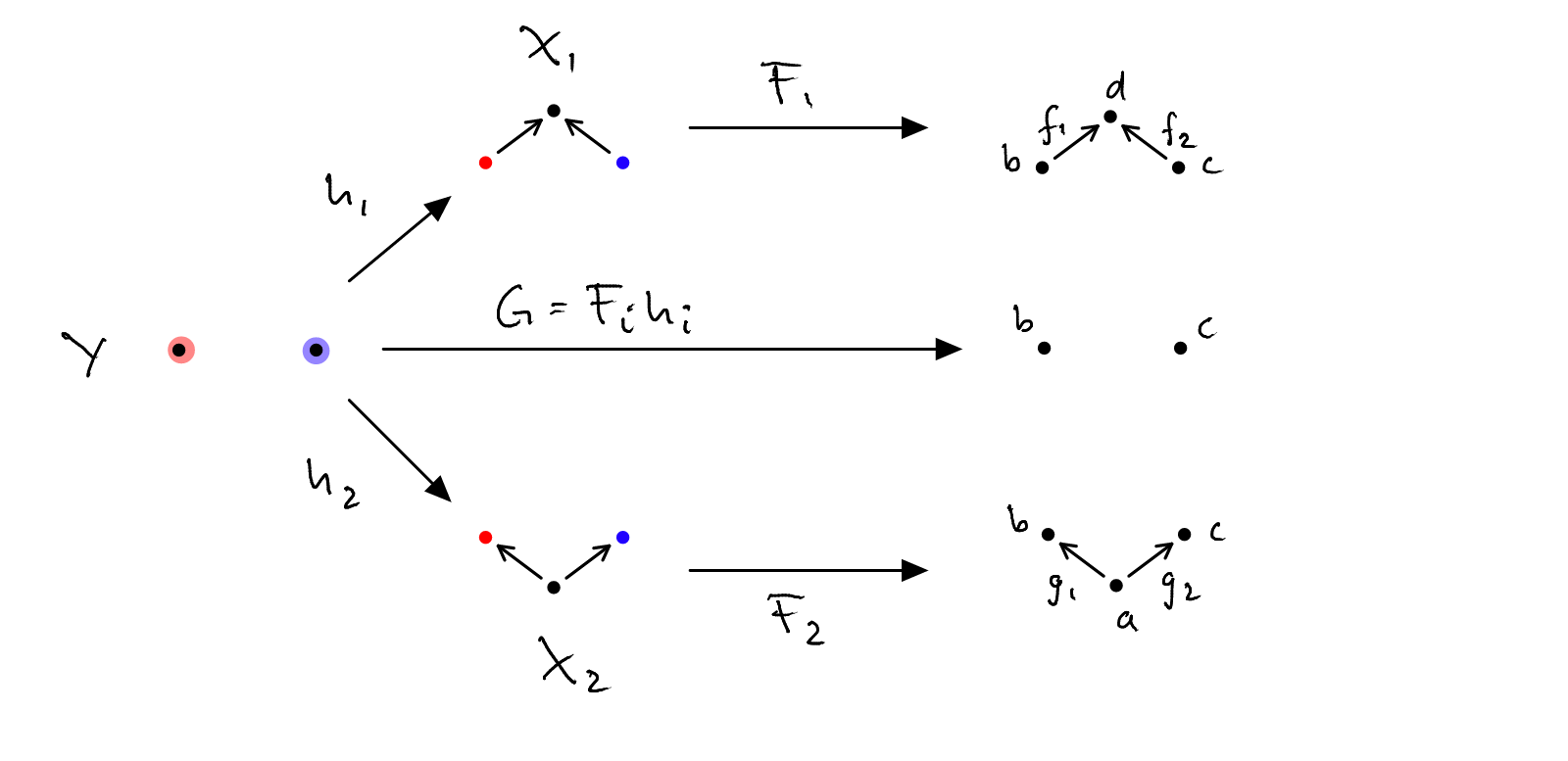}
\endgroup\end{restoretext}
such that $f_1 g_1 \neq f_2 g_2$. Then the previous construction fails, and we see that pushouts in $\Poss \cC$ do not always exist.
\end{eg}

\chapter{Intervals} \label{ch:intervals}

In this chapter we lay the foundations of singular interval theory. Singular intervals are posetal models of (finite) stratifications of the interval $[0,1]$ with many surprising combinatorial properties. They are naturally organised into a category $\SI$, which admits a functor $\SiR$ into $\PRel$. This functor associates to each morphism in $\SI$ a profunctorial relation, whose total poset can be interpreted as a triangulation of the square $[0,1]^2$. These definitions are discussed in  \autoref{sec:SI}. In the next chapter, $n$-dimensional space will be build by inductively bundling intervals over $(n-1)$-dimensional space. The ``interval bundles" used in each such step will be defined in \autoref{sec:SI_fam}.

\section{Singular intervals, their morphisms, and associated relations} \label{sec:SI}

\subsection{The definition of singular intervals}

A singular interval $I$ is a finite ``zig-zag" poset together with a second order on its objects called the ``direction order". We will define singular intervals by an explicit representation in terms of sets of natural numbers. The ordering of numbers will represent the direction order. We make this choice because, firstly, it is amenable to calculation and computer implementation. Secondly, fixing a representation by numbers allows us to understand equivalence of singular intervals as strict equality. And finally, using numbers to represent the second order avoids having to define a new categorical structure with two orders, and  allows us to keep working with posets instead.

\begin{defn}[Singular intervals]\label{defn:singular_intervals}
A \textit{singular interval} $I$ is a poset with objects being a set of natural numbers
\begin{equation}
\obj(I) = \Set{0,1,2, \ldots, 2\iH_I-1, 2\iH_I}
\end{equation}
for some $\iH_I \in \lN$. $\iH_I$ is called the \textit{height} of $I$. The poset order $\to$ of $I$, also called the \textit{singular order} on $I$, is generated by the morphisms $a \to a+1$ if $a \in I$ is even, and $a \ot a+1$ if $a \in I$ is odd. The unique interval of height $n \in \lN$ is denoted by $\singint n$. $\singint 0$ is called the \textit{initial singular interval}, $\singint 1$ is called the \textit{terminal singular interval}.

We further define the subset
\begin{equation}
\singcont(I) := \Set{i\in I ~|~ i \text{~odd}} \subset \obj(I)
\end{equation}
whose elements are called \textit{singular heights}, and the subset
\begin{equation}
\regcont(I) :=\Set{i\in I ~|~ i \text{~even}} \subset \obj(I)
\end{equation}
whose elements are called \textit{regular segments}. For technical convenience, we also define
\begin{equation}
\extsing(I) := \Set{-1} \cup \singcont(I) \cup \Set{2\iH_I + 1} \subset \lZ
\end{equation}
\end{defn}

\begin{rmk}[Singular order and direction order]\label{notn:monotonicity1} Note that the objects of $I$ are integers and thus determine a full subposet of $\lZ$, which (abusing notation when no confusion arises) we again denote by $I$. The subposet order $\leq$ on $I \subset \lZ$ is called the \textit{direction (order)} on $I$ (cf. \autoref{eg:posets}). This is different from the above poset order $\to$ on $I$, which we call the singular order.
\end{rmk}

Examples illustrating the terminology and notation for singular intervals include the following.

\begin{egs}[Singular intervals]The following depicts the posetal structure and object naming convention  of $\singint n$ for $0 \leq n \leq 4$
\begin{restoretext}
\begingroup\sbox0{\includegraphics{test/page1.png}}\includegraphics[clip,trim=0 {.05\ht0} 0 {.1\ht0} ,width=\textwidth]{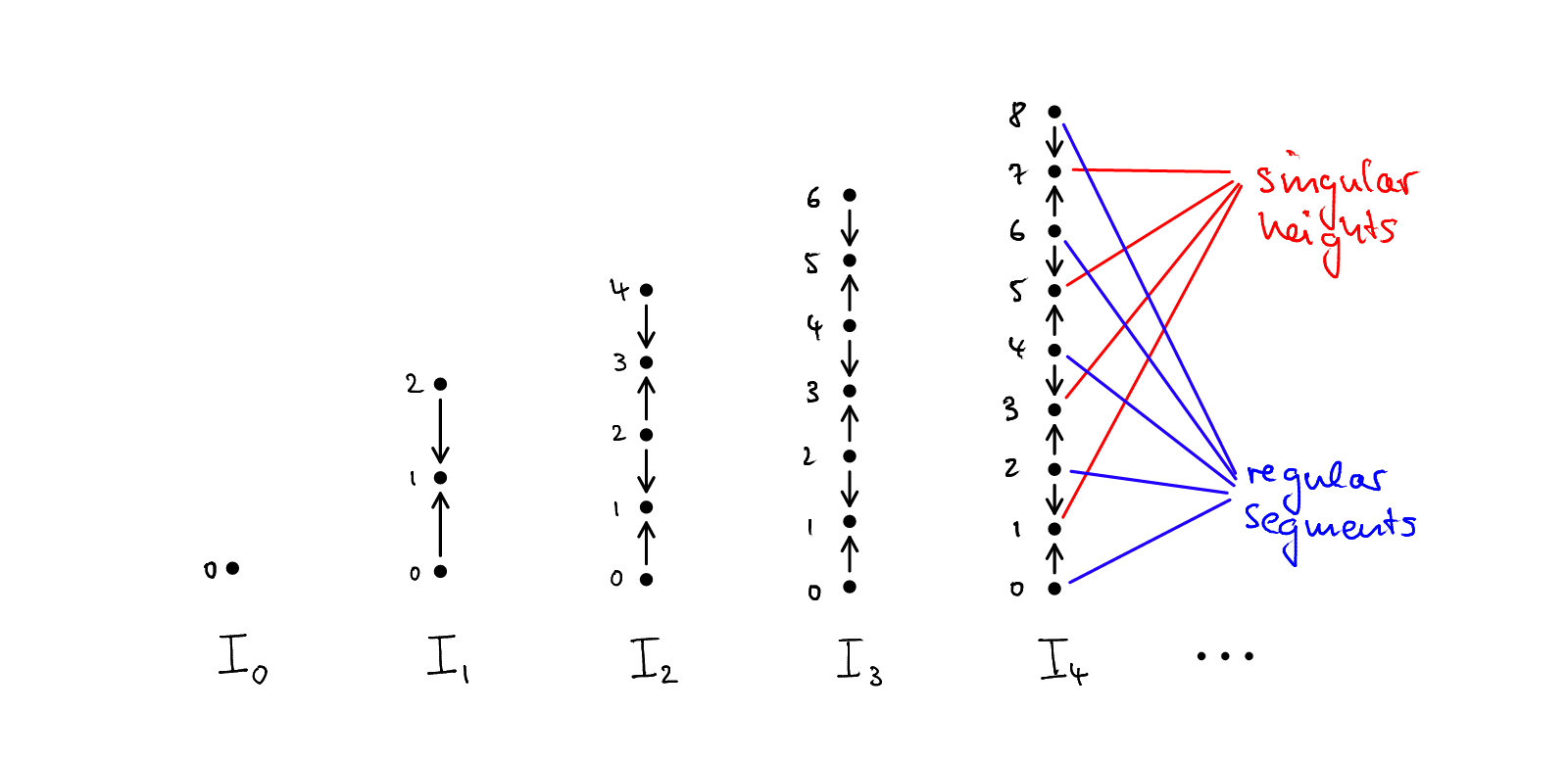}
\endgroup\end{restoretext}
\end{egs}

\begin{notn}[Depicting the direction order of singular intervals]  If $X$ is a poset and $i : X \iso \singint k$ is an isomorphism for some $k$, then $X$ has a single (involutive) non-identity automorphism $\alpha$. Thus there are exactly two isomorphisms $i, i\alpha : X \iso \singint k$. Fixing a choice of either isomorphism, means choosing a ``direction order" on the objects of $X$; that is, objects of $X$ get identified with numbers of $\singint k$. Such a direction order can be fully determined by ``embedding" a depiction of $X$ (cf. \autoref{notn:depicting_posets}) in 1-dimensional space and then giving direction to that space---this then is the direction of increasing numbers for objects. By convention this direction will be from bottom to top and left to right, unless otherwise indicated. Based on this, the following depictions of $\singint 3$ have identical meaning
\begin{restoretext}
\begingroup\sbox0{\includegraphics{test/page1.png}}\includegraphics[clip,trim=0 {.15\ht0} 0 {.2\ht0} ,width=\textwidth]{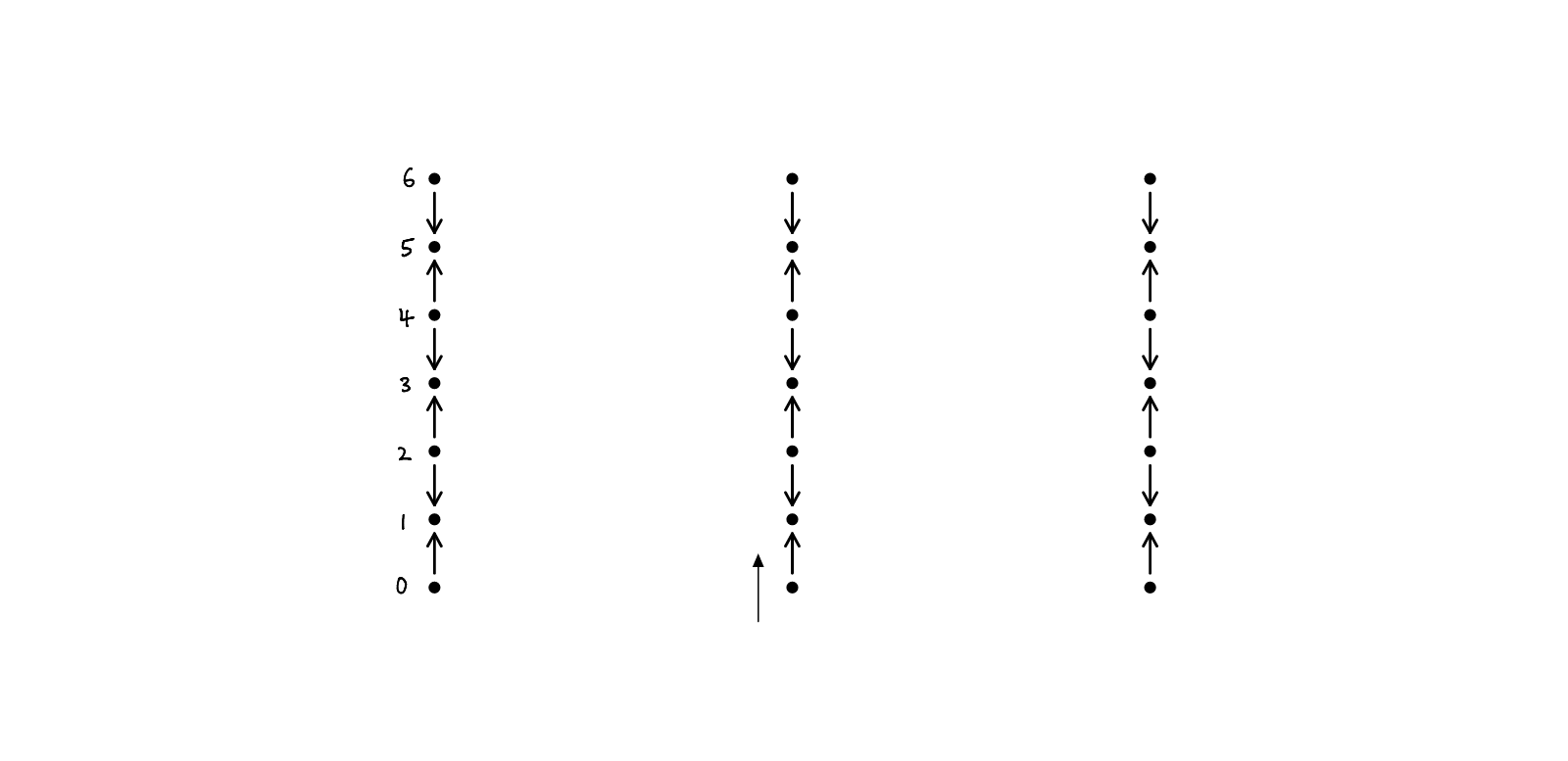}
\endgroup\end{restoretext}
On the left, objects are given by their numbers. In this middle these numbers can be derived from linear arrangement of the given posets together with the direction indicated by a ``coordinate arrow". On the right, these numbers can be derived from linear arrangement and the bottom-to-top/left-to-right convention.
\end{notn}

\begin{defn}[Morphisms of singular intervals] \label{defn:singular_intervals_morphism} Let $I$, $J$ be singular intervals. 
\begin{enumerate}
\item A \textit{(singular-height) morphism} $f$ between singular intervals $I$ and $J$, also called \SI-morphism and denoted by $f: I \to_{\SI} J$, is a monotone map $f: \singcont(I) \to \singcont(J)$. Morphisms of singular intervals compose as functions of sets. Singular intervals and their morphisms form the category $\SI$. 

\item A \textit{regular-segment morphism} $g$ between singular intervals $I$ and $J$, also called an $\SI\regop$-morphism and denoted by $g: I \to\regop_{\SI} J$, is monotone map $g : \regcont(I) \to \regcont(J)$ that also preserves endpoints, i.e. $g(0) = 0$ and $g(2\iH_I) = 2\iH_J$. Regular-segment morphisms of singular intervals compose as functions of sets. Singular intervals and their regular-segment morphisms form a category $\SI^{\mathrm{reg}}$.
\end{enumerate}
\end{defn}

\begin{rmk}[\SI{} is equivalent to finite totally ordered sets] \label{rmk:SI_is_cat_of_tot_ord} There is an embedding $S : \SI \to \Pos$, mapping objects $\singint n$ to $(\bnum{n+1})$ and morphisms $f : \singint n \to\oSI \singint m$ to the poset functor $S(I) : (\bnum{n+1}) \to (\bnum {m+1})$, defined to map objects as
\begin{equation}
i \mapsto \frac{f(2i+1)-1}{2}
\end{equation}
The image of the embedding $S$ is the full subcategory $\BDelta_+$ of finite totally ordered posets in $\Pos$.
\end{rmk}

\begin{egs}[Maps between singular intervals] \label{eg:singint_map} Examples of maps between singular intervals (and their notation) are the following: the maps $f_1 : \singint 2 \to_{\SI} \singint 2$, $f_2 : \singint 3 \to_{\SI} \singint 4$ can be depicted as follows
\begin{restoretext}
\begingroup\sbox0{\includegraphics{test/page1.png}}\includegraphics[clip,trim=0 {.0\ht0} 0 {.13\ht0} ,width=\textwidth]{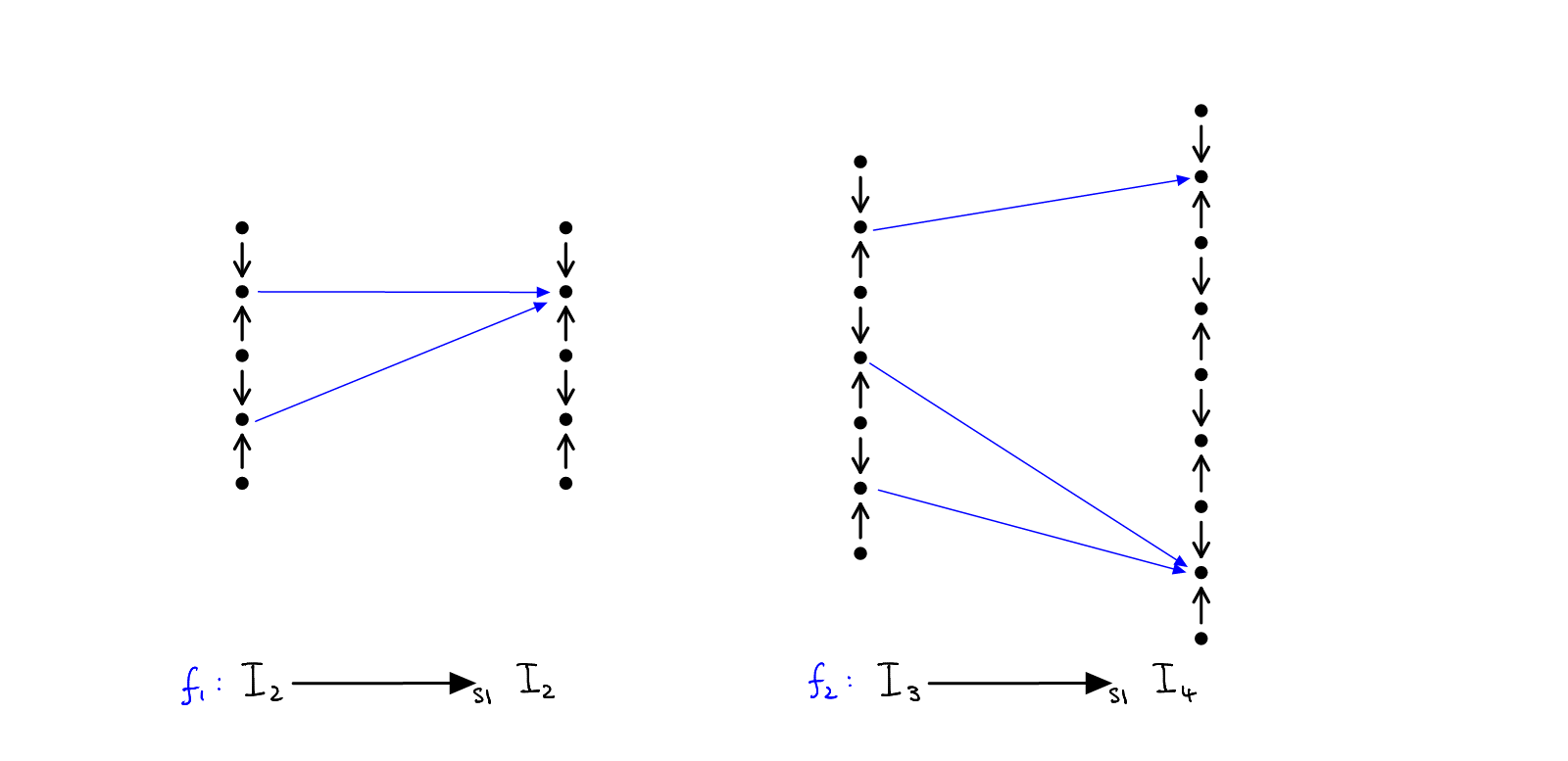}
\endgroup\end{restoretext}
Here, we depicted mappings $f : X \to Y$ between sets $X$, $Y$ by arrows from $a \in X$ to $f(a) \in Y$. Note that $f_2$ cannot be obtained by restricting a functor of posets to the singular heights of its domain and codomain (while $f_1$ can be obtain that way). We further define $f_3 : \singint 1 \to_{\SI} \singint 2$ and $f_4 : \singint 0 \to_{\SI} \singint 2$ as follows
\begin{restoretext}
\begingroup\sbox0{\includegraphics{test/page1.png}}\includegraphics[clip,trim=0 {.2\ht0} 0 {.25\ht0} ,width=\textwidth]{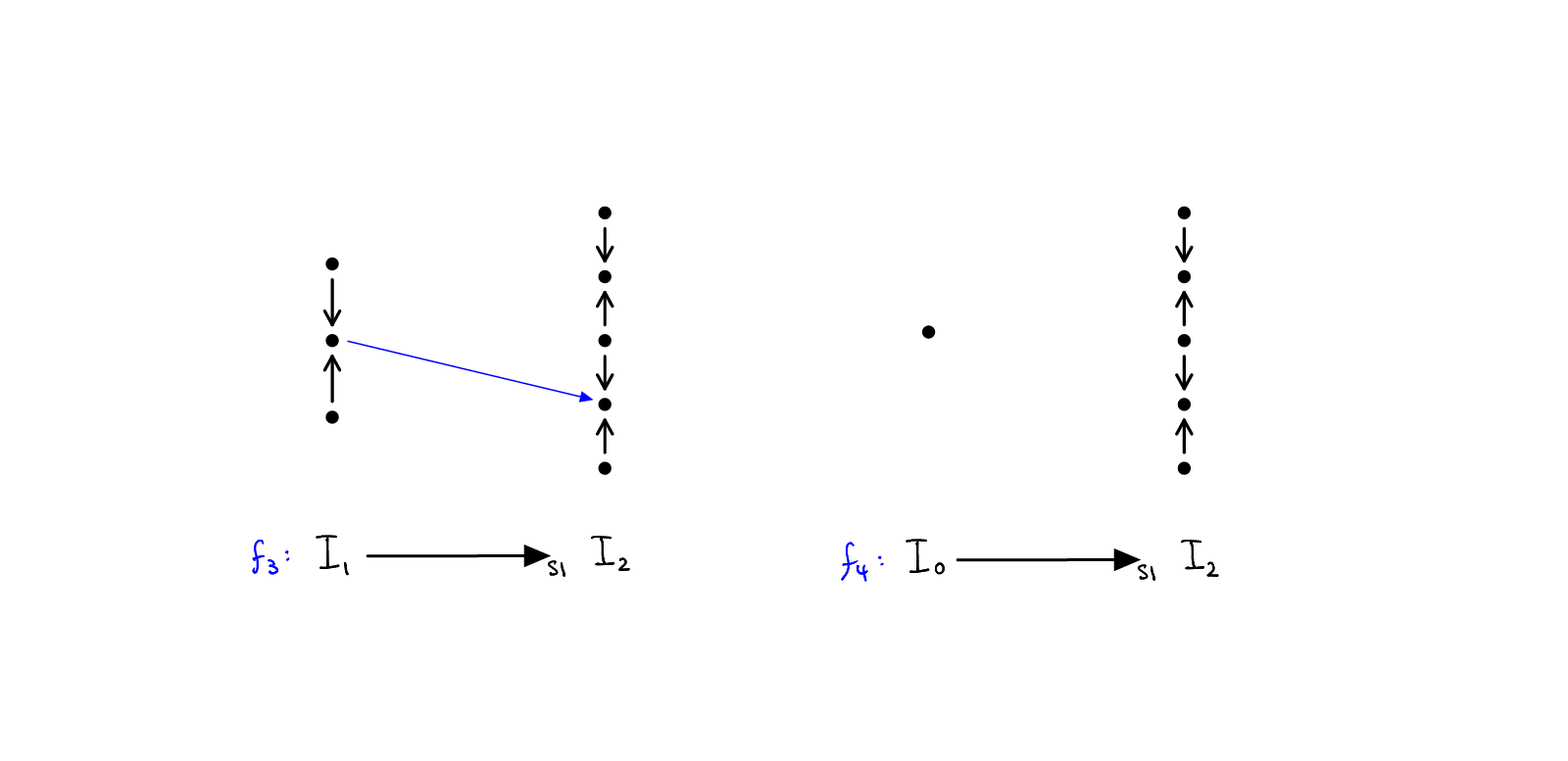}
\endgroup\end{restoretext}
In particular, note that the initial singular interval does not contain any singular heights. On the other hand,
\begin{restoretext}
\begingroup\sbox0{\includegraphics{test/page1.png}}\includegraphics[clip,trim=0 {.3\ht0} 0 {.2\ht0} ,width=\textwidth]{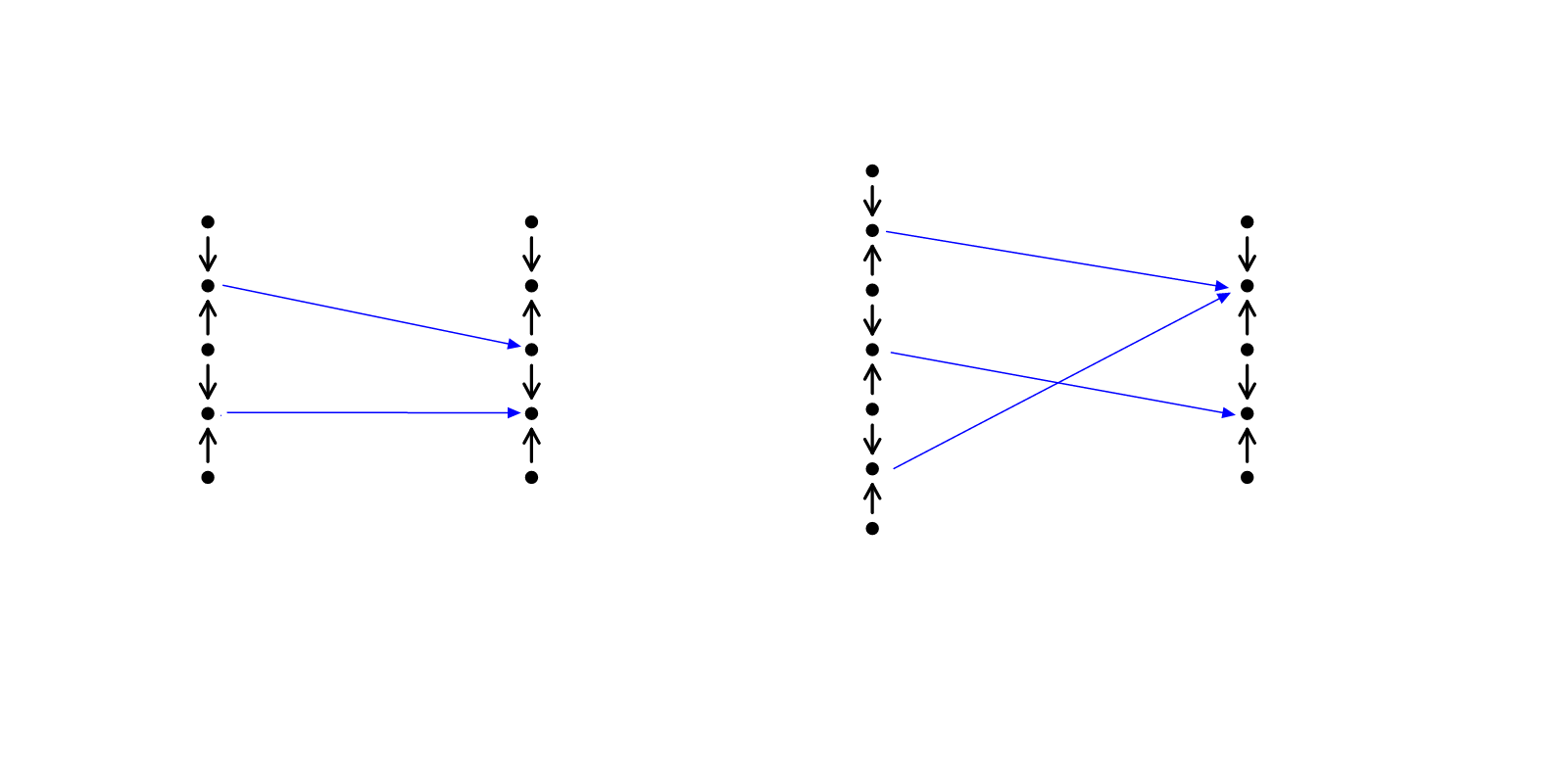}
\endgroup\end{restoretext}
are not examples of singular interval morphisms: the left map maps singular heights to regular segments, while the right map is not monotone.

Similarly, the following
\begin{restoretext}
\begingroup\sbox0{\includegraphics{test/page1.png}}\includegraphics[clip,trim=0 {.05\ht0} 0 {.24\ht0} ,width=\textwidth]{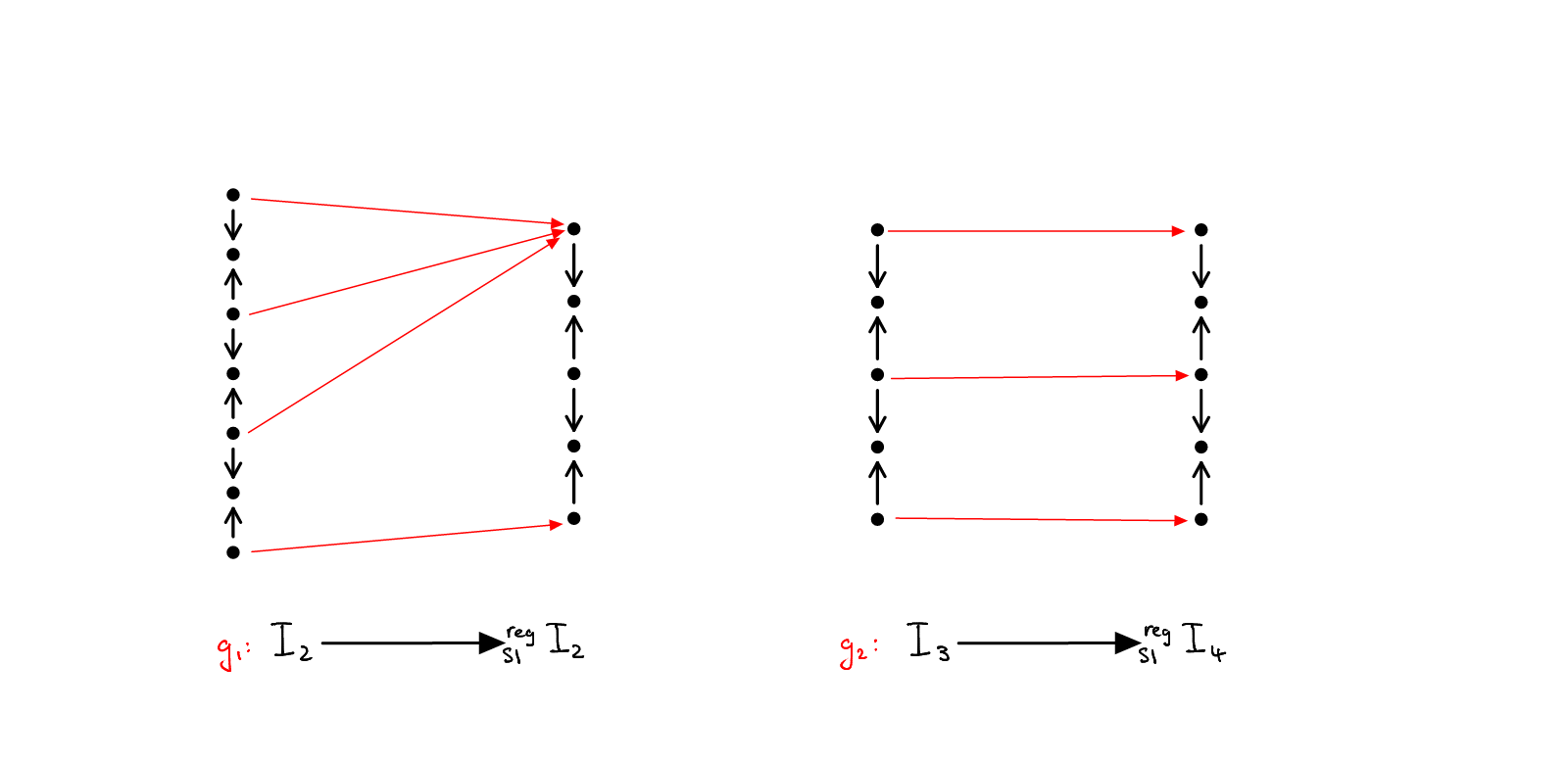}
\endgroup\end{restoretext}
are examples of regular height morphisms.
\end{egs}

\begin{rmk}[Equivalence is strict equality for elements of $\SI$]
\label{rmk:direction_preserving_iso_is_eq} Let $I$, $J$ be singular intervals. We make two observations relating to notions of equivalence that singular intervals could admit.
\begin{itemize}
\item Let $F : I \to J$ be a functor of posets. If $F$ is an isomorphism of posets and monotone (that is, it preserves direction order) then we must have $F = \id$ and $I = J$ by \autoref{defn:singular_intervals}. 
\item Let $f : I \to\oSI J$ be a morphism in $\SI$. If $f$ is an isomorphism in $\SI$, then $f$ is the identity morphism.
\end{itemize}
\end{rmk}

\begin{defn}[Extensions of singular morphisms] \label{notn:singular_morphism_boundary_cases} Given $f: I \to_{\SI} J$, for technical convenience we also define $\wwidehat f : \extsing (I) \to \extsing (J)$, called the \textit{extension} of $f$, by setting $\wwidehat f(-1) = -1$, $\wwidehat f(2\iH_I +1) = 2\iH_J + 1$ and $\wwidehat f (a) = f(a)$ if $a \in \singcont(I)$. 
\end{defn}

\begin{eg} Depicting mappings $f : X \to Y$ between sets $X$, $Y$ as before (by arrows from $a \in X$ to $f(a) \in Y$) we can depict extensions of singular interval morphisms $f_2$ and $f_4$ defined in the previous example as follows
\begin{restoretext}
\begingroup\sbox0{\includegraphics{test/page1.png}}\includegraphics[clip,trim=0 {.1\ht0} 0 {.1\ht0} ,width=\textwidth]{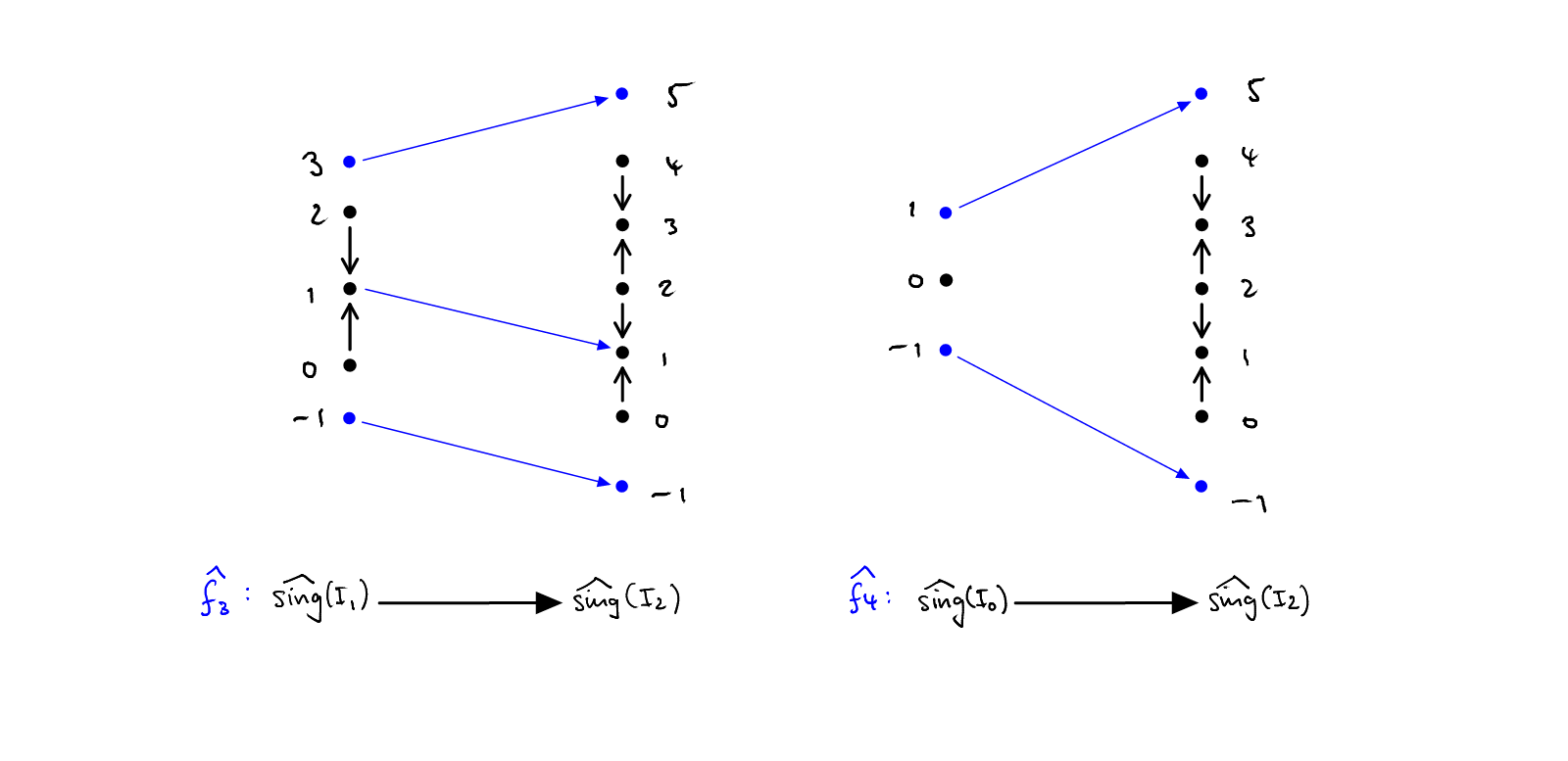}
\endgroup\end{restoretext}
\end{eg}

\subsection{Duality of singular-height and regular-segment morphisms}

Singular intervals have many interesting combinatorial properties. In this section we shine light on one of their most fundamental properties, namely a ``duality relation" between their two types of morphisms.

\begin{constr}[Regular and singular \rsdual{}s{}]\label{defn:regop_singop} To every singular-height morphism $f : I \to_{\SI} J$ we associate a regular-segment morphism $f\regop : J \to\regop_{\SI} I$ as follows: let $b \in \regcont(I)$ and $a \in \regcont(J)$, then
\begin{align} \label{eq:fop_from_f}
f\regop (a) = b \quad \iff \quad \big(\wwidehat f(b - 1) < a < \wwidehat f(b + 1)\big)
\end{align}
$f\regop$ is called the \textit{regular \rsdual{}} of $f$. The fact that this construction defines a monotone function of singular heights as required in \autoref{defn:singular_intervals_morphism} follows from \autoref{rmk:regop_singop_welldef}.

Conversely, to every regular-segment morphism $g : J \to\regop_{\SI} I$ we associate a singular-height morphism $g\singop : I \to_{\SI} J$ as follows. Let $d \in \singcont(J)$ and $c \in \singcont(I)$, then
\begin{align} \label{eq:f_from_fop}
g\singop (c) = d \quad \iff \quad \big(g(d - 1) < c < g(d + 1)\big)
\end{align}
$g\singop$ is called the \textit{singular \rsdual{}} of $g$. The fact that this construction defines a monotone function of regular segments as required in \autoref{defn:singular_intervals_morphism} follows from \autoref{rmk:regop_singop_welldef}.
\end{constr}

\begin{eg}[Regular and singular \rsdual{}s{}] Based on the singular interval morphisms $f_1, f_2$ defined in \autoref{eg:singint_map} we can construct the following regular \rsdual{}s{}
\begin{restoretext}
\begingroup\sbox0{\includegraphics{test/page1.png}}\includegraphics[clip,trim=0 {.0\ht0} 0 {.0\ht0} ,width=\textwidth]{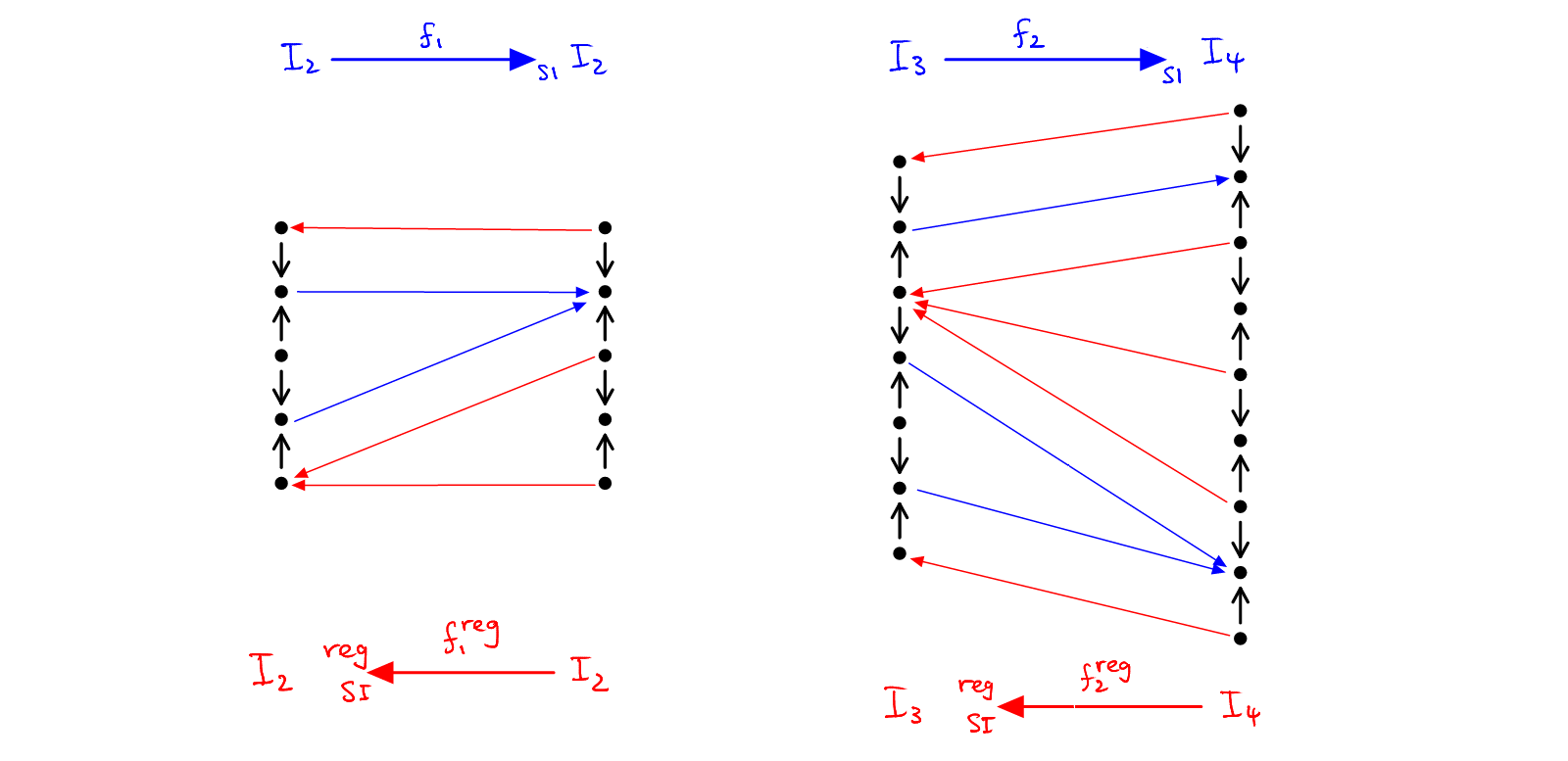}
\endgroup\end{restoretext}
Conversely, using $g_1, g_2$ as defined in \autoref{eg:singint_map} we construct the following singular \rsdual{}s{}
 \begin{restoretext}
\begingroup\sbox0{\includegraphics{test/page1.png}}\includegraphics[clip,trim=0 {.0\ht0} 0 {.0\ht0} ,width=\textwidth]{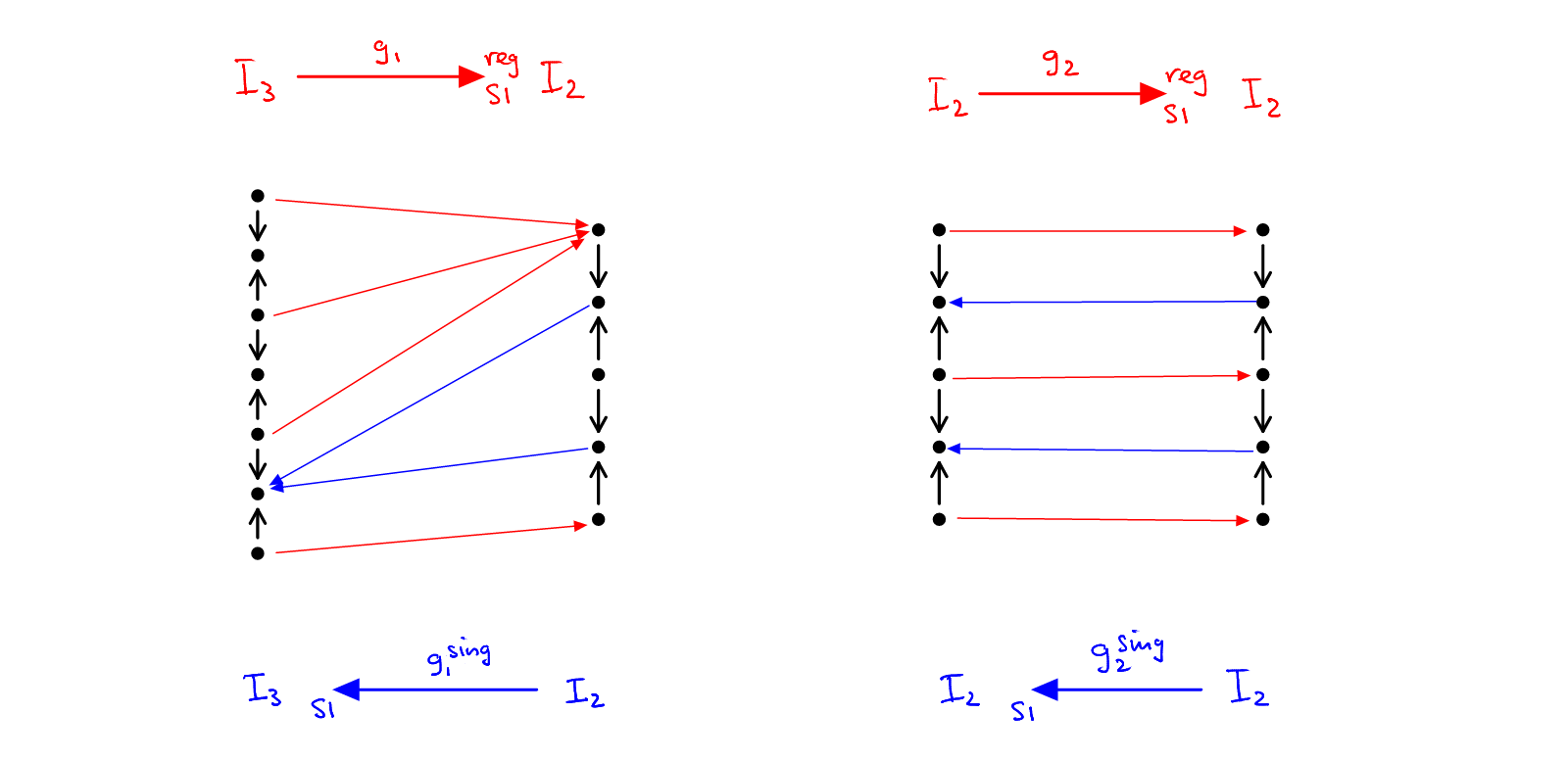}
\endgroup\end{restoretext}
\end{eg}

\begin{claim}[Well-definedness of \rsdual{}s{}] \label{rmk:regop_singop_welldef} The given definition \eqref{eq:fop_from_f} for $f\regop$ (resp. \eqref{eq:f_from_fop} for $g\singop$) yields a monotone function of regular segments (resp. singular heights). 
\proof  The proof is \stfwd{}. We argue in detail for $f \regop$ being a monotone function (the case of $g\singop$ follows similarly). Note that
\begin{equation} \label{eq:regular_covering}
[-1,2\iH_J + 1] = \left(\bigsqcup_{b \in \regcont(I)}  \hyperref[glsentry-openint]{]}\wwidehat f(b-1),\wwidehat f(b+1)\hyperref[glsentry-openint]{[}\right) ~ \sqcup ~ \im(\wwidehat f)
\end{equation}
where $\bigsqcup$ and $\sqcup$ denote disjoint unions of sets. This follows from $\wwidehat f(-1) = -1$,  $\wwidehat f(2\iH_I + 1) = 2\iH_J + 1$ and monotonicity of $\wwidehat{f}$ (since $f$ is monotone) which implies the claimed disjointness of the unions on the right.

We deduce from \eqref{eq:regular_covering} that $f\regop$ is well-defined as a function. To prove monotonicity of $f\regop$ we let $a < a'$, $a, a' \in \regcont(J)$. By \eqref{eq:regular_covering} we have $a \in ]f(b-1),f(b+1)[$, $a' \in ]f(b'-1),f(b'+1)[$ for unique $b, b' \in \regcont(J)$, and thus $f\regop(a) = b$, $f\regop(a') = b'$ by definition. Assume by contradiction $b' < b$. Then $f(b'+1) \leq f(b-1)$ by monotonicity, contradicting $a < a'$ and the choice of $b, b'$. \qed
\end{claim}

\begin{claim}[Ambidexterity of \rsdual{}s{}] \label{claim:involutivity_of_op} Using notation from \autoref{defn:regop_singop} the two pairs of maps
\begin{equation}
(s,r) \in \Set{(f, f\regop), (g\singop, g)}
\end{equation}
satisfy that for all $a \in \regcont(J)$ and $c \in \extsing (I)$
\begin{align} \label{eq:ambidexterity_condition_extended}
 \big(\wwidehat s(c) <  a &\iff c <  r(a)\big) \\
\text{~and~} ~\big( a <  \wwidehat s(c) &\iff r(a) <  c \big)
\end{align}
As a consequence we obtain $(f\regop)\singop = f$ and $(g \singop)\regop = g$.

\begin{rmk} \label{rmk:ambidext_cond}
\eqref{eq:ambidexterity_condition_extended} is called the \textit{ambidexterity condition} for the pair $(s,r)$.
\end{rmk}

\proof The proof is \stfwd{}. We argue for the first pair $(s,r) = (f, f\regop)$, the second argument for $(s,r) = (g \singop, g)$ is similar. For the first statement of the condition, assume $\wwidehat f (c) < a$. By \eqref{eq:regular_covering} there is $b$ such that $\wwidehat f (b - 1) < a < \wwidehat f(b + 1)$. Since $\wwidehat f$ is monotone we must have $c \leq b - 1$ and thus $c < b = f\regop(a)$. Conversely assume $c < f\regop(a)$. Again using \eqref{eq:regular_covering} we find $b$ such that $\wwidehat f (b - 1) < a < \wwidehat f (b + 1)$. By \eqref{eq:fop_from_f} we have $c < b = f\regop(a)$ and deduce $c \leq \wwidehat f (b - 1) < a$ as required. Finally, to show that $(f\regop)\singop = f$ and $(g \singop)\regop = g$ it suffices to note that pairs $(s,r)$ of a singular-height and a regular-segment morphisms which satisfy the ambidexterity condition \eqref{eq:ambidexterity_condition_extended} mutually determine each other uniquely. \qed
\end{claim}

\begin{lem}[Taking \rsdual{}s{} induces an equivalence] \label{lem:singular_reg_is_singular_op} There is an isomorphism of categories $\SI\op$ and $\SI\regop$ given by the assignment $I \mapsto I$ on objects, and $(f : I \to_{\SI} J) \mapsto (f\regop : J \to\regop_{\SI} I)$ on morphisms.

\proof It remains to prove functoriality of the assignment; that is, we need to prove that $\id\regop_I = \id_I$ and $(f_2 f_1) \regop = f\regop_1 f\regop_2$ for $f_1 : I_1 \to I_2$, $f_2 : I_2 \to I_3$. For the first statement note that $\wwidehat{\id_I} = \id_{\extsing(I)}$ by \autoref{notn:singular_morphism_boundary_cases} and thus \eqref{eq:fop_from_f} becomes
\begin{align}
\id_I\regop (a) = b \quad &\iff \quad  \wwidehat{\id_I}(b - 1) < a < \wwidehat{\id_I}(b+ 1) \\
\quad &\iff \quad b - 1 < a < b+ 1
\end{align}
implying $\id_I\regop (b) = b$. For the second statement, note that $\wwidehat{f_2 f_1} = \wwidehat{f_2} \wwidehat{f_1}$ by \autoref{notn:singular_morphism_boundary_cases} and thus, using \eqref{eq:fop_from_f} we find
\begin{align}
(f_2 f_1) \regop(a) = b &\iff \wwidehat{f_2 f_1}(b-1) < a < \wwidehat{f_2 f_1}(b + 1) \\
&\iff \wwidehat{f_2} \wwidehat{f_1}(b-1) < a < \wwidehat{f_2} \wwidehat{f_1} (b + 1)\\
&\iff  \wwidehat{f_1}(b-1) < f_2\regop (a) < \wwidehat{f_1} (b + 1) \\
&\iff  f_1\regop f_2\regop (a) = b 
\end{align}
where for the third implication we used the ambidexterity condition \eqref{eq:ambidexterity_condition_extended}. \qed
\end{lem}

\subsection{The profunctorial realisation of a map of singular intervals}

We now define an important functor on the category $\SI$, which produces a ``triangulation" of space based on a map of singular intervals. Technically, it is the unique ``minimal completion" of the graphs of a singular-height morphism and its regular \rsdual{} to a profunctorial relation.

\begin{defn}[Profunctorial realisation] \label{defn:order_realisation} Given $f : I \to_{\SI} J$ we define a relation $\SiR(f) : I \xslashedrightarrow{} J$, called the \textit{profunctorial realisation} of $f$, by a case distinction in the first argument as follows
\begin{equation} \tag{\textsc{Sing1}} \label{eq:defn_order_realisation_1}
\SiR(f)(a \in \singcont(I),b \in J)  \quad \iff \quad f(a) = b
\end{equation}
and
\begin{equation} \tag{\textsc{Sing2}} \label{eq:defn_order_realisation_2}
\SiR(f)(a \in \regcont(I), b \in J)\quad \iff \quad \wwidehat f(a-1) \leq b \leq \wwidehat f(a+1)
\end{equation}
Relations of the form $\SiR(f) : I \xslashedrightarrow{} J$ are called \textit{\SI-relations}, and $f$ is called the \textit{underlying \SI-morphism} of $\SiR(f)$.
\end{defn}

\begin{eg}[Profunctorial realisation] \label{eg:prel_real} Using $f_1$ and $f_2$ as defined in \autoref{eg:singint_map} we can construct $\SiR(f_1)$ and $\SiR(f_2)$ to be the following relations
\begin{restoretext}
\begingroup\sbox0{\includegraphics{test/page1.png}}\includegraphics[clip,trim=0 {.05\ht0} 0 {.0\ht0} ,width=\textwidth]{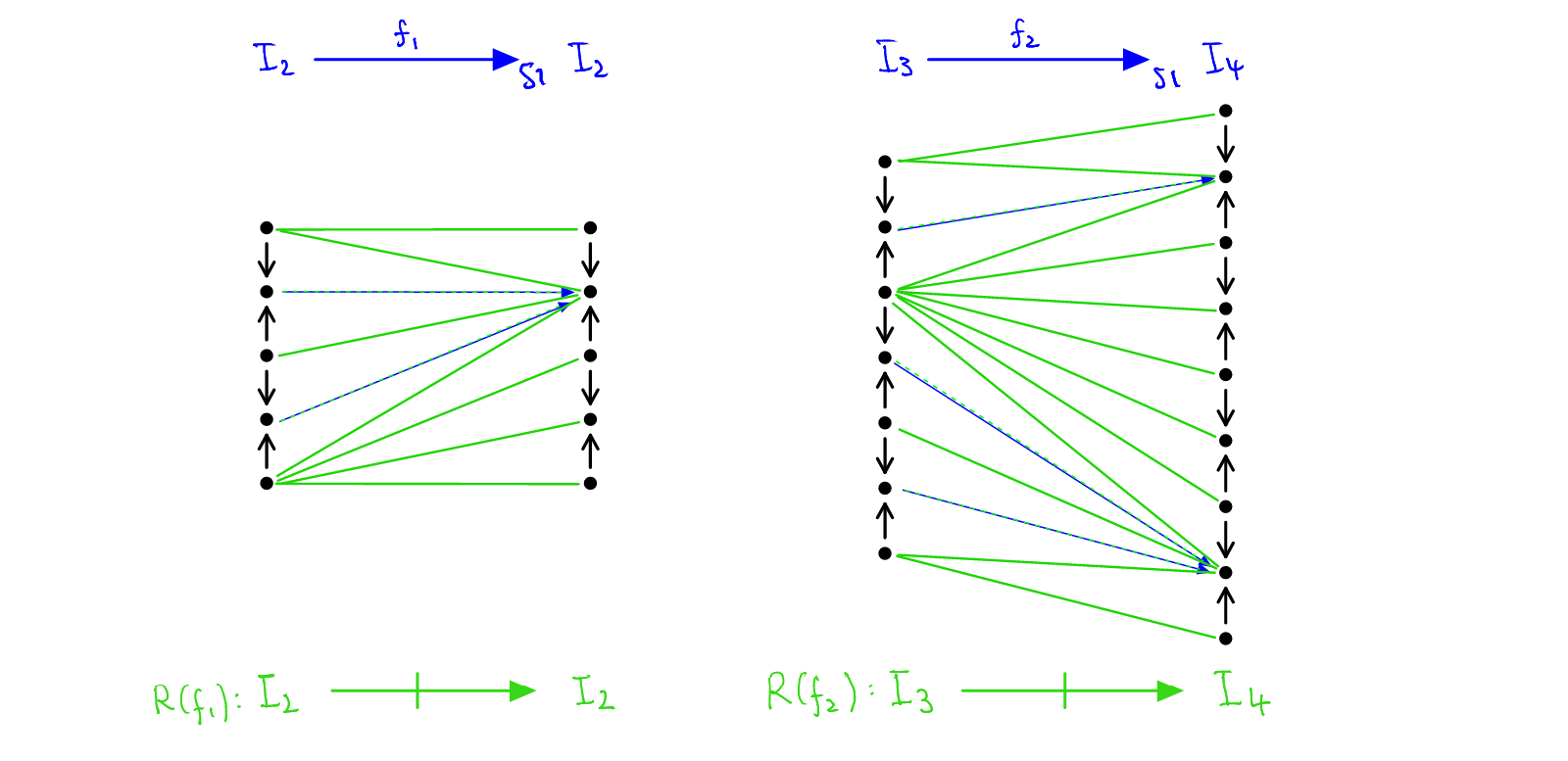}
\endgroup\end{restoretext}
Here, elements of the relation $\SiR(f_i)$ are indicated as union of the \cblue{} edges and \cgreen{} arrows.
\end{eg}

The preceding definition of the profunctorial realisation uses a \textit{singular}-height morphism for the defining inequalities \eqref{eq:defn_order_realisation_1} and \eqref{eq:defn_order_realisation_2}. The next claim shows that their is an equivalent characterisation using the \textit{regular} \rsdual{}.

\begin{claim}[\rsDual{} definition of profunctorial realisation] \label{claim:order_realisation_symmetric_def} Given $g : J \to\regop_{\SI} I$, set $f := g\singop$ and define a relation $\SiR(g) : I \xslashedrightarrow{} J$ by setting
\begin{equation} \tag{\textsc{Reg1}} \label{eq:defn_order_realisation_3}
\SiR(g)(a \in I,b \in \regcont(J))  \quad \iff \quad g(b) = a
\end{equation}
and
\begin{equation} \tag{\textsc{Reg2}} \label{eq:defn_order_realisation_4}
\SiR(g)(a \in I, b \in \singcont(J))\quad \iff \quad g(b-1) \leq a \leq g(b+1)
\end{equation}
Then we claim $\SiR(f) = \SiR(g)$ where $\SiR(f)$ was defined in \autoref{defn:order_realisation}.
\proof The proof is \stfwd{}. We argue by case distinction:
\begin{enumerate}
\item $a \in \singcont(I), b \in \singcont(J)$. Then $\SiR(f)(a,b)$ if and only if $f(a)=b$ and $\SiR(g)(a,b)$ if and only if $g(b-1) < a < g(b+1)$. By \eqref{eq:f_from_fop} both are equivalent.
\item $a \in \singcont(I), b \in \regcont(J)$. Then $\SiR(f)(a,b)$ if and only if $f(a)=b$ and $\SiR(g)(a,b)$ if and only if $g(b) = a$. Both are impossible by \autoref{defn:singular_intervals_morphism}.
\item $a \in \regcont(J), b \in \regcont(J)$. Then $\SiR(f)(a,b)$ if and only if $\wwidehat f (a - 1) < b < \wwidehat f (a + 1)$ and $\SiR(g)(a,b)$ if and only if $g(b) = a$. By \eqref{eq:fop_from_f} both are equivalent.
\item $a \in \regcont(J), b \in \singcont(J)$. Then $\SiR(f)(a,b)$ if and only if $\wwidehat f(a-1) \leq b \leq \wwidehat f(a+1)$ and $\SiR(g)(a,b)$ if and only if $g(b-1) \leq a \leq g(b+1)$. Equivalence follows since firstly $\wwidehat f(a-1) \leq b \iff a \leq g(b+1)$ and secondly $b \leq \wwidehat f(a+1) \iff g(b-1) \leq a$. We argue for the first statement (the second follows similarly). Assume $\wwidehat f(a-1) \leq b < b + 1$. Then ambidexterity  \eqref{eq:ambidexterity_condition_extended} yields $a -1 < g(b+1)$ and thus $a \leq g(b+1)$. Conversely, assume $a - 1 < a \leq g(b+1)$. Then ambidexterity \eqref{eq:ambidexterity_condition_extended} yields $\wwidehat{f} (a - 1) < b + 1$ and thus $\wwidehat f (a-1) \leq b$ as required. \qed
\end{enumerate}
\end{claim}

\begin{eg}[\rsDual{} definition of profunctorial realisation] Using $f_1$ and $f_2$, two examples of the construction in the claim's statement are the following
\begin{restoretext}
\begingroup\sbox0{\includegraphics{test/page1.png}}\includegraphics[clip,trim=0 {.0\ht0} 0 {.0\ht0} ,width=\textwidth]{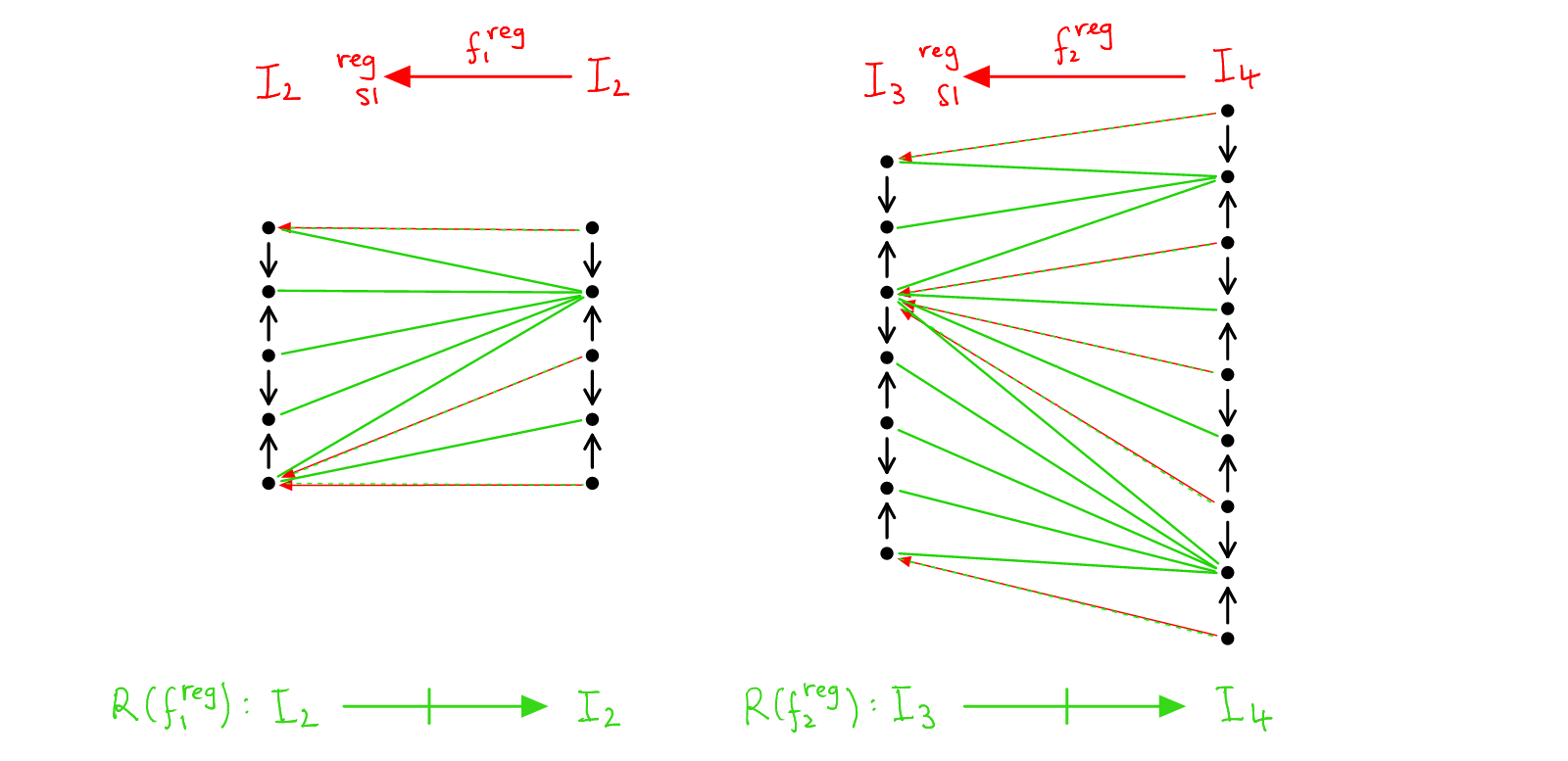}
\endgroup\end{restoretext}
Comparing to \autoref{eg:prel_real}, we find $\SiR(f_1) = \SiR(f\regop_1)$ and $\SiR(f_2) = \SiR(f\regop_2)$ as claimed.
\end{eg}

As suggested by its name, an \SI-relation $\SiR(f) : I \xslashedrightarrow{} J$ with underlying \SI-morphism $f : I \to_{\SI} J$ is not only a relation but a profunctorial relation (cf. \autoref{defn:prerequisites}). Thus $\SiR(f)$ is a morphism in $\PRel$. In fact, the assignment $f \mapsto \SiR(f)$ is functorial. To motivate both statements we first give an example.

\begin{eg} Using \autoref{notn:depicting_prel}, the following depicts a composition of \SI-relations (whose underlying \SI-morphisms are indicated by \cblue{} edges from left to right)
\begin{restoretext}
\begingroup\sbox0{\includegraphics{test/page1.png}}\includegraphics[clip,trim=0 {.29\ht0} 0 {.15\ht0} ,width=\textwidth]{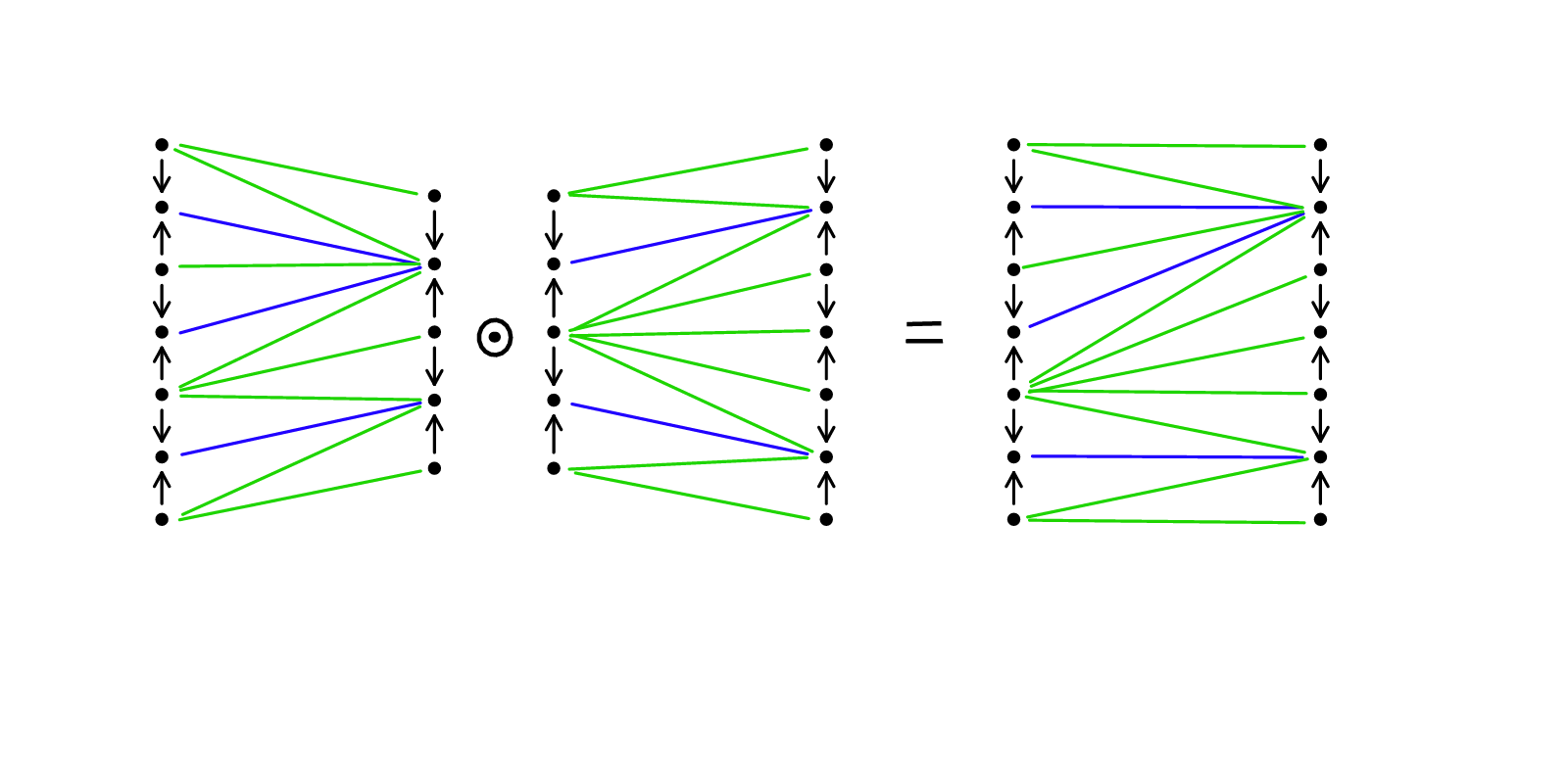}
\endgroup\end{restoretext}
We observe that in this case, profunctorial realisation of the composition of the two singular interval morphisms is the composition (as relations) of the profunctorial realisation of each individual singular interval morphism.

Similarly, we find that profunctorial realisations of the identity yield the identity profunctorial realisation, which is illustrated by the following example
\begin{restoretext}
\begingroup\sbox0{\includegraphics{test/page1.png}}\includegraphics[clip,trim=0 {.1\ht0} 0 {.2\ht0} ,width=\textwidth]{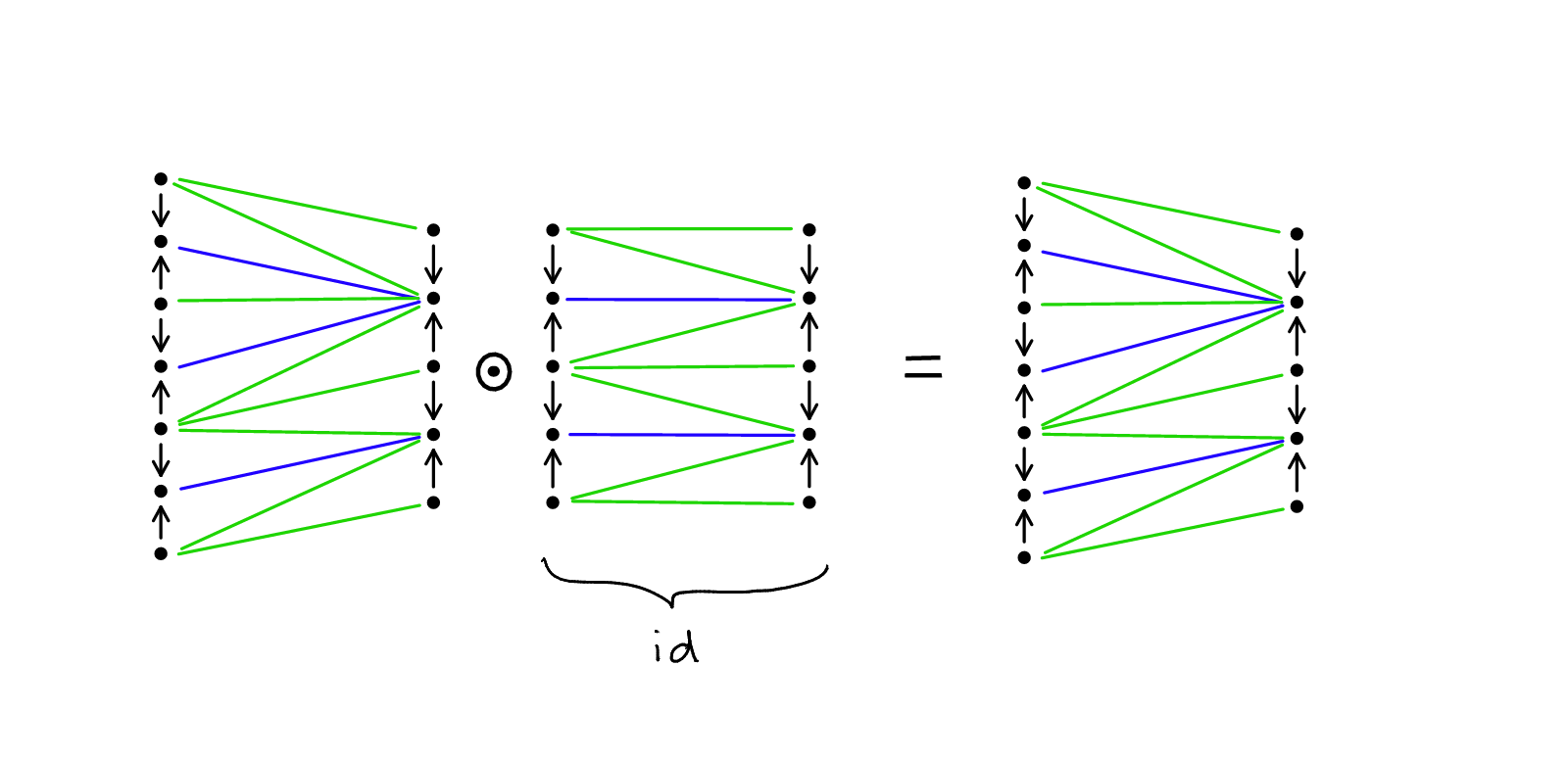}
\endgroup\end{restoretext}
\end{eg}

\begin{lem}[Profunctorial realisation functor] \label{lem:order_realisation_functorial} \hfill
\begin{enumerate}
\item The relation $\SiR(f) : I \xslashedrightarrow{} J$ as defined in \autoref{defn:order_realisation} is a profunctorial relation between the posets $I$ and $J$ (in singular order)
\item The assignment $f \mapsto \SiR(f)$ is functorial, giving rise to a functor $\SiR : \SI \to \PRel$ called the \emph{profunctorial realisation}.
\end{enumerate}

\proof 
\begin{enumerate}
\item The proof is \stfwd{}. For profunctoriality of $\SiR(f)$ we need to verify \eqref{eq:profunctorial_relation}. Let $a' \to a$ in $I$ and $\SiR(f)(a,b)$. We need to show $\SiR(f)(a',b)$. Excluding the case where $a' = a$ (in which the statement holds trivially), by \autoref{defn:singular_intervals} we either have $a' = a - 1$ or $a' = a + 1$, and $a' \in \regcont(I), a \in \singcont(I)$ in both cases. We assume $a' = a - 1 \to a$ (and argue similarly in the second case). Since $a \in \singcont(I)$ we have $b = f(a) = f(a' + 1)$ by \eqref{eq:defn_order_realisation_1}. Thus $ f(a' - 1) \leq  b \leq f(a' + 1)$ by monotonicity of $f$ and we deduce $\SiR(f)(a',b)$ by \eqref{eq:defn_order_realisation_2} as required. Next assume $b \to b'$ in $J$ and $\SiR(f)(a,b)$. We need to show $\SiR(f)(a,b')$. Excluding again $b = b'$ we either have $b' = b + 1$ or $b' = b - 1$ and argue in the former case. Since $b \in \regcont(J)$ we have $a = f \regop (b)  = f \regop (b' + 1)$ by \eqref{eq:defn_order_realisation_3}. Thus $f \regop(b' - 1) \leq a \leq f \regop(b' + 1)$ by monotonicity of $f\regop$ and we deduce $\SiR(f)(a,b')$ by \eqref{eq:defn_order_realisation_4} as required.
\item The proof is \stfwd{}. For functoriality we need to check firstly, that 
\begin{equation}
\SiR(\id : I \to_{\SI} I) = \Hom_I(-,-)
\end{equation}
and secondly, that
\begin{equation}
\SiR(f_2 f_1) = \SiR(f_1) \odot \SiR(f_2)
\end{equation}
for any $f_1 : I_1 \to_{\SI} I_2$, $f_2 : I_2 \to_{\SI} I_3$. 

For the first statement recall \autoref{constr:PRel} where we defined the Hom relation. The statement then follows by comparison of \eqref{eq:defn_order_realisation_1} and \eqref{eq:defn_order_realisation_2} (in the case that $f = \id$ and thus $\wwidehat f = \id_{\extsing(I)}$) with \autoref{defn:singular_intervals} defining the singular order on $I$. 

For the second statement we show for $a \in I_1$, $c \in I_3$ that
\begin{equation}
\SiR(f_2 f_1)(a,c) \quad \iff \quad \big(\SiR(f_1) \odot \SiR(f_2)\big)(a,c)
\end{equation}
To prove ``$\Leftarrow$", assume $\SiR(f_1)(a,b)$ and $\SiR(f_2)(b,c)$ for some $b \in I_2$. We want to show $\SiR(f_2f_1)(a,c)$. There are the following cases
\begin{enumerate}
\item $a \in \singcont(I_1)$. Then $b = f_1(a) \in \singcont(I_2)$ and $c = f_2(b)$ by \eqref{eq:defn_order_realisation_1} and thus $\SiR(f_2f_1)(a,c)$.
\item $a \in \regcont(I_1)$ and $b \in \singcont(I_2)$. Then $c = f_2(b)$ by \eqref{eq:defn_order_realisation_1} and $\wwidehat {f_1} (a -1) \leq b \leq \wwidehat {f_1} (a + 1)$ by \eqref{eq:defn_order_realisation_2}. Monotonicity of $\wwidehat {f_2}$ as well as $\wwidehat{f_2}\wwidehat {f_1} =\wwidehat {f_2 f_1}$ implies that $\wwidehat {f_2 f_1}(a-1) \leq c \leq \wwidehat {f_2 f_1}(a+1)$ and thus $\SiR(f_2f_1)(a,c)$ as required.
\item $a \in \regcont(I_1)$ and $b \in \regcont(I_2)$. Then $\wwidehat {f_1} (a -1) \leq b-1 < b < b+1\leq \wwidehat f (a + 1)$ and $\wwidehat {f_2} (b -1) \leq c \leq \wwidehat {f_2} (b + 1)$ by \eqref{eq:defn_order_realisation_2}. Monotonicity of $\wwidehat {f_2}$ as well as $\wwidehat{f_2}\wwidehat {f_1} =\wwidehat {f_2 f_1}$ implies that $\wwidehat {f_2 f_1}(a-1) \leq c \leq \wwidehat {f_2 f_1}(a+1)$ and thus $\SiR(f_2f_1)(a,c)$ as required.
\end{enumerate}
To prove ``$\imp$" assume $\SiR(f_2f_1)(a,c)$. We want to show that there is $b \in I_2$ such that $\SiR(f_1)(a,b)$ and $\SiR(f_2)(b,c)$. There are the following cases.
\begin{enumerate}
\item $a \in \singcont(I_1)$. Then $f_2f_1(a) = c$ and we can take $b = f_1(a)$  and \eqref{eq:defn_order_realisation_1} yields $\SiR(f_1)(a,b)$ and $\SiR(f_2)(b,c)$.
\item $a \in \regcont(I_1)$, $c \in \regcont(I_3)$. Then $a = (f_2f_1)\regop(c)$ and we can take $b = f_2 \regop (c)$ and \eqref{eq:defn_order_realisation_3} yields $\SiR(f_1)(a,b)$ and $\SiR(f_2)(b,c)$.
\item $a \in \regcont(I_1)$. $c \in \singcont(I_3)$. Then $\wwidehat {f_2 f_1}(a-1) \leq c \leq \wwidehat {f_2 f_1}(a+1)$. Since $\wwidehat {f_2 f_1} = \wwidehat{f_2}\wwidehat {f_1}$ this implies $\wwidehat{f_2}\wwidehat {f_1} (a-1) \leq c \leq \wwidehat{f_2}\wwidehat {f_1} (a+1)$. Using the proof of \autoref{claim:order_realisation_symmetric_def} case (iv) this implies $\wwidehat {f_1} (a-1) + 1 \leq (f_2)\regop (c+1)$ and $(f_2)\regop (c-1) \leq \wwidehat {f_1} (a+1) - 1$. Thus, neither $\wwidehat{f_2}(a + 1) < f_2\regop(c-1)$ nor $f_2\regop(c+1) < \wwidehat{f_1}(a - 1)$ is possible. In other words, the intersection of the intervals $[\wwidehat{f_1}(a - 1), \wwidehat{f_2}(a + 1)]$ and $[f_2\regop(c-1), f_2\regop(c+1)]$ is non-trivial, and there is a $b\in I_2$ which simultaneously satisfies $\wwidehat{f_1}(a - 1) \leq b \leq \wwidehat{f_2}(a + 1)$ and $f_2\regop(c-1) \leq b \leq f_2\regop(c+1)$. Then \eqref{eq:defn_order_realisation_2} and  \eqref{eq:defn_order_realisation_4} yield $\SiR(f_1)(a,b)$ and $\SiR(f_2)(b,c)$ respectively.
\end{enumerate}
\end{enumerate}
\qed
\end{lem}

\begin{rmk}[$\SiR$ is an embedding of categories] \label{rmk:profunctorial_real_inj} Note that the profunctorial realisation $R : \SI \to \PRel$ is an embedding of $\SI$ into $\PRel$: It is injective on objects and faithful on morphism since for $f_1, f_2 : I \to_{\SI} J$ we have
\begin{align}
\SiR(f_1) = \SiR(f_2) ~ &\imp ~ \big(\SiR(f_1)(a \in \sing (I), b \in J) \iff \SiR(f_2)(a \in \singcont(I), b \in J) \big) \\
&\imp ~ f_1 = f_2
\end{align}
\end{rmk}

\subsection{Total order on the edges of profunctorial realisations}

In this section we develop an important tool for later on proofs, called \textit{edge induction}. If one visualises a profunctorial realisations of singular height morphisms one can observe that edges can always be assumed to have no intersections. This allows to put a total order on edges by traversing them one by one. We will establish this result by first defining successors and predecessors, and then observing that these are mutually invertible constructions.

\begin{defn}[Edges of singular interval morphisms] \label{defn:edge_sets} Given a morphism $f : I \to_{\SI} J$ of singular intervals, define the \textit{set of edges} of $f$ to be the set
\begin{equation}
\edgeset(f) = \Set{ (a,a') ~|~ a \in I, a' \in J, \SiR(f)(a,a') }
\end{equation}
We also set $\edgeset(I) := \edgeset(\id_I)$ for a singular interval $I$. Given an edge $\eda = (a,a') \in \edgeset(f)$, we denote the \textit{source} $a \in I$ of $\eda$ by a subscript $\soe$ as follows
\begin{equation}
\eda\ssoe  := a
\end{equation}
and the \textit{target} $a' \in J$ of $\eda$ by a subscript $\tae$ as follows
\begin{equation}
\eda\ttae  := a'
\end{equation}
Further define the \textit{norm} of an edge $\eda \in \edgeset(f)$ by setting
\begin{equation}
\avg{\eda} = \eda\ssoe  + \eda\ttae 
\end{equation} 

\end{defn}

\begin{eg}[Set of edges] Consider the following \SI-relation (its underlying \SI-morphism is indicated by \cblue{} edges)
\begin{restoretext}
\begingroup\sbox0{\includegraphics{test/page1.png}}\includegraphics[clip,trim=0 {.05\ht0} 0 {.05\ht0} ,width=\textwidth]{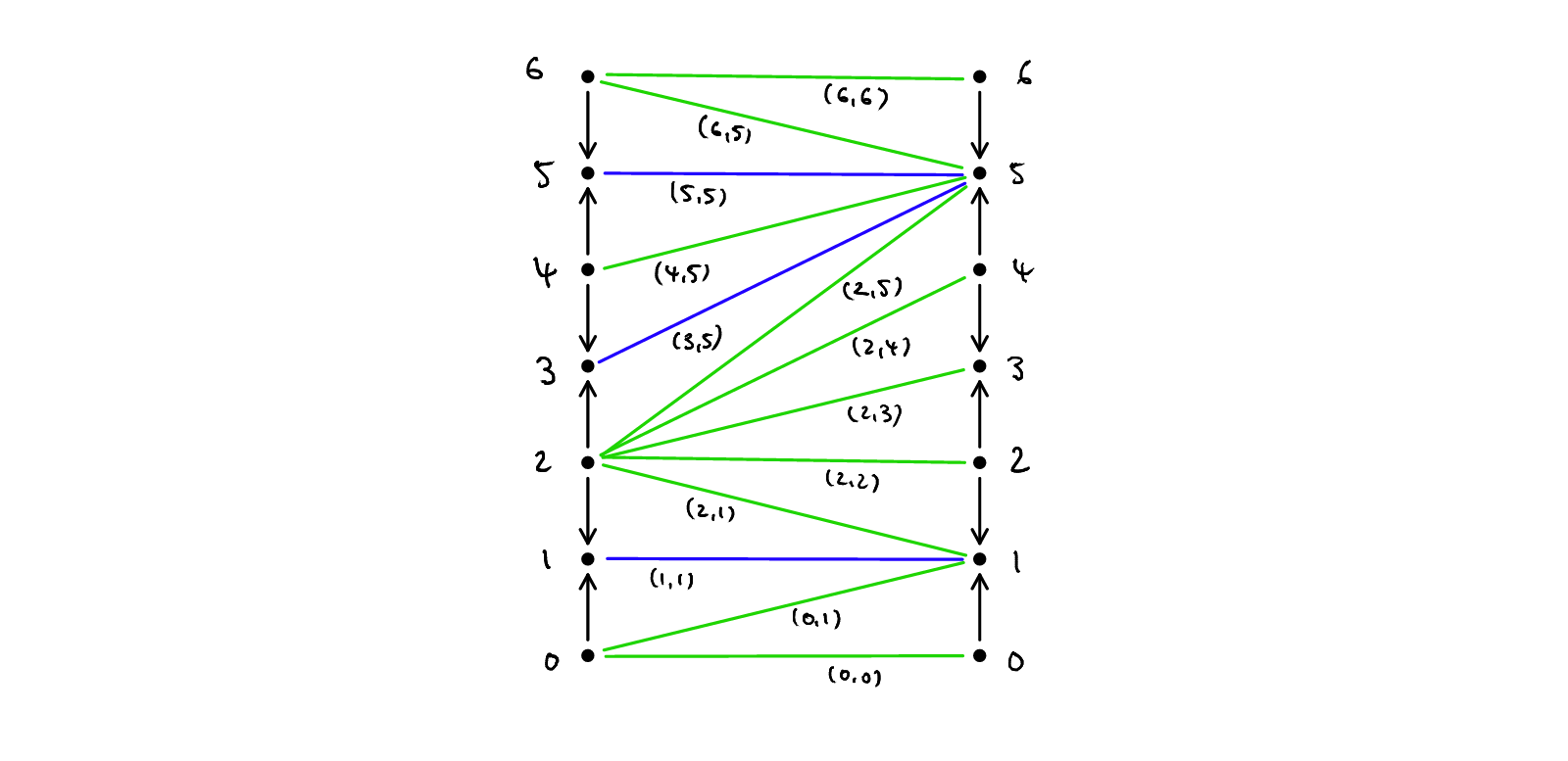}
\endgroup\end{restoretext}
The set of edges of this \SI-relation is the set of tuples of number written next to each edge.
\end{eg}

\begin{constr}[Successor edge] \label{defn:successor}  Let $f : I_\sop \to_{\SI} {I_\tap}$ and $\eda \in \edgeset(f)$ such that $\avg{\eda} < 2(\iH_{I_\sop} + \iH_{I_\tap})$. We define the \textit{successor edge} $\succ\eda = (\succ\eda\ssoe , \succ\eda\ttae ) \in \edgeset(f)$ of $x$, and its \textit{filler edge} $\succfill\eda \in \edgeset(I_{\succindex\eda})$, where $\succindex\eda \in \Set{ \sop, \tap }$ is the \textit{filler's index}. The definition is by a case distinction of four cases that can be visualised as follows
\begin{restoretext}
\begingroup\sbox0{\includegraphics{test/page1.png}}\includegraphics[clip,trim=0 {.0\ht0} 0 {.0\ht0} ,width=\textwidth]{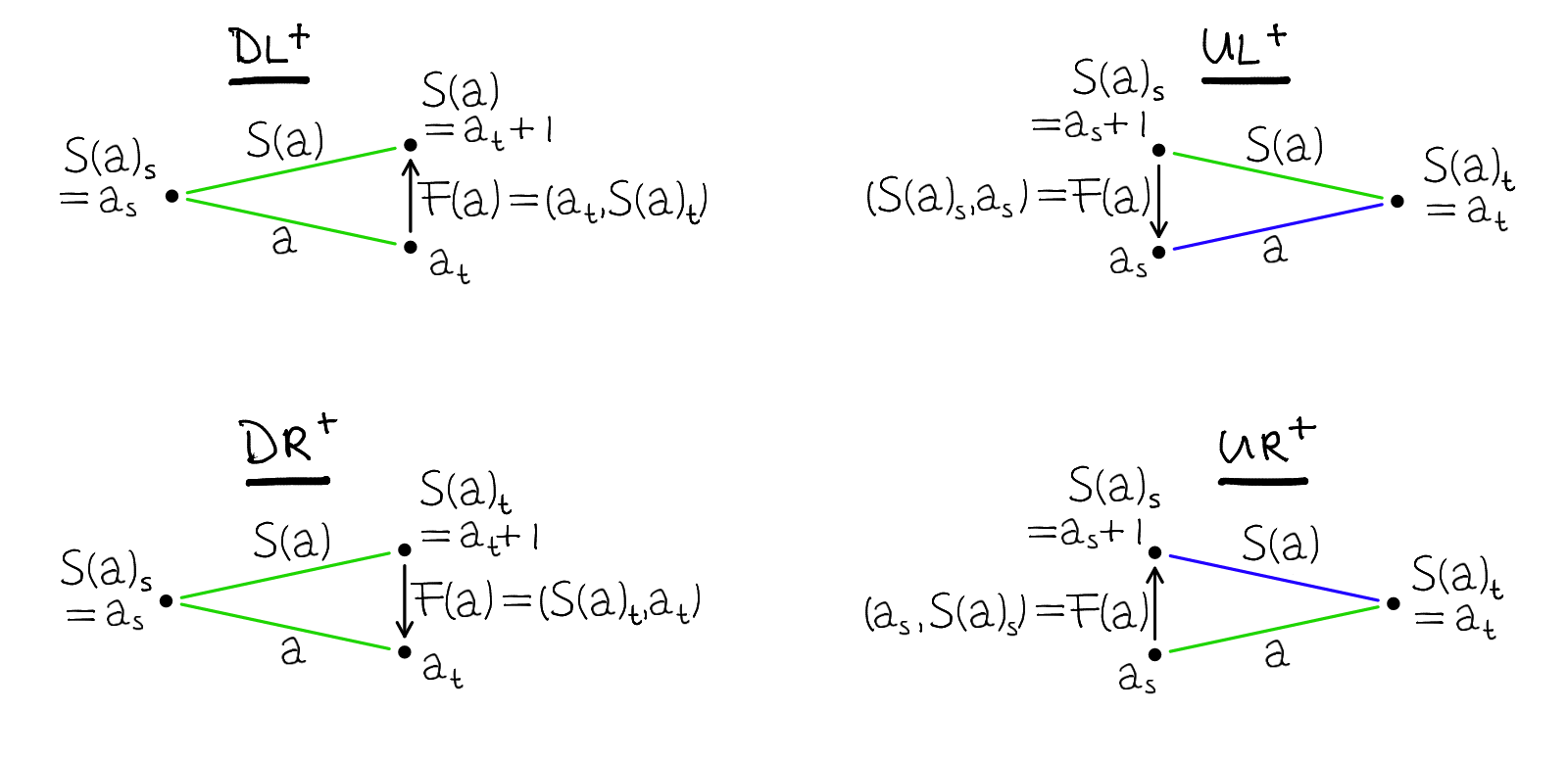}
\endgroup\end{restoretext}
Here, \cblue{} edges are edges between singular heights. Comparing to the previous example, the reader should observe that all triangles appearing therein are covered by the above four cases. Below we will see that this holds in general.

Recall that for $(\succ\eda\ssoe , \succ\eda\ttae )$ to be indeed an element of $\edgeset(f)$ we need $\SiR(f)(\succ\eda\ssoe , \succ\eda\ttae )$, which we verify in each case.

\begin{enumerate}
\item[(DL$^+$)] $\eda\ssoe  \in \regcont(I_\sop)$, $\eda\ttae  \in \regcont({I_\tap})$. Then we can set $\succ\eda\ssoe  = \eda\ssoe $ and $\succ\eda\ttae  = \eda\ttae  + 1$ ($\succ\eda\ttae  \in {I_\tap}$ by \autoref{rmk:successor_well_def}). Since $\eda\ttae  \to \succ\eda\ttae $ in ${I_\tap}$, profunctoriality of $\SiR(f)$ implies $\SiR(f)(\succ\eda\ssoe , \succ\eda\ttae )$ as required for $(\succ\eda\ssoe , \succ\eda\ttae ) \in \edgeset(f)$. We define the filler edge as $\succfill\eda := (\eda\ttae ,\succ\eda\ttae ) \in \edgeset({I_{\succindex\eda}})$ where $\succindex\eda = \tap$. Note that the pair $(\edb\ssoe ,\edb\ttae ) := (\succ\eda\ssoe ,\succ\eda\ttae )$ satisfies case (DL$^-$) below.

\item[(UL$^+$)] $\eda\ssoe  \in \singcont({I_\tap})$, $\eda\ttae \in \singcont(I_\sop)$. Then we can set $\succ\eda\ssoe  = \eda\ssoe  + 1$ and $\succ\eda\ttae  = \eda\ttae $. Since $\succ\eda\ssoe  \to \eda\ssoe $ in $I_\sop$, profunctoriality of $\SiR(f)$ implies $\SiR(f)(\succ\eda\ssoe , \succ\eda\ttae )$ as required for $(\succ\eda\ssoe , \succ\eda\ttae ) \in \edgeset(f)$. We define the filler edge as $\succfill\eda := (\succ\eda\ssoe ,\eda\ssoe ) \in \edgeset(I_{\succindex\eda})$ where $\succindex\eda = \sop$. Note that the pair $(\edb\ssoe ,\edb\ttae ) := (\succ\eda\ssoe ,\succ\eda\ttae )$ satisfies case (UL$^-$) below.

\item[(DR$^+$)] $\eda\ssoe  \in \regcont(I_\sop)$, $\eda\ttae  \in \singcont({I_\tap})$ and $\wwidehat f(\eda\ssoe +1) > \eda\ttae $. Then we can set $\succ\eda\ssoe  = \eda\ssoe $ and $\succ\eda\ttae  = \eda\ttae  + 1$. Using \eqref{eq:defn_order_realisation_2} we obtain $\SiR(f)(\succ\eda\ssoe ,\succ\eda\ttae )$ as required for $(\succ\eda\ssoe , \succ\eda\ttae ) \in \edgeset(f)$. We  define the filler edge as $\succfill\eda := (\succ\eda\ttae ,\eda\ttae ) \in \edgeset({I_{\succindex\eda}})$ where $\succindex\eda = \tap$. Note that the pair $(\edb\ssoe ,\edb\ttae ) := (\succ\eda\ssoe ,\succ\eda\ttae )$ satisfies case (DR$^-$) below.

\item[(UR$^+$)] $\eda\ssoe  \in \regcont(I_\sop)$, $\eda\ttae  \in \singcont({I_\tap})$ and $\wwidehat f(\eda\ssoe +1) = \eda\ttae $. Then we can set $\succ\eda\ssoe  = \eda\ssoe  + 1$ and $\succ\eda\ttae  = \eda\ssoe $ ($\succ\eda\ssoe  \in I_\sop$ by \autoref{rmk:successor_well_def_2}). Using \eqref{eq:defn_order_realisation_1} we obtain $\SiR(f)(\succ\eda\ssoe ,\succ\eda\ttae )$ as required for $(\succ\eda\ssoe , \succ\eda\ttae ) \in \edgeset(f)$. We define the filler edge as $\succfill\eda := (\eda\ssoe ,\succ\eda\ssoe ) \in \edgeset(I_{\succindex\eda})$ where $\succindex\eda = \sop$. Note that the pair $(\edb\ssoe ,\edb\ttae ) := (\succ\eda\ssoe ,\succ\eda\ttae )$ satisfies case (UR$^-$) below.
\end{enumerate}
\end{constr}

\begin{rmk} \label{rmk:successor_well_def} In case (DL$^+$) our choice $\succ\eda\ttae  = \eda\ttae  + 1$ is valid (that is, $\succ\eda\ttae  \in I_\tap$) by the following argument: first note $\eda\ssoe  = f\regop(\eda\ttae )$ by \eqref{eq:defn_order_realisation_3}. Assume $\succ\eda\ttae  = \eda\ttae  + 1 \notin {I_\tap}$, i.e. $\eda\ttae  = 2\iH_{I_\tap}$. Then $\eda\ssoe  = 2\iH_{I_\sop}$ since $f\regop$ is end-point preserving (cf. \autoref{defn:singular_intervals_morphism}). This contradicts our assumption on $x$ that $\avg{x} < 2(\iH_{I_\sop} + \iH_{I_\tap})$, and thus we must have $\eda\ttae  < 2\iH_{I_\tap}$. 
\end{rmk}

\begin{rmk} \label{rmk:successor_well_def_2}
In case (UR$^+$) our choice $\succ\eda\ssoe  = \eda\ssoe  + 1$ is valid (that is, $\succ\eda\ssoe  \in I_\sop$) since by assumption in that case we have $\wwidehat f (\eda\ssoe  + 1) = \eda\ttae $. Then, using \autoref{defn:singular_intervals_morphism} together with $\eda\ttae  < 2\iH_{I_\tap}$ we must have $\eda\ssoe  < 2\iH_{I_\sop}$ (otherwise $\wwidehat f(\eda\ssoe +1) = \wwidehat f(2\iH_{I_\sop}+1) = 2\iH_{I_\tap} +1$ which cannot equal $\eda\ttae $ and thus contradicts the assumptions of that case). 
\end{rmk}

\begin{constr}[Predecessor edge] \label{defn:predecessor} Let $f : I_\sop \to_{\SI} {I_\tap}$ and $\edb \in \edgeset(f)$ such that $\avg{\edb} > 0$. We define the \textit{predecessor} $\pred\edb \equiv (\pred\edb\ssoe , \pred\edb\ttae ) \in \edgeset(f)$ of $y$ by the following case distinction. The cases can be visualised as follows
\begin{restoretext}
\begingroup\sbox0{\includegraphics{test/page1.png}}\includegraphics[clip,trim=0 {.0\ht0} 0 {.0\ht0} ,width=\textwidth]{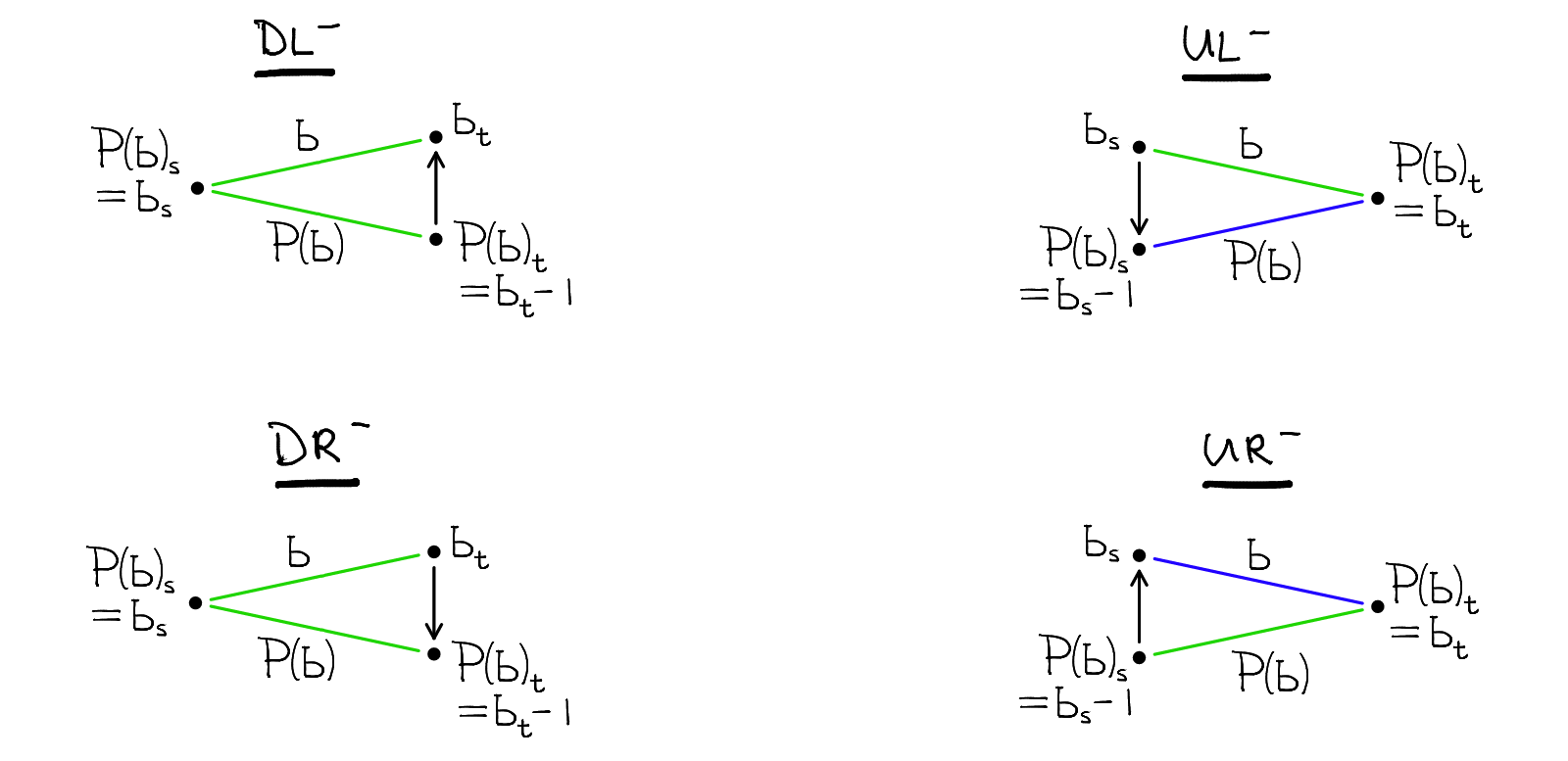}
\endgroup\end{restoretext}
Here, \cblue{} edges are edges between singular heights. Note the correspondence of these four cases to the cases in the successor construction.

Similar to the successor construction, in each of the below cases we need to show $\SiR(f)(\pred\eda\ssoe , \pred\eda\ttae )$.
\begin{enumerate}
\item[(DL$^-$)] $\edb\ssoe  \in \regcont(I_\sop)$, $\edb\ttae  \in \singcont({I_\tap})$ and $f\regop(\edb\ttae  - 1) = \edb\ssoe $. Then we can set $\pred\edb\ssoe  = \edb\ssoe $ and $\pred\edb\ttae  = \edb\ttae  - 1$. Using \eqref{eq:defn_order_realisation_3} we obtain $\SiR(f)(\pred\edb\ssoe ,\pred\edb\ttae )$ as required. Note that the pair $(\eda\ssoe ,\eda\ttae ) := (\pred\edb\ssoe ,\pred\edb\ttae )$ satisfies (DL$^+$) above.

\item[(UL$^-$)] $\edb\ssoe  \in \regcont(I_\sop)$, $\edb\ttae  \in \singcont({I_\tap})$ and $f\regop(\edb\ttae  - 1) < \edb\ssoe $. Then we can set $\pred\edb\ssoe  = \edb\ssoe  - 1$ and $\pred\edb\ttae  = \edb\ttae $ ($\pred\edb\ssoe  \in I_\sop$ by \autoref{rmk:predecessor_well_def_1}). Using \eqref{eq:defn_order_realisation_4} we obtain $\SiR(f)(\pred\edb\ssoe ,\pred\edb\ttae )$ as required. Note that the pair $(\eda\ssoe ,\eda\ttae ) := (\pred\edb\ssoe ,\pred\edb\ttae )$ satisfies case  (UL$^+$) above.

\item[(DR$^-$)] $\edb\ssoe  \in \regcont(I_\sop)$, $\edb\ttae  \in \regcont({I_\tap})$. Then we can set $\pred\edb\ssoe  = \edb\ssoe $ and $\pred\edb\ttae  = \edb\ttae  - 1$ ($\pred\edb\ttae  \in {I_\tap}$ by \autoref{rmk:predecessor_well_def_2}). Since $\edb\ttae  \to \pred\edb\ttae $ in ${I_\tap}$, profunctoriality of $\SiR(f)$ implies $\SiR(f)(\pred\edb\ssoe ,\pred\edb\ttae )$ as required. Note that the pair $(\eda\ssoe ,\eda\ttae ) := (\pred\edb\ssoe ,\pred\edb\ttae )$ satisfies case (DR$^+$) above.

\item[(UR$^-$)] $\edb\ssoe  \in \singcont(I_\sop)$, $\edb\ttae  \in \singcont(I_\sop)$. Then we can set $\pred\edb\ssoe  = \edb\ssoe  - 1$ and $\pred\edb\ttae  = \edb\ttae $. Since $\pred\edb\ssoe  \to \edb\ssoe $ in $I_\sop$, profunctoriality of $\SiR(f)$ implies $\SiR(f)(\pred\edb\ssoe , \pred\edb\ttae )$ as required. Note that the pair $(\eda\ssoe ,\eda\ttae ) := (\pred\edb\ssoe , \pred\edb\ttae )$ satisfies case (UR$^+$) above.
\end{enumerate}
\end{constr}

\begin{rmk} \label{rmk:predecessor_well_def_1} In case (UL$^-$) our choice $\pred\edb\ssoe  = \edb\ssoe  - 1$ is valid since by assumption in that case we have $f\regop(\edb\ttae  -1) < \edb\ssoe $. But $0\leq f\regop(\edb\ttae  -1)$ and thus $\edb\ssoe  > 0$ from which we deduce $\pred\edb\ssoe  \in I_\sop$.
\end{rmk}

\begin{rmk} \label{rmk:predecessor_well_def_2} In case (DR$^-$) our choice $\pred\edb\ttae  = \edb\ttae  - 1$ is valid by the following argument: first not that $\edb\ssoe  = f\regop (\edb\ttae )$. Now assume $\pred\edb\ttae  \notin {I_\tap}$ by contradiction, i.e. equivalently that $\edb\ttae  = 0$. Since $f \regop$ is endpoint-preserving we deduce $\edb\ssoe  = 0$. But this would mean $\avg{y} = 0$ contradicting our assumption on $y$.
\end{rmk}

\begin{notn}
Given $f : I_\sop \to_{\SI} {I_\tap}$ and $\eda,\edb \in \edgeset(f)$, if the respective successor and predecessor edges exist, we inductively set $\succn \eda 1 = \succ\eda$, $\predn \edb 1 = \pred \edb$ and $\succn \eda n = \succ{\succn \eda {n-1}}$, $\predn \edb n = \pred{\predn \edb {n-1}}$, for $n>1$.
\end{notn}

\begin{lem}[Edge induction] \label{claim:edge_set_properties} Let $f : I_\sop \to_{\SI} {I_\tap}$, and assume $\eda,\edb \in \edgeset(f)$. 
\begin{enumerate}
\item If $\edb = \succ\eda$ then $\avg{\edb} = \avg{\eda} + 1$, $\edb\ssoe  \geq \eda\ssoe $, $\edb\ttae  \geq \eda\ttae $ and $\pred\edb = \eda$. 

\item If $\eda = \pred\edb$ then $\avg{\eda} = \avg{\edb} - 1$, $\eda\ssoe  \leq \edb\ssoe $, $\eda\ttae  \leq \edb\ttae $ and $\succ\eda = \edb$

\item $\edgeset(f)$ has a unique $\eda$ with $\avg{\eda} = 0$ and this is the element with minimal norm. Similarly, $\edgeset(f)$ has unique $\edb$ with $\avg{\edb} = 2\iH_{I_\sop} + 2\iH_{I_\tap}$ and this is the element with maximal norm.

\item If $\avg{\eda} = \avg{\edb}$, then $\eda = \edb$
\end{enumerate}
This establishes that $\avg{-} : \edgeset(f) \to \lN$ is a bijection with its image $[0, 2(\iH_{I_\sop} + \iH_{I_\tap})]$. Thus $\edgeset(f)$ inherits a linear order $<\oedgeset$ on its elements by setting $\eda <\oedgeset \edb \iff \avg{\eda} < \avg{\edb}$ for $\eda,\edb \in \edgeset(f)$. We deduce the principle of \emph{edge induction}: Every element $\eda \in \edgeset(f)$ is of the form $\eda = \succn {0,0} m$ for $m = \avg{\eda}$.

\proof The proof for each statement is \stfwd{}.
\begin{enumerate} 

\item The first three statements follow from the observation that, for $\mathtt{r} \in \Set{\soe,\tae}$, we have firstly, $\succ \eda_{\mathtt{r}} = \eda_{\mathtt{r}} + 1$ if $\mathtt{r} = \succindex \eda$ and secondly, $\succ \eda_{\mathtt{r}} = \eda_{\mathtt{r}}$ if $\mathtt{r} \neq \succindex \eda$. This follows by separately inspecting each case $C^+$, $C \in \Set{\text{DL, UL, DR, UR}}$ of \autoref{defn:successor}. The last statement follows since, as we observed in \autoref{defn:successor}, in each case $C^+$ the successor edge $\edb = \succ\eda$ satisfies the conditions of case $C^-$ of \autoref{defn:predecessor}. Then, arguing separately in each case $C^-$, we can the see that $\pred \edb = \pred {\succ \eda}= \eda$ .

\item The argument is similar to the previous item, with the roles of successor and predecessor reversed.

\item Let $\eda \in \edgeset(f)$ and $\eda \equiv (\eda\ssoe ,\eda\ttae )$. Since $\eda\ssoe ,\eda\ttae $ are positive numbers we have $\avg{\eda} \geq 0$ and $\avg{\eda} = 0$ if and only if $\eda\ssoe  = \eda\ttae  = 0$. Note that $(0,0) \in \edgeset(f)$ by \eqref{eq:defn_order_realisation_3} together with $f\regop(0) = 0$ (which holds by \autoref{defn:singular_intervals_morphism}). Similarly, $\eda\ssoe  \leq 2\iH_{I_\sop}$ and $\eda\ttae  \leq 2\iH_{I_\tap}$. Thus $\avg{\eda} \leq 2\iH_{I_\sop} + 2\iH_{I_\tap}$ and $\avg{\eda} = 2\iH_{I_\sop} + 2\iH_{I_\tap}$ iff $\eda\ssoe  = 2\iH_{I_\sop}$ and $\eda\ttae  = 2\iH_{I_\tap}$. Now, note that $(2\iH_{I_\sop}, 2\iH_{I_\tap}) \in \edgeset(f)$ by \eqref{eq:defn_order_realisation_3} together with $f\regop(2\iH_{I_\tap}) = 2\iH_{I_\sop}$ (which holds by \autoref{defn:singular_intervals_morphism}).

\item Assume $\edb \in \edgeset(f)$ and $m = \avg \edb = \avg \eda$. Then $\predn \edb m = \predn \eda m = (0,0)$ by part (iii). Using that successor and predecessor are mutually inverse, we then calculate
\begin{align*}
\edb &= \succn{\predn \edb m} m \\
&= \succn{0,0} m \\
&= \succn{\predn \eda m} m \\
&= \eda
\end{align*}
as required. \qed
\end{enumerate}

\end{lem}

\begin{eg}[Ordering on edges] In the following example, the ordering on the set of edges is indicated by \cred{} arrows
\begin{restoretext}
\begingroup\sbox0{\includegraphics{test/page1.png}}\includegraphics[clip,trim=0 {.1\ht0} 0 {.05\ht0} ,width=\textwidth]{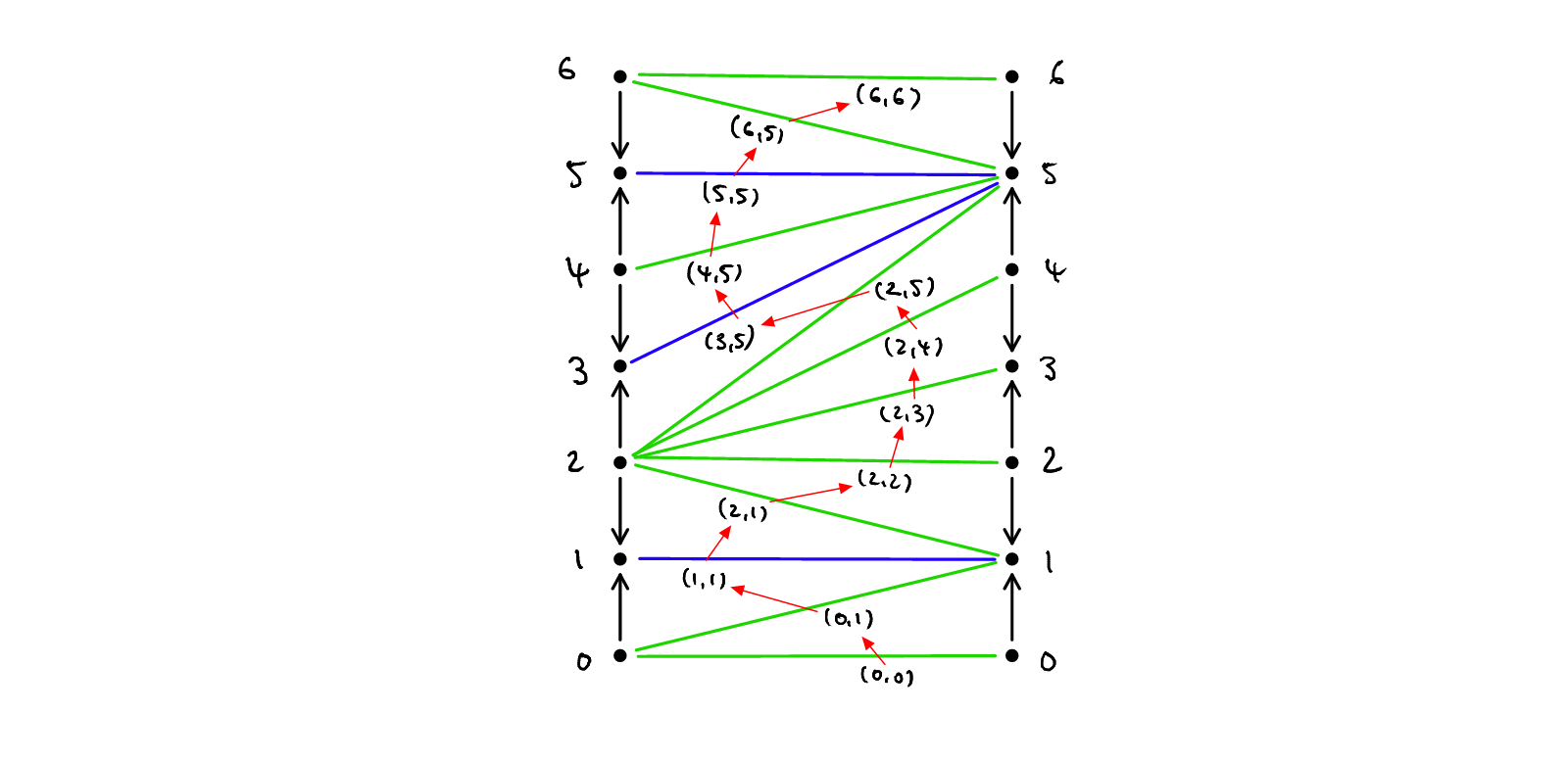}
\endgroup\end{restoretext}
\end{eg}

\begin{cor}[Bimonotonicity of profunctorial realisation]
\label{claim:order_realisations_monotone} Whenever $\eda,\edb \in \edgeset(f)$ then
\begin{gather}
(\eda\ssoe  <  \edb\ssoe ) \quad \imp \quad (\eda\ttae  \leq  \edb\ttae )
\end{gather}
and 
\begin{gather}
(\eda\ttae  <  \edb\ttae ) \quad \imp \quad (\eda\ssoe  \leq  \edb\ssoe )
\end{gather}
We refer to this property by saying that the relation $\SiR(f)$ is \emph{bimonotone}.
\proof We prove the first statement, and the second follows similarly. Assume $\eda\ssoe  < \edb\ssoe $. Then $\eda \neq \edb$ since equality of edges would imply $\eda\ssoe  = \edb\ssoe $. By \autoref{claim:edge_set_properties} we either have $\eda = \succn \edb m$ or $\eda = \succn \edb m$ for some $m \in \lN_{>0}$. The former case implies $\eda\ssoe  \geq \edb\ssoe $ by \autoref{claim:edge_set_properties} contradicting our assumption $\eda\ssoe  < \edb\ssoe $ and thus we must have $\eda = \predn \edb m$. Using \autoref{claim:edge_set_properties} again, this yields $\eda\ttae  \leq \edb\ttae $ as required. \qed
\end{cor}

\begin{cor}[Bisurjectivity of profunctorial realisation] \label{cor:relation_fullness} Let $f : I_\sop \to_{\SI} {I_\tap}$. Then for any $a \in {I_\tap}$ there is at least one $\eda \in \edgeset(f)$ with $\eda\ttae  = a$ and for any $b \in I_\sop$ there is at least one $\edb \in \edgeset(f)$ with $\edb\ssoe  = b$. We refer to this property by saying that the relation $\SiR(f)$ is \emph{bisurjective}.
\proof Applying the successor construction to $\edc \in \edgeset(f)$ increments either $\edc\ssoe $ or $\edc\ttae $ by $1$ while leaving the other constant (cf. proof of \autoref{claim:edge_set_properties}). Starting at $\edc = (0,0) \in \edgeset(f)$, since we will eventually reach the last edge $(2\iH_{I_\sop}, 2\iH_{I_\tap})$ the claimed edges must exist. \qed
\end{cor}

We also record the following statement, that will become useful later on.

\begin{cor}[Filler edge bounds] \label{cor:filler_edge} Let $\eda,\edb \in \edgeset(f)$. If $\avg{\eda} < \avg{\edb}$ then for $a \in \Set{\succfill \eda\ssoe , \succfill \eda\ttae }$ we have
\begin{equation}
\eda_{\succindex \eda} \leq a \leq \edb_{\succindex \eda}
\end{equation}

\proof First note that since $\avg{\eda} < \avg{\edb}$ we have $\avg{\eda} < 2\iH_{I_\sop} + 2\iH_{I_{\tap}}$, and thus the successor of $\eda$ exists by \autoref{defn:successor}. Next, $\avg{\eda} < \avg{\edb}$ implies $\avg{\succ\eda} \leq \avg{\edb}$ by \autoref{claim:edge_set_properties}, and the latter then also yields $\eda_{\succindex \eda} \leq \succ \eda_{\succindex \eda} \leq \edb_{\succindex \eda}$. By inspecting each case $C^+$, $C \in \Set{\text{DL, UL, DR, UR}}$ of \autoref{defn:successor} we see that $\succfill \eda$ is a pair $(\succfill \eda\ssoe , \succfill \eda\ttae ) \in E_{\succindex \eda}$ of numbers such that $\succfill \eda\ssoe , \succfill \eda\ttae  \in \Set{\eda_{\succindex \eda}, \succ\eda_{\succindex \eda}}$ which thus implies statement. \qed
\end{cor}

\section{Singular interval families} \label{sec:SI_fam}

Having defined $\SI$ and $\SiR : \SI \to \PRel$, in this section we put both ideas together to define a notion of singular interval family, and singular interval bundle. These bundles of intervals will later on form individual maps in towers of bundles, projecting down from $n$-dimensional space, one dimension at a time.

\subsection{Total posets of singular interval families}

We discuss total posets $\sG(\SiR\scA)$ of profunctorial realisations of $\SI$-families $\scA$, and introduce the definition of singular interval families as well as relevant notation.

\begin{constr}[Total posets of $\SI$-families] \label{constr:SI_families}
Given a poset $X$ and a functor $\scA : X \to \SI$, then by postcomposing $\scA$ with $\SiR : \SI \to \PRel$ we obtain a functor $\SiR\scA : X \to \PRel$ to which we can apply the Grothendieck construction given in \autoref{defn:grothendieck_construction}. We obtain a bundle $\pi_{\SiR\scA} : \sG(\SiR\scA) \to X$. Note that for $x \in X$, the fibre $\pi\inv_{\SiR\scA}(x)$ is a singular interval, or more precisely, we have
\begin{equation}
\pi\inv_{\SiR\scA}(x) = \Set{x} \times \SiR\scA(x)
\end{equation}
and can thus thus be identified with the singular interval $\scA(x)$. Following \autoref{defn:grothendieck_construction}, we will represent objects and morphisms as follows
\begin{itemize}
\item \textit{Objects}: An object $w \in \pi\inv_{\SiR\scA}(x) \subset \sG(\SiR\scA)$ is a tuple $(x,a)$, where $x \in X$ and $a \in \scA(x) = \Set{0, 1, \ldots, 2\iH_{\scA(x)}}$ is a natural number (cf. \autoref{defn:singular_intervals})
\item \textit{Morphisms}: $\big((x,a) \to (y,b)\big) \in \mor(\sG(\SiR\scA))$ if and only if $(x \to y) \in \mor(X)$ and $\SiR\scA(x \to y)(a,b)$ holds. Thus (cf. \autoref{defn:edge_sets}) we have 
\begin{align} \label{eq:grothendieck_edges}
(x,a) \to (y,b) \in &\mor(\sG(\SiR\scA))  \\
&\iff  (x \to y) \in \mor X \text{~and~} (a,b) \in \edgeset(\scA(x \to y))
\end{align}
A morphism in $\sG(\SiR\scA)$ is therefore a tuple $(r, \eda)$ where $r \in \mor(X)$ and $\eda \in \edgeset(\scA(r))$. 
\end{itemize}
\end{constr}

\begin{defn}[$\SI$-families]

For $X$ a poset, a functor $\scA : X \to \SI$ is called a \textit{singular interval family}, or more concisely, an \textit{$\SI$-family}, with base poset $X$. $\sG(\SiR\scA)$ is called its \textit{total poset} and $\pi_{\SiR\scA}$ is called an \SI-bundle.

Given two \SI-families $\scA_1, \scA_2 : X \to \SI$, a \textit{map of \SI-bundles} $G : \pi_{\SiR\scA_1} \to \pi_{\SiR\scA_2}$ is a map of bundles $G: \pi_{\SiR\scA_1} \to \pi_{\SiR\scA_2}$ that preserves direction order: in other words, we require for each $x \in X$ that $\rest G x$ is monotone.
\end{defn}

\begin{eg}[\SI-families] Using \autoref{eg:singint_map} consider the \SI-family $\scA : \singint 1 \to \SI$ defined by setting $\scA(0 \to 1) = f_1$ and $\scA(2 \to 1) = f_2$. In this case the Grothendieck construction yields the bundle $\pi_{\SiR\scA} : \sG(\SiR\scA) \to \singint 1$ which maps
\begin{restoretext}
\begingroup\sbox0{\includegraphics{test/page1.png}}\includegraphics[clip,trim=0 {.1\ht0} 0 {.1\ht0} ,width=0.8\textwidth]{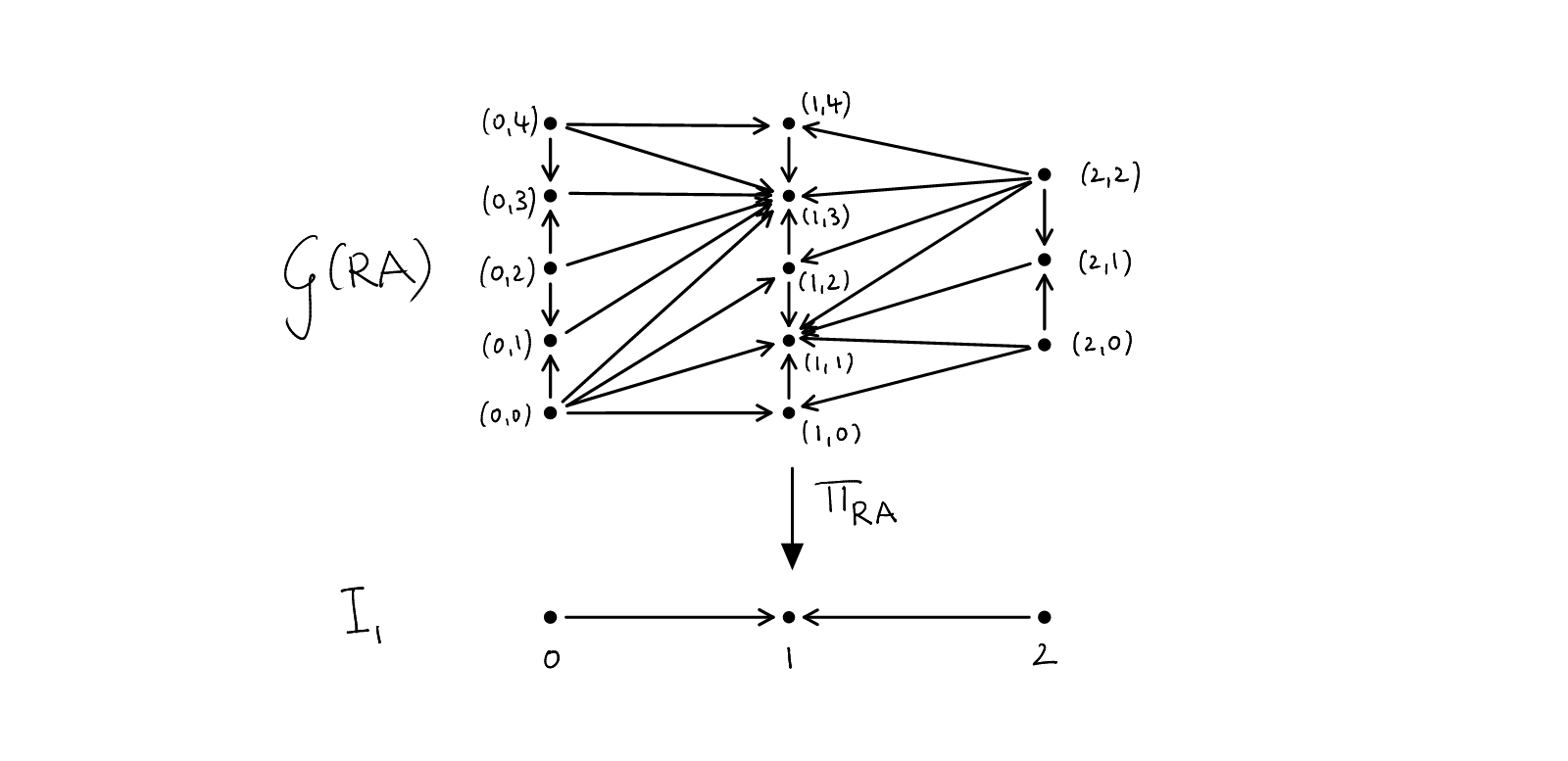}
\endgroup\end{restoretext}
\end{eg}

Note that since $\sG(\SiR\scA)$ is a poset, it comes with a notion of composition $(r,\eda) \circ (s,\edb)$ for compatible morphisms $(r,\eda)$, $(s,\edb)$ (and compatibility means $s = (w \to x)$, $r = (x \to y)$ for some $w,x,y \in X$ as well as $\eda\ssoe  = \edb\ttae $). In later proofs we will be particularly interested in the case of the successor construction. For this, we explicitly state the following observation
\begin{rmk}[Compositions with filler edges] \label{rmk:successor_compositionality} Let $(r,\eda) \in \mor(\sG(\SiR\scA))$ such that $\eda$ has a successor $\succ\eda$ in $\edgeset(\scA(r))$. Then,  writing $r = (x \to y)$, the four cases of \autoref{defn:successor} respectively yield the following compositions
\begin{enumerate}
\item[(DL$^+$)] $(\id_y,\succfill \eda) \circ (r,\eda) = (r,\succ \eda)$
\item[(UL$^+$)] $(r,\eda) \circ  (\id_x,\succfill \eda) = (r,\succ \eda)$
\item[(DR$^+$)] $(r,\eda) = (\id_y,\succfill \eda) \circ (r,\succ \eda)$
\item[(UR$^+$)] $(r,\eda) = (r,\succ \eda) \circ (\id_x,\succfill \eda)$
\end{enumerate}
which can be visualised for the cases (DL$^+$) and (UL$^+$) as
\begin{restoretext}
\begingroup\sbox0{\includegraphics{test/page1.png}}\includegraphics[clip,trim=0 {.15\ht0} 0 {.05\ht0} ,width=\textwidth]{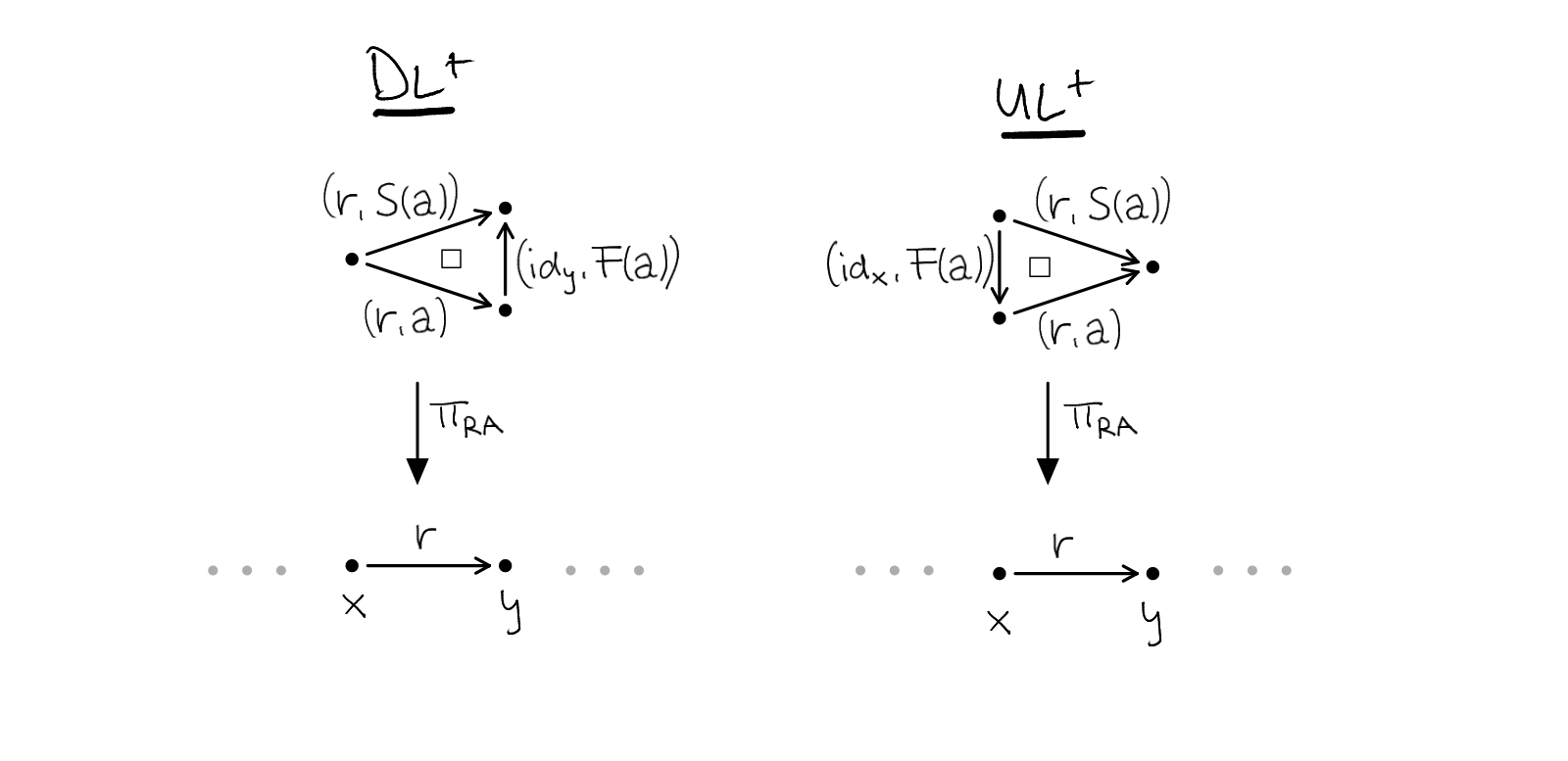}
\endgroup\end{restoretext}
and for the cases (DR$^+$) and (UR$^+$) as
\begin{restoretext}
\begingroup\sbox0{\includegraphics{test/page1.png}}\includegraphics[clip,trim=0 {.15\ht0} 0 {.05\ht0} ,width=\textwidth]{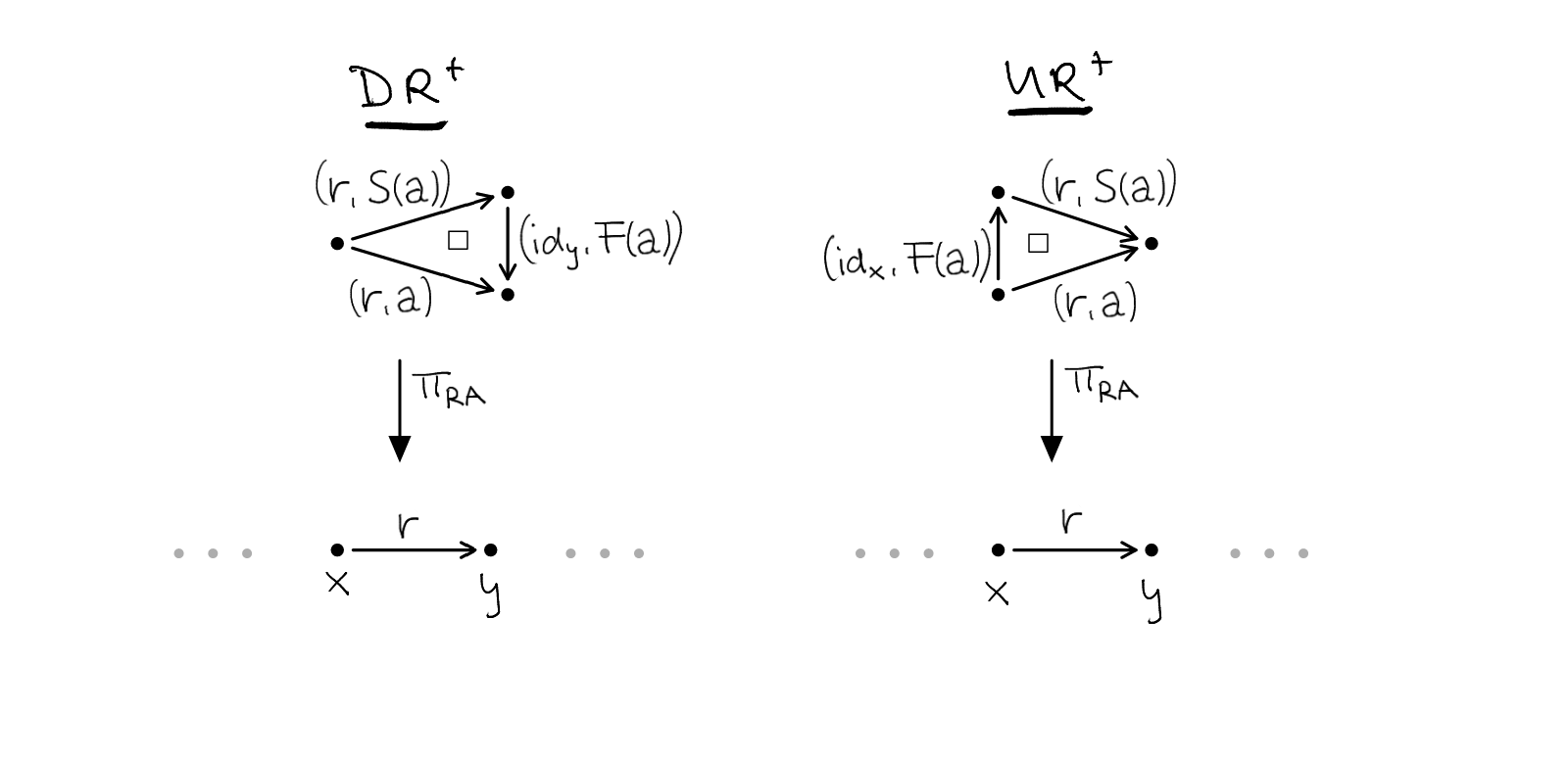}
\endgroup\end{restoretext}
\end{rmk}

The following are important simplifications of notation for the rest of the document.

\begin{notn}[Implicit profunctorial realisation] \label{notn:SI_families} To ease readability, throughout the rest of document we will denote the \SI-bundle 
\begin{itemize}
\item[] $\pi_{\SiR\scA} : \sG(\SiR\scA) \to X$
\end{itemize}
obtained from a singular interval family $\scA : X \to \SI$ by
\begin{itemize}
\item[] $\pi_{\scA} : \sG(\scA) \to X$
\end{itemize}
thus keeping post-composition with $\SiR$ implicit.
\end{notn}

\begin{notn}[Extension notation for $\SI$-families] \label{notn:SI_families2} We further extend \autoref{notn:singular_morphism_boundary_cases} as follows: if $X$ is a poset, $\scA : X \to \SI$ a singular interval family and $x \to y$ a morphism in $X$, then for readability we denote 
\begin{itemize}
\item[] $\wwidehat{\scA(x \to y)}$
\end{itemize}
by
\begin{itemize}
\item[] $\wwidehat{\scA} (x \to y)$
\end{itemize}
\end{notn}

We end this section with a first (simple) observation about $\sG(\scA)$.

\begin{rmk}[Connectedness of $\sG(\scA)$] \label{rmk:connectedness} If $X$ is a connected poset then for $\scA : X \to \SI$, $\sG(\scA)$ is a connected poset. This is because each $\pi\inv_\scA(x)$ is non-empty and connected, and $\edgeset(\scA(x \to y))$ is never empty.
\end{rmk}

\subsection{Strict uniqueness of pullbacks}

The following remark is an important observation which in addition to previous remarks highlights that many concepts evolving from singular intervals naturally have strict equality taking the place of a notion of isomorphism. 

\begin{rmk}[Uniqueness of pullbacks for $\SI$-families] \label{rmk:grothendieck_span_construction_basechange} 
In the setting of \SI-families (and unlike $\PRel$-families) pullbacks are unique on-the-nose. This is due to equivalence being strict equality as remarked in \autoref{rmk:direction_preserving_iso_is_eq}. 

More concretely, we first define what we mean by a pullback of \SI-bundles: let $H : X\to Y$ be a functor of posets and $B : Y \to \SI$ an $\SI$-family over $Y$. Then an  $\SI$-family pullback of $B$ along $H$ is a pullback of posets
\begin{align} \label{eq:SI_family_pullback}
\xymatrix@R=0.4cm{ \sG(\scA) \ar[r]^{K} \ar[d]_{\pi_{\scA}} \pullback & \sG(\scB) \ar[d]^{\pi_\scB} \\
X \ar[r]_{H} & Y }
\end{align}
such that $\pi_A$ is a bundle arising from an $\SI$-family $\scA : X \to \SI$ and $\rest K x$ is monotone for each $x \in X$ (cf. \autoref{notn:fiber_restrictions}).

Now, we show these pullbacks are unique up to strict equality. First, given a map $G : \pi_{\scA_1} \to \pi_{\scA_2}$ of \SI-bundles $\scA_1$, $\scA_2$ over $X$, we claim that if it is an isomorphism (that is, it has an inverse $G\inv : \pi_{\scA_2} \to \pi_{\scA_2}$) then we have $G = \sG(\id_X)$ as defined in \autoref{defn:grothendieck_base_change}. From \autoref{rmk:direction_preserving_iso_is_eq} we infer that $\scA_1(x) = \scA_2(x)$ and $\rest G x = \id$ for all $x \in X$. Note that functoriality and invertibility of $G$ then imply for all $(x \to y) \in \mor(X)$ that
\begin{equation}
\SiR\scA_1(x \to y)(a,b) \iff \SiR\scA_2(x \to y)(a,b)
\end{equation}
and thus by $\SiR$ being a faithful embedding (see \autoref{rmk:profunctorial_real_inj}) we find
\begin{equation}
\scA_1(x \to y) = \scA_2(x \to y)
\end{equation}
Therefore $G = \sG(\id_X)$ as claimed.

Next, assume the pullback \eqref{eq:SI_family_pullback}. We claim $\scA = \scB H$ showing that $\scA$ is uniquely determined by $\scB$ and $H$. First, by \autoref{defn:grothendieck_base_change} we have a pullback of \SI-bundles
\begin{align}
\xymatrix@R=0.4cm{ \sG(\scB H) \ar[r]^{\sG(H)} \ar[d]_{\pi_{\scB H}} \pullback & \sG(\scB) \ar[d]^{\pi_\scB} \\
X \ar[r]_{H} & Y }
\end{align}
By the pullback property there is a poset isomorphism $F : \sG(\scA) \to \sG(\scB H)$ such that $\pi_{\scB H} F = \pi_\scA$ and $\sG(H) F = K$. Since $K$ and $\sG(H)$ preserve direction order (that is, they are fibrewise monotone) so does $F$. This makes $F$ into a map of \SI-bundles. We deduce $F = \sG(\id_X)$ by our previous argument and thus pullbacks are unique.
\end{rmk}

The preceding remark establishes a one-to-one correspondence of $\SI$-bundle pullbacks \eqref{eq:SI_family_pullback} and base change functors $H : X \to Y$.

\chapter{Cubes} \label{ch:cubes}

In this chapter we build a notion of labelled singular $n$-cubes from the previously defined notion of intervals (that is, $1$-cubes). These $n$-cubes will be towers (of height $n$) of bundles of intervals, together with a functor (called ``labelling") on their total space. Interestingly, these cubes organise into a category whose morphism are bundles of cubes themselves. \autoref{sec:labels} will be dedicated to the inductive step in building this category. \autoref{sec:cube_fam} will then use this step to build the category of labelled singular $n$-cubes and discuss their representation as towers of bundles, as well as their most important properties.

\section{Labels} \label{sec:labels}

The goal of this section will be to define the category $\SIvertone  \cC$ of $\cC$-labelled singular intervals, and then prove that there is an one-to-one correspondence between functors $F$ from a poset $X$ into $\SIvertone  \cC$
\begin{equation}
R : X \to \SIvertone  \cC
\end{equation}
and tuples of functors $(V,U)$
\begin{equation}
\xymatrix{ E \ar[d]_{V} \ar[r]^{U} & \cC \\ X & }
\end{equation}
where $V$ is an $\SI$-bundle, and $U$ is any functor (which will be called the labelling of $R$). Using the examples provided early onwards in the section, the reader might want to attempt to spot this correspondence as soon as possible, as this will allow us to gloss over some of the tedious notation involved in the fully formal construction of it.

\subsection{Labelled singular intervals}

\begin{notn} \label{notn:functor_morphism_action} Given a functor $F : \cC \to \cD$ denote by
\begin{equation}
\nu_F : \Hom_\cC(-,-) \imp \Hom_\cD(F-,F-)
\end{equation}
the natural transformation with components
\begin{equation}
{(\nu_F)}_{(x,y)} : (g : x\to y) \mapsto (F g : F x \to F y)
\end{equation}
\end{notn}

\begin{constr}[$\cC$-labelled singular intervals] \label{defn:SIvert} Given a category $\cC$, define the category $\SIvertone  {\cC}$, called the category of \textit{$\cC$-labelled singular intervals}, as follows: 
\begin{itemize}
\item Objects $I$ are pairs $(\SIf I \in \SI, \SIs I : \SIf I \to \cC)$ consisting of a singular interval $\SIf I$ and a functor $\SIs I : \SIf I \to \cC$. 
\item Morphisms $f : I \to J$ are pairs $(\SIf f, \SIs f) : (\SIf I,\SIs I) \to (\SIf J, \SIs J)$ consisting of an \SI-morphism
\begin{equation}
\SIf f : \SIf I \to_{\SI} \SIf J
\end{equation}
of singular intervals together with a natural transformation
\begin{equation} \label{eq:mor_in_SIvert}
\SIs f : \iP\SiR(\SIf f)(-,-) \to \Hom_\cC(\SIs I-,\SIs J-)
\end{equation}
Note that $\iP$ will be usually kept implicit from now on.
\end{itemize}
Note that components of $\SIs f$ are maps
\begin{equation}
\SIs f _{(a,b)} : \iP \SiR(\SIf f)(a,b) \to \Hom_\cC(\SIs I a,\SIs J b)
\end{equation}
If $\SiR(\SIf f)(a,b)$ then we write $\SIs f_{(a,b)}$ for the morphism $\SIs f_{(a,b)}(*)$ (which is an element in $\Hom_\cC (\SIs I (a), \SIs J (b))$). On the other hand, if $(a,b) \notin \edgeset(\SIf f)$ then $\SIs f _{(a,b)}$ is trivial as it has empty domain and image.
\begin{itemize}
\item The composition of two morphisms $f_1 = (\SIf {f_1},\SIs {f_1}) : I_1 \to I_2$, $f_2 = (\SIf {f_2},\SIs {f_2}) : I_2 \to I_3$ is given by a pair
\begin{equation}
(\SIf {f_2} \circ \SIf {f_1}, \SIs {f_1} \odot \SIs {f_2}) : I_1 \to I_3
\end{equation}
where $\SIs {f_1} \odot \SIs {f_2}$ is the ``horizontal composition" of $\SIs {f_1}$ and $\SIs {f_2}$ defined as follows: for any $(a_1,a_3)$ satisfying the left hand side of \eqref{eq:horizontal_comp_SIC} we need to define $(\SIs {f_1} \odot \SIs {f_2})_{(a_1, a_3)}$. Recall that $\SiR(\SIf {f_2}\circ \SIf {f_1}) = \SiR(\SIf {f_1}) \odot \SiR(\SIf {f_2})$ by \autoref{lem:order_realisation_functorial}, and thus by \autoref{defn:prerequisites} we find
\begin{equation} \label{eq:horizontal_comp_SIC}
(a_1,a_3) \in \edgeset(\SIf {f_2}\SIf {f_1}) \quad \iff \quad \exists a_2. (a_1,a_2) \in \edgeset(\SIf {f_1}) \text{~and~} (a_2,a_3) \in \edgeset(\SIf {f_2}) 
\end{equation}
Choose any $a_2$ for which the right hand side holds, and then set
\begin{equation} \label{eq:horizontal_composition}
(\SIs {f_1} \odot \SIs {f_2})_{(a_1, a_3)} := \SIs {f_2}_{(a_2, a_3)} \circ \SIs {f_1}_{(a_1, a_2)}
\end{equation}
In \autoref{claim:horizontal_composition_well_def} below we show that this definition is well-defined; that is, it doesn't depend on the choice of the integer $a_2 \in I_2$. 
\end{itemize}
This constructs the category $\SIvertone  \cC$. We note that identities in $\SIvertone {\cC}$ are given by $\id_{I} = (\id_{\SIf I}, \nu_{\SIs I})$ (cf. \autoref{notn:functor_morphism_action}). 
\end{constr}

\begin{rmk}[Labelling terminology] $\SIvertone {\cC}$ is called the category of ``$\cC$-labelled singular intervals" based on the following idea: if $I \in \SIvertone {\cC}$ then $\SIs I : \SIf I \to \cC$ should be thought of as a singular interval whose objects $a$ and morphisms $a_1\to a_2$ are labelled by objects $\SIs I (a)$ and morphisms $\SIs I (a_1\to a_2)$ in $\cC$ respectively. 

Similarly, if $(f : I \to J) \in \mor(\SIvertone {\cC})$ then $\SIs f_{(a,b)} \in \mor(\cC)$ should be thought of as a label for the edge $(a,b) \in \edgeset(\SIf f)$ of $\SIf f$. By naturality of $\SIs f$ these labels are ``compatible" with the labels $\SIs I$, $\SIs J$ of $\SIf I$, $\SIf J$ in the natural way as we will see in the next example.
\end{rmk}

\begin{eg}[A morphism in $\SIvertone  \cC$] \label{eg:SIvert} We define the objects $(\SIf I, \SIs I), (\SIf J, \SIs J) \in \SIvertone  \cC$ by setting $\SIf I = \singint 2$, $\SIf J = \singint 2$ and defining $\SIs I, \SIs J : \singint 2 \to \cC$ to be the functors
\begin{restoretext}
\begingroup\sbox0{\includegraphics{test/page1.png}}\includegraphics[clip,trim=0 {.3\ht0} 0 {.25\ht0} ,width=\textwidth]{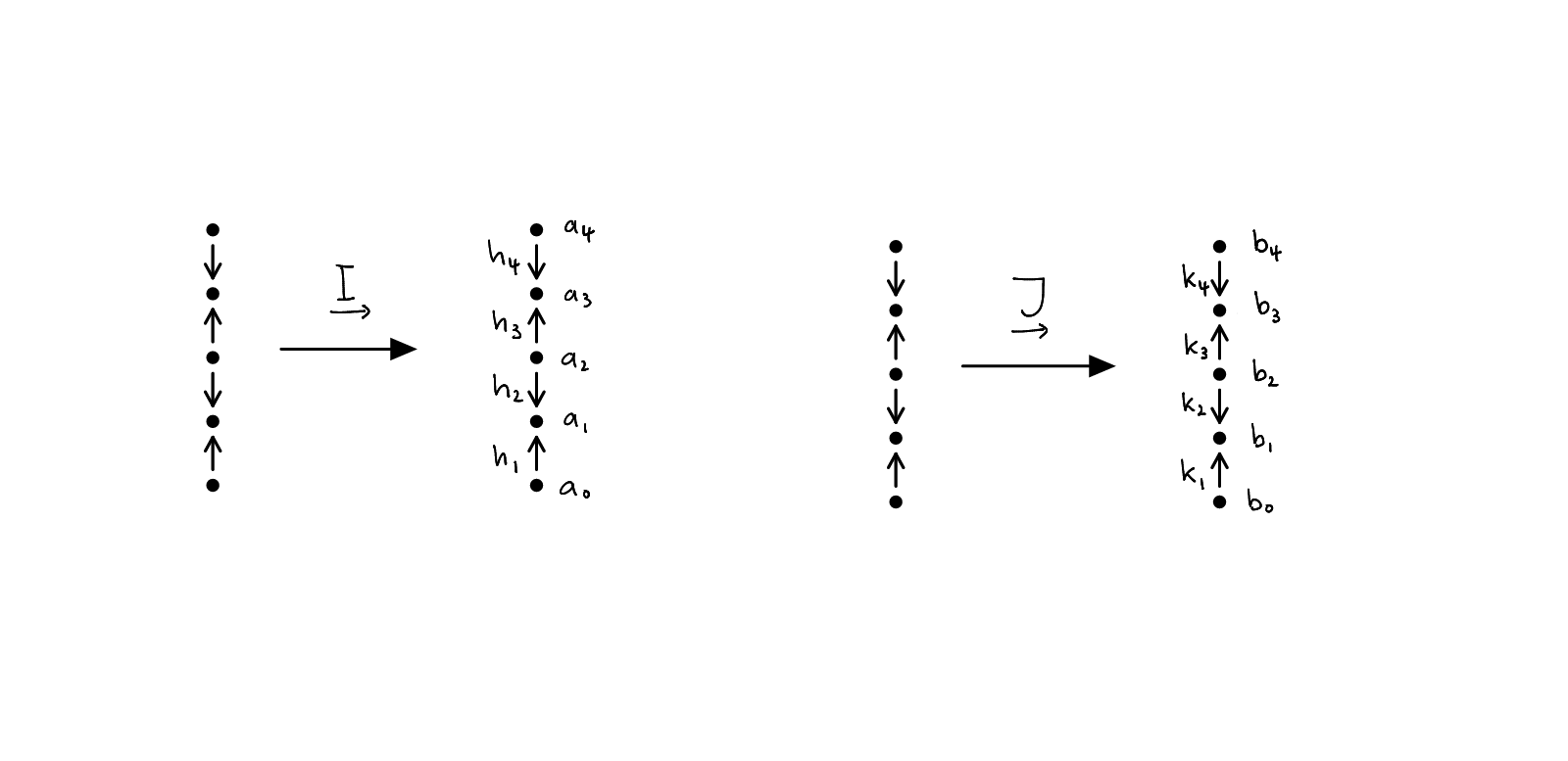}
\endgroup\end{restoretext}
where $a_i, b_i$ are objects, and $h_i, k_i$ morphisms in $\cC$. We define $(\SIf f , \SIs f) : (\SIf I, \SIs I) \to (\SIf J, \SIs J)$ by first setting $\SIf f : \SIf I \to_{\SI} \SI f$ to be the \SI-morphism with profunctorial realisation 
\begin{restoretext}
\begingroup\sbox0{\includegraphics{test/page1.png}}\includegraphics[clip,trim=0 {.35\ht0} 0 {.25\ht0} ,width=\textwidth]{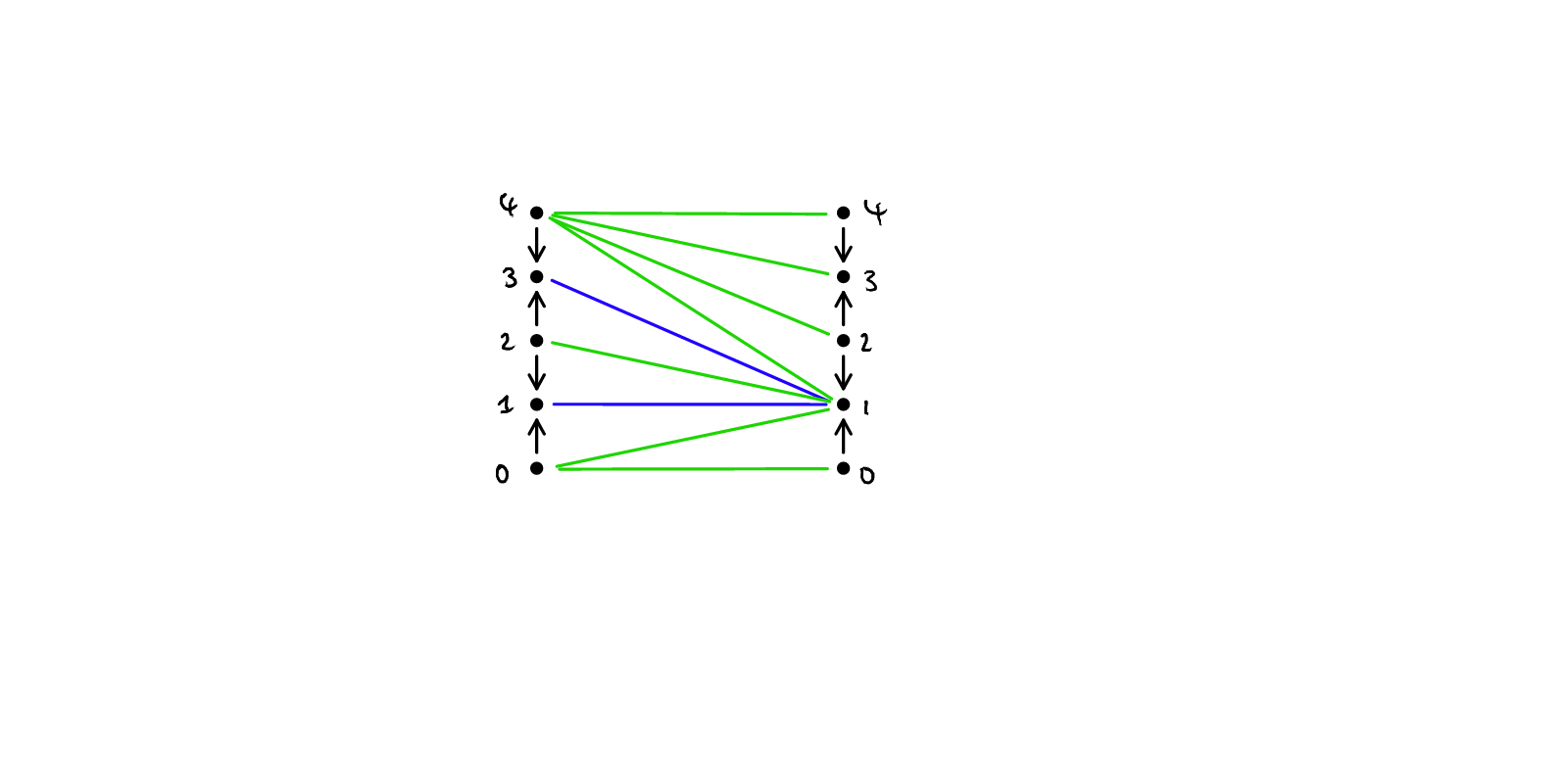}
\endgroup\end{restoretext}
(the mapping of $\SIf f$ itself is indicated by the \cblue{} edges). Now, $\SIs f$ is a natural transformation as defined in \eqref{eq:mor_in_SIvert} with non-trivial components $(\SIs f)_{(a,b)}$ whenever $(a,b) \in \edgeset (\SIf f)$ is an edge of $f$, in which case we have $(\SIs f)_{(a,b)} \in \Hom_\cC (\SIs I (a), \SIs J (b))$ (cf. \autoref{defn:SIvert}). Naturality (in the first component) of $\SIs f$ means that for all $(a' \to a) \in \mor(SIf I)$ we have
\begin{align}
 (\SIs f)_{(a',b)} &= \Hom_\cC(\SIs I (a' \to a),\SIs J (b)) \big((\SIs f)_{(a,b)}\big) \\
 &= (\SIs f)_{(a,b)}) \SIs I (a' \to a)
\end{align}
and naturality in the second component means for all $(b \to b') \in \mor(\SIf J)$ we have
\begin{align}
 (\SIs f)_{(a,b')} &= \Hom_\cC(\SIs I (a),\SIs J (b \to b')) \big((\SIs f)_{(a,b)}\big) \\
 &= \SIs J (b \to b') (\SIs f)_{(a,b)}) 
\end{align}
This means, that data for $\SIf f$ is an assignment $(a,b) \in \edgeset(f) \mapsto (\SIs f)_{(a,b)}$ such that 
\begin{restoretext}
\begingroup\sbox0{\includegraphics{test/page1.png}}\includegraphics[clip,trim=0 {.1\ht0} 0 {.1\ht0} ,width=\textwidth]{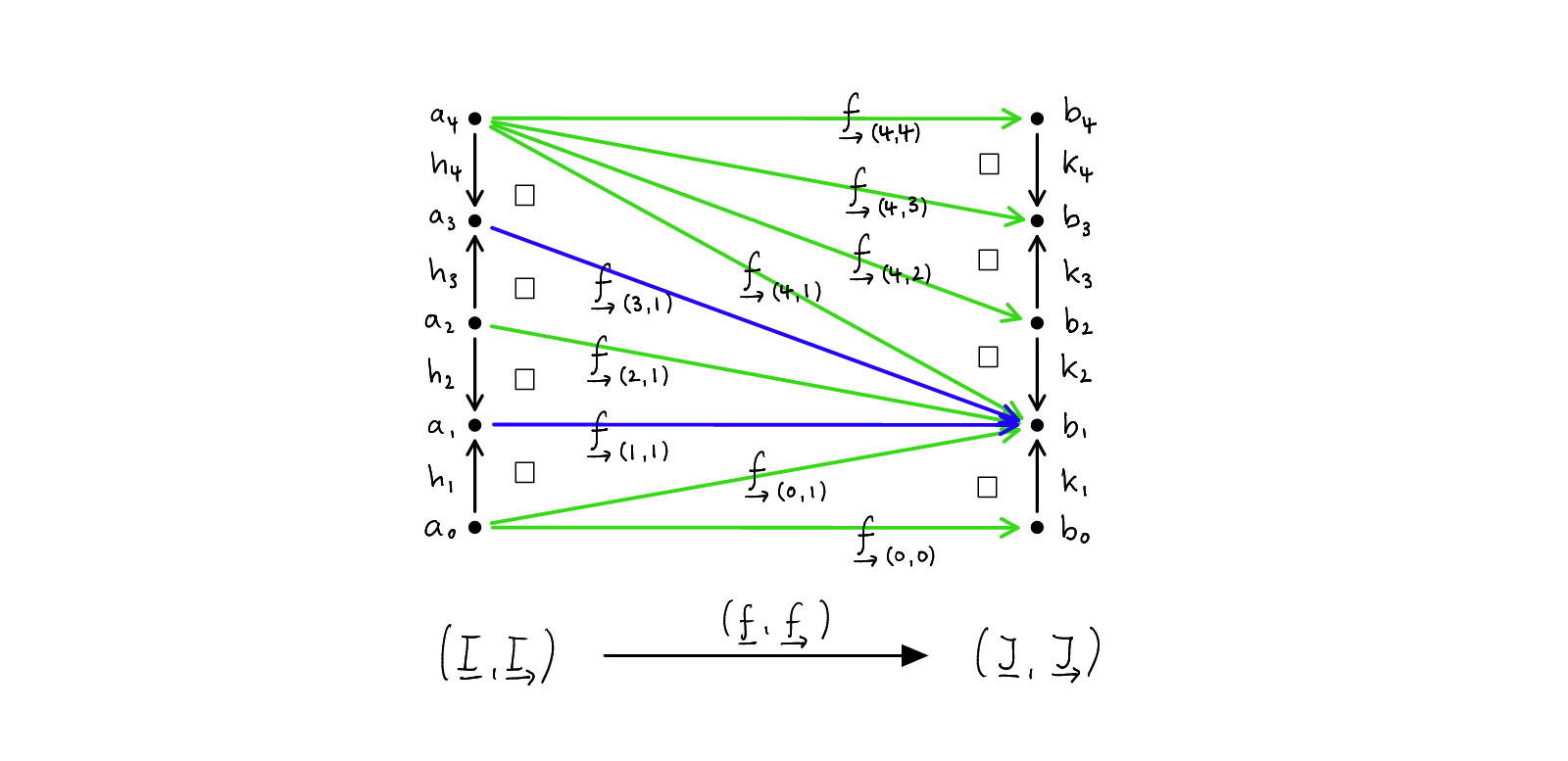}
\endgroup\end{restoretext}
commutes. But the commuting diagram above (without the labels by morphisms in $\cC$) is just $\sG(\Delta_{\SIf f})$, where $\Delta_{\SIf f} : \bnum 2 \to \SI$ has image $\und f$. This means that $(\SIf f , \SIs f) : (\SIf I, \SIs I) \to (\SIf J, \SIs J)$ gives a functor
\begin{equation} \label{eq:labelled_total_poset_of_single_mor}
\sG(\Delta_{\SIf f}) \to \cC
\end{equation}
And conversely, any such functor determines $\SIs I$, $\SIs J$ and $\SIs f$.
\end{eg}

\begin{eg}[Composition of morphisms in $\SIvertone  \cC$] \label{eg:SIvert_hor_comp} Given two composable morphisms $(f_1 : I_1 \to I_2), (f_2 : I_2 \to I_3) \in \mor(\SIvertone  \cC)$ then the components of $(\SIs {f_1 \odot f_2})_{(a,c)}$ are compositions of components of $\SIs {f_1}$ and $\SIs {f_2}$ (cf. \eqref{eq:horizontal_composition}) but the choice of this composition is non-unique as the following example shows
\begin{restoretext} 
\begingroup\sbox0{\includegraphics{test/page1.png}}\includegraphics[clip,trim=0 {.07\ht0} 0 {.15\ht0} ,width=\textwidth]{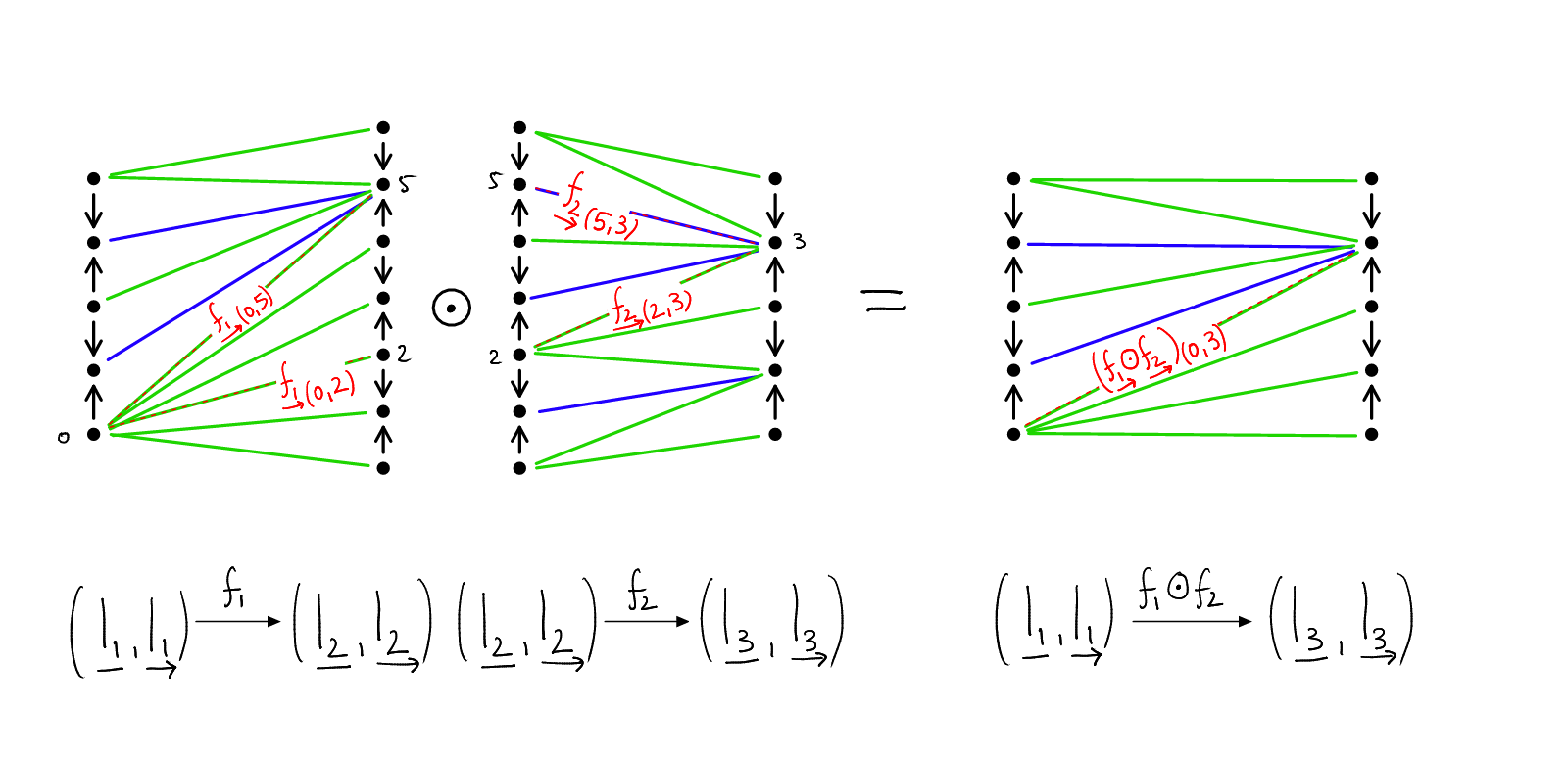}
\endgroup\end{restoretext}
The fact that horizontal composition is still well-defined is the content of the next lemma.
\end{eg}

\begin{claim}[Well-definedness of compostion in $\SIvertone  {\cC}$] \label{claim:horizontal_composition_well_def} $\SIs {f_1} \odot \SIs {f_2}$ as defined in \eqref{eq:horizontal_composition} satisfies that
\begin{enumerate}
\item it is single-valued. In other words, $(\SIs {f_1} \odot \SIs {f_2})_{(a_1,a_3)}$ is independent of the choice of $a_2$. 
\item it is a natural transformation
\begin{equation}
\SiR(\SIf {f_2 f_1})(-,-) \to \Hom_\cC(\SIs {I_1}-,\SIs {I_3}-)
\end{equation}
\end{enumerate}

\proof 
\begin{enumerate}
\item To show that the definition given in \eqref{eq:horizontal_composition} is single-valued, assume $\SiR(\SIf {f_1})(a_1,a_2)$ and $\SiR(\SIf {f_2})(a_2,a_3)$ as well as $\SiR(\SIf {f_1})(a_1,a'_2)$ and $\SiR(\SIf {f_2})(a'_2,a_3)$ for $a_2, a'_2 \in I_2$. In other words, both $a_1$ and $a_2$ are witnesses for the existential quantification in \eqref{eq:horizontal_comp_SIC}. We need to show that
\begin{equation} \label{eq:horizontal_composition_well_def_1}
(\SIs {f_2})_{(a_2, a_3)} \circ (\SIs {f_1})_{(a_1, a_2)} = (\SIs {f_2})_{(a'_2, a_3)} \circ (\SIs {f_1})_{(a_1, a'_2)}
\end{equation}
Assume $a_2 < a'_2$ or switch the roles of $a_2$, $a'_2$ in the following. Note that for any $a_2 < a < a'_2$ we have $\SiR(\SIf {f_1})(a_1,a)$ and $\SiR(\SIf {f_2})(a,a_3)$ by monotonicity (cf. \autoref{claim:order_realisations_monotone}). Thus, arguing inductively, it is enough to show \eqref{eq:horizontal_composition_well_def_1} for $a'_2 = a_2 + 1$. There are two cases.
\begin{enumerate}
\item $a_2$ is even: in this case $a_2 \to a'_2$ and we find
\begin{align}
(\SIs {f_2})_{(a_2, a_3)} \circ (\SIs {f_1})_{(a_1, a_2)} &= (\SIs {f_2})_{(a'_2, a_3)} \circ 
\SIs {I_2} (a_2 \to a'_2) \circ (\SIs {f_1})_{(a_1, a_2)} \\
&= (\SIs {f_2})_{(a'_2, a_3)} \circ (\SIs {f_1})_{(a_1, a'_2)}
\end{align}
where we used naturality of $\SIs {f_2}$ in the first step, and naturality of $\SIs {f_1}$ in the second.
\item $a_2$ is odd: in this case $a'_2 \to a_2$ and we find
\begin{align}
(\SIs {f_2})_{(a_2, a_3)} \circ (\SIs {f_1})_{(a_1, a_2)} &= (\SIs {f_2})_{(a_2, a_3)} \circ \SIs {I_2} (a'_2 \to a_2) \circ (\SIs {f_1})_{(a_1, a'_2)}  \\
&= (\SIs {f_2})_{(a'_2, a_3)} \circ (\SIs {f_1})_{(a_1, a'_2)}
\end{align}
where we used naturality of $\SIs {f_1}$ in the first step, and naturality of $\SIs {f_2}$ in the second.
\end{enumerate}

\item Finally, the fact that $\SIs {f_1} \odot \SIs {f_2}$ is a natural transformation now follows from \eqref{eq:horizontal_composition} and  $\SIs {f_1}$, $\SIs {f_2}$ being natural transformations individually.
\end{enumerate}
 \qed
\end{claim}

\begin{constr}[Relabelling functor] \label{defn:transfer_of_coloring} Given a functor $F : \cC \to \cD$, we construct the functor $\SIvertone F : \SIvertone {\cC} \to \SIvertone \cD$, called \textit{relabelling} by $F$, as follows
\begin{itemize}
\item Objects $I = (\SIf I, \SIs I) $ are mapped to $(\SIf I,F\SIs I)$
\item Morphisms $f = (\SIf f, \SIs f): I \to J$ are mapped to $(\SIf f, \nu_F(\SIs I\op \times \SIs J) \circ \SIs f)$, where $\nu_F(\SIs I\op \times \SIs J) :  \Hom_\cC(\SIs I-,\SIs J-) \to \Hom(F\SIs I-,F\SIs J-)$ is $\nu_F$ (cf. \autoref{notn:functor_morphism_action}) pre-whiskered with $\SIs I\op \times \SIs J : \SIf I \op \times \SIf J \to \cC\op \times \cC$.
\end{itemize}
For functoriality of $\SIvertone  {F}$ we need to verify
\begin{equation} \label{eq:label_transfer_functoriality}
\SIs {\SIvertone  {F} f_1} \odot \SIs {\SIvertone  {F} f_2} = \SIs {\SIvertone  {F} (f_1 \odot f_2)}
\end{equation}
Note that by definition of $Ff$ and $\nu_F$ (cf. \autoref{notn:functor_morphism_action}), for $f : I \to J \in \mor(\SIvertone  \cC)$, $a \in \SIf I, b \in \SIf J$, we find
\begin{align} \label{eq:label_transfer_action}
(\SIvertone  {F} f)_{(a,b)} &:= \big(\nu_F(\SIs {I}\op \times \SIs {J}) \circ \SIs {f}\big)_{(a,b)} \\ &= F (\SIs f_{(a,b)}) ~ : ~ FI(a) \to FJ(b)
\end{align}
and thus \eqref{eq:label_transfer_functoriality} follows from functoriality of $F$ and \eqref{eq:horizontal_composition}.
\end{constr}

\begin{rmk}[Compositionality of relabelling for \SI-families] \label{rmk:transfer_compositional} The preceding definition entails that for $F : \cC \to \cD$, $G : \cD \to \cE$ we have
\begin{gather}
\SIvertone  G \SIvertone F = \SIvertone  {GF}
\end{gather}
\end{rmk}

\begin{eg}[Relabelling] Using \autoref{eg:SIvert}, the action of the relabelling functor can be illustrated by the following mapping
\begin{restoretext}
\begingroup\sbox0{\includegraphics{test/page1.png}}\includegraphics[clip,trim=0 {.25\ht0} 0 {.1\ht0} ,width=\textwidth]{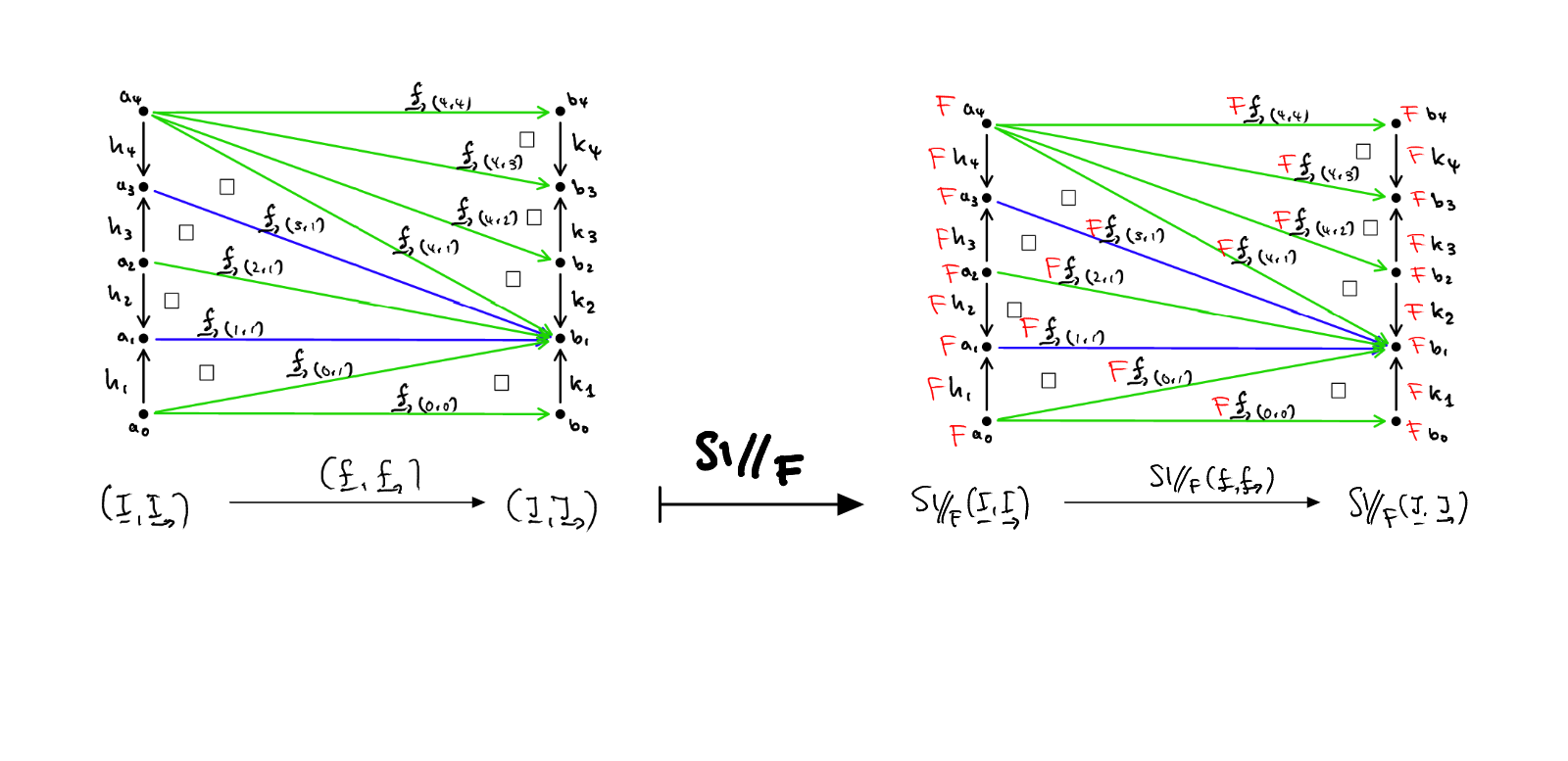}
\endgroup\end{restoretext}
\end{eg}

\begin{defn}[Forgetting labels] \label{defn:label_forgetting} The \textit{label-forgetting functor} $\und{(-)} : \SIvertone {\cC} \to \SI$ is defined to map $I \in \obj(\SIvertone \cC)$ to $\SIf I$ and $(f : I \to J) \in \mor(\SIvertone \cC)$ to $\SIf f : \SIf I \to_{\SI} \SIf J$.
\end{defn}

\begin{eg}[Forgetting labels] Using \autoref{eg:SIvert}, the action of the label-forgetting functor can be illustrated by the following mapping
\begin{restoretext}
\begingroup\sbox0{\includegraphics{test/page1.png}}\includegraphics[clip,trim=0 {.2\ht0} 0 {.15\ht0} ,width=\textwidth]{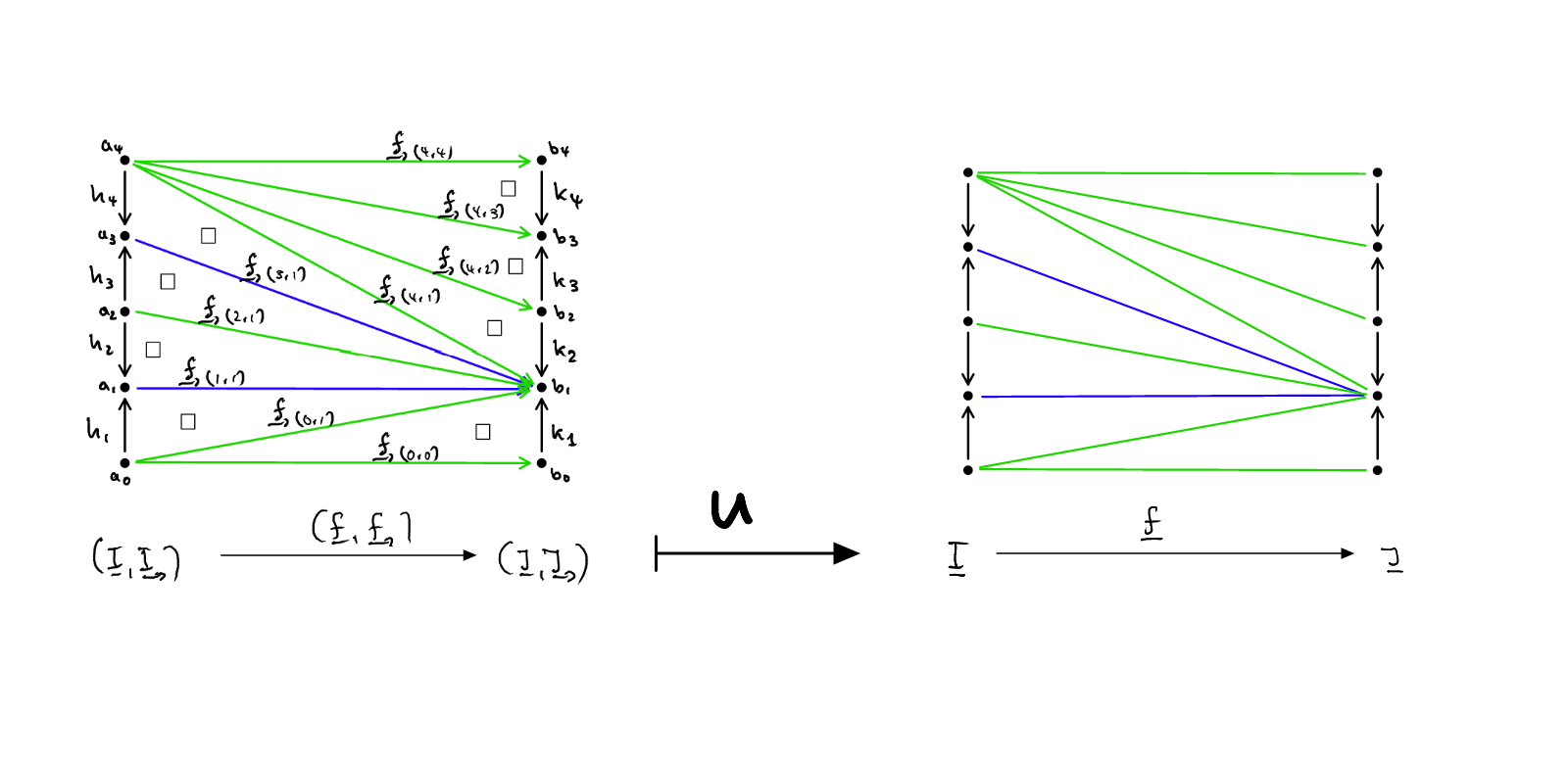}
\endgroup\end{restoretext}
\end{eg}

The following remark gives a broader context in which to see the preceding definitions.

\begin{rmk}[Relation to double category $\PProf$] The ``double-slash" notation for $\SIvertone \cC$ is chosen to indicate that its construction is analogous to over-categories in the double category $\PProf$. $\PProf$ has categories (as objects), functors (as vertical arrows), profunctors (as horizontal arrows) and natural transformations (as squares). Morphisms $f : I \to J$ in $\SIvertone  {\cC}$ then correspond to squares in $\PProf$ of the form
\begin{equation}
\xymatrix@C=1.5cm@R=0.8cm{ \SIf I \ar[r]^-{\SiR(\SIf f)}|-*=0@{|} \ar[d]_{\SIs I} \ar@{}[dr]|{ \Downarrow\SIs f} & \SIf J \ar[d]^{\SIs J} \\
\cC \ar[r]_{\Hom_{\cC}}|-*=0@{|} & \cC }
\end{equation}
Since $\Hom_{\cC}$ however plays the role of an (horizontal) identity, such a square can be thought of as a ``square over $\cC$". However, note that $\SIvertone \cC$ (unlike $\PProf$ or, for instance, $\Bool$-$\PProf$) does not have interesting $2$-categorical structure, since $\SI$ does not admit such structure either.
\end{rmk}

\subsection{Labelled singular intervals functor}

Using \autoref{rmk:transfer_compositional}, we can now define

\begin{defn}[Labelled singular intervals functor] \label{defn:SIvertone_functor} The labelled \textit{singular interval functor}
\begin{equation}
\SIvertone {-} : \Cat \to \Cat
\end{equation}
is defined to map $\cC$ to $\SIvertone \cC$ and $F : \cC \to \cD$ to $\SIvertone F$.
\end{defn}

We also introduce the following terminology (effectively, these are objects in the comma category $\Pos \slash (\SIvertone {-})$).

\begin{defn}[$\SIvertone  \cC$-family] \label{defn:SIvert_fam} A \textit{$\cC$-labelled singular interval family} $\scA$ over a poset $X$, also called a $\SIvertone  \cC$-family over $X$, is a functor $\scA : X \to \SIvertone {\cC}$.
\end{defn}

\begin{notn}[Label and label-forgetting notation of $\cC$-labelled $\SI$-families] Given a poset $X$, a morphisms $(x \to y)$ in $X$, and a $\cC$-labelled singular interval family $\scA : X \to \SIvertone  \cC$ over $X$, then (in order to simplify notation) we denote
\begin{itemize}
\item[] $\SIs{\scA(x)}, \SIf{\scA(x)}, \SIs{\scA(x \to y)}$ and $\SIf {\scA(x \to y)}$
\end{itemize}
by
\begin{itemize}
\item[] $\SIs \scA(x), \SIf \scA(x), \SIs \scA(x \to y)$ and $\SIf \scA(x \to y)$
\end{itemize}
respectively.
\end{notn}

\subsection{Unpacking and repacking labels}

In this section we discuss a ``generalised Grothendieck construction" which applies to $\SIvertone  \cC$. That is, $\SIvertone  \cC$ turns out to be the classifying space for $\SI$-bundles $\pi_V$ together with a labelling functor $U : \sG(V) \to \cC$  on its total space. We start by illustrating this behaviour in the following example. 

\begin{eg}[Unpacking $\SIvertone \cC$-families] \label{eg:unpacking_SIvert} Using \autoref{eg:SIvert_hor_comp} define an $\SIvertone \cC$-family $\scA : \bnum{3} \to \SIvertone \cC$ as follows
\begin{restoretext}
\begingroup\sbox0{\includegraphics{test/page1.png}}\includegraphics[clip,trim=0 {.3\ht0} 0 {.2\ht0} ,width=\textwidth]{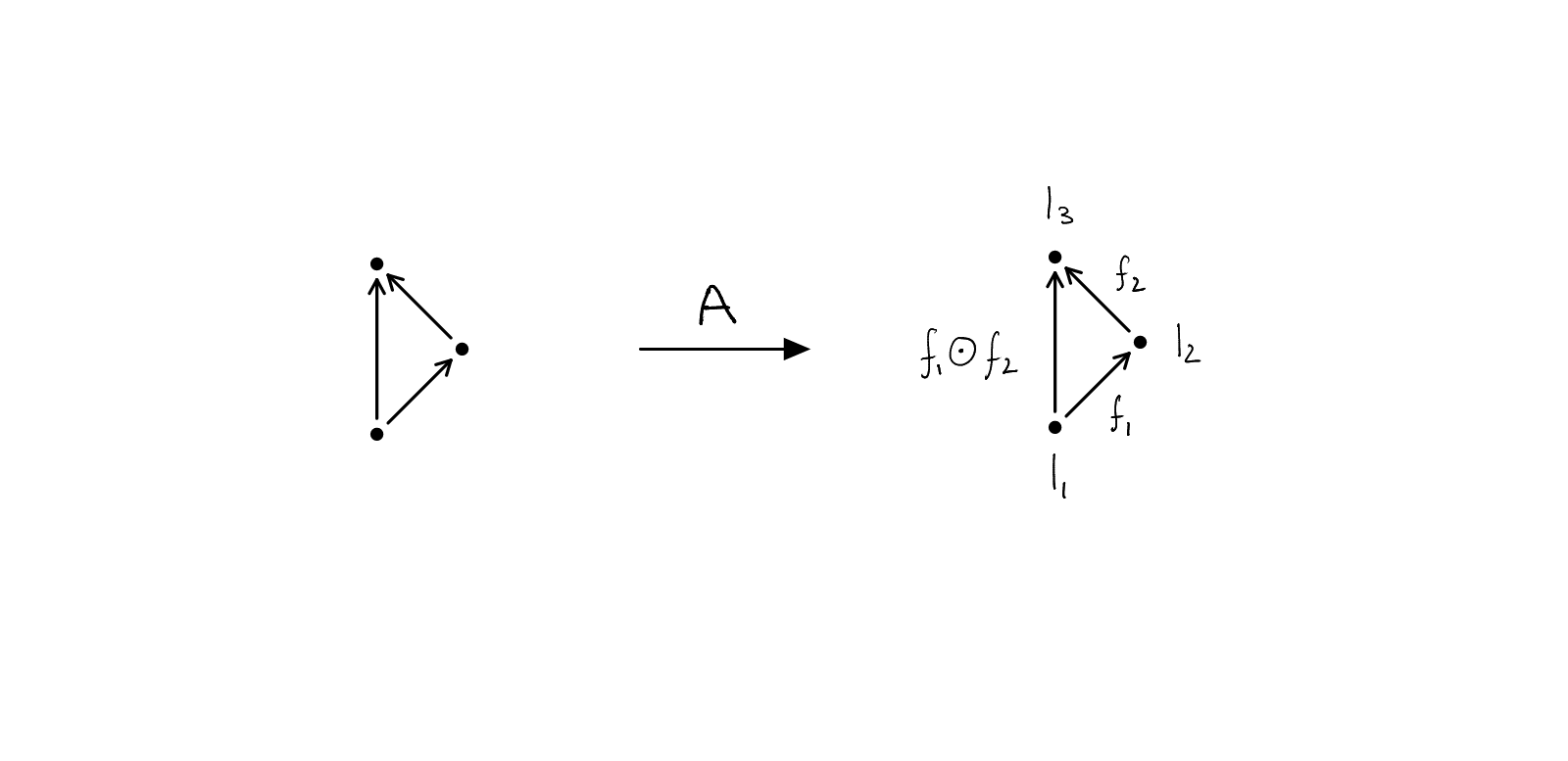}
\endgroup\end{restoretext}
Then, for each $i,j \in \bnum{3}$, $i < j$, using the argument of \autoref{eg:SIvert} for the morphism $\scA(i \to j)$, we find a functor analogous to  \eqref{eq:labelled_total_poset_of_single_mor} which maps
\begin{equation}
\sG(\Delta_{\und\scA(i \to j)}) \to \cC
\end{equation}
and which ``contains" the labelling of $\scA(i \to j)$. These functors glue together to give a functor $\sU_\scA : \sG(\und \scA) \to \scC$
\begin{restoretext}
\begingroup\sbox0{\includegraphics{test/page1.png}}\includegraphics[clip,trim=0 {.0\ht0} 0 {.0\ht0} ,width=\textwidth]{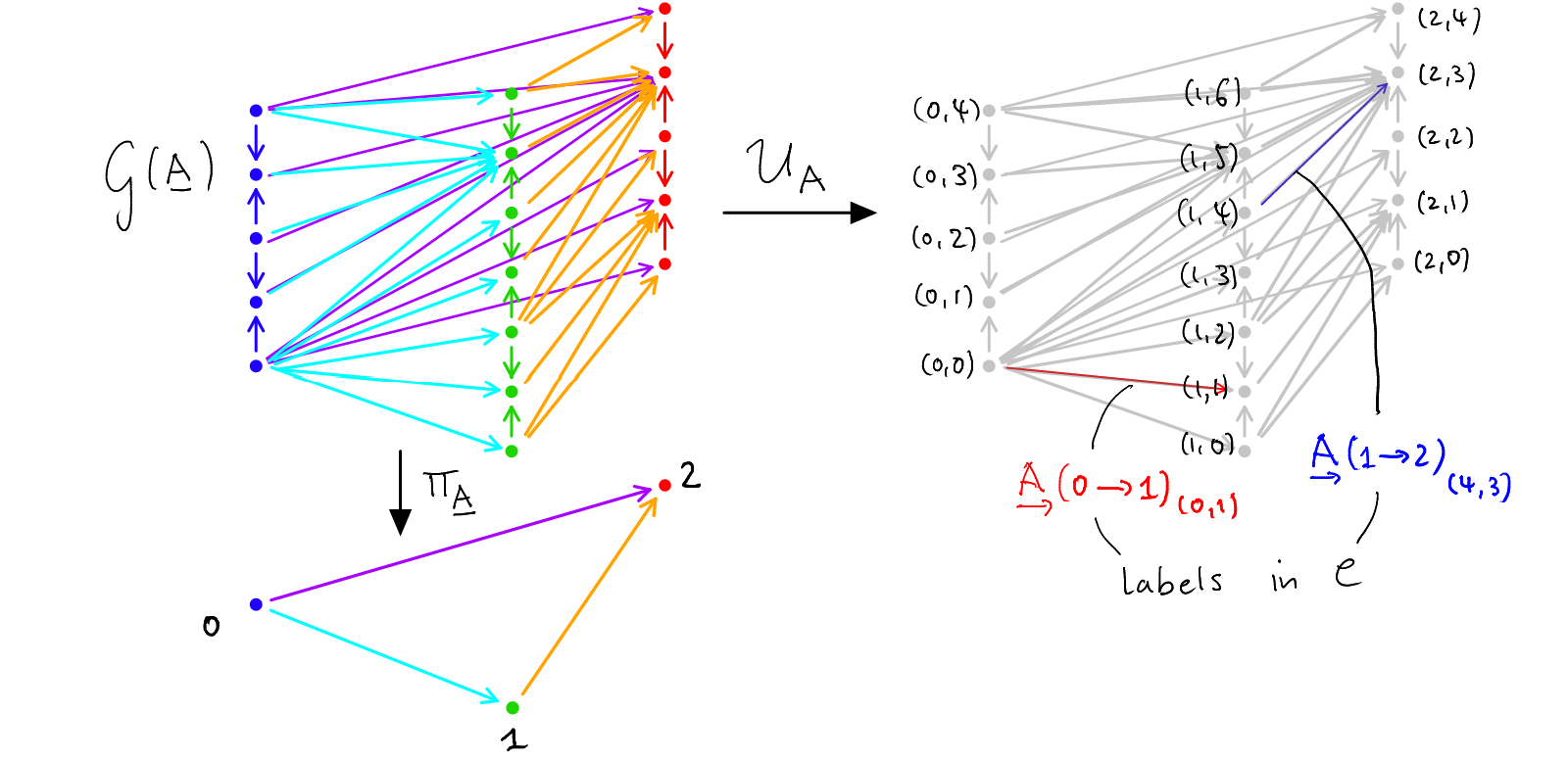}
\endgroup\end{restoretext}
Here we indicated the mapping of $\sU_A$ only in two places: on the arrow $0 \to 1$ over $0 \to 1$ marked in \cred{} where $\sU_A$ has value $\SIs \scA (0 \to 1)_{(0,1)}$, and on the  arrow $4 \to 3$ over $1 \to 2$ marked in \cblue{} where $\sU_A$ has value $\SIs \scA (1 \to 2)_{(4,3)}$.
\end{eg}

This correspondence of $\SIvertone \cC$-families $\scA$ to ``labelling functors" $\sU_\scA$ on total posets $\sG(\und\scA)$ is the content of the following construction.

\begin{constr}[Unpacking and repacking $\SIvertone \cC$-families] \label{defn:unpacking_and_repacking} We define unpacking and repacking operations for labelled singular interval families. Let $X$ be a poset, and $\cC$ be a category. Unpacking takes a functor
\begin{equation}
R : X \to \SIvertone  \cC
\end{equation}
and produces a tuples
\begin{equation}
(\sV_R : X \to \SI, \sU_R : \sG(\sV_R) \to \cC)
\end{equation}
Repacking takes a tuple
\begin{equation}
(V : X \to \SI, U : \sG(V) \to \cC)
\end{equation}
and produces a functor
\begin{equation}
\sR_{V,U} : X \to \SIvertone  \cC
\end{equation}
Formally, the construction does the following.
\begin{enumerate}
\item \textit{Unpacking}: Assume a labelled singular interval family $R : X \to \SIvertone {\cC}$. We define the \textit{unpacking} of $R$ to be the tuple of functors $(\sV_R : X \to \SI, \sU_R : \sG(\sV_R) \to \cC)$ given as follows. We set $\sV_R = \SIf R$. We define $\sU_R$ by $\sU_R(x,a) = \SIs R(x)(a)$ on $(x,a)\in \sG(\sV_R)$, and for $(x,a) \to (y,b)\in \mor(\sG(\sV_R))$ we set 
\begin{equation}\label{eq:chi_G_defn}
\sU_R((x,a) \to (y,b)) := {\SIs R}(x \to y)_{(a,b)} : \SIs R(x)(a) \to \SIs R(y)(b)
\end{equation} 
Note that functoriality of $\sU_R$ follows from functoriality of $R$: Indeed, Given $(x,a) \to (y,b) \to (z,c)$ in $\sG(\sV_R)$ then ${\SIs R}(x \to z) = {\SIs R}(x \to y) \odot {\SIs R} (y \to z)$ implies by \eqref{eq:horizontal_composition} and \eqref{eq:grothendieck_edges} that
\begin{equation}
{\SIs R}(x \to z)_{(a,c)} = {\SIs R}(y \to z)_{(b,c)} \circ {\SIs R}(x \to y)_{(a,b)}
\end{equation}
and thus
\begin{equation}
\sU_R((x,a) \to (z,c)) = \sU_R((y,b) \to (z,c)) \circ \sU_R((x,a) \to (y,b))
\end{equation}
as required.

\item \textit{Repacking}: Let $(V : X \to \SI, U : \sG(V) \to \cC)$ by a tuple of functors. We define the \textit{repacking} of $(V,U)$ to be the functor $\sR_{V,U}: X \to \SIvertone {\cC}$ given as follows. On $x$ in $X$ we set $\sR_{V,U}(x) := (V(x),\rest U x)$, where $\rest U x : V(x) \to \cC$ denotes the restriction of $U$ to the fiber $V(x)$ (cf. \autoref{notn:fiber_restrictions}). Given $(x \to y)$ in $X$, we further define $\sR_{V,U}$ by setting
\begin{equation}
\sR_{V,U}(x \to y) := (V(x\to y), {\SIs {\sR_{V,U}}}(x\to y)) : (V(x), \rest U x ) \to (V(y), \rest U y)
\end{equation}
where $\SIs {\sR_{V,U}} : R(V(x \to y))(-,-) \to \Hom_{\cC}(\rest U x -, \rest U y -)$ is the natural transformation with components
\begin{equation}
\SIs {\sR_{V,U}}(x \to y)_{(a,b)} := U((x,a) \to (y,b))
\end{equation}
Note that naturality of ${\SIs {\sR_{V,U}}}$ as defined above then follows from functoriality of $U$: If $a' \to a$ in $V(x)$ then 
\begin{align*}
{\SIs {\sR_{V,U}}}(x \to y)_{(a,b)}U((x,a') \to (y,a)) &= U((x,a) \to (y,b))U((x,a') \to (y,a)) \\
&= U((x,a') \to (y,b)) \\
&= {\SIs {\sR_{V,U}}}(x \to y)_{(a',b)}
\end{align*}
and similarly one verifies naturality in the second component.
\end{enumerate}
\end{constr}

\begin{lem}[Unpacking and repacking are inverse operations] \label{claim:unpacking_and_repacking} 
The unpacking operation mapping $R \mapsto (\sV_R, \sU_R)$ and the repacking operation mapping $(V,U) \to \sR_{V,U}$ defined in \autoref{defn:unpacking_and_repacking} are mutually inverse to each other.

\proof The proof is very \stfwd{}. 

\qed
\end{lem} 

As an application of the unpacking construction, the next lemma shows that relabelling is acting by post-composition with the labelling functor.

\begin{lem}[Unpacking a relabelling functor] \label{lem:unpacking_coloring_transfer} Let $\scA : X\to \SIvertone \cC$ and $F : \cC \to \cD$. Then 
\begin{equation}
\sU_{\SIvertone F \scA} = F \sU_{\scA}
\end{equation}
\proof Using \autoref{defn:unpacking_and_repacking} and \autoref{defn:transfer_of_coloring} (in particular \eqref{eq:label_transfer_action}) we find
\begin{align}
\sU_{\SIvertone F \scA}((x,a) \to (y,b)) &= \SIs {\SIvertone F \scA} (x \to y) _{(a,b)} \\
&= F (\SIs {\scA} (x \to y) _{(a,b)} ) \\
&= F \sU_{\scA}((x,a) \to (y,b)) 
\end{align}
as required.
\end{lem}

\subsection{Base change}

Just as the Grothendieck construction in \autoref{ch:prerequisites}, the generalised Grothendieck construction introduced in the previous section allows to form new families by pulling back along base change functors $H$. In fact, this is such that labelling functor commute with the total base change $\sG(H)$. We start with the following example.

\begin{eg}[Base change for $\SIvertone \cC$-families] \label{eg:basechange_SIvert} Using $\scA$ as defined in \autoref{eg:unpacking_SIvert} we consider the $\SIvertone \cC$-family $\scA \delta^2_0 : \bnum{2} \to \SIvertone \cC$. The following diagram commutes
\begin{restoretext} 
\begingroup\sbox0{\includegraphics{test/page1.png}}\includegraphics[clip,trim=0 {.0\ht0} 0 {.0\ht0} ,width=\textwidth]{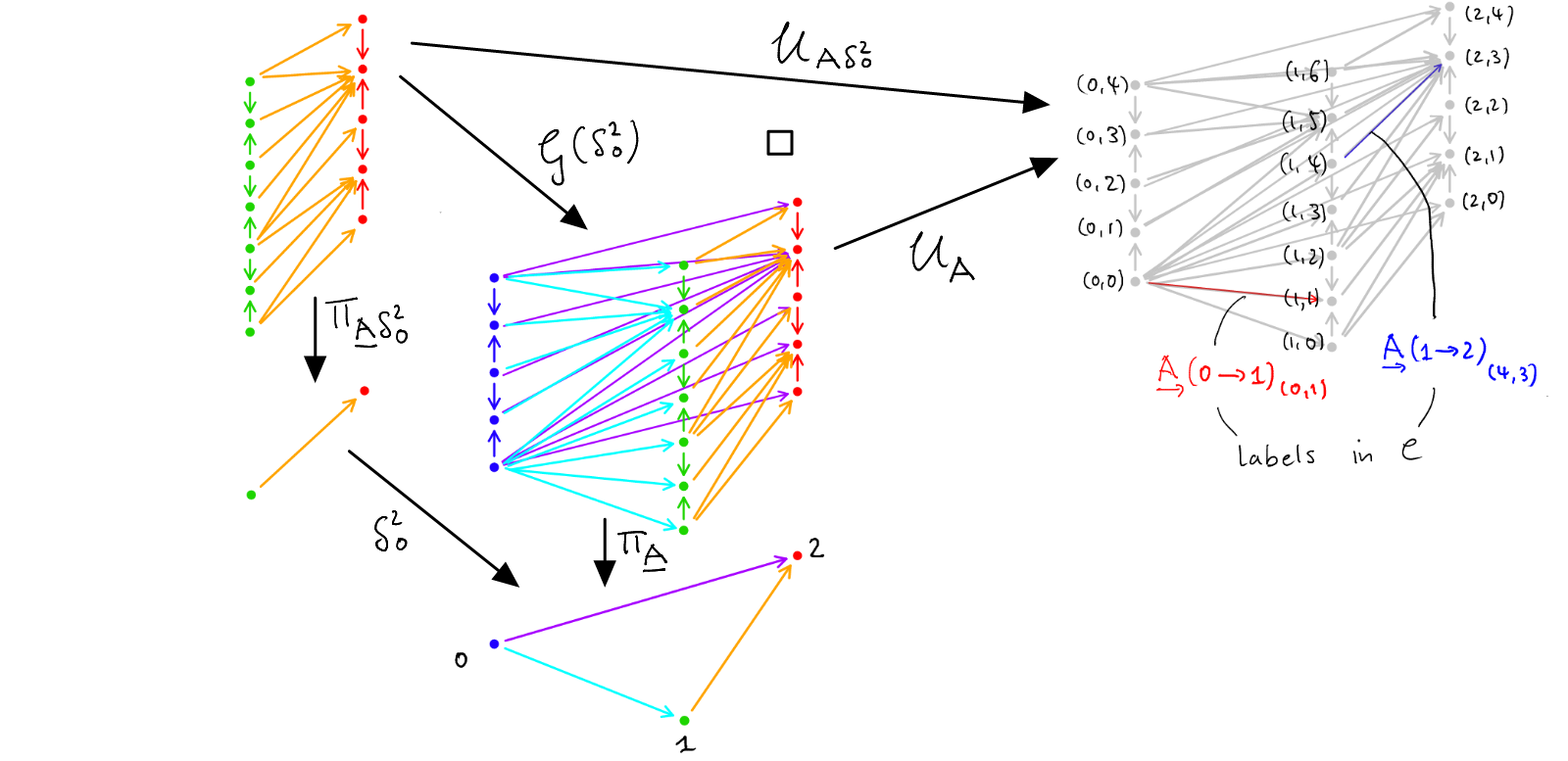}
\endgroup\end{restoretext}
Indeed, the labelling of $\scA \delta^2_0$ is determined by its image morphism
\begin{equation}
\scA \delta^2_0 (0 \to 1) = \scA(1 \to 2)
\end{equation}
and thus the labellings of $\scA\delta^2_0$ agree with those of $\scA$ on the image of $\sG(\delta^2_0)$ (cf. \autoref{defn:grothendieck_base_change}). 
\end{eg}

This observation is more generally recorded in the next Lemma.

\begin{lem}[Base change of $\cC$-labelled $\SI$-families] \label{claim:grothendieck_span_construction_basechange}
Given a family $\scB : Y \to \SIvertone {\cC}$ and $H : X \to Y$. Then,  
\begin{equation}
\sR_{\sV_\scB, \sU_\scB} H = \sR_{\sV_\scB H, \sU_\scB \sG(H)}
\end{equation}
In other words (using \autoref{claim:unpacking_and_repacking}), given $\scA : X \to \SIvertone {\cC}$, then $\scA = \scB H$ if and only if both $\sV_\scA = \sV_\scB H$ and $\sU_\scA = \sU_\scB \sG(H)$.

\proof The proof is \stfwd{}. Set $\scC = \scB H$. Then $\SIf {\scC} = \SIf{\scB H}$ and thus $\sV_\scC = \sV_\scB H$ by \autoref{defn:unpacking_and_repacking}. Recall from \autoref{defn:grothendieck_base_change} that $\sG(H): \sG(\sV_\scC) \to \sG(\sV_\scB)$ is the functor induced by pullback along $H$ (for the factorisation $\sV_\scC = \sV_\scB H$). Reusing notation from \autoref{defn:unpacking_and_repacking} we also have $\SIs \scC(x \to y) = \SIs \scB(H(x) \to H(y))$. We find
\begin{align}
\sU_\scC((x,a) \to (y,b)) &= \SIs \scC(x \to y)_{(a,b)} \\
&= \SIs \scB(H(x) \to H(y))_{(a,b)} \\
&= \sU_\scB((H(x),a) \to (H(y),b)) \\
&= \sU_\scB\big(\sG(H)((x,a) \to (y,b))\big)
\end{align}
where, in the first and third line we used \eqref{eq:chi_G_defn}, in the second our assumption on $\scC$, in the last line we used the definition of $\sG(H)$ in \autoref{defn:grothendieck_base_change}. This establishes $\sU_\scC= \sU_\scB\sG(H)$ as required. \qed

\end{lem}

To summarise the preceding two sections, we have a bijective correspondence
\begin{restoretext}
\begingroup\sbox0{\includegraphics{test/page1.png}}\includegraphics[clip,trim=0 {.35\ht0} 0 {.25\ht0} ,width=\textwidth]{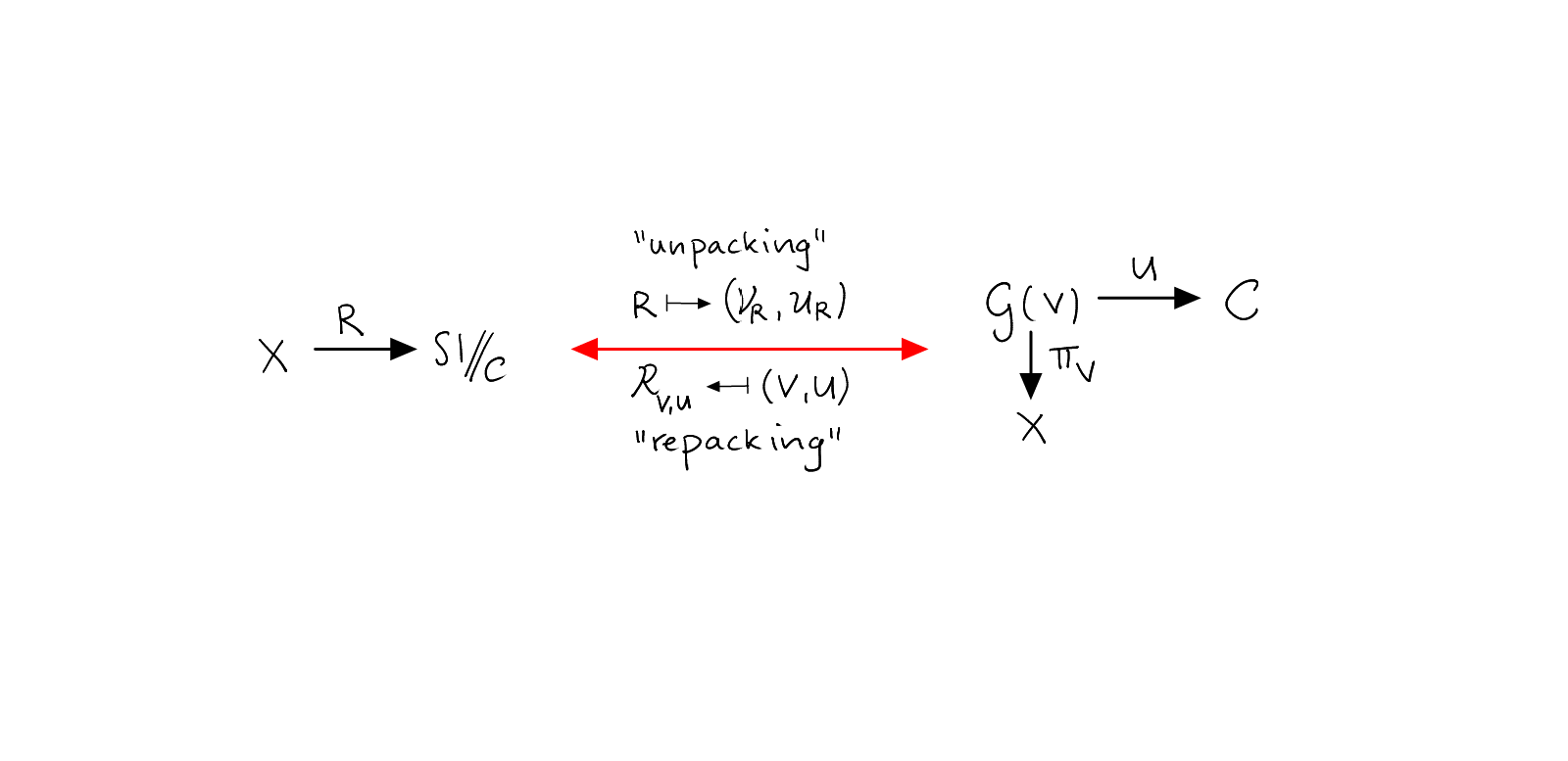}
\endgroup\end{restoretext}
and this correspondence is compatible with pullbacks in the sense that
\begin{restoretext}
\begingroup\sbox0{\includegraphics{test/page1.png}}\includegraphics[clip,trim=0 {.2\ht0} 0 {.25\ht0} ,width=\textwidth]{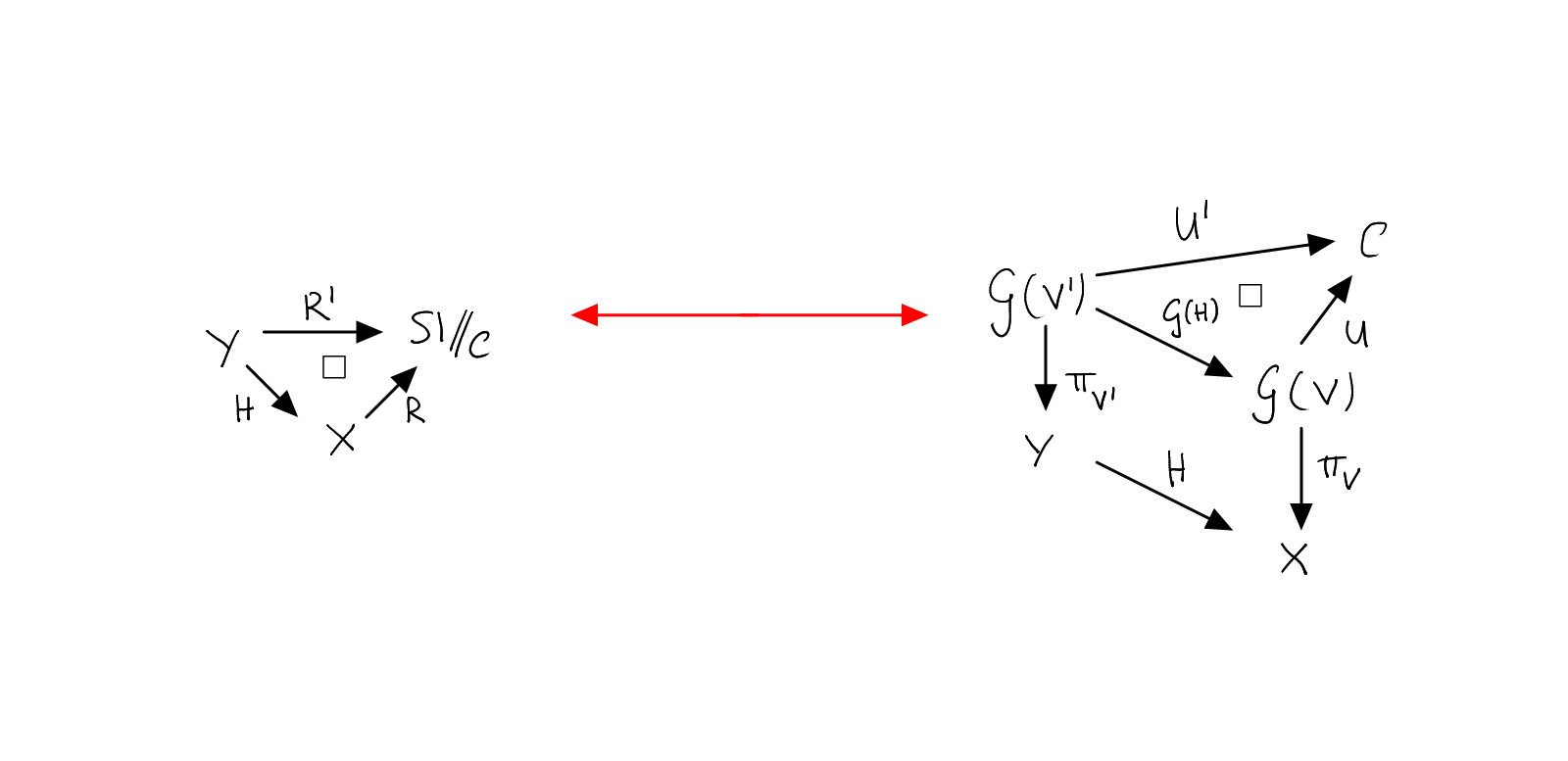}
\endgroup\end{restoretext}
(Here, $R$ corresponds to $(V,U)$ and $R'$ corresponds to $(V',U')$). Further, we showed that the correspondence also satisfies
\begin{restoretext}
\begingroup\sbox0{\includegraphics{test/page1.png}}\includegraphics[clip,trim=0 {.35\ht0} 0 {.35\ht0} ,width=\textwidth]{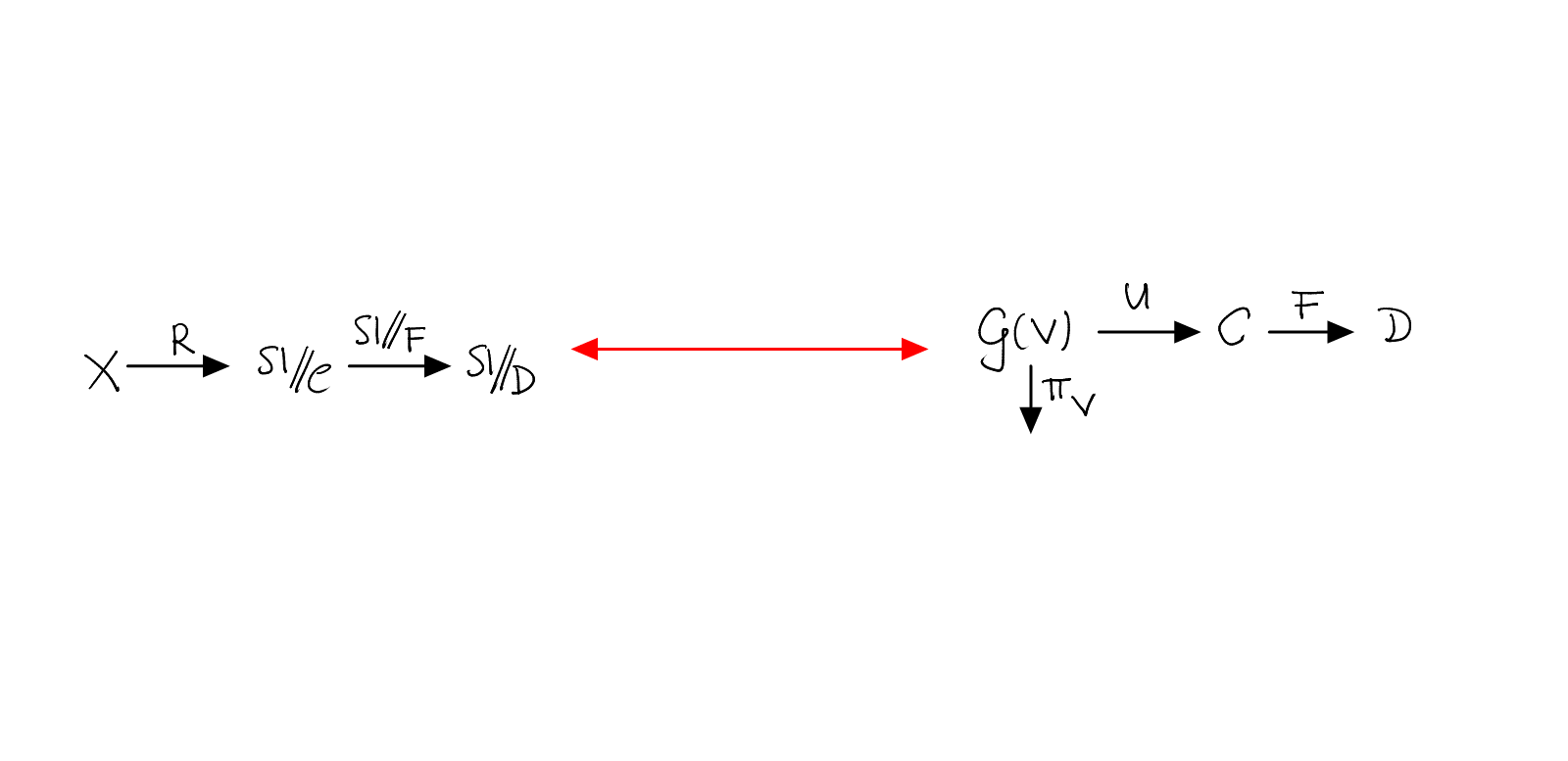}
\endgroup\end{restoretext}
(Here, the left hand side corresponds to $(V,FU)$ and also $R$ corresponds to $(V,U)$).

\section{Singular cube families} \label{sec:cube_fam}

We are now in the position to define (labelled) singular $n$-cubes. We will see that there are two equivalent approaches to the definition: one starts with ``towers" of $\SI$-bundles together with a labelling, and the other inductively uses the $\SIvertone  \cC$ construction.

\subsection{Towers of singular interval families}

We start with the definition of towers of $\SI$-bundles.

\begin{defn}[Towers of $\SI$-families] \label{defn:towers_of_bundles} A \textit{(height $n$) tower $T$ of $\SI$-families} over a poset $X$ is an ordered sequence of \SI-families $\Set{U^{n-1}, U^{n-2}, ... , U^0}$ of length $n$ such that $\dom(U^k) = \cod(\pi_{U^{k-1}})$ ($0 < k < n$) and $\dom(U^0) = X$.
\begin{itemize}
\item $U^k$ is called the \textit{$k$-level labelling} and $\pi_{U^{k}}$ the \textit{$k$-level bundle} of $T$
\item $\sG(U^{k-1})$ is called the \textit{$k$-level total poset} (of $T$), $X$ is called the \textit{base poset} or \textit{$0$-level total poset} of $T$
\item The height $k$ tower $T^k := \Set{U^{k-1}, U^{k-2}, ... , U^0}$ is called the \textit{$k$-level truncation} of $T$
\end{itemize}
\end{defn}

\begin{eg}[Towers of \SI-families] \hfill
\begin{enumerate}
\item As a first example, we construct a tower $T_0 = \Set{U^2,U^1,U^0}$ over $\bnum{1}$ as follows. Set $U^0 : X \to \SI$ to be $\Delta_{\singint 1}$, that is, $X = \bnum{1}$ and $U(0) = \singint 1$. Computing $\sG(U^0)$, we then define $U^1 : \sG(U^0) \to \SI$ by 
\begin{restoretext}
\begingroup\sbox0{\includegraphics{test/page1.png}}\includegraphics[clip,trim=0 {.25\ht0} 0 {.22\ht0} ,width=\textwidth]{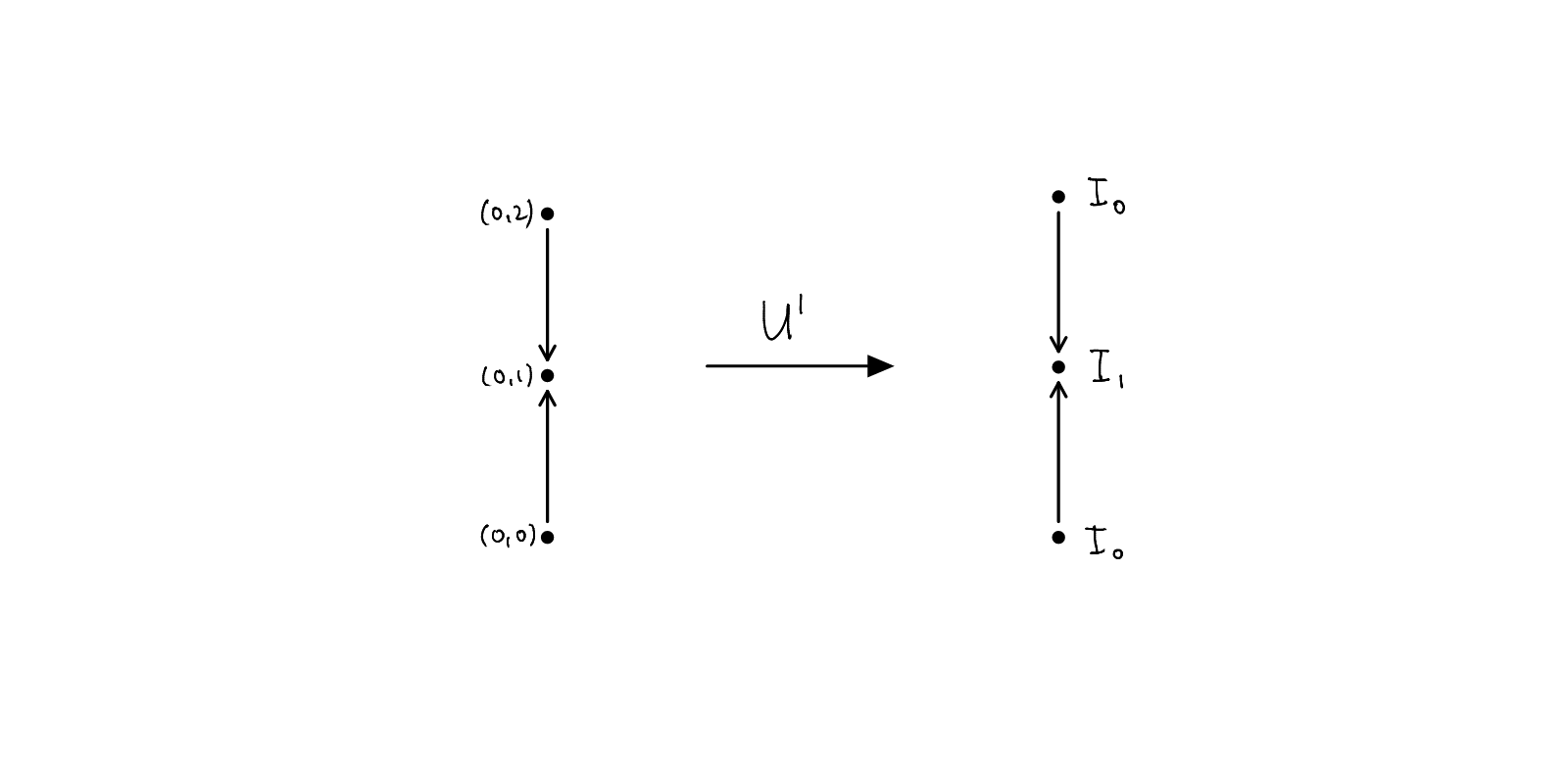}
\endgroup\end{restoretext}
Note that since $\singint 1$ is the terminal singular interval no morphisms need to be specified. Computing the total poset of $U^1$, and setting $f, f' :  \singint 1 \to_{\SI} \singint 2$ to map $f(1) = 3$ and $f'(1) = 1$ respectively, we can then define $U^2 : \sG(U^1) \to \SI$ by
\begin{restoretext}
\begingroup\sbox0{\includegraphics{test/page1.png}}\includegraphics[clip,trim=0 {.1\ht0} 0 {.1\ht0} ,width=\textwidth]{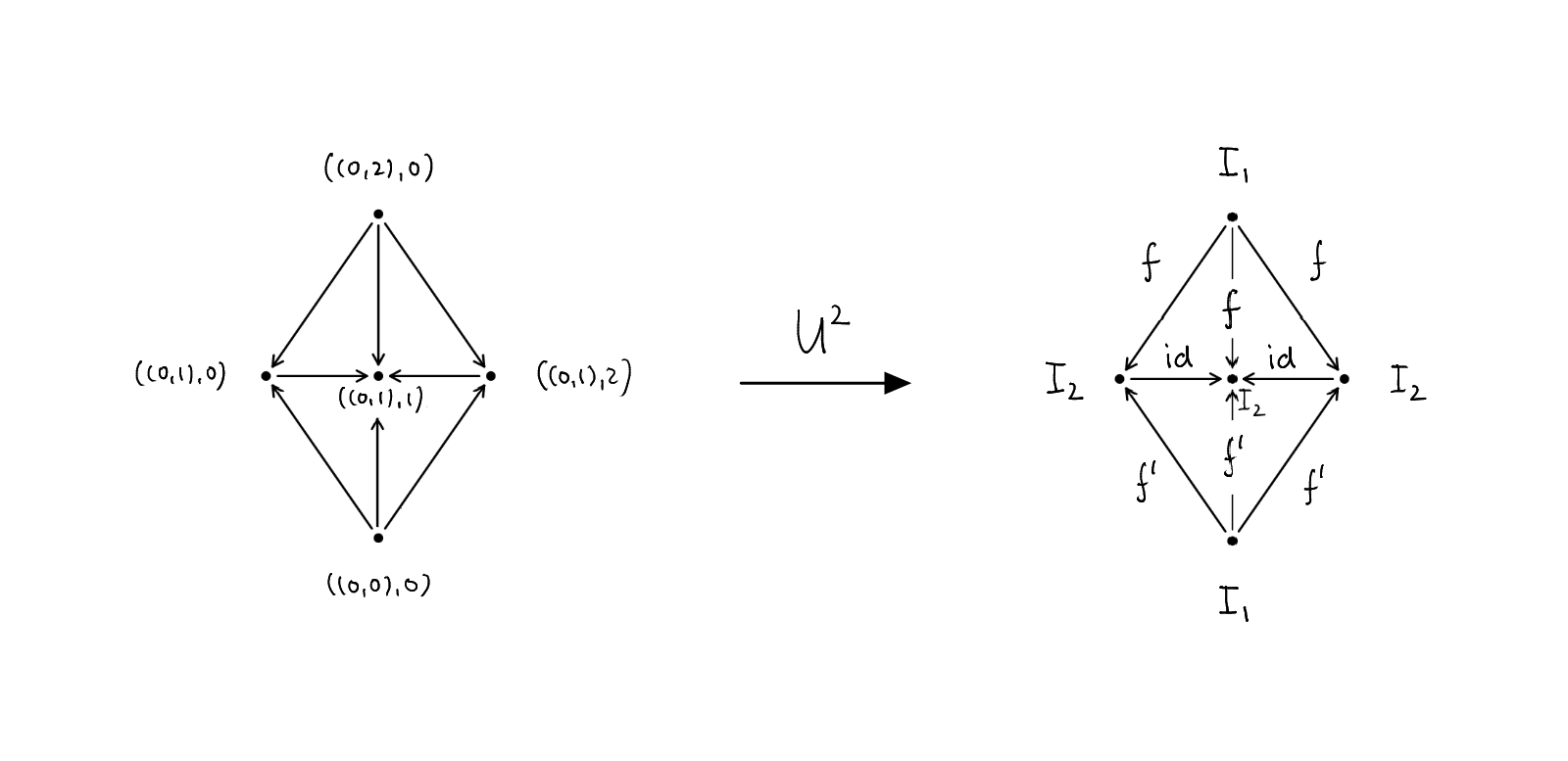}
\endgroup\end{restoretext}
The entire tower $T_0$ can then be visualised by the \SI-bundles $\pi_{U^i}$ (for $i \in \bnum 3$) as a sequence of functors as follows
\begin{restoretext} 
\begingroup\sbox0{\includegraphics{test/page1.png}}\includegraphics[angle=90, origin=c,clip,trim=0 {.0\ht0} 0 {.0\ht0} ,width=7cm]{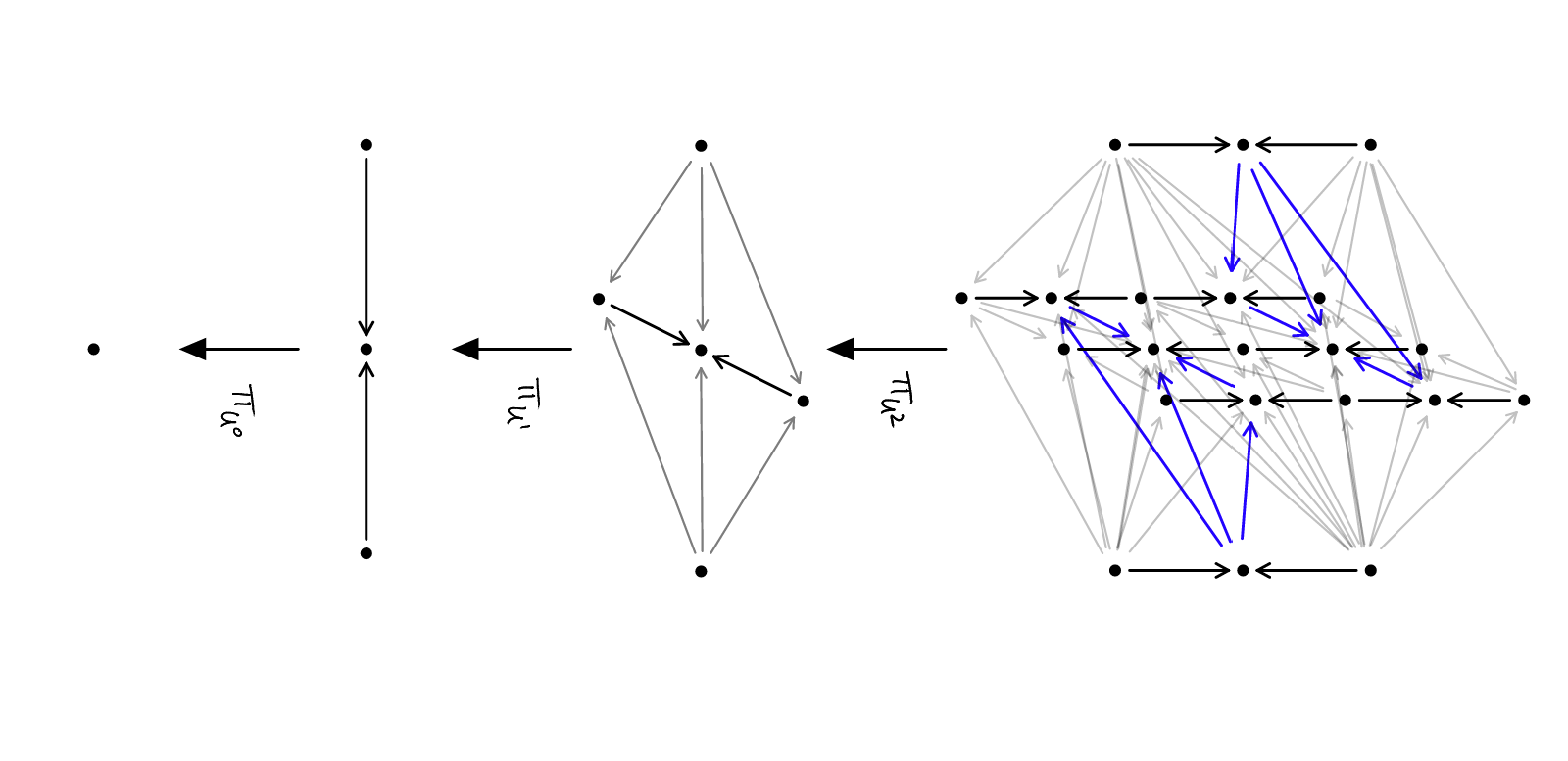}
\endgroup\end{restoretext}
Since a profunctorial realisation $\SiR(f)$ is determined by a function of singular heights (namely, $f$) it is in general not necessary to depict all arrows in $\sG(U^k) = \sG(\SiR U^k)$. Instead, it suffices to depict arrows between singular heights which then determine morphisms in $\SI$. In the above these arrows are marked in \cblue{} (note in particular that $\dom(U^2)$ does not contain any \cblue{} arrows). While we will often adopt this notational simplification, we will however always depict all arrows in the fiber $\pi\inv_{U^k}(x)$ over a single object $x \in \dom(U^k)$. These arrows are marked in \cblack{} above.

\item Depicting only arrows between singular heights as just explained, we define a second tower $T_1 = \Set{U^1, U^0}$ over $\bnum{2}$ by
\begin{restoretext} 
\begingroup\sbox0{\includegraphics{test/page1.png}}\includegraphics[angle=90, origin=c,clip,trim=0 {.0\ht0} 0 {.0\ht0} ,width=7cm]{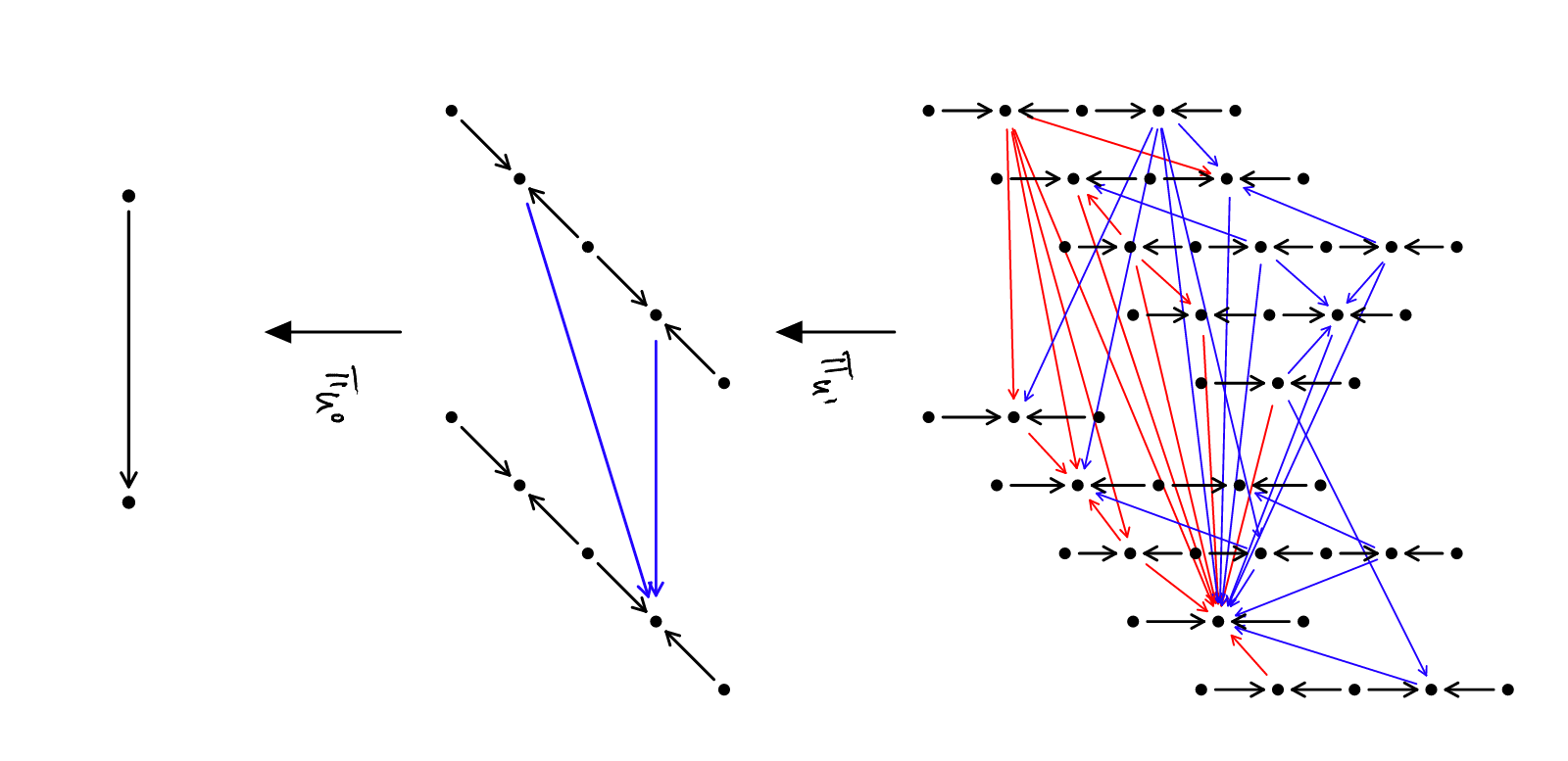}
\endgroup\end{restoretext}
Note that in the top total poset we colored arrows in different colors (in particular distinguishing two different ``layers of arrows" in between fibers) for the sole purpose of giving a clearer picture of the data we are defining.
\end{enumerate}
\end{eg}

\subsection{Category of singular $n$-cubes}
The repacking operation from \autoref{defn:unpacking_and_repacking}, takes a \SI-family $V : X \to \SI$ and labelling functor $U : \sG(V) \to \cC$ and turns it into an $\SIvertone \cC$-family $\sR_{V,U} : X \to \SIvertone \cC$. Using this inductively, we can encode a tower of \SI-families $T = \Set{U^{n-1}, U^{n-2}, \dots, U^1, U^0}$ together with a \textit{labelling functor} $U^n : \sG(U^{n-1}) \to \cC$ in some category $\cC$ as a single functor into a special category, namely, the category of singular cubes. Explicitly, the inductive use of the repacking operation (for $n$ times) takes the following form
\begin{restoretext} 
\begin{noverticalspace}
\begingroup\sbox0{\includegraphics{test/page1.png}}\includegraphics[clip,trim={.25\ht0} {.0\ht0} {.3\ht0} {.15\ht0} ,width=\textwidth]{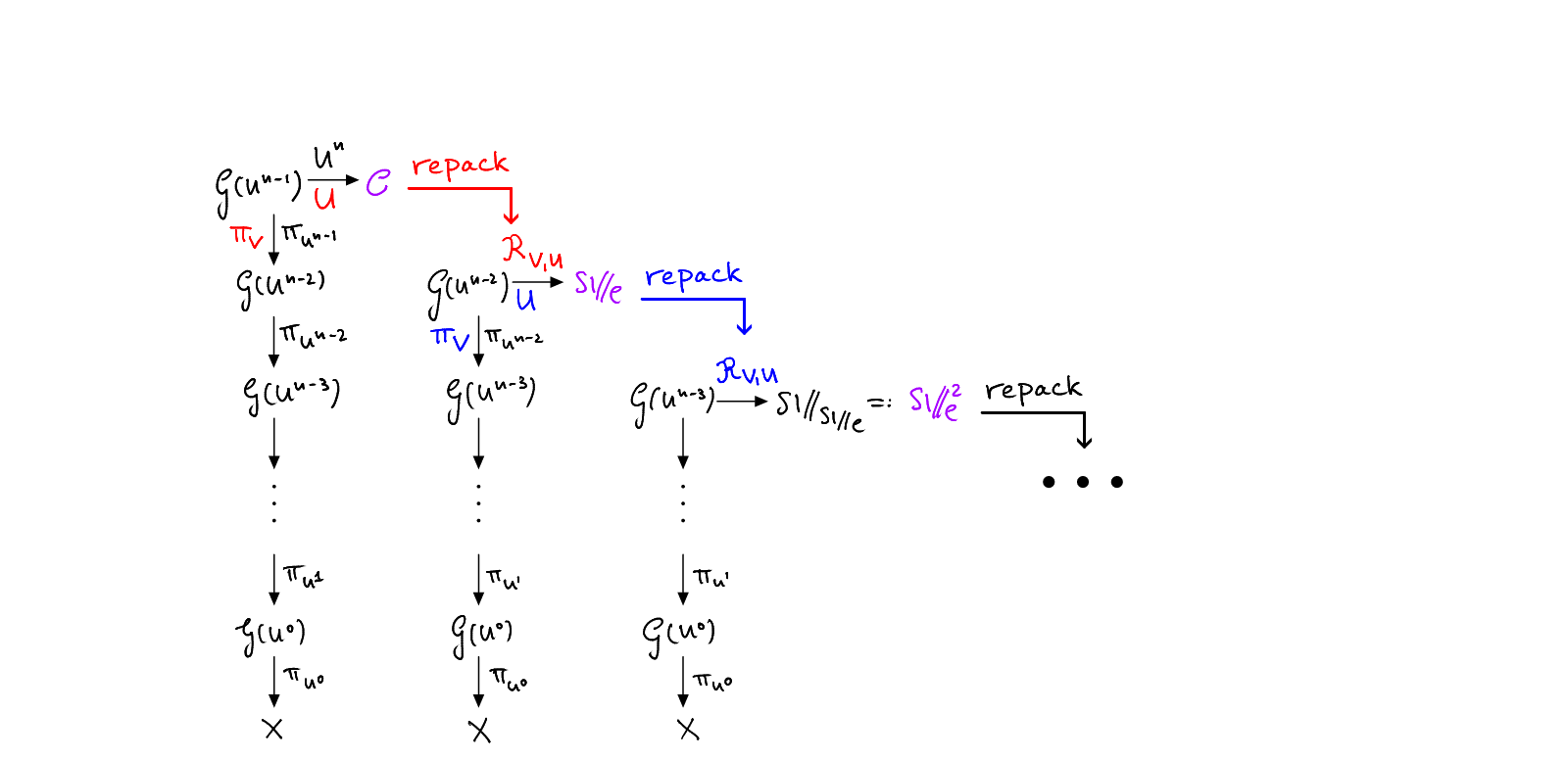}
\endgroup \\
\begingroup\sbox0{\includegraphics{test/page1.png}}\includegraphics[clip,trim={.25\ht0} {.0\ht0} {.3\ht0} {.55\ht0} ,width=\textwidth]{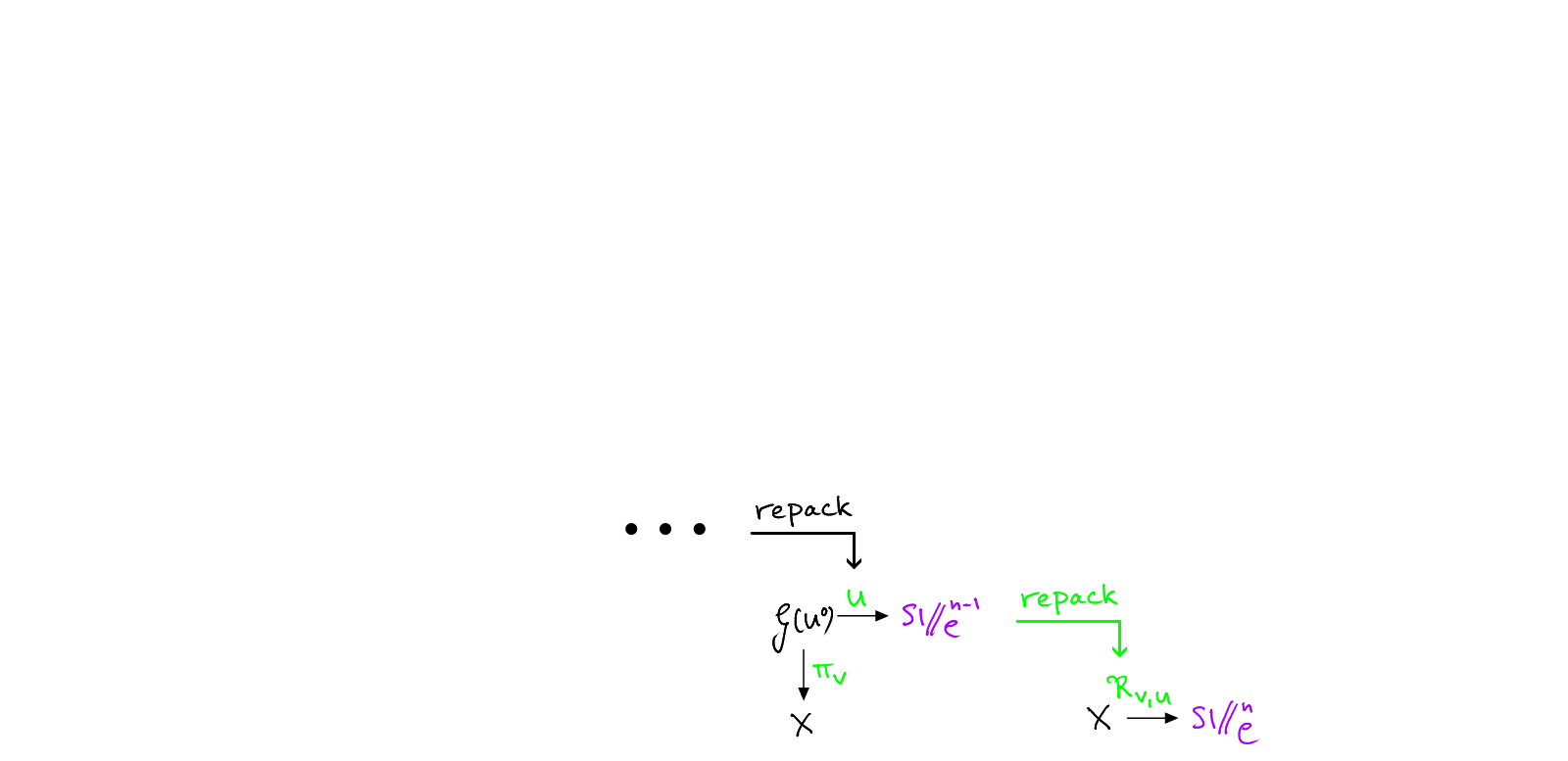}
\endgroup
\end{noverticalspace}
\end{restoretext}
This inductive construction, which will be fledged out in the next section, motivates the following definition.

\begin{defn}[$\cC$-labelled $\SI^n$-families] \label{defn:grothendieck_cubes} Let $\cC$ be a category, $L : \cC \to \cD$ a functor, and $X$ a finite poset. 
\begin{itemize}
\item We inductively define $\SIvert {0} {\cC} = \cC$ and for $n > 0$ we set
\begin{equation}
\SIvert {n} {\cC} = \SIvertone {(\SIvert {n-1}\cC)}
\end{equation}
A functor $\scA : X \to \SIvert n \cC$ is called a \textit{$\cC$-labelled singular $n$-cube family} (or a  \textit{$\SIvert n \cC$-family}) indexed by $X$. %

\item Using \autoref{defn:transfer_of_coloring}, we inductively set $\SIvert {0} {L} = L$ and for $n > 0$
\begin{equation}
\SIvert n L = \SIvertone  {(\SIvert {n-1} L)} : \SIvert n \cC \to \SIvert n \cD
\end{equation}
The functor $\SIvert n L$ is called the \textit{relabelling} induced by $L$ on $n$-cube families.
\end{itemize}
Further, a $\cC$-labelled singular $n$-cube family $\scA : \bnum{1} \to \SIvert n \cC$ indexed by $\bnum{1}$ is also called a \textit{$\cC$-labelled singular $n$-cube} (or a $\SIvert n \cC$-cube).
\end{defn}

\begin{rmk}[Compositionality of relabelling for $\SI^n$-families] \label{rmk:SIvert_endofunctor} Using \autoref{defn:SIvertone_functor} (and \autoref{rmk:transfer_compositional}) we see that the above constructions are values of the functor
\begin{equation}
\SIvert n {-} = \underbrace{\SIvertone {} \circ ... \circ \SIvertone {-}}_{\text{$n$ times}}: \Cat \to \Cat
\end{equation}
In particular, for later use we record the fact, that for $F : \cC \to \cD$, $G : \cD \to \cE$, using for the preceding definition we find that
\begin{gather} \label{rmk:transfer_compositional_for_cubes} 
\SIvert n G \SIvert n F = \SIvert n {GF}
\end{gather}
\end{rmk}

We defer discussing examples of objects in $\SIvert n \cC$ until the next section, which will provide us tools to easily represent such objects.

\subsection{$n$-Unpacking and $n$-repacking of labels}

We will now formally define the inductive repacking operation (and its inverse) which was mentioned in the previous section.

\begin{constr}[$n$-unpacking and $n$-repacking of $\cC$-labelled $\SI^n$-families] \label{defn:complete_unpacking_and_repacking} We define the operations of \textit{$n$-unpacking} and \textit{$n$-repacking} which translate between $\cC$-labelled singular $n$-cube families $\scA : X \to \SIvert n \cC$, and height $n$ towers of singular interval families $\Set{U^{n-1}, U^{n-2}, ... , U^0}$ over $X$ together with a labelling $U^n : \sG(U^{n-1}) \to \cC$.

\begin{enumerate}
\item \textit{$n$-Unpacking}: Let $\scA : X \to \SIvert n \cC$ be a $\cC$-labelled singular $n$-cube family over $X$. We construct its \textit{$n$-unpacking} which consists of the tower of \SI-families $\sT^n_\scA$ (called the  \textit{associated tower} of $\scA$) and a functor $\tsU n_\scA$ (called the \textit{associated labelling} of $\cA$). These are defined based on the following inductive construction of functors
\begin{equation}
\tsU {k}_\scA : \tsG {k}(\scA) \to \SIvert {n-k} \cC
\end{equation}
Set $\tsG 0(\scA) := X$ as well as $\tsU 0_\scA := \scA : \tsG 0(\scA) \to \SIvert n \cC$. For $0 < k \leq n$, note $\tusU {k-1}_\scA : \tsG {k-1}(\scA) \to \SI$. Set 
\begin{equation}
\tsG k(\scA) := \sG(\tusU {k-1}_\scA)
\end{equation}
and define
\begin{equation}
\tpi k_\scA := \pi_{\tusU {k-1}_\scA} : \tsG k(\scA)  \to \tsG {k-1}(\scA)
\end{equation}  
and further (using \autoref{defn:unpacking_and_repacking})
\begin{equation}
\tsU k_\scA := \sU_{\tsU {k-1}_\scA} : \tsG k(A) \to \SIvert {n-k} \cC
\end{equation} 
To complete our construction we now set $\tsU n_\scA : \tsG n(\scA) \to \cC$ to be the associated of $\scA$. And the associated tower of $\scA$ is defined by the sequence
\begin{align} \label{eq:associated_tower}
\sT^n_\scA = \Set{\tusU {n-1}_{\scA}, \tusU {n-2}_{\scA}..., \tusU {0}_{\scA}}
\end{align}

\item \textit{$n$-Repacking}: Conversely, assume a height $n$ tower of \SI-families 
\begin{equation}
T = \Set{U^{n-1}, U^{n-2}, ... , U^0}
\end{equation}
over base space $X$ together with a functor $U^n : \sG(U^{n-1}) \to \cC$. We define its \textit{$n$-repacking} to be a $\cC$-labelled singular $n$-cube family $\tsR n_{T,U^n} : X \to \SIvert n \cC$ which is constructed as follows. Inductively set $\tsR 0_{T,U^n} = L$ and
\begin{equation}
\tsR k_{T,U^n} = \sR_{U^{k-1},\tsR {k-1}_{T,U^n}} : \dom(U^{k-1}) \to \SIvert k \cC
\end{equation}
for $0 < k \leq n$ (using \autoref{defn:unpacking_and_repacking}). $\tsR n_{T,U^n} : X \to \SIvert n \cC$ is then the claimed $n$-repacking of $(T,U^n)$.
\end{enumerate}
Note that by using \autoref{claim:unpacking_and_repacking} inductively we deduce that $n$-unpacking and $n$-repacking are mutually inverse, in the sense that on one hand
\begin{equation}
\Set{\tusU {n-1}_{\tsR n_{T,U^n}}, \tusU {n-2}_{\tsR n_{T,U^n}}..., \tusU {0}_{\tsR n_{T,U^n}}} = T
\end{equation}
and
\begin{equation}
\tsU n_{\tsR n_{T,U^n}} = U^n
\end{equation}
On the other hand,
\begin{equation} \label{eq:complete_unpack_repack_inverse}
\tsR n_{\sT^n_\scA,\tsU n_\scA} = A
\end{equation}
\end{constr}

\begin{rmk}
We also remark that $(\sT^n_\scA)^k = \sT^k_\scA$ by the above construction (cf. \autoref{defn:towers_of_bundles}), where $k \leq n$. The second expression $\sT^k_\scA$ makes sense since every $n$-cube family $\scA : X \to \SIvert n \cC$ is also an $k$-cube family since $\SIvert n \cC = \SIvert k {\SIvert {n-k} \cC}$.
\end{rmk}

\begin{notn} Further to \autoref{defn:towers_of_bundles}, given $A : X \to \SIvert n \cC$ we introduce the following terminology (for $0 \leq k \leq n$)
\begin{itemize}
\item $\tpi k_A$ is the \textit{$k$-level bundle} (of $A$)
\item $\tsU k_A$ is the \textit{$k$-level labelling} (of $A$)
\item $\tsG k(A)$ is the \textit{$k$-level total poset} (of $A$)
\end{itemize}
\end{notn}

\begin{rmk}[Relabelling cubes] \label{rmk:recoloring_cubes} For functors $\scA : X \to \SIvert n \cC$ and $F : \cC \to \cD$, by applying \autoref{lem:unpacking_coloring_transfer} inductively we find that
\begin{equation}
\SIvert n F \scA = \tsR n_{\sT_\scA, F\tsU n_\scA}
\end{equation}
\end{rmk}

\subsection{Examples}

We give two examples of $\SIvert n \cC$-families using the ``$n$-unpacking" terminology developed in the previous section.

\begin{eg}[$\SIvert n \cC$-families] \label{eq:SIvert_n_C_families} \hfill
\begin{enumerate}
\item We define a $\SIvert 3 \cC$-family $\scA$ over $\bnum{1}$ where $\cC$ is a poset as depicted below (concretely, it is the opposite of the category of elements of the $3$-globe as will be defined later on). We define $\scA$ by giving its tower of \SI-families $\sT_\scA$ and its labelling $\tsU n_\scA$ as follows
\begin{restoretext} 
\begin{noverticalspace}
\begingroup\sbox0{\includegraphics{test/page1.png}}\includegraphics[clip,trim=0 {.0\ht0} 0 {.0\ht0} ,width=.8\textwidth]{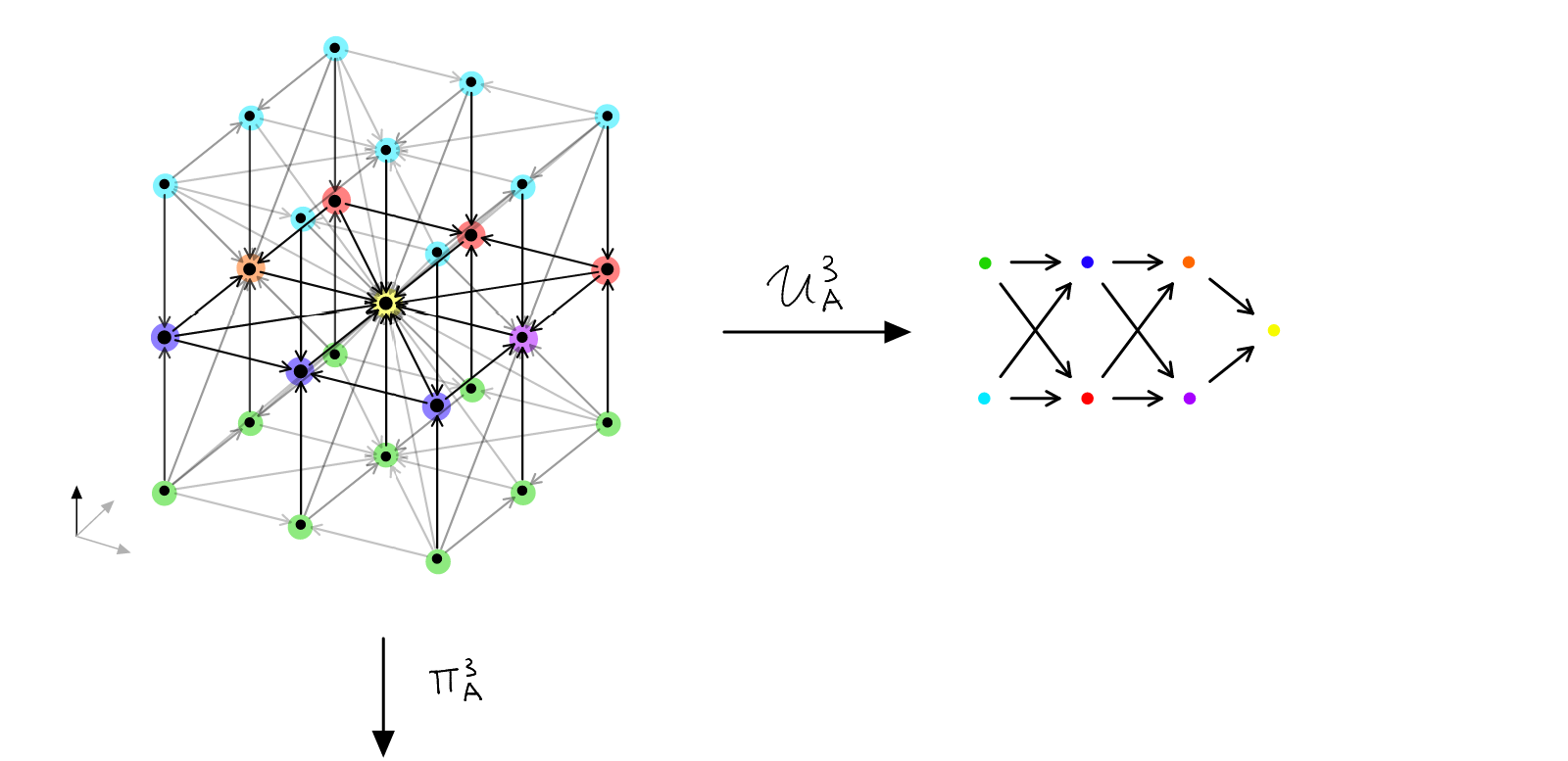}
\endgroup \\*
\begingroup\sbox0{\includegraphics{test/page1.png}}\includegraphics[clip,trim=0 {.0\ht0} 0 {.0\ht0} ,width=.8\textwidth]{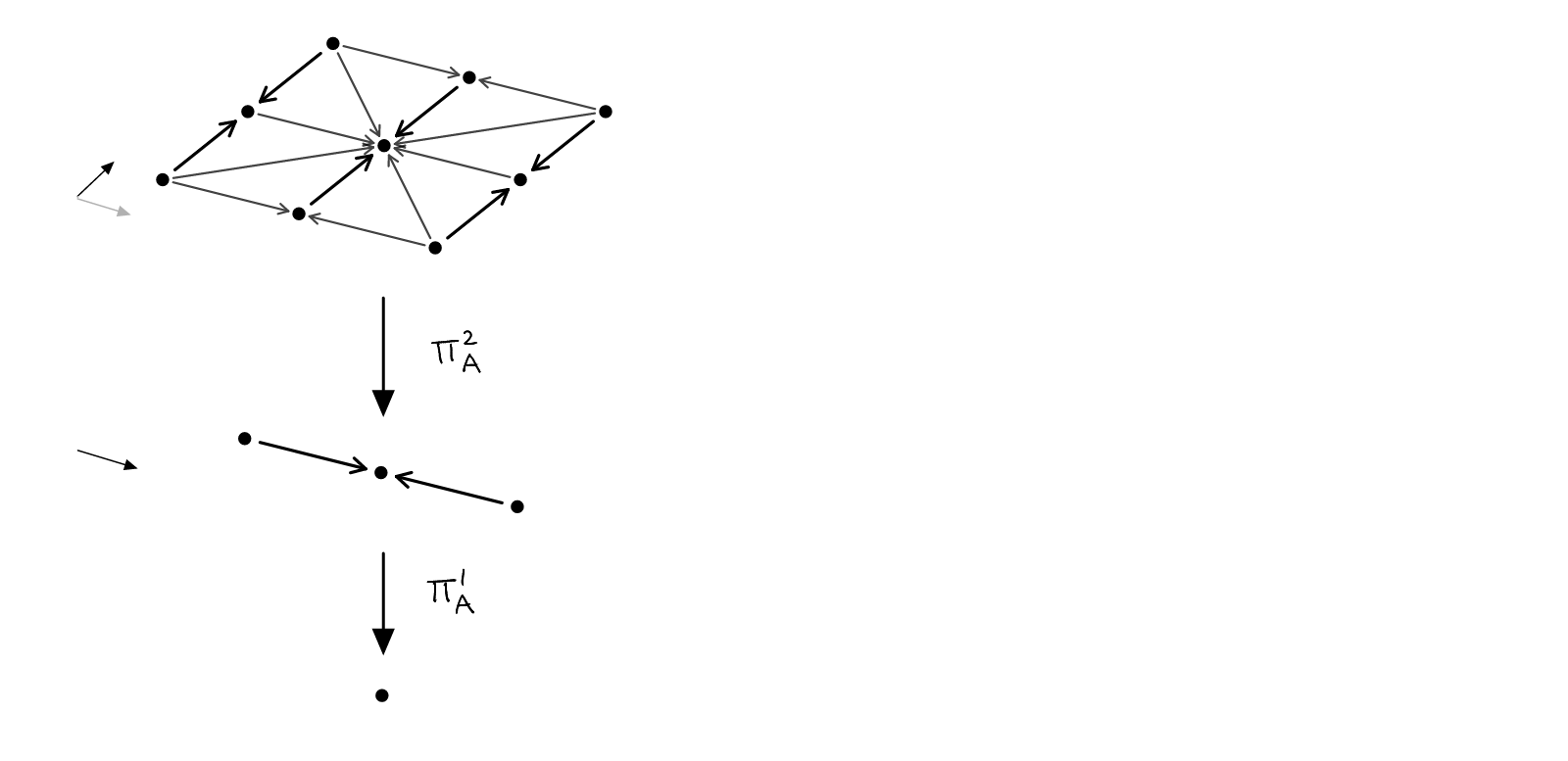}
\endgroup
\end{noverticalspace}
\end{restoretext}
In other words, $\tusU 0_\scA := \const_{\singint 1} : \bnum{1} \to \SI$, $\tusU 1_\scA := \const_{\singint 1} : \sG(U^0) \to \SI$, $\tusU 2_\scA := \const_{\singint 1} : \sG(U^1) \to \SI$ and the functor of posets $\tusU 3_\scA$ is depicted by coloring its preimages.

\item We next define a $\SIvert 2 \cC$-family $\scC : \bnum{3} \to \SIvert 2 \cC$ where $\cC$ is an arbitrary category with some objects $a \in \obj(\cC)$. We define $\scC$ by setting its $2$-unpacked data to be
\begin{restoretext}
\begingroup\sbox0{\includegraphics{test/page1.png}}\includegraphics[clip,trim={.15\ht0} {.0\ht0} {.15\ht0} {.0\ht0} ,width=.8\textwidth]{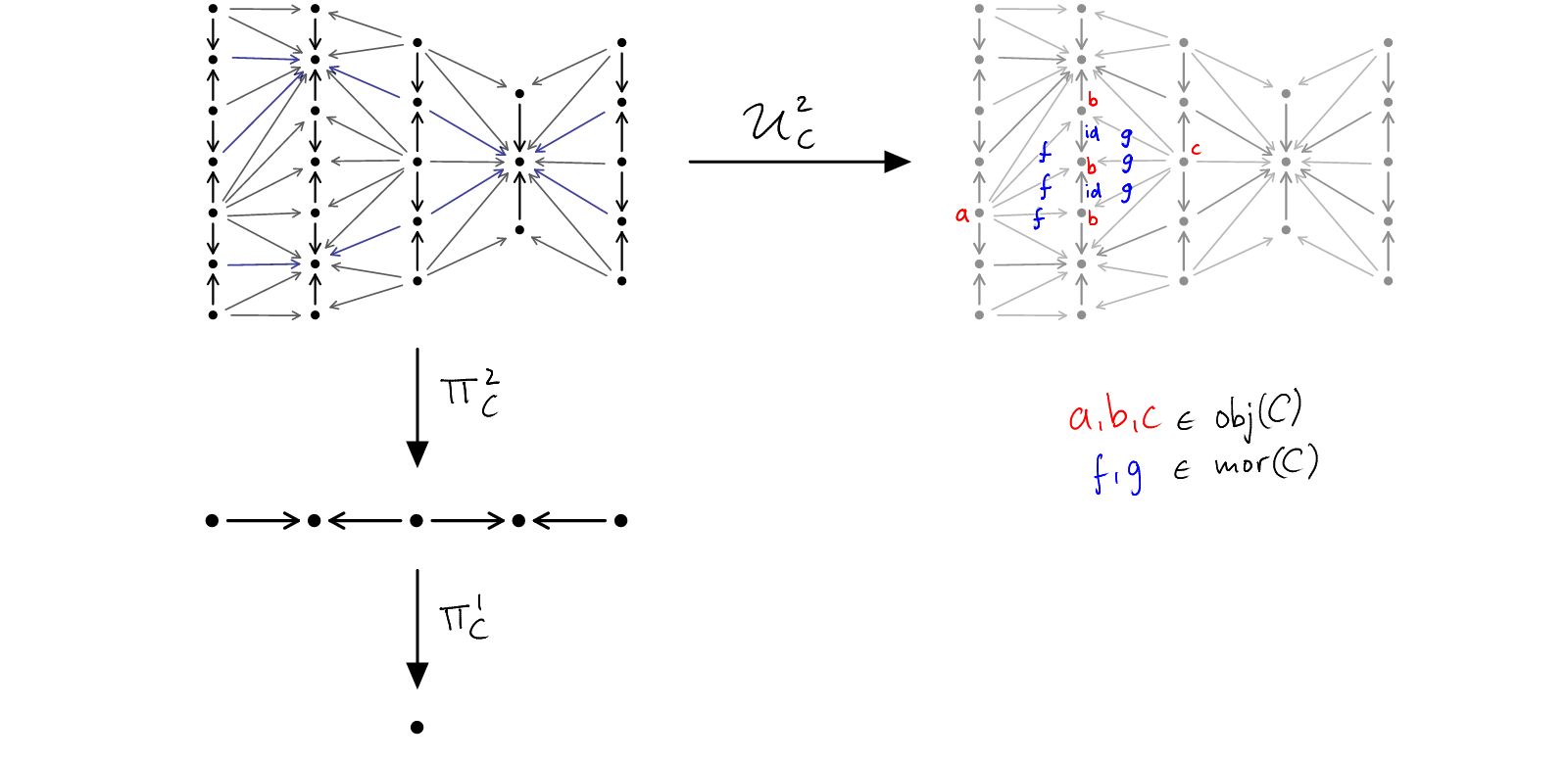}
\endgroup\end{restoretext}
Here we left the labelling of $\scC$ in $\cC$ undefined, apart from 5 objects (which are labelled by $a$, $b$ and $c$) and 8 morphisms (which are labelled by the identity on $f$, $\id_b$ and $g$). The importance of this example lies in the fact that, as we will understand in the coming chapters, it represents a generic situation in which singular heights can be collapsed into regular segments. Concretely, we will introduce an operation of ``collapse" that will act on the above by merging all three copies of $f$, $b$ and $g$ into a single copy (of $f$, $b$ and $g$ respectively).
\end{enumerate}
\end{eg}

We summarise the core terminology, of $k$-level labellings/bundle/total spaces, repacking and unpacking, introduced in this section and the previous section, as follows
\begin{restoretext}
\begingroup\sbox0{\includegraphics{test/page1.png}}\includegraphics[clip,trim={.4\ht0} {.1\ht0} {.4\ht0} {.05\ht0} ,width=0.8\textwidth]{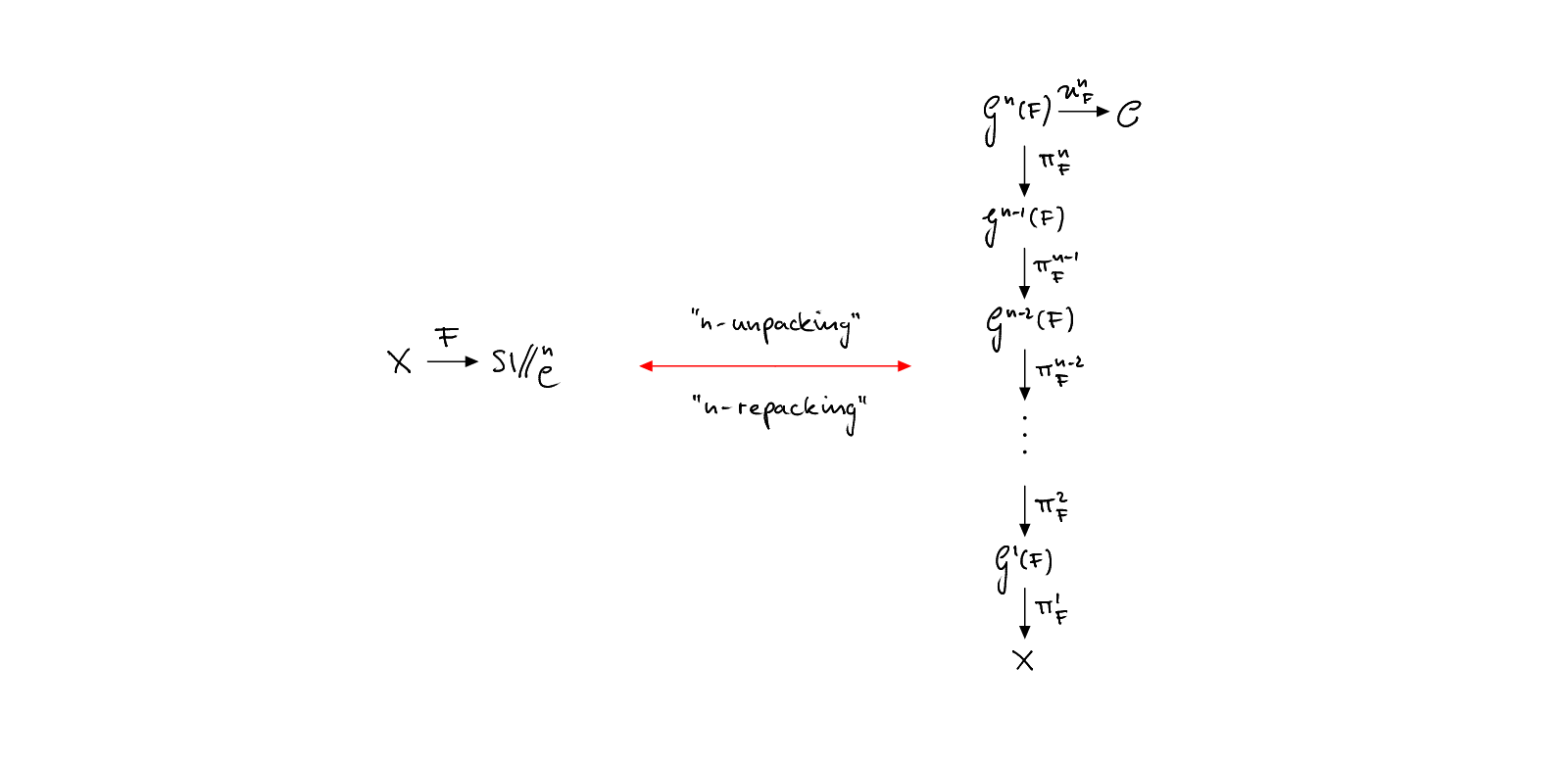}
\endgroup\end{restoretext}

\subsection{Identities}

As a first application of repacking, the next definition gives a notion of identities. As we will later on see, ``geometrically" these are degenerate $(n+1)$-cubes obtained from $n$-cubes.

\begin{defn}[Identities]\label{defn:identities} Let $\scA : X \to \SIvert n \cC$ be a $\cC$-labelled singular $n$-cube family over $\cC$. We define the \textit{identity on $\scA$} $\Id_\scA : X \to \SIvert {n+1} \cC$ by setting\footnote{Note that in this definition we implicitly identify $\sG(\const_{\singint 0}) = X \times\bnum 1 \iso X$, cf. \autoref{notn:product_categories}.}
\begin{equation}
\Id_\scA := \sR_{\const_{\singint 0}, \scA}
\end{equation}
We inductively define $\Id^0_\scA = \scA$ and $\Id^m_\scA = \Id_{\Id^{m-1}_\scA}$.
\end{defn}

\begin{eg}[Identities] Starting from \autoref{eq:SIvert_n_C_families}, we obtain an $\SIvert 1 \cC$-family $\sU^1_\scC$. $\Id_{\sU^1_\scC}$ is then the $\SIvert 2 \cC$-family constructed from the following tower of \SI-families and labelling
\begin{restoretext}
\begingroup\sbox0{\includegraphics{test/page1.png}}\includegraphics[clip,trim={.15\ht0} {.0\ht0} {.15\ht0} {.0\ht0} ,width=\textwidth]{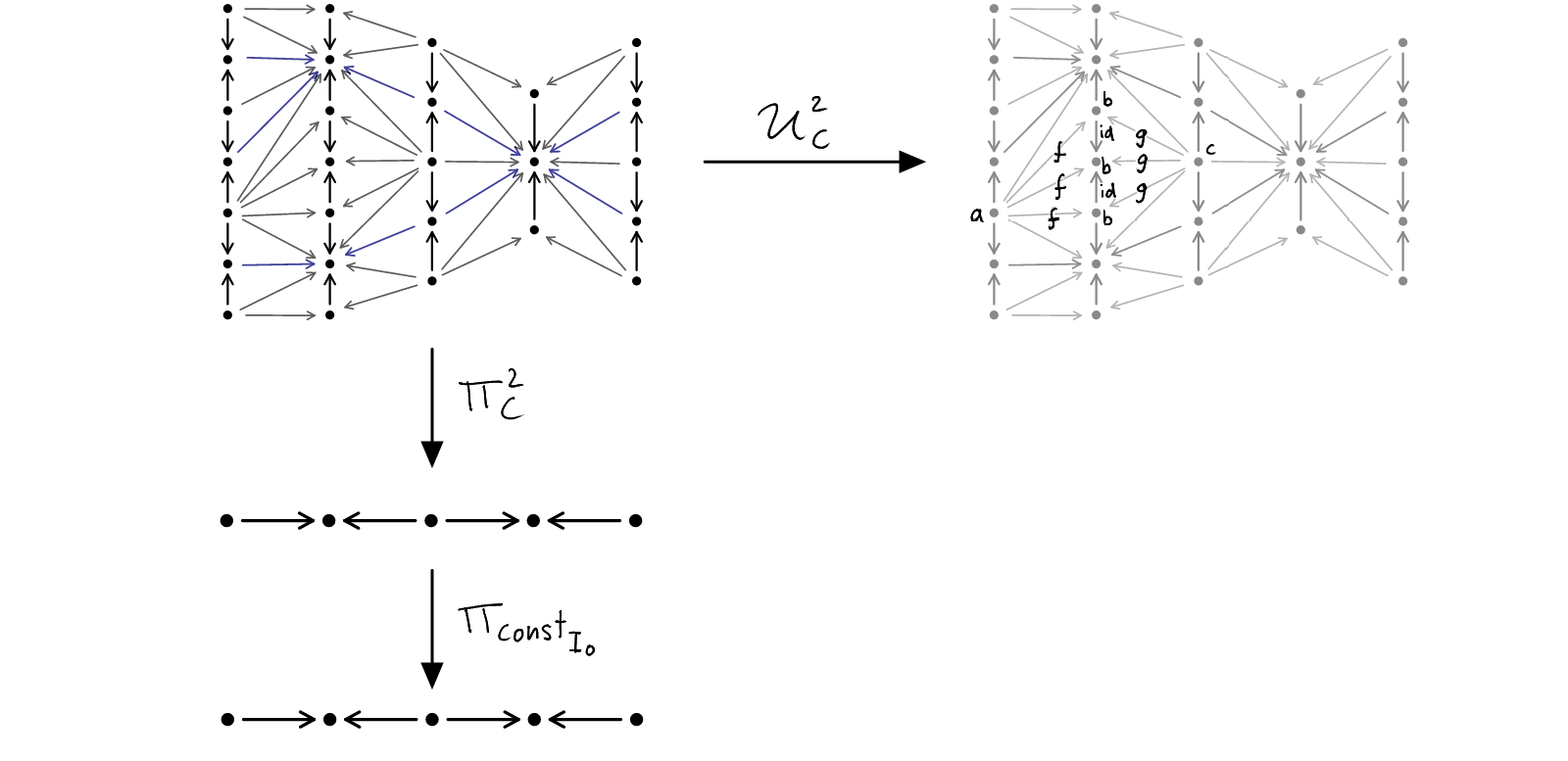}
\endgroup\end{restoretext}
\end{eg}

\subsection{Combinatorial region types} \label{ssec:regions}

Having seen the definition of labelled singular cubes, we will now have a closer look at the structure of the posets involved. In particular, we will give a classification of points of total posets by ``\stratatype " and region dimension. Roughly speaking the dimension gives information about the the maximal dimension of simplices starting at a given point, while the type additionally gives information about the direction of these simplices. The reader might find it helpful to recall \autoref{ssec:coloring}. 

\begin{constr}[Regions, dimension and projection] \label{defn:regions} Let $\scA : X \to \SIvert n \cC$ be a $\cC$-labelled singular $n$-cube family over $\cC$. Objects $p \in \tsG n(\scA)$ will also be called a \textit{region} in $A$. 
\begin{itemize}
\item Regions are classified by a \textit{\stratatype } functor
\begin{equation}
\ctyp^n_\scA : \tsG n(\scA) \to (X \times \bnum{2}^n)
\end{equation}
where $\bnum{2}^n$ denotes the $n$-fold product of the single arrow category $\bnum{2}$. $\ctyp^n_\scA(p)$ is called the \textit{\stratatype } of $p$, and constructed inductively as follows: set $\ctyp^0_\scA : \tsG 0(\scA) \to (X \times \bnum{2}^0)$ to be $\id_X$ (up to the canonical identification $X \iso X \times \bnum{2}^0$). Next, assume to have defined $\ctyp^{m-1}_\scA : \tsG {m-1}(\scA) \to (X \times \bnum{2}^{m-1})$, $0 < m \leq n$. Let $(q,a) \in \tsG m(\scA)$. We define
\begin{equation}
\ctyp^m_\scA(q,a) = (\ctyp^{m-1}_\scA(p), a \mathrm{~mod~} 2) \quad \in \quad \obj\big((X\times \bnum{2}^{m-1}) \times \bnum{2}\big) 
\end{equation}
This extends to a functor of posets
\begin{equation}
\ctyp^m_\scA :  \tsG m(\scA) \to (X \times \bnum{2}^m)
\end{equation}
since if $(q,a) \to (q',a')$ then $a \in \singcont(\tsU {m-1}_\scA(q))$ implies $a' \in \singcont(\tsU {m-1}_\scA(q'))$ by \eqref{eq:defn_order_realisation_1}.
\item Now, given a region $p \in \tsG n(\scA)$, such that $\ctyp^n_\scA(p) = (x,a_1, a_2, ..., a_n)$ for integers $a_i \in \Set{0,1}$ we define the \textit{region co-dimension} of $p$ to be (cf. \autoref{defn:regions})
\begin{equation}
\ctypsum^n_\scA(p) = \sum_{1 \leq i \leq n} a_i
\end{equation}
\item Finally, given a region $p \in \tsG n(\scA)$ we define its \textit{$k$-level projection} $p^k \in \tsG k(\scA)$ inductively by setting $p^n = p$ and $p^k := \tpi {k+1} p^{k+1}$ for $n > k \geq 0$. 
\end{itemize}
\end{constr}

\begin{rmk}[Reverse sequence] Note that the above region types will have the same meaning as the region types $\bnum 2 ^n$ which were already defined in \autoref{ssec:po_mfld_diag} up to reversing the sequence associated to any given type.
\end{rmk}

The following examples illustrates these definitions. It further introduces a way of visualising and coloring posets of $\SIvert n \cC$-cubes. 

\begin{eg}[Regions, projections and coloring] Starting from \autoref{eq:SIvert_n_C_families} \label{eg:regions} consider the following points $p,q,r\in \tsG 2(\scC)$
\begin{restoretext}
\begingroup\sbox0{\includegraphics{test/page1.png}}\includegraphics[clip,trim={.15\ht0} {.0\ht0} {.15\ht0} {.0\ht0} ,width=\textwidth]{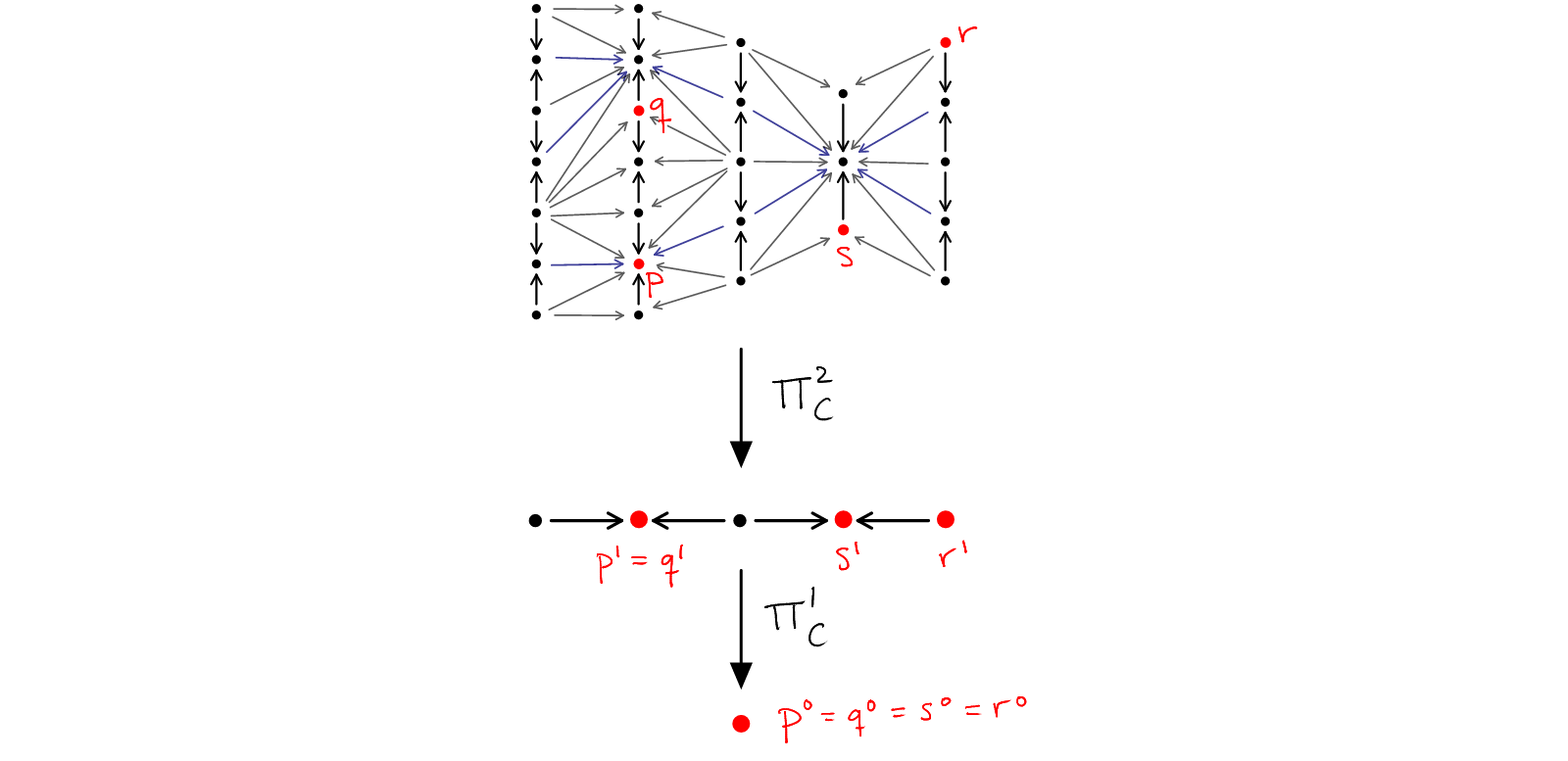}
\endgroup\end{restoretext}
Then we find $\ctyp^2_\scC(p) = (0,1,1) \in \bnum{1} \times \bnum{2}^2$, $\ctyp^2_\scC(q) = (0,1,0)$, $\ctyp^2_\scC(r) = (0,0,0)$ and $\ctyp^2_\scC(s) = (0,1,0)$. Thus $p$ is of region co-dimension $2$, $s$ and $q$ are of region co-dimension $1$ and $r$ is of region co-dimension $0$. We further marked the $k$-level projections for all four points in the tower of \SI-families for $\scC$. In the following illustration we marked all different \stratatype s by colors (\cdarkgreen{}, \cblue{}, \cred{}, \cyellow{}) on the left, and illustrated their ``geometric dimension" (which is $n=2$ minus their region co-dimension) on the right
\begin{restoretext}
\begingroup\sbox0{\includegraphics{test/page1.png}}\includegraphics[clip,trim=0 {.25\ht0} 0 {.25\ht0} ,width=\textwidth]{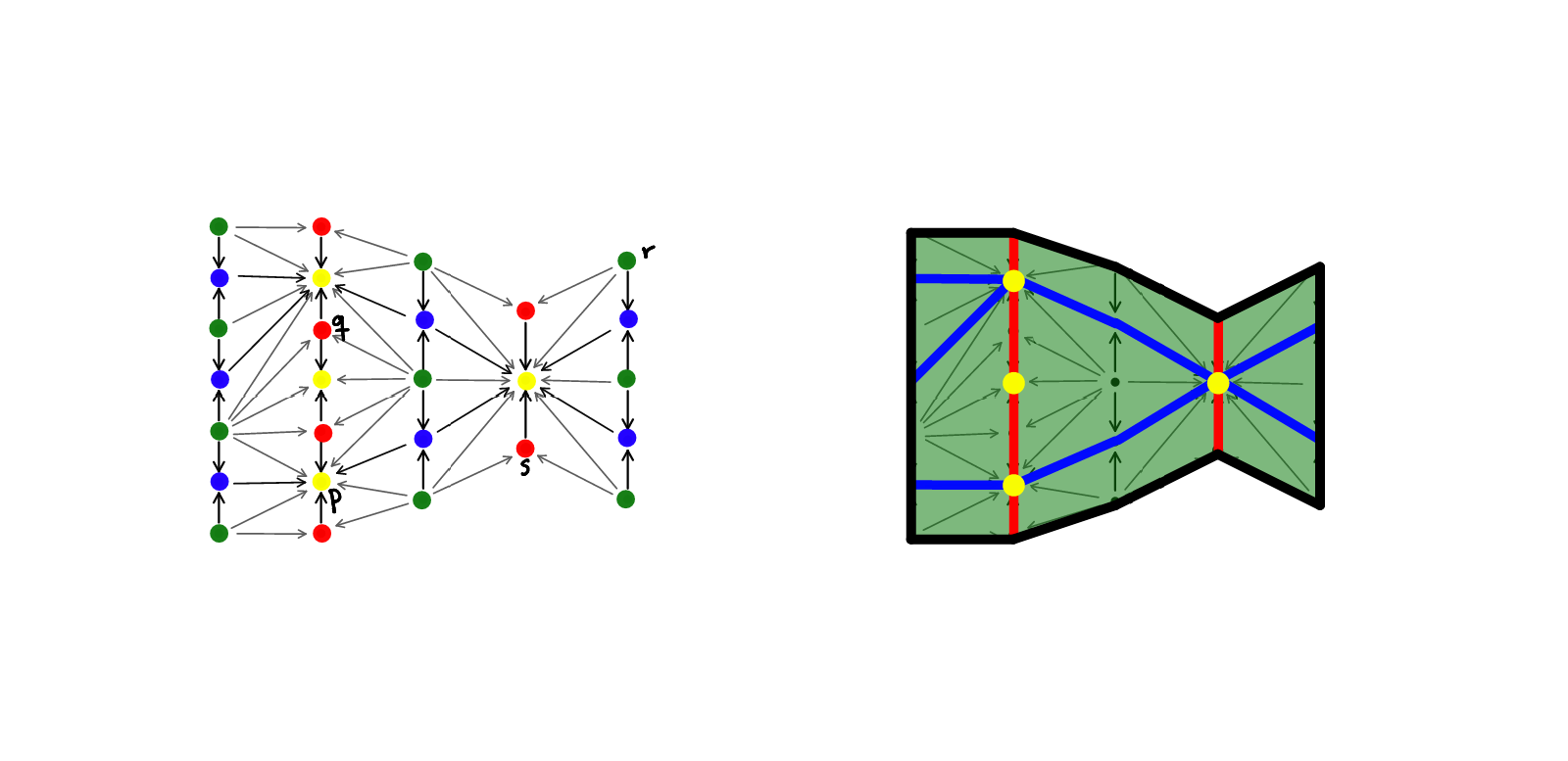}
\endgroup\end{restoretext}
Note that there is a 1-to-1 correspondence between colored points on the left, and colored connected regions on the right. How was the right picture obtained from the left? Any $k$-simplex in $\sG^2(\scC)$ starting at a region $p \in \sG^2(\scC)$ was colored in that region's color from the left (which in turn corresponds to its type): Here, by a \textit{$k$-simplex in a poset $Y$} we mean an injective-on-objects map $f : (\bnum {k+1}) \to Y$. By a $k$-simplex \textit{starting} at $p \in Y$ we mean that $f(0) = p$. Finally by \textit{coloring} a $k$-simplex in $Y$, we mean coloring the interior of the convex closure of its points $f(0), f(1), ..., f(k-1)$ in a visualisation of $Y$.

Further, comparing to \autoref{ssec:coloring}, we now note that 2-simplices \textit{only} start in \cdarkgreen{} regions (which are exactly those of co-dimension $0$), $1$-simplices \textit{also} start in \cblue{} and \cred{} regions (which are of co-dimension $1$) and $0$-simplices start in all regions, in particular in \cyellow{} ones (which are of co-dimension $2$). The resulting connected colored regions are exactly prestrata from \autoref{ssec:coloring}.
\end{eg}

\begin{rmk}[Geometric realisation]
We re-emphasize that the geometric realisation from \autoref{ssec:coloring} (either as strata, or as prestrata like above) is the easiest way to understand the $n$-dimensional geometric content of $\SIvert n \cC$-cubes. This geometric intuition will be helpful for all subsequent definitions.
\end{rmk}

\subsection{$k$-level base change}

We will now introduce the appropriate notion of pullback (or base change) for towers of bundles. The central result of this section is \autoref{cor:basechange_compact}. We start by constructing pullbacks of towers.

\begin{constr}[Base change for towers of \SI-families] \label{constr:basechange_towers} Let 
\begin{equation}
T = \Set{U^{n-1}, U^{n-2}, ... , U^0}
\end{equation}
be a tower of \SI-families, such that $\dom(U^1) = X$. Let $H : Y \to X$ be a functor of posets. The \textit{base change of $T$ along $H$ is a tower of \SI-families} $TH = \Set{V^{n-1}, V^{n-2}, ... , V^0}$ which, together with the $k$-level base change maps $\upi H k : \dom(V^k) \to \dom(U^k)$, is inductively defined  for $k = 0,1, ... (n-1)$ below. Inductively we claim that
\begin{equation} \label{eq:basechange_towers}
V^l = U^l \upi H l
\end{equation}
for all $0 \leq l < n$. For the base case of the induction we set $\tsG 0 (H) = H$ and $V^0 = U^0 H$. For $k > 0$ we use \autoref{defn:grothendieck_base_change} on the factorisation \eqref{eq:basechange_towers} for $l = (k-1)$. This allows us to set
\begin{equation}
\upi H k := \sG(\upi H {k-1})
\end{equation}
and further define
\begin{equation}
V^k := U^k \upi H k
\end{equation}
This satisfies our inductive assumption \eqref{eq:basechange_towers} for $l = k$, and thus completes the inductive construction. Note that $T$ is implicit in the notation for $\tsG k (H)$.
\end{constr}

The preceding construction can be visualised as follows
\begin{restoretext}
\begingroup\sbox0{\includegraphics{test/page1.png}}\includegraphics[clip,trim={.43\ht0} {.15\ht0} {.57\ht0} {.1\ht0} ,width=.8\textwidth]{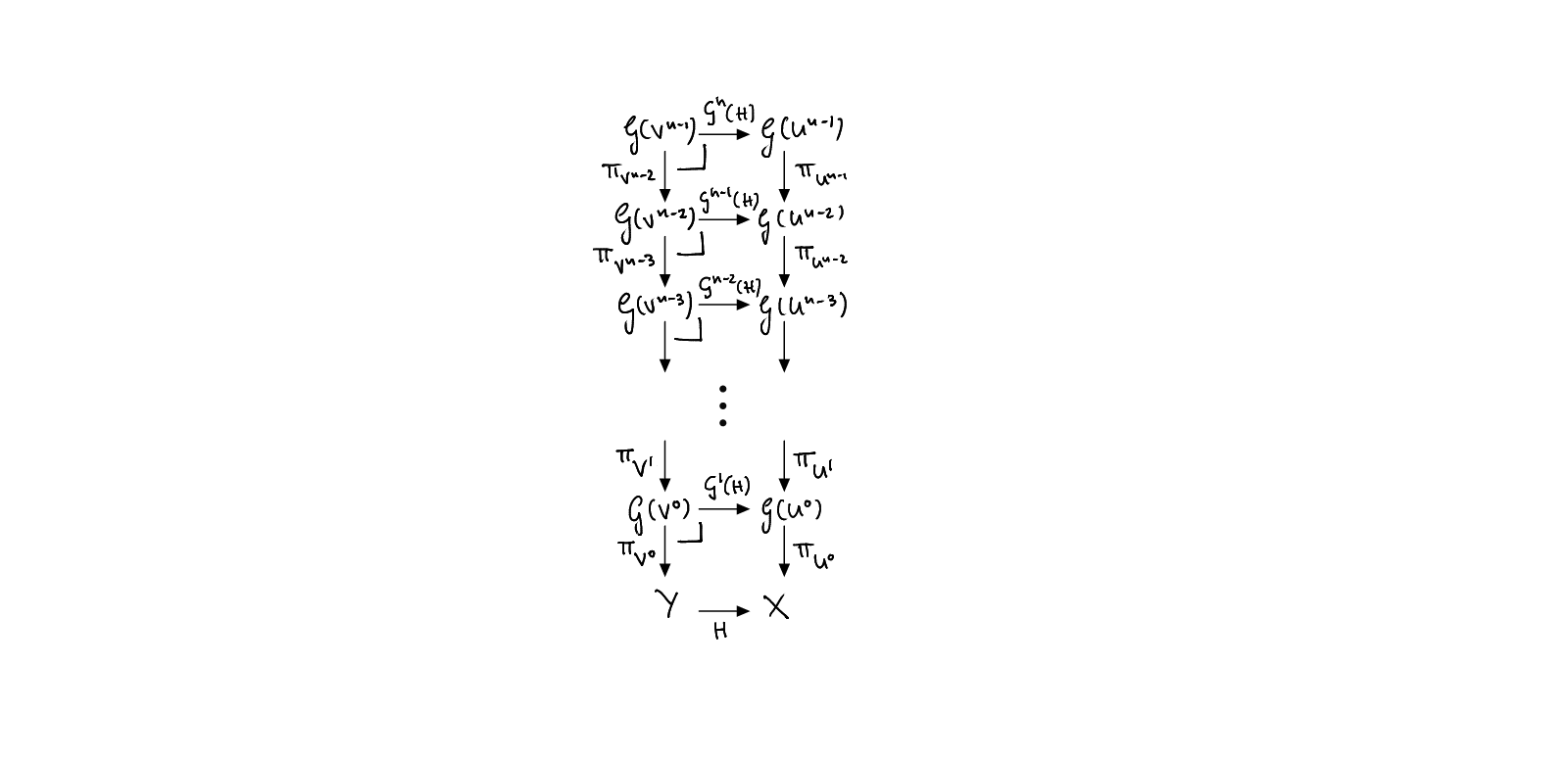}
\endgroup\end{restoretext}
We will further often work with two towers $T_1$, $T_2$ of \SI-families, that coincide up to level $k$, that is, $T^k_1 = T^k_2$. This situation will often be depicted as
\begin{restoretext}
\begingroup\sbox0{\includegraphics{test/page1.png}}\includegraphics[clip,trim={.4\ht0} {.0\ht0} {.6\ht0} {.0\ht0} ,width=.8\textwidth]{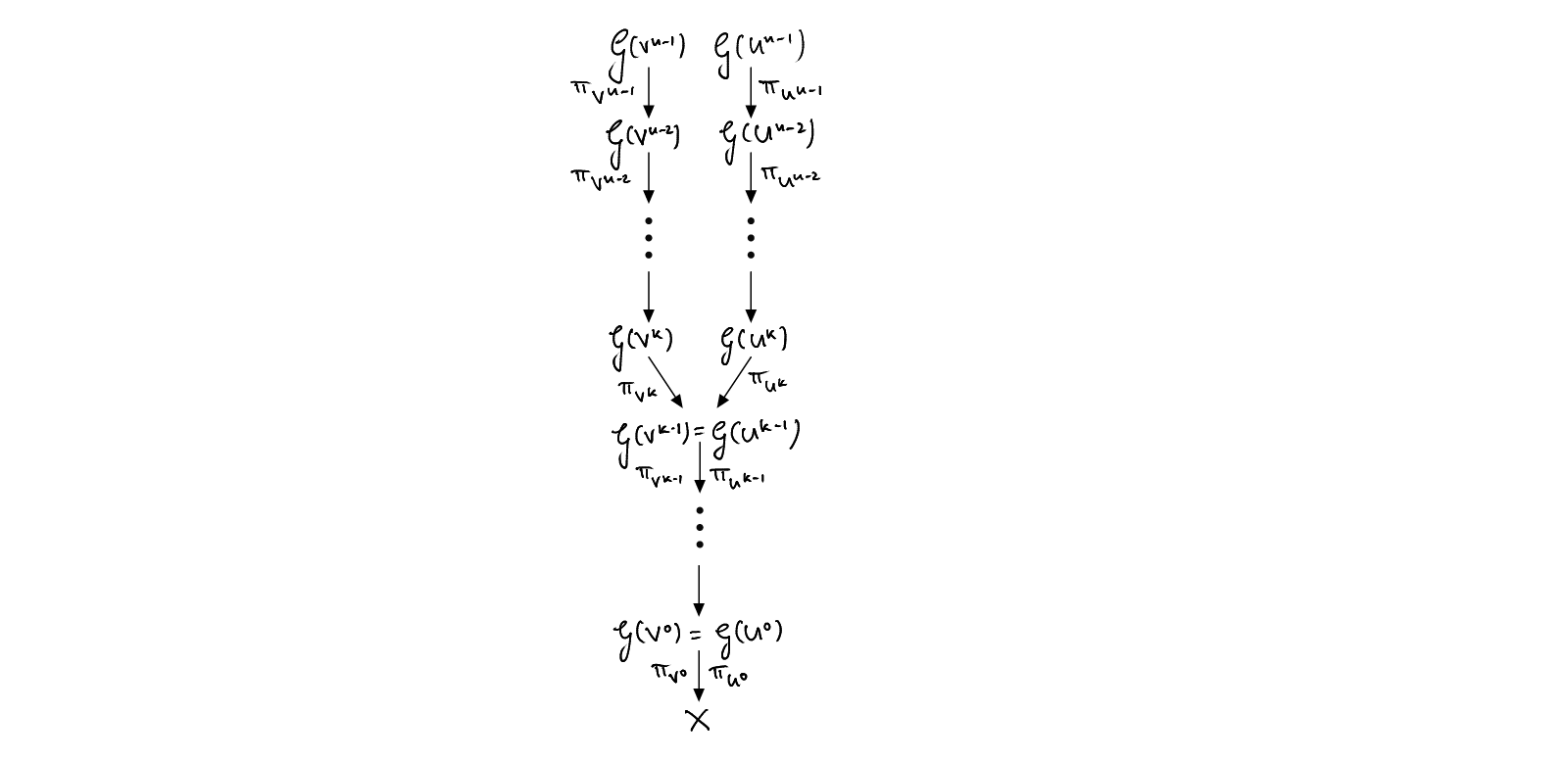}
\endgroup\end{restoretext}
In such a situation we can ask under what conditions the towers above level $k$ are obtained by a base change; that is, we are asking for a condition such that
\begin{restoretext}
\begingroup\sbox0{\includegraphics{test/page1.png}}\includegraphics[clip,trim={.38\ht0} {.0\ht0} {.62\ht0} {.0\ht0},width=.8\textwidth]{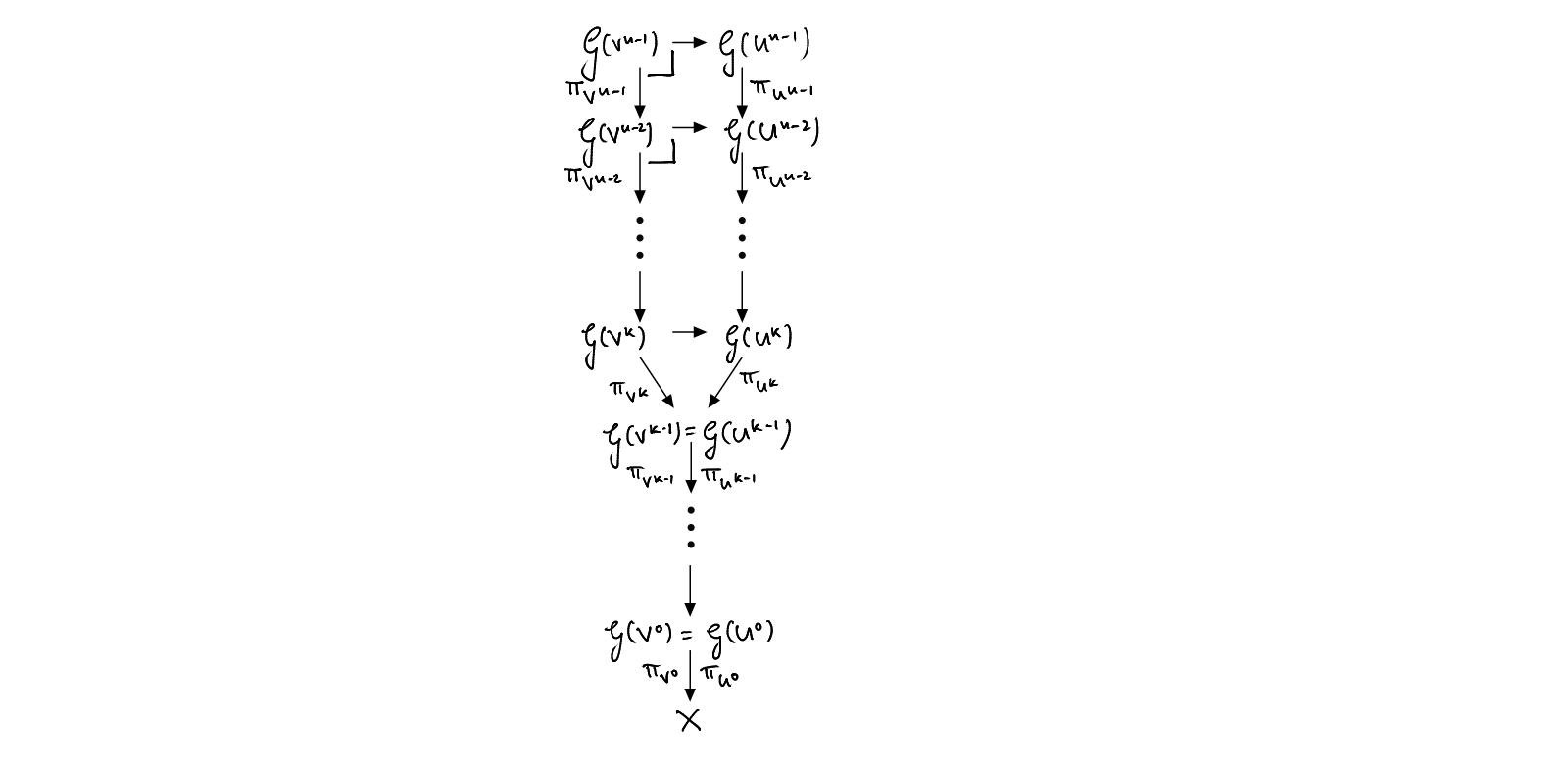}
\endgroup\end{restoretext}
In the case of \textit{labelled} $n$-cubes, a sufficient condition is given in \autoref{constr:unpacking_collapse} below, which in \autoref{rmk:unpacking_collapse} and \autoref{lem:unpacking_collapse} we will then see to be a necessary condition as well.

\begin{defn}[$k$-level base change] \label{defn:klvl_base_change} Let $\scA : X \to \SIvert n \cC$, $\scB : Y \to \SIvert n \cC$. Assume that $\sT^{k-1}_\scA =\sT^{k-1}_\scB$ for some $0 \leq k \leq n$. Let $H : \tpi k_\scA \to \tpi k_\scB$ be a map of bundles such that $\tsU {k}_{\scA} = \tsU {k}_{\scB} H$. In this case we say that $\scA$ is obtained as the \textit{$k$-level base change} of $\scB$ along $H$.
\end{defn}

\begin{claim} [$k$-level base change for labelled cube families] \label{constr:unpacking_collapse} Let $\scA : X \to \SIvert n \cC$, $\scB : Y \to \SIvert n \cC$ and assume $\scA$ is obtained as the $k$-level base change of $\scB$ along $H$ as just defined. Then
\begin{equation}
\sT^{n-k}_{\tsU {k}_{\scA}} = \sT^{n-k}_{\tsU {k}_{\scB}} H
\end{equation}
and for each $i$ with $k \leq k + i \leq n$ we have 
\begin{equation} \label{eq:complete_unpacking_ind}
\tsU {k+i}_{\scA} = \tsU {k+i}_{\scB} \upi H i
\end{equation}
\proof The proof is \stfwd{} and visualised below. The base case $(i = 0)$ of \eqref{eq:complete_unpacking_ind} follows by assumption in the claim since $\upi H 0 = H$. Now, for $i > 0$ we find by 
\autoref{claim:grothendieck_span_construction_basechange} applied to the factorisation 
\begin{equation}
\tsU {k+i-1}_{\scA} = \tsU {k+i-1}_{\scB} \upi H i : \tsG {k+i-1}(\scA) \to \SIvert {n-k-i+1} \cC = \SIvertone  {\SIvert {n-k-i} \cC}
\end{equation}
that
\begin{equation}
\tsU {k+i}_{\scA} = \tsU {k+i}_{\scB} \sG(\upi H {i-1})
\end{equation}
which completes the inductive proof of \eqref{eq:complete_unpacking_ind}. From the latter together with the definitions in \autoref{defn:complete_unpacking_and_repacking} and \autoref{constr:basechange_towers} we then deduce that 
\begin{equation}
\sT^{n-k}_{\tsU {k}_{\scA}} = \sT^{n-k}_{\tsU {k}_{\scB}} H
\end{equation}
as claimed. \qed
\end{claim}

The inductive unpacking procedure used in the proof of the previous claim can be visualised as follows
\begin{restoretext} 
\begin{noverticalspace}
\begingroup\sbox0{\includegraphics{test/page1.png}}\includegraphics[clip,trim={.3\ht0} {.0\ht0} {.4\ht0} {.15\ht0} ,width=.8\textwidth]{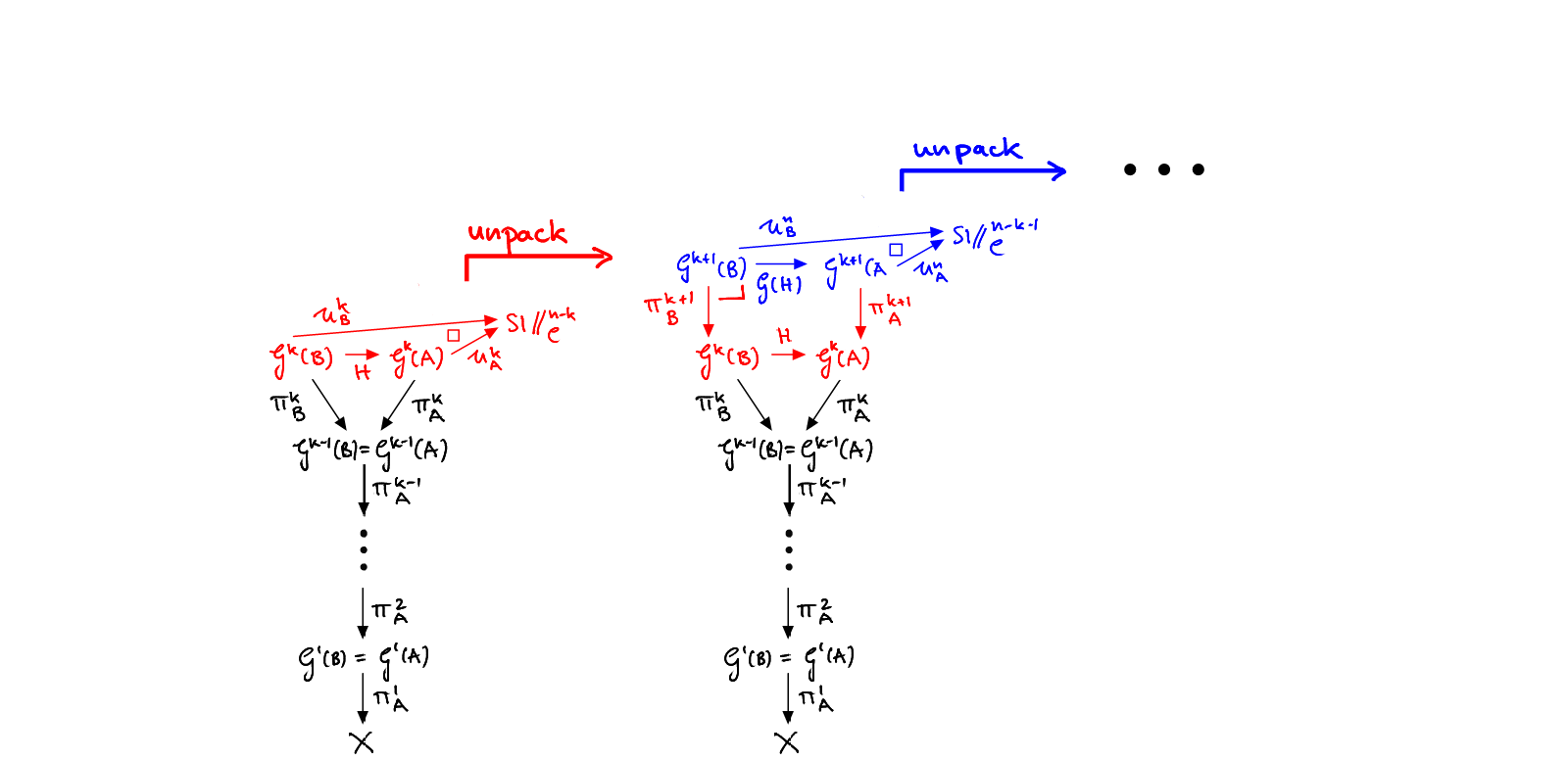}
\endgroup \\*
\begingroup\sbox0{\includegraphics{test/page1.png}}\includegraphics[clip,trim={.3\ht0} {.0\ht0} {.4\ht0} {.0\ht0} ,width=.8\textwidth]{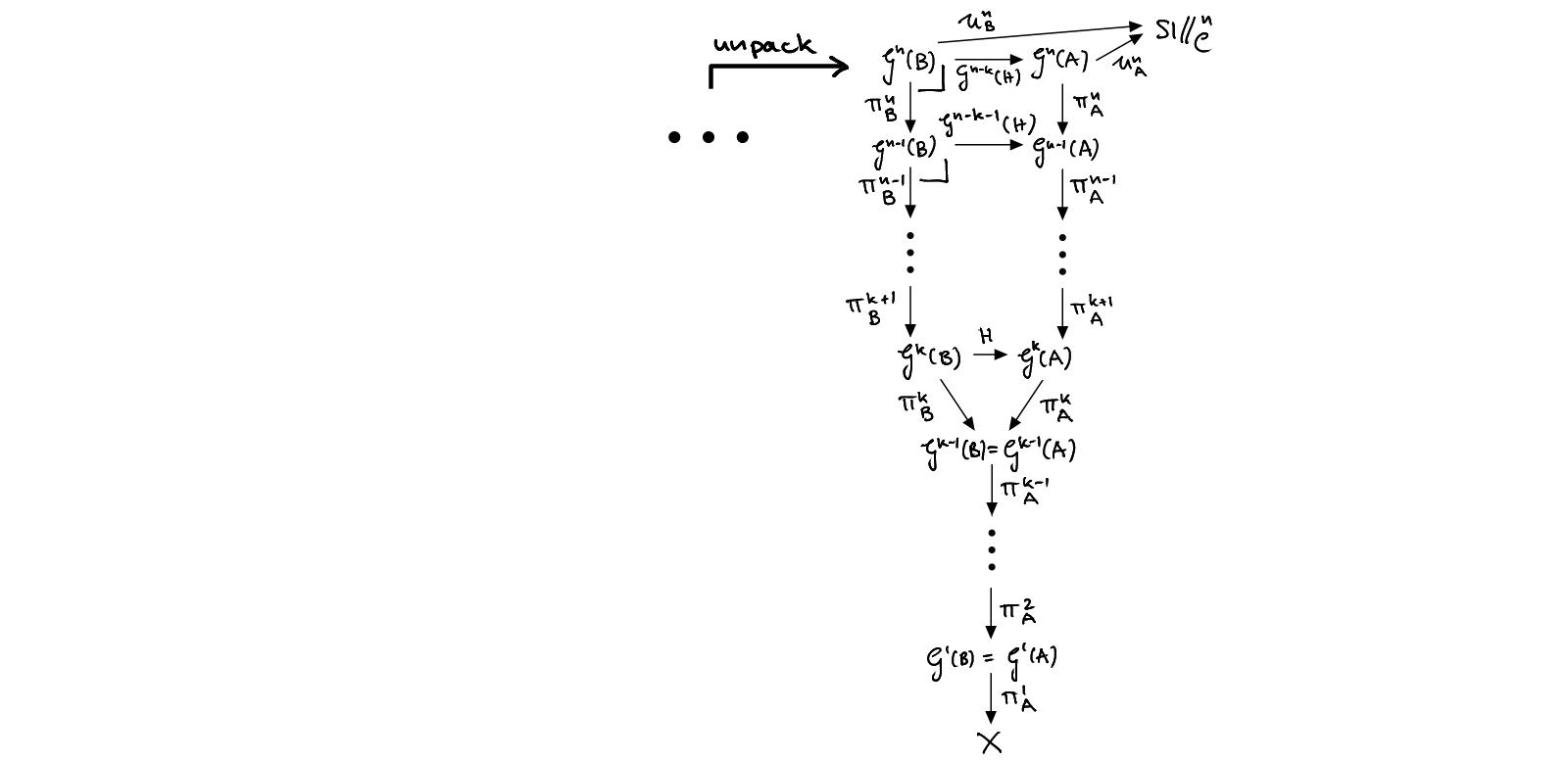}
\endgroup
\end{noverticalspace}
\end{restoretext}

As a corollary we find the following.

\begin{cor}[Base change for labelled cube families] \label{cor:basechange_compact} Let $\scB : X \to \SIvert n \cC$, and $H : Y \to X$. Then
\begin{equation}
\tsR n_{\sT^n_\scB, \tsU n_\scB} H = \tsR n_{\sT^n_\scB H, \tsU n_\scB \upi H n}
\end{equation}
\proof This follows from the previous \autoref{constr:unpacking_collapse} (with $k = 0$) and using \eqref{eq:complete_unpack_repack_inverse}.  \qed
\end{cor}

The preceding corollary is the central result of the section. However, we will need a more hands-on description of $k$-level base change in later proofs. For this, note that the claim gives a sufficient condition for a ``tower base change above level $k$" which we previously asked for. However, we know that all of the above unpacking operations are in fact invertible. This is made precise in the next lemma, to which the following remark provides the inductive step.

\begin{rmk}[Repacking base change for labelled cube families] \label{rmk:unpacking_collapse} Let $\scA$ and $\scB$ as before. If for some $i$, $0 \leq i < n-k$ there is a pullback of the $\SI$-family bundle $\tpi {k+i}_\scB$ along a map $H^i : \tsG {k+i} (\scA) \to \tsG {k+i}(\scB)$ of the form
\begin{equation} 
\xymatrix{ \tsG {k + i + 1}(\scA) \ar[r]^{H^{i+1}} \ar[d]_{\tpi {k+i +1}_{\scA}} \pullback & \tsG {k + i + 1}(\scB)\ar[d]^{\tpi {k+i+1}_{\scB}} \\
\tsG {k + i}(\scA) \ar[r]^{H^{i}}  & \tsG {k + i}(\scB) }
\end{equation}
such that 
\begin{equation} \label{eq:unpacking_collapse}
\tsU {k+i+1}_{\scA} = \tsU {k+i+1}_{\scB} H^{i+1}
\end{equation}
then, firstly using \autoref{rmk:grothendieck_span_construction_basechange}, the pullback forces
\begin{equation}
\tusU {k+i}_{\scA} = \tusU {k+i}_{\scB} H^i
\end{equation}
and, together with our assumption \eqref{eq:unpacking_collapse}, \autoref{claim:grothendieck_span_construction_basechange} then implies
\begin{equation}
\tsU {k+i}_{\scA} = \tsU {k+i}_{\scB} H^i
\end{equation}
\end{rmk}

The preceding \autoref{constr:unpacking_collapse} and remark together establish the following

\begin{lem}[Iterated pullbacks characterise $k$-level base change] \label{lem:unpacking_collapse} Let $\scA, \scB : X \to \SIvert n \cC$ and $0 < k \leq n$, and assume $\sT^{k-1}_\scA = \sT^{k-1}_\scB$ as well as a map of bundles $H : \tpi k_\scA \to \tpi k_\scB$. Then the following are equivalent
\begin{enumerate}
\item $\tsU {k}_{\scA} = \tsU {k}_{\scB} H$, that is $H$ is a $k$-level base change
\item There are $H^i : \tsG {k + i}(\scA) = \tsG {k + i}(\scB)$ for $0 \leq i \leq n-k$ with $H^0 = H$ such that firstly,
\begin{equation} \label{eq:collapse_square_unpack}
\xymatrix{ \tsG {k + i + 1}(\scA) \ar[r]^{H^{i+1}} \ar[d]_{\tpi {k+i +1}_{\scA}} \pullback & \tsG {k + i + 1}(\scB)\ar[d]^{\tpi {k+i+1}_{\scB}} \\
\tsG {k + i}(\scA) \ar[r]^{H^{i}}  & \tsG {k + i}(\scB) }
\end{equation}
for $0 \leq i < n-k$ and secondly,
\begin{equation}
\tsU {n}_{\scA} = \tsU {n}_{\scB} H^{n-k}
\end{equation}
\end{enumerate}
\proof $(i) \imp (ii)$ has been shown in \autoref{constr:unpacking_collapse}, by setting $H^i = \upi H i$. Conversely, $(ii) \imp (i)$ follows from inductively applying \autoref{rmk:unpacking_collapse}. \qed
\end{lem}

\begin{notn}[Base change for negative indices] \label{notn:k_lvl_basechange} Let $\scA, \scB : X \to \SIvert n \cC$ and $0 < k \leq n$, and assume $\sT^{k-1}_\scA = \sT^{k-1}_\scB$ as well as a functor of bundles $H : \tpi k_{\scA} \to \tpi k_{\scB}$ such that $\tsU {k}_{\scA} = \tsU {k}_{\scB} H$. \autoref{constr:unpacking_collapse} defines the notation $\upi H i$ for $0 \leq i \leq (n-k)$. We extend the notation to negative indices by setting
\begin{equation}
\upi H i := \id : \tsG {k + i} (\scA) \to \tsG {k + i} (\scA)
\end{equation}
for $-k \leq i < 0$.  Note that then, all squares of the form \eqref{eq:collapse_square_unpack} commute for $-k \leq i \leq n-k$, and they are pullbacks at every level apart from level $k$. This is summarised in the following illustration
\begin{restoretext}
\begingroup\sbox0{\includegraphics{test/page1.png}}\includegraphics[clip,trim={.3\ht0} {.0\ht0} {.3\ht0} {.0\ht0},width=\textwidth]{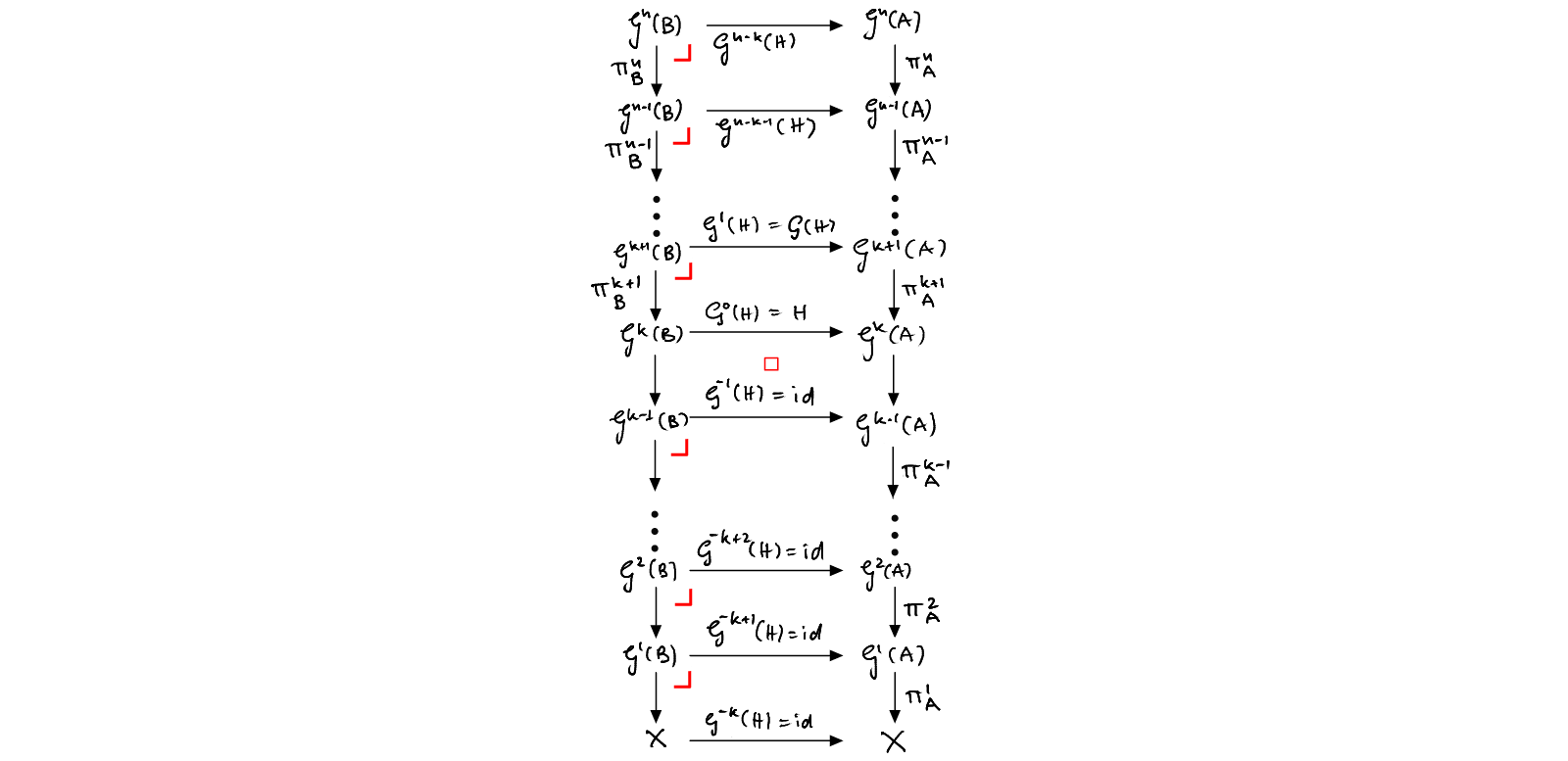}
\endgroup\end{restoretext}
\end{notn}

\subsection{Multi-level base change} \label{ssec:multi_bc}

We generalise the definition of the previous section from ``$k$-level" to ``multi-level" base change. While this will not be of relevance for our later proofs, it provides a unifying perspective of our subsequent theoretical developments.

\begin{defn}[Multi-level base change and the category $\Bunbc^n_\cC$] \label{defn:multilevel_base_change} Let $\scA : X \to \SIvert n \cC$ and $\scB : Y \to \SIvert n \cC$. A multi-level base change $\vvec H : \scB \to \scA$ (from $\scA$ to $\scB$) consists of functors of posets $\vvec H^i : \sG^i(\scB) \to \sG^i(\scA)$ for each $0 \leq i \leq n$ such that
\begin{equation}
\xymatrix{ & \cC & \\
\tsG n(\scA) \ar[ur]^{\tsU n_\scA} \ar[rr]_{\vvec H^n} && \tsG n(\scB) \ar[ul]_{\tsU n_\scB} }
\end{equation}
and for $0 \leq i < n$
\begin{equation}
\xymatrix{ \tsG {i+1}(\scA) \ar[r]^{\vvec H^{i+1}} \ar[d]_{\tpi {i+1}_\scA} & \tsG {i+1}(\scB) \ar[d]^{\tpi {i+1}_\scB} \\
\tsG i(\scA) \ar[r]^{\vvec H^i} & \tsG i(\scB)  }
\end{equation}
Multi-level base changes compose by component-wise composition. That is, if $\vvec H : \scA \to \scB$, and $\vvec K : \scB \to \scC$, then there is $\vvec K \vvec H : \scA \to \scC$ with components
\begin{equation}
(\vvec K \vvec H)^i = \vvec K^i \vvec H^i
\end{equation}
$\SIvert n \cC$-families and multi-level base change then form a category denoted by $\Bunbc^n_\cC$.
\end{defn}

Note that any $k$-level basechange $H$ of $\scA : X \to \SIvert n \cC$ to $\scB : Y \to \SIvert n \cC$ gives rise to a multi-level base change $\vec \cM^k_{H} : \scA \to \scB$ by defining (using \autoref{notn:k_lvl_basechange})
\begin{equation} \label{eq:klvl_bc_abuse}
(\vvec \cM^k_{H})^i := \sG^{k-i}(H)
\end{equation}
Abusing terminology, a multi-level base change $\vvec \cM_H$ derived from a $k$-level base change $H$ will itself also be called a \textit{$k$-level base change} in many cases.

Conversely, given a multi-level collapse $\vvec H$, there is a \textit{unique} decomposition
\begin{equation}
\vvec H = \vec \cM^n_{K^n} \vec \cM^{n-1}_{K^{n-1}} ... \vec \cM^1_{K^1} \vec \cM^0_{K^0} 
\end{equation}
where $\vec \cM^i_{K^i}$ is the multi-level base change associated to a $i$-level base change $K^i$. The uniqueness of this decomposition follows inductively using the definition of $\vec \cM$ and the universal property of pullbacks.

In the next three chapters, we will in much detail study multi-level base changes and their decompositions into $k$-level base changes in two special cases. The first special case concerns fibrewise ``open and surjective" base changes (the $k$-level construction is in \autoref{ch:collapse_intervals} and the multi-level construction in \autoref{ch:collapse_cubes}). The second special case concerns fibrewise ``open and injective" base changes (and is discussed in \autoref{ch:emb}). The decomposition into and classification of $k$-level base changes in both cases is \textit{crucial} to prove important results about the multi-level counterparts.

\chapter{Collapse of intervals} \label{ch:collapse_intervals}

While there is no non-trivial notion of natural isomorphism of singular interval families (cf. \autoref{rmk:grothendieck_span_construction_basechange}), there is a non-trivial notion of natural monomorphism. Understanding this notion of monomorphism is the topic of this chapter, and will be the inductive basis for a theory of normalisation for singular $n$-cube families. In \autoref{sec:injections} we start by discussing the notion of \textit{injections} which will later determine a map of bundles called \textit{collapse}. An injection is the ``minimal data" determining a collapse in the same sense that the profunctorial relation $\SiR(f)$ is fully determined by a function of sets $f$. In \autoref{sec:interval_collapse} we then introduce the notion of collapse for labelled singular interval families before giving certain universal constructions for it in \autoref{sec:properties_collapse}. These constructions will ultimately be leading up to a proof of the normal form theorem (\autoref{thm:normal_forms_unique}) in the next chapter.

\begin{rmk}[Injection terminology] To distinguish the ensuing discussion of natural monomorphisms for singular interval families from what is usually meant by ``bundle monomorphism", or ``bundle embedding", we will adapt the terminology of ``(family) injections" for them.
\end{rmk}

\section{Injections of interval families} \label{sec:injections}

\subsection{Injections as stable singular subset sections} 

We start with definitions of injections and (stable) singular subset sections, showing that there is a straight-forward correspondence between them (namely, the image of an injection is a stable singular subset section).

\begin{defn}[Injections] \label{defn:natural_injections}  Let $X$ be a poset, and let $\scA, \scB : X \to \SI$ be singular interval families. An \textit{injection} $\lambda : \scA \into \scB$ of singular interval families $\scA, \scB : X \to \SI$ is a natural transformation $\scA \to \scB$ such that each component $\lambda_x : \scA(x) \to_{\SI} \scB(x)$ is monic, that is, each component $\lambda_x$, $x \in X$, is an injective map
\begin{equation}
\lambda_x : \singcont(\scA(x)) \into \singcont(\scB(x))
\end{equation}
\end{defn}

We give three examples of injections.

\begin{egs}[Injections] \label{eg:injections} We first define two \SI-morphisms, $f_0 : \singint 2 \to\oSI \singint 3$ and $g_0 : \singint 3 \to\oSI \singint 3$ by setting
\begin{restoretext}
\begingroup\sbox0{\includegraphics{test/page1.png}}\includegraphics[clip,trim=0 {.2\ht0} 0 {.15\ht0} ,width=\textwidth]{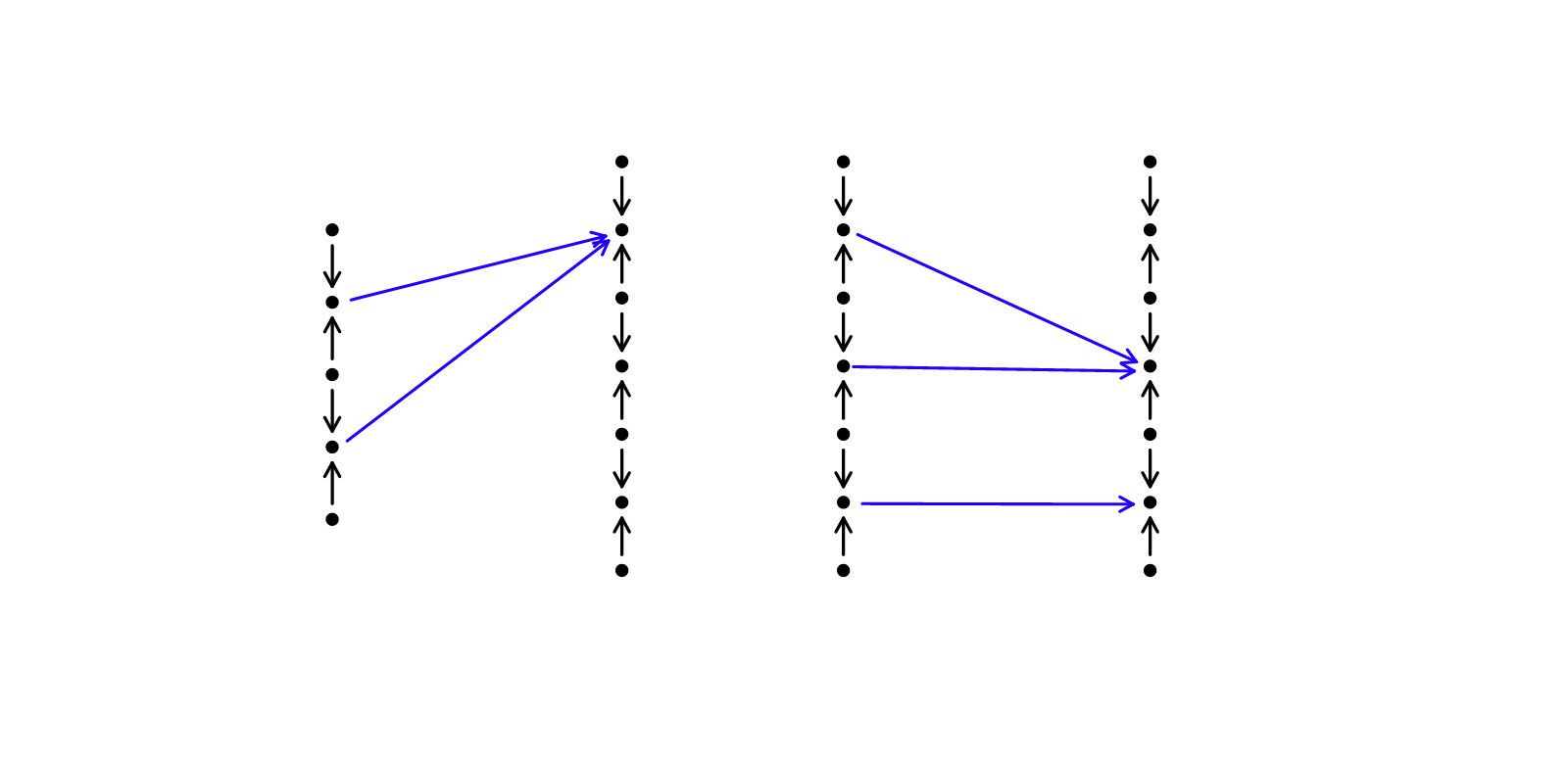}
\endgroup\end{restoretext}
Here, \cblue{} arrows indicate the mapping of $f_0$ on the left and similarly on the right for $g_0$. This allows us to define a functor $\scC : \bnum{3} \to \SI$ by setting $\scC(0 \to 1) = f_0$ and $\scC(1 \to 2) = g_0$. We now give three examples of injections into $\scC$.
\begin{enumerate}
\item We construct an inclusion $\lambda_1 : \scB_1 \into \scC$ where $\scB_1 : \bnum{3} \to \SI$ is defined by $\scB_1(0 \to 1) = f_1$, $\scB_1(1 \to 2) = g_1$ with $f_1$ and $g_1$ being the morphisms whose mappings are indicated by \cpurple{} arrows as
\begin{restoretext}
\begingroup\sbox0{\includegraphics{test/page1.png}}\includegraphics[clip,trim=0 {.25\ht0} 0 {.15\ht0} ,width=\textwidth]{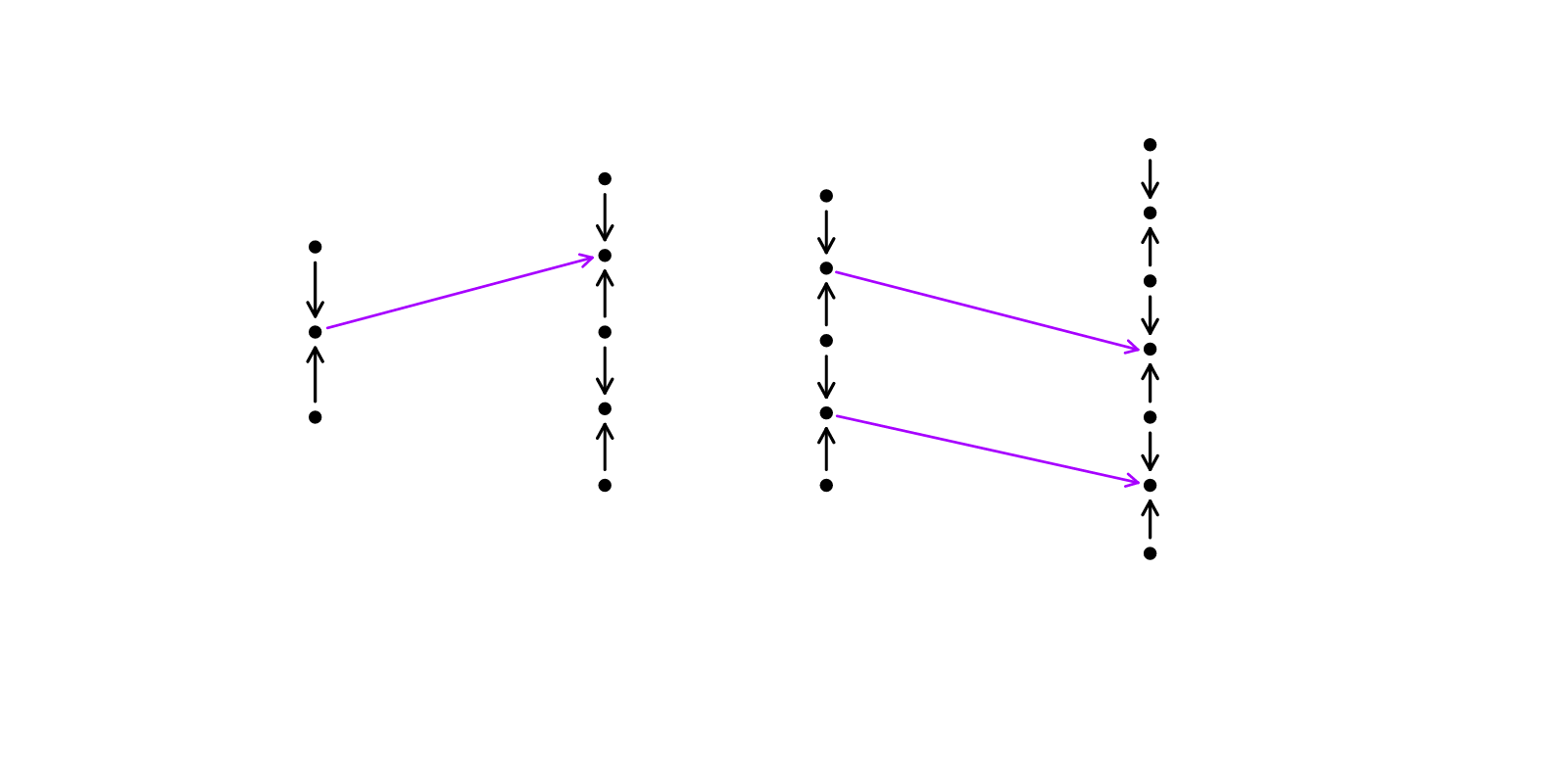}
\endgroup\end{restoretext}
on the left and right respectively.
Then we can define the components $(\lambda_1)_0$, $(\lambda_1)_1$, $(\lambda_1)_2$ by the following \SI-morphisms
\begin{restoretext}
\begingroup\sbox0{\includegraphics{test/page1.png}}\includegraphics[clip,trim=0 {.2\ht0} 0 {.1\ht0} ,width=\textwidth]{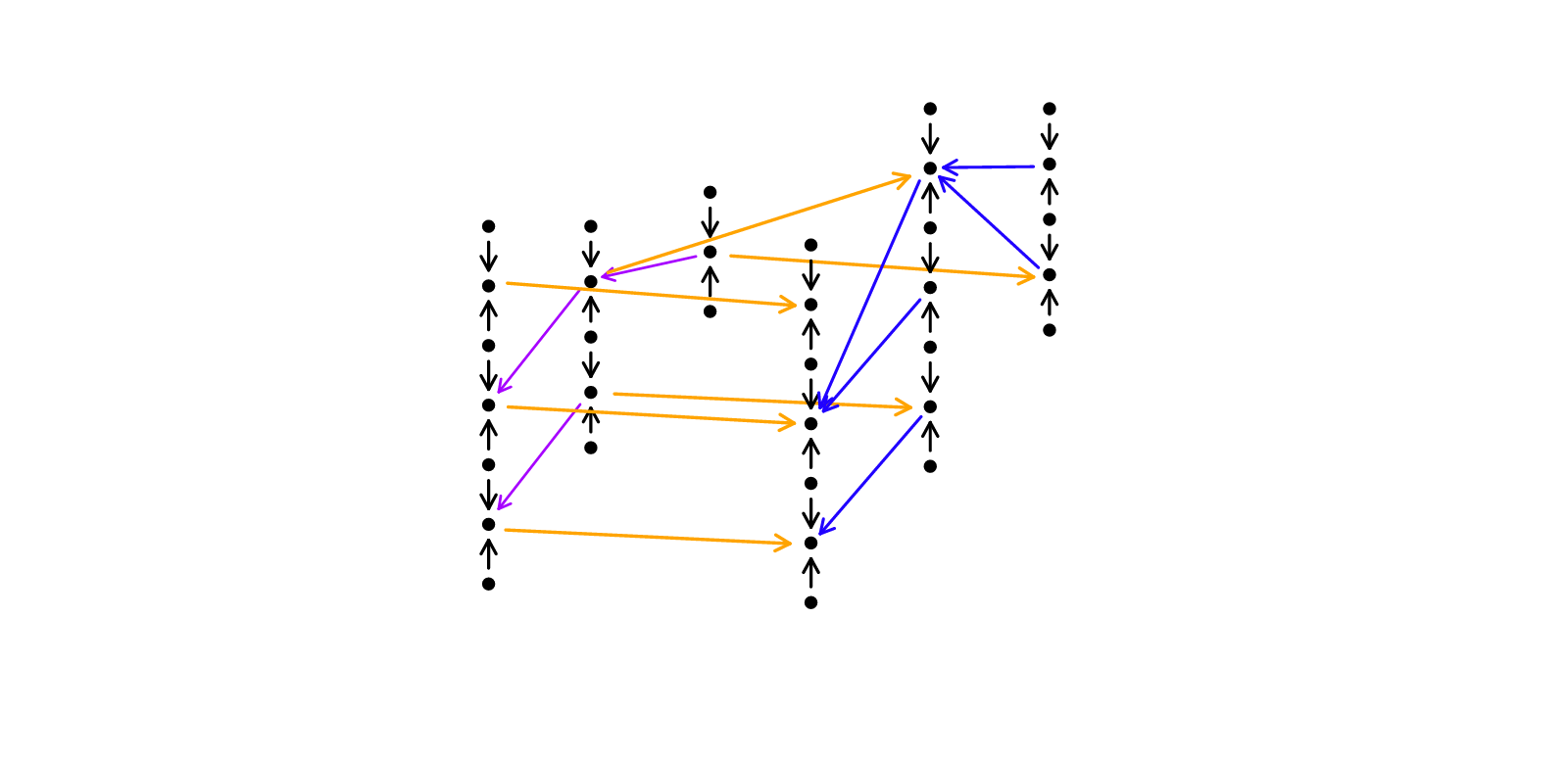}
\endgroup\end{restoretext}
whose respective mappings are marked in \corange{}. The reader can verify that this satisfies injectivity and naturality.

\item Next, we construct an inclusion $\lambda_2 : \scB_2 \into \scC$ where $\scB_2 : \bnum{3} \to \SI$ is defined by $\scB_2(0 \to 1) = f_2$, $\scB_2(1 \to 2) = g_2$ with $f_2$ and $g_2$ being the morphisms whose mappings are given by \cpink{} arrows
\begin{restoretext}
\begingroup\sbox0{\includegraphics{test/page1.png}}\includegraphics[clip,trim=0 {.3\ht0} 0 {.25\ht0} ,width=\textwidth]{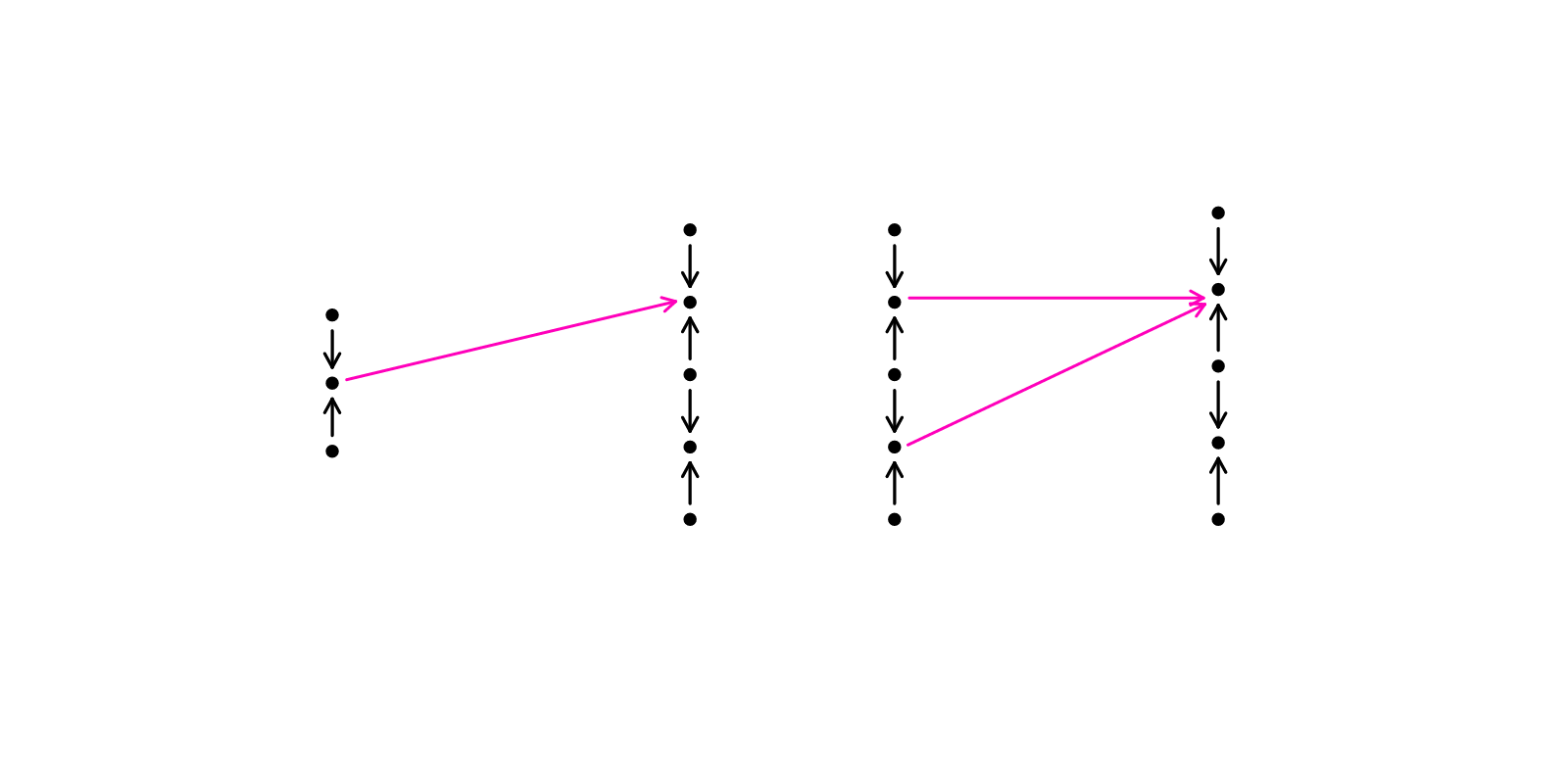}
\endgroup\end{restoretext}
on the left and right respectively.
Then we can define the components $(\lambda_2)_0$, $(\lambda_2)_1$, $(\lambda_2)_2$ by the following \SI-morphisms
\begin{restoretext}
\begingroup\sbox0{\includegraphics{test/page1.png}}\includegraphics[clip,trim=0 {.15\ht0} 0 {.15\ht0} ,width=\textwidth]{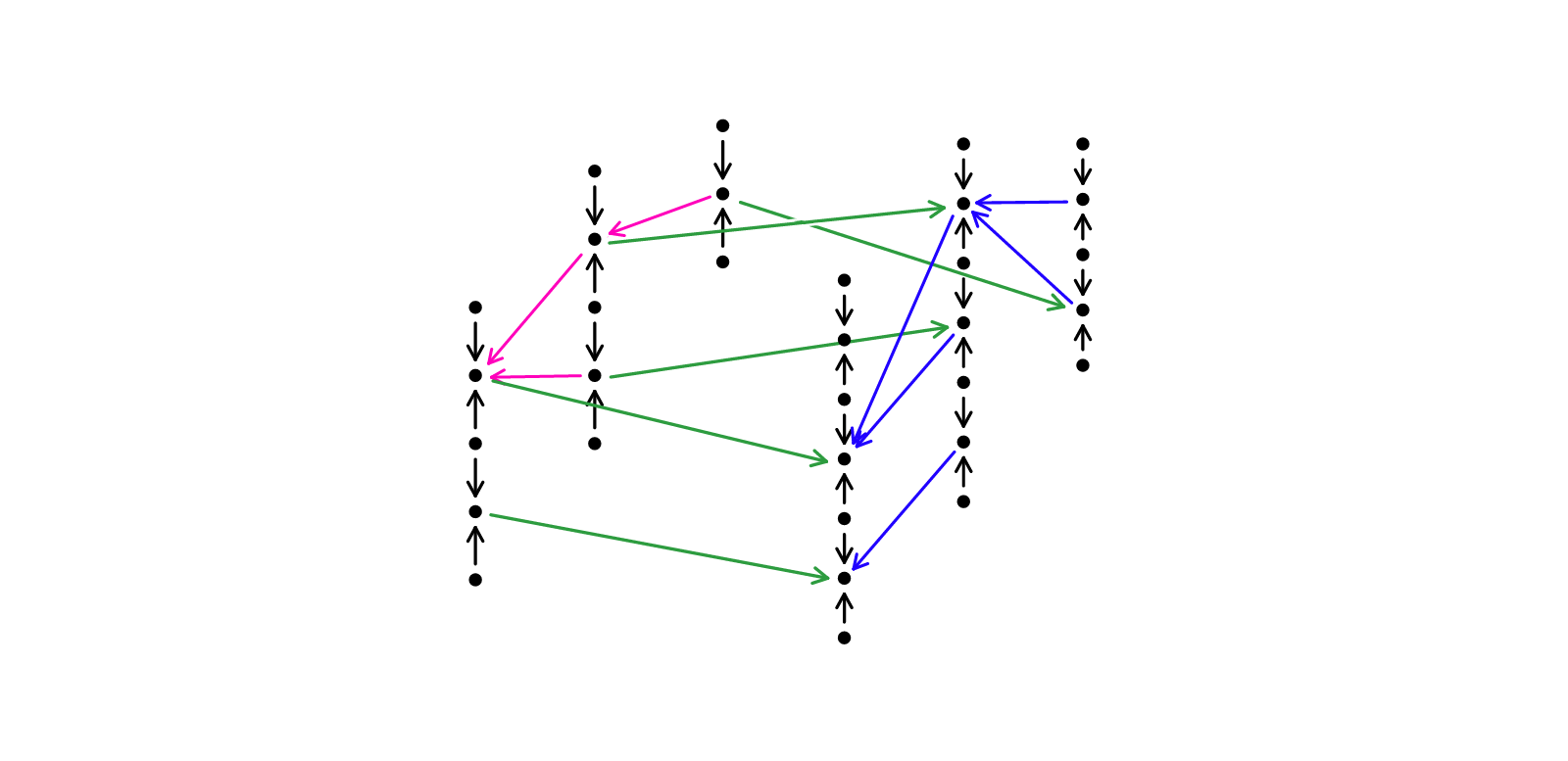}
\endgroup\end{restoretext}
marked in \cdarkgreen{}.

\item Finally, we construct an inclusion $\eps : \scA \into \scC$ where $\scA : \bnum{3} \to \SI$ is defined by $\scA(0 \to 1) = f_1$, $\scA(1 \to 2) = g_1$ with $f_1$ and $g_1$ being the morphisms whose mappings are given by (lighter)  \cgreen{} arrows
\begin{restoretext}
\begingroup\sbox0{\includegraphics{test/page1.png}}\includegraphics[clip,trim=0 {.3\ht0} 0 {.3\ht0} ,width=\textwidth]{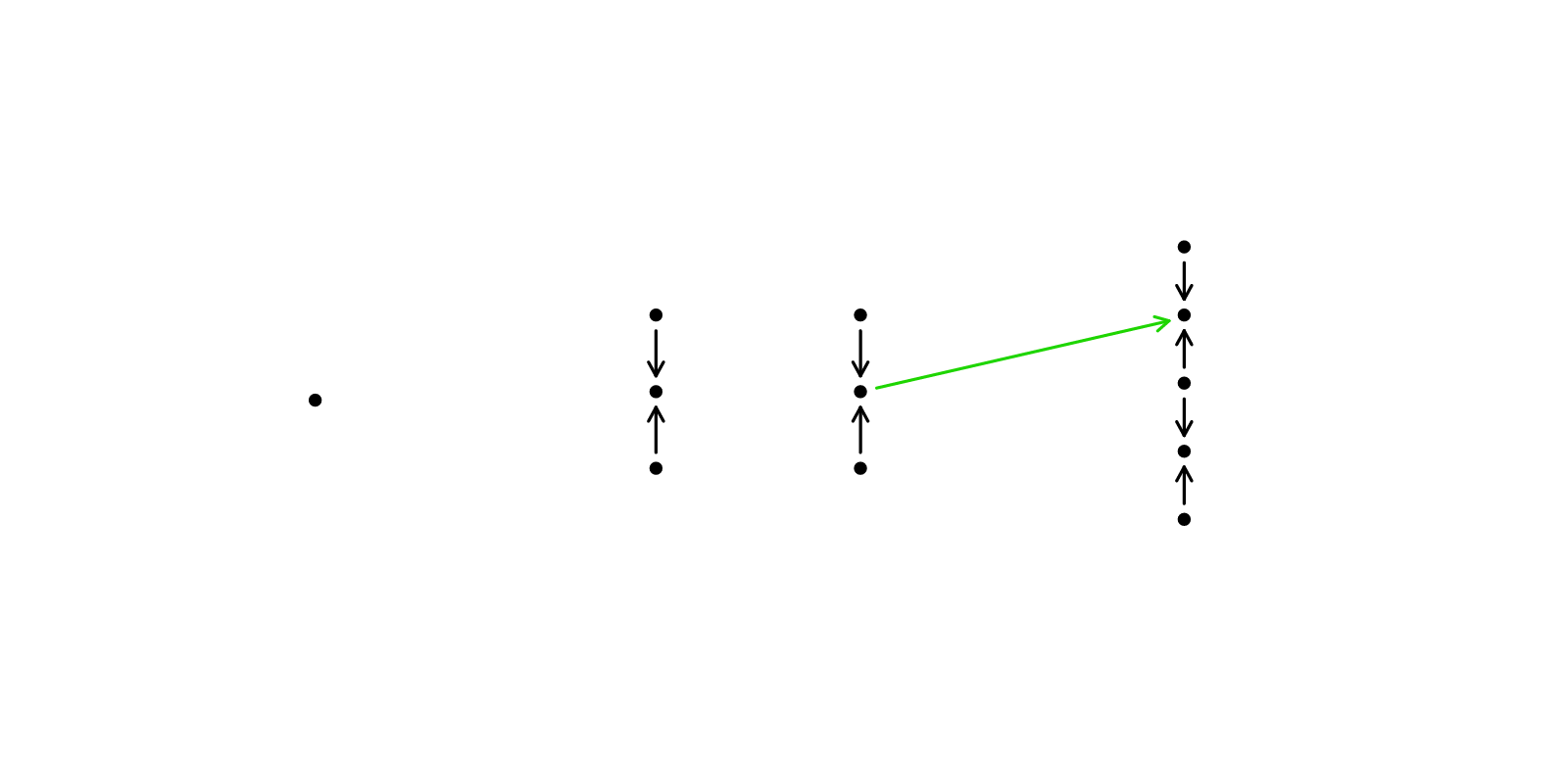}
\endgroup\end{restoretext}
on the left and right respectively. Note that on the left we are giving a mapping from the initial singular interval, thus no arrows were drawn.
Then we can define the components $\eps_0$, $\eps_1$, $\eps_2$ by the following \SI-morphisms
\begin{restoretext}
\begingroup\sbox0{\includegraphics{test/page1.png}}\includegraphics[clip,trim=0 {.15\ht0} 0 {.15\ht0} ,width=\textwidth]{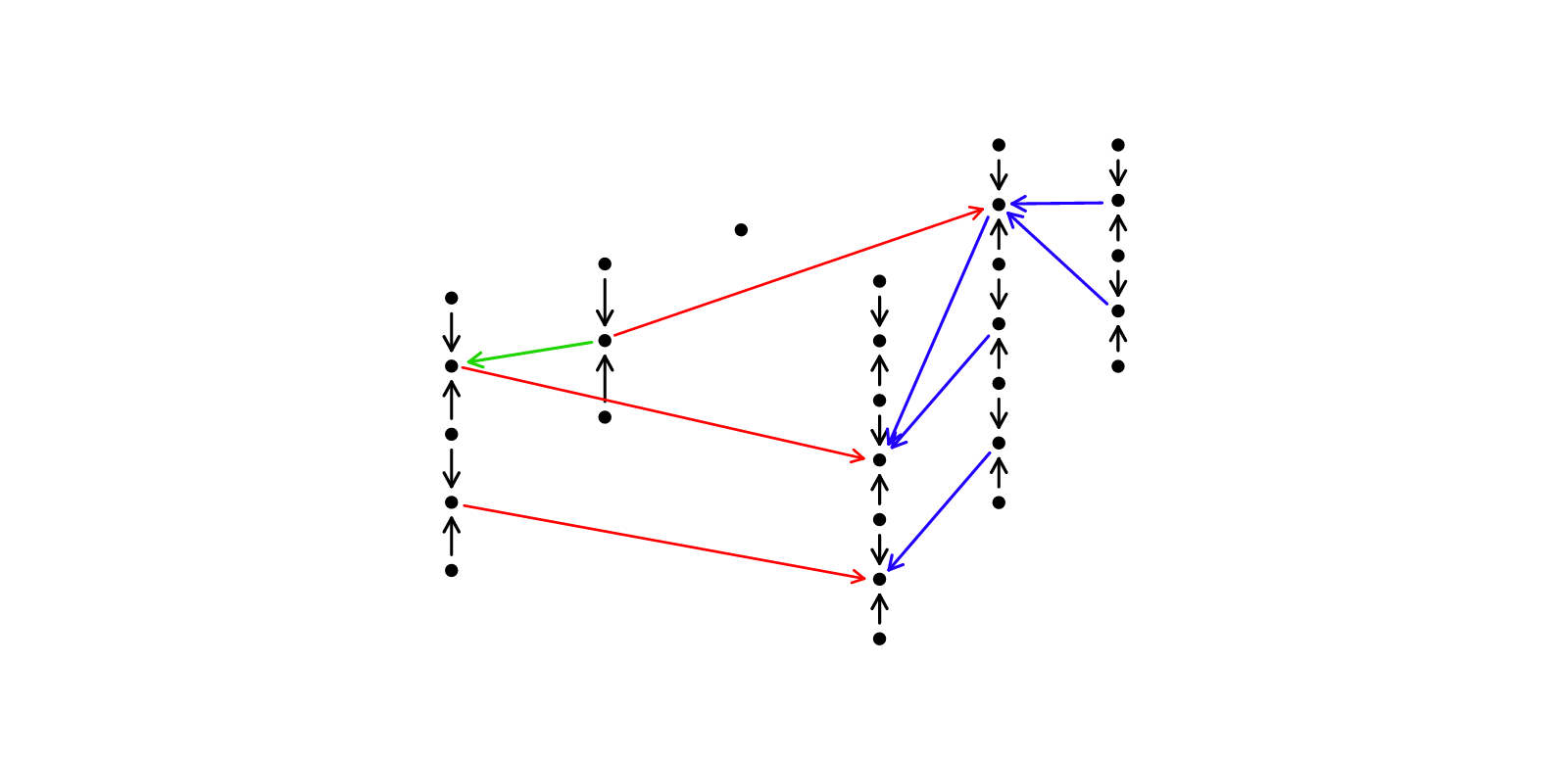}
\endgroup\end{restoretext}
marked in \cred{}.
\end{enumerate}
\end{egs}

\begin{constr}[Injections of subsets of singular heights] \label{defn:eta_inclusion} Given a singular interval $I \in \SI$, and a subset $S \subset \singcont(I)$ of singular heights, denote by $\intrel{S}$ the singular interval of height $\iH_{\intrel{S}} = \#S$ (where $\#S$ denotes the number of elements in $S$) and let
\begin{equation} 
\eta_{S} : \singcont(\intrel{S}) \into \singcont(I)
\end{equation}
be the unique monotone injection with image $S \subset \singcont(I)$. Explicitly, $\eta_{S}$ maps $(2j - 1) \in \intrel{S}$ to $a_j$, where $a_j \in S = \Set{a_1 < a_2 < ... < a_{\#S}}$ is the $j$th element of $S$ in $\lZ$-order. $\eta_{S}$ is called the \textit{injection of the subset $S$} of singular heights.
\end{constr}

\begin{eg}[Injections of subset of singular heights] For the concrete choice of $S = \Set{1,5,7} \subset \singcont(\singint 5)$, the injection $\eta_S$ is given by the \SI-morphism
\begin{restoretext}
\begingroup\sbox0{\includegraphics{test/page1.png}}\includegraphics[clip,trim=0 {.1\ht0} 0 {.1\ht0} ,width=\textwidth]{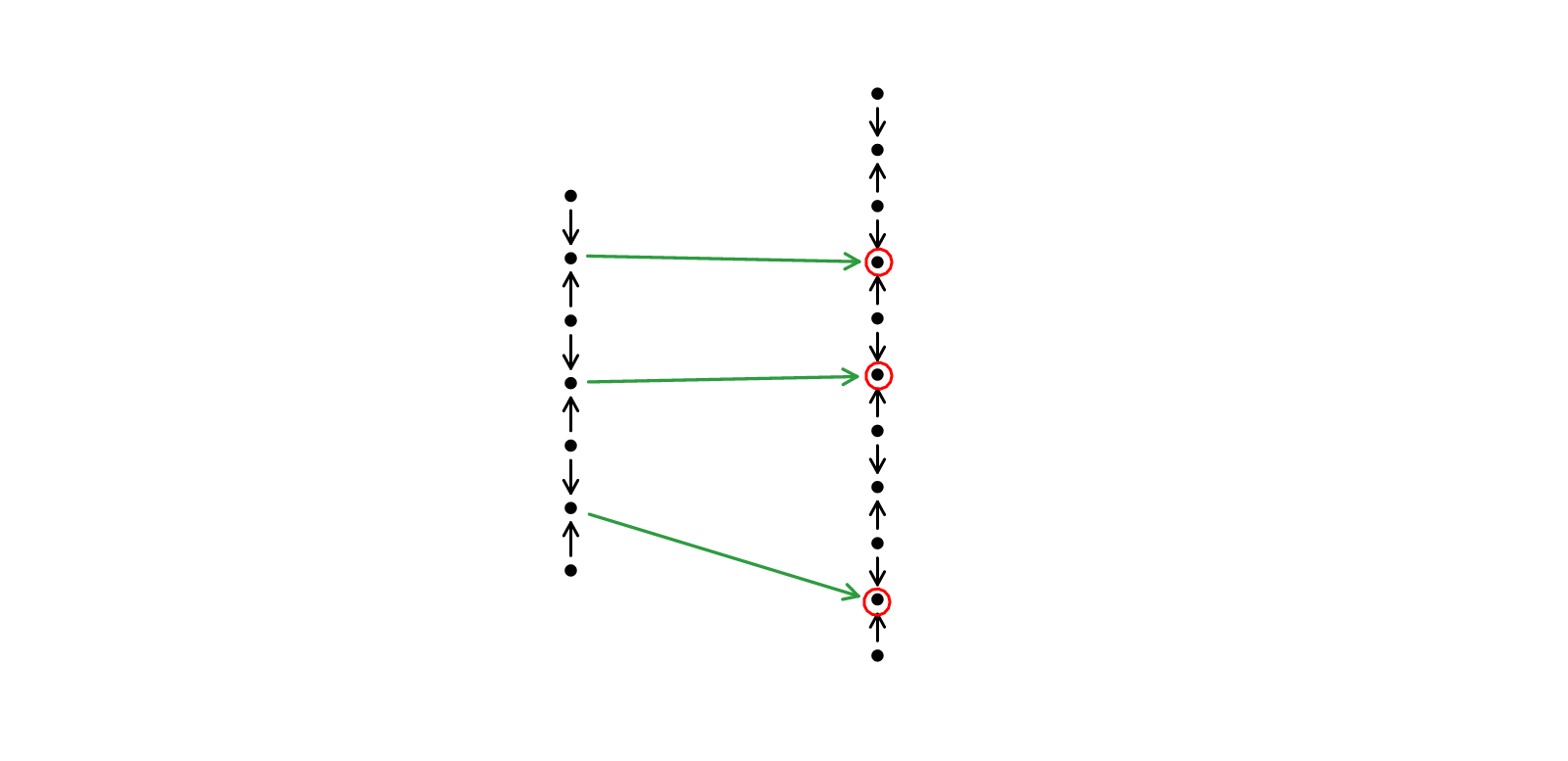}
\endgroup\end{restoretext}
($S$ is marked by \cred{} circles).
\end{eg}

\begin{defn}[Singular subset sections] \label{defn:subset_sections}
Given a singular interval family $\scB : X \to \SI$, a \textit{singular subset section $\cF$ of $\scB$} is a family of subsets $\cF_x \subset \singcont(\scB(x))$ indexed by objects $x\in X$. We say $\cF$ is a \textit{stable} singular subset section of $\scB$ if for all $(x \to y) \in \mor(X)$ we have 
\begin{equation} \label{eq:natural_injection_condition}
\scB(x \to y) (\cF_x) \subset \cF_y
\end{equation}
\end{defn}

\begin{eg}[Singular subset sections] A singular subset section $\cF$ for $\scC$ as defined in \autoref{eg:injections} consists of three sets $\cF_i \subset \singcont(\scC(i))$, $i \in \obj(\bnum{3})$. For instance we can set $\cF_0 = \Set{3}$, $\cF_1 = \Set{1,3,5}$ and $\cF_2 = \Set{3}$ which can be visualised as
\begin{restoretext}
\begingroup\sbox0{\includegraphics{test/page1.png}}\includegraphics[clip,trim=0 {.15\ht0} 0 {.15\ht0} ,width=\textwidth]{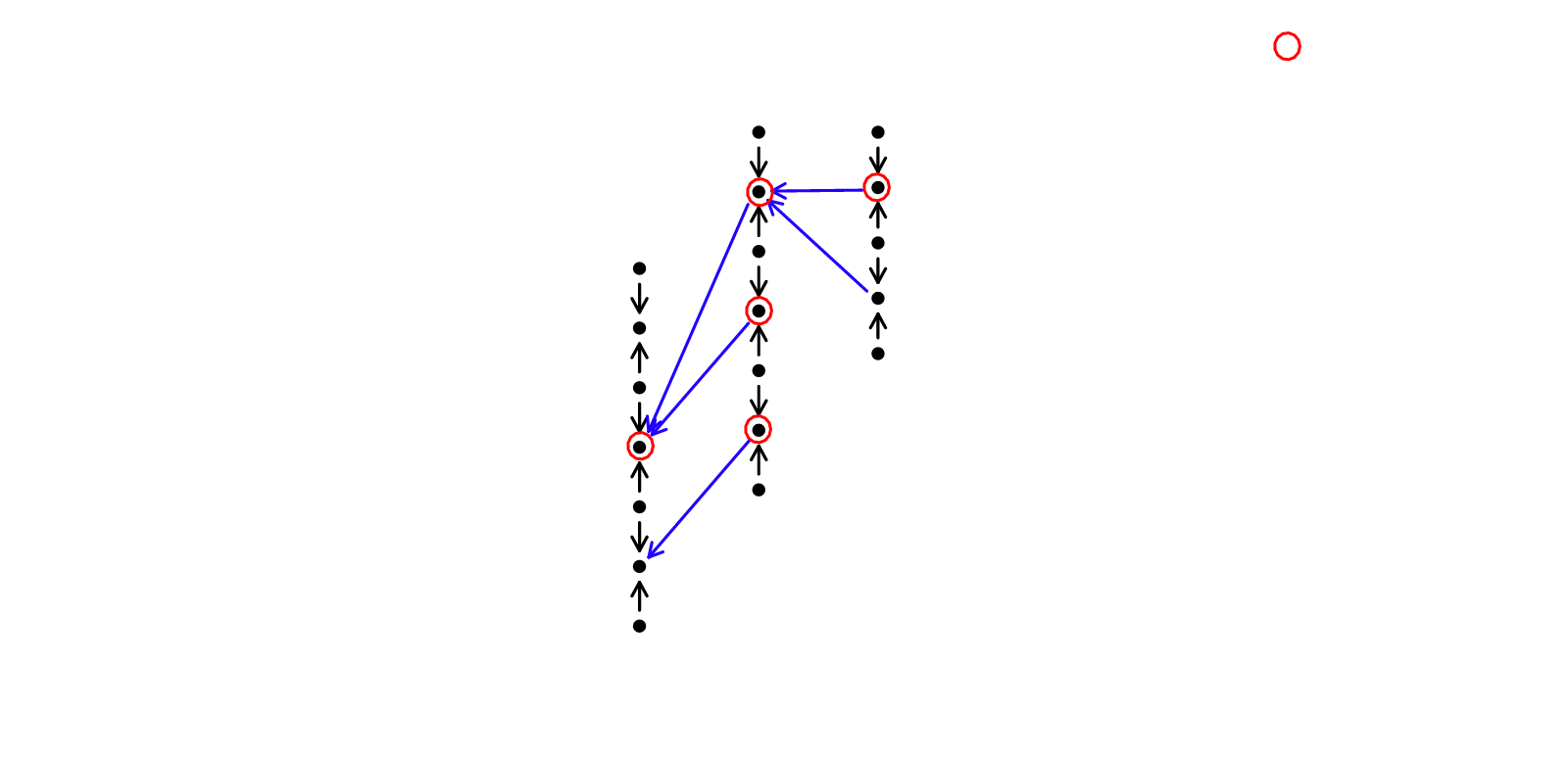}
\endgroup\end{restoretext}
(the elements of the singular subset section $\cF$ are marked by \cred{} circles). We remark that the above subset section is not stable. However the following subset section $\cF$, $\cF'_0 = \Set{3}$, $\cF'_1 = \Set{3,5}$ and $\cF'_2 = \Set{3,5}$ is stable
\begin{restoretext}
\begingroup\sbox0{\includegraphics{test/page1.png}}\includegraphics[clip,trim=0 {.15\ht0} 0 {.15\ht0} ,width=\textwidth]{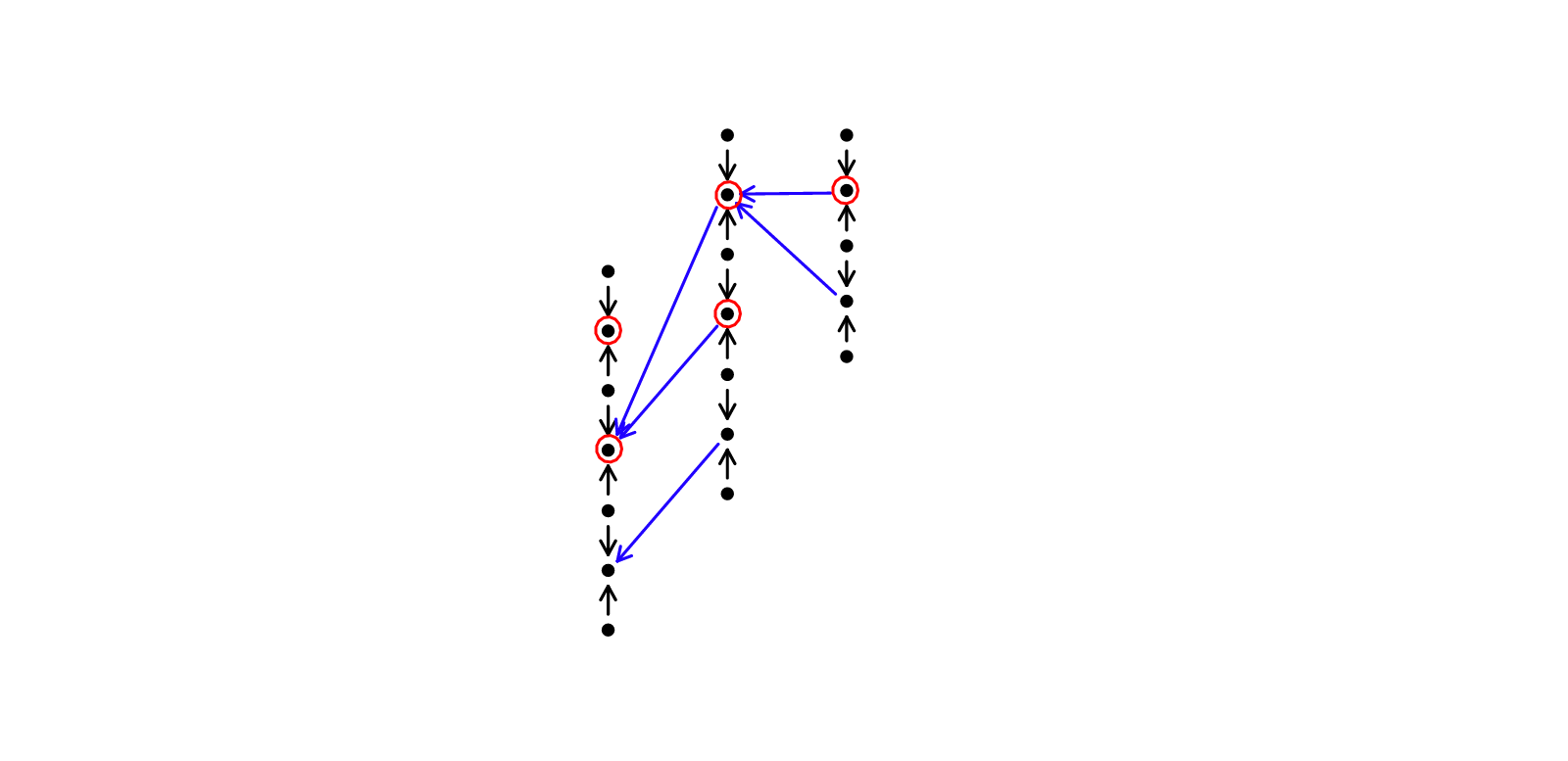}
\endgroup\end{restoretext}
(the elements of the singular subset section $\cF'$ are marked by \cred{} circles).
\end{eg}

\begin{rmk}[Images of injections are stable] \label{rmk:im_of_inj_stable}
Note that any injection $\lambda : \scA \into \scB$ defines a singular subset section, namely the one that is given by its (componentwise) image. For instance $\eps$ from \autoref{eg:injections} define the singular subset section
\begin{restoretext}
\begingroup\sbox0{\includegraphics{test/page1.png}}\includegraphics[clip,trim=0 {.15\ht0} 0 {.15\ht0} ,width=\textwidth]{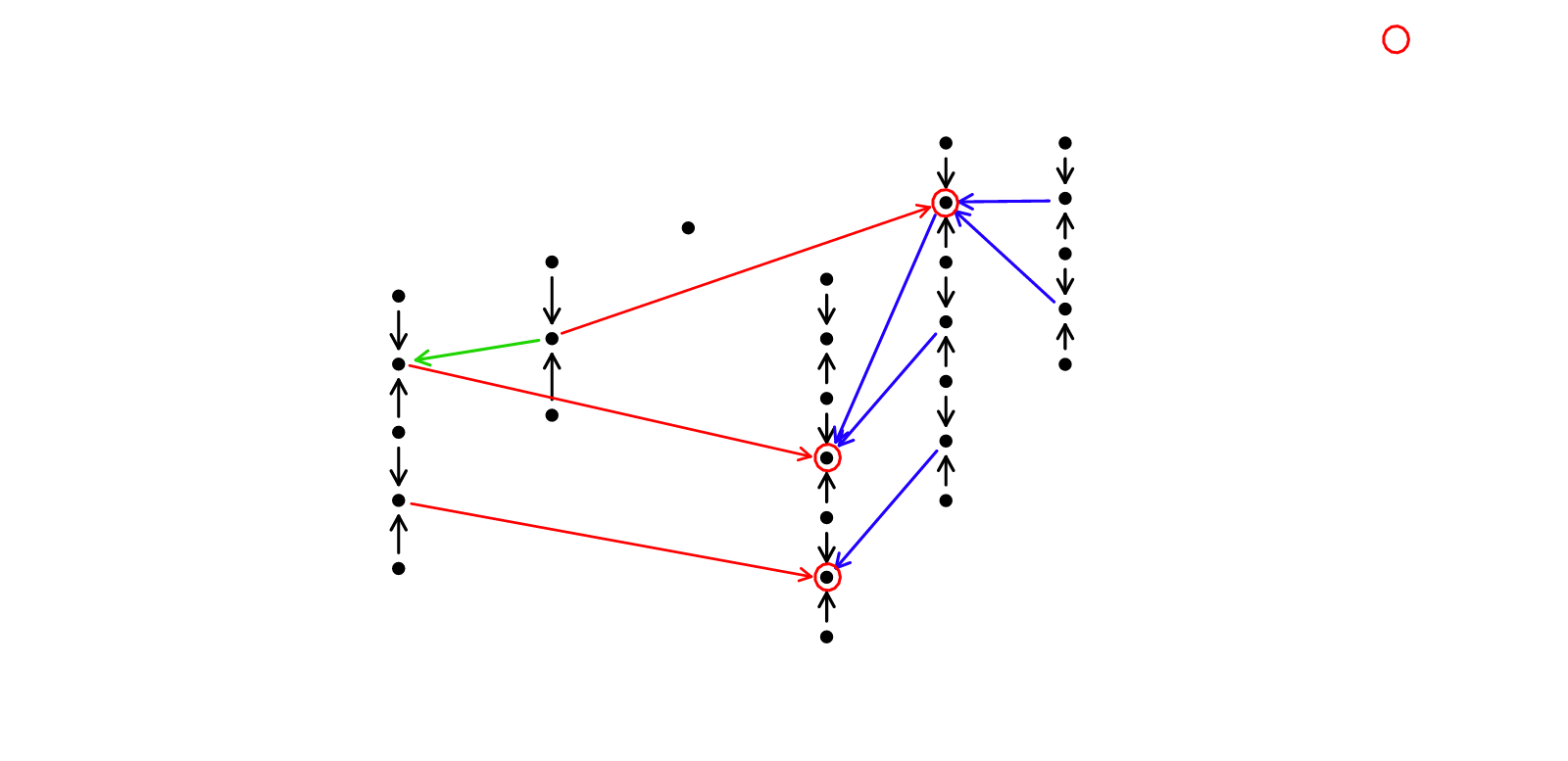}
\endgroup\end{restoretext}
We observe that this singular subset section is stable since $\eps$ is natural. The correspondence and notation of stable singular subset sections and injections is given by the next construction.
\end{rmk}

\begin{constr}[Injections associated to stable singular subset sections and stable singular subset sections associated to injections] \label{defn:stability_vs_injections} Let $\scB : X \to \SI$ be a singular interval family.
\begin{enumerate}
\item Let $\cF$ be a stable singular subset section of $\scB$. We define the \textit{$\cF$-associated} singular interval family $\intrel{\cF} : X \to \SI$ by setting $\intrel{\cF}(x) = \intrel{\cF_x}$, for $x \in X$, and for $(x \to y) \in \mor(X)$ we have
\begin{equation} \label{eq:natural_injection_naturality}
\intrel{\cF} (x \to y) := \eta_{\cF_y} \inv \scB(x\to y) \eta_{\cF_x}
\end{equation}
Note that this is well-defined by \eqref{eq:natural_injection_condition}, which guarantees that $\im(\scB(x\to y) \eta_{\cF_x}) \subset \im (\eta_{\cF_y})$. Functoriality of $\intrel{\cF}$ follows from functoriality of $\scB$. Next, we define the \textit{$\cF$-associated} injection $\eta_{\cF} : \intrel{\cF} \into \scB$ componentwise by setting
\begin{equation}
(\eta_\cF)_x := \eta_{\cF_x}
\end{equation}
Note that by \eqref{eq:natural_injection_naturality} this is indeed a natural transformation. By \autoref{defn:natural_injections} and \autoref{defn:eta_inclusion} it is thus an injection $\eta_\cF : \intrel{\cF} \into \scB$.
\item Let $\lambda : \scA \into \scB$ be an injection. We define the \textit{$\lambda$-associated} stable singular subset section $\cS^\lambda$ of $\scB$ by setting
\begin{equation}
\cS^\lambda_x := \im(\lambda_x) \subset \singcont(\scB(x))
\end{equation}
Note that this defines a stable singular subset section $\cS^\lambda$ because naturality of $\lambda$ means $\lambda_y \scA (x \to y) =  \scB(x\to y) \lambda_x$ which implies $\im( \lambda_y \scA (x \to y) ) \subset \im (\lambda_y) = \cS^\lambda_y$, and thus $\cS^\lambda$ satisfies  \eqref{eq:natural_injection_condition}.
\end{enumerate}
\end{constr}

As motivated earlier, stable singular subset sections and injections are in one-to-one correspondence, which in case it is not yet obvious, is recorded by the following remark.

\begin{rmk}[Association of sub-families with injections is bijective] \label{lem:stability_vs_injections} The constructions $\lambda \mapsto \cS^\lambda$ and $\cF \mapsto \eta_\cF$ Given in \autoref{defn:stability_vs_injections} are mutually inverse to each other. That is, given $\scB : X \to \SI$, $\cF$ a stable family of $\scB$ and $\lambda : \scA \into \scB$ an injection then $\cS^{\eta_\cF} = \cF$, and $\eta_{\cS^\lambda} = \lambda$.

\end{rmk}

\subsection{Pullbacks of injections}

The goal of this section is to give a construction for pullbacks of injections, which will ultimately be relevant for our discussion of normal forms of cube families.

\begin{claim}[Factorisation of injections] \label{claim:injection_factorisation} Given injections $\eps : \scA \into \scC$, $\lambda : \scB \into \scC$ for $A,B,C : X \to \SI$, then there is a (necessarily unique) factorising injection $\mu : A \into B$ such that
\begin{equation}
\xymatrix@R=0.5cm@C=1cm{\scA \ar@{-->}[dr]_{\mu} \ar[rr]^\eps & & \scC \\
& \scB \ar[ur]_\lambda & }
\end{equation}
if and only if $\cS^\eps \subset \cS^\lambda$ (by which we mean $\cS^\eps_x \subset \cS^\lambda_x$ for all $x \in X$)
\proof The proof is \stfwd{}. If $\eps = \lambda\mu$ then we have 
\begin{align*} \label{eq:stable_families_of_composites}
\cS^{\eps}_x  & = \cS^{\lambda\mu}_x \\
&= \im\big(\lambda_x \mu_x\big) \\
&= \lambda_x\big(\im(\mu_x)\big) \\
&= \lambda_x(\cS^{\mu}_x) \subset \cS^\lambda_x
\end{align*}
where for the second line we used \autoref{defn:stability_vs_injections}.

Conversely, assume $\cS^\eps \subset \cS^\lambda$. Then, we write $\lambda\inv(\cS^\eps)$ for the singular subset section of $\scB$ with members $\lambda\inv(\cS^\eps)_x := \lambda\inv_x(\cS^\eps_x)$. We claim this singular subset section is stable: for $(x \to y) \in \mor(X)$ we verify
\begin{align*}
\scB(x \to y) (\lambda\inv_x(\cS^\eps_x)) &= \lambda\inv_y(\scC(x \to y)(\cS^\eps_x)) \\
&\subset \lambda\inv_y(\cS^\eps_y) 
\end{align*}
where we used naturality of $\lambda$ in the first step and stability of $\cS^\eps$ in the second step. Since $\lambda\inv(\cS^\lambda)$ is stable, we can use \autoref{defn:stability_vs_injections} to construct $\mu := \eta_{\lambda\inv(\cS^\eps)}$. We then compute
\begin{align}
\cS^{\lambda  \mu } &= \lambda (\cS^{\mu }) \hspace{7em}\\
&= \lambda  (\lambda\inv  (\cS^{\eps}))\\
&= \cS^{\eps}
\end{align}
In the first step we used \eqref{eq:stable_families_of_composites}, then applied the definition of $\mu$ and finally computed componentwise $(\lambda )_x(\lambda )\inv_x(\cS^\eps_x) = \cS^\eps_x$ which uses $\cS^\eps_x \subset \cS^{\lambda }_x$. By \autoref{lem:stability_vs_injections} we infer $\eps = \lambda\mu$ as required. Note that $\lambda (\cS^{\mu}) = \cS^\eps$ now also implies that our choice for $\cS^\mu$ was unique. \qed
\end{claim}

\begin{eg}[Factorisation of injections] \label{eg:factorisation_of_injections} Using the definitions in \autoref{eg:injections} we see that $\cS^{\eps} \subset \cS^{\lambda_1}$ and we thus obtain $\mu_1 : \scA \into \scB_1$ which (together with $\lambda_1$ in \cdarkgreen{} and $\eps$ in \cred{}) can be depicted as follows
\begin{restoretext}
\begingroup\sbox0{\includegraphics{test/page1.png}}\includegraphics[clip,trim=0 {.15\ht0} 0 {.1\ht0} ,width=\textwidth]{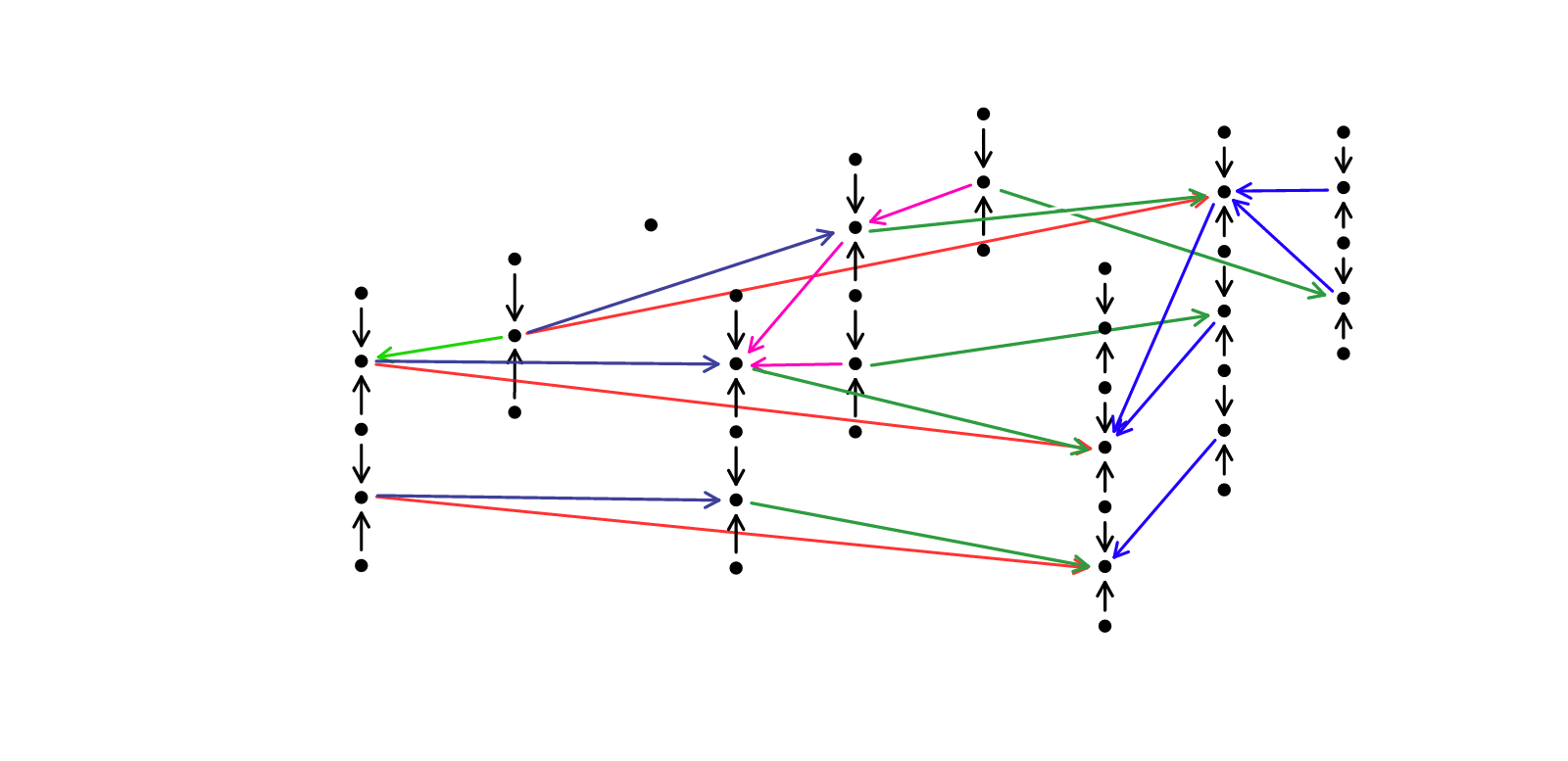}
\endgroup\end{restoretext}
(Here, the components of $\mu_1$ are marked in \cdarkblue{}). 

Similarly, we find from the definition in \autoref{eg:injections} that $\cS^{\eps} \subset \cS^{\lambda_2}$ and we thus obtain $\mu_2 : \scA \into \scB_2$ which (together with $\lambda_2$ in \corange{} and $\eps$ in \cred{}) can be depicted as follows
\begin{restoretext}
\begingroup\sbox0{\includegraphics{test/page1.png}}\includegraphics[clip,trim=0 {.1\ht0} 0 {.2\ht0} ,width=\textwidth]{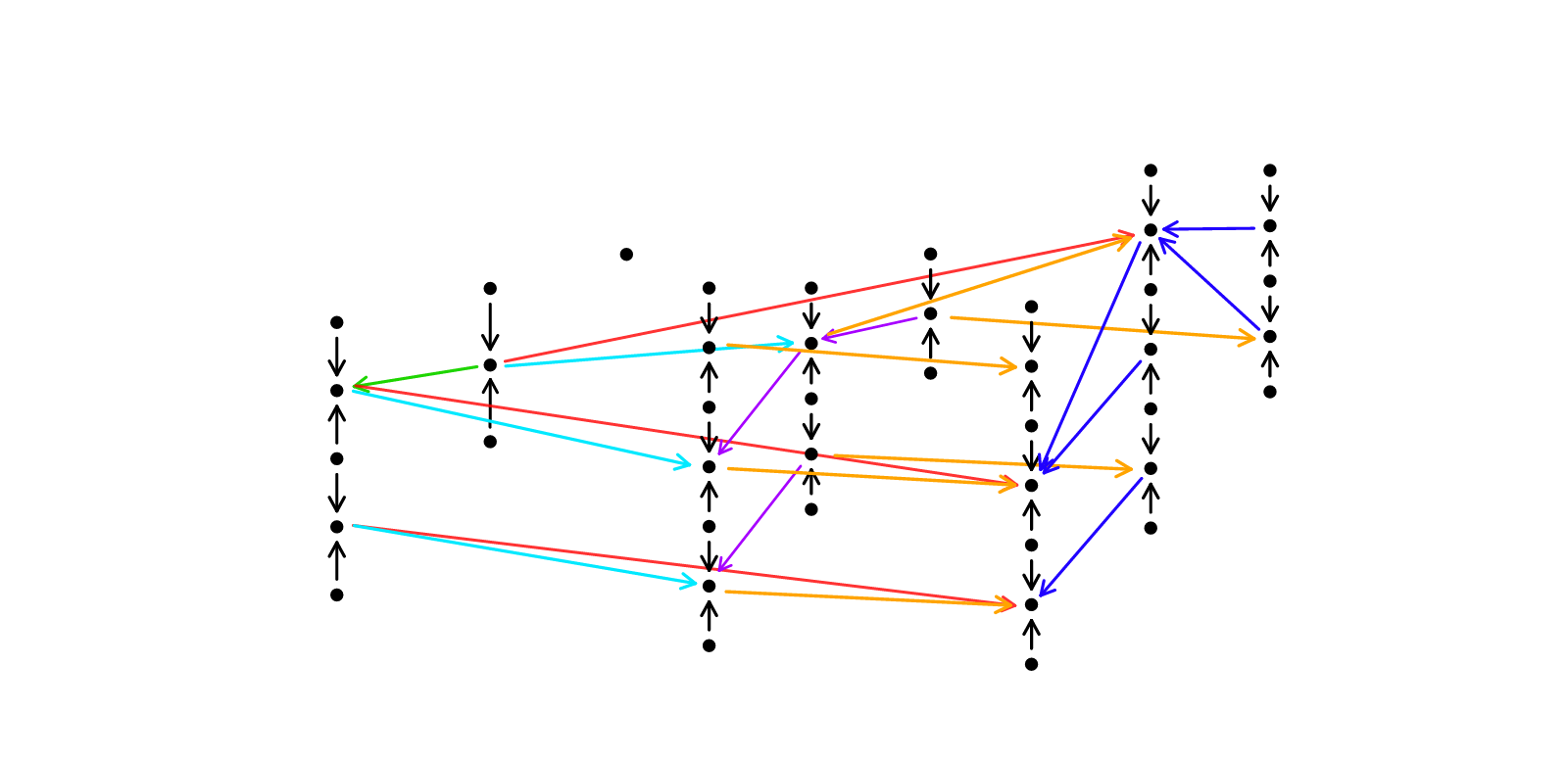}
\endgroup\end{restoretext}
(Here, the components of $\mu_2$ are marked in \cturquoise{}).
\end{eg}

\begin{lem}[Pullbacks of injections] \label{claim:natural_injection_church_rosser} For $i \in \Set{1,2}$, let $\lambda_i : \scB_i \into \scC : X \into \SI$ be injections. Then, in the category of singular interval families (as objects) and injections (as morphisms), the pullback
\begin{equation} 
\xymatrix{ \scA \pullback \ar[r]^{\mu_1} \ar[d]_{\mu_2} & \scB_1 \ar[d]^{\lambda_1} \\
\scB_2 \ar[r]_{\lambda_2} & \scC}
\end{equation}
exists. The injection $\eps : \scA \into \scC$, defined by 
\begin{equation}
\eps = \lambda_1 \mu_1 = \lambda_2 \mu_2
\end{equation}
is called the \emph{intersection} of $\lambda_1, \lambda_2$, with \emph{factorisations} $\mu_1, \mu_2$.

\proof[Proof of \autoref{claim:natural_injection_church_rosser}] Write $\cS^{\lambda_1} \cap \cS^{\lambda_1}$ for the  singular subset section of $\scC$ with members $\cS^{\lambda_1}_x \cap \cS^{\lambda_1}_x$, $x \in X$. We claim this singular subset section is stable. For $(x \to y) \in \mor(X)$ we need to verify \eqref{eq:natural_injection_condition}, that is
\begin{equation} \label{eq:natural_injection_eps_condition}
\scC (x \to y)(\cS^{\lambda_1}_x \cap \cS^{\lambda_2}_x) \subset (\cS^{\lambda_1}_y \cap \cS^{\lambda_2}_y)
\end{equation}
But this follows since $\scC (x \to y)(\cS^{\lambda_i}_x) \subset \cS^{\lambda_i}_y$ by \eqref{eq:natural_injection_condition} individually for each $\lambda_i$, $i \in \Set{1,2}$. Using \autoref{defn:stability_vs_injections}, we can then define an injection
\begin{equation}
\eps := \eta_{\cS^{\lambda_1} \cap \cS^{\lambda_2}} : \intrel{\cS^{\lambda_1} \cap \cS^{\lambda_2}} \into \scC 
\end{equation}

Next, we define the legs $\mu_i$ as the unique factorising injections of $\eps$ and $\lambda_i$, which were constructed in claim \autoref{claim:injection_factorisation}. That is, we define $\cS^{\mu_i} := \lambda\inv_i (\cS^\eps)$.

Finally, Given any other commutative square
\begin{equation}
\xymatrix{ \scA' \ar[r]^{\mu'_1} \ar[d]_{\mu'_2} \ar[dr]^{\eps'} & \scB_1 \ar[d]^{\lambda_1} \\
\scB_2 \ar[r]_{\lambda_2} & \scC}
\end{equation}
then we know by \autoref{claim:injection_factorisation} that $\cS^{\eps'} \subset \cS^{\lambda_i}$ and thus $\cS^{\eps'} \subset \cS^{\lambda_1} \cap \cS^{\lambda_2}$. In particular, \autoref{claim:injection_factorisation} then implies that $\eps$ uniquely factors through $\eps'$ as $\eps = \eps' \alpha$ for some $\alpha : \scA \into \scA'$. In turn, we obtain $\mu'_i \alpha = \mu_i$ since both are the unique factorisation of $\eps$ through $\lambda_i$. This gives the universal property of the pullback as claimed. \qed
\end{lem}

\begin{eg}[Pullback of injections] \label{eg:pullback_of_injections} Using the definitions in \autoref{eg:injections} we find the stable singular subset section $\cS^{\lambda_1}$ associate to $\lambda_1$ to be
\begin{restoretext}
\begingroup\sbox0{\includegraphics{test/page1.png}}\includegraphics[clip,trim=0 {.15\ht0} 0 {.15\ht0} ,width=\textwidth]{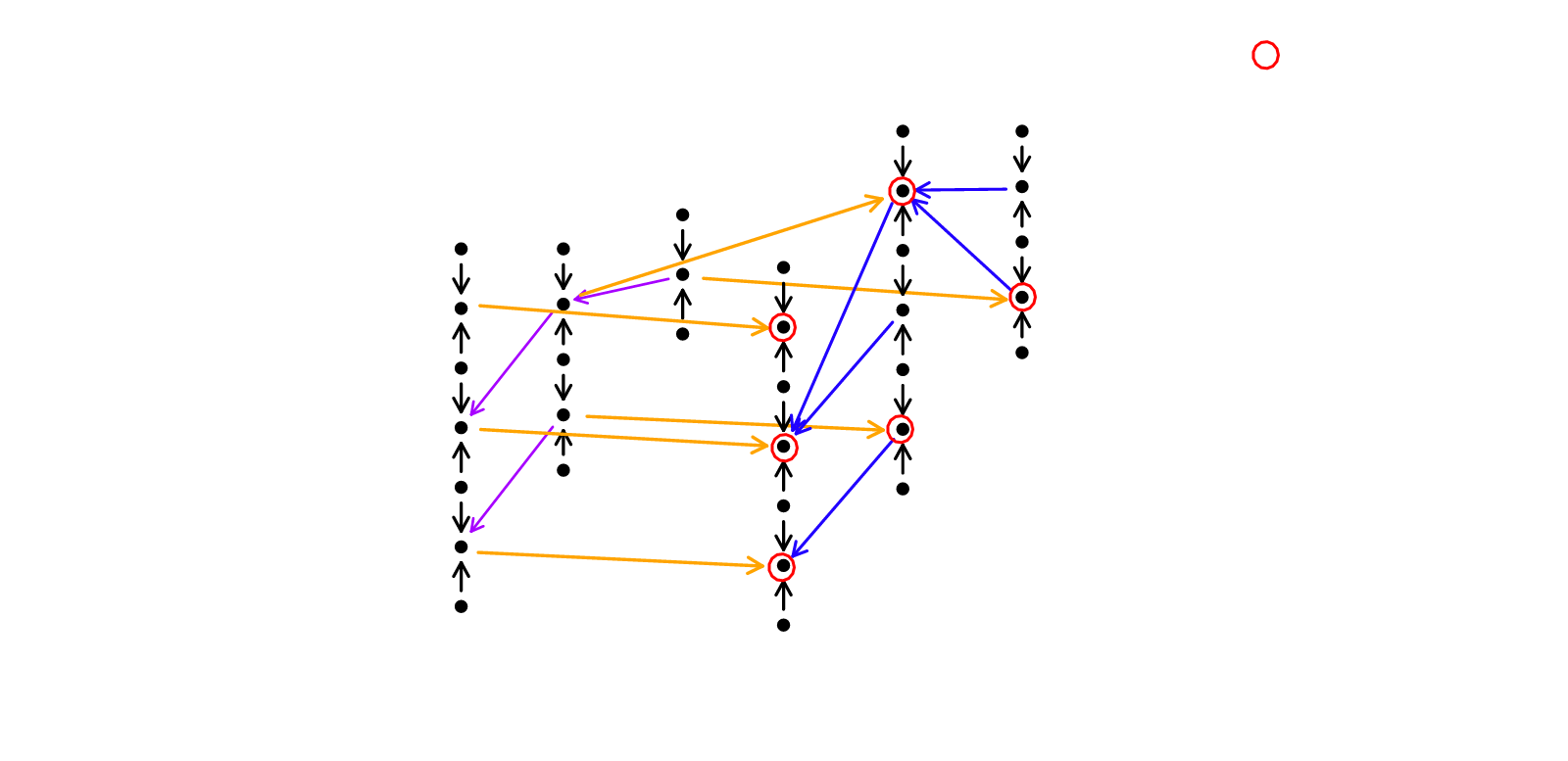}
\endgroup\end{restoretext}
(Elements of $\cS^{\lambda_1}$ are marked by \cred{} circles), and the stable singular subset section $\cS^{\lambda_2}$ associate to $\lambda_2$ to be
\begin{restoretext}
\begingroup\sbox0{\includegraphics{test/page1.png}}\includegraphics[clip,trim=0 {.15\ht0} 0 {.13\ht0} ,width=\textwidth]{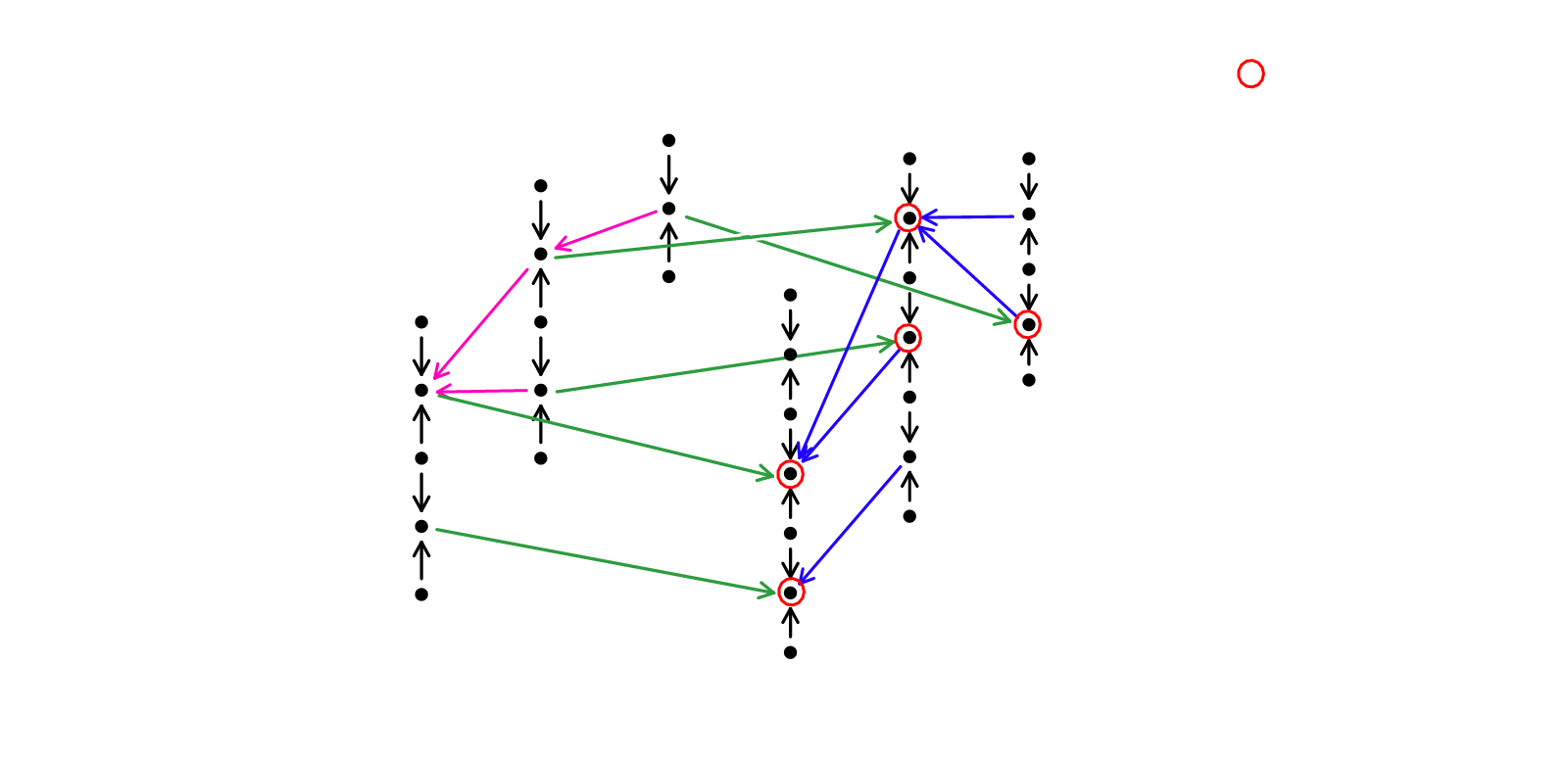}
\endgroup\end{restoretext}
We observe that their intersection is exactly $\cS^{\eps}$ (cf. \autoref{rmk:im_of_inj_stable}, or the \cred{} circles below). Further, using \autoref{eg:factorisation_of_injections} $\eps$ factors as $\lambda_i \mu_i$, and thus we obtain 
\begin{restoretext}
\begingroup\sbox0{\includegraphics{test/page1.png}}\includegraphics[clip,trim=0 {.0\ht0} 0 {.0\ht0} ,width=\textwidth]{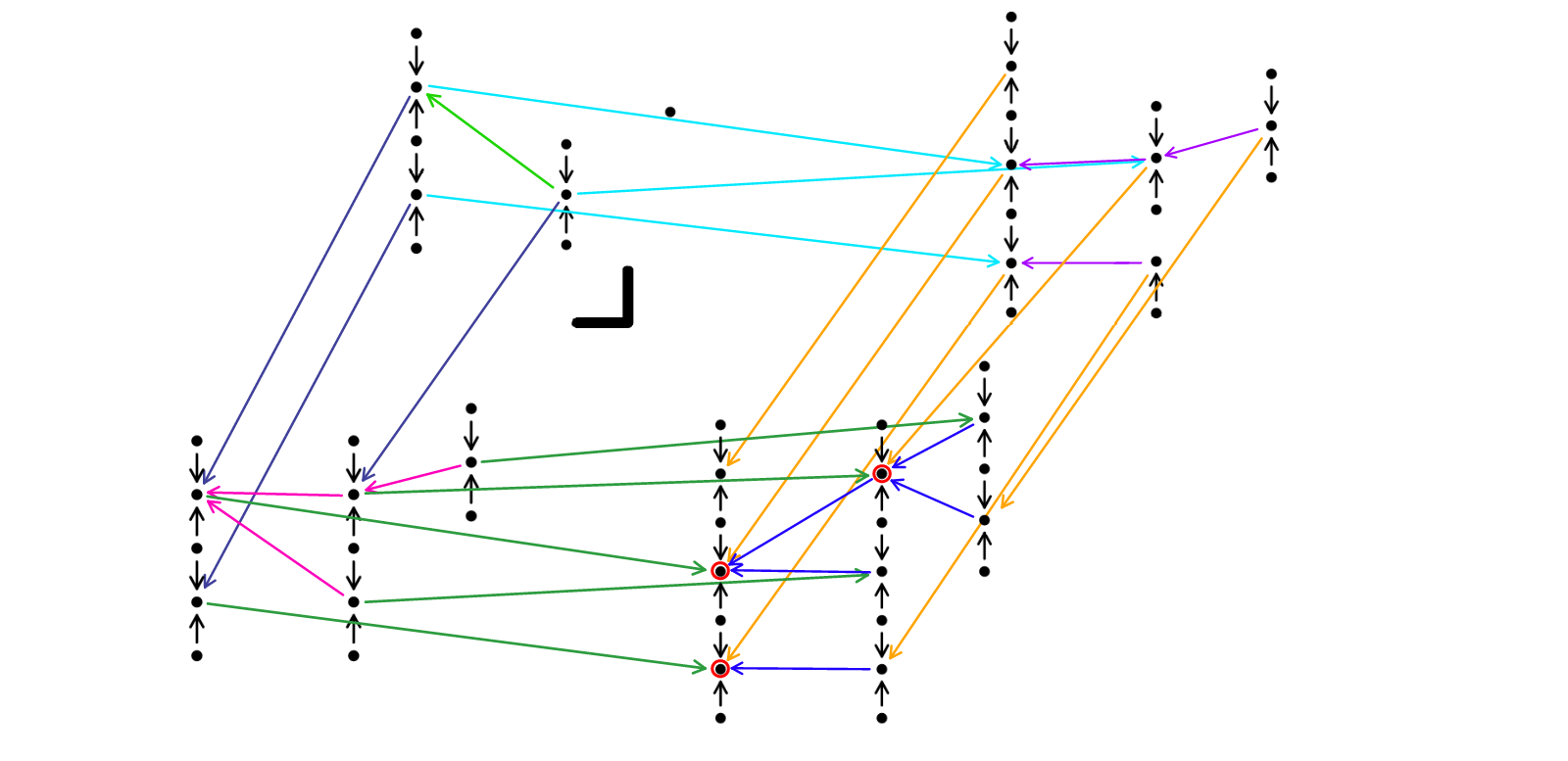}
\endgroup\end{restoretext}
This is an example of a pullback as constructed in the previous Lemma.
\end{eg}

\subsection{Open functors and collapse functors}

In this section we will turn our attention to a simple class of functors on singular intervals called ``open" functors. Topologically these can be interpreted as \textit{open} maps of intervals. Open functors specialise to collapse functors (and in \autoref{ch:emb}, to embedding functors). They will be used as a fibrewise description for family collapse functors in the next section. The goal of this section is to further establish a correspondence of collapse functors and their ``underlying monomorphisms".

\begin{defn}[Open functors of singular intervals] \label{defn:open_maps} A monotone functor $f : I \to J$ between singular intervals $I,J \in \SI$ is said to be \textit{open} if it preserves regular segments. That is,
\begin{equation}
a \in \regcont(I) \quad \im \quad f(a) \in \regcont(J)
\end{equation}
\end{defn}

\begin{eg}[Open functors] \label{eg:open_maps} \hfill
\begin{enumerate}
\item An example of an open functor $f_0 : \singint 5 \to \singint 4$ is the functor of posets defined by
\begin{restoretext}
\begingroup\sbox0{\includegraphics{test/page1.png}}\includegraphics[clip,trim=0 {.1\ht0} 0 {.1\ht0} ,width=\textwidth]{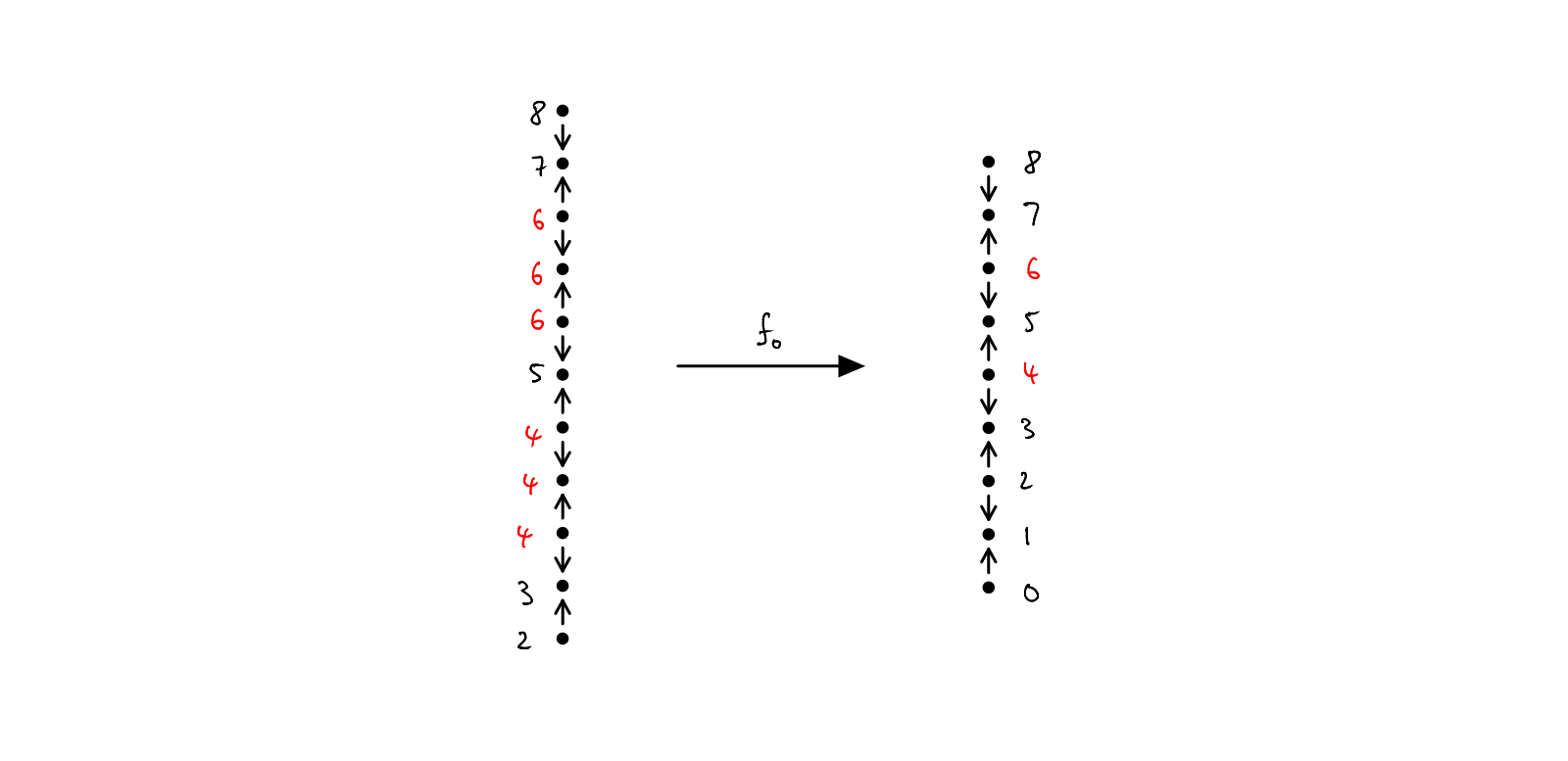}
\endgroup\end{restoretext} 
Here, we labelled the preimages of $f_0$ with the same labels as its image points. We highlight non-injective image points (and their preimages) in \cred{}. Note that labelling preimages is in line with \autoref{notn:depicting_posets} (while this depiction might seem cumbersome at this stage, it will allow for more readable notation later on).
\item An example of a monotone ``open" function $f_2 : \singint 2 \to \singint 1$ is the function defined by
\begin{restoretext}
\begingroup\sbox0{\includegraphics{test/page1.png}}\includegraphics[clip,trim=0 {.1\ht0} 0 {.1\ht0} ,width=\textwidth]{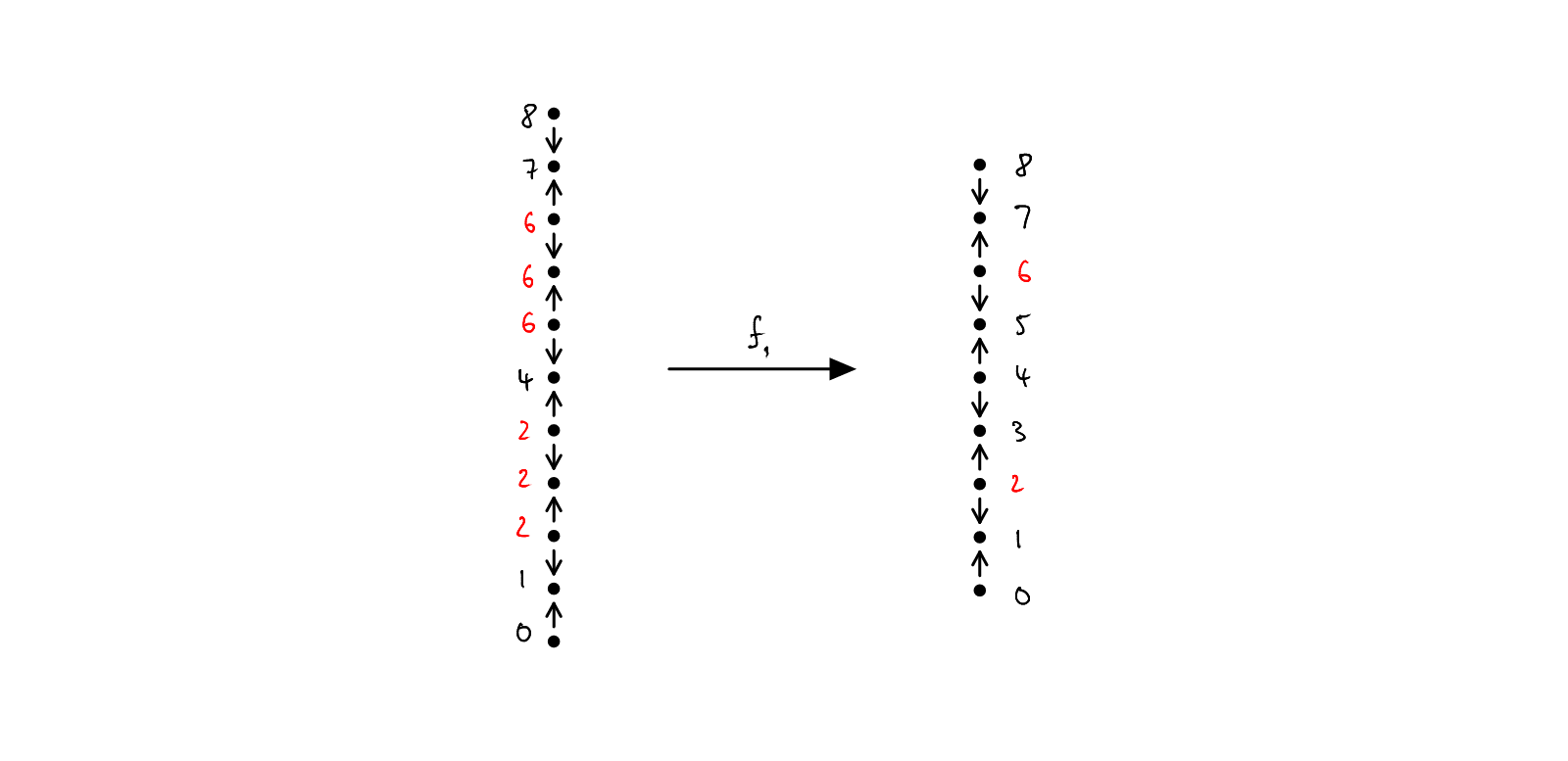}
\endgroup\end{restoretext} 
$f_1$ preserves regular segments as a function, but is not functorial. 
\item An example of an monotone functor $f_2 : \singint 2 \to \singint 1$ which is not open, is the functor of posets defined by
\begin{restoretext}
\begingroup\sbox0{\includegraphics{test/page1.png}}\includegraphics[clip,trim=0 {.3\ht0} 0 {.25\ht0} ,width=\textwidth]{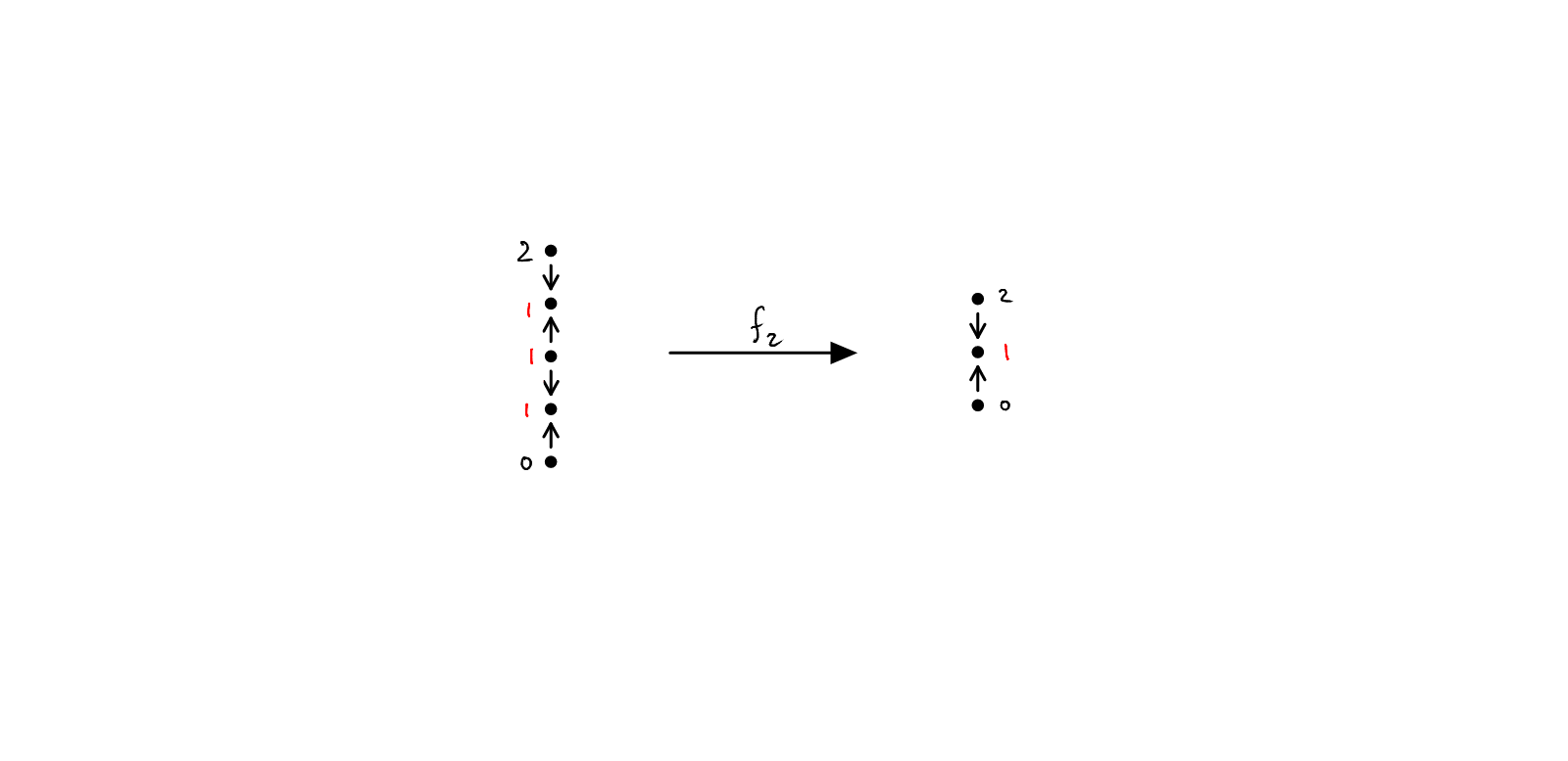}
\endgroup\end{restoretext} 
$f_2$ is functorial but doesn't preserve regular segments.
\end{enumerate}
\end{eg}

\begin{defn}[Collapse maps of singular intervals] \label{defn:singular_collapse_map} Given two singular intervals $I, J \in \SI$, a \textit{collapse functor} $s$ is a functor $s : J \to I$ that is open and surjective. 
\end{defn}

We remark that dually to the preceding definition, for the theory of embeddings, we will later on also study functors which are open an injective. These functors will be called embedding functors.

\begin{eg}[Collapse functors] \label{eg:collapse_map} $f_0$ from \autoref{eg:open_maps} is not a collapse functor since it is not surjective. However, the map $f'_0 : \singint 4 \to \singint 3$ defined by setting 
\begin{restoretext}
\begingroup\sbox0{\includegraphics{test/page1.png}}\includegraphics[clip,trim=0 {.0\ht0} 0 {.1\ht0} ,width=\textwidth]{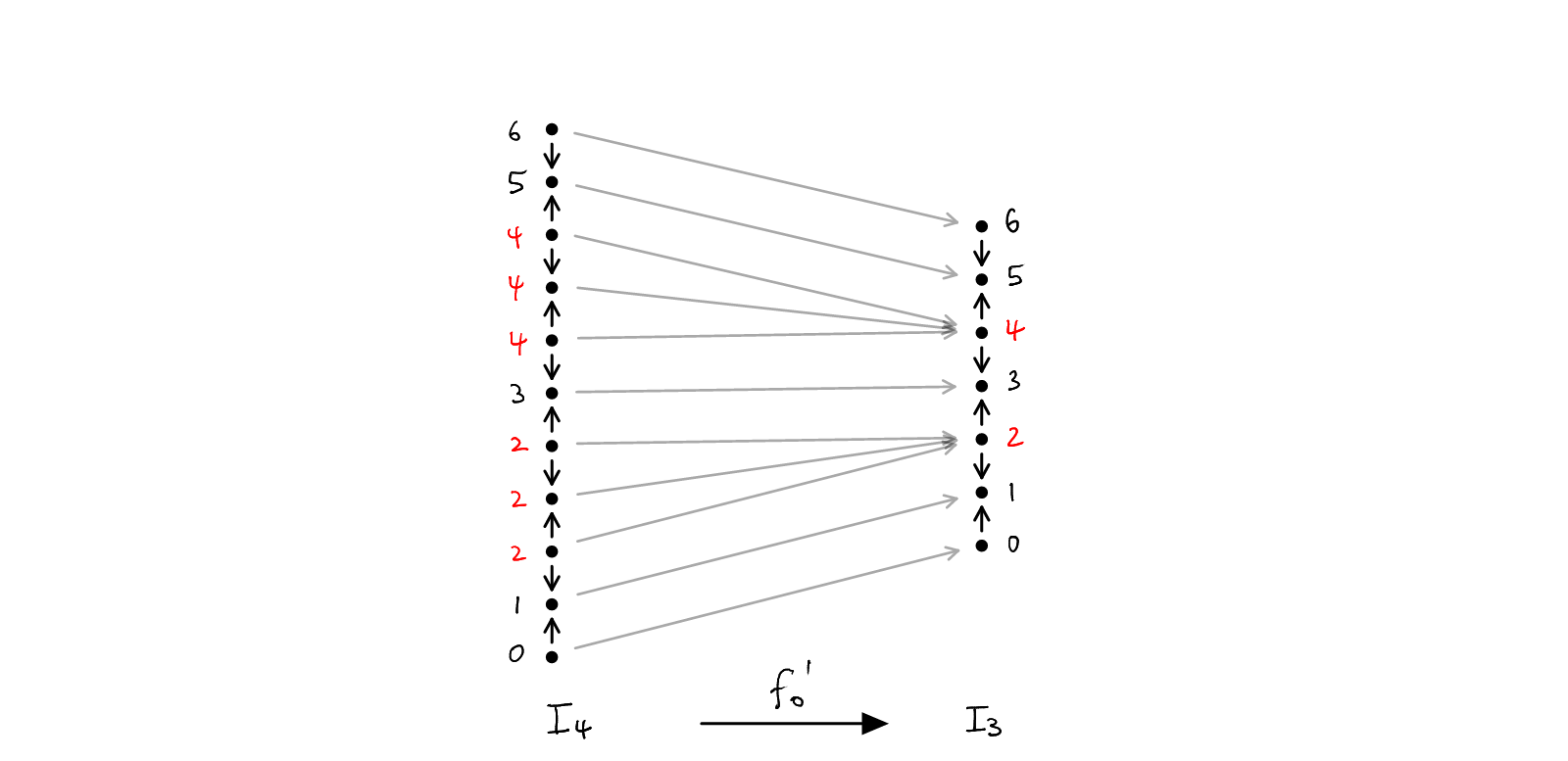}
\endgroup\end{restoretext}
is a collapse functor. For clarity, in addition to labelling preimages, we have also depicted the graph $\grph_{f'_0}$ of $f'_0$ (cf. \autoref{defn:graph_of_functors})  by \cgray{} arrows.
\end{eg}

As stated in the beginning of the section, to each collapse functor we will associate an ``underlying monomorphism". We start with the following example  of this fact.

\begin{eg}[Underlying monomorphism] The underlying monomorphism $(f'_0)\inv\monosing : \singcont(\singint 3) \to \singcont(\singint 4)$ for $f'_0$ from \autoref{eg:collapse_map} is given by
\begin{restoretext}
\begingroup\sbox0{\includegraphics{test/page1.png}}\includegraphics[clip,trim=0 {.0\ht0} 0 {.1\ht0} ,width=\textwidth]{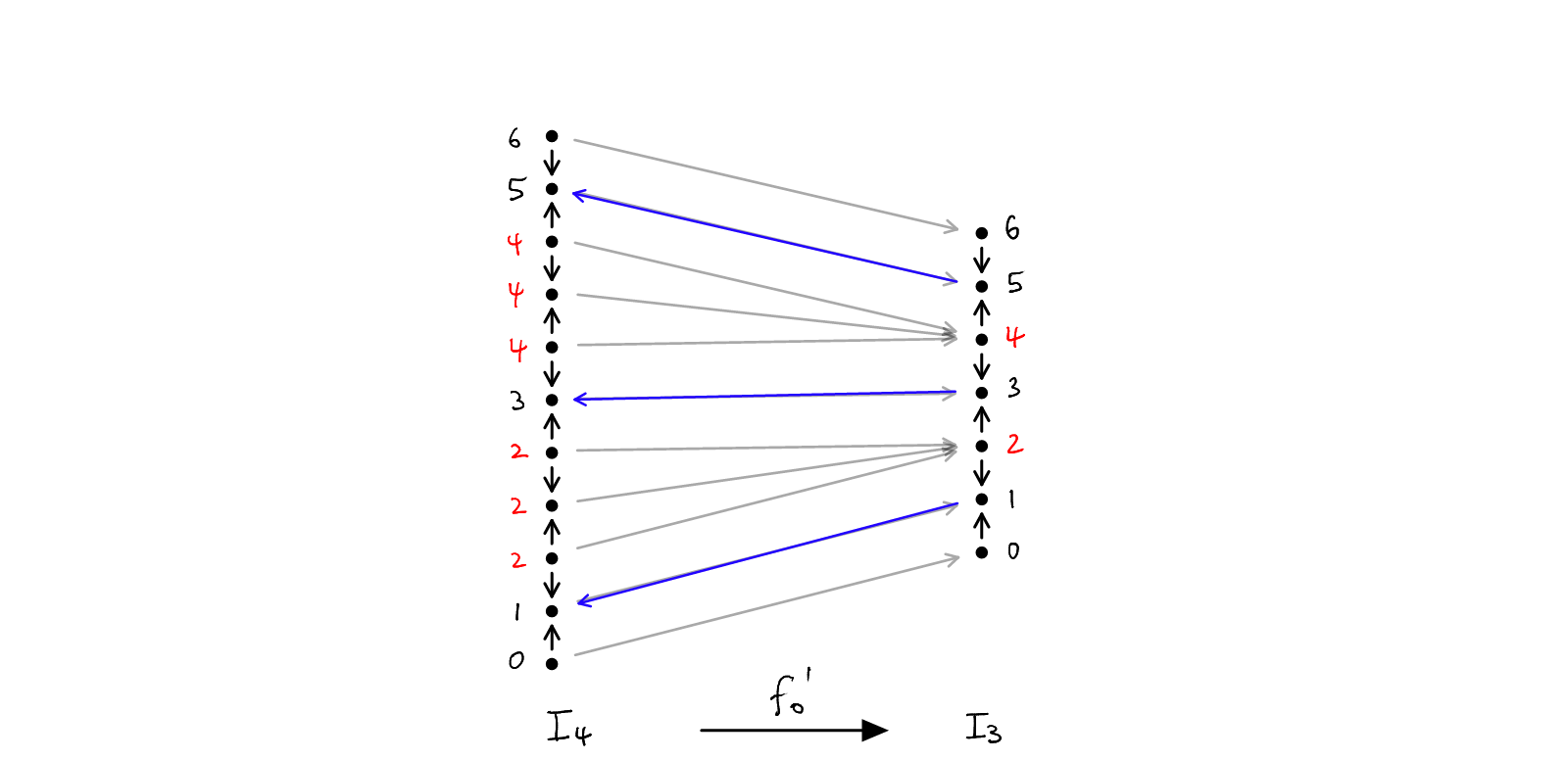}
\endgroup\end{restoretext}
Here, \cblue{} arrows indicate the mapping of the underlying monomorphism, which in particular goes in the opposite direction than $f'_0$.
\end{eg}

The next construction defines underlying monomorphisms more formally.

\begin{constr}[Underlying morphism of collapse functors] \label{constr:underlying_monomorphism_fctr} Let $s$ be a collapse functor. This implies that the graph $\grph_s$ of $s$ (cf. \autoref{defn:graph_of_functors}) has an converse relation $\grph\converse_s$ that restricts to give the graph of a function $s\inv\monosing$ on singular heights. In other words,
\begin{equation}
\grph_{s\inv\monosing} = \rest {\grph\converse_s} {\singcont(I) \times \singcont{J}}
\end{equation}
Indeed arguing by contradiction first assume that $\rest {\grph\converse_s} {\singcont(I) \times J}$ was a non-functional relation. By surjectivity of $s$, this means there is $b, b' \in I$ such that $b < b'$ and $s(b) =s(b') = a \in \singcont(J)$. There is at least one $c \in [b,b']$ such that $c \in \regcont(I)$. Monotonicity of $s$ implies $s(c) = a$ which contradicts openness of $s$. Openness similarly implies that the functional relation yields a function with image in $\singcont(J)$, which then proves the equality stated above.

This guarantees existence of $s\inv$, which is necessarily monomorphic and called the \emph{underlying monomorphism of $s$}. Note that $s\inv\monosing : I \to_{\SI} J$. By openness we further have $\rest s {\regcont(J)} : J \to\regop_{\SI} I$.
\end{constr}

The next claim establishes that any open functor with ``underlying monomorphism" must be a collapse functor, and that the association of underlying monomorphisms to collapse functors is injective.

\begin{claim}[Underlying monomorphisms determine collapse functors] \label{claim:underlying_mono} Let $s : J \to I$ be an open functor such that
\begin{equation}
\grph_{t} = \rest {\grph\converse_s} {\singcont(I) \times \singcont{J}}
\end{equation}
for some (necessarily monotone)
\begin{equation}
t  : \singcont(I) \to \singcont(J)
\end{equation}
Then $s$ is a collapse functor, and determined by $t$ as follows
\begin{alignat}{3} \label{eq:glambda_defn}
&b \in \singcont(I) &&\text{~then}\quad & \forall a \in J. ~\big( s(a) = b &\iff t(b) = a\big) \\
&b \in \regcont(I) &&\text{~then}&  \forall a \in J.  ~\big(s(a) = b  &\iff (\wwidehat t(b - 1) < a < \wwidehat t(b+1))\big)
\end{alignat}

\proof The proofs are \stfwd{}. The first statement of \eqref{eq:glambda_defn} follows from the assumed property of $t$ and openness of $s$. 

For the second statement of \eqref{eq:glambda_defn} we first proof $\imp$. Assume $s(a) = b$ and $b \in \regcont(J)$, $a \in I$. We prove $\wwidehat{t}(b-1) < a$. If $b = 0$, then $\wwidehat{t}(b-1) < a$ follows from the boundary cases in \autoref{notn:singular_morphism_boundary_cases}. If $b > 0$, then $s(\wwidehat{t}(b-1)) = b-1$. If $a < \wwidehat{t}(b-1)$ this would contradict monotonicity of $s$. Thus we must have $\wwidehat{t}(b-1) < a$. Similarly, we find $a < \wwidehat{t}(b+1)$ and thus the implication $\imp$ follows. 

Next we prove $\limp$. Assume $(\wwidehat t(b - 1) < a < \wwidehat t(b+1))$ for $b \in \regcont(J)$, $a \in I$. Using monotonicity of $s$, the first statement of \eqref{eq:glambda_defn} as well as the boundary cases in \autoref{notn:singular_morphism_boundary_cases} we find $b-1 \leq s(a) \leq b + 1$. Since $b\pm 1$ are singular heights, $s(a) = b\pm 1$ the first statement of \eqref{eq:glambda_defn} would then imply $a = t(b\pm 1)$ contradicting our assumption. Thus we must have $s(a) = b$.

Finally, note that $t \regop = \rest s {\regcont(J)}$. This follows from comparing \eqref{eq:glambda_defn} with \eqref{eq:fop_from_f}. But $t\regop$ is endpoint-preserving. Thus $0$ and $2\iH_J$ are in the image of $s$. Since the domain of $s$ is connected, so is its image and thus $s$ must be surjective, it is a collapse functor. \qed
\end{claim}

The previous claim shows that there is a (possibly partial) left inverse to the underlying monomorphism construction. The next claim shows that this inverse is in fact total, and two-sided.

\begin{claim}[Singular collapse functors from monomorphisms] \label{constr:sing_coll_map_from_mono} Given two singular intervals $I, J \in \SI$, and a monomorphism $t : I \to_{\SI} J$, we claim that $s$ defined from $t$ by \eqref{eq:glambda_defn} is a collapse functor. This collapse functor will be denote by
\begin{equation}
\sS^t := s
\end{equation}
The underlying monomorphism of $\sS^t$ is $t$.
\proof The proof is \stfwd{}. We check well-definedness, monotonicity, preservation of regular segments, surjectivity and functoriality.
\begin{itemize}
\item Firstly, by \eqref{eq:regular_covering} the definition \eqref{eq:glambda_defn} indeed covers all choices for $a \in J$ unambiguously and thus defines a function of sets $\sS^t : J \to I$. 
\item Next we check monotonicity of $\sS^t$. Note that both statements of \eqref{eq:glambda_defn} can be seen to imply one of \eqref{eq:defn_order_realisation_1}, \eqref{eq:defn_order_realisation_2} defining $R(t)(a,b)$. In other words, we observe that $\sS^t (a) = b$ implies $R(t)(a,b)$. Now, assume $a_1 < a_2 \in J$ and $\sS^t (a_i) = b_i$, $i \in \Set{1,2}$. Then $R(t)(a_i,b_i)$ holds and thus bimonotonicity of $R(t)$ (cf. \autoref{claim:order_realisations_monotone}) implies $b_1 \leq b_2$ as required for monotonicity of $\sS^t$. 
\item We check that $\sS^t$ preserves regular segments. Since $\sS^t(a) = b$ has $b \in \singcont(I)$ only if $a = t(b) \in \singcont(J)$ (by the first condition in \eqref{eq:glambda_defn}) we find that $\sS^t$ preserves regular segments. %
\item Surjectivity of $s$ follows from the conditions in \eqref{eq:glambda_defn} and $t$ being monomorphic.
\item Functoriality now follows from combining the previous three items.
\end{itemize}
Thus, we deduce that $\sS^t$ is indeed a collapse functor, and that (by the first condition in \eqref{eq:glambda_defn}) it has underlying monomorphism $t$. \qed
\end{claim}

We record the result of the preceeding two claims and one construction in the following remark.

\begin{rmk}[Collapse functors are in bijection with monomorphisms] \label{rmk:collapse_vs_mono} The construction $(t \mapsto \sS^t)$ establishes a bijection from monomorphisms $I \to_{\SI} J$ to singular collapse functors $J \to I$. This follows from \autoref{claim:underlying_mono} and \autoref{constr:sing_coll_map_from_mono}. The inverse operation to $(t \mapsto \sS^t)$ is taking underlying monomorphisms (cf. \autoref{constr:underlying_monomorphism_fctr}).
\end{rmk}

\subsection{Collapse of \SI-families}

We now define the notion of family collapse functor on total spaces of singular interval families. Later in this chapter this will be used to define collapse of \textit{labelled} singular interval families. Apart from defining collapse functors on families and giving examples, the goal of this section is to understand that collapse maps are in correspondence with injections. 

The reason for introducing this correspondence is the following: injections are the ``minimal data" for collapse functors. This makes injections sometimes easier to work with in proofs.

\begin{defn}[Family collapse functors] \label{defn:family_collapse_map} Given two singular interval families $\scA, \scB : X \to \SI$, a \textit{(family) collapse functor} $S$ of $\scB$ to $\scA$ is a map of bundles $S : \pi_\scB \to \pi_\scA$ such that that $\rest S x : \scB(x) \to \scA(x)$ is a collapse functor of intervals for each $x \in X$.
\end{defn}

The following two constructions show that family collapse functors are fully and faithfully described by injections.

\begin{constr}[Underlying injections associated to family collapse functors] \label{defn:underlying_injections}
Every collapse functor $S : \pi_\scB \to \pi_\scA$ has an \textit{underlying injection} $\lambda : \scA \into \scB$ defined componentwise by
\begin{equation}
\lambda_x = (\rest S x)\inv\monosing
\end{equation}
that is, $\lambda_x$ is the underlying monomorphism of $\rest S x$. To check naturality of $\lambda$ we need to verify that for any morphism $x_1 \to x_2$ in $X$ and $a \in \singcont(\scA(x_1))$, we have
\begin{equation}
\lambda_{x_2}\scA(x_1 \to x_2)(a_1) = \scB(x_1 \to x_2)\lambda_{x_1}(a_1)
\end{equation}
Setting $a_2 := \scA(x_1 \to x_2)(a_1)$, this follows from functoriality of $S$ together with definition of $\lambda$ (which implies $S((x_1,b_1) \to (x_2,b_2)) = (x_1,a_1) \to (x_1,a_2)$ where $b_i := \lambda_{x_i}(a_i)$, $i \in \Set{1,2}$) and from \eqref{eq:defn_order_realisation_1}.
\end{constr}

\begin{constr}[Family collapse functors associated to injections]\label{defn:glambda} Let $\lambda : \scB \to \scA$ be an injection for singular interval families $\scA, \scB : X \to \SI$. We define a collapse functor of $\scB$ to $\scA$ 
\begin{equation}
\sS^\lambda : \pi_\scB \to \pi_\scA
\end{equation}
called the family collapse functor associated to $\lambda$. $\sS^\lambda$ is defined by setting its fibrewise maps (over $x \in X$) to equal
\begin{equation}
\rest {\sS^\lambda} x  := \sS^{\lambda_x} : \scB(x) \to \scA(x)
\end{equation}
where we used \autoref{constr:sing_coll_map_from_mono} for the monomorphism $\lambda_x$. This defines a function of sets $\sG(\scA) \to \sG(\scB)$. To show that $\sS^\lambda$ given in \autoref{defn:glambda} defines a map of bundles we need to check that $\sS^\lambda$ functorial and preserves fibers. First note that $\sS^\lambda$ as given in \eqref{eq:glambda_defn} is defined fibrewise for each $x \in X$ and thus satisfies $\pi_\scB \sS^\lambda = \pi_\scA$. We are left with showing functoriality: let $(x_1, a_1) \to (x_2, a_2)$ be a morphism in $\sG(\scA)$. We need to show there is a morphism $\sS^\lambda (x_1, a_1) \to \sS^\lambda (x_2, a_2)$ in $\sG(B)$. We distinguish two cases:
\begin{enumerate}
\item In the first case assume $a_1 \notin \cS^\lambda_{x_1}$. Then by \eqref{eq:glambda_defn} there is $b_1 \in \regcont(\scB(x_1))$ such that 
\begin{equation} \label{eq:Glambda_case_1}
\wwidehat\lambda_{x_1}(b_1 - 1) < a_1 < \wwidehat\lambda_{x_1}(b_1 + 1)
\end{equation} 
and thus $\sS^\lambda(x_1,a_1) = (x_1, b_1)$ by \eqref{eq:glambda_defn}. By naturality of $\lambda$ we find
\begin{equation} \label{eq:Glambda_case_1b}
\wwidehat \scA(x_1 \to x_2)\wwidehat\lambda_{x_1}(b_1 \pm 1) = \wwidehat\lambda_{x_2} \wwidehat \scB(x_1 \to x_2) (b_1 \pm 1)
\end{equation} 
By \autoref{claim:order_realisations_monotone}, i.e. monotonicity of $R \scA(x_1 \to x_2)$, applied to \eqref{eq:Glambda_case_1} and using our assumption $(a_1,a_2) \in \edgeset(\scA(x_1 \to x_2))$, we deduce that
\begin{equation}
\wwidehat \scA(x_1 \to x_2) \wwidehat\lambda_{x_1}(b_1 - 1)\leq a_2 \leq \wwidehat \scA(x_1 \to x_2) \wwidehat \lambda_{x_1}(b_1 + 1)
\end{equation}
and thus by \eqref{eq:Glambda_case_1b}
\begin{equation}
\wwidehat\lambda_{x_2} \wwidehat \scB(x_1 \to x_2) (b_1 - 1) \leq a_2 \leq \wwidehat\lambda_{x_2} \wwidehat \scB(x_1 \to x_2) (b_1 + 1)
\end{equation}
Setting $(x_2, b_2) := \sS^\lambda(x_2, a_2)$ and inspecting both cases of \eqref{eq:glambda_defn} for $a_2$ (namely, either $\lambda_{x_2}(b_2) = a_2$ or $\lambda_{x_2}(b_2 -1) < a_2 < \lambda_{x_2}(b_2 + 1)$), by monotonicity of $\lambda_{x_2}$ this necessarily implies
\begin{equation}
\wwidehat \scB(x_1 \to x_2) (b_1 - 1) \leq b_2 \leq \wwidehat \scB(x_1 \to x_2) (b_1 + 1)
\end{equation}
By \eqref{eq:defn_order_realisation_2} we infer $R \scB(x_1 \to x_2) (b_1, b_2)$, i.e. $(b_1,b_2) \in \edgeset(\scB(x_1 \to x_2))$. Since $(x_i,b_i) = \sS^\lambda(x_i,a_i)$ this implies $\sS^\lambda (x_1, a_1) \to \sS^\lambda (x_2, a_2)$ as required.

\item In the second case assume $a_1 \in \cS^\lambda_{x_1}$. Then there is $b_1 \in \singcont(\scB(x_1))$ such that $\lambda_{x_1} (b_1) = a_1$ and thus $\sS^\lambda(x_1,a_1) = (x_1, b_1)$ by \eqref{eq:glambda_defn}. Since $\lambda$ is natural, we have
\begin{equation}
a_2 = \scA(x_1\to x_2)(\lambda_{x_1} (b_1)) = \lambda_{x_2}(\scB(x_1 \to x_2)(b_1))
\end{equation}
And thus setting $b_2 :=\scB(x_1 \to x_2)(b_1)$, by \eqref{eq:glambda_defn}, we have $\sS^\lambda(x_2,a_2) = (x_2,b_2)$. By \eqref{eq:defn_order_realisation_1} we have $R \scB(x_1 \to x_2) (b_1, b_2)$ as required.
\end{enumerate}
This completes the construction of the map of bundles $\sS^\lambda$ as claimed.
\end{constr}

\begin{eg}[family collapse functors associated to injections] \label{eq:family_collapse} \hfill
\begin{enumerate}
\item Recall $\lambda_1$ in \autoref{eg:injections} was defined to have components (indicated by \corange{} arrows)
\begin{restoretext}
\begingroup\sbox0{\includegraphics{test/page1.png}}\includegraphics[clip,trim=0 {.2\ht0} 0 {.1\ht0} ,width=\textwidth]{ANCimg/page71.png}
\endgroup\end{restoretext}
It's associated family collapse $\sS^{\lambda_1}$ is the mapping
\begin{restoretext}
\begingroup\sbox0{\includegraphics{test/page1.png}}\includegraphics[clip,trim={.1\ht0} {.0\ht0} {.3\ht0} {.0\ht0} ,width=.7\textwidth]{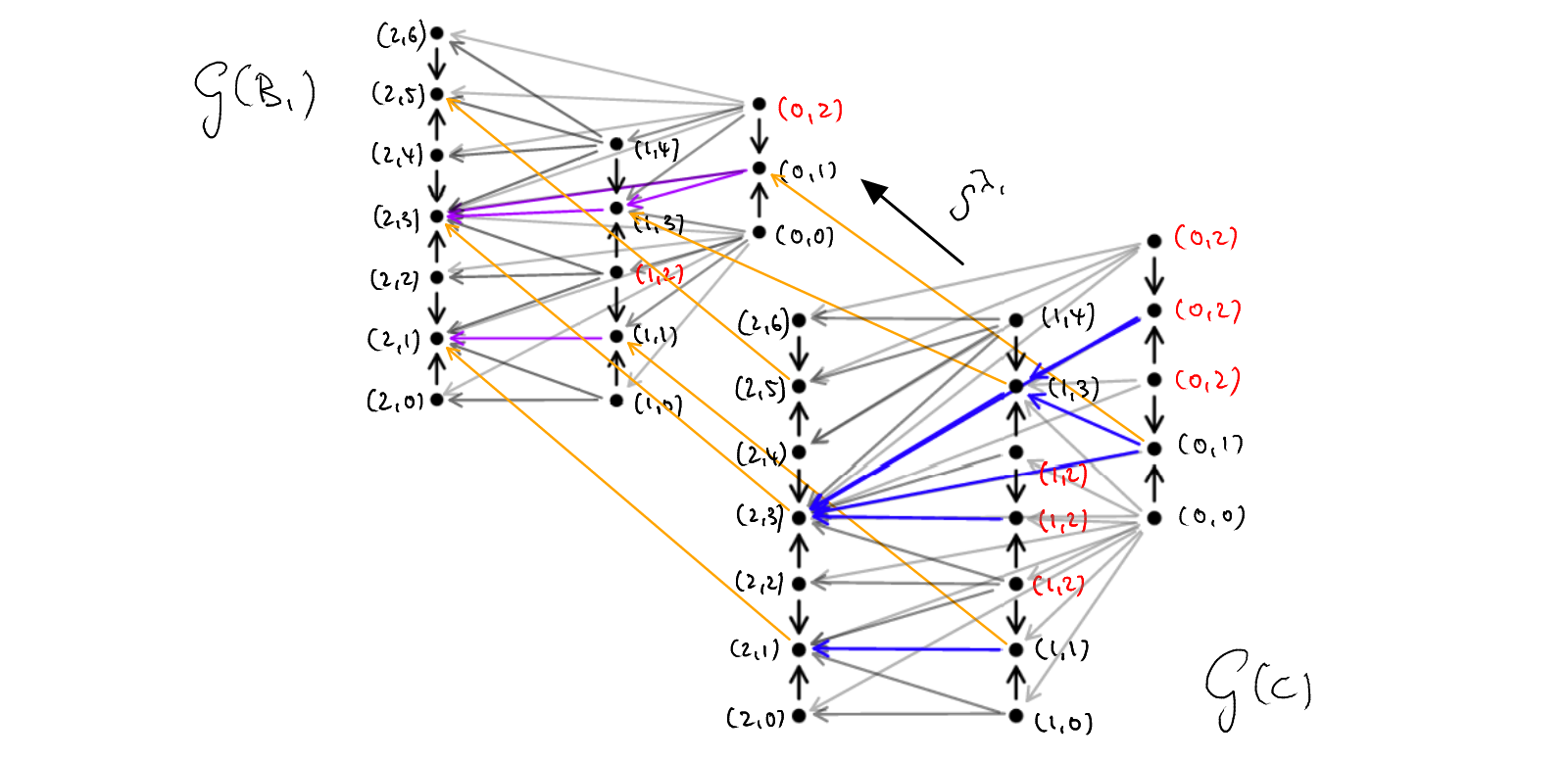}
\endgroup\end{restoretext}
Here, we labelled the preimages of $\sS^{\lambda_1}$ with the same labels as its image points. We highlight non-injective image points (and their preimages) in \cred{}. We also indicate by \cdarkgreen{} arrows the mapping of $\sS^{\lambda_1}$ on singular heights, which can be seen to invert the mapping of $\lambda_1$ in the previous picture. Arrows in \cblue{} and \cpurple{} correspond to \cblue{} and \cpurple{} arrows in the previous picture (the above picture has more arrows however, since it depicts the entire total space).

\item Recall $\lambda_2$ in \autoref{eg:injections} was defined to have components (indicated by \corange{} arrows)
\begin{restoretext}
\begingroup\sbox0{\includegraphics{test/page1.png}}\includegraphics[clip,trim=0 {.15\ht0} 0 {.15\ht0} ,width=\textwidth]{ANCimg/page73.png}
\endgroup\end{restoretext}
It's associated family collapse $\sS^{\lambda_2}$ is the mapping
\begin{restoretext}
\begingroup\sbox0{\includegraphics{test/page1.png}}\includegraphics[clip,trim={.1\ht0} {.0\ht0} {.3\ht0} {.0\ht0} ,width=.7\textwidth]{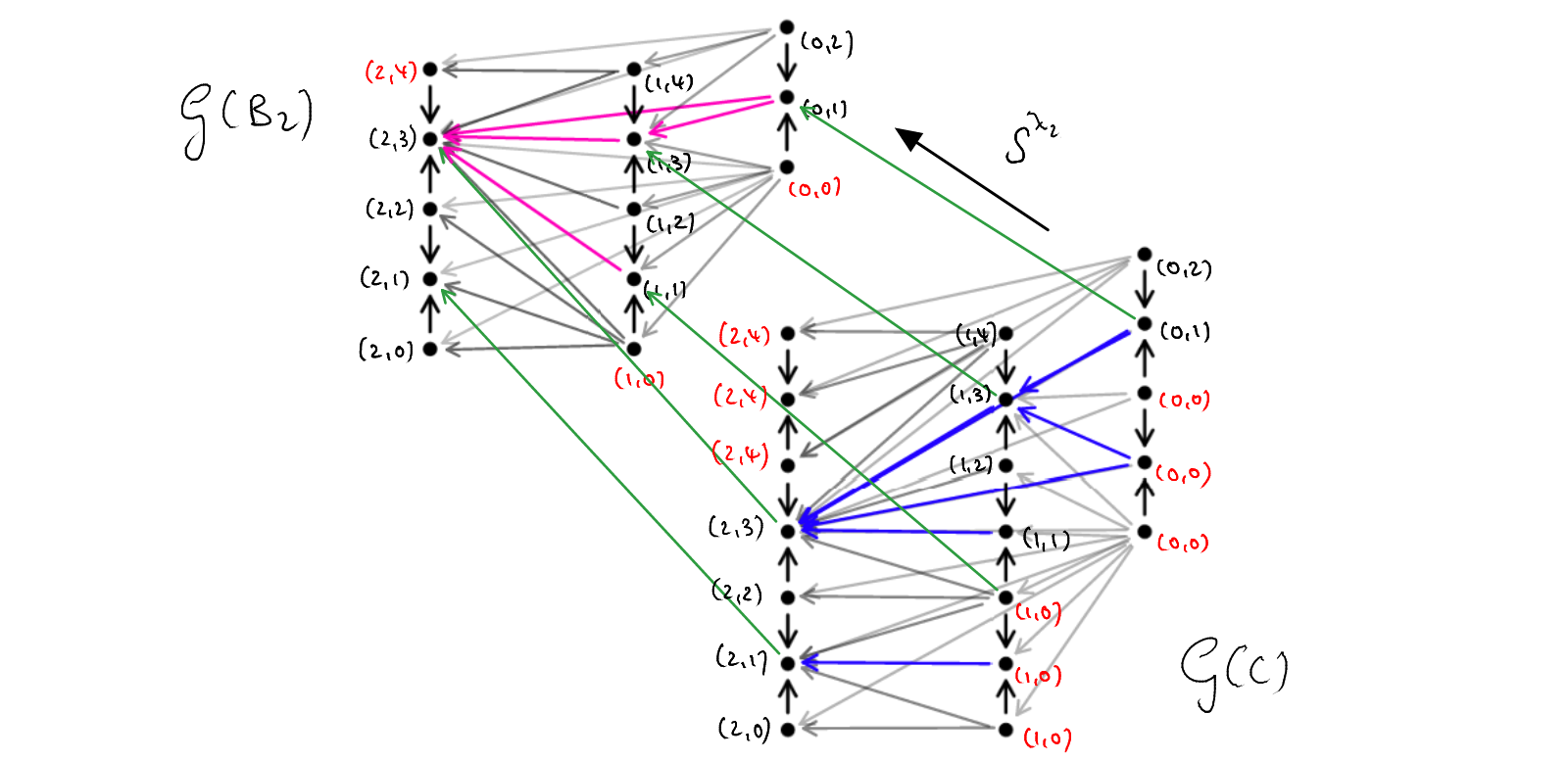}
\endgroup\end{restoretext}
Here, again we labelled the preimages of $\sS^{\lambda_2}$ with the same labels as its image points. We highlight non-injective image points (and their preimages) in \cred{}. We also indicate by \corange{} arrows the mapping of $\sS^{\lambda_2}$ on singular heights, which can be seen to invert the mapping of $\lambda_2$ in the previous picture. Arrows in \cblue{} and \cpink{} correspond to \cblue{} and \cpink{} arrows in the previous picture (the above picture has more arrows however, since it depicts the entire total space).
\end{enumerate}
\end{eg}

\noindent We record the correspondence of injections and family collapse functors in the following theorem.

\begin{thm}[Injections correspond to family collapse functors] \label{thm:collapse_maps_vs_injections} Let $\scA, \scB : X \to \SI$. The construction $(\lambda \mapsto \sS^\lambda)$ establishes a bijection from injections $\scA \into \scB$ to family collapse functors $\pi_\scB \to \pi_\scA$.
\proof The inverse mapping is given by taking underlying injections as constructed in \autoref{defn:underlying_injections}. The fact that this is an inverse follows (arguing fibrewise) from \autoref{rmk:collapse_vs_mono}. \qed
\end{thm}

\begin{cor}[Uniqueness of collapse functors] \label{lem:glambda} Let $\lambda : \scA \into \scB$. If $F : \pi_\scB \to \pi_\scA$ satisfies (for all $x \in X$)
\begin{enumerate}
\item $\rest {F} x$ is open %
\item $\rest {F} x$ inverts $\lambda_x$, i.e.\ $\rest {F} x (\lambda_x(a)) = a$ ($a \in \singcont(\scB(x))$)
\end{enumerate}
then $F = \sS^\lambda$. 
\proof The condition on $\rest {F} x$ imply that it is a collapse functor of singular intervals $\scB(x) \to \scA(x)$ with underlying monomorphism $\lambda_x$ (cf. \autoref{claim:underlying_mono}), and thus equals $\rest {\sS^\lambda} x$. \qed

\end{cor}

\begin{cor}[Functoriality of singular collapse functor assocation] \label{cor:composition_of_Glambda} Given injections $\mu : \scB_0 \into \scB_1$, $\lambda : \scB_1 \into \scB_2$ of singular interval families $\scB_0, \scB_1,\scB_2 : X \into \SI$, then we have
\begin{equation}
\sS^{\lambda\mu} = \sS^\mu\sS^{\lambda}
\end{equation}
Further, let $\id : \scB_0 \to \scB_0$ the identity transformation. Then 
\begin{equation}
\sS^{\id} = \id
\end{equation}
\proof Since $\sS^\mu\sS^{\lambda}$ (resp. $\id$) satisfies the properties stated in \autoref{lem:glambda} for the injections $\lambda\mu$ (resp. for $\id$), by \autoref{lem:glambda} they must equal $\sS^{\lambda\mu}$ (resp. $\sS^{\id}$).  \qed
\end{cor}

\subsection{Collapse has lifts}

The fact that collapse bundle maps have \gls{lifts} is an important observation, albeit with a somewhat tedious proof.

\begin{lem}[Collapse has lifts] \label{lem:collapse_has_lifts} The functor of posets $\sS^\lambda$ given in \autoref{defn:glambda} has \gls{lifts} (cf. \autoref{defn:having_lifts}). 

\proof To see that $\sS^\lambda$ has \gls{lifts} we need to verify that whenever $(r : x_\sop \to x_\tap)\in \mor(X)$, $(r,\edb) \in \mor(\sG(\scB))$ and $\sS^\lambda (x_\tap, a) = (x_\tap, \edb\ttae )$ then there is $(r,\eda) \in \mor(\sG(\scA))$ with $\eda\ttae  = a$ and $\sS^\lambda(r,\eda) = (r,\edb)$. %
We distinguish the following cases

\begin{enumerate}
\item Assume $\edb\ssoe  \in \regcont(\scB(x_\sop)), \edb\ttae  \in \regcont(\scB(x_\tap))$. Using \autoref{cor:relation_fullness} for $\scA(r)$, we can find $\eda \in \edgeset(\scA(r))$ with $\eda\ttae  = a$. Define $(r,\edb') = \sS^\lambda (r,\eda)$. We have $\edb'\ttae  = \edb\ttae $ since $\eda\ttae  = a$. We claim $\edb'\ssoe  = \edb\ssoe $: Indeed, by \eqref{eq:defn_order_realisation_3} both must equal $\edb'\ssoe  = \edb\ssoe  = \scB(r)\regop(\edb\ttae )$. Thus $(r,\eda)$ gets mapped to $(r,\edb)$ by $\sS^\lambda$ as required.

\item Assume $\edb\ssoe  \in \singcont(\scB(x_\sop))$, $\edb\ttae  \in \singcont(\scB(x_\tap))$. The first assumption implies $\scA(r)(\edb\ssoe ) = \edb\ttae $ by \eqref{eq:defn_order_realisation_1}. The last assumption implies $\lambda_{x_\tap}(\edb\ttae ) = a$ by \eqref{eq:glambda_defn}. We define $a_\sop = \lambda_{x_\sop}(\edb\ssoe )$. Then naturality of $\lambda$ implies 
\begin{align}
a &= \lambda_{x_\tap}(\scB(r)(\edb\ssoe )) \\
&= \scA(r)(\lambda_{x_\sop}(\edb\ssoe )) \\
&= \scA(r)(a_\sop)
\end{align}
We deduce $\scA(r)(a_\sop) = a$ and $\eda := (a_\sop,a) \in \edgeset(\scA(r))$ by \eqref{eq:defn_order_realisation_1}. Since  $\sS^\lambda$ inverts $\lambda$ we have $\sS^\lambda(x_\sop,\eda\ssoe ) = (x_\sop,\edb\ssoe )$. Thus $(r,\eda)$ gets mapped onto $(r,\edb)$ by $\sS^\lambda$ as required.

\item %
Assume $\edb\ssoe  \in \regcont(\scB(x_\sop)), \edb\ttae  \in \singcont(\scB(x_\tap))$ and $\edb\ttae  \in \im(\scB(r))$. The second assumption implies $\lambda_{x_\tap}(\edb\ttae ) = a$. The last assumptions implies there is $b_\sop \in \singcont(\scB(x_\sop))$ such that $\scB(r)(b_\sop) = \edb\ttae $. Assume $b_\sop < \edb\ssoe $ (in the case $b_\sop > \edb\ssoe $ the argument is similar). In particular, $\edb\ssoe  > 0$. Set $b_\tap := \scB(r)(\edb\ssoe -1)$. $b_\sop < \edb\ssoe $ implies $b_\sop \leq \edb\ssoe  - 1$. If $b_\sop = \edb\ssoe  - 1$ then $b_\tap = \edb\ttae $. Otherwise $x < \edb\ssoe  - 1 < \edb\ssoe $, and thus $(b_\sop,\edb\ttae ), (\edb\ssoe -1,b_\tap), (\edb\ssoe ,\edb\ttae ) \in \edgeset(\scB(r))$ together with monotonicity of $R\scB(r)$ (cf. \autoref{claim:order_realisations_monotone}) implies $\edb\ttae  \leq b_\tap \leq \edb\ttae $. That is, we again have $b_\tap = \edb\ttae $. We infer $\scB(r)(\edb\ssoe -1) = \edb\ttae $ in either case. Now, define $a_\sop = \lambda_{x_\sop} (\edb\ssoe  - 1) + 1$ (note $\edb\ssoe  - 1 \in \singcont(\scB(x_\sop))$ since $\edb\ssoe  > 0$). By naturality of $\lambda$ we have
\begin{align}
a &= \lambda_{x_\tap}(\scB(r)(\edb\ssoe  - 1)) \\
&= \scA(r)(\lambda_{x_\sop} (\edb\ssoe  - 1)) \\
&= \scA(r)(a_\sop -1)
\end{align}
and thus by \eqref{eq:defn_order_realisation_1} there is an edge $(a_\sop - 1, a)\in \edgeset(\scA(r))$. Further, $a_\sop  = \lambda_{x_\sop} (\edb\ssoe  - 1) + 1$ implies $\lambda_{x_\sop} (\edb\ssoe  - 1)  < a_\sop < \wwidehat \lambda_{x_\sop} (\edb\ssoe  +1)$ by injectivity of $\lambda_{x_\sop}$ (which implies injectivity of $\wwidehat \lambda_{x_\sop}$). Then, \eqref{eq:glambda_defn} implies that $(x_\sop,a_\sop)$ gets mapped to $(x_\sop,\edb\ssoe )$ by $\sS^\lambda$. Further, profunctoriality of $\scA(r)$ applied to $(a_\sop, a_\sop) \in \edgeset(\scA(x_\sop))$, $(a_\sop - 1, a) \in \edgeset(\scA(r))$ implies $\eda := (a_\sop,a) \in \edgeset(\scA(r))$. We deduce that $(r,\eda)$ is the required lift of $(r,\edb)$ under $\sS^\lambda$.

\item Assume $\edb\ssoe  \in \regcont(\scB(x_\sop)), \edb\ttae  \in \singcont(\scB(x_\tap))$ and $\edb\ttae  \notin \im(\scB(r))$. The second assumption implies $\lambda_{x_\tap}(\edb\ttae ) = a$. Set $b^\pm_\sop = \scB(r)\regop (\edb\ttae  \pm 1)$. Since $\edb\ttae  - 1 < \edb\ttae  < \edb\ttae  + 1$ and $(b^-_\sop,\edb\ttae -1), (\edb\ssoe ,\edb\ttae ), (b^+_\sop,\edb\ttae  + 1) \in \edgeset(\scB(r)$ we find $b^-_\sop \leq \edb\ssoe  \leq b^+_\sop$ by monotonicity of $R\scB(r)$ (cf. \autoref{claim:order_realisations_monotone}). We claim $b^\pm_\sop = \edb\ssoe $. Assume by contradiction $b^-_\sop < b^+_\sop$ and thus there is $b_\sop \in \singcont(\scB(x_\sop))$, $b^-_\sop < b_\sop < b^+_\sop$. We infer $(b_\sop,\edb\ttae ) \in \edgeset(\scB(r))$ by \eqref{eq:defn_order_realisation_4}. But then $\scB(r)(b_\sop) = \edb\ttae $ by \eqref{eq:defn_order_realisation_1}, contradicting our assumptions. We infer $\scB(r)\regop (\edb\ttae  \pm 1) = \edb\ssoe $. This implies $\edb^\pm = (\edb\ssoe ,\edb\ttae  \pm 1) \in \edgeset(\scB(r))$. We want to construct a lift for the morphism $(r,\edb^-)$ (the choice of sign is arbitrary here). Indeed, note 
\begin{equation}
\wwidehat\lambda_{x_\tap}((\edb\ttae  -1) - 1) < a - 1 < \lambda_{x_\tap}((\edb\ttae  -1) + 1) = (a - 1) +1
\end{equation}
(where we used injectivity of $\wwidehat \lambda_{x_\sop}$ for the first inequality). Thus by \eqref{eq:glambda_defn} we see $(x_\tap,a - 1)$ gets mapped to $(x_\tap,\edb\ttae  - 1)$ by $\sS^\lambda$. We set $a_\sop = \scA(r)\regop(a - 1)$ and find $\eda^- := (a_\sop,a-1) \in \edgeset(\scA(r))$ by \eqref{eq:defn_order_realisation_3}. Then $(r,\eda^-)$ maps onto $(r,\edb^-)$ by the argument used in part (a). Using profunctoriality of $\scA(r)$, $\eda^- \in \edgeset(\scA(r))$ and $(\eda^-\ttae ,a) \in \edgeset(\scA(x_\tap))$ imply $\eda := (\eda^-\ssoe , a) \in \edgeset(\scA(r))$. We deduce that $(r,\eda)$ is the required lift of $(r,\edb)$ under $\sS^\lambda$. \qed
\end{enumerate}

\end{lem}

\begin{eg}[Collapse has lifts] Using the family collapse $\sS^{\lambda_2}$ constructed in \autoref{eq:family_collapse}, the following depicts multiple liftings in $\sG(\scC)$ (marked as \corange{} arrows) of of the same arrow in $\sG(\scB_2)$ (also marked in \corange{})
\begin{restoretext}
\begingroup\sbox0{\includegraphics{test/page1.png}}\includegraphics[clip,trim={.2\ht0} {.0\ht0} {.2\ht0} {.0\ht0} ,width=.7\textwidth]{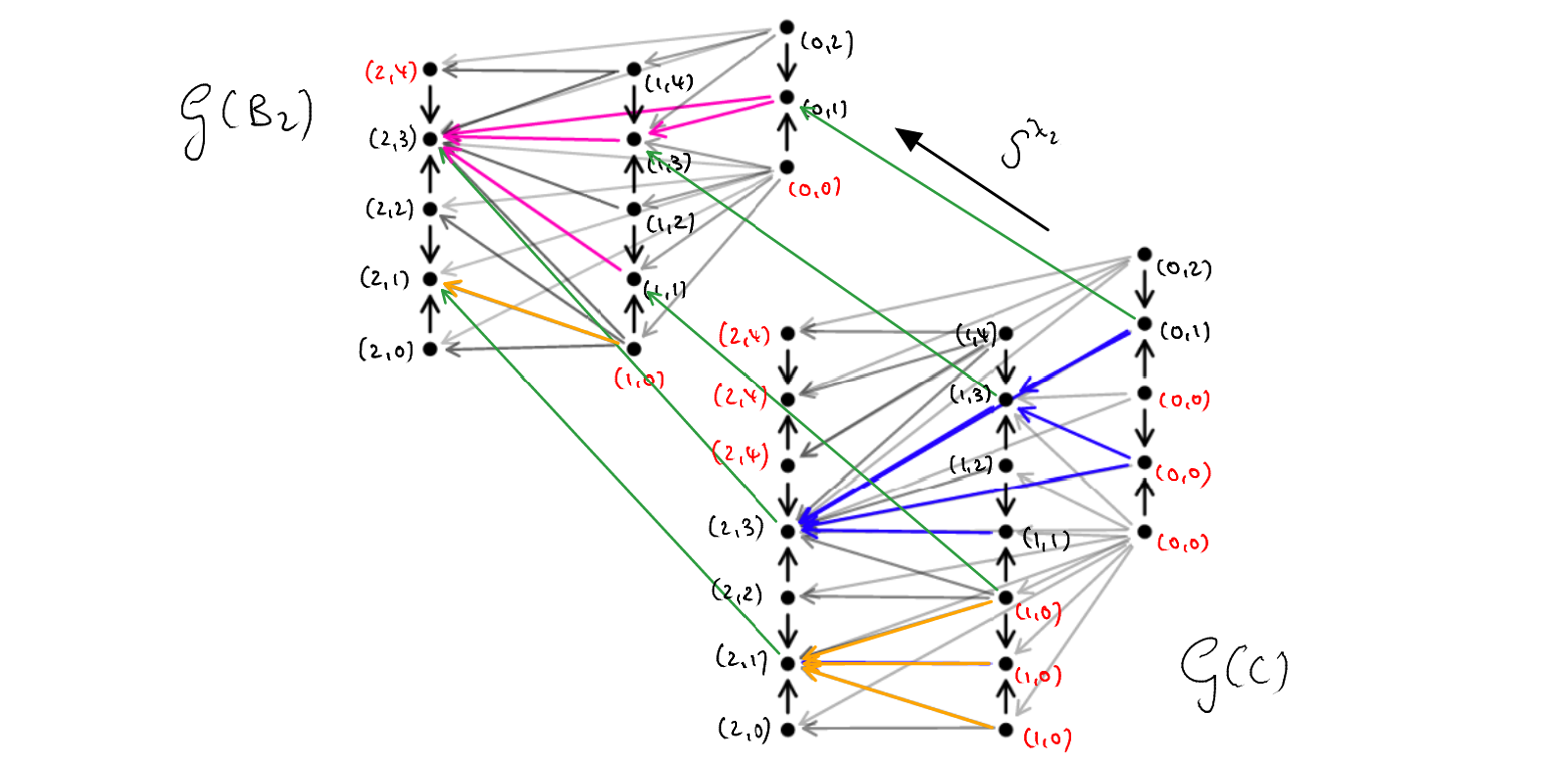}
\endgroup\end{restoretext}
\end{eg}

\begin{cor}[Collapse surjective on objects and morphisms] \label{cor:collapse_surjective} The functor of posets $\sS^\lambda$ given in \autoref{defn:glambda} is surjective on objects and on morphisms, i.e. for all $g \in \mor(\sG(\scB))$ there is $f \in \mor(\sG(\scA))$ such that $\sS^\lambda (f) = g$.

\proof We first show that $\sS^\lambda$ is surjective on objects. Let $(x,b) \in \sG(\scB)$. If $b \in \singcont(\scB(x))$ we have $\sS^\lambda (x,\lambda_x (b)) = (x,b)$ by \eqref{eq:glambda_defn}. Otherwise, if $b \in \regcont(\scB(x))$ we have at least one $a$ such that $\wwidehat \lambda_x(b-1) <  a < \wwidehat \lambda_x(b+1)$ since $\lambda_x$ is injective (which implies injectivity of $\wwidehat \lambda_x$). Then $\sS^\lambda(x,a) = (x,b)$  by \eqref{eq:glambda_defn}. Thus $\sS^\lambda$ is surjective on objects. 

Surjectivity of morphism now follows as corollary to surjectivity on objects together with $\sS^\lambda$ having \gls{lifts} as proved in the previous Lemma. \qed
\end{cor}

\subsection{Base change for collapse}

Base change interacts with collapse as follows.

\begin{cor}[Base change for collapse functors] \label{claim:precomposition_of_collapse_map} Let $\lambda : \scB \into \scA$ for $\scA, \scB : X \to \SI$, and let $H : Y \to X$ be a functor of posets. Then $\sS^{\lambda H}$ is described by the pullback (using \autoref{defn:grothendieck_base_change})
\begin{equation}
\xymatrix{ \sG(\scA H) \ar[r]^{\sG(H)} \ar[d]_{\sS^{\lambda H}} \pullbackfar & \sG(\scA) \ar[d]^{\sS^\lambda} \\
\sG(\scB H) \ar[r]_{\sG(H)} & \sG(\scB) }
\end{equation}
\proof 
We claim to have the following diagram
\begin{equation}
\xymatrix{ \sG(\scA H) \ar@/_2pc/[dd]_{\pi_{\scA H}} \ar[d]^{S} \ar[r]^{\sG(H)} \pullback & \sG(\scA) \ar[d]_{\sS^\lambda}  \ar@/^2pc/[dd]^{\pi_{\scA}}\\
{\sG(\scB H)} \ar[r]^{\sG(H)} \ar[d]^{{\pi_{\scB H}}} \pullback & \sG(\scB) \ar[d]_{\pi_{\scB}} \\
Y \ar[r]^{H} & X }
\end{equation}
The outer square and lower square are the pullbacks obtained by \autoref{defn:grothendieck_base_change}. Note the right outer triangle commutes since $\sS^\lambda$ is a map of bundles. Then $S$ is the universal factorisation of the outer pullback through the lower pullback. Now the upper square is also pullback by ``pullback cancellation on the right": This fact means that in any diagram of the above from, whenever the lower and outer square are pullbacks (which they are by \autoref{defn:grothendieck_base_change}) then this implies that the upper square is also a pullback. 

We now claim that $S = \sS^{{\lambda H}}$. This follows as a corollary to the uniqueness statement of \autoref{lem:glambda}: Indeed, $S$ is a map of bundles since the left triangle commutes by $S$ being a universal factorisation. It preserves regular segments on each fibre since both maps labelled $\sG(H)$ restrict to isomorphisms on each fibre (cf. \autoref{defn:grothendieck_base_change}) and $\sS^\lambda$ preserves regular segments. The fibrewise isomorphisms also imply that $S$ inverts $({\lambda H})_x = \lambda_{H(x)}$ on each fibre since $\sS^\lambda$ does. Thus, $S$ satisfies the conditions in \autoref{lem:glambda} and must equal $\sS^{{\lambda H}}$. \qed
\end{cor}

\begin{eg}[Base change for collapse functors] \label{eg:basechange_collapse} Using $\sS^{\lambda_2}$ as constructed in \autoref{eg:collapse_map}, we can apply the preceding lemma for the case of $H = \delta^2_0$ to obtain
\begin{restoretext}
\begingroup\sbox0{\includegraphics{test/page1.png}}\includegraphics[clip,trim={.2\ht0} {.0\ht0} {.3\ht0} {.0\ht0} ,width=.8\textwidth]{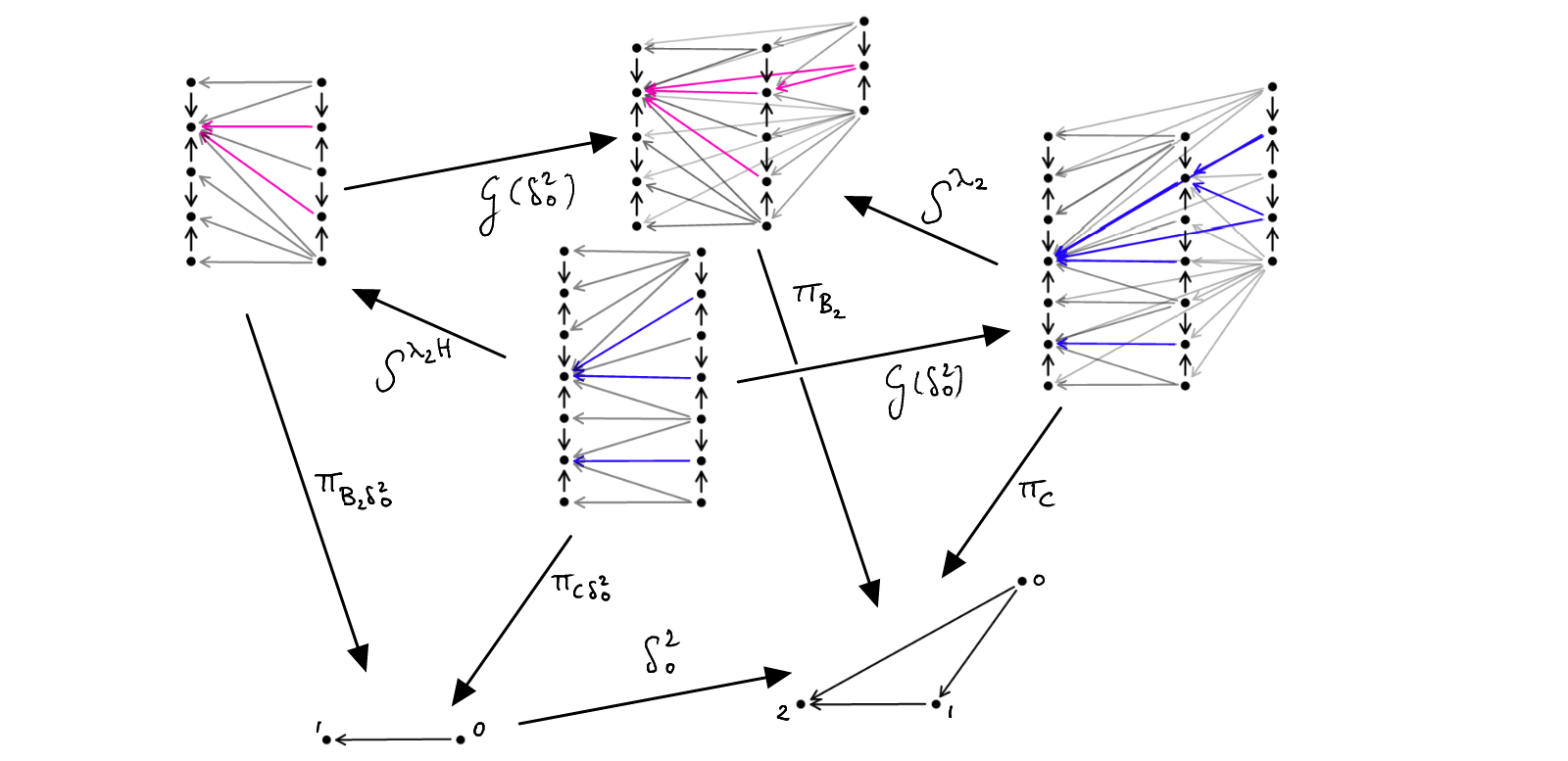}
\endgroup\end{restoretext}
\end{eg}

\begin{rmk}[Unwinding $\sS^{\lambda H}$] We remark that by \autoref{defn:stability_vs_injections}
\begin{equation} \label{eq:lambda_pullback}
\cS^{{\lambda H}}_x = \cS^{{\lambda}}_{H(x)}
\end{equation}
is a subset of $\singcont(\SIf{\scB}(H(x))) = \singcont(\SIf{\scA}(x))$ (since $\scA =  \scB H$). 
\end{rmk}

\subsection{What does collapse do?}

A collapse bundle map erases certain singular heights and absorbs (or ``collapses") them into a their respective neighbouring  regular segments. The following is a technical claim that illustrates this point, and will become useful later on.

\begin{claim}[Singular heights collapse to regular neighbours] \label{claim:collapse_sandwich} Let $\lambda : \scB \into \scA$ be an injection of singular interval families $\scA, \scB : X \to \SI$. Let $a \in \singcont(\scA(x))$ but $a \notin \cS^\lambda_x$. Then
\begin{equation}
\rest {\sS^\lambda} x (a\pm 1) = \rest {\sS^\lambda} x (a)
\end{equation}
\proof Set $b^\pm = \rest {\sS^\lambda} x (a\pm 1)$ and $b = \rest {\sS^\lambda} x (a)$. Since $\rest {\sS^\lambda} x$ preserves regular segments, $b^\pm$ must be regular segments. By \eqref{eq:glambda_defn} we have 
\begin{equation}
\wwidehat\lambda_x(b^\pm-1) < a\pm 1< \wwidehat\lambda_x(b^\pm + 1)
\end{equation}
Since $a \notin \cS^\lambda_x$ we have $\abs{\wwidehat\lambda_x(b^\pm - 1) - a} \geq 2$ and $\abs{\wwidehat\lambda_x(b^\pm + 1) - a} \geq 2$. On the other hand, we have $\abs{a - (a \pm 1)} \leq 1$ and thus deduce
\begin{equation}
\wwidehat\lambda_x(b^\pm-1) < a- 1<a<a+1< \wwidehat\lambda_x(b^\pm + 1)
\end{equation}
Using \eqref{eq:glambda_defn} again we therefore obtain $\sS^\lambda(x,a\pm 1) = (x,b)$ as claimed. \qed
\end{claim}

\begin{eg}[Singular heights collapse to regular neighbours] \label{eg:sing_height_coll}  Starting from \autoref{eq:SIvert_n_C_families} define the \SI-family $\und\scA : \singint 2 \to \mathbf{SI}$ by
\begin{equation}
\und\scA :=\und\sU^1_\scC
\end{equation}
In other words, $\und\scA$ is defined by the \SI-bundle
\begin{restoretext}
\begingroup\sbox0{\includegraphics{test/page1.png}}\includegraphics[clip,trim={.0\ht0} {.1\ht0} 0 {.1\ht0} ,width=\textwidth]{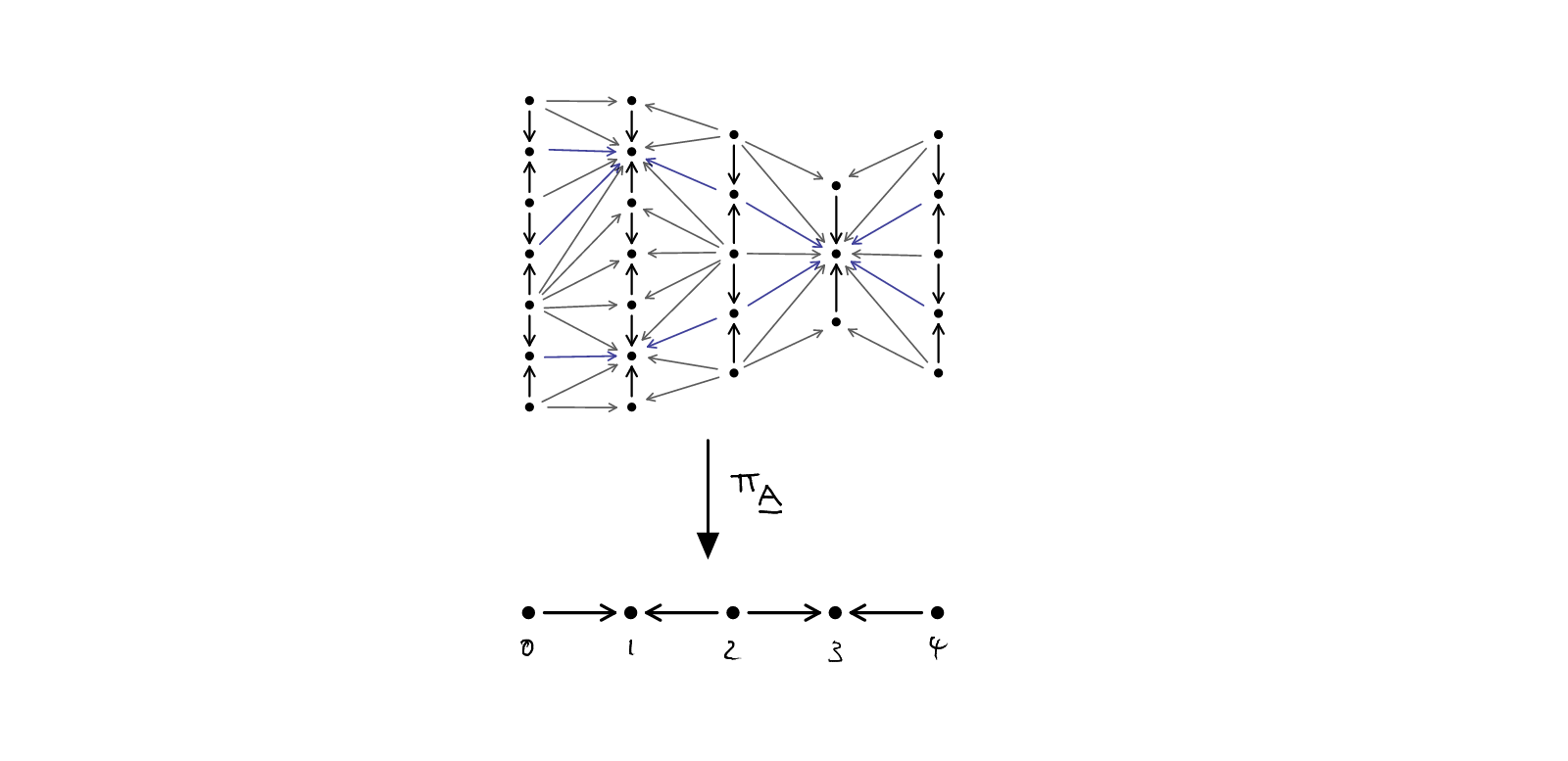}
\endgroup\end{restoretext}
Define a stable singular subset section $\cF$ of $\und\scA$ by the following subset of singular heights (marked in \cred{})
\begin{restoretext}
\begingroup\sbox0{\includegraphics{test/page1.png}}\includegraphics[clip,trim=0 {.2\ht0} 0 {.25\ht0} ,width=\textwidth]{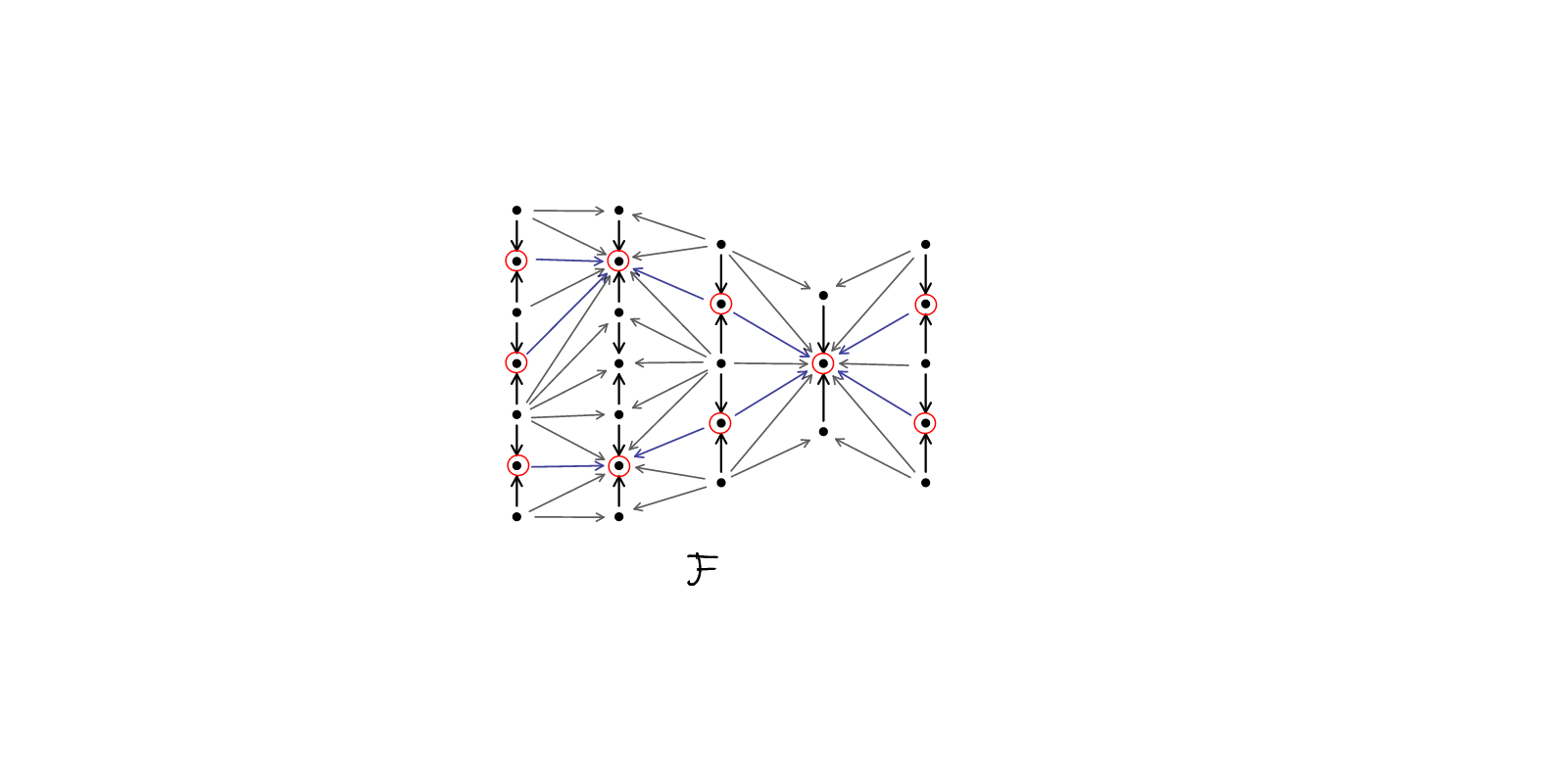}
\endgroup\end{restoretext}
We observe that $\cF$ is stable. Then, using \autoref{defn:stability_vs_injections} we find a family $\intrel{\cF} : \singint 2 \to \mathbf{SI}$ with \SI-bundle
\begin{restoretext}
\begingroup\sbox0{\includegraphics{test/page1.png}}\includegraphics[clip,trim=0 {.1\ht0} 0 {.1\ht0} ,width=\textwidth]{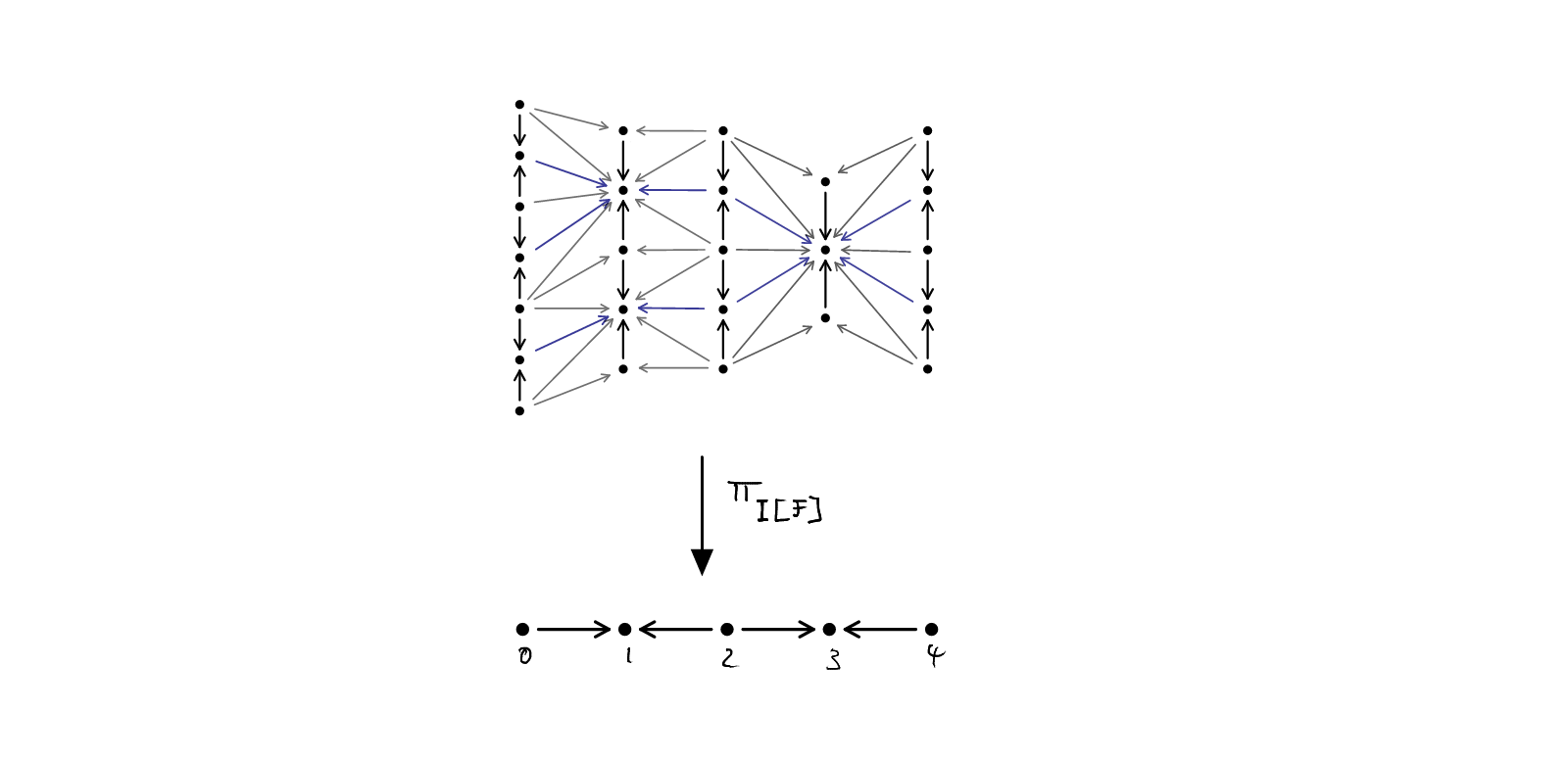}
\endgroup\end{restoretext}
and a natural injection $\mu := \eta_{\cF} : \intrel{\cF} \into \und\scA$. Note that only one singular height is missing from $\cF$, namely $3 \in \singcont(\und\scA(1))$ but $3 \in \cS^\mu_1$. Then the collapse functor $\sS^\mu$ is injective on morphisms (i.e. preimages are singletons) apart from on three morphisms which (together with their preimages) are marked in  \cred{}, \cpurple{} and \cdarkgreen{} in the following illustration
\begin{restoretext}
\begingroup\sbox0{\includegraphics{test/page1.png}}\includegraphics[clip,trim={.2\ht0} {.0\ht0} {.2\ht0} {.0\ht0} ,width=\textwidth]{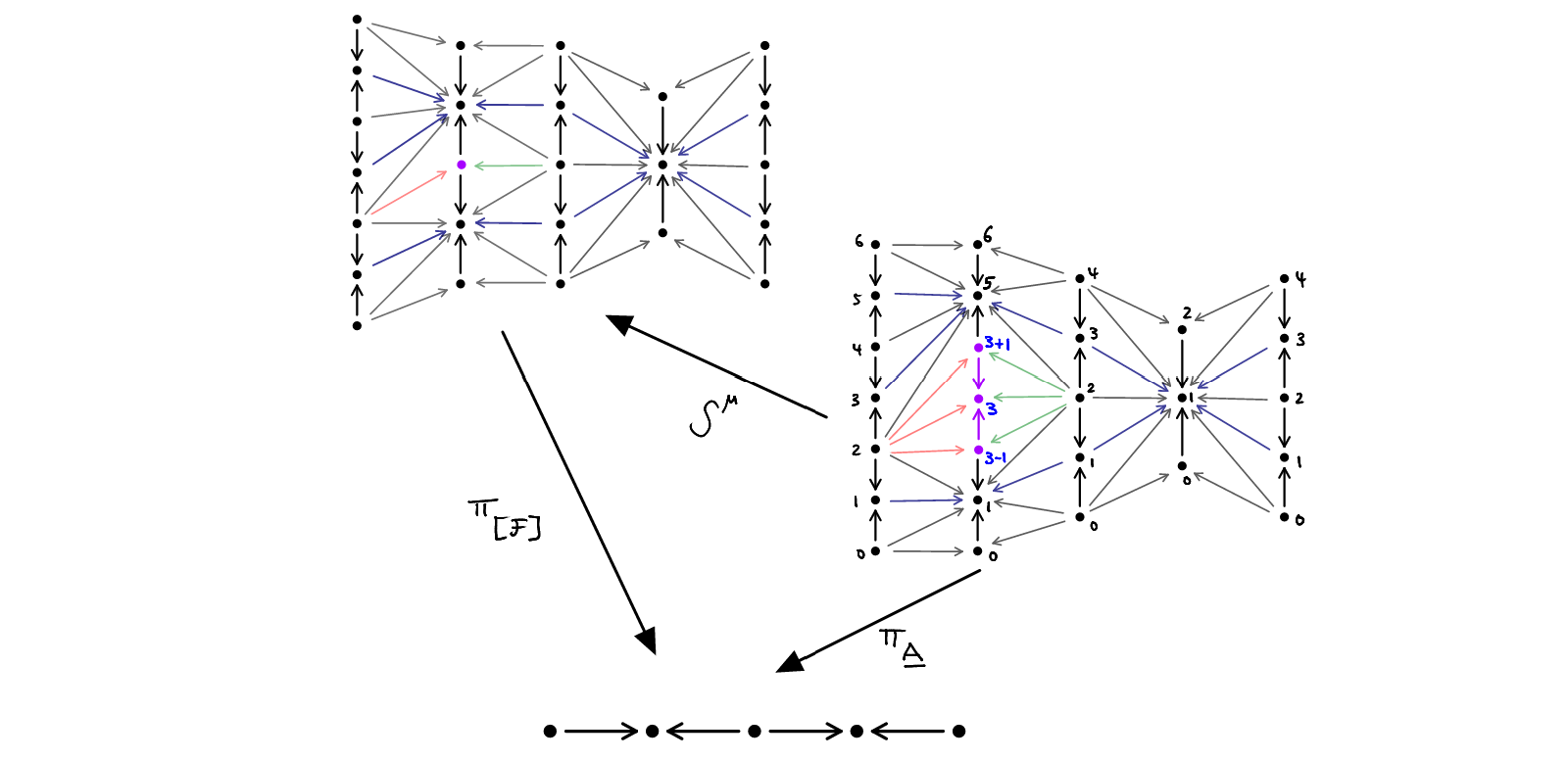}
\endgroup\end{restoretext}
In particular, we see that $a-1, a, a+1$ for $a = 3$ in the fibre over $1$ (of $\pi_{\und\scA}$) all get mapped by $\rest {\sS^\mu} 1$ to the same regular segment 
\begin{equation}
\rest {\sS^\mu} {1} (a \pm 1) = \rest {\sS^\mu} {1} (a) = 2
\end{equation}
in the fibre over $1$ (of $\pi_{\intrel{\cF}}$) as shown in \autoref{claim:collapse_sandwich}.
\end{eg}

\section{Collapse of labelled \SI-families} \label{sec:interval_collapse}

We now extend our discussion of collapse from singular interval families to labelled singular interval families. In the next chapter, this will be used to establish a notion of collapse of labelled singular cube families (since these are inductively build from interval families). In particular, we will see that collapse of cube families is a rewriting relation with unique normal forms. To prove this, this and the next section will establish the necessary definitions and results in the simpler context of interval families.

\subsection{Definition and examples}

\begin{defn}[Collapse of labelled \SI-families] \label{defn:collapse_of_SI_families} Let $\scA, \scB : X \to \SIvertone \cC$ be $\cC$-labelled $\SI$-families, and $\lambda : \SIf \scB \into \SIf \scA$ an injection of $\SI$-families. Then we say $\lambda$ is a witness for a \textit{labelled family collapse of $\scA$ to $\scB$}, written $\lambda : \scA \mcoll \scB$ or $\xymatrix@1{\scA\ar@{~>}[r]^{\lambda} & \scB}$, if
\begin{equation}
\sU_\scA = \sU_\scB \sS^\lambda
\end{equation}
We say \textit{$\scA$ collapses to $\scB$} (or alternatively, $\scB$ is a collapse of $\scA$), written $\scA \mcoll \scB$, if there is some $\lambda$ such that $\lambda :  \scA \mcoll \scB$.
\end{defn}

\begin{eg}[Collapse of labelled singular interval families] \label{eg:collapse_of_SIC_family} \hfill
\begin{enumerate}
\item Recall $\lambda_1$, $\scB_2$, $\scC$ and $H = \delta^2_0$ from \autoref{eg:basechange_collapse}. We construct a $\SIvertone \cC$-family $\scA : \bnum 2 \to \SIvertone \cC$ by setting
\begin{restoretext} %
\begingroup\sbox0{\includegraphics{test/page1.png}}\includegraphics[clip,trim={.1\ht0} {.1\ht0} {.3\ht0} {.15\ht0} ,width=\textwidth]{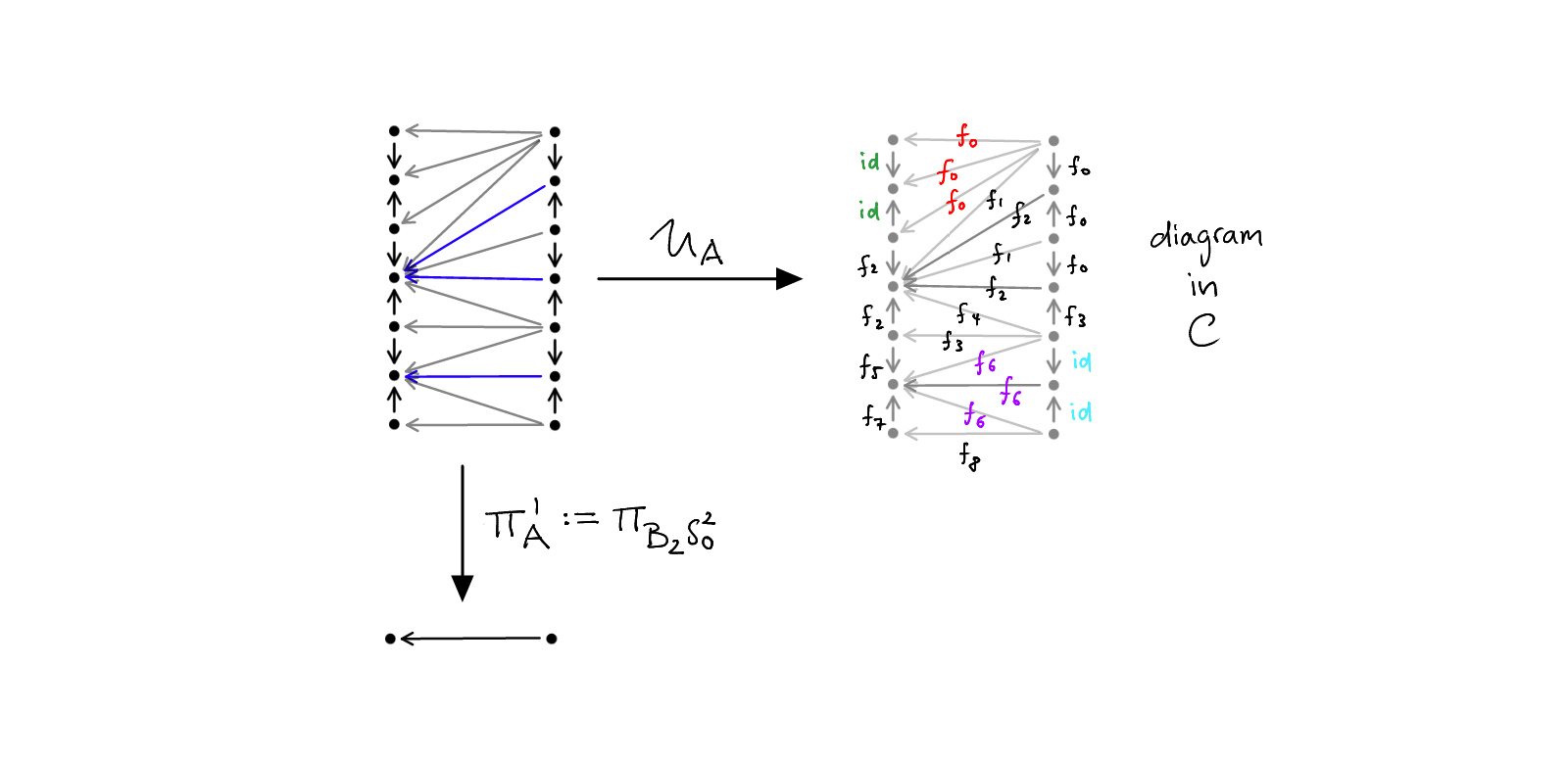}
\endgroup\end{restoretext}
and we further construct an $\SIvertone \cC$-family $\scB : \bnum 2 \to \SIvertone \cC$ by setting
\begin{restoretext} %
\begingroup\sbox0{\includegraphics{test/page1.png}}\includegraphics[clip,trim={.2\ht0} {.15\ht0} {.2\ht0}  {.15\ht0} ,width=\textwidth]{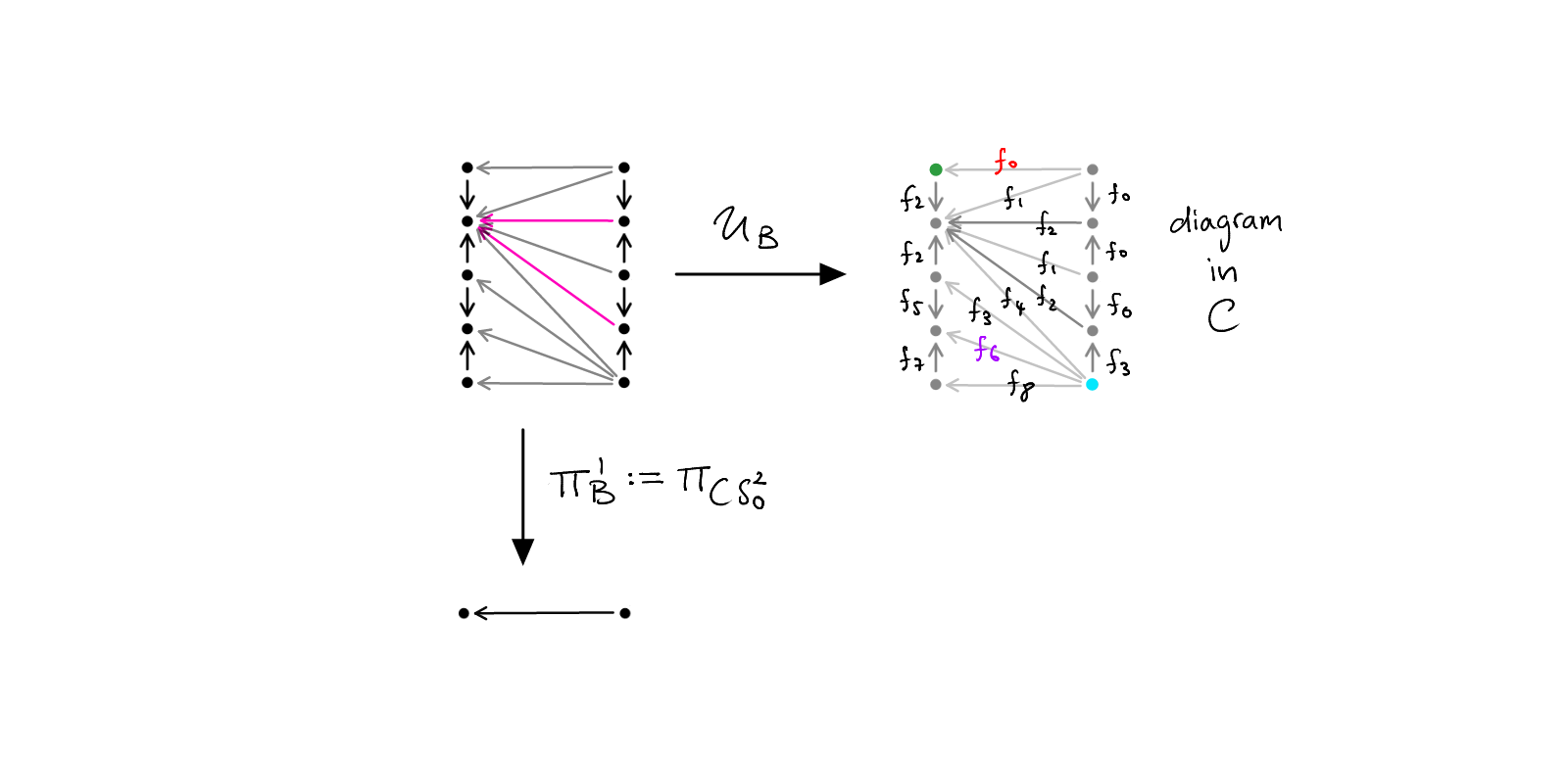}
\endgroup\end{restoretext}
Now, note that the morphism highlighted in \cdarkgreen{}, \cred{}, \cpurple{} and \cturquoise{} are those which are non-injectively mapped by $\sS^{\lambda_1 H}$ (this is illustrated below). In particular, the \cdarkgreen{} and \cturquoise{} points in the labelling of $\sG(\scB)$ by $\sU_\scB$ represent the identity morphism A necessary condition for $\sU_\scA$ to factor as $\sU_\scB \sS^{\lambda_1 H}$, is that morphisms in $\sG(\scA)$ which collapse to these identities by $\sS^{\lambda_1 H}$ need to be labelled by the identity themselves under $\sU_\scA$. Since this is the case for our choice of $\sU_\scA$, and by choice of all other labels above, we find that
\begin{restoretext}
\begingroup\sbox0{\includegraphics{test/page1.png}}\includegraphics[clip,trim=0 {.0\ht0} {.4\ht0} {.1\ht0} ,width=\textwidth]{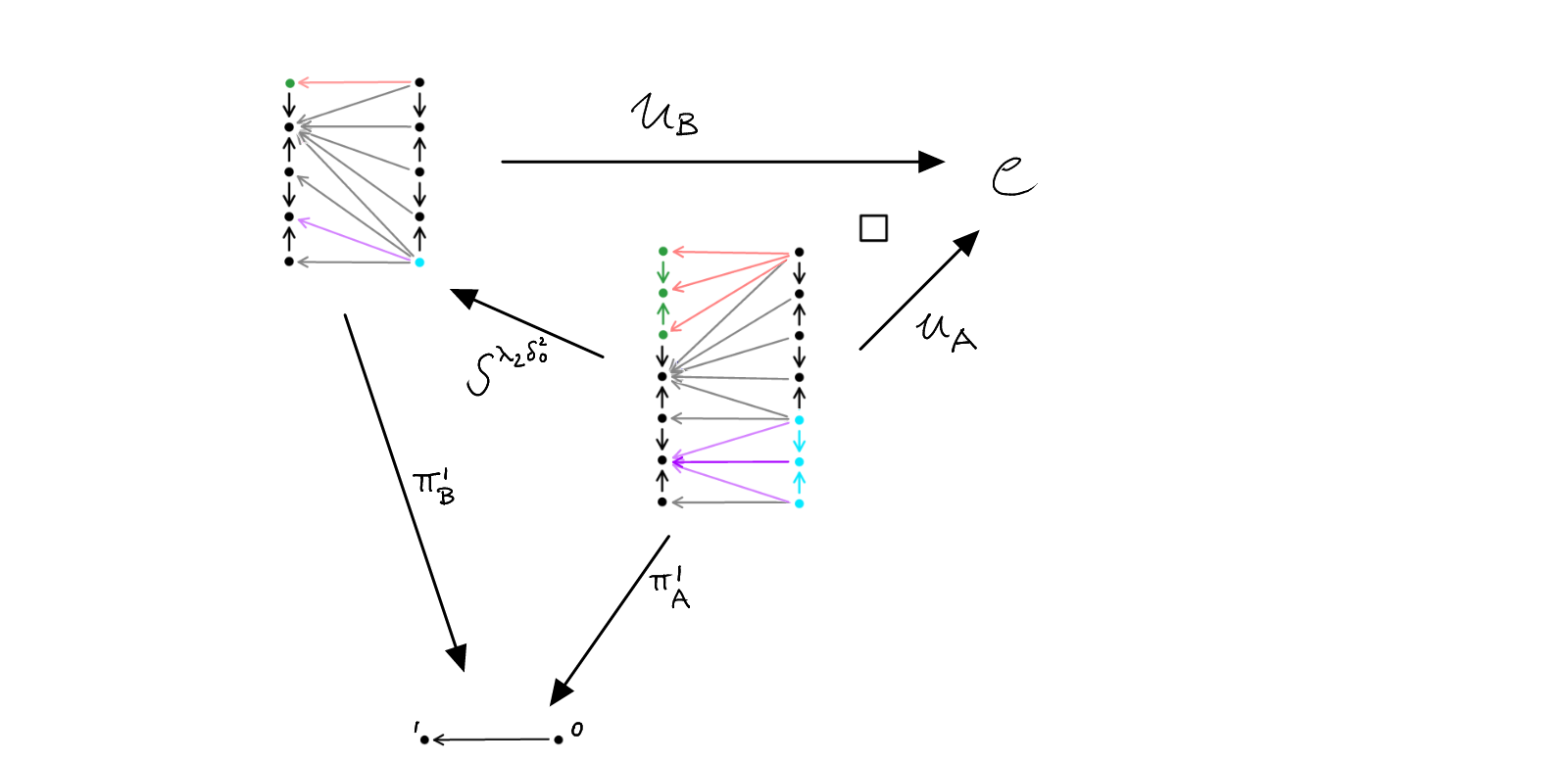}
\endgroup\end{restoretext}
This is precisely the condition for $\scA \mcoll \scB$.

\item From \autoref{eq:SIvert_n_C_families} recall the family $\scC : \bnum 1 \to \SIvert 2 \cC$.
The collapse $\mu : \intrel{\cF} \mcoll \tusU 1_\scC$ derived in  \autoref{eg:sing_height_coll} is non-injective only on three imagine points, and by choice of $\scC$, $\tsU 2_\scC$ is constant on their preimages (marked in \cred{}, \cpurple{} and \cdarkgreen{} below). Thus $\tsU 2_\scC$ factors through $\sS^\mu$ by some $U'$
\begin{restoretext}
\begingroup\sbox0{\includegraphics{test/page1.png}}\includegraphics[clip,trim={.2\ht0}  {.0\ht0} {.2\ht0}  {.0\ht0} ,width=.8\textwidth]{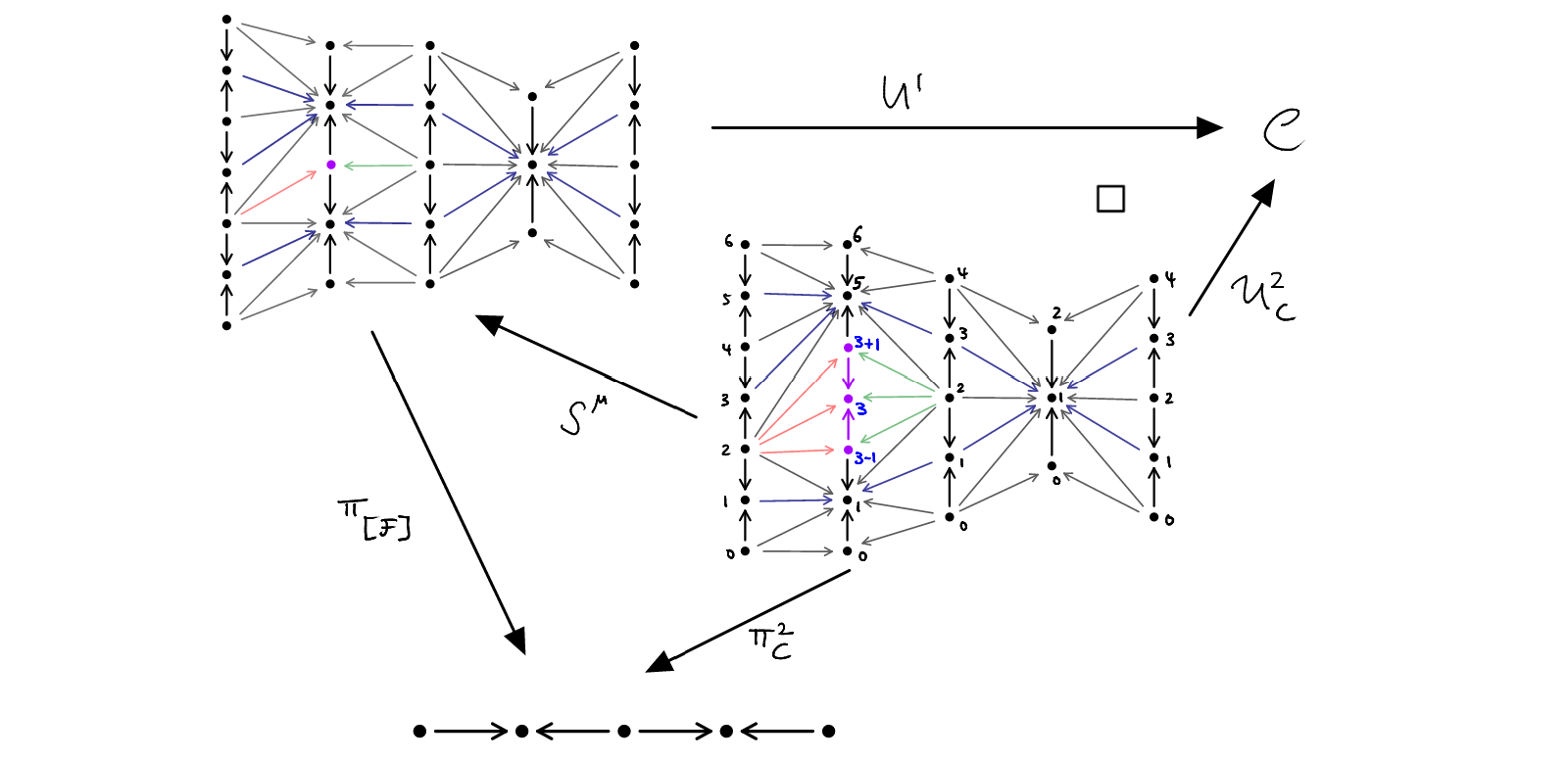}
\endgroup\end{restoretext}
In particular, we find
\begin{equation}
\mu : \tsU 1_\scC \mcoll \sR_{\intrel{\cF},U'}
\end{equation}
\end{enumerate} 
\end{eg}

\subsection{Normal forms of interval families}

We give the definition of normal forms of a labelled singular interval family $\scA$. These are families obtainable by collapse from $\scA$ to which no non-trivial collapse applies. We will later see that $\scA$ has exactly one normal form.

\begin{defn}[Normal forms] Given $\scA : X \to \SIvertone  \cC$ we say $\scA$ is in \textit{normal form}, if no non-identity collapse applies to $\scA$. That is, there is no $\lambda : \scA \mcoll \scB$ for any choice of $\scB, \lambda$ but the trivial one (namely, $\scB = \scA$ and $\lambda = \id_{\SIf \scA}$). The set of normal forms of $\scA$ is denoted by $\NF{\scA}$:
\begin{equation}
\NF{\scA} = \Set{\scB : X \to \SIvertone  \cC ~|~ \scA \mcoll \scB \text{~and~} \scB \text{~is in normal form}}
\end{equation}
\end{defn}

\begin{eg}[Normal forms] If all $f_i$ are different morphisms in the definition of $\scB$ in \autoref{eg:collapse_of_SIC_family}, then $\scB$ is in normal form. If for instance, $f_1 = f_2$ and $f_0  = \id$ then $\scB$ is not in normal form.
\end{eg}

\section{Properties of collapse} \label{sec:properties_collapse}

\subsection{Pushouts of collapses}

In this section we construct the category $\Buncoll \cC$, the category of labelled singular interval families and collapses, and show that it has pushouts. We will later on use this result established for interval families and apply it to cube families. This will be part of the proof of collapse normal forms being unique.

\begin{constr}[Category of interval families and collapses] \label{constr:cat_of_bun_and_coll} We construct $\Buncoll \cC$ as the category consisting of $\cC$-labelled $\SI$-families $\scA : X \to \SIvertone \cC$ (as objects) and witnesses of collapses $\lambda : \scA \mcoll \scB$ (as morphisms). Composition is given by composition of the underlying natural transformations.

To see that this is well-defined, note the following: if $\mu : \scA \mcoll \scB$, $\lambda : \scB \mcoll \scC$ for $\cC$-labelled $\SI$-families $\scA$, $\scB$ and $\scC$, then $\lambda\mu : \scA \mcoll \scC$: Indeed, using \autoref{cor:composition_of_Glambda} we verify
\begin{align}
\sU_\scA &= \sU_\scB \sS^\mu \\
&= \sU_\scC \sS^\lambda \sS^\mu \\
&= \sU_\scC \sS^{\lambda\mu}
\end{align}
which shows that $\lambda\mu$ witnesses the collapse $\scA \mcoll \scC$. 
\end{constr}

\noindent We will show that this category has pushouts by using a previous result from \autoref{claim:natural_injection_church_rosser}. To understand what having pushouts entails for $\Buncoll \cC$ we first discuss the following example.

\begin{eg}[Pushout of collapse functors] \label{eg:pushout_of_coll_funct} Using the definitions in \autoref{eg:pullback_of_injections}, and turning injections into collapse functors by $\sS$, then we find that the following commuting square
\begin{restoretext}
\begingroup\sbox0{\includegraphics{test/page1.png}}\includegraphics[clip,trim=0 {.0\ht0} 0 {.0\ht0} ,width=\textwidth]{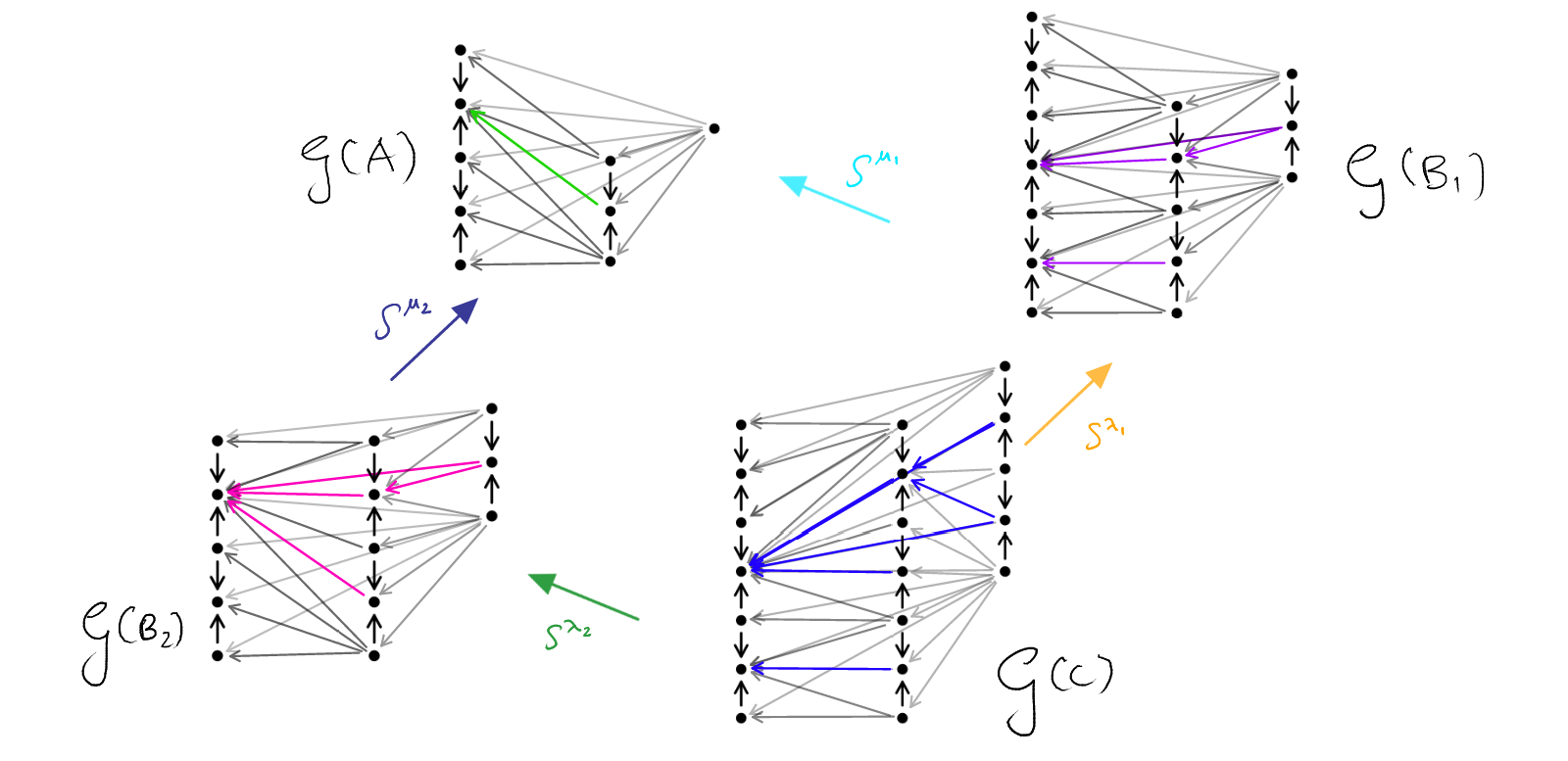}
\endgroup\end{restoretext}
Note that the diagonal map of this square (using the names of \autoref{eg:pullback_of_injections}) is given by $\sS^\eps$.

We claim that the above square is a pushout in $\Cat$. For this to be true we need to show that whenever we have a cocone on the above diagram (marked in \cblue{} below) there is a unique factorisation $F$ through $\sG(A)$ (marked in \cred{} below), that is
\begin{restoretext}
\begingroup\sbox0{\includegraphics{test/page1.png}}\includegraphics[clip,trim=0 {.2\ht0} 0 {.0\ht0} ,width=\textwidth]{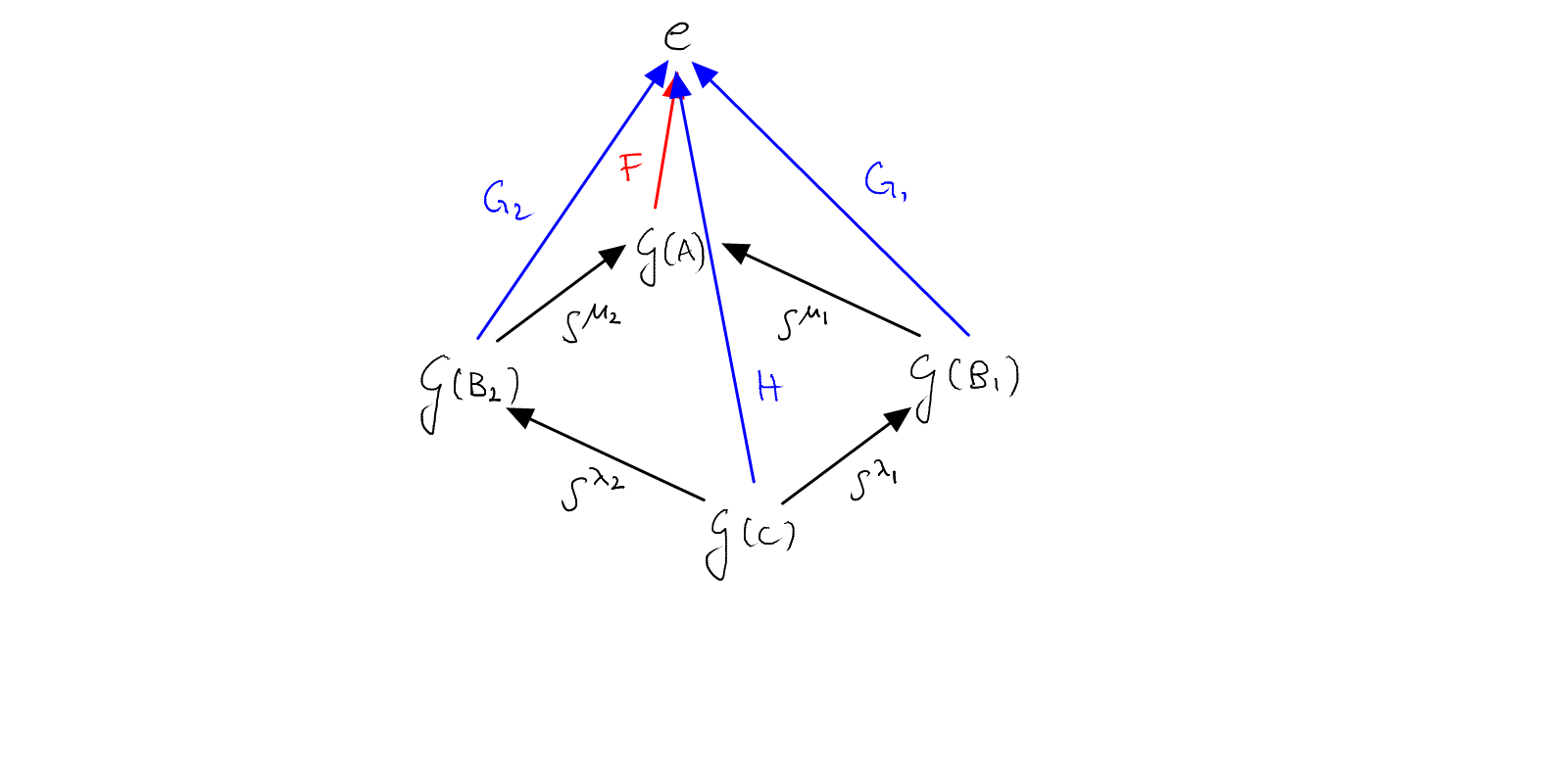}
\endgroup\end{restoretext}
Setting $\eps= \lambda_i\mu_i$ ($i \in \Set{1,2}$) we deduce that we must have $H = F \sS^\eps$. For such $F$ to exist we must show that $H$ is constant on preimages of $\sS^\eps$. Two arrows in $\sG(\scC)$ (marked in \cred{} below) in the preimage of an arrow in $\sG(\scA)$ (marked in \cred{} below) are shown in the following
\begin{restoretext}
\begingroup\sbox0{\includegraphics{test/page1.png}}\includegraphics[clip,trim=0 {.0\ht0} 0 {.0\ht0} ,width=\textwidth]{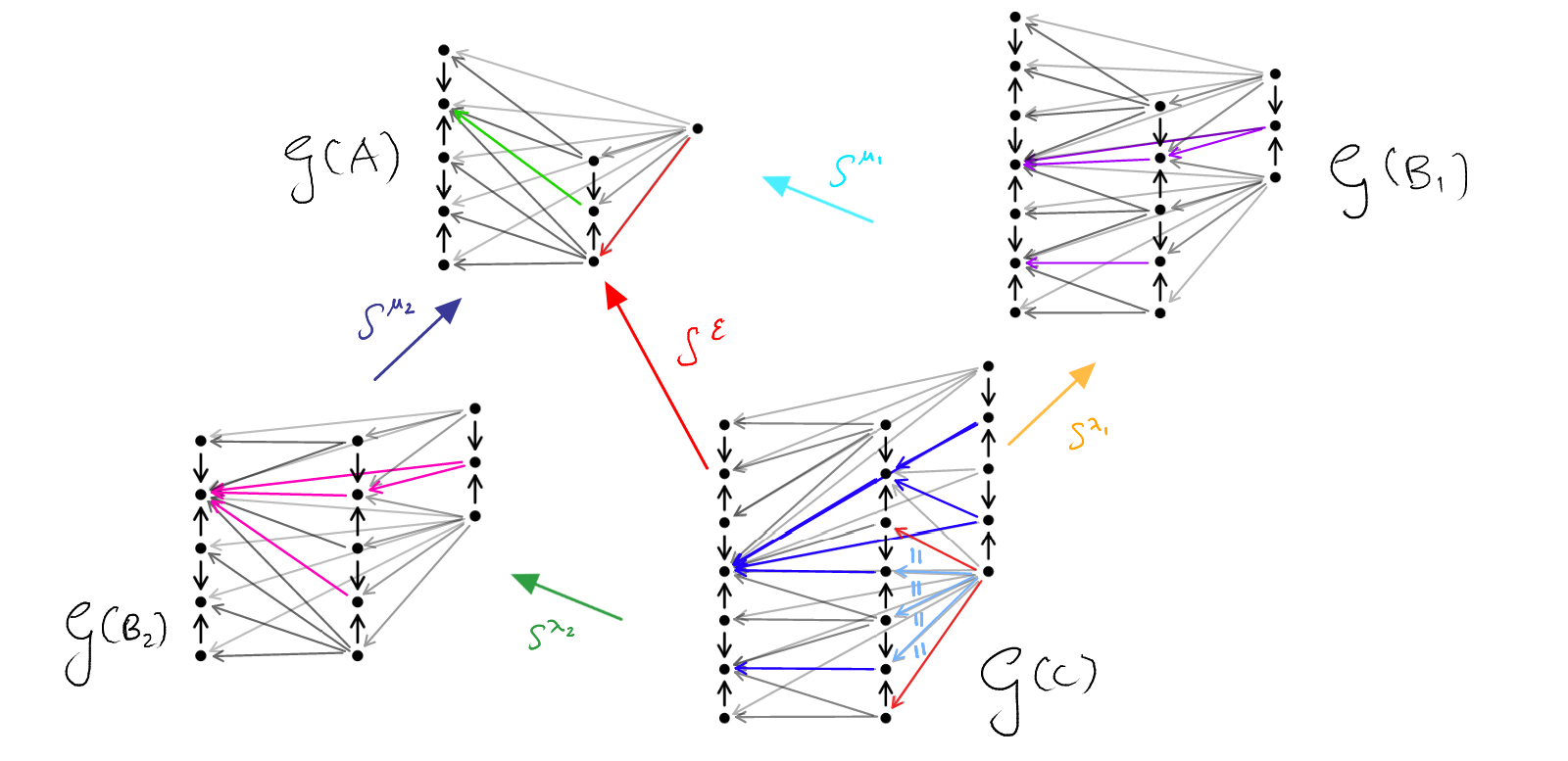}
\endgroup\end{restoretext}
To show that $H$ agrees on both arrows, we will use edge induction and ``induct over" the arrows in between the \cred{} arrows in $\sG(\scC)$, which are marked in \clightblue{}, showing that $H$ agrees on any adjacent two of them in each step. This will in turn follow since $\eps$ is the intersection of $\lambda_1$ and $\lambda_2$ (cf. \autoref{claim:natural_injection_church_rosser}), and thus, as shown below, each inductively claimed equality follows from inspection of one of $H = G_i \sS^{\lambda_i}$, $i \in \Set{1,2}$
\begin{restoretext}
\begingroup\sbox0{\includegraphics{test/page1.png}}\includegraphics[clip,trim=0 {.0\ht0} 0 {.0\ht0} ,width=\textwidth]{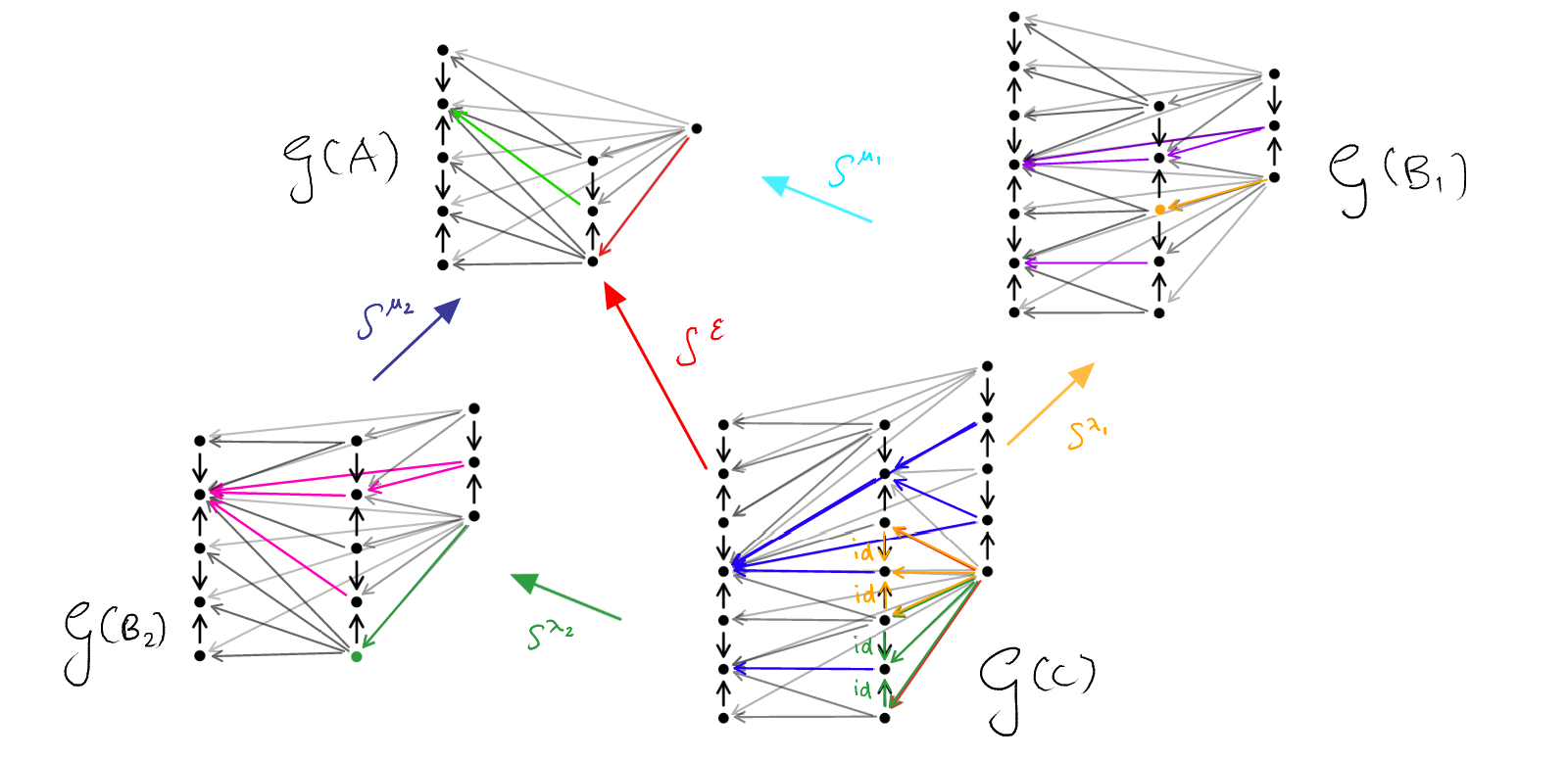}
\endgroup\end{restoretext}
Here, the edges marked by $\id$ in $\sG(\scC)$ are mapped by $H$ to $\id$ (for those colored in \corange{} this follows from  $H = G_1 \sS^{\lambda_1}$, and for those marked in \cdarkgreen{} this follows from $H = G_2 \sS^{\lambda_2}$).
\end{eg}

\begin{lem}[Intersection of collapse gives $\Cat$ pushouts] \label{claim:natural_injection_pushout} Let $\scC, \scB_1, \scB_2 : X \to \SIvertone \cC$ are $\cC$-labelled $\SI$-families over a poset $X$. Let $\lambda_i : \scB_i \into \SIf \scC$, $i \in \Set{1,2}$, be injections and let
\begin{equation}
\xymatrix{ \scA \ar[r]^{\mu_1} \ar[d]_{\mu_2} \pullback & \scB_1 \ar[d]^{\lambda_1} \\
\scB_2 \ar[r]_{\lambda_2} & \scC}
\end{equation}
be the pullback of injections as constructed in \autoref{claim:natural_injection_church_rosser} and write $\eps = \lambda_i\mu_i$ (note $\scA : X \to \SI$). Then,
\begin{equation} \label{eq:pushout_sing_coll_map}
\xymatrix{ \sG(\scC) \ar[r]^{\sS^{\lambda_1}} \ar[d]_{\sS^{\lambda_2}} & \sG(\scB_1) \ar[d]^{\sS^{\mu_1}} \\
\sG(\scB_2) \ar[r]_{\sS^{\mu_2}} &  \sG(\scA) \pushoutfar}
\end{equation}
is a pushout in $\Cat$, the category of categories and functors.

\proof[Proof of \autoref{claim:natural_injection_pushout}]  Let $\cC \in \Cat$ be a category and assume we are given functors $G_j : \sG(\scB_j) \to \cC$, $H : \sG(\scC) \to \cC$ such that $H = G_j \sS^{\lambda_j}$. We need to show that there is a unique $F : \sG(\scA) \to \cC$ such that $G_j = F \sS^{\mu_j}$, $j \in \Set{1,2}$. We construct $F$ explicitly as follows.

Let $(r : x_\sop \to x_\tap) \in \mor(X)$ and $\eda \in \edgeset(\scA(r)))$, that is, let $(r,\eda)$ be a morphism in $\sG(\scA)$ (cf.  \eqref{eq:grothendieck_edges}). Since $\sS^\eps$ by \autoref{cor:collapse_surjective} is surjective on objects and morphisms we can find $(r,\edc)\in \mor(\sG(\scC))$ which is mapped to $(r,\eda)$ by $\sS^\eps$. We define $F$ (on morphism \textit{and} objects) by setting
\begin{equation} \label{eq:intersection_functor}
F(r,\eda) := H(r,\edc)
\end{equation}
and need to check that this is well-defined: that is, if $(r,\edc') \in \mor(\sG(\scC))$ also satisfies
\begin{equation}
\sS^\eps(r,\edc') = (r,\eda)
\end{equation}
then for \eqref{eq:intersection_functor} to be well-defined we need to show
\begin{equation} \label{eq:pushout_welldef}
H(r,\edc) = H(r,\edc')
\end{equation}
We argue by edge induction. That is, we argue by induction on the distance 
\begin{equation}
n = \abs{\avg{\edc'} - \avg{\edc}}
\end{equation}
where $\avg{-}$ is the norm defined in \autoref{defn:edge_sets}. Note that by \autoref{claim:edge_set_properties}, if $n = 0$ we must have $\edc = \edc'$ and thus \eqref{eq:pushout_welldef} is satisfied. Without loss of generality assume $\avg{\edc} < \avg{\edc'}$ (otherwise switch their roles in the following). Recall \autoref{defn:successor} and the construction of the filler index $\succindex \edb$. By assumption we have
\begin{equation}
\sS^\eps(x_{\succindex \edc},\edc_{\succindex \edc}) =\sS^\eps(x_{\succindex \edc},\edc'_{\succindex \edc}) = (x_{\succindex \edc},\eda_{\succindex \edc})
\end{equation}
Then, by applying monotonicity of $\rest {\sS^\eps} {x_{\succindex \edc}}$ to the inequalities of \autoref{cor:filler_edge} (obtained in the case of $\avg{\edc} < \avg{\edc'}$), we can infer
\begin{equation}
\sS^\eps(x_{\succindex \edc}, \succfill \edc\ssoe ) = \sS^\eps(x_{\succindex \edc}, \succfill \edc\ttae ) = (x_{\succindex \edc}, \eda_{\succindex \edc})
\end{equation}
Since $\succfill \edc\ssoe $ is a regular segment, and $\rest {\sS^\eps} {x_{\succindex \edc}}$ preserves regular segments this forces $\eda_{\succindex \edc}$ to be a regular segment. Thus $\succfill \edc\ttae  \notin \cS^{\eps}_{x_{\succindex \edc}} = \im(\eps_{x_{\succindex \edc}})$, since otherwise $\eda_{\succindex \edc}$ would be a singular height by \eqref{eq:glambda_defn}. 

By definition of the injection $\eps$ as intersection of the injections $\lambda_1, \lambda_2$ (cf. \autoref{claim:natural_injection_church_rosser}) there is $k \in \Set{1,2}$ such that $\succfill \edc\ttae  \notin \cS^{\lambda_k}_{x_{\succindex \edc}}$. Using \autoref{claim:collapse_sandwich} together with $\succfill \edc\ttae  \notin \cS^{\lambda_k}_{x_{\succindex \edc}}$ and $\succfill \edc\ssoe  = \succfill \edc\ttae  \pm 1$ (for one choice of $\Set{+,-}$) we find
\begin{equation}
\sS^{\lambda_k}(x_{\succindex \edc}, \succfill \edc\ssoe ) = \sS^{\lambda_k}(x_{\succindex \edc}, \succfill \edc\ttae )
\end{equation} 
and thus
\begin{equation}
\sS^{\lambda_k}(\id_{x_{\succindex \edc}}, \succfill \edc) = \id
\end{equation} 
Using \autoref{rmk:successor_compositionality} and functoriality of $\sS^{\lambda_k}$ we find
\begin{equation}
\sS^{\lambda_k}(r,\edc) = \sS^{\lambda_k}(r,\succ \edc)
\end{equation}
Since further $H = G_k\sS^{\lambda_k}$ we deduce
\begin{equation}
H(r,\edc) = H(r,\succ \edc)
\end{equation}
Since $\abs{\avg{\succ \edc} - \avg{\edc'}} = n - 1$ this completes the inductive proof of well-definedness of \eqref{eq:intersection_functor}. 

We still need to verify functoriality of $F$ as defined via \eqref{eq:intersection_functor}: It suffices to note that since $\sS^\eps$ has \gls{lifts} and is surjective on objects (by \autoref{lem:glambda} ) then for any chain $(x,a) \to (x', a') \to (x'', a'')$ in $\sG(\scA)$ we can find a chain $(x,b) \to (x', b') \to (x'', b'')$ in $\sG(\scC)$ that is mapped onto it by $\sS^\eps$. Using \eqref{eq:intersection_functor}, functoriality of $F$ thus follows from functoriality from $H$. This finishes the construction of $F$, shows that it is necessarily unique (since $\sS^\eps$ is surjective on objects and morphism) and by its definition \eqref{eq:intersection_functor} we obtain $H = F \sS^\eps$.

Finally, we also obtain (for $j \in \Set{1,2}$)
\begin{align*}
G_j \sS^{\lambda_j} &= H \\
&= F \sS^\eps \\
&= F \sS^{\mu_j} \sS^{\lambda_j}
\end{align*} 
where, in the last step we used \autoref{cor:composition_of_Glambda}. Now, since $\sS^{\lambda_j}$ is surjective on objects and morphisms, we deduce $G_j = F \sS^{\mu_j}$ as required. \qed
\end{lem}

\begin{thm}[Pushout of collapse] \label{thm:family_pushouts} The category $\Buncoll \cC$ has pushouts. 

\proof Assume
\begin{equation}
\xymatrix{ \scC \ar@{~>}[r]^{ {\lambda_1}} \ar@{~>}[d]_{ {\lambda_2}} & \scB_1  \\
\scB_2  &  }
\end{equation}
where $\scC, \scB_1, \scB_2 : X \to \SIvertone \cC$ are $\cC$-labelled $\SI$-families over some poset $X$. Define $\mu_i$, $\scA : X \to \SI$ via the pushout \eqref{eq:pushout_sing_coll_map} in $\Cat$. Note that the maps $\sU_\scC$ and $\sU_{\scB_i}$ give a cocone on this pushout (they play the role of $H$ and $G_i$ respectively in \autoref{claim:natural_injection_pushout}). The pushout property thus yields a factorising map $F : \sG(\scA) \to \cC$. We obtain
\begin{equation}
\xymatrix{ \scC \ar@{~>}[r]^{ {\lambda_1}} \ar@{~>}[d]_{ {\lambda_2}} & \scB_1 \ar@{~>}[d]^{ {\mu_1}} \\
\scB_2 \ar@{~>}[r]_{ {\mu_2}} &  \widecheck\scA \pushoutfar}
\end{equation}
by setting
\begin{equation}
\widecheck\scA = \sR_{\scA, F}
\end{equation}
The fact that this is a pushout in $\Buncoll \cC$ follows from the pullback property proven in \autoref{claim:natural_injection_church_rosser}, and the pushout property proven in \autoref{claim:natural_injection_pushout}. \qed

\end{thm}

\begin{eg} By adding labelling $\sU_\scC$ and $\sU_{\scB_i}$ (for $i \in \Set{1,2}$) such that $\sU_\scC = \sU_{\scB_i} \sS^{\lambda_i}$ to the pushout of total posets $\sG(\scC)$ and $\sG(\scB_i)$ in \autoref{eg:pushout_of_coll_funct}, we obtain $\SIvertone \cC$ families $\scC$, $\scB_1$ and $\scB_2$ such that $\lambda_i : \scC \mcoll \scB_i$. The pushout $\scA$ of collapses can then be visualised as follows
\begin{restoretext}
\begingroup\sbox0{\includegraphics{test/page1.png}}\includegraphics[clip,trim={.3\ht0} {.0\ht0} {.3\ht0} {.0\ht0} ,width=\textwidth]{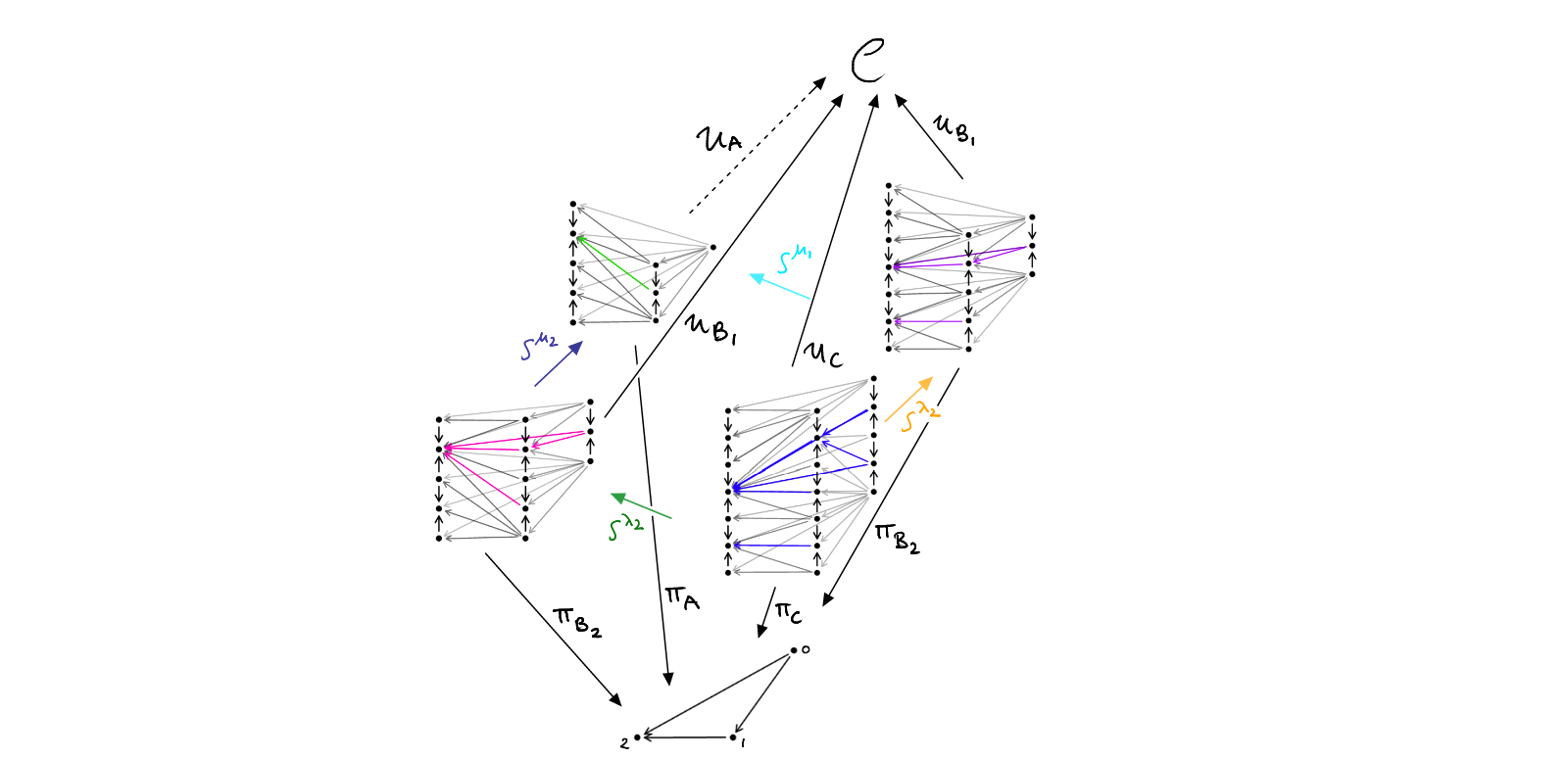}
\endgroup\end{restoretext}
\end{eg}

\subsection{Pullbacks of collapses along base change}

We will now investigate the relation of collapse and base change. We will find that collapse can be ``pulled back" along a base change functor. We will later on use this result established for interval families and apply it to cube families. This will be part of the proof of collapse normal forms being unique.

\begin{lem}[Base change for collapse] \label{rmk:precomposition_of_collapse} If $\lambda : \scA \mcoll \scB$ for $\scA, \scB : X \to \SIvertone \cC$, and $H : Y \to X$ is a functor of posets, then $\lambda H : \scA H \mcoll \scB H$ (cf. \autoref{notn:basic_category_theory}). 

\proof This follows since \autoref{claim:grothendieck_span_construction_basechange} together with \autoref{claim:precomposition_of_collapse_map} allows to compute
\begin{align}
\sU_{\scA H} &= \sU_\scA \sG(H) \\
&=\sU_\scB \sS^\lambda \sG(H) \\
&=\sU_\scB  \sG(H)  \sS^{\lambda H}\\
&=\sU_{\scB H}  \sS^{\lambda H}
\end{align}
as required for $\lambda H$ to witness a collapse. \qed
\end{lem}

\begin{eg}[Base change for collapse] The preceding proof can be illustrated by adding labels to \autoref{eg:basechange_collapse} (where $H = \delta^2_0$) as follows
\begin{restoretext}
\begingroup\sbox0{\includegraphics{test/page1.png}}\includegraphics[clip,trim={.3\ht0} {.0\ht0} {.3\ht0} {.0\ht0} ,width=\textwidth]{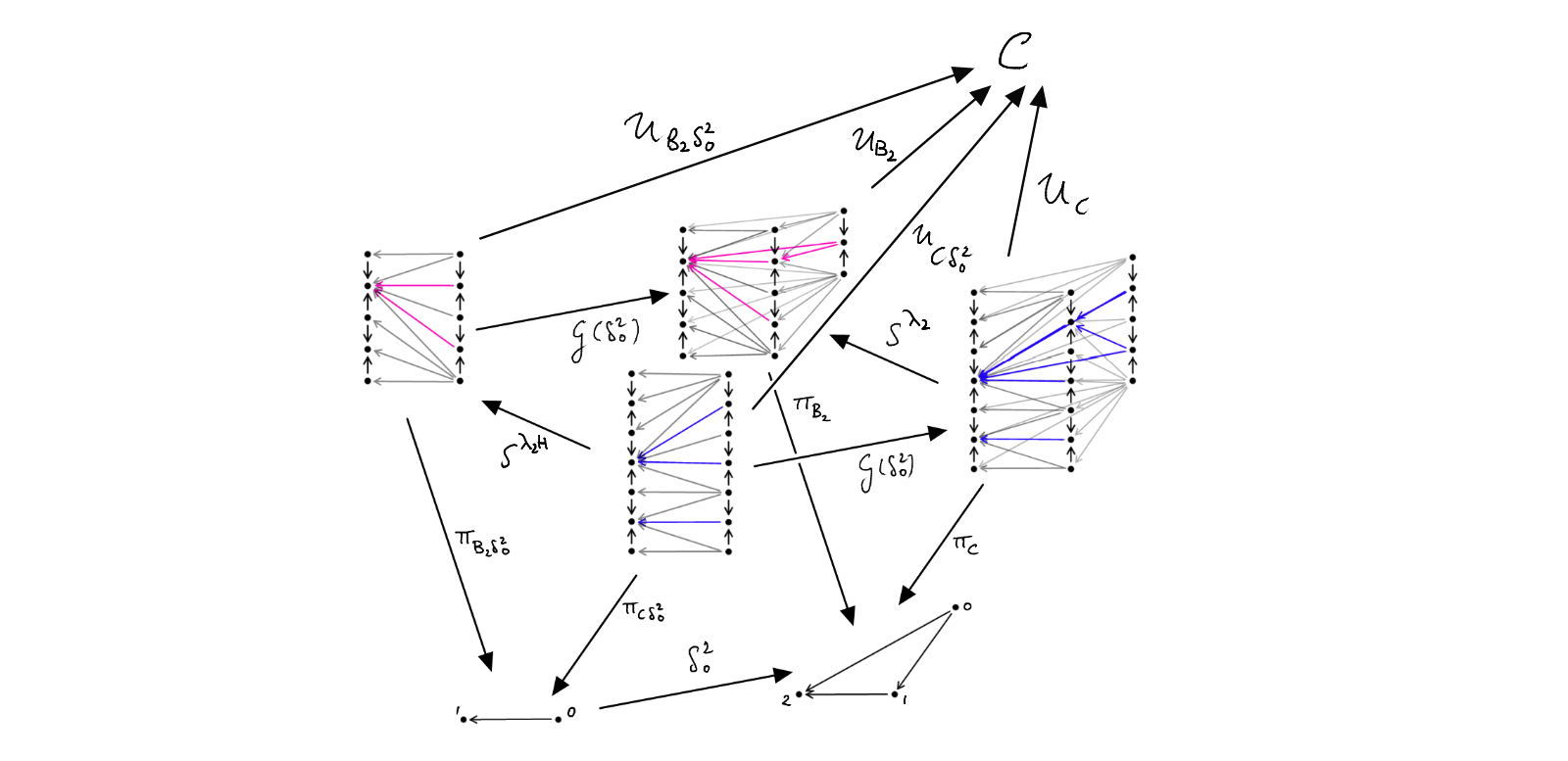}
\endgroup\end{restoretext}
Thus, we see that under the assumption that $\lambda_2 : \scC \mcoll \scB_2$, we can deduce that $\lambda_2 \delta^2_0 : \scC \delta^2_0 \mcoll B_2 \delta^2_0$ as proven in the previous lemma.
\end{eg}

We reformulate the above as follows

\begin{cor}[Pullback of collapses] \label{claim:collapse_dimension_interaction1}
Consider 
\begin{restoretext}
\begingroup\sbox0{\includegraphics{test/page1.png}}\includegraphics[clip,trim=0 {.0\ht0} 0 {.1\ht0} ,width=\textwidth]{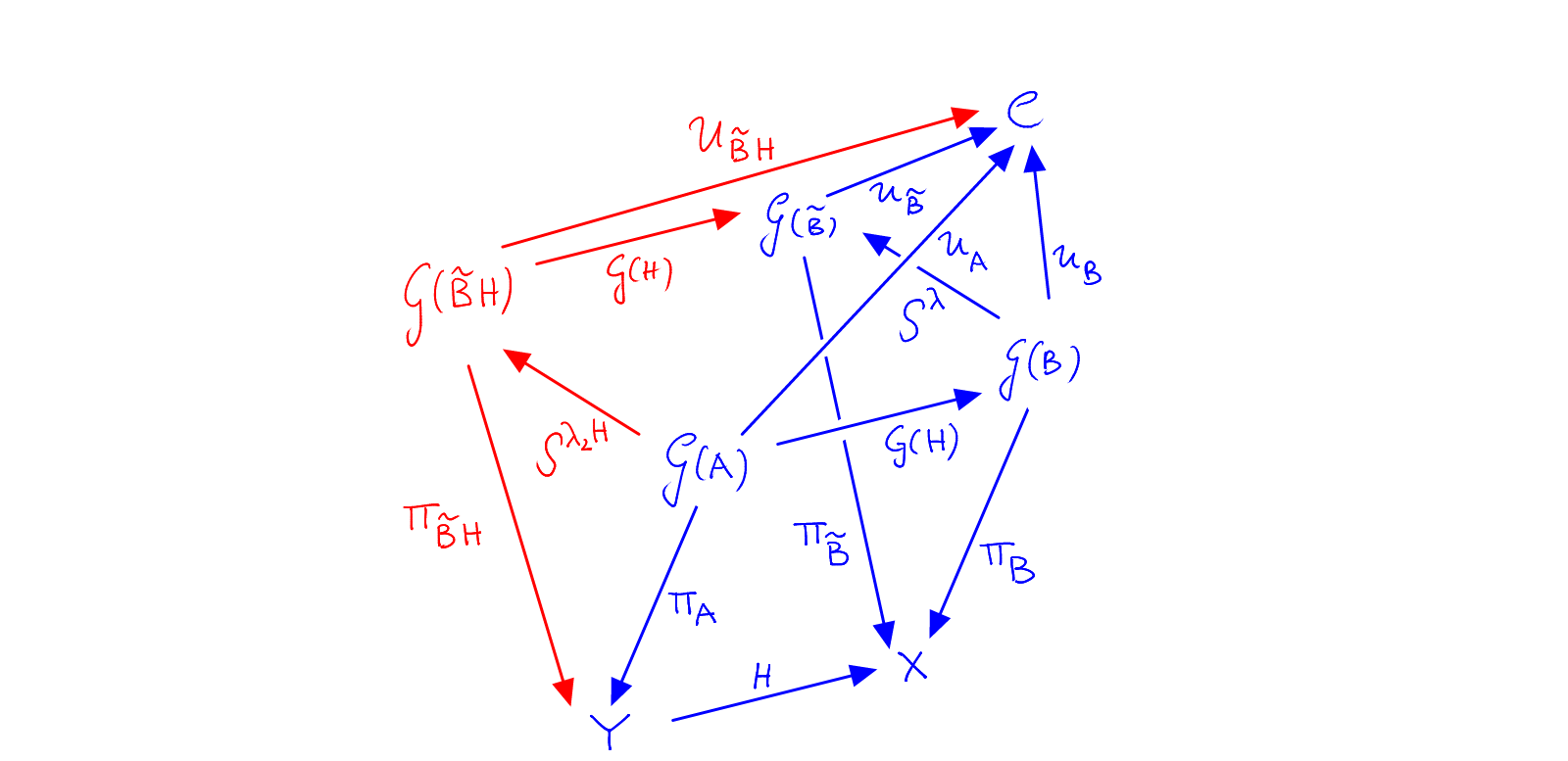}
\endgroup\end{restoretext}
Assume existence of the \cblue{} part of the diagram (which is equivalent to the assumption of $\scA : X \to \SIvertone {\cC}$, $\scB : Y \to \SIvertone {\cC}$, $H: X \to Y$ such that $\scA = \scB H$ and $\lambda : \scB \mcoll\widetilde \scB$). Then, the rest of the diagram commutes.
\end{cor}

\subsection{Push-forwards of collapses along base change}

In this final section, we will construct an analogue of the previous section's result for push-forwards along base changes. Again, we will later on apply this result to cube families, and this will be part of the proof of collapse normal forms being unique.

More precisely, consider the following diagram
\begin{restoretext}
\begingroup\sbox0{\includegraphics{test/page1.png}}\includegraphics[clip,trim=0 {.0\ht0} 0 {.1\ht0} ,width=\textwidth]{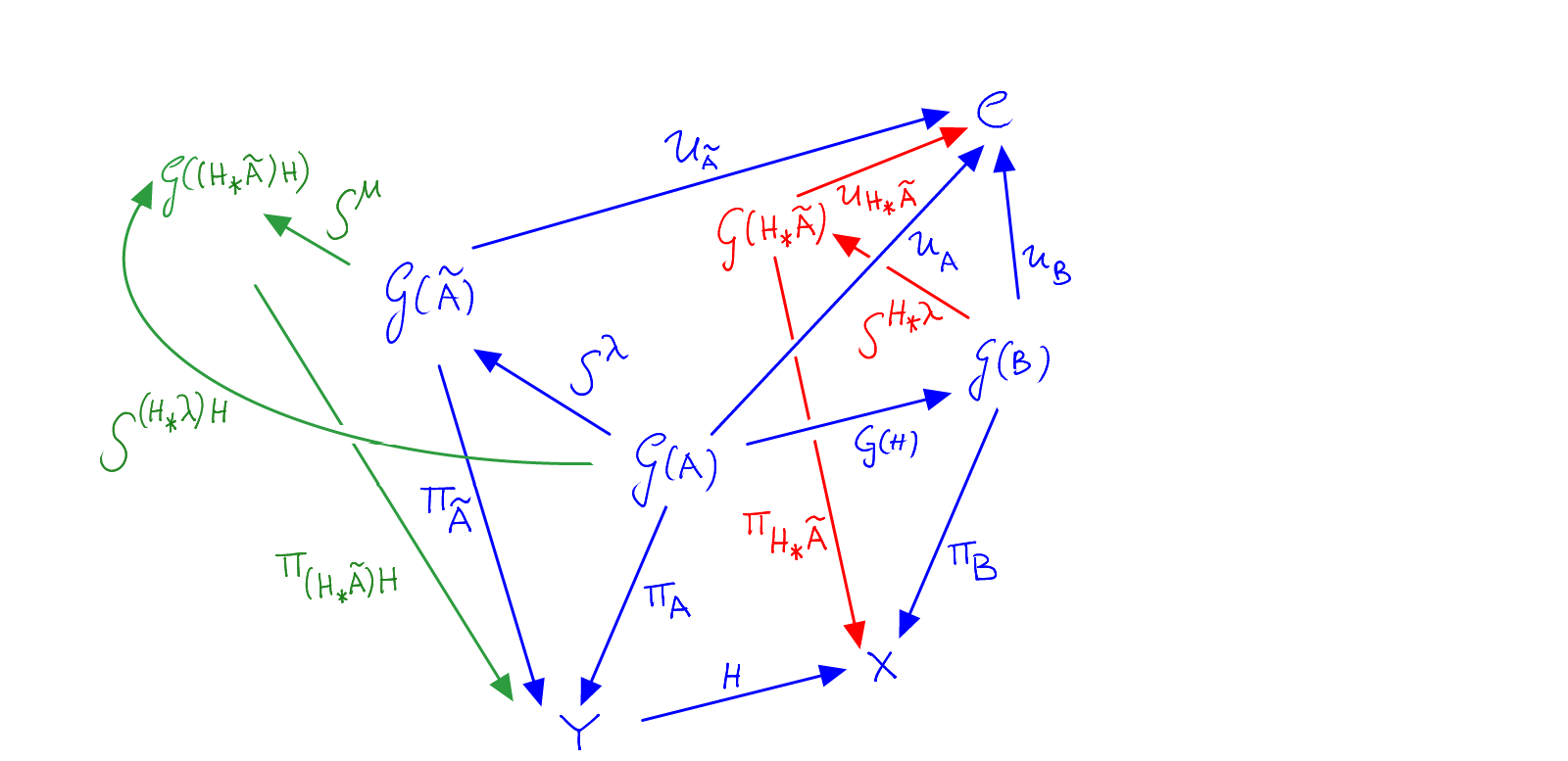}
\endgroup\end{restoretext}
We will show that under the assumption of the \cblue{} part of the diagram (subject to the condition that $H$ has \gls{lifts}), the \cred{} part can be constructed. In fact, as illustrated by the \cdarkgreen{} part, the constructed collapse $H\psstar \lambda  : \scB \to \widetilde\scB$ will be such that $(H\psstar \lambda)H$ factors through $\lambda$ by some $\mu$, and is in fact the ``minimal" collapse with that property (that is any other $\lambda'$ such that $\lambda' H$ factors through $\lambda$, already factors itself through $H\psstar  \lambda$). 

To illustrate the argument which will be used for this construction, we start with an example: consider the following pullback along $H$ of \SI-families $\und \scA$ and $\und \scB$
\begin{restoretext}
\begingroup\sbox0{\includegraphics{test/page1.png}}\includegraphics[clip,trim=0 {.0\ht0} 0 {.0\ht0} ,width=\textwidth]{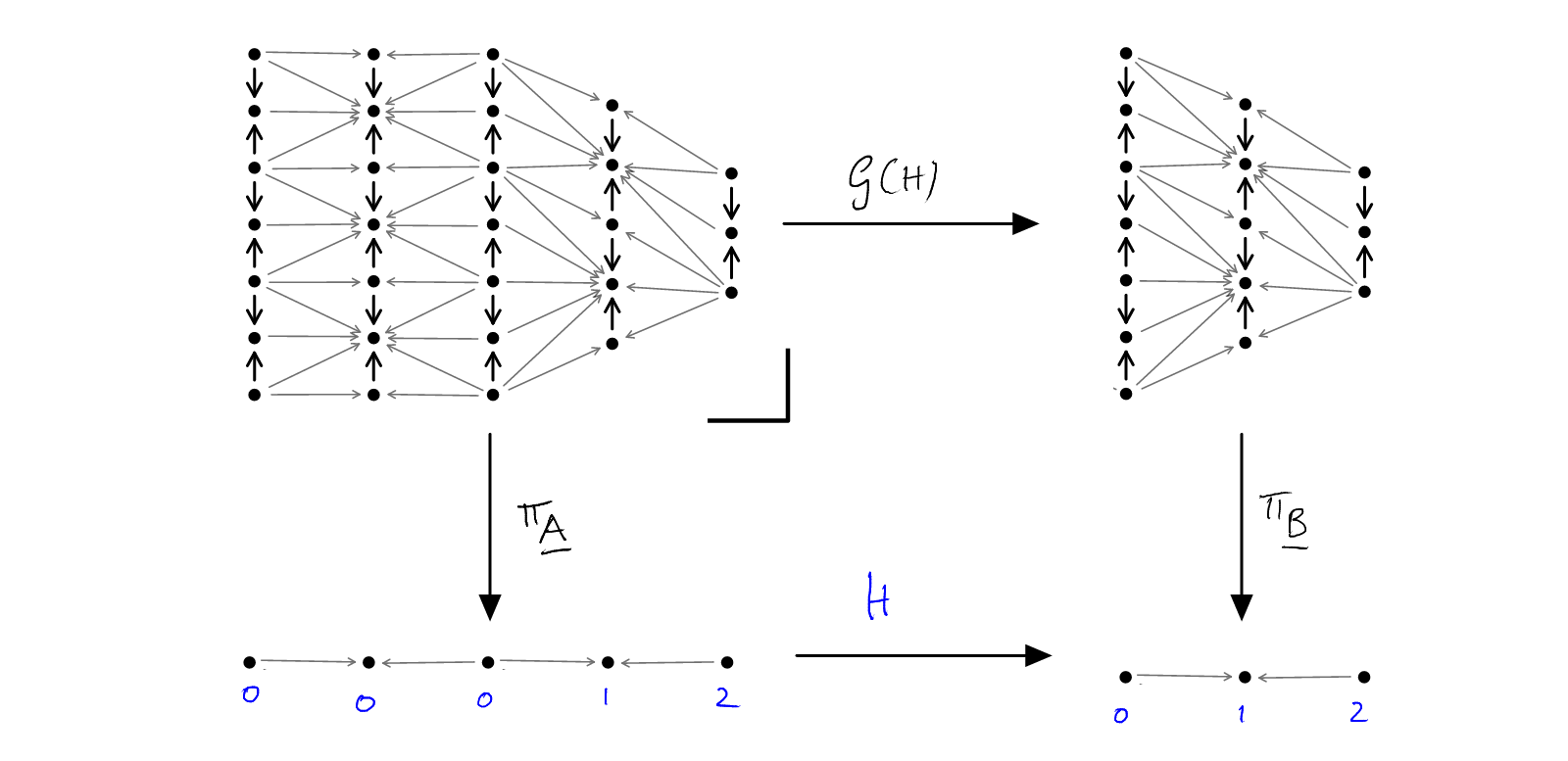}
\endgroup\end{restoretext}
Here $H$ is given by defining its preimages of image points (both marked in \cblue{}). Now we construct an injection $\lambda$ into $\und\scA$ by defining its stable singular subset section $\cF := \cS^\lambda$ to contain the following singular heights (marked by \cred{} circle)
\begin{restoretext}
\begingroup\sbox0{\includegraphics{test/page1.png}}\includegraphics[clip,trim=0 {.25\ht0} 0 {.25\ht0} ,width=\textwidth]{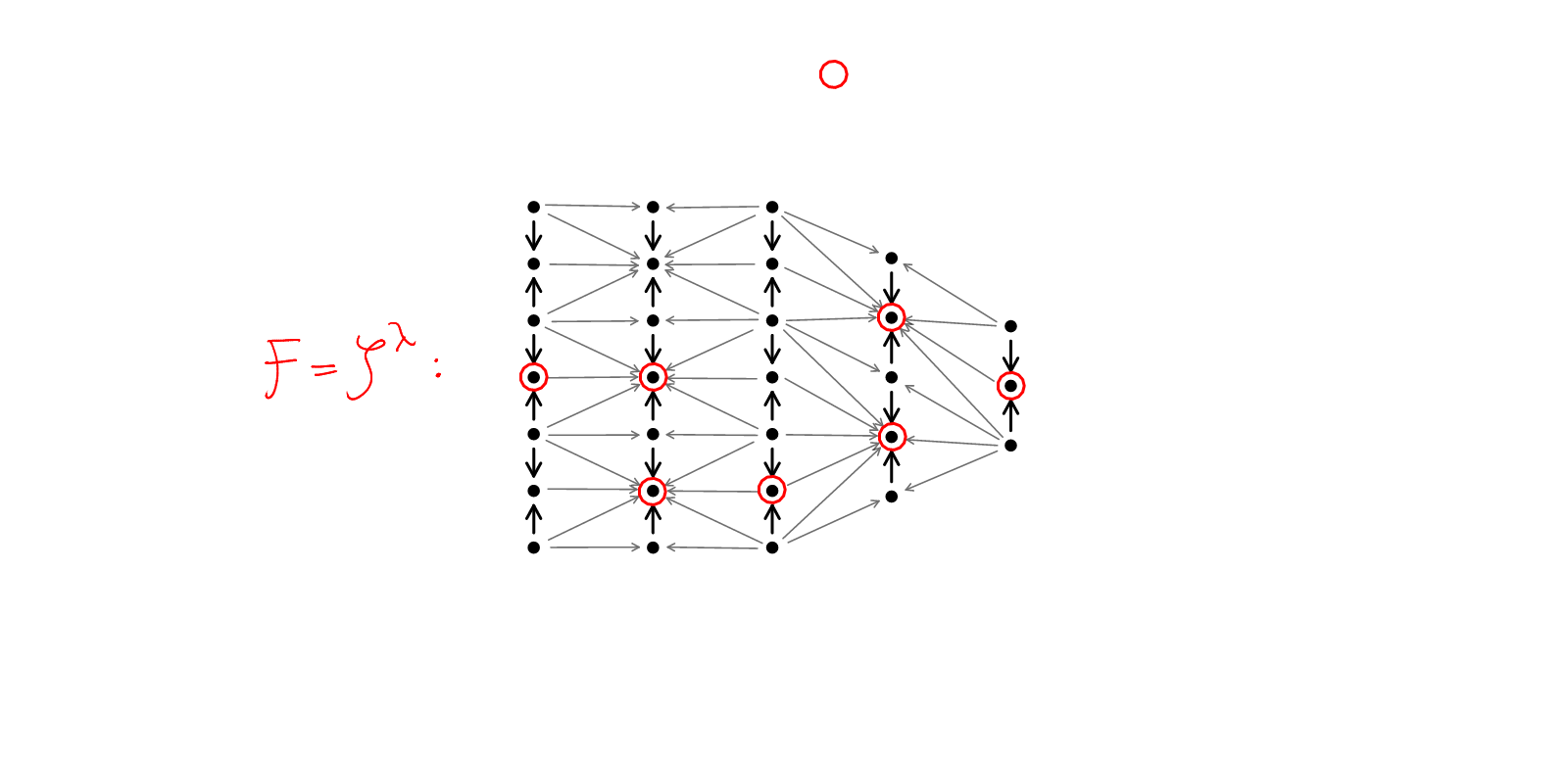}
\endgroup\end{restoretext}
We thus obtain a collapse functor $\sS^\lambda : \pi_{\und \scA} \to \pi_{\intrel{\cF}}$. We assume families $\scA, \widetilde \scA$ such that $\lambda : \scA \mcoll \widetilde\scA$. This can be visualised as
\begin{restoretext}
\begingroup\sbox0{\includegraphics{test/page1.png}}\includegraphics[clip,trim={.2\ht0} {.0\ht0} {.2\ht0} {.0\ht0} ,width=\textwidth]{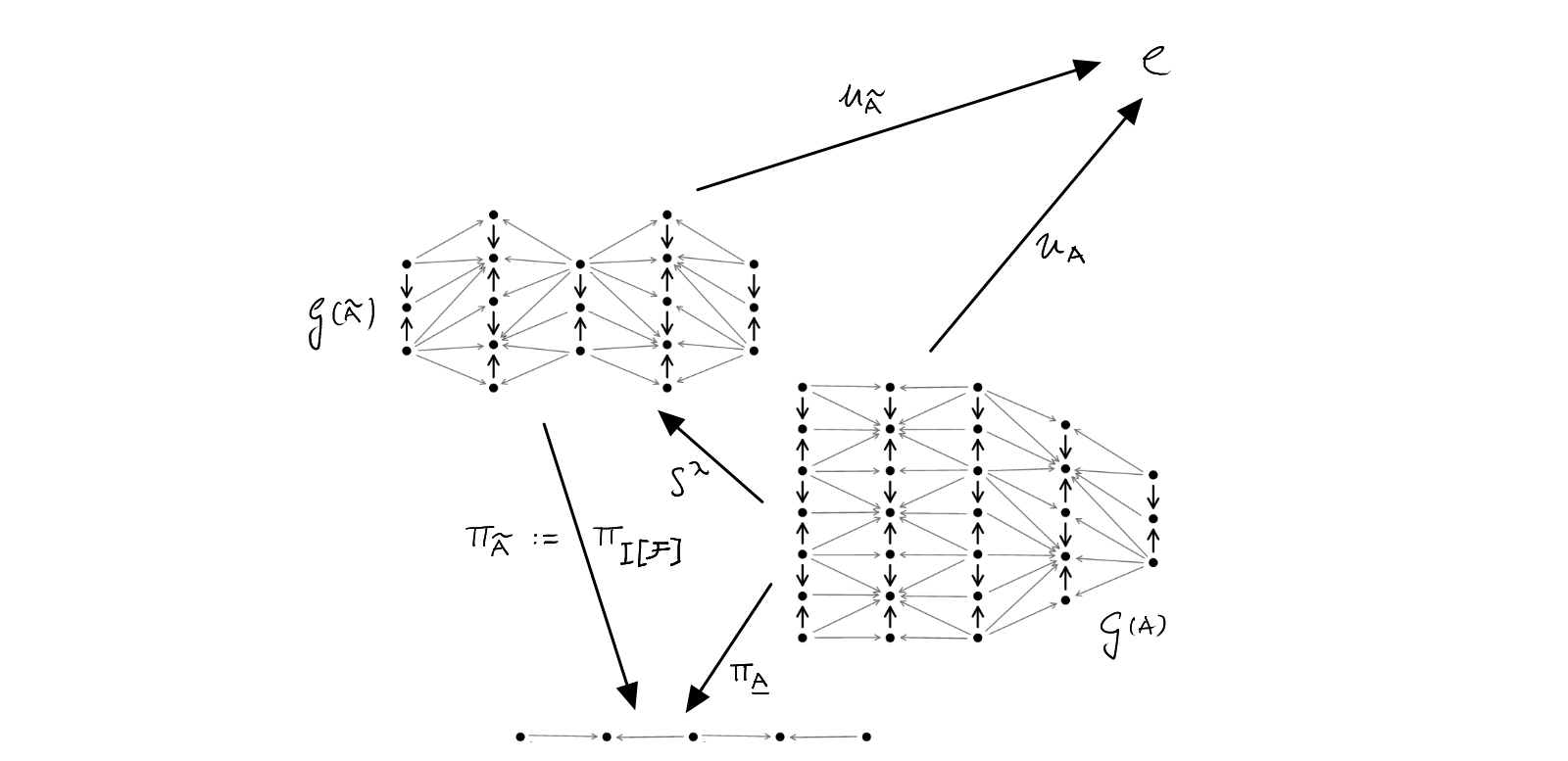}
\endgroup\end{restoretext}
We also assume a $\SIvertone \cC$-family $\scB$ such that $\scB H = \scA$ (in particular, this implies $\sU_\scB \sG(H) = \sU_\scA$ by \autoref{claim:grothendieck_span_construction_basechange}). 

The crucial step of the construction is now the following: $H\psstar  \lambda$ is constructed by defining its stable singular subset section $\cF' := \sS^{H\psstar  \lambda}$ as follows: $\cF'_y$ is the intersection of sets $\cF_x$ with $x \in H\inv (y)$. In our case this means $\cF'$ is given by the following singular heights (marked by \cblue{} circles)
\begin{restoretext}
\begingroup\sbox0{\includegraphics{test/page1.png}}\includegraphics[clip,trim=0 {.0\ht0} 0 {.0\ht0} ,width=\textwidth]{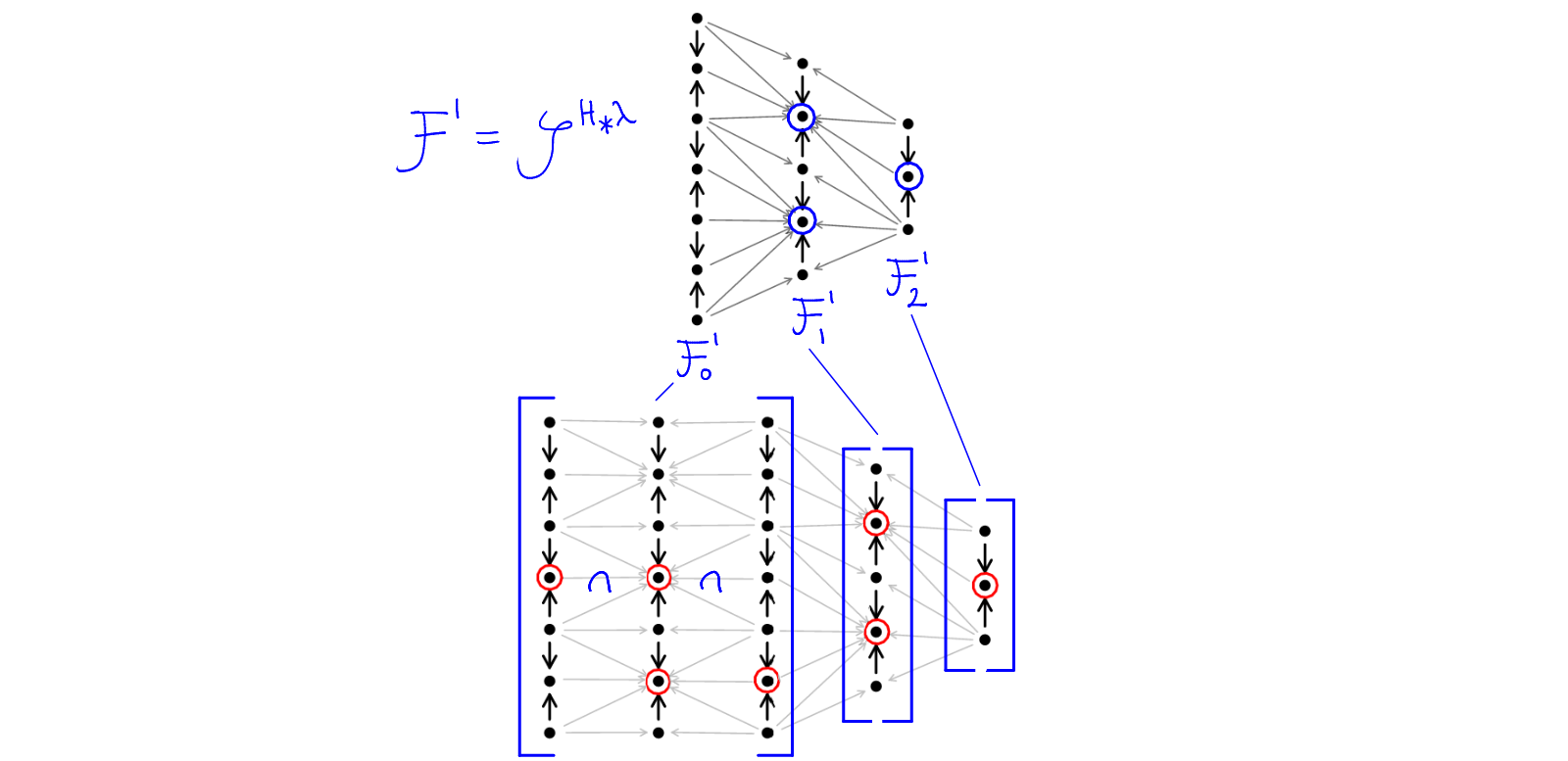}
\endgroup\end{restoretext}
Note that for instance $\cF'_0$ is the empty subset, since $H\inv(0) = \Set{0,1,2}$ and $\cF_0$, $\cF_1$ and $\cF_2$ intersect to give the empty subset. $\cF'$ can be shown to be stable, and thus we constructed an injection $H\psstar \lambda : \intrel{\cF'} \into \und \scB$. We now want to construct $H\psstar\widetilde\scA$ satisfying (colored morphisms are to be ignored for now)
\begin{restoretext}
\begingroup\sbox0{\includegraphics{test/page1.png}}\includegraphics[clip,trim={.1\ht0} {.0\ht0} {.1\ht0} {.0\ht0} ,width=\textwidth]{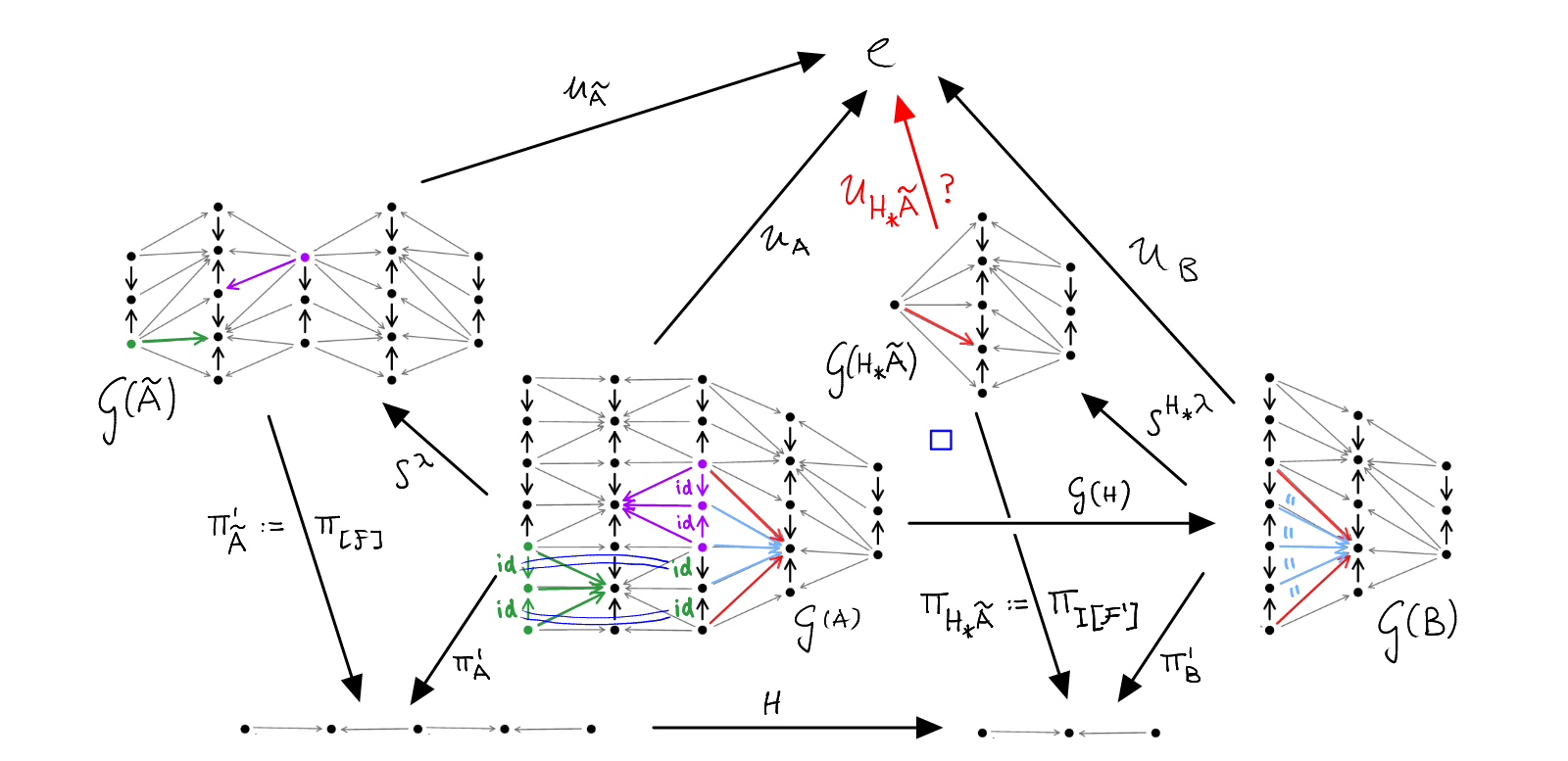}
\endgroup\end{restoretext}
And thus we need to define $\sU_{H\psstar\widetilde\scA}$ such that $\sU_{H\psstar\widetilde\scA} \sS^{H\psstar \lambda} = \sU_\scB$. The obstruction to such a definition is that morphisms in $\sG(H\psstar\widetilde\scA)$ might have multiple preimages under $\sS^{H\psstar \lambda}$. An example of this is depicted above: the \cred{} arrow in $\sG(H\psstar\widetilde\scA)$ has multiple preimages in $\sG(\scB)$, two of which are marked in \cred{}. For $\sU_{H\psstar\widetilde\scA}$ (such that $\sU_{H\psstar\widetilde\scA} \sS^{H\psstar \lambda} = \sU_\scB$) to exist we require that the value of $\sU_\scB$ agrees on both of these preimages. To show this, we once more argue by edge induction arguing that each arrow in between the \cred{} arrows (marked in \clightblue{}) agrees in value (of $\sU_\scB$) with its neighbouring arrow. This follows since it holds for corresponding arrows in $\sU_{\scA}$: Namely, preimages of the arrows in $\sG(\scB)$ along $\sG(H)$ (marked again in \cred{} and \clightblue{} in $\sG(\scA)$) agree in value (of $\sU_\scA)$). Indeed, we find that each pair of neighbouring (\clightblue{}/\cred{}) arrows agrees in value because their filler edge is the identity (marked by \id labels in \cdarkgreen{} an \cpurple{}). The latter statement follows in general from the choice of $H\psstar  \lambda$. For the above, note that the \cdarkgreen{} and \cpurple{} arrows are being collapsed into a single arrow by $\sS^\lambda$ which then forces the $\id$-labelled arrows to be mapped to $\id$ by $\sU_\scA$. Note that for the (\cblue{}) equality signs connecting \cdarkgreen{} identities we used that $\scA = \scB H$. 

This argument is generalised and made precise in the proof of the following construction.

\begin{constr}[Pushforward of injections] \label{claim:collapse_dimension_interaction2}
Let $\scA : X \to \SIvertone {\cC}$, $\scB : Y \to \SIvertone {\cC}$ be $\cC$-labelled singular $n$-cube families and let $H: X \to Y$ be a surjective functor which has \gls{lifts} and such that $\scA = \scB H$ (cf. \autoref{defn:having_lifts}). Assume $\widetilde{\scA} : X \to \SIvertone {\cC}$ and $\lambda : \scA \mcoll \widetilde \scA$. We will construct ${H\psstar \widetilde{\scA}} : Y \to {\SIvertone {\cC}}$ together with $H\psstar \lambda : \SIf{{H\psstar \widetilde{\scA}}} \into \SIf{\scB}$ such that $H\psstar  \lambda : \scB \mcoll \widetilde \scB$. Note that $\lambda$ is implicit in the notation of $H\psstar \widetilde{\scA}$.

By \autoref{constr:unpacking_collapse} we can define a $\cC$-labelled singular interval family ${H\psstar \widetilde{\scA}} : Y \to \SIvertone {\cC}$ by giving a pair consisting of a functor $\SIf{{H\psstar \widetilde{\scA}}} : Y \to \SI$ and a functor $\sU_{{H\psstar \widetilde{\scA}}} : \sG(\SIf{{H\psstar \widetilde{\scA}}}) \to \cC$.

We first define $\SIf{{H\psstar \widetilde{\scA}}}$. By assumption on $H$, we note that $\SIf{\scB}(y) = \SIf{\scA}(x)$ for any $x\in H\inv(y)$. Then define a singular subset section $\cF$ as follows
\begin{equation}
\cF_y := \bigcap_{x \in H\inv(y)} \cS^{\lambda}_x \subset \SIf{\scB}(y)
\end{equation}
In order to define an injection $H\psstar \lambda := \eta_\cF$ we use \autoref{lem:stability_vs_injections}: We need to verify that $\cF$ is a stable singular subset section of $\SIf{\scB}$. That is, whenever $y \to y'$ in $Y$ we need to show
\begin{equation} \label{eq:collapse_pushforward_defn}
\SIf{\scB}(y \to y')(\cF_y) \subset \cF_{y'}
\end{equation}
Assume by contradiction that there is $a \in \cF_y$ such that $a' := \SIf{\scB}(y \to y')(a) \notin \cF_{y'}$. By definition of $\cF$ there is $x' \in H\inv(y')$ such that $a' \notin \cS^\lambda_{x'}$. Since $H(x') = y'$, and by the lifting property of $H$, we find $x \to x'$ in $X$ such that $H(x \to x') = y \to y'$. Since $a \in \cF_y$ by assumption, we also have $a \in \cS^\lambda_x$ by definition of $\cF$. But then
\begin{align*}
\SIf{\scA}(x \to x')(a) &= \SIf{\scB} H (x \to x')(a) \\
&= \SIf{\scB} (y \to y')(a) \\
&=a' \notin \cS^\lambda_{x'}
\end{align*}
which means
\begin{equation}
\SIf{\scA}(x \to x')(\cS^\lambda_x) \not\subset \cS^\lambda_{x'}
\end{equation}
and thus contradicts stability of $\cS^\lambda$. Therefore, $\cF$ must be stable too and using \autoref{lem:stability_vs_injections} we can set $H\psstar  \lambda := \eta_\cF$ as well as $\SIf{{H\psstar \widetilde{\scA}}} := \intrel{\cS^{H\psstar \lambda}} = \intrel{\cF}$.

Next, we define $\sU_{{H\psstar \widetilde{\scA}}}$. Recall \autoref{defn:edge_sets} of edge sets and the notation of \eqref{eq:grothendieck_edges}. Given $(r_y : y_\sop \to y_\tap) \in \mor(Y)$ and  $(r_y,\widetilde{\edb}) \in \mor(\sG(\SIf{{H\psstar \widetilde{\scA}}}))$ we set
\begin{equation} \label{eq:collapse_pushforward_defn_2}
\sU_{{H\psstar \widetilde{\scA}}}(r_y,\widetilde{\edb}) := \sU_{\scB}(r_y,\edb)
\end{equation}
for any $(r_y,\edb)\in \mor(\sG(\SIf{\scB}))$ with $\sS^{H\psstar \lambda}(r_y,\edb) = (r_y,\widetilde{\edb})$ (note that at least one such preimage exists since $\sS^{H\psstar \lambda}$ is surjective on objects and morphisms by \autoref{lem:glambda}).

We need to show that this is well-defined, i.e. for $(r_y,\edb')\in \mor(\sG(\SIf{\scB}))$ with $\sS^{H\psstar \lambda}(r_y,\edb') = (r_y,\widetilde{\edb})$ we must have
\begin{equation} \label{eq:chi_G2_welldef}
\sU_{\scB}(r_y,\edb') = \sU_{\scB}(r_y,\edb)
\end{equation}
We argue similar to the proof of \autoref{claim:natural_injection_pushout}, namely, by induction on the distance 
\begin{equation}
n = \abs{\avg{\edb'} - \avg{\edb}}
\end{equation}
where $\avg{-}$ is the norm defined in \autoref{defn:edge_sets}. Note that by \autoref{claim:edge_set_properties}, if $n = 0$ we must have $\edb = \edb'$ and thus \eqref{eq:pushout_welldef} is satisfied. Without loss of generality assume $\avg{\edb} < \avg{\edb'}$ (otherwise switch their roles in the following). Recall \autoref{defn:successor} and the construction of the filler index $\succindex \edb$. Using our assumption that
\begin{equation}
\sS^{H\psstar \lambda}(y_{\succindex \edb},\edb_{\succindex \edb}) =\sS^{H\psstar \lambda}(y_{\succindex \edb},\edb'_{\succindex \edb}) = (y_{\succindex \edb},\widetilde{\edb}_{\succindex \edb})
\end{equation}
then, by applying monotonicity of $\rest {\sS^{H\psstar \lambda}} {y_{\succindex \edb}}$ to the inequalities of \autoref{cor:filler_edge} (obtained in the case of $\avg{\edb} < \avg{\edb'}$), we can infer
\begin{equation}
\sS^{H\psstar \lambda}(x_{\succindex \edb}, \succfill \edb\ssoe ) = \sS^{H\psstar \lambda}(x_{\succindex \edb}, \succfill \edb\ttae ) = (x_{\succindex \edb}, \widetilde{\edb}_{\succindex \edb})
\end{equation}
Since $\succfill \edb\ssoe $ is a regular segment, and $\rest {\sS^{H\psstar \lambda}} {x_{\succindex \edb}}$ preserves regular segments this forces $\widetilde{\edb}_{\succindex \edb}$ to be a regular segment. Thus $\succfill \edb\ttae  \notin \cS^{{H\psstar \lambda}}_{x_{\succindex \edb}} = \im({H\psstar \lambda}_{x_{\succindex \edb}})$, since otherwise $\widetilde{\edb}_{\succindex \edb}$ would have to be a singular height by \eqref{eq:glambda_defn}. 

Since $\succfill \edb\ttae  \notin \cS^{H\psstar  \lambda}_{y_{\succindex \edb}}$, by definition of $\cS^{H\psstar  \lambda}_{y_{\succindex \edb}}$ as the intersection of $\cS^{\lambda}_{x}$, $x \in H\inv(y_{\succindex \edb})$, there is $x_{\succindex \edb} \in H\inv(y_{\succindex \edb})$ such that $\succfill \edb\ttae  \notin\cS^{\lambda}_{x_{\succindex \edb}}$. Using \autoref{claim:collapse_sandwich} together with $\succfill \edb\ttae  \notin \cS^{\lambda}_{x_{\succindex \edb}}$ and $\succfill \edb\ssoe  = \succfill \edb\ttae  \pm 1$ (for one choice of $\Set{+,-}$) we find
\begin{equation}
\sS^{\lambda}(x_{\succindex \edb}, \succfill \edb\ssoe ) = \sS^{\lambda}(x_{\succindex \edb}, \succfill \edb\ttae )
\end{equation} 
and thus
\begin{equation}
\sS^{\lambda}(\id_{x_{\succindex \edb}}, \succfill \edb) = \id
\end{equation}
Since $\sU_{\scA} = \sU_{\widetilde{\scA}}\sS^\lambda$ we then deduce
\begin{equation}
\sU_{\scA}(\id_{x_{\succindex \edb}}, \succfill \edb) = \id
\end{equation}
Since $\sU_{\scA} = \sU_{\scB} \sG(H)$ and $H(\widetilde x) = y_{\succindex \edb}$ this implies 
\begin{equation}
\sU_{\scB}(\id_{x_{\succindex \edb}}, \succfill \edb) = \id
\end{equation}
Using \autoref{rmk:successor_compositionality} and functoriality of $\sU_{\scB}$ we find
\begin{equation}
\sU_{\scB}(r,\edb) = \sU_{\scB}(r,\succ \edb)
\end{equation}
Since $\abs{\avg{\succ \edb} - \avg{\edb'}} = n - 1$ this completes the inductive proof of well-definedness of \eqref{eq:chi_G2_welldef}. 

$\SIf{{H\psstar \widetilde{\scA}}}$ and $\sU_{{H\psstar \widetilde{\scA}}}$ together provide us with a $\cC$-labelled singular interval family ${H\psstar \widetilde{\scA}} : Y \to \SIvertone {\cC}$ that satisfies $\sU_{\scB} = \sU_{{H\psstar \widetilde{\scA}}} \sS^{H\psstar \lambda}$ by construction (cf. \eqref{eq:collapse_pushforward_defn_2}).

This completes the construction of $H\psstar\lambda$ and $H\psstar\scA$.
\end{constr}

The following remark explain the \cdarkgreen{} arrows in the diagram at the beginning of the section. 

\begin{rmk}\label{rmk:factorising_of_pull_push} Following the previous construction, we remark that $H\psstar \lambda H$ factors through $\lambda$ by a natural injection $\mu$: That is, there is $\mu :  \SIf {H\psstar \widetilde{\scA} H} \into \SIf {\widetilde{\scA}}$ such that $H\psstar \lambda H = \lambda \mu$: Indeed, recalling \autoref{claim:injection_factorisation} the statement that $H\psstar \lambda H$ factors through $\lambda$ by some $\mu$ follows from definition of $\cS^{H\psstar \lambda}$ together with \eqref{eq:lambda_pullback}. We then find
\begin{align} \sU_{\widetilde\scA} \sS^\lambda &= \sU_{\scA} \\
&= \sU_{H\psstar \widetilde \scA H} \sS^{H\psstar \lambda H} \\
&= \sU_{H\psstar \widetilde \scA H} \sS^{\lambda \mu} \\
&= \sU_{H\psstar \widetilde \scA H} \sS^\mu \sS^\lambda
\end{align}
By surjectivity of $ \sS^\lambda$ we deduce
\begin{equation}
\sU_{\widetilde\scA} = \sU_{H\psstar \widetilde \scA H} \sS^\mu
\end{equation}
making $\mu :\widetilde \scA \mcoll H\psstar \widetilde \scA H$ a collapse, and thus
\begin{equation}
\xymatrix{\scA \ar@{~>}@/^2pc/[rr]^-{H\psstar \lambda H} \ar@{~>}[r]^-\lambda & \widetilde \scA \ar@{~>}[r]^-\mu & H\psstar \widetilde \scA H}
\end{equation}
is a commuting triangle in $\Buncoll \cC$. The minimality property of $H\psstar\lambda$ (claimed in the beginning of the section) is now straight-forward and left to the reader.
\end{rmk}

\chapter{Collapse of cubes} \label{ch:collapse_cubes}

This chapter introduces collapse for labelled singular $n$-cube families, and proves that this gives rise to a notion of normalisation. Many definitions and proofs will be direct analogues of statements from the previous chapters. In \autoref{sec:coll_cube} we give the definitions of collapse, and its reformulation as a base change. We also define normal forms. In \autoref{sec:coll_norm} we will then proceed to show that normal forms exist uniquely. The central result of this (and the previous) chapter will be \autoref{thm:normal_forms_unique}. Finally, in \autoref{sec:coll_multi} we generalise ($k$-level) collapse to multi-level collapse. This will then finally recover the notion of collapse of labelled singular $n$-cubes that was introduced in the summary in \autoref{sec:sum_norm}.

\section{Collapse of labelled singular cube families} \label{sec:coll_cube}

\subsection{$k$-Level collapse}

Using \autoref{defn:collapse_of_SI_families}, we can make the following definition.

\begin{defn}[$k$-level collapse] \label{defn:singular_cube_families} Let $\scA, \scB : X \to \SIvert n \cC$ be $\cC$-labelled singular $n$-cube families. We say \textit{$\lambda$ is a witness of $\scA$ collapsing to $\scB$ at level $k$}, written $\lambda : \scA \kcoll k \scB$ or $\xymatrix@1{\scA\ar@{~>}[r]^{\lambda}_k & \scB}$, if
\begin{equation}
 \text{$\sT_\scA^{k-1} = \sT^{k-1}_\scB$ and $\lambda : \tsU {k-1}_\scA \mcoll \tsU {k-1}_\scB$}
 \end{equation} 
We say \textit{$\scA$ collapses to $\scB$ at level} $k$ (or alternatively, that $\scB$ is a $k$-level collapse of $\scA$), written $\scA \kcoll k \scB$, if there is some $\lambda$ such that $\lambda : \scA \kcoll k \scB$.
\end{defn}

\begin{rmk}[Unwinding the definition of $k$-level collapse]\label{rmk:k_collapse_defn} Note that in the preceding definition $\tsU {k-1}_{\scA}$ and $\tsU {k-1}_{\scB}$ are $\SIvert{n-k} \cC$-families 
\begin{equation}
\tsU {k-1}_{\scA}, \tsU {k-1}_{\scB} : \tsG {k-1}(\scA) = \tsG {k-1}(\scB) \to \SIvert {n-k} \cC
\end{equation}
and thus they are in particular $\SIvertone  {\SIvert {n-k-1} \cC}$-families. From \autoref{defn:complete_unpacking_and_repacking} recall that we then have $\sU_{\tsU {k-1}_\scA} = \tsU k_{\scA}$ and $\sU_{\tsU {k-1}_\scB} = \tsU k_{\scB}$. Thus the condition $\lambda : \tsU {k-1}_\scA \mcoll \tsU {k-1}_\scB$ means that the following equation holds
\begin{equation} \label{eq:k_collapse_factorisation}
\big(\tsG {k}(\scA) \xto {\tsU {k}_{\scA}} \SIvert {n-k} \cC \big) \quad = \quad \big( \tsG {k}(\scA) \xto {\sS^{\lambda}} \tsG {k}(\scB) \xto {\tsU {k}_{\scB}} \SIvert {n-k} \cC \big)
\end{equation}
Here, $\sS^{\lambda} : \sG(\tusU {k-1}_{\scA}) \to \sG(\tusU {k-1}_{\scB})$ was constructed in \autoref{defn:glambda}.
\end{rmk}

\begin{rmk}[Relabelling for $k$-level collapse] \label{rmk:k_coll_label_transfer} Using \autoref{rmk:recoloring_cubes}, we find that whenever $\lambda : \scA \kcoll k \scB$ for $\scA, \scB : X \to \SIvert n \cC$ and $F : \cC \to \cD$, then also
\begin{equation}
\lambda : \SIvert n F \scA \kcoll k \SIvert n F \scB
\end{equation}
Conversely, if $F$ is injective on objects and morphisms then we find
\begin{equation}
\big(\lambda : \scA \kcoll k \scB\big) ~\iff~ \big( \lambda : \SIvert n F \scA \kcoll k \SIvert n F \scB \big)
\end{equation}
with the $\limp$ direction now following from relabelling along a partial inverse $F\inv$ of $F$ (defined on its image).
\end{rmk}

\subsection{$k$-Level collapse is a $k$-level base change}

We make the following observation.

\begin{rmk}[$k$-level collapse is a $k$-level base change] \label{rmk:coll_is_basechange} Comparing the previous remark to \autoref{constr:unpacking_collapse} we see that $\lambda : \scA \kcoll k \scB$ is a special case of a $k$-level base change: the foregoing definition of $k$-level collapse $\scA \kcoll k \scB$ can equivalently be stated as a $k$-level base change along a map of the form $\sS^\lambda$ for some $\lambda : \tusU {k-1}_{\scB} \into \tusU {k-1}_\scA$. This in turn has an equivalent formulation in terms of pullbacks by \autoref{lem:unpacking_collapse}. We illustrate the latter as follows: for $\scA, \scB : X \to \SIvert n \cC$, the statement $\lambda : \scA \kcoll k \scB$ is equivalent to the existence of the following diagram
\begin{restoretext}
\begingroup\sbox0{\includegraphics{test/page1.png}}\includegraphics[angle=90, origin=c,clip,trim=0 {.0\ht0} 0 {.0\ht0} ,width=7cm]{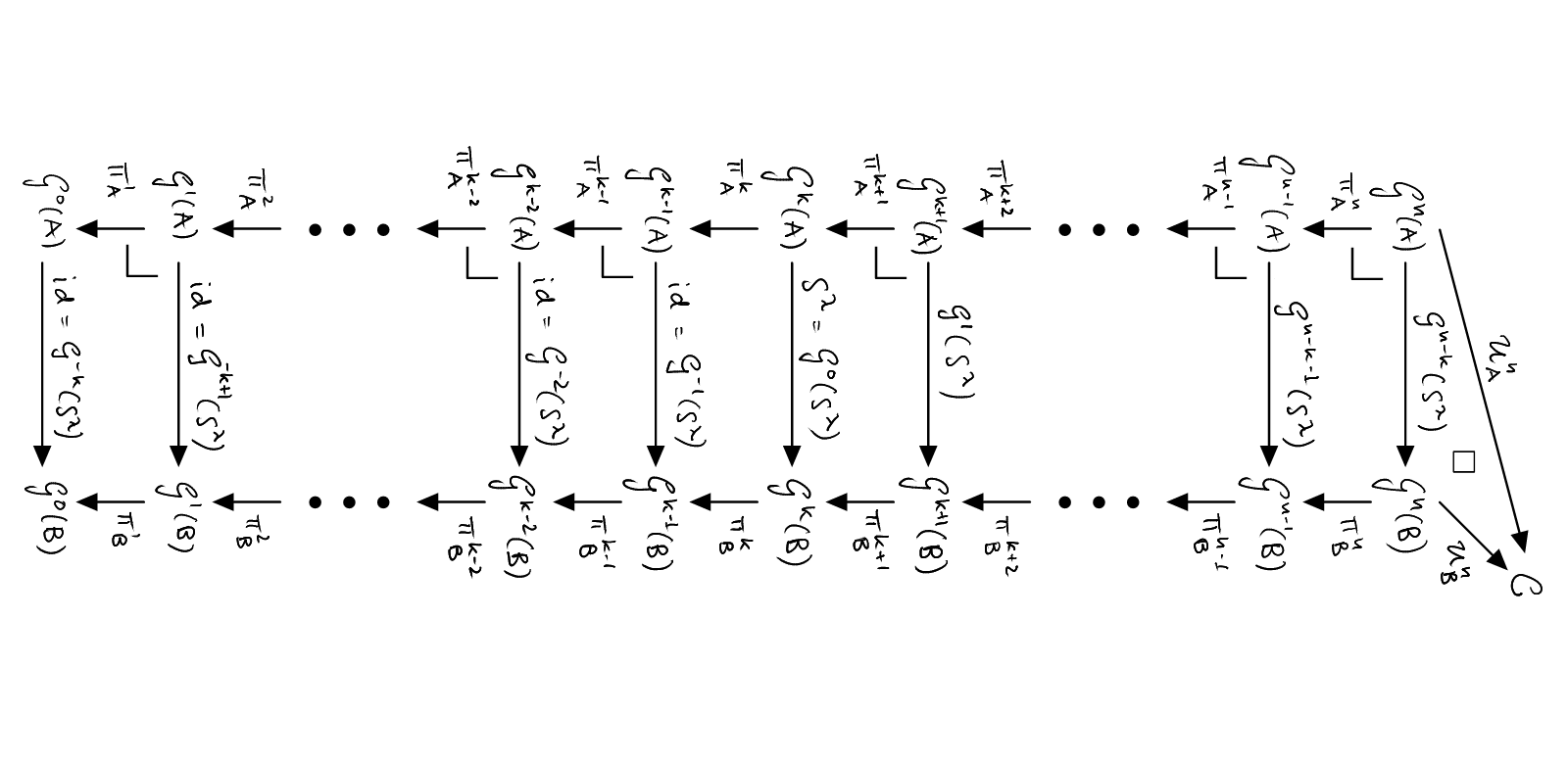}
\endgroup\end{restoretext}
\end{rmk}

\subsection{Normal forms}

We now define a notion of normal forms of a cube family. These will be the end terms of maximal sequences of collapses of that family.

\begin{defn} \label{defn:collapse_normal_form} Let $\scA, \scB: X \to \SIvert n \cC$ be $\cC$-labelled singular $n$-cube families. 
\begin{itemize}
\item We say $\scA$ \textit{collapses} to $\scB$, written $\scA \starcoll \scB$, if there is a finite sequence 
\begin{equation}
\scA = \scA_0 \kcoll {k_1} \scA_1 \kcoll {k_2} \scA_2 \kcoll {k_2} ... \kcoll {k_m} \scA_m = \scB
\end{equation}
of $k_i$-collapses, $0 < k_i \leq n$, which starts at $\scA$ and ends in $\scB$.
\item For a given $k$, $0 < k \leq n$, we further say $\scB$ is in \textit{$k$-level normal form}, if no non-identity $k$-level collapse applies to $\scB$. That is, there is no $\lambda : \scB \kcoll k \scB'$, for any choice of $\scB'$, $\lambda$, but the trivial one: namely, $\scB' = \scB$ and $\lambda = \id_{{\tusU {k-1}_{\scB}}}$.
\item Finally, we say $\scB$ is in \textit{normal form (up to level $n$)} if no non-identity $k$-level collapse for any $k$, $0 < k \leq n$ applies to $\scB$. The set of of normal forms that a given $\scA$ collapses to is denoted by $\NF{\scA}^n$:
\begin{equation}
\NF{\scA}^n = \Set{\scB : X \to \SIvert n \cC ~|~ \scA \starcoll \scB \text{~and $\scB$ is in normal form up to level $n$}}
\end{equation}
If $n = 1$, then $\NF{\scA}^1$ is sometimes just denoted by $\NF{\scA}$.
\end{itemize}  
\end{defn}

\begin{eg}[Normal forms] \label{eg:normal_forms} Let $\cC$ be the poset
\begin{restoretext}
\begingroup\sbox0{\includegraphics{test/page1.png}}\includegraphics[clip,trim=0 {.45\ht0} 0 {.3\ht0} ,width=\textwidth]{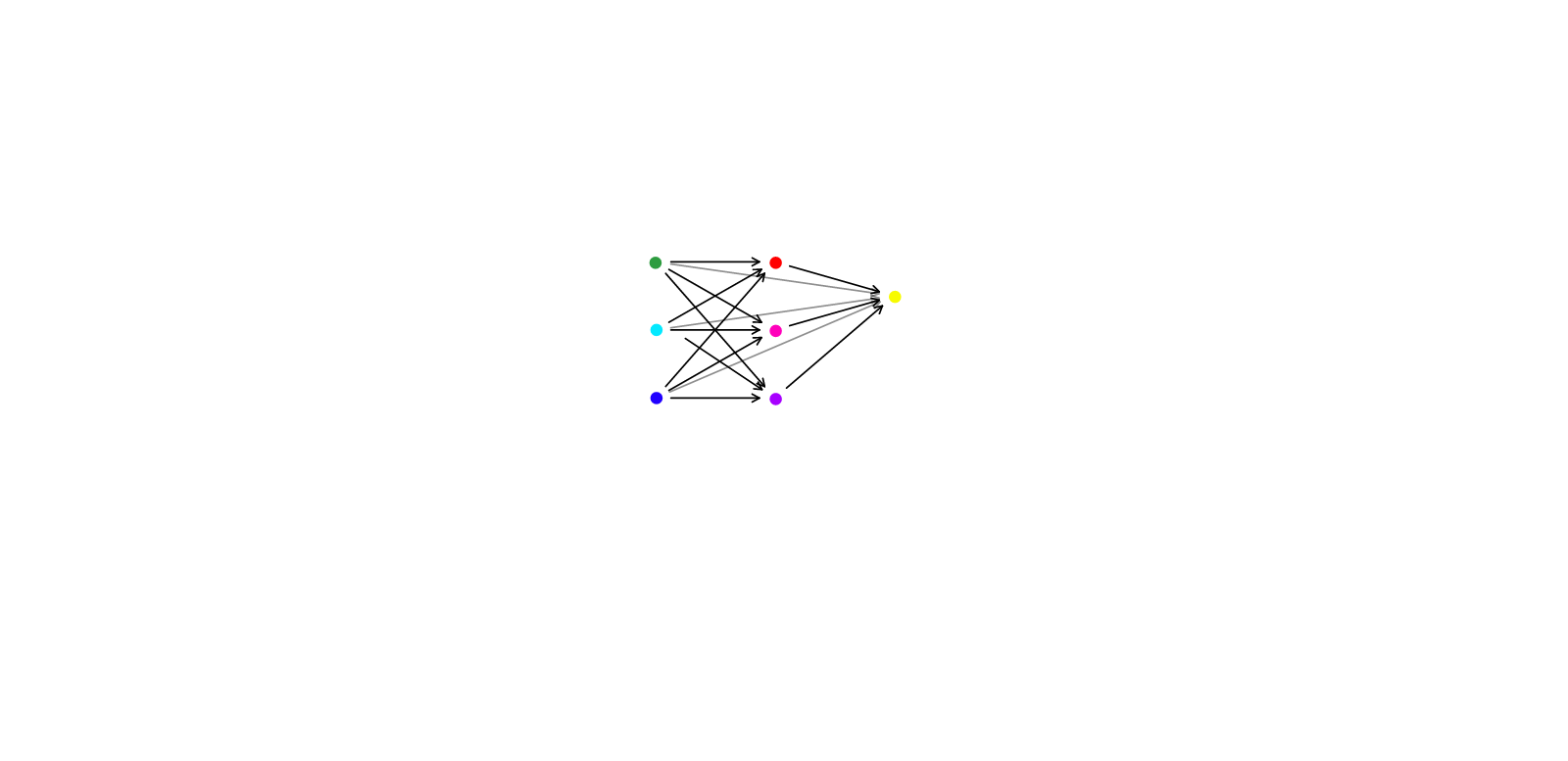}
\endgroup\end{restoretext}
Let $\scA$ be the $\SIvert 2 \cC$ family defined by
\begin{restoretext}
\begingroup\sbox0{\includegraphics{test/page1.png}}\includegraphics[clip,trim=0 {.0\ht0} 0 {.0\ht0} ,width=\textwidth]{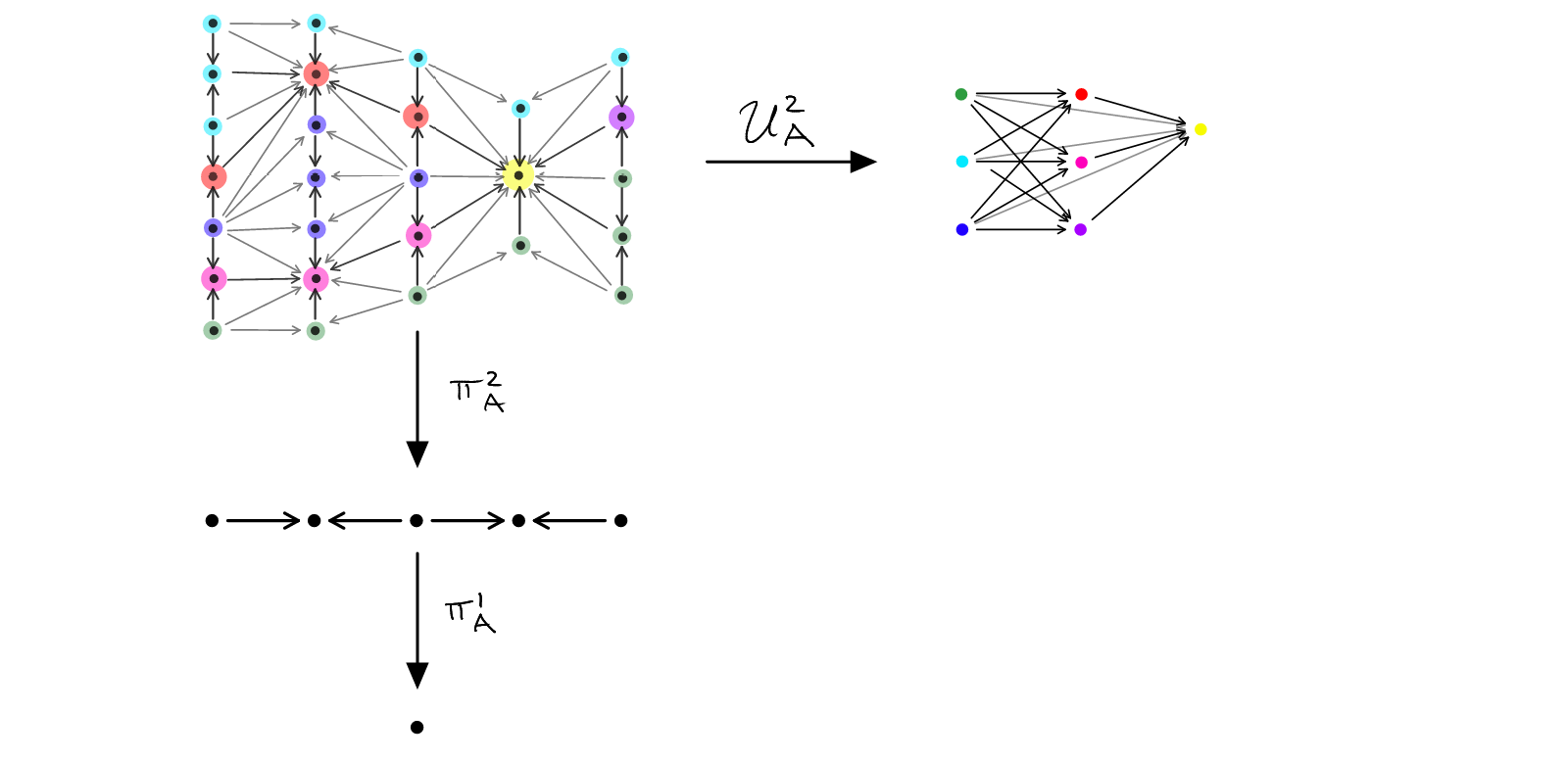}
\endgroup\end{restoretext}
Consider the singular subset section $\cF$ of $\tusU 1_\scA$ defined to include the following singular heights marked in \cblue{}
\begin{restoretext}
\begingroup\sbox0{\includegraphics{test/page1.png}}\includegraphics[clip,trim=0 {.1\ht0} 0 {.1\ht0} ,width=\textwidth]{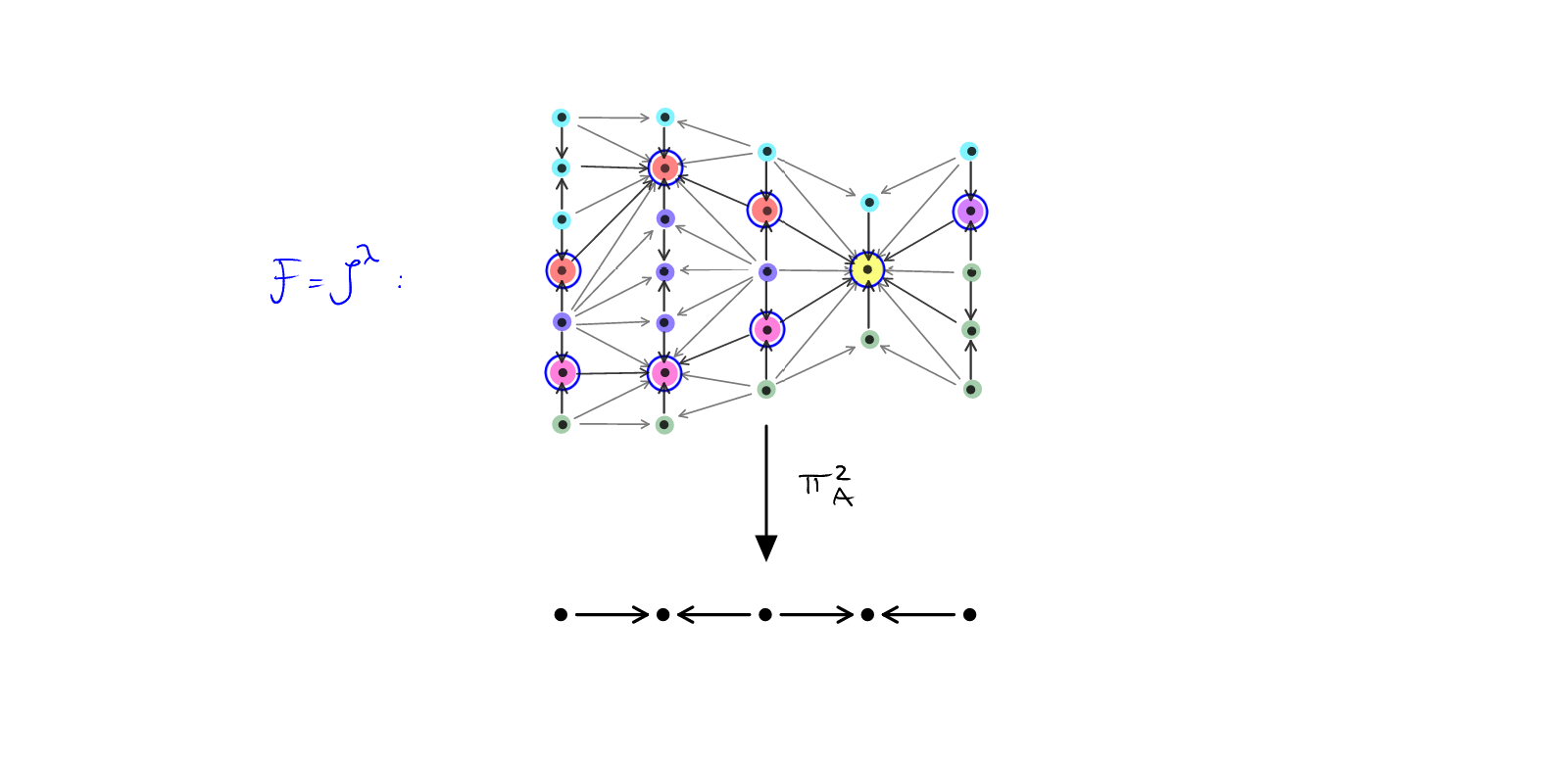}
\endgroup\end{restoretext}
$\cF$ is stable and thus there is an injection $\lambda : \intrel{\cF} \into \tusU 1_\scA$. By choice of $\lambda$ and $\scA$ there is $\scB$ such that $\lambda : \scA \kcoll 2 \scB$ defined as follows (cf. \autoref{rmk:coll_is_basechange})
\begin{restoretext}
\begin{noverticalspace}
\begingroup\sbox0{\includegraphics{test/page1.png}}\includegraphics[clip,trim=0 {.0\ht0} 0 {.4\ht0} ,width=\textwidth]{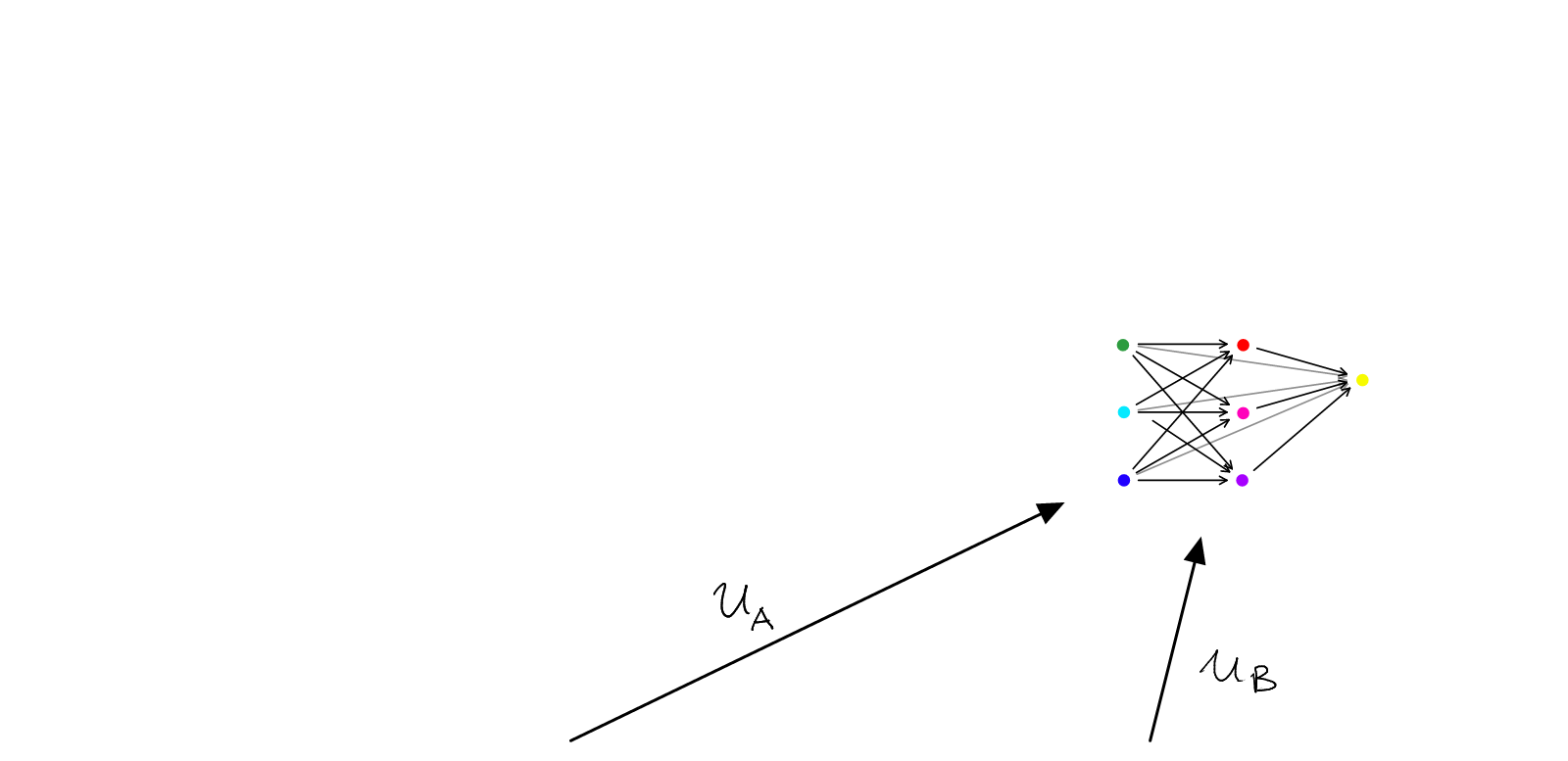}
\endgroup \\*
\begingroup\sbox0{\includegraphics{test/page1.png}}\includegraphics[clip,trim=0 {.0\ht0} 0 {.0\ht0} ,width=\textwidth]{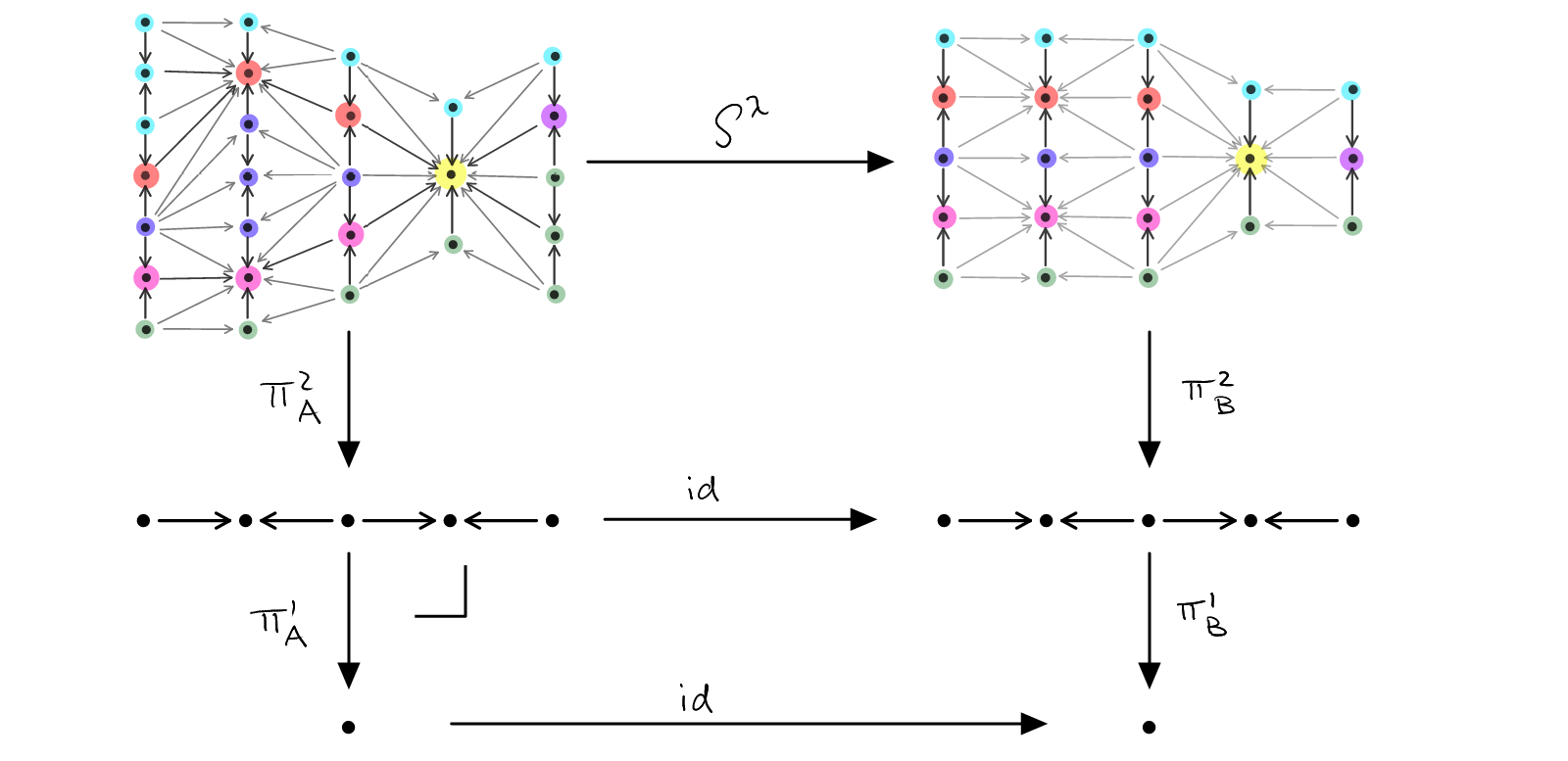}
\endgroup
\end{noverticalspace}
\end{restoretext}
Now $\scB$ can be seen to be in $2$-level normal form. However, it is not in $1$-level normal form, as we will next construct a non-identity $\mu : \scB \kcoll 1 \scC$ as follows. $S^\mu = \cF'$ is the stable singular subset section of $\tusU 0_\scB$ with the single element
\begin{restoretext}
\begingroup\sbox0{\includegraphics{test/page1.png}}\includegraphics[clip,trim=0 {.3\ht0} 0 {.35\ht0} ,width=\textwidth]{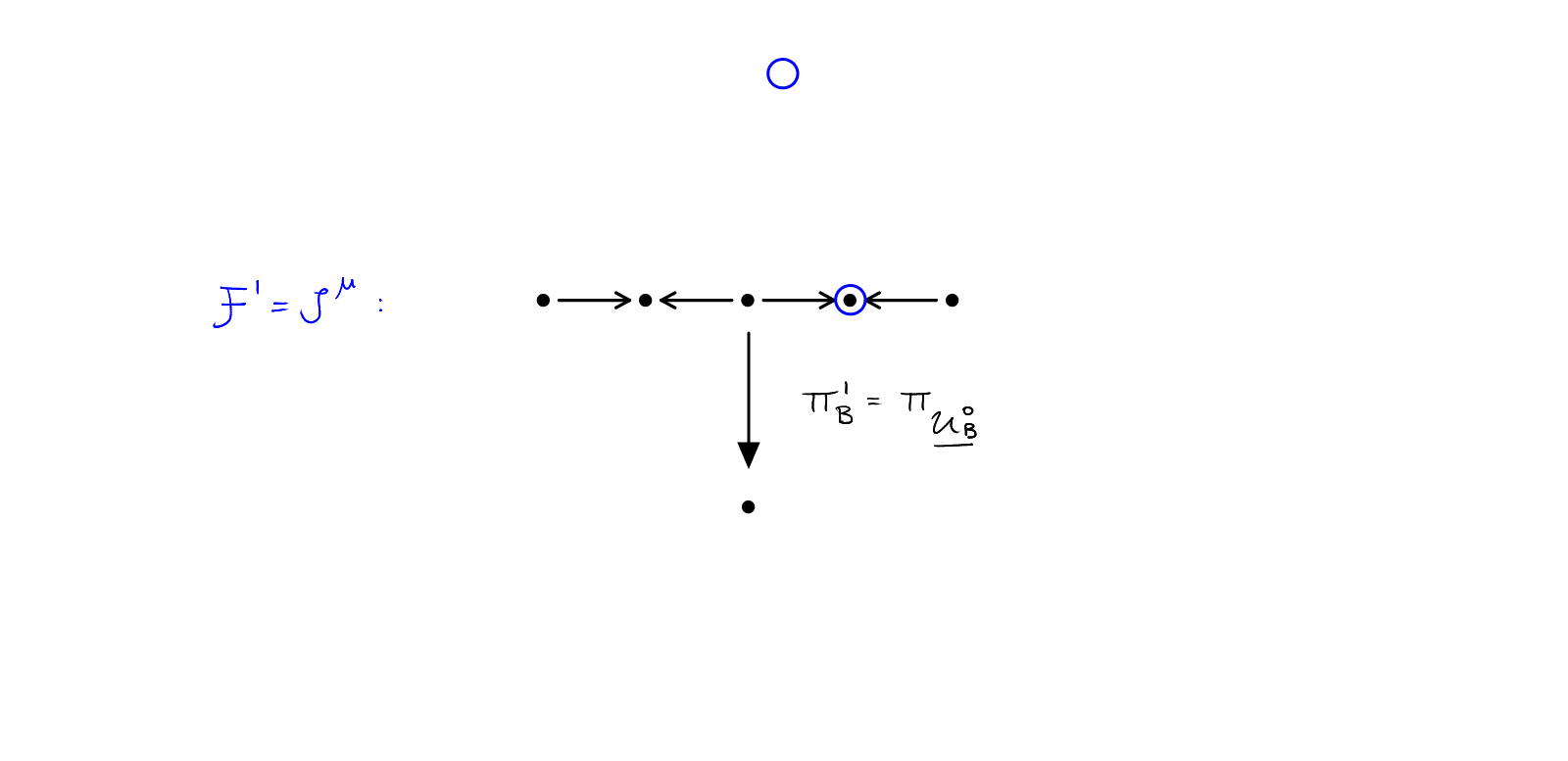}
\endgroup\end{restoretext}
This can then be seen to induce a $1$-level collapse  $\mu : \scB \kcoll 1 \scC$ with the following data
\begin{restoretext}
\begin{noverticalspace}
\begingroup\sbox0{\includegraphics{test/page1.png}}\includegraphics[clip,trim=0 {.0\ht0} 0 {.4\ht0} ,width=\textwidth]{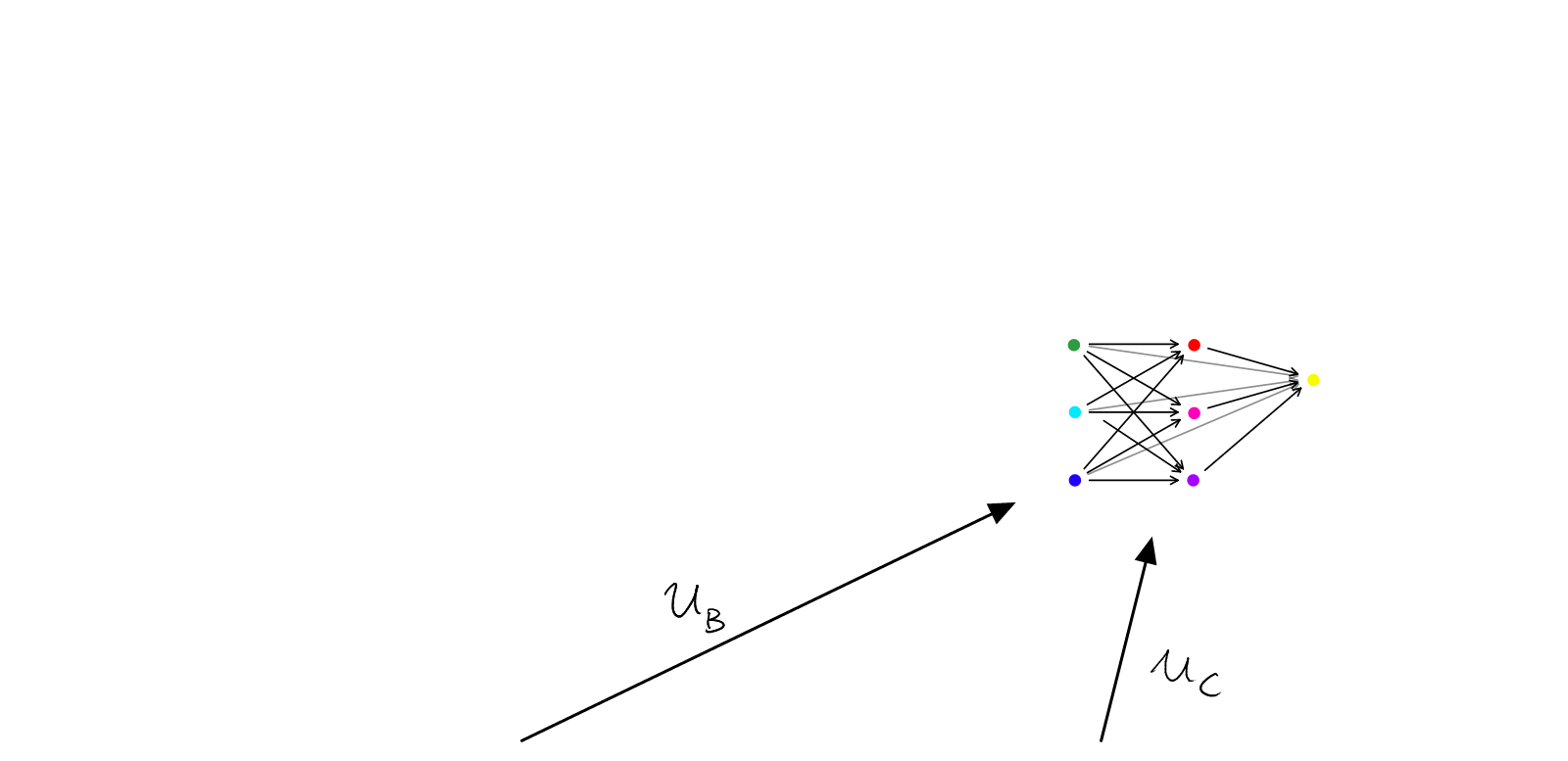}
\endgroup \\*
\begingroup\sbox0{\includegraphics{test/page1.png}}\includegraphics[clip,trim=0 {.0\ht0} 0 {.0\ht0} ,width=\textwidth]{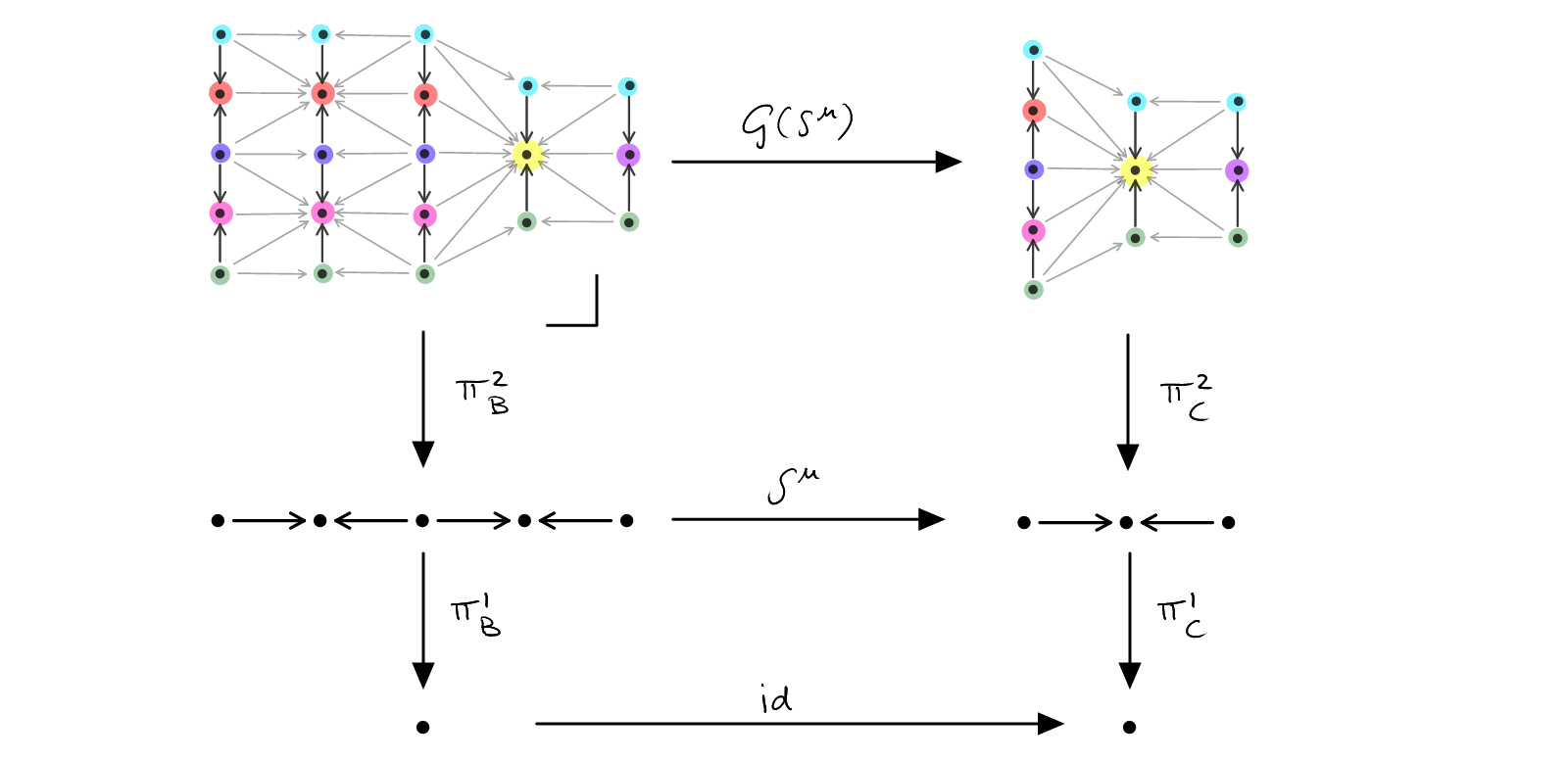}
\endgroup
\end{noverticalspace}
\end{restoretext}
We observe that now $\scC$ is in normal, as no other $1$-level or $2$-level collapse applies to it.
\end{eg}

\section{Properties of normalisation} \label{sec:coll_norm}

The goal of this section is to show existence and uniqueness of normal forms.

\subsection{Existence of collapse normal forms}

Note that for a given $\scA : X \to \SIvert n \cC$ the set $\NF{\scA}^n$ of normal forms is non-empty and finite, since each non-identity $k$-level collapse $\lambda : \scA \kcoll k \scB$ properly reduces the number of objects in $\tsG n(\scA)$, that is $\#\obj(\tsG n(\scA)) > \#\obj(\tsG n(\scB))$ and $\#\obj(\tsG n(\scA))$ is finite.

\subsection{Uniqueness of collapse normal form} We will see now that in fact $\NF{\scA}^n$ is a singleton set. We will thus refer to $\NF{\scA}^n : X \to \SIvert n \cC$ as $\cC$-labelled singular $n$-cube over $X$ itself, called the (unique) \textit{normal form of $\scA$}. Conversely, we say $\scA$ normalises to $\NF{\scA}^n$.

For the following lemmas, let $n \in \lN$ and $1 \leq k < l \leq n$. Let $\scC, \scC_i$ ($i \in \lN$) denote $\cC$-labelled singular $n$-cube families over a poset $X$.

\begin{lem}[$l$-level collapses compose] \label{lem:coll_comp} If
\begin{equation}
\xymatrix{ C_1 \ar@{~>}[r]^{\lambda_1}_l & C_2 \ar@{~>}[r]^{\lambda_2}_l & C_3}
\end{equation}
then
\begin{equation}
\xymatrix{ C_1 \ar@{~>}[r]^{\lambda_1 \lambda_2}_l & C_3}
\end{equation}
\proof \autoref{cor:composition_of_Glambda} implies $\sS^{\lambda_1 \lambda_2} = \sS^{\lambda_2} \sS^{\lambda_1}$. Thus we have
\begin{equation}
\sT^{l-1}_{\scC_1} = \sT^{l-1}_{\scC_2} = \sT^{l-1}_{\scC_3} 
\end{equation}
and (cf. \autoref{constr:cat_of_bun_and_coll})
\begin{equation}
\sU_{\scC_1} = \sU_{\scC_2} \sS^{\lambda_1} = \sU_{\scC_3} \sS^{\lambda_2} \sS^{\lambda_1} = \sU_{\scC_3} \sS^{\lambda_1 \lambda_2}
\end{equation}
which together imply $\lambda_1 \lambda_2 : \scC_1 \kcoll l \scC_3$ as required.
\qed
\end{lem}

\begin{lem}[$k$-level pullback of $l$-level collapse] \label{lem:coll_pull}
If
\begin{equation}
\xymatrix{ & C_3 \\
C_1  \ar@{~>}[r]^{\mu}_k & C_2 \ar@{~>}[u]^{\lambda}_l}
\end{equation}
then we have
\begin{equation}
\xymatrix{ C_4 \ar@{~>}[r]^{\mu}_k & C_3 \\
C_1 \ar@{~>}[u]^{\mu\pbstar \lambda}_l \ar@{~>}[r]^{\mu}_k & C_2 \ar@{~>}[u]^{\lambda}_l}
\end{equation}
where 
\begin{equation}
\mu\pbstar \lambda = \lambda \sG^{l-k-1}(\sS^\mu)
\end{equation}
$\mu\pbstar \lambda$ is called the \emph{$k$-level pullback} of $\lambda$ along $\mu$.

\proof The \textit{assumptions} in the claim are equivalent to the existence of the \textit{blue} part of the following diagram
\begin{restoretext}
\begingroup\sbox0{\includegraphics{test/page1.png}}\includegraphics[angle=90, origin=c,clip,trim=0 {.0\ht0} 0 {.0\ht0} ,width=9cm]{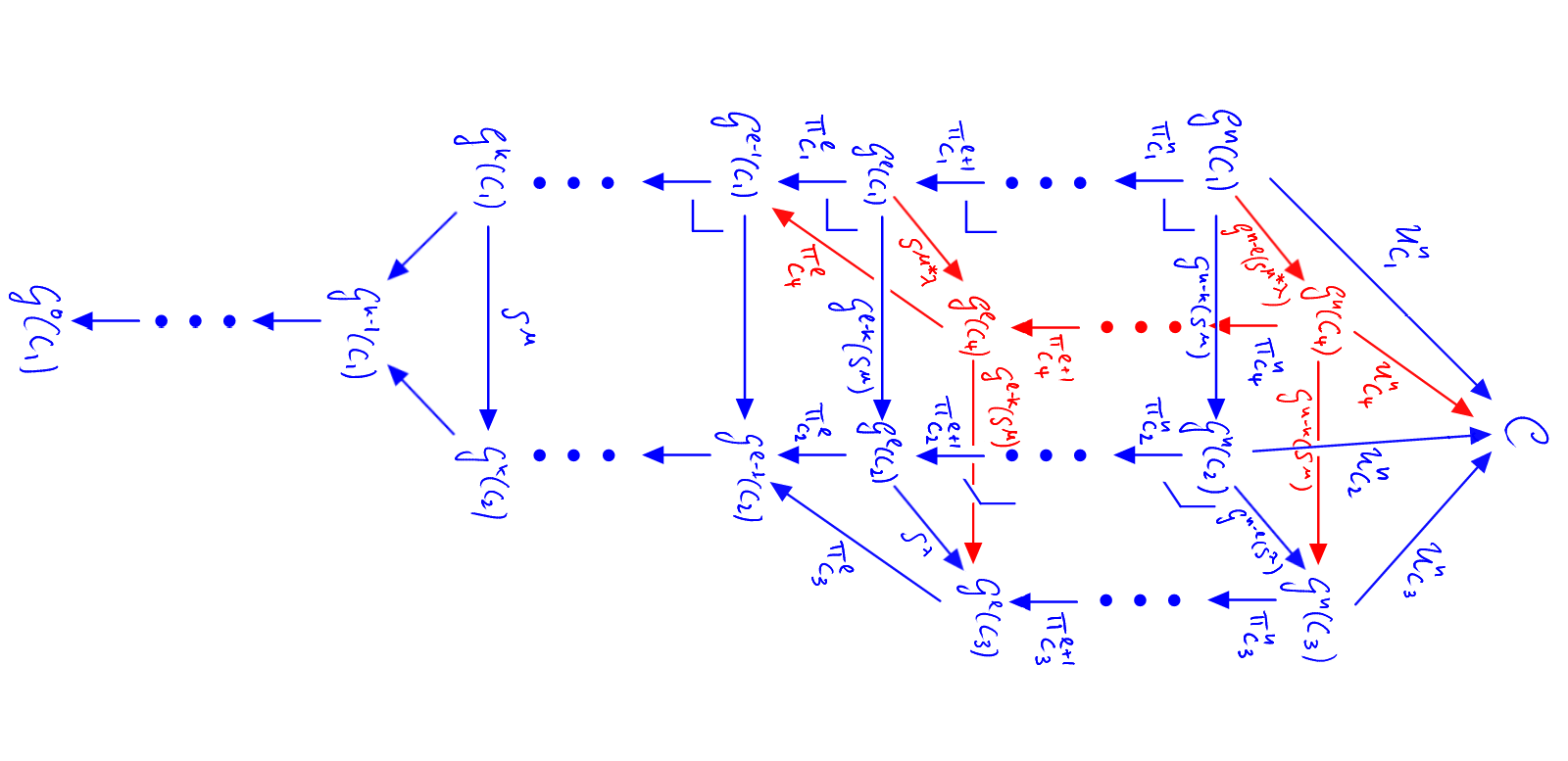}
\endgroup\end{restoretext}
The statement of the claim is equivalent to the existence of the \cred{} part of the diagram (such that the above diagram commutes). Using \autoref{lem:unpacking_collapse} the above diagram is equivalent to the following diagram
\begin{restoretext}
\begingroup\sbox0{\includegraphics{test/page1.png}}\includegraphics[angle=90, origin=c,clip,trim= {.0\ht0} {.0\ht0} {.4\ht0} {.0\ht0} ,width=9cm]{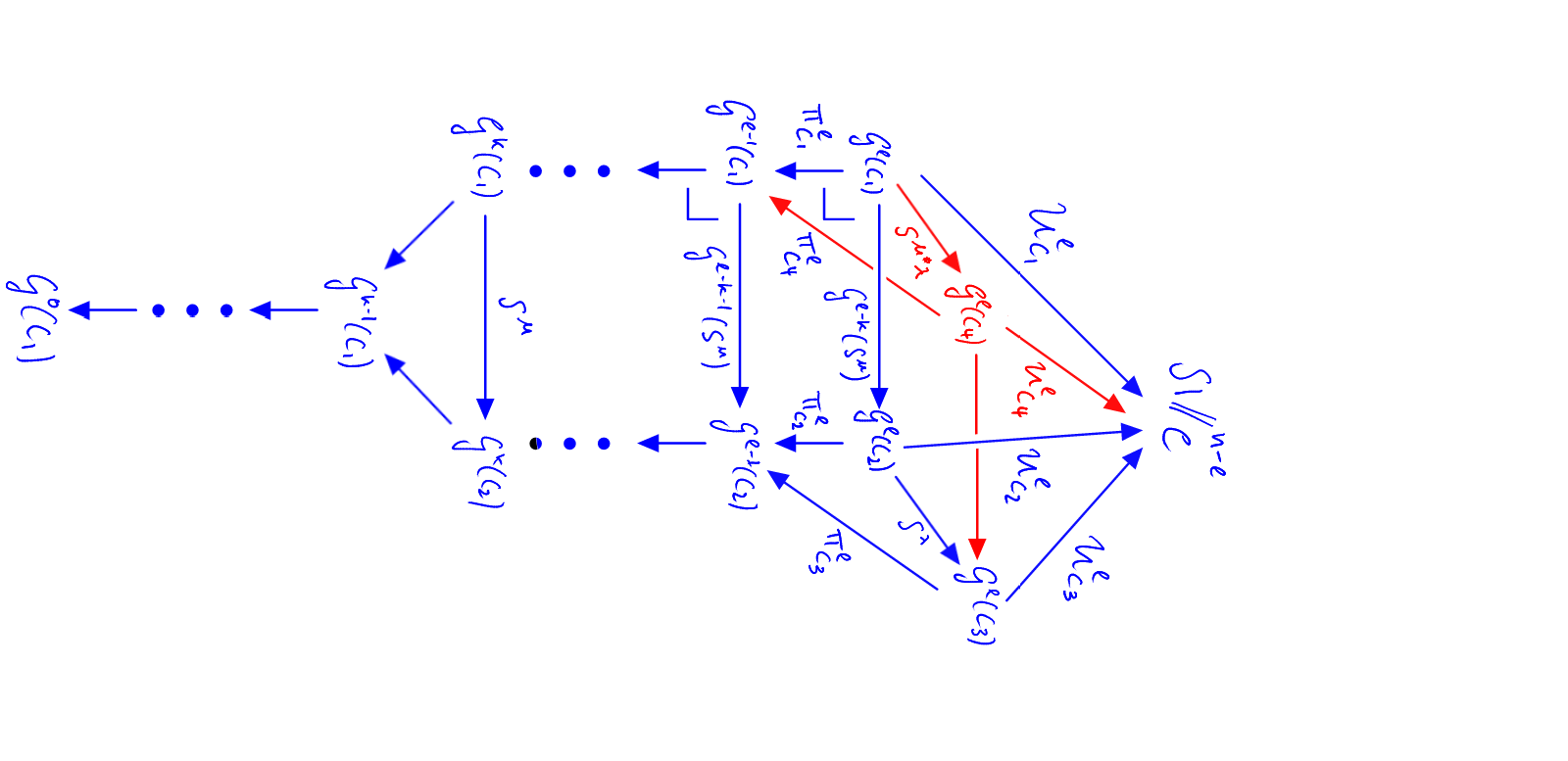}
\endgroup\end{restoretext}
But now we can see that the claim's statement is just an instance of \autoref{claim:collapse_dimension_interaction1}. More formally, we have the following: $\mu : \scC_1 \kcoll k \scC_2$ implies $\tsU k_{\scC_1} = \tsU k_{\scC_2} \sS^\mu$. Then note $\tsU {l-1}_{\scC_1} = \tsU {l-1}_{\scC_2} \upi {\sS^\mu} {l-1-k}$ by \autoref{constr:unpacking_collapse}. Given $\lambda : \scC_2 \kcoll l \scC_3$ we can now apply the statement of \autoref{claim:collapse_dimension_interaction1} with the following notational identification: $\scA := \tsU {l-1}_{\scC_1}$, $\scB := \tsU {l-1}_{\scC_2}$, ${\widetilde{\scB}} := \tsU {l-1}_{\scC_3}$, $\lambda := \lambda$ and $H := \upi {\sS^\mu} {l-k-1}$.  \autoref{claim:collapse_dimension_interaction1} then allows us to obtain ${{\widetilde{\scB}H}} = \tsU {l-1}_{\scC_3} H$ together with $\lambda H: \und{\widetilde{\scB}} H \into \tusU l_{\scC_1}$ and $\tsU l_{\scC_1} = \sU_{{{\widetilde{\scB} H}}} \sS^{\lambda H}$. 
Using \autoref{defn:complete_unpacking_and_repacking}, we define 
\begin{equation}
\scC_4 := \tsR {l-1}_{\sT^{l-1}_{\scC_1},{\widetilde{\scB} H}}
\end{equation}
Note that since ${\widetilde{\scB} H} : \tsG {l-1}(\scC_1) \to \SIvert {n-(l-1)} \cC$, we find $\scC_4$ to have the type of a $\cC$-labelled singular $n$-cube family, that is, $\scC_4 : X \to \SIvert n \cC$, which by its construction satisfies $\tsU {l-1}_{\scC_4} = {{\widetilde{\scB} H}}$ and $\sT^{l-1}_{\scC_4} = \sT^{l-1}_{\scC_1}$. Further, we have $\tsU {l}_{\scC_4} = \sU_{{{\widetilde{\scB} H}}}$ and setting
\begin{equation}
\mu\pbstar  \lambda := \lambda H
\end{equation}
this  $\tsU {l}_{\scC_1} = \tsU {l}_{\scC_4} \sS^{\lambda H}$. Together this means that $\mu\pbstar \lambda : \scC_1 \kcoll l \scC_4$.

Now, we claim $\mu : \scC_4 \kcoll k \scC_3$
which can be seen as follows: since $\mu : \scC_1 \kcoll k \scC_2$ we have $\tsU {k}_{\scC_1} = \tsU {k}_{\scC_2} \sS^\mu$.  Applying \autoref{rmk:unpacking_collapse}(i) we find for $i$ satisfying $k < k + i \leq n$
\begin{equation} 
\xymatrix@C=2cm{ \tsG {k + i}(\scC_1) \ar[r]^{\upi {\sS^\mu} i} \ar[d]_{\tpi {k+i}_{\scC_1}} \pullback & \tsG {k + i}(\scC_2)\ar[d]^{\tpi {k+i}_{\scC_2}} \\
\tsG {k + i -1}(\scC_1) \ar[r]^{\upi {\sS^\mu} {i-1}}  & \tsG {k + i - 1}(\scC_2) }
\end{equation}
Since $\sT^{l-1}_{\scC_2} = \sT^{l-1}_{\scC_3}$ and $\sT^{l-1}_{\scC_1} =\sT^{l-1}_{\scC_4}$ this implies for $i$ satisfying $k < k + i < l$
\begin{equation} \label{eq:packing_collapses}
\xymatrix@C=2cm{ \tsG {k + i}(\scC_4) \ar[r]^{\upi {\sS^\mu} i} \ar[d]_{\tpi {k+i}_{\scC_4}} \pullback & \tsG {k + i}(\scC_3)\ar[d]^{\tpi {k+i}_{\scC_3}} \\
\tsG {k + i -1}(\scC_4) \ar[r]^{\upi {\sS^\mu} {i-1}}  & \tsG {k + i - 1}(\scC_3) }
\end{equation}
However, $\tsU {l-1}_{\scC_4} = {\widetilde{\scB} H}$ and ${{\widetilde{\scB}H}}$ denotes $\tsU {l-1}_{\scC_3} \upi {\sS^\mu} {l-k-1}$. This means $\tsU {l-1}_{\scC_4} = \tsU {l-1}_{\scC_3} \upi {\sS^\mu} {l-k-1}$ which together with the pullbacks \eqref{eq:packing_collapses} then allows us to apply the implication $(ii) \imp (i)$ of \autoref{lem:unpacking_collapse} to obtain $\tsU {k}_{\scC_4} = \tsU {k}_{\scC_3} \upi {\sS^\mu} 0$. Since also $\sT^{k-1}_{\scC_4} = \sT^{k-1}_{\scC_3}$ (which follows from $\sT^{l-1}_{\scC_3} = \sT^{l-1}_{\scC_2}$, $\sT^{k-1}_{\scC_2} = \sT^{k-1}_{\scC_1}$, $\sT^{l-1}_{\scC_1} =\sT^{l-1}_{\scC_4}$ and $l > k$) the conditions for $\mu : \scC_4 \kcoll k \scC_3$ have been verified.
\qed
\end{lem}

\begin{lem}[$k$-level push-forward of $l$-level collapse] \label{lem:coll_push}
If
\begin{equation}
\xymatrix{ C_3 \\
C_1 \ar@{~>}[u]^{\lambda}_l \ar@{~>}[r]^{\mu}_k & C_2}
\end{equation}
then we have
\begin{equation} \label{eq:coll_push_full}
\xymatrix{ C'_3 \ar@{~>}[r]^{\mu}_k & C_4 \\
 C_3 \ar@{~>}[u]^\rho_l & \\
C_1 \ar@{~>}@/^2pc/[uu]^{\mu\pbstar \mu\psstar  \lambda}_l \ar@{~>}[u]^{\lambda}_l \ar@{~>}[r]^{\mu}_k & C_2 \ar@{~>}[uu]^{\mu\psstar  \lambda}_l}
\end{equation}
where $\mu\psstar \lambda$ is called the \emph{$k$-level pushforward} of $\lambda$ along $\mu$, and $\rho$ is the factorisation of $\mu\pbstar \mu\psstar \lambda$ through $\lambda$ (cf. \autoref{claim:injection_factorisation}).
\proof The assumptions in the claim are equivalent to the existence of the \cblue{} part of the following diagram
\begin{restoretext}
\begingroup\sbox0{\includegraphics{test/page1.png}}\includegraphics[angle=90, origin=c,clip,trim=0 {.0\ht0} 0 {.0\ht0} ,width=9cm]{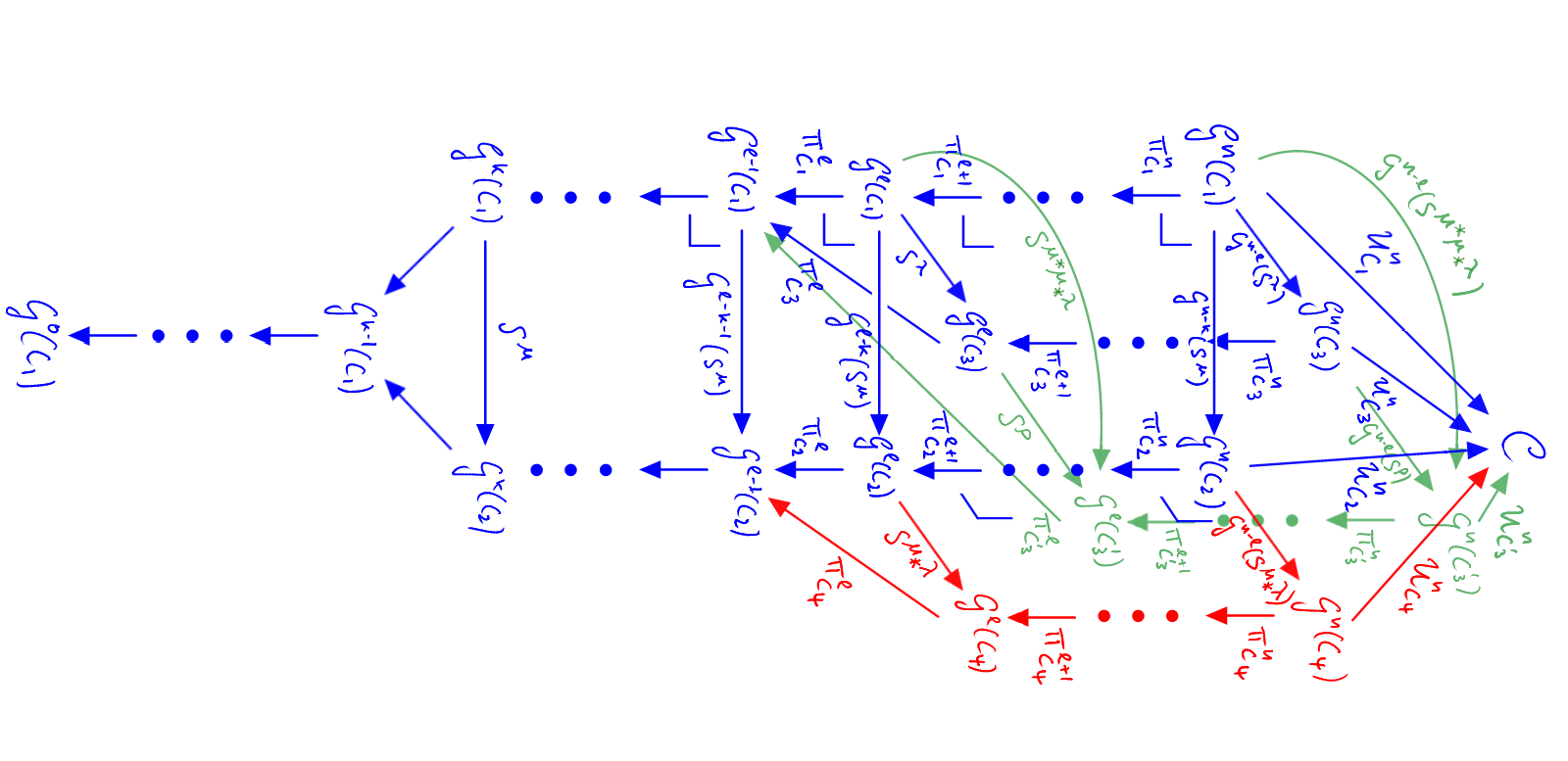}
\endgroup\end{restoretext}
The statement of the claim is equivalent to the existence of the \cred{} and \cdarkgreen{} part of the diagram. Using \autoref{lem:unpacking_collapse} the above diagram is equivalent to the following diagram
\begin{restoretext}
\begingroup\sbox0{\includegraphics{test/page1.png}}\includegraphics[angle=90, origin=c,clip,trim=0 {.0\ht0} {.4\ht0} {.0\ht0} ,width=9cm]{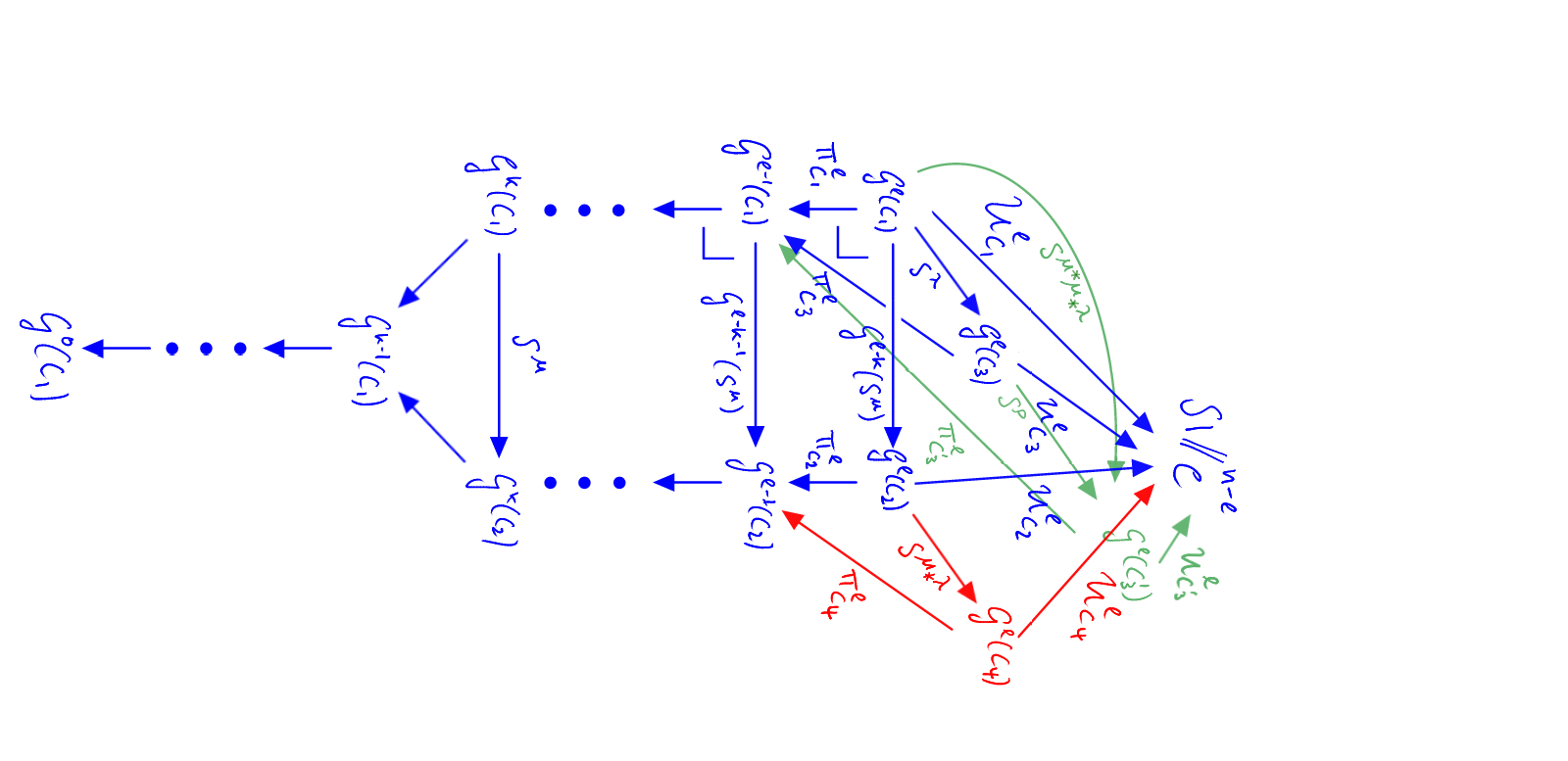}
\endgroup\end{restoretext}
But now we can see that the claim's statement is just an instance of \autoref{claim:collapse_dimension_interaction2}. More formally, we have the following: note $\sT^{l-1}_{\scC_1} = \sT^{l-1}_{\scC_3}$ and $\tsU {l}_{\scC_1} = \tsU {l}_{\scC_2} \sS^\lambda$. We can apply \autoref{claim:collapse_dimension_interaction2} with the following identification of notation: $ \scA := \tsU {l-1}_{\scC_1}$, $\widetilde{\scA} := \tsU {l-1}_{\scC_3}$, $\scB := \tsU {l-1}_{\scC_2}$, $\lambda := \lambda$ and $H := \upi {\sS^{\mu}} {l-k-1}$ (cf. \autoref{constr:unpacking_collapse}). From \autoref{claim:collapse_dimension_interaction2} we then obtain $H\psstar \widetilde{\scA}$ with $(\upi {\sS^{\mu}} {l-k-1})\psstar \lambda : \und{H\psstar \widetilde{\scA}} \into \und{\scB}$ and such that $\tsU {l}_{\scC_2} = \sU_{H\psstar \widetilde{\scA}} \sS^{(\upi {\sS^{\mu}} {l-k-1})\psstar \lambda}$ (Note, that by construction in \autoref{claim:collapse_dimension_interaction2} again $(\upi {\sS^{\mu}} {l-k-1})\psstar \lambda$ is non-identity if $\lambda$ is). Using \autoref{constr:unpacking_collapse}, we can define
\begin{equation}
\scC_4 := \tsR {l-1}_{\sT^{l-1}_{\scC_1}, {H\psstar \widetilde{\scA}}}
\end{equation}
Note that since ${H\psstar \widetilde{\scA}} : \tsG {l-1}(\scC_2) \to \SIvert {n-(l-1)} \cC$, we find $\scC_4$ to have the type of a $\cC$-labelled singular $n$-cube family, that is, $\scC_4 : X \to \SIvert n \cC$, which by its definition satisfies $\tsU {l-1}_{\scC_2} = {{H\psstar \widetilde{\scA}}}$ and $\tpi i_{\scC_4} = \tpi i_{\scC_2}$ for $i < l$. Thus we have $\sT^{l-1}_{\scC_2} = \sT^{l-1}_{\scC_4}$ and $\tsU {l}_{\scC_4} = \sU_{H\psstar \widetilde{\scA}}$ implying $\tsU {l}_{\scC_2} = \tsU {l}_{\scC_4} \sS^{H\psstar \lambda}$. This implies $H\psstar \lambda : \scC_2 \kcoll l \scC_4$ and thus we can set
\begin{equation}
\mu\pbstar \lambda := H\psstar \lambda
\end{equation}
The existence of $\rho$ and the outer square in \eqref{eq:coll_push_full} then follows from both \autoref{rmk:factorising_of_pull_push} and \autoref{lem:coll_pull}.
\qed
\end{lem}

\begin{cor}[$l$-level normalisation independent of $k$-level collapse] \label{cor:coll_norm_coll_ind}
If $\mu : \scC_1 \kcoll k \scC_2$ then $\scC_1$ is in $l$-normal form if and only if $\scC_2$ is.
\proof First, assume $\scC_2$ is not in $l$-level normal form, i.e. there is a non-identity $\lambda : \scC_2 \kcoll l \scC_3$. By \autoref{lem:coll_pull} we can then construct $\mu\pbstar  \lambda : \scC_1 \kcoll l \scC_4$ (which by the construction in \autoref{claim:collapse_dimension_interaction1} can be seen to be a non-identity collapse if $\lambda$ is not the identity). Thus $\scC_1$ cannot be in $l$-level normal form.

Conversely, assume that $\scC_1$ is not in $l$-level normal form, i.e. there is a non-identity $\lambda : \scC_1 \kcoll l \scC_3$. By \autoref{lem:coll_push} we can then construct $\mu\psstar  \lambda : \scC_2 \kcoll l \scC_4$ (which by \autoref{claim:collapse_dimension_interaction2} can be seen to be a non-identity collapse if $\lambda$ is not the identity). Thus $\scC_2$ cannot be in $l$-level normal form.
\qed
\end{cor}

\begin{defn}[Ordered collapse sequence] \label{defn:ordered_coll_seq} An ordered ($n$-level) collapse sequence $\vvec \lambda : \scC \starcoll \scC_1$ (notationally indicated by a vector sign) is a sequence of collapse of the form  
\begin{equation}
 \xymatrix{\scC \ar@{~>}[r]^-{\vvec\lambda ^n}_-n & \scC_{n} \ar@{~>}[r]^-{\vvec\lambda ^{n-1}}_-{n-1} & \scC_{n-1} \dots \ar@{~>}[r]^-{\vvec\lambda ^1}_-1 & \scC_1}
\end{equation}
\end{defn}

\begin{lem}[Ordered collapse sequences for normal forms] \label{lem:coll_seq}
If $\scC_1$ is a normal form of $\scC$ then there is an ordered collapse sequence
\begin{equation}
\scC \kcoll n \scC_n \kcoll {n-1} ... \kcoll 1 \scC_1
\end{equation}
ending in $\scC_1$. Further, for any such sequence $\scC_l$ is in $l$-level collapse normal form.
\proof Assume $\scC_1$ is a normal form of $\scC$ obtained by a sequence (cf. \autoref{defn:collapse_normal_form})
\begin{equation}
\scC \kcoll {k_0} \scA_0 \kcoll {k_1} \scA_1 \kcoll {k_2} ... \kcoll {k_m} \scA_m = \scC_1
\end{equation}
From this, using the commutative squares from \autoref{lem:coll_pull}, and the compositionality of collapses from \autoref{lem:coll_comp}, we can obtain a sequence
\begin{equation}
\scC \kcoll {n} \scC_n \kcoll {n-1} \scC_{n-1} \kcoll {n-2} ... \kcoll {1}  \scC_1
\end{equation}
Note that some of the collapses in this sequence might be trivial.

Now, for any such sequence the statement of \autoref{cor:coll_norm_coll_ind} implies that $\scC_l$ is in $l$-level normal form if and only if $\scC_{l-1}$ is if and only if $\scC_{l-2}$ is in $l$-level normal form, and inductively, if and only if $\scC_1$ is $l$-level normal form. Since $\scC_1$ is in normal form, and thus in $l$-level normal form, we thus infer that $\scC_l$ is in $l$-normal form as claimed.
\qed
\end{lem}

\begin{defn}[Normal form collapse] If $\lambda : C_0 \kcoll l C_1$ and $C_1$ is in $l$-level collapse normal form, then $\lambda : C_0 \kcoll l C_1$ is called a $l$-level \textit{normal form collapse} of $C_0$.
\end{defn}

\begin{lem}[Uniqueness of $l$-level normal form collapse] \label{lem:coll_uniq} If $\lambda_1 : \scC_0 \kcoll l \scC_1$, $\lambda_2 : \scC_0 \kcoll l \scC_2$ such that both $\scC_1$, $\scC_2$ are in $l$-level collapse normal form, then $\lambda_1 = \lambda_2$ and $\scC_1 = \scC_2$. 
\proof $\lambda_i : \scC_0 \kcoll l \scC_i$ by definition is equivalent to $\lambda_i : \tsU {l-1}_{\scC_0} \mcoll \tsU {l-1}_{\scC_i}$. Then \autoref{thm:family_pushouts} applies and we find
\begin{equation}
\xymatrix{ \tsU {l-1}_{\scC_0} \ar@{~>}[r]^{ {\lambda_1}} \ar@{~>}[d]_{ {\lambda_2}} & \tsU {l-1}_{\scC_1} \ar@{~>}[d]^{ {\mu_1}} \\
\tsU {l-1}_{\scC_2} \ar@{~>}[r]_{ {\mu_2}} &  \widecheck{\scC} \pushoutfar}
\end{equation}
We can then define
\begin{equation}
\scC_3 := \tsR {l-1}_{\sT^{l-1}_{\scC_0}, \widecheck{\scC}}
\end{equation}
which implies $\mu_i : \scC_i \kcoll l \scC_3$. Since $\scC_i$ are assumed to be in normal form we must have $\mu_i = \id$. We deduce $\lambda_1 = \lambda_2$ and consequently $\scC_1 = \scC_2$ as claimed. \qed
\end{lem}

The section culminates in the following result.

\begin{thm}[Normal forms are unique] \label{thm:normal_forms_unique} Given $\scC : X \to \SIvertone \cC$, there is a unique normal form $\NF{\scC}^n$ and a unique ordered collapse sequence
\begin{equation}
 \xymatrix{\scC = \NF{\scC}^n_{n+1} \ar@{~>}[r]^-{\nfc_\scC^n}_-n & \NF{\scC}^n_{n} \ar@{~>}[r]^-{\nfc_\scC^{n-1}}_-{n-1} & \NF{\scC}^n_{n-1} \dots \ar@{~>}[r]^-{\nfc_\scC^1}_-1 & \NF{\scC}^n_1 = \NF{\scC}^n}
\end{equation}
which we denote by $\vec \nfc_\scC$.
\proof In the previous section we showed that there is at least one normal form $\scA_1$ of $\scC$. Thus \autoref{lem:coll_seq} implies there is at least one ordered collapse sequence
\begin{equation}
\xymatrix{\scC \ar@{~>}[r]^-{\nfc^n}_-n & \scA_n \ar@{~>}[r]^-{\nfc^{n-1}}_-{n-1} & \scA_{n-1} \dots \ar@{~>}[r]^-{\nfc^1}_-1 & \scA_1}
\end{equation}
Assume a different normal form $\scB_1$ and an ordered collapse sequence
\begin{equation}
\xymatrix{\scC \ar@{~>}[r]^-{\mu^n}_-n & \scB_n \ar@{~>}[r]^-{\mu^{n-1}}_-{n-1} & \scB_{n-1} \dots \ar@{~>}[r]^-{\mu^1}_-1 & \scB_1}
\end{equation}
Using \autoref{lem:coll_uniq} inductively for $l = n, (n-1), ... , 1$, we find $\nfc^l = \mu^l$ and $\scA_l = \scB_l$ thus proving the theorem. 
\qed
\end{thm}

\begin{eg}[Ordered collapse sequence to normal form] In \autoref{eg:normal_forms}, we already constructed an example of an ordered collapse sequence to normal form 
\begin{equation}
\Big(\xymatrix{\scA = \NF{\scA}^2_{3} \ar@{~>}[r]^-{\lambda_\scA^2}_-2 & \NF{\scA}^2_{2} \ar@{~>}[r]^-{\lambda_\scA^{1}}_-{1} & \NF{\scA}^2_{1}}\Big)  \quad =\quad \Big( \xymatrix{\scA \ar@{~>}[r]^-{\lambda}_-2 & \scB \ar@{~>}[r]^-{\mu}_-{1} & \scC}\Big)
\end{equation}
with the following data
\begin{restoretext}
\begingroup\sbox0{\includegraphics{test/page1.png}}\includegraphics[clip,trim={.15\ht0} {.0\ht0} {.15\ht0} {.0\ht0} ,width=\textwidth]{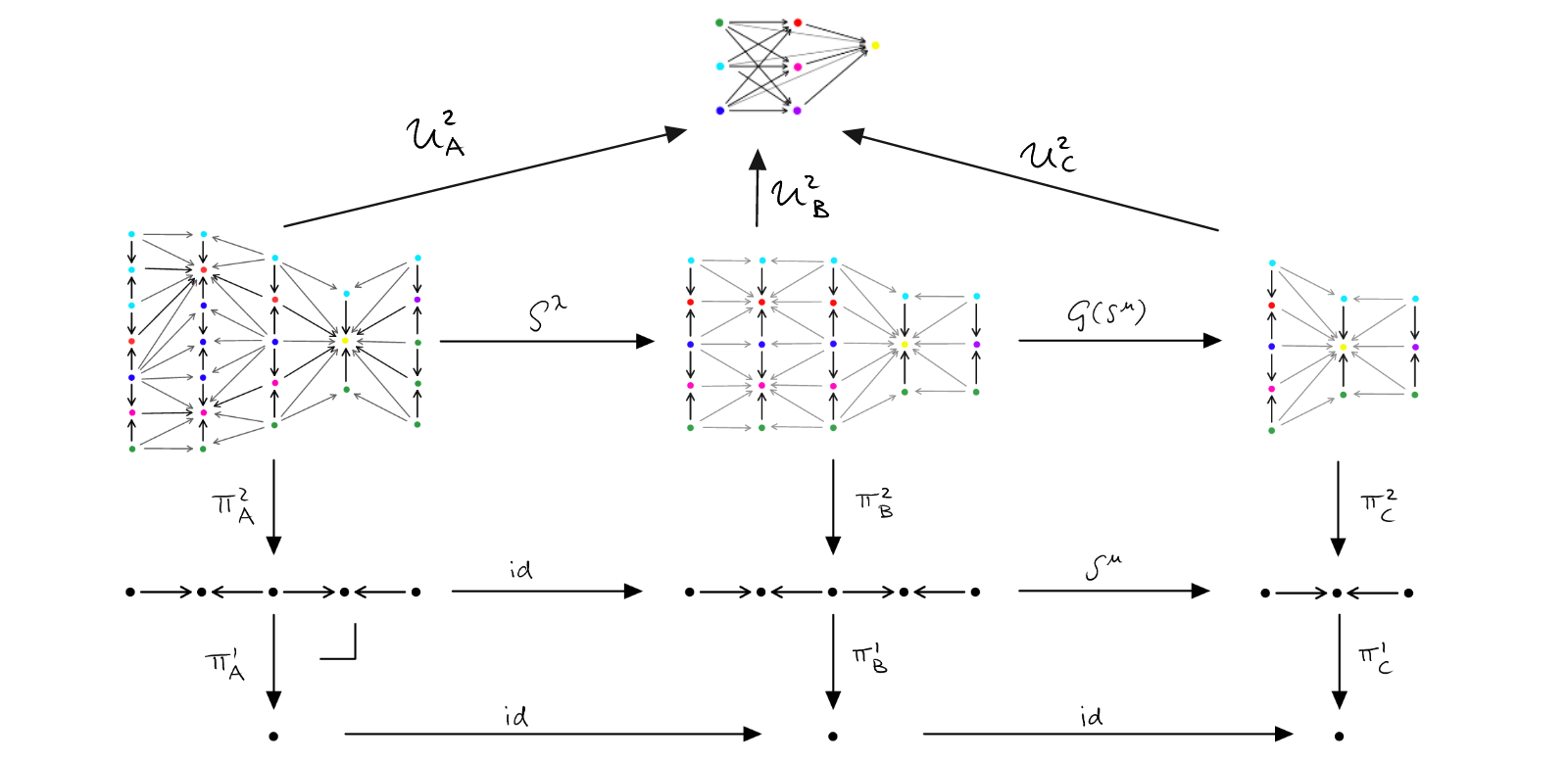}
\endgroup\end{restoretext}
\end{eg}

\section{Multi-level collapse} \label{sec:coll_multi}

\subsection{Definition}

In this section, we connect our discussion of collapse with the ideas presented in the \hyperref[ch:summary]{Summary} (\autoref{sec:sum_norm}). We will also understand that the results of the previous sections can be phrased more abstractly as the construction of terminal objects in certain subcategories of $\Bunbc^n_\cC$ (cf. \autoref{defn:multilevel_base_change}).

We start with the following definition:

\begin{defn}[Multi-level collapse] \label{defn:multilevel_collapse} A \textit{(multi-level family) collapse} $\vvec S : \scA \mcoll \scB$  of $X$-indexed $\SIvert n \cC$-families $\scA$, $\scB$ is a ``fibrewise open and surjective" multi-level base change $\scA \to \scB$ (cf. \autoref{defn:multilevel_base_change}). More explicitly, it consists of functors
\begin{equation}
\vvec S^k : \tsG k(\scA) \to \tsG k(\scB)
\end{equation}
such that 
\begin{enumerate}
\item Firstly, 
\begin{equation}
\xymatrix{ & \cC & \\
\tsG n(\scA) \ar[ur]^{\tsU n_\scA} \ar[rr]_{\vvec S^n} && \tsG n(\scB) \ar[ul]_{\tsU n_\scB} }
\end{equation}
and for $0 \leq k < n$
\begin{equation}
\xymatrix{ \tsG {k+1}(\scA) \ar[r]^{\vvec S^{k+1}} \ar[d]_{\tpi {k+1}_\scA} & \tsG {k+1}(\scB) \ar[d]^{\tpi {k+1}_\scB} \\
\tsG k(\scA) \ar[r]^{\vvec S^k} & \tsG k(\scB)  }
\end{equation}
\item Secondly, we require for each $\vvec S^k(x) = y$ that the restriction of $\vvec  S^{k+1}$ to the fiber 
\begin{equation}
\rest {\vvec S^{k+1}} {x} : \tusU k_\scA (x) \to \tusU k_\scB(y)
\end{equation}
is a collapse functor $\tusU k_\scA(x) \to \tusU k_\scB(y)$.
\item And finally, we also require $\vvec S^0 = \id_X$.
\end{enumerate}
Multi-level family collapse functors are stable under composition, that is for multi-level family collapse functors $\vvec S : \scB \to \scC$, $\vvec T : \scA \to \scB$, $(\vvec S \vvec T)^k = \vvec S^k \vvec T^k$ defines a multi-level collapse   $\vvec S \vvec T$.
\end{defn}

We remark that the last condition of the definition (required $\vvec S^0 = \id_X$)  was not present in \autoref{sec:sum_norm}. As we will be ultimately working with cube families over a point ($X = \bnum 1$) this difference will not play a role. However, for general $X$, the condition is needed for the next result, which establishes a 1-to-1 correspondence of multi-level collapse an ordered collapse sequences. 

\subsection{Decomposing multi-level collapses into $k$-level collapse}

\begin{constr}[1-to-1 correspondence of multi-level collapse functors and ordered collapse sequences] \label{constr:multilevel_collapse} Let $\vvec \lambda$ be an ordered collapse sequence
\begin{equation}
 \xymatrix{\scC_{n+1} \ar@{~>}[r]^-{\vvec\lambda ^n}_-n & \scC_{n} \ar@{~>}[r]^-{\vvec\lambda ^{n-1}}_-{n-1} & \scC_{n-1} \dots \ar@{~>}[r]^-{\vvec\lambda ^1}_-1 & \scC_1}
\end{equation}
We define a multi-level family collapse functor $\vsS{\vvec \lambda} : \scC_{n+1} \to \scC_1$ as follows. Given $\mu :  \scA \kcoll l \scB$ we define a multi-level family collapse functor $\vsS\mu : \scA \to \scB$ by setting (cf. \autoref{notn:k_lvl_basechange})
\begin{equation}
(\vsS{\mu})^k = \tsG {k-l}(\sS^\mu)
\end{equation}
We then set
\begin{equation}
\vsS{\vvec \lambda} := \vsS{\vvec\lambda^1} \vsS{\vvec\lambda^2} \dots \vsS{\vvec\lambda^n}
\end{equation}

Conversely, assume a multi-level collapse $\vvec S : \scC_{n+1} \to \scC_1$. We define an ordered collapse sequence $\vvec \lambda$ such that $\vvec S = \vsS{\vvec \lambda}$ as follows. We inductively define multi-level collapse functors $\vvec S_k : \scC_{k+1} \to \scC_k$ and $\vvec S_{\und k} : \scC_{n+1} \to \scC_{k+1}$ such that
\begin{equation}
\vvec S = \vvec S_1 \vvec S_2 ... \vvec S_k \vvec S_{\und k}
\end{equation}
such that $\vvec S_{\und k}^l = \id$ for $l \leq k$.

For $k = 0$, set $\vvec S_0 = \id$ and $\vvec S_{\und 0} = \vvec S$. 

For $k > 0$, define $\vvec S_k : \scC_{k+1} \to \scC_k$ to be the $k$-level base change (cf. \autoref{notn:k_lvl_basechange} and recall \autoref{ssec:multi_bc})
\begin{equation}
\vvec S^l_k = \tsG {k-l}(\vvec S^k_{\und {k-1}}) 
\end{equation}
(note that this defines $\scC_{k+1}$ from $\scC_k$ by a $k$-level base change). Using the universal property of pullbacks, we find a unique factorisation $\vvec S_{\und k}$ such that
\begin{equation}
\vvec S_{\und {k-1}} = \vvec S_{k} \vvec S_{\und k}
\end{equation}
Note that $\vvec S_{\und k}$ is again a multi-level collapse functor as claimed inductively. This completes the inductive decomposition. Note that using \autoref{thm:collapse_maps_vs_injections}, each  $\vvec S_k$ can by definition be written as $\vsS{\vvec\lambda^k}$ for $\vvec \lambda^k : \scC_{k+1} \kcoll k \scC_{k}$ (such that $\vvec S^k_{\und {k-1}} = \sS^{\vvec \lambda^k}$). We conclude $\vvec S = \vsS{\vvec \lambda}$ as claimed.

Also note that the composition of $\vvec S$ into functors of the form of a $k$-level base change is unique (which can be seen inductively from the above argument). Thus the above constructions are inverse to each other.
\end{constr}

\begin{rmk}[Uniqueness of normal form collapse for multi-level functors] As a corollary of \autoref{thm:normal_forms_unique} and the preceding construction we obtain the following. Let 
\begin{equation}
\vsS{\vvec \mu} : \scA \to \scB
\end{equation}
be a multi-level family collapse functor of families $\scA, \scB : X \to \SIvert n \cC$ associated to an ordered collapse sequence $\vvec \mu : \scA \starcoll \scB$. Then we have
\begin{equation} \label{eq:uniqueness_for_multi-level_coll}
\sS^{\vec \nfc_\scA} = \vsS{\vec \nfc\scB} \vsS{\vvec \mu}
\end{equation}
\end{rmk}

We end with the remark that using \autoref{defn:multilevel_base_change}, the central results of this chapter can be summarised by saying that any connected component of the subcategory generated by multi-level collapses in $\Bunbc^n_\cC$ has a terminal object.

\chapter{Embedding of cubes} \label{ch:emb}

In this chapter we discuss how to embed a cube family into another cube family. The definition of embedding is ``dual" to the definition of collapse, in the sense that the former is based on injective open maps of fibers, while the latter is based on surjective open maps. 

In \autoref{sec:emb_def} we will give the definition of embeddings of cube families. We will also observe some basic properties of this definition, and give examples that will be used throughout the rest of the chapter. In \autoref{sec:emb_end}, we will decompose embeddings into their \textit{level-wise operations} which are determined by ``endpoint sections" and called \textit{family embedding functors} (these are analogous to family collapse functors: the level-wise operations of multi-level collapse). We warn the reader, that this section will introduce quite tedious notation to classify embeddings, and recommend to focus on examples when reading the chapter for the first time. In \autoref{sec:emb_min} we will then use family embedding functors determined by ``minimal" endpoint sections to construct a minimal embedding around a given point of the total space. Finally, in \autoref{sec:emb_norm} we discuss the natural interplay of embedding and collapse, showing (again by a rather tedious, but straight-forward argument) that a collapse of a parent family induces a collapse on families embedded in it.

\section{Embeddings of singular cube families} \label{sec:emb_def}

\subsection{Definition of embeddings}

We start by defining embedding functors. These will be the fibrewise mappings of the components of embeddings.

\begin{defn}[Embedding functors] \label{defn:embedding_fctr_of_int} An open functor $f : I \to J$ of singular intervals $I, J$ (cf. \autoref{defn:open_maps}) is called an \textit{(singular) embedding functor} if it is injective on objects.
\end{defn}

\begin{claim}[Embedding functors are linear] \label{claim:linearity_subfamilies} Let $f : I \to J$ be an embedding functor. Then $f$ maps
\begin{equation}
d \mapsto d + f(0)
\end{equation}

\proof Note that since $f$ preserves regular segments, we have $f (0) \in \regcont(J)$. The proof of the statement is by induction (it holds for $d = 0$), in each step showing $f(d+1) = f (d) + 1$. Assume $d$ is a regular segment, $d+1$ is a singular height and thus $d \to d+1$ in $I)$ (argue similarly if $d$ is singular). Then, since $f$ is monotone and injective we must have $f (d+1) > f (d)$. But since it is a functor of posets, we also have $f (d) \to f (d+1)$. By \autoref{defn:singular_intervals} this implies $f (d+1) = f (d) + 1$ as required.  \qed
\end{claim}

Recall that open surjective functors (that is, collapse functors) are classified by their underlying monomorphisms which record the mappings between singular heights. We will see that open injective functors (that is, embedding functors) are classified by their mapping of (regular) interval \textit{endpoints}. This is foreshadowed by the following remark.

\begin{constr}[Endpoints determine embedding functors] \label{rmk:endpoints_of_embedding_fctr}
Given $J \in \SI$ and $q_-, q_+ \in \regcont(J)$, $q_- < q_+$ then the functor $\sJ\restsec{[q_-,q_+]} : \singint {\frac{q_+-q_-}{2}} \to J$ defined by
\begin{equation}
d \mapsto d + q_-
\end{equation}
is an embedding functor. Conversely, using \autoref{claim:linearity_subfamilies} every embedding functor $f : I \to J$ equals $\sJ\restsec{[f(0), f(2\iH_I)]}$. Thus embedding functors $f : I \to J$ (with fixed $J$ and variable $I$) correspond to choices of \textit{endpoints} $q_-, q_+ \in \regcont(J)$, $q_- < q_+$. This statement will be revisited in a more general context later on (cf. \autoref{constr:endpoint_inclusions}).
\end{constr}

\begin{egs}[Embedding functors] \hfill
\begin{itemize}
\item Consider the function $f_0 : \singint 2 \to \singint 5$ defined by
\begin{restoretext}
\begingroup\sbox0{\includegraphics{test/page1.png}}\includegraphics[clip,trim=0 {.1\ht0} 0 {.05\ht0} ,width=\textwidth]{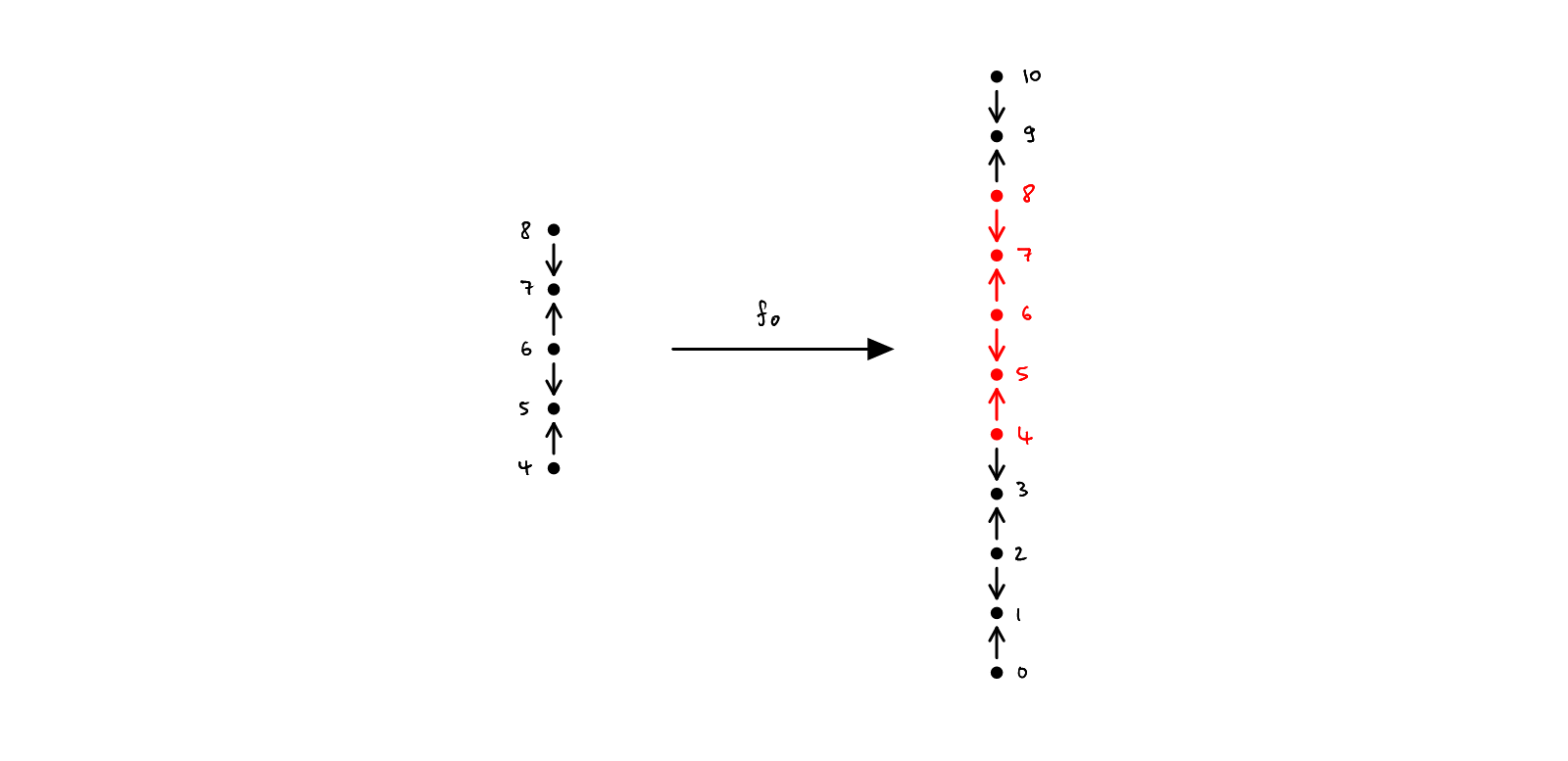}
\endgroup\end{restoretext}
$f_0$ is an injective functor of posets, monotone and preserves regular segments. Thus $f_0$ is an embedding functor. Note that the mapping of $f_0$ is indeed linear in the sense of \autoref{claim:linearity_subfamilies}. Further, in the above we coloured the image of $f_0$ in \cred{}. The image, and in fact the ``endpoints" of the image, determine $f_0$ fully in the sense of \autoref{rmk:endpoints_of_embedding_fctr}.
\item  Consider the function $f_1 : \singint 0 \to \singint 3$ defined by
\begin{restoretext}
\begingroup\sbox0{\includegraphics{test/page1.png}}\includegraphics[clip,trim=0 {.25\ht0} 0 {.2\ht0} ,width=\textwidth]{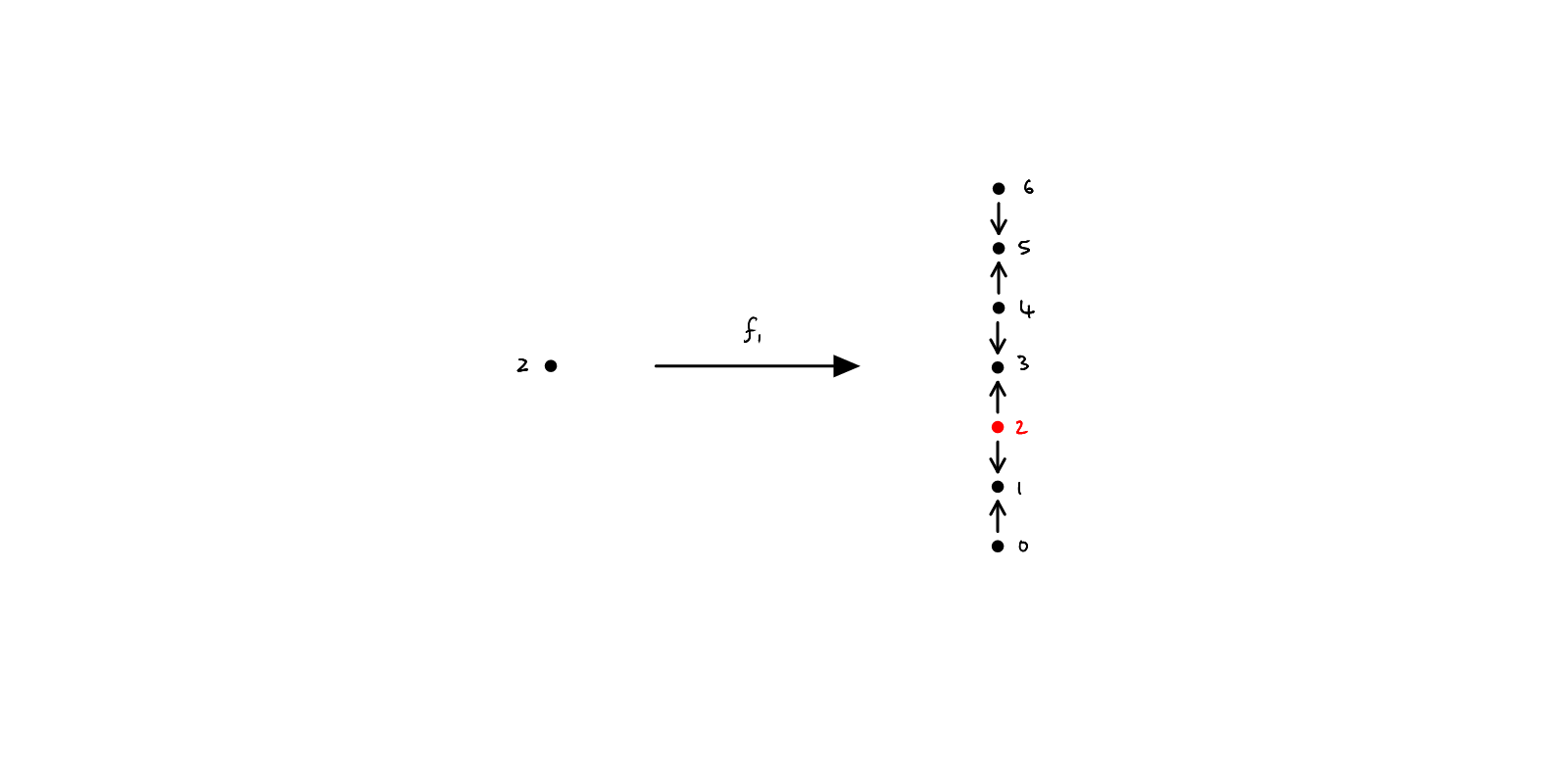}
\endgroup\end{restoretext}
$f_1$ is an embedding functor. Note that the endpoints of $f_1$ in the sense of  \autoref{rmk:endpoints_of_embedding_fctr} coincide.
\item Consider the function $f_2 : \singint 0 \to \singint 3$ defined by
\begin{restoretext}
\begingroup\sbox0{\includegraphics{test/page1.png}}\includegraphics[clip,trim=0 {.25\ht0} 0 {.2\ht0} ,width=\textwidth]{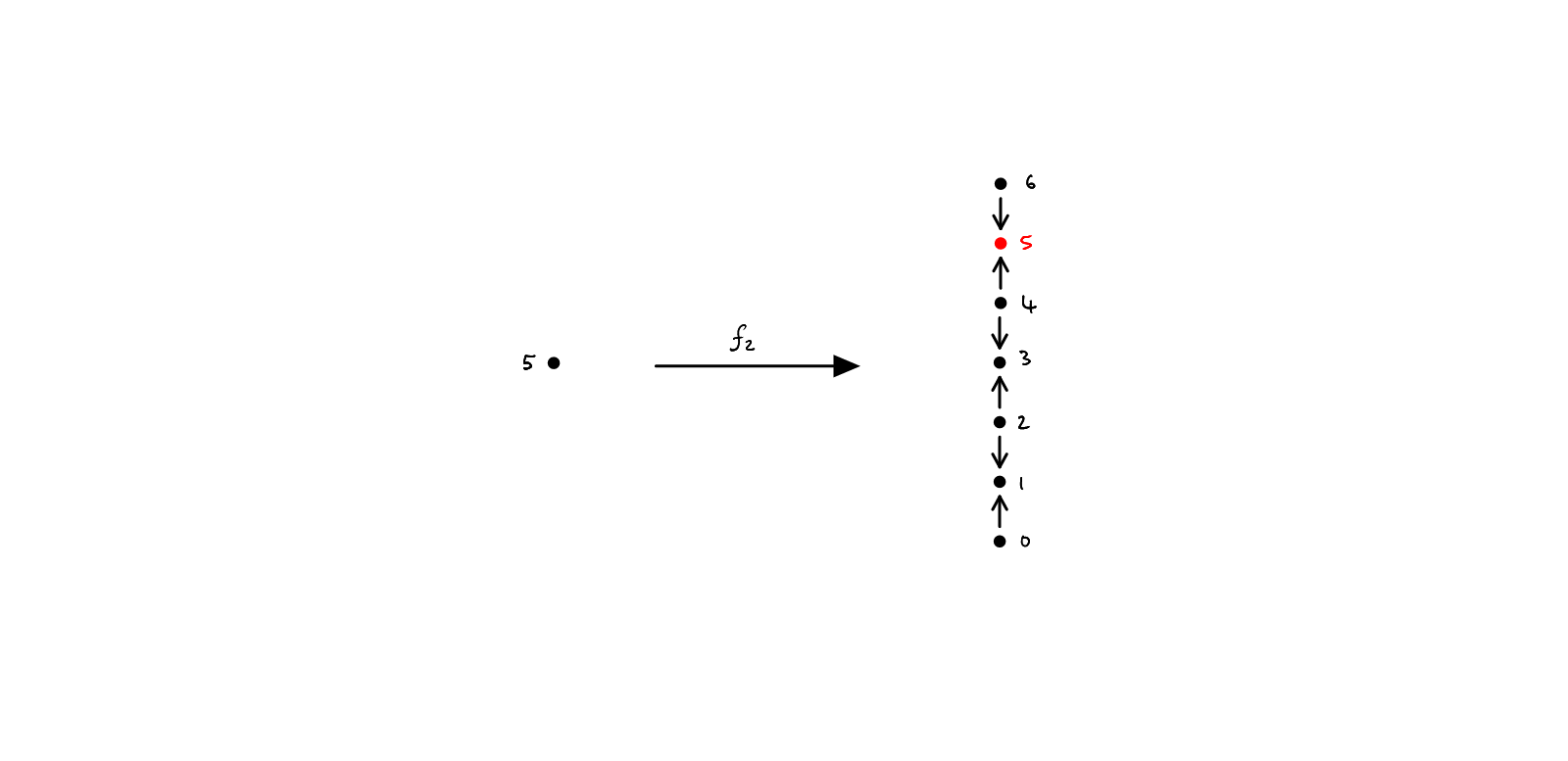}
\endgroup\end{restoretext}
$f_2$ is an injective functor of posets and monotone. However, $f_2$ does not preserve regular segments and is thus not open, as would be required for being an embedding functor.
\end{itemize}
\end{egs}

\begin{defn}[Multi-level embeddings] \label{defn:subfamilies}  Let $\scA : X \to \SIvert n \cC$ be a $\cC$-labelled singular $n$-cube family over $X$. An \textit{(multi-level) embedding} $\theta : \scB \mono \scA$ into $\scA$ is a ``fibrewise open and injective" multi-level base change $\scB \to \scA$ (cf. \autoref{defn:multilevel_base_change}). More explicitly, consists of a $\cC$-labelled singular $n$-cube family $\scB$ over $Y$ together with poset inclusions $\theta^k : \tsG k(\scB) \into \tsG k(\scA)$ which satisfy
\begin{enumerate}
\item The following commutes
\begin{equation} \label{eq:subbund_commute_1}
\xymatrix{ & \cC & \\
\tsG n(\scB) \ar[ur]^{\tsU n_\scB} \ar[rr]_{\theta^n} && \tsG n(\scA) \ar[ul]_{\tsU n_\scA} }
\end{equation}
\item For $0 < k \leq n$, the following commutes
\begin{equation} \label{eq:subbund_commute_2}
\xymatrix{ \tsG {k}(\scB) \ar[d]_{\tpi k_\scB} \ar[r]^{\theta^k} &\tsG k(\scA) \ar[d]^{\tpi k_\scA} \\
\tsG {k-1}(\scB) \ar[r]_{\theta^{k-1}} & \tsG {k-1}(\scA) }
\end{equation}
\item For $0 < k \leq n$ and $y \in \tsG {k-1}(\scB)$, $\rest {\theta^k} y$ is an embedding functor
\end{enumerate}
With reference to a specific embedding $\theta : \scB \mono \scA$ we also call $\scB$ a \textit{subfamily} (of $\scA$) and $\scA$ a \textit{parent family} (of $\scB$). Note that multi-level embeddings are stable under composition (as multi-level base changes).
\end{defn}

Note that for embeddings, we usually drop the predicate ``multi-level" as this will be our default assumption. Note also that we do not use the vector notation as we did for multi-level base change and multi-level collapse. 

We remark two immediate facts about embeddings based on its definition.

\begin{rmk}[Fibrewise endpoints determine embedding components] \label{rmk:fibrewise_endpoints_of_embedding_fctr} Let $\theta : \scB \mono \scA$ be an embedding. Given $\theta^{k-1}$ then $\theta^k$ is fully specified by giving endpoints for each $\rest {\theta^k} y$ in the sense of \autoref{rmk:endpoints_of_embedding_fctr}. In particular $\theta^k$ is fully specified by its image. 
\end{rmk}

\begin{rmk}[Subfamily components are fully faithful] \label{rmk:subfamily_components_ff} The components $\theta^l$ of an embedding $\theta : \scB \mono \scA$ are necessarily fully faithful. Indeed, \autoref{claim:linearity_subfamilies} and \autoref{claim:edge_set_properties} together imply that for $(a,b) \in \edgeset(\tusU k_\scB(x \to y))$ we must have
\begin{equation}
\succ {\rest {\theta^{k+1}} x a, \rest {\theta^{k+1}} y b} = (\rest {\theta^{k+1}} x \succ {a,b}\ssoe , \rest {\theta^{k+1}} y \succ {a,b}\ttae )
\end{equation}
that is, $\theta^{k+1}$ commutes with successors. Thus, any edge with source and target in the image of $\theta^{k+1}$ must be in the image of $\theta^{k+1}$. This means, $\theta^{k+1}$ is full (and as a functor of posets necessarily faithful).
\end{rmk}

\subsection{Examples of embeddings}

\begin{eg}[Embeddings] \label{eg:subfamilies} \hfill
\begin{enumerate}
\item Define $\scA_a$ to be the $\SIvert 2 \cC$-family
\begin{restoretext}
\begingroup\sbox0{\includegraphics{test/page1.png}}\includegraphics[clip,trim=0 {.0\ht0} 0 {.0\ht0} ,width=\textwidth]{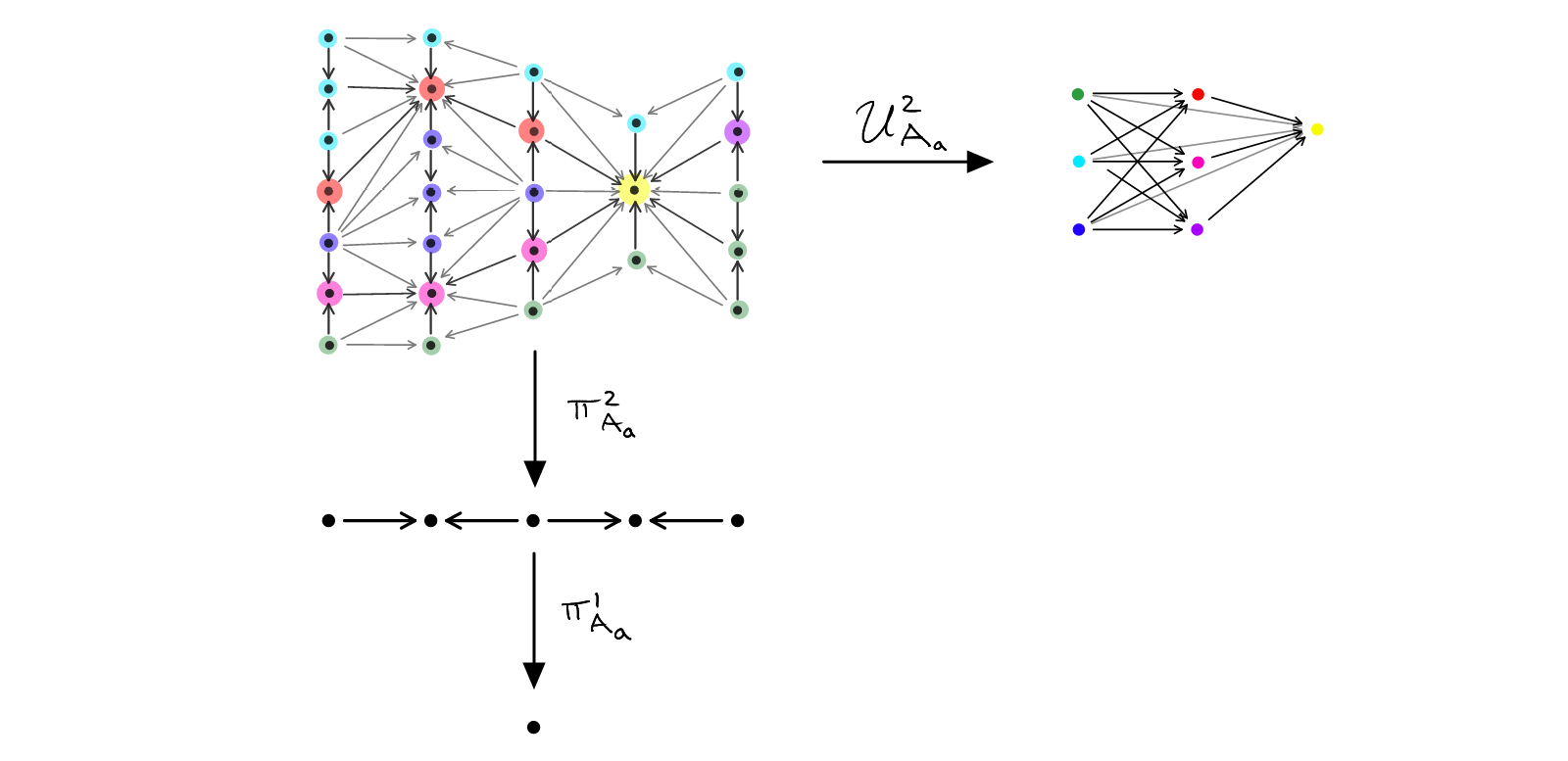}
\endgroup\end{restoretext}
Then, we can define an embedding $\theta_a : \scB_a \mono \scA_a$, for $\scB_a : \bnum{1} \to \SIvert 2 \cC$ by the following data
\begin{restoretext}
\begin{noverticalspace}
\begingroup\sbox0{\includegraphics{test/page1.png}}\includegraphics[clip,trim=0 {.0\ht0} 0 {.55\ht0} ,width=\textwidth]{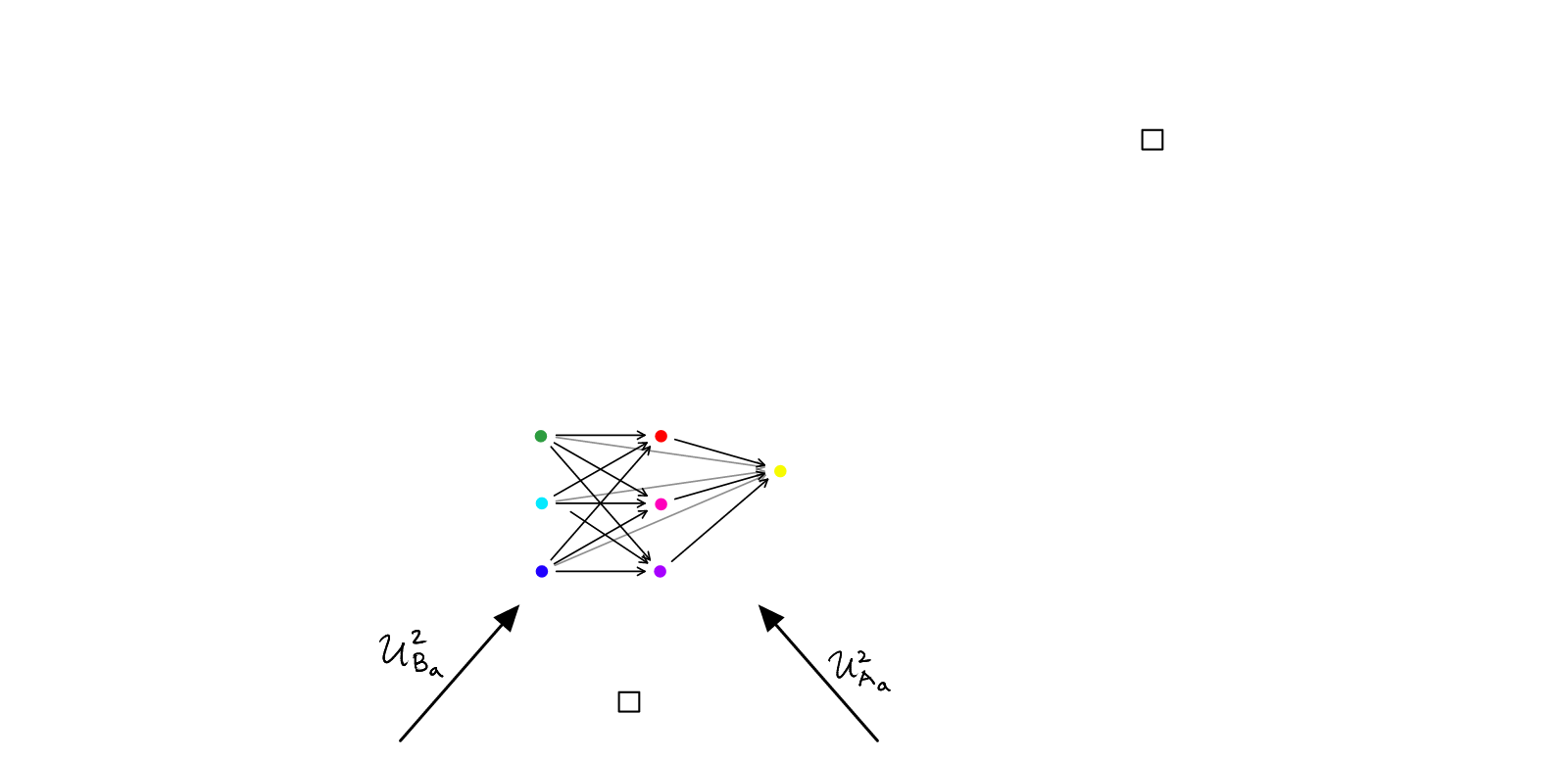}
\endgroup \\*
\begingroup\sbox0{\includegraphics{test/page1.png}}\includegraphics[clip,trim=0 {.0\ht0} 0 {.0\ht0} ,width=\textwidth]{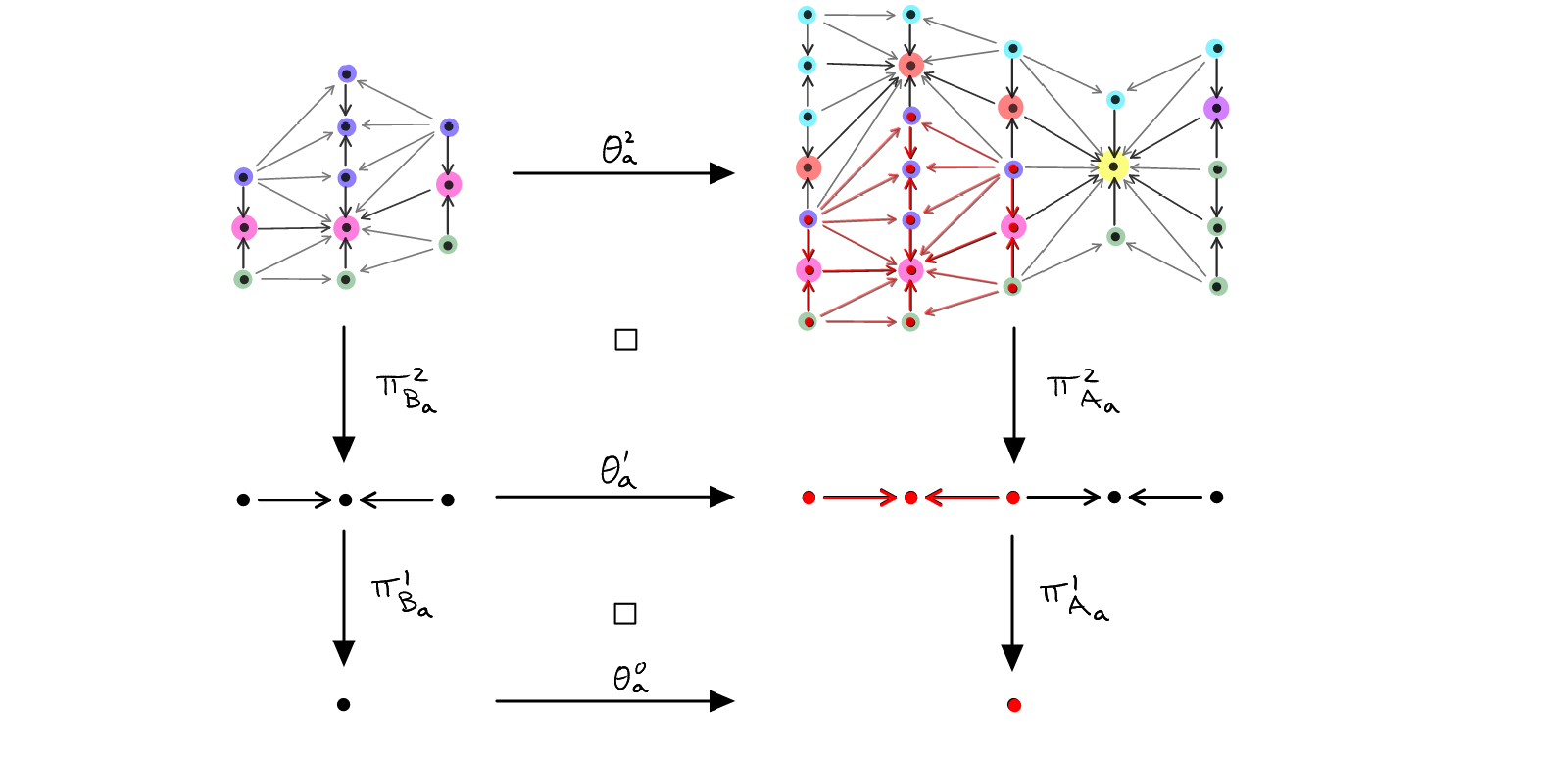}
\endgroup
\end{noverticalspace}
\end{restoretext}
Here, using \autoref{rmk:fibrewise_endpoints_of_embedding_fctr} we specified  the components $\theta^k_a$ by giving their image, whose objects and morphisms are marked in \cred{} in the above (note that $\theta^0_a$ must always be the terminal map for subcubes). Note that however the image of $\theta^{k-1}_a$ is determined by $\tpi k_\scA$ together with the image of $\theta^k_a$. This leads to the following observation: in general, given a family $\scA : X \to \SIvert n \cC$, then the up to base change by an isomorphism, embeddings $\theta : \scB \mono \scA$ (into $\scA$) are determined by the image of $\theta^n$. This will be illustrated by the next examples.
\item Let $\cC$ be the poset 
\begin{restoretext}
\begingroup\sbox0{\includegraphics{test/page1.png}}\includegraphics[clip,trim=0 {.4\ht0} 0 {.35\ht0} ,width=\textwidth]{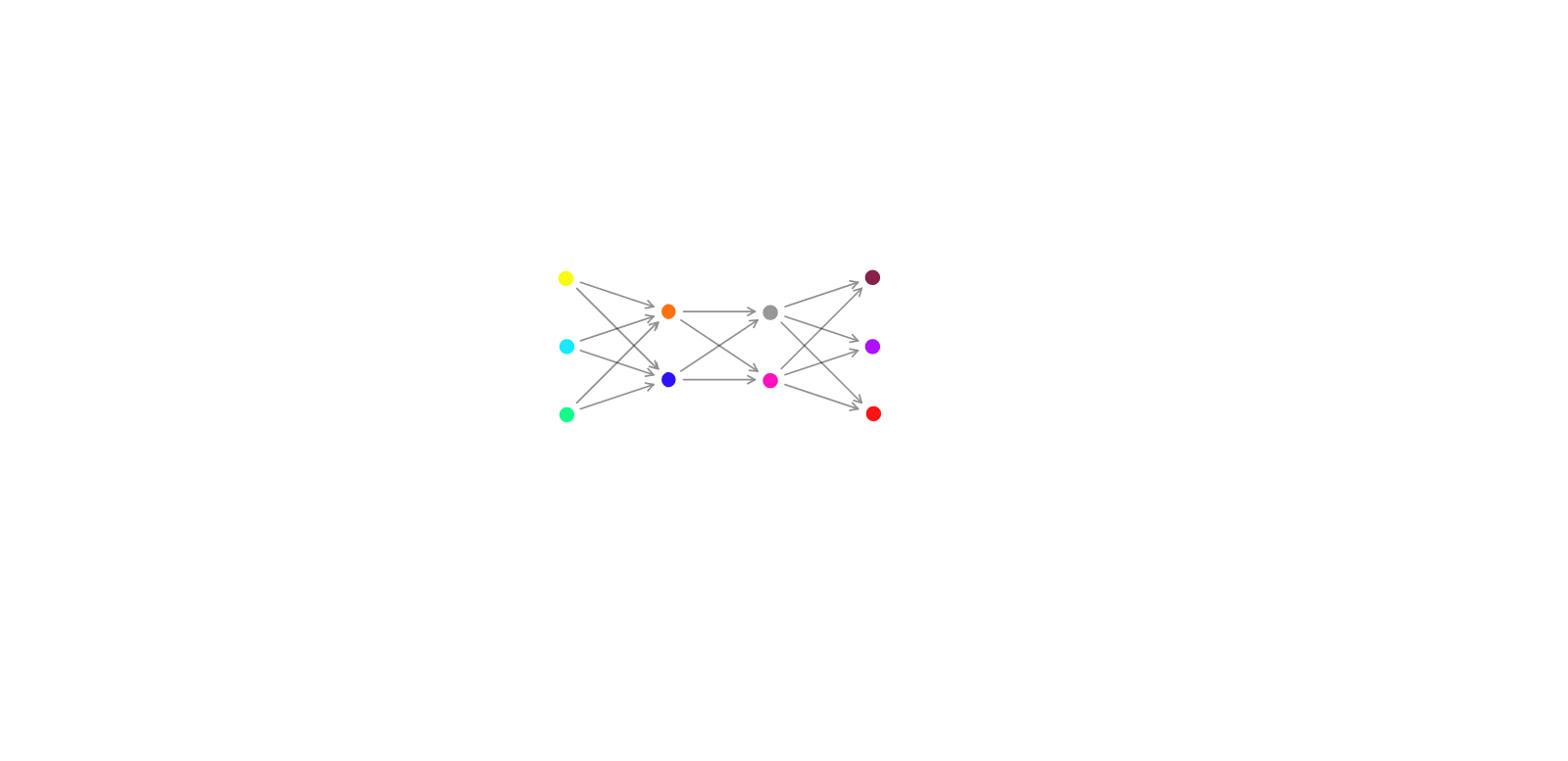}
\endgroup\end{restoretext}
To given an example in dimension 3, we define the family $\scA_b : \bnum{1} \to \SIvert 3 \cC$ as follows
\begin{restoretext}
\begin{noverticalspace}
\begingroup\sbox0{\includegraphics{test/page1.png}}\includegraphics[clip,trim=0 {.0\ht0} 0 {.0\ht0} ,width=.95\textwidth]{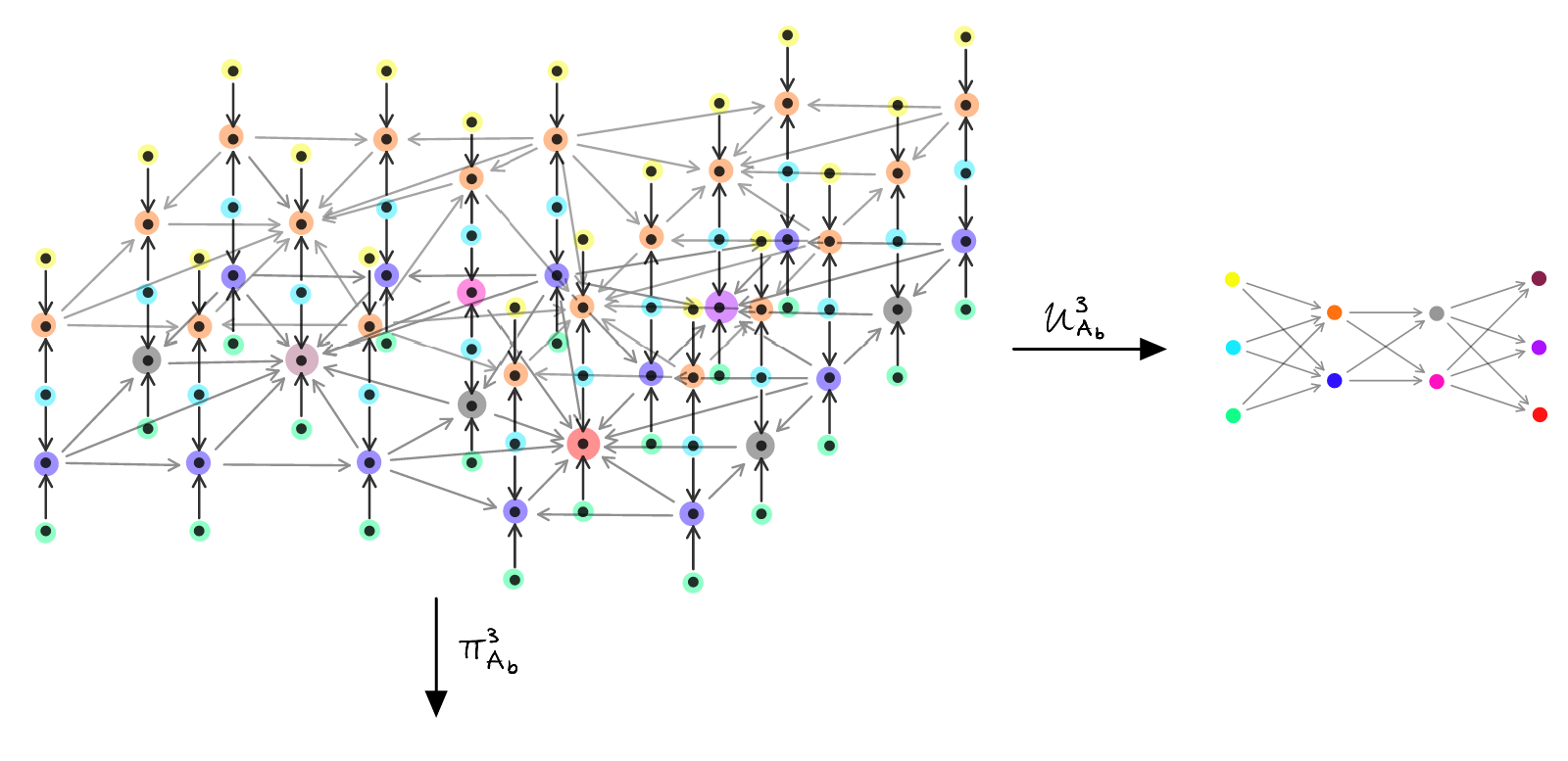}
\endgroup \\*
\begingroup\sbox0{\includegraphics{test/page1.png}}\includegraphics[clip,trim=0 {.07\ht0} 0 {.0\ht0} ,width=.95\textwidth]{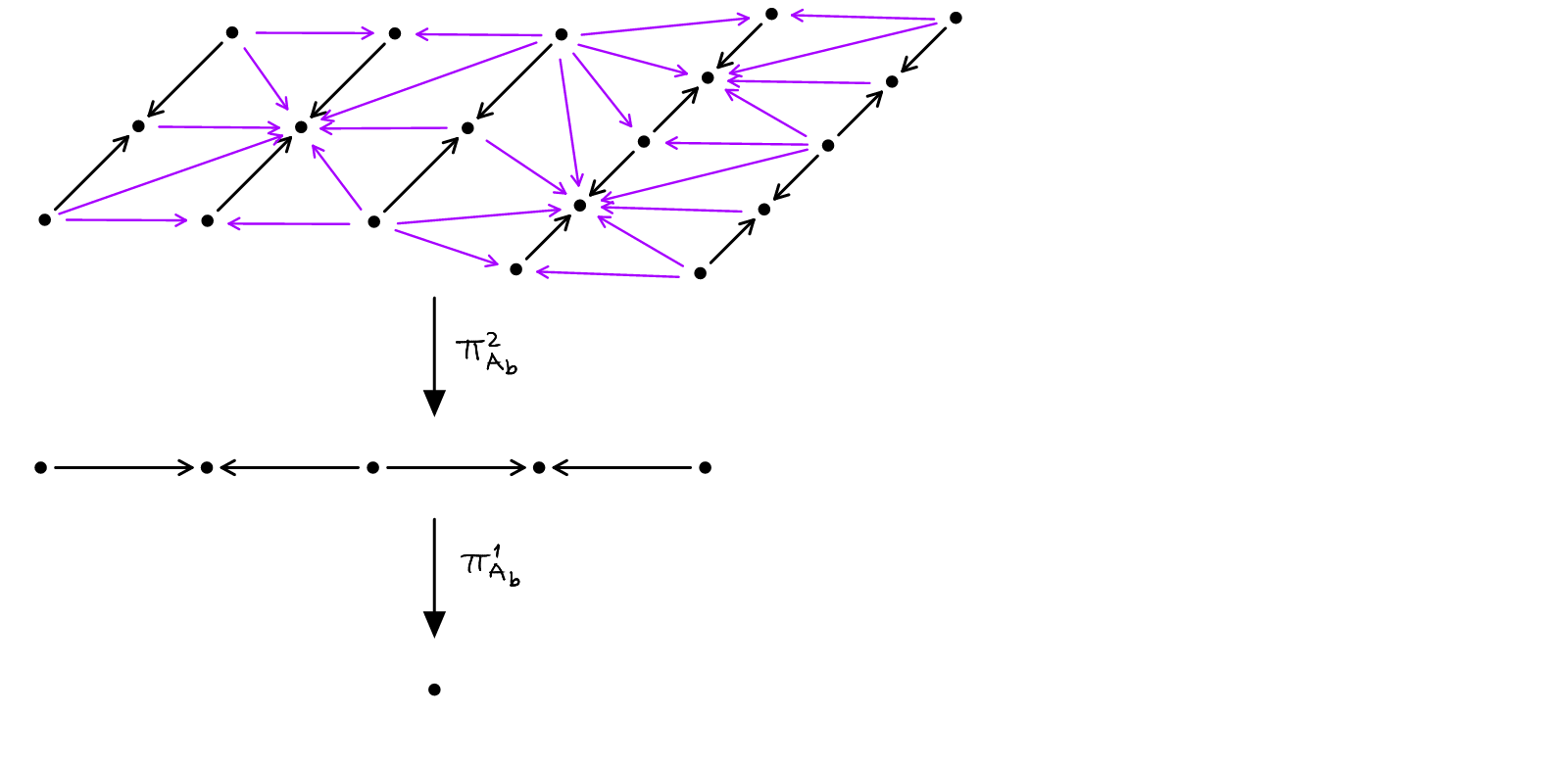}
\endgroup
\end{noverticalspace}
\end{restoretext}
Note that geometrically (cf. \autoref{ssec:coloring}) the situation that the example depicts is rather simple, and the reader is invited to sketch the corresponding colored $3$-cube (which will contain two $2$-dimensional, two $1$-dimensional and three $0$-dimensional singularities). For clarity, we illustrate $\tsG 3(\scA_b)$ again without the labelling colors
\begin{restoretext}
\begingroup\sbox0{\includegraphics{test/page1.png}}\includegraphics[clip,trim=0 {.15\ht0} 0 {.1\ht0} ,width=\textwidth]{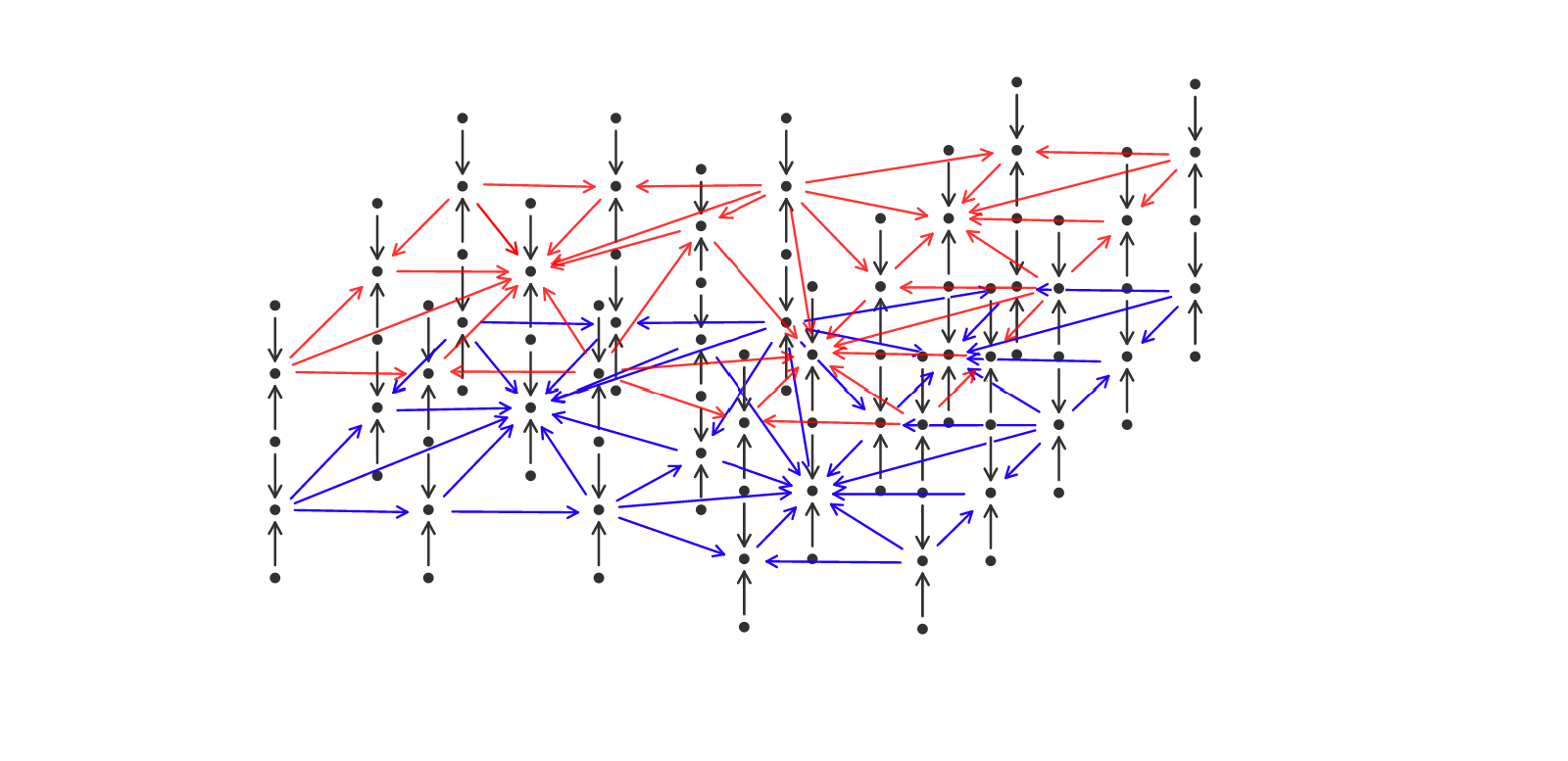}
\endgroup\end{restoretext}
Here, both \cred{} and \cblue{} arrows record arrows between singular heights in $\tsG 3(\scA_b)$.

Now, the following subposet of $\tsG 3(\scA_b)$ (marked by highlighting objects and arrows of the subposet in \cred{})
\begin{restoretext}
\begingroup\sbox0{\includegraphics{test/page1.png}}\includegraphics[clip,trim=0 {.15\ht0} 0 {.1\ht0} ,width=\textwidth]{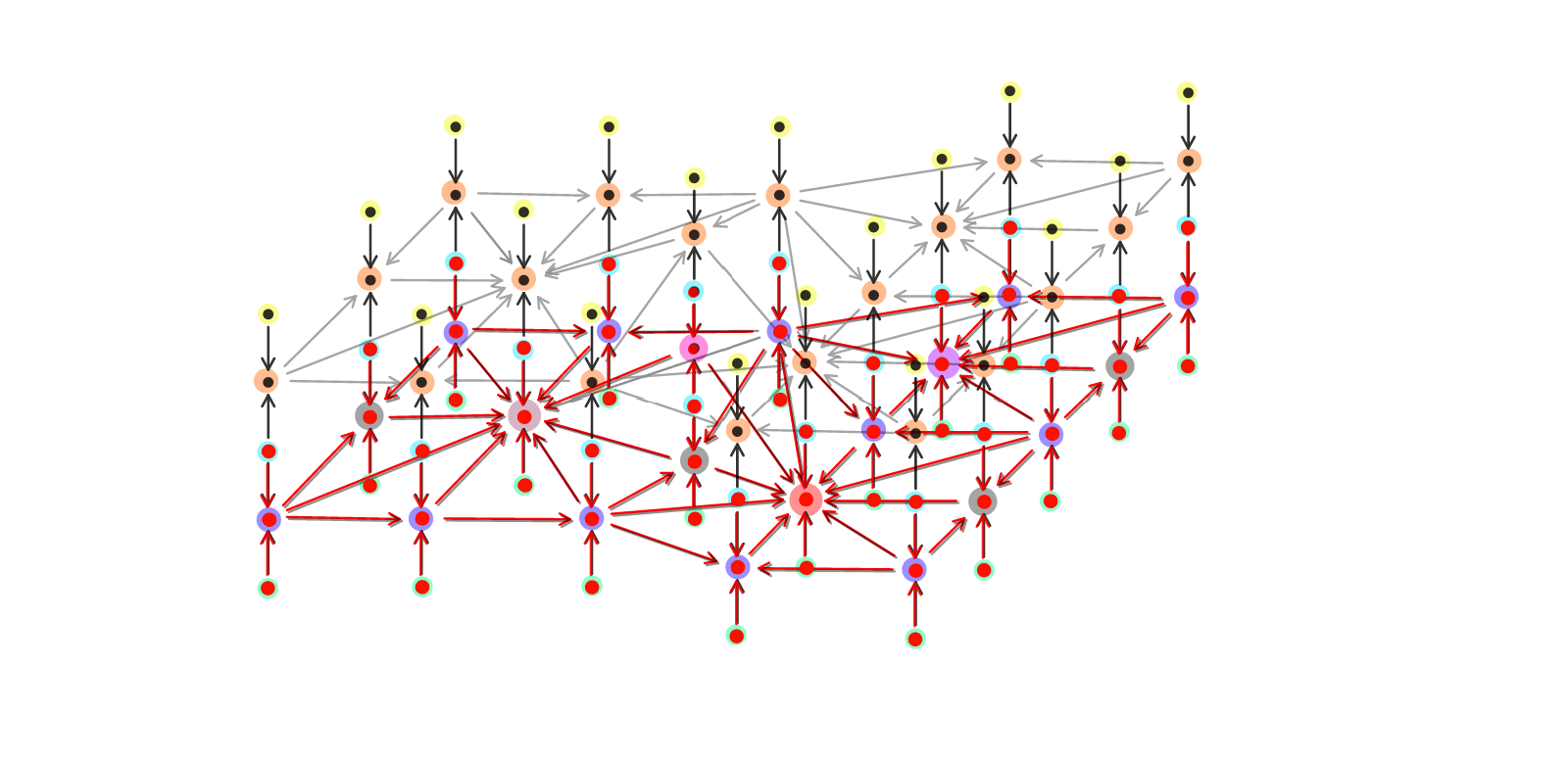}
\endgroup\end{restoretext}
determines an embedding by projecting it along $\tpi 3_{\scA_b}$, $\tpi 2_{\scA_b}$ and $\tpi 1_{\scA_b}$ as follows
\begin{restoretext}
\begin{noverticalspace}
\begingroup\sbox0{\includegraphics{test/page1.png}}\includegraphics[clip,trim=0 {.0\ht0} 0 {.0\ht0} ,width=\textwidth]{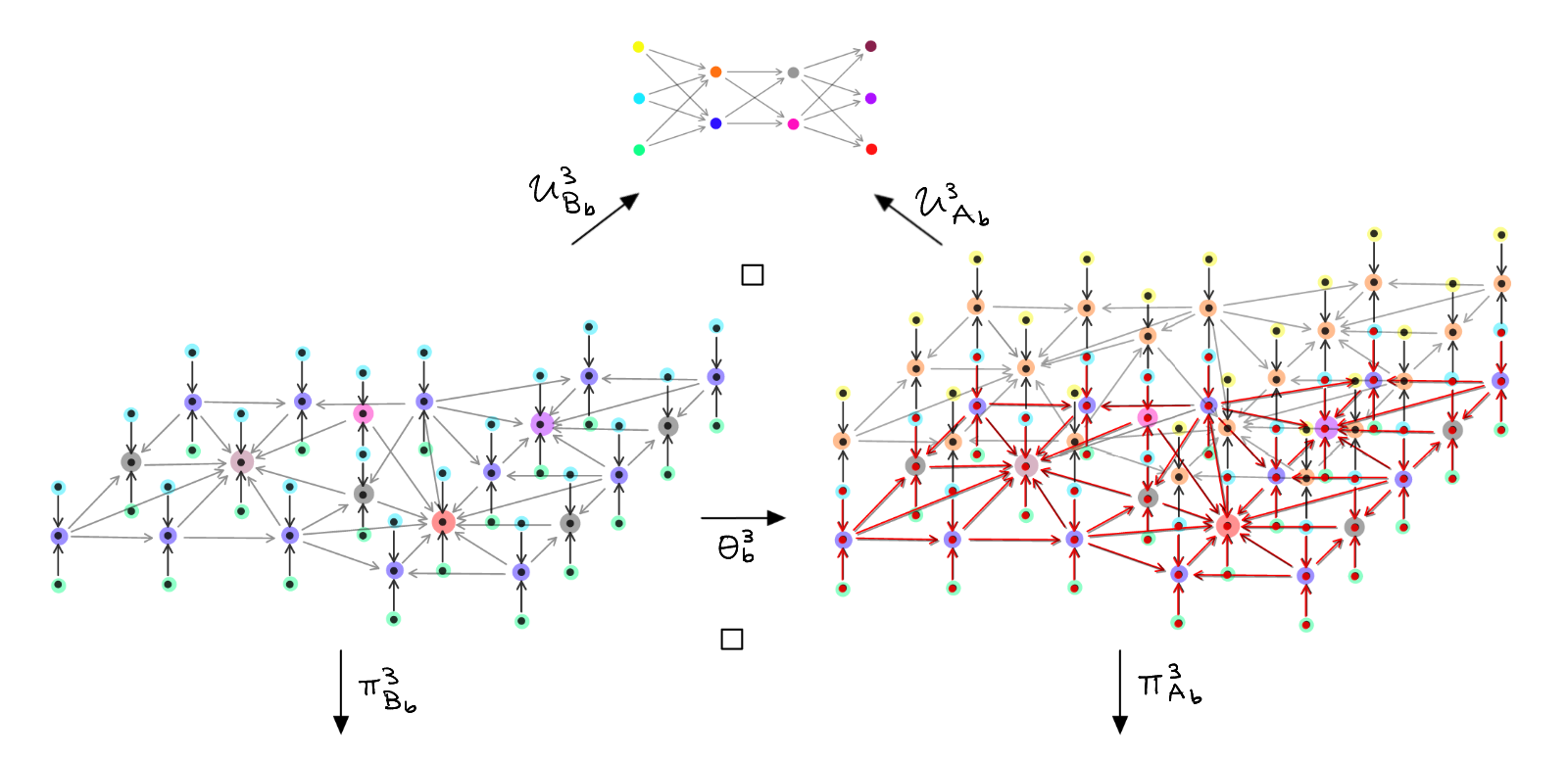}
\endgroup \\*
\begingroup\sbox0{\includegraphics{test/page1.png}}\includegraphics[clip,trim=0 {.25\ht0} 0 {.0\ht0} ,width=\textwidth]{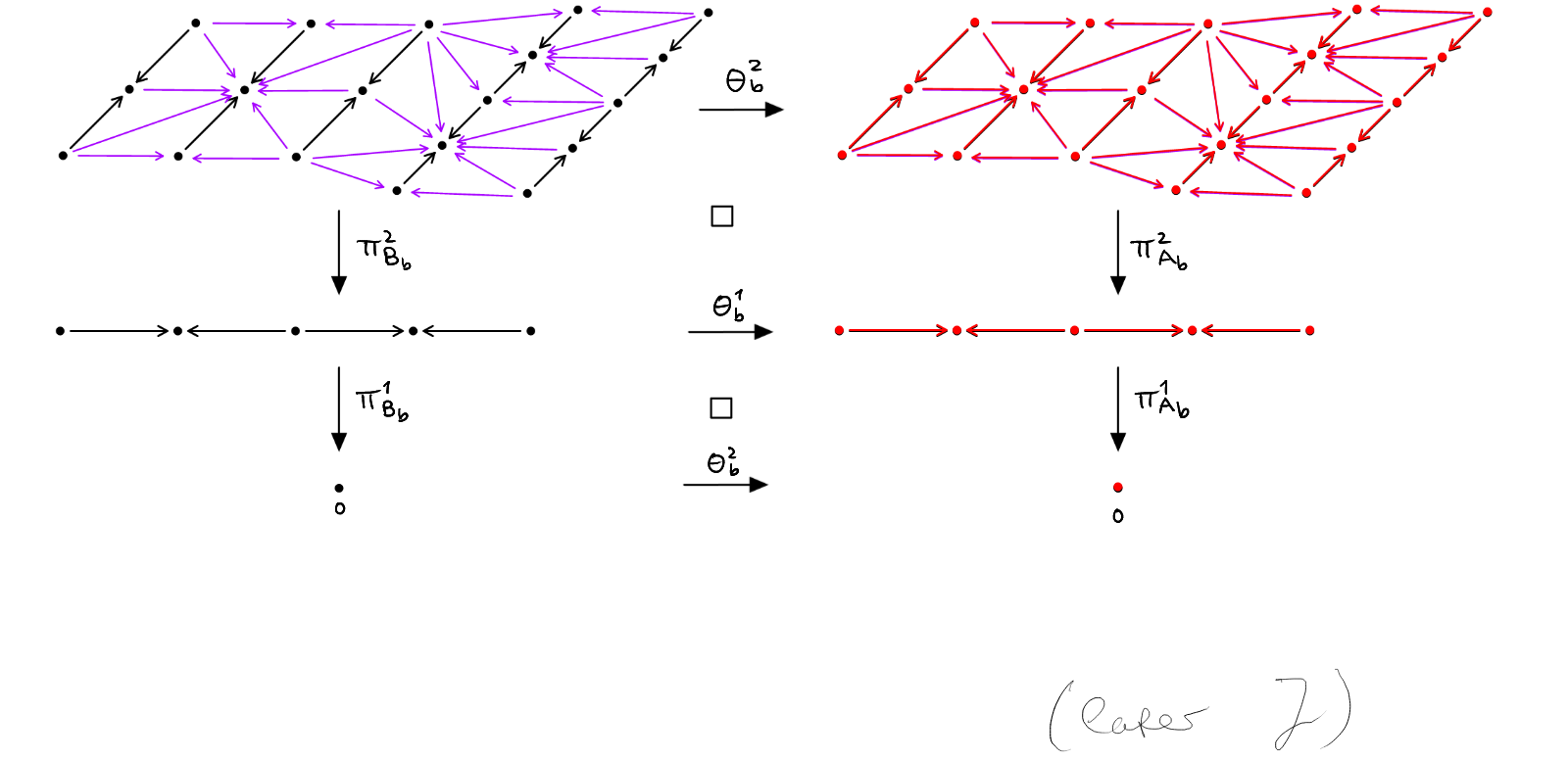}
\endgroup
\end{noverticalspace}
\end{restoretext}
up to base change isomorphism (in the above we chose the base space to be $\bnum{1}$ making $\theta^0$ a subposet inclusion. But we could have chosen it to be any other one-object poset, cf. \autoref{rmk:material_set_theory}). By default, if we define $\theta : \scB \mono \scA$ by $\im(\theta^n)$, then we will assume $\theta^0$ to be a subposet inclusion thereby fixing the base space $\tsG 0(\scB)$.

\item Similarly, the following
\begin{restoretext}
\begingroup\sbox0{\includegraphics{test/page1.png}}\includegraphics[clip,trim=0 {.15\ht0} 0 {.1\ht0} ,width=\textwidth]{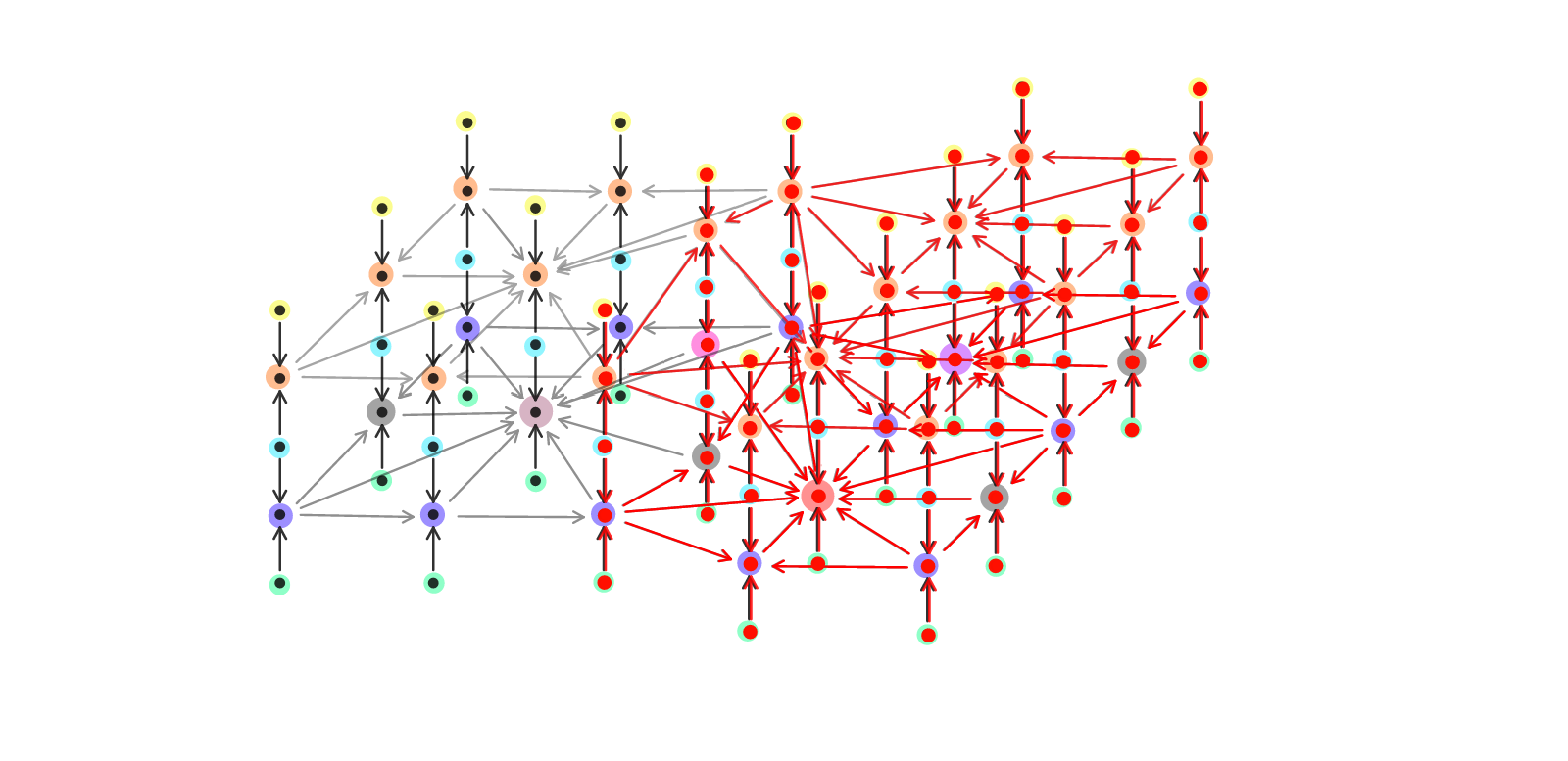}
\endgroup\end{restoretext}
denotes an embedding $\phi_b : \scC_b \mono \scA_b$, when interpreted as subposet of $\sG^3(\scA_b)$. But note that the latter equals $\sG^2(\tsU 1_{\scA_b})$. Relative to $\tsU 1_{\scA_b}$, it thus denotes an embedding $\widetilde\phi_b : \tsU 1_{\scC_b} \mono \tsU 1_{\scA_b}$ with components $\widetilde\phi_b^k = \phi^{k+1}_b$. Namely, we obtain the following embedding
\begin{restoretext}
\begin{noverticalspace}
\begingroup\sbox0{\includegraphics{test/page1.png}}\includegraphics[clip,trim=0 {.0\ht0} 0 {.0\ht0} ,width=\textwidth]{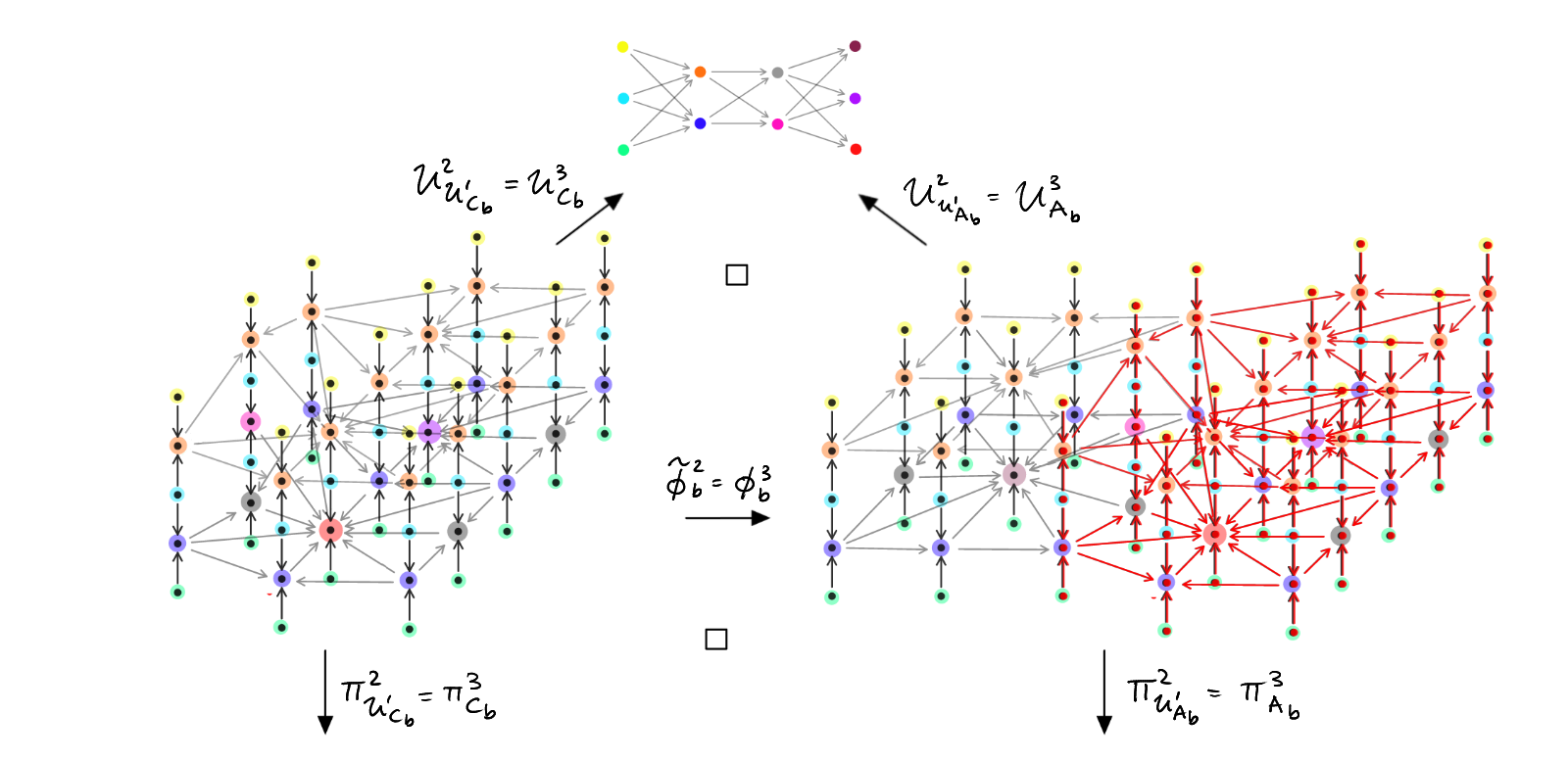}
\endgroup\\*
\begingroup\sbox0{\includegraphics{test/page1.png}}\includegraphics[clip,trim=0 {.5\ht0} 0 {.0\ht0} ,width=\textwidth]{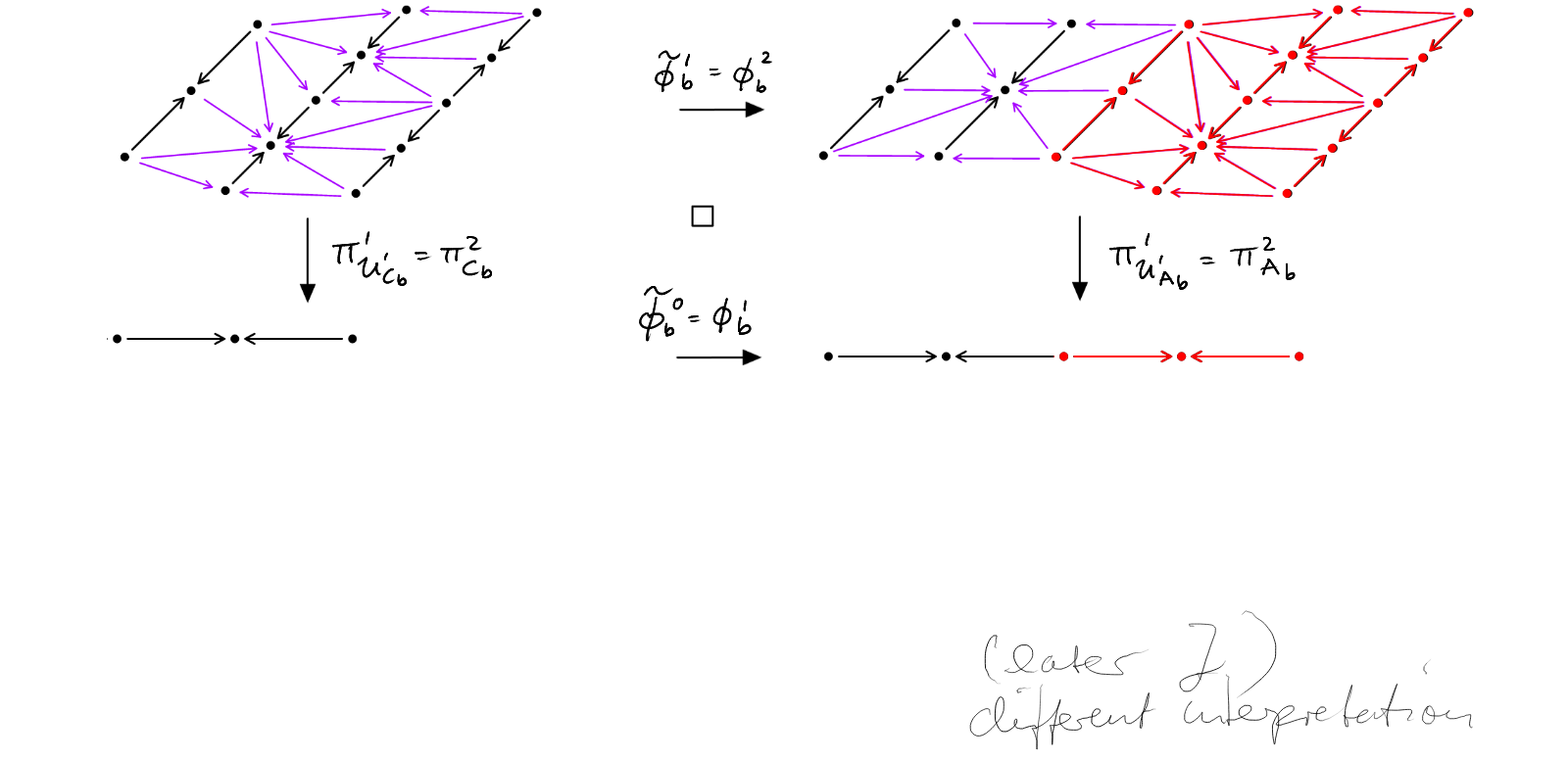}
\end{noverticalspace}
\endgroup\end{restoretext}
\end{enumerate}
\end{eg}

Recall the notion of composition for multi-level base changes (cf. \autoref{defn:multilevel_base_change}). From the definition of embeddings it can be seen that when $\theta : \scB \mono \scA$ and $\phi: \scC \mono \scB$ we obtain $\theta\phi : \scC \mono \scA$ which has components
\begin{equation} \label{rmk:subfamilies_compose}
(\theta\phi)^l = \theta^l \phi^l
\end{equation}
We give an example.

\begin{eg}[Composition of embeddings] \label{eg:composition_of_subbund} Consider the embedding $\phi_a : \scC_a \mono \scB_a$ (where $\scB_a$ was defined in \autoref{eg:subfamilies}) defined by
\begin{restoretext}
\begingroup\sbox0{\includegraphics{test/page1.png}}\includegraphics[clip,trim=0 {.35\ht0} 0 {.3\ht0} ,width=\textwidth]{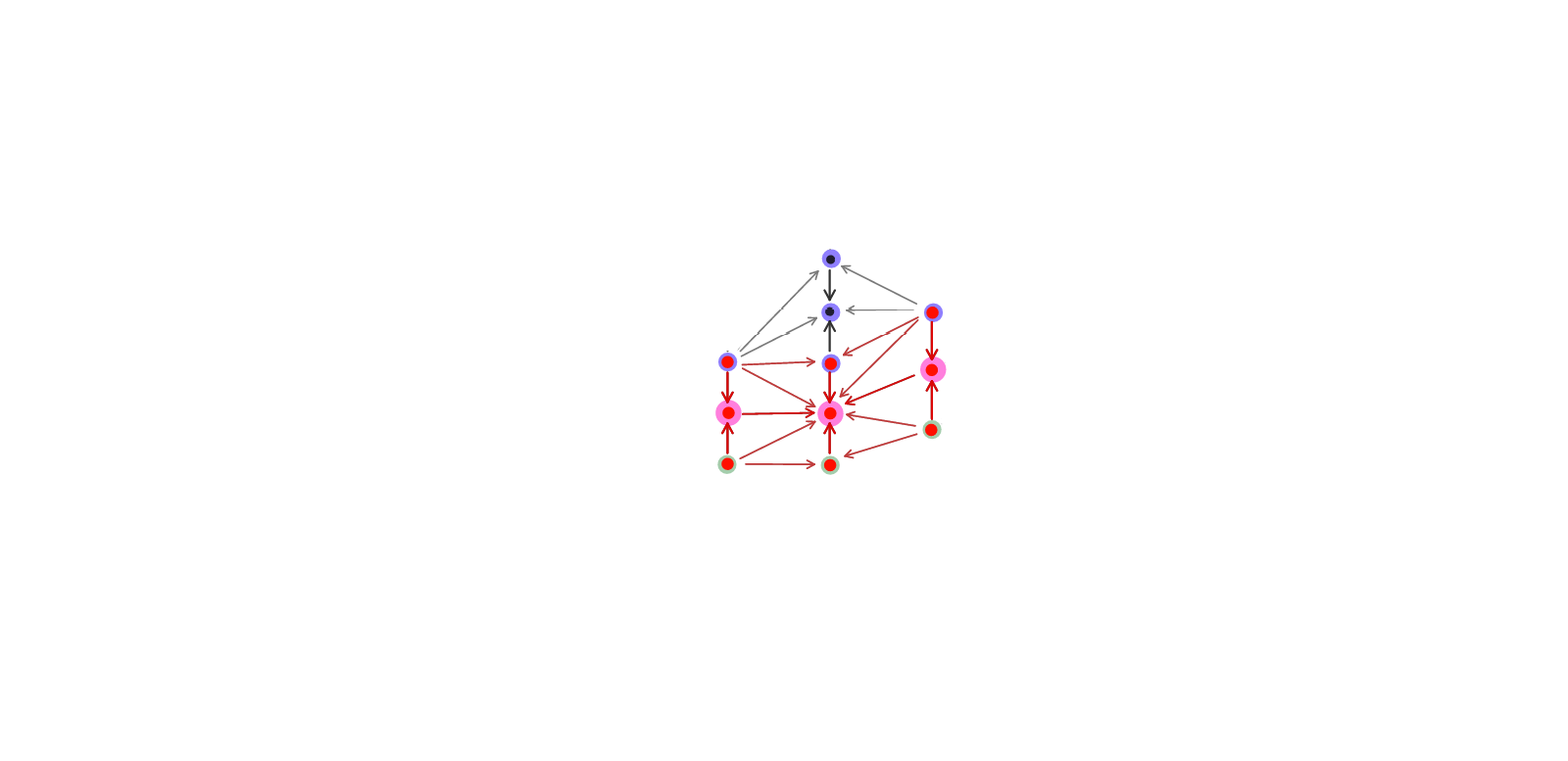}
\endgroup\end{restoretext}
Then the composition $\chi_a := \theta_a\phi_a : \scC_a \mono \scA_a$ is the embedding
\begin{restoretext}
\begingroup\sbox0{\includegraphics{test/page1.png}}\includegraphics[clip,trim=0 {.25\ht0} 0 {.25\ht0} ,width=\textwidth]{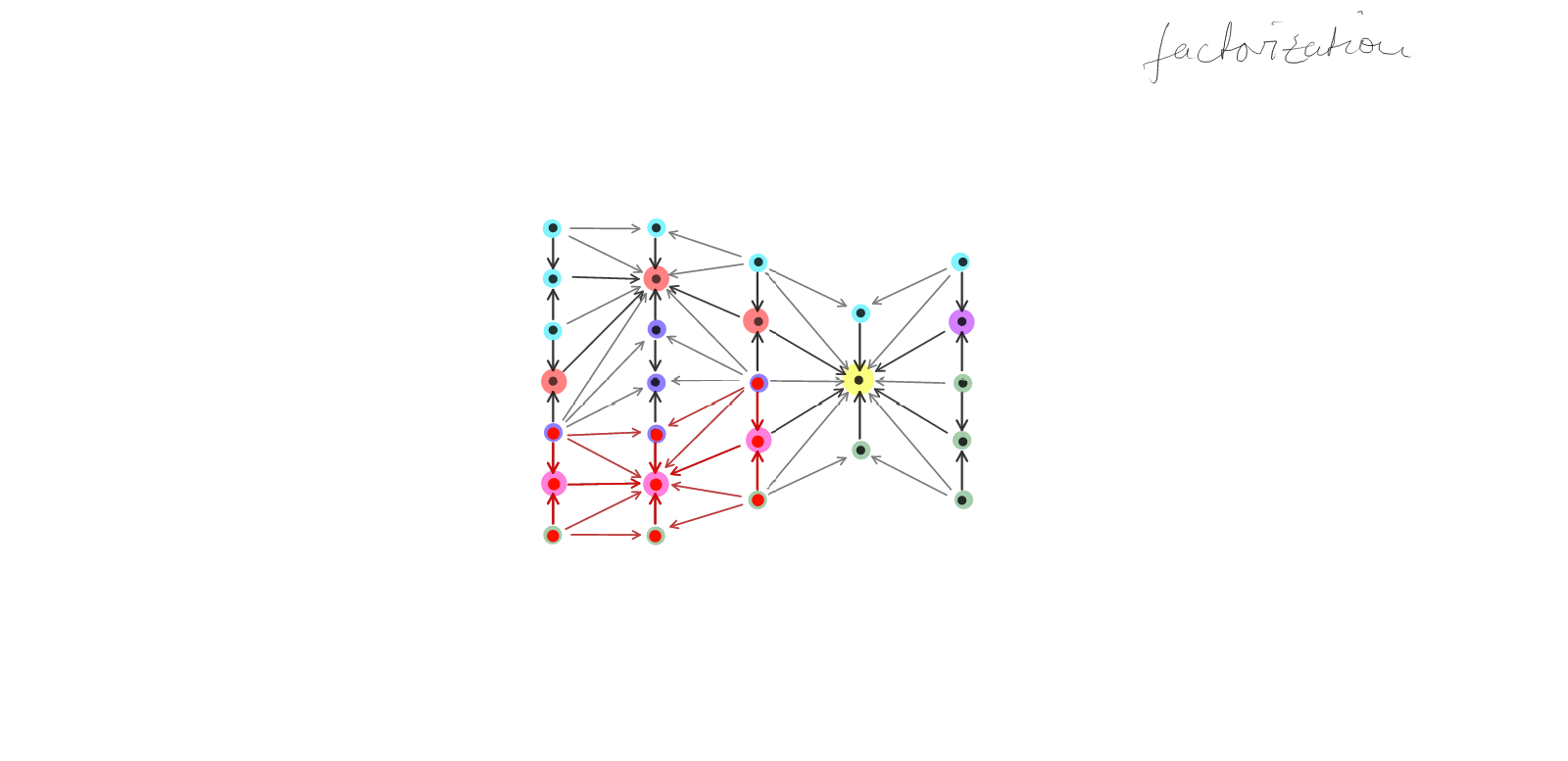}
\endgroup\end{restoretext}
\end{eg}

\subsection{Restriction embeddings}

We discuss embeddings obtained by restriction of the base space.

\begin{constr}[Subfamilies by restriction] \label{eg:subfamily_by_restriction} Let $\scA : X \to \SIvert n \cC$ and let $H : Y \into X$ be a subposet inclusion). Then we can obtain an embedding $\restemb_H : \rest \scA Y \mono  \scA$, called the \textit{restriction embedding} of $\scA$ to $Y$ by setting $\restemb^k_H = \tsG k(H)$ (cf. \autoref{constr:unpacking_collapse}) for $0 \leq k \leq n$.
\end{constr}

\begin{eg}[Subfamilies by restriction] The embedding $\widetilde\phi_b : \tsU 1_{\scC_b} \mono \tsU 1_{\scA_b}$ constructed in \autoref{eg:subfamilies} is an example of an embedding by restriction, namely, setting $H = \phi^1_b$ we find
\begin{equation}
\widetilde\phi_b = \restemb_{H}
\end{equation}
\end{eg}

\subsection{Factorisation of embeddings}

We discuss the (straight-forward) necessary and sufficient conditions for one embedding to factor through another, and more importantly, introduce notation for this case.

\begin{claim}[Factorisation of embeddings] \label{claim:factorisation_subfamilies} Let $\scA : X \to \SIvert n \cC$, $\scB : Y \to \SIvert n \cC$, $\scC : Z \to \SIvert n \cC$ and $\theta : B \mono A$, $\phi : C \mono A$ be embeddings. There is an embedding $\chi : C \mono B$ such that $\theta^k \chi^k = \phi^k$ if and only if $\im(\phi^k) \subset \im(\theta^k)$ for all $0\leq k \leq n$. The latter holds if and only if $\im(\phi^n) \subset \im(\theta^n)$. In this case we write $\phi \mono \theta$ and denote $\chi$ by $\theta\inv\phi$.

\proof If there is $\chi : C \mono B$ such that $\theta^k \chi^k = \phi^k$ then $\im(\phi^k) \subset \im(\theta^k)$. Conversely, assume $\im(\phi^k) \subset \im(\theta^k)$ and define $\chi : C \mono B$ componentwise by setting
\begin{equation}
\chi^k := (\theta^k)\inv\phi^k : \tsG k(\cC) \to \tsG k(\cB)
\end{equation}
This satisfies the conditions of an embedding since both $\phi$ and $\theta$ do.
\qed
\end{claim}

\begin{eg}[Factorisation of embeddings] \hfill
\begin{enumerate}
\item The embedding $\chi_a$ and $\theta_a$ defined in \autoref{eg:composition_of_subbund} satisfy that $\im(\chi^2_a)$ is contained in $\im(\theta^2_a)$. Thus, $\chi_a$ factors through $\theta_a$ and in this case we have
\begin{equation}
(\theta_a)\inv \chi_a = \phi_a
\end{equation}
as defined in \autoref{eg:composition_of_subbund}.
\item The embedding $\chi_b : \scD_b \mono \scA_b$ (with $\scA_b$ defined in \autoref{eg:subfamilies}) is defined by 
\begin{restoretext}
\begingroup\sbox0{\includegraphics{test/page1.png}}\includegraphics[clip,trim=0 {.15\ht0} 0 {.05\ht0} ,width=\textwidth]{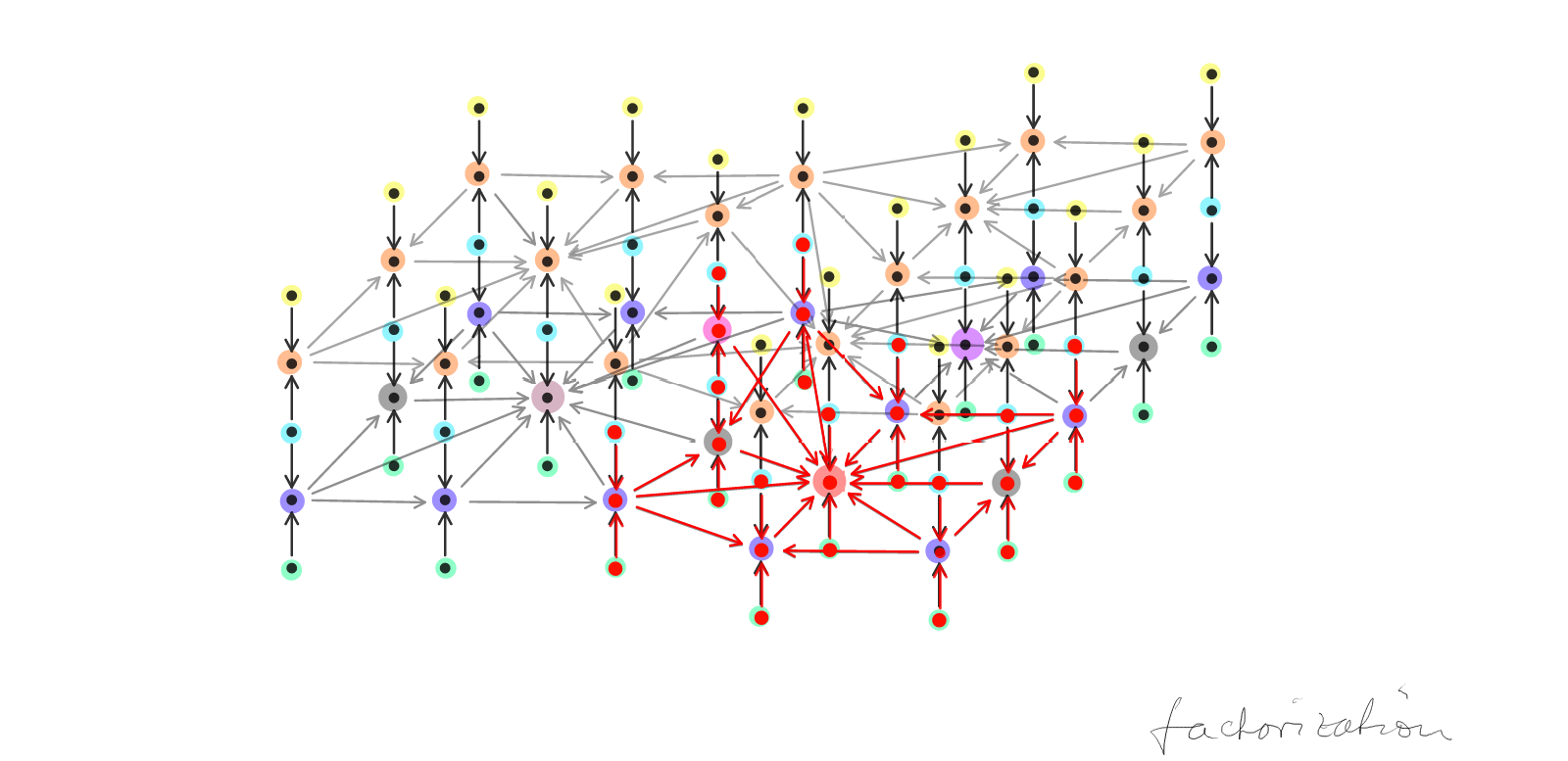}
\endgroup\end{restoretext}
Thus the image of $\chi^3_b$ is contained both in the image of $\phi^3_b$ and $\theta^3_b$. Thus $\chi_b \mono \phi_b$ and $\chi_b \mono \theta_b$.
\end{enumerate}
\end{eg}

\section{Endpoints} \label{sec:emb_end}

We will now classify embeddings as defined in the previous section. By the end of this section we will see that any (multi-level) embedding can be decomposed into a composition of ($k$-level) embeddings determined by so-called $k$-level endpoint sections.

\subsection{Endpoint sections}

We start with the definition of (open) sections.

\begin{defn}[Sections and open sections] \label{defn:sections} Let $\scA : X \to \SI$ be an \SI-family. A \textit{section} $q$ \textit{in $\scA$} is a functor $q : X \to \sG(\scA)$ such that 
\begin{equation}
\pi_\scA q = \id_X
\end{equation}
Note that to any such section we can associate a function $\secp q : X \to \lN$ (cf. \autoref{notn:total_poset_projs}) defined to map
\begin{equation}
x \to \secp q (x)
\end{equation}
A section $q$ is called \textit{open} if for all $x \in X$, $\secp q (x) \in \regcont(\scA(x))$.
\end{defn}

\begin{notn}[Depicting sections] \label{notn:depicting_sections} Since a section is determined by its value in each fibre, we will depict sections by marking points in each fibre of an \SI-bundle. This is illustrated by the next section.
\end{notn}

\begin{rmk}[Open sections have lifts] \label{rmk:src_and_tgt_have_lifts} Open sections $q : X \to \sG(\scA)$ have \gls{lifts} (cf. \autoref{defn:having_lifts}). Indeed, if $(y,a) \to (x,q(x))$ then $y \to x$ and $(a,q(x)) \in E(\scA(y\to x))$. By \eqref{eq:defn_order_realisation_3} $a$ is the unique element of $\scA(y)$ with that property. By functoriality of $q$ we must thus have $q(y) = a$.
\end{rmk}

\subsection{Source and target section}

Source and target section are special open section that any $\SI$-family admits. We define them, and give examples.

\begin{constr}[Source and target endpoints] \label{constr:source_and_target_inclusion} Let $\scA : X \to \SI$ be an $\SI$-family over a poset $X$. We define an open section $\msrc_A : X \to \sG(A)$ (called the \textit{source section} of $\scA$) to be the map
\begin{equation}
x \in X \quad \mapsto \quad (x,0) \in \sG(A)
\end{equation}
and an open section $\mtgt_A : X \to \sG(A)$ (called the \textit{target section} of $\scA$) by mapping
\begin{equation}
x \in X \quad \mapsto \quad (x,\iH_{A(x)}) \in \sG(A)
\end{equation}
Note that both maps are indeed functors of posets: if $(x \to y) \in \mor(X)$ then both $(x,0) \to (y,0)$ and $(x,\iH_{A(x)})\to (y,\iH_{A(y)})$ in $\sG(\scA)$. To see this, by \autoref{constr:SI_families} it suffices to remind ourselves that  $(0,0) \in E(\scA(x\to y))$ and $(\iH_{A(x)}, \iH_{A(y)}) \in E(\scA(x \to y))$ (for instance by arguments in the proof of \autoref{claim:edge_set_properties}). 
\end{constr}

\begin{eg}[Sections, open sections, source and target] \label{eg:sections} \hfill
\begin{enumerate}
\item Recall $\tpi 2_{\scA_a}$ from \autoref{eg:subfamilies}. We define two sections $q_1 , q_2 : X \to \tsG 2(\scA_a)$ in the \SI-family $\tusU 1_{\scA_a}$ as follows
\begin{restoretext}
\begingroup\sbox0{\includegraphics{test/page1.png}}\includegraphics[clip,trim=0 {.15\ht0} 0 {.15\ht0} ,width=\textwidth]{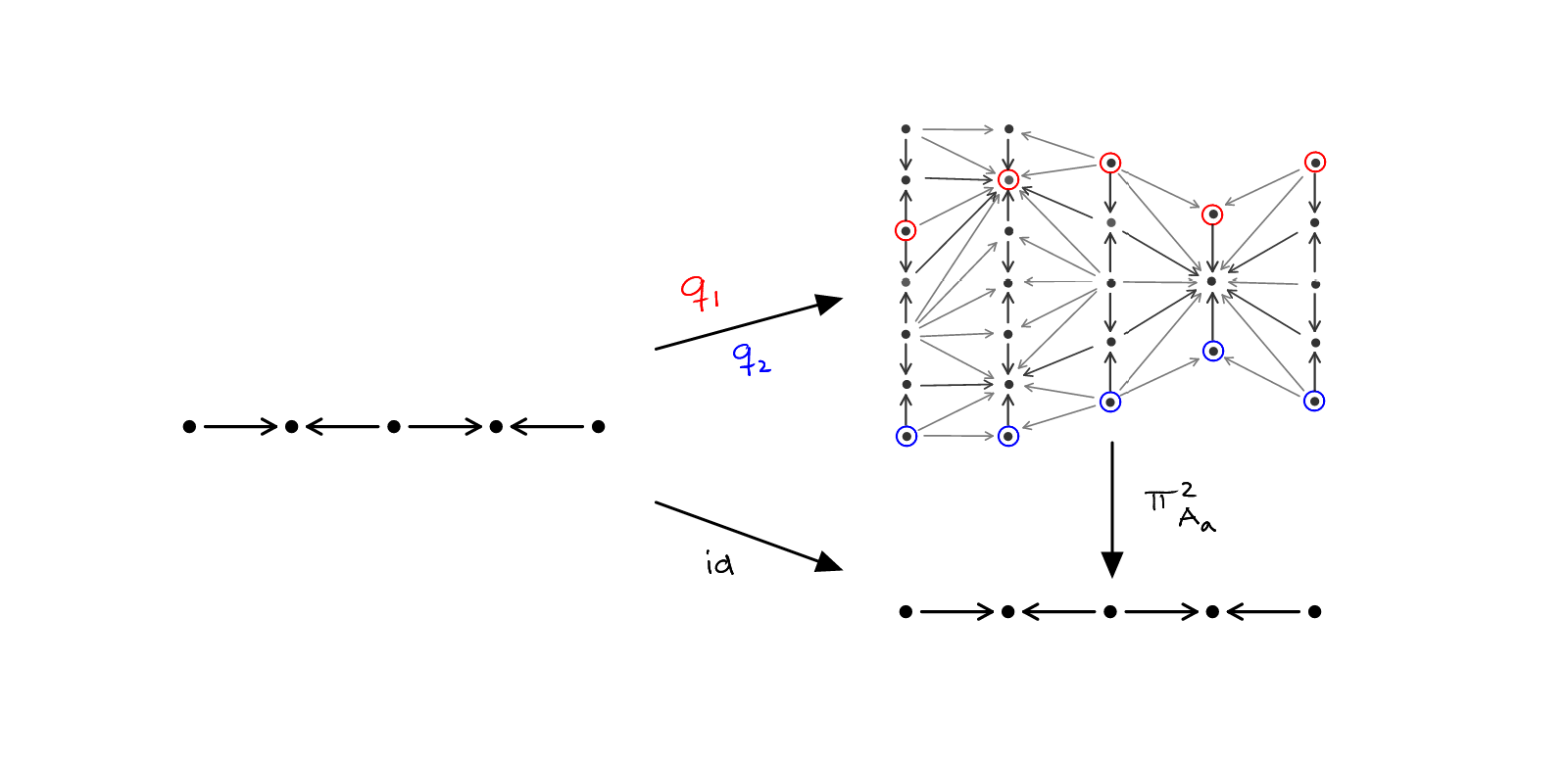}
\endgroup\end{restoretext}
The image of $q_1$ is marked by \cred{} circles, the image of $q_2$ is marked by \cblue{} circles (cf. \autoref{notn:depicting_sections}). $q_1$ is a section, but it is not an open section. $q_2$ is an open section. In fact, $q_2$ is the source section of $\tsU 1_{\scA_a}$.

\item Next, recalling $\scA_b$ from \autoref{eg:subfamilies}, consider the two sections $q^b_1 , q^b_2 : X \to \tsG 3(\scA_b)$ in in the \SI-family $\tusU 2_{\scA_b}$ as follows
\begin{restoretext}
\begingroup\sbox0{\includegraphics{test/page1.png}}\includegraphics[clip,trim=0 {.0\ht0} 0 {.0\ht0} ,width=\textwidth]{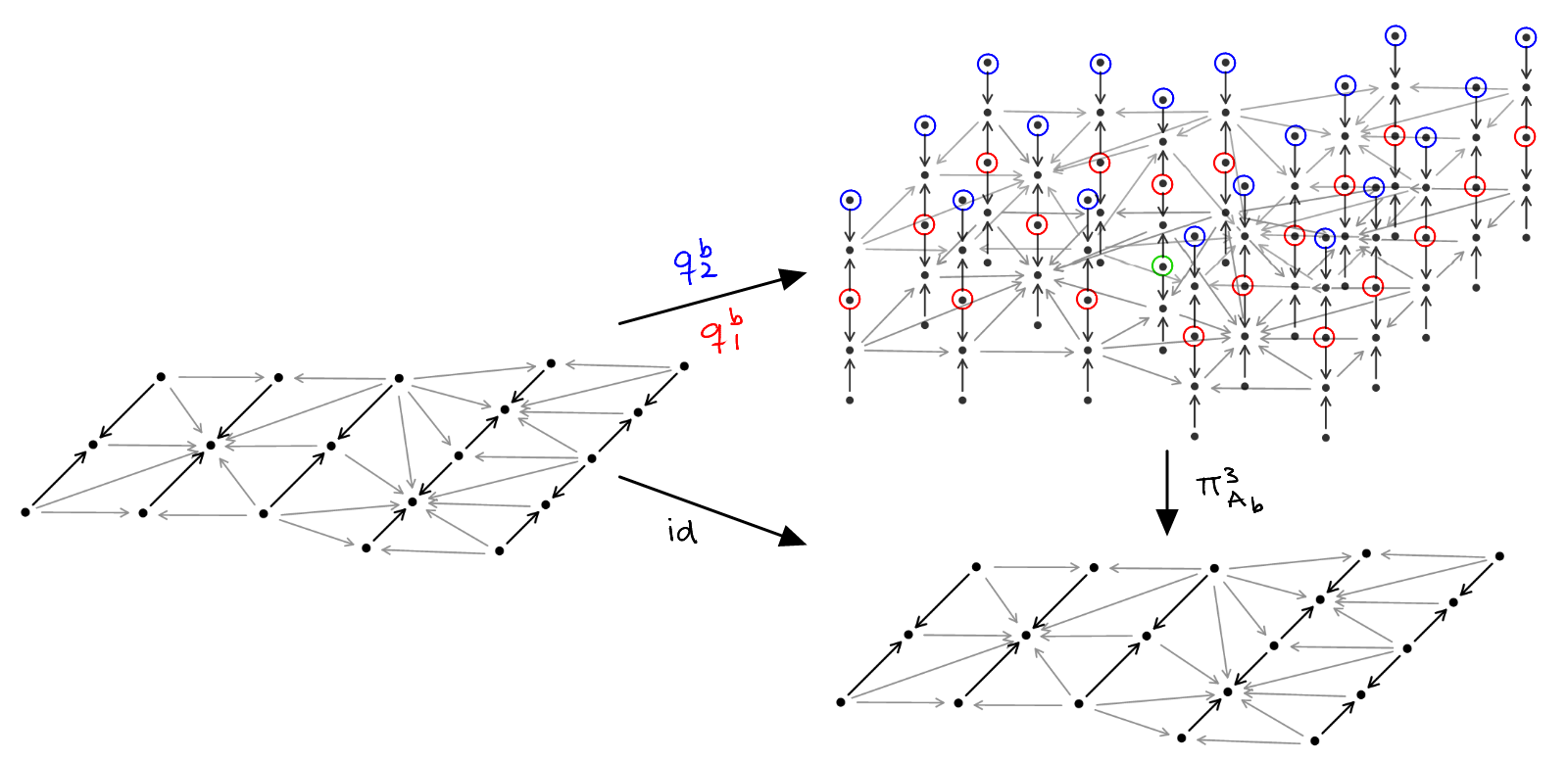}
\endgroup\end{restoretext}
The image of $q^b_1$ is marked by \cred{} circles, the image of $q^b_2$ is marked by \cblue{} circles. Both sections are open. In fact, $q^b_2$ is the target section of $\tsU 2_{\scA_b}$. Note that for instance replacing the value of $q^b_1$ in the fibre containing the \cdarkgreen{} circle by the point circled in \cdarkgreen{}, would not yield a valid section.
\end{enumerate}
\end{eg}

\begin{rmk}[Base change for source and target] \label{rmk:basechange_for_src_and_tgt} Let $\scA : X \to \SI$, and $H : Y \to X$. Then, since $\sG(H)$ is a fibrewise isomorphism, we find that
\begin{align}
\msrc_\scA H &= \sG(H) \msrc_{\scA H} \\
\mtgt_\scA H &= \sG(H) \mtgt_{\scA H}
\end{align}
This observation is generalised by \autoref{claim:factorising_subbund_through_endpoints}.
\end{rmk}

\subsection{Family embedding functors}

We define family embedding functors $\sJ^\scA\restsec{[q_-,q_+]}$. We start with the following observation.

\begin{rmk}[Sections determining embeddings] \label{rmk:sections_determine_subbund} Further to \autoref{eg:sections}, note that two open sections such as $q^b_2$ and $q^b_1$ of $\tusU 2_{\scA_b}$ fibrewise determine endpoints in the sense of \autoref{rmk:endpoints_of_embedding_fctr}, and thus an embedding functor for that fibre. Combining these embeddings gives an embedding $\scB \mono \tusU 2_{\scA_b}$ such that $\scB$ has the same base space as $\scA$. This is illustrated in the following picture in the case of the source section $\msrc_{\tusU 2_{\scA_b}}$ and the open section $q^b_1$:
\begin{restoretext}
\begingroup\sbox0{\includegraphics{test/page1.png}}\includegraphics[clip,trim=0 {.0\ht0} 0 {.0\ht0} ,width=\textwidth]{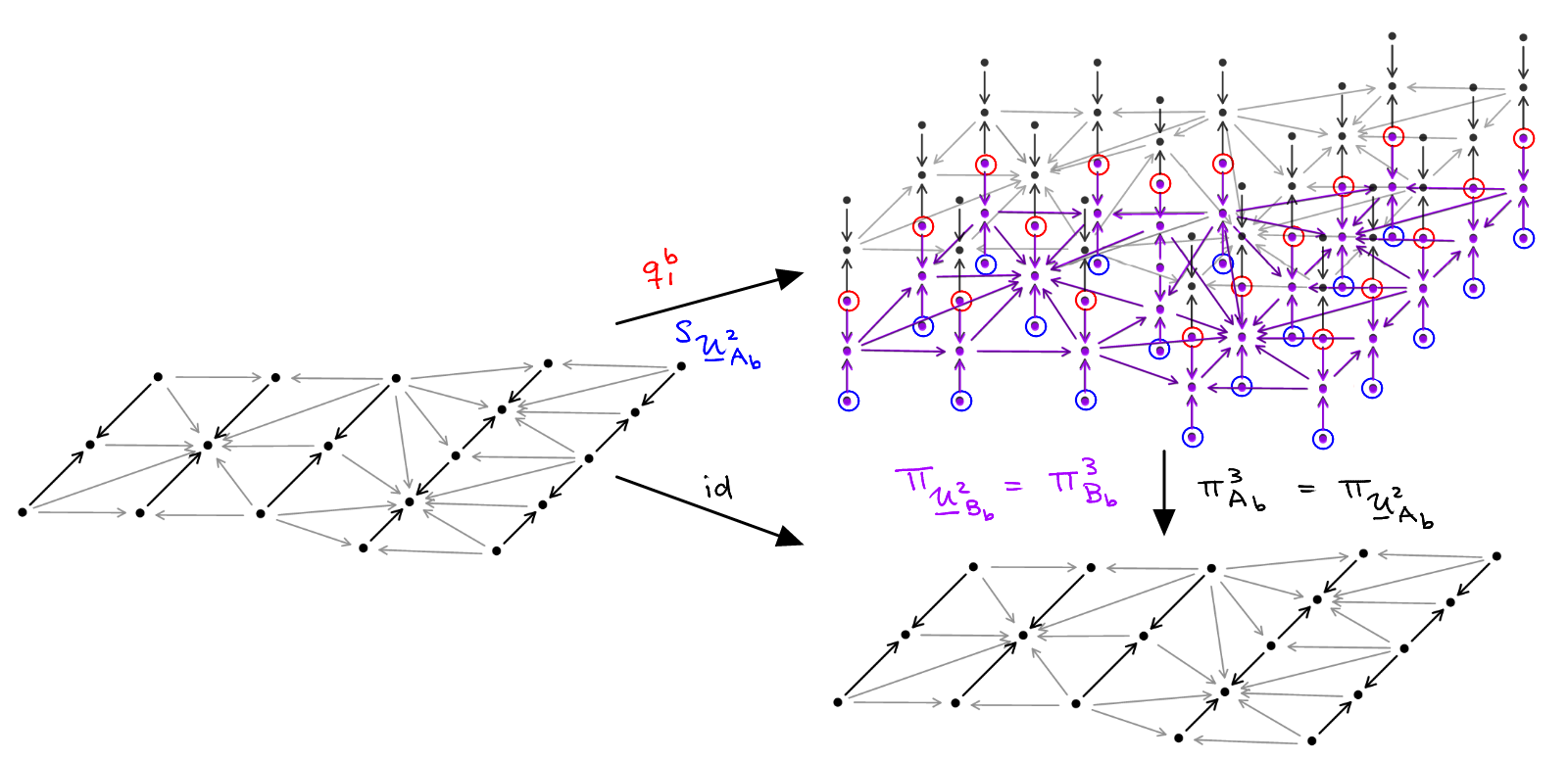}
\endgroup\end{restoretext}
In the above, the \SI-family marked in \cpurple{} is $\tusU 2_{\scB_b}$ from \autoref{eg:sections}. It is the subfamily of the \SI-family $\tusU 2_{\scA_b}$ which is delimited by $\msrc_{\tusU 2_{\scA_b}}$ and $q^b_1$. Note that it has indeed the same base space as its parent family $\tusU 2_{\scA_b}$. This observation is generalised in the next construction.
\end{rmk}

\begin{constr}[Family embedding functors] \label{constr:endpoint_inclusions}
Let $\scA : X \to \SI$. A continuous choice of endpoints in $\scA$ is a tuple of maps of open sections $q_+, q_- : X \to \sG(\scA)$ such that for all $x \in X$
\begin{equation}
\secp q_- (x) \leq \secp q_+(x)
\end{equation}
We will define an \SI-family $\scA\restsec{[q_-,q_+]}$ together with an injective $\sJ^\scA\restsec{[q_-,q_+]} : \sG(\scA\restsec{[q_-,q_+]}) \to \sG(\scA)$, called the \textit{family embedding functor for the endpoints $q_-,q_+$}, whose image is delimited by $q_-, q_+$ in the sense of the previous remark.

The construction is \stfwd{}. For $x \in X$ denote $E_x := [\secp q_-(x),\secp q_+(x)] \subset \scA(x)$. Note that for $(x \to y) \in \mor(X)$ by  $(\secp q_\pm(x), \secp q_\pm(y)) \in \edgeset(\scA(x\to y))$ and bimonotonicity of $\scA(x\to y)$ we have
\begin{equation} \label{eq:endpoint_boundaries}
(b,b') \in \edgeset(\scA(x\to y)) \quad \imp \quad \big(b \in E_x \iff b' \in E_y\big)
\end{equation}
We first construct a family $\scA\restsec{[q_-,q_+]}: X \to \SI$ called the (unlabelled) \textit{subfamily of $\scA$ with endpoints $q_-$, $q_+$}. We define $\scA\restsec{[q_-,q_+]}$ on objects $x \in X$ by setting
\begin{equation}
\iH_{\scA\restsec{[q_-,q_+]} (x)} = \frac{ q_+ (x) - q_- (x)}{2}
\end{equation}
and for $(x \to y)$ in $X$
\begin{equation}
\scA\restsec{[q_-,q_+]}(x \to y) : \scA\restsec{[q_-,q_+]}(x) \to \scA\restsec{[q_-,q_+]}(y)
\end{equation}
is the morphism of singular intervals defined by the mapping
\begin{equation}
b \in \singcont(\scA\restsec{[q_-,q_+]}(y)) \quad \mapsto \quad \scA(x \to y)(b + q_-(x)) - q_-(y)
\end{equation}
It is straight-forward to see that this is well-defined. Explicitly, we need to check the right hand side lives in $\singcont(\scA\restsec{[q_-,q_+]}(z))$ for all $b$. This is equivalent to requiring that for all $b$ we have
\begin{equation}
\scA(x \to y)(b + q_-(x)) \in E_y
\end{equation}
which in turn is equivalent to 
\begin{equation}
\scA(x \to y)(E_x) \subset E_y
\end{equation}
(note that here $\scA(y \to z)$ is only defined on the singular heights of $E_y$). The latter holds by \eqref{eq:defn_order_realisation_1} and \eqref{eq:endpoint_boundaries}.

Next, we verify that $\scA\restsec{[q_-,q_+]}$ as defined above is functorial. Assume $x \to y \to z  \in \mor(X)$. Then 
\begin{align*}
\scA\restsec{[q_-,q_+]}(x \to z)(b) &= \scA(x \to z)(b + q_-(x)) - q_-(z) \\
&= \scA(y \to z)\Big( \big(\scA(x \to y)(b + q_-(x)) - q_-(y)\big) + q_-(y)\Big) - q_-(z) \\
&= \scA\restsec{[q_-,q_+]} (y \to z)\big( \scA\restsec{[q_-,q_+]}(x \to y) (b)\big)
\end{align*}
as required. We have thus constructed $\scA\restsec{[q_-,q_+]} : X \to \SI$. 

We now define $\sJ^\scA\restsec{[q_-,q_+]} : \sG(\scA\restsec{[q_-,q_+]}) \to \sG(\scA)$ by setting
\begin{equation} \label{eq:endpoint_inclusion}
\sJ^\scA\restsec{[q_-,q_+]}(x,b) = (x,b + q_-(y))
\end{equation}
Functoriality of this assignment will be proven in the next claim. We note that (on objects) we have
\begin{equation}\label{eq:endpoint_inclusion_image}
\im(\rest {\sJ^\scA\restsec{[q_-,q_+]}} x) = [\secp{q}_-(x),\secp {q}_+(x)]
\end{equation}
as claimed in the beginning of the section.
\end{constr}

\begin{eg}[Family embeddings] \label{eg:family_embeddings} \hfill
\begin{enumerate} 
\item In \autoref{rmk:sections_determine_subbund} we have seen the preceding construction in a special case, namely setting $q^b_- = \msrc_{\tusU 2_{\scA_b}}$ and $q^b_+ = q^b_2$ we have
\begin{equation}
(\tusU 2_{\scA_b})\restsec{[q^b_-,q^b_+]} = \tusU 2_{\scB_b}
\end{equation}
and the family embedding for these endpoint equals
\begin{equation}
\sJ^{\tusU 2_{\scA_b}}\restsec{[q_-,q_+]} = \theta^3_b
\end{equation}

\item As a simpler example, consider the open sections $r^b_-, r^b_+$ in $\und {\scA_b}$ (from \autoref{defn:sections}) defined by
\begin{restoretext}
\begingroup\sbox0{\includegraphics{test/page1.png}}\includegraphics[clip,trim=0 {.25\ht0} 0 {.2\ht0} ,width=\textwidth]{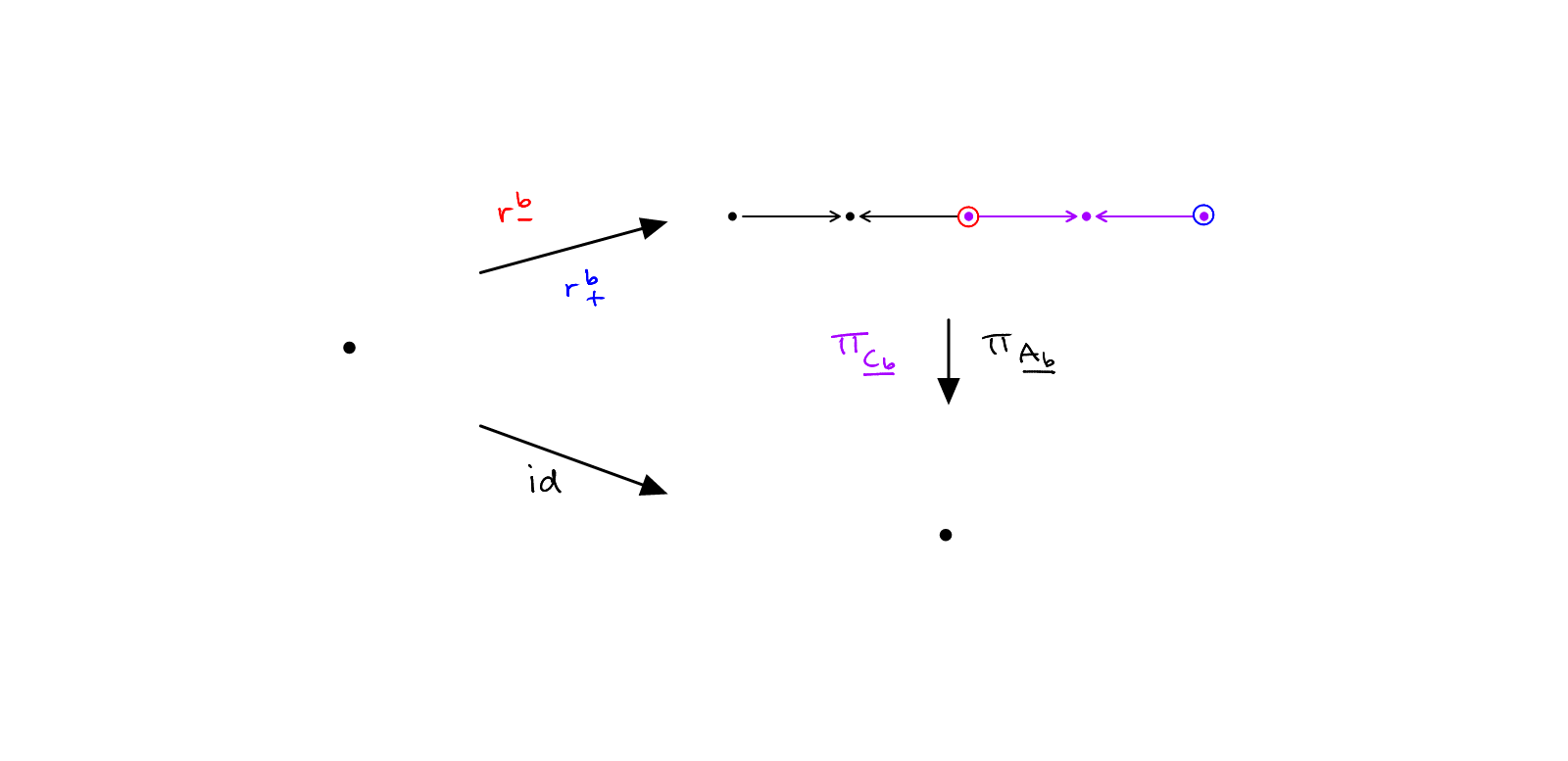}
\endgroup\end{restoretext}
Then we have
\begin{equation}
(\und {\scA_b})\restsec{[r^b_-,r^b_+]} = \und{\scC_b}
\end{equation}
namely, $\und{\scC_b}(0) = \singint 1$ (the interval with one singular height, highlighted in \cpurple{} above). And further, the family embedding for these endpoint equals
\begin{equation}
\sJ^{\und {\scA_b}}\restsec{[r^b_-,r^b_+]} = \phi^1_b
\end{equation}
which includes $\singint 1$ as the \cpurple{} subposet into $\singint 2$ in the illustration above.
\end{enumerate}
\end{eg}

\begin{claim}[Endpoint embedding is functorial] \label{claim:endpoint_embedding_functorial}
$\sJ^\scA\restsec{[q_-,q_+]}$ is functorial.
\proof The proof is \stfwd{}. We need to check
\begin{equation} \label{eq:pullback_injection_ff}
(b,c) \in \edgeset((\scA\restsec{[q_-,q_+]})(y \to z)) \iff (b+q_-(y), c + q_-(z)) \in \edgeset(\scA(y \to z)) 
\end{equation}
We do so by a case distinction
\begin{enumerate}

\item If $b$ is a singular height (if and only if $b+q_-(y)$ is a singular height) then by \eqref{eq:defn_order_realisation_1} the left hand side of   \eqref{eq:pullback_injection_ff} holds iff $(\scA\restsec{[q_-,q_+]})(y \to z)(b) = c$ which holds iff $\scA(y \to z)(b + q_-(y)) - q_-(z) = c$. The latter again by \eqref{eq:defn_order_realisation_1} is equivalent to the right hand side of \eqref{eq:pullback_injection_ff}.

\item If $b$ is a regular segment (if and only if $b+q_-(y)$ is a regular segment), then by \eqref{eq:defn_order_realisation_2} the left hand side of \eqref{eq:pullback_injection_ff} holds iff
\begin{equation} \label{eq:even_edge_condition_pullback_0}
\wwidehat{(\scA\restsec{[q_-,q_+]})}(y \to z)(b - 1) \leq c \leq \wwidehat{(\scA\restsec{[q_-,q_+]})}(y \to z)(b + 1)
\end{equation}
which we claim (and prove below) to hold iff
\begin{equation} \label{eq:even_edge_condition_pullback}
\wwidehat{\scA}(y \to z)(b - 1 + q_-(y)) - q_-(z) \leq c \leq \wwidehat{\scA}(y \to z)(b + 1+ q_-(y)) - q_-(z)
\end{equation}
Then, \eqref{eq:even_edge_condition_pullback} again by \eqref{eq:defn_order_realisation_2} is equivalent to the right hand side of  \eqref{eq:pullback_injection_ff} proving the statement.

To see the mutual implication \eqref{eq:even_edge_condition_pullback_0} $\iff$ \eqref{eq:even_edge_condition_pullback}, first note that if $0 < b < q_+(y) - q_-(y)$ then the claim follows from unpacking the definition of $\scA\restsec{[q_-,q_+]}$. In the boundary cases, if $b = 0$ then we have $\wwidehat{(\scA\restsec{[q_-,q_+]})}(y \to z)(b - 1)  = -1 < c$ by \autoref{notn:singular_morphism_boundary_cases} as well as $\wwidehat{\scA}(y \to z)(b - 1 + q_-(y)) - q_-(z) \leq c$ by the following argument: for any singular height $d \in \extsing(\scA(y))$ with $d < q_-(y)$ we must have $\wwidehat{\scA}(y \to z)(d) \leq q_-(z)$, since either $(d,\wwidehat{\scA}(y \to z)(d)) = (-1,-1)$ or $(d,\wwidehat{\scA}(y \to z)(d)) <\oedgeset (q_-(y),q_-(z)) \in \edgeset(\scA(y \to z))$. By a similarly argument, if $b = q_+(y) - q_-(y)$ we find $c < q_+(y) - q_-(y) + 1 = \wwidehat{(\scA\restsec{[q_-,q_+]})}(y \to z)(b + 1)$ as well as $c \leq  \wwidehat{\scA}(y \to z)(b + 1 + q_-(y)) - q_-(z)$.
\end{enumerate}
\end{claim}

\subsection{$k$-level embeddings}

In this section we extend the construction of the previous section from $\SI$-families to $\SIvert n \cC$-families. Explicitly, we will construct cube embeddings denoted by $\sJ^{\scA,k}\restsec{[q_-,q_+]}$ from endpoints $(q_-, q_+)$ ``at level $k$", and give examples thereof.

\begin{constr}[Constructing cube family embeddings by specifying endpoints] \label{constr:subfamilies_from_endpoints}
Now let $\scA : X \to \SIvert n \cC$ and $0 \leq k < n$, and let $(q_+, q_-)$ be a continuous choice of endpoints for $\scB := \tusU k_\scA : \tsG k(\scA) \to \SI$. Using \autoref{constr:endpoint_inclusions}, we define a family $\scA^k\restsec{[q_-,q_+]} : X \to \SIvert n \cC$ called the (labelled) \textit{subfamily of $\scA$ determined by $k$-level endpoints $q_-$, $q_+$} by setting
\begin{equation}
T = \Set{\pi_{\scB\restsec{[q_-,q_+]}}, \tpi k_\scA, \dots , \tpi 1_\scA}
\end{equation}
and
\begin{equation}
\scA^k\restsec{[q_-,q_+]} := \tsR k_{T,\tsU {k+1}_\scA \sJ^{\scB}\restsec{[q_-,q_+]}}
\end{equation}
This can be made into an embedding $\sJ^{\scA,k}\restsec{[q_-,q_+]} : \scA^k\restsec{[q_-,q_+]} \mono \scA$ by setting (cf. \autoref{notn:k_lvl_basechange})
\begin{equation}
(\sJ^{\scA,k}\restsec{[q_-,q_+])})^l = \tsG {l-k-1}(\sJ^{\scB}\restsec{[q_-,q_+]}) 
\end{equation}
Using \eqref{eq:endpoint_inclusion} and \autoref{constr:unpacking_collapse} this can be seen to give an embedding as required. It is called the \textit{$k$-level embedding into $\scA$ determined by endpoints $q_-$, $q_+$}. 
\end{constr}
Note that a $k$-level embedding is in particular a $k$-level base change (in the sense that was discussed in \autoref{ssec:multi_bc}).

We already saw examples of subfamilies of $\scA$ with $k$-level endpoints and we recall them in the following

\begin{eg}[Subfamilies from endpoints] \hfill
\begin{enumerate}
\item Using definitions in \autoref{eg:family_embeddings} and \autoref{eg:sections}, setting $\scB := \tusU 2_{\scA_b}$ we have
\begin{align}
\theta^3_b &= \sJ^{\scB}\restsec{[q^b_-,q^b_+]} \\
\theta^2_b &= \id \\
\theta^1_b &= \id \\
\theta^0_b &= \id
\end{align}
and thus
\begin{equation}
\big(~\theta_b : \scB_b \mono \scA_b~\big) \quad = \quad \big(~\sJ^{\scA_b,3}\restsec{[q^b_-,q^b_+]} : (\scA_b)^k\restsec{[q^b_-,q^b_+]} \mono \scA_b ~\big)
\end{equation}
\item Using definitions in \autoref{eg:family_embeddings} and \autoref{eg:sections}, setting $\scB := \und{\scA_b}$ we have
\begin{align}
\phi^3_b &=  \tsG 2(\sJ^{\scB}\restsec{[r^b_-,r^b_+]})\\
\phi^2_b &= \sG(\sJ^{\scB}\restsec{[r^b_-,r^b_+]}) \\
\phi^1_b &= \sJ^{\scB}\restsec{[r^b_-,r^b_+]} \\
\phi^0_b &= \id
\end{align}
and thus
\begin{equation}
\big(~\phi_b : \scC_b \mono \scA_b~\big) \quad = \quad \big(~\sJ^{\scA_b,1}\restsec{[r^b_-,r^b_+]} : (\scA_b)^k\restsec{[r^b_-,r^b_+]} \mono \scA_b ~\big)
\end{equation}
\end{enumerate}
\end{eg}

\begin{rmk}[Source-target endpoint inclusion] Given $\scA : X \to \SIvertone  \cC$ $0 \leq k < n$, and setting $\scB := \tusU k_\scA$, the previous \autoref{constr:source_and_target_inclusion} provides us with the \textit{trivial $k$-level endpoints} $(\msrc_{\scB},\mtgt_{\scB})$ for $\scB$. Using 
\autoref{constr:subfamilies_from_endpoints} we find that
\begin{equation}
\scA^k\restsec{[\msrc_{\scB},\mtgt_{\scB}]} = \scA
\end{equation}
This can be deduced as a corollary to \autoref{claim:factorising_subbund_through_endpoints}, which we will discuss in the next section.
\end{rmk}

\subsection{Decomposing (multi-level) embeddings into $k$-level embeddings}

We will now show that any embedding can be decomposed as the composition of embeddings determined by endpoints. We first show how individual squares of the form \eqref{eq:subbund_commute_2} decompose. 

\begin{eg}[Factoring embedding squares by endpoint family embeddings] \label{eg:factoring_subfamilies} For instance, using \autoref{eg:sections} consider the square 
\begin{equation}
\xymatrix{ \tsG {2}(\scB_a) \ar[r]^{\theta^{2}_a} \ar[d]_{\tpi 2_{\scB_a}} & \tsG {2}(\scA_a) \ar[d]^{\tpi 2_{\scA_a}} \\
\tsG {1}(\scB_a) \ar[r]_{\theta^1_a} & \tsG {1}(\scA_a) }
\end{equation}
This can be decomposed into a pullback and a family embedding functor as follows
\begin{restoretext}
\begingroup\sbox0{\includegraphics{test/page1.png}}\includegraphics[clip,trim=0 {.15\ht0} 0 {.1\ht0} ,width=\textwidth]{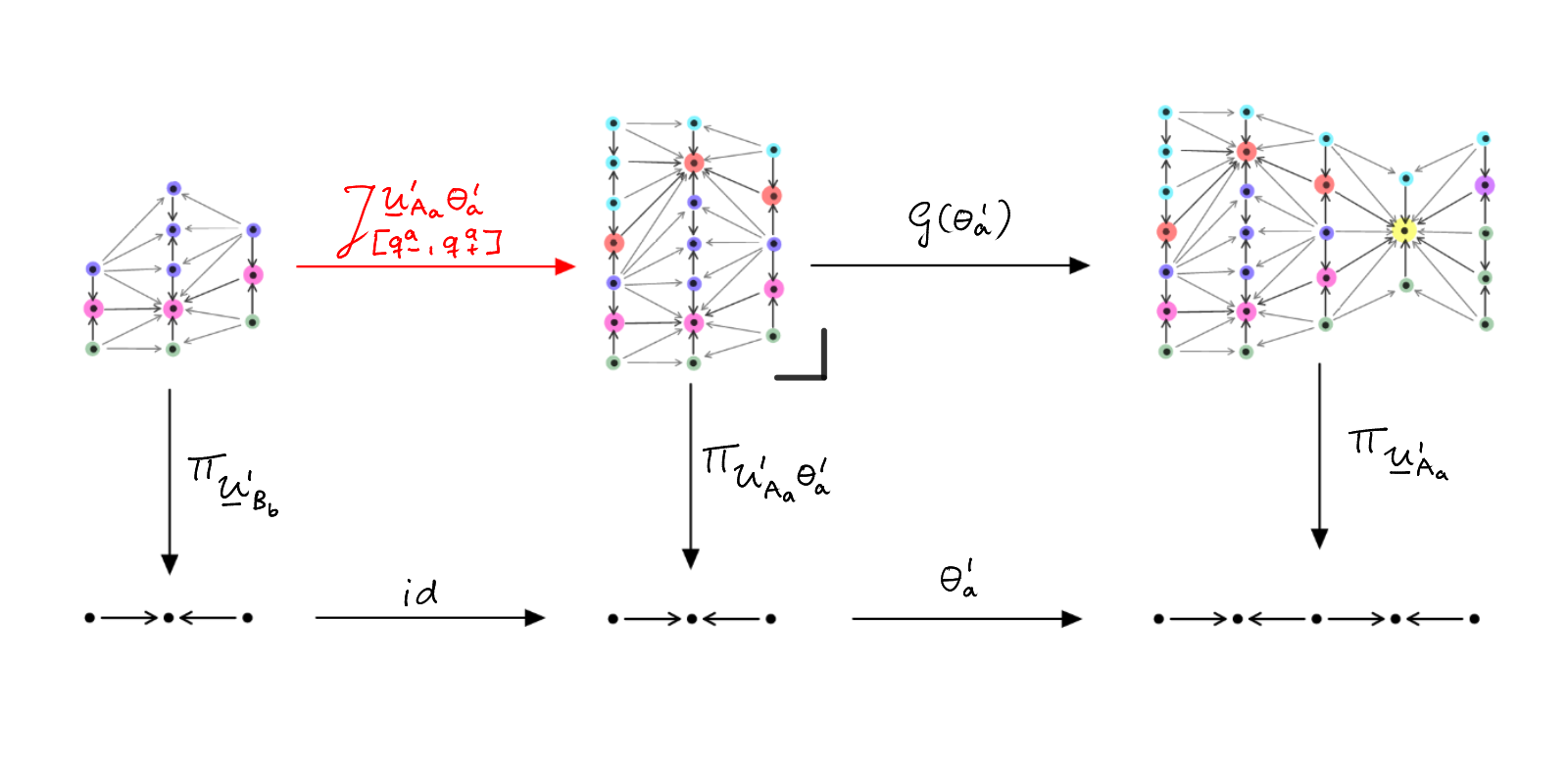}
\endgroup\end{restoretext}
where $q^a_-$ and $q^a_+$ are defined as
\begin{restoretext}
\begingroup\sbox0{\includegraphics{test/page1.png}}\includegraphics[clip,trim=0 {.2\ht0} 0 {.15\ht0} ,width=\textwidth]{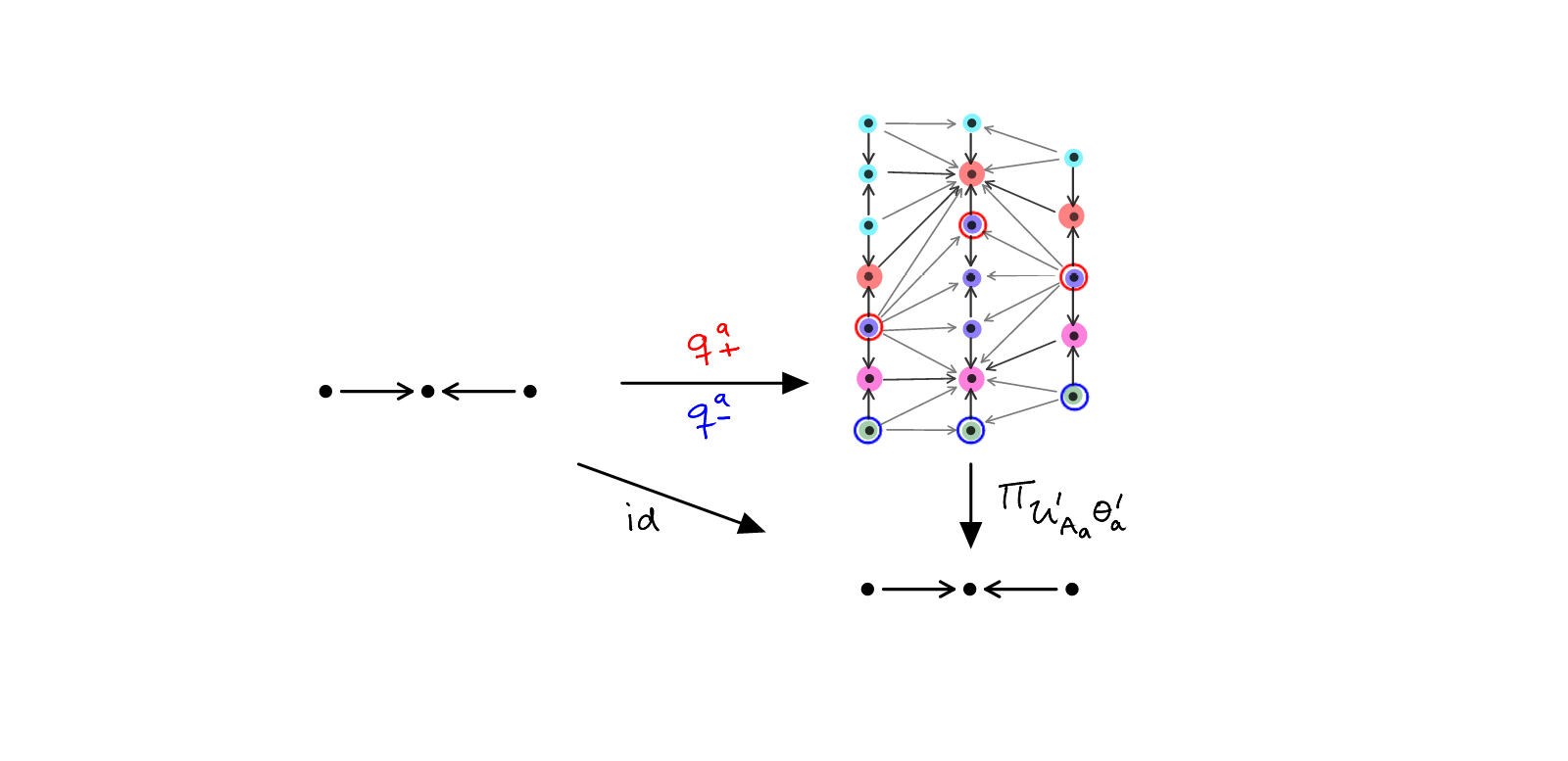}
\endgroup\end{restoretext}
\end{eg} 

This observation is generalised in the following claim.

\begin{claim}[Factoring embedding squares by endpoint family embeddings] \label{claim:factorising_subbund_through_endpoints} Let $\theta : \scB \mono \scA$ for $\scA : X \to \SIvert n \cC$, $\scB: Y \to \SIvert n \cC$. The square
\begin{equation} \label{eq:factoring_subbun_through_rest2}
\xymatrix{ \tsG {k+1}(\scB) \ar[r]^{\theta^{k+1}} \ar[d]_{\tpi k_\scB} & \tsG {k+1}(\scA) \ar[d]^{\tpi k_\scA} \\
\tsG {k}(\scB) \ar[r]_{\theta^k} & \tsG {k}(\scA) }
\end{equation}
is the composite of the squares
\begin{equation} \label{eq:factoring_subbun_through_rest}
\xymatrix@C=4cm{ \tsG {k+1}(\scB) \ar[r]^-{\sJ^{\tusU k_\scA \theta^k}_{\left[u^{k+1}_\theta\msrc_{\tusU k_\scB},u^{k+1}\mtgt_{\tusU k_\scB}\right]}} \ar[d]_{\tpi k_\scB} & \sG(\tusU k_\scA \theta^k) \ar[r]^{\sG(\theta^k)} \ar[d]_{\pi_{\tusU k_\scA \theta^k}} \pullbackfar & \tsG {k+1} (\scA) \ar[d]^{\tpi k_\scA} \\
\tsG {k}(\scB) \ar[r]_{\id} & \tsG {k}(\scB) \ar[r]_{\theta^k} & \tsG {k}(\scA) }
\end{equation}
where $u^{k+1}_\theta$ is the factorisation of \eqref{eq:factoring_subbun_through_rest2} through the right pullback square in \eqref{eq:factoring_subbun_through_rest}, that is
\begin{equation} \label{eq:factoring_subbun_through_rest_u}
\xymatrix{ & \tsG {k+1}(\scB) \ar[dd]|(0.5){\txt{\color{white}HACK}}^(0.7){\tpi k_\scB} \ar[dl]_{u^{k+1}_\theta} \ar[dr]^{\theta^{k+1}} & \\
\sG(\tusU k_\scA \theta^k) \pullbackfar \ar[rr]^(0.4){\sG(\theta^k)} \ar[dd]_{\pi_{\tusU k_\scA \theta^k }} && \tsG {k+1}(\scA) \ar[dd]^{\tpi k_\scA} \\
& \tsG k(\scB) \ar[dl]_{\id} \ar[dr]^{\theta^k} & \\
\tsG k(\scB)\ar[rr]^{\theta^k} && \tsG k(\scA) }
\end{equation}
\proof 
The proof is \stfwd{}. Note that since $\sG(\theta^k)$ is a fibrewise identity, we have ($x \in \tsG {k}(\scB)$)
\begin{equation} \label{eq:factoring_subbun_through_rest3}
\rest {u^{k+1}_\theta} x = \rest {\theta^{k+1}} x
\end{equation}
Abbreviate $\scC := \tusU k_\scA$ and $\scD := \tusU k_\scB$. To see that the upper left arrow of \eqref{eq:factoring_subbun_through_rest} is well-typed, we first show that
\begin{equation}
 \scD = ( {\scC} \theta^k )\restsec{[u^{k+1}_\theta\msrc_{ \scD},u^{k+1}_\theta\mtgt_{ \scD}]}
\end{equation}
The right hand side was defined in \autoref{constr:endpoint_inclusions}. Since $\tsG 1(\scD) = \sG( \scD)$ this shows well-definedness of the the squares in \eqref{eq:factoring_subbun_through_rest}. 

On objects $y \in Y$, \autoref{constr:endpoint_inclusions} gives
\begin{align}
(\scC \theta^k )\restsec{[u^{k+1}_\theta\msrc_{ \scD},u^{k+1}_\theta\mtgt_{ \scD}]}(y) &=\singint {\frac{\secp{u^{k+1}_\theta}\mtgt_{ \scD}(y) - \secp{u^{k+1}_\theta}\msrc_{ \scD}(y)}{2}} \\
&= \singint {\frac{\secp{\mtgt_{ \scD}}(y) - \secp{\msrc_{ \scD}}(y)}{2}} \\
&= \scD(y)
\end{align}
where, in the second step we used linearity of $u^{k+1}_\theta$ (cf. \autoref{claim:linearity_subfamilies}) and in the third step we used the definition of $\msrc$ and $\mtgt$ in \autoref{constr:source_and_target_inclusion}.

Similarly, on morphisms $(y \to z)$ in $Y$ we find by \autoref{constr:subfamilies_from_endpoints}
\begin{align}
(a,b) \in &\edgeset\big( ( \scC \theta^k )\restsec{[u^{k+1}_\theta\msrc_{ \scD},u^{k+1}_\theta\mtgt_{ \scD}]} (y \to z) \big) \\
& \iff (a + \secp{u^{k+1}_\theta}\msrc_{ \scD}(y), b + \secp{u^{k+1}_\theta}\msrc_{ \scD}(z)) \in \edgeset\big(\rest { \scC} Y (y \to z) \big) \\
& \iff (\rest {\theta^{k+1}} y a, \rest {\theta^{k+1}} z a) \in \edgeset\big( \scC (y \to z) \big) \\
& \iff (a,b) \in \edgeset\big( \scD (y \to z) \big)
\end{align}
where, in the second step we used \autoref{claim:linearity_subfamilies} and \eqref{eq:factoring_subbun_through_rest3}, and in the last step used that $\theta^{k+1}$ is fully faithful by \autoref{rmk:subfamily_components_ff}.

Now, the equation 
\begin{equation}
\theta^{k+1} = \sG(\theta^k) \sJ^{\scC \theta^k}\restsec{[u^{k+1}_\theta\msrc_{ \scD},u^{k+1}_\theta\mtgt_{ \scD}]}
\end{equation}
follows fibrewise. Indeed, for $y \in Y$, $\rest {\sG(\theta^k)} y$ is just the identity, and $\rest {(\sJ^{\rest \scC Y}\restsec{[u^{k+1}_\theta\msrc_{ \scD},u^{k+1}_\theta\mtgt_{ \scD}]})} y$ can be seen to give the same mapping as $\rest {u^{k+1}_\theta} y$ by using \eqref{eq:factoring_subbun_through_rest3} and comparing \eqref{eq:endpoint_inclusion} with \autoref{claim:linearity_subfamilies}. \qed
\end{claim}

The preceding claim now allows us to decompose any embedding as follows.

\begin{constr}[Decomposing embeddings into embeddings determined by endpoints] \label{constr:subfamilies_by_specifying_endpoints} We show that all embeddings $\theta : \scB \mono \scA$ can be inductively (and canonically) decomposed as a composite of embeddings with endpoints (starting with an initial restriction embedding). This is based on \autoref{eg:subfamily_by_restriction}, \autoref{claim:linearity_subfamilies} and \autoref{constr:subfamilies_from_endpoints}. 

Let $\theta : \scB \mono \scA$ for $\scA : X \to \SIvert n \cC$, $\scB: Y \to \SIvert n \cC$ and $Y \in X$. For $k = 0, 1, \dots , n$, we inductively construct $\theta_k : \scC^\theta_k \mono \scA$ such that $\theta \mono \theta_k$ and further, for $l \leq k$
\begin{equation} \label{eq:endpoints_decomposition}
\theta^l_k = \theta^l : \tsG l(\scB) \to \tsG l(\scA)
\end{equation}
and for $l \geq k$
\begin{equation} \label{eq:endpoints_decomposition2}
\theta^l_k = \tsG {l-k}(\theta^k) : \tsG l(\scB) \to \tsG l(\scA)
\end{equation}
\begin{itemize}
\item
For $k = 0$ we set 
\begin{equation} \label{eq:subfamily_decomp_base}
\theta_0 = \restemb_{\theta^0} :  \scA \theta^0 \mono \scA
\end{equation}
This satisfies \eqref{eq:endpoints_decomposition}, \eqref{eq:endpoints_decomposition2} and $\theta \mono \theta_0$ by definition of $\restemb_{\theta^0}$ in \autoref{eg:subfamily_by_restriction} (in particular $\theta \mono \theta_0$ follows factoring over the pullbacks that define $\theta_k^l$, $l > 0$).
\item
For $k > 0$, note that by inductive assumption we have $\tusU {k-1}_{\scC^\theta_{k-1}} = \tusU {k-1}_{\scA} \theta^{k-1}$. Borrowing notation from \autoref{claim:factorising_subbund_through_endpoints}, we define
\begin{equation} \label{eq:subfamily_decomp_ind}
\theta_k = \theta_{k-1} \sJ^{\tusU {k-1}_{\scA}\theta^{k-1},k}_{\left[u^{k}_\theta\msrc_{\tusU k_\scB},u^{k}_\theta\mtgt_{\tusU k_\scB}\right]}
\end{equation}
This satisfies \eqref{eq:endpoints_decomposition} by inductive assumption, \autoref{constr:subfamilies_from_endpoints} and \autoref{claim:factorising_subbund_through_endpoints}. It satisfies \eqref{eq:endpoints_decomposition2} by inductive assumption and \autoref{constr:subfamilies_from_endpoints}. It also satisfies $\theta \mono \theta_0$  by inductive assumption and \autoref{constr:subfamilies_from_endpoints} (in particular $\theta \mono \theta_0$ follows factoring over the pullbacks that define $\theta_k^l$, $l > k$).
\end{itemize}
Since $\theta_n = \theta$, by unwinding the inductive step \eqref{eq:subfamily_decomp_ind} (and base case \eqref{eq:subfamily_decomp_base}) we find a decomposition of the form
\begin{equation}
\theta =  \sJ^{\cdot~,~ n}\restsec{[\cdot~,~ \cdot]} \dots \sJ^{\cdot~,~ 2}\restsec{[\cdot~,~ \cdot]} \sJ^{\cdot~,~ 1}\restsec{[\cdot~,~ \cdot]} \restemb_{\theta^0}
\end{equation}
(where dots represent the expressions derived in \eqref{eq:subfamily_decomp_ind}).
\end{constr}

\begin{rmk}[Analogy with levelwise decomposition of collapse] The decomposition of $\theta$ into ``levelwise embeddings" determined by $k$-level endpoints, should be regarded analogous to the levelwise decomposition of multi-level collapse in \autoref{constr:multilevel_collapse}. Indeed both decomposition are instances of decompositions into $k$-level basechanges as discussed in \autoref{ssec:multi_bc}.
\end{rmk}

\begin{eg}[Decomposing embeddings into embeddings determined from endpoints] Recall $\theta_a : \scB_a \mono \scA_a$ from \autoref{eg:sections}. The preceding construction gives a decomposition of $\theta_a$ of the form
\begin{equation}
\theta_a = \sJ^{\cdot~,~ 2}\restsec{[\cdot~,~ \cdot]} \sJ^{\cdot~,~ 1}\restsec{[\cdot~,~ \cdot]} \restemb_{\theta_0}
\end{equation}
In case of $\theta_a$, $\theta^0_a = \id : \bnum{1} \to \bnum{1}$ and thus $\restemb_{\theta_0} = \id$. The remaining two embeddings of the decomposition, can be visualised as follows
\begin{restoretext}
\begin{noverticalspace}
\begingroup\sbox0{\includegraphics{test/page1.png}}\includegraphics[clip,trim=0 {.0\ht0} 0 {.55\ht0} ,width=\textwidth]{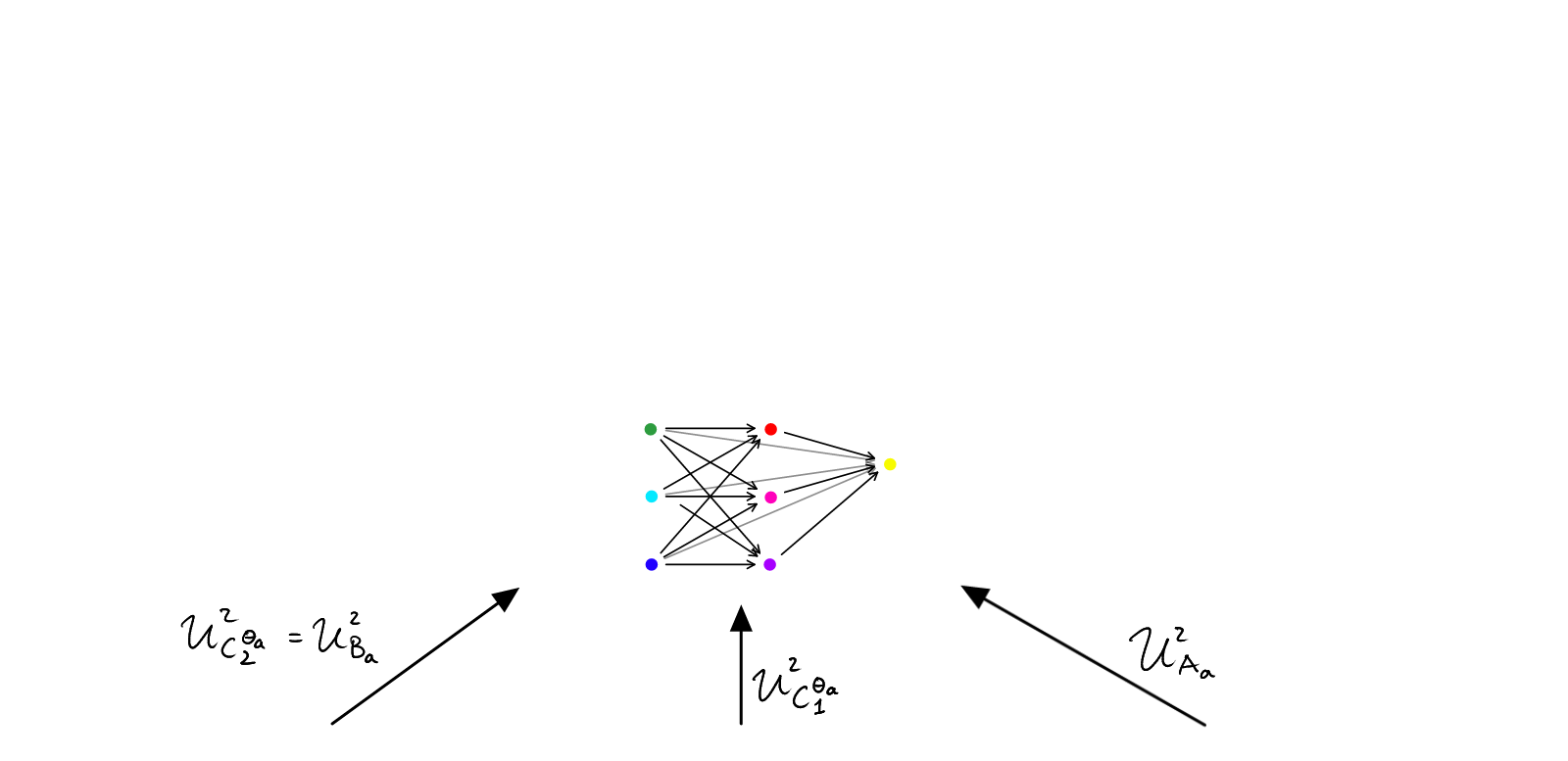}
\endgroup \\*
\begingroup\sbox0{\includegraphics{test/page1.png}}\includegraphics[clip,trim=0 {.0\ht0} 0 {.0\ht0} ,width=\textwidth]{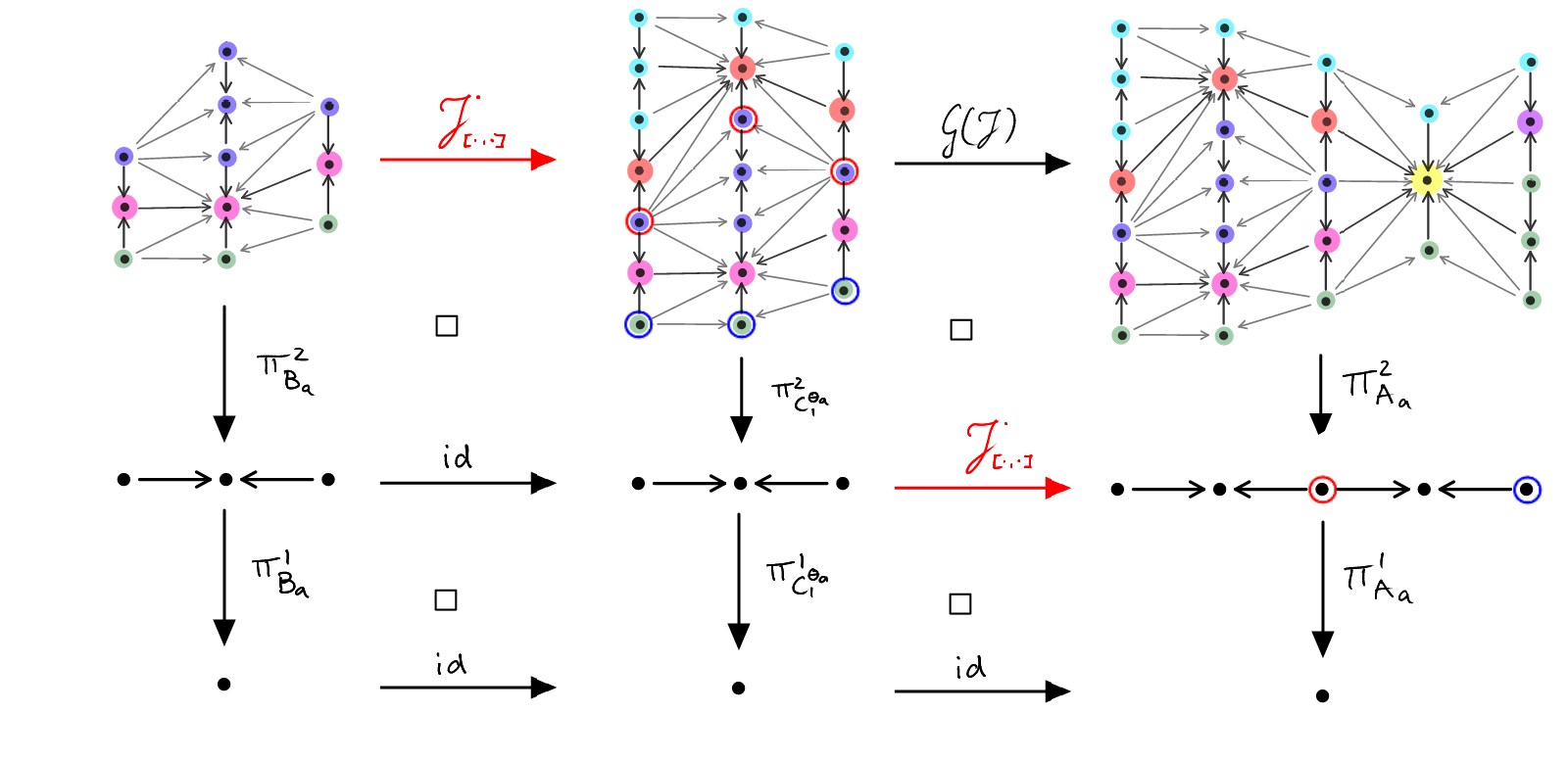}
\endgroup
\end{noverticalspace}
\end{restoretext}
Where in each of the two embeddings we marked the endpoints for the family embeddings (marked by \cred{} arrows) by \cred{} and \cblue{} circles.
\end{eg}

More examples of such a decompositions will be given in the next section. 

\section{Minimal embeddings} \label{sec:emb_min}

Using techniques from the previous section we will now construct a specific embedding, namely the minimal one ``around a point".

\subsection{Minimal endpoints}

We start with a discussion of certain endpoint sections that will be involved in the construction of minimal embedding. These will be called minimal endpoints. First, recall the following

\begin{notn}[Overcategories] \label{notn:overcats} \hfill
\begin{enumerate}
\item Given a category $\cC$, and an object $c \in \cC$, denote by $\cC \sslash c$ the overcategory of $\cC$ over $c$. That is, the objects $f$ of $\cC \sslash c$ are morphisms $(f : d \to c)$ in $\cC$, and morphisms $g : (f : d \to c) \to (f' : d' \to c)$ of $\cC \sslash c$ are morphisms $(g : d \to d')$ such that $f = f'g$. We have a forgetful functor $U : \cC \sslash c \to \cC$ that maps $(f : d \to c)$ to $d$, and $g : (f : d \to c) \to (f' : d' \to c)$ to $g : d \to d'$. 
\item Given a functor $F : \cC \to \cD$, and $c \in \cC$, we have an induced functor $F\sslash c : \cC \sslash c \to \cD \sslash Fc$ that maps  $(f : d \to c)$ to $(Ff : Fd \to Fc)$ and $g : (f : d \to c) \to (f' : d' \to c)$ to $Fg : (Ff : Fd \to Fc) \to (Ff' : Fd' \to Fc)$. 
\item Given a poset $X$, and $x \in X$, then $X \sslash x$ is a poset and the forgetful functor $U : X \sslash x \to X$ is fully-faithful and injective on objects. Thus we will identify $X \sslash x$ with the full subcategory of $X$ consisting of objects $y$ such that $y \to x$, i.e. we denote objects of $X \sslash x$ by $y$ instead of $y \to x$ (where $y \to x \in \mor(X)$).
\end{enumerate} 
\end{notn}

\noindent In particular, given a functor $\scA : X \to \SI$ and $p = (x,a) \in \sG(\scA)$, then from $\pi_{\scA} : \sG(\scA) \to X$ we obtain
\begin{equation}
\pi_{\scA} \sslash p : \sG(\scA) \sslash p \to X \sslash x
\end{equation}
which equals the restriction of $\pi_{\scA}$ to $\sG(\scA) \sslash p$ on its domain, and $X \sslash x$ on its codomain (cf. \autoref{notn:subsets_and_restrictions}).

\begin{constr}[Minimal endpoints] \label{constr:minimal_endpoints} Let $\scA : X \to \SI$, $p \in \sG(\scA)$. Assume $X \sslash p^0 = X$. We define endpoints $(\qmin p_+, \qmin p_-)$ for $\scA$ called minimal endpoints around $p$ by the following case distinction. 
\begin{enumerate}
\item If $\secp p \in \singcont(I)$ we set
\begin{equation}
\qmin p_\pm (y) = \big(y,\und\scA(y \to p^0)\regop (\secp p \pm 1)\big)
\end{equation}
This is functorial by \eqref{eq:defn_order_realisation_3}.
\item If $\secp p \in \regcont(I)$ we set
\begin{equation}
\qmin p_\pm (y) = \big(y,\und\scA(y \to p^0)\regop (\secp p)\big)
\end{equation}
This is functorial by \eqref{eq:defn_order_realisation_3}.
\end{enumerate}
\end{constr}

\begin{eg}[Minimal endpoints] \hfill
\begin{enumerate}
\item In the following we define an \SI-bundle $\pi_\scD : \sG(\scD) \to \singint 1$ and select a point $p \in \sG(\scD)$ which is marked in \cgreen{} below. Note that $\pi_\scD (p)$ indeed satisfies
\begin{equation}
(\singint 1 \sslash \pi_\scD (p)) = \singint 1
\end{equation}
Thus we can use the above construction to obtain $\qmin p_\pm$ as follows 
\begin{restoretext}
\begingroup\sbox0{\includegraphics{test/page1.png}}\includegraphics[clip,trim=0 {.0\ht0} 0 {.0\ht0} ,width=\textwidth]{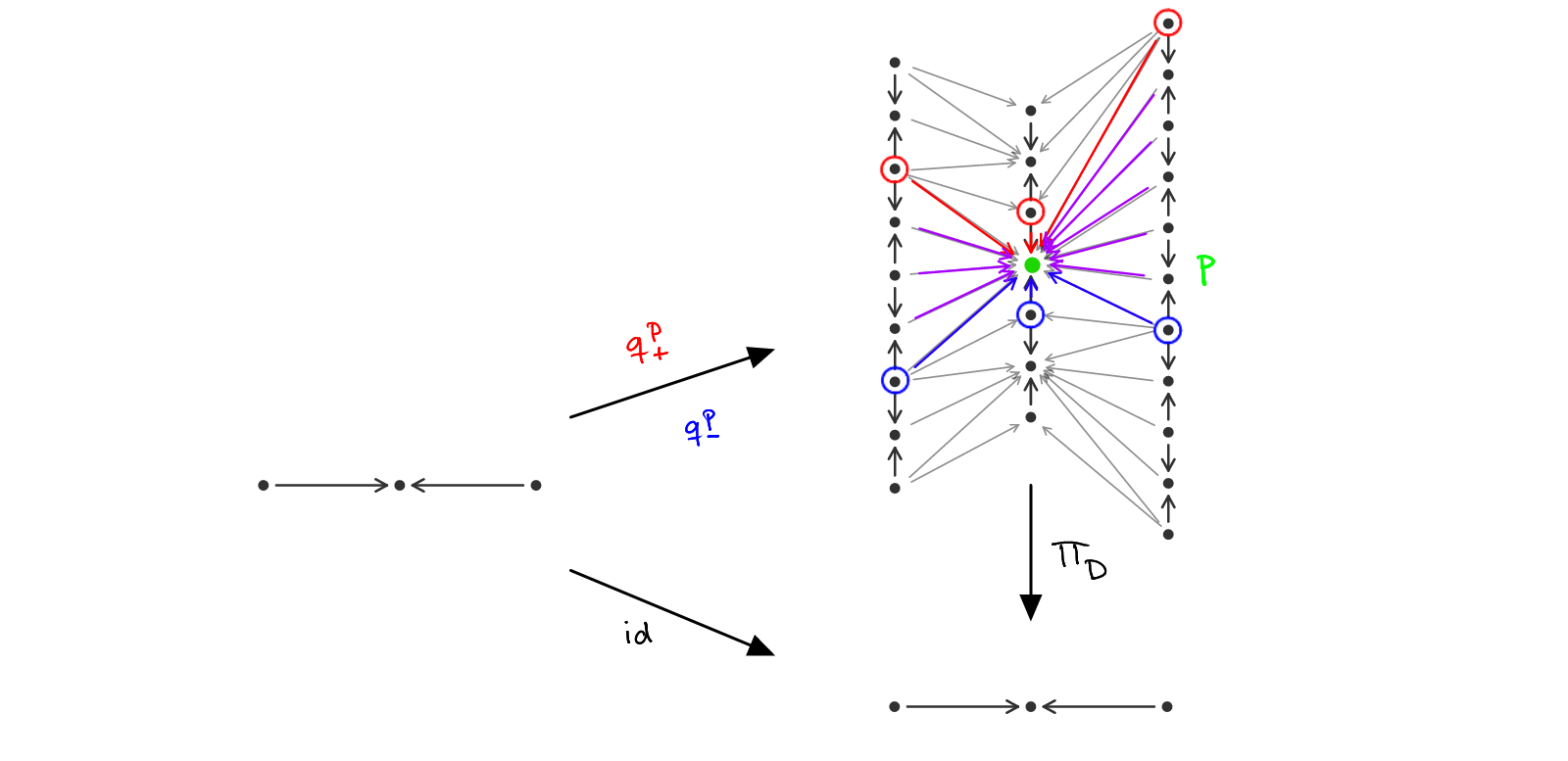}
\endgroup\end{restoretext}
Note that the endpoints $\qmin p_\pm(x)$ for $x \in \singint 1$ are the outer boundaries for regions $r \in \pi\inv_\scD(x)$ which have an arrow $(r \to p)$ to $p$ in $\sG(\scD)$ (the arrows $\qmin p_\pm(x) \to p$ are marked in \cred{} and \cblue{} respectively, the rest of such arrows is marked in \cpurple{} in the above). This observation is made precise in the next claim.
\item The following defines an \SI-bundle $\pi_{\widetilde \scD} : \sG(\widetilde \scD) \to X$ (where $X$ is $\sG(\const_{\singint 1} : \singint 1 \to \SI)$). We select a point $p \in \sG(\widetilde\scD)$ by marking it in \cgreen{}. Once more we have $(X \sslash \pi_{\widetilde \scD} (p)) = X$ and thus we can construct minimal endpoints $\qmin p_\pm$ which are given as follows
\begin{restoretext}
\begingroup\sbox0{\includegraphics{test/page1.png}}\includegraphics[clip,trim=0 {.0\ht0} 0 {.0\ht0} ,width=\textwidth]{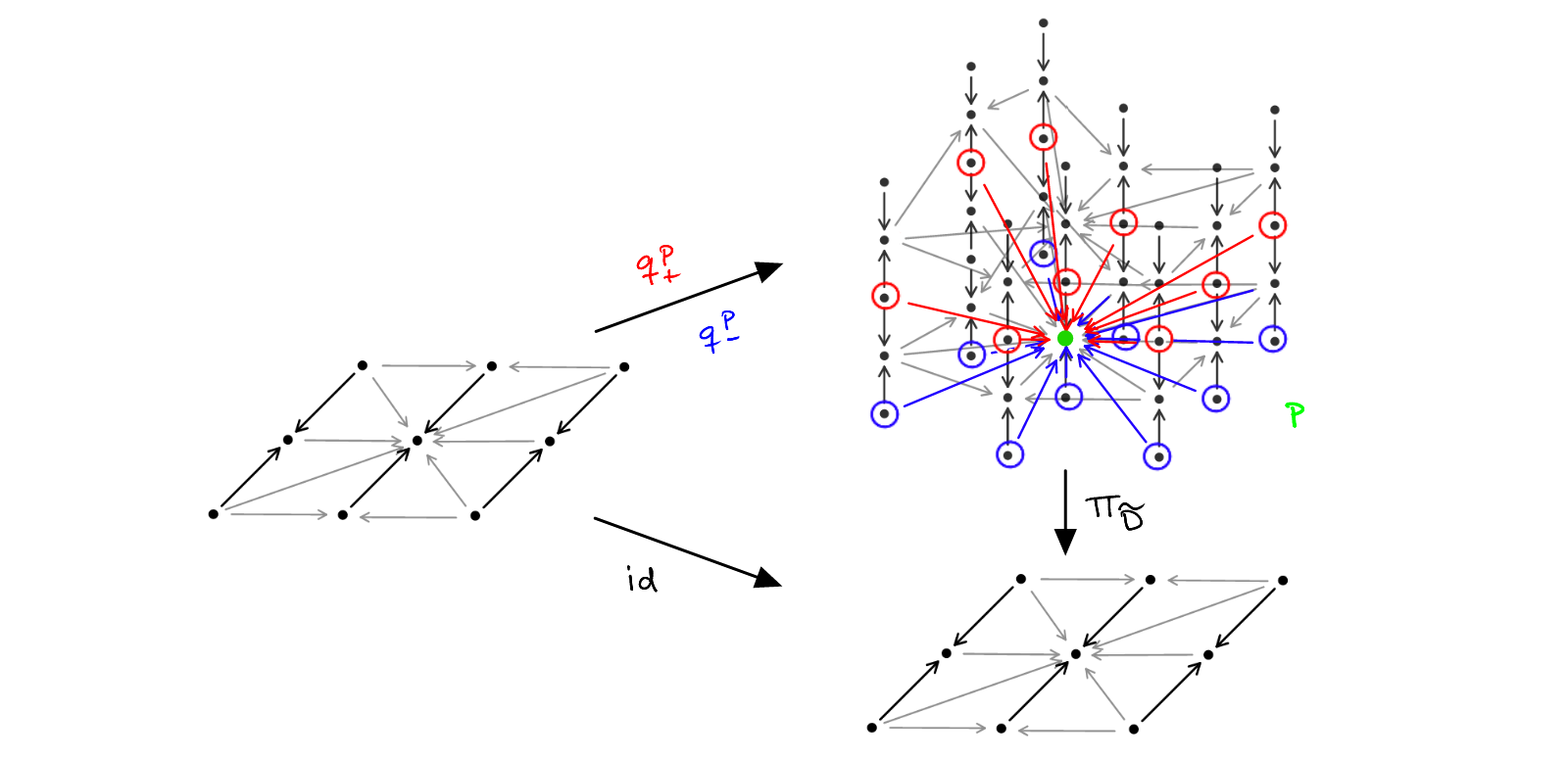}
\endgroup\end{restoretext}
\end{enumerate}
\end{eg}

\begin{claim}[Minimal endpoints delimit over-posets] \label{claim:minimal_endpoints_delimit_neighbourhoods} Let $\scA : X \to \SI$, $p \in \sG(\scA)$. Then $\sJ^{\scA}\restsec{[\qmin p_-,\qmin p_+]}$ restricts (on its codomain) to give an isomorphism
\begin{equation}
\sJ^{\scA}\restsec{[\qmin p_-,\qmin p_+]} : \sG(\scA\restsec{[\qmin p_-,\qmin p_+]}) ~ \iso ~ \sG(\scA) \sslash p
\end{equation}
\proof We now that $\sJ^{\scA}\restsec{[\qmin p_-,\qmin p_+]}$ is injective and has fibrewise image
\begin{equation}
\im\big(\rest {\sJ^{\scA}\restsec{[\qmin p_-,\qmin p_+]}} y\big) = E_y := [\secp {\qmin p_-},\secp {\qmin p_+}]
\end{equation}
it thus remains to show that 
\begin{equation}
(q \to p) \in \mor(\sG(\scA)) \quad \iff \quad \secp q \in E_{q^0}
\end{equation}
The implication from left to right follows from \eqref{eq:endpoint_boundaries}. The implication from right to left follows by a case distinction on $\secp p$. 
\begin{enumerate}
\item If $\secp p \in \singcont(\scA(p^0))$ then $(\secp q, \secp p) \in \edgeset(\scA(q^0 \to p^0))$ by \eqref{eq:defn_order_realisation_4} and \autoref{constr:minimal_endpoints}
\item If $\secp p \in \regcont(\scA(p^0))$ then $(\secp q, \secp p) \in \edgeset(\scA(q^0 \to p^0))$ by \eqref{eq:defn_order_realisation_3} and \autoref{constr:minimal_endpoints} \qed
\end{enumerate}
\end{claim}

\subsection{Construction}

Recall from \autoref{defn:regions} that $p \in \tsG n(\scA)$ for $\scA : X \to \SIvert n \cC$ has projection $p^k \in \tsG k(\scA)$ for each $0 \leq k \leq n$. We can use our minimal endpoint construction for each of these projections, and let the endpoints determine embeddings. If we compose these embeddings (after possible starting with a restriction embedding) we will obtain the ``minimal subfamily around $p$". For instance, in the case of a point $p \in \sG(\scA_a)$, where the $\SIvert 2 \cC$-family $\scA_a$ was defined in \autoref{eg:subfamilies}, we obtain the following (both $p^1$ and $p^2$ are marked in \cgreen{} below)
\begin{restoretext}
\begin{noverticalspace}
\begingroup\sbox0{\includegraphics{test/page1.png}}\includegraphics[clip,trim=0 {.55\ht0} 0 {.0\ht0} ,width=\textwidth]{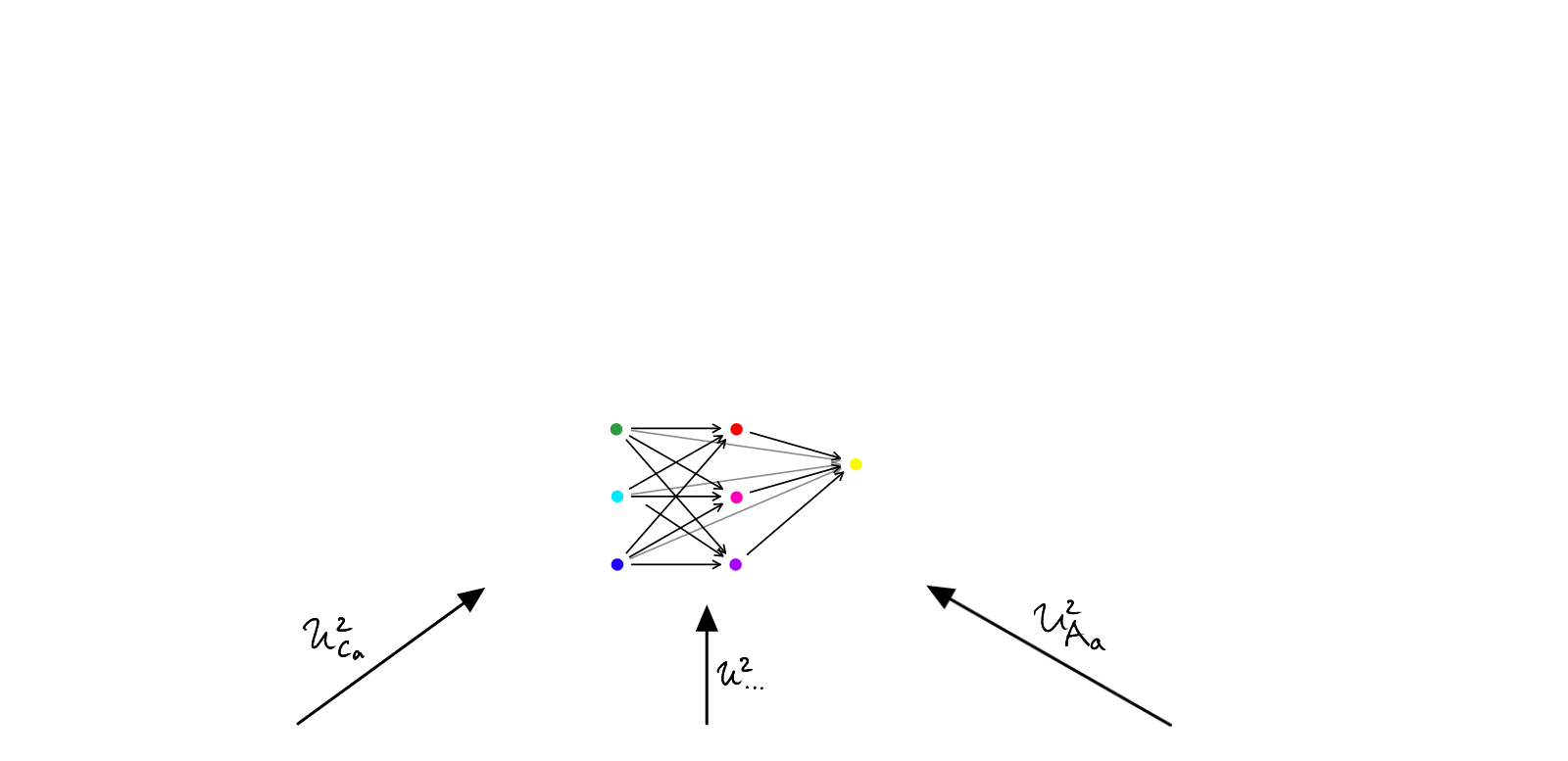}
\endgroup \\*
\begingroup\sbox0{\includegraphics{test/page1.png}}\includegraphics[clip,trim=0 {.0\ht0} 0 {.0\ht0} ,width=\textwidth]{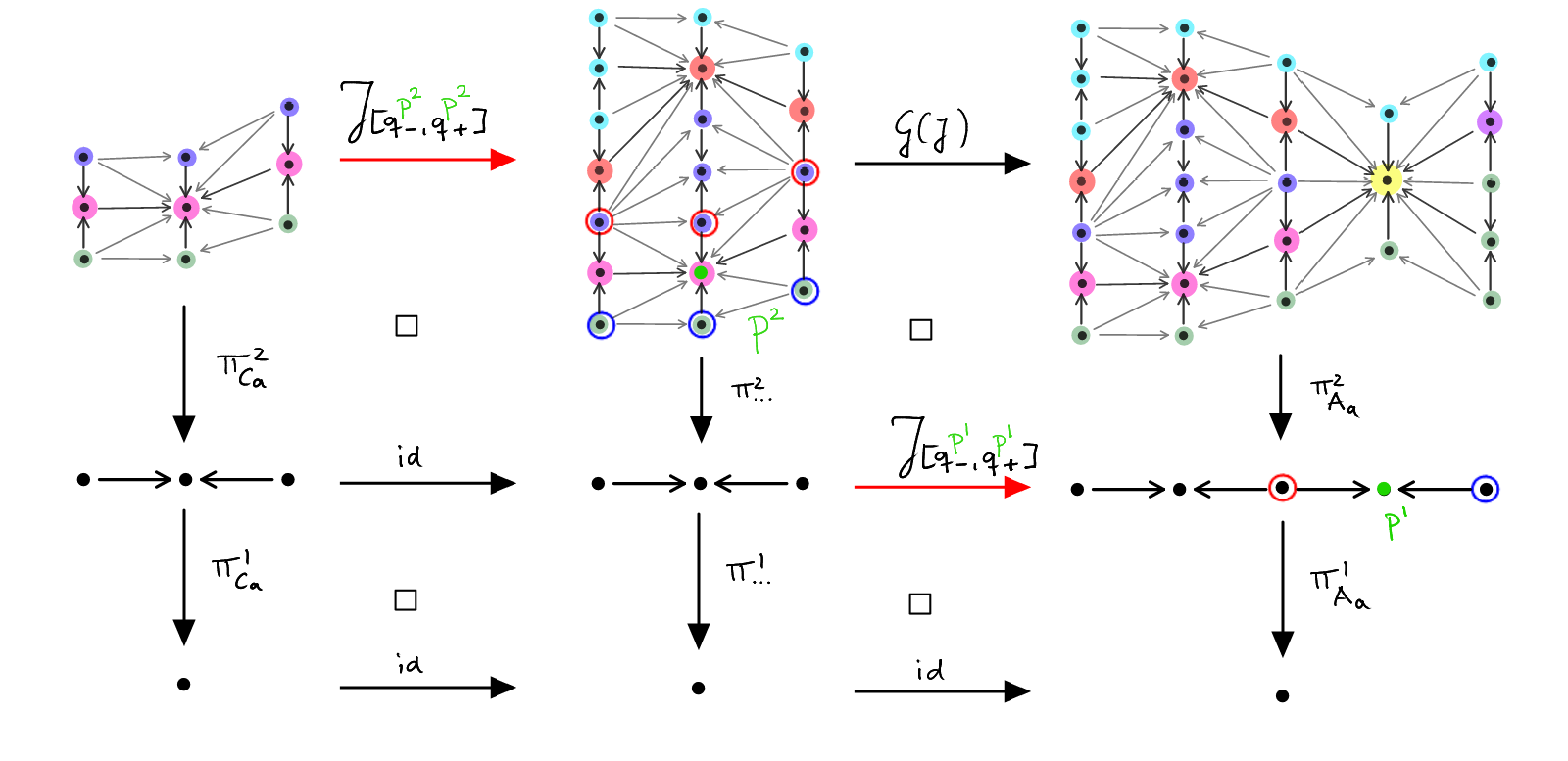}
\endgroup
\end{noverticalspace}
\end{restoretext}
Here, we named the $\SIvert 2 \cC$-family on the left $\scC_a$. Now $\phi_a : \scC_a \mono \scA_a$ is the embedding obtained from a composition of embeddings determined by endpoints (in fact, minimal endpoint). The resulting $\scC_a$ satisfies that $\tsG 2(\scC_a)$ is exactly the over-poset of $\tsG n(\scA_a)$ over $p$. This will be called the minimal embedding of $p$. We generalise the idea in the following construction.

\begin{constr}[Minimal embeddings] \label{constr:minimal_subfamilies} Let $\scA : X \to \SIvert n \cC$ and $p \in \tsG n(\scA)$. For $k = 0,1,2, ... n$, we inductively define embeddings $\iota^p_k : \scA^{p}_k \mono \scA$ such that there is $p_k \in \tsG n(\scA^{p}_k)$ with
\begin{equation} \label{eq:minimal_subfamily_1}
(\iota^p_k)^n(p_k) = p
\end{equation}
Further, for all $l \leq k$, $(\iota^p_k)^l$ will restrict (on its codomain) to an  isomorphism
\begin{equation} \label{eq:minimal_subfamily_2}
(\iota^p_k)^l = (\iota^p_l)^l \quad : \quad  \tsG l(\scA^{p}_k) =  \tsG l(\scA^{p}_l) ~ \xiso {} ~~\tsG l(\scA) \sslash p^l
\end{equation}
and for $l \geq k$, we have
\begin{equation} \label{eq:minimal_subfamily_3}
(\iota^p_k)^l = \tsG {l-k}((\iota^p_k)^k) \quad : \quad  \tsG l(\scA^{p}_k) ~\to ~\tsG l(\scA)
\end{equation}

\begin{itemize}
\item For $k = 0$ we set
\begin{equation}
\iota^p_0 = \restemb_{X \sslash p^0} : \scA(X \sslash p^0) \mono \scA
\end{equation}
Using the universal property of the pullbacks defining the components of $\restemb_{X \sslash p^0}$, we find a unique $p_0 \in \tsG n(\scA_0)$ satisfying the inductive conditions \eqref{eq:minimal_subfamily_1}. From \autoref{eg:subfamily_by_restriction} we deduce that $\iota^p_0$ also satisfies \eqref{eq:minimal_subfamily_2} and \eqref{eq:minimal_subfamily_3}.
\item For $k > 0$, the inductive assumptions allow us to define (using \autoref{constr:subfamilies_from_endpoints} and \autoref{constr:minimal_endpoints})
\begin{equation} \label{eq:minimal_subfamily_ind}
\iota^p_k = \iota^p_{k-1} \sJ^{\scA^{p}_{k-1},k}_{\left[\qmin {p^k_{k-1}}_-, \qmin {p^k_{k-1}}_+\right]} ~:~ \scA_k \mono \scA
\end{equation}
Using \autoref{claim:minimal_endpoints_delimit_neighbourhoods} and the universal property of the pullbacks defining the components (above level $k$) of $\sJ^{\scA^{p}_{k-1},k}_{\left[\qmin {p^k_{k-1}}_-, \qmin {p^k_{k-1}}_+\right]}$, we find a unique $p_k \in \tsG n(\scA_k)$ satisfying the inductive conditions \eqref{eq:minimal_subfamily_1}. Using \autoref{constr:subfamilies_from_endpoints} and the inductive assumptions, we find that  $\iota^p_k$ also satisfies conditions \eqref{eq:minimal_subfamily_2} and \eqref{eq:minimal_subfamily_3}.
\end{itemize}

\noindent We define
\begin{equation}
\iota^p_\scA := \iota^p_n : \scA^{p}_n \mono \scA
\end{equation}
which is called the \textit{minimal embedding} around $p$ in $\scA$, and denote
\begin{equation}
\scA \sslash p := \scA^{p}_n
\end{equation}
called the \textit{minimal subfamily} around $p$ in $\scA$.
\end{constr}

\begin{eg}[Minimal embeddings] We have already seen one example of the preceding construction. The following example follows through the construction of the minimal embedding for a point $p \in \tsG 3(\scA_b)$ (marked in \cgreen{} below, together with its projections) where $\scA_b$ was defined in \autoref{eg:subfamilies}
\begin{restoretext}
\begin{noverticalspace}
\begingroup\sbox0{\includegraphics{test/page1.png}}\includegraphics[clip,trim=0 {.0\ht0} 0 {.0\ht0} ,width=\textwidth]{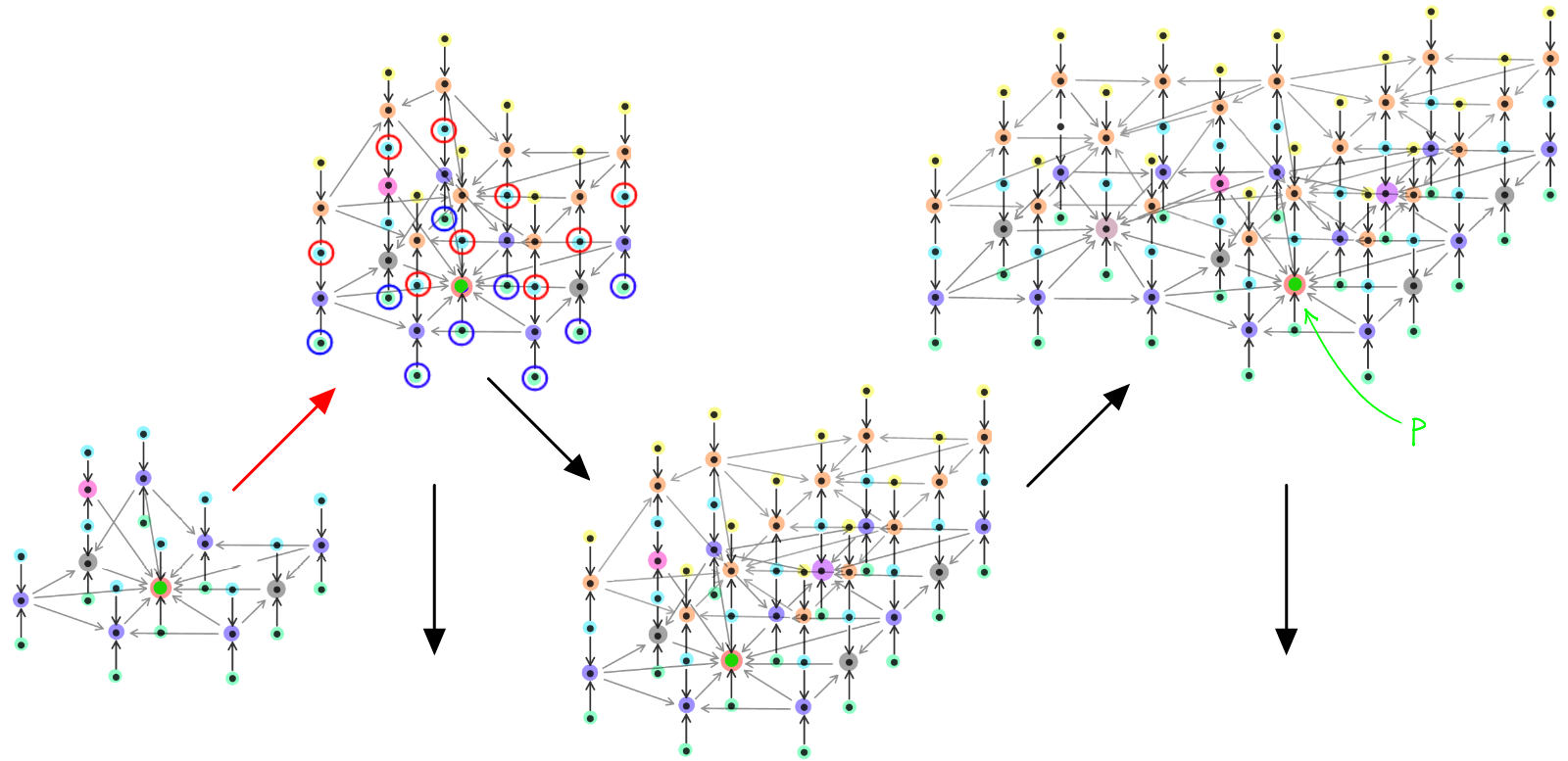}
\endgroup \\*
\begingroup\sbox0{\includegraphics{test/page1.png}}\includegraphics[clip,trim=0 {.0\ht0} 0 {.0\ht0} ,width=\textwidth]{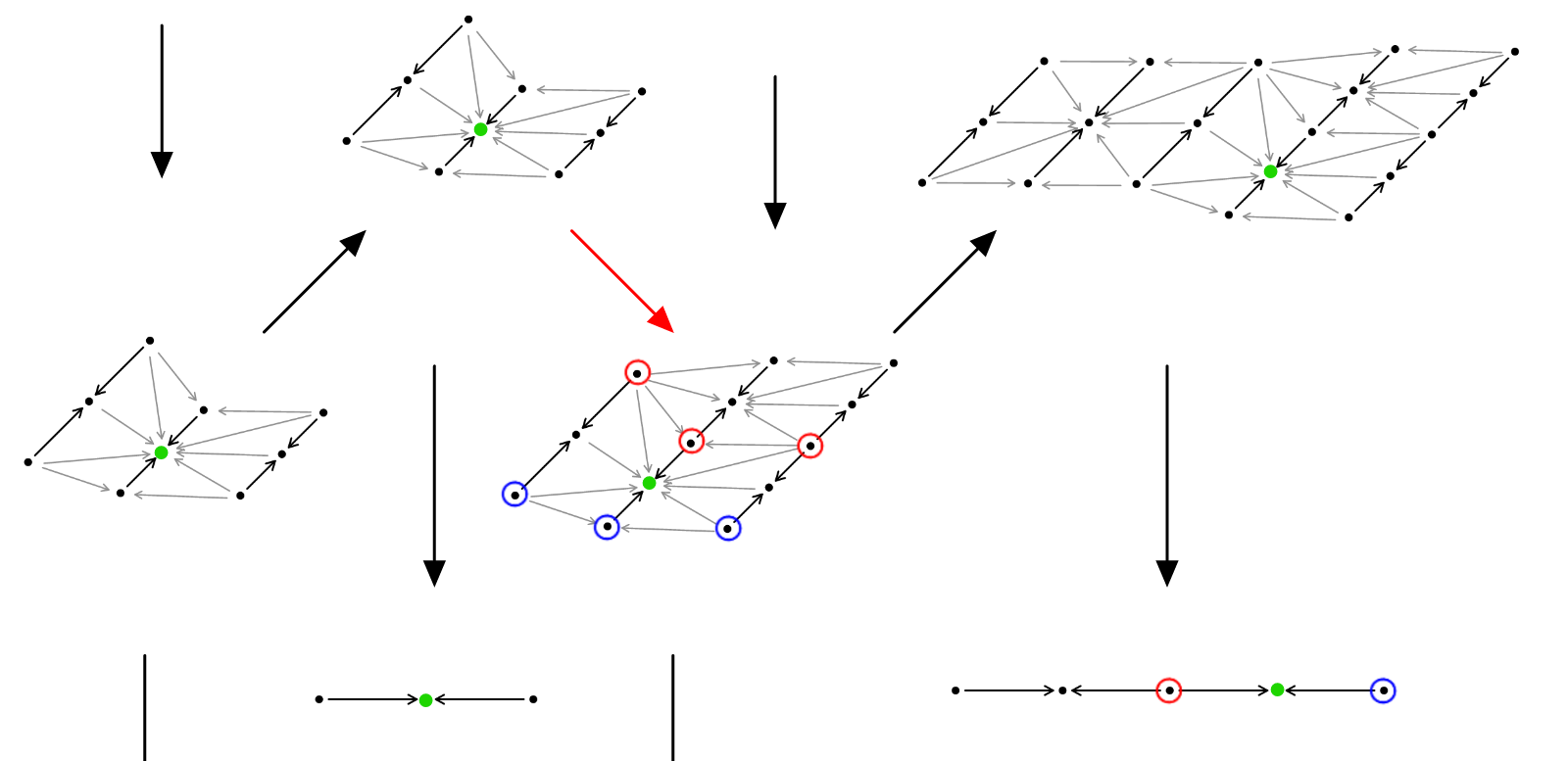}
\endgroup \\*
\begingroup\sbox0{\includegraphics{test/page1.png}}\includegraphics[clip,trim=0 {.25\ht0} 0 {.0\ht0} ,width=\textwidth]{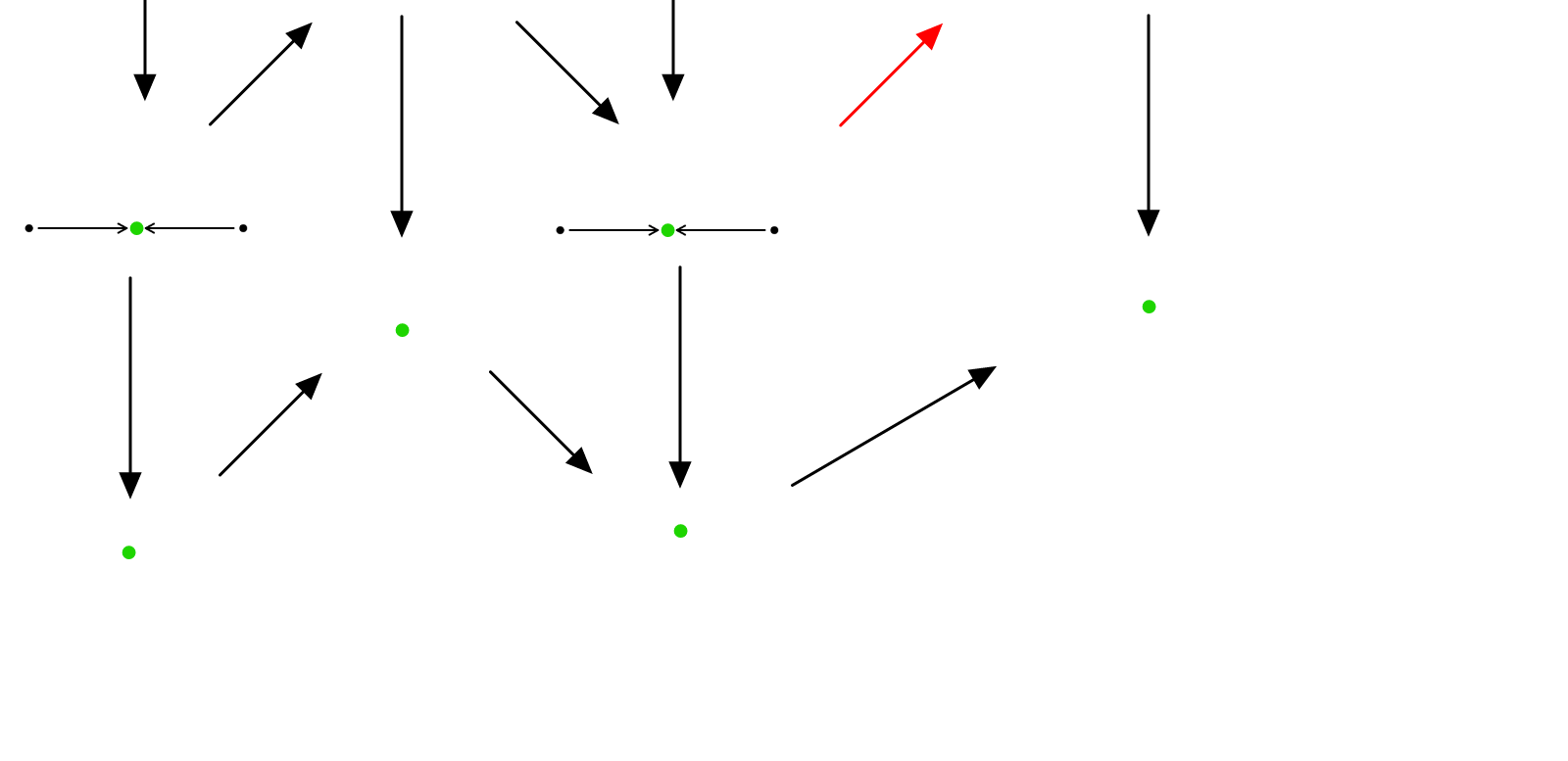}
\endgroup
\end{noverticalspace}
\end{restoretext}
Family embedding functors $\sJ^{~\cdot~}\restsec{[\cdot~,~ \cdot]}$ are highlighted by \cred{} arrows. Their (minimal) endpoints are given by \cred{} and \cblue{} circles.
\end{eg}

\subsection{Properties}

Connecting back to \autoref{sec:sum_norm}, we will now discuss the ``universal properties" of the minimal embedding.

\begin{thm}[Properties of minimal embeddings] \label{claim:minimal_subfamily_is_minimal} Let $\scA : X \to \SIvert n \cC$ be a family indexed by $X$, and $p \in \tsG n(\scA)$ a region in $\scA$.
\begin{enumerate}
\item For all $0 \leq k \leq n$ we have $\tsG k(\scA) \sslash p^k \subset \im(\iota^p_\scA)$
\item Let $\theta : \scB \subset \scA$ be an embedding of $\scA$. If $p \in \im(\theta^n)$ and $\im(\theta^0)$ is \gls{downwardclosed} in $X$ (in particular, $X \sslash p^0 \subset \im(\theta^0)$) then we have
\begin{equation}
\iota^p_\scA \mono \theta
\end{equation}
\end{enumerate}

\proof The proof is \stfwd{}. (i) follows from \autoref{constr:minimal_subfamilies} (in particular, c.f. \eqref{eq:minimal_subfamily_2}). We proceed to show (ii). Comparing \autoref{constr:minimal_subfamilies} and \autoref{constr:subfamilies_by_specifying_endpoints} it suffices to inductively show for $k = 0,1,2, ... n$ that
\begin{equation}
\iota^p_k  \mono \theta_k
\end{equation}
We have $\iota^p_0 \mono \theta_0$ by assumption that $\theta$ contains the neighbourhood of $p$ in the base $X$. Assume $\iota^p_{k-1}  \mono \theta_{k-1}$,  and set
\begin{equation}
\phi_{k-1} := (\iota^p_{k-1})\inv \theta_{k-1}
\end{equation}
The components $(\phi_{k-1})^l$ for $l \geq k$ are family pullbacks (since this holds for both $\iota^p_{k-1}$ and $\theta_{k-1}$) and thus restrict to isomorphisms of fibers. Comparing \eqref{eq:minimal_subfamily_ind} and \eqref{eq:subfamily_decomp_ind} it now suffices to show that
\begin{equation}
\phi_{k-1} \sJ^{\scA^{p}_{k-1},k}_{\left[\qmin {p^k_{k-1}}_-, \qmin {p^k_{k-1}}_+\right]} \quad \mono \quad  \sJ^{\tusU {k-1}_{\scC}\theta^{k-1},k}_{\left[u^{k}_\theta\msrc_{\tusU k_\scB},u^{k}_\theta\mtgt_{\tusU k_\scB}\right]}
\end{equation}
By \autoref{claim:factorisation_subfamilies} this means we need to check that images of the left hand side are contain in images of the right hand side. We can check this fibrewise. Namely, by \autoref{constr:subfamilies_from_endpoints} we need to check that for all $x \in \tsG {k-1}(\scA^p_{k-1})$ we have
\begin{equation} \label{eq:min_subbund_inc1}
\left[\secp {\phi^k_{k-1}} \qmin {p^k_{k-1}}_-(x), \secp {\phi^k_{k-1}} \qmin {p^k_{k-1}}_+(x)\right] \subset \left[\secp {u^{k}_\theta}\msrc_{\tusU k_\scB}(\phi^{k-1}_{k-1}(x)), \secp {u^{k}_\theta}\mtgt_{\tusU k_\scB}(\phi^{k-1}_{k-1}(x))\right]
\end{equation}
Note, that for $x = p^{k-1}_{k-1}$ we have 
\begin{equation} \label{eq:min_subbund_inc2}
\secp {u^{k}_\theta}\msrc_{\tusU k_\scB}(\phi^{k-1}_{k-1}(x)) \leq \secp {
\phi^k_{k-1}}(p^k_{k-1}) - 1 = \secp {\phi^k_{k-1}} \qmin {p^k_{k-1}}_-(x)
\end{equation}
where the last equation follows from \autoref{constr:minimal_endpoints} and linearity of $\phi^k_{k-1}$ (cf. \autoref{claim:linearity_subfamilies}). Using bimonotonicity and \eqref{eq:defn_order_realisation_3} we thus deduce that for any $(y \to p^{k-1}_{k-1}) \in \tsG {k-1}(\scA^p_{k-1})$, we have
\begin{equation} \label{eq:min_subbund_inc3}
\secp {u^{k}_\theta}\msrc_{\tusU k_\scB}(\phi^{k-1}_{k-1}(y)) \leq \secp {\phi^k_{k-1}} \qmin {p^k_{k-1}}_-(y)
\end{equation}
This shows the lower endpoint of the left interval in \eqref{eq:min_subbund_inc1} is greater or equal that the lower endpoint of the right interval. Similarly, we find the upper endpoint of the left interval in \eqref{eq:min_subbund_inc1} is less or equal that the upper endpoint of the right interval. This completes the proof. 
\qed
\end{thm}

\begin{cor}[Downward closed embeddings] \label{cor:downward_closed_subbund} Let $\scA : X \to \SIvert n \cC$, $\scB : Y \to \SIvert n \cC$, $p \in \tsG n(\scA)$ and let $\theta : \scB \subset \scA$ be an embedding. If $\im(\theta^0)$ is \gls{downwardclosed} in $X$ then $\theta^k(\tsG k(\scB))$ is \gls{downwardclosed} in $\tsG k(\scA)$ for $0 \leq k \leq n$. In this case $\scB$ is called a \gls{downwardclosed}.
\proof This follows from \autoref{claim:minimal_endpoints_delimit_neighbourhoods} and \autoref{claim:minimal_subfamily_is_minimal}. \qed
\end{cor}

\section{Normalisation on embeddings} \label{sec:emb_norm}

In this section we show that given $\theta : \scB \mono \scA$, then collapse of a parent family $\scA$ induces a collapse on its subfamily $\scB$. In more abstract terms we will show that the subcategory of $\Bunbc^n_\cC$ (cf. \autoref{defn:multilevel_base_change}) generated by multi-level collapses and embeddings admits an ``(epi, mono) factorization": Every morphism in this subcategory can be uniquely factored in a multi-level collapse followed by an embedding.

\subsection{Restricting collapse along embeddings}

We start with two examples for this section.

\begin{eg}[Collapse on subfamilies] \label{eg:subbund_collapse} \hfill
\begin{enumerate}
\item Recall $\theta_a : \scB_a \mono \scA_a$ from \autoref{eg:subfamily_by_restriction} and from \autoref{eg:normal_forms} recall the collapse $\lambda : \scA \kcoll 2 \widetilde\scB$. Set $\widetilde\scA_a := \scB$ and note $\scA_a$ equals $\scA$ of the latter example. Thus we have $\lambda : \scA_a \to \widetilde\scA_a$. We can then find $\lambda\postar \theta_a : \widetilde \scB_a \mono \widetilde \scA_a$ and $(\theta_a)\pbstar \lambda : \scB_a \kcoll 2 \widetilde\scB_a$ with the following data (such that the following commutes)
\begin{restoretext}
\begin{noverticalspace}
\begingroup\sbox0{\includegraphics{test/page1.png}}\includegraphics[clip,trim=0 {.0\ht0} 0 {.65\ht0} ,width=\textwidth]{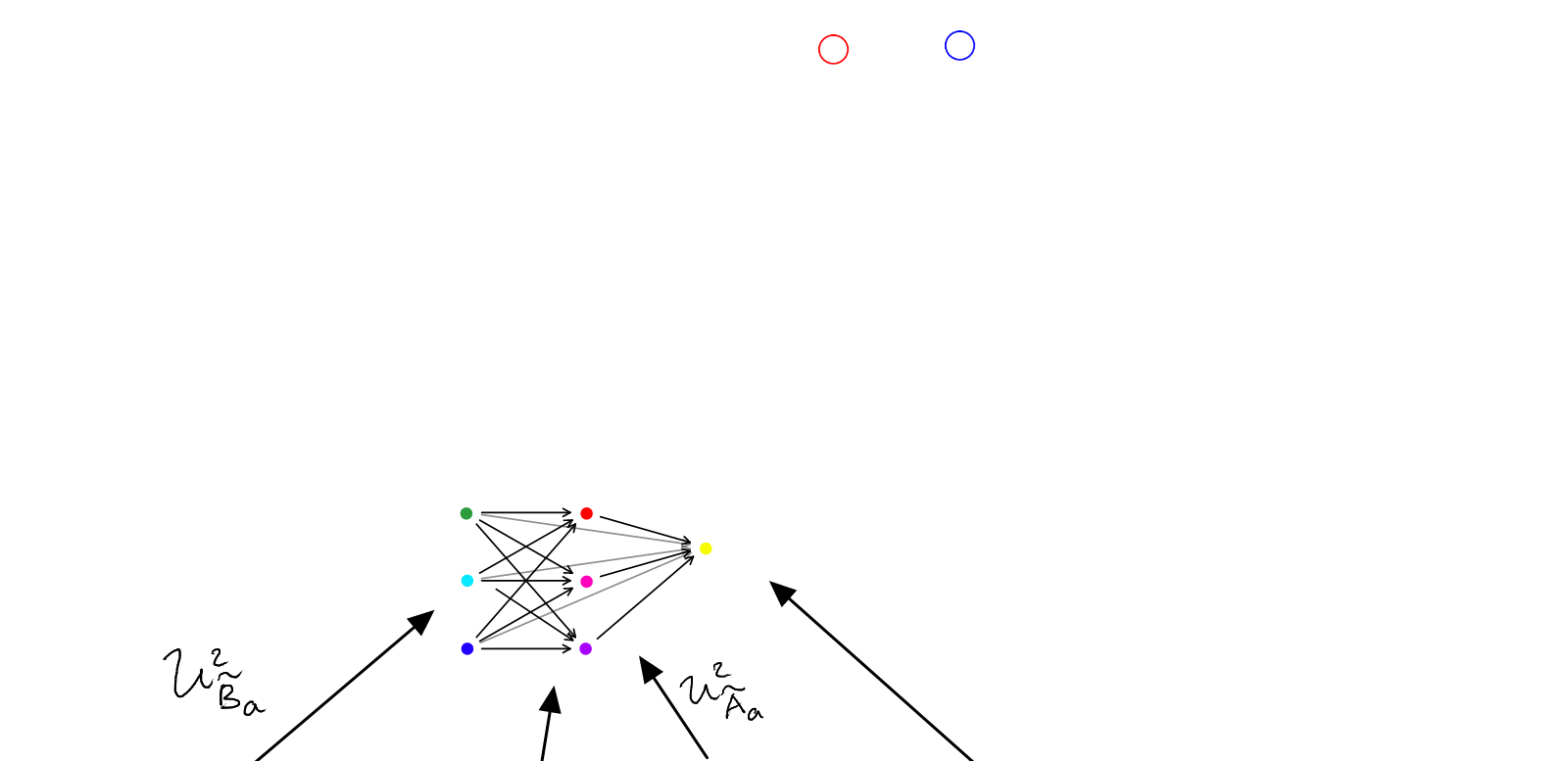}
\endgroup \\*
\begingroup\sbox0{\includegraphics{test/page1.png}}\includegraphics[clip,trim=0 {.0\ht0} 0 {.0\ht0} ,width=\textwidth]{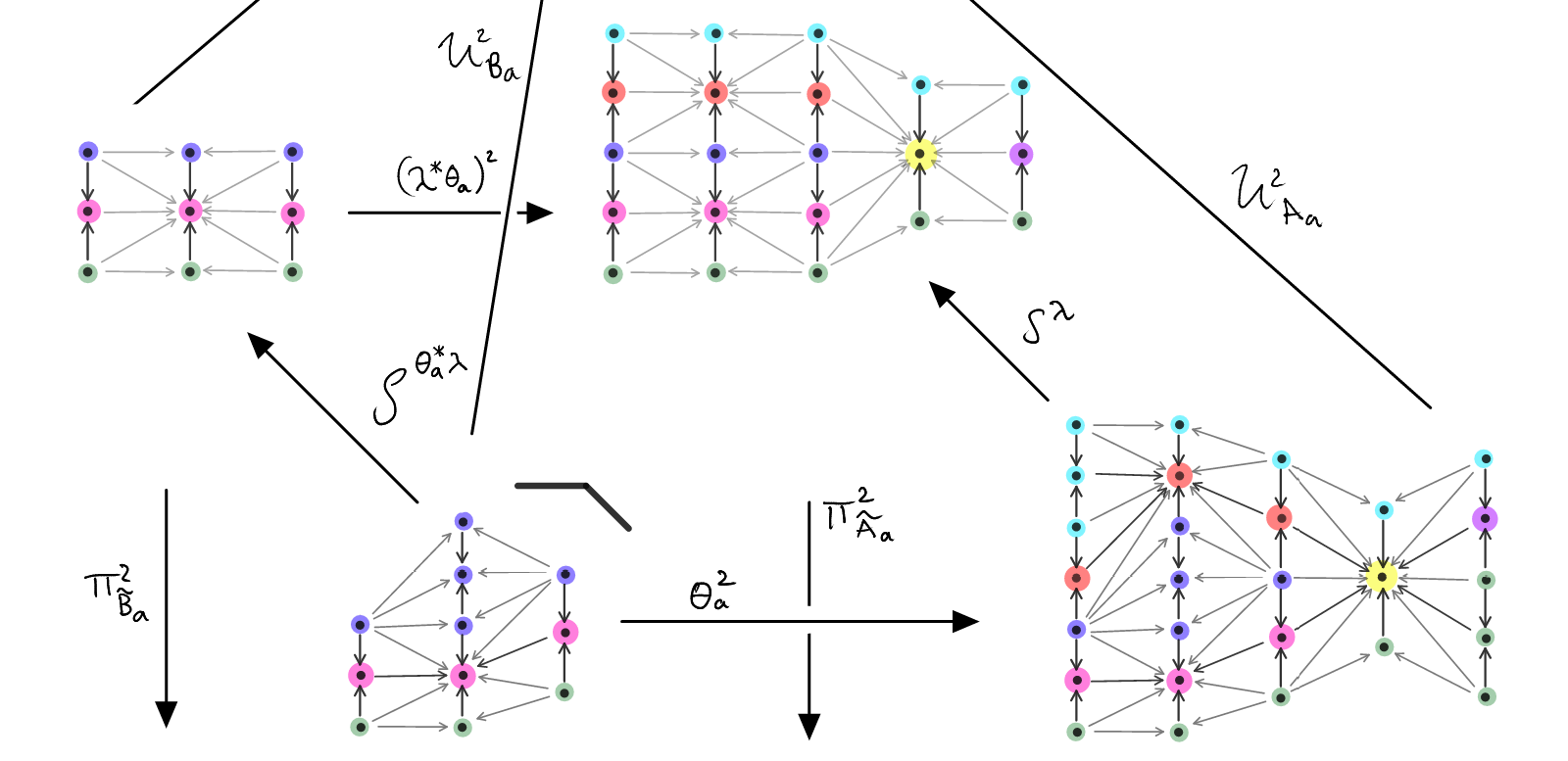}
\endgroup \\*
\begingroup\sbox0{\includegraphics{test/page1.png}}\includegraphics[clip,trim=0 {.0\ht0} 0 {.0\ht0} ,width=\textwidth]{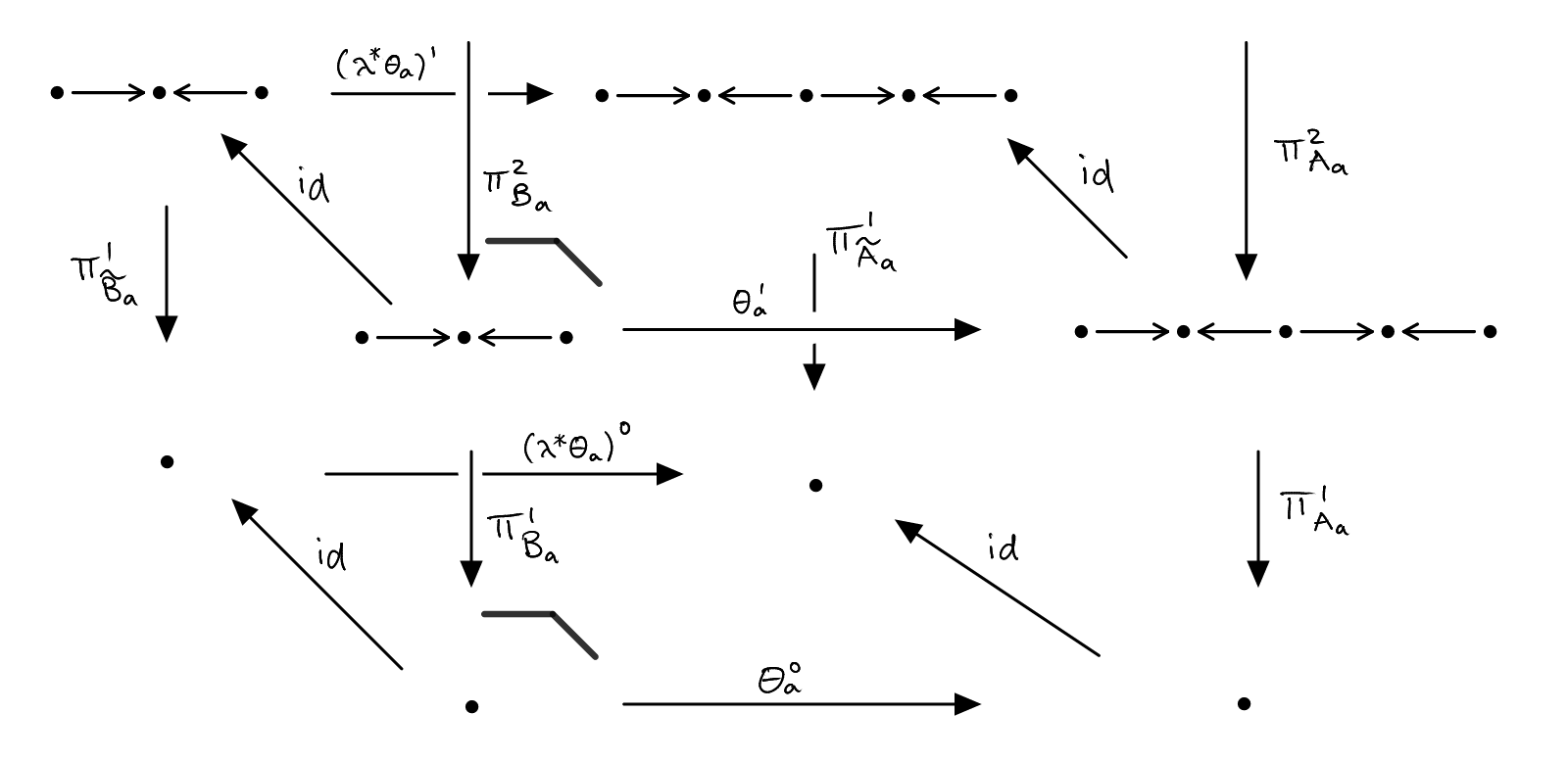}
\endgroup
\end{noverticalspace}
\end{restoretext}

\item Further to the above, set $\chi_a = \lambda\postar \theta : \widetilde \scB_a$,  from \autoref{eg:normal_forms} recall the collapse $\mu : \scB \kcoll 1 \scC$. Set $\widetilde{\widetilde\scA}_a := \scC$ and note $\widetilde \scA_a$ equals $\scB$ from the latter example.
We can then find $\lambda\postar \theta_a : \widetilde{\widetilde \scB}_a \kcoll 1 \widetilde{\widetilde\scA}_a$ and $(\theta_a)\pbstar \lambda : \widetilde\scB_a \kcoll 1 \widetilde{\widetilde\scB}_a$ such that the following commutes
\begin{restoretext}
\begin{noverticalspace}
\begingroup\sbox0{\includegraphics{test/page1.png}}\includegraphics[clip,trim=0 {.0\ht0} 0 {.6\ht0} ,width=\textwidth]{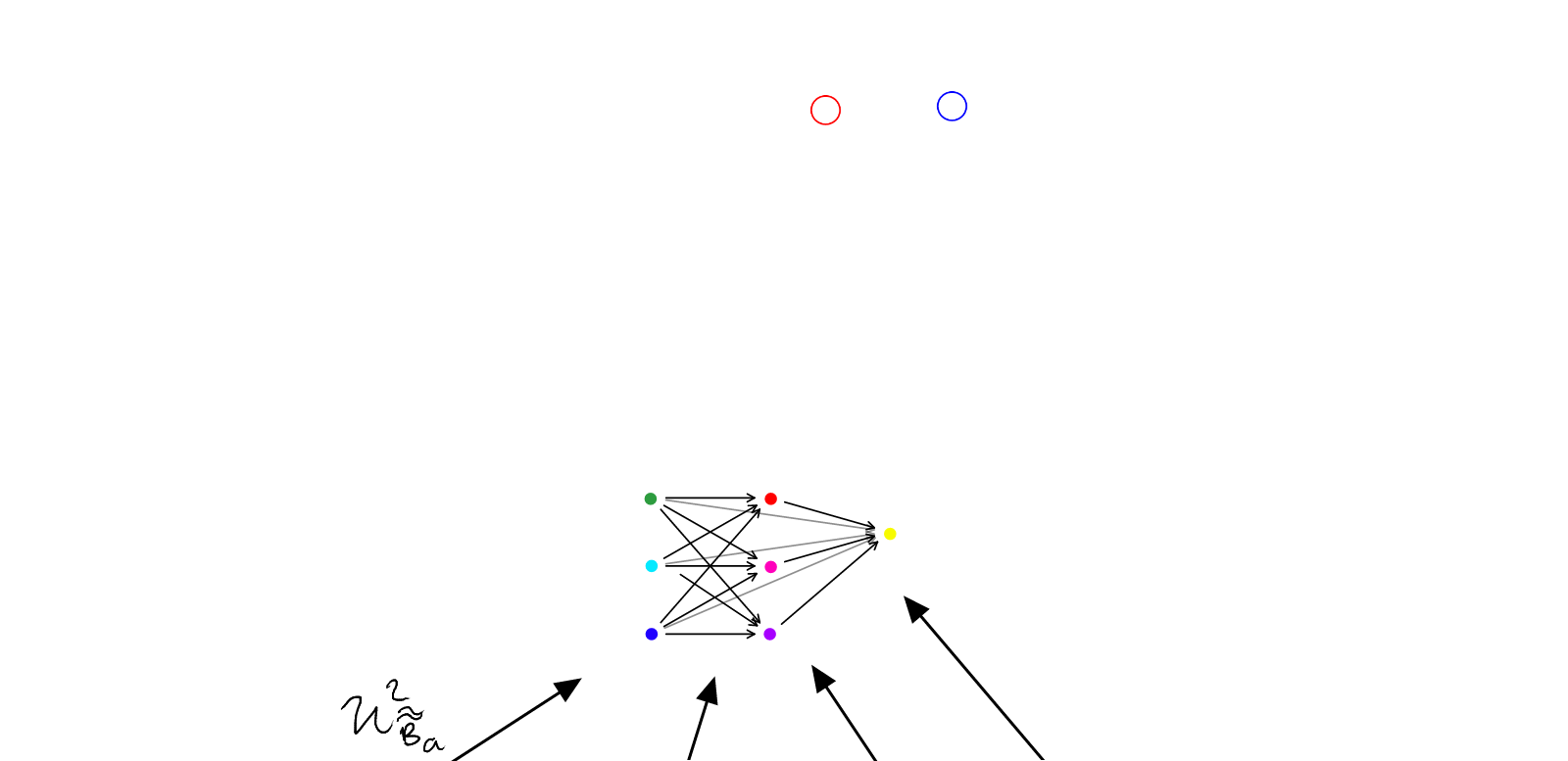}
\endgroup \\*
\begingroup\sbox0{\includegraphics{test/page1.png}}\includegraphics[clip,trim=0 {.0\ht0} 0 {.0\ht0} ,width=\textwidth]{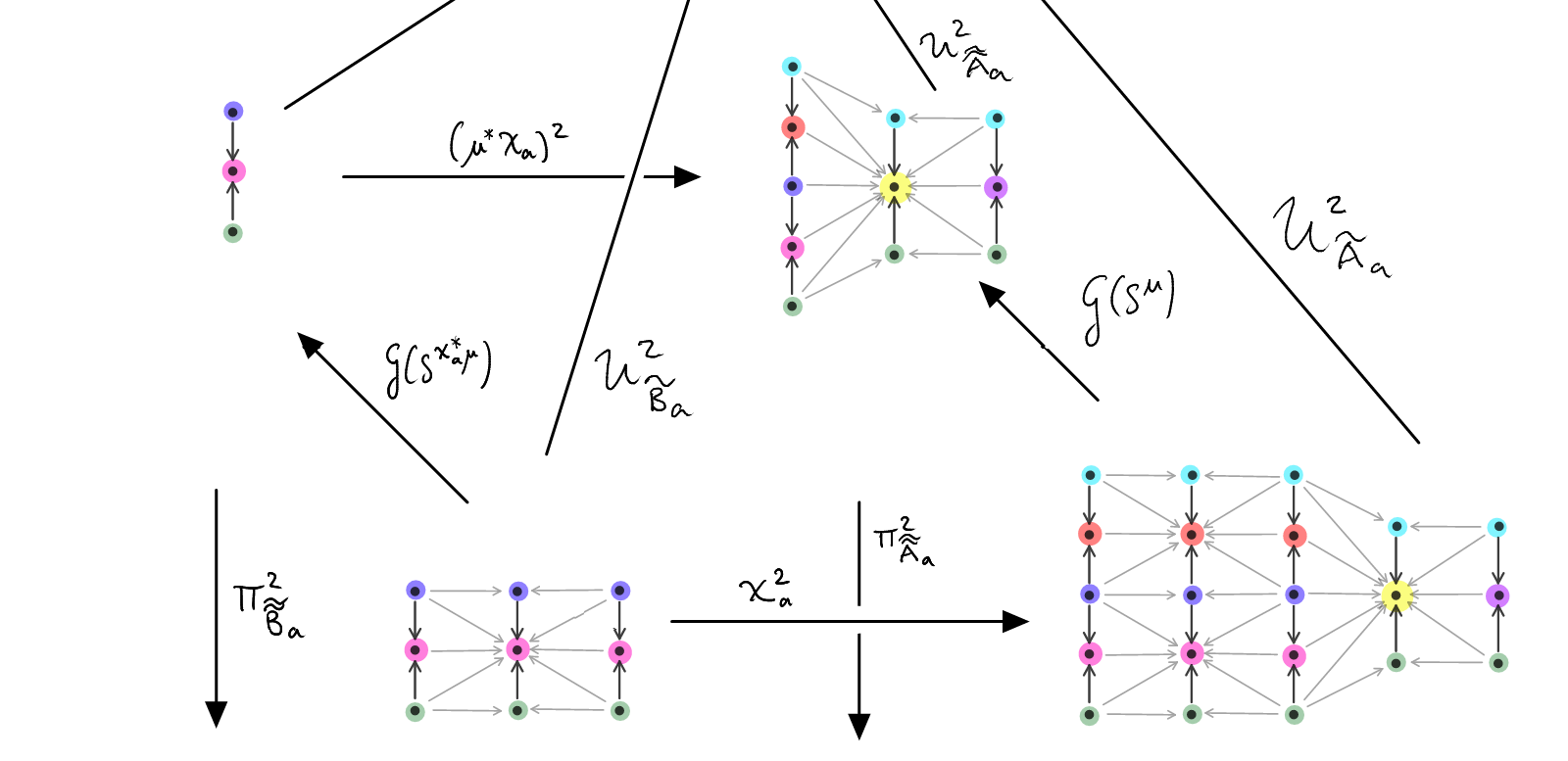}
\endgroup \\*
\begingroup\sbox0{\includegraphics{test/page1.png}}\includegraphics[clip,trim=0 {.0\ht0} 0 {.0\ht0} ,width=\textwidth]{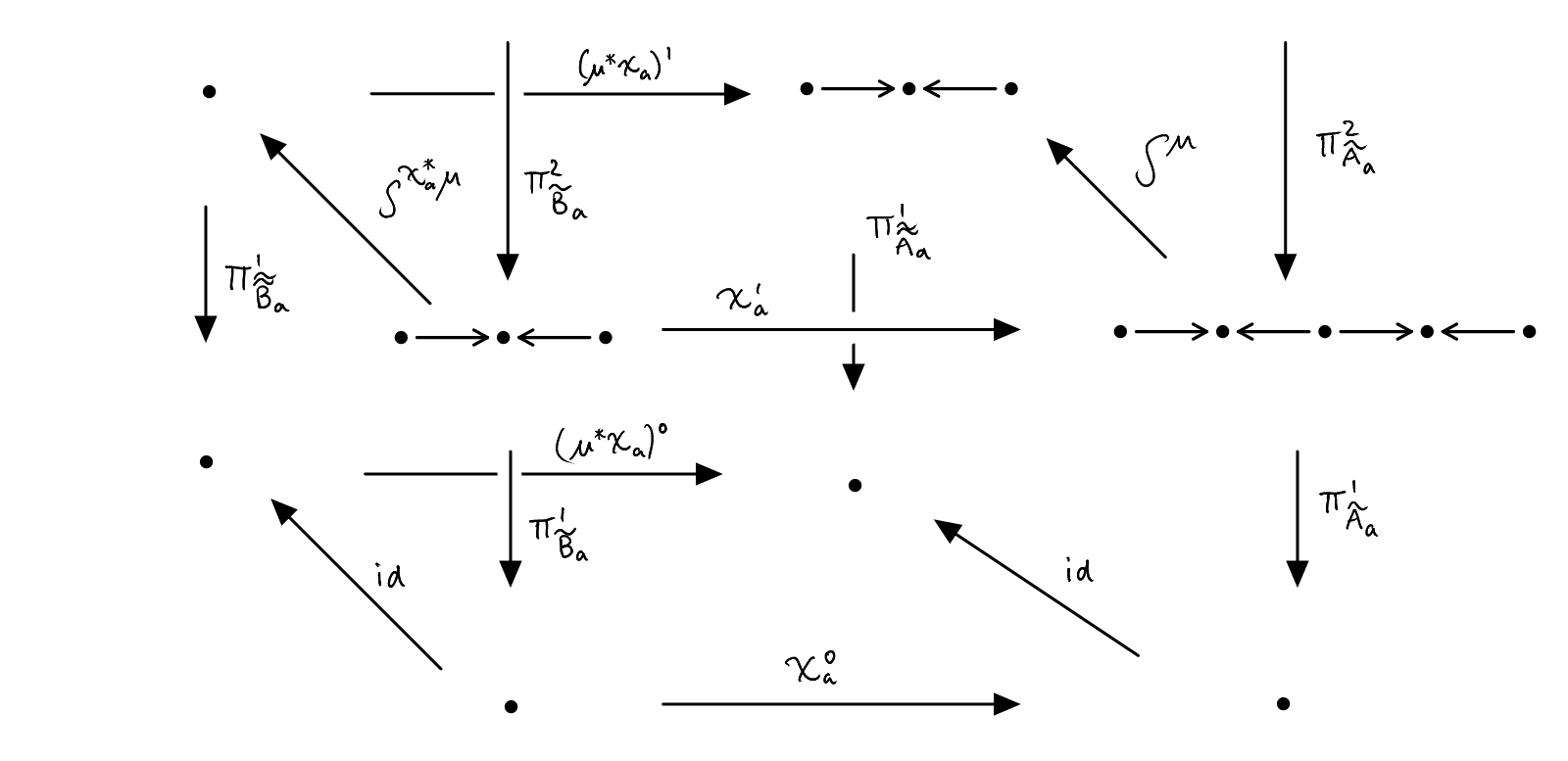}
\endgroup
\end{noverticalspace}
\end{restoretext}
\end{enumerate}
\end{eg}

The ``induced" embeddings and collapses which were denoted using pullback notation ($^*$) in the previous example, are part of a more general construction as the following will show.

\begin{constr}[Restricting a collapse to an embedding] \label{constr:collapse_on_subfamilies}  Let $\scA : X\to \SIvert n \cC$, $\scB : Y \to \SIvert n \cC$ and $\theta : \scB \subset \scA :  X \to \SIvert n \cC$ be a $\cC$-labelled $n$-cube embedding (cf. \autoref{defn:subfamilies}). Assume $\lambda : \scA \kcoll k \widetilde\scA$ ($1 \leq k \leq n$). We will construct $\theta\pbstar \lambda : \scB \kcoll k \widetilde\scB$ (called the \textit{restriction of $\lambda$ to $\theta$}) and $\lambda\postar \theta : \widetilde\scB \subset \widetilde\scA$ (called the \textit{collapse of $\theta$ by $\lambda$}) such that the following commutes
\begin{equation} \label{eq:subfamily_and_collapse}
\xymatrix@C=2cm{ \tsG l(\scB) \ar[d]_{\tsG {l-k}(\sS^{\theta\pbstar \lambda})} \ar[r]^{\theta^l} \pullbackfar & \tsG l(\scA) \ar[d]^{\tsG {l-k}(\sS^\lambda)} \\
\tsG {l}(\widetilde\scB) \ar[r]_{(\lambda\postar \theta)^l}  & \tsG {l}(\widetilde\scA) }
\end{equation}
and is a pullback. 

Using \autoref{constr:subfamilies_by_specifying_endpoints} and the fact that pullbacks compose we can assume $\theta$ to be a $k$-level base change, that is, it is either of the form
\begin{equation}
\theta = \restemb_H
\end{equation}
for some $H : X \to Y$, or it is of the form
\begin{equation}
\theta = \sJ^{\scA,l}\restsec{[q_-,q_+]}
\end{equation}
for some $l$ and $l$-level endpoints $(q_-,q_+)$. For notational convenience in this second case, we will assume $\theta$ to be the $l$-level basechange by a functor denoted by $H$. In other words
\begin{equation}
H = (\sJ^{\scA,l}\restsec{[q_-,q_+]})^l
\end{equation}

\begin{enumerate}
\item In the first case, we can set 
\begin{equation}
\lambda\postar\theta := \restemb_H : \widetilde\scB \mono \widetilde\scA
\end{equation}
We argue about a functor $S$ in the following diagram
\begin{equation}
\xymatrix{
\sG^k(\widetilde\scB)
	\ar[rr]^{\tsG k(H)}
	\ar[dd]^{\tpi k_{\widetilde\scB}}
&& \tsG k(\widetilde\scA)
	\ar[dd]|(0.45){\color{white}\txt{H\\A}}^(0.7){\tpi k_{\widetilde\scA}}
\\
& \tsG k(\scB)
	\ar[rr]^{\tsG k(H)}
	\ar[dd]^(0.3){\tpi k_\scB}
	\ar[ul]_{S}
&& \tsG k (\scA)
	\ar[ul]_{\sS^{\lambda}}
	\ar[dd]^(0.5){\tpi k_{\scA}}
\\
\tsG {k-1}(\widetilde\scB)
	\ar[rr]|(0.51){\hspace{9pt}}_(0.7){\tsG {k-1}(H)}
&& \tsG {k-1}(\widetilde\scA)
\\
& \tsG {k-1}(\scB)
	\ar[rr]_{\tsG {k-1}(H)}
	\ar[ul]^-{\id}
&& \tsG {k+i-1}(\scB)
	\ar[ul]^-{\id}
}
\end{equation}
Note $S$ exists since the back square is a pullback (of posets) and $S$ is a factorisation through it. The front square is also a pullback and thus the top square is a pullback by pullback cancellation on the right. Since $\sS^\lambda$ is surjective this means $S$ is surjective. It is also open since $\sG(H)$ are open and injective and $\sS^\lambda$ is open. Thus $S$ is a family collapse functor, and we deduce 
\begin{equation}
S = \sS^{\theta^*\lambda}
\end{equation}
for some $\theta^*\lambda : \sU^{k-1}_{\widetilde\scB} \into \sU^{k-1}_\scB$. Now arguing inductively in $i = 1,2,...,n-k$ we find diagrams and factorisations $S^i$
\begin{equation}
\xymatrix{
\tsG{k+i}(\widetilde\scB)
	\ar[rr]^{\tsG k(H)}
	\ar[dd]^{\tpi {k+i}(_{\widetilde\scB}}
&& \tsG{k+i}(\widetilde\scA)
	\ar[dd]|(0.45){\color{white}\txt{H\\A}}^(0.7){\tpi {k+i}(_{\widetilde\scA}}
\\
& \tsG {k+i}(\scB)
	\ar[rr]^{\tsG k(H)}
	\ar[dd]^(0.3){\tpi {k+i}_\scB}
	\ar[ul]_{S^i}
&& \tsG {k+i}( (\scA)
	\ar[ul]_{\tsG {i+1}(\sS^{\lambda})}
	\ar[dd]^(0.5){\tpi {k+i}_{\scA}}
\\
\tsG {k+i-1}(\widetilde\scB)
	\ar[rr]|(0.51){\hspace{9pt}}_(0.7){\tsG {k-1}(H)}
&& \tsG {k+i-1}(\widetilde\scA)
\\
& \tsG {k+i-1}(\scB)
	\ar[rr]_{\tsG {k-1}(H)}
	\ar[ul]^-{\tsG {i}(\sS^{\theta^*\lambda})}
&& \tsG {k+i-1}(\scB)
	\ar[ul]^-{\tsG {i}(\sS^{\lambda})}
}
\end{equation}
such that all sides but the left one are pullbacks (the bottom one is by inductive assumption). By pullback cancellation the left side is a pullback as well, and in particular this forces
\begin{equation}
\tsG {i+1}(\sS^{\theta^*\lambda}) = S^i
\end{equation}
We finally compute
\begin{align}
\tsU n_{\widetilde\scB} \tsG {n-k}(\sS^{\theta\pbstar \lambda}) &=   \tsU n_{\widetilde \scA} (\lambda\postar \theta)^n \tsG {n-k}(\sS^{\theta\pbstar \lambda}) \\
&=   \tsU n_{\widetilde \scA} \tsG {n-k}(\sS^{\lambda}) \theta^n \\
&=   \tsU n_{\scA} \theta^n \\
&=   \tsU n_{\scB}
\end{align}
And this implies $\theta^*\lambda :  \scB \kcoll k \widetilde\scB$ as required.

\item For the second case we set
\begin{equation}
\lambda\postar\theta := \sJ^{\widetilde\scB,k}\restsec{[\sG^{l-k}(\sS^\lambda)q_-,\sG^{l-k}(\sS^\lambda)q_+]} : \widetilde\scB \mono \widetilde\scA
\end{equation}
Assume that the $l$-level base change $\lambda\postar\theta$ to be given by a functor $\widetilde H$. In other words,
\begin{equation}
\widetilde H = (\sJ^{\widetilde\scB,k}\restsec{[\sG^{l-k}(\sS^\lambda)q_-,\sG^{l-k}(\sS^\lambda)q_+]})^l
\end{equation}

We distinguish subcases $l < k$, $l = k$ and $l > k$. First, if $l < k$ (and with our notational convention for $\theta$ using $H$) then the argument is similar to the previous item. If $l = k$, then we will find a diagram of the following form
\begin{equation}
\xymatrix{
\tsG k(\widetilde\scB)
	\ar[rr]^{\widetilde H}
	\ar[dd]^{\tpi k_{\widetilde\scB}}
&& \tsG k(\widetilde\scA)
	\ar[dd]|(0.45){\color{white}\txt{H\\A}}^(0.7){\tpi k_{\widetilde\scA}}
\\
& \tsG k(\scB)
	\ar[rr]^{H}
	\ar[dd]^(0.3){\tpi k_\scB}
	\ar[ul]_{S}
&& \tsG k (\scA)
	\ar[ul]_{\sS^{\lambda}}
	\ar[dd]^(0.5){\tpi k_{\scA}}
\\
\tsG {k-1}(\widetilde\scB)
	\ar[rr]|(0.51){\hspace{9pt}}_(0.7){\id}
&& \tsG {k-1}(\widetilde\scA)
\\
& \tsG {k-1}(\scB)
	\ar[rr]_{\id}
	\ar[ul]^-{\id}
&& \tsG {k-1}(\scB)
	\ar[ul]^-{\id}
}
\end{equation}
We argue about $S$. Note that by definition of $\widetilde H$ and \eqref{eq:endpoint_inclusion_image}, the objects in the fibrewise image of $\widetilde H$ are exactly the objects in the fibrewise image under $\sS^{\lambda} H$. Thus there is a factorising functor $S$ making the top left square commute. Since $\widetilde H$ is an embedding (cf. \autoref{rmk:subfamily_components_ff}) $S$ must be surjective. Arguing again fibrewise using the definition of $\widetilde H$, \eqref{eq:endpoint_inclusion_image} and monotonicity of $\sS^\lambda$, the image of $H$ is the preimage of the image of $\widetilde H$ under $\sS^\lambda$. This makes the top square a pullback. The rest of the argument is similar to the previous item.

Finally, we consider the case $l > k$. In this case we set
\begin{equation}
\theta^*\lambda := \lambda
\end{equation}
We find a diagram of the form
\begin{equation}
\xymatrix{
\tsG l(\widetilde\scB)
	\ar[rr]^{\widetilde H}
	\ar[dd]^{\tpi l_{\widetilde\scB}}
&& \tsG l(\widetilde\scA)
	\ar[dd]|(0.45){\color{white}\txt{H\\A}}^(0.7){\tpi l_{\widetilde\scA}}
\\
& \tsG l(\scB)
	\ar[rr]^{H}
	\ar[dd]^(0.3){\tpi l_\scB}
	\ar[ul]_{\tsG {l-k}(\sS^\lambda)}
&& \tsG l (\scA)
	\ar[ul]_{\tsG {l-k}(\sS^\lambda)}
	\ar[dd]^(0.5){\tpi l_{\scA}}
\\
\tsG {l-1}(\widetilde\scB)
	\ar[rr]|(0.51){\hspace{9pt}}_(0.7){\id}
&& \tsG {l-1}(\widetilde\scA)
\\
& \tsG {l-1}(\scB)
	\ar[rr]_{\id}
	\ar[ul]^-{\tsG {l-k-1}(\sS^\lambda)}
&& \tsG {l-1}(\scB)
	\ar[ul]^-{\tsG {l-k-1}(\sS^\lambda)}
}
\end{equation}
The top square commutes by definition of $\widetilde H$ (this can be argued for on each fibre). The left and the right square of this diagram are pullbacks and thus the top square is a pullback by pullback cancellation on the right. The rest of the argument is then similar to the previous item.
\end{enumerate}
This completes the construction.
\end{constr}

\subsection{Collapse sequences}

We end this chapter with a discussion of the interaction of embedding and collapse \textit{sequences}. Recall that collapse sequences are in correspondence with multilevel collapse. The decomposition of multilevel collapse into $k$-level base changes is analogous to the decomposition of embeddings into $k$-level base changes. Just as the previous construction generalised from the case of $k$-level embedding to multilevel embedding, we will now generalise to multilevel collapse. This will recover the result claimed in \autoref{thm:restriction_of_collapse}.

\begin{constr}[Restricting a collapse sequence to a subfamily] \label{constr:restricting_collapse_seq} Assume we are given an ordered collapse sequence $\vvec \lambda$
\begin{equation}
 \xymatrix{\scC_{n+1} \ar@{~>}[r]^-{\vvec\lambda ^n}_-n & \scC_{n} \ar@{~>}[r]^-{\vvec\lambda ^{n-1}}_-{n-1} & \scC_{n-1} \dots \ar@{~>}[r]^-{\vvec\lambda ^1}_-1 & \scC_1}
\end{equation}
and an embedding $\theta : \scD_{n+1} \mono \scC_{n+1}$. Using \autoref{constr:collapse_on_subfamilies} and setting $(\vvec \lambda\postar  \theta)_{n+1} := \theta$, we inductively set
\begin{align}
(\vvec \lambda\postar  \theta)_k &:= (\vvec\lambda^k)\postar  (\vvec \lambda\postar  \theta)_{k+1} \\
(\theta\pbstar \vvec\lambda)^k &:= (\vvec \lambda\postar  \theta)\pbstar _{k+1} \vvec\lambda^k
\end{align}
Putting the defined maps together we obtain
\begin{equation}
 \xymatrix@C=1.5cm{\scC_{n+1} \ar@{~>}[r]^-{\vvec\lambda ^n}_-n & \scC_{n} \ar@{~>}[r]^-{\vvec\lambda ^{n-1}}_-{n-1} & \scC_{n-1} \dots \ar@{~>}[r]^-{\vvec\lambda ^1}_-1 & \scC_1 \\
 \scD_{n+1} \ar@{>->}[u]^{\theta = (\vvec \lambda\postar  \theta)_{n+1}} \ar@{~>}[r]^-{(\theta\pbstar \vvec\lambda)^n}_-n & \scD_{n} \ar@{>->}[u]^{(\vvec \lambda\postar  \theta)_n} \ar@{~>}[r]^-{(\theta\pbstar \vvec\lambda)^{n-1}}_-{n-1} & \scD_{n-1} \ar@{>->}[u]_{(\vvec \lambda\postar  \theta)_{n-1}} \dots \ar@{~>}[r]^-{(\theta\pbstar \vvec\lambda)^1}_-1 & \scD_1 \ar@{>->}[u]_{(\vvec \lambda\postar  \theta)_1} }
\end{equation}
The bottom collapse sequence is denoted by $\theta\pbstar \vvec\lambda$. This is called the \textit{restriction of the collapse sequence $\vvec \lambda$ to the embedding $\theta$}. We also denote $\vvec\lambda\postar \theta := (\vvec\lambda\postar \theta)_1$. This is called the \textit{collapsed embedding $\theta$ witness by $\vvec \lambda$}. Using \autoref{eq:subfamily_and_collapse} and \autoref{constr:multilevel_collapse} we find that for all $k$ we have pullback squares
\begin{equation} \label{eq:restriction_of_multi-level_coll}
\xymatrix{\tsG k(\scD_{n+1} \ar[r]^{\theta^k} \ar[d]_{(\vsS{\theta\pbstar \vvec\lambda})^k} \pullbackfar & \tsG k(\scC_1) \ar[d]^{(\vsS{\vvec\lambda})^k} \\
\tsG k(\widetilde \scD_1) \ar[r]_{(\vvec\lambda\postar \theta)^k} & \tsG k{\widetilde \scC_1} }
\end{equation}
as claimed in \autoref{thm:restriction_of_collapse}. 
\end{constr}

\begin{rmk}[Restriction of collapse to minimal subbundle] \label{rmk:restrction_of_coll_on_minimal} If we restrict a collapse $\vvec\lambda$ to a minimal embedding $\iota^p_\scA$, then the embedding resulting from pullback along the collapse is also minimal. This can be seen level-wise from \autoref{constr:collapse_on_subfamilies}. As a consequence of this and \eqref{eq:restriction_of_multi-level_coll}, writing $\vsS{\vvec\lambda} : \scA \to \scB$ (for $\scA,\scB : X \to \SIvert n \cC$) and $q = (\vsS{\vvec\lambda})^n(p)$, we have for all $k$
\begin{equation} \label{eq:restrction_of_coll_on_minimal}
(\vsS{\vvec\lambda})^k (\iota^p_\scA)^k = (\iota^q_\scB)^k (\vsS{(\iota^p_\scA)\pbstar \vvec\lambda})^k
\end{equation}
\end{rmk}

We end with the remark that using \autoref{defn:multilevel_base_change}, the central result of this section can be summarised by saying that a composition of multi-level embedding and multi-level collapse in $\Bunbc^n_\cC$ has a unique ``(epi,mono) factorisation" into a composition of multi-level collapse followed by a multi-level embedding.

\chapter{Globes} \label{ch:globes}

The purpose of this chapter is two-fold: the first part will culminate in \autoref{lem:local_trivial_to_global_trivial}, which discusses a normalisation property of so-called locally trivial cube families. We will not need the result until the next chapter. To prove it we need to introduce the theory of \textit{sequenced collapse limits}, which will be a description of multi-level collapse maps of $n$-cubes families over $X$ in terms of families over $X \times (\bnum {n+1})$. This theory will yield the tools to make several important observation about normal forms (in particular, apart from \autoref{lem:local_trivial_to_global_trivial} it will also provide  insights into algorithms for computing normal forms).

The notion of local triviality is the basis for the second part of this chapter, where we will then enter the world of so-called \textit{globular} $n$-cubes (such cubes will also be called \textit{$n$-globes}). We will discuss interactions of collapse and embeddings with globularity and we will also show that $n$-globes admit a so-called ``double cone" construction. This in turn will be based on a natural generalisation of the ``adjoin-a-terminal-object"-monad to $n$-cubes, which we call the $\top$-monad.

More in detail we will see the following. Sequenced collapse limits are a reformulation of ordered collapse sequences as chains of morphisms of $\SIvert n \cC$ (at least in the case when the base space is trivial). In particular, multilevel collapses of cubes correspond to a certain class of morphisms in $\SIvert n \cC$. This is a special observation about collapse, which does not hold for its ``dual" notion of embedding. The discussion of this observation and its repercussions will take place in \autoref{sec:globe_lim}. In \autoref{sec:globe_globular} we will define $n$-globes to be cubes that are (at any projection level) locally trivial on regular segments. As we will see this will play a central role in capturing our usual intuition about manifold diagrams as semantics for $n$-categories (and sets them apart from $n$-fold categories). In \autoref{sec:globe_emb_coll}, we show that embeddings inherit globularity and globularity is preserved and reflected by collapse. In the final section \autoref{sec:globe_cone}, we construct $n$-globes called ``double cones" which are determined by their source and target: morally (and degenerating to the case of disks), this takes two $k$-cells which agree on their boundary, and then produces a new $(k+1)$-cell by taking the cone of the $k$-sphere formed by the two $k$-disks. This construction is only possible under the condition of globularity. As an application we will  define a special class of globes called \textit{terminal} $n$-globes, and these then correspond to classical $n$-globes (cf. \cite{leinster-operads}).

\section{Encoding collapse as limits} \label{sec:globe_lim}

In this section we will show that it is possible to encode the $k$-level collapse between two $\SIvert n \cC$-families over $X$ as a $\SIvert n \cC$-family over $X \times \bnum{2}$. In particular, if we restrict to working with families over $\bnum{1}$ this means that $k$-level collapses actually corresponds to certain morphisms in the category $\SIvert n \cC$. This reformulation will enable us to (concisely) write down a central proof later on, however, it will not play an important role otherwise in this thesis. We give two guiding examples for the reader, which she might refer to when trying to make sense of the constructions in this section.

\begin{eg}[Encoding $k$-level collapse as families]
\begin{enumerate}
\item Recall from \autoref{eg:normal_forms} the collapse $\lambda : \scA \kcoll 2 \scB$ where $\scA, \scB : \bnum{1} \to \SIvert 2 \cC$, which has the data
\begin{restoretext}
\begin{noverticalspace}
\begingroup\sbox0{\includegraphics{test/page1.png}}\includegraphics[clip,trim=0 {.0\ht0} 0 {.4\ht0} ,width=\textwidth]{ANCimg/page126.png}
\endgroup \\*
\begingroup\sbox0{\includegraphics{test/page1.png}}\includegraphics[clip,trim=0 {.0\ht0} 0 {.0\ht0} ,width=\textwidth]{ANCimg/page127.png}
\endgroup
\end{noverticalspace}
\end{restoretext}

This collapse corresponds to a family $L_\lambda : \bnum{2} \to \SIvert 2 \cC$ with the following data
\begin{restoretext}
\begin{noverticalspace}
\begingroup\sbox0{\includegraphics{test/page1.png}}\includegraphics[clip,trim=0 {.0\ht0} 0 {.6\ht0} ,width=\textwidth]{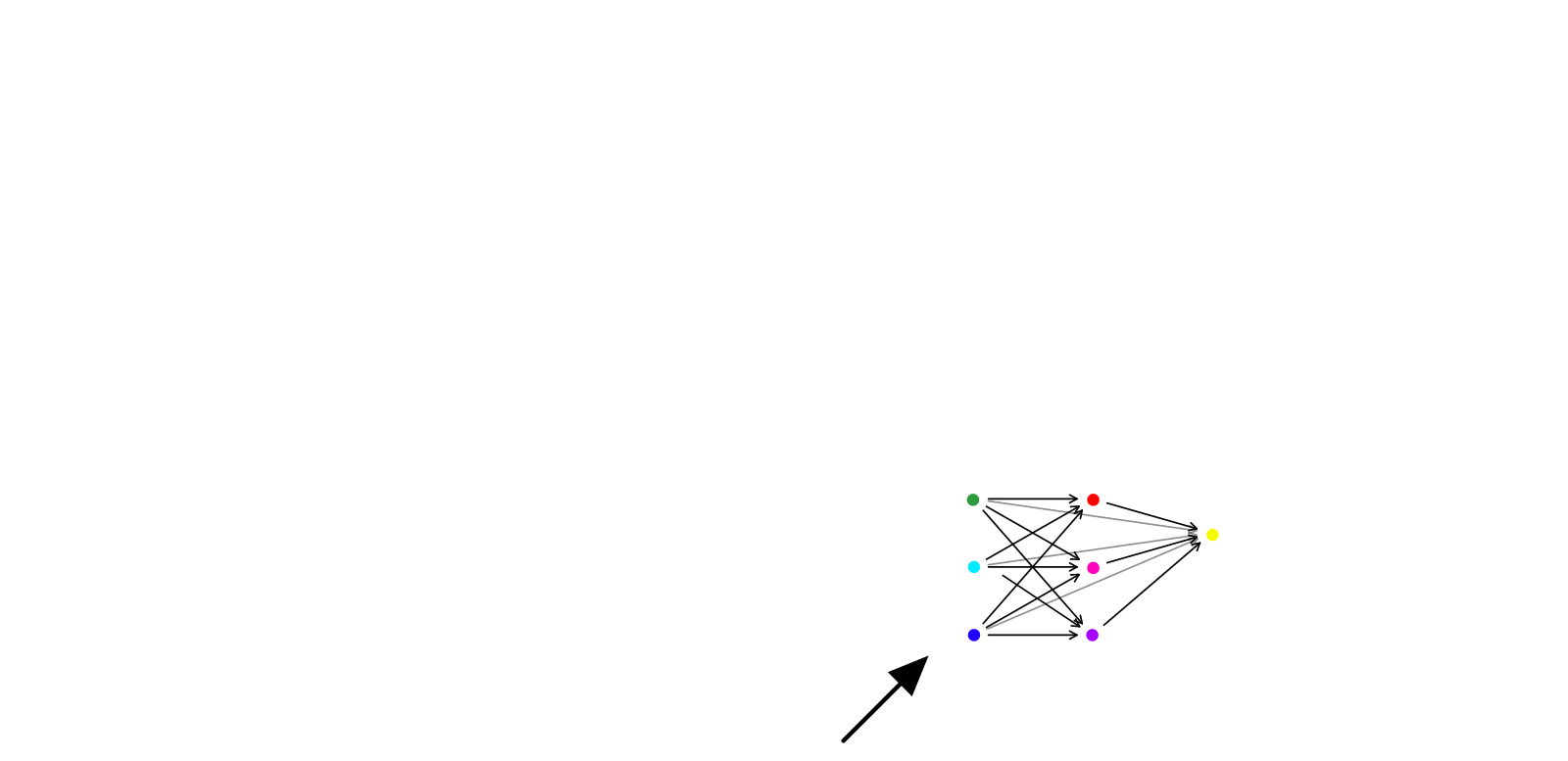}
\endgroup \\*
\begingroup\sbox0{\includegraphics{test/page1.png}}\includegraphics[clip,trim=0 {.0\ht0} 0 {.0\ht0} ,width=\textwidth]{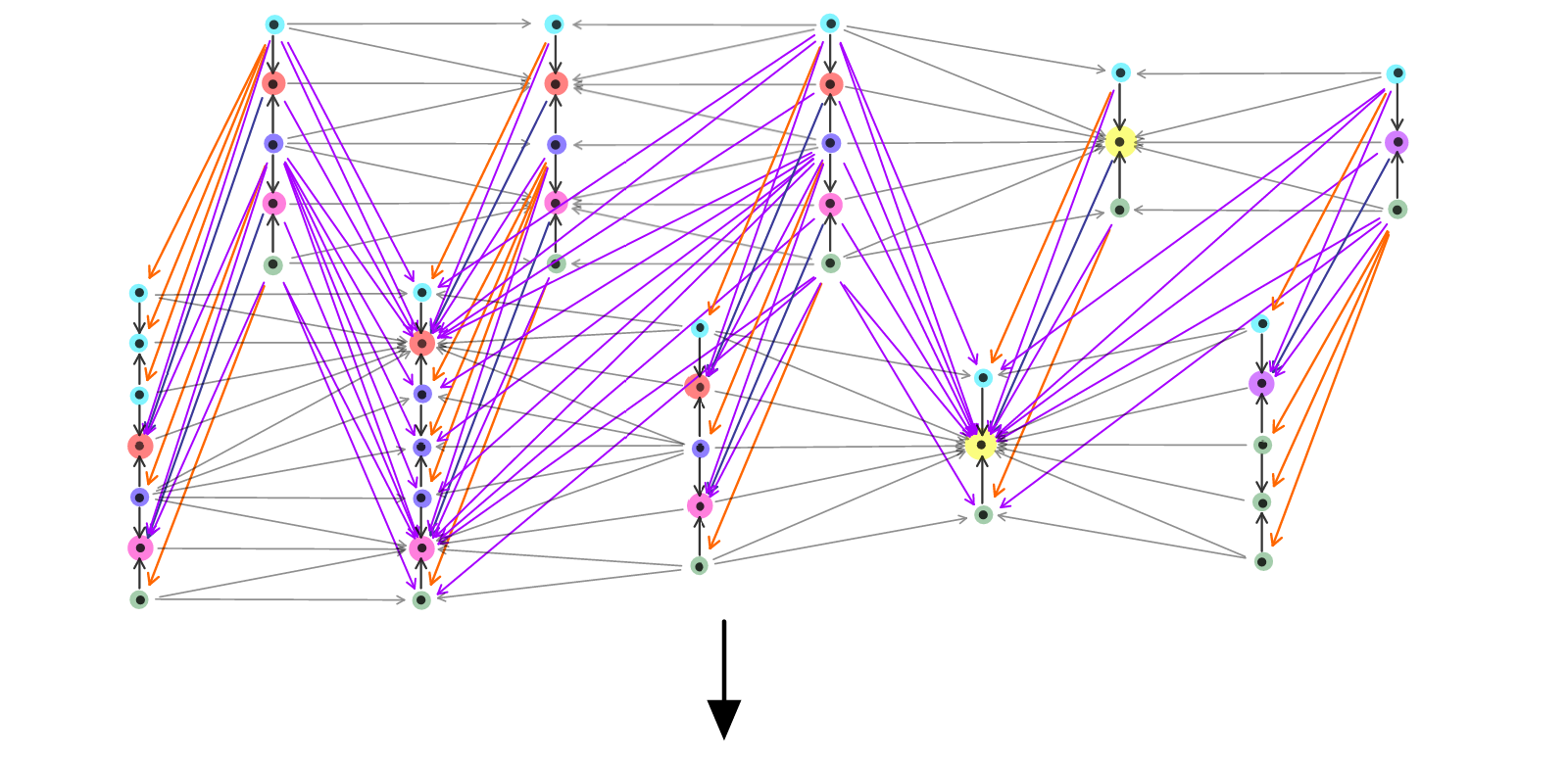}
\endgroup \\*
\begingroup\sbox0{\includegraphics{test/page1.png}}\includegraphics[clip,trim=0 {.0\ht0} 0 {.0\ht0} ,width=\textwidth]{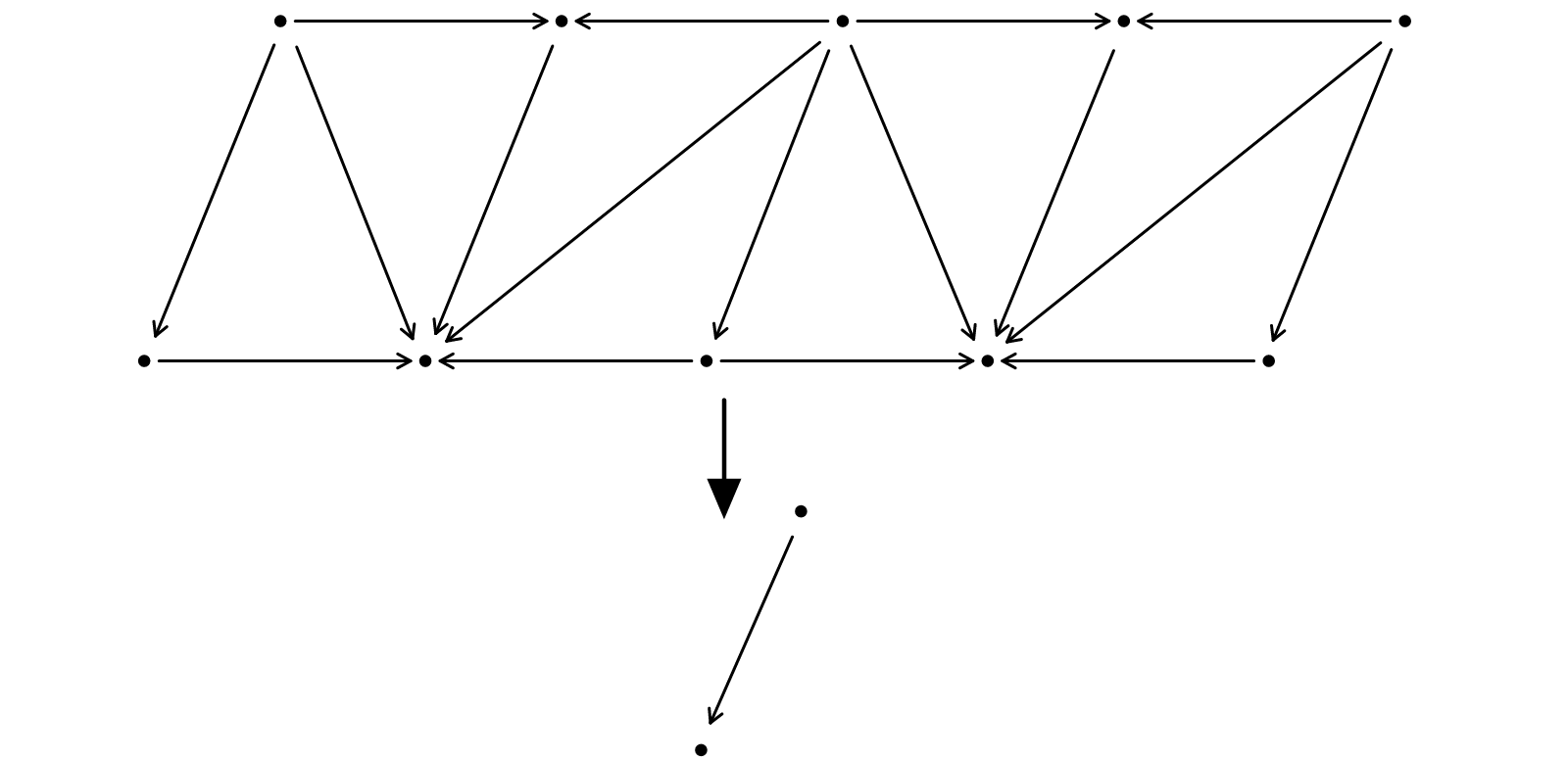}
\endgroup
\end{noverticalspace}
\end{restoretext}
Note that the restriction $L\Delta_0 : \bnum{1} \to \SIvert 2 \cC$ equals $\scB$ (the ``back part" of the above diagrams) and $L\Delta_1 : \bnum{1} \to \SIvert 2 \cC$ equals $\scA$ (the ``front part" of the above diagrams). Thus there are embeddings $\restemb^i_{\Delta_0} : \tsG i(\scC) \to \tsG i(L_\lambda)$ and $\restemb^i_{\Delta_1} : \tsG i(\scB) \to \tsG i(L_\lambda)$. Also note that \corange{} and \cdarkblue{} arrows in $\tsG 2(L_\lambda)$ exactly connect images $\sS^\lambda(x)$ (in $\im(\restemb^2_{\Delta_0})$) of $\sS^\lambda$ with their (possibly multiple) preimages $x$ (in $\im(\restemb^2_{\Delta_1})$). The \cdarkblue{} arrows by themselves connect preimages $a$ (in the fibre over $\restemb^1_{\Delta_0}(y)$, $y \in \tsG 1(\scB) = \tsG 1(\scA)$) with images $\lambda_y(a)$ (in the fibre over $\restemb^1_{\Delta_1}(y)$). As we will see, the \cpurple{} arrows can then be filled in uniquely under the condition of $\tusU 1_{L_\lambda}$ being an \SI-family. 

\item Recall from \autoref{eg:normal_forms} the collapse $\mu : \scB \kcoll 1 \scA$, which has the data
\begin{restoretext}
\begin{noverticalspace}
\begingroup\sbox0{\includegraphics{test/page1.png}}\includegraphics[clip,trim=0 {.0\ht0} 0 {.4\ht0} ,width=\textwidth]{ANCimg/page129.png}
\endgroup \\*
\begingroup\sbox0{\includegraphics{test/page1.png}}\includegraphics[clip,trim=0 {.0\ht0} 0 {.0\ht0} ,width=\textwidth]{ANCimg/page130.png}
\endgroup
\end{noverticalspace}
\end{restoretext}
This collapse corresponds to a family $L_\mu : \bnum{2} \to \SIvert 2 \cC$ with the following data
\begin{restoretext}
\begin{noverticalspace}
\begingroup\sbox0{\includegraphics{test/page1.png}}\includegraphics[clip,trim=0 {.0\ht0} 0 {.6\ht0} ,width=\textwidth]{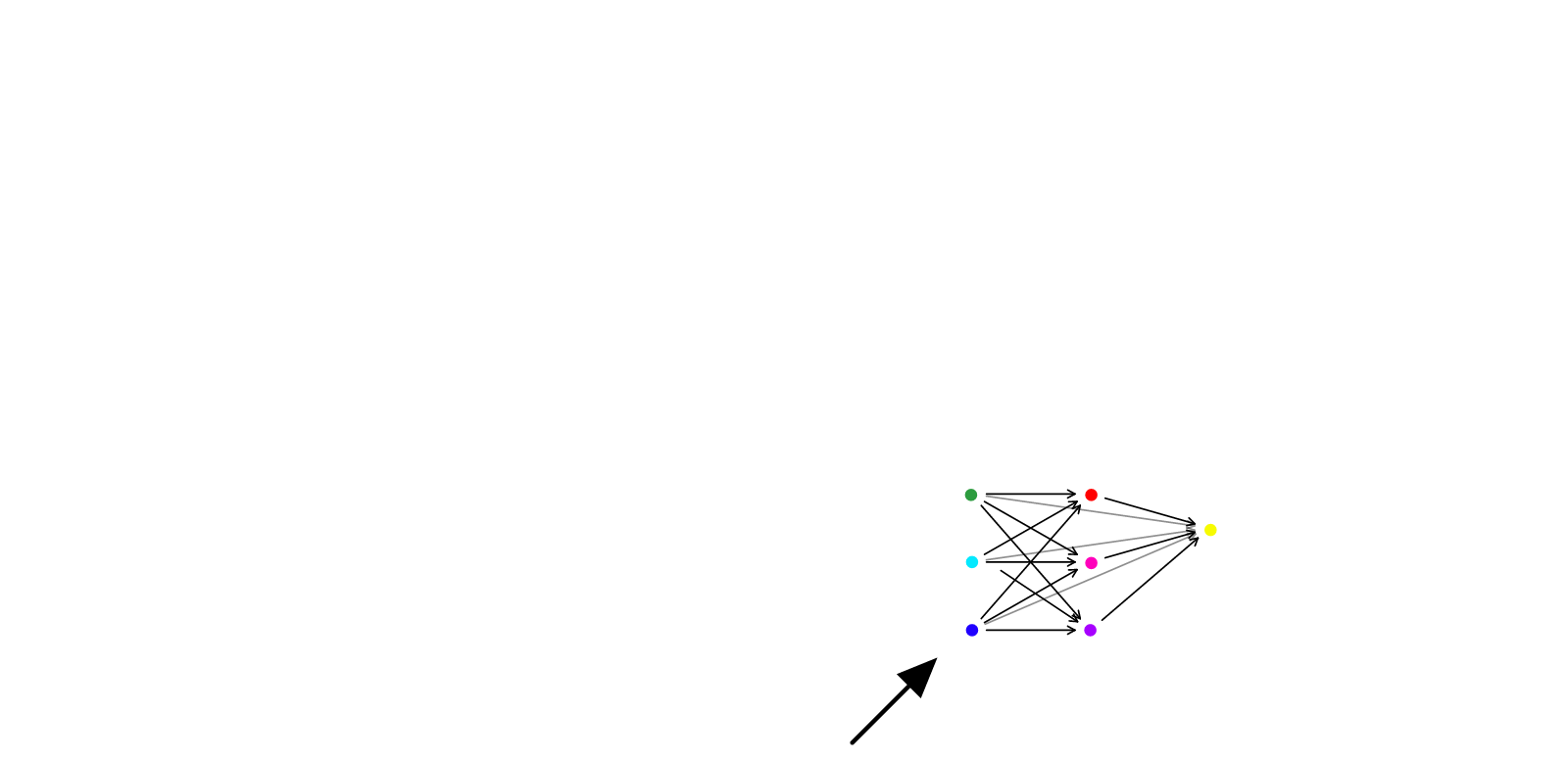}
\endgroup \\*
\begingroup\sbox0{\includegraphics{test/page1.png}}\includegraphics[clip,trim=0 {.0\ht0} 0 {.0\ht0} ,width=\textwidth]{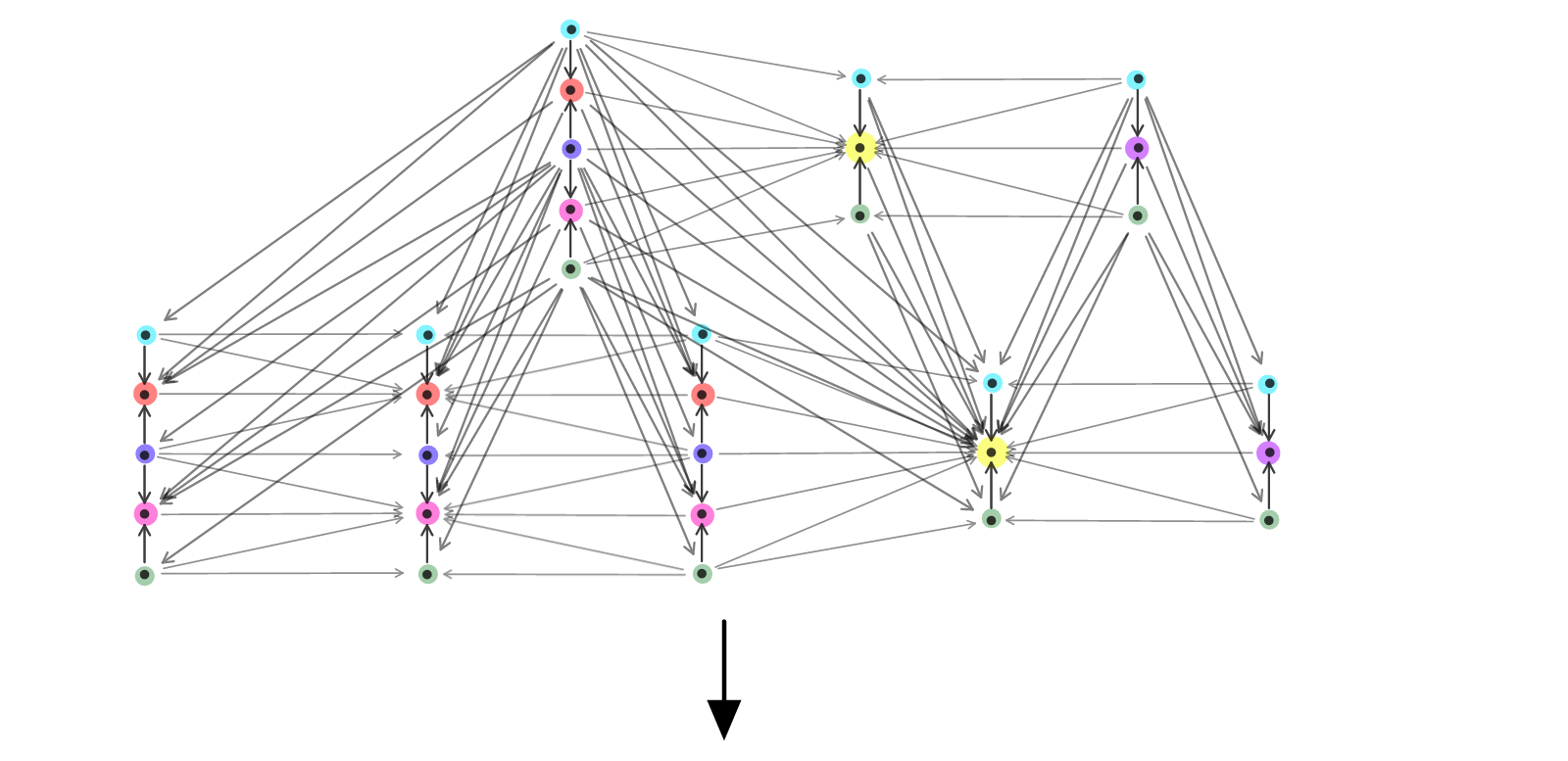}
\endgroup \\*
\begingroup\sbox0{\includegraphics{test/page1.png}}\includegraphics[clip,trim=0 {.0\ht0} 0 {.0\ht0} ,width=\textwidth]{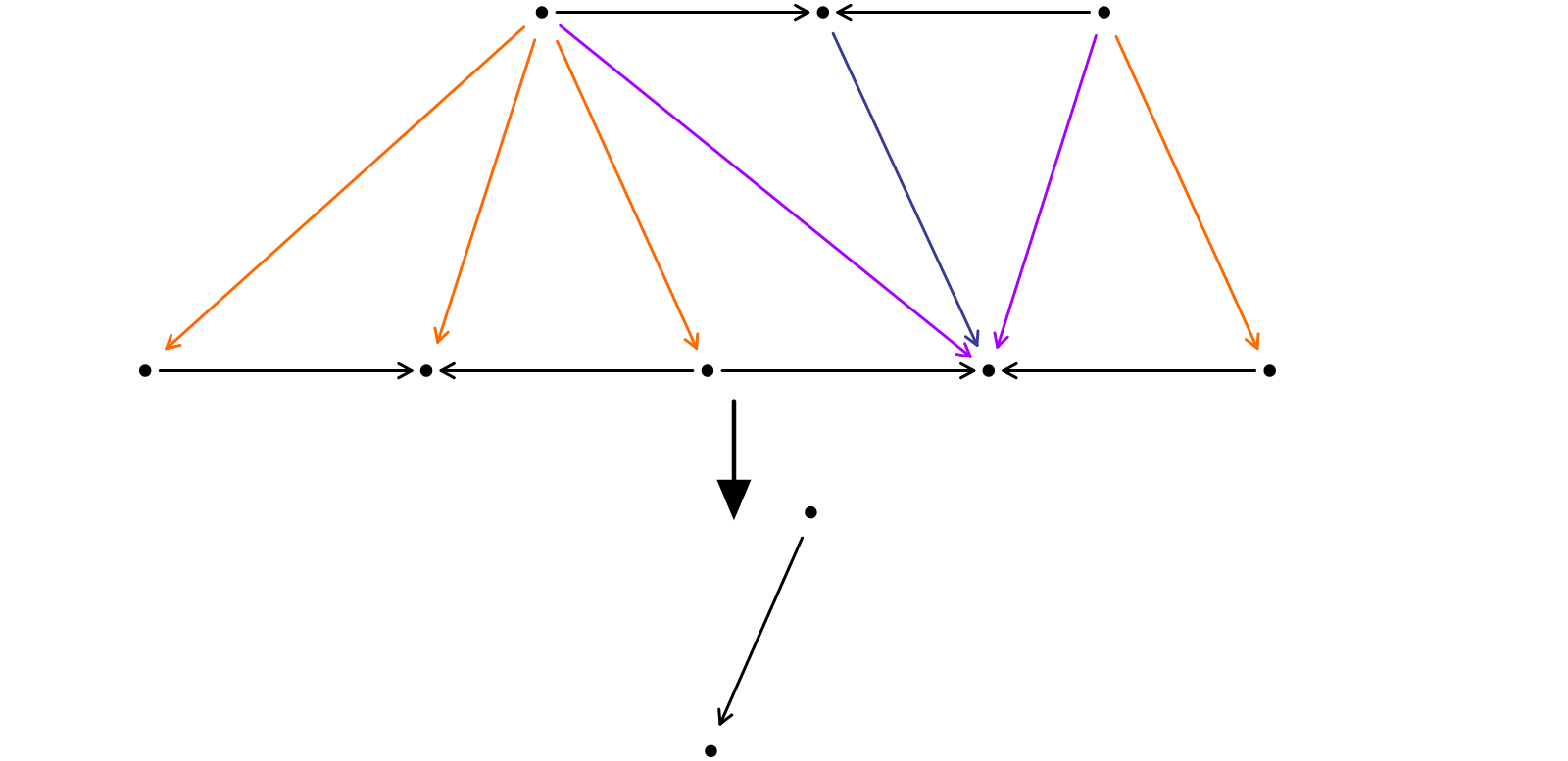}
\endgroup
\end{noverticalspace}
\end{restoretext}
Note that the restriction $L\Delta_0 : \bnum{1} \to \SIvert 2 \cC$ equals $\scC$ (the ``back part" of the above diagrams)and $L\Delta_1 : \bnum{1} \to \SIvert 2 \cC$ equals $\scB$ (the ``front part" of the above diagrams). thus there are inclusions $\restemb^i_{\Delta_0} : \tsG i(\scC) \to \tsG i(L_\mu)$ and $\restemb^i_{\Delta_1} : \tsG i(\scB) \to \tsG i(L_\mu)$. Also note that \corange{} and \cdarkblue{} arrows in $\tsG 1(L_\mu)$ exactly connect images $\sS^\mu(x)$  (in $\im(\restemb^1_{\Delta_0})$) of $\sS^\mu$ with their (possibly multiple) preimages $x$ (in $\im(\restemb^1_{\Delta_1})$). The \cdarkblue{} arrows in turn connect preimages $a$ (in the fibre over $\restemb^0_{\Delta_0}(y)$, $y \in \tsG 0(\scC) = \tsG 0(\scB)$) with images $\mu_y(a)$ (in the fibre over $\restemb^0_{\Delta_1}(y)$). The \cpurple{} arrows can be then filled in uniquely be the condition for $\tusU 1_{L_\mu}$ to be an \SI-family. 
\end{enumerate}
\end{eg}

Note that in the above example $L_\mu(1) = L_\lambda(0)$ (namely $L_\lambda\Delta_0 = L_\mu\Delta_1 = \scB$). Thus the morphisms $L_\mu(0 \to 1)$ and $L_\lambda(0 \to 1)$ compose. This means the sequence of collapses $(\lambda, \mu) : \scA \kcoll 2 \scB \kcoll 1 \scC$ actually corresponds to a family
\begin{equation}
L_{(\lambda, \mu)} : \bnum 3 \to \SIvert 2 \cC
\end{equation}
More generally, a central goal of this section will be to construct a correspondence between ordered collapse sequences and collapse limits as follows
\begin{equation}
\big( \xymatrix{\scA = \scA_{n+1} \ar@{~>}[r]^-{\lambda^n}_-n & \scA_{n} \ar@{~>}[r]^-{\lambda^{n-1}}_-{n-1} & \scA_{n-1} \dots \ar@{~>}[r]^-{\lambda^1}_-1 & \scA_1}\big)\quad \mapsto \quad (\sL^n_{\vvec \lambda} : X \times (\bnum {n+1}) \to \SIvert n \cC)
\end{equation}
where $\scA : X \to \SIvert n \cC$. This will be proven in \autoref{thm:n_lvl_collapse_limits_represent_collapses_seq}, which this section will culminate in.

\subsection{Interval family collapse}

We start by recalling the following two pieces of notation.

\begin{notn}[Two-variable functors] \label{notn:two_variable_functors} A functor $F : \cC \times \cD \to \cE$ is called a \textit{two-variable functor}, and has \textit{partial functors} $F(c,-) : \cD \to \cE$ and $F(-,d) : \cC \to \cE$ for $c \in \cC$, $d \in \cD$, defined by setting $F(c,-)(d') := F(c,d')$ on objects $d'$, $F(c,-)(f) := F(\id_c,f)$ on morphisms $f$, respectively $F(-,d)(c') := F(c',d)$ on objects $c'$, $F(-,d)(g) = F(g,\id_d)$ on morphisms $g$. We remark that alternatively $F(-,d)$ can be defined as the composite
\begin{equation} \label{eq:partial_functor_characterisation}
F(-,d) := \big(\cC \iso \cC \times \bnum{1} \xto {\cC \times \Delta_d} \cC \times \cD \xto F \cE\big)
\end{equation}
A similar reformulation holds for $F(c,-)$. We will often use the isomorphism $\cC \iso \cC \times \bnum{1}$ implicitly.

Conversely, given families of functors $F(c,-) : \cD \to \cE$, $F(-,d) : \cC \to \cE$ for $c \in \cC$, $d \in \cD$ then there is a unique two-variable functor $F : \cC \times \cD \to \cE$ which recovers $F(c,-)$, $F(-,d)$ as partial functors if and only if $F(c,-)(d) = F(-,d)(c)$ and the following \textit{commutativity condition} holds: all squares of the form
\begin{equation} 
\xymatrix@C=2cm{F(c,d) \ar[r]^{F(f : c \to c',d)} \ar[d]_{F(c, g : d \to d')} & F(c',d) \ar[d]^{F(c', g: d \to d')} \\ F(c,d') \ar[r]_{F(f : c \to c',d')} & F(c',d') }
\end{equation}
commute. Here, and in the following, we denote $F(c,-)(x) = F(c,x)$ and $F(-,d)(x) = F(x,d)$ for appropriate objects $x$ or morphisms $x$ (this notation is single-valued by assumption that $F(c,-)(d) = F(-,d)(c)$, and consistent with the definition of $F$ which shares the same notation).
\end{notn}

\begin{notn}[Uncurrying natural transformations] \label{notn:cartesian_closed}
$\Cat$, the category of categories and functors is \textit{cartesian closed}: Given categories $\cC, \cD$ we denote by $[\cC,\cD]$ the internal hom (which is the category of functor and natural transformations between them). Thus, a natural transformation $\lambda : F \to G$ is thus equivalently a functor $\bnum{2} \to [\cC, \cD]$ mapping $(0 \to 1)$ to $\lambda$. Using cartesian closedness we further find
\begin{equation} \label{eq:cartesian_closedness}
\Fun(\bnum{2},[\cC,\cD]) \iso \Fun(\cC \times \bnum{2}, \cD)
\end{equation}
We refer to the natural isomorphism from right to left as \textit{currying}, and from left to right as \textit{uncurrying}. Thus a natural transformation $\lambda : F \to G$ is equivalently a functor $\cY_\lambda : \cC \times \bnum{2} \to \cD$. $\cY_\lambda$ is called the \textit{uncurrying of $\lambda$}. Explicitly, $\cY_\lambda$ has partial functors $\cY_\lambda(-,0) = F$, $\cY_\lambda(-,1) = G$ and $\cY_\lambda(c,-) = \lambda_c$ (for $c \in \cC$). Conversely, a functor $F : \cC \times \bnum{2} \to \cD$ is equivalently a natural transformation $\cN_F : F(-,0) \to F(-,1)$ with components $(\cN_F)_c = F(c,-)$ 
\end{notn}

We now introduce a notion of ``extended" singular collapse (this is unrelated to extension $\wwidehat f$ of $\SI$-morphisms $f$, but instead extends the domain of collapse to include its codomain in a natural way).

\begin{constr}[Extended singular collapse] \label{constr:uncurrying} Let $\lambda : \scA \into \scB$ be in injection of $\SI$-families $\scA, \scB : X \to \SI$ over a poset $X$. Using cartesian closedness of $\Cat$, we construct $\cY_\lambda : X \times \bnum{2} \to \SI$ by uncurrying. $\cY_\lambda$ is thus a $\SI$-family over the poset $X \times \bnum{2}$.

We now construct a functor of posets
\begin{equation}
\sE^\lambda : \sG(\cY_\lambda) \to \sG(\scA) = \sG(\cY_\lambda(-,0))
\end{equation}
called the \textit{extended singular collapse} by setting
\begin{equation}
\sE^\lambda ((x,0),a) = (x,a)
\end{equation}
and
\begin{equation}
\sE^\lambda ((x,1),b) = \sS^\lambda (x, b)
\end{equation}
It is \stfwd{} to verify that this defines a functor of posets. Functoriality of $\sE^\lambda$ can be inferred from functoriality of $\sS^\lambda$ (which was proven in \autoref{lem:glambda}) and its lifting property (which was proven in \autoref{lem:collapse_has_lifts}) as follows: let
\begin{equation}
f \equiv ((z_1,j_1),c_1) \to ((z_2,j_2),c_2)
\end{equation}
be a morphisms in the poset $\sG(\cY_\lambda)$. To prove that $\sE^\lambda$ is a functor of posets we need to show
\begin{equation}
\sE^\lambda ((z_1,j_1),c_1) \to \sE^\lambda ((z_2,j_2),c_2)
\end{equation}
Here, $j_i \in \Set{0,1}$ and $j_1 \leq j_2$. We distinguish three cases.
\begin{enumerate}
\item If $j_1 = j_2 = 0$, then $(c_1,c_2) \in \edgeset(\scA(z_1 \to z_2))$, and $\sE^\lambda(f) = (z_1,c_2) \to (z_2,c_2)$ is a morphism in $\sG(A)$
\item If $0 = j_1 < j_2 = 1$, then we prove
\begin{equation}
\big((z_1,c_1) \to \sS^\lambda(z_2,c_2)\big) \in \mor(\sG(\scA))
\end{equation}
Set $\sS^\lambda(z_2,c_2) \equiv (z_2,c_3)$. First note
\begin{align}
\cY_\lambda ((z_1,j_1) \to (z_2,j_2)) &= \cY_\lambda ((z_2,0) \to (z_2,1)) \cY_\lambda ((z_1,0) \to (z_2,0)) \\
&= \lambda_{z_2} \scA(z_1 \to z_2)
\end{align}
Since $(c_1, c_2) \in \edgeset( \lambda_{z_2} \scA(z_1 \to z_2) )$ we have one of the following
\begin{enumerate}
\item $c_1 \in \singcont(\scA(z_1))$. Then by \eqref{eq:defn_order_realisation_1}
\begin{equation}
 \lambda_{z_2} \scA(z_1 \to z_2) (c_1) = c_2
\end{equation}
and since $\rest {\sS^\lambda} {z_2}$ in inverts $\lambda_{z_2}$ we find 
\begin{equation}
\scA(z_1 \to z_2) (c_2) = c_3
\end{equation}
and thus $(c_1, c_3) \in \edgeset(\scA(z_1 \to z_2))$ as required.
\item $c_1 \in \regcont(\scA(z_1))$. Then by \eqref{eq:defn_order_realisation_1}
\begin{equation}
 \lambda_{z_2} \scA(z_1 \to z_2) (c_1 - 1) \leq c_2 \leq  \lambda_{z_2} \scA(z_1 \to z_2) (c_1 + 1)
\end{equation}
Since $\rest {\sS^\lambda} {z_2}$ in inverts $\lambda_{z_2}$ and is monotone we find 
\begin{equation}
\scA(z_1 \to z_2) (c_1 - 1) \leq c_3 \leq  \scA(z_1 \to z_2) (c_1 + 1)
\end{equation}
and thus $(c_1, c_3) \in \edgeset(\scA(z_1 \to z_2))$ as required.
\end{enumerate}
\item If $j_1 = j_2 = 1$, then $(c_1,c_2) \in \edgeset(\scB(z_1 \to z_2)$, and $\sE^\lambda(f) = \sS^\lambda((z_1,c_2) \to (z_2,c_2))$ is a morphism in $\sG(A)$ by functoriality of $\sS^\lambda$.
\end{enumerate}
This finishes the verification of $\sE^\lambda$ being a functor of posets.
\end{constr}

\begin{rmk}[Restriction of extended collapses] \label{rmk:extended_collapse_props} Let $\lambda : \scA \into \scB$ be in injection of $\SI$-families $\scA, \scB : X \to \SI$ over a poset $X$. We remark that the above constructed map $\sE^\lambda$ satisfies
\begin{equation} \label{eq:extended_collapse_props1}
\Big( \sG(\scA)  \xinto {\sG(X \times \delta^1_1)} \sG(\cY_\lambda) \xto {\sE^\lambda} \sG(\scA) \Big) = \Big( \sG(\scA) \xto {\id} \sG(\scA) \Big)
\end{equation}
as well as
\begin{equation} \label{eq:extended_collapse_props2}
\Big( \sG(\scB)  \xinto {\sG(X \times \delta^1_0)} \sG(\cY_\lambda) \xto {\sE^\lambda} \sG(\scA) \Big) = \Big( \sG(\scB) \xto {\sS^\lambda} \sG(\scA) \Big)
\end{equation}
where we used \autoref{defn:grothendieck_base_change} for the commuting triangle ($i \in \Set{0,1}$)
\begin{equation}
\xymatrix{ X \iso X \times \bnum{1} \ar[rr]^{X \times \delta^1_i} \ar[dr]_{{\cY_\lambda}(-,i)} && X \times \bnum{2} \ar[dl]^{{\cY_\lambda}} \\
&\SIvertone \cC& }
\end{equation}
(which commutes by definition, cf. \eqref{eq:partial_functor_characterisation}).
\end{rmk}

\begin{rmk}[Base change of extended collapse] \label{rmk:precomposition_extended_collapse} Let $\lambda : \scA \into \scB$ be in injection of $\SI$-families $\scA, \scB : X \to \SI$ over a poset $X$. Let $H : Y \to X$ be a functor of posets. Using cartesian closedness of $\Cat$, together with the definition of $\sE^\lambda$ and the commutation relation derived in \autoref{claim:precomposition_of_collapse_map}, we find that
\begin{equation}
\xymatrix@C=1.7cm@R=0.6cm{ \sG(\cY_\lambda (H\times \bnum{2})) \ar[r]^-{\sG(H \times \bnum{2})} \ar[d]_{\sE^{ \lambda H}} & \sG(\cY_\lambda) \ar[d]^{\sE^\lambda} \\
\sG(\scA H) \ar[r]_{\sG(H)} & \sG(\scA) }
\end{equation}
is well-defined and commutes.
\end{rmk}

\begin{defn}[Collapse limits] Let $X$ be a poset. A \textit{collapse limit} over $X$  is a $\cC$-labelled $\SI$-family $\scC : X \times \bnum{2} \to \SIvertone \cC$, such that the following holds
\begin{itemize}
\item Firstly, we require that $\und \scC (x,0 \to 1)$ is injective for all $x \in X$. This means, currying $\und \scC$ we obtain an injection of $\SI$-families
\begin{equation} \label{eq:collapse_limit1}
\cN_{\und \scC} : \und \scC (-,0) \into \und \scC (-,1)
\end{equation}
\item Secondly, we then require
\begin{equation} \label{eq:collapse_limit2}
\sU_\scC = \sU_{\scC(-,0)} \sE^{\cN_{\und \scC}}
\end{equation}
\end{itemize}
\end{defn}

\begin{constr}[Correspondence of collapse limits to collapses] \label{constr:collapse_limits} 
Given a collapse limit $\scC : X \times \bnum{2} \to \SIvertone \cC$, set $\scB := \scC (-,0)$, $\scA := \scC (-,1)$. Then we claim that $\cN_{\und \scC}$ witnesses a collapse
\begin{equation}
\cN_{\und \scC} : \scA \mcoll \scB
\end{equation}
To see that this is indeed a collapse, by definition we need to verify that $\cN_{\und \scC} : \und \scB \into \und \scA$
is a natural transformation satisfying
\begin{equation} 
\sU_ \scA = \sU_\scB \sS^{\cN_{\und \scC}}
\end{equation}
The latter holds by precomposing \eqref{eq:collapse_limit2} with $\sG(X \times \delta^1_0)$ and then using \eqref{eq:extended_collapse_props2} together with \autoref{claim:grothendieck_span_construction_basechange}.

Conversely, let $\lambda : \scA \mcoll \scB$ be a collapse of $\cC$-labelled $\SI$-families $\scA,\scB : X \to \SIvertone \cC$ over a poset $X$, in particular we now have
\begin{equation} \label{eq:factorization_condition_limit}
\sU_ \scA = \sU_\scB \sS^{\lambda}
\end{equation}
We construct a family $\sL_\lambda : X \times \bnum{2} \to \SIvertone \cC$, called the \textit{collapse limit of $\lambda$}, by setting
\begin{equation} \label{eq:collapse_lim_of_injection}
\sL_\lambda = \sR_{\cY_\lambda, \sU_\scA \sE^{\lambda}}
\end{equation}
As is proven explicitly in \autoref{lem:collapse_lim_equiv_to_collapse} below, this satisfies $\sL_\lambda(-,0) = \scA$ and $\sL_\lambda(-,1) = \scB$ (the proof uses \autoref{rmk:extended_collapse_props} and \autoref{claim:grothendieck_span_construction_basechange}). We further need to check $\sL_\lambda$ satisfies condition \eqref{eq:collapse_limit2} for collapse limits. But this now holds by its definition \eqref{eq:collapse_lim_of_injection}:
\begin{align}
\sU_{\sL_\lambda} &:= \sU_\scA \sE^{\lambda} \\&=  \sU_{\sL_\lambda(-,0)} \sE^{\lambda}
\end{align}
\end{constr}

\begin{lem} \label{lem:collapse_lim_equiv_to_collapse} The constructions $(\lambda : \scA \mcoll \scB) \mapsto \sL_\lambda$ (turning a collapse into a collapse limit) and $\scC \mapsto \big(\cN_{\und \scC} : \scC(-,1) \mcoll \scC(-,0)\big)$ (turning a collapse limit into a collapse) are mutually inverse.
\proof The proof is \stfwd{}. We first show $\scC = \sL_{\cN_{\und \scC}}$. We need to show $\und \scC = \cY_{\cN_{\und \scC}}$ and $\sU_\scC = \sU_{\scC\delta^1_1} \sE^{\cN_{\und \scC}}$. The former follows from the isomorphism \eqref{eq:cartesian_closedness}, the latter follows from assumption \eqref{eq:collapse_limit2}.

Conversely, we want to show $(\lambda : \scA \mcoll \scB) = \cN_{\und {\sL_{\lambda}}}$. For this we first need to show $\scA = \sL_{\lambda} (-,0)$ and $\scB = \sL_{\lambda} (-,1)$. This follows by \autoref{claim:unpacking_and_repacking}, since firstly, both the underlying $\SI$-families ($\und {\sL_{\lambda}(-,0)} := \und \scA$ and $\und {\sL_{\lambda} (-,1)} := \und \scB$) coincide and secondly, the labelling functors coincide---the latter by the following arguments:
\begin{align}
\sU_{\sL_{\lambda}(-,0)} &= \sU_{\sL_{\lambda}} \sG(X \times  \delta^1_1) \\
&= \sU_{\scA} \sE^\lambda \sG(X \times  \delta^1_1) \\
&= \sU_{\scA} \id
\end{align}
and
\begin{align}
\sU_{\sL_{\lambda} (-,1)} &= \sU_{\sL_{\lambda}} \sG(X \times  \delta^1_0) \\
&= \sU_{\scA} \sE^\lambda \sG(X \times  \delta^1_0) \\
&= \sU_{\scA} \sS^\lambda \\
&= \sU_\scB
\end{align}
where we used \autoref{claim:grothendieck_span_construction_basechange} (and \eqref{eq:partial_functor_characterisation}) in the first step of both derivations, the definition of $\sL_\lambda$ in the second step, and further that $\sE^\lambda \sG(X \times  \delta^1_1) = \id$ and $\sE^\lambda \sG(X \times  \delta^1_0) = \sS^\lambda$ by \autoref{rmk:extended_collapse_props}. Finally, 
\begin{equation}
\cN_{\und {\sL_\lambda}} = \cN_{\cY_\lambda} = \lambda : \und {\sL_{\lambda}}(-,0) \into \und {\sL_{\lambda}}(-,1)
\end{equation}
by the isomorphism \eqref{eq:cartesian_closedness}. \qed
\end{lem}

\begin{rmk}[Base change of collapse limits] \label{rmk:precomposition_of_collapse_limits} Further to \autoref{rmk:precomposition_of_collapse} and \autoref{rmk:precomposition_extended_collapse} we compute
\begin{align}
\sL_{\lambda H} &= \sR_{\cY_{ \lambda H}, \sU_{\scA H} \sE^{ \lambda H}} \\
&=  \sR_{\cY_{\lambda}(H \times \bnum{2}), \sU_{\scA} \sG(H)\sE^{ \lambda H}} \\
&=  \sR_{\cY_{\lambda}(H \times \bnum{2}), \sU_{\scA} \sE^{\lambda} \sG(H\times \bnum{2})} \\
&= \sL_\lambda (H\times \bnum{2})
\end{align}
where we used \autoref{claim:grothendieck_span_construction_basechange} as well as the definition of $\sL_\lambda$ in \eqref{eq:collapse_lim_of_injection} in the last step.
\end{rmk}

\begin{claim}[Characterising collapse limits on single morphism] \label{claim:collapse_limit_characterisation_individual} Let $\scC : \bnum{2} \times \bnum{2} \to  \SIvertone \cC$ such that $\scC(i,-)$, $i \in \bnum{2}$, are collapse limits. Then $\scC$ is a collapse limit. 
\proof The proof is \stfwd{}. We need to show 
\begin{equation} \label{eq:coll_lim_loc1}
\sU_{\scC} = \sU_{\scC(-,0)} \sE^{\cN_{\und \cC}}
\end{equation}
Define 
\begin{equation}
\lambda := \cN_{\und \cC}
\end{equation}
that is, $\lambda$ has components $\lambda_i := \SIf \scC (i,0\to 1)$, $i \in \bnum 2$. To show \eqref{eq:coll_lim_loc1}, we evaluate both sides on an arbitrary $f \in \mor(\sG(\scC))$. We need to show
\begin{equation} \label{eq:wts_collapse_limits_condition}
\sU_{\scC}(f) = \sU_{\scC(-,0)} \sE^{\lambda}(f)
\end{equation}
Note that $f$ is of the form 
\begin{equation}
f = ((i,j),a) \to ((k,l),b)
\end{equation}
where $i,j,k,l \in \bnum 2$, $i \leq k$, $j \leq l$. We argue by case distinction. Since
\begin{align}
\sU_\scC \sG(X \times \delta^1_1) &= \sU_{\scC(-,0)} \\
&= \sU_{\scC(-,0)} \sE^{\lambda} \sG(X \times \delta^1_1)
\end{align}
(where we first used \autoref{claim:grothendieck_span_construction_basechange} and then \autoref{rmk:extended_collapse_props}), the statement of \eqref{eq:wts_collapse_limits_condition} follows whenever $j = l = 0$. Next, since each $\scC(i,-)$ individually is a collapse limit the statement of \eqref{eq:wts_collapse_limits_condition} also follows for $i = k$. The only morphisms $f$ left to consider are then of the form $f = ((0,1),a_0) \to ((1,1),a_1) \in \mor(\sG(\scC))$ or $f = ((0,0),a_0) \to ((1,1),a_1) \in \mor(\sG(\scC))$. 
\begin{itemize}
\item We first consider the former case of $f = ((0,1),a_0) \to ((1,1),a_1) \in \mor(\sG(\scC))$. For $i \in \Set{1,2}$ set 
\begin{equation}
(i, a'_i) := \sS^\lambda (i,a_i)
\end{equation}
such that by definition of $\sE^\lambda$ we have
\begin{equation}
\sE^\lambda(f) = (0,a'_0) \to (1,a'_1) \in \mor(\sG(\scC(-,0)))
\end{equation}
Since $\sS^\lambda$ (or equivalently $\sE^\lambda$) is functorial by \autoref{lem:glambda}, we have $(a_0',a_1') \in \edgeset(\scC((0,0) \to (1,0)))$. Since the graph of $\rest {\sS^\lambda} i$ is a subrelation of $R(\lambda_i)$ (cf. \autoref{lem:glambda}) we also have $(a_i,a'_i) \in \edgeset(\scC((0,0) \to (0,1)))$, $i \in \Set{0,1}$. Since $\scC(i,-)$ are collapse morphisms, and using the definition of $\sE^{\lambda_i}$, we find
\begin{equation} \label{eq:collapse_morphism_id}
\sU_\scC(((i,0),a'_i) \to ((i,1),a_i)) = \id
\end{equation}
Since $\scC$ is functorial we also have
\begin{equation}
 \scC((0,0) \to (0,1)) \odot \scC((0,1) \to (1,1)) = \scC((0,0)) \to (1,0)) \odot \scC((1,0) \to (1,1)) 
\end{equation}
This is an equation of horizontal compositions of natural transformations. The $(a'_0, a_1)$ component of this transformation can be derived in two ways as
\begin{align}
&\scC((0,1) \to (1,1))_{(a_0,a_1)} \circ \scC((0,0) \to (0,1))_{(a'_0, a_0)} \\ &= \scC((1,0) \to (1,1))_{(a'_1,a_1)} \circ \scC((0,0)) \to (1,0))_{(a'_0, a'_1)}
\end{align}
Using the construction of $\sU_\scC$ given in \autoref{defn:unpacking_and_repacking} and using \eqref{eq:collapse_morphism_id} this implies
\begin{equation}
\scC((0,1) \to (1,1))_{(a_0,a_1)} = \scC((0,0) \to (1,0))_{(a'_0, a'_1)}
\end{equation}
Using \autoref{defn:unpacking_and_repacking} again we in turn infer
\begin{equation}
\sU_\scC (((0,1),a_0) \to ((1,1),a_1)) = \sU_{\scC(-,0)} ((0,a'_0) \to (1,a'_1))
\end{equation}
This proves \eqref{eq:wts_collapse_limits_condition} when evaluated on $f$ as required.

\item In the second case, we consider morphisms of the form $f = ((0,0),a_0) \to ((1,1),a_1) \in \mor(\sG(\scC))$. Setting
\begin{equation}
(1,a'_1) := \sE^\lambda((1,1),a_1) = \sS^\lambda(1,a_1)
\end{equation}
Since $(0,a_0) = \sE^\lambda ((0,0),a_0)$, functoriality of $\sE^\lambda$ (as proven in \autoref{constr:uncurrying}) implies
\begin{equation}
\sE^\lambda(f) = (0,a_0) \to (1,a_1) \in \mor(\sG(\scC(-,0)))
\end{equation}
and as before we have $((1,0),a'_1) \to ((1,1),a_1) \in \mor(\sG(\scC))$ such that, since $\scC(1,0 \to 1)$ is a collapse morphism, we have
\begin{equation}
\sU_\scC ((1,0),a'_1) \to ((1,1),a_1) = \id
\end{equation}
Using functoriality of $\scC$ and the definition of horizontal composition as before, we then deduce that
\begin{align*}
\scC((0,1) \to (1,1))_{(a_0,a_1)}  &= \scC((1,0) \to (1,1))_{(a'_1,a_1)} \circ \scC((0,0) \to (1,0))_{(a_0,a'_1)} \\
&= \id \circ \scC((0,0) \to (1,0))_{(a_0,a'_1)}
\end{align*}
Using \autoref{defn:unpacking_and_repacking} again we in turn infer
\begin{equation}
\sU_\scC (((0,0),a_0) \to ((1,1),a_1)) = \sU_{\scC(-,0)} ((0,a_0) \to (1,a'_1))
\end{equation}
which proves \eqref{eq:wts_collapse_limits_condition} when evaluated on $f$ as required.
\end{itemize}
\qed
\end{claim}

\begin{thm} \label{claim:collapse_limit_characterisation_collective} Let $\scC : X \times \bnum{2} \to  \SIvertone \cC$. Then $\scC$ is a collapse limit if and only if it is a pointwise collapse limit: that is, $\scC(x,-)$ is a collapse limits for all $x \in X$.
\proof  The proof is \stfwd{}. Clearly, if $\scC$ is a collapse limit then for all $x \in X$, by precomposition (cf. \autoref{rmk:precomposition_of_collapse_limits}) with $\Delta_x \times \bnum{2}$, all $\scC(\Delta_x \times \bnum{2}) : \bnum{1} \times \bnum{2} \to \SIvertone \cC$ are collapse limits.

We now prove the converse. Assume $\scC(x,-)$ are collapse limits for all $x \in X$. We want to show 
\begin{equation}
\sU_{\scC} = \sU_{\scC(-,0)} \sE^{\lambda}
\end{equation}
where $\lambda = \cN_{\und \scC}$, that is, it has components $\lambda_x := \und \scC (i,0 \to 1)$, $x \in X$. Like in the previous proof, this equality can be shown by evaluating the left and right hand side on all morphisms
\begin{equation}
f = ((x,i),a_0) \to ((y,j),a_1) \in \mor(\sG(\scC))
\end{equation}
Let $F : \bnum{2} \to X$ map $(0 \to 1)$ to $(x \to y)$. Then we can restrict to the image of $F$ to complete the proof as follows. Setting $f' = ((0,i),a_0) \to ((1,j),a_1)$ (thus $f = \sG(F \times \bnum 2)(f')$)  we find
\begin{alignat}{3}
&& ~\sU_{\scC}(f) &= \sU_{\scC(-,0)} \sE^{\lambda}(f) \\
&\iff &~ \sU_{\scC}\big(\sG(F \times \bnum{2})(f')\big) &= \sU_{\scC(-,0)} \sE^{\lambda}\big(\sG(F \times \bnum{2})(f')\big) \\
&\iff &~ \sU_{\scC}\big(\sG(F \times \bnum{2})(f')\big) &= \sU_{\scC(-,0)} \sG(F) \sE^{ \lambda F}(f') \\
&\iff &~ \sU_{\scC(F \times \bnum{2})}(f') &= \sU_{\scC(-,0)F} \sE^{ \lambda F}(f')
\end{alignat}
Here, in the second step we used \autoref{rmk:precomposition_extended_collapse} and in the last step we used \autoref{claim:grothendieck_span_construction_basechange} (on both sides). Now, the last of the above equations follows from \autoref{claim:collapse_limit_characterisation_individual} applied to $\scC(F \times \bnum{2}) : \bnum{2} \times \bnum{2} \to \SIvertone \cC$. \qed 
\end{thm}

We record the results of the preceding discussion in the following theorem.

\begin{thm}[Collapse limits encode collapses] \label{thm:collapse_limits_represent_collapses} The construction
\begin{equation}
(\lambda : \scA \mcoll \scB) \mapsto \sL_\lambda
\end{equation}
is a bijection of collapse witnesses $\lambda$ of $\cC$-labelled $\SI$-families over $X$ with $\SI$-families over $X \times \bnum{2}$ that are pointwise collapse limits. Its inverse is given by
\begin{equation}
 \scC \mapsto \big(\cN_{\und \scC} : \scC(-,1) \mcoll \scC(-,0)\big)
\end{equation} 
\qed
\end{thm}

\subsection{Trivial product bundles}

\begin{constr}[Trivial product interval families] \label{rmk:product_bundles} Let $Z \neq \emptyset$ be a non-trivial poset. Let $X = Y \times Z$ be a product of posets, and denote by $\ppi_Y : Y \times Z \to Y$ the projection on the first component $Y$. An \SI-family $\scA : X \to \SI$ which factors as
\begin{equation}
\scA = \scB \ppi_Y
\end{equation}
for some $\scB : Y \to \SI$, is called a \textit{trivial product interval family} (of $\scB$ and $Z$). Note that this implies that for any $z_0 \in Z$ we have
\begin{equation}
\scB = \scA(-,z_0)
\end{equation}
There is an isomorphism
\begin{align}
\sW^\scB_Z : \sG(\scB) \times Z ~&\iso~ \sG(\scA) \\
((y,a), z) ~&\mapsto~ ((y,z),a) 
\end{align}
(The notation $\sW^\scB_Z$ indicates that $\scB$ is ``multiplied" with $Z$ to obtain $\scA$). We remark that these definitions imply
\begin{equation}
\pi_{\scA} \sW^\scB_Z = \pi_{\scB} \times Z
\end{equation}
\end{constr}

\begin{constr}[Trivial product towers of interval families] \label{constr:proj_towers} Let $\ppi_Y : (Y \times Z) \to Y$ be a projection of a product of posets. Let $T = \Set{ U^{n-2}, U^{n-2}, \dots, U^0}$ be a tower of \SI-families of height $n$ over $Y$. Then we construct a tower $T \ttimes Z = \Set{ V^{n-1}, V^{n-2}, \dots, V^0}$ of \SI-families of height $n$ over $X := (Y \times Z)$, together with isomorphisms of posets  (for $1 \leq k \leq n$)
\begin{equation}
\sW^{T,k}_Z :  \sG(U^{k-1})\times Z \iso \sG(V^{k-1})
\end{equation}
$T \ttimes Z$ is called the \textit{trivial product tower of $T$ and $Z$}, $T$ is the \textit{projected tower} of $T$ along $Z$. $\sW^{T,k}_Z $ is called the \textit{$k$-level trivial product isomorphism}. The construction of $V^{k}$ and $\sW^{T,k}_Z$ is inductive in $k = 0,1, \dots, (n-1)$. We inductively claim that
\begin{equation} \label{eq:tower_proj_ind_assumption}
U^k \ppi_{\sG(U^{k-1})} = V^k \sW^{T,k}_Z
\end{equation}
that is, $V^k \sW^{T,k}_Z$ is a trivial product \SI-family with $Z$ (as was defined in \autoref{rmk:product_bundles}).

For $k = 0$ we set $V^{0} = U^0 \ppi_Y : X \to \SI$ and $\sW^{T,0}_Z = \id_X$ which satisfies the inductive assumption. 

For $k > 0$, using \autoref{rmk:product_bundles} on the inductive assumption \eqref{eq:tower_proj_ind_assumption} (for $k-1$) we define
\begin{equation}
W^{T,k}_Z := \sW^{U^{k-1}}_{Z} ~:~ \sG(U^{k-1}) \times Z \iso \sG(V^{k-1} \sW^{T,k-1}_Z)
\end{equation}
Note that since $\sW^{T,k-1}_Z$ is an isomorphism, we obtain an isomorphism
\begin{equation}
\sG(\sW^{T,k-1}_Z) : \sG(V^{k-1} \sW^{T,k-1}_Z) \iso \sG(V^{k-1})
\end{equation}
We can then define
\begin{equation}
\sW^{T,k}_Z := \sG(\sW^{T,k-1}_Z) W^{T,k}_Z ~:~ \sG(U^{k-1}) \times Z \iso \sG(V^{k-1})
\end{equation}
and set
\begin{equation}
V^k := U^k \ppi_{\sG(U^{k-1})} (\sW^{T,k}_Z)\inv
\end{equation}
which completes the inductive construction.
\end{constr}

\begin{rmk}[Correspondence to pullback of towers] Comparing the previous construction to \autoref{constr:basechange_towers} one can inductively see that (for $k \geq 1$)
\begin{equation}
\ppi_{\sG(U^{k-1})} (\sW^{T,k}_Z)\inv = \sG^k(\ppi_Y)
\end{equation}
and thus in fact (using pullback notation for towers)
\begin{equation}
T \times Z = T\ppi_Y
\end{equation}
We will however usually use the notation on the left hand side.
\end{rmk}

\begin{claim}[Repacking trivial product interval families] \label{claim:repacking_prod_bun} Let $\ppi_Y : (Y \times Z) \to Y$ be a projection of a product of posets and $U : Y \to \SI$ an \SI-family. For $L : \sG(U) \to \cC$ and the projection $\ppi_{\sG(U)} : \sG(U) \times Z \to \sG(U)$ we claim that
\begin{equation}
\sR_{U\ppi_Y, L \ppi_{\sG(U)} (\sW^U_Z) \inv} = (\sR_{U, L}) \ppi_Y
\end{equation}
\proof  The proof is \stfwd{}. First note
\begin{align}
\sV_{(\sR_{U, L}) \ppi_Y} &= \und{(\sR_{U, L})} \ppi_Y \\
&= U \ppi_Y 
\end{align}
Then, for $((x,z),a) \to ((y,w),b)$ we compute
\begin{align}
\sU_{(\sR_{U, L}) \ppi_Y}((x,z),a) \to ((y,w),b) &= \SIs {(\sR_{U, L})} \ppi_Y ((x,z) \to (y,w))_{(a,b)} \\
&=\SIs {(\sR_{U, L})} (x \to y)_{(a,b)} \\
&=L((x,a) \to (y,b))
\end{align}
which equals $L \ppi_{\sG(U)} (\sW^\scA_Z) \inv ((x,z),a) \to ((y,w),b)$. Thus, by \autoref{claim:unpacking_and_repacking} the claim follows.
\qed 
\end{claim}

\begin{cor}[Repacking trivial product cube families] \label{claim:repacking_prod_cube_bun}  Let $\ppi_Y : (X := Y \times Z) \to Y$ be a projection of a product of posets and $\scB : Y \to \SIvert n \cC$. Then 
\begin{equation}
\tsR n_{\sT_\scB \ttimes Z, \tsU n_\scB \ppi_{\tsG n(\scB)} (\sW^{\sT_\scB,n}_Z)\inv} = \scB \ppi_Y
\end{equation}
\proof Inductively apply \autoref{claim:repacking_prod_bun}. \qed
\end{cor}

\subsection{Cube family collapse}

\begin{defn}[$n$-sequenced collapse limits] \label{defn:n_lvl_collapse_limit} The definition is inductive.
\begin{enumerate}
\item A $1$-sequenced collapse limit is a map $\scC : X \times \bnum{2} \to \SIvertone  \cC$ which is a collapse limit.
\item An $n$-sequenced collapse limit is a map $\scC : X \times (\bnum{n+1}) \to \SIvert n \cC$ such that the following two conditions are satisfied
\begin{itemize}
\item $\rest \scC {0 \to 1} : X \times \bnum{2} \to \SIvert n \cC = \SIvertone  {\SIvert {n-1} \cC}$ is a collapse limit.
\item $\scC (X \times \delta^n_0) : X \times \bnum{n} \to \SIvert n \cC$ satisfies
\begin{equation} \label{eq:k_lvl_limit_proj}
\und \scC (X \times \delta^n_0) = \und \scC(-,1) \ppi_X
\end{equation}
and
\begin{equation} \label{eq:k_lvl_limit_trunc}
\scC_{\bnum{n}} := \sU_{\scC(X \times \delta^n_0)} \sW^{\und \scC(-,1)}_{\bnum{n}} \quad : \quad \sG(\und{\scC}(-,1)) \times \bnum{n} ~\to~ \SIvert {n-1} \cC
\end{equation}
is an $(n-1)$-sequenced collapse limit.
\end{itemize}
\end{enumerate}
\end{defn}

\begin{rmk}[Truncations of $n$-sequenced collapse limits] \label{rmk:trunc_of_n_lvl_coll_lim} Given an  $n$-sequenced collapse limit $\scC : X \times (\bnum{n+1}) \to \SIvertone  n$ we first want to understand the definition of $\scC_{\bnum{n}}(-,k) : \sG(\und\scC(-,1)) \to \SIvert {n-1} \scC$ for $0 \leq k \leq (n-1)$. Recall from \autoref{rmk:product_bundles} that $\sG(\und\scC(-,1)) = \sG(\und\scC(-,k+1))$. We can then unwind the definition of $\scC_{\bnum{n}}(-,k)$ as follows: using \eqref{eq:partial_functor_characterisation} and \eqref{eq:k_lvl_limit_trunc} and $\scC_{\bnum{n}}(-,1)$ equals
\begin{align}
&\big(\sG(\und\scC(-,1)) \xto {(\sG(\und{\scC}(-,1)) \times \Delta_k)} \sG(\und\scC(-,1)) \times \bnum{n} \xto {\sW^{\und \scC (-,1)}_{\bnum{n}}} \sG(\und{\scC}(X \times \delta^n_0))\\
& \hspace{8.5cm} \xto {\tsU 1_{\scC (X \times \delta^n_0)}} \SIvert {n-1} \scC\big)\\
=~ &\big(\sG(\und\scC(-,1)) \iso \pi_{\und \scC(X \times \delta^n_0)}\inv(X \times \Set{k})  \xto {\tsU 1_{\scC(X \times \delta^n_0)}}  \SIvert {n-1} \scC\big)\\
=~ &\big(\sG(\und\scC(-,1)) \iso \pi_\scC\inv (X \times \Set{k+1}) \xto {\tsU 1_\scC} \SIvert {n-1} \scC\big) \\
=~ &\big(\sG(\und\scC(-,k+1)) \xto {\tsU 1_{\scC(-,k+1)}}  \SIvert {n-1} \scC\big)
\end{align}
Here, in the first step we computed the image of the first line, and in the second and third step we used \autoref{claim:grothendieck_span_construction_basechange} (together with \eqref{eq:partial_functor_characterisation} in the last step). In summary, we have shown
\begin{equation}
\scC_{\bnum{n}}(-,k) = \sU_{\scC(-,k+1)}
\end{equation}
In particular, this shows that we can infer the following representation of the domain of $(\scC_{\bnum{n}})_{\bnum{n-1}}$
\begin{equation}
(\scC_{\bnum{n}})_{\bnum{n-1}} : \tsG 2(\scC(-,2)) \times (\bnum{n-1}) \to \SIvert {n-2} \cC
\end{equation}
And inductively for $1 \leq k < n$ we find that
\begin{equation}
\big(\dots(\scC_{\bnum{n}})_{\bnum{n-1}} \dots\big)_{\mathbf{k}}  : \tsG k(\scC(-,n-k)) \times \mathbf{k} \to \SIvert {k} \cC
\end{equation}
We denote
\begin{equation}
\scC_{\mathbf{k+1}} := \big(\dots(\scC_{\bnum{n}})_{\bnum{n-1}} \dots\big)_{\mathbf{k+1}}
\end{equation}
This is a $k$-sequenced collapse limit called the $k$-level truncation of $\scC$.
\end{rmk}

\begin{thm}[$n$-sequenced collapse limits can be defined pointwise] \label{thm:pointwise_n_lvl_coll_lim} A family $\scC : X \times (\bnum{n+1}) \to \SIvert n \cC$ is a $n$-sequenced collapse limit if and only if it is pointwise an $n$-sequenced collapse limit: that is, for each $x \in X$, $\scC(x,-)$ is a $n$-sequenced collapse limit.
\proof The proof is \stfwd{}. It is by induction on $n$. The base case was proven in the previous section. Since both conditions in \autoref{defn:n_lvl_collapse_limit} for being an $n$-sequenced collapse limit can be restricted to subposets of $X$, the direction $\imp$ of the statement follows. Conversely, assume that for each $x \in X$, $\scC(x,-)$ is a $n$-sequenced collapse limit. We first want to show that \eqref{eq:k_lvl_limit_proj} holds, that is
\begin{equation} \label{eq:proj_condition}
\und \scC (X \times \delta^n_0) = \und \scC(-,1) \ppi_X
\end{equation}
Assume $(x \to y) \in \mor(X)$, and $(i \to j) \in \mor(\bnum{n})$. Then we find
\begin{align}
\und \scC (X \times \delta^n_0)((x,i) \to (y,j)) &= \und \scC ((x,i+1) \to (y,j+1)) \\
&=  \und \scC ((x,i+1) \to (y,j+1)) \und \scC ((x,1) \to (x,i+1)) \\
&=\und \scC ((x,1) \to (y,j+1))  \und \scC ((y,j+1) \to (y,1)) \\
&= \und \scC ((x,1) \to (y,1)) 
\end{align}
which proves \eqref{eq:proj_condition}. Here we used that both $\und \scC ((x,1) \to (x,i+1))  = \id$ and $\und \scC (y,j+1) \to (y,1)) = \id$ (since by assumption $\scC(x,-)$ and $\scC(y,-)$ are $n$-sequenced collapse limits) as well as functoriality of $\scC$. Next we need to verify that \eqref{eq:k_lvl_limit_trunc} is a $(n-1)$-sequenced collapse limit, that is 
\begin{equation}
\scC_{\bnum{n}} = \sU_{\scC(X \times \delta^n_0)} \sW^{\und \scC (-.1)}_{\bnum{n}} \quad : \quad \sG(\und{\scC}(-,1)) \times \bnum{n} ~\to~ \SIvert {n-1} \cC
\end{equation}
is an $(n-1)$-sequenced collapse limit. For each $x \in X$ the restriction of this map to $\sG(\und{\scC}(x,1)) \subset \sG(\und{\scC}(-,1))$ equals
\begin{equation}
\scC(x,-)_{\bnum{n}} = \sU_{\scC(\Delta_x \times \delta^n_0)} \sW^{\und \scC(-,1)}_{\bnum{n}} \quad : \quad \sG(\und{\scC}(x,1)) \times \bnum{n} ~\to~ \SIvert {n-1} \cC
\end{equation}
By inductive assumption, this map is an $(n-1)$-sequenced collapse limit. Thus, $\scC_{\bnum{n}}$ is a pointwise $(n-1)$-sequenced collapse limit. Using the statement of the theorem inductively we deduce that $\scC_{\bnum{n}}$ is an $(n-1)$-sequenced collapse limit. This verifies the conditions that $\scC$ is $n$-sequenced collapse limit and thus proves the theorem. \qed
\end{thm}

\begin{constr}[$n$-sequenced collapse limits from ordered collapse sequences] \label{constr:n_lvl_coll_lim_from_coll} Let $\scA : X \to \SIvert n \cC$. Let
\begin{equation}
 \xymatrix{\scA = \scA_{n+1} \ar@{~>}[r]^-{\lambda^n}_-n & \scA_{n} \ar@{~>}[r]^-{\lambda^{n-1}}_-{n-1} & \scA_{n-1} \dots \ar@{~>}[r]^-{\lambda^1}_-1 & \scA_1}
\end{equation}
be an ordered collapse sequence and denote this sequence by $\vvec\lambda$. By induction on $n$, we construct an $n$-sequenced collapse limit
\begin{equation}
\sL^n_{\vvec \lambda} : X \times (\bnum{n+1}) \to \SIvert n \cC
\end{equation}
such that (for $i \in (\bnum{n+1})$)
\begin{equation}
\sL^n_{\vvec \lambda}(-,i) = \scA_{i + 1}
\end{equation}
For $n = 1$, the construction is given by
\begin{equation}
\sL^1_{\vvec \lambda} = \sL_{\lambda^1}
\end{equation}
For $n > 1$ we proceed as follows. Note that
\begin{equation}
 \xymatrix{\tsU 1_{\scA_{n+1}} \ar@{~>}[r]^-{\lambda^n}_-{n-1} & \tsU 1_{\scA_{n}} \ar@{~>}[r]^-{\lambda^{n-1}}_-{n-2} & \tsU 1_{\scA_{n-1}} \dots \ar@{~>}[r]^-{\lambda^2}_-1 & \tsU 1_{\scA_1}}
\end{equation}
is an ordered $(n-1)$-level collapse sequence. Denote this sequence by $\vvec \mu$.  We inductively obtain
\begin{equation}
\sL^{n-1}_{\vvec \mu} : \tsG 1(\scA) \times \bnum{n} \to \SIvert {n-1} \cC
\end{equation}
and define 
\begin{equation}
L^{n-1} := \sR_{\und\scA\ppi_X,\sL^{n-1}_{\vvec \mu} (\sW^{\und\scA}_{\bnum{n}})\inv}
\end{equation}
where $\und\scA\ppi_X : X \times  \bnum{n} \to \SI$. Arguing inductively we know that (for $2 \leq i \leq n+1$)
\begin{equation}
\sL^{n-1}_{\vvec \mu}(-,i) = \tsU 1_{\scA_{i+1}}
\end{equation}
We thus find that
\begin{align}
L^{n-1}(-,i) &= L^{n-1}(X \times \Delta_i) \\
&= \sR_{\und\scA\ppi_X (X \times \Delta_i),\sL^{n-1}_{\vvec \mu} (\sW^{\und\scA}_{\bnum{n}})\inv \sG(X \times \Delta_i)} \\
&= \sR_{\und\scA,\sL^{n-1}_{\vvec \mu}(-,i)} \\
&=\scA_{i+1}
\end{align}
Now, on the other hand $\lambda^1 : \scA_2 \mcoll \scA_1$ yields
\begin{equation}
\sL_{\lambda^1} : X \times \bnum{2} \to \SIvertone  {\SIvert {n-1} \cC} = \SIvert n \cC
\end{equation}
such that $\sL_{\lambda^1}(-,1) = \scA_2$ and $\sL_{\lambda^1}(-,0) = \scA_1$. This means, $\sL_{\lambda^1}$ and $L^{n-1}$ agree on a poset $X$ when included as $X \iso X \times \Set{1}$ into $X \times \bnum{2}$ and as $X \iso X \times \Set{0}$ into $X \times \bnum{n}$. Thus, using the gluing construction from \autoref{constr:pushouts_in_labelled_posets} we obtain
\begin{equation}
\sL^n_{\vvec \lambda} = (\sL_{\lambda^1} \cup_X L^{n-1})\alpha\inv_X
\end{equation}
where we identified
\begin{equation}
\alpha_X \quad : \quad (X \times \bnum{2}) \cup_X (X \times \bnum{n}) \quad \iso \quad X \times (\bnum{n+1}) 
\end{equation}
by the map mapping $(x,i) \in X \times \bnum{2}$ to $(x,i) \in X \times (\bnum{n+1}) $ and $(x,i) \in X \times \bnum{n}$ to $(x,i+1) \in X \times (\bnum{n+1})$. By construction and inductive assumption, this satisfies the conditions for an $n$-sequenced collapse limit.
\end{constr}

\begin{thm}[$n$-sequenced collapse limits encode ordered collapse sequences]   \label{thm:n_lvl_collapse_limits_represent_collapses_seq} The mapping defined in  \autoref{constr:n_lvl_coll_lim_from_coll}
\begin{equation}
\big( \xymatrix{\scA = \scA_{n+1} \ar@{~>}[r]^-{\lambda^n}_-n & \scA_{n} \ar@{~>}[r]^-{\lambda^{n-1}}_-{n-1} & \scA_{n-1} \dots \ar@{~>}[r]^-{\lambda^1}_-1 & \scA_1}\big)\quad \mapsto \quad \big( \sL^n_{\vvec \lambda} : X \times (\bnum {n+1}) \to \SIvert n \cC\big)
\end{equation}
is a bijection of ordered $n$-level collapse sequences of $\cC$-labelled $\SI$-families over $X$ to $\SI$-families over $X \times (\bnum{n+1})$ that are pointwise $n$-sequenced collapse limits.

\proof  The proof is \stfwd{}. We construct an inverse to the above map. The construction is inductive in $n$. The base case ($n = 1$) follows by \autoref{thm:collapse_limits_represent_collapses} applied to the base case of \autoref{constr:n_lvl_coll_lim_from_coll}.

Let $n > 1$. First note, by \autoref{thm:pointwise_n_lvl_coll_lim} $\SI$-families over $X \times (\bnum{n+1})$ that are pointwise $n$-sequenced collapse limits are precisely $n$-sequenced collapse limits. Assume an $n$-sequenced collapse limit $\scC : X \times (\bnum{n+1}) \to \SIvert n \cC$. We construct an ordered collapse sequence $\vvec \cN_\scC$ consisting of collapses
\begin{equation} \label{eq:coll_seq_from_lim}
\xymatrix{\scC(-,n) \ar@{~>}[r]^-{\cN_\scC^n}_-n & \scC(-,n-1) \ar@{~>}[r]^-{\cN_\scC^{n-1}}_-{n-1} & \scC(-,n-2) \dots \ar@{~>}[r]^-{\cN_\scC^1}_-1 & \scC(-,0)}
\end{equation}
Arguing inductively, we find $\vvec \cN_{\scC_{\bnum{n}}}$ to yield a collapse sequence
\begin{equation}
\xymatrix{{\scC_{\bnum{n}}}(-,n-1) \ar@{~>}[r]^-{\cN_{\scC_{\bnum{n}}}^{n-1}}_-{n-1} & {\scC_{\bnum{n}}}(-,n-2) \ar@{~>}[r]^-{\cN_{\scC_{\bnum{n}}}^{n-2}}_-{n-2} & {\scC_{\bnum{n}}}(-,n-3) \dots \ar@{~>}[r]^-{\cN_{\scC_{\bnum{n}}}^1}_-1 & {\scC_{\bnum{n}}}(-,0)}
\end{equation}
We set (for $0 \leq k < n$)
\begin{equation}
\cN^{k+1}_\scC := \cN^k_{\scC_{\bnum{n}}}
\end{equation}
and using \autoref{rmk:trunc_of_n_lvl_coll_lim} we find a sequence
\begin{equation}
\xymatrix{\sU_{\scC(-,n)} \ar@{~>}[r]^-{\cN_\scC^n}_-n & \sU_{\scC(-,n-1)} \ar@{~>}[r]^-{\cN_\scC^{n-1}}_-{n-1} & \sU_{\scC(-,n-2)}\dots \ar@{~>}[r]^-{\cN_\scC^1}_-1 & \sU_{\scC(-,1)}}
\end{equation}
Finally, by the first condition in \autoref{defn:n_lvl_collapse_limit} we find using \autoref{constr:collapse_limits} 
\begin{equation}
\cN^1_\scC := \cN_{\und{\rest \scC {0 \to 1}}} : \scC(-,1) \mcoll \scC(-,0)
\end{equation}
This constructs a collapse sequence \eqref{eq:coll_seq_from_lim} as required. The claim that this construction is left and right inverse to \autoref{constr:n_lvl_coll_lim_from_coll} now follows from the inductive structure of \autoref{constr:n_lvl_coll_lim_from_coll} and inductive assumption of the present construction. \qed
\end{thm}

\subsection{Normalisation on restrictions}

We now turn to applications of the theory developed in the previous sections.

\begin{lem}[Normalisation on restriction of interval families] \label{lem:normalisation_on_restrictions} Let $\scA : X \to \SIvertone  \cC$ for a poset $X$ and a category $\cC$. Let $H: Y \into X$ be a \gls{downwardclosed} inclusion of posets. Then $\scA H$ is normalised if $\scA$ is normalised.

\proof Note that by assumption $H$ restricts to an isomorphisms $Y \iso \im(H)$ where the image of $H$, $\im(H) \subset X$ is a \gls{downwardclosed} subposet of $X$. Arguing by contradiction, assume $\scA$ is normalised, but $\scA H$ is not normalised. In particular we have a non-identity $\lambda : \scA H \mcoll \scB$ where $\scB :  Y \to \SIvertone  \cC$. Using \autoref{constr:collapse_limits} we find $\sL_\lambda : Y \times \bnum{2} \to \SIvertone  \cC$ and denote its collapse morphisms associated to each $\lambda_y$, $y \in Y$, by
\begin{equation}
f_y := \sL_\lambda (y,0\to 1)
\end{equation}
We will now construct a non-identity collapse limit
\begin{equation}
\scC : X \times \bnum{2} \to \SIvertone  \cC
\end{equation}
such that $\scC(-,1) = \scA$ (and also $\rest {\scC(-,0)} {\im(H)} = \scB H\inv$) which contradicts the assumption that $\scA$ is normalised. 

Since $\scC$ is a two-variable functor, we can define it by defining the family of partial functors $\scC(-,i)$ ($i \in \obj(\bnum{2}) = \Set{0,1}$), $\scC(x,-)$ ($x \in X$) and show their commutativity (cf. \autoref{notn:two_variable_functors}). Explicitly we define
\begin{itemize}
\item $\scC(-,1) := \scA$
\item $\scC(-,0)$ is defined by setting $\rest {\scC(-,0)} {\im(H)} = \scB$ and $\rest {\scC(-,0)} {X \setminus \im(H)} = \scA$. This leaves us with defining $\scC(-,0)$ on morphisms $(y \to x)$ for $y \in \im(H)$ and $x \in X \setminus \im(H)$, for which we set
\begin{equation}
\scC(x \to y,0) = \scA(x \to y) f_y
\end{equation}
\item $\scC(x,-)$ ($x \in X$) are defined, if $x \in \im(H)$, by
\begin{equation}
\scC(x, 0 \to 1) = f_x
\end{equation}
and otherwise, if $x \in X \setminus \im(H)$, by
\begin{equation}
\scC(x, 0 \to 1) = \id
\end{equation}
\end{itemize}
Note that $\scC(x,-)$ and $\scC(-,1)$ are functorial (the latter since $\scA$ is functorial). We check functoriality of $\scC(-,0)$ as follows: let $x_1 \to x_2 \to x_3$ be a chain of two non-identity morphisms in the poset $X$, that is, such that $x_l \neq x_k$ for $l < k \in \Set{1,2,3}$. We want to show
\begin{equation}
\scC(x_1 \to x_2 \to x_3,0) = \scC(x_2 \to x_3,0) \circ \scC(x_1 \to x_2,0)
\end{equation}
This holds if for all $k$ we have $x_k \in \im(H)$ or $x_k \in X \setminus \im(H)$ by definition of $\scC(-,0)$. We further distinguish the following remaining two cases
\begin{itemize}
\item $x_1 \in \im(H)$, $x_2,x_3 \notin \im(H)$ then by definition of $\scC(-,0)$ we have
\begin{align} 
\scC(x_1 \to x_2 \to x_3,0) &=  \scA(x_1 \to  x_3)\circ f \\
&= \scA(x_2 \to  x_3) \circ \scA(x_1 \to  x_2) \circ f \\
&= \scC(x_2 \to x_3,0) \circ \scC(x_1 \to x_2,0) 
\end{align}
\item $x_1, x_2 \in \im(H)$, $x_3 \notin \im(H)$ then by definition of $\scC(-,0)$ we have
\begin{align}
\scC(x_1 \to x_2 \to x_3,0) &=  \scA(x_1 \to  x_3) \circ f_{x_1} \\
&= \scA (x_2 \to  x_3) \circ \scA(x_1 \to  x_2) \circ f_{x_1} \\
&= \scA (x_2 \to  x_3) \circ f_{x_2} \circ \scB(x_1 \to  x_2) \\
&= \scC(x_2 \to x_3,0) \scC(x_1 \to x_2,0) 
\end{align}
where we used $\scA(x_1 \to  x_2) \circ f_{x_1} = f_{x_2} \circ \scB(x_1 \to  x_2)$ which follows from $\sL_\lambda$ being functorial.
\end{itemize}
To show that these partial functors define $\scC$ we now need to verify their commutativity, which amounts to checking
\begin{equation} \label{eq:downward_closed_comm}
\xymatrix@C=2cm{\scC(x_1,0) \ar[r]^{\scC(x_1 \to x_2,0)} \ar[d]_{\scC(x_1, 0 \to 1)} & \scC(x_2,0) \ar[d]^{\scC(x_2, 0 \to 1)} \\ \scC(x_1,1) \ar[r]_{\scC(x_1 \to x_2,1)} & \scC(x_2,1) }
\end{equation}
for $x_1, x_2 \in X$. This holds if $x_1 = x_2$ so once more we can assume $x_1 \neq x_2$. If $x_1, x_2 \in \im(H)$ then \eqref{eq:downward_closed_comm} holds since it then is equivalent to
\begin{equation}
\scA(x_1 \to  x_2) \circ f_{x_1} = f_{x_2} \circ \scB(x_1 \to  x_2)
\end{equation}
(which is true by functoriality of $\sL_\lambda$). Further, if $x_1, x_2 \notin Y$ then \eqref{eq:downward_closed_comm} just witnesses $\scA(x_1 \to  x_2)  = \scA(x_1 \to  x_2)$ (and the left and right sides are identities). In the remaining case of $x_1 \in \im(H)$ and $x_2 \notin \im(H)$ we find \eqref{eq:downward_closed_comm} to equal
\begin{equation} 
\xymatrix@C=2cm{\scB(x_1) \ar[r]^{\scA(x_1 \to x_2) f_{x_1}} \ar[d]_{f_{x_1}} & \scA(x_2) \ar[d]^{\id} \\ \scA(x_1) \ar[r]_{\scA(x_1 \to x_2)} & \scA(x_2) }
\end{equation}
which also holds. This completes the construction of $\scC$ and thus the proof.
\qed
\end{lem}

\begin{thm}[Normalisation on restrictions of n-cube families] \label{thm:normalisation_on_dc_restrictions} Let $\scA : X \to \SIvert n \cC$ for a poset $X$ and category $\cC$. Let $H : Y \to X$ be a \gls{downwardclosed} subposet inclusion. Then $\scA H$ is normalised (up to level $n$) if $\scA$ is normalised (up to level $n$).
\proof Using \autoref{eg:subfamily_by_restriction} we have inclusions $\restemb^k_H = \tsG k(H) : \tsG k(\scA H) \to \tsG k(\scA)$. By \autoref{claim:pullback_preserves_having_lifts} we know that these inclusions are \gls{downwardclosed}. Further, from \autoref{constr:unpacking_collapse} we obtain an equality
\begin{equation} \label{eq:restriction_colorings}
\tsU k_\scA \restemb^k_H = \tsU k_{\scA H}
\end{equation}
Now $\scA$ is normalised up to level $n$ if and only if all $\tsU k_\scA$ are normalised (as labelled singular interval families). By \autoref{lem:normalisation_on_restrictions} we deduce that  $\tsU k_\scA \restemb^k_H$ is normalised. By \eqref{eq:restriction_colorings} we then infer that $ \scA H$ is normalised up to level $n$ as claimed. \qed
\end{thm}

The condition of \gls{downwardclosed}ness cannot be dropped as the following example shows.

\begin{eg}[Non-normalised subfamily in normalised parent family] Let $\cC$ be the groupoid with two objects $\Set{0,1}$ and two non-identity morphisms $f : 0 \to 1$ and $f\inv : 1 \to 0$. Define $\scA : \singint 1 \to \SIvert 1 \cC$ by
\begin{restoretext}
\begingroup\sbox0{\includegraphics{test/page1.png}}\includegraphics[clip,trim=0 {.25\ht0} 0 {.15\ht0} ,width=\textwidth]{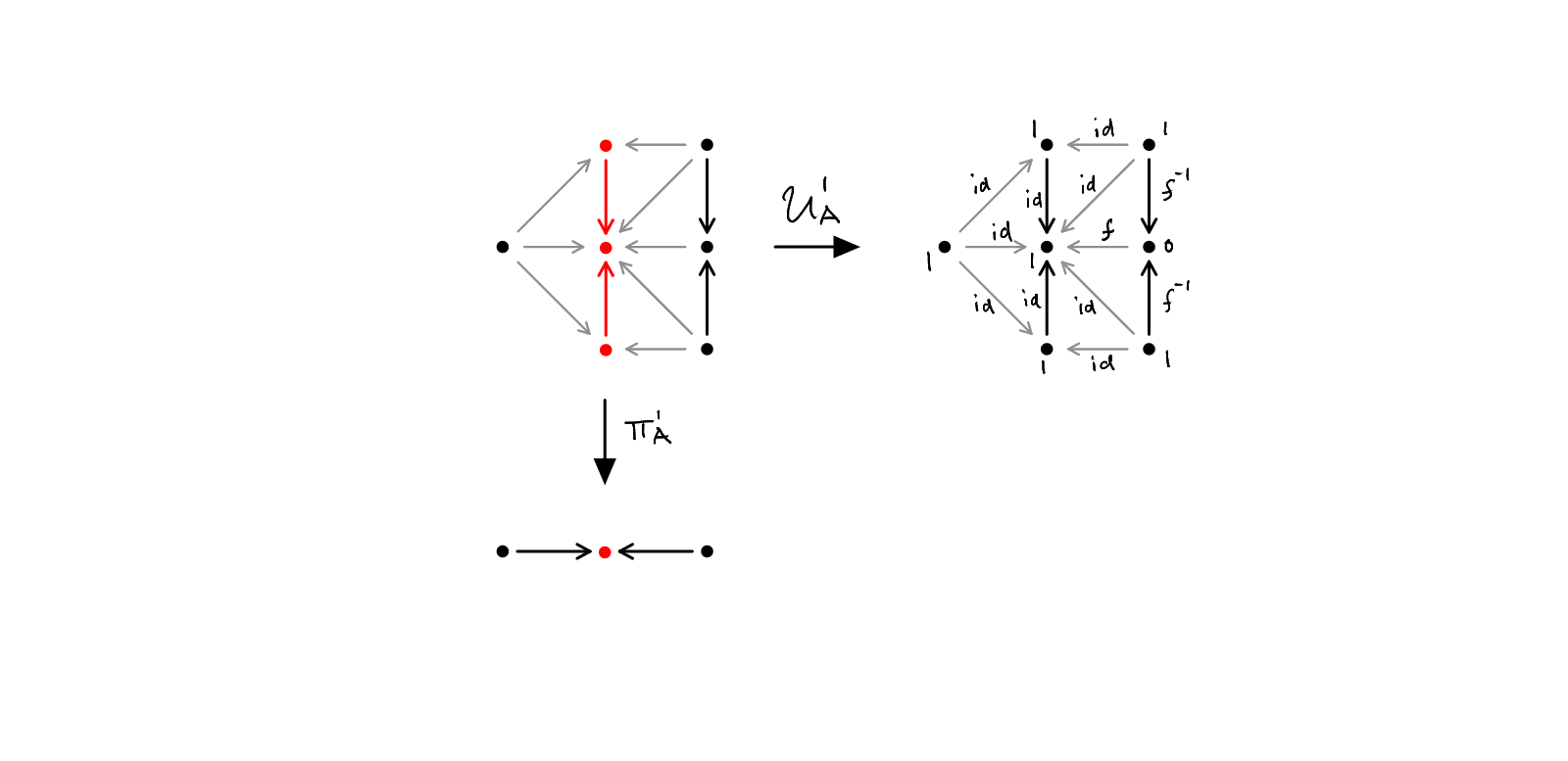}
\endgroup\end{restoretext}
Then $\scA$ is in normal form, but its subfamily marked in \cred{} is not. Note that this does not contradict the preceding theorem since the poset inclusion $\Set{1} \subset \singint 1$ is not \gls{downwardclosed}. 
\end{eg}

\subsection{Locally trivial cubes}

\begin{defn}[Locally and globally trivial families] \label{defn:loc_triv_fam} Let $\scA : X \to \SIvert n \cC$. We say $\scA$ is \textit{locally trivial} if for any $(x \to y) \in \mor(X)$, we have
\begin{equation}
\NF{\rest \scA {x \to y}}^n = \const
\end{equation}
We say $\scA$ is \textit{globally trivial} if it is the constant functor on all connected components of $X$.
\end{defn}

\begin{eg}[Locally and globally trivial families] \label{eg:locally_trivial} \hfill
\begin{enumerate}
\item Let $\cC$ be the terminal category $\bnum{1}$, and define the $\SIvert 1 \cC$ family $\scA$ by
\begin{restoretext}
\begingroup\sbox0{\includegraphics{test/page1.png}}\includegraphics[clip,trim=0 {.15\ht0} 0 {.2\ht0} ,width=\textwidth]{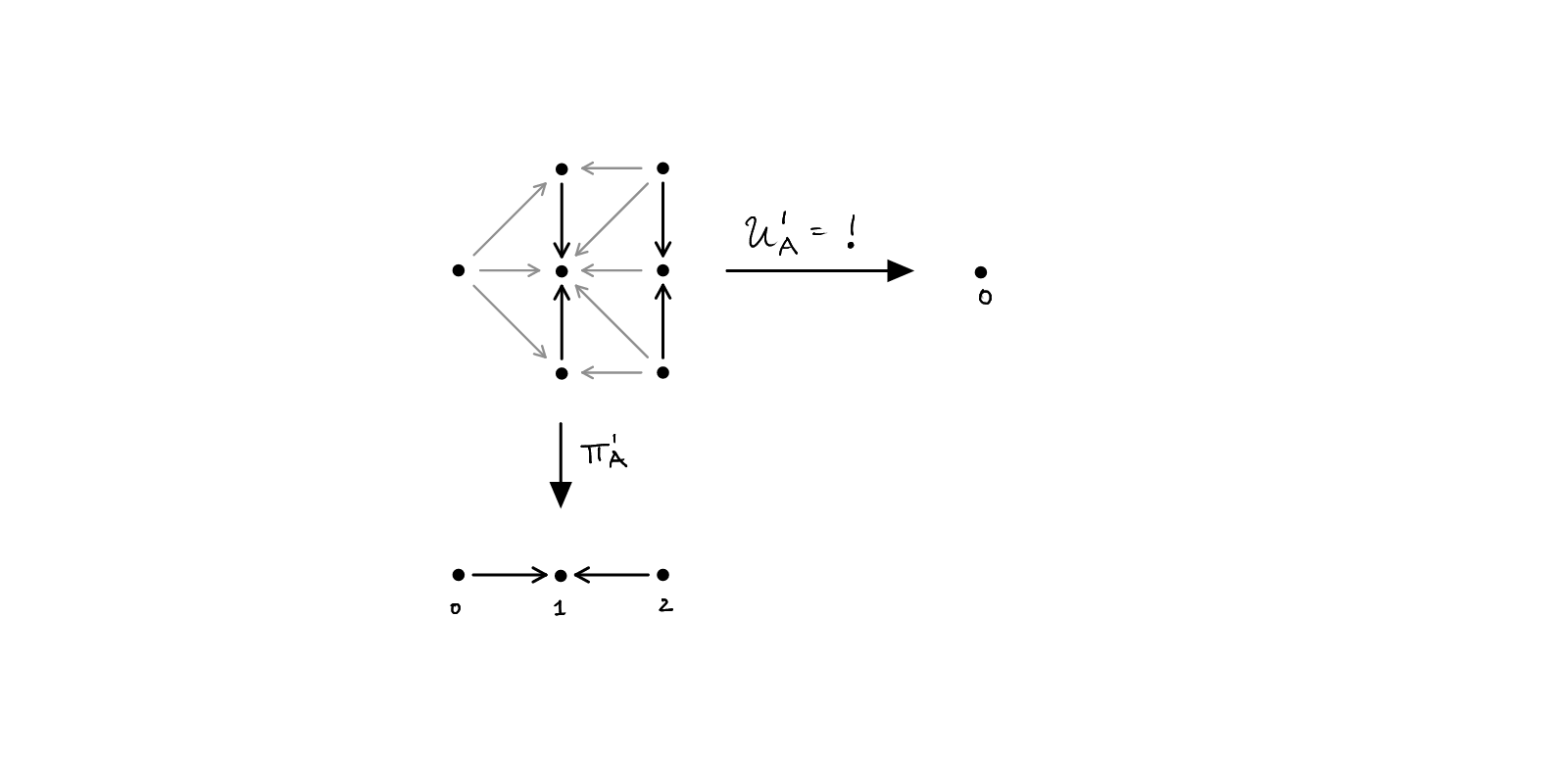}
\endgroup\end{restoretext}
Here, recall that $\bang$ is the unique functor to terminal category. This family is locally trivial. Indeed we verify that $\rest \scA {0 \to 1}$ normalises to the constant family, namely to
\begin{equation}
\NF{\rest \scA {0 \to 1}}^1 = \const_{\tsR 1_{\Delta_{\singint 0}, \bang}}
\end{equation}
which can be visualised 
\begin{restoretext}
\begingroup\sbox0{\includegraphics{test/page1.png}}\includegraphics[clip,trim=0 {.2\ht0} 0 {.15\ht0} ,width=\textwidth]{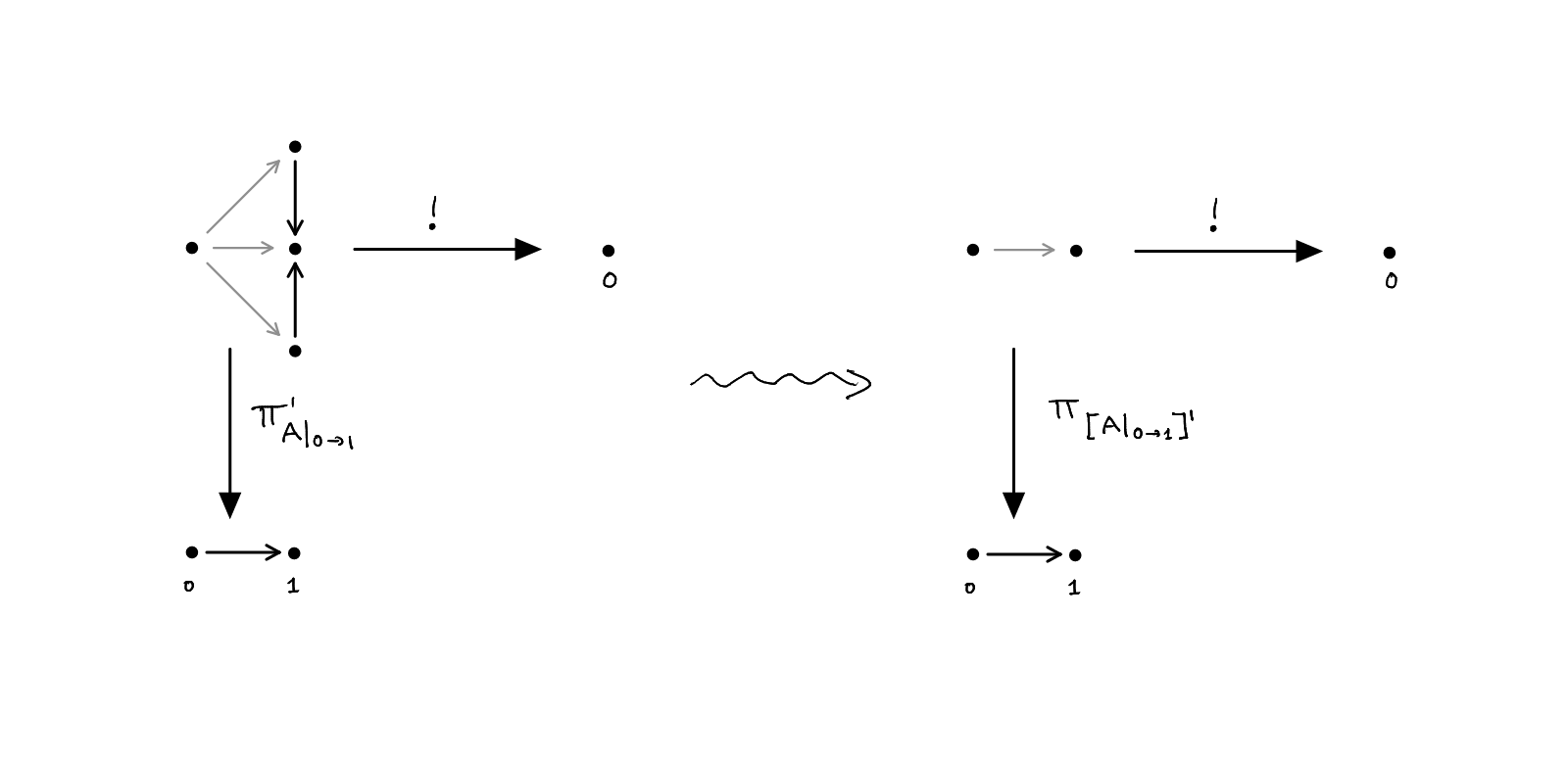}
\endgroup\end{restoretext}
And further, we verify that $\rest \scA {2 \to 1}$ normalises to the constant family, namely to
\begin{equation}
\NF{\rest \scA {2 \to 1}}^1 = \const_{\tsR 1_{\Delta_{\singint 0}, \bang}}
\end{equation}
which can be visualised 
\begin{restoretext}
\begingroup\sbox0{\includegraphics{test/page1.png}}\includegraphics[clip,trim=0 {.2\ht0} 0 {.1\ht0} ,width=\textwidth]{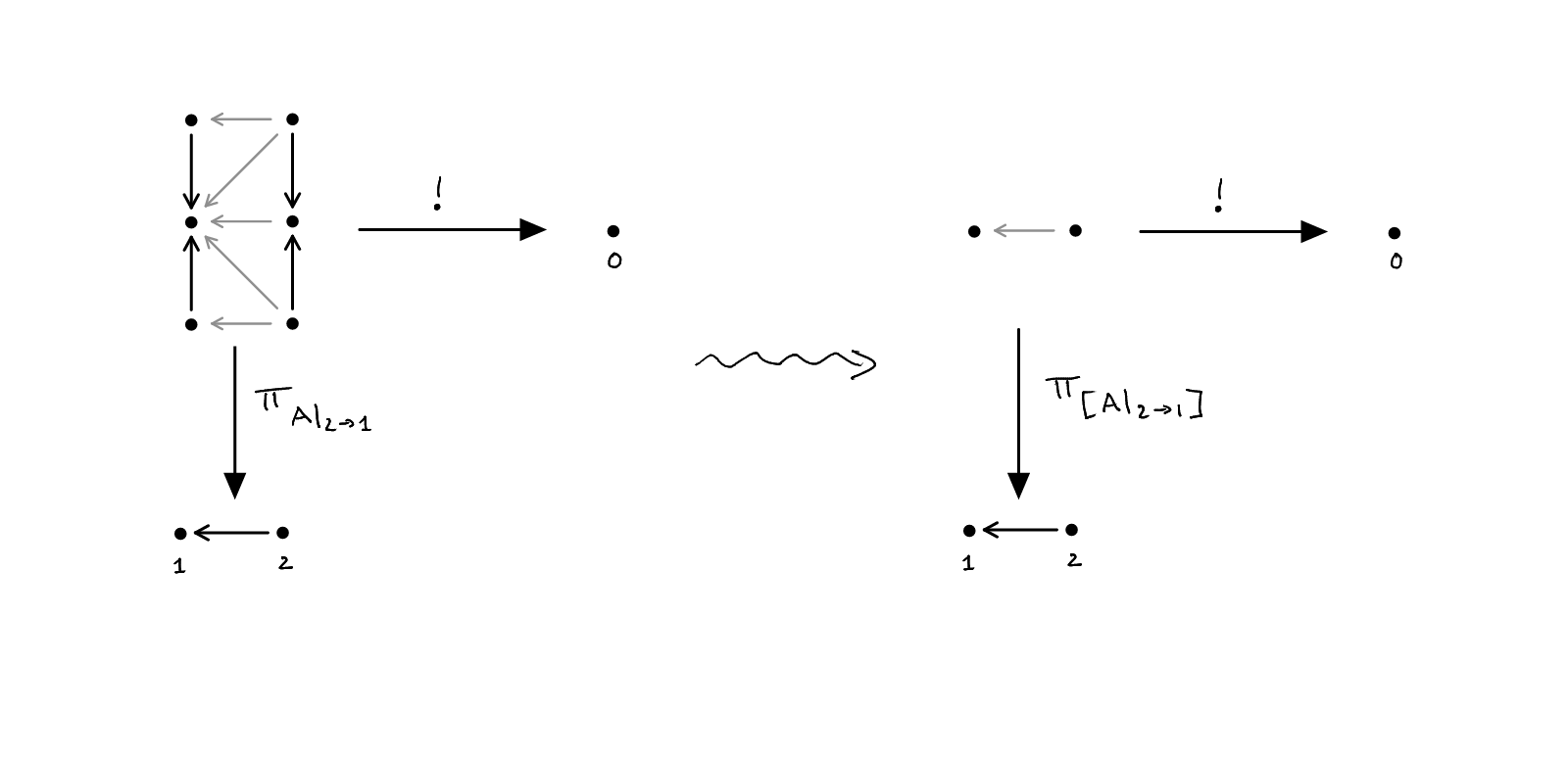}
\endgroup\end{restoretext}
Note that in this second case, $\rest \scA {2 \to 1}$ is in fact globally trivial. In this case  \autoref{claim:norm_pres_glob_triv} below applies, showing that $\rest \scA {2 \to 1}$ will normalise to the constant family. In particular, we do not need to check $\rest \scA {i \to i}$ normalises to the constant family either, and have thus verified $\scA$ to be locally trivial.

\item An example of a globally trivial family is $\scB : \singint 1 \to \SIvert 1 {\bnum{1}}$ is
\begin{equation}
\scB := \const_{\tsR 1_{\Delta_{\singint 1}, \bang}}
\end{equation}
which can be visualised as
\begin{restoretext}
\begingroup\sbox0{\includegraphics{test/page1.png}}\includegraphics[clip,trim=0 {.15\ht0} 0 {.15\ht0} ,width=\textwidth]{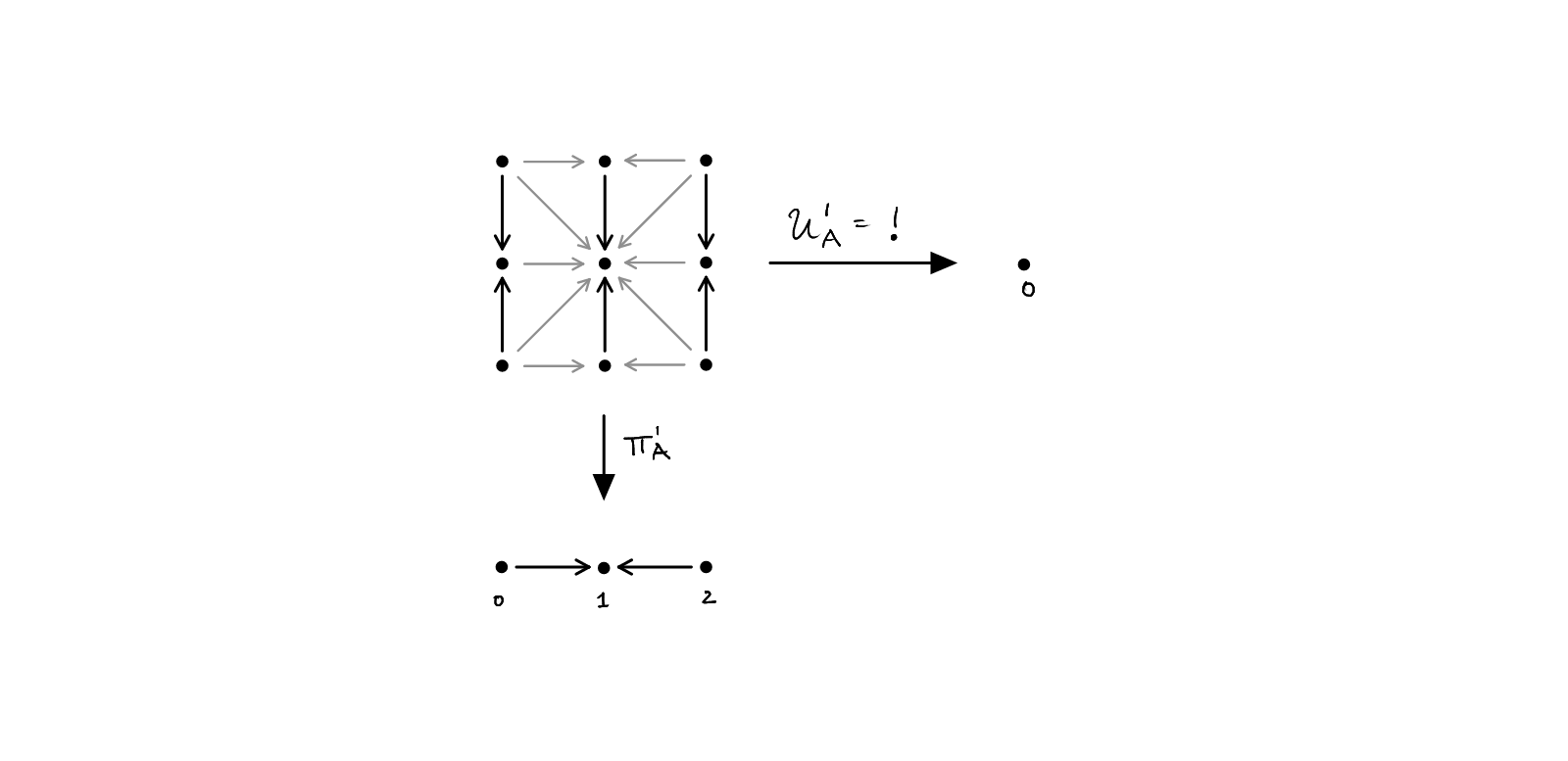}
\endgroup\end{restoretext}
The reader can verify that this normalised to a globally trivial family, as guaranteed more generally by the next claim.
\end{enumerate} 
\end{eg}

\begin{claim}[Normalisation preserves global triviality] \label{claim:norm_pres_glob_triv} If $\scA : X \to \SIvert n \cC$ is a globally trivial family, then $\scA$ normalises to a globally trivial family.
\proof We can assume $X$ to be connected (the argument holds for each connected component individually). Let $x_0 \in X$ and let $\bang : X \to \bnum{1}$ be the terminal functor. Then global triviality of $\scA$ means
\begin{equation}
\scA = \scA \Delta_{x_0} \bang
\end{equation}
Let $\scC : (\bnum{n+1}) \to \SIvert n \cC$ be the collapse limit of $\scA \Delta_{x_0}$ to its normal form $\NF{\scA \Delta_{x_0}}$. Then $\scC (\bang \times (\bnum{n+1}))$ is a collapse limit collapsing $\scA$ to
\begin{equation}
\scA' = \NF{\scA \Delta_{x_0}}\bang
\end{equation}
$\scA'$ is in normal form since any collapse $\lambda : \scA' \kcoll k \scA''$ will have a non-trivial restriction $\restemb\pbstar _x \lambda : \rest {\scA'} x \kcoll k \rest {\scA''} x$ for at least one $x \in X$ (cf. \autoref{constr:collapse_on_subfamilies}), which then contradicts $\NF{\scA \Delta_{x_0}} = \rest {\scA'} x$ being in normal form. \qed
\end{claim}

\begin{defn}[Locally normalised families] Let $\scA : X \to \SIvert n\cC$. We say $\scA$ is \textit{locally normalised} if for any $x \subset X$, we have
\begin{equation}
\NF{\rest \scA {x}}^n = \rest \scA {x}
\end{equation}
\end{defn}

\begin{lem}[Locally normalised implies normalised] Let $\scA : X \to \SIvertone \cC$. If $\scA$ is locally normalised, then it is normalised.
\proof Assume by contradiction a non-identity $\lambda : \scA \kcoll k \scB$ for some $0 < k \leq n$. Then there is $x \in X$ such that for some $y \in \tsG k(\rest \scA x) \subset \sG(\scA)$ we have $\lambda_x \neq \id$. Using \autoref{constr:collapse_on_subfamilies}, we find an induced collapse $\restemb_x \pbstar  \lambda : \rest \scA x \to \scA'$ on the embedding $\restemb_x : \rest \scA x \subset \scA$, and by choice of $x$ (cf. \eqref{eq:subfamily_and_collapse}) this collapse is non-trivial. This contradicts $\scA$ being locally normalised. \qed
\end{lem}

\begin{defn}[Downward distance] Let $X$ be a finite poset. On objects $x \in X$ we define the \textit{downward distance} $\mathrm{ddist}(x) = m \in \lN$ if $m$ is the maximum length of a chain 
\begin{equation}
x_1 \to x_2 \to ... \to x_{m-1} \to x_m = x
\end{equation}
with no identity arrows (that is, $x_i \neq x_{i+1}$). Define for $k \geq 1$ the full subposets $X_{\mathrm{dd}\leq k} \subset X$ to have objects
\begin{equation}
\obj(X_{\mathrm{dd}\leq k}) = \Set{ x\in X ~|~\mathrm{ddist}(x) \leq k }
\end{equation}
Note that this is \gls{downwardclosed} by the following remark.
\end{defn}

\begin{eg}[Downward distance] In the following poset, the downward distance of each point is indicated by a number next to it
\begin{restoretext}
\begingroup\sbox0{\includegraphics{test/page1.png}}\includegraphics[clip,trim={.1\ht0} {.3\ht0} {.1\ht0} {.25\ht0} ,width=\textwidth]{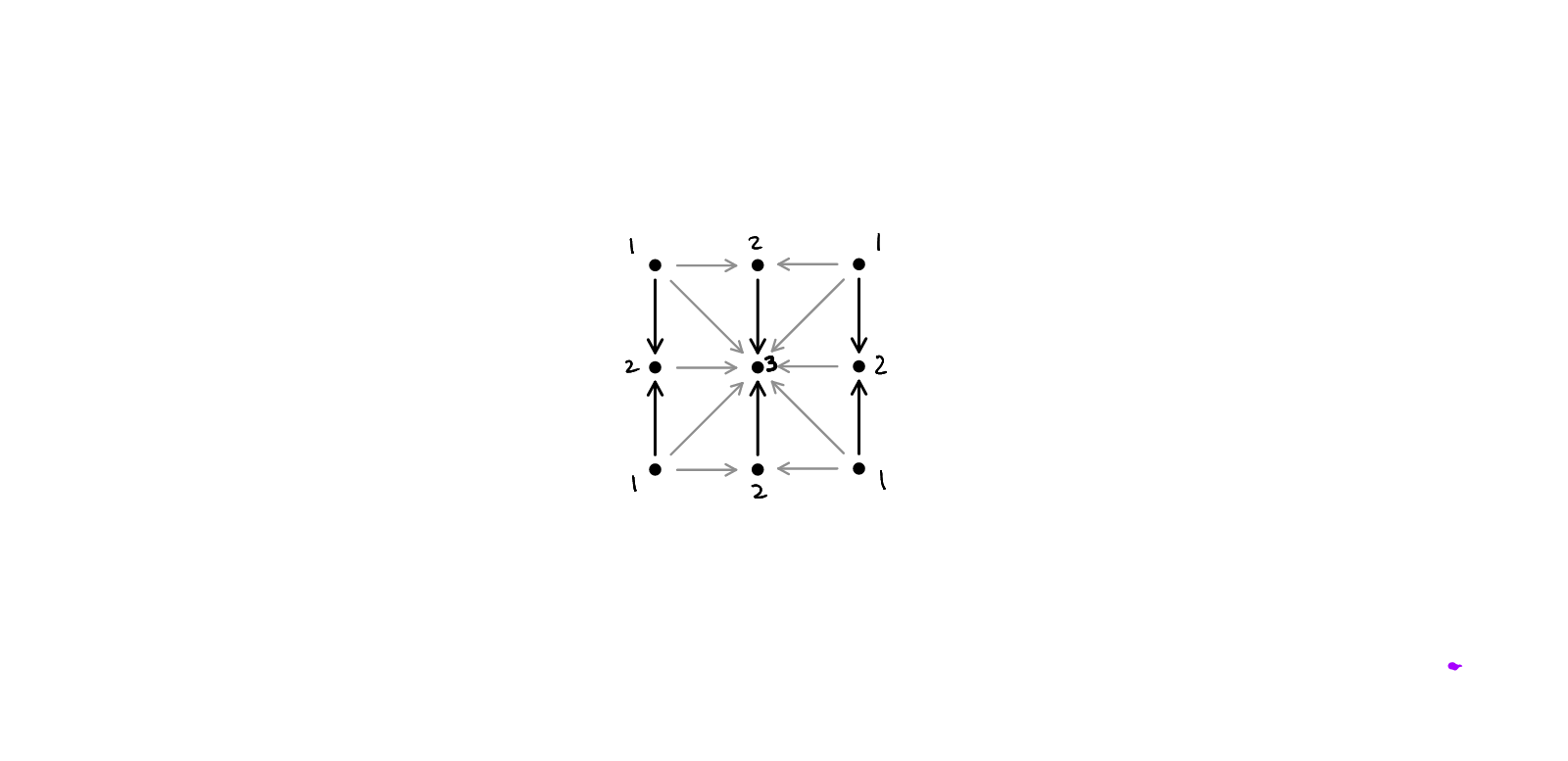}
\endgroup\end{restoretext}
\end{eg}

\begin{rmk}[Inductive property of downward distance] \label{rmk:ddist_prop} We observe that if $\mathrm{ddist}(y) = m$ and $y' \to y$ in $X$, then $\mathrm{ddist}(y) < m$.
\end{rmk}

We now prove that there is no normalised and locally trivial family which is not also globally trivial. For instance, $\scA$ from \autoref{eg:locally_trivial} is locally trivial but not globally trivial, so it cannot be normalised, which is indeed the case.

\begin{thm}[Normalised and locally trivial implies locally normalised and globally trivial] \label{lem:local_trivial_to_global_trivial} Let $\scA : X \to \SIvert n \cC$. If $\scA$ is normalised and locally trivial, then it is globally trivial and locally normalised.
\proof We inductively prove that $Y_k := \rest \scA {X_{\mathrm{dd}\leq k}}$ is globally trivial. Note that $\rest \scA {Y_1}$ is globally trivial (by \autoref{lem:normalisation_on_restrictions}) and locally normalised (since $Y_1$ is a discrete poset). Now assume that the statement is proven for $\rest \scA {Y_{k-1}}$. First note, by \autoref{lem:normalisation_on_restrictions} we infer that $\rest \scA {Y_k}$ is in normal form. Note that if it is globally trivial then it must also be locally normalised since any connected component of $Y_k$ contains points of $Y_{k-1}$ at which $\scA$ is normalised by inductive assumption. Thus, arguing by contradiction assume that $\rest \scA {Y_k}$ is not globally trivial.

If $\rest \scA {Y_k}$ is not globally trivial we can find $(y_0 \to y_1) \in \mor(Y_k)$ such that
\begin{equation}
\rest \scA {y_0 \to y_1} \neq \id
\end{equation}
Note that we must have $y_0 \in Y_{k-1}$ (cf. \autoref{rmk:ddist_prop}). Define $Y \subset X$ to be the full (and \gls{downwardclosed}) subset of $X$ with objects
\begin{equation}
\obj(Y) = \obj(Y_{k-1}) \cup \Set{y_1}
\end{equation}
By \autoref{lem:normalisation_on_restrictions} we know that $\rest \scA Y$ is normalised. 

Now, assume $(x \to y_1) \in \mor(Y)$ and define
\begin{gather}
{\scB_x} := \rest \scA {x \to y_1}
\end{gather}
Let $\vec\nfc_{\scB_x}$ be the ordered collapse sequence collapsing ${\scB_x}$ into normal form $\NF{\scB_x}^n$ (cf. \autoref{thm:normal_forms_unique}). By functoriality of $\sL^n_{\vec \nfc_{\scB_x}}$ we find for any $x \in Y_{k-1}$ that
\begin{equation}
\xymatrix@C=3cm{ 
\NF{\scB_x}^n(x) 
		\ar[r]^{\sL^n_{\vec \nfc_{\scB_x}}(x \to y_1, 0)} \ar[d]_{\sL^n_{\vec \nfc_{\scB_x}}(x, 0 \to n)} 
& \NF{\scB_x}^n(y_1) 
		\ar[d]^{\sL^n_{\vec \nfc_{\scB_x}}(y_1, 0 \to n)} 
\\ 
{\scB_x}(x)  
		\ar[r]_{\sL^n_{\vec \nfc_{\scB_x}}(x \to y_1,n)} 
& {\scB_x}(y_1) }
\end{equation}
Now we must have $\sL^n_{\vec \nfc_{\scB_x}}(x \to y_1, 0) = \id$ since ${\scB_x} = \rest \scA {x \to y_1}$ normalises to the identity by assumption on $\scA$. But we also must have $\sL^n_{\vec\nfc_{\scB_x}}(x, 0 \to n) = \id$ since $\rest {\scB_x} x = \rest \scA x$ is already normalised by inductive assumption, and thus now non-trivial ordered collapse chain applies to it. We deduce that
\begin{align}
\sL^n_{\vec\nfc_{\scB_x}}(y_1, 0 \to n) &= \sL^n_{\vec\nfc_{\scB_x}}(x \to y_1,n) \\
&= {\scB_x}(x \to y_1)
\end{align}
where, in the last step we used the properties of the $\sL$ construction (cf. \autoref{constr:n_lvl_coll_lim_from_coll}). Since ${\scB_x}(x)$ is in normal form we deduce that $\sL^n_{\vec\nfc_{\scB_x}}(y_1,0)$ is in normal form. Thus by  \autoref{thm:normal_forms_unique} and \autoref{thm:n_lvl_collapse_limits_represent_collapses_seq}
\begin{equation}
\vvec \cN_{\sL^n_{\vec\nfc_{\scB_x}}(y_1,-)} = \vvec\nfc_{\rest \scA {y_1}}
\end{equation}
is the unique ordered collapse sequence to normal form for $\rest \scA {y_1}$ with $k$-level normal forms 
\begin{equation}
\sL^n_{\vec\nfc_{\scB_x}}(y_1,k) = \NF{\rest \scA {y_1}}^n_{k+1}
\end{equation}
On the other hand, choosing $x = y_0$ note that since $\scA(y_0 \to y_1)$ was assumed to not equal the identity, neither can $\sL^n_{\vec\nfc_{\scB_x}}(y_1, 0 \to n)$ equal the identity (that is, $\rest \scA {y_1}$ is not in normal form).  

We now construct a collapse limit $\scC : Y \times (\bnum{n+1}) \to \SIvert n \cC$ corresponding to a non-identity collapse and thus contradicting that $\rest \scA Y$ is in normal form. Since $\scC$ is a two-variable functor we can define it by defining the family of partial functors $\scC(-,i)$ ($i \in \bnum{2}$), $\scC(x,-)$ ($x \in Y$) and show their commutativity (cf. \autoref{notn:two_variable_functors}). Explicitly we define
\begin{itemize}
\item For $k \in (\bnum{n+1})$ we define $\scC(-,k)$ by setting $\rest {\scC(-,k)} {Y_{k-1}} = \rest \scA {Y_{k-1}}$ and
\begin{equation}
\scC(y_1,k) = \NF{\rest \scA {y_1}}^n_{k+1}
\end{equation}
It remains to define $\scC(-,0)$ on morphisms $(x \to y_1) \in \mor(Y)$ for which we set
\begin{equation}
\scC(x \to y_1,k) = \sL^n_{\vec\nfc_{\scB_x}}(y_1,0 \to k)
\end{equation}
We check functoriality of this definition: let $x_1 \to x_2 \to x_3$ be a chain of two non-identity morphisms in the poset $X$ (that is, such that $x_l \neq x_j$ for $l < j \in \Set{1,2,3}$). We want to show
\begin{equation}
\scC(x_1 \to x_2 \to x_3,k) = \scC(x_2 \to x_3,0) \circ \scC(x_1 \to x_2,k)
\end{equation}
This holds if $x_j \in Y_{k-1}$. We need to consider the remaining case that $x_1, x_2 \in Y_{k-1}$ and $x_3 = y_1$. In this case we calculate
\begin{align}
\scC(x_1 \to x_2 \to x_3,0) &=  \sL^n_{\vec\nfc_{\scB_x}}(y_1,0 \to k) \\
&= \sL^n_{\vec\nfc_{\scB_x}}(y_1,0 \to k) \circ \id \\
&= \scC(x_2 \to y_1,0) \circ \scC(x_1 \to x_2,0) 
\end{align}
where we use $\scC(x_1 \to x_2,0) = \id$ by global triviality of $\rest \scA {Y_{k-1}}$ (which follows from our inductive assumption).

\item For $x \in Y$ we define $\scC(x,-)$ as the partial functor
\begin{equation}
\scC(x,-) = \sL^n_{\vec\nfc_{\scB_x}}(x,-)
\end{equation}
Note that if $x \neq y_1$, then by the above remarks $\scC(x,-)$ is the constant functor.
\end{itemize}

We are left which showing commutativity of partial functors (cf. \autoref{notn:two_variable_functors}). Specifically, we need to verify (for $x_1, x_2 \in Y$ and $l,k \in (\bnum{n+1})$)
\begin{equation} 
\xymatrix@C=2cm{\scC(x_1,l) \ar[r]^{\scC(x_1 \to x_2,l)} \ar[d]_{\scC(x_1, l \to k)} & \scC(x_2,l) \ar[d]^{\scC(x_2, l \to k)} \\ \scC(x_1,k) \ar[r]_{\scC(x_1 \to x_2,k)} & \scC(x_2,k) }
\end{equation}
The only interesting case is $x_1 = x \in Y_{k-1}$ and $x_2 = y_1$, for which we obtain
\begin{equation} 
\xymatrix@C=3cm{
\scA(x) 
		\ar[r]^-{\sL^n_{\vec\nfc_{\scB_x}}(y_1,0 \to l)} \ar[d]_-{\id} 
& \NF{\rest \scA {y_1}}^n_l(y_1) 
		\ar[d]^-{\sL^n_{\vec\nfc_{\scB_{y_1}}}(y_1,l \to k)} 
\\ 
\scA(x) 
		\ar[r]_-{\sL^n_{\vec\nfc_{\scB_x}}(y_1,0 \to k)} 
& \NF{\rest \scA {y_1}}^n_k(y_1) }
\end{equation}
which commutes by functoriality of $\sL^n_{\vec\nfc_{\scB_x}}$ after observing that
\begin{equation}
\sL^n_{\vec\nfc_{\scB_x}}(y_1,-) = \sL^n_{\vec\nfc_{\scB_{y_1}}}(y_1,-)
\end{equation}
by uniqueness of the normal form collapse sequence for $\rest \scA {y_1}$.

We have thus constructed a non-trivial collapse limit $\scC$ which translates into a non-trivial collapse of $\cN_\scC : \rest \scA Y \to \scC(-,0)$ (indeed, since $\rest \scA {y_0} \neq \rest \scA {y_1}$ but $\scC(y_0 \to y_1,0) = \id$, $\cN_\scC$ must be non-trivial). This contradicts that $\rest \scA Y$ is in normal form and thus finishes the proof. 
\qed
\end{thm}

\section{Globularity} \label{sec:globe_globular}

In this section we introduce globular cubes. Globularity distinguishes the case of $n$-fold categories from $n$-categories (the latter containing morphisms which satisfy globularity). We will formulate two versions of globularity
\begin{itemize}
\item $\partial$-globularity: A cube is $\partial$-globular (pronounced ``boundary globular") if it is constant along its sides in a certain way.
\item Strict globularity: A cube is strictly globular if itself and all its sub-cubes are $\partial$-globular.
\end{itemize}
By default ``globular" will refer to ``strictly globular". The advantage of strict globularity over $\partial$-globularity is that it is stable under taking subcubes. Geometrically, it also disallows certain ``highly singular" situations from arising, by preventing singular behaviour to happen on regular segments. 

\begin{rmk}[Without loss of generality, assume strict globularity] \label{rmk:strict_glob} Importantly, all relevant results up to \autoref{ch:groupoids} (except \autoref{thm:subfamilies_inherit_globularity}) which are proven for strict globularity in this thesis, can be easily amended to statements about $\partial$-globularity. More specifically, the exchangeability of ``strictly globular" and ``$\partial$-globular" holds for definitions of presented associative $n$-categories (in \autoref{ch:presented}) and $n$-groupoids (in \autoref{ch:groupoids}) but not our discussion of generic composites (in \autoref{ch:composition} and \autoref{ch:associative}). Ultimately, being stable under taking subcubes appears to make strict globularity the more canonical choice.
\end{rmk}

Both versions of globularity will be formulated as a local triviality condition of $k$-level labellings. Globular cubes will also be called \textit{$n$-globes}, but they can have a lot of non-trivial structure when compared to the simple cell complexes that are usually called ``$n$-globes" in the literature. However, a class of globes which directly represent elements of \textit{globular sets} and thus ``classical $n$-globes", will be introduced at the end of the chapter, and will be called \textit{terminal $n$-globes}.

\subsection{Globular cubes}

We start with a definition of $\partial$-globularity.

\begin{defn}[$\partial$-globularity] \label{defn:partial_glob} Let $\scA : X \to \SIvert n \cC$ and $1 \leq k \leq n$. We say that $\scA$ is $k$-level $\partial$-globular if $\sU^k_\scA$ is locally trivial on both the image of $\msrc_{\tusU {k-1}_\scA}$ and the image of $\mtgt_{\tusU {k-1}_\scA}$. 

We say $\scA$ is $\partial$-globular if it is $k$-level $\partial$-globular for all $k$.
\end{defn}

We now turn to a stricter version of globularity. We first give the definition of singular and regular content which will enable us to formulate this stricter definition.

\begin{defn}[Singular/regular content of $\SI$-families] \label{defn:singular_content} Let $X$ be a poset, and $\scA: X \to \SI$ a singular interval family. 
\begin{itemize}
\item We define the \textit{singular content} $\singcont(\sG(\scA))$ of $\scA$ to be the full subposet of $\sG(\scA)$ generated by the objects $(x,a) \in \sG(\scA)$ such that $a \in \singcont(\scA(x))$.
\item Similarly, we define the \textit{regular content} $\regcont(\sG(\scA))$ of $\scA$ to be the full subposet of $\sG(\scA)$ generated by the objects in $(x,b) \in \sG(\scA)$ such that $b \in \regcont(\scA(x))$.
\end{itemize} 
\end{defn}

\begin{eg}[Singular/regular content] Using previous examples we give two \SI-bundles and in each case circle points in the regular content in \cblue{}, and points in the singular content in \cred{}. First consider the following \SI-bundle 
\begin{restoretext}
\begingroup\sbox0{\includegraphics{test/page1.png}}\includegraphics[clip,trim=0 {.3\ht0} 0 {.3\ht0} ,width=\textwidth]{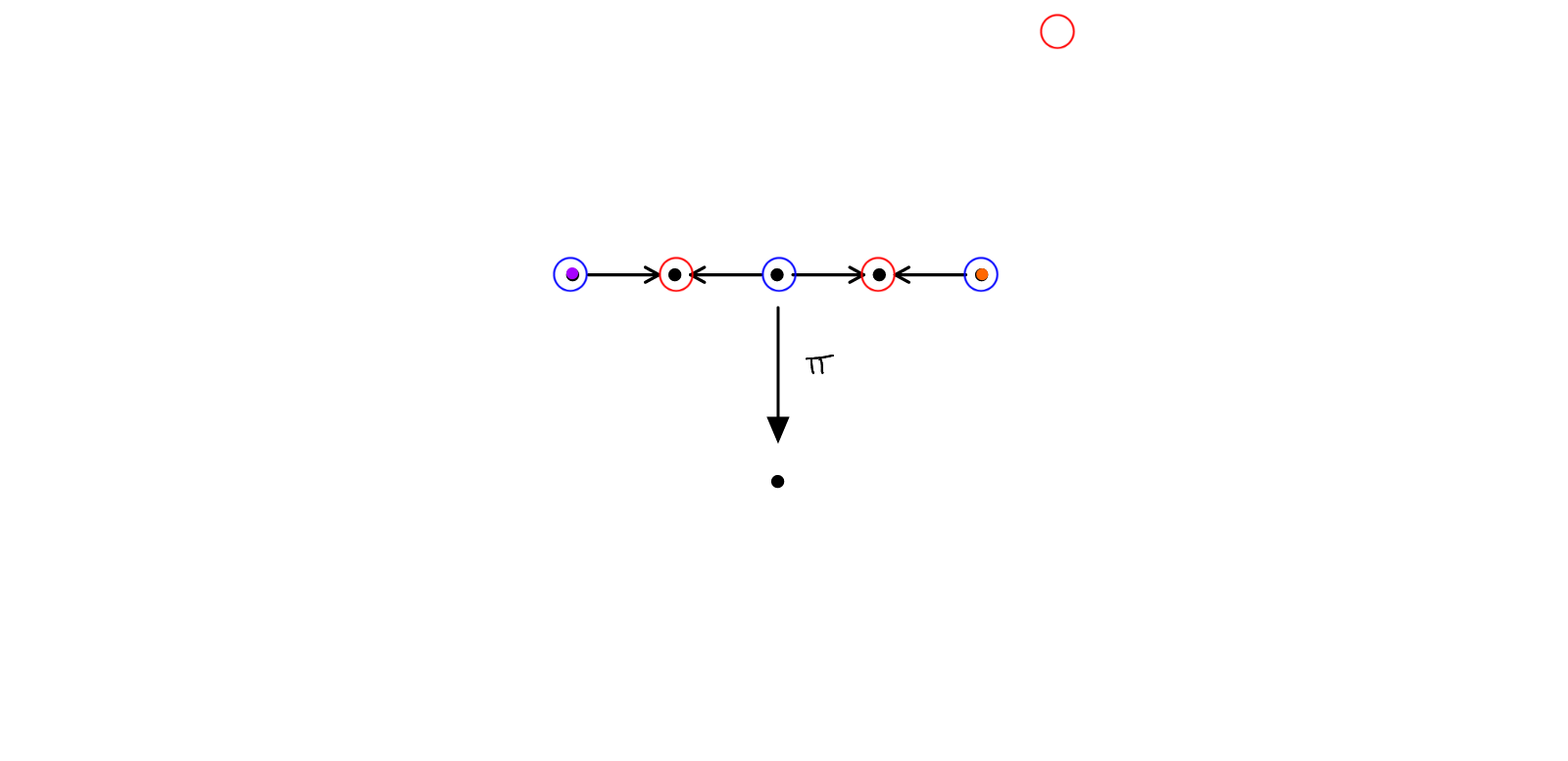}
\endgroup\end{restoretext}
Secondly consider
\begin{restoretext}
\begingroup\sbox0{\includegraphics{test/page1.png}}\includegraphics[clip,trim=0 {.2\ht0} 0 {.1\ht0} ,width=\textwidth]{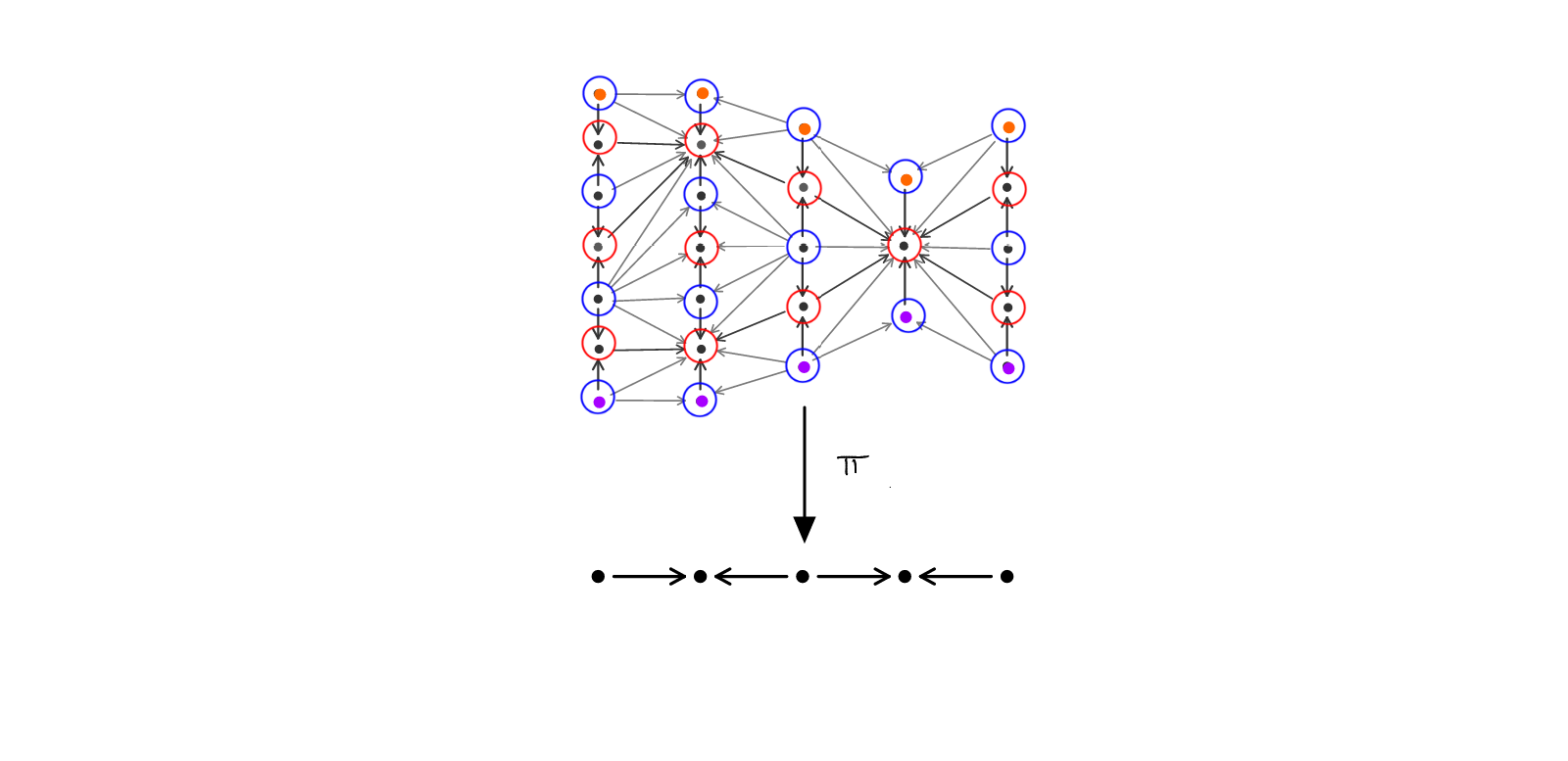}
\endgroup\end{restoretext}
Note that in each case we also marked the source and target sections by \cpurple{} and \corange{} dots respectively.
\end{eg}

\begin{rmk}[Regular content is downward closed] \label{rmk:regular_content} Note that $\regcont(\sG(\scA)) \subset \sG(\scA)$ is a \gls{downwardclosed} subset: indeed if $(x,a) \to (y,b)$, $(y,b) \in \regcont(\sG(\scA))$, then $(a,b) \in \edgeset(\scA(x \to y))$ and $b \in \regcont(\scA(y))$ and thus by \eqref{eq:defn_order_realisation_3} we must have $a = \scA(x \to y)\regop (b)$ implying $a \in \regcont(\scA(x))$ as required for $(x,a) \in \regcont(\sG(\scA))$.
\end{rmk}

\begin{defn}[Globular families] \label{defn:globular_families} Let $\scA : X \to \SIvert n \cC$ and $1 \leq k \leq n$. We say $\scA$ is \textit{$k$-level globular} if $\tsU k_\scA$ is locally trivial on $\regcont(\tsG k(\scA))$. 

If $X$ itself is explicitly assumed to be of the form $\sG(F)$ for some $F : Y \to \SI$, then we further say $\scA$ is $0$-level globular if $\scA$ is locally trivial on  $\regcont(\sG(F))$. Otherwise, $\scA$ is $0$-level globular by default. 

We say $\scA$ is \textit{globular} if it is $k$-level globular for all $k$, $0 \leq k \leq n$. 
\end{defn}

A globular $\cC$-labelled $n$-cube $\scA : \bnum{1} \to \SIvert n \cC$ is also called a ($\cC$-labelled) $n$-globe.

\begin{rmk}[Source, target inclusions are in regular content] \label{rmk:src_tgt_are_regular_content} To compare \autoref{defn:partial_glob} and \autoref{defn:globular_families}, recall that $\msrc_\scA, \mtgt_\scA$ are open sections and thus $\im(\msrc_\scA), \im(\mtgt_\scA) \subset \regcont(\sG(\scA))$. Further note that both $\im(\msrc_\scA), \im(\mtgt_\scA)$ are \gls{downwardclosed} as well as connected if the base space $X$ of $\scA : X \to \SI$ is connected (cf. \autoref{rmk:connectedness}).
\end{rmk}

\begin{egs}[Globular families] The reader is invited to verify that all families in \autoref{ch:emb} are globular. Recall the example of a cube $\scA_b : \bnum{1} \to \SIvert 3 \cC$ from \autoref{eg:subfamilies}. We want to see that $\scA_b$ is $2$-level globular and thus pick some $(x \to y) \in \regcont(\tsG 2(\scA_b))$ and need to check that $\tsU 2_{\scA_b}$ restricted to $(x \to y)$ normalises to the constant family: indeed, $\rest {\tsU 2_{\scA_b}} {x \to y}$ for a choice of $(x \to y)$ is highlighted in the below illustration (the rest of the parent family is greyed out)
\begin{restoretext}
\begin{noverticalspace}
\begingroup\sbox0{\includegraphics{test/page1.png}}\includegraphics[clip,trim=0 {.0\ht0} 0 {.0\ht0} ,width=\textwidth]{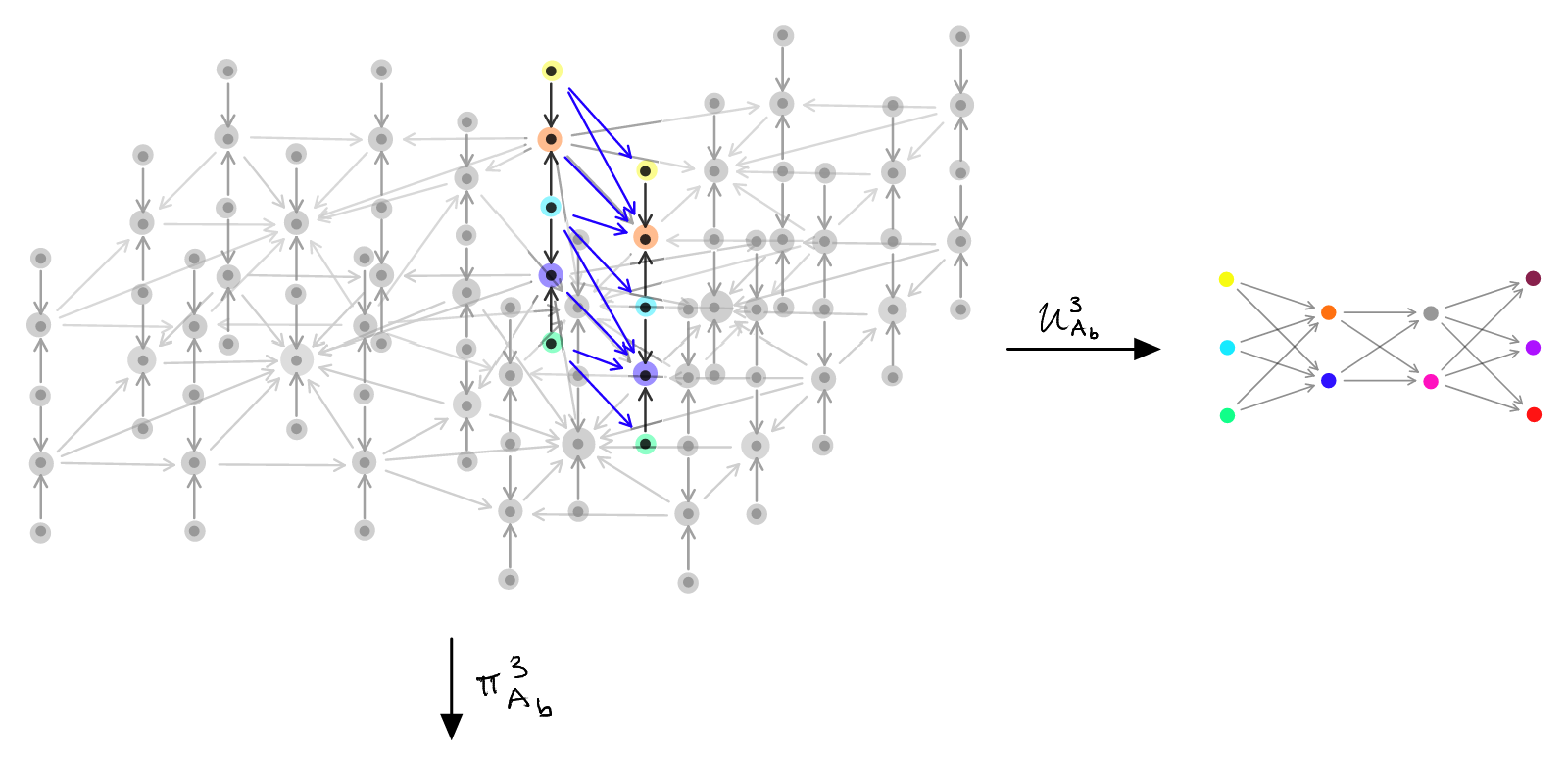}
\endgroup \\*
\begingroup\sbox0{\includegraphics{test/page1.png}}\includegraphics[clip,trim=0 {.6\ht0} 0 {.0\ht0} ,width=\textwidth]{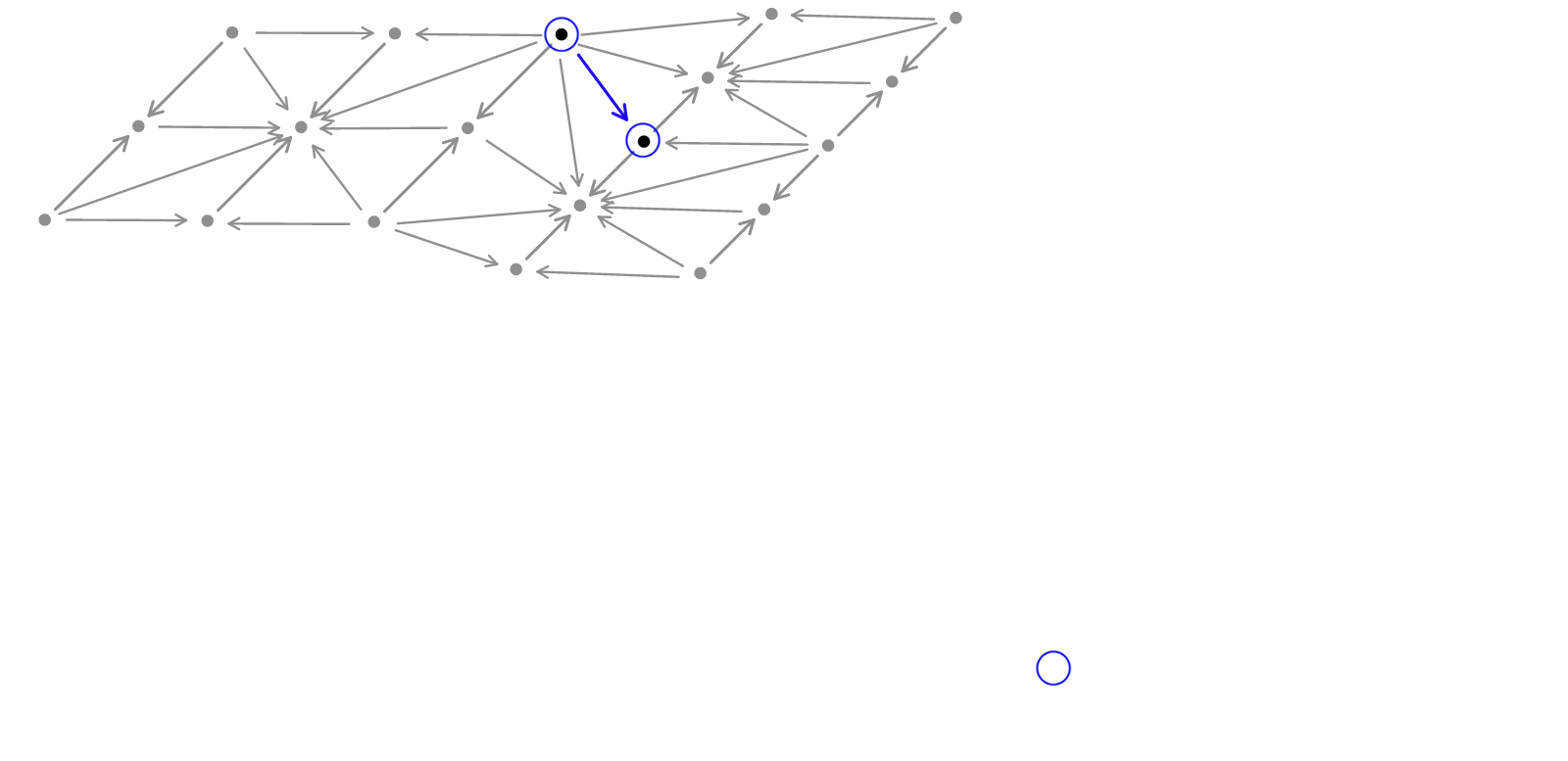}
\endgroup
\end{noverticalspace}
\end{restoretext}
Thus, $\rest {\tsU 2_{\scA_b}} {x \to y}$ normalises to (in fact is already) a constant family, and the same holds for all other choices of $(x \to y) \in \regcont(\tsG 2(\scA_b))$.

Next, we want to see $\scA_b$ is $3$-level globular. Thus again, we pick $(x \to y) \in \regcont(\tsG 3(\scA_b))$ and now need to check that $\rest {\tsU 2_{\scA_b}} {x \to y}$ is constant. Note that $\rest {\tsU 2_{\scA_b}} {x \to y}$ is a $\SIvert 0 \cC$-family, that is, a functor into $\cC$, and thus no normalisation applies to it. We check $3$-level globularity for the three highlighted choice of $(x \to y)$ below
\begin{restoretext}
\begingroup\sbox0{\includegraphics{test/page1.png}}\includegraphics[clip,trim=0 {.2\ht0} 0 {.0\ht0} ,width=\textwidth]{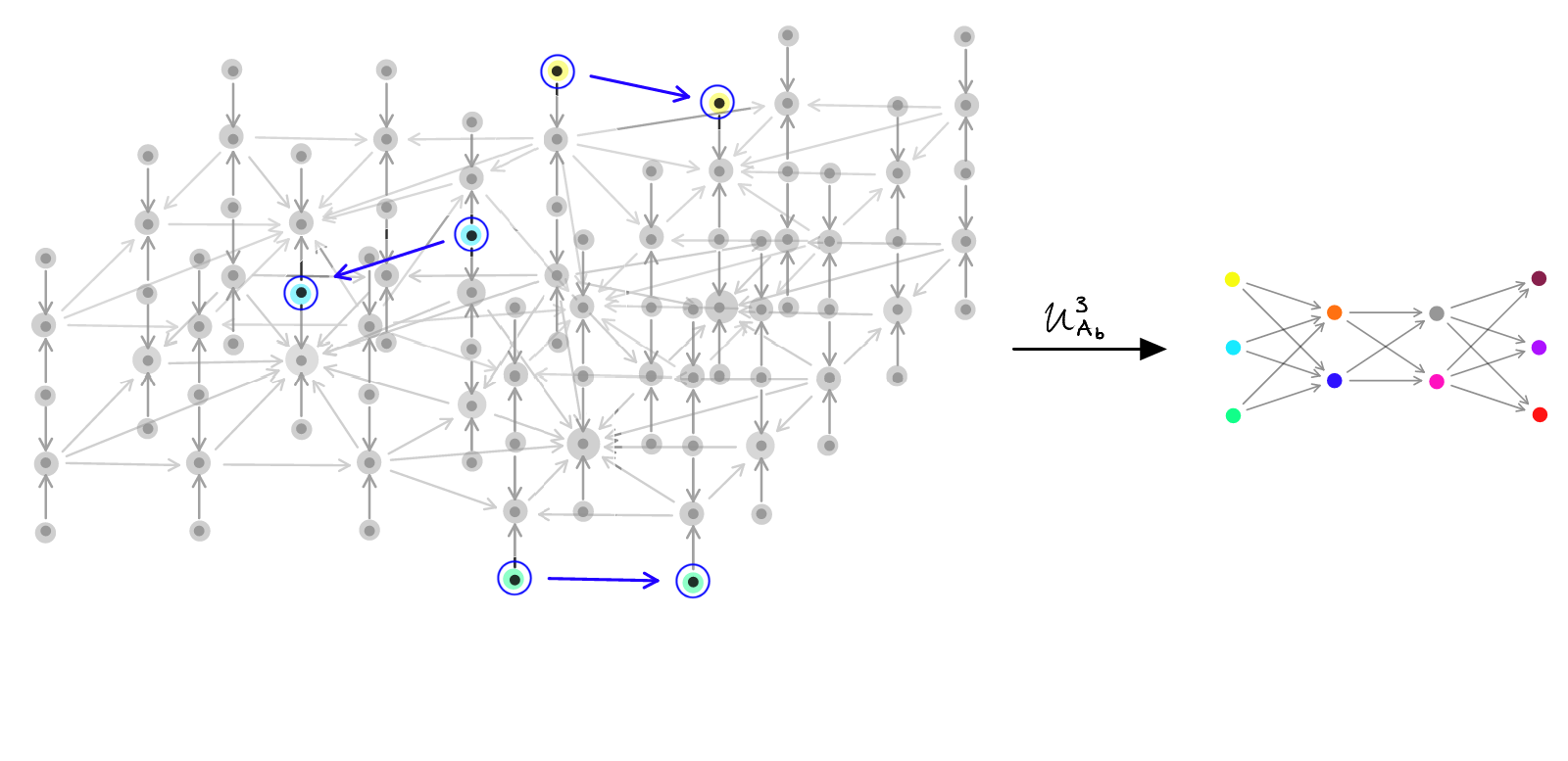}
\endgroup\end{restoretext}
Note that $\scA_b$ is also $1$-level globular since there are no non-trivial $(x \to y) \in \regcont(\tsG 1(\scA))$. Thus we conclude that $\scA_b$ is globular.
\end{egs}

\begin{claim} \label{claim:globularity_inherited_by_subposets} Let $\scA : X \to \SIvert n \cC$ be $k$-level globular and $H : Y \to X$ be an inclusions of posets. Then $\scA H$ is $k$-level globular too.
\proof  The proof is \stfwd{}. Since the components of $\restemb_H : \scA H \subset \scA$ (cf. \autoref{eg:subfamily_by_restriction}) fibrewise preserve regular segments, we have
\begin{equation}
\restemb_H^k(\reg (\tsG k(\scA H))) \subset \regcont(\tsG k(\scA))
\end{equation}
Note that (cf. \eqref{eq:restriction_colorings}) 
\begin{equation}
{\tsU k_{\scA}} \restemb_H^k = \tsU k_{\scA H}
\end{equation}
Thus, $\rest {\tsU k_{\scA H}} {\reg (\tsG k(\scA H)) }$ is locally trivial since $\rest {\tsU k_{\scA}} {\regcont(\tsG k(\scA))}$ is locally trivial, which gives $k$-level globularity of $\scA H$ as required.
\end{claim}

\begin{cor}[Restrictions inherit globularity] \label{cor:globularity_inherited_by_subposets} If $\scA : X \to \SIvert n \cC$ is globular and $H : Y \to X$ is an inclusions of subposet, then $\scA H$ is globular too. \qed
\end{cor}

\subsection{Globular source and target}

In the next chapter, the following definition will play an important role when gluing globes together.

\begin{defn} \label{defn:globular_src_tgt} Let $\scA : \bnum{1} \to \SIvert n \cC$ be globular, $n > 0$. In this case, we define
\begin{align}
\gsrc(\scA) &:=  \sU^1_\scA \msrc_{\und\scA}: \bnum{1} \to \SIvert {n-1} \cC\\
\gtgt(\scA) &:=  \sU^1_\scA \mtgt_{\und\scA} : \bnum{1} \to \SIvert {n-1} \cC
\end{align}
called the \textit{globular source} and \textit{globular target} (of $\scA$) respectively.
\end{defn}

As an application of the results in the previous section we remark the following.

\begin{rmk}[Source, target inherit globularity] \label{cor:vertical_degeneracy_src_tgt} Let $\scA :  \bnum{1} \to \SIvert n \cC$ be globular. Using \autoref{cor:globularity_inherited_by_subposets} we remark that then $\gsrc(\scA), \gtgt(\scA) : \bnum{1} \to \SIvert {n-1} \cC$ are globular too.
\end{rmk}

\begin{rmk}[Boundary of a point is empty] \label{notn:boundary_of_a_point} The above defines sources and targets for $n$-cubes, $n > 0$. Recall \autoref{conv:minus_one_cubes}. We extend this convention by defining that (cubical and globular) sources and targets of $0$-cubes will always be understood as the empty $(-1)$-cube $\emptyset : \emptyset \to \cC$. In particular, all sources and targets of $0$-cubes are equal. The source and target of $\emptyset$ are in turn defined to be again $\emptyset$.
\end{rmk}

Examples of globular sources and target will be given in \autoref{ssec:glob_comp}.

\section{Embedding and collapse} \label{sec:globe_emb_coll}

In this short section, we discuss the interaction of globularity with embedding and collapse.

\subsection{Subfamilies inherit globularity}

As the following result is not critical for the rest of the development in this thesis, some details in its proof are left to the reader.

\begin{thm}[Subfamilies inherit globularity] \label{thm:subfamilies_inherit_globularity} Let $\theta : \scB \mono \scA$ where $\scA : X \to \SIvert n \cC$ is a globular $\SIvert n \cC$-family indexed by $X$. Then $\scB$ is a globular too.
\proof Using \autoref{constr:subfamilies_by_specifying_endpoints} is enough to show this in the case of $\sJ^{\scA,k-1}\restsec{[q_-,q_+]} : \scB \mono \scA$ for $0 < k \leq n$.  We verify $\scB$ is $l$ globular. 
\begin{enumerate}
\item If $l \geq k$, this follows immediately from $l$-level globularity of $\scA$, since $\tsU l_\scB \Delta_{x \to y} = \tsU l_\scA \Delta_{\theta^l(x \to y)}$.
\item If $l < k$ without loss of generality assume $l = 0$ (otherwise redefine $\scA$ by $\tsU l_\scA$ and $\scB$ by $\tsU l_\scB$). By restriction we then obtain an embedding
\begin{equation}
\chi := \sJ^{\scA, k-1}\restsec{[q_-,q_+]} \restemb_{\Delta_{x\to y}} : \scB \Delta_{x\to y} \mono \scA \Delta_{x\to y}
\end{equation}
By globularity of $\scA$ we find
\begin{equation}
\vvec \lambda := \vvec \nfc_{ \scA \Delta_{x\to y} } :  \scA \Delta_{x\to y} \starcoll \scC \ppi_{\bnum 0}
\end{equation}
where $\ppi_{\bnum 0} : \bnum 1 \to \bnum 0$ and $\scC : \bnum 0 \to \SIvert l {\cC}$.

Then, using \autoref{constr:restricting_collapse_seq} we find
\begin{equation}
\vvec \lambda\postar  \chi : \scD \mono \const_\scC
\end{equation}
where $\chi\pbstar \vvec\lambda : \scB \Delta_{x\to y} \starcoll \scD$. We need to show $\scD$ is of the form $E \ppi_{\bnum 1}$ for some family $E$. This follows from \autoref{constr:collapse_on_subfamilies} together with the trivial product bundle construction (cf. \autoref{claim:repacking_prod_cube_bun}). The former implies $\vvec \lambda\postar  \chi : \scD \mono \const_\scC$ is of the form 
\begin{equation}
\vvec \lambda\postar  \chi = \sJ^{\scC\ppi_{\bnum 0},k-1}\restsec{[r_-,r_+]} 
\end{equation}
for some open sections $r_-$, $r_+$. The latter implies, that 
\begin{equation} 
\xymatrix@C=2cm{ \tsG {k}(\scC) \times \bnum 1 \ar[r]^{\sW^{\sT_\scC,k}_{\bnum 1}} \ar[d]_{\tpi {k}_\scC \times \bnum 1} & \tsG {k}(\scC\ppi_{\bnum 0})  \ar[d]^{\tpi {k}_{\scC\ppi_{\bnum 0}}} \\
\tsG {k-1}(\scC) \times \bnum 1\ar[r]_{\sW^{\sT_\scC,k-1}_{\bnum 1}} &  \tsG {k-1}(\scC\ppi_{\bnum 0}) }
\end{equation}
Openness of the sections $r_-$, $r_+$ then forces 
\begin{equation}
r_\pm(\sW^{\sT_\scC,k-1}_{\bnum 1}(x,0)) = r_\pm(\sW^{\sT_\scC,k-1}_{\bnum 1}(x,1)) =: s_\pm(x)
\end{equation}
Using \autoref{constr:subfamilies_from_endpoints}, and noting $\scD =  (\scC\ppi_{\bnum 0})^{k-1}\restsec{[r_-,r_+]}$ (by definition of $r_\pm$), we find
\begin{equation}
\scD = (\scC^{k-1}\restsec{[s_-,s_+]}) \ppi_{\bnum 0}
\end{equation}
as required. \qed
\end{enumerate}
\end{thm}

\subsection{Collapse preserves and reflects globularity}

\begin{thm}[Collapse preserves globularity] \label{thm:collapse_preserves_globularity} Let $\vsS{\vvec \lambda} : \scA \to \scB$. Then $\scA$ is a globular $\SIvert n \cC$-cube if and only if  $\scB$ is a globular $\SIvert n \cC$-cube.
\proof  The proof is \stfwd{}. Assume $\scA$ is globular. Let $(x \to y) \in \regcont(\tsG l(\scB)))$. By surjectivity and openness of $(\vsS{\vvec \lambda})^l$ we find $(u \to v)\in \regcont(\tsG l(\scA)))$ such that $(\vsS{\vvec \lambda})^l(u \to v) = (x \to y)$. The multi-level collapse $\vsS{\vvec \lambda} : \scA \to \scB$ then restricts to give multi-level collapse $\tsU l_\scA \Delta_{u \to v} \to \tsU l_\scB \Delta_{x \to y}$. Since the former normalises to the identity (by globularity of $\scA$) the latter does too. This shows globularity of $\scB$. 

Conversely, assume $\scB$ is globular. Let $(x \to y) \in \regcont(\tsG l(\scA))$. Set $(u \to v) := (\vsS{\vvec \lambda})^l(x \to y) \in \regcont(\tsG l(\scA))$. The multi-level collapse $\vsS{\vvec \lambda} : \scA \to \scB$ then restricts to give multi-level collapse $\tsU l_\scA \Delta_{x \to y} \to \tsU l_\scB \Delta_{u \to v}$. Since the latter normalises to the identity (by globularity of $\scA$) the former does too. This shows globularity of $\scA$. \qed
\end{thm}

\section{Globular cones} \label{sec:globe_cone}

The goal of this section is to first show that a certain monad (called the $\top$-monad and acting on categories by adjoining a terminal object) can be generalised to normalised globular cube families. This in turn can be applied to construct \textit{double cones} which are globular cubes of conical shapes, and which will be the content of the second part of the section. More precisely, double cones are $(k+1)$-cubes obtained from ``coning" two $k$-cubes and joining them on the side of their vertex point (this requires their boundaries to agree). In the last part of the section, we will then focus on special double cones, which are called \textit{terminal $n$-globes}.

\subsection{The $\top$-monad}

We start by recalling the definition of the usual $\top$-monad.

\begin{notn}[$\top$-Monad] \label{notn:top_monad} There is a monad $(-)\topmon : \Cat \to \Cat$, called the \textit{top monad}, which when applied to a category $\cC$, yields the category $\cC\topmon$ which is obtained from $\cC$ by adjoining a single terminal object called $\top_\cC$. It maps a functor $F : \cC \to \cD$ to a functor $F\topmon : \cC\topmon \to \cD\topmon$ which maps $\top_\cC$ to $\top_\cD$ and otherwise acts as $F$. The unit is given by the canonical embedding $\eta^\top_\cC : \cC\to \cC\topmon$, and the multiplication $\mu^\top_\cC : (\cC\topmon)\topmon \to \cC\topmon$ maps $\top_{\cC\topmon}$ and $\top_\cC$ to $\top_\cC$ and acts otherwise as the identity on $\cC$.
\end{notn}

The $\SIvert n \empty$-construction is compatible with the $\top$-monad as the following construction shows.

\begin{constr}[$\top$-Monad for $\SIvert  n \cC$-familes]  \label{constr:top_monad}  Let $\scA : X \to \SIvert n \cC$ be a normalised and globular $\SIvert n \cC$-family such that $\scA$ is constant on the images of $\msrc_{\und\scA}$ and $\mtgt_{\und\scA}$ (if $X$ is connected this condition is always fulfilled). We construct $\scA\topmon : X\topmon \to \SIvert n {\cC\topmon}$, called the \textit{cone of $\scA$} for which we make the following claims
\begin{enumerate}
\item $\scA\topmon \eta^\top_X = \SIvert n {\eta^\top_\cC} \scA$
\item $\scA\topmon$ is globular and normalised
\item There is a unique $\top_\scA \in \tsG n(\scA\topmon)$ with $\tsU n_{\scA\topmon}(\top_\scA) = \top_\scC$ and this $\top_\scA$ satisfies $\ctypsum^n_{\scA\topmon}(\top_\scA) = n$ as well as (cf. \autoref{constr:minimal_subfamilies})
\begin{equation}
\scA\topmon \sslash \top_\scA = \scA\topmon
\end{equation}
\item $\scA\topmon$ is the unique $\SIvert n {\cC\topmon}$-family over $X\topmon$ with the above three properties
\end{enumerate}

The construction of $\scA\topmon$ and the proof of the claimed properties are inductive in $n$. For $n = 0$, the construction equals the usual $\top$-monad and verification of the properties is trivial (in this case $\top_\scA = \top_X$).

Let $n > 0$. We construct $\scA\topmon$ by first setting $\tusU 0_{\scA\topmon} : X\topmon \to \SI$ to map $\top_X$ to $\singint 1$ (the terminal object of $\SI$) and otherwise act as $\tusU 0_\scA$, that is
\begin{equation}
\tusU 0_{\scA\topmon} \eta^\top_X = \tusU 0_\scA
\end{equation}
Setting $Y$ to be the full subposet of $Z := \sG(\tusU 0_{\scA\topmon})$ consisting both of objects $(\top_X, 0)$, $(\top_X,2)$ and objects in the image of $\sG(\eta^\top_X)$, we note that there is an isomorphism
\begin{equation}
\alpha : Y\topmon \iso Z
\end{equation}
which includes $Y$ as a subset and maps $\top_Y$ to $(\top_X,1)$. Now, define $\tsU 1_{\scA\topmon} : Z \to \SIvert {n-1} {\cC\topmon}$ as follows. We first define $U^1_\scA : Y \to \SIvert {n-1} {\cC}$ on $Y$, by making it agree with $\tsU 1_{\scA}$ on the image of $\sG(\eta^\top_X)$ and further requiring it to be constant on the images of $\msrc_{\tusU 0_{\scA\topmon}}$ and $\mtgt_{\tusU 0_{\scA\topmon}}$ (this is possible by our assumption on $\scA$). Arguing inductively, we can now define
\begin{equation}
\tsU 1_{\scA\topmon} := (U^1_\scA)\topmon \alpha\inv
\end{equation}
This completes the construction of $\scA\topmon$. We now prove its claimed properties, which is \stfwd{}.
\begin{enumerate}
\item $\scA\topmon\eta^\top_X = \scA$ follows since $\tusU 0_{\scA\topmon}\eta^\top_X = \tusU 0_\scA$ and $\tsU 1_{\scA\topmon} \tsG 1(\eta^\top_X) = \SIvert {n-1} {\eta^\top_\cC} \tsU 1_\scA$.
\item First note that $U^1_\scA$ is globular and normalised, since it is globular and normalised on the \gls{downwardclosed} subposets $\tsG 1(\scA)$, $Y \sslash (\top_X,0)$ and $Y \sslash (\top_X,2)$ (in fact, it was chosen to be constant on the latter two subposets) and these three subposets cover $Y$. Thus arguing inductively $(U^1_\scA)\topmon$ is globular and normalised. This implies $\scA\topmon$ is globular and normalised: globularity follows since firstly, $\tusU 1_{\scA\topmon}$ was chosen to be constant on the images of $\msrc_{\tusU 0_{\scA\topmon}}$ and $\mtgt_{\tusU 0_{\scA\topmon}}$ and secondly, since $\tsU 1_{\scA\topmon}$ is globular. Normalisation follows since $\tsU 1_{\scA\topmon}$ is normalised and no $1$-level collapse can apply to $\scA\topmon$. To see the latter note that since $\scA$ is normalised and thus the only non-trivial component of a $1$-level collapse $\lambda$ can be $\lambda_{\top_X}$. However, in the next item we will find $\top^1_\scA = (\top_X,1)$ and thus $\tsU 1_{\scA\topmon}(\top_X,-)$ cannot be constant (which would be required for such a non-trivial $1$-level collapse $\lambda$).

\item Define $\top_\scA := \tsG {n-1}(\alpha)(\top_{U^1_\scA}) \in \tsG n (\scA)$, which by induction applied to $U^1_\scA$ is the unique region labelled by $\top_\scC$. Note that
\begin{align}
\top^1_\scA &= \alpha(\top^0_{U^1_\scA}) \\
&= \alpha(\top_Y) \\
&= (\top_X,1)
\end{align}
Inductively we also find
\begin{align}
\ctypsum^n_{\scA\topmon}(\top_\scA) &= 1 + \ctypsum^{n-1}_{U^1_\scA} (\top_{U^1_\scA})\\
&= 1 + (n-1) \\
&= n
\end{align}
where, in the first step we used that $\alpha(\top_Y) = (\top_X,1)$. Arguing inductively, we know
\begin{equation}
(\tsU 1_\scA)_\bot \sslash \top_{U^1_\scA} = \tsU 1_\scA
\end{equation}
We deduce that for $1 \leq k \leq n$, $\tsG k(\scA) \sslash \top^k_\scA = \tsG k(\scA)$. We further find that
\begin{align}
\tsG 0(\scA) \sslash \top^0_\scA &= X\topmon \sslash \top_X \\
&= X\topmon
\end{align}
Thus \autoref{claim:minimal_subfamily_is_minimal} implies that $\scA\topmon \sslash \top_\scA = \scA\topmon$ as claimed.
\item Finally, we prove that $\scA\topmon$ is the unique bundle with the above properties. Assume there was a different bundle $\scB$ over $X\topmon$ satisfying firstly that $\scB \eta^\top_X = \SIvert n {\eta^\top_\cC} \scA$, secondly that $\scB$ is globular and normalised, and finally that there is a unique $t \in \scB$ with $\tsU n_\scB(t) = \top_\cC$ such that $\scB \sslash t = \scB$. First note we must have $t^0 = \top_X$. Next we must have $\und \scB(\top_X) = \singint 1$ and $t^1 = (\top_X, 1)$ since otherwise $\tsG 1(\scB) \sslash t^1 \subsetneq \tsG 1(\scB)$ (which would contradict $\scB \sslash t = \scB$). Then we can argue by induction for $Y$-indexed $(n-1)$-cube family $U^1_\scA$, to find that $(U^1_\scA)\topmon = \tsU 1_\scB \alpha$ and consequently we have $\tsU 1_{\scA\topmon} = \tsU 1_\scB$. Since $\und{\scA\topmon} = \und\scB$ (by assumption $\scB\eta^\top_X = \scA$ and our observation that $\scB(\top_X) = \singint 1$) we infer $\scA\topmon = \scB$ as claimed.
\end{enumerate}
\end{constr}

\begin{rmk}[$\top$-restriction is normalised] \label{rmk:top_restriction_normalised} Further to the previous construction, note that $\scA\Delta\topmon$ is normalised.
\end{rmk}

\begin{eg}[$\top$-Monad for families] Consider the family $\scA : X \to \SIvert 1 \cC$ given by
\begin{restoretext}
\begingroup\sbox0{\includegraphics{test/page1.png}}\includegraphics[clip,trim=0 {.05\ht0} 0 {.1\ht0} ,width=\textwidth]{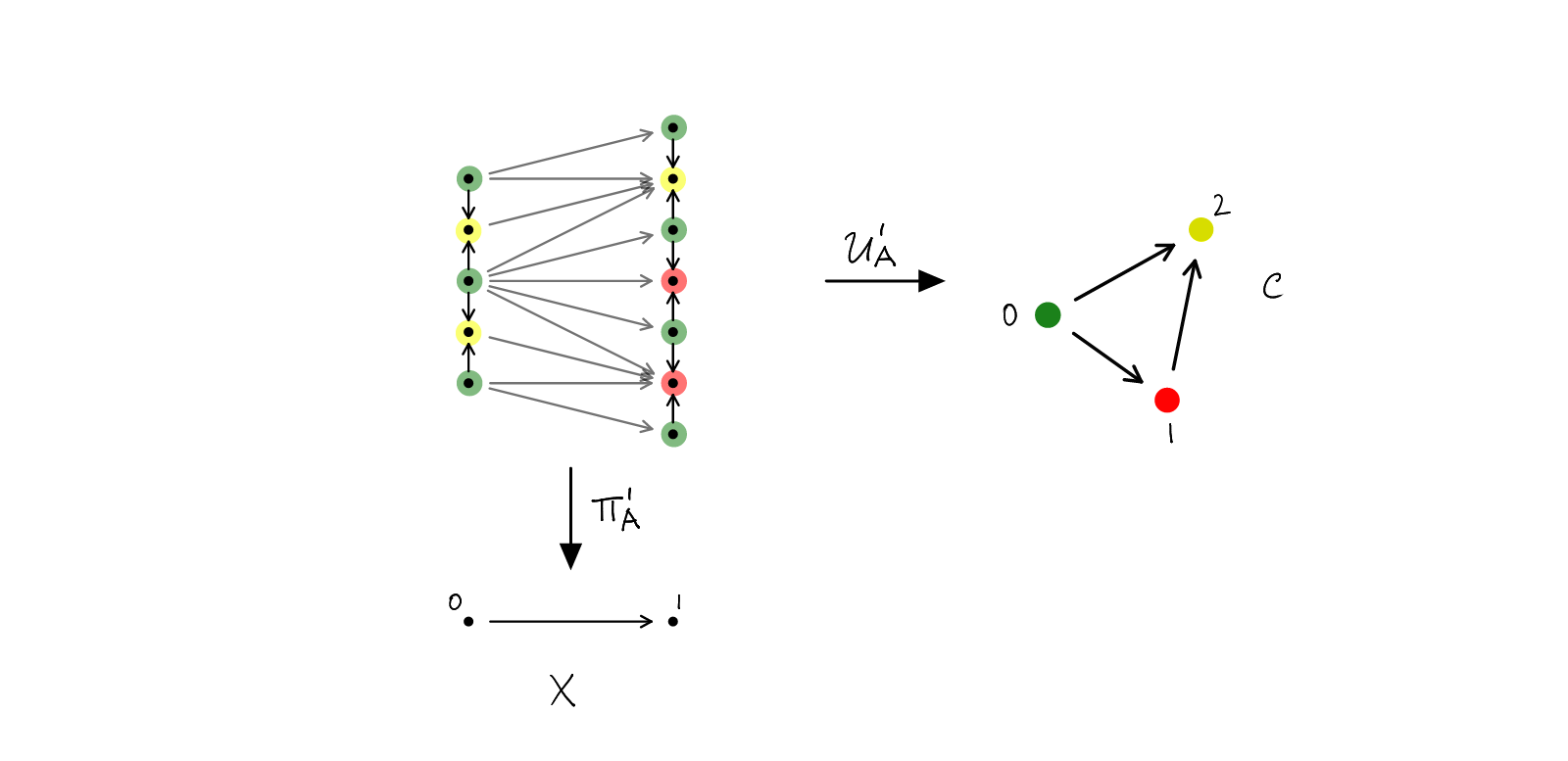}
\endgroup\end{restoretext}
Then $\scA\topmon : X\topmon \to \SIvert 1 {\cC\topmon}$ is given by
\begin{restoretext}
\begingroup\sbox0{\includegraphics{test/page1.png}}\includegraphics[clip,trim=0 {.0\ht0} 0 {.1\ht0} ,width=\textwidth]{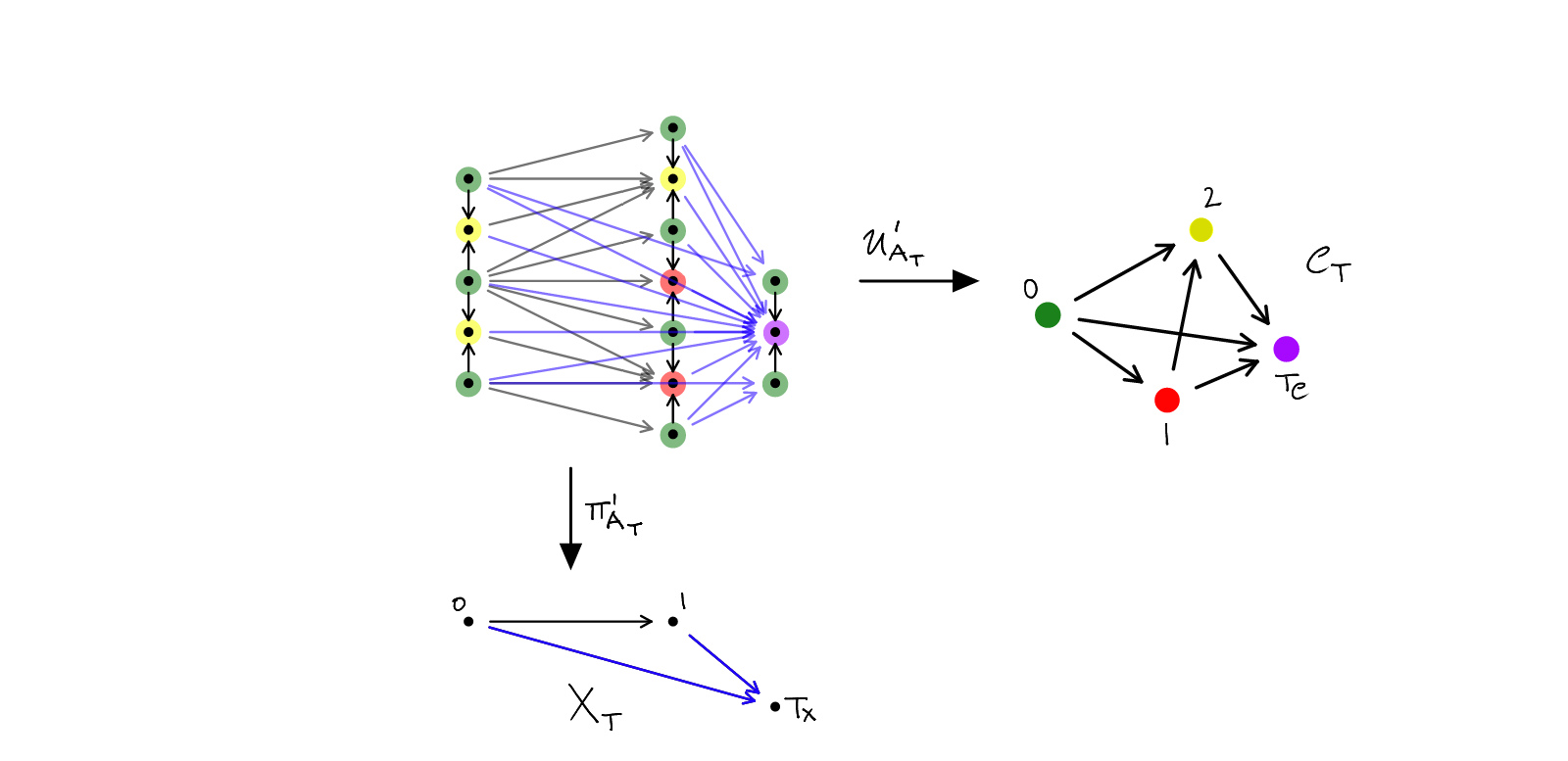}
\endgroup\end{restoretext}
\end{eg}

\subsection{Double cone construction}

Using the $\top$-monad, it is possible to construct an $(n+1)$-globe with a single ``singular vertex point" from two given $n$-globes that coincide on their boundary. This $(n+1)$-globe is called a double cone. Before we give its construction, we introduce the following.

\begin{notn}[Factorisation of labels] \label{notn:restrictions_of_labels} Let $\scA : X \to \SIvert n \cC$ and $F : \cD \to \cC$ be some functor such that there is $U$ with 
\begin{equation}
U^n_\scA = F U
\end{equation}
Then we denote
\begin{equation}
F\pbstar  \scA := \tsR n_{\sT^n_\scA, U} : X \to \SIvert n {\cC_0}
\end{equation}
called the \textit{factorisation of labels of $\scA$ through $F$} (if $F = \cC_0 : \cC_0 \subset \cC$ is a subcategory of $\cC$ then we speak of a \textit{restriction of labels}). Note that
\begin{equation}
\SIvert n {F} (F\pbstar  \scA) = \scA
\end{equation}
If $F$ is injective and faithful then there is a $1$-to-$1$ correspondence of $X$-indexed $\SIvert n \cC$-families with labels factoring through $F$ and $X$-indexed $\SIvert n {\cD}$-families. Also note that factorisation of labels along compositions of functors is (contravariantly) functorial.
\end{notn}

\begin{constr}[Double cones of sources and targets] \label{constr:double_cones_of_src_and_tgt} Let $S, T$ be normalised globular $\SIvert {n-1} \cC$-cubes whose globular sources and targets coincide, that is 
\begin{align}
\gsrc(S) &=\gsrc(T)\\
\gtgt(S) &= \gtgt(T)
\end{align}
Let $g \in \cC$ be such that for each $s \in \im(\tsU {n-1}_S)$ and $t \in \im(\tsU {n-1}_T)$ there are unique morphisms $s \to g$ and $t \to g$. There is a unique $\SIvert n \cC$-cube satisfying
\begin{enumerate}
\item $\abss{S \xto g T} = S$ is globular and normalised
\item $\gsrc\abss{S \xto g T} = S$ and $\gtgt\abss{S \xto g T} = T$
\item There is a unique $p_g \in \tsG n\abss{S \xto g T}$ such that $\tsU n_{\abss{S \xto g T}}(p_g) = g$
\item $\ctypsum^n_{\abss{g}} (p_g) = n$ and
\begin{equation}
\abss{S \xto g T} \sslash p_g = \abss{S \xto g T}
\end{equation}
\end{enumerate}
We construct $\abss{S \xto g T}$ as follows. Note that the family $\scA := (S,T) : \bnum{1} + \bnum{1} \to \SIvert {n-1} \cC$ (here $\bnum{1}+\bnum{1}$ is the coproduct in $\Cat$ and $(S,T)$ the unique factorising map) satisfies the conditions for \autoref{constr:top_monad}, and we can form
\begin{equation}
\scA_\bot : (\bnum{1} + \bnum{1})_\bot \to \SIvert {n-1} {\cC_\bot}
\end{equation}
Note that there is an isomorphism $\alpha : \singint 1 \iso (\bnum{1} + \bnum{1})_\bot$ such that $\scA_\bot \alpha \Delta_0 = S$ and $\scA_\bot \alpha \Delta_2 = T$. We can thus define
\begin{equation}
\abss{S \to T} := \sR_{\Delta_{\singint 1}, \scA_\bot\alpha} : \bnum{1} \to \SIvert n {\cC_\bot}
\end{equation}
This is called the \textit{double cone with source $S$ and target $T$}. Now, by assumption on $g$, $S$ and $T$ there is an canonical inclusion 
\begin{equation}
\beta : \cD := \im(\tsU n_{\abss{S \to T}}) \into \cC
\end{equation}
such that $\beta(\top_\cC) = g$ and $\beta(c) = c$ otherwise. Using \autoref{notn:restrictions_of_labels} we then define
\begin{equation}
\abss{S \xto g T} := \SIvert n \beta (\cD\pbstar  \abss{S \to T})
\end{equation}
This has the claimed properties, and the uniqueness claim follows from the uniqueness in \autoref{constr:top_monad}. $\abss{S \xto g T}$ is called the \textit{double cone with source $S$,  target $T$ and vertex point $g$}.
\end{constr}

\begin{eg}[Double cones] The following bundle has globular source $\scA$ and target $\scB$, and is itself equal to the double cone $\abss{\scA \to \scB}$
\begin{restoretext}
\begingroup\sbox0{\includegraphics{test/page1.png}}\includegraphics[clip,trim=0 {.0\ht0} 0 {.1\ht0} ,width=\textwidth]{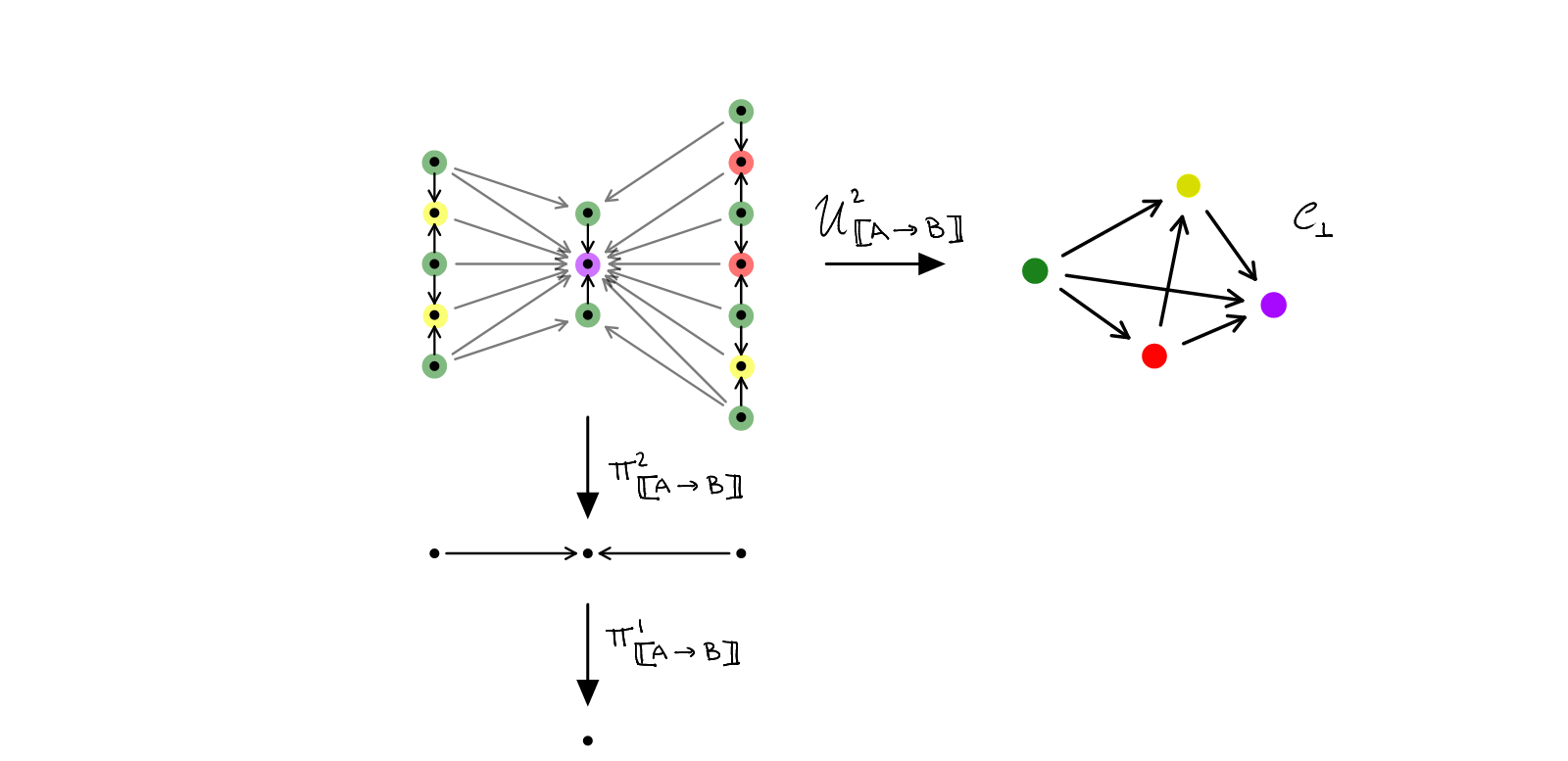}
\endgroup\end{restoretext}
\end{eg}

\subsection{Terminal $n$-globes}

We will now discuss a special class of double cones, called \textit{terminal $n$-globes}. These globular $n$-cubes will be a cubical representation of classical $n$-globes from a globular set. We start by recalling facts about the latter.

\begin{notn}[Globular sets] \label{notn:globeset} The globe category $\lG$ has as objects natural number $n \in \lN$, and morphisms are generated by 
\begin{equation}
\sigma_n, \tau_n : n \to (n+1)
\end{equation}
subject to the conditions
\begin{align}
\sigma_{n+1} \tau_n &= \sigma_{n+1} \sigma_n \\ 
\tau_{n+1} \tau_n &= \tau_{n+1} \sigma_n  
\end{align}
The category $\lG$ can be visualized (for its first 5 objects $0,1,2,3,4$ and morphisms between them) as follows
\begin{restoretext}
\begingroup\sbox0{\includegraphics{test/page1.png}}\includegraphics[clip,trim=0 {.35\ht0} 0 {.0\ht0} ,width=\textwidth]{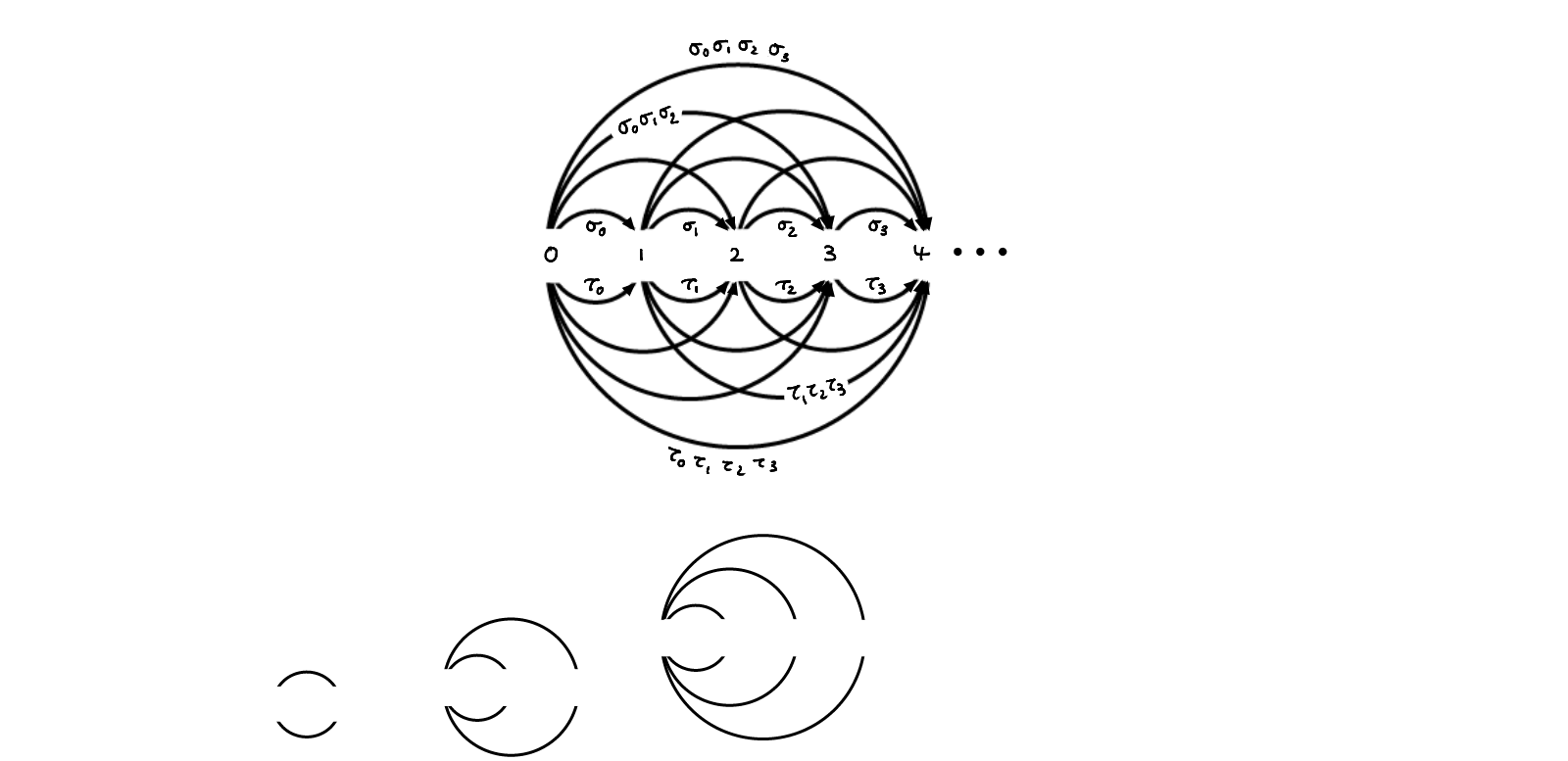}
\endgroup\end{restoretext}
We define the category of globular sets $\globset$ to be the presheaf category on $\lG$, that is,
\begin{equation}
\globset = [\lG\op ,\SetCat]
\end{equation}
A globular set $S$ is a contravariant functor on $\lG$
\begin{equation}
S : \lG\op \to \SetCat
\end{equation}
A map $\alpha : S_1 \to S_2$ of globular sets $S_1$, $S_2$ is a natural transformation $S_1 \to S_2$. 

As an example consider the globular set $S_a$ defined on objects and generating morphisms
\begin{restoretext}
\begingroup\sbox0{\includegraphics{test/page1.png}}\includegraphics[clip,trim=0 {.1\ht0} 0 {.0\ht0} ,width=\textwidth]{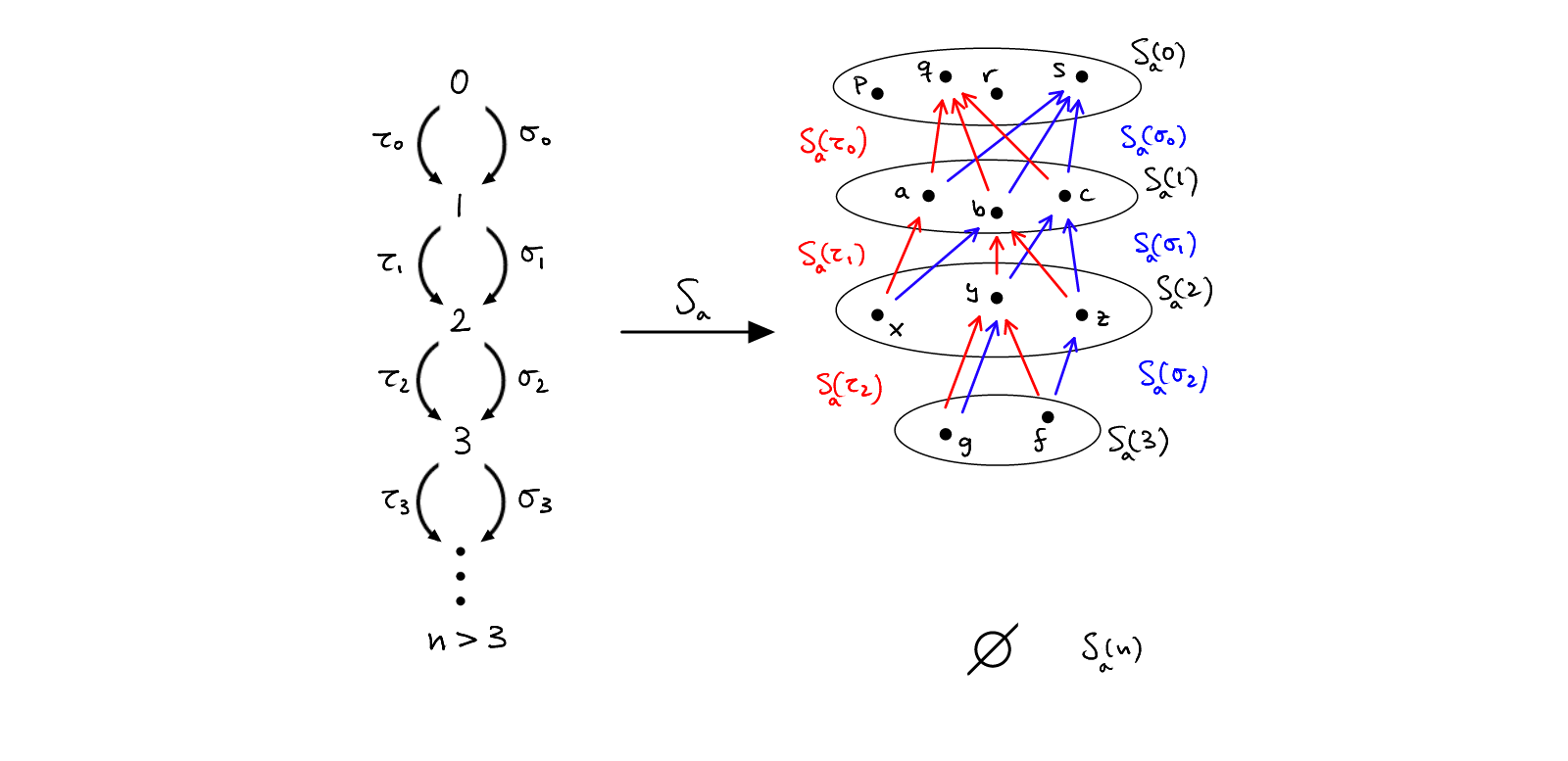}
\endgroup\end{restoretext}
Note that $S_a$ satisfies $S_a(n) = \emptyset$ for $n > 3$. Later such a globular set will be called $3$-truncated. 

The $n$-globes $\lGR n$ are special globular sets given by the hom-functors
\begin{equation}
\lGR n := \lG(-,n) : \lG\op \to \SetCat
\end{equation}
Note that $\lG(k,n)$ has exactly two elements if $k < n$: Namely, $\sigma_{k,n}$ (the equivalence class spanned by $\sigma_k \sigma_{k+1} ... \sigma_{n-1}$) and $\tau_{k,n}$ (the equivalence class spanned by $\tau_k \tau_{k+1} ... \tau_{n-1}$). If $k = n$, it has exactly one element, namely $\id_n$. If $k > n$ it has zero elements. The elements are often thought of us the components of the $n$-globe, for instance elements of $\lG(k,3)$ corresponds to components of the $3$-globe 
\begin{restoretext}
\begingroup\sbox0{\includegraphics{test/page1.png}}\includegraphics[clip,trim=0 {.25\ht0} 0 {.2\ht0} ,width=\textwidth]{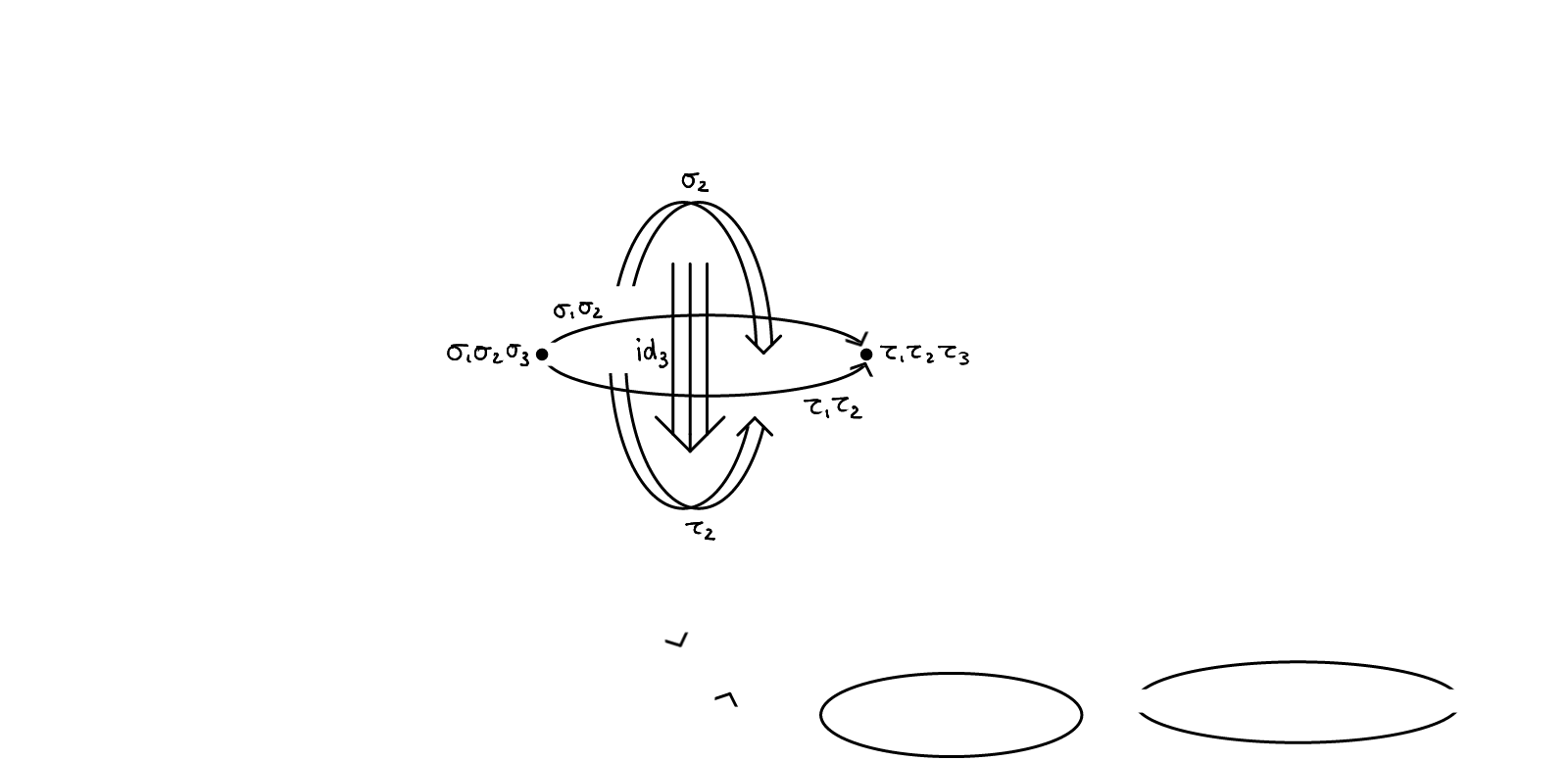}
\endgroup\end{restoretext}
\end{notn}

\begin{notn}[Sources and targets in globular sets] \label{notn:glob_set_src_tgt} Let $S$ be a globular set as defined above. We usually write $S_n = S(n)$, as well as $s_n = S(\sigma_n)$, $s_{n,k} = S(\sigma_{n,k})$ and $t_n = S(\tau_n)$, $t_{n,k} = S(\tau_{n,k})$. 
\end{notn}

\begin{constr}[Categories of elements of globular sets] \label{constr:elcat} Let $S : \lG \to \SetCat$ be a globular set. Using \autoref{rmk:set_valued_grothendieck_construction} we define the \textit{category of elements} $\elcat(S)$ by
\begin{equation}
\elcat(S) := \sG(\Discr S)
\end{equation} 
Explicitly, $\elcat(S)$ has elements of the form $(n \in \lN, g \in S_n)$ and morphism $(n,g) \to (k,h)$ consist of $f \in \lG(k,n)$ such that $S(f)(g) = h$. Using $S_a$ as defined above, its category of elements can be given by
\begin{restoretext}
\begingroup\sbox0{\includegraphics{test/page1.png}}\includegraphics[clip,trim=0 {.25\ht0} 0 {.2\ht0} ,width=\textwidth]{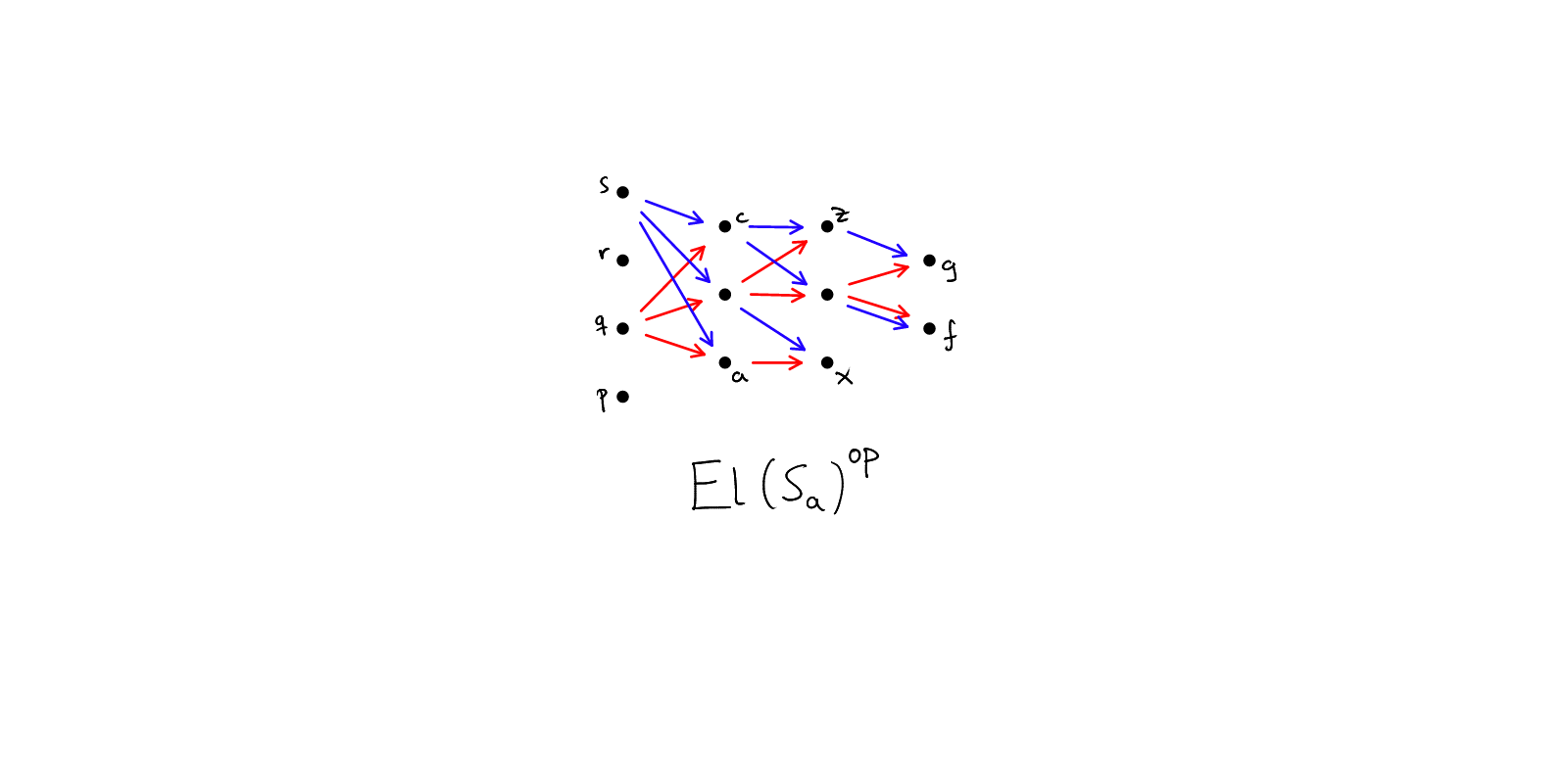}
\endgroup\end{restoretext}
where we only recorded the generating morphisms.

Note that in the case of $\lGR n$, since $\lG(k,n)$ has two elements, there are at most two morphism between any two objects in $\elcat(S)$. In particular we define
\begin{equation}
\cG^n := \elcat(\lGR n)\op
\end{equation}
and note that $\cG^n$ is a poset with objects $\cup_k \lG(k,n)$. For instance, the category $\cG^n$ is the poset
\begin{restoretext}
\begingroup\sbox0{\includegraphics{test/page1.png}}\includegraphics[clip,trim=0 {.3\ht0} 0 {.25\ht0} ,width=\textwidth]{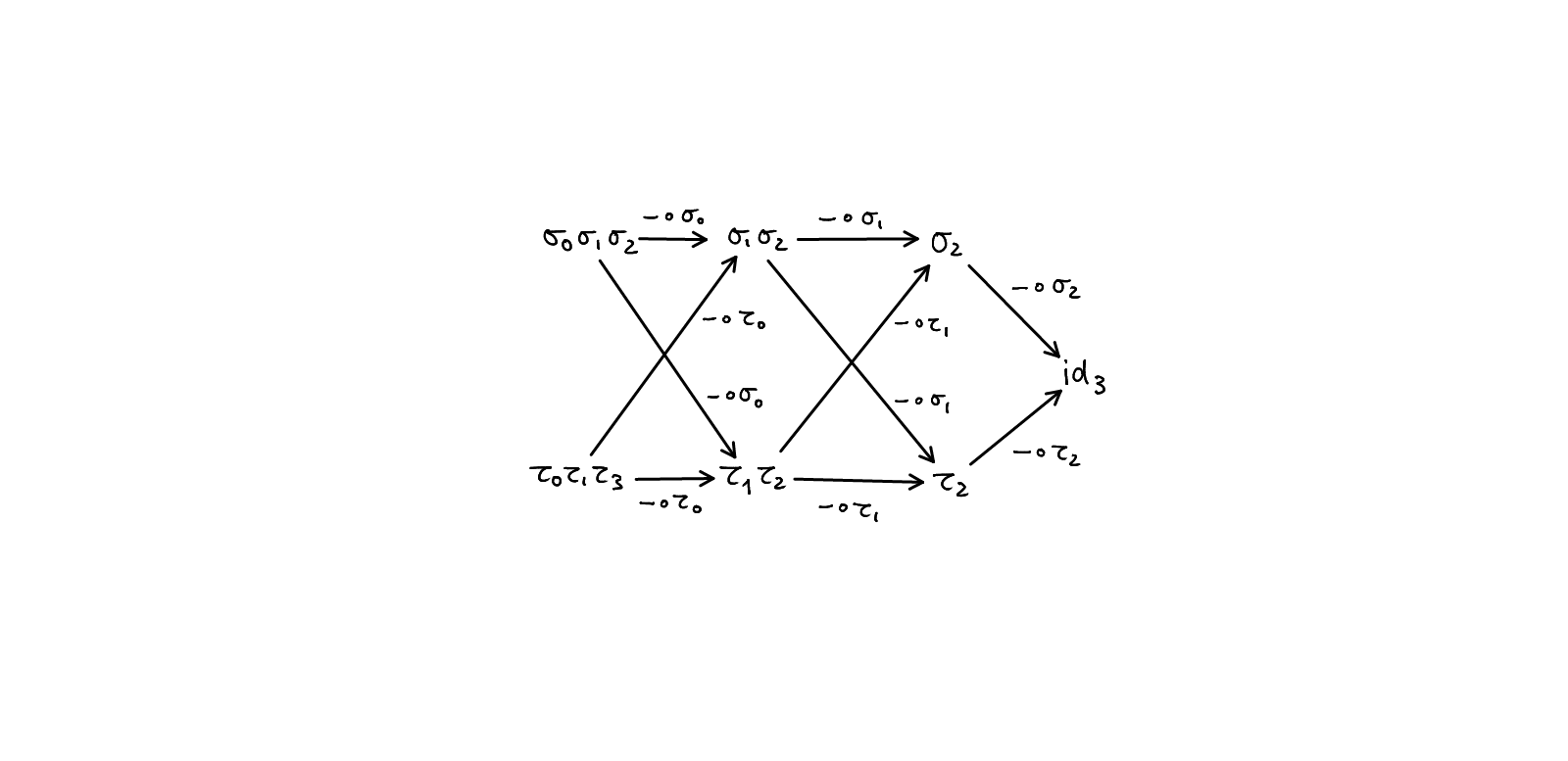}
\endgroup\end{restoretext}
which we will represent by colors as
\begin{restoretext}
\begingroup\sbox0{\includegraphics{test/page1.png}}\includegraphics[clip,trim=0 {.25\ht0} 0 {.35\ht0} ,width=\textwidth]{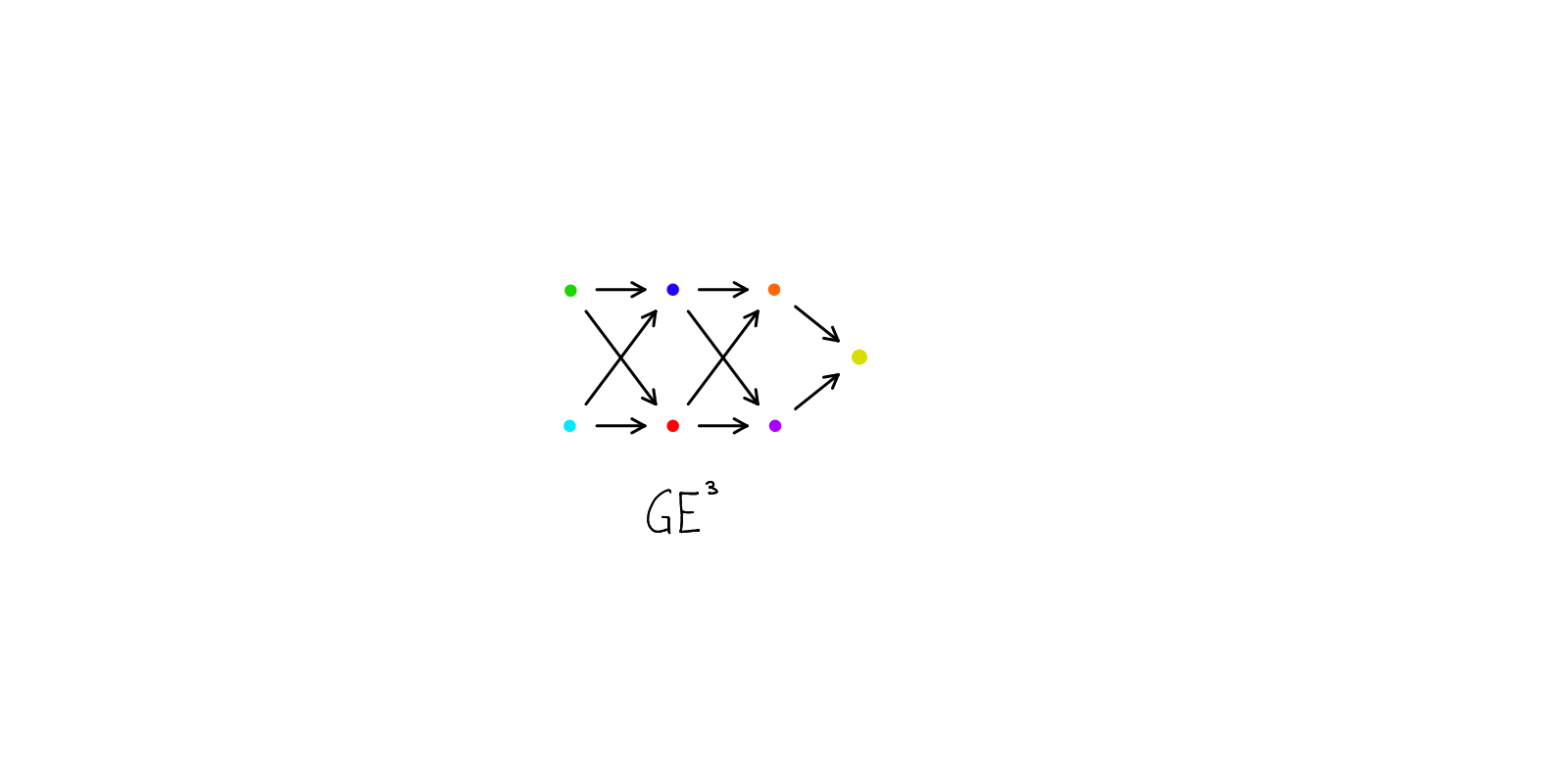}
\endgroup\end{restoretext}
Similarly we find for instance $\cG^2$, $\cG^1$ and $\cG^0$ to equal the posets
\begin{restoretext}
\begingroup\sbox0{\includegraphics{test/page1.png}}\includegraphics[clip,trim=0 {.1\ht0} 0 {.1\ht0} ,width=\textwidth]{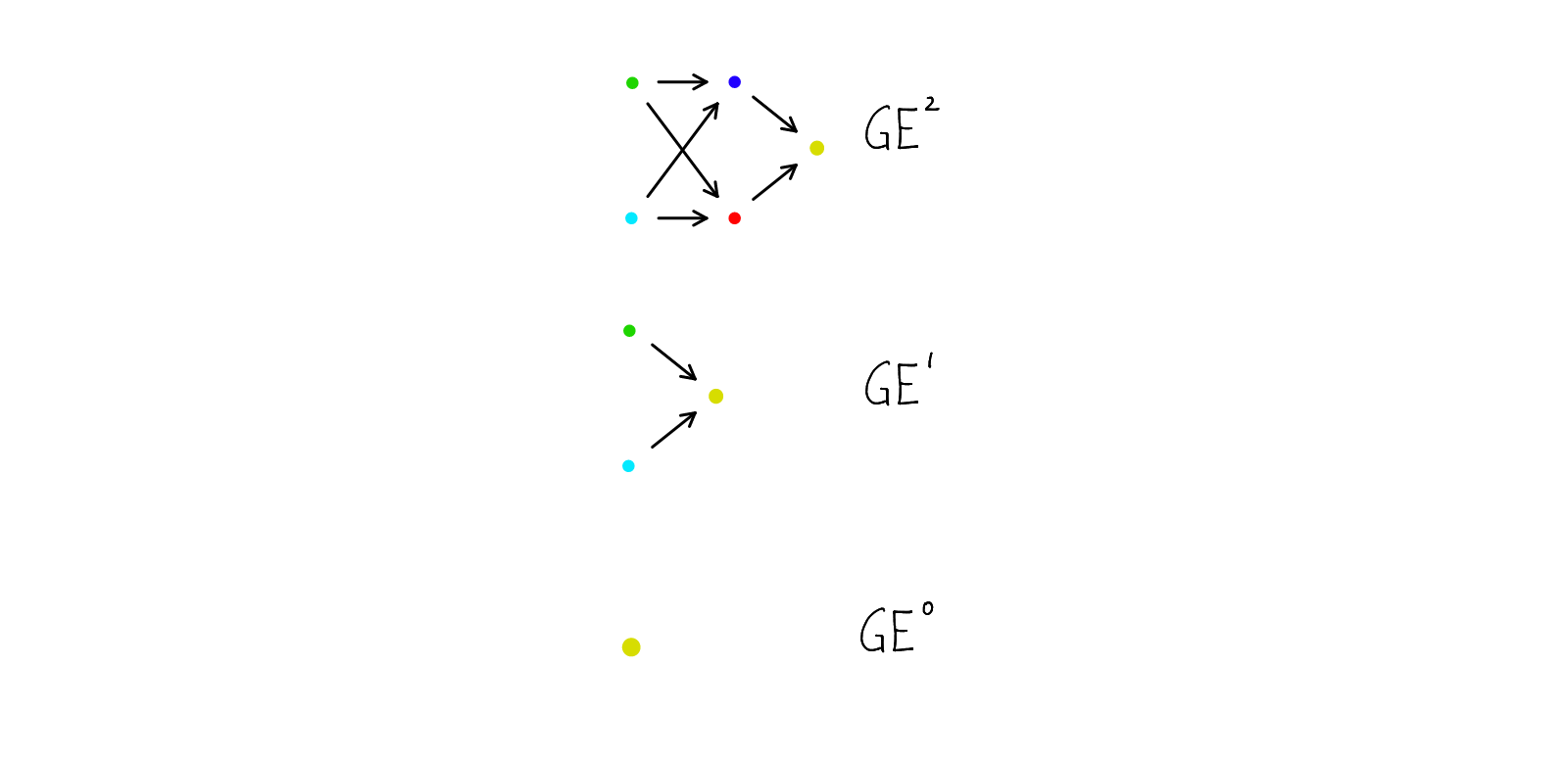}
\endgroup\end{restoretext}

Finally, given a map $\alpha : S_1 \to S_2$ of globular sets, this induces functor
\begin{equation}
\elcat(\alpha) : S_1 \to S_2
\end{equation}
which on objects maps
\begin{equation}
(n,g) \quad \mapsto \quad (n, \alpha_n(g))
\end{equation}
and on morphisms maps
\begin{align}
(f \in \lG(k,n)&, S_1(f)(g) = h) ~:~ (n,g) \to (k,h) \\
&\mapsto \quad (f \in \lG(k,n), S_2(f)(\alpha_n(g)) = \alpha_k(h)) ~:~ (n,\alpha_n(g)) \to (k,\alpha_k(h))
\end{align}
where we used naturality of $\alpha$ to infer $S_2(f)(\alpha_n(g)) = \alpha_k(h)$ from $S_1(f)(g) = h$. Note that for a map of globular sets $\beta : S_2 \to S_3$ we have
\begin{equation}\label{eq:elcat_compositional}
\elcat(\beta)\elcat(\alpha) = \elcat(\beta\alpha)
\end{equation}
\end{constr}

There is a canonical way of labelling certain cubes by $\cG^n$. Namely, these are $n$-cubes $\scA$ in which every fiber of all of its \SI-bundles is the terminal singular interval $\singint 1$, that is, $\tusU k_{\scA} = \const_{\singint 1}$ for $0 \leq k < n$. These cubes will be called terminal $n$-cubes, and we give the following examples. The terminal $0$-cube $\tgl^0$ is the $\SIvert 0 {\cG^0}$-cube with data
\begin{restoretext}
\begingroup\sbox0{\includegraphics{test/page1.png}}\includegraphics[clip,trim=0 {.45\ht0} 0 {.35\ht0} ,width=\textwidth]{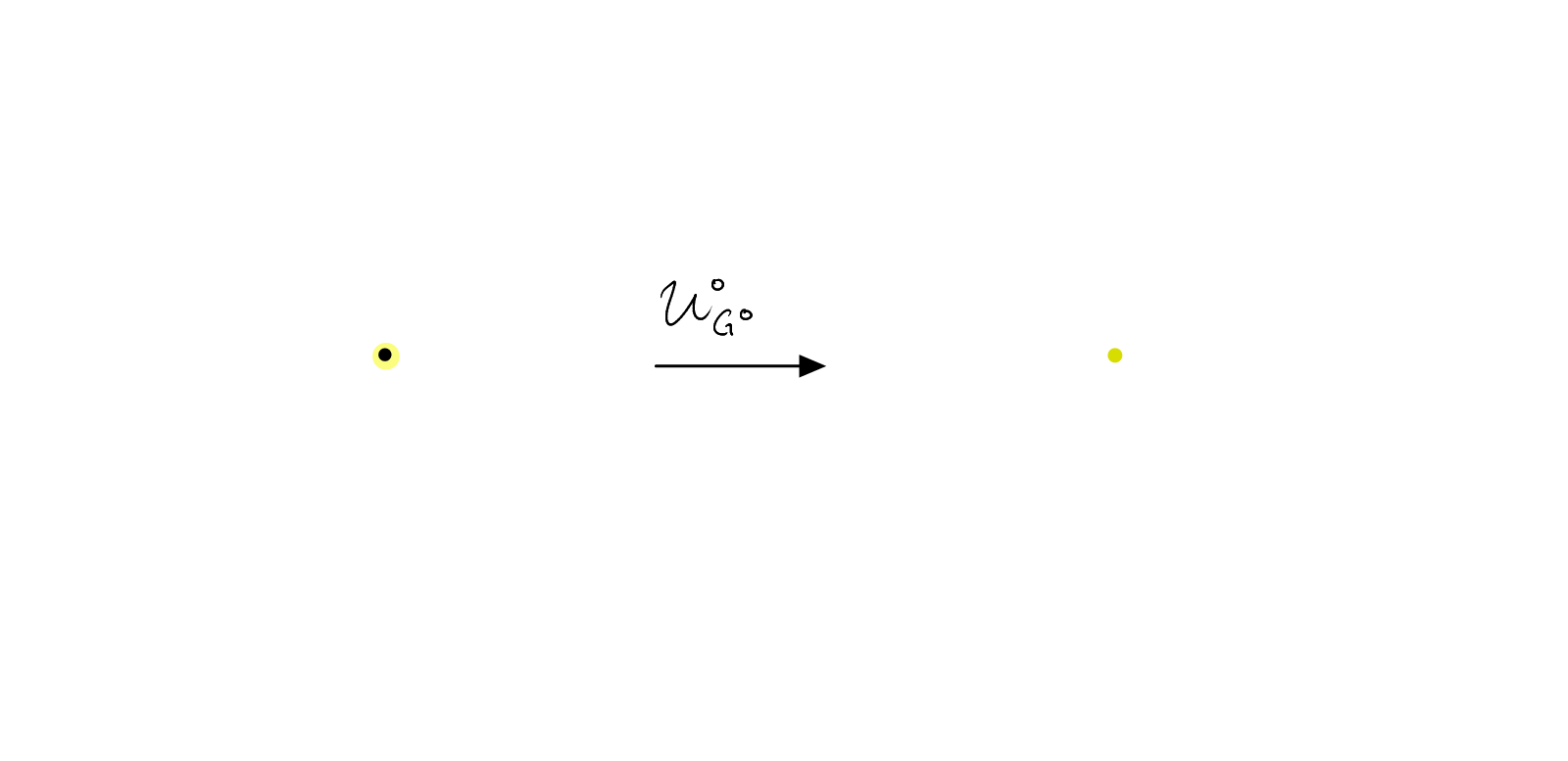}
\endgroup\end{restoretext}
The terminal $1$-cube $\tgl^1$ is the $\SIvert 1 {\cG^1}$-cube with data
\begin{restoretext}
\begingroup\sbox0{\includegraphics{test/page1.png}}\includegraphics[clip,trim=0 {.0\ht0} 0 {.2\ht0} ,width=\textwidth]{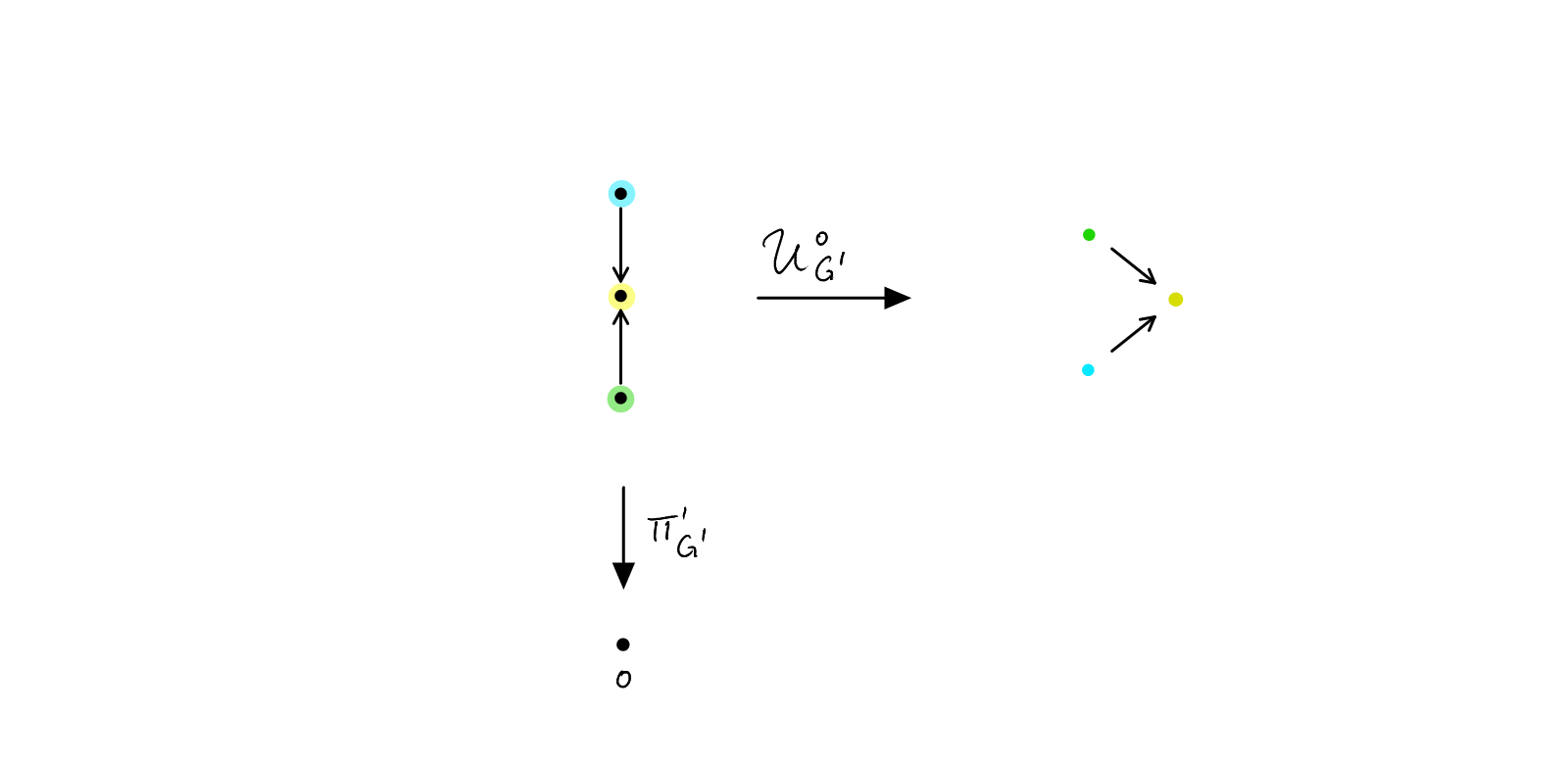}
\endgroup\end{restoretext}
The terminal $2$-cube $\tgl^2$ is the $\SIvert 2 {\cG^2}$-cube with data
\begin{restoretext}
\begingroup\sbox0{\includegraphics{test/page1.png}}\includegraphics[clip,trim=0 {.0\ht0} 0 {.0\ht0} ,width=\textwidth]{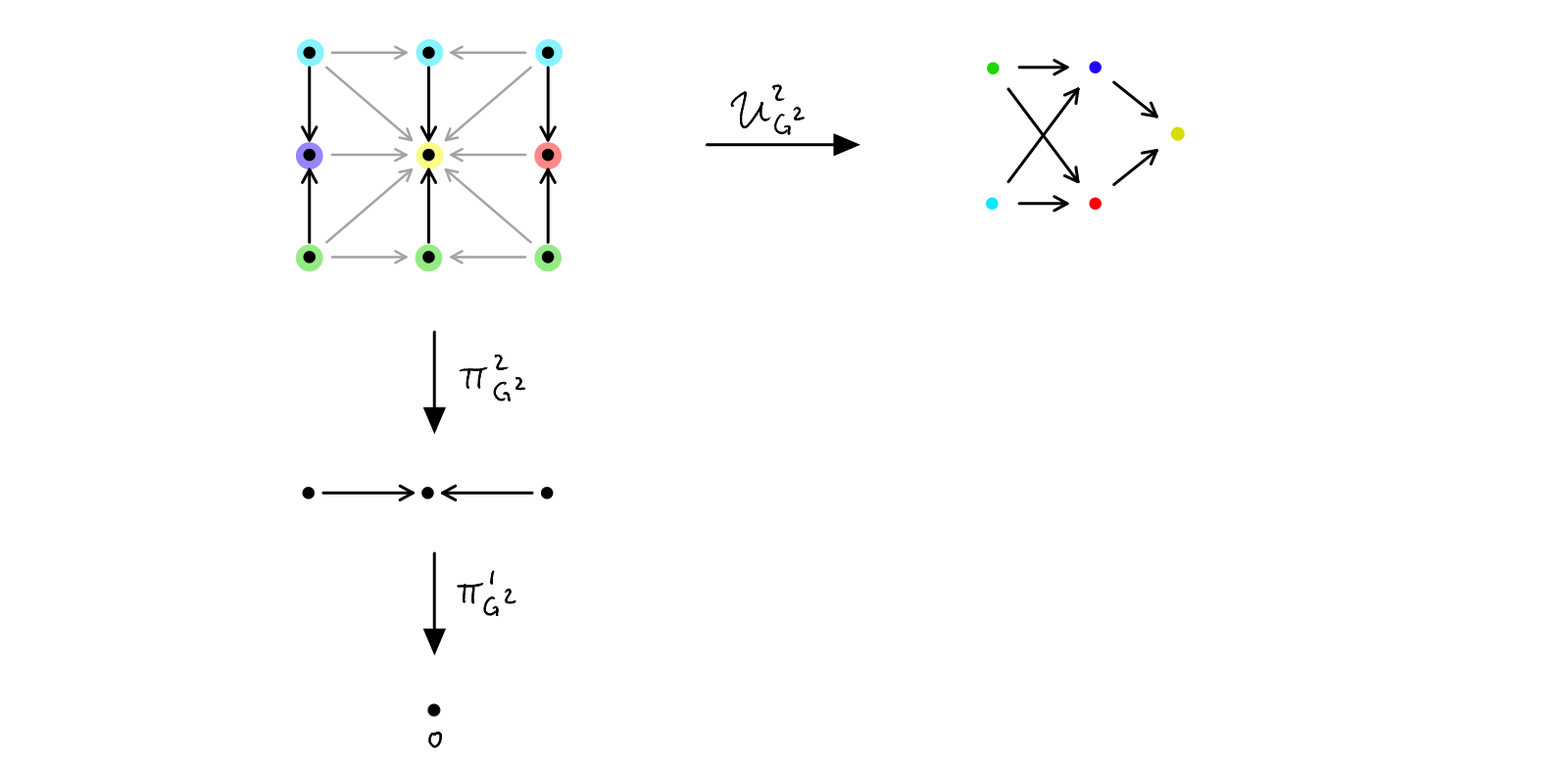}
\endgroup\end{restoretext}
and as a final example, note that the terminal $3$-cube $\tgl^3$ was in fact defined already as the cube $\scA$ in \autoref{eq:SIvert_n_C_families}.

The notion of terminal $n$-globe is formalised in the next construction.

\begin{constr}[Terminal $n$-cubes] \label{constr:terminal_n_globe} We construct the \textit{terminal $n$-globe} $G^n$, which is a globular and normalised $\cG^n$-labelled $n$-cube $\tgl^n : \bnum{1} \to \SIvert n {\cG^n}$. We proceed by induction on $n$

For $n=0$ we need to define $G^0 : \bnum{1} \to \cG^0$. Note that both $\bnum{1}$ and $\cG^0$ are the (terminal) one-object category ($\cG^0$'s only object is $\id_{0} \in \lG(0,0)$) and thus we are left with a unique choice for $\tgl^0$ which maps
\begin{equation}
\tgl^0 (0) = \id_{0}
\end{equation}
This definition makes $\tgl^0$ both globular and normalised.

For $n > 0$, we first define
\begin{align}
s &:= \elcat(\sigma_{n-1,n})\op : \cG^{n-1} \to \cG^n\\
t &:= \elcat(\tau_{n-1,n})\op : \cG^{n-1} \to \cG^n
\end{align}
where both $\sigma_{n-1,n}$ and $\tau_{n-1,n}$ are natural transformations by the Yoneda lemma
\begin{equation}
\lG(n-1,n) \iso \mathrm{Nat}(\lGR {n-1}, \lGR n)
\end{equation}
We then define $\tgl^{n} : \bnum 1 \to \SIvert {n+1} {\cG^n}$ by setting
\begin{equation} \label{claim:src_tgt_of_terminal_cube}
\abss{ \SIvert n s \tgl^{n-1} \xto {\id_n} \SIvert n t \tgl^{n-1} }
\end{equation}
where we used the double cone construction from the previous section, and thus in particular globular and normalised.
\end{constr}

\chapter{Foundation-independent higher categories} \label{ch:presented}

We will now define a first notion of $n$-categories, called presented associative $n$-categories. This is not a notion of higher categories in the classical sense (of morphisms sets with operations). Instead it captures the idea of \textit{presentations} of higher categories consisting of ``generators and relations" (as we will see generators and relations become part of the same structure in higher categories). The definition of presented associative $n$-categories is much more elementary and clear-cut than that of associative $n$-categories in the final \autoref{ch:associative}. In particular, no unbounded constructions will be used and (if we restrict ourselves to finite presentations) then presented associative $n$-categories are fully computer implementable and ``foundation-independent". We remark, that the same is to be expected from a ``full theory" of presented associative $n$-categories, involving functors, natural transformations, etc., which are the next steps in the research programme described in the introduction to this thesis, but lie beyond the scope of the present work.

In \autoref{sec:pres_def} we will define presented associative $n$-categories, and discuss a natural notion of morphisms (associated to such presentations). Further, we discuss how to \textit{adjoin} a generating morphisms to a given presentation which will make use of the previously defined double cone construction. 

In \autoref{sec:panc_ex} we will then discuss examples of associative $n$-categories for $n = -1,0,1,2$ and $3$. We will use these example to foreshadow ideas of the next chapter (discussing $n$-groupoids). Namely, we will include examples of truncated presentations of the $2$-sphere $S^2$ and perform elementary computations in its homotopy groups.

\section{\Free{} associative $n$-categories} \label{sec:pres_def}

\subsection{Definition}

We recall the following definitions and notation from \autoref{ssec:sum_namescopes}.

\begin{notn}[Full profunctorial relation] \label{rmk:full_pro_rel} Given two posets $X$ and $Y$ there is a profunctorial relation $R^{\mathrm{full}}_{X,Y}$, called the \textit{full profunctorial relation} between $X$ and $Y$, defined by setting $R^{\mathrm{full}}_{X,Y}(x,y)$ to be true for all $x \in X$ and $y \in Y$.
\end{notn}

\begin{defn}[Namescope and realisation] \label{defn:namescope}  A \textit{namescope} $\sN$ of dimension $n$, where $n \in \Set{-1} \cup \lN \cup \Set{\infty}$, is a list of sets $\sN_0, \sN_1, \sN_2, ... $ of length $(n+1)$. Elements of $\sN_k$ are called names or labels. Its  \textit{$k$-truncated full realisation}, for $1 \leq k \leq n+1$, is given by the functor
\begin{equation}
\Gamma^{k}_{\sN} : \bnum{k} \to \PRel
\end{equation}
mapping $i \mapsto \Discr \sN_i$ (cf. \autoref{rmk:set_valued_grothendieck_construction}), and $(i \to j)$ to $R^{\mathrm{full}}_{\Discr \sN_i,\Discr \sN_j}$ (cf. \autoref{rmk:full_pro_rel}).  Note that for $\bnum m : \bnum m \subset \bnum k$ we have $\Gamma^{k}_{\sN} \bnum m = \Gamma^{m}_{\sN}$. We write $\Gamma_{\sN} = \Gamma^{n+1}_{\sN}$. We also write $\GGamma{}\sN = \sG(\Gamma_{\sN})$ and $\GGamma k \sN = \sG(\Gamma^k_{\sN})$.

If $g \in \GGamma{}{\sN})$ then (identifying $g$ and $\secp g$ for total spaces of namescopes) there is $k \leq n$ with $g \in \sN_k$. We write $k = \dim(g)$, which is called the dimension of the label $g$.
\end{defn}

We now turn to the central definition of this chapter, that of presented associative $n$-categories. A concise version of this was given in \autoref{defn:sum_pres_anc}, in which morphisms and generator types where defined in a mutually inductive manner. Here we will take an alternative, less concise approach, which allows us to define presented associative $n$-categories only in terms of generators and separately from their morphisms.

\begin{defn}[Presented associative $n$-categories] \label{defn:pres_ANC} A \textit{\free{} associative $n$-category} $\sC$ is a namescope $\sC$ of dimension $(n+1)$ together with an assignment of \textit{types}. This means each $g \in \sC_m$ gets assigned 
\begin{equation}
\abss{g} : \bnum{1} \to \SIvert {\dim(g)} {\GGamma{}\sC}
\end{equation}
called its type, and which satisfies the following
\begin{itemize}
\item \textit{Minimality around $g$}: There is $\ip_g \in \tsG m(\abss{g})$, called the \textit{vertex point} of $\abss{g}$, which satisfies $\tsU m_{\abss{g}}(\ip_g) = g$ and
\begin{equation}
\abss{g} = \abss{g} \sslash \ip_g
\end{equation}
\item \textit{Correctness of dimension}: $\ip_g$ is of region co-dimension $m$ in $\abss{g}$, that is (cf. \autoref{defn:regions})
\begin{equation}
\ctypsum^m_{\abss{g}} (\ip_g) = m
\end{equation}
\item \textit{Well-typedness}: For all $p \in \tsG m(\abss{g})$ with $\tsU m_{\abss{g}}(p) = f \in \sC_l$, we have
\begin{equation}
\NF{~\abss{g} \sslash \ip_g~}^m = \Id^{m - l}_{\abss{f}}
\end{equation}
\end{itemize}
If $n < \infty$ and $m = n+1$ (that is, $g$ is a generating equality relation), then we also add the following condition.
\begin{itemize}
\item \textit{Proof irrelevance and invertibility of equalities}: For any $g, h \in \sC_{n+1}$ we have
\begin{equation}
\gsrc (\abss{h}) = \gsrc (\abss{g}) \text{~and~} \gtgt (\abss{h}) = \gtgt (\abss{g}) \quad \imp \quad h = g
\end{equation}
referred to as \textit{proof irrelevance} of equality, and further there must be $g\inv \in \sC_{n+1}$ with
\begin{equation}
\gsrc (\abss{g\inv}) = \gtgt (\abss{g}) \text{~and~} \gtgt (\abss{g\inv}) = \gsrc (\abss{g})
\end{equation}
referred to as \textit{invertibility} of equality.
\end{itemize}
Note that, if $m = 0$, then the above minimality condition implies $\abss{g} : \mathbf{1} \to \GGamma{}\sC$ has image $g$, that is, $\abss{g} = \Delta_g$. Note also that we sometimes write $g$ in place of $\abss{g}$ to lighten notation (especially if $m = 0$).
\end{defn}

\begin{rmk} We remark that this definition is well-founded in the sense that conditions for each $g \in \sC_i$ may only refer to $f \neq g$ which live in $\sC_j$ for $j < i$, i.e. $f$ of properly smaller dimension than $g$. In other words, the typability condition for $\abss{g}$ refers to $\abss{f}$ where $f$ is the label of some region in $\abss{g}$. But by the minimality condition we must have $f \in \GGamma{}\sC \sslash g$, and thus $\dim(f) < \dim(g)$ or $f = g$ (in the latter case the typability condition is trivial). Practically speaking, this means that \free{} associative $n$-categories can be build inductively by dimension as we will see in later examples.
\end{rmk}

\subsection{Category of $n$-category presentations}

The data necessary for the definition of a presented associative $n$-category is also called a \textit{presentation}. We now consider the question of how to map between presentations. The notion we define below is the simplest such notion of maps of presentations. However, it will become useful for defining presentations by colimits as discussed in the subsequent \autoref{rmk:colimit_of_presentations}. Note that our notion of maps of presentations should generalise to some notion of functors, but is strictly more special than that.

\begin{defn}[Maps of presentations] \label{defn:map_of_presentations} Let $\sC, \sD$ be \free{} associative $n$-categories. A \textit{map of presentations} $\alpha : \sC \to \sD$ is a bundle map $\alpha : \pi_{\Gamma_\sC} \to \pi_{\Gamma_\sD}$ (that is, a functor $\alpha : \GGamma{}{\sC} \to \GGamma{}{\sD}$ commuting with the bundles $\pi_{\Gamma_\sC}$ and $\pi_{\Gamma_\sD}$) such that for each $g \in \sC_k$ we have
\begin{equation} \label{eq:type_compatibility}
\SIvert k \alpha \abss{g} = \abss{\alpha(g)}
\end{equation}
If $\alpha$ is injective we call it an \textit{inclusion of presentations}. Note that $\alpha$ is determined by its components $\rest \alpha k : \sC_k \to \sD_k$ (cf. \autoref{notn:fiber_restrictions}) which are functions of sets.
\end{defn}

\begin{defn}[Category of $n$-category presentations] \label{defn:cat_of_n_pres}The category of $n$-category presentations $\pCat_n$ has as objects presented associative $n$-categories and as morphisms maps of presentations. Note that we canonically have $\pCat_n \into\pCat_m$ for $n < m$.
\end{defn}

\begin{rmk}[Colimits of presentations] \label{rmk:colimit_of_presentations} Given a (possibly transfinite) sequence of inclusions of presentations
\begin{equation}
\sC^0 \xinto {\alpha^1} \sC^1 \xinto {\alpha^2} \sC^2 \into ...
\end{equation}
we can form a \free{} associative $n$-category $\colim_i (\sC^i)$ such that the inclusions of presentations $\alpha^\infty_j : \sC^j \into \colim_i(\sC^i)$ have components $\rest {\alpha^\infty_j} k : \sC^j_k \into \colim_i (\sC^i)_k$ defined as inclusions into a colimit taken in \SetCat{} given by
\begin{equation}
\sC^0_k \xinto {\rest {\alpha^1} k} \sC^1_k \xinto {\rest {\alpha^2} k } \sC^2_k \into ...
\end{equation}
For types, note that given $g \in  \colim_i (\sC^i)_k$ then $g = \alpha^\infty_j (h)$ for some $j$ and $h$. We define
\begin{equation}
\abss{g} := \SIvert k {\alpha^\infty_j } \abss{h} 
\end{equation}
which is well-defined and satisfies \eqref{eq:type_compatibility} for all $\alpha^\infty_i$ since $\alpha^\infty_i \alpha^i = \alpha^\infty_{i-1}$.
\end{rmk}

\subsection{Composites of a \free{} associative category}

We now introduce a notion of morphisms, synonymously called composites\footnote{The terminology of ``composites" highlights that this is notion is associated to presentations and not classical categories---this is analogous to speaking of ``words" for the elements of a presented group.}, associated to the presentations of the previous sections.

\begin{constr}[Globular set of morphisms of \free{} associative $n$-categories] \label{defn:PANC_mor} Let $\sC$ be a \free{} associative $n$-category. We define an $(n+1)$-truncated globular set $\Comp(\sC)$ of \textit{morphisms} (or, synonymously, of \textit{composites}) as follows. Elements $\scD \in \Comp(\sC)_m$, $0 \leq m \leq (n+1)$, in degree $m$ are $\GGamma{}\sC$-labelled $m$-cubes, that is $\scD : \bnum{1} \to \SIvert m {\GGamma{}\sC}$, satisfying the following conditions
\begin{itemize}
\item $\scD$ is normalised
\item $\scD$ is globular
\item $\scD$ is well-typed: for all regions $p \in \tsG m(\scD)$ with $\tsU m_\scD (p)=g \in \sC_k$ we have
\begin{equation} \label{eq:typability}
\NF{\scD \sslash p}^m = \Id^{m-k}_{\abss{g}}
\end{equation}
\end{itemize}
The source and target maps $s_{m-1}, t_{m-1}$ are defined by
\begin{equation}
s_{m-1} \scD = \gsrc (\scD)
\end{equation}
and
\begin{equation}
t_{m-1} \scD = \gtgt (\scD)
\end{equation}
To see this is well-defined we need to check that $\gsrc (\scD)$ and $\gtgt (\scD)$ are globular, normalised and well-typed. The first two statements follow from \autoref{cor:vertical_degeneracy_src_tgt} and \autoref{lem:src_tgt_normalisation}. The last statement follows since both $\gsrc (\scD)$ and $\gtgt (\scD)$ are subfamilies of $\tsU 1_\scD$ and since minimal neighbourhoods $\scD \sslash p$ are minimal (cf. \autoref{claim:minimal_subfamily_is_minimal}).
\end{constr}

\begin{notn}[$(-1)$-morphisms] \label{notn:minus_one_mor} Recall \autoref{conv:minus_one_cubes}. As a consequence, it makes sense to also define a set of ``$(-1)$-morphisms" by setting
\begin{equation}
\Comp(\sC)_{-1} := \Set{\emptyset}
\end{equation}
This together with \autoref{notn:boundary_of_a_point} will facilitate inductive arguments later on.
\end{notn}

\begin{rmk}[Maps of presentations induce maps of morphisms] \label{rmk:mor_functor} Let $\alpha : \sC \to \sD$ be a map of presentations associative $n$-categories, and let $f \in \Comp(\sC)_m$. Then we find 
\begin{equation}
\SIvert m \alpha f \in \Comp(\sD)_m
\end{equation}
For this to hold we need to check that the latter is globular, normalised and well-typed. Being globular and normalised follow by \autoref{rmk:k_coll_label_transfer}. Well-typedness follows from \autoref{rmk:k_coll_label_transfer} and the observation that minimal subcubes commute with relabelling. This gives rise to a functor from presentation to globular sets:
\begin{equation}
\Comp : \pCat_n \to [\lG\op,\SetCat]
\end{equation}
\end{rmk}

The following construction formally points out the obvious fact that $m$-morphisms can be defined only knowing about $k$-generators for $k \leq m$.

\begin{constr}[Restricting $m$-morphism labels] \label{rmk:inductive_PANC_def} Let $\sC$  be a presented associative $n$-category. There is a natural correspondence of $m$-morphisms with globular, normalised, ``well-typed" $\GGamma m \sC$-labelled $m$-cubes.

We have a base change (recall the notation $\bnum m  : \bnum m \subset (\bnum {n+2})$ from \autoref{notn:subsets_and_restrictions})
\begin{equation}
\Gamma^m_\sC = \Gamma_\sC \bnum m
\end{equation}
which induces a map of total posets
\begin{equation}
\sG(\bnum m) : \GGamma m \sC \into \GGamma{}\sC
\end{equation}

Let $\scD \in \Comp(\sC)_m$ as defined above. Then $\scD$ being well-typed implies that for any $x \in \tsG m(D)$ we have $\tsU m_\scD(x) \in \sC_k$ for $k \leq m$. Thus, using \autoref{notn:restrictions_of_labels} we find $\scD$ factors as
\begin{equation}
\scD = \SIvert m {\sG(\bnum m)}  \sG(\bnum m)\pbstar  \scD
\end{equation}
Note that the cube
\begin{equation}
\sG(\bnum m)\pbstar  D : \bnum 1 \to \SIvert m {\GGamma m \sC}
\end{equation}
is globular and normalised. 

Conversely, we say that a globular, normalised $\scD : \bnum 1 \to \SIvert m {\GGamma m \sC}$ is \textit{well-typed} if $\SIvert n {\sG(\bnum m)} \scD$ is well-typed. The latter is the case if and only if for all regions $p \in \tsG m(\scD)$ with $\tsU m_\scD (p)=g \in \sC_k$ we have
\begin{equation}
\NF{\scD \sslash p}^m = \Id^{m-k}_{\sG(\bnum m)\pbstar \abss{g}}
\end{equation}
Using label restriction by $\sG(\bnum m)\pbstar $ and relabelling by $\SIvert n {\sG(\bnum m)} $, we thus find a 1-to-1 correspondence of globular, normalised, well-typed $\scD : \bnum 1 \to \SIvert m {\GGamma m \sC}$ and $m$-morphisms. 
\end{constr}

In other words, the definition of $m$-morphisms $\Comp(\sC)_m$ only depends on generators (and their types) of dimension $k \leq m$. We will make use of this fact in situations when $\Comp(\sC)_m$ is needed but $\sC_k$ (and $\abss{g}, g \in \sC_k$) has only been defined for $k \leq m$.

\subsection{Adjoining generators}

The goal of this section is to give a construction of new presentations obtained by adding, or ``adjoining", generators and relations to given presentations.

Recall \autoref{constr:double_cones_of_src_and_tgt}. Comparing with \autoref{defn:pres_ANC} we observe that types are double cones, that is,
\begin{equation}
\abss{g} = \abss{\gsrc{\abss{g}} \xto g \gtgt{\abss{g}}}
\end{equation}
The following construction shows the converse. Given well-typed sources and targets, we can define a well-typed type for a new generator using the double cone construction. The crucial technical difficulty is to show that the double cone is well-typed, which we will do in the following.

\begin{constr}[Adjoining a generating $k$-morphism] \label{constr:adjoining_gen} Let $\sC$ be an \free{} associative $n$-category, let $\ig$ be a new name, and assume $x, y  \in \Comp(\sC)_{k-1}$ for $k < n$, such that $x$ and $y$ have the same sources and targets. In this case we define $(\sC \gadd{x,y} \ig)$, \textit{the category obtained by adjoining a generator $\ig$ between $x$ and $y$}, as follows. We set
\begin{equation}
(\sC \gadd{x,y} \ig)_i = \begin{cases}
\sC_k \cup \Set{\ig} & \text{~if~}  i = k \\
\sC_k & \text{~if~}  i \neq k 
\end{cases}
\end{equation}
There is a canonical injective bundle map $\alpha: \pi_{\Gamma_{\sC}} \into \pi_{\Gamma_{\sC \gadd{x,y} \ig}}$. For all $g \in (\sC \gadd{x,y} \ig)_i$ with $g \neq \ig$ types $\abss{g}$ are then inherited by types of generators $g \in \sC_i$ in $\sC$ such that this inclusion becomes an inclusion of presentations, that is, such that \eqref{eq:type_compatibility} is satisfied.

On the other hand, if $g = \ig$ then we set
\begin{equation}
\abss{g} = \abss{x \xto \ig y}
\end{equation}
As shown in \autoref{constr:top_monad} this is globular and normalised. It also satisfies minimality and correctness of dimension by \autoref{constr:double_cones_of_src_and_tgt} as required in \autoref{defn:pres_ANC} for notation. It remains to verify $\abss{g}$ is well-typed, which we prove as follows. Note that the proof is \stfwd{} albeit uses quite complex notation.

Let $p \in \tsG k(\abss{g})$ with $p \neq \ip_g$, and set $f = \tsU k_{\abss{g}}(p) \in \sC_l$. We need to show that
\begin{equation}
\vec\nfc{\abss{g}\sslash p} : \abss{g} \sslash p \starcoll \Id^{k-l}_{\abss{f}}
\end{equation}

Let $j \geq 1$ be minimal such that $p^j$ is a regular height. From this choice of $j$ and by inspecting \autoref{constr:top_monad}, we infer that $\tsU {j-1} (p^{j-1}) = \singint 1$, and thus we have either $\secp p^j = 0$ or $\secp p^j = 2$. We argue in the former case (the latter is similar, replacing sources by targets in the following). Recall the embedding 
\begin{equation}
\restemb_{\msrc_{\tusU {j-1}_{\abss{g}}}} : \csrc^j\abss{g} \mono \tsU j_{\abss{g}}
\end{equation}
By assumption on $p$ we have $p = \restemb^{k-j}_{\msrc_{\tsU {j-1}_{\abss{g}}}} (q)$ for some $q \in \tsG {k-j}(\csrc^j\abss{g})$. We define
\begin{equation}
\scA_j = \tsR {j-1}_{\sT^{j-1}_{\abss{g}}, \Id_{ \csrc^j \abss{g} }}
\end{equation}
and note that
\begin{equation}
\theta : \scA_j \mono \abss{g}
\end{equation}
where $\theta^i = \id$ for $i < j$ and $\theta^i = \restemb^{i-j}_{\msrc_{\tusU {j-1}_{\abss{g}}}}$ otherwise.

Note that by globularity and normalisation, $\csrc^j \abss{g}$ is constant, namely we have $\csrc^j \abss{g} = \const_{\gsrc^j\abss{g}}$. Thus there is an ordered collapse sequence
\begin{equation}
\vvec\mu : \scA_j \starcoll \Id^{j}_{\gsrc^j \abss{g}}
\end{equation}
where $\vvec\mu^i = \id$ for $i \geq (j -1)$ and $\vvec\mu^i : \const_{\singint 0} \to \tusU {i}_{\abss{g}}$ for $i < (j - 1)$ (recall that $\singint 0$ is initial in \SI).

Now, using minimality of $\iota^p_{\abss{g}}$ we first find
\begin{equation}
\iota^p_{\abss{g}} = \theta \iota^q_{\scA_j}
\end{equation}
and in particular
\begin{equation}
\abss{g} \sslash p = \scA_j \sslash q
\end{equation}
Next, setting $r = (\vsS{\vvec\mu})^k(q)$ and using \autoref{rmk:restrction_of_coll_on_minimal} we find
\begin{equation}
(\iota^q_{\scA_j})\pbstar \vvec\mu : \scA_j \sslash q \starcoll (\Id^j_{\gsrc^j \abss{g}} \sslash r)
\end{equation}
Finally, using well-typedness of $\gsrc^j \abss{g}$ we find
\begin{equation}
(\Id^j_{\gsrc^j \abss{g}} \sslash r) \starcoll \Id^{j + (k - j) - l}_{\abss{f}}
\end{equation}
Combining these facts proves the statement that we set out to prove.
\end{constr}

\begin{constr}[Adjoining sets of generators] \label{constr:adjoining_sets_of_gen} Let $\sC$ be a \free{} associative $n$-category and recall \autoref{rmk:colimit_of_presentations} about colimits of presentations. For $i\in \Set{1,2}$ let $\ig_i$ be names and let $x_i, y_i \in \Comp(\sC)_{k_i}$ ($k_i < n$) such that $x_i$ and $y_i$ have the same source and target. Note that in this case the previous construction yields
\begin{equation}
((\sC \gadd{x_1,y_1} \ig_1) \gadd{x_2,y_2} \ig_2) = ((\sC \gadd{x_2,y_2} \ig_2) \gadd{x_1,y_1} \ig_1)
\end{equation}
Consequently, given a label set $I$ and pairs $(x_\ig,y_\ig)$, $\ig \in I$ of $k_\ig$-morphisms in $\sC$ ($k_\ig < n$) such that $x_\ig$ and $y_\ig$ have the same source and target, then the colimit of the (possibly transfinite) sequence
\begin{equation}
\sC \into (\sC \gadd{x_{\ig_1},y_{\ig_1}} \ig_{1})\into ((\sC \gadd{x_{\ig_1},y_{\ig_1}} \ig_{1}) \gadd{x_{\ig_2},y_{\ig_2}} \ig_{2}) \into \dots
\end{equation}
is independent of any ordering $\ig_1, \ig_2, ...$ of elements in $I$. We will denote this colimit by
\begin{equation}
\sC \gadd{x_\ig,y_\ig} \Set{\ig \in I}
\end{equation}
\end{constr}

The above two constructions only concern the addition of generating morphisms. To adjoin a generating equality we proceed similarly, but need to consider proof irrelevance and symmetry. In particular, with the additional convention that \textit{all} generating equalities of a \free{} associative $n$-category $\sC$ have names of the form $\ie_{x = y} \in \sC_{n+1}$ for $x,y \in \Comp(\sC)_n$ we can give the following construction.

\begin{constr}[Adjoining generating equalities] \label{constr:adjoin_eq}  Fix $x,y \in \Comp(\sC)_n$ with the same globular source and target. Then $\sC \gadd{x,y} \ie_{x = y}$ equals $\sC$ if we already have $\ie_{x = y} \in \sC$. Otherwise, $\sC \gadd{x,y} \ie_{x = y}$ is defined by
\begin{equation}
(\sC \gadd{x,y} \ie_{x = y})_i = \begin{cases}
\sC_k \cup \Set{\ie_{x = y}, \ie_{y = x}} & \text{~if~}  i = n+1 \\
\sC_k & \text{~if~}  i < n + 1 
\end{cases}
\end{equation}
The rest of the discussion is analogous to \autoref{constr:adjoining_gen} and \autoref{constr:adjoining_sets_of_gen} can then be extended to adjoining sets of pairs of both morphisms and equalities.
\end{constr}

\subsection{Presented associative $\infty$-categories from globular sets} \label{ssec:pres_to_globe}

Recall from \autoref{ssec:sum_panc_from_gset} that presented associative $\infty$-categories can be canonically constructed from globular sets. In this section we will discuss this statement in more detail.

\begin{constr}[\Free{} associative $n$-categories from globular sets] \label{constr:panc_from_glob} Let $S$ be a globular set. We define a presented associative $\infty$-category $\kC(S)$ from $S$ by first setting $\kC(S)_k = S_k$ for $k \geq 0$. Using \autoref{notn:globeset} we find a functor
\begin{equation}
\elpsk_S : \elcat(S)\op \to \GGamma{} {\kC(S)}
\end{equation}
which is defined by acting as the identity on objects\footnote{There is a different more abstract (and much less readable) way to define this auxiliary functor. Recall (e.g. from \autoref{sec:poset}) that the inclusion $\Pos \into \Cat$ has a left adjoint $\PSk : \Cat \to \Pos$ called the posetal skeleton (which maps a Hom-sets to the truth-value of their non-emptiness and identifies objects antisymmetrically), and let $\mathsf{\mu}$ be the unit of this adjunction. Recall \autoref{defn:sum_red_lab_cat} (the ``minimal labelling category" $i : \redGamma{\kC(S)} \into \GGamma{}{\kC(S)}$ containing only the ``necessary" arrows). We now note that $\PSk(\elcat(S)\op))$ is precisely $\redGamma{\kC(S)}$. Then we can write
\begin{equation}
\elpsk_S = \PSk(i) \mathsf{\mu}_{\elcat(S)\op} : \elcat(S)\op \to \GGamma{}{\kC(S)}
\end{equation}}. We then set
\begin{equation}
\abss{g} = \SIvert n {\elpsk_S} \SIvert n {\elcat(g)\op} \tgl^n : \bnum 1 \to \SIvert k {\GGamma{} \sC}
\end{equation}
where we employed the Yoneda lemma to find $g : \lG(-,n) \to S$ and then used functoriality of $\elcat(-)$ and $\SIvert n {-}$ (cf. \autoref{rmk:SIvert_endofunctor}) to find
\begin{equation}
\SIvert k {\elcat(g)\op} : \SIvert k {\cG^k} \to \SIvert k {\elcat(S)\op}
\end{equation}
We need to check that $\abss{g}$ satisfies the three conditions in \autoref{defn:pres_ANC}. The first two conditions are preserved under relabelling and thus follow from the properties of $\tgl^k$ (see \autoref{constr:terminal_n_globe}). The final condition of well-typedness now follows inductively using \autoref{constr:adjoining_gen}, since the source (resp. target) of $\abss{g}$ are the source (resp. target) of $\tgl^k$ up to relabelling and thus inductively equal $\abss{s(g)}$ (resp. $\abss{t(g)}$).

This completes the construction of $\kC(S)$.
\end{constr}

In words, $\kC(S)$ is the presented associative $\infty$-category whose generating $k$-morphisms are the elements of $S_k$ with types determined by the sources and targets in $S$. 

\begin{rmk}[$n$-truncated case] Note that if we start with an $n$-truncated globular set $S$ then the preceding constructions yield a presented associative $n$-category (without equality relations).
\end{rmk}

\begin{constr}[The functor $\kC$] Note that the previous construction canonically gives rise to a functor
\begin{equation}
\kC : [\lG\op,\SetCat] \to \pCat_\infty
\end{equation}
which maps $S \mapsto \kC(S)$ on objects, and for $\alpha : S_1 \to S_2$ produces a map of presentations $\kC(\alpha) : \kC(S_1) \to \kC(S_2)$ defined on objects to be the function
\begin{equation}
\kC(\alpha) := \elcat(\alpha) : \GGamma{}{\kC(S_1)} \to \GGamma{}{\kC(S_1)}
\end{equation}
This satisfies \eqref{eq:type_compatibility} as verified by direct computation: namely, (on objects)
\begin{equation}
\elcat(\alpha) \elpsk_S \elcat(g) = \elcat(\alpha(g))
\end{equation}
This completes the construction of the functor $\kC$.
\end{constr}

\begin{eg}[Presented associative $\infty$-categories from globular sets] \label{eg:panc_from_gset}
Take the globular set $S$
\begin{restoretext}
\begingroup\sbox0{\includegraphics{test/page1.png}}\includegraphics[clip,trim=0 {.2\ht0} 0 {.14\ht0} ,width=\textwidth]{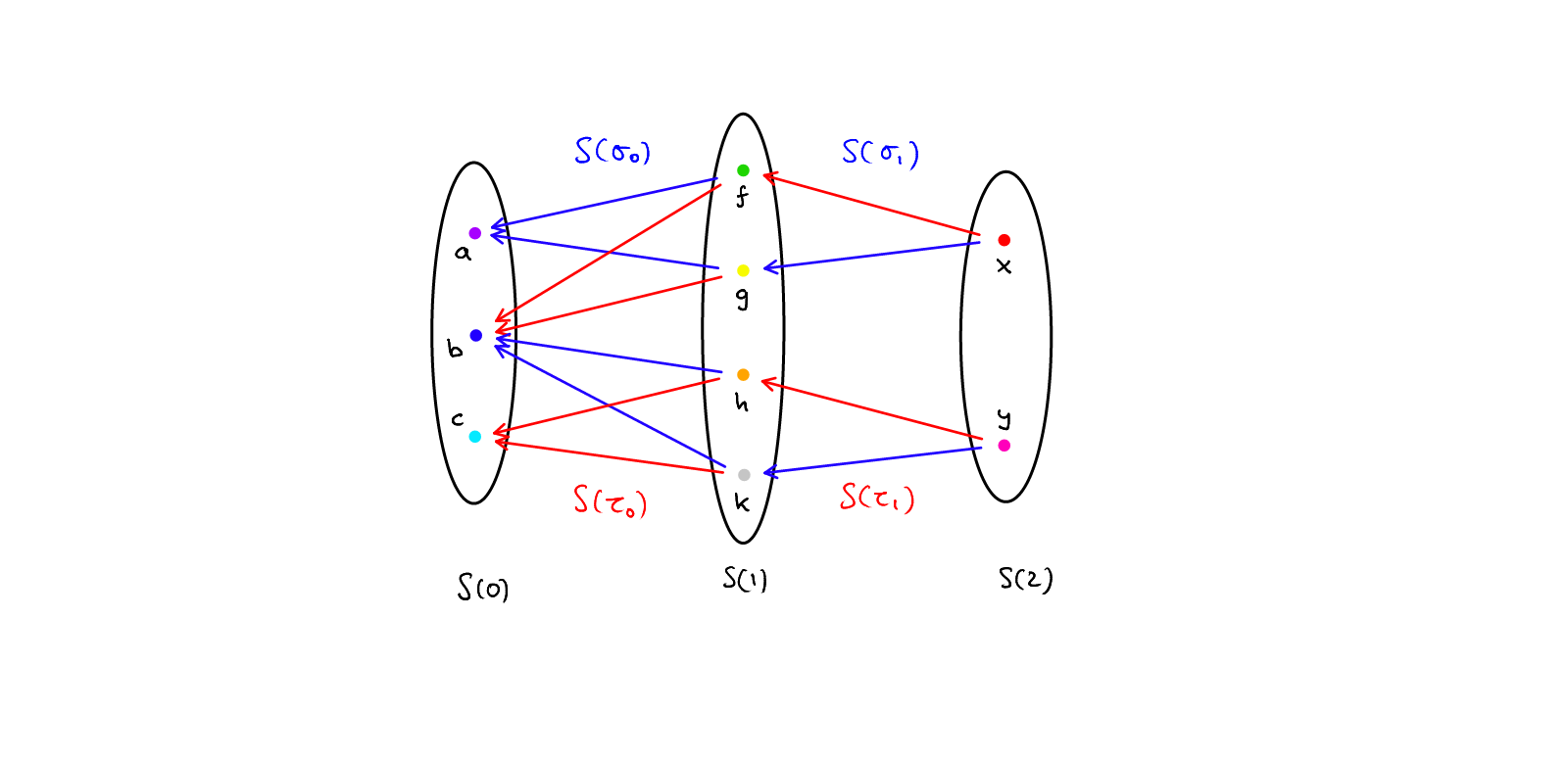}
\endgroup\end{restoretext}
which thus has category of elements
\begin{restoretext}
\begingroup\sbox0{\includegraphics{test/page1.png}}\includegraphics[clip,trim=0 {.18\ht0} 0 {.2\ht0} ,width=\textwidth]{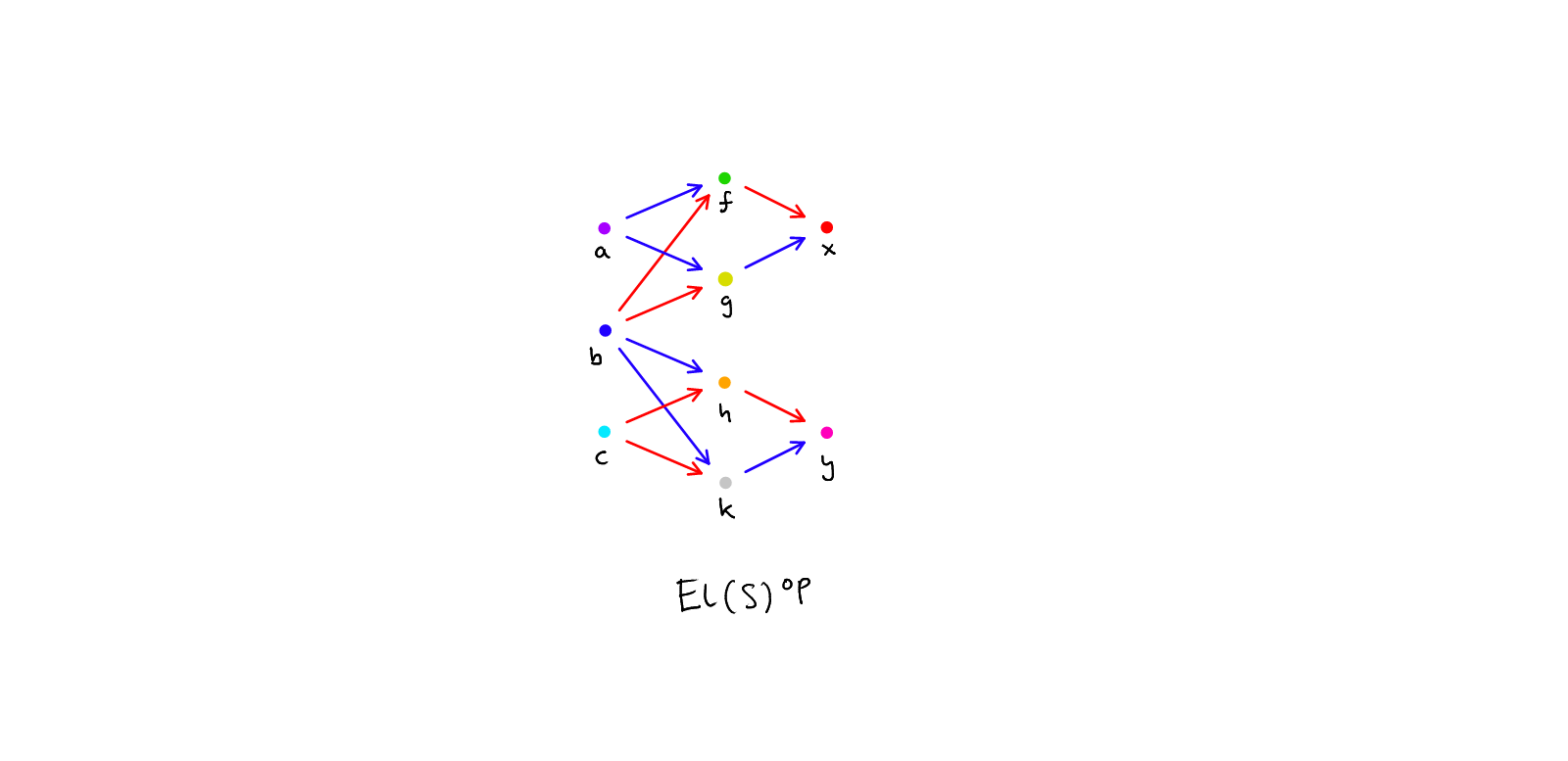}
\endgroup\end{restoretext}
We want to construct $\kC(S)$. First note, $\GGamma{}{\kC(S)}$ is the poset generated by
\begin{restoretext}
\begingroup\sbox0{\includegraphics{test/page1.png}}\includegraphics[clip,trim=0 {.18\ht0} 0 {.22\ht0} ,width=\textwidth]{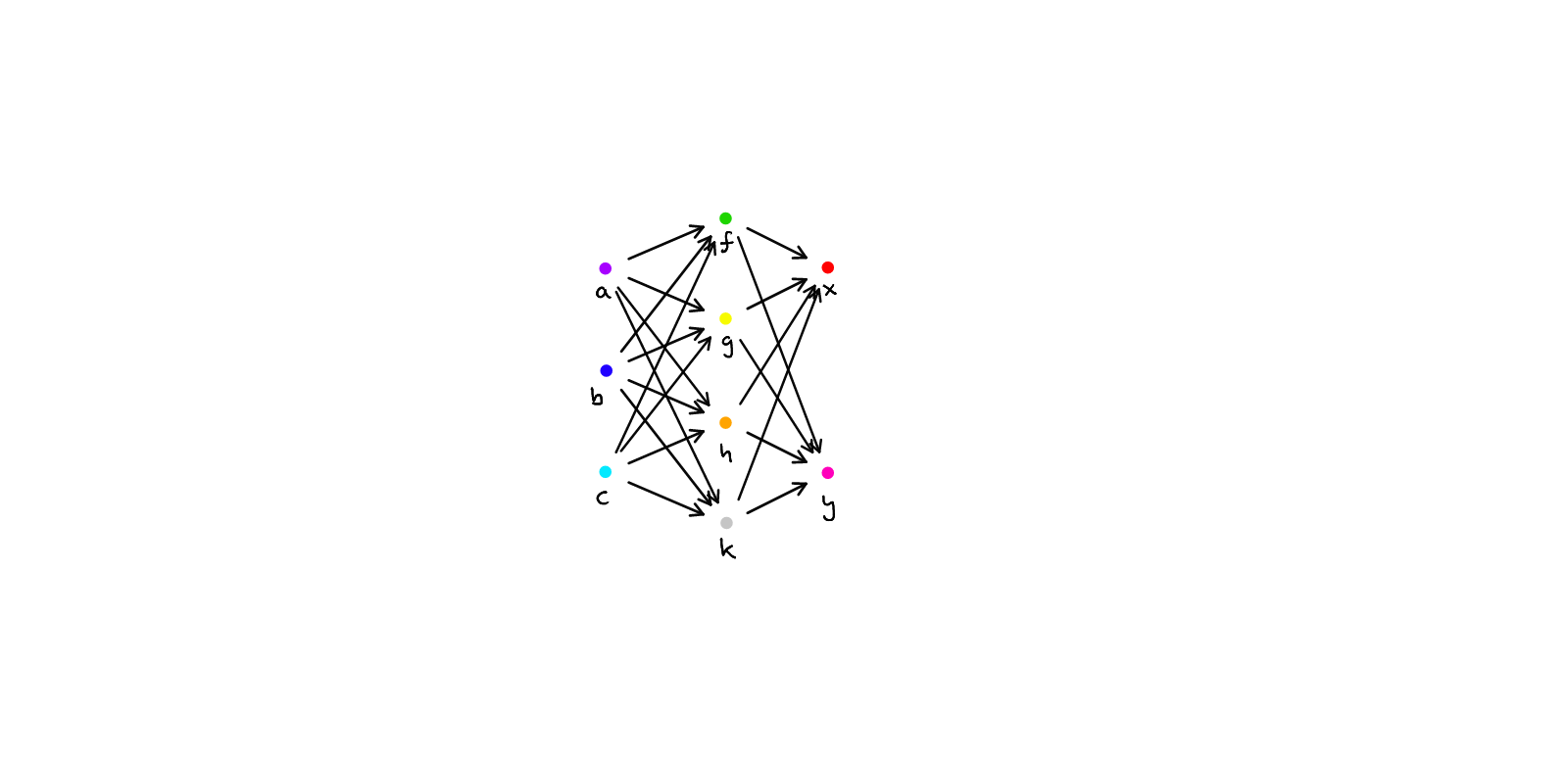}
\endgroup\end{restoretext}
$\elpsk_S$ simply forgets about the distinction of source and target maps, and acts as the identity on objects.

Now, by the Yoneda lemma we have an identification
\begin{equation}
S_n \iso \mathrm{Nat}(\lG(-,n), S)
\end{equation}
and thus (abusing notation) $g$ corresponds to a natural transformation $g \in \mathrm{Nat}(\lG(-,n), S)$. Using \autoref{constr:elcat} we find
\begin{equation}
\elcat (g) : \cG^n \to \elcat(S)
\end{equation}
We then define $\abss g$ to be the composite
\begin{equation} \label{eq:abss_g}
\abss{  g} \quad := \quad \big(\bnum{1} \xto {\tgl^n } \SIvert{n}{\cG^n} \xto {\SIvert{n}{\elcat (g)\op}} \SIvert{n}{\elcat(S)\op} \xto {\SIvert n {\elpsk_S}} \SIvert n {\GGamma{} {\kC(S)}} \big)
\end{equation}

Concretely, for our choice of $S$ we find for instance $\abss{g}$ to be the $\SIvert 1 {\kC(S)}$-family
\begin{restoretext}
\begingroup\sbox0{\includegraphics{test/page1.png}}\includegraphics[clip,trim=0 {.2\ht0} 0 {.08\ht0} ,width=\textwidth]{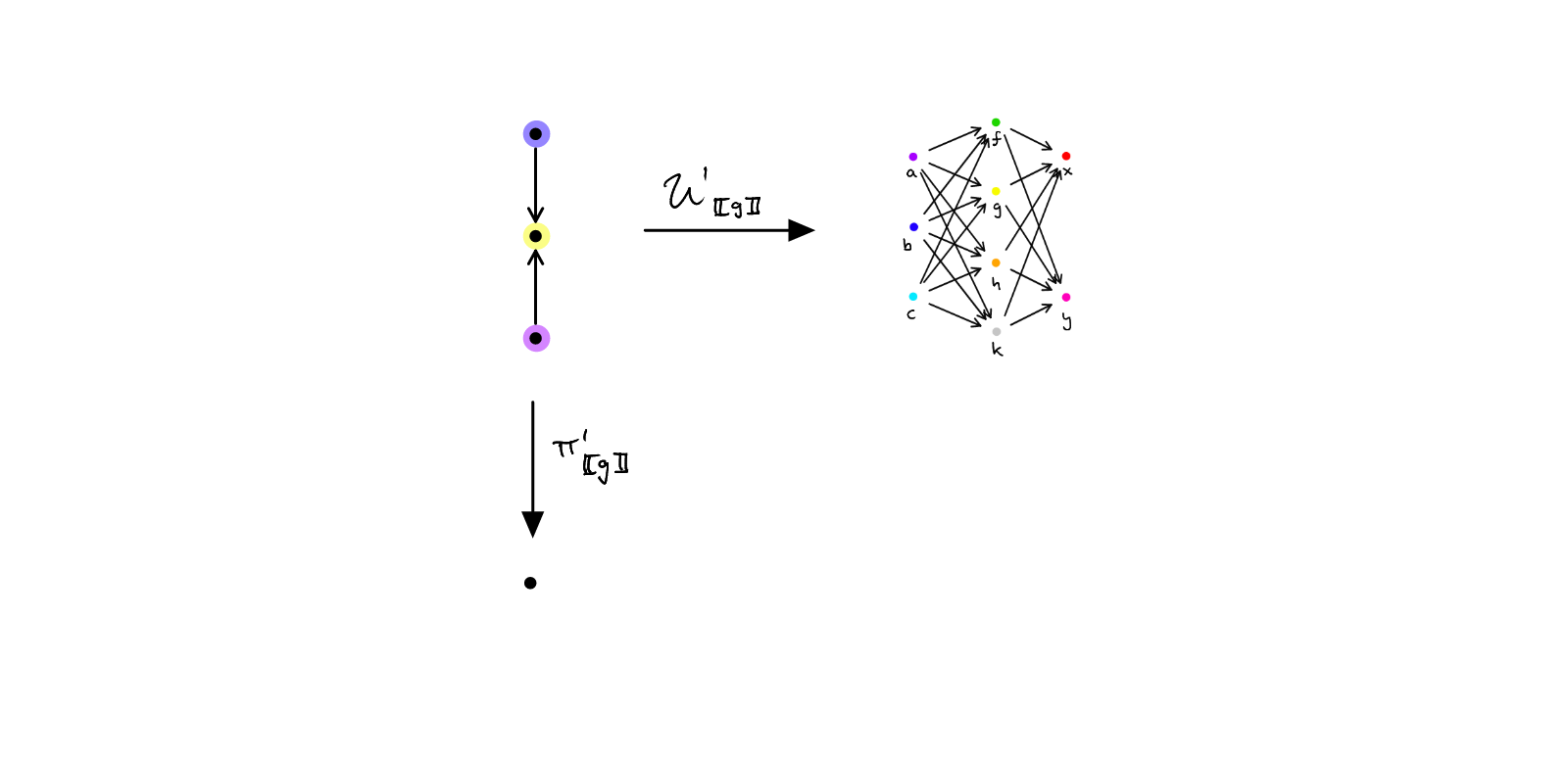}
\endgroup\end{restoretext}
$\abss{h}$ is the $\SIvert 1 {\kC(S)}$- family
\begin{restoretext}
\begingroup\sbox0{\includegraphics{test/page1.png}}\includegraphics[clip,trim=0 {.2\ht0} 0 {.09\ht0} ,width=\textwidth]{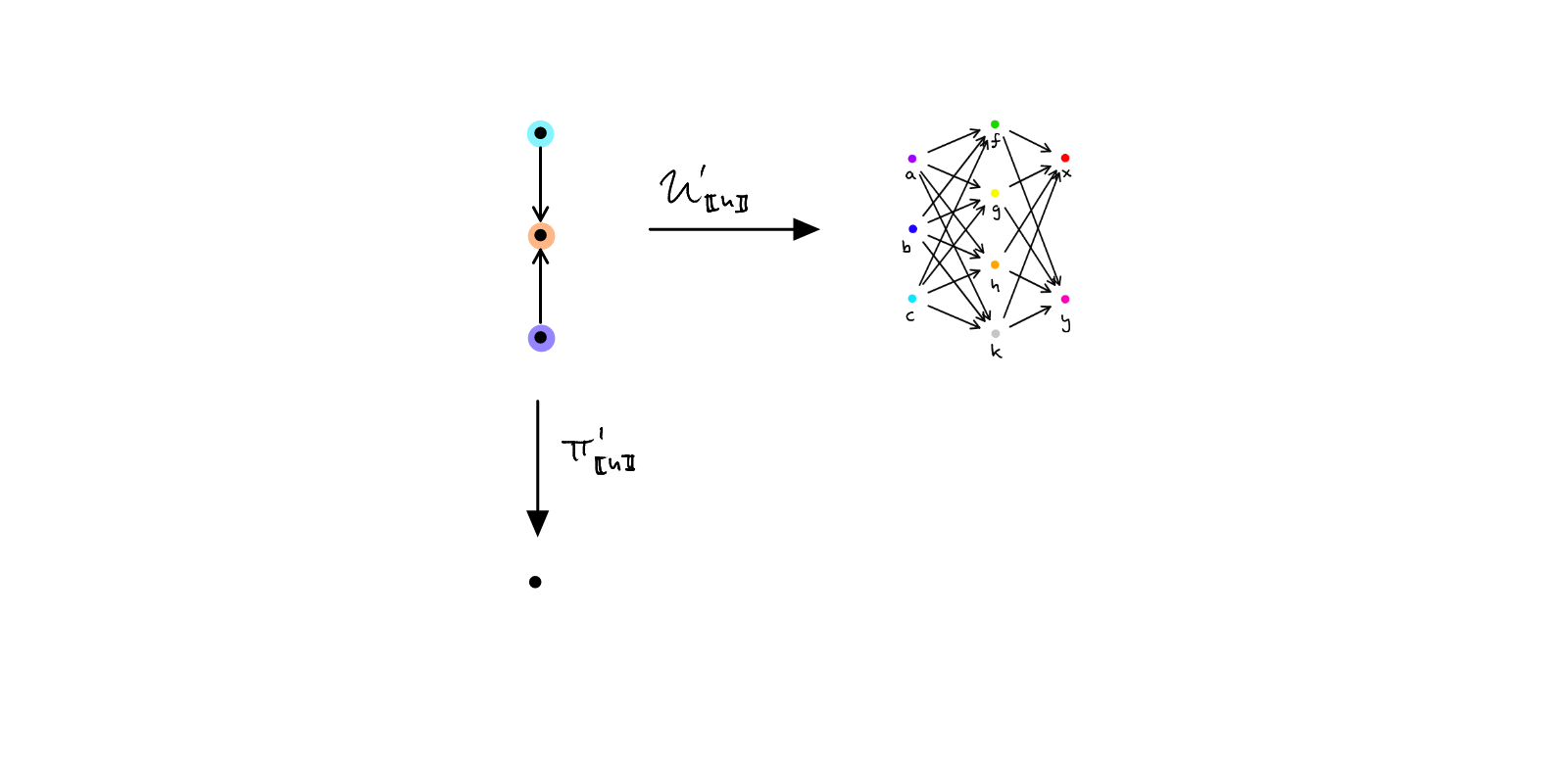}
\endgroup\end{restoretext}
$\abss{  x}$ is the $\SIvert 2 {\kC(S)}$-family
\begin{restoretext}
\begingroup\sbox0{\includegraphics{test/page1.png}}\includegraphics[clip,trim=0 {.0\ht0} 0 {.0\ht0} ,width=\textwidth]{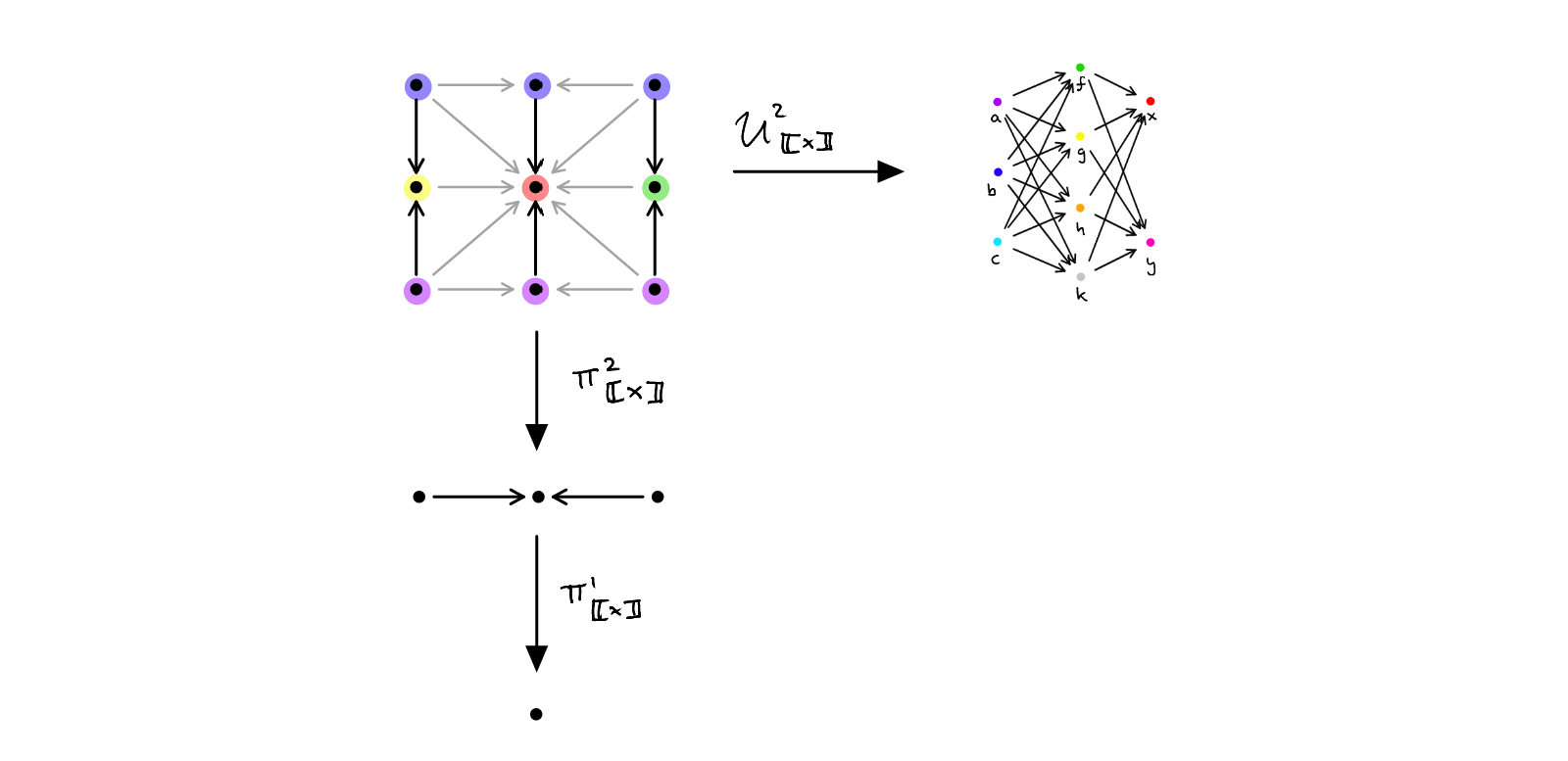}
\endgroup\end{restoretext}
Similarly, $\abss{y}$ is the $\SIvert 2 {\kC(S)}$-family
\begin{restoretext}
\begingroup\sbox0{\includegraphics{test/page1.png}}\includegraphics[clip,trim=0 {.0\ht0} 0 {.0\ht0} ,width=\textwidth]{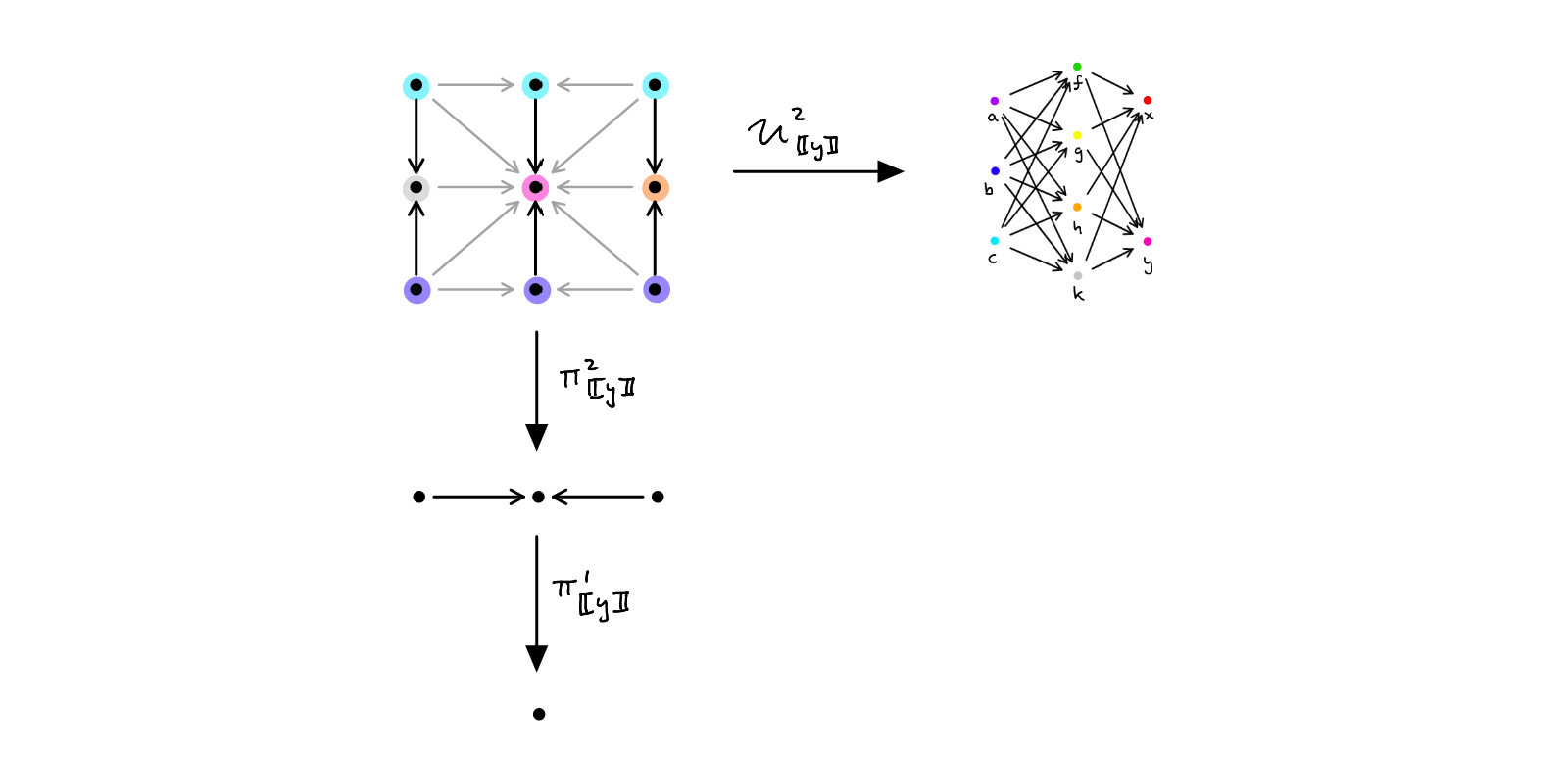}
\endgroup\end{restoretext}
\end{eg}

\section{Examples} \label{sec:panc_ex}

We give detailed examples of the definition of presented associative $n$-categories in the previous section. We first introduce the following

\begin{notn}[$n$-quotients of globular sets] \label{constr:globular_quotients} Let $n \in \lN \cup \Set{\infty}$ and $S$ be a globular set. We define the \textit{globular $n$-quotient} $S^{\quotg}_{\truncleq n}$ to be the $n$-truncated globular set with elements
\begin{equation}
(S^{\quotg}_{\truncleq n})_i := S_i
\end{equation}
for $i < n$, and (cf. \autoref{notn:glob_set_src_tgt})
\begin{equation} \label{eq:quotient_glob}
(S^{\quotg}_{\truncleq n})_n := \coeq(s_n,t_n)
\end{equation}
In other words, $(S^{\quotg}_{\truncleq n})_n$ is the quotient of $S_n$ by the equivalence relation generated by $g_1 \sim g_2 \text{~iff~} \exists h \in S_{n+1}. s_n h = g_1 \wedge t_n h = g_2$.
\end{notn}

Note that since both $s_{n-1}$ and $t_{n-1}$ coequalise $(s_n,t_n)$ (of $S$) by the globularity condition for $S$, we can now chose $s_{n-1}$, $t_{n-1}$ for $S^{\quotg}_{\truncleq n}$ to be the maps induced by $s_{n-1}$, $t_{n-1}$ for $S$ on the quotient \eqref{eq:quotient_glob}. For  $i < n-1$ we further chose $s_i$ for $S^{\quotg}_{\truncleq n}$ to equal $s_i$ from $S_i$.

\subsection{Dimension -1 and 0} 

\Free{} associative $(-1)$-categories consist of a single set $\scC_0$ which by proof irrelevance can be either empty or contain a single element. Thus, \free{} associative $(-1)$-categories correspond to \textit{truth values}.

\Free{} associative $0$-categories consist of a set $\scC_0$ (of objects) and a set $\scC_1$ (of generating equalities). Consider for instance the sets
\begin{restoretext}
\begingroup\sbox0{\includegraphics{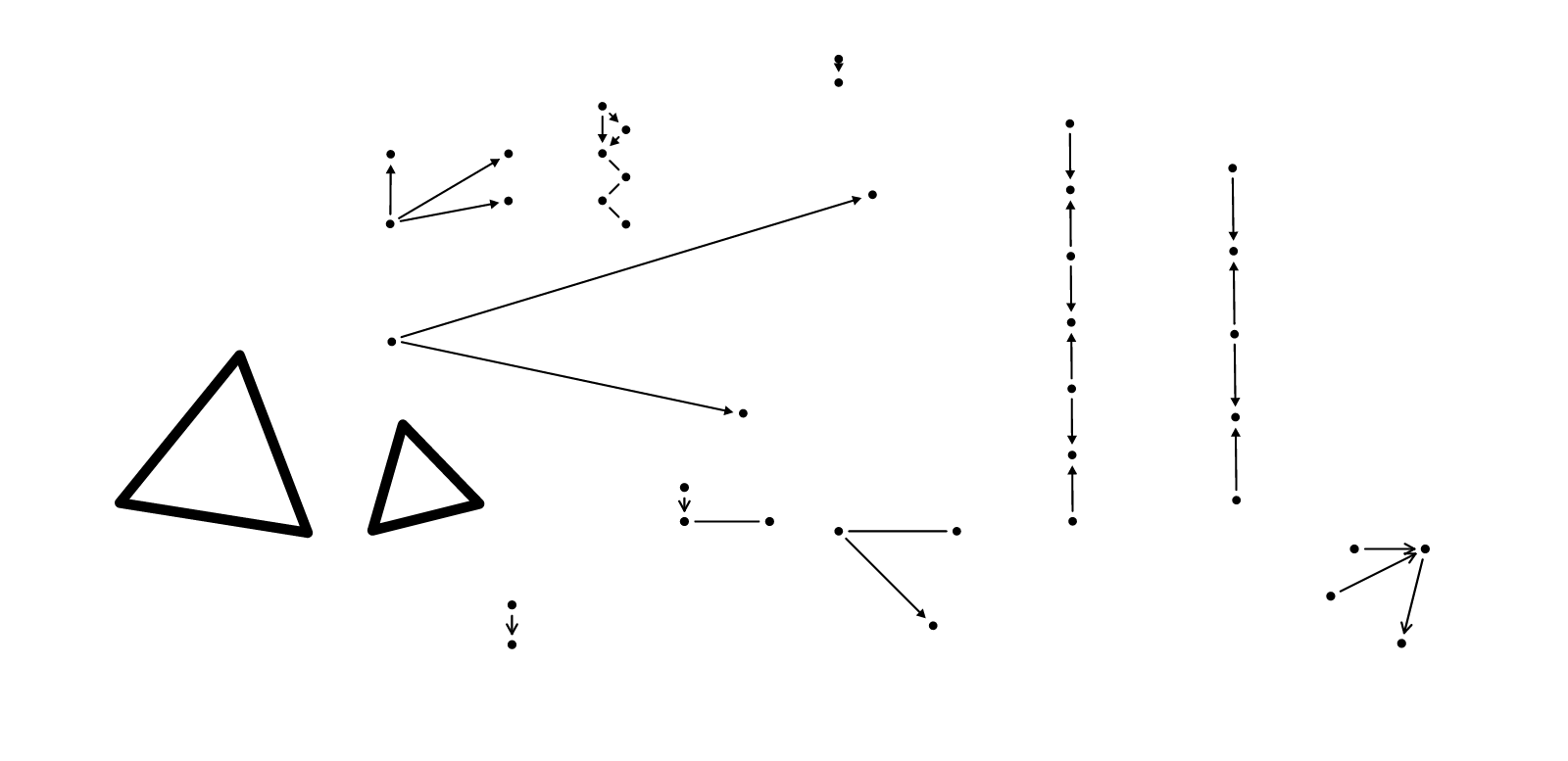}}\includegraphics[clip,trim=0 {.4\ht0} 0 {.1\ht0} ,width=\textwidth]{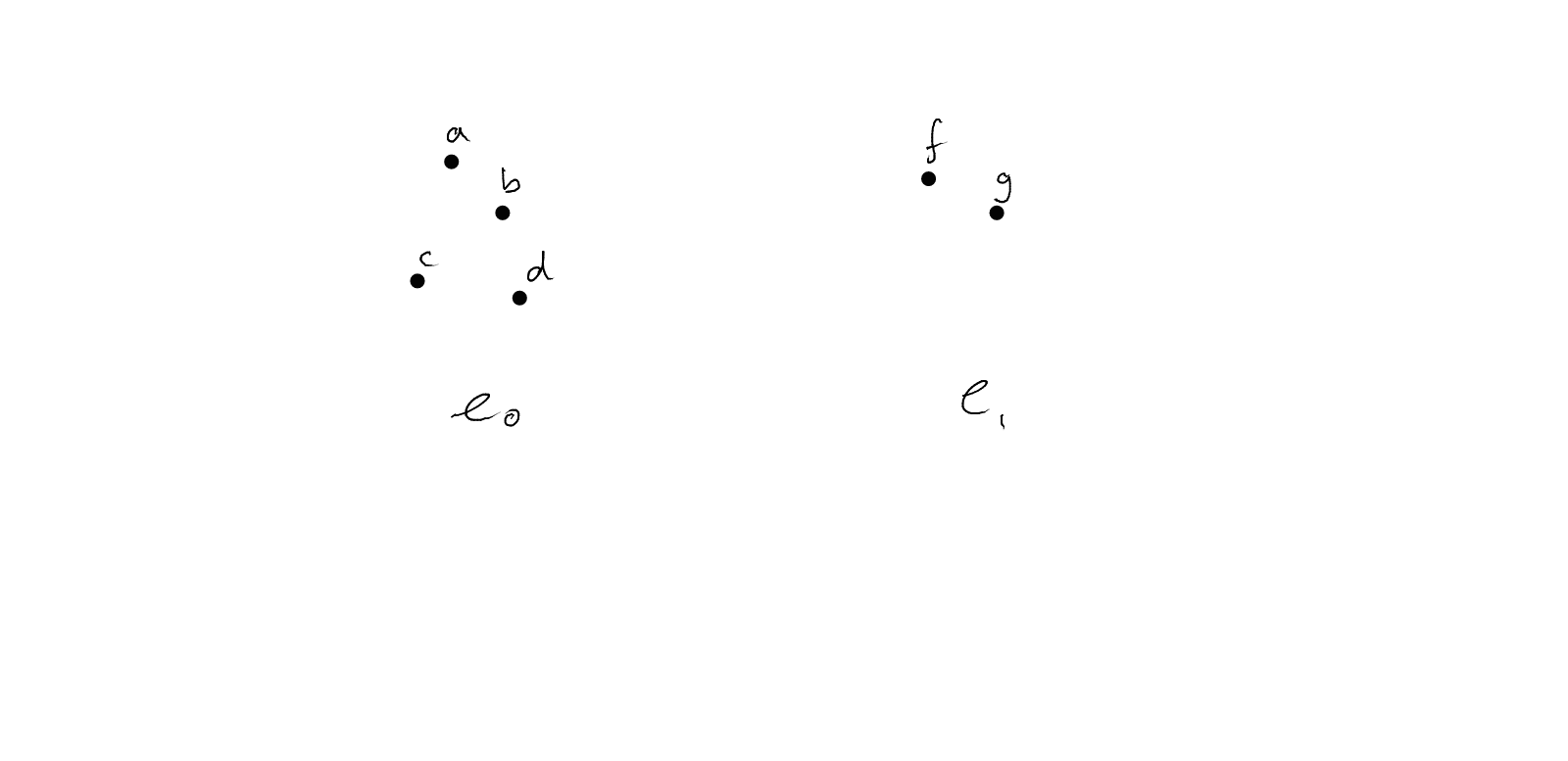}
\endgroup\end{restoretext}
Then $\GGamma{}\sC$ is poset of the form
\begin{restoretext}
\begingroup\sbox0{\includegraphics{ANCimg/page1.png}}\includegraphics[clip,trim=0 {.2\ht0} 0 {.2\ht0} ,width=\textwidth]{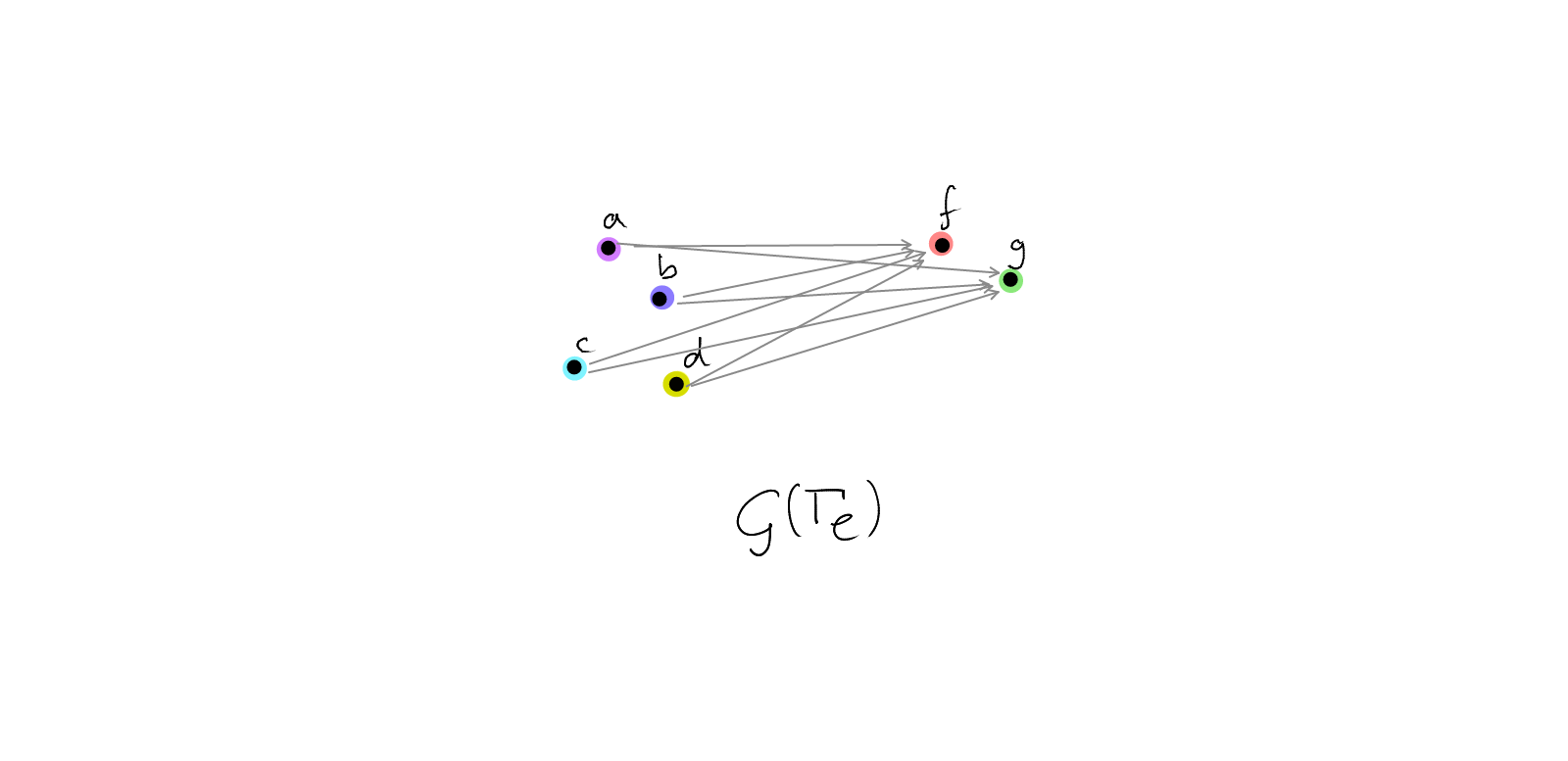}
\endgroup\end{restoretext}
To define an assignment of types in dimension $0$, each $x\in \sC_0$ needs an associated $0$-dimensional cube $\abss{x}$. Concretely, identifying labels and colors correctly, we define
\begin{restoretext}
\begingroup\sbox0{\includegraphics{ANCimg/page1.png}}\includegraphics[clip,trim=0 {.0\ht0} 0 {.0\ht0} ,width=\textwidth]{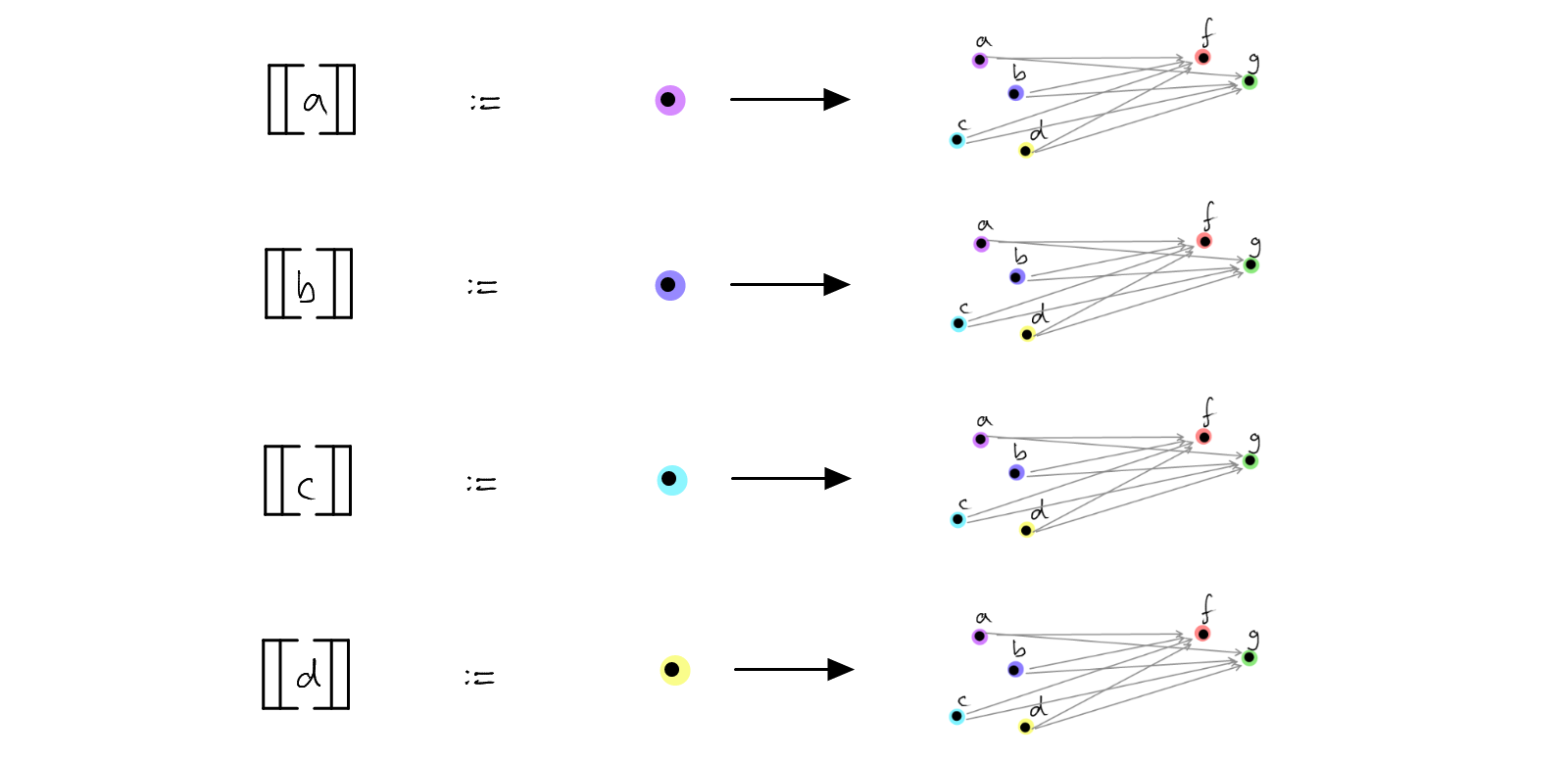}
\endgroup\end{restoretext}
These are the unique choices of $0$-types satisfying minimality as required in \autoref{defn:pres_ANC}.

To complete the definition of the assignment of types we further need to give, for every $x \in \scC_1$, an associated $1$-dimensional cube $\abss{x}$. Due to minimality, this must in fact be of the form of a terminal 1-cube. For our example we set
\begin{restoretext}
\begingroup\sbox0{\includegraphics{ANCimg/page1.png}}\includegraphics[clip,trim=0 {.0\ht0} 0 {.0\ht0} ,width=\textwidth]{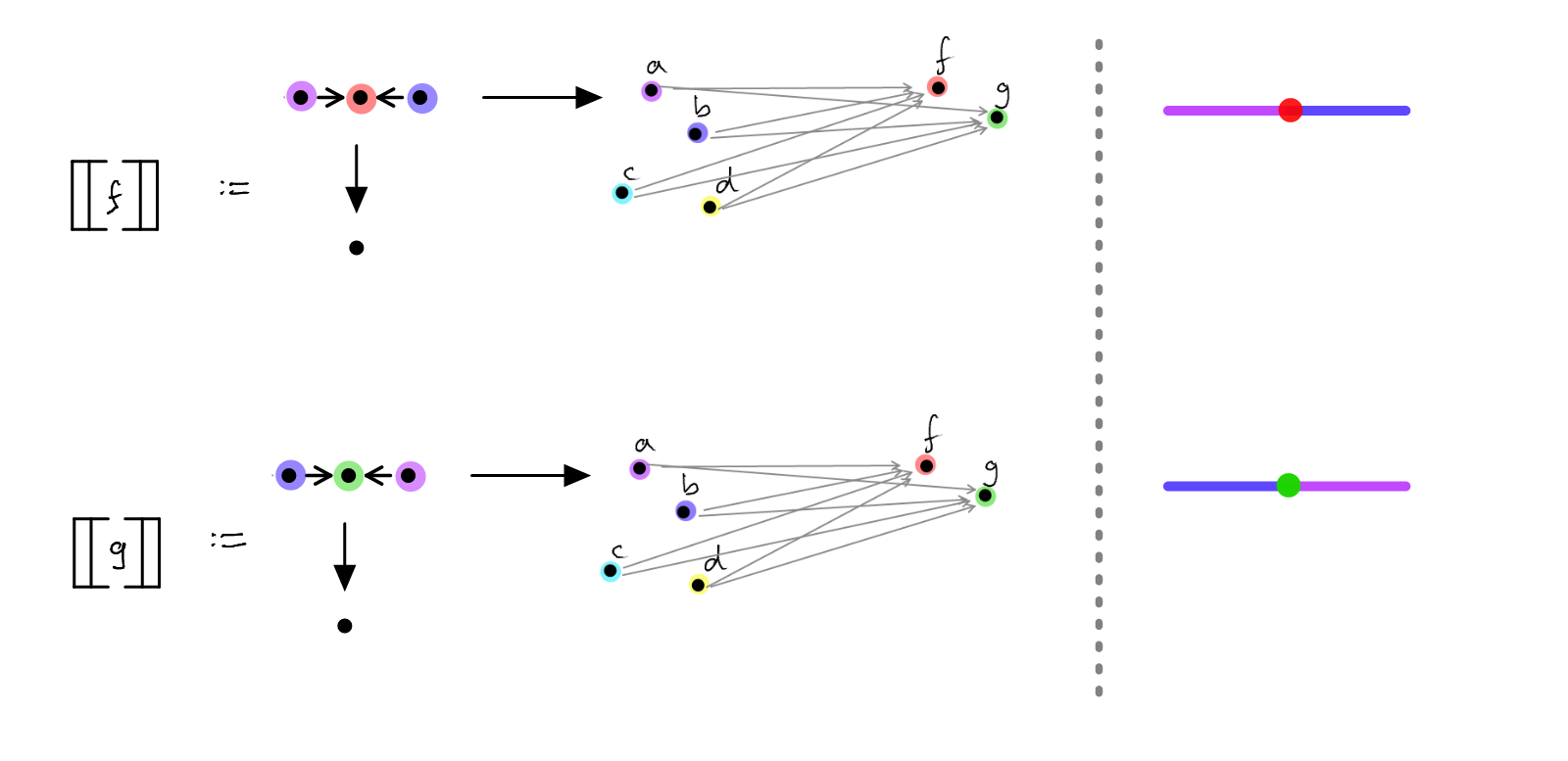}
\endgroup\end{restoretext}
Here, on the right we depict manifold diagrams corresponding to algebraic data on the left (cf. \autoref{ssec:coloring}).

Note that for instance the following would have been an invalid choice
\begin{restoretext}
\begingroup\sbox0{\includegraphics{ANCimg/page1.png}}\includegraphics[clip,trim=0 {.35\ht0} 0 {.3\ht0} ,width=\textwidth]{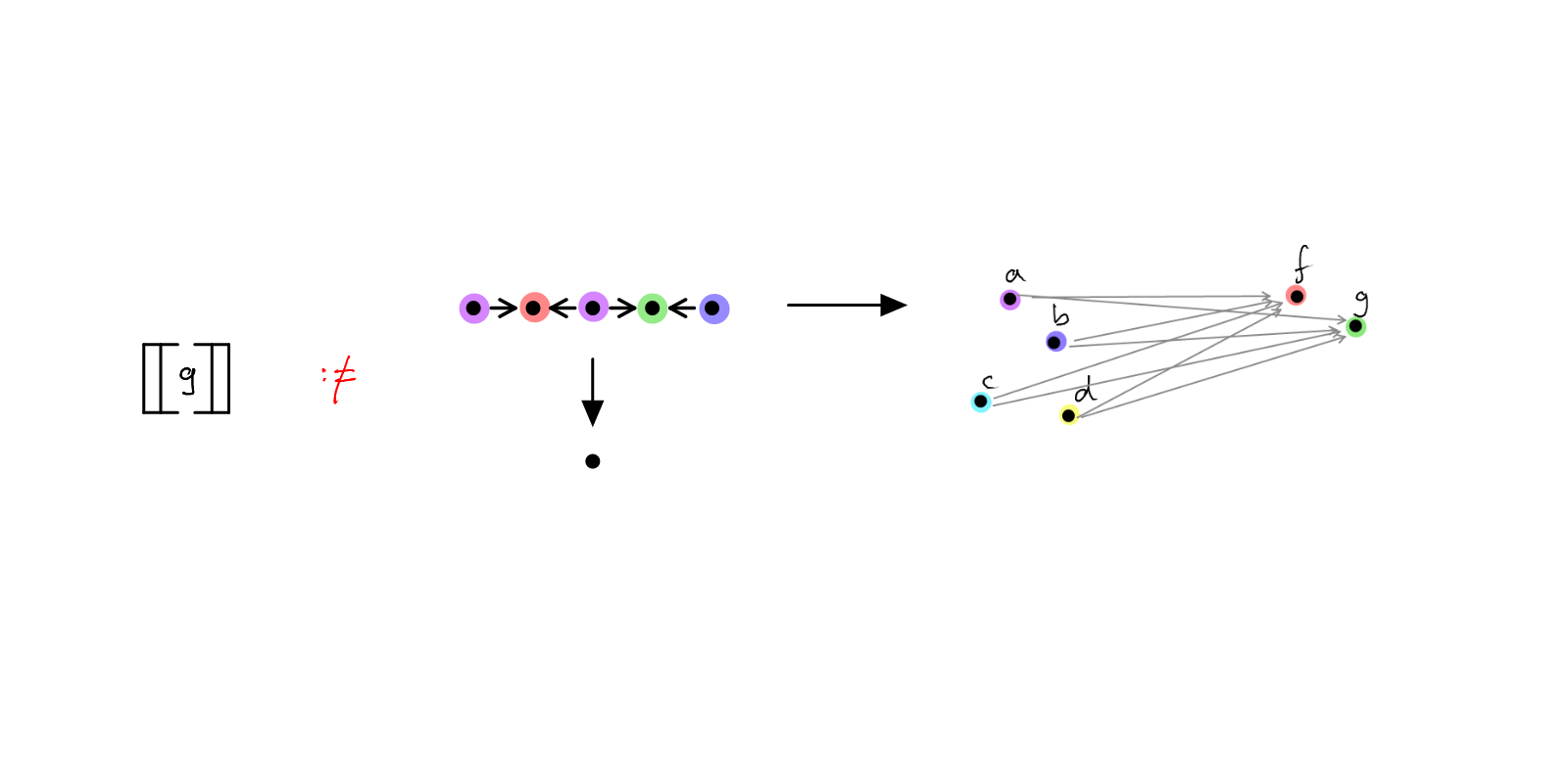}
\endgroup\end{restoretext}
as this doesn't satisfy the minimality condition. 

Correctness of dimensions can be seen to hold in the above definitions. The typability condition is trivial in this dimension, for instance, for $p = (0,0) \in \tsG 1(\abss{g})$ we find
\begin{restoretext}
\begingroup\sbox0{\includegraphics{ANCimg/page1.png}}\includegraphics[clip,trim=0 {.15\ht0} 0 {.25\ht0} ,width=\textwidth]{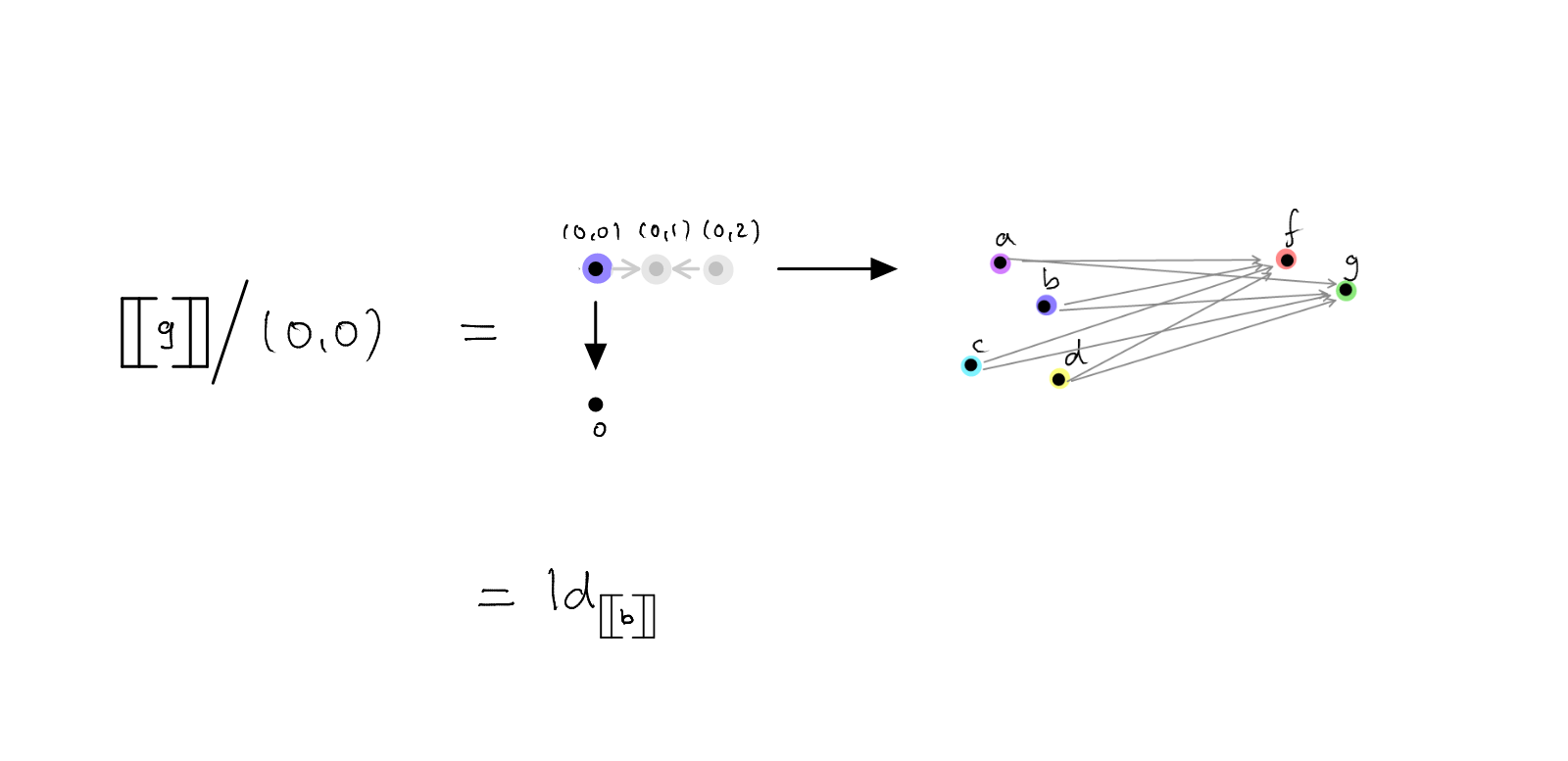}
\endgroup\end{restoretext}
as required for typability. Further, note that our choice of generating relations $f,g \in \sC_1$ satisfies the symmetry and and proof irrelevance conditions. 

Elements in $\Comp(\sC)_1$ (cf. \autoref{defn:PANC_mor}) are of the form
\begin{restoretext}
\begingroup\sbox0{\includegraphics{ANCimg/page1.png}}\includegraphics[clip,trim=0 {.1\ht0} 0 {.15\ht0} ,width=\textwidth]{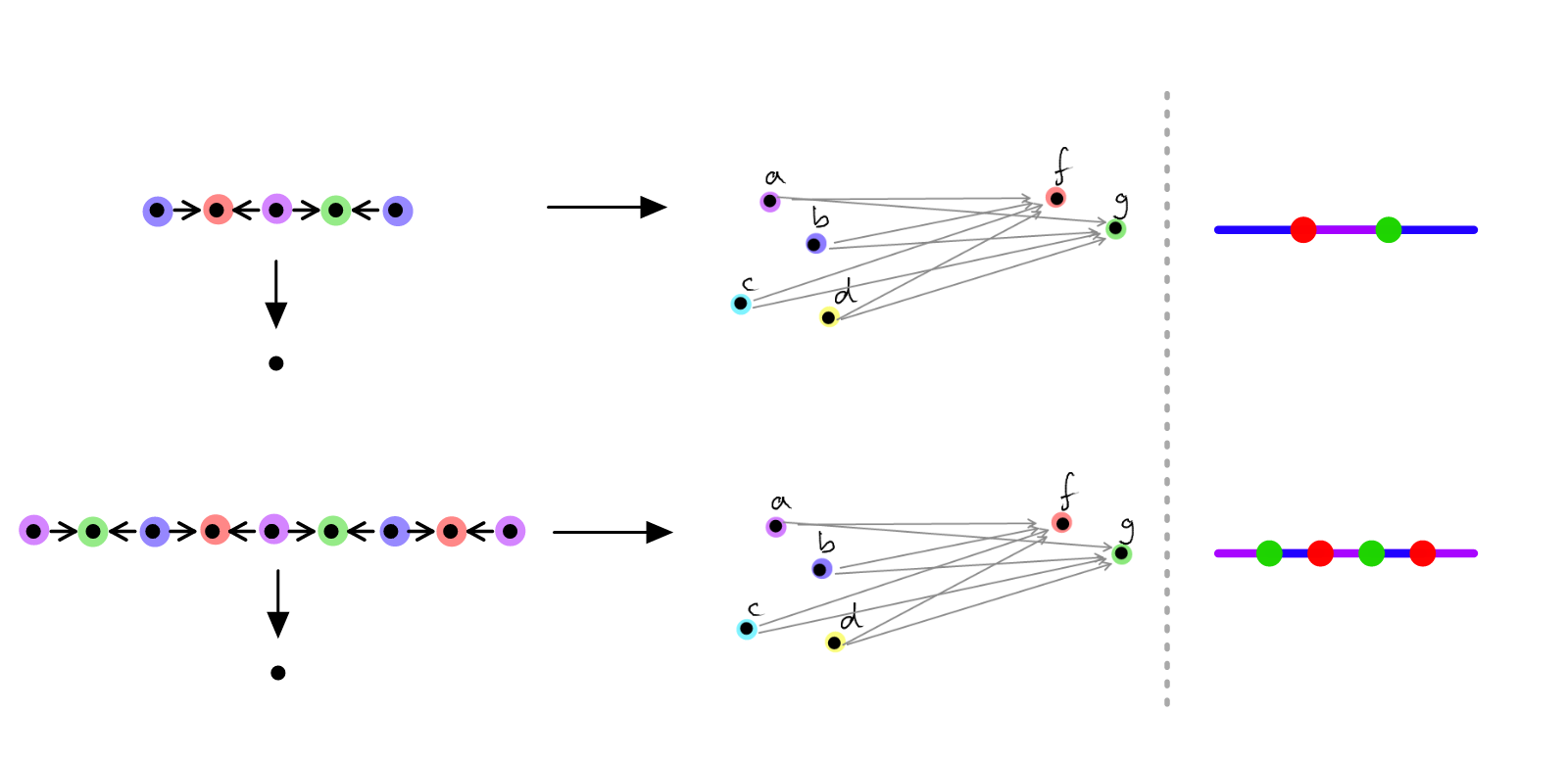}
\endgroup\end{restoretext}
$\Comp(\sC)_1$ induces a relation on the set $\Comp(\sC)_0 = \sC_0$, by relating its source and target of each elements $x \in \Comp(\sC)_1$ by
\begin{equation}
s_0 x \leq t_0 x
\end{equation}
The relation given by all elements in $\Comp(\sC)_1$ is transitive since $1$-cube can be glued along their endpoints (see \autoref{ch:composition}), it is symmetric by the symmetry assumption on the generating relations and it is reflexive since $\Comp(\sC)_1$ also contains identity cubes. 

Thus, a \free{} associative $0$-category is a set $\sC_0 = \Comp(\sC)_0$ with an equality relation induced by $\Comp(\sC)_1$. Such a set with equality relation is also called a \textit{setoid} in the literature. Note that the set resulting from quotienting $\sC_0$ by the equality relation induced by $\Comp(\sC)_1$ precisely gives $(\Comp(\sC)^{\quotg}_{\truncleq 0})_0$ (cf. \autoref{constr:globular_quotients}). For our specific example above this is the set $\Set{[a],c,d}$ where $[a]$ is the equivalence class containing $\Set{a,b}$.

\subsection{Dimension 1} 

\Free{} associative $1$-categories consist of a set $\scC_0$ (of objects), a set $\scC_1$ (of generating $1$-morphisms), and a set $\scC_2$ (of generating equalities between $1$-morphisms) together with an assignment of types. For instance consider $\sC$ with objects
\begin{restoretext}
\begingroup\sbox0{\includegraphics{ANCimg/page1.png}}\includegraphics[clip,trim=0 {.25\ht0} 0 {.25\ht0} ,width=\textwidth]{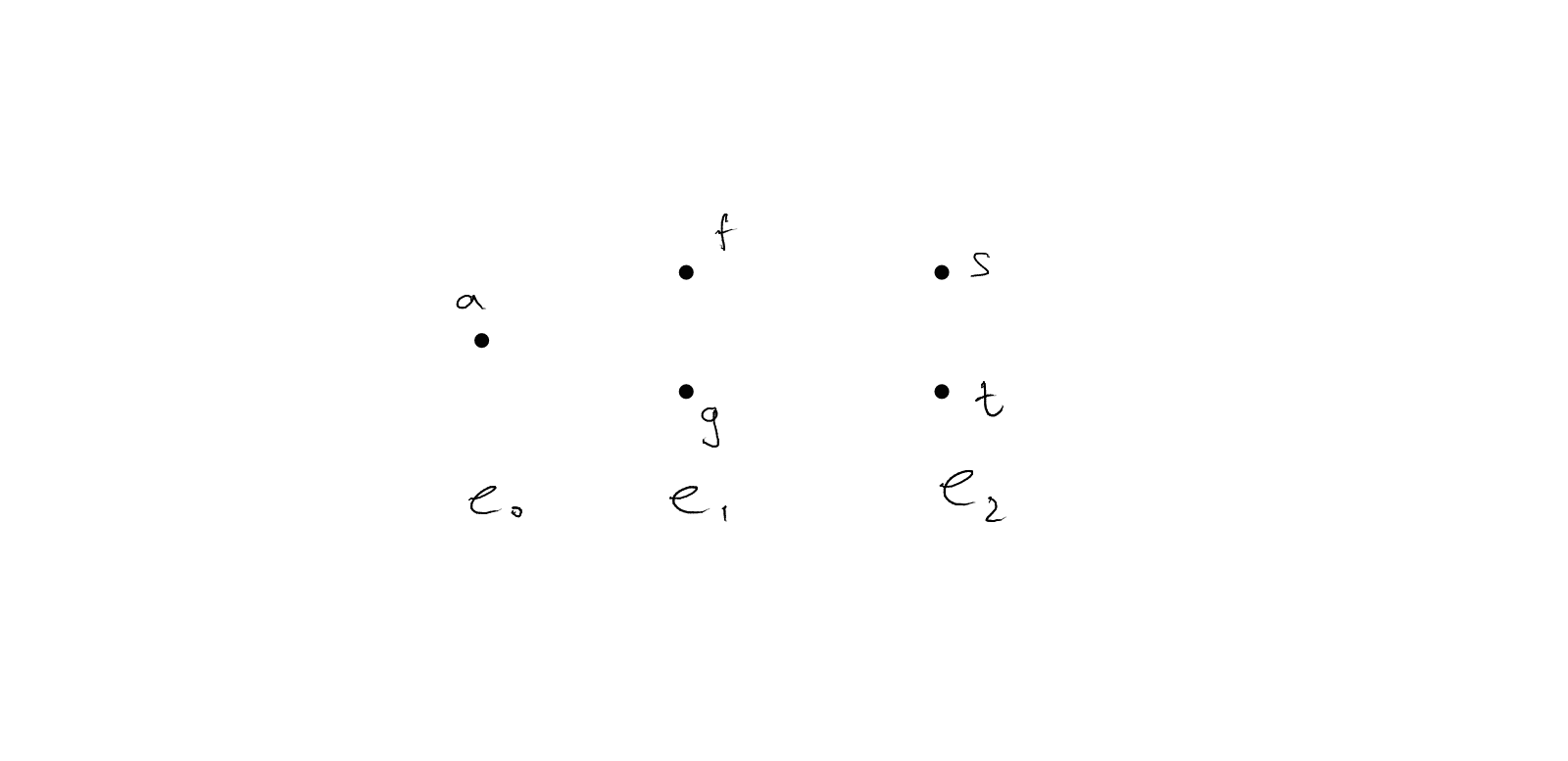}
\endgroup\end{restoretext}
In this case $\GGamma{}\sC$ takes the form
\begin{restoretext}
\begingroup\sbox0{\includegraphics{ANCimg/page1.png}}\includegraphics[clip,trim=0 {.2\ht0} 0 {.25\ht0} ,width=\textwidth]{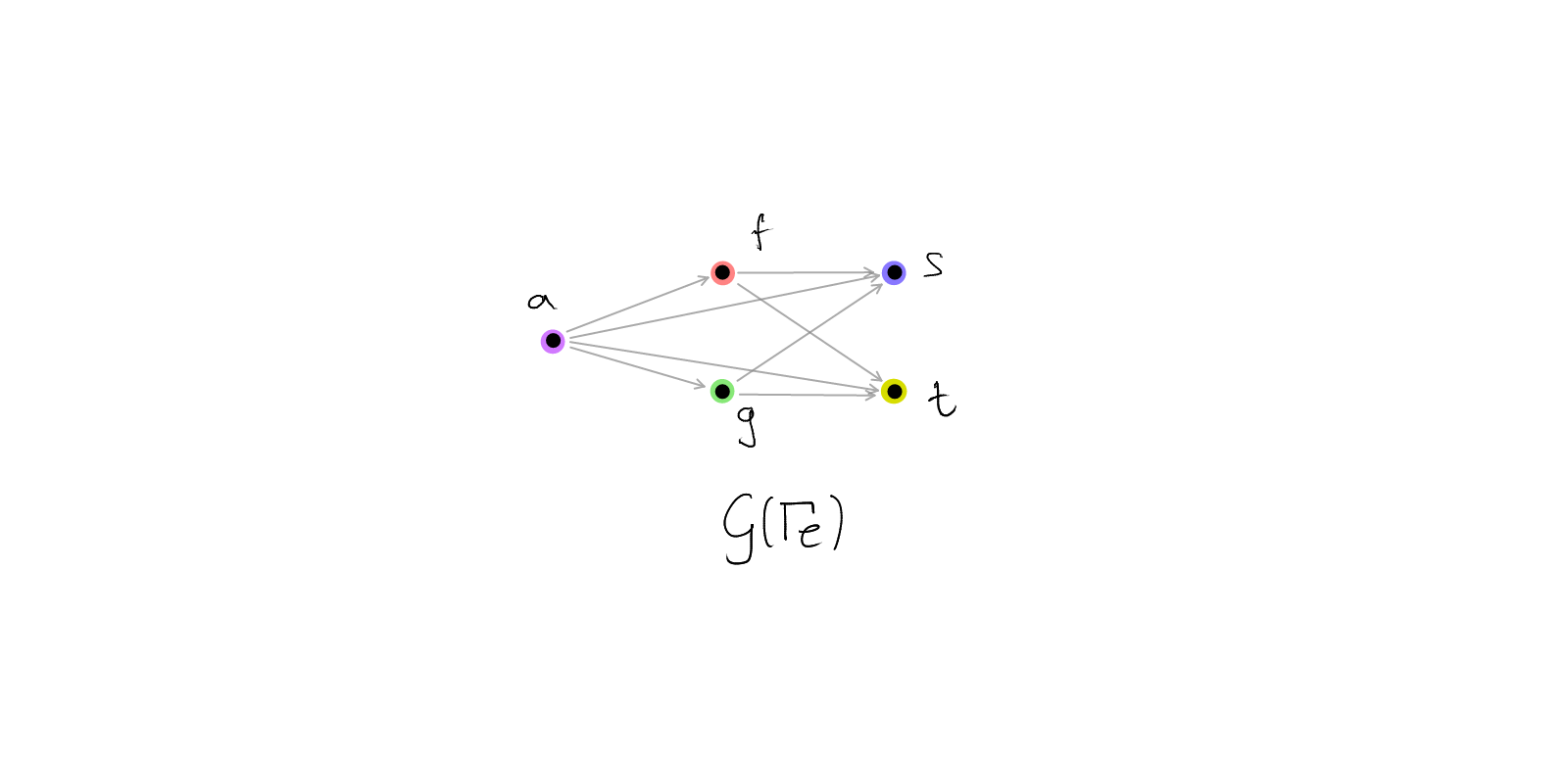}
\endgroup\end{restoretext}
For the assignment of types, for each $x \in \sC_m$ we need to define $\abss{x}$. For $m = 0$, $\abss{a}$ is defined as before (the functor from $\bnum{1} \to \GGamma{}\sC$ mapping $0$ to $a$). For $m>0$, we will depict both the $\SIvert m {\GGamma{}\sC}$-cube (on the left) as well as its associated manifold diagram (on the right). We make the following definition.
\begin{restoretext}
\begingroup\sbox0{\includegraphics{ANCimg/page1.png}}\includegraphics[clip,trim=0 {.15\ht0} 0 {.0\ht0} ,width=\textwidth]{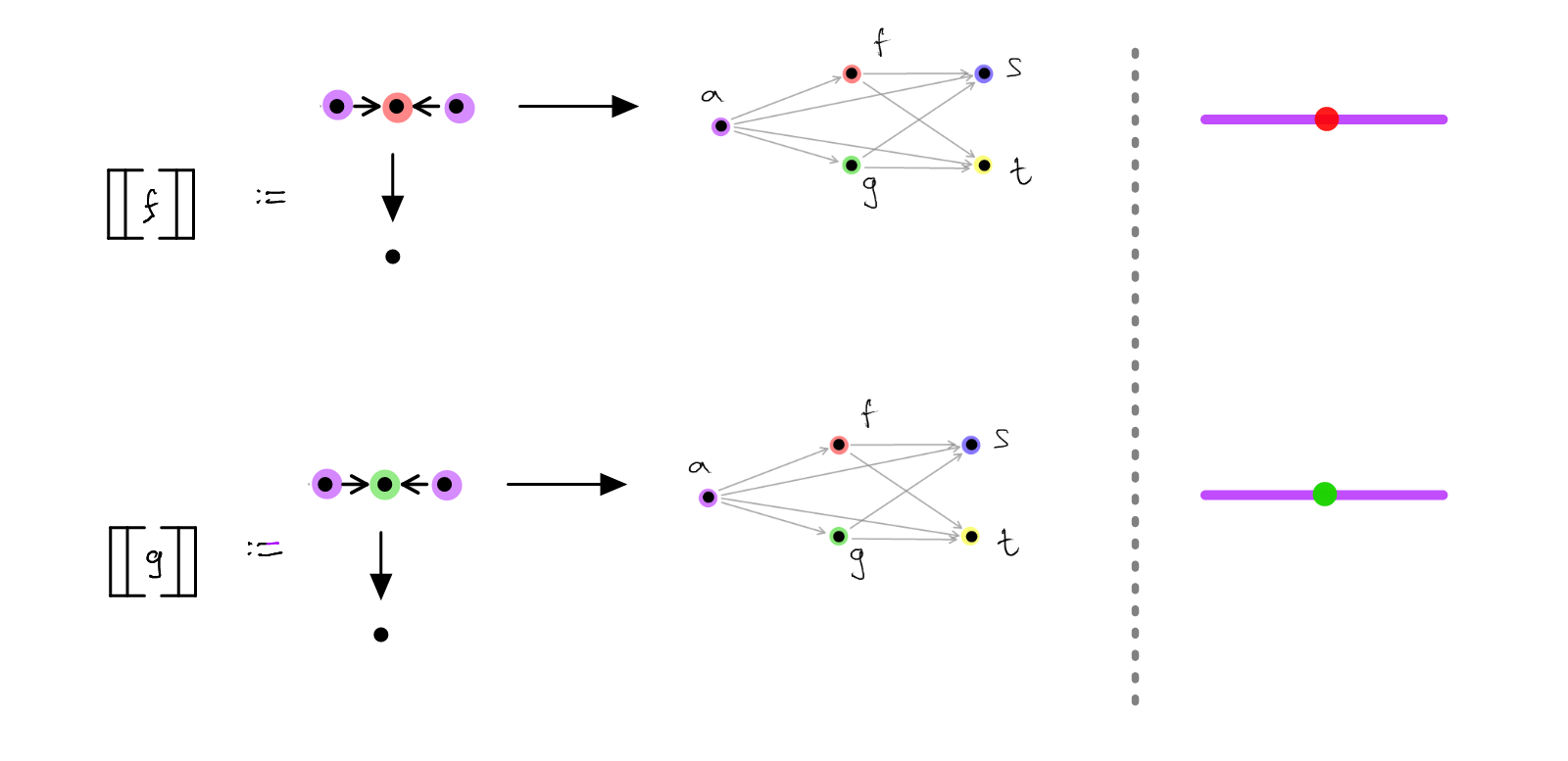}
\endgroup\end{restoretext}
\begin{restoretext}
\begingroup\sbox0{\includegraphics{ANCimg/page1.png}}\includegraphics[clip,trim=0 {.15\ht0} 0 {.15\ht0} ,width=\textwidth]{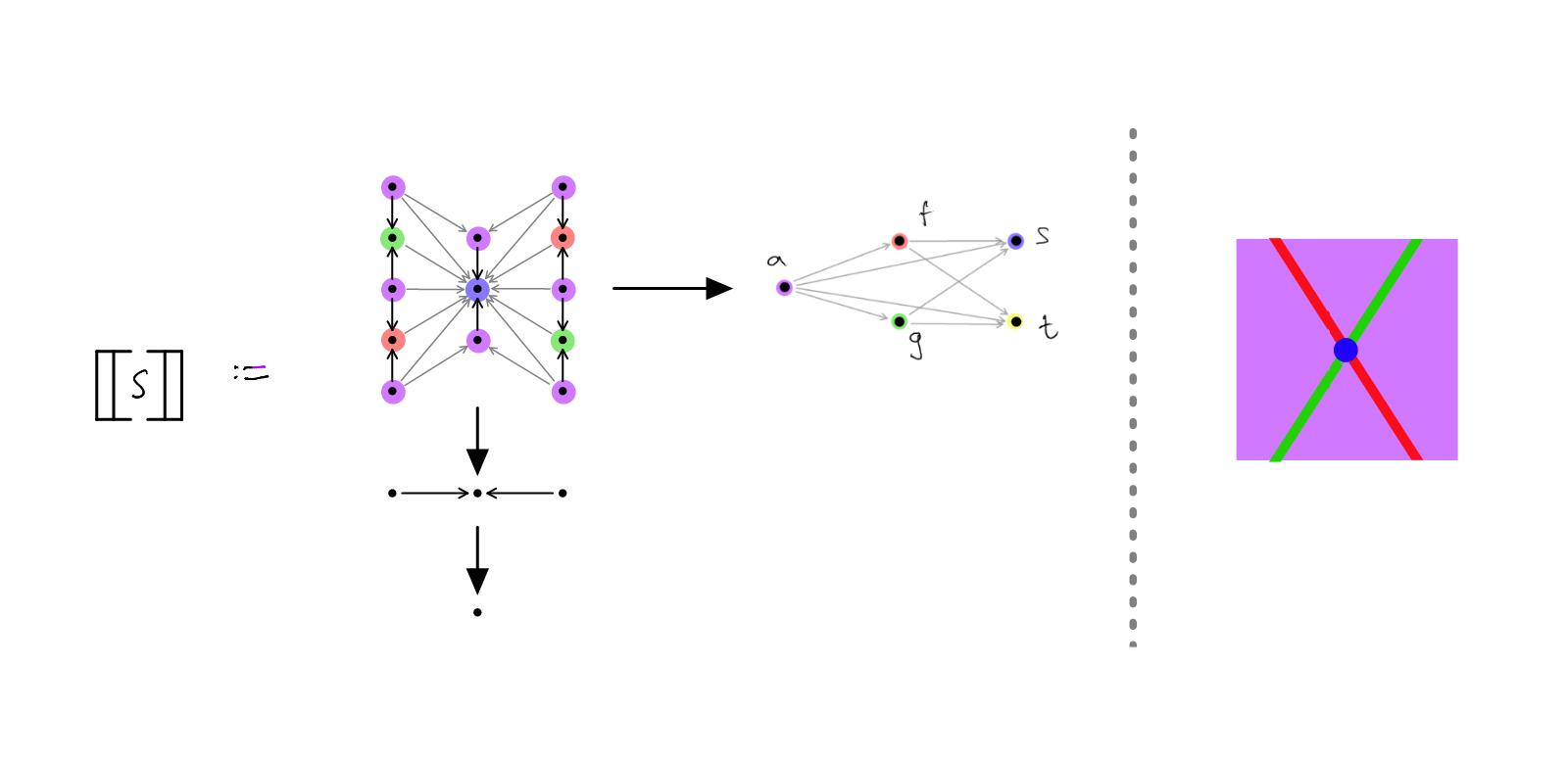}
\endgroup\end{restoretext}
and
\begin{restoretext}
\begingroup\sbox0{\includegraphics{ANCimg/page1.png}}\includegraphics[clip,trim=0 {.15\ht0} 0 {.15\ht0} ,width=\textwidth]{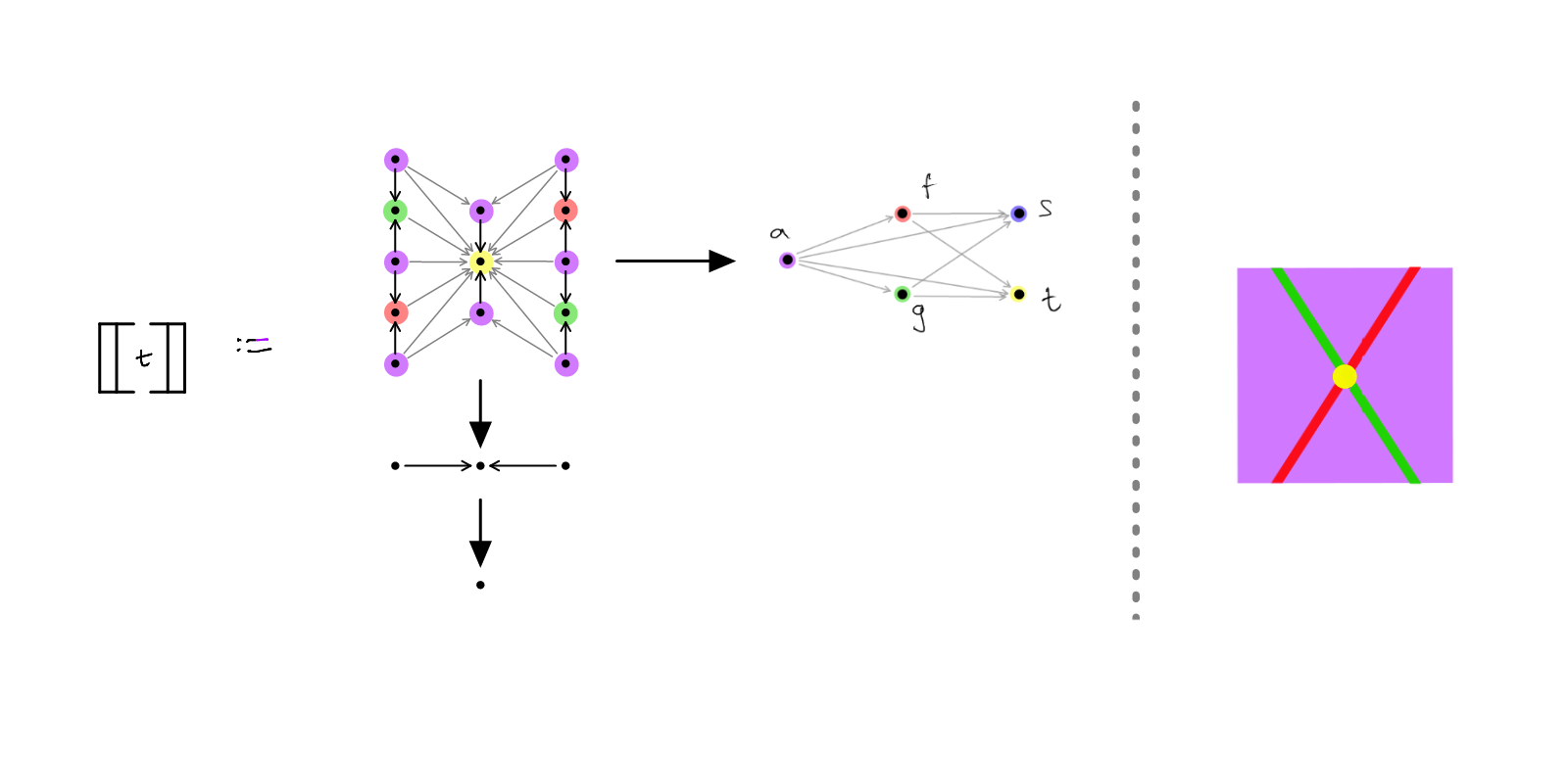}
\endgroup\end{restoretext}
Recall these need to satisfy the minimality, dimensionality and typability condition. The first two are easily checked, for the last we give the following example in the case of $\abss{t}$ and $((0,2),3) \in \tsG 2(\abss{t})$. 
\begin{restoretext}
\begingroup\sbox0{\includegraphics{ANCimg/page1.png}}\includegraphics[clip,trim=0 {.15\ht0} 0 {.0\ht0} ,width=\textwidth]{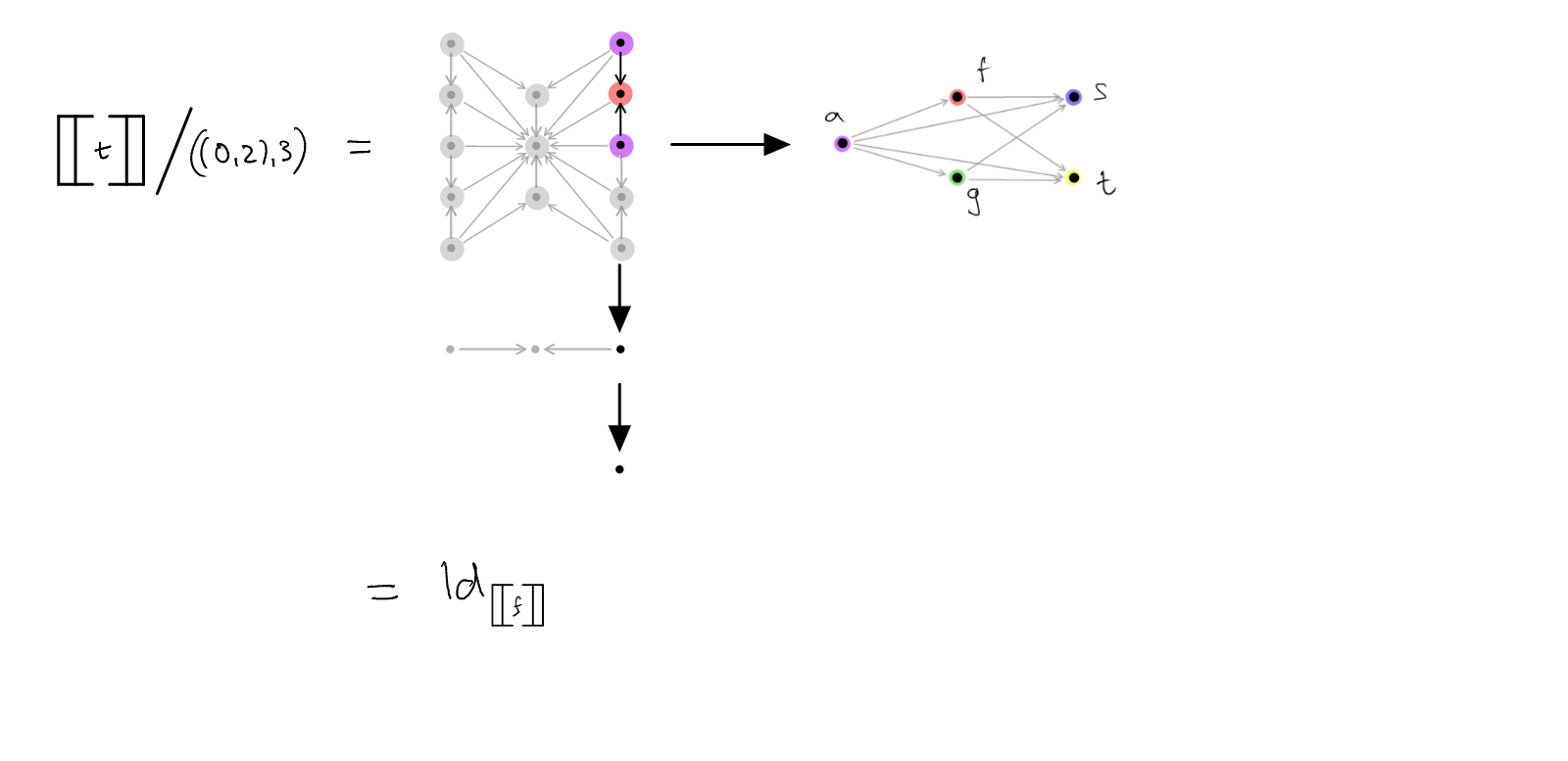}
\endgroup\end{restoretext}
In terms of manifold diagrams, this type check can be understood as looking at the neighbourhood (marked in \cblue{}) of the \cred{} wire, and checking that it is consisted with the given type of the \cred{} wire:
\begin{restoretext}
\begingroup\sbox0{\includegraphics{ANCimg/page1.png}}\includegraphics[clip,trim=0 {.25\ht0} 0 {.35\ht0} ,width=\textwidth]{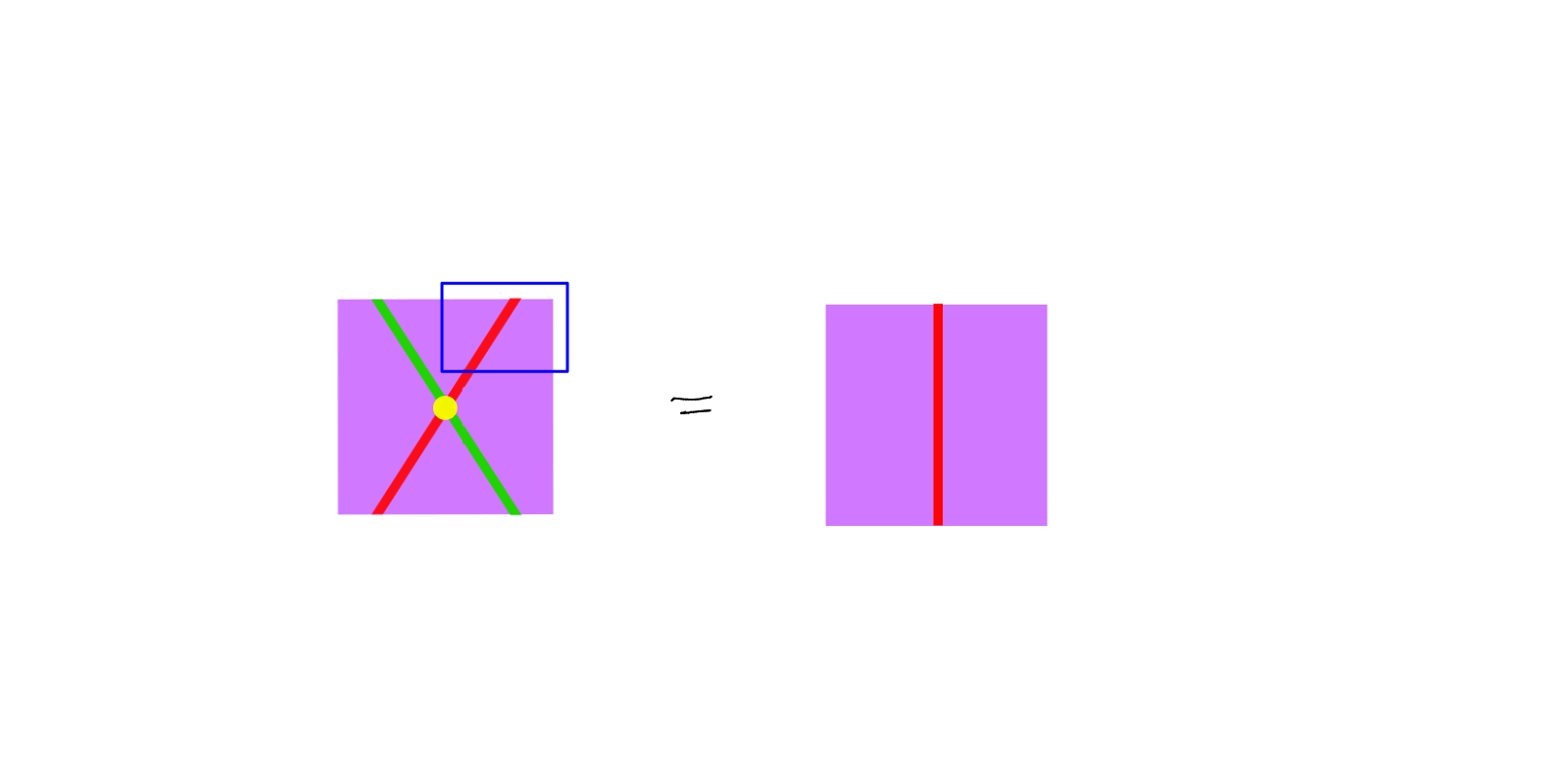}
\endgroup\end{restoretext}
Similar to the case of $0$-categories, we interpret $2$-dimensional (normalised, globular, \typed) cubes in $\sC$ (that is, elements of $\Comp(\sC)_2$) as equality relations on the set of  $1$-dimensional (normalised, globular, \typed) cubes in $\sC$. Namely, each $\scA 
\in \Comp(\sC)_2$ gives a relation
\begin{equation}
\gsrc(\scA) \leq \gtgt(\scA)
\end{equation}
By the symmetry condition we know this is an equivalence relation (it is transitive since cubes can be glued together, see \autoref{ch:composition}). Quotienting the set of $\Comp(\sC)_1$ by this relation gives $(\Comp(\sC)^{\quotg}_{\truncleq 1})_1$. As an example, the cube $\abss{t} \in \Comp(\sC)_2$ witnesses the equality
\begin{equation}
\gsrc \abss{t} = \gtgt \abss{t}
\end{equation}
which using manifold diagrams can be expressed
\begin{restoretext}
\begingroup\sbox0{\includegraphics{ANCimg/page1.png}}\includegraphics[clip,trim=0 {.45\ht0} 0 {.35\ht0} ,width=\textwidth]{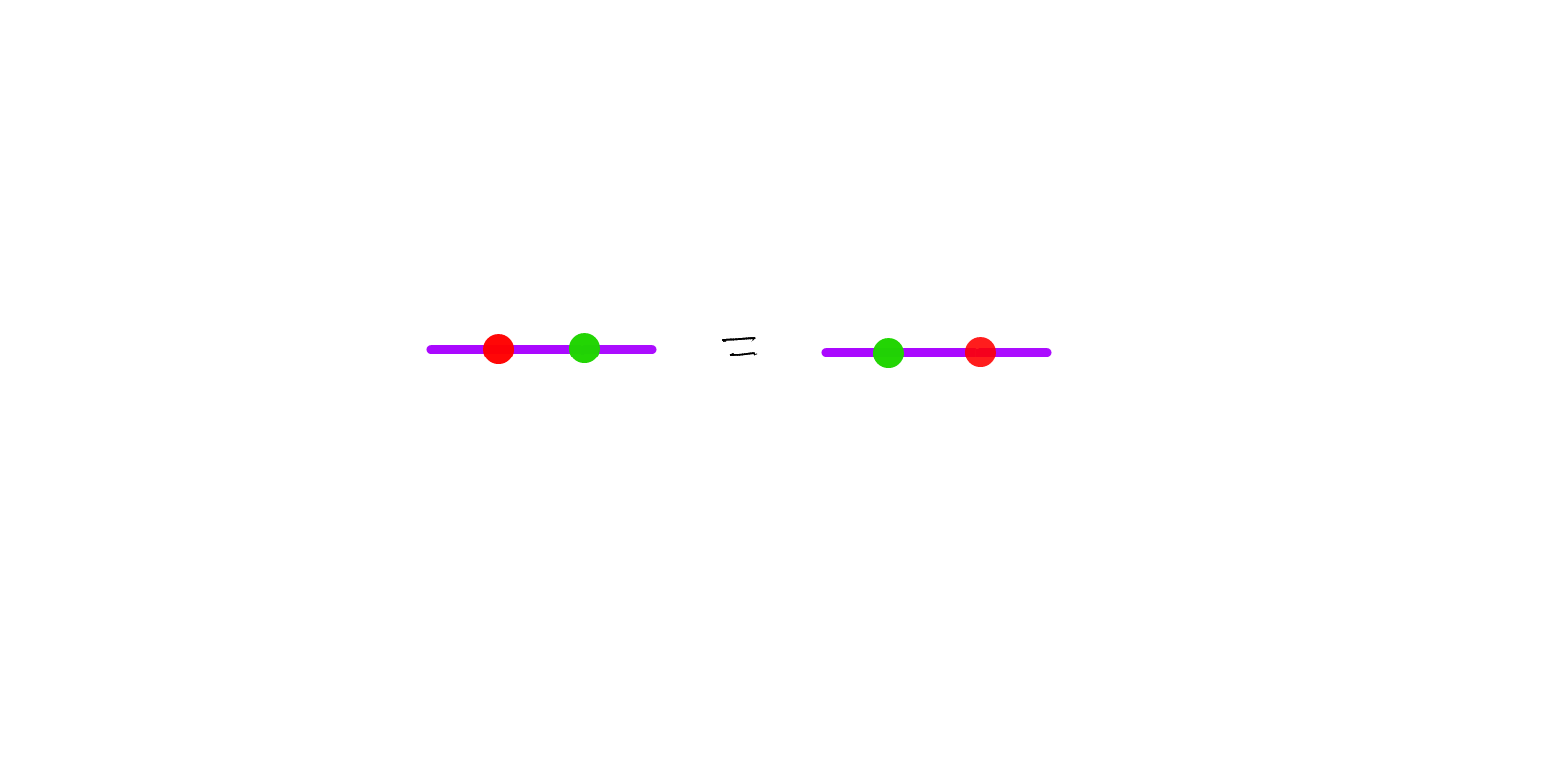}
\endgroup\end{restoretext}
Similarly, we have the equalities
\begin{restoretext}
\begingroup\sbox0{\includegraphics{ANCimg/page1.png}}\includegraphics[clip,trim=0 {.1\ht0} 0 {.25\ht0} ,width=\textwidth]{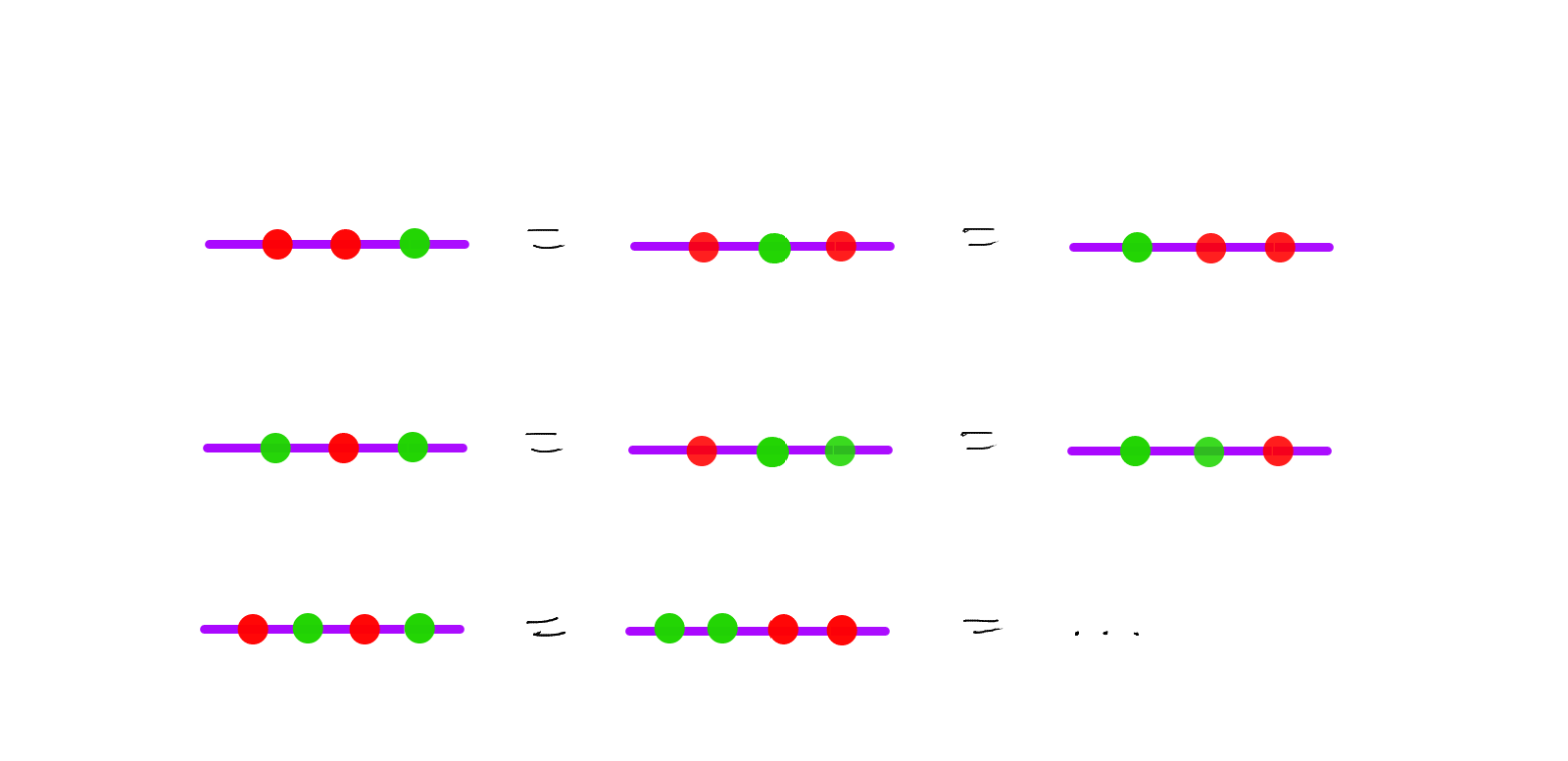}
\endgroup\end{restoretext}
some of which are witnessed for instance by the following elements of $\Comp(\sC)_2$
\begin{restoretext}
\begingroup\sbox0{\includegraphics{ANCimg/page1.png}}\includegraphics[clip,trim=0 {.3\ht0} 0 {.3\ht0} ,width=\textwidth]{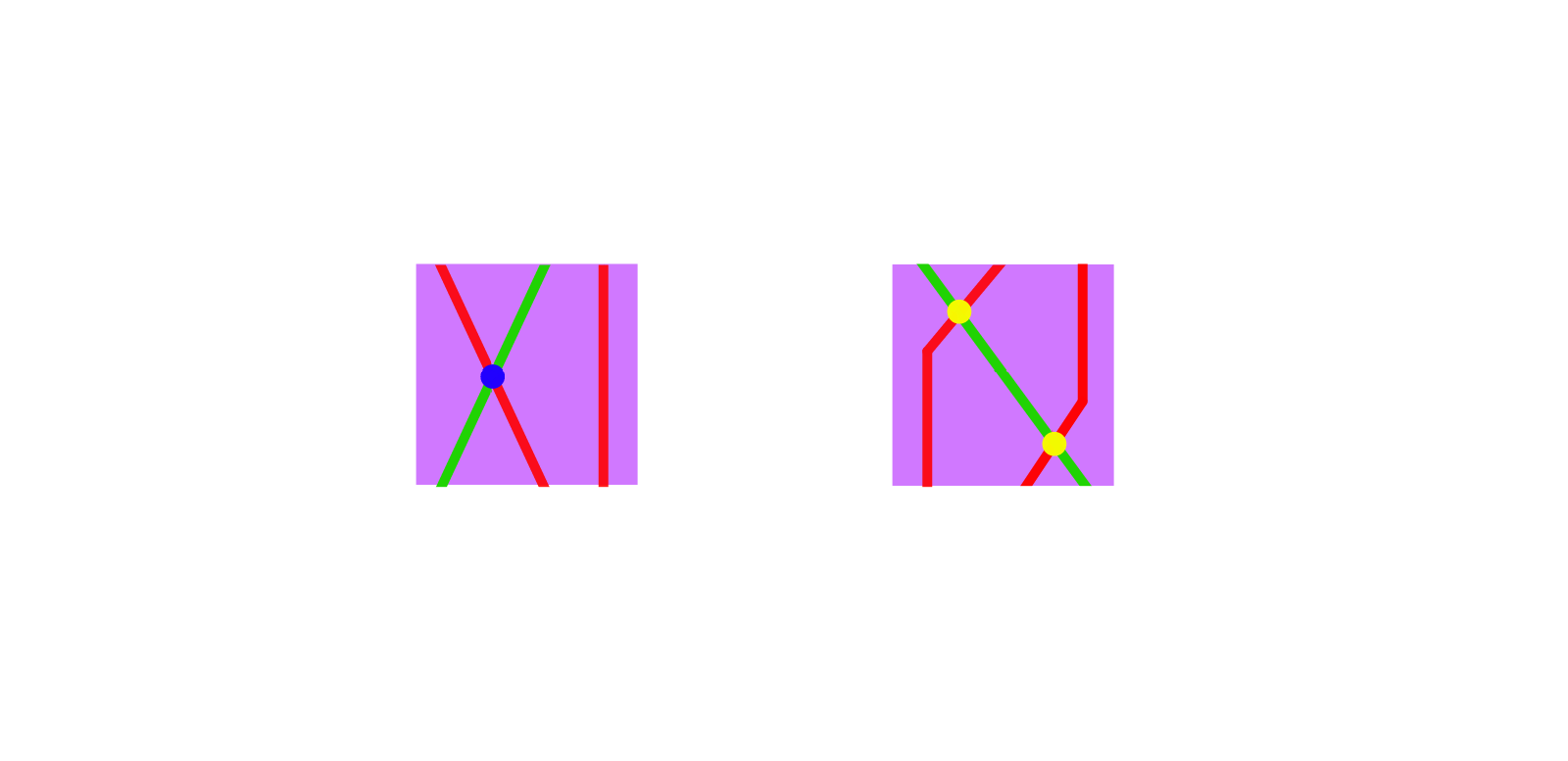}
\endgroup\end{restoretext}
For these cubes (represented by the above manifold diagrams) to be elements of $\Comp(\sC)_2$, they need to be normalised, globular and well-typed. We leave the verification of these conditions to the reader, but illustrate the typability condition as before using manifold diagrams in the two following special cases
\begin{restoretext}
\begingroup\sbox0{\includegraphics{ANCimg/page1.png}}\includegraphics[clip,trim=0 {.1\ht0} 0 {.1\ht0} ,width=\textwidth]{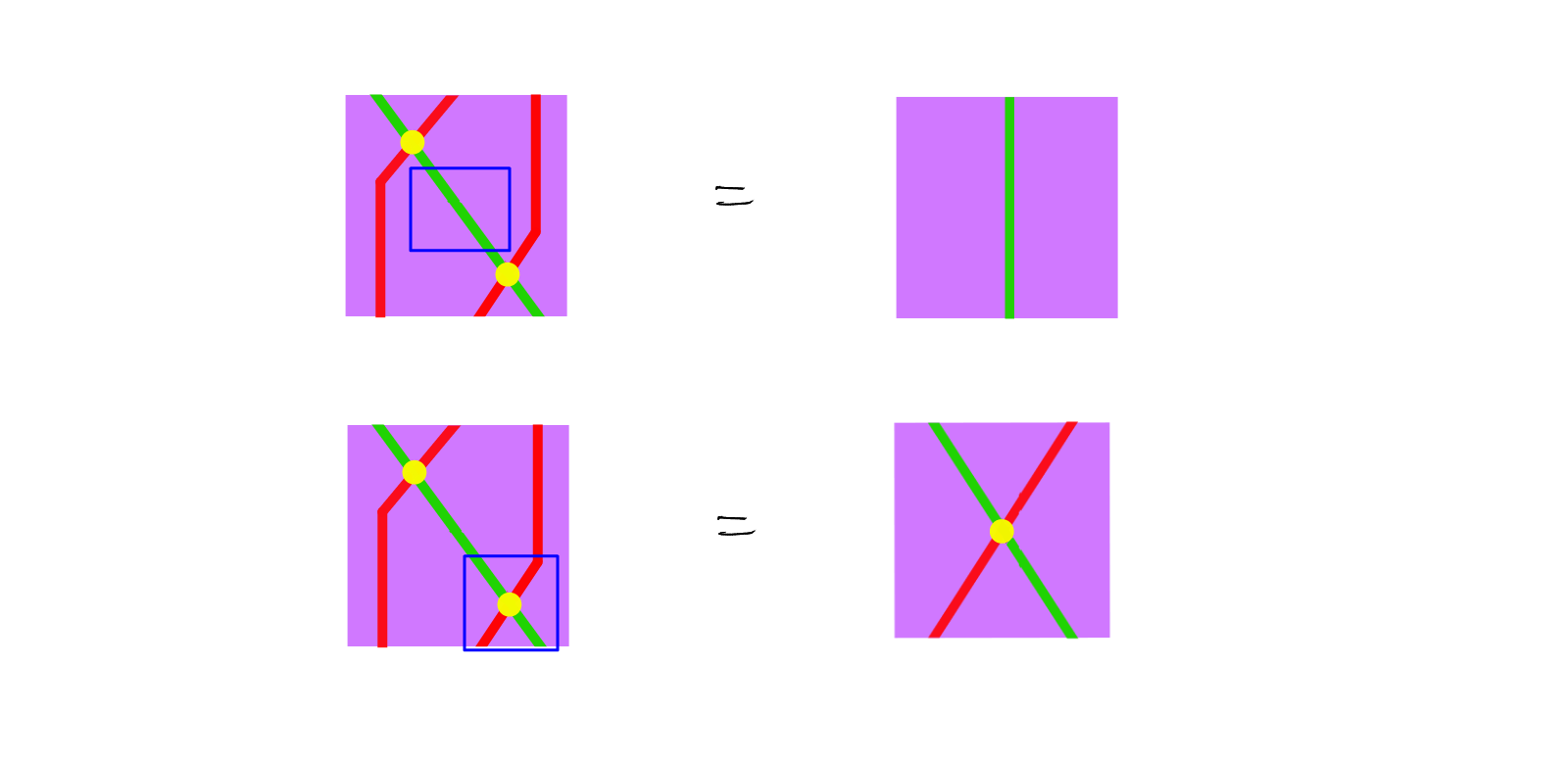}
\endgroup\end{restoretext}

In general, a \free{} associative $1$-category is a free category, with objects $\sC_0$, generating morphisms $\sC_1$ and generating equality relations $\sC_2$ between morphisms. As a category, its morphism set is given by $(\Comp(\sC)^{\quotg}_{\truncleq 1})_1$ (cf. \autoref{constr:globular_quotients}). For our specific example of $\sC$ above, we compute $\sC$ to be the free commutative monoid on two generators (that is, $\lN \times \lN$).

\subsection{Dimension 2} \label{sec:pres_2}

\Free{} associative $2$-categories consist of a set $\scC_0$ (of objects), a set $\scC_1$ (of generating $1$-morphisms), a set $\scC_2$ (of generating $2$-morphisms), and a set $\scC_3$ (of generating equalities between $2$-morphisms) together with an assignment of types. For instance, consider $\sC$ with objects
\begin{restoretext}
\begingroup\sbox0{\includegraphics{ANCimg/page1.png}}\includegraphics[clip,trim={.2\ht0} {.2\ht0} {.2\ht0} {.25\ht0} ,width=.9\textwidth]{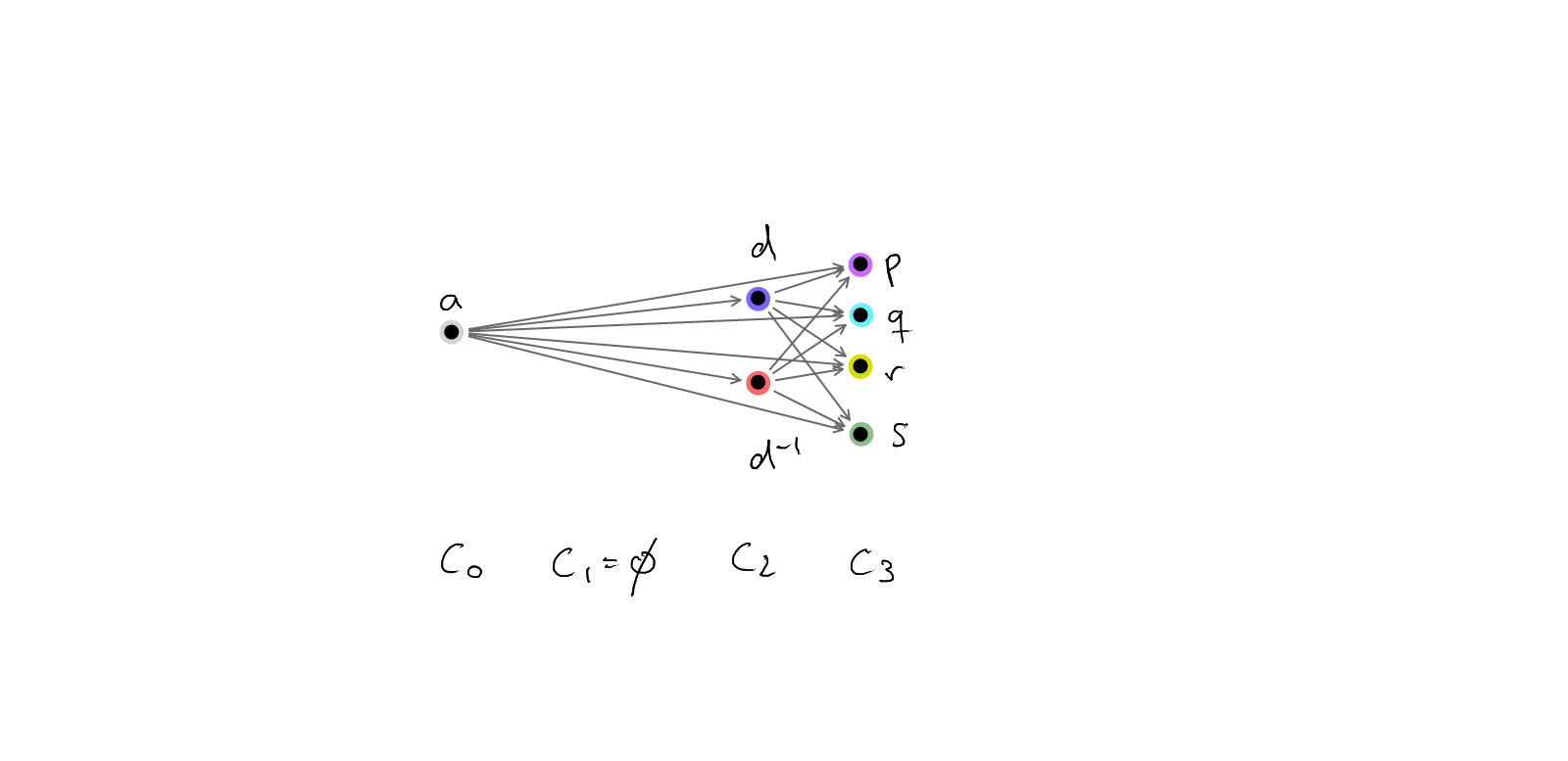}
\endgroup\end{restoretext}
in particular, note that here we chose $\sC_1 = \emptyset$.  We also depicted the poset $\GGamma{}\sC$ by its generating arrows (in \cgray{}). 

We give the assignment of types in dimension $2$ as follows (corresponding manifold diagrams are depicted on the right as before)
\begin{restoretext}
\begingroup\sbox0{\includegraphics{ANCimg/page1.png}}\includegraphics[clip,trim=0 {.15\ht0} 0 {.15\ht0} ,width=.9\textwidth]{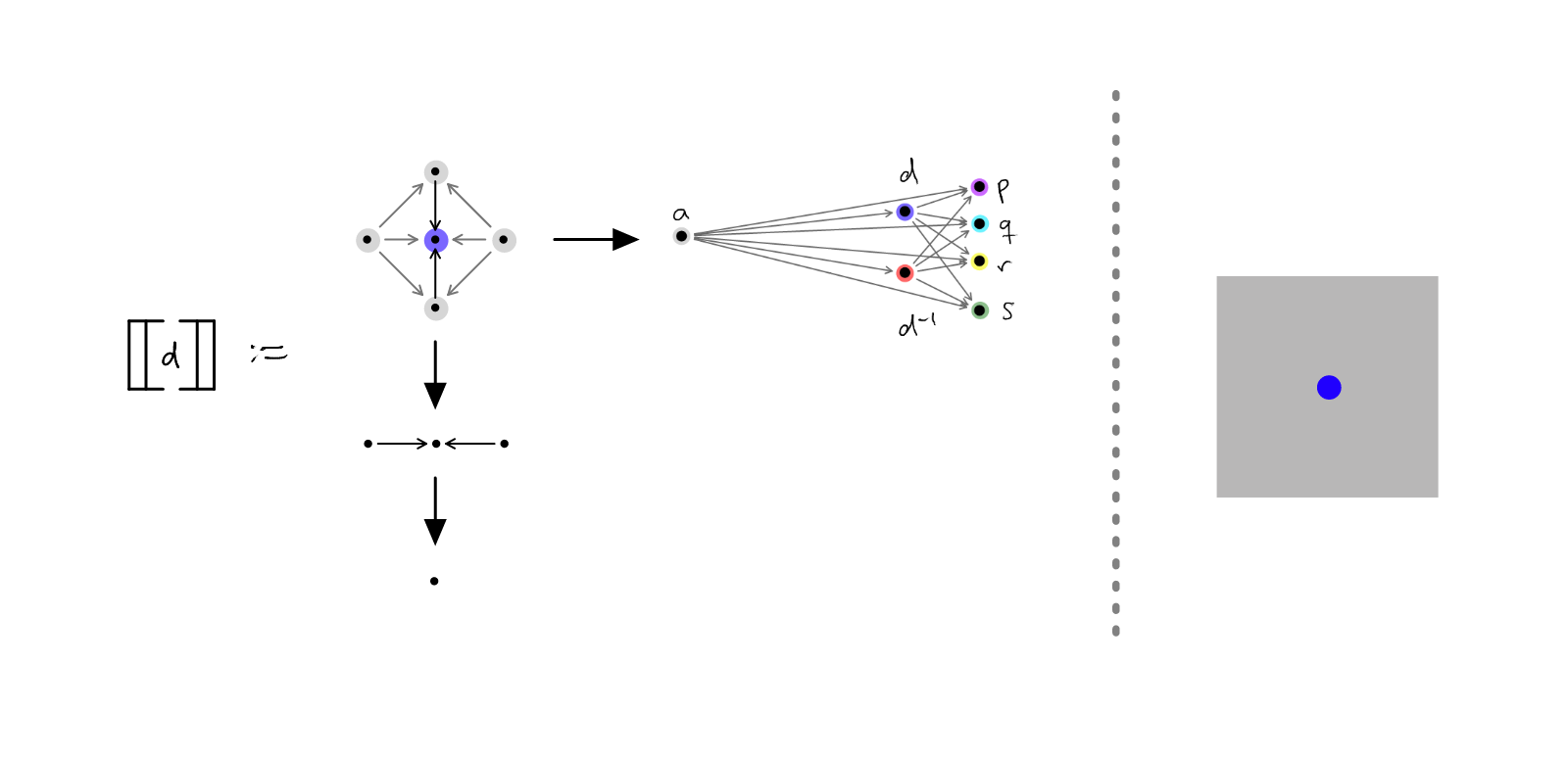}
\endgroup\end{restoretext}
\begin{restoretext}
\begingroup\sbox0{\includegraphics{ANCimg/page1.png}}\includegraphics[clip,trim=0 {.15\ht0} 0 {.15\ht0} ,width=.9\textwidth]{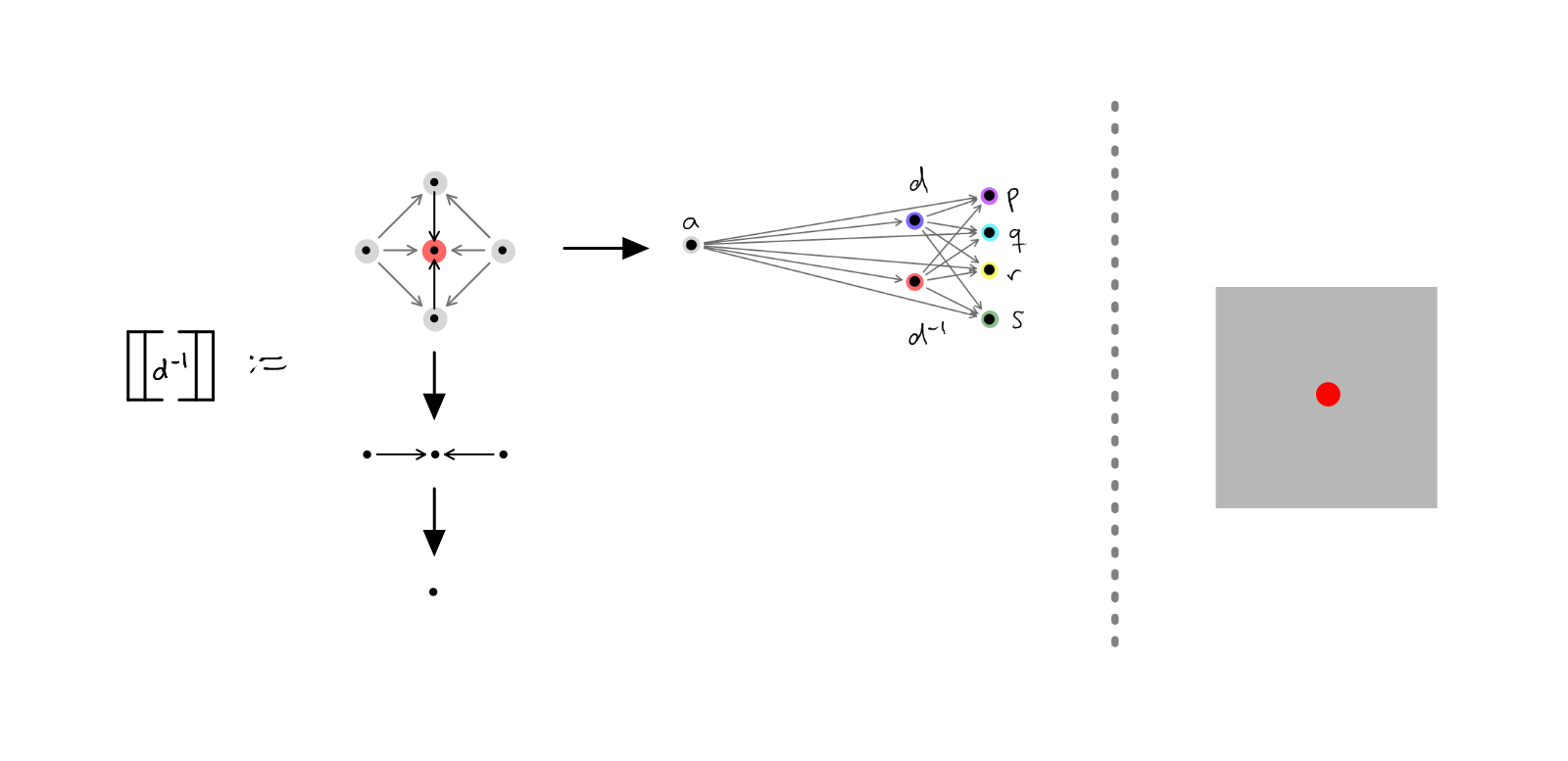}
\endgroup\end{restoretext}
Here, the cube data is defined on the left, and its associated manifold diagram is given on the right. Similarly, we define the the assignment of types in dimension $3$ by setting
\begin{restoretext}
\begingroup\sbox0{\includegraphics{ANCimg/page1.png}}\includegraphics[clip,trim=0 {.0\ht0} 0 {.0\ht0} ,width=.9\textwidth]{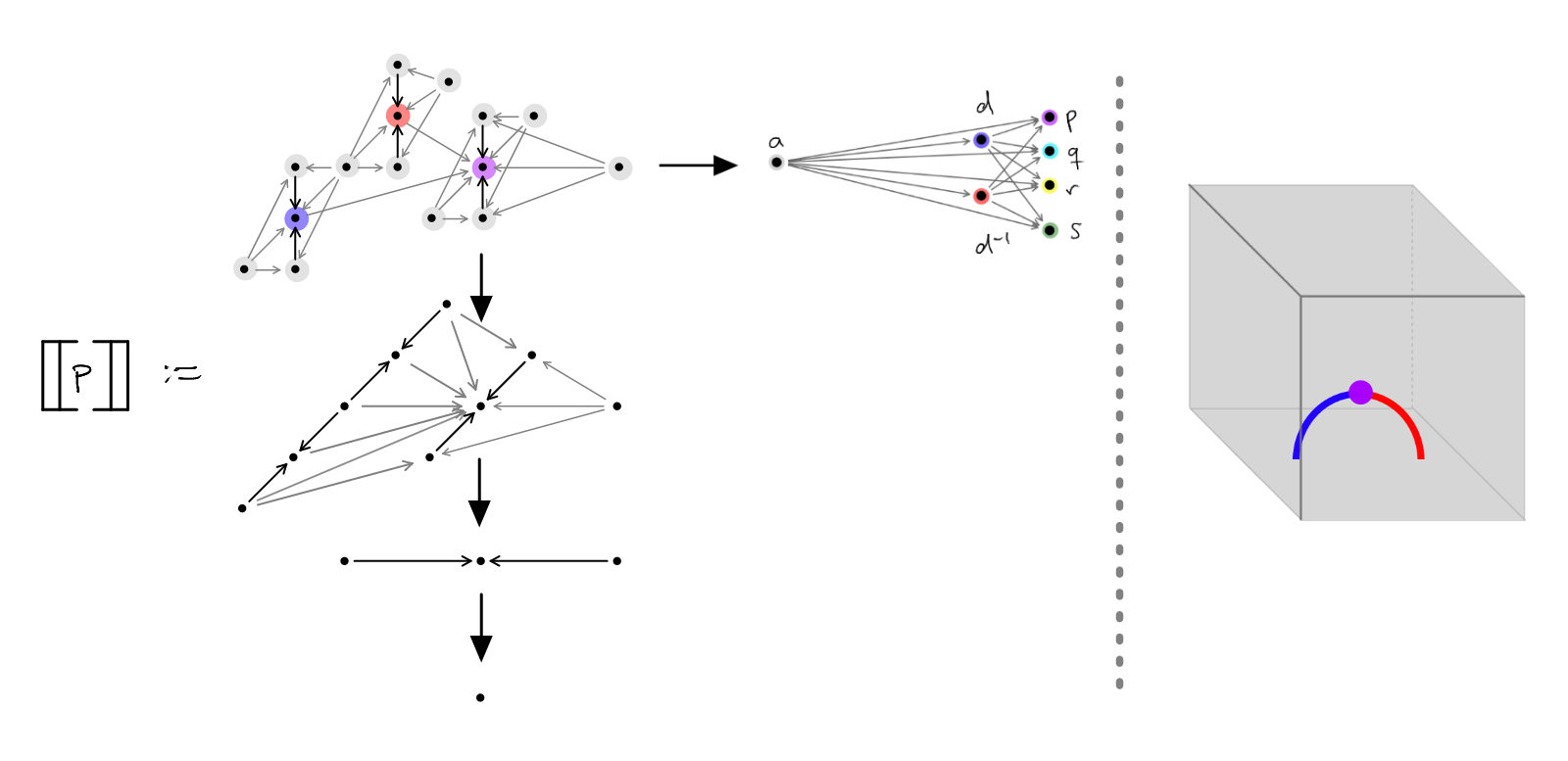}
\endgroup\end{restoretext}
and further (only depicting the manifold diagram and omitting cube data)
\begin{restoretext}
\begingroup\sbox0{\includegraphics{ANCimg/page1.png}}\includegraphics[clip,trim=0 {.3\ht0} 0 {.15\ht0} ,width=.9\textwidth]{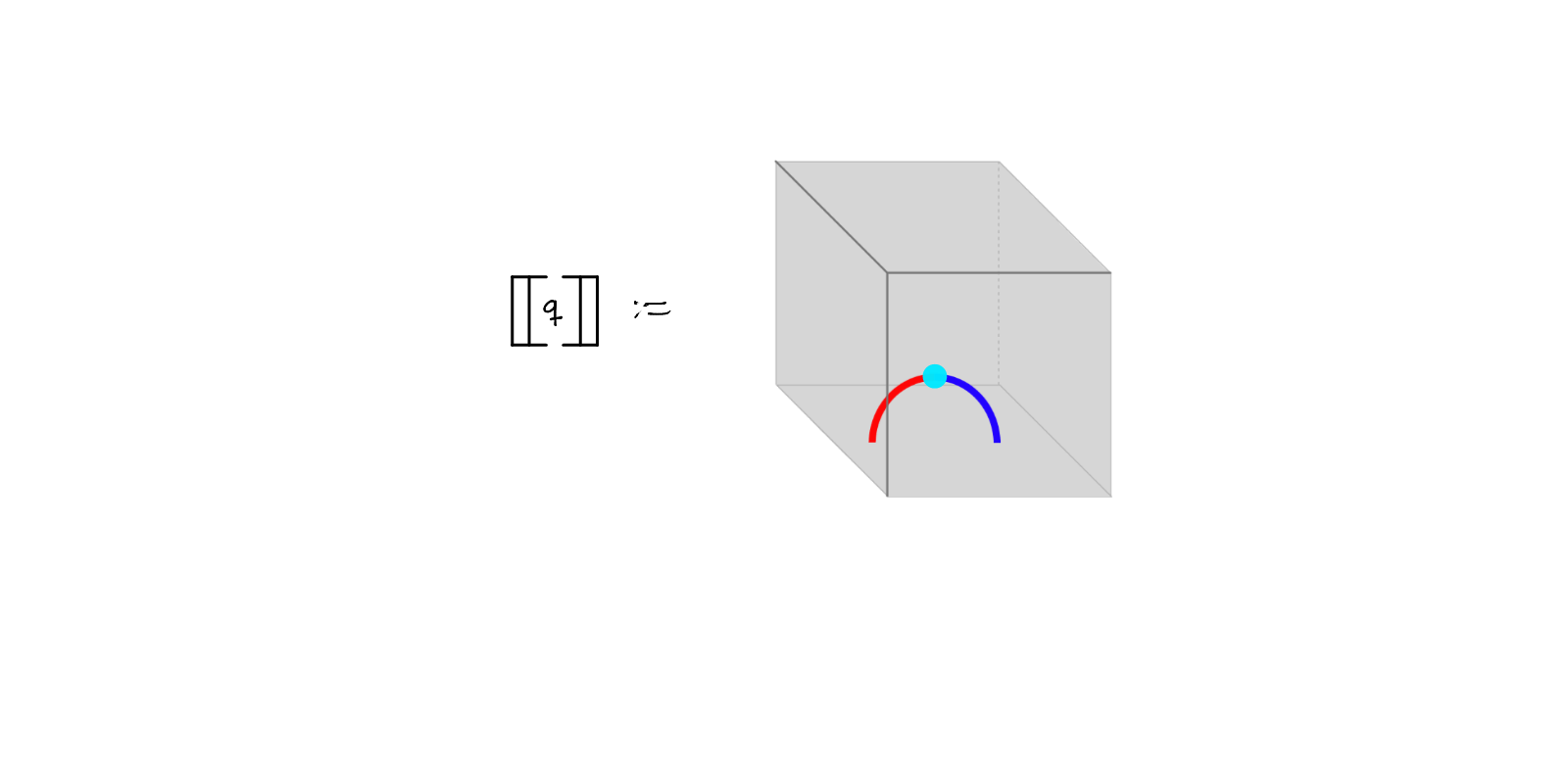}
\endgroup\end{restoretext}
\begin{restoretext}
\begingroup\sbox0{\includegraphics{ANCimg/page1.png}}\includegraphics[clip,trim=0 {.0\ht0} 0 {.0\ht0} ,width=.9\textwidth]{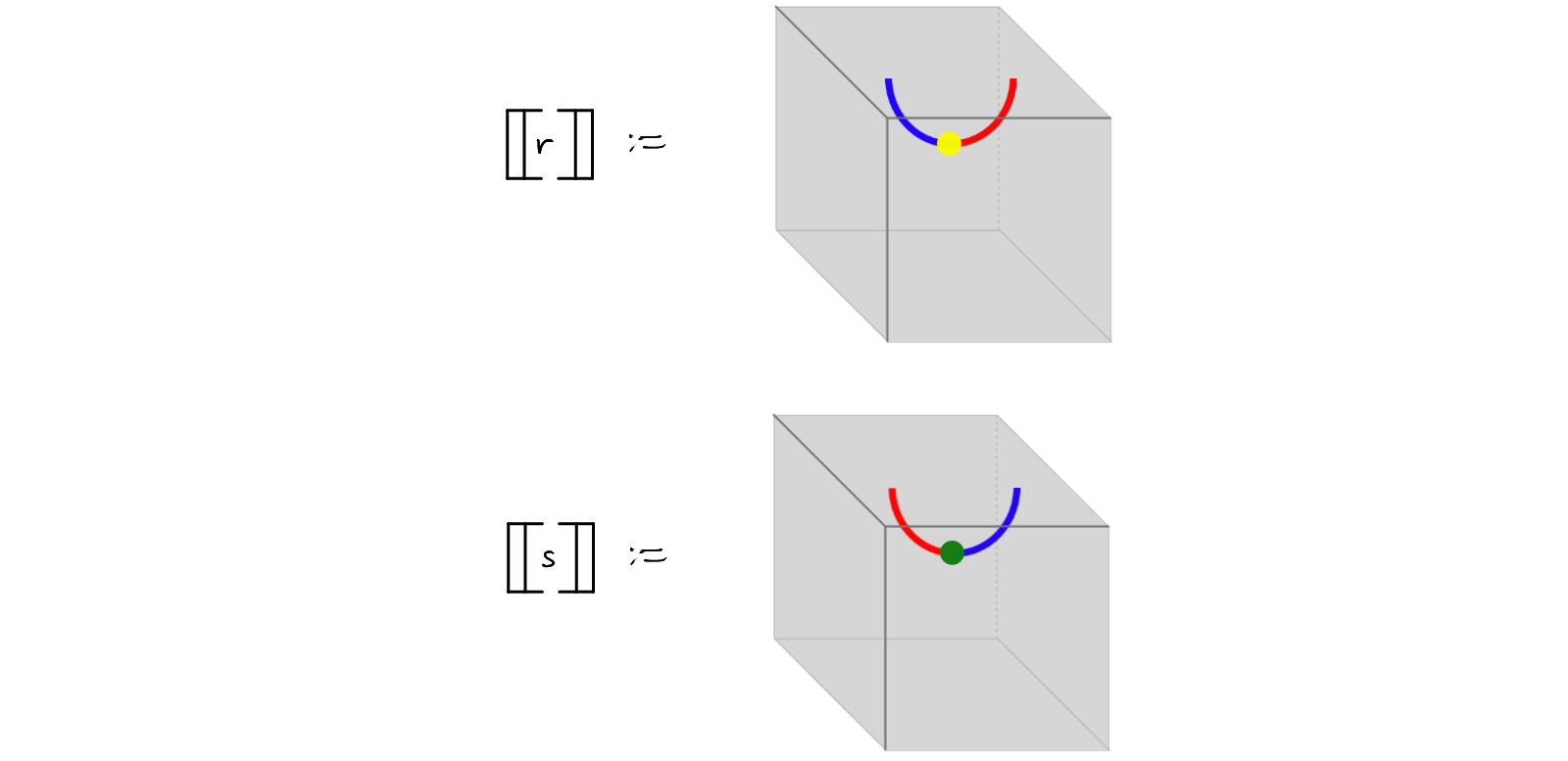}
\endgroup\end{restoretext}
As before elements of $\Comp(\sC)_3$ (cf. \autoref{defn:PANC_mor}) give an equivalence relation on $\Comp(\sC)_2$ by relating their sources and targets. However, unlike before $\Comp(\sC)_3$ not only contains composites of generating relations, but it now also contains ``homotopies". For instance, the following is a valid element $\scD$ of $\Comp(\sC)_3$
\begin{restoretext}
\begingroup\sbox0{\includegraphics{ANCimg/page1.png}}\includegraphics[clip,trim=0 {.0\ht0} 0 {.0\ht0} ,width=.9\textwidth]{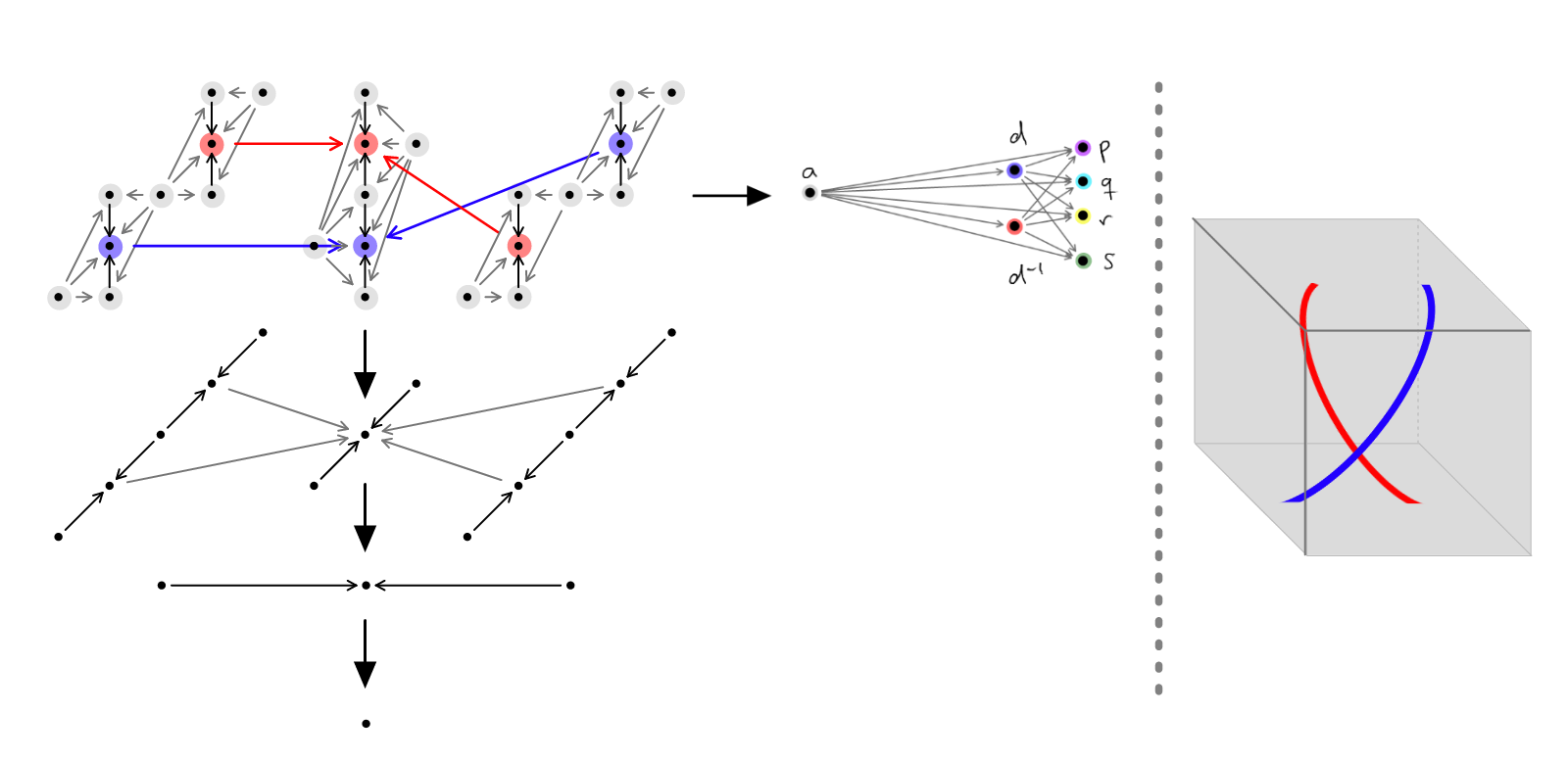}
\endgroup\end{restoretext}
To see this we need to check typability. Consider for instance the following minimal neighbourhood (note that the small cube on the right is supposed to only intersect the \cred{} wire, but not the \cblue{} one)
\begin{restoretext}
\begingroup\sbox0{\includegraphics{ANCimg/page1.png}}\includegraphics[clip,trim=0 {.0\ht0} 0 {.0\ht0} ,width=.9\textwidth]{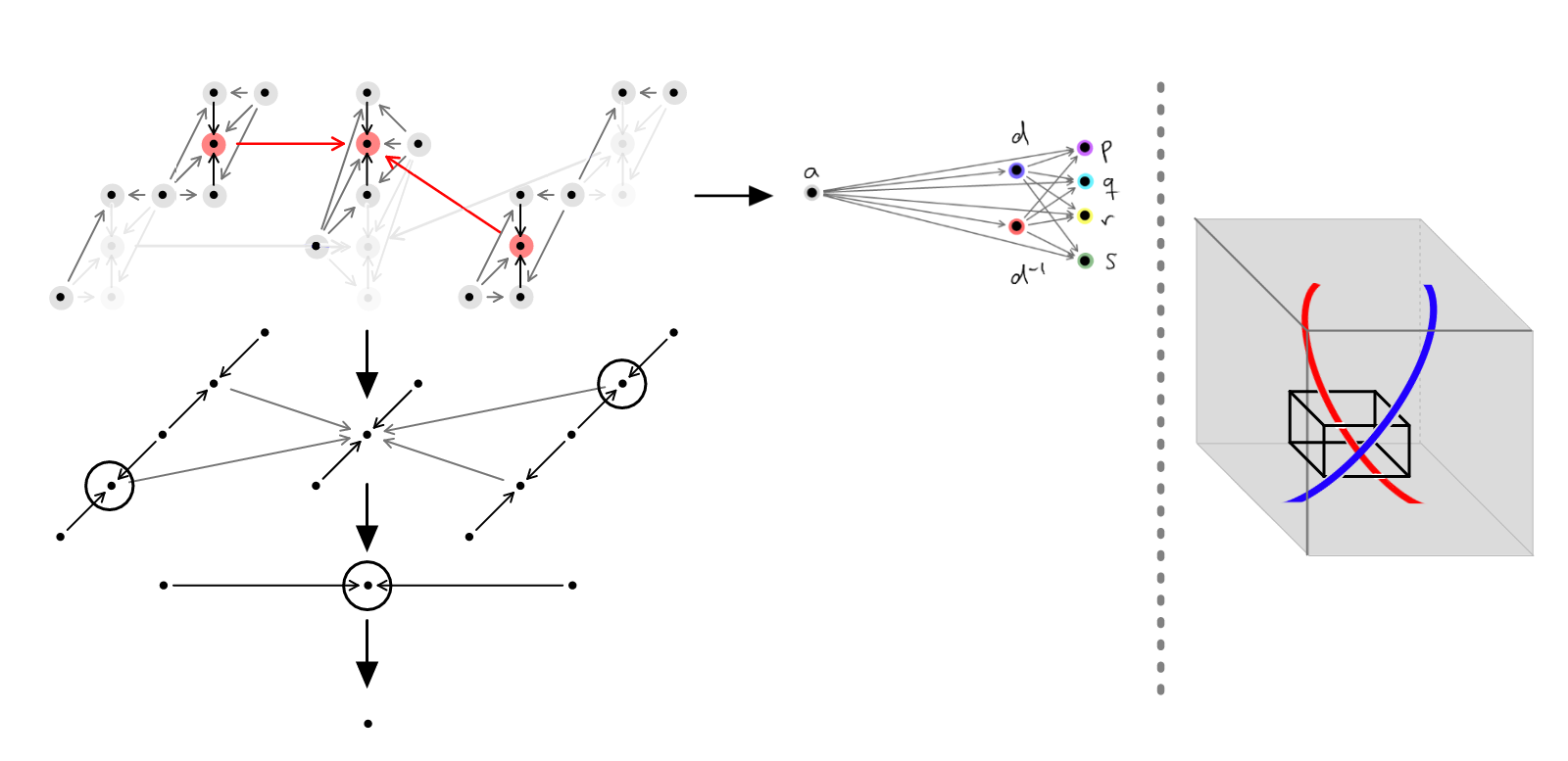}
\endgroup\end{restoretext}
Using a $2$-level collapse which collapses the singular heights circled in \cblack{} in $\tsG 2(\scD)$ and then a $1$-level collapse which collapses the singular height circled in \cblack{} in $\tsG 1(\scD)$ we find this subcube collapses to
\begin{restoretext}
\begingroup\sbox0{\includegraphics{ANCimg/page1.png}}\includegraphics[clip,trim=0 {.05\ht0} 0 {.1\ht0} ,width=.9\textwidth]{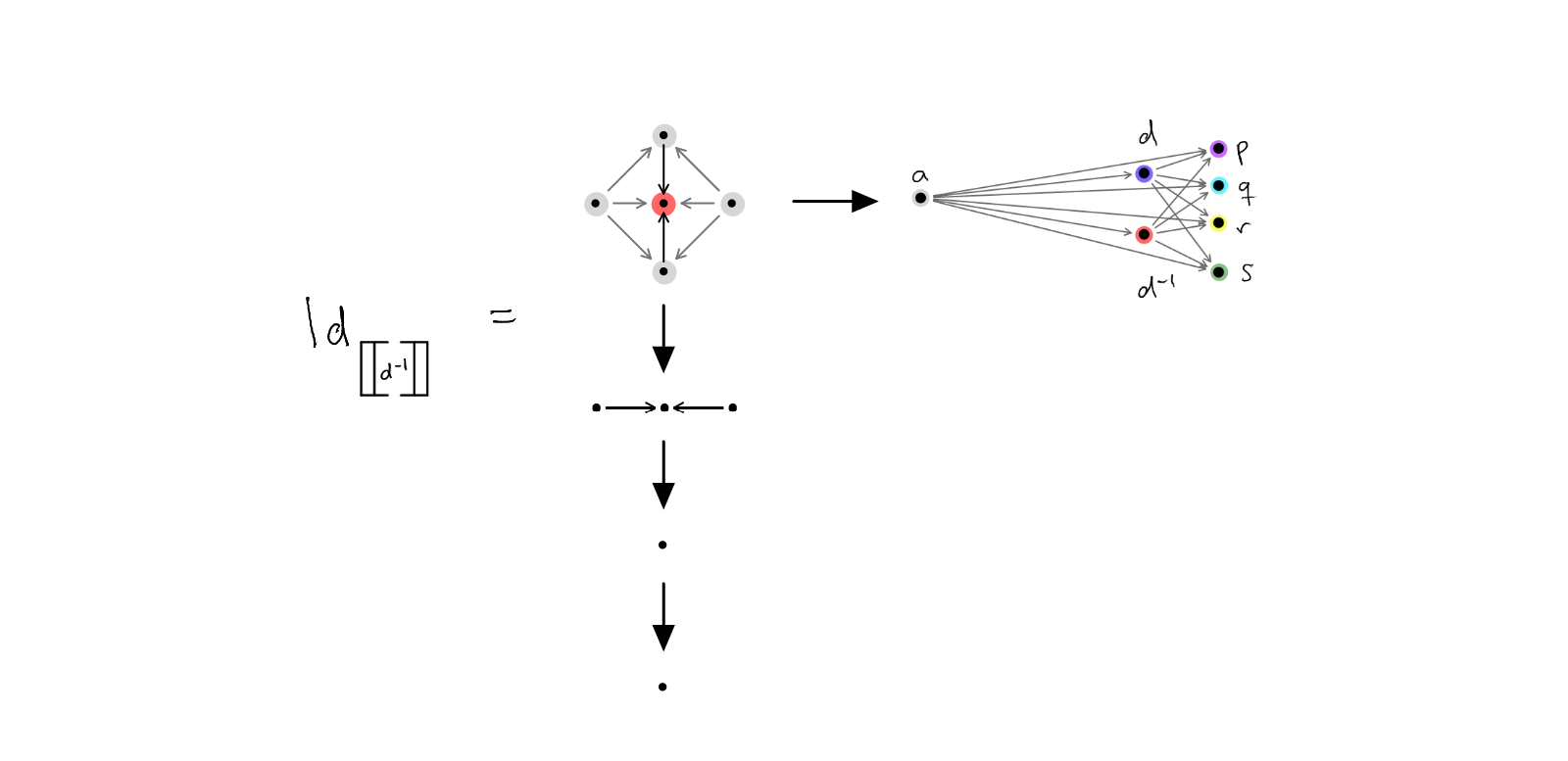}
\endgroup\end{restoretext}
as required for typability. We will call $\scD$ an interchanger homotopy. The equality relation that $\scD$ gives, permutes the order of $d$ and $d\inv$. There are similar homotopies permuting more copies of $d$ and $d\inv$. We obtain equalities such as
\begin{restoretext}
\begingroup\sbox0{\includegraphics{ANCimg/page1.png}}\includegraphics[clip,trim=0 {.1\ht0} 0 {.1\ht0} ,width=.9\textwidth]{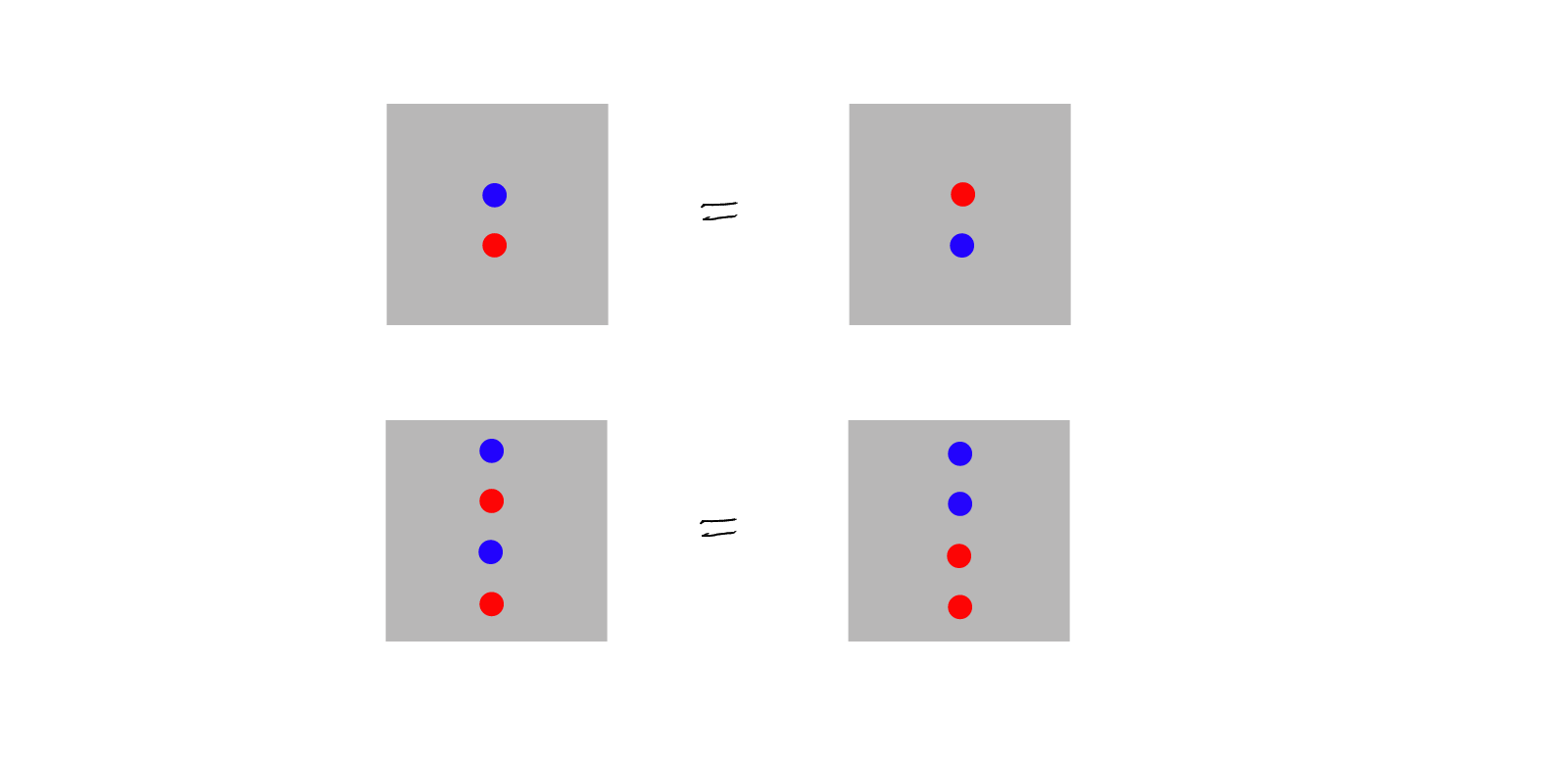}
\endgroup\end{restoretext}
On the other hand, the generating relations (that is, elements of $\sC_3$), 
 allow us to ``cancel" adjacent occurrences of $\alpha$ and $\beta$, and yield equalities such as
\begin{restoretext}
\begingroup\sbox0{\includegraphics{ANCimg/page1.png}}\includegraphics[clip,trim=0 {.35\ht0} 0 {.25\ht0} ,width=.9\textwidth]{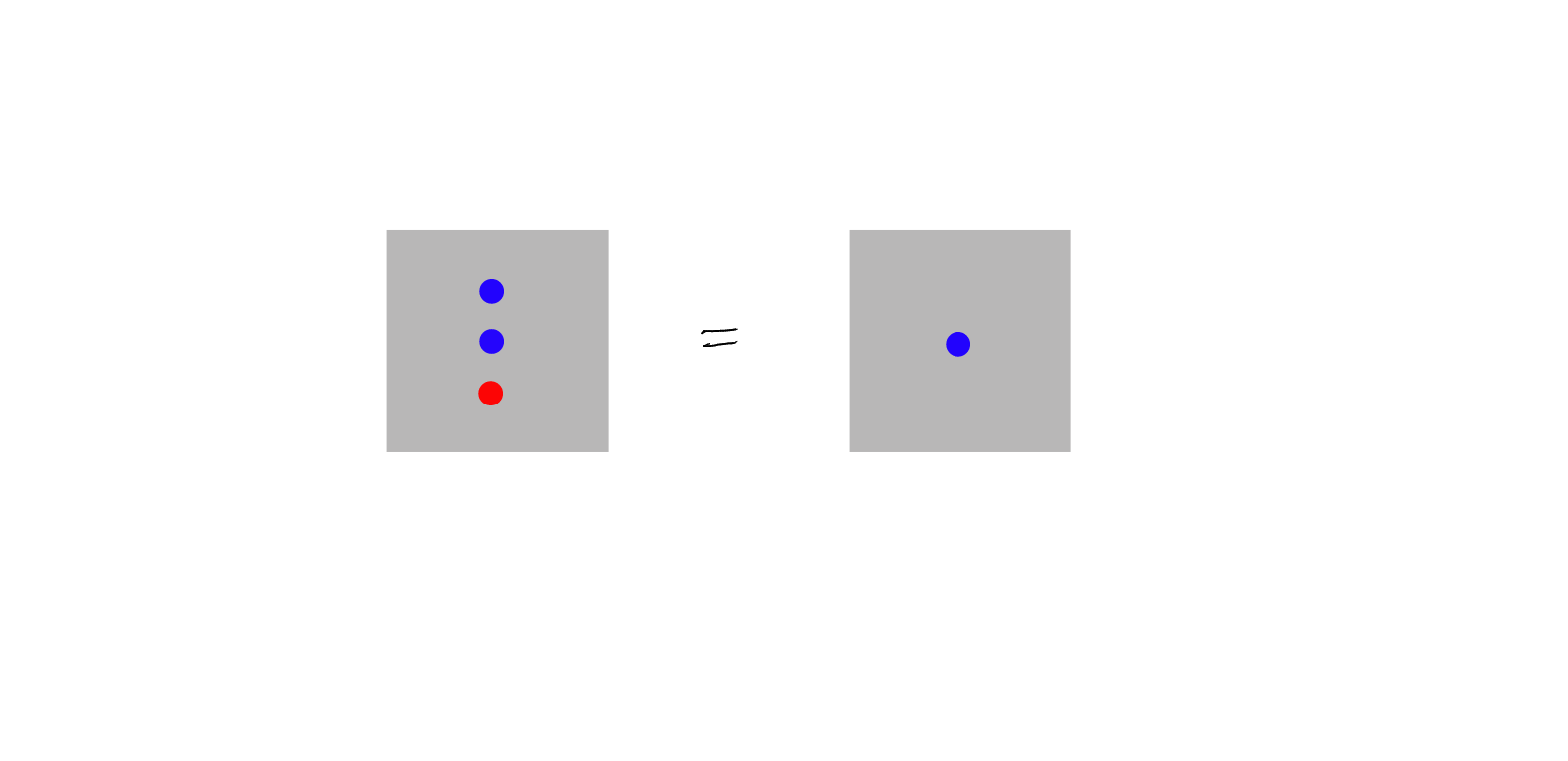}
\endgroup\end{restoretext}
Recall $(\Comp(\sC)^{\quotg}_{\truncleq 2})_2$ is $\Comp(\sC)_2$ quotiented by the relation relating sources and targets of elements in $\Comp(\sC)_3$ (cf. \autoref{constr:globular_quotients}). We find that $(\Comp(\sC)^{\quotg}_{\truncleq 2})_2$ is given by
\begin{restoretext}
\begingroup\sbox0{\includegraphics{ANCimg/page1.png}}\includegraphics[clip,trim=0 {.3\ht0} 0 {.3\ht0} ,width=.9\textwidth]{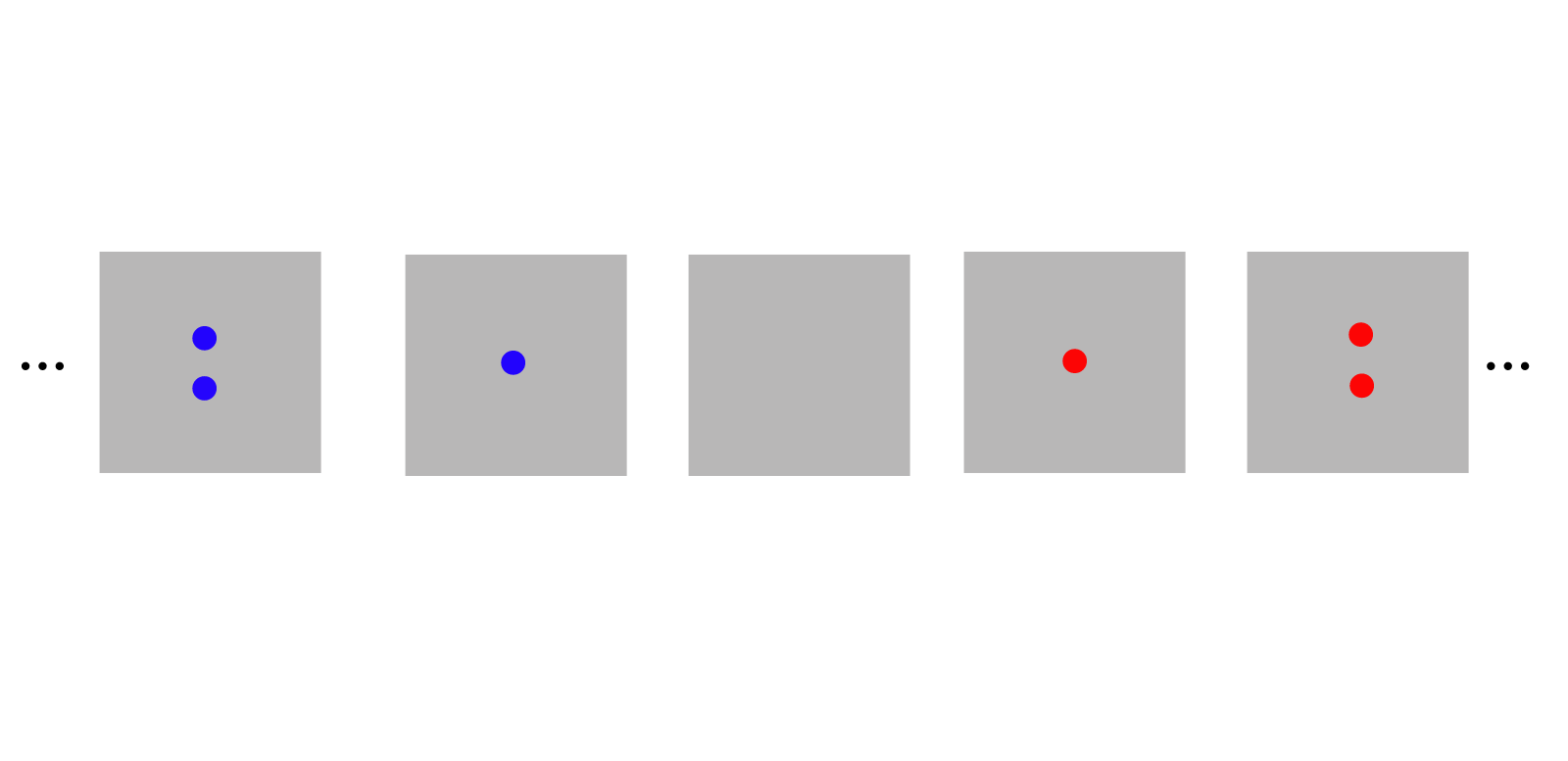}
\endgroup\end{restoretext}
By gluing $2$-globes along their $1$-boundary (formally, this is $\whisker 1 2$ composition described in \autoref{ch:composition}) this set $(\Comp(\sC)^{\quotg}_{\truncleq 2})_2$ obtains a group structure, which can be seen to be isomorphic to $(\lZ,+)$. It is no coincidence that this is the structure of the second homotopy group of the $2$-sphere $S^2$
\begin{equation}
\pi_2(S^2) = \lZ
\end{equation}
In fact, our category $\sC$ is a representation of the $2$-truncated homotopy type of $S^2$. We will learn about its $3$-truncation in the next dimension.

\subsection{Dimension 3} \label{sec:pres_3}

\Free{} associative $3$-categories consist of a set $\scC_0$ (of objects), for $1 \leq i \leq 3$, a set $\scC_i$ (of generating $i$-morphisms), a set $\scC_2$ (of generating $2$-morphisms), and a set $\scC_4$ (of generating equalities between $3$-morphisms). Consider for instance the following choice of $\sC_i$
\begin{restoretext}
\begingroup\sbox0{\includegraphics{ANCimg/page1.png}}\includegraphics[clip,trim={.25\ht0} {.2\ht0} {.25\ht0} {.25\ht0} ,width=\textwidth]{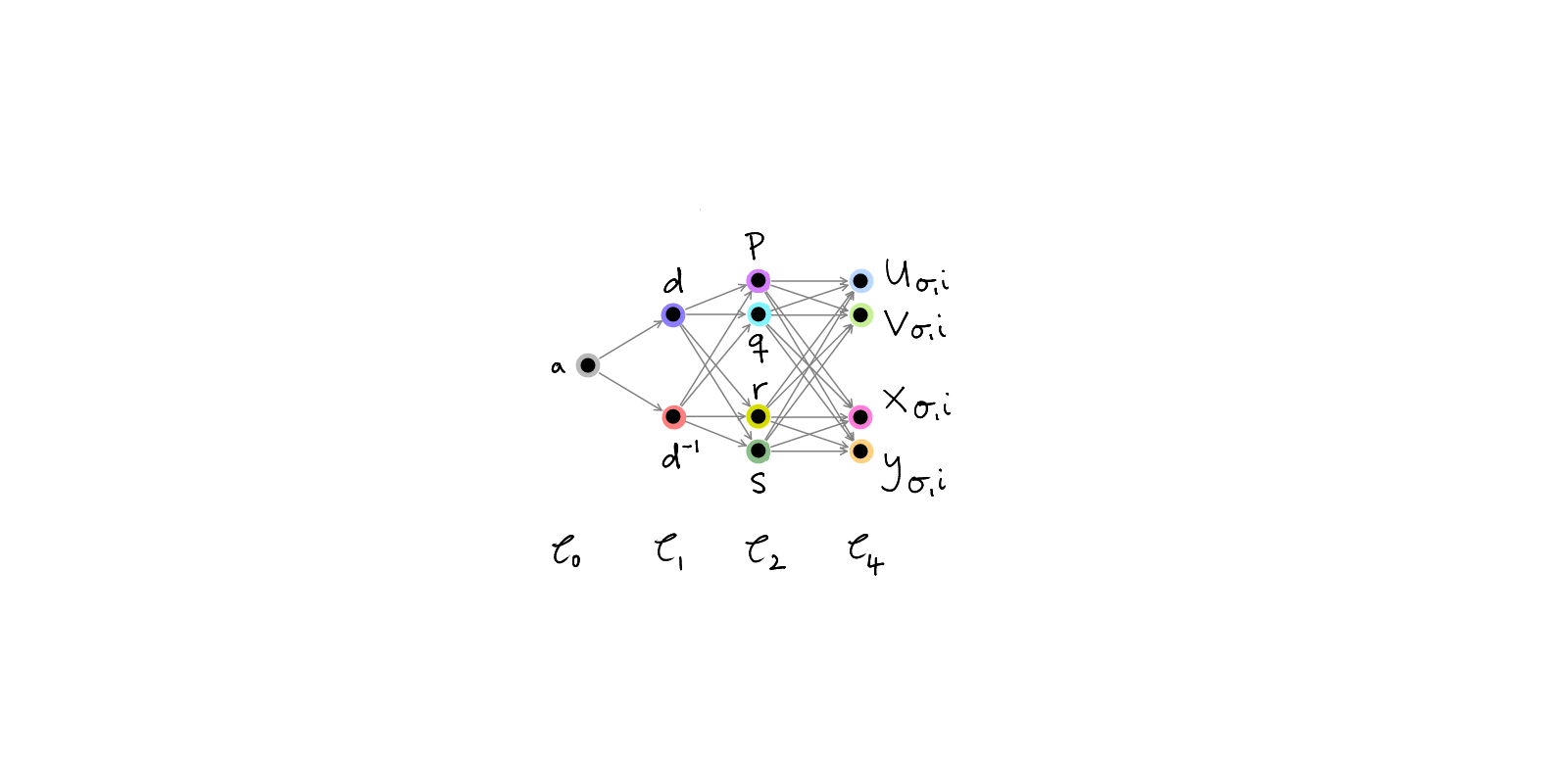}
\endgroup\end{restoretext}
Here, the indices $\sigma$ and $i$ range over $\sigma \in \Set{+,-}$ and $i \in \Set{1,-1}$ respectively (thus, $\sC_4$ contains a total of 16 elements). We also depicted the poset $\GGamma{}\sC$ by its generating arrows (in \cgray{}).

The assignment of types for $\sC$ in dimension 0,1,2 and 3 is the same as for our example in dimension $2$. It remains to assign types for the elements of $\sC_4$. The type of $u_{+,1}$ is given by the following $4$-cube (a manifold diagram representation is given after the formal algebraic definition below)
\begin{restoretext}
\begin{noverticalspace}
\begingroup\sbox0{\includegraphics{ANCimg/page1.png}}\includegraphics[clip,trim=0 {.0\ht0} 0 {.1\ht0} ,width=.9\textwidth]{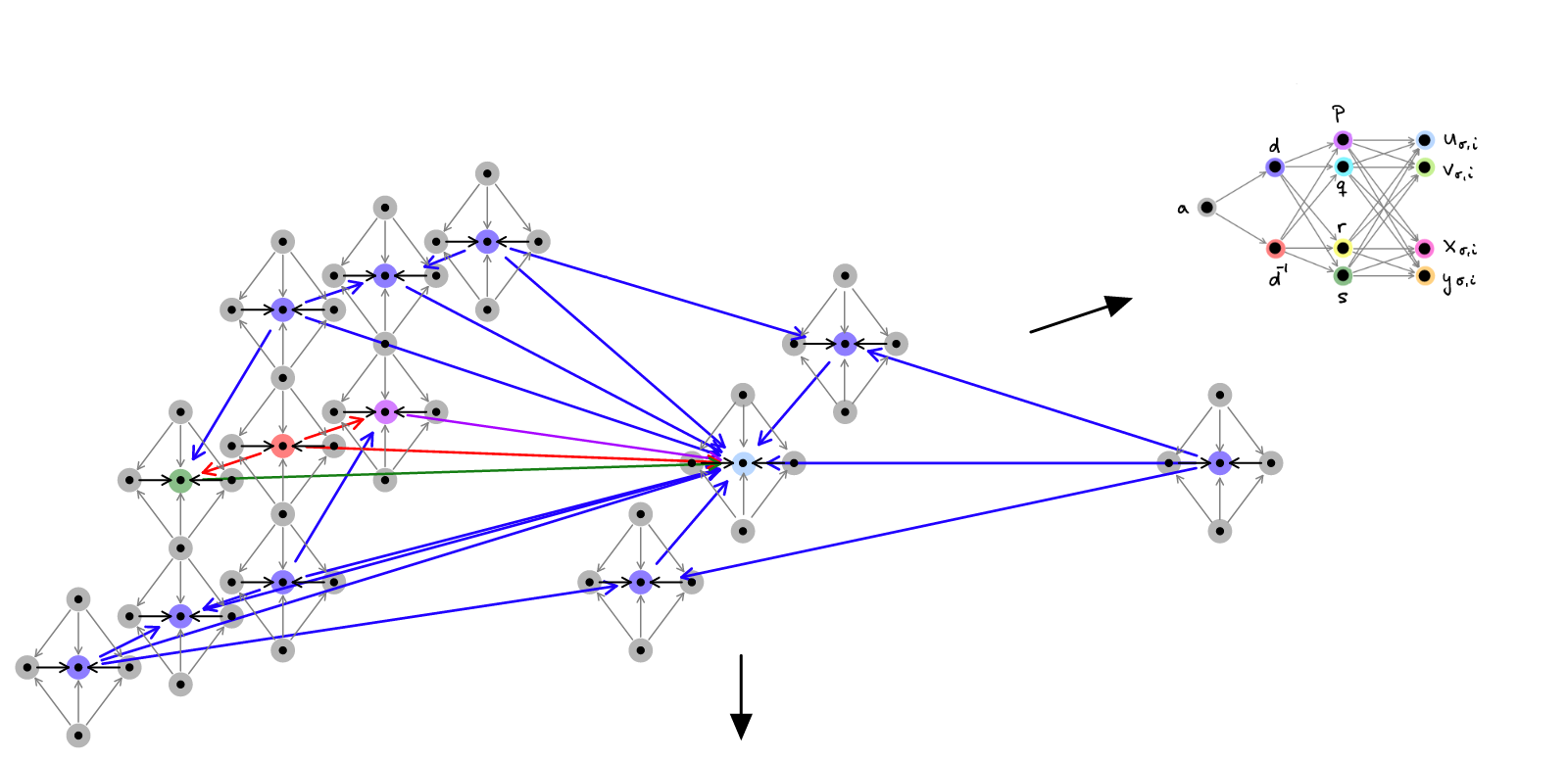}
\endgroup \\*
\begingroup\sbox0{\includegraphics{ANCimg/page1.png}}\includegraphics[clip,trim=0 {.2\ht0} 0 {.0\ht0} ,width=.9\textwidth]{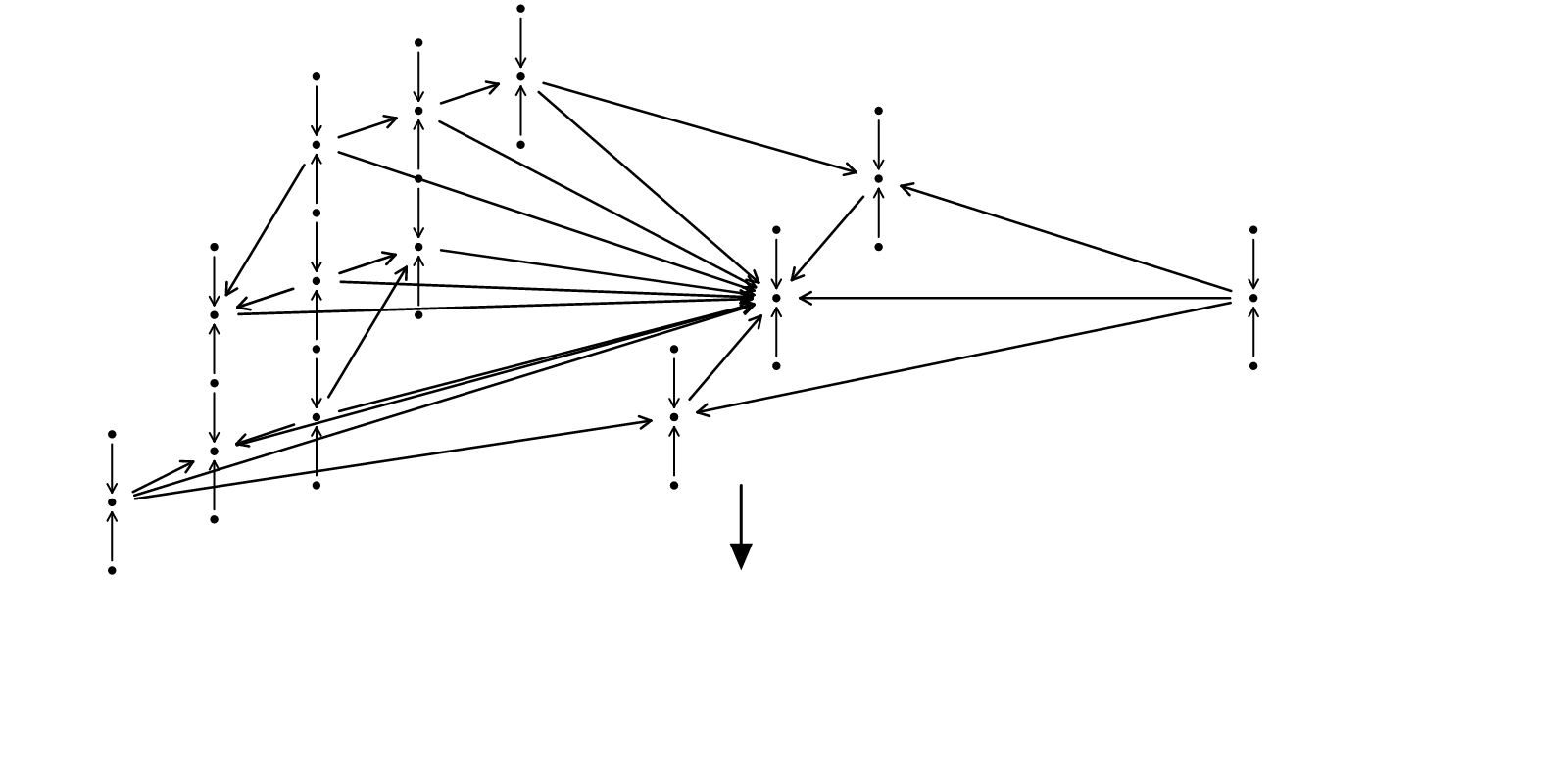}
\endgroup \\*
\begingroup\sbox0{\includegraphics{ANCimg/page1.png}}\includegraphics[clip,trim=0 {.0\ht0} 0 {.0\ht0} ,width=.9\textwidth]{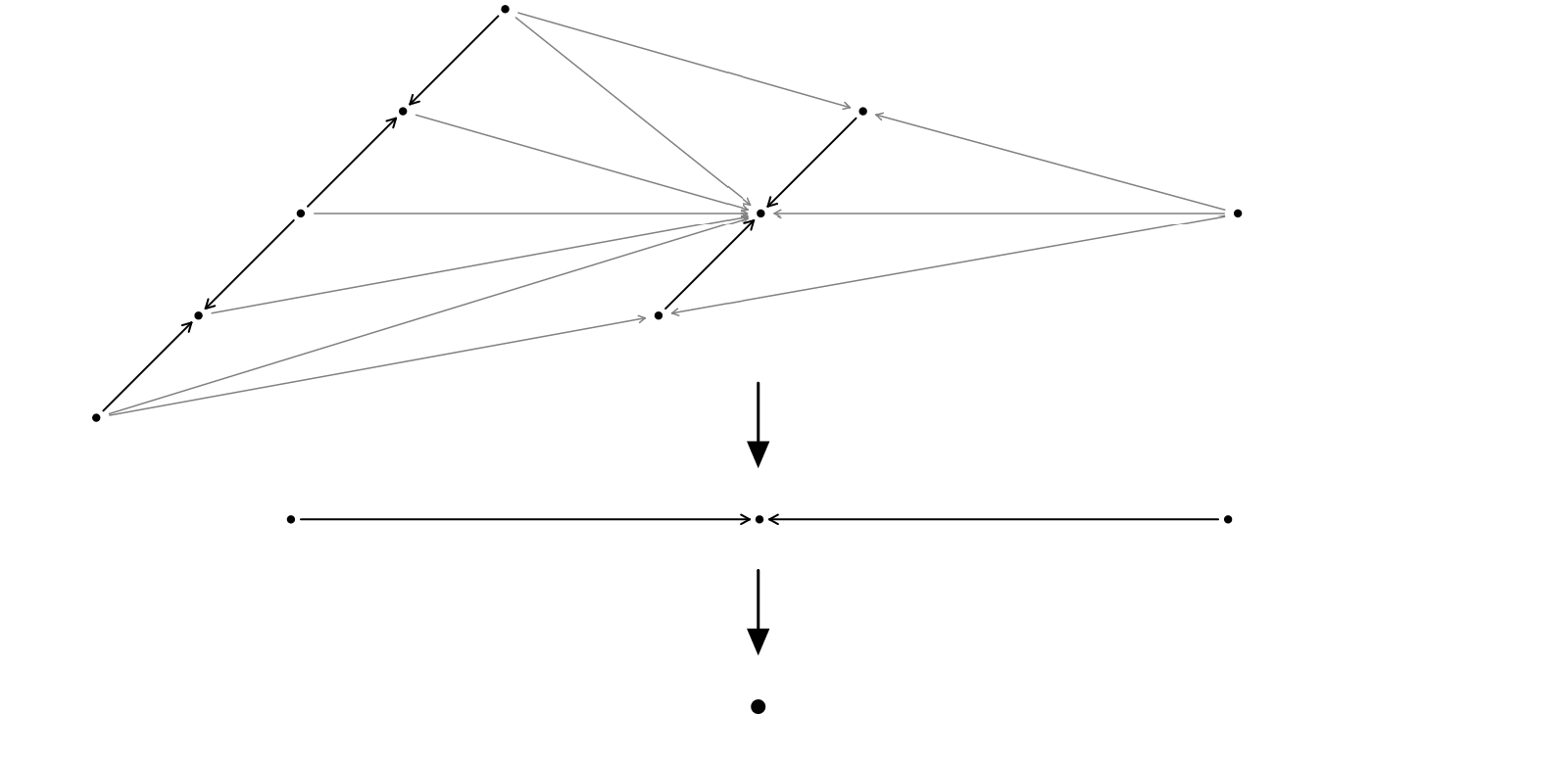}
\endgroup
\end{noverticalspace}
\end{restoretext}
Using manifold diagrams, this can be represented as follows
\begin{restoretext}
\begingroup\sbox0{\includegraphics{ANCimg/page1.png}}\includegraphics[clip,trim=0 {.0\ht0} 0 {.2\ht0} ,width=\textwidth]{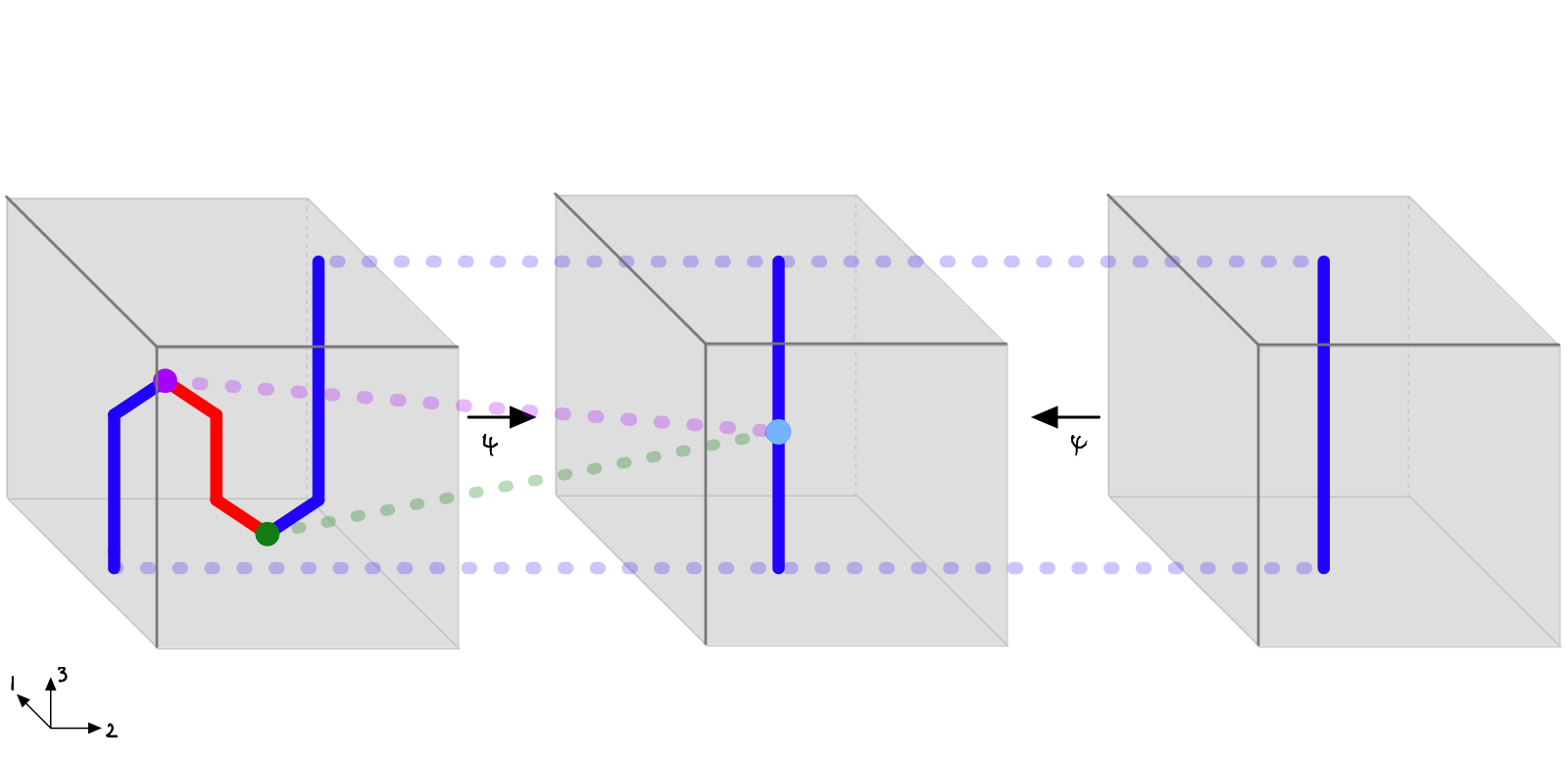}
\endgroup\end{restoretext}
Geometrically speaking, along the left arrow the \cpurple{} and \cdarkgreen{} points converge closer and closer together until they meet in the \clightblue{} point as represented in the middle picture. 

$\abss{u_{-,1}}$ is defined by the same type, but with \cblue{} and \cred{} colors interchanged (note that this also switches \cpurple{} to \cturquoise{} and \cdarkgreen{} to \cyellow{}). $\abss{u_{\sigma,-1}}$ is the ``inverse" of $\abss{u_{\sigma,1}}$, and is defined by ``reading $\abss{u_{\sigma,1}}$ from right to left" (formally, this can be achieved by using the automorphism $\singint 1 \to \singint 1$ that maps $0 \mapsto 2$, $1 \mapsto 1$, $2 \mapsto 0$, and precomposing $\tsU 1_{\abss{u_{\sigma,1}}}$ with it).

Similarly, $\abss{v_{+,1}}$ is defined by
\begin{restoretext}
\begingroup\sbox0{\includegraphics{ANCimg/page1.png}}\includegraphics[clip,trim=0 {.0\ht0} 0 {.2\ht0} ,width=\textwidth]{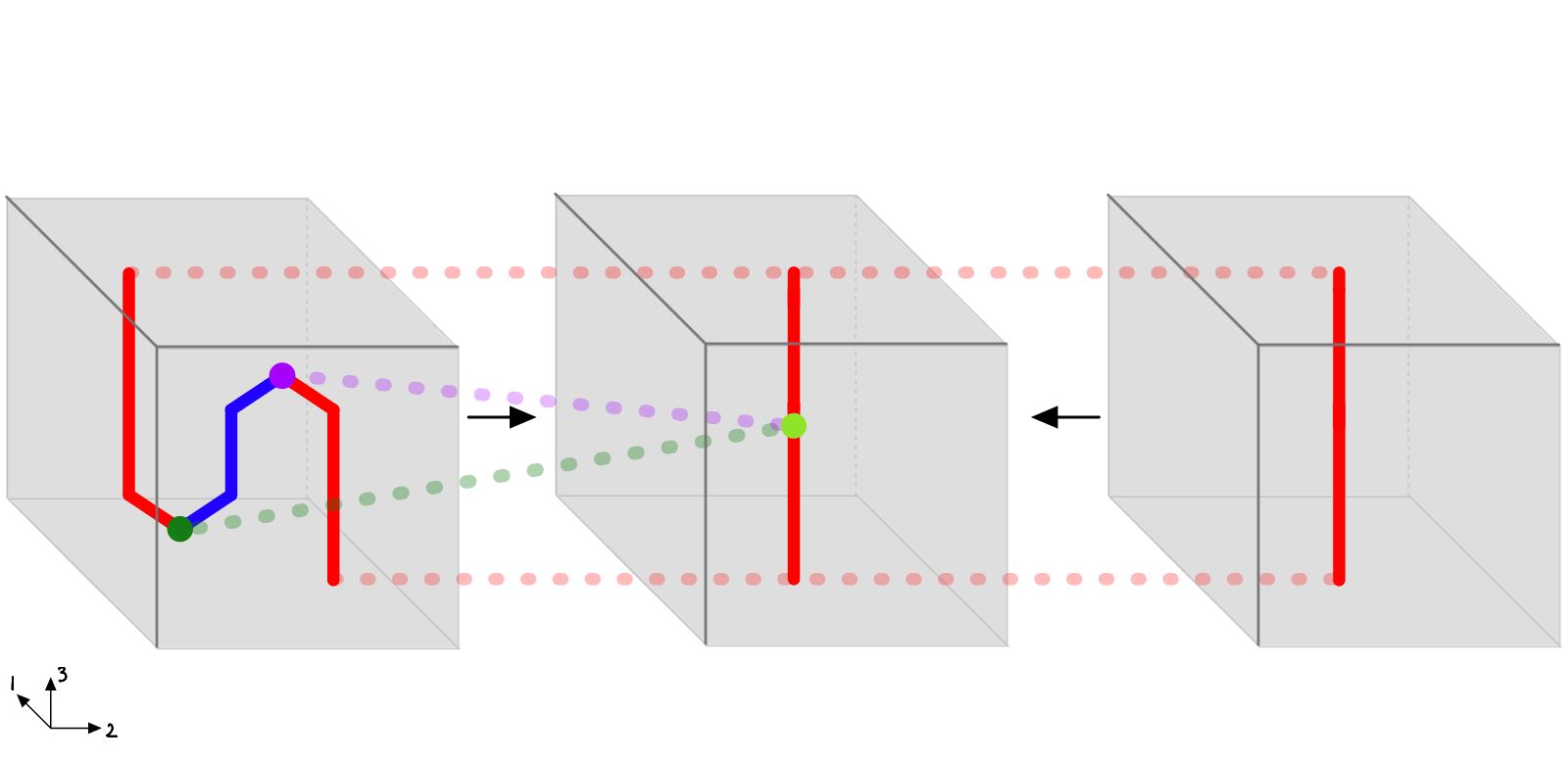}
\endgroup\end{restoretext}
and $\abss{v_{-,1}}$ is obtained by exchanging \cblue{} and \cred{} colors. $\abss{v_{\sigma,-1}}$ are the inverses of $\abss{v_{\sigma,1}}$ as before.

We remark that the types of $u_{\sigma,i}$ and $v_{\sigma,i}$ are usually referred to as ``snake singularities" and this is usually a $3$-dimensional singularity whereas here we meet the snake in dimension $4$. The fact that the singularity already lives in dimension $3$ can be seen by looking at the $3$-level total poset of the above tower of \SI-bundles. More formally, in the next chapter we will see the snake to be just one of many singularities in the ``theory of an invertible \textit{$1$-morphism}" $\TI$, where it will appear as 3-dimensional singularity. We will then understand, that the reason why the snake appears in dimension 4 here, is that we applied the theory of invertibility to a \textit{$2$-morphism} $d$.

Importantly, we will now see that the snake in dimension 4 is more interesting than the snake in dimension 3: it will allow for interesting interactions with homotopies such as the interchange. For instance, we have a $4$-dimensional cube 
\begin{restoretext}
\begingroup\sbox0{\includegraphics{ANCimg/page1.png}}\includegraphics[clip,trim=0 {.0\ht0} 0 {.2\ht0} ,width=\textwidth]{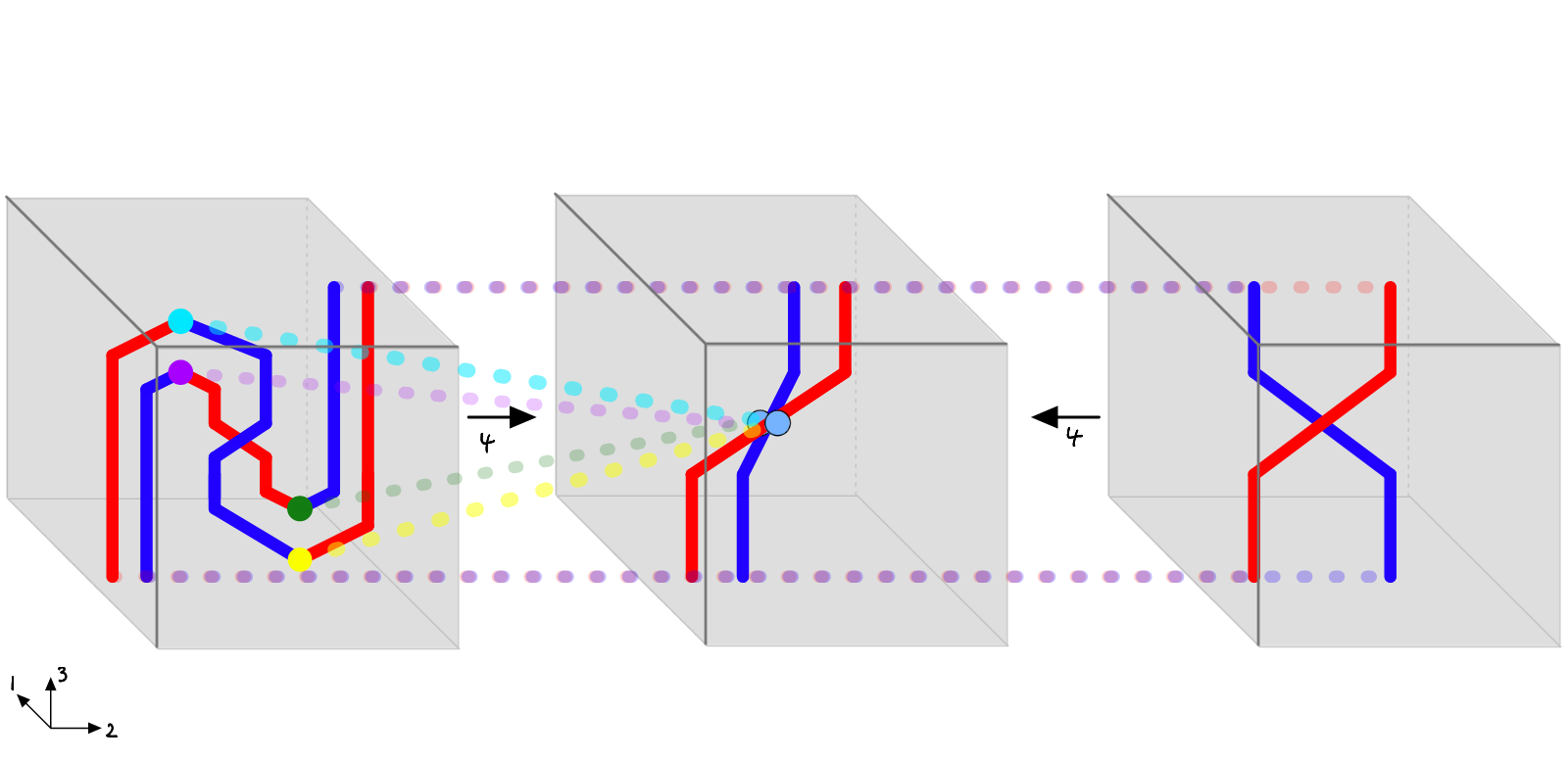}
\endgroup\end{restoretext}

Next, we assign $x_{+,1}$ the type $\abss{x_{+,1}}$ given by the $4$-dimensional manifold diagram
\begin{restoretext}
\begingroup\sbox0{\includegraphics{ANCimg/page1.png}}\includegraphics[clip,trim=0 {.0\ht0} 0 {.2\ht0} ,width=\textwidth]{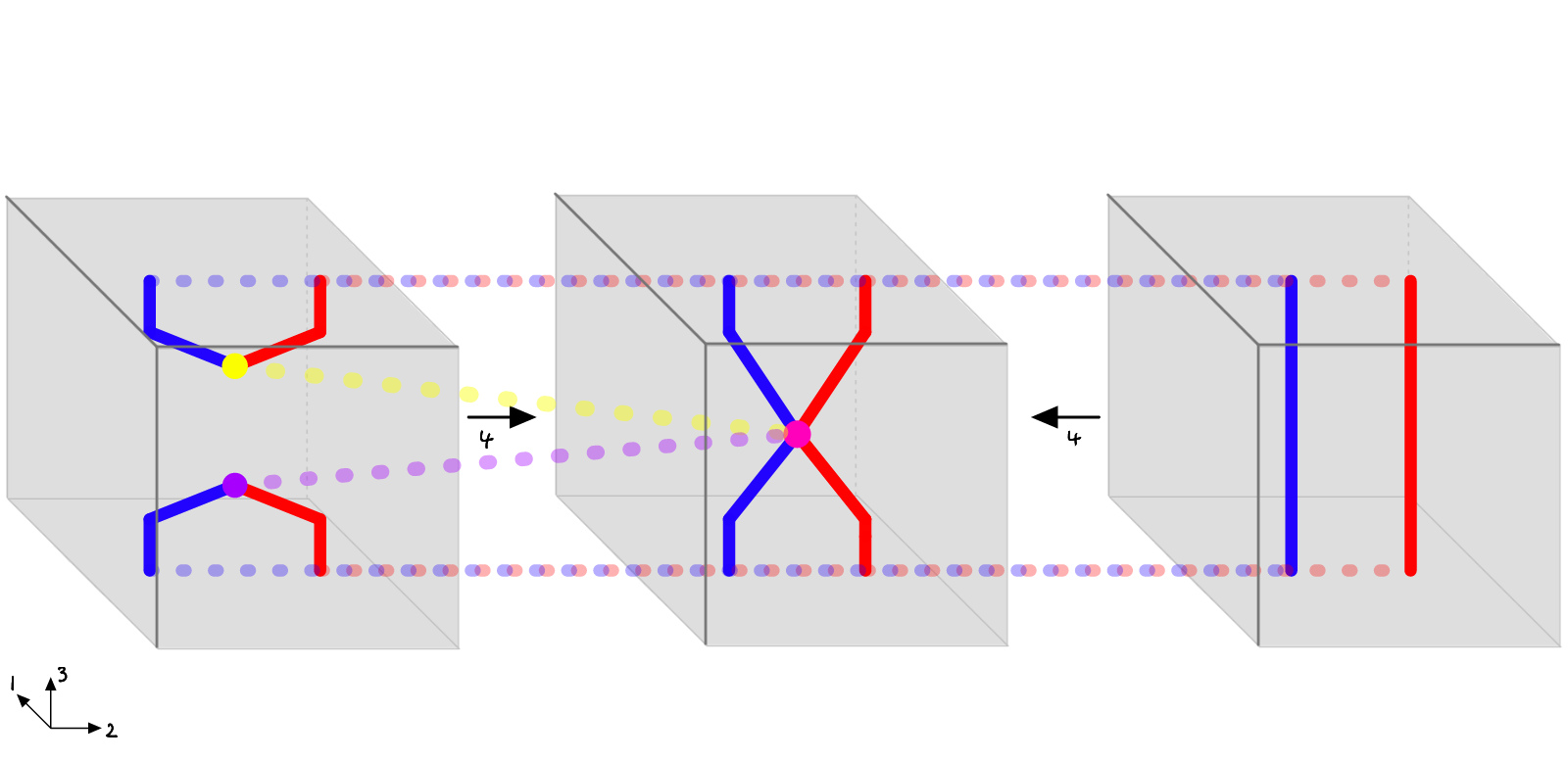}
\endgroup\end{restoretext}
Geometrically, along the left arrow the two arcs converge until they touch at their higher (purple) point and lowest (yellow) point. Along the right arrow, two identities converge closer and closer together until they touch. This type is also known as a ``crotch singularity" (and as the snake singularity it can already live in dimension 3). Formally, this is given by the following data
\begin{restoretext}
\begin{noverticalspace}
\begingroup\sbox0{\includegraphics{ANCimg/page1.png}}\includegraphics[clip,trim=0 {.0\ht0} 0 {.0\ht0} ,width=.9\textwidth]{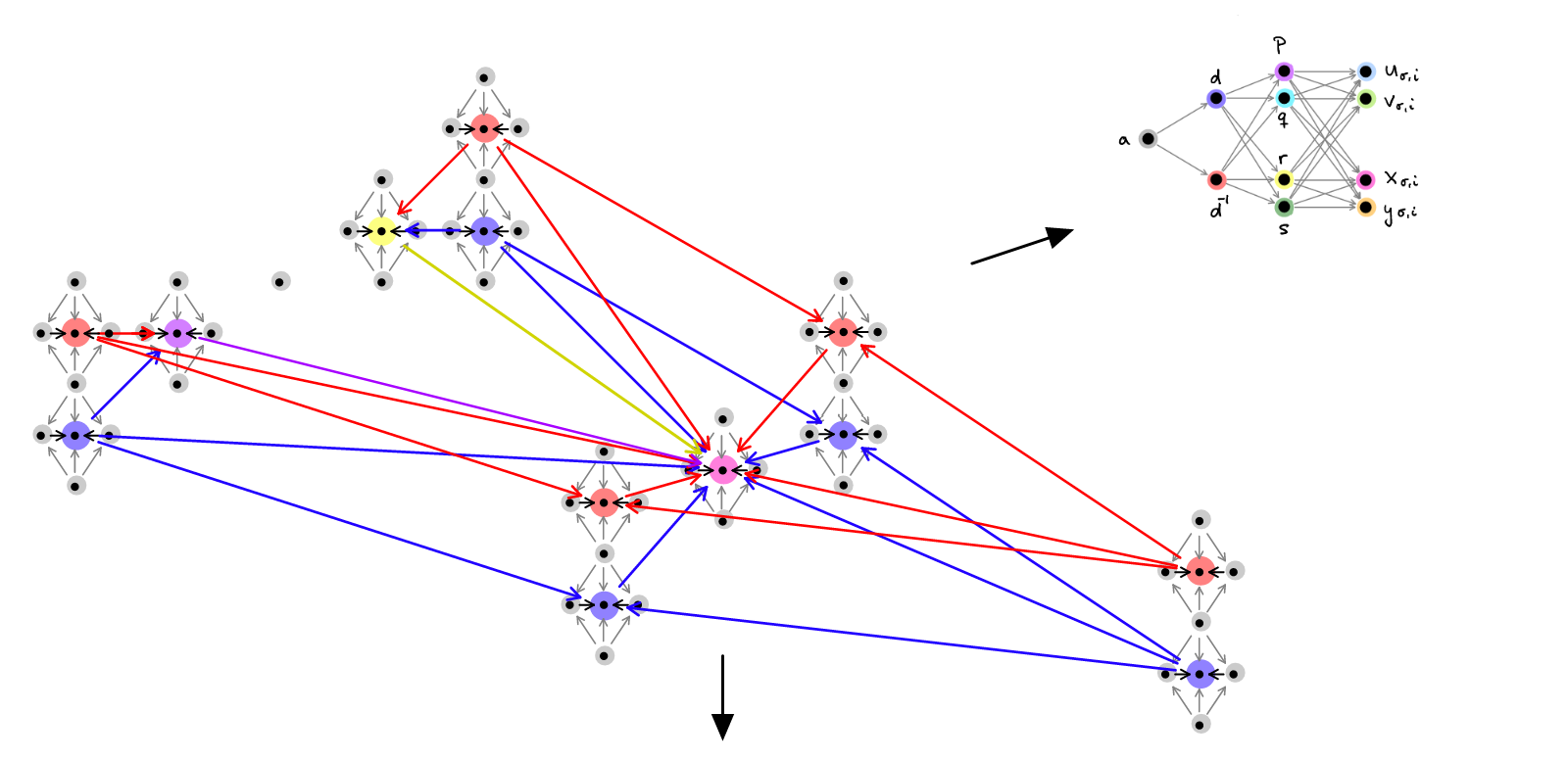}
\endgroup \\*
\begingroup\sbox0{\includegraphics{ANCimg/page1.png}}\includegraphics[clip,trim=0 {.15\ht0} 0 {.0\ht0} ,width=.9\textwidth]{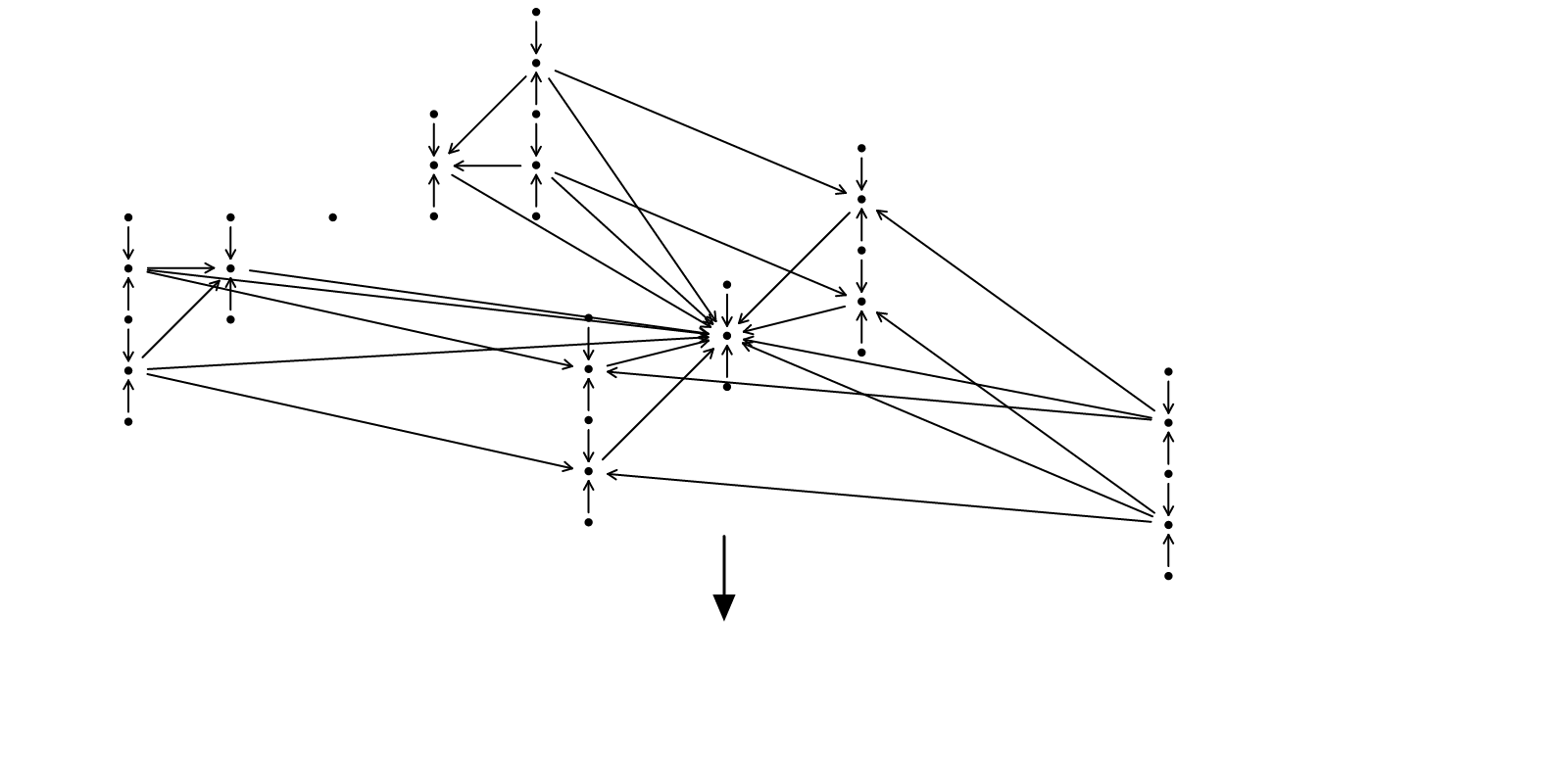}
\endgroup \\*
\begingroup\sbox0{\includegraphics{ANCimg/page1.png}}\includegraphics[clip,trim=0 {.0\ht0} 0 {.0\ht0} ,width=.9\textwidth]{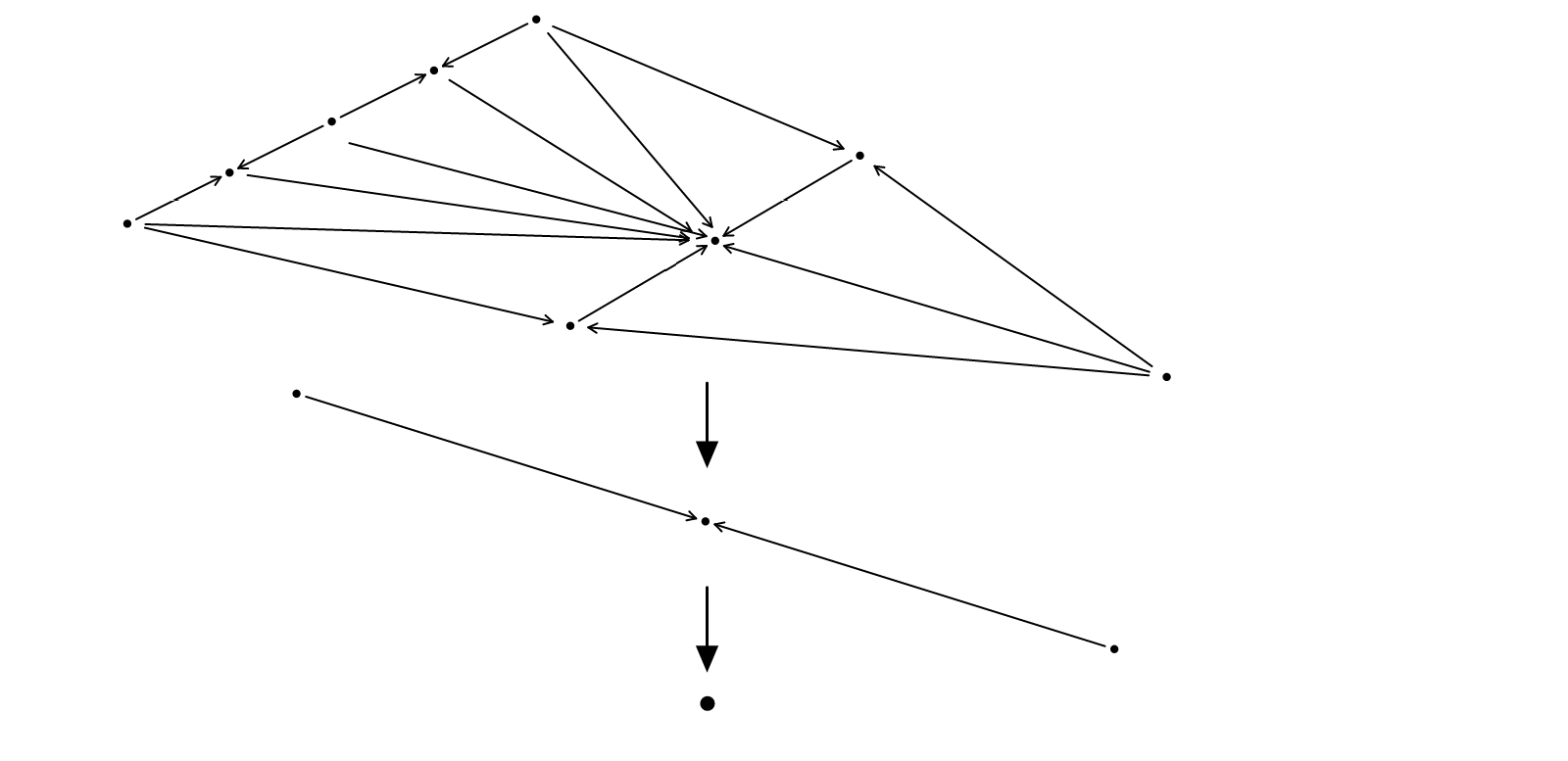}
\endgroup
\end{noverticalspace}
\end{restoretext}
As before $\abss{x_{-,1}}$ is obtained from $\abss{x_{+,1}}$ by exchanging \cred{} and \cblue{} colors, and $\abss{x_{\sigma,-1}}$ is obtained as the inverse of $\abss{x_{\sigma,1}}$ (the latter two are called ``saddle singularities").

Finally, the type $\abss{y_{+,1}}$ is given by
\begin{restoretext}
\begingroup\sbox0{\includegraphics{ANCimg/page1.png}}\includegraphics[clip,trim=0 {.0\ht0} 0 {.2\ht0} ,width=\textwidth]{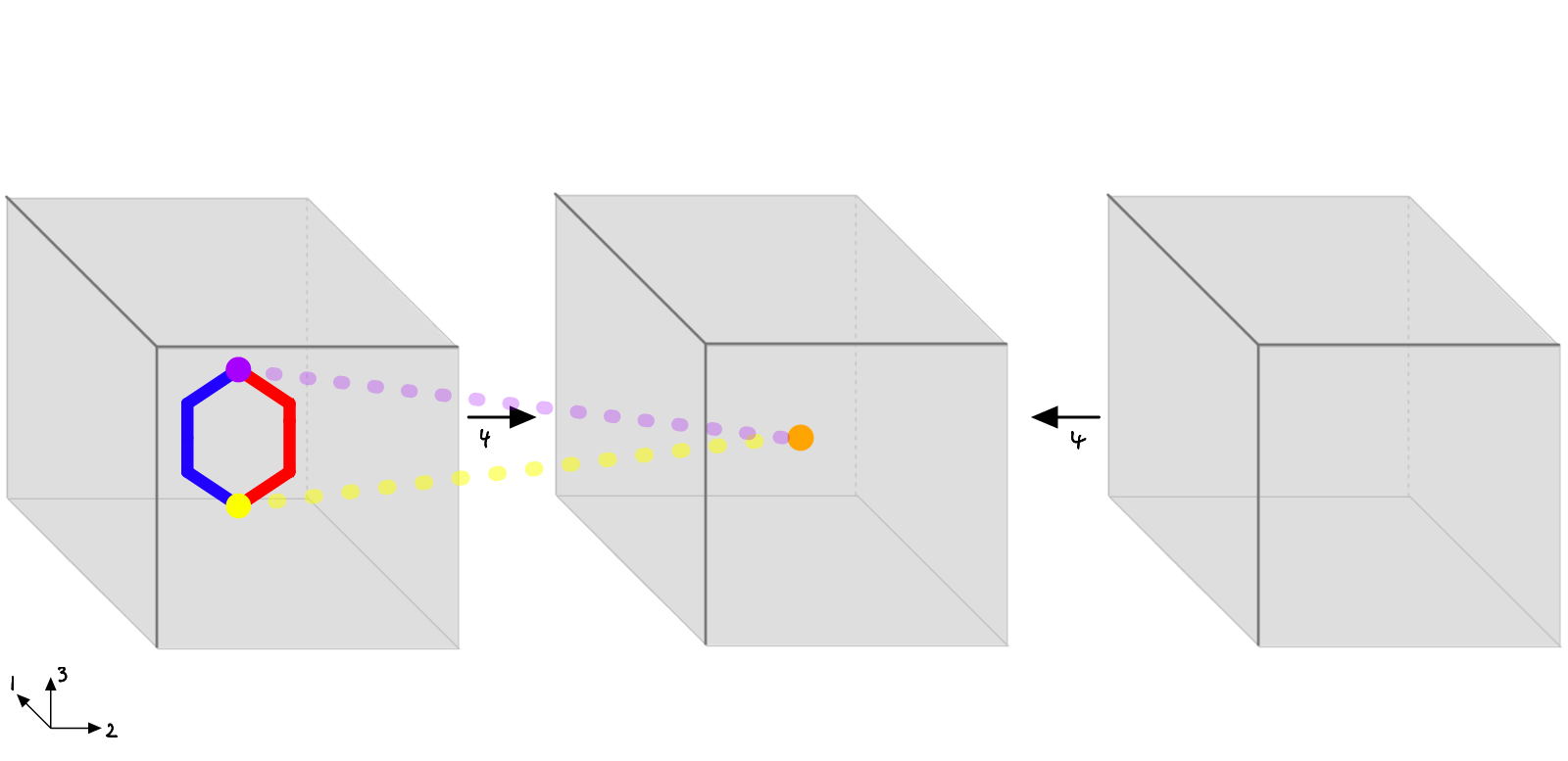}
\endgroup\end{restoretext}
and $\abss{y_{-,1}}$ is obtained from $\abss{y_{+,1}}$ by exchanging \cred{} and \cblue{} colors, and $\abss{y_{\sigma,-1}}$ is obtained as the inverse of $\abss{y_{\sigma,1}}$. These singularities are called ``birth of circle" and ``death of circle" singularities.

To better understand the presentable associative $3$-category defined by the above assignments, we will now discuss the set (cf. \autoref{constr:globular_quotients})
\begin{equation}
\sC(\Id^2_{\abss{a}},\Id^2_{\abss{a}}) \subset (\Comp(\sC)^{\quotg}_{\truncleq 3})_3
\end{equation}
consisting only of those morphisms whose source and target equals $\Id^2_{\abss{a}}$. 

We first discuss the following ``chain of elements in $\Comp(\sC)_4$" which will induce a chain of equalities on elements in $\Comp(\sC)_3$. We start with
\begin{restoretext}
\begingroup\sbox0{\includegraphics{ANCimg/page1.png}}\includegraphics[clip,trim=0 {.0\ht0} 0 {.2\ht0} ,width=\textwidth]{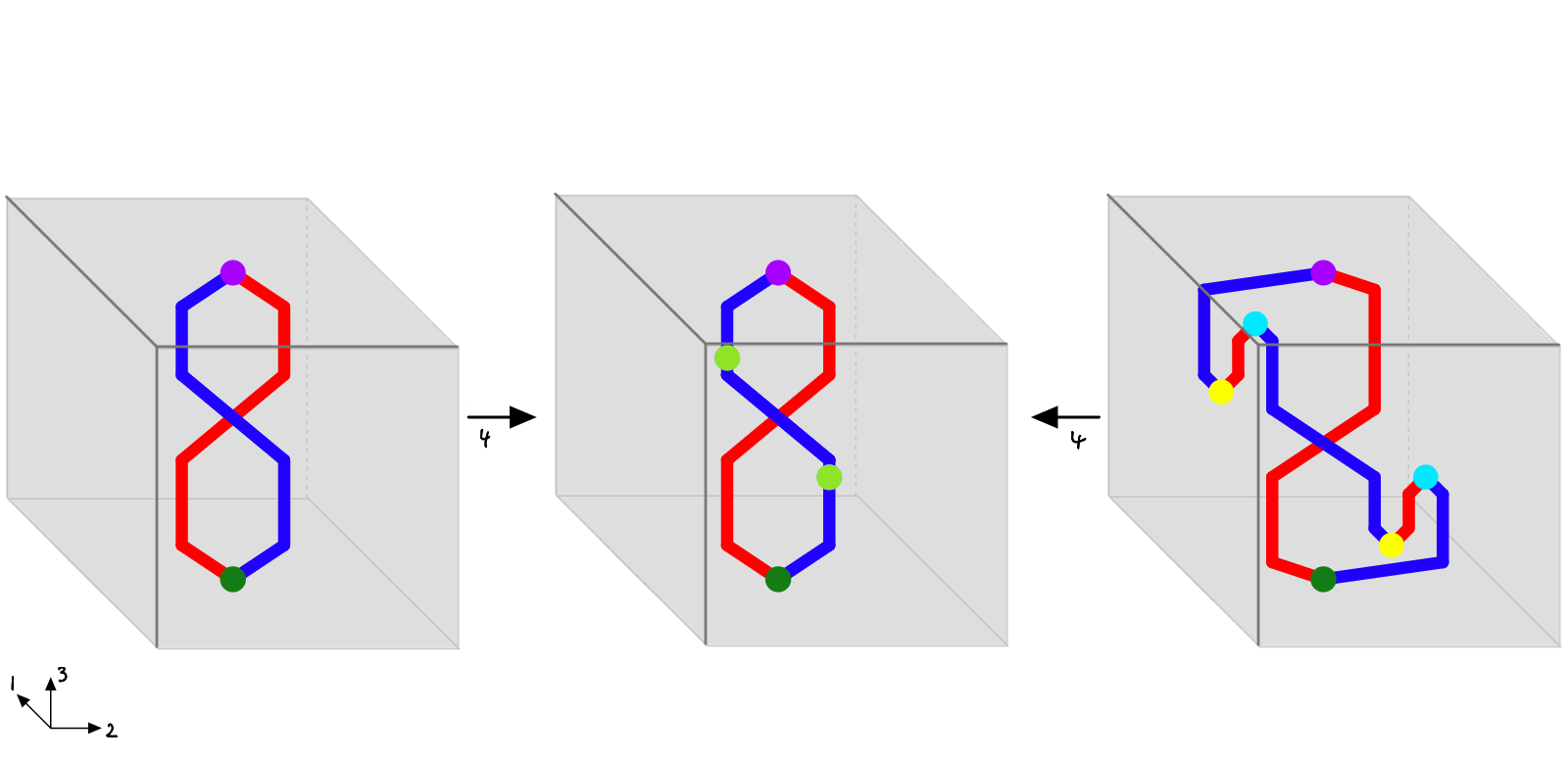}
\endgroup\end{restoretext}
The closed shape on the left we will called the ``closed over-braid". Going from left to right we deform the closed over-braid by applying two snake equations to it. Next consider
\begin{restoretext}
\begingroup\sbox0{\includegraphics{ANCimg/page1.png}}\includegraphics[clip,trim=0 {.0\ht0} 0 {.2\ht0} ,width=\textwidth]{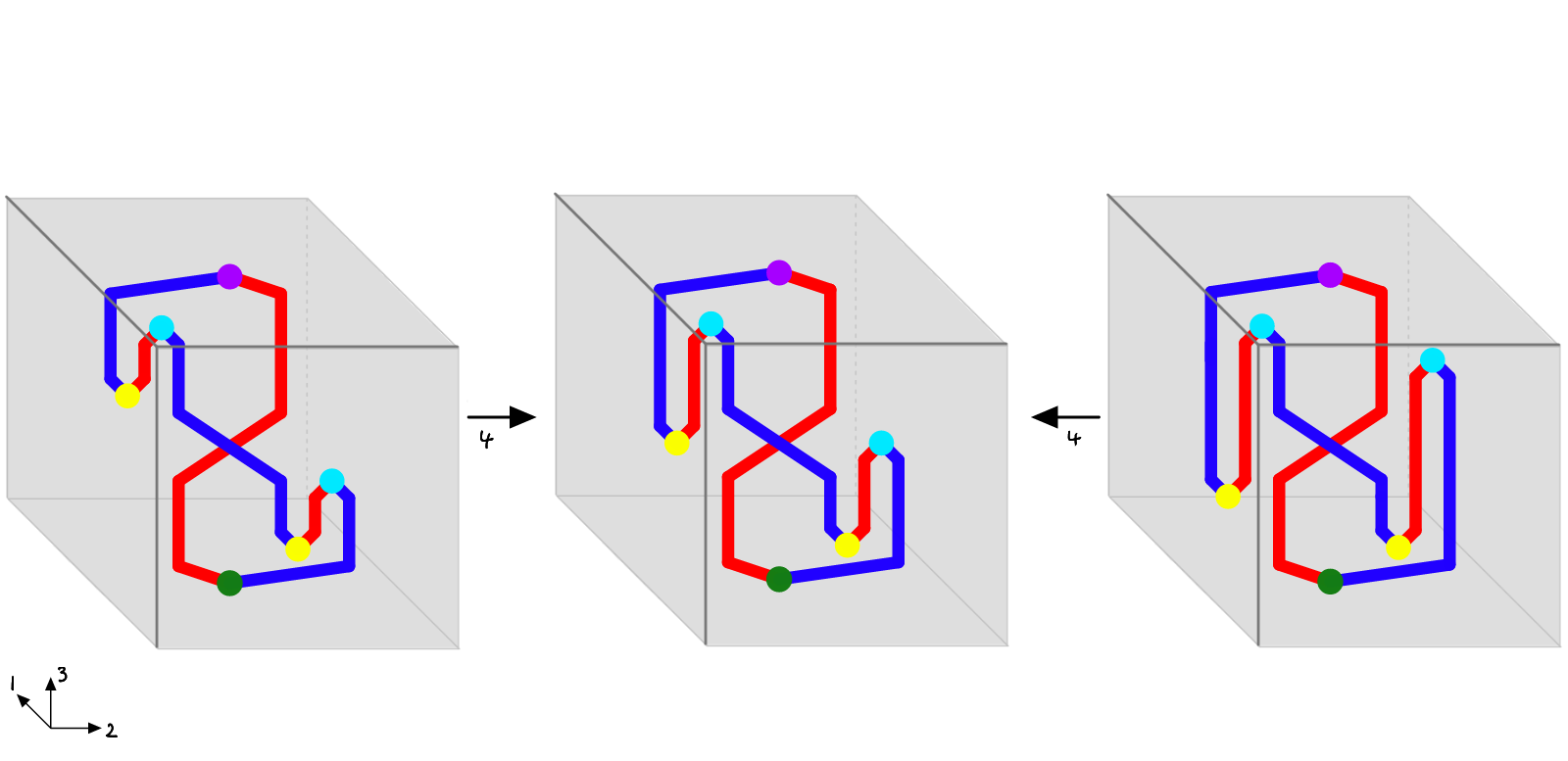}
\endgroup\end{restoretext}
\begin{restoretext}
\begingroup\sbox0{\includegraphics{ANCimg/page1.png}}\includegraphics[clip,trim=0 {.0\ht0} 0 {.2\ht0} ,width=\textwidth]{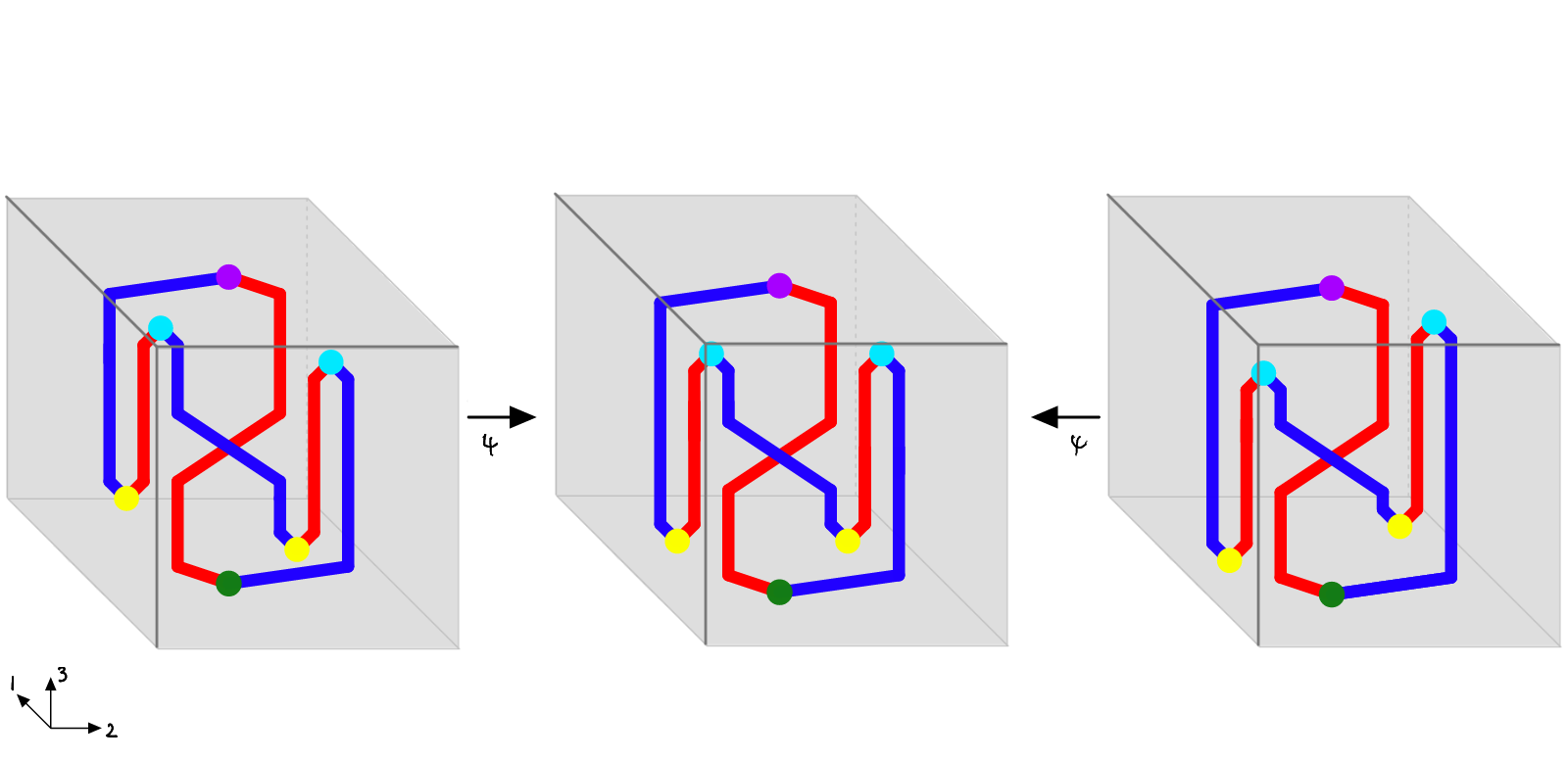}
\endgroup\end{restoretext}
\begin{restoretext}
\begingroup\sbox0{\includegraphics{ANCimg/page1.png}}\includegraphics[clip,trim=0 {.0\ht0} 0 {.2\ht0} ,width=\textwidth]{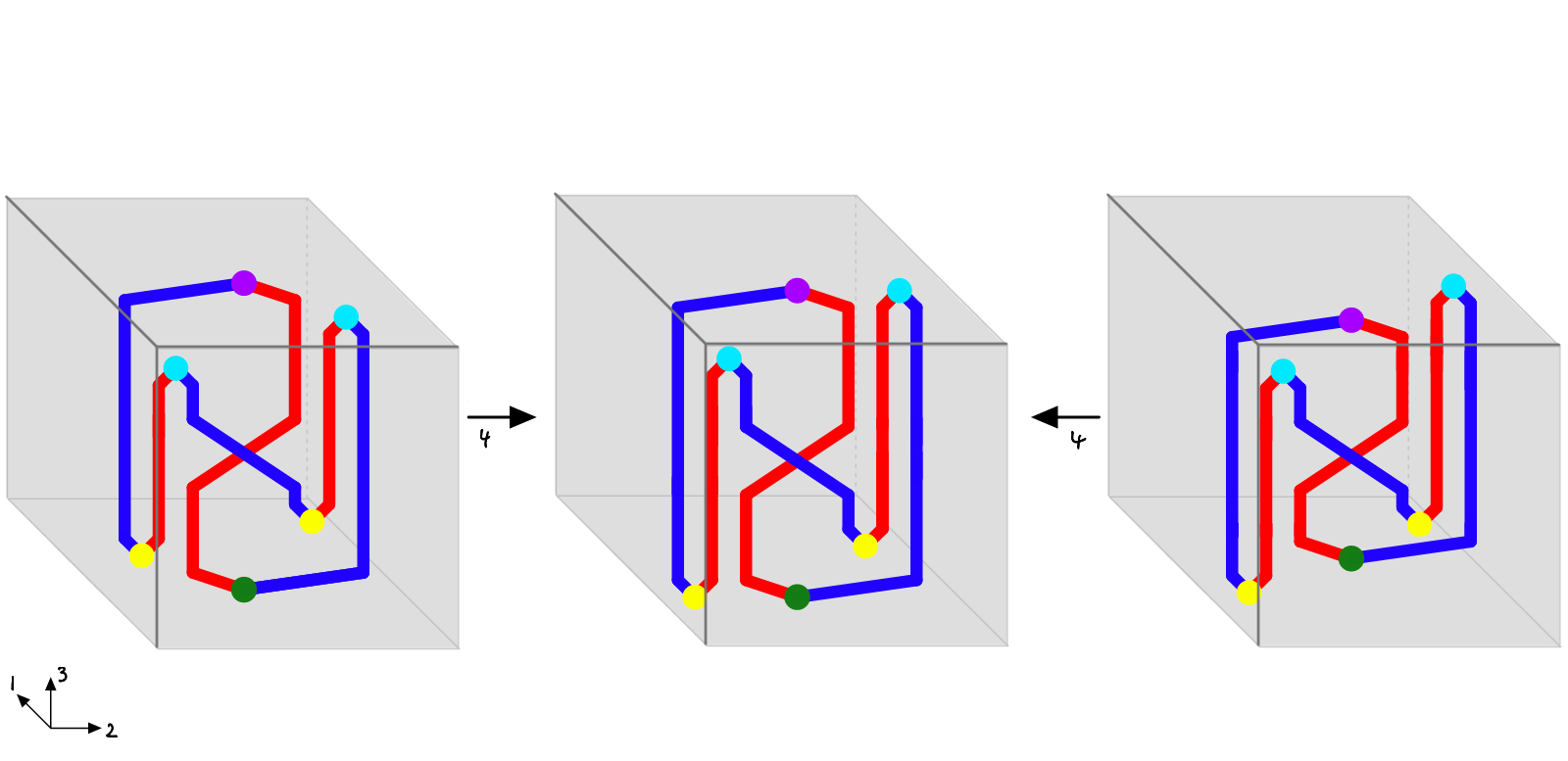}
\endgroup\end{restoretext}
In each of the above we are ``interchange" the heights of $3$-morphisms using (simultaneous) interchangers. Now we can apply (as previously discussed) two snakes in the middle of the last $3$-cube as follows
\begin{restoretext}
\begingroup\sbox0{\includegraphics{ANCimg/page1.png}}\includegraphics[clip,trim=0 {.0\ht0} 0 {.2\ht0} ,width=\textwidth]{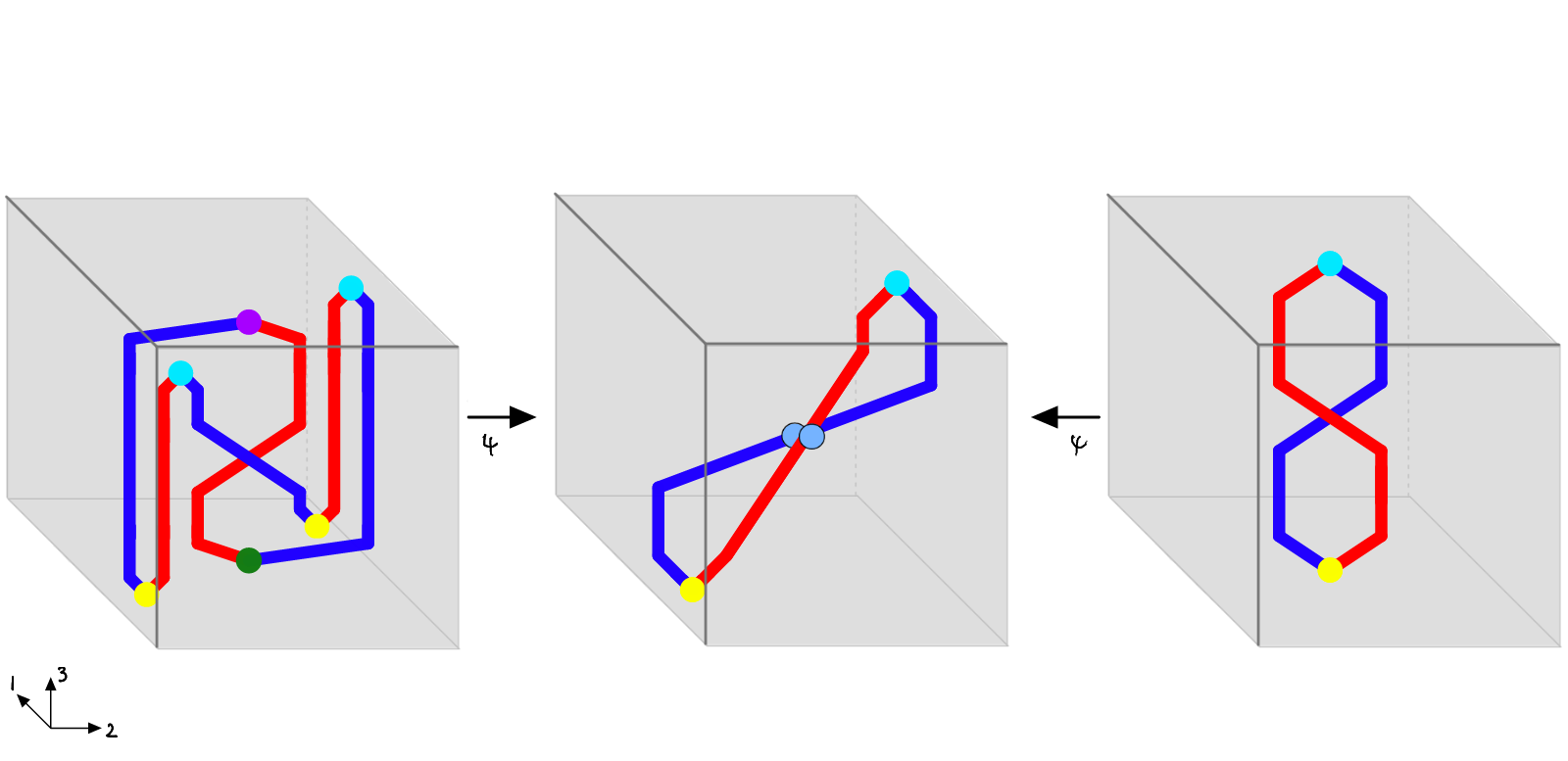}
\endgroup\end{restoretext}
The chain of equations of the above elements of $\Comp(\sC)_4$ (which can also be glued together yielding a single cube) allows as to derive that
\begin{restoretext}
\begingroup\sbox0{\includegraphics{ANCimg/page1.png}}\includegraphics[clip,trim=0 {.1\ht0} 0 {.2\ht0} ,width=\textwidth]{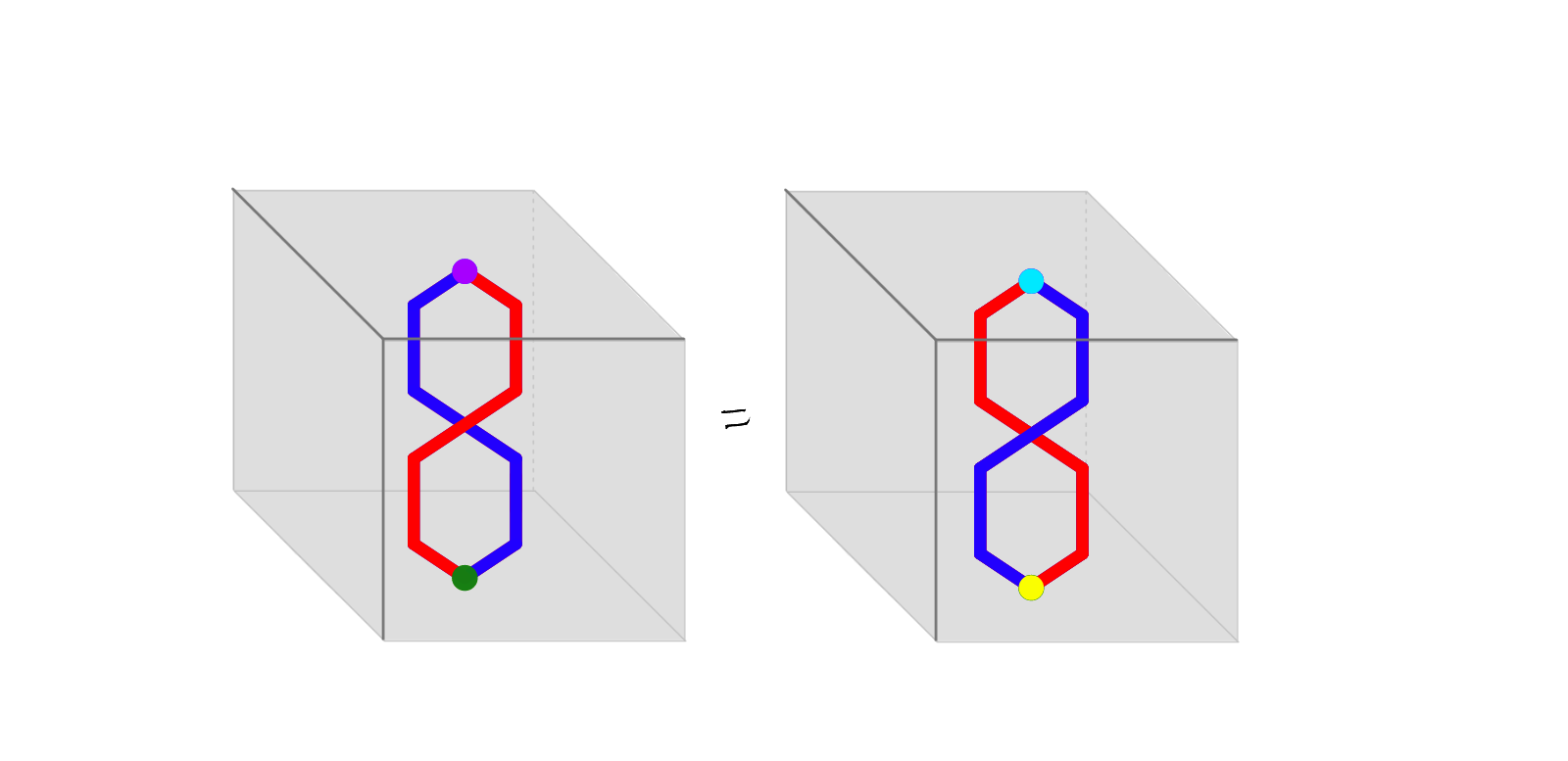}
\endgroup\end{restoretext}
Further using the saddle singularity we find
\begin{restoretext}
\begingroup\sbox0{\includegraphics{ANCimg/page1.png}}\includegraphics[clip,trim=0 {.0\ht0} 0 {.2\ht0} ,width=\textwidth]{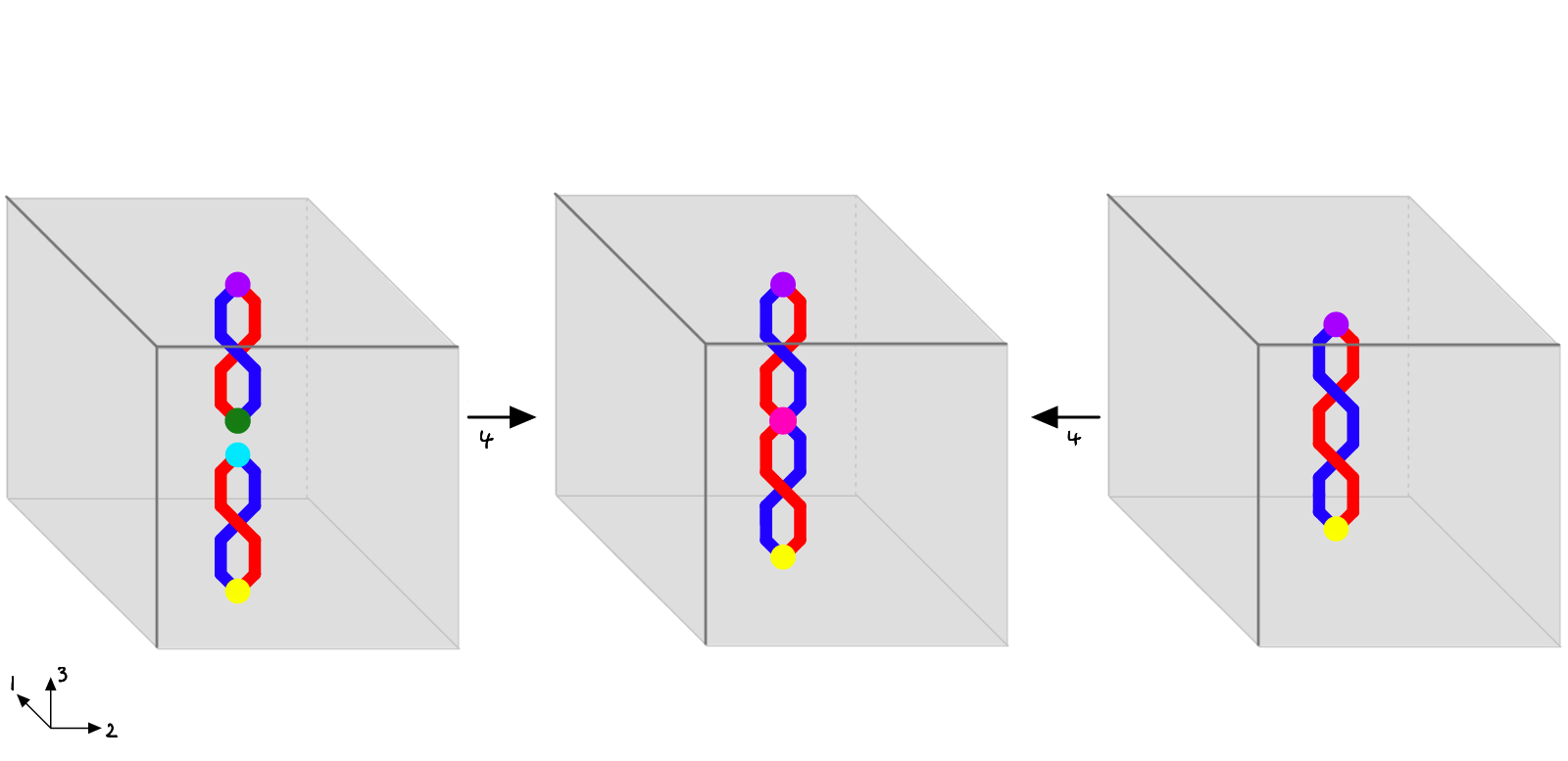}
\endgroup\end{restoretext}
Combining these results we find the equation
\begin{restoretext}
\begingroup\sbox0{\includegraphics{ANCimg/page1.png}}\includegraphics[clip,trim=0 {.1\ht0} 0 {.2\ht0} ,width=\textwidth]{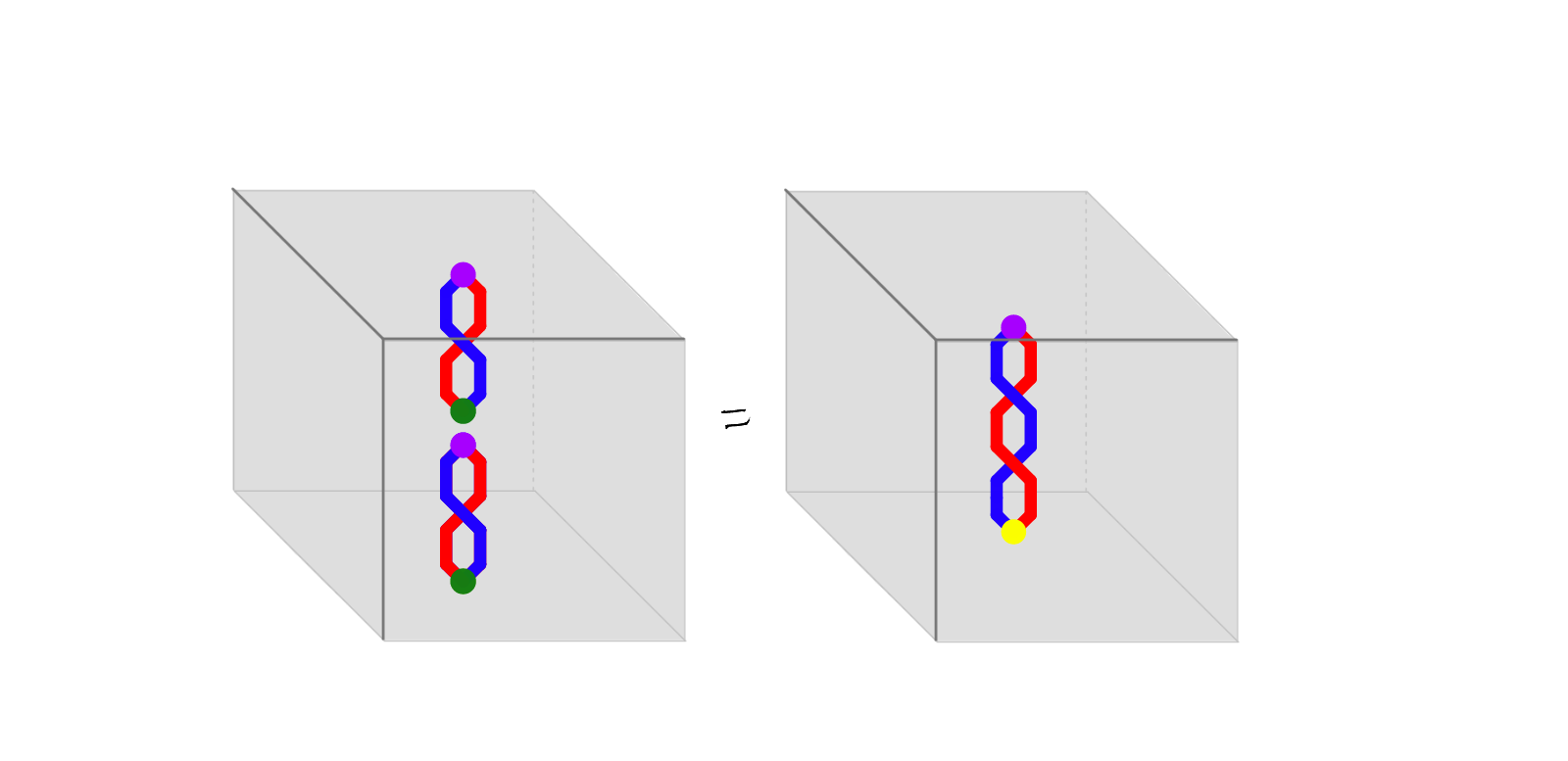}
\endgroup\end{restoretext}
A fully analogous argument allows us to derive that
\begin{restoretext}
\begingroup\sbox0{\includegraphics{ANCimg/page1.png}}\includegraphics[clip,trim=0 {.1\ht0} 0 {.2\ht0} ,width=\textwidth]{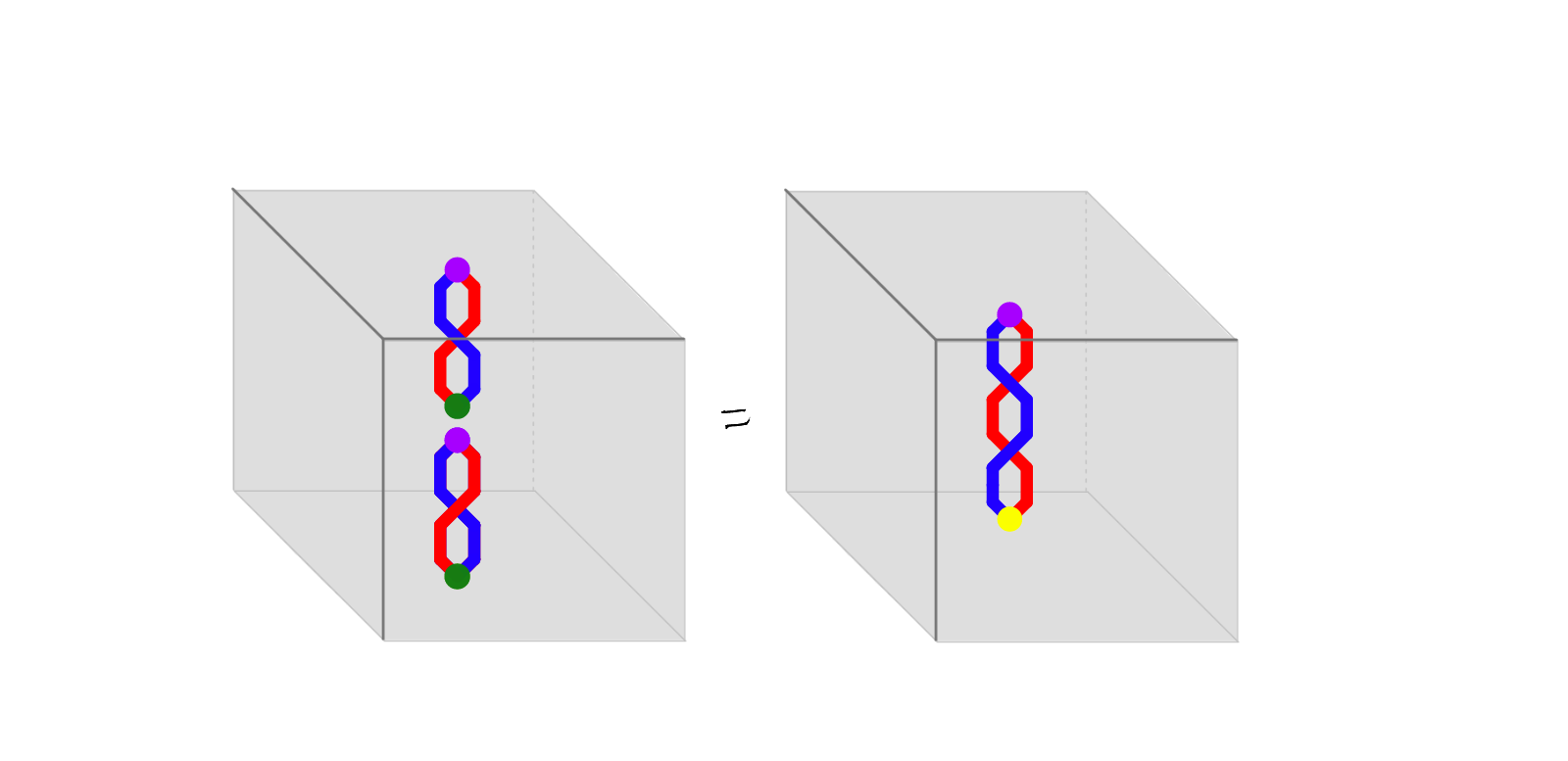}
\endgroup\end{restoretext}
Note that here in the left singular $3$-cube the closed braid on the bottom is an under-braid and not an over-braid. Now we also have the singular $4$-cube
\begin{restoretext}
\begingroup\sbox0{\includegraphics{ANCimg/page1.png}}\includegraphics[clip,trim=0 {.0\ht0} 0 {.2\ht0} ,width=\textwidth]{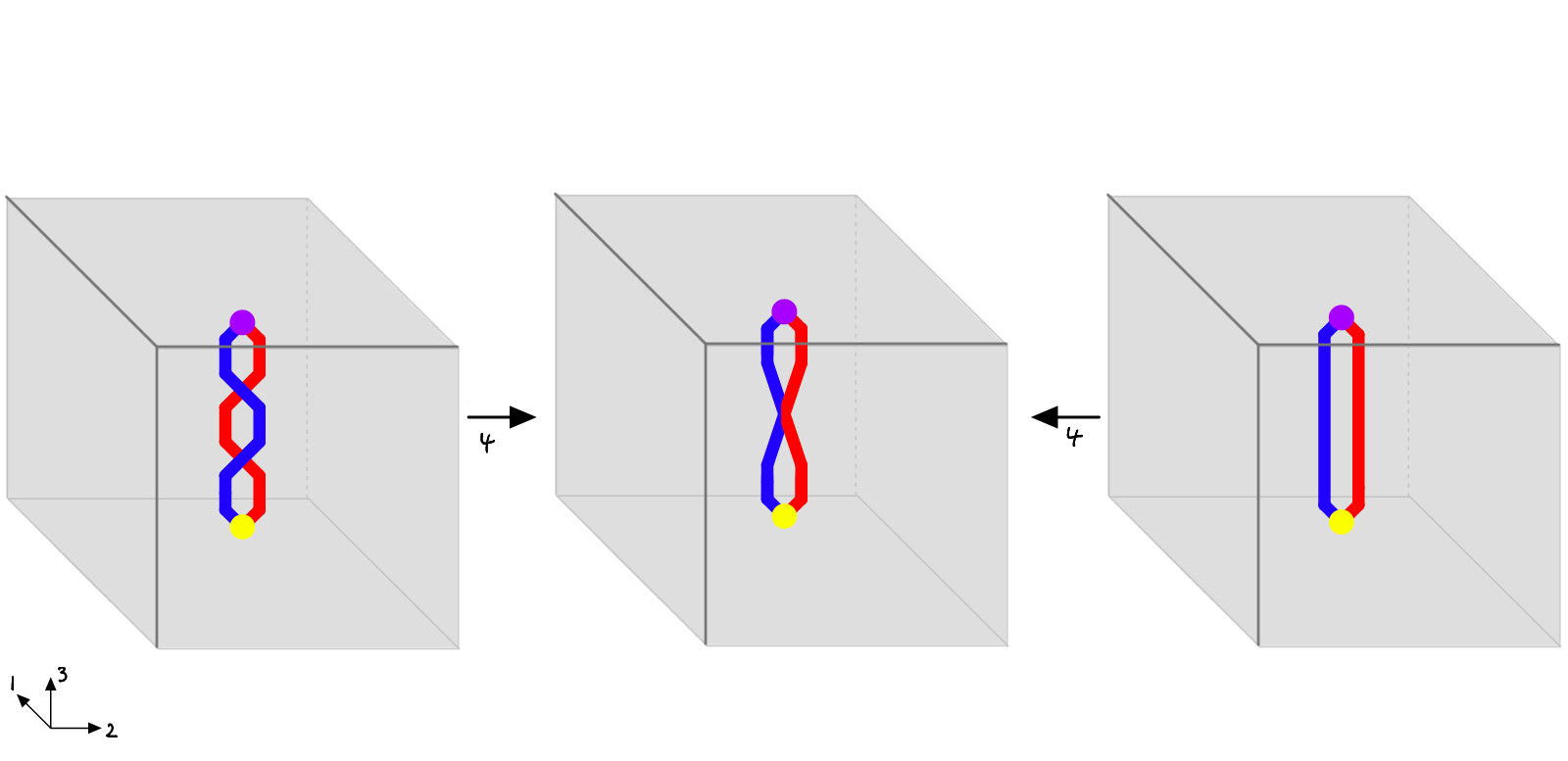}
\endgroup\end{restoretext}
And (using the type of $y_{+,1}$) this means we have the equality
\begin{restoretext}
\begingroup\sbox0{\includegraphics{ANCimg/page1.png}}\includegraphics[clip,trim=0 {.2\ht0} 0 {.15\ht0} ,width=\textwidth]{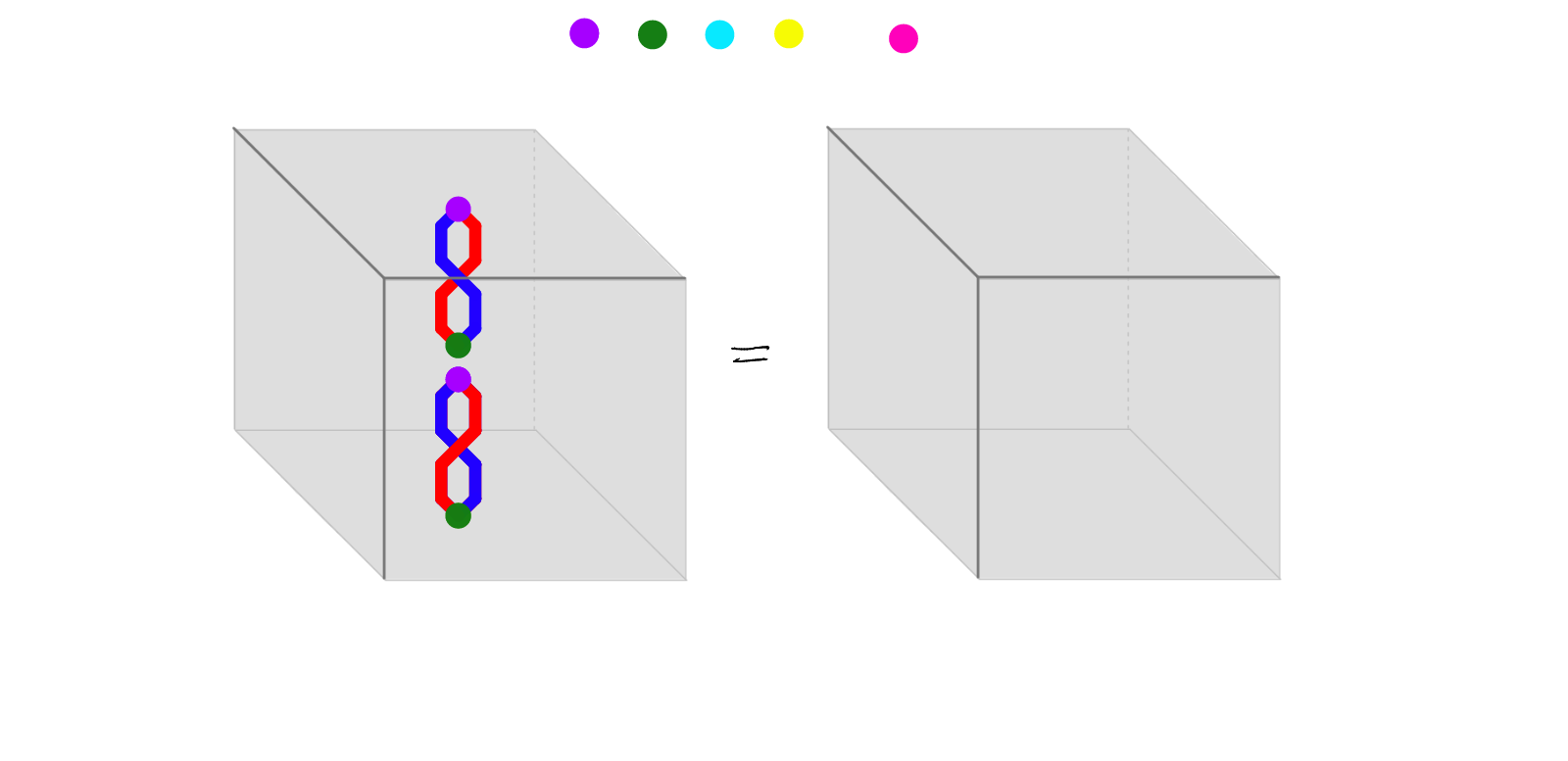}
\endgroup\end{restoretext}
In fact, any element in $\sC(\Id^2_{\abss{a}},\Id^2_{\abss{a}})$ can be seen to lie in the equivalence class of one of
\begin{restoretext}
\begingroup\sbox0{\includegraphics{ANCimg/page1.png}}\includegraphics[clip,trim=0 {.25\ht0} 0 {.25\ht0} ,width=\textwidth]{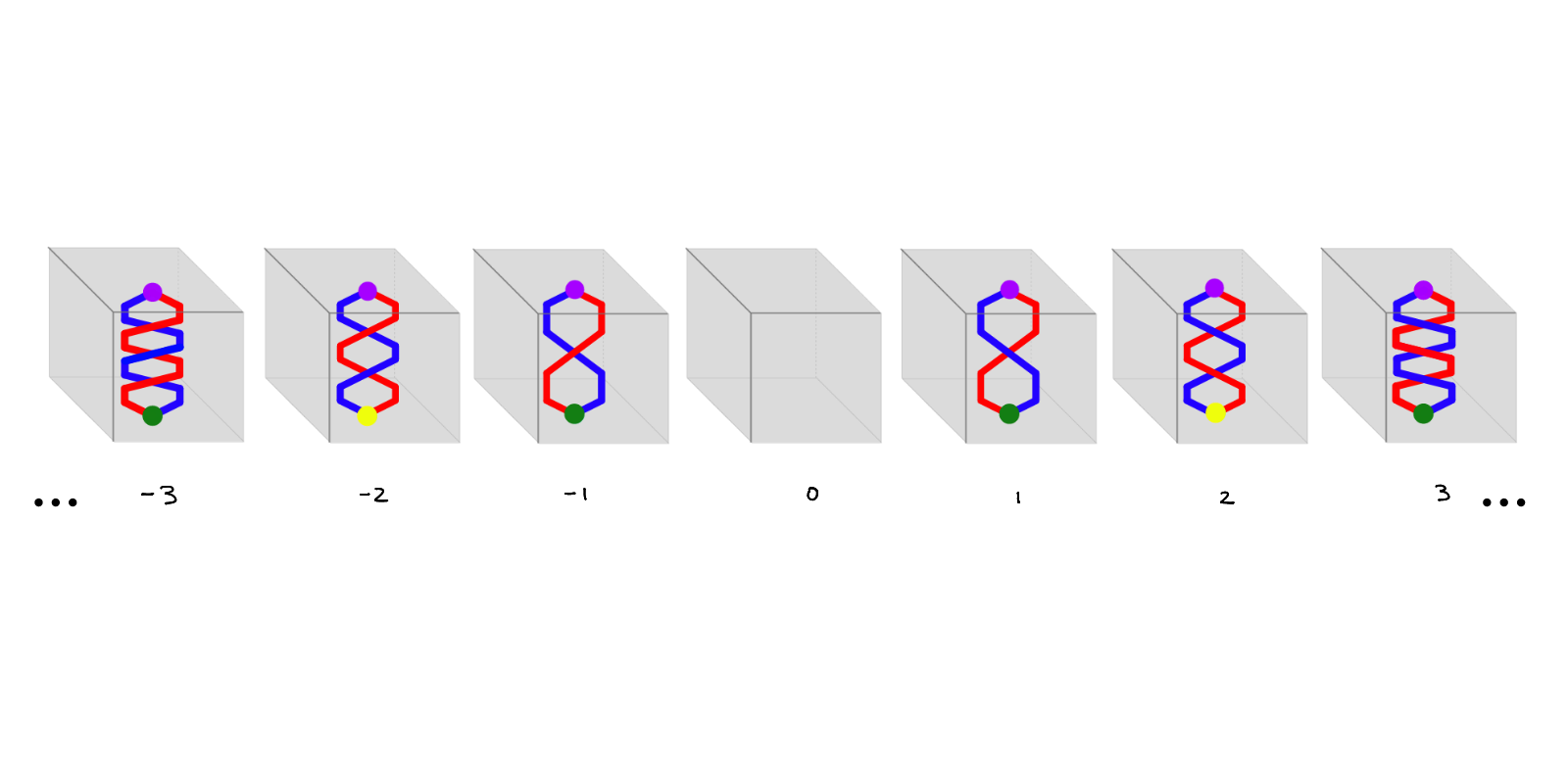}
\endgroup\end{restoretext}
Furthermore, this set obtains additive structure given by gluing cubes, or more formally by $\whisker 1 3$ composition (which means glueing along $2$-boundaries of $3$-globes as we will learn in \autoref{ch:composition}). This allows us to identity 
\begin{equation}
\sC(\Id^2_{\abss{a}},\Id^2_{\abss{a}}) \iso \lZ
\end{equation}
The closed over- and under-braid behave like $1$ and $-1$.

As we will understand later on, the chosen $\sC$ was a $3$-truncated representation of the homotopy type of $S^2$, and $\sC(\Id^2_{\abss{a}},\Id^2_{\abss{a}})$ represent homotopy classes of maps $D^3 \to S^2$ such that the boundary $\partial D^3$ lands in the basepoint $a$. In other words, the above equality states that 
\begin{equation}
\lZ \iso \pi_3(S^2)
\end{equation}
An analogous discussion can be given for an $(n+1)$-presentation $\sC$ of $S^n$, $n > 2$, which shifts all $\sC_k$ for $k > 0$ by $n-2$ when compared to our $S^2$ presentation above). It is then easy to see that 
\begin{equation}
\lZ_2 \iso \pi_{n+1}(S^n)
\end{equation}

\chapter{Coherent invertibility and $n$-groupoids} \label{ch:groupoids}

In \autoref{sec:pres_3} we saw a \free{} associative $3$-category which we claimed modelled the homotopy $3$-type of the $2$-sphere $S^2$. In this chapter we will  understand that this $3$-presentation can be interpreted as a ``higher groupoid" and that its generators then fall into one of two categories: They are either ``fundamental generators" (for our example, those living in dimensions 0 and 2) or derived coherences witnessing invertibility of fundamental generators (such as the cup, cap, snake, saddle, crotch, death and birth of the circle in dimension 3 and 4).

The goal of this chapter is to give a candidate definition for the theory of coherent invertibility. There are multiple natural ways of defining this theory (all of which, however, are conjecturally equivalent). We will here present a spectrum of choices ranging from ``manifold coherent" to ``full coherent". We pick ``full coherence" as our default definition (as this is the most general definition, and for instance satisfied by homotopies). We will then define general invertible elements and how to adjoin them to existing presented associative $n$-categories. This will lay the foundation for defining presented associative $n$-groupoids in \autoref{sec:group_def}. 

Later on, in \autoref{ch:geom}, we will explore further the theory of ``manifold coherent" invertibility: the relevance of this theory derives from the fact, that manifold-like invertibility witnesses look like (framed) $k$-tangles, and appear naturally in the context of the generalised Thom-Pontryagin construction, suggesting a bridge from CW-complexes to presented associative (manifold coherent) $\infty$-groupoids.

\section{The $r$-connected theory of invertibility} \label{sec:group_TI}

\subsection{Definition}

In order to understand general invertible $k$-morphisms, we will need to understand the \textit{theory of invertibility} which is an associative higher category that should be thought of as the ``prototypical presented associative $n$-category on a single invertible $1$-morphism". As it turns out, we are confronted with a spectrum of choices for such a theory.

We will suggest to understand (points of) this spectrum using a characterisation by the following predicates $P^k_r(X)$ (for $-1 \leq r \leq \infty + 1$, $k \in \lN$, and $X$ a space)
\begin{equation}
P^k_r(X) := \begin{cases} \text{$X$ is $r$-connected} & -1 \leq r \leq \infty \\
\text{$X$ is homeomorphic to a $k$-ball} & r = \infty + 1
\end{cases} 
\end{equation}
We remark that our conditions $P^k_\infty(X)$ and $P^k_{\infty+1}(X)$ are possibly equivalent subject to the additional assumptions that $X$ is a compact $k$-dimensional manifold embedded in $\lR^k$ (which will hold in all cases considered below). However, as a precautionary measure we list these two predicates separately.

\begin{constr}[Manifold coherent theory of invertibility] \label{constr:TI}
The \textit{$r$-connected theory of $n$-invertibility} $\TI^{n,r}$ is a \free{} associative $n$-category. The construction of $\TI^{n,r}$ is by induction. $\TI^{n,r}_0$ contains two objects
\begin{equation}
\TI^{n,r}_{0} = \Set{-,+}
\end{equation}
Now let $0 \leq k \leq n$ and inductively assume $\TI^{n,r}_m$ (and thus $\Comp(\TI^{n,r})_m$ by \autoref{rmk:inductive_PANC_def}) have been defined for $m \leq k$. $\TI^{n,r}_{k+1}$ is then defined as follows. Let $S,T \in \Comp(\TI^{n,r})_{k}$, $S \neq T$, such that globular sources and targets of $S,T$ coincide. Using \autoref{constr:double_cones_of_src_and_tgt}, define $\scD := \abss{S \to T}$ as a double cone with center $p$ (i.e. $\sU^{k+1}_\scD(p) =\top$).

We now formulate the following condition on $\scD$: 
\begin{itemize}
\item Let $\norm{\scD}$ be the manifold diagram associated to $\scD$ and $\partial \norm{\scD}$ the labelled stratification of the $(k+1)$-cube's boundary $\partial [0,1]^{k+1}$ (obtained by appropriate closure). Then the regions labelled by $-$ and $+$ in $\partial \norm{\scD}$ satisfy the predicate $P^k_r$.
\end{itemize}
Formally, this condition can be phrased as follows. Let $\abs{\scD}^{k+1}_\bullet$ be the triangulation constructed from $\scD$ as discussed in \autoref{ssec:coloring} and remove from it the center vertex $p$ to obtain $\abs{\partial\scD}^{k+1}_{\bullet}$. For $x \in \Set{-,+}$ define
\begin{equation}
\abs{x} = \bigcup \Set{ f \in \abs{\partial\scD}^{k+1}_{k} ~|~\tsU {k+1}_\scD(f(0))=x}
\end{equation}
Now the above condition formally ask that for both $x$
\begin{equation}
P^k_r(\abs{x}) \text{ is true}
\end{equation}
If this holds for the chosen $S,T$ we define $\TI^{n,r}_k$ to contain an element
\begin{equation}
\ic_{S\equiv T} \in \TI^{n,r}_{k+1}
\end{equation}
called the \textit{coherence between $S$ and $T$}. To complete the definition of $\TI^{n,r}$ we also need to assign types. For this we set 
\begin{equation} \label{eq:cobordism_minimal_def}
\abss{\ic_{S \equiv T}} =  \abss{S \xto {\ic_{S \equiv T}} T}
\end{equation}
where the right hand side is defined using \autoref{constr:double_cones_of_src_and_tgt}. The latter construction (together with the argument in \autoref{constr:adjoining_gen}) implies that this satisfies the conditions for types in \autoref{defn:pres_ANC}.
\end{constr}

\begin{rmk}[Proof irrelevance and symmetry] 
For $k = n$, we note that $\ic_{S \equiv T} \in \TI^{n,r}_{k+1}$ is uniquely identified by its type (i.e. by its source and target), and that $\ic_{S \equiv T} \in \TI^{n,r}_{k+1}$ implies $\ic_{T \equiv S} \in \TI^{n,r}_{k+1}$. The former gives proof-irrelevance, the latter gives symmetry as required in \autoref{defn:pres_ANC} of a \free{} associative $n$-category.
\end{rmk}

\begin{rmk}[Computability of $\TI$] For general $k$, the condition $P^k_r(\abs{x})$ is decidable in general only if $r = -1$, in which case it reduces to checking \textit{non-emptiness} of $(\sU^k_\scD)\inv(x)$. For all other $r$ the condition becomes undecidable. However, it remains semi-decidable (which would be sufficient e.g. for the purposes of a proof assistant).
\end{rmk}

Note that there is a (transfinite) chain of inclusions of presentations
\begin{equation}
\TI^{n,\infty+1} \into \TI^{n,\infty} \into \quad  ... \quad \into \TI^{n,2} \into \TI^{n,1} \into \TI^{n,0} \into \TI^{n,-1} 
\end{equation}
$\TI^{n,\infty+1}$ will also be called the \textit{manifold coherent} theory of invertibility, whereas $\TI^{n,-1}$ is called the \textit{full coherent} theory of invertibility.

By default we will work with the full coherent theory of invertibility, which in several ways is the ``safest" choice. We make the following notational conventions.

\begin{notn} $\TI^n$ denotes $\TI^{n,-1}$ and we usually denote $\TI^\infty$ simply by $\TI$.
\end{notn}

\subsection{Examples} \label{ssec:TI_eg}

We draw examples of elements in $\TI^{n,r}_k$ (for $k \leq n$).

\begin{notn} For $\cC$-labelled singular $1$-cubes $A, B$ such that the labelling of the last regular segment of $A$ coincides with the labelling of the first regular segment of $B$, denote by $A \whisker 1 1 B$ the $1$-cube obtained by ``gluing" $A$ and $B$ together along these segments. This notation will be generalised and and formalised in the next chapter.
\end{notn}

 In dimension $1$, $\TI^{n,k}_1$ for each $k$ has exactly two elements, namely $\ic_{- \equiv +}$ and $\ic_{+\equiv -}$, whose types are given by
\begin{restoretext}
\begingroup\sbox0{\includegraphics{ANCimg/page1.png}}\includegraphics[clip,trim=0 {.35\ht0} 0 {.25\ht0} ,width=\textwidth]{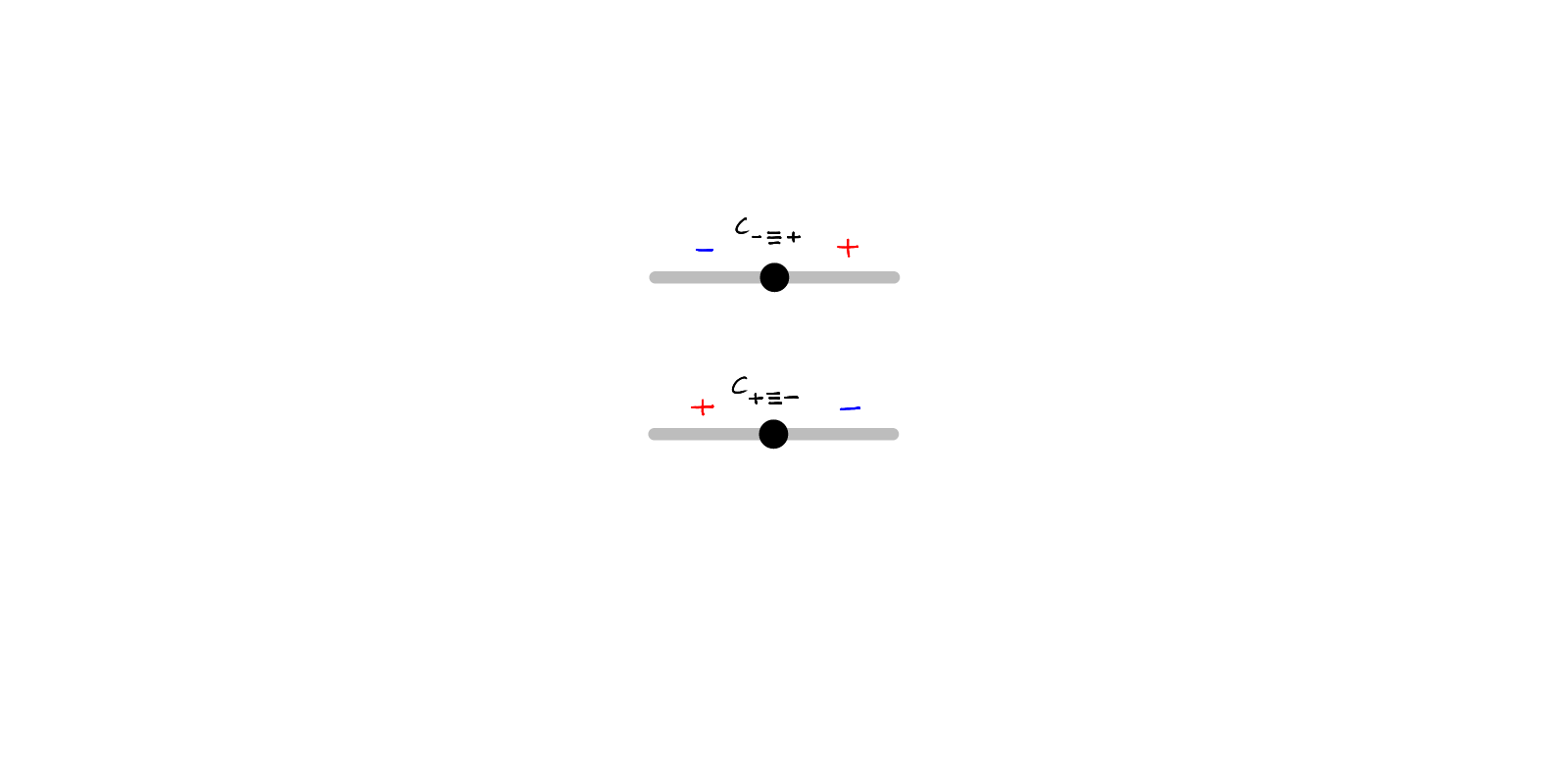}
\endgroup\end{restoretext}
Note that we cannot have generating $1$-morphisms whose source and target are both $+$ (or both $-$) as this would contradict the non-emptiness condition in \autoref{constr:TI}. General $1$-morphisms in $\Comp(\TI)_1$ are then composites of these two generating $1$-morphisms.

In dimension $2$, however, differences between theories start to appear: $\TI^{n,-1}_2$ has elements such as (also see notation above)
\begin{equation}
\ic_{\abss{\ic_{-\equiv +}} \whisker 1 1 \abss{\ic_{+ \equiv -} \whisker 1 1 \abss{\ic_{-\equiv +}}} \equiv \abss{\ic_{-\equiv +}}} \\
\end{equation}
which has the following corresponding manifold diagram
\begin{restoretext}
\begingroup\sbox0{\includegraphics{ANCimg/page1.png}}\includegraphics[clip,trim=0 {.3\ht0} 0 {.3\ht0} ,width=\textwidth]{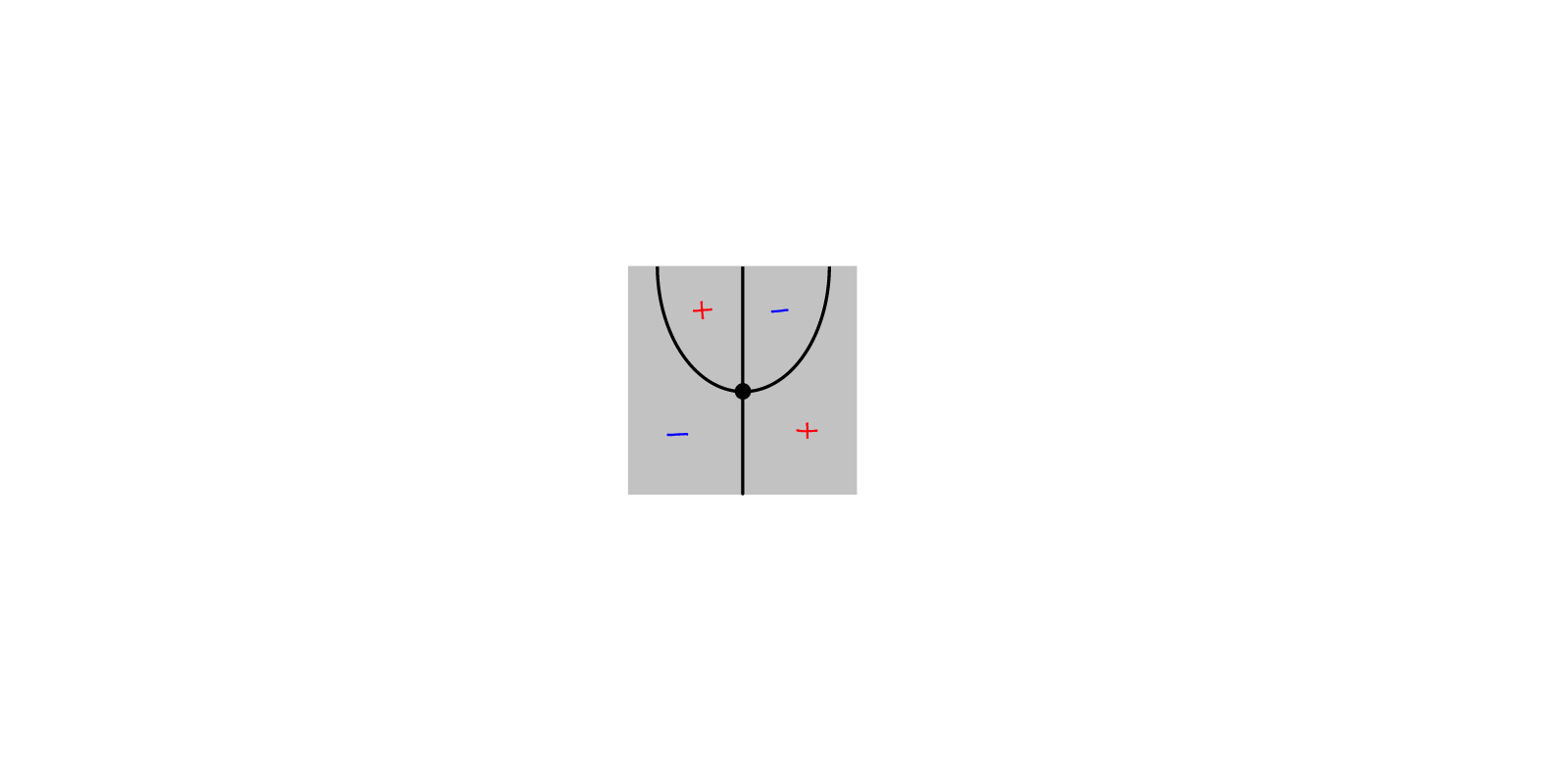}
\endgroup\end{restoretext}
But this is not an element of $\TI^{n,k}_2$ for $k \geq 0$ since boundary regions labelled by $-$ or $+$ are not $0$-connected (i.e. connected) but they are $(-1)$-connected (i.e. non-empty). In fact, for $k \geq 0$, $\TI^{n,k}_2$ only contain exactly four elements that look like $1$-dimensional manifolds, namely
\begin{align}
\ic_{\abss{\ic_{-\equiv +}} \whisker 1 1 \abss{\ic_{+ \equiv -}} \equiv \Id_{-}} \\
\ic_{\abss{\ic_{+\equiv -}} \whisker 1 1 \abss{\ic_{- \equiv +}} \equiv \Id_{+}} \\
\ic_{\Id_{+}\equiv\abss{\ic_{+\equiv -}} \whisker 1 1 \abss{\ic_{- \equiv +}}} \\
\ic_{\Id_{-} \equiv \abss{\ic_{-\equiv +}} \whisker 1 1 \abss{\ic_{+ \equiv -}}}
\end{align}
whose types are of the form
\begin{restoretext}
\begingroup\sbox0{\includegraphics{ANCimg/page1.png}}\includegraphics[clip,trim=0 {.0\ht0} 0 {.1\ht0} ,width=\textwidth]{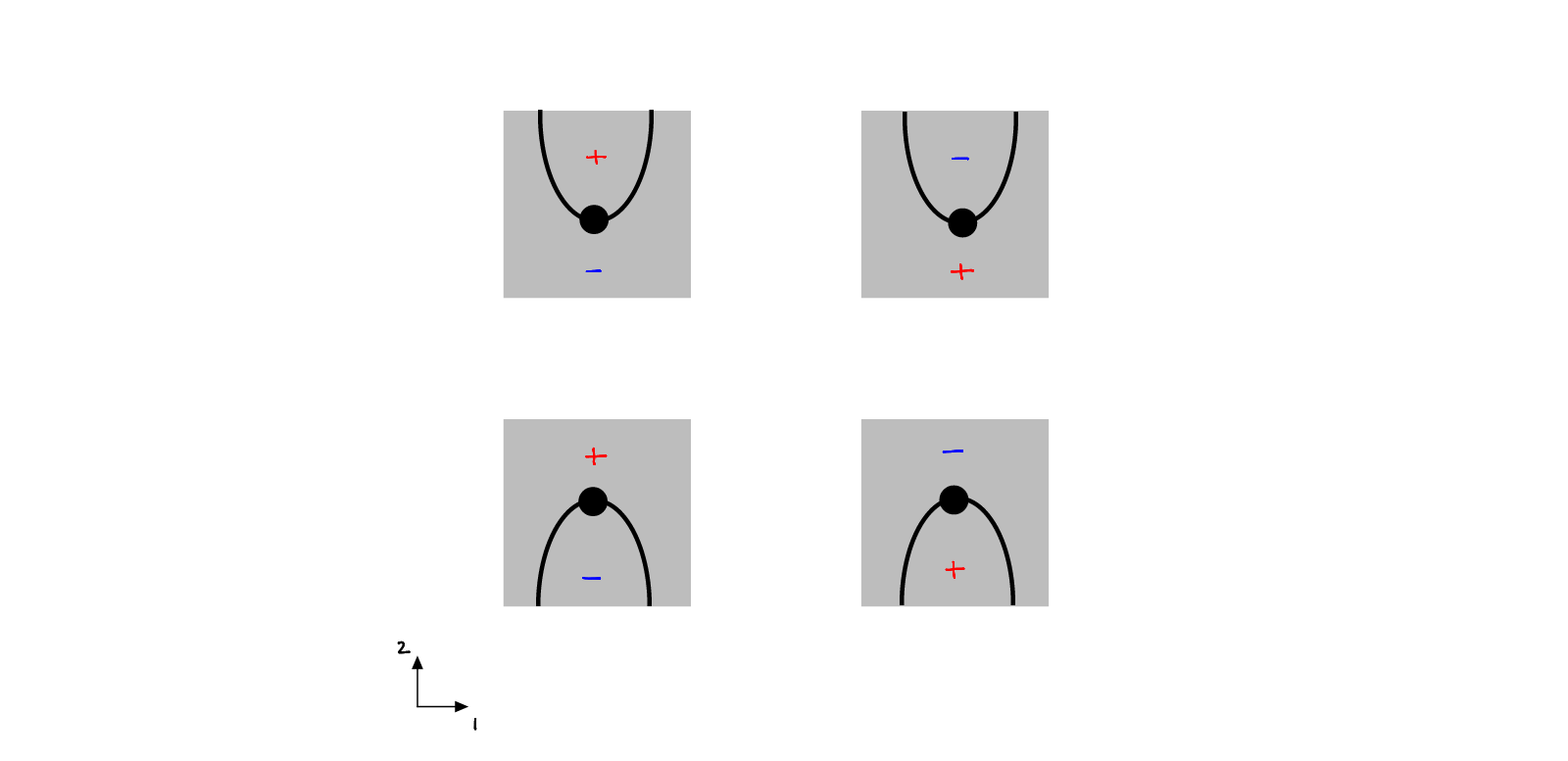}
\endgroup\end{restoretext}
Types in the upper row are called cap singularities (or ``caps" for short) and types in the lower row are called cup singularities (or ``cups" for short).

In higher dimensions $k$, we see more differences between $\TI^{n,k}_k$ appearing. For instance, for $k = 4$, the following generator
\begin{restoretext}
\begingroup\sbox0{\includegraphics{ANCimg/page1.png}}\includegraphics[clip,trim=0 {.0\ht0} 0 {.1\ht0} ,width=\textwidth]{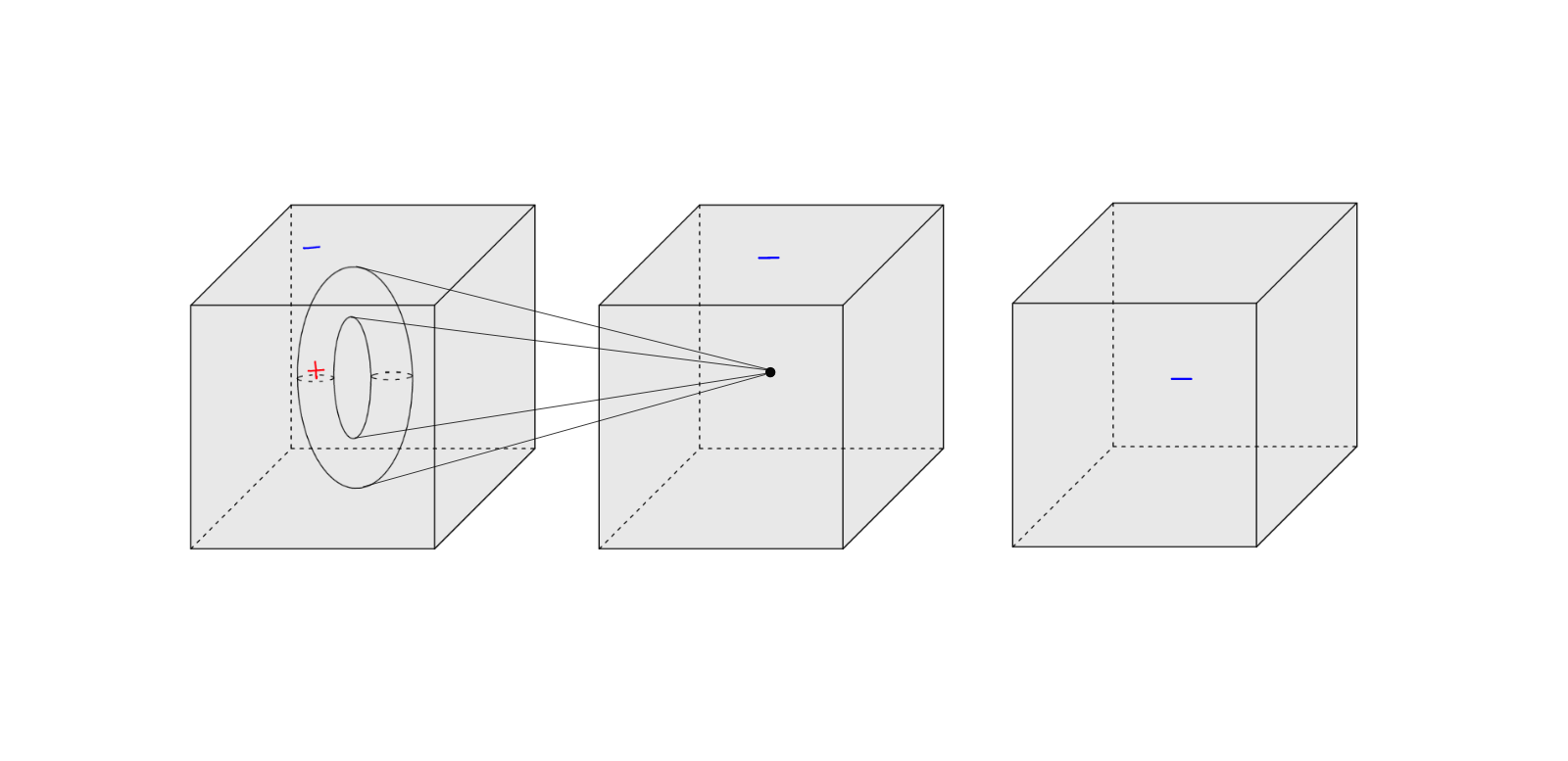}
\endgroup\end{restoretext}
is part of $\TI^{n,0}_4$ (and thus of $\TI^{n,-1}_4$) but not part of $\TI^{n,1}_3$, since the solid torus and its complement in $\partial [0,1]^4$ ($\iso S^3$) are 0-connected (i.e. path-connected) but not 1-connected. Note that for the above generator, the complement of the union of regions $-$ and $+$ does not form a manifold. Only for generators in $\TI^{n,\infty + 1}_k$ (and, based on our initial remark, possibly also $\TI^{n,\infty}_k$) can we guarantee that this complement is always a manifold.

\subsection{Adjoining invertible generators}

We now discuss how to adjoin an invertible generator to a presented associative $n$-category $\sC$. The construction will be inductive, and use two auxiliary maps $\tinv$ and $\ttinv$. These maps translate morphisms respectively labels from $\TI$ into the category obtained after adjoining the new generator to $\sC$ (more precisely, $\ttinv$ translates labels to certain \textit{fibers} in types, according to the dimension of the new generator).

\begin{constr}[Freely adjoining an invertible generator] \label{constr:adjoin_inv_gen} Let $\sC$ be an associative $n$-category. Let $\ig$ be a new label (cf. \autoref{rmk:material_set_theory}), and $x,y \in \Comp(\sC)_m$, $m < n$, such that $\abss{x}$ and $\abss{y}$ have the same globular source and target. We want to define the category $\sC \igadd{x,y} \ig$ obtained from $\sC$ by freely adjoining an invertible element $\ig$ between $x$ and $y$. We will define $\sC \igadd{x,y} \ig$ and its types inductively.

First, for $k > 0$, we define (note that in the following, $\icg{\ig}_{S\equiv T}$ is nothing but the name of a new element, cf. \autoref{rmk:material_set_theory})
\begin{equation}
\TI^{\ig}_k = \Set{ \icg{\ig}_{S\equiv T} ~|~ \ic_{S \equiv T} \in \TI_{k}}
\end{equation}
We extend this to non-positive indices by setting $\TI^{\ig}_k = \emptyset$ for $-\infty < k \leq 0$. We then set (for $0 \leq l \leq n+1$)
\begin{equation}
(\sC \igadd{x,y} \ig)_l = \sC_l \cup \TI^{\ig}_{l-m}
\end{equation}
Types for $g \in \sC_l \subset (\sC \igadd{x,y} \ig)_l $ are given by types of $\sC$. We will inductively (for $k = 0,1, ..., n - m + 1$) define the remaining types together with a function
\begin{equation}
\tinv^{\ig}_k : \Comp(\TI)_{k} \to \Comp(\sC \igadd{x,y} \ig)_{k + m}
\end{equation}
and a functor
\begin{equation}
\ttinv^{\ig}_k : \redGamma{\TI^k} \to \SIvert m {\GGamma {k + m} {\sC \igadd{x,y} \ig}}
\end{equation}
which, for any $\scD \in \Comp(\TI)_k$, satisfy the following 
\begin{enumerate}
\item \textit{Compatibility with sources and targets}:
\begin{align} \label{eq:inv_prop1}
\tinv^{\ig}_{k-1}(\gsrc(\scD)) &= \gsrc (\tinv^{\ig}_k(\scD)) \\
\tinv^{\ig}_{k-1}(\gtgt(\scD)) &= \gtgt (\tinv^{\ig}_k(\scD)) 
\end{align}

\item \textit{Compatibility with ($k$-truncated) towers:}
\begin{equation} \label{eq:inv_prop2}
\sT^k_\scD = \sT^k_{\tinv^{\ig}_k(\scD)}
\end{equation}

\item \textit{Compatibility with ($k$-level) labelling}:
\begin{equation} \label{eq:inv_prop3}
\ttinv^{\ig}_k \tsU k_\scD = \tsU k_{\tinv^{\ig}_k(\scD)}
\end{equation}

\item For each $b \in \redGamma{\TI^k}_l$, the $l$-cube with image $\ttinv^{\ig}_k(b)$ is globular and normalised.

\end{enumerate}

\noindent We remark that, using \autoref{notn:restrictions_of_labels}, conditions (ii) and (iii) can be summarised as
\begin{equation}
\SIvert k {\ttinv^{\ig}_k} \scD = \tinv^{\ig}_k(\scD)
\end{equation}

In the base case of $k = 0$ we define
\begin{align}
\tinv^{\ig}_k(-) = \ttinv^{\ig}_k(-) &= x \\
\tinv^{\ig}_k(+) = \ttinv^{\ig}_k(+) &= y
\end{align}
This trivially satisfies all inductive hypotheses (for the first condition recall our \autoref{conv:minus_one_cubes} for $(-1)$-cubes).

For $k > 0$, inductively assume that $\tinv^{\ig}_l$ and types of elements in $\TI^{\ig}_l$ have been defined for $l < k$. For $g = \icg{\ig}_{S \equiv T} \in \TI^{\ig}_k$ we define its type by
\begin{equation}
\abss{g} := \abss{\tinv^{\ig}_k(S) \xto g \tinv^{\ig}_k(T)}
\end{equation}
Writing $c = \ic_{S \equiv T}$, we first set 
\begin{equation}
\tinv^{\ig}_k(\abss{c}) := \abss{g}
\end{equation}
By \autoref{constr:double_cones_of_src_and_tgt} we know this is normalised and globular. From  \autoref{constr:adjoining_gen} we know that this is well-typed. Thus $\abss{g} \in \Comp(\sC \igadd{x,y} \ig)_{k+m}$ as required. We further set (for the notation $\ip_g$ see \autoref{defn:pres_ANC})
\begin{equation}
\ttinv^{\ig}_k(c) := \tsU k_{\abss{g}}\Delta_{\ip^k_g}
\end{equation}
We check this is normalised and globular. Normalisation follows from  \autoref{rmk:top_restriction_normalised} together with the inductive argument of \autoref{constr:top_monad}. Globularity follows since subfamilies (in particular those obtained by restrictions, cf. \autoref{eg:subfamily_by_restriction}) inherit globularity.

We next verify the other inductive conditions for these definitions. Comparing \eqref{eq:cobordism_minimal_def} to the definition of $\abss{g}$ we observe that \eqref{eq:inv_prop1} and \eqref{eq:inv_prop2} are satisfied. To see that \eqref{eq:inv_prop3} holds recall that \autoref{constr:double_cones_of_src_and_tgt} refers to \autoref{constr:top_monad} which defines $\abss{\tinv^{\ig}_k(S) \to \tinv^{\ig}_k(T)}$ in $(k + m -1)$ inductive steps. Unpacking the first $(k-1)$ steps, and using \eqref{eq:inv_prop3} inductively for $S$ and $T$, we find that for any $(p \to q) \in \tsG k(\abss{c}) = \tsG k(\abss{g})$ with $q \neq \ip^k_g$ we have
\begin{equation}
\tsU k_{\abss{g}}(p \to q) = \ttinv^{\ig}_{k-1} \tsU k_{\abss{c}}(p \to q)
\end{equation}
Unpacking the last $m$ steps, then for any $p,q \in \tsG k(\abss{c}) = \tsG k(\abss{g})$ with $\tsU k_{\abss{c}}(p) = \tsU k_{\abss{c}}(q) = b \in \TI_l$ we have $f_b$ such that
\begin{equation}
\tsU k_{\abss{g}} (q \to \ip^k_g) = f_b = \tsU k_{\abss{g}} (p \to \ip^k_g)
\end{equation}
We then define for $b \in \TI_l$, $b' \in \TI_j$, $j \leq l < k$ and $b' \to b$
\begin{align}
\ttinv^{\ig}_k(b' \to b) &:= \ttinv^{\ig}_{k-1}(b' \to b) \\
\ttinv^{\ig}_k(b \to c) &:= f_b 
\end{align}
which implies \eqref{eq:inv_prop3} for $\abss{g}$ as required. This completes the definition of $\ttinv^{\ig}_k$ (after repeating the above for any $g = \icg{\ig}_{S \equiv T} \in \TI^\ig_k$). Note that functoriality follows from functoriality of $\sU^k_{\abss{g}}$. The above also defines $\tinv^{\ig}_k$ on the image of $\abss{-} : \TI_k \to \Comp(\TI)_k$ such that all inductive conditions are satisfied.

It remains to define $\tinv^{\ig}_k$ on general morphisms. Let $\scD \in \Comp(\TI)_k$ be an arbitrary $k$-morphisms in $\TI$. We define $\scC := \tinv^{\ig}_k(\scD)$ by setting
\begin{equation}
\scC = \tsR k_{\sT_\scD, \ttinv^{\ig}_k \tsU k_\scD}
\end{equation}
This satisfies \eqref{eq:inv_prop1} \eqref{eq:inv_prop2} and \eqref{eq:inv_prop3} (the first of these conditions follows by using the last inductively). It remains to show that $\scC$ is normalised, globular and well-typed. In other words, we need to check that
\begin{equation}
\scC \in \Comp(\sC \igadd{x,y} \ig)_{m+k}
\end{equation}
The proof is \stfwd{} and detailed below.

\begin{enumerate}
\item\textit{Normalisation}: By contradiction, assume $\scC$ is not normalised and that an $l$-collapse applies to it, $\lambda : \scC \kcoll l \scB$. We argue in the following cases
\begin{enumerate}
\item If $l \leq k$, then
\begin{equation}
\tsU k_\scC = \tsU k_\scB (\vsS\lambda)^k
\end{equation}
By \eqref{eq:inv_prop3} we find factorisation through some $F$ as
\begin{equation}
\tsU k_\scD = F (\vsS\lambda)^k
\end{equation}
which together with \eqref{eq:inv_prop2} implies
\begin{equation}
\vsS\lambda : \scD \to \tsR k_{\sT_\scD, F}
\end{equation}
This contradicts $\scD$ being normalised.
\item If $l > k$, then there is $x \in \tsG {l-1}(\scC)$ such that $\lambda_x \neq \id$. Let $x^k \in \tsG k(\scD) = \tsG k(\scC)$ be its projection to level $k$. Set $b = \sU^k_\scD(x^k)\in \TI_l$. Then
\begin{equation}
\tsU k_\scC \Delta_{x^k} = \Delta_{\tinv^{\ig}_l(b)}
\end{equation}
is a subbundle of $\tsU k_\scC$, which is normalised and globular by the properties of $\ttinv^{\ig}_l$, $l \leq k$. However, by choice of $x$, $\lambda$ induces a non-trivial $(l- k)$-level collapse which yields a contradiction.
\end{enumerate}

\item \textit{Globularity}: Assume $(x \to y) \in \regcont(\tsG l(\scC))$. We need to show $\tsU l_{\scC} \Delta_{x \to y}$ normalises to the constant functor (in fact, by normalisation of $\scC$ we can just require it to be constant from the beginning). We argue in the following cases.
\begin{enumerate}
\item  If $l \geq k$, then consider the projection to level $k$, $(x^k \to y^k) \in \tsG k(\scD) = \tsG k(\scC)$. Setting $b = \tsU k_\scD(y^k) \in \TI_l$, using well-typedness of $\scD$ we find
\begin{equation}
\vvec \lambda_{\scD\sslash y^k} : \scD \sslash y^k \starcoll \Id^{k-l}_{\abss{b}}
\end{equation}
the components of this collapse (padded by $m$ identity collapses in dimension $m, m-1, ..., k+1$), induce a collapse
\begin{equation}
\vvec\mu : \scC_{y^k} \starcoll \Id^{k-l}_{\tinv^{\ig}_l\abss{b}}
\end{equation}
where the left-hand side is the subcube
\begin{equation}
\theta_{y^k} : \scC_{y^k} \mono \scC
\end{equation}
given by $\theta_{y^k}^i = (\iota^{y^k}_\scD)^i$ for $i\leq k$ and $\theta_{y^k}^i = \restemb^{i-k}_{(\iota^p_\scD)^k}$ for $i \geq k$. Now, $\Id^{k-l}_{\tinv^{\ig}_k\abss{b}}$ is globular and thus by \autoref{thm:collapse_preserves_globularity} so is $\scC_{y^k}$. By choice of $y^k$ this implies $\tsU l_{\scC} \Delta_{x \to y}$ normalises to the constant functor as required.

\item If $l < k$, then by globularity of $\scD$ we have
\begin{equation}
\vvec \lambda : \tsU l_\scD \Delta_{x \to y} \starcoll \const
\end{equation}
and as before by \eqref{eq:inv_prop3} this induces
\begin{equation}
\vvec \mu: \tsU l_\scC \Delta_{x \to y} \starcoll \const
\end{equation}
as required.
\end{enumerate}

\item \textit{Well-typedness}: Let $p \in \tsG {k+m}(\scC)$, and set $b = \tsU k_\scD(p^k) \in \TI_l$. Using the definition of $\theta_p$ in the previous item (replacing $y^k$ with $p$), we find 
\begin{equation}
\scC_p \starcoll \Id^{k-l}_{\tinv^{\ig}_l\abss{b}}
\end{equation}
Note that $\scC \sslash p$ is a subbundle of $\scC_p$ by minimality, and thus collapses to a minimal sub-bundle of $\Id^{k-l}_{\tinv^{\ig}_l\abss{b}}$. But since $\tinv^{\ig}_l\abss{b}$ is well-typed, $\Id^{k-l}_{\tinv^{\ig}_l\abss{b}}$ is well-typed too. From this well-typedness of $\scC$ follows.
\end{enumerate}

This completes the inductive construction of $\sC \igadd{x,y} \ig$.
\end{constr}

\begin{constr}[Adjoining sets of invertible generators] \label{constr:adjoining_sets_of_inv_generators} Analogous to \autoref{constr:adjoining_sets_of_gen} we can adjoin \textit{sets} of invertible generators. That is, given a label set $I$ and pairs $(x_\ig,y_\ig)$, $\ig \in I$ of $k_\ig$-morphisms in $\sC$ ($k_\ig < n$) such that $x_\ig$ and $y_\ig$ have the same source and target, then (for any ordering $\ig_1,\ig_2, ...$ of $I$) the colimit (cf. \autoref{rmk:colimit_of_presentations}) of the possibly transfinite sequence
\begin{equation}
\sC \into (\sC \igadd{x_{\ig_1},y_{\ig_1}} \ig_{1})\into ((\sC \igadd{x_{\ig_1},y_{\ig_1}} \ig_{1}) \igadd{x_{\ig_2},y_{\ig_2}} \ig_{2}) \into \dots
\end{equation}
will be denoted by
\begin{equation}
\sC \igadd{x_{\ig},y_\ig} \Set{\ig \in I}
\end{equation}
\end{constr}

\begin{notn}[Invertible generators] We usually write $\ig$ for $\icg{\ig}_{-\equiv +}$ and $\ig\inv$ for $\icg\ig_{+\equiv -}$.
\end{notn}

\begin{eg}[Adjoining an invertible generator] \label{eg:adjoining_inv_gen} Let $\sC$ be the \free{} associative $n$-category with $\sC_0 = \Set{a}$ and $\sC_1 = \Set{f}$ (in particular, the source and target of $\abss{f}$ must be $a$). Then \free{} associative $n$-category
\begin{equation}
\sD := (\sC \igadd{\abss{f}, \abss{f} \whisker 1 1 \abss{f}} \ig)
\end{equation}
obtained by adjoining an invertible generator $\ig$ between $\abss{f}$ and $\abss{f} \whisker 1 1 \abss{f}$, has the following generating morphisms. In dimension $0$ and $1$ the presentations coincide, that is, $\sD_0 = \sC_0$ and $\sD_1 = \sC_1$. $\sD_2$ contains two elements, $\ig$ and $\ig\inv$ which have types
\begin{restoretext}
\begingroup\sbox0{\includegraphics{ANCimg/page1.png}}\includegraphics[clip,trim={.2\ht0} {.2\ht0} {.2\ht0} {.3\ht0} ,width=\textwidth]{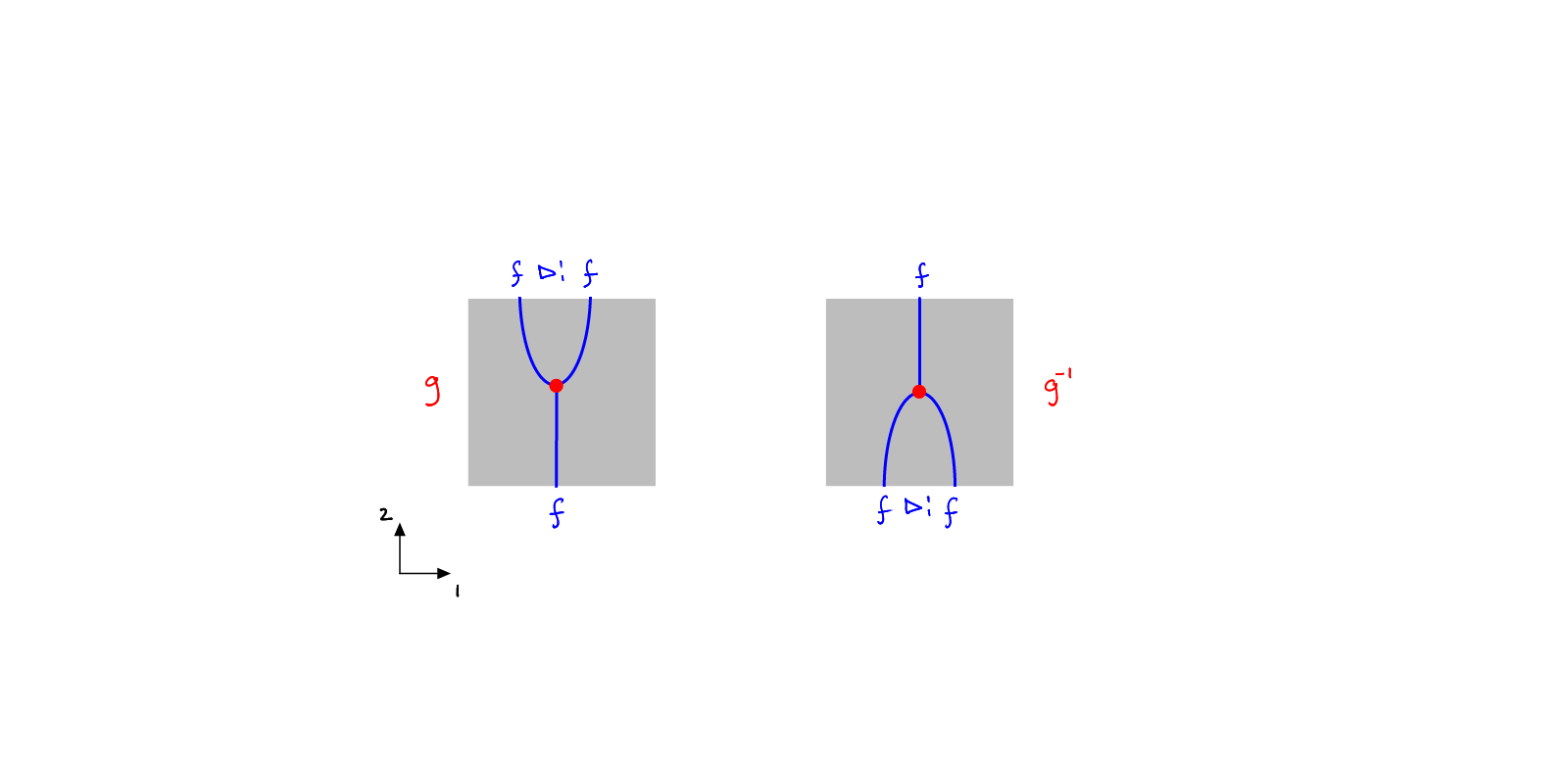}
\endgroup\end{restoretext}
These correspond to the two elements in $\TI_1$. Similarly, $\sD_3$ has elements with types such as
\begin{restoretext}
\begingroup\sbox0{\includegraphics{ANCimg/page1.png}}\includegraphics[clip,trim=0 {.0\ht0} 0 {.1\ht0} ,width=\textwidth]{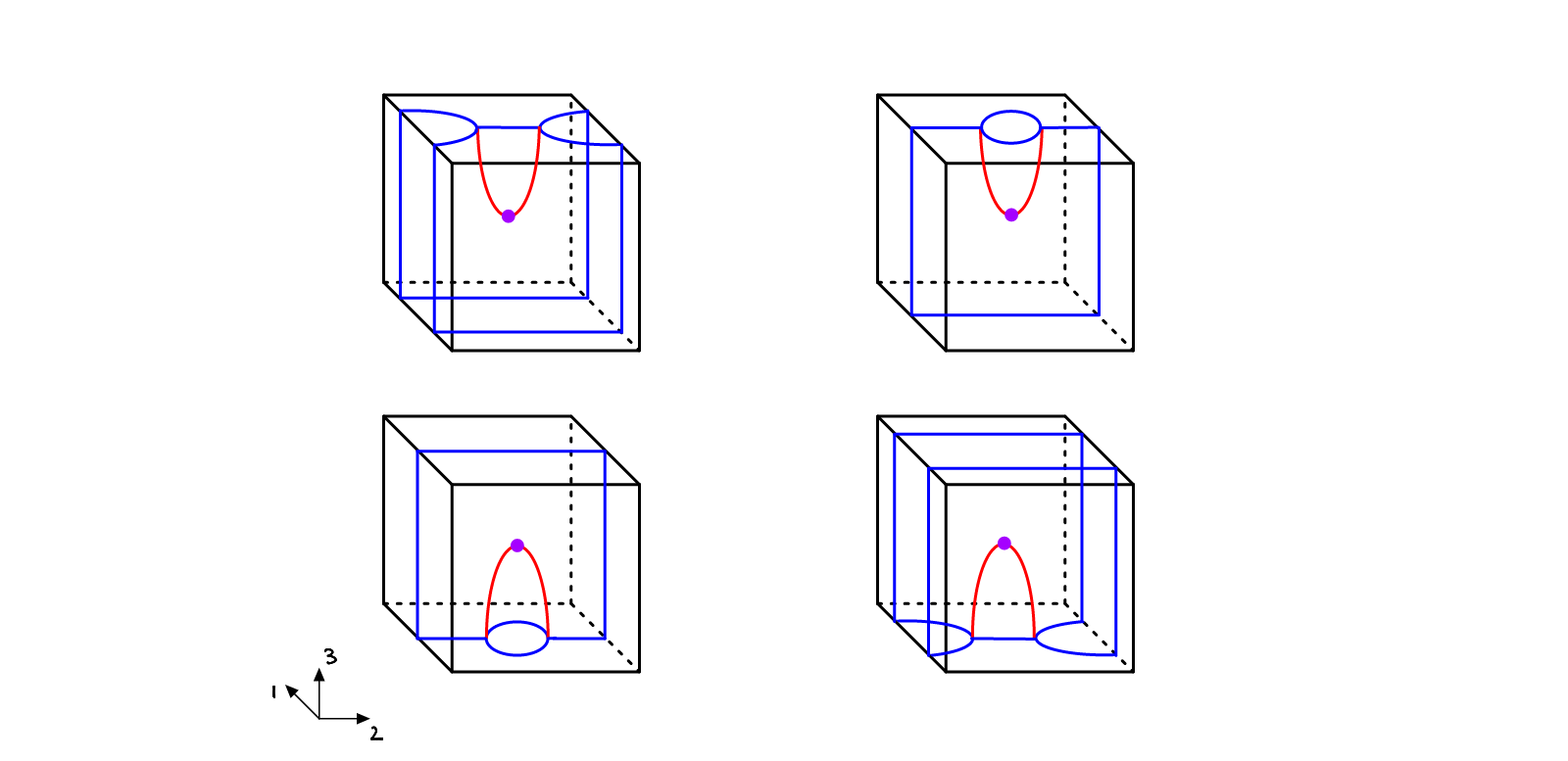}
\endgroup\end{restoretext}
These correspond to the cups and caps in $\TI_2$ which we previously discussed in \autoref{ssec:TI_eg}. More generally, $\sD_k$ with have elements that correspond to the elements of $\TI_{k-1}$.
\end{eg}

\section{Presented associative $n$-groupoids}

\subsection{Definition} \label{sec:group_def}

We are now ready to give the definition of presented associative $n$-groupoids.

\begin{constr}[\Free{} associative $n$-groupoids] \label{constr:pres_gpd} A \free{} associative $n$-groupoid $\sX$ is a \free{} associative $n$-category inductively constructed as follows. In the base case, we define $\sX^{(0)}\in \pCat_n$ by setting  $\sX^{(0)}_0 = \sX^{0}$ for some set $\sX^{0}$ called the \textit{set of $0$-cells} and $\sX^{(0)}_k$, $k > 0$ is empty.

In the $k$th step ($1 \leq k \leq n$), we then take a set $\sX^{k}$ called the \textit{set of $k$-cells}, together with morphisms $s_c, t_c \in \Comp(\sX^{(k-1)})_{k-1}$ for each $c \in \sX^{k}$, such that $s_c$ and $t_c$ agree on their sources and targets. We inductively define
\begin{equation}
\sX^{(k)} = \sX^{(k-1)} \igadd{s_c,t_c} \Set{c \in \sX^{k}}
\end{equation}
Finally, if $n < \infty$ we set $\sX = \sX^{(n)}$ and otherwise define $\sX$ as the colimit
\begin{equation}
\sX^{(0)} \into \sX^{(1)} \into \sX^{(2)} \into \dots
\end{equation}
\end{constr}

\subsection{Examples}

\begin{egs}[Groupoids] \hfill
\begin{enumerate}
\item Our previous example of a \free{} associative $2$-category $\sC$ in \autoref{sec:pres_2} is also a \free{} associative $2$-groupoid: it is obtained by adjoining a single invertible $2$-generator $d$ between $\Id_{a}$ and $\Id_{a}$, the identity on the single generating $0$-morphism in $\sC_0 = \Set{a}$. The elements of $\sC_3$ are cups and caps for $d$.

\item Our previous examples of a \free{} associative $3$-category $\sC$ in \autoref{sec:pres_2} is not an example of a \free{} associative $3$-groupoid, since $\sC_4$ only includes a biased selection of ``binary" generators for the invertibility of $d$. Namely, it contains only the snake, saddle/crotch and circle death/birth singularities, which are all binary generators (see \autoref{rmk:finite_cob_rep} for further discussion). Other singularities of (the infinitely many) singularities in $\TI_3$ have not been added to $\sC_4$. However, conjecturally there is a sub-theory $\TI^{\infty,\infty+1,\mathrm{bin}}\subset \TI$ of (certain) binary generators which is a minimal sub-theory that still equivalent to $\TI$. For such a choice of the theory of invertibility, our example of a $3$-presentation of $S^2$ would then satisfy the definition of $3$-groupoids.
\end{enumerate} 
\end{egs}

\chapter{``Perturbation-stable" generic composites} \label{ch:composition}

In this chapter we will study generic composites. Intuitively speaking, generic composites are composites (i.e. morphisms) of presented associative $n$-categories whose manifold diagrams have strata in ``generic position"---that is, slight directed perturbations of parts of the manifold diagrams will not change its equivalence class. In this sense the composite is ``perturbation-stable".

A few examples and non-examples of generic composites are the following
\begin{restoretext}
\begingroup\sbox0{\includegraphics{ANCimg/page1.png}}\includegraphics[clip,trim=0 {.05\ht0} 0 {.05\ht0} ,width=\textwidth]{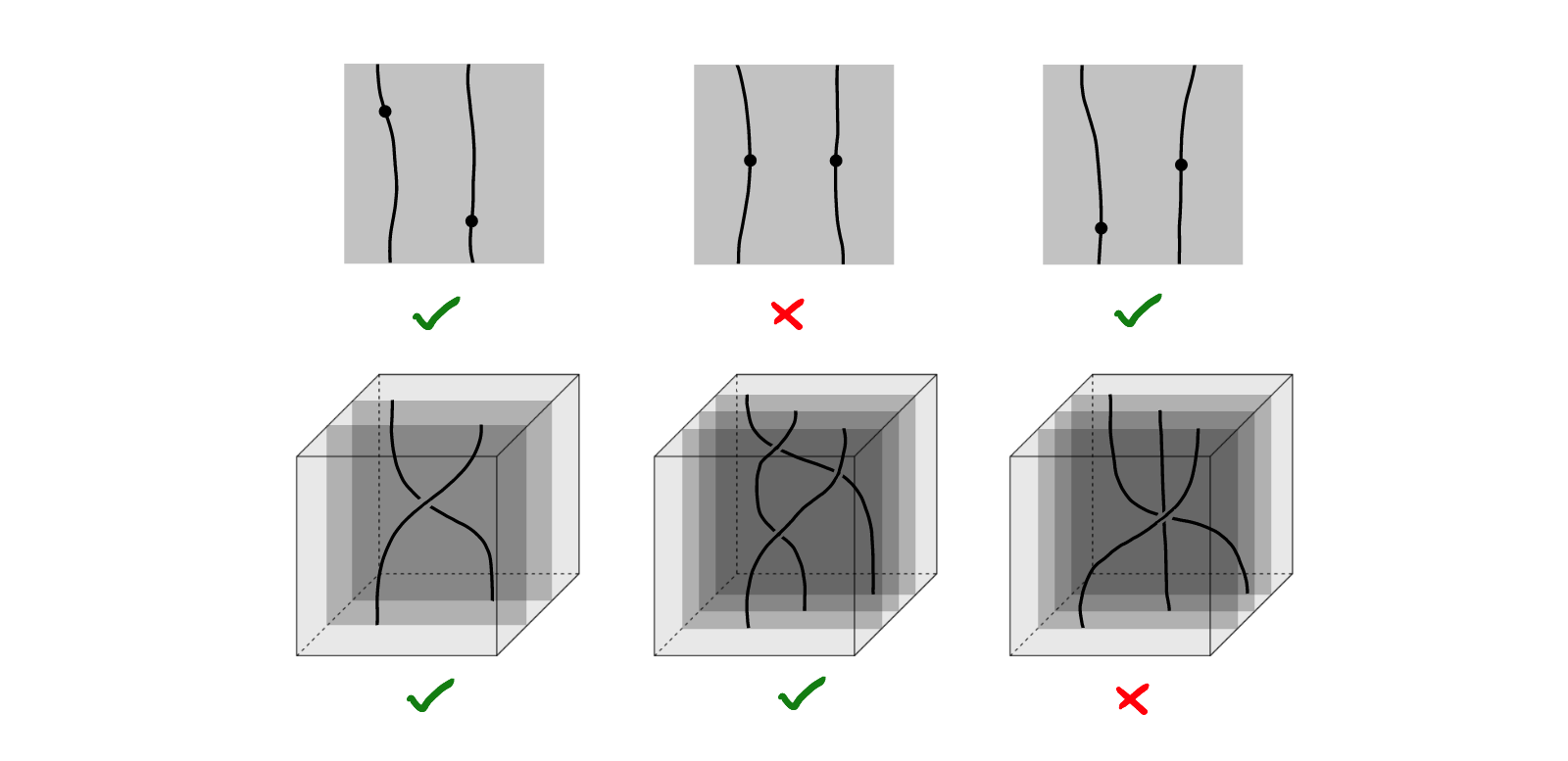}
\endgroup\end{restoretext}
Diagrams which are ``perturbation-stable" (that is, generic) have green check-marks. For the other two we remark: the second diagram can be perturbed to yield the first (or the third) diagram, and is thus not perturbation-stable; the sixth diagram can be perturbed to yield the fifth diagram (or another similar diagram, not depicted here) and is thus not perturbation-stable either.

We will suggest two characterisations of generic composites (and sketch a proof that they coincide). A first characterisation, given in \autoref{sec:top_down_gcomp}, will characterise generic composites as composites satisfying a certain genericity condition. An important definition that will be introduced for that purpose, is the (overdue and very central) notion of \textit{homotopies}, which describes manifold diagrams without $0$-dimensional singularities.

In our second characterisation, given in \autoref{sec:bottom_up_gcomp}, we  characterise generic composites as a set built inductively starting from elementary composites (namely, generators and so-called elementary homotopies) and using a ``generic composition" operation called whiskering. For this, we will first discuss the two basic modes of composition. The first one, called stacking, applies to general singular cubes families. It should be thought of as taking two cubes and gluing them along one of their sides. The second one, called whiskering, is based on the first but only applies to globular cubes. This should be thought of taking two globes and gluing them a long a mutual boundary (which is of codimension 1 to at least one of the globes, in particular, whiskering includes the ``usual" composition of two $n$-globes along a mutual $(n-1)$-boundary). We will prove that these composition operations preserve properties such as normalisation or globularity. We will then characterise generic composites of an presented associative $n$-category as being built inductively from generators and elementary homotopies using only whiskering operation.

Looking ahead to \autoref{ch:associative}, generic composites will be crucial to our definition of associative $n$-categories: since they are characterised as a set build inductively by operations from elementary composites, it will be possible to fully define functions (such as, for instance,``composition data") on generic composites by first assigning values to elementary composites and then requiring compatibility with the composition operations.

\section{``Top-down" characterisation as subset} \label{sec:top_down_gcomp}

We will define a genericity condition yielding a definition of generic composites as composites (of a presented associative $n$-category) satisfying that condition. Our discussion will rephrase (and slightly extend) the discussion in \autoref{ssec:sum_gcomps}.

\subsection{Homotopies}

As discussed in the introduction of this thesis the central difficulty of higher category theory lies in determining which coherences a higher category ought to have. In our framework coherences will be deformations of manifold diagrams which are not equivalent to the identity, but (by virtue of being deformations) they are ``locally the identity". Coherences will also be called homotopies and can be defined as follows.

\begin{defn}[Homotopies] \label{defn:homotopies} Let $\sC \in \pCat_m$. A \textit{homotopy} $H$ is a globular normalised cube $H : \bnum 1 \to \SIvert {n+1} {\GGamma{}{\cC}}$ satisfying the following condition
\begin{itemize}
\item \textit{Local triviality}: $\tsU 1_{\scD \sslash p}$ is locally trivial for all $p \in \tsG {n+1}(\scD)$. 
\end{itemize}
We note that the condition can be shown (see \autoref{lem:ctrees_well_typed}) to be equivalent to well-typedness and thus we have $H \in \Comp(\sC)_n$.
\end{defn}

Note also that the local triviality condition makes sense for general cubes $\bnum 1 \to \SIvert {n+1} \cC$, and, dropping globularity and normalisation, this lets us define homotopies for such cubes as well (however, we will not need this level of generality).

\begin{notn}[Homotopies] If $\tsG 1 (H) = \singint k$ then $H$ is called a \textit{height $k$} homotopy\footnote{A similar terminology could be introduced for all morphisms}. Setting (for $0 \leq i \leq 2k$)
\begin{equation}
f_i = \sU^1_H \Delta_{2i}
\end{equation}
then the fact that $H$ is a homotopy is also denoted by
\begin{equation}
H : f_0 \xiso {} f_1 \xiso{} ... \xiso{} f_{2k}
\end{equation}
If $k = 1$, then we also write $f_0 \xiso {H} f_1$, and if $k = 0$ then note that $H = \Id(f_0)$.
\end{notn}

Choosing the ``local triviality" condition over the usual ``well-typedness" condition was merely a cosmetic choice. We provide the following proof of their equivalence.

\begin{lem}[Technical lemma for equivalence of local-triviality and well-typedness] \label{lem:ctrees_well_typed} Let $\sC \in \pCat_\infty$ and let $\scD : \mathbf{1} \to \SIvert n {\GGamma{} {\sC}}$ be globular and normalised such that for all $0 \leq i \leq \iH_{\sU^0_\scD}$ we have
\begin{equation}
\scA_i := \sU^1_\scD \Delta_{2i} \in \Comp(\sC)_{n-1}
\end{equation}
The conditions
\begin{enumerate}
\item For all $p \in \tsG n (\scA)$ we have $\sU^1_{\scD \sslash p}$ is locally trivial
\item For all $p \in \tsG n (\scA)$ we have $\sU^n_{[\scD]} = g \in \sC_k$ satisfies $k < n$ and $[\scD \sslash p] = \Id^{n-k}_{\abss{g}}$.
\end{enumerate}
are equivalent
\proof The proof is \stfwd{}. Showing that the first condition implies the second on is harder, and will be done here. Let $p \in \tsG n (\scA)$. By minimality of $\scD \sslash p$ without loss of generality we can assume the $\iH_{\sU^0_\scD} = 1$.

We first consider the case $p^1 = 0$ or $p^1 = 2$. In this case $p$ lies in the embedding image of one of $\scA_i$. Then $(\tsU 1_\scA \sslash p) \mono \abss{T_i}$ by minimality of the minimal neighbourhood $(\tsU 1_\scA \sslash p)$ (cf. \autoref{claim:minimal_subfamily_is_minimal}). By assumption we know $\scA_0$ and $\scA_1$ are well-typed. Also note that for a general $n$-cube $\scB$ we have
\begin{equation}
\scB \sslash p = \Id_{\tsU 1_\scB \sslash p}
\end{equation}
Thus well-typedness follows from that of $\scA_i$. 

It remains to consider the case where $p^1 = 1$. Let
\begin{equation}
\vec\nfc _{\scD\sslash p} : \scD \sslash p \starcoll \NF{\scD \sslash p}^n
\end{equation}
be the sequence of collapses leading to normal form of $\scD \sslash p$. Note $\tsU 1_{\scD \sslash p}$ has globally trivial normal form by definition of a homotopy. Note that this normal form is given by $\tsU 1_{\NF{\scD\sslash p}^n_2}$ (by the above chain of collapse except the last one, cf. \autoref{thm:normal_forms_unique} where we constructed ``$k$-level normal forms" $\NF{\scA}^n_k$). Thus $\tsU 1_{\NF{\scD\sslash p}^n_2}$ is normalised and globally trivial
\begin{equation}
\tsU 1_{\NF{\scD\sslash p}^n_2} = \const
\end{equation}
Note that $\restemb_{\Delta_1} : \sU^1_{\scD\sslash p}  \Delta_1 \mono \sU^1_{\scD\sslash p}$ contains $\widetilde p := ((\iota^p_\scD)^n)\inv(p)$ in its image, by assumption on $p$. Thus $(\sU^1_{\scD\sslash p}  \Delta_1 \sslash p) = \sU^1_{\scD\sslash \widetilde p}  \Delta_1$. Setting
\begin{equation}
q = \vsS{\vec\nfc _ {\sU^1_{\scD\sslash p} }} (p)
\end{equation}
by surjectivity of multilevel collapse we find $(\sU^1_{\NF{\scD\sslash p}^n_2}  \Delta_1 \sslash q) = \sU^1_{\NF{\scD\sslash p}^n_2} \Delta_1$. By definition of multilevel collapse we also have $\sU^1_{\NF{\scD\sslash p}^n_2}(q) = g$. In words, $\sU^1_{\NF{\scD\sslash p}^n_2}  \Delta_1$ is a minimal neighbourhood around $q$, which is labelled by $g$.

By constancy of $\tsU 1_{\NF{\scD\sslash p}^n_2}$, the same must hold for $\tsU 1_{\NF{\scD\sslash p}^n_2} \Delta_0$ and $\tsU 1_{\NF{\scD\sslash p}^n_2}\Delta_2$, that is they are minimal neighbourhoods of a point labelled by $g$. But $\tsU 1_{\NF{\scD\sslash p}^n_2} \Delta_0$ and $\tsU 1_{\NF{\scD\sslash p}^n_2}\Delta_2$ are normalised (by \autoref{thm:normalisation_on_dc_restrictions}). Thus they must equal $\Id^{n-1-m}_{\abss{g}}$ by $\scA_1$ resp. $\scA_2$ being well-typed. We infer that $\tsU 1_{\NF{\scD\sslash p}^n_2} \Delta_1$ also equals $\Id^{n-1-m}_{\abss{g}}$, and thus for all $i \in \singint 1$
\begin{equation}
\tsU 1_{\NF{\scD\sslash p}^n_2} = \const_{\Id^{n-1-m}_{\abss{g}}(0)}
\end{equation}
The final $1$-level collapse $\nfc_{\scD\sslash p}^1 : \NF{\scD\sslash p}^n_2 \to \NF{\scD\sslash p}^n$, by virtue of being the maximal possible collapse, must then be witnessed by single-component natural transformation given by $(\nfc_{\scD\sslash p}^1)_{0}: \singint 1 \to \singint 0$ for the object $0 \in \bnum{1}$, and we find
\begin{equation}
\NF{\scD \sslash p}^{n+1} = \Id^{n-m}_{\abss{g}}
\end{equation}
as required. \qed
\end{lem}

\subsection{Definition of $\GComps$}

We start with two small pieces of terminology.

\begin{notn}[Homotopically non-trivial sub-morphisms] \label{notn:non_identity_subcubes} Let $\sC \in \pCat_\infty$. Let $f \in \Comp(\sC)_n$, and assume we are given a subcube $\theta : g \mono f$. We say $\theta$ is a \textit{homotopically non-trivial} subcube if there does not exists a $h \in \Comp(\sC)_{n-1}$ and a homotopy $H \in \SIvert {n+1} \cC$ such that
\begin{align}
\gsrc (H) = \NF{g}
\end{align}
and
\begin{equation}
\gtgt (H) = \Id(h)
\end{equation}
\end{notn}

We remark that for $f \in \Comp(\sC)_n$ and $g \mono f$, we always have $\NF{g} \in \Comp(\sC)_n$. Note that this would not be the case if we had chosen $\partial$-globularity over strict globularity (cf. \autoref{rmk:strict_glob}).

\begin{notn}[Subcubes concentrated around a singular height] Let $\sC \in \pCat$. Let $f \in \Comp(\sC)_n$, and assume we are given a subcube $\theta : g \mono f$, and let $a \in \sing(\sG^1(f))$. We say $\theta$ is \textit{concentrated at the singular height $a$} of $\scA$ if $\im(\theta^1) = \left.\sG^1(\scA) \sslash a\right.$. 
\end{notn}

Now we are in the position to formalise the idea of generic composites as follows.

\begin{defn}[Generic composites] \label{defn:top_down_gen_comp} Let $\sC$ be a presented associative $n$-category. We define a globular subset $\GComps(\sC) \subset \Comp(\sC)$ of \textit{generic composites} inductively in each degree as follows. In degree $0$ we set
\begin{equation}
\GComps(\sC)_0 := \Comp(\sC)_0
\end{equation}
Now let $f \in \Comp(\sC)_k$, $k > 0$. We say $f$ is a \textit{generic $k$-composite} and declare $f \in \GComp(S)$ if the following two conditions hold.
\begin{enumerate}
\item Firstly, for each $a \in \reg(\sG^1(f))$ we have that
\begin{equation}
\sU^1_f \Delta_a 
\end{equation}
is a generic $(k-1)$-composite.
\item Secondly, for each $a \in \sing(\sG^1(f))$, there is a unique homotopically non-trivial subcube $\theta_a : f_a \mono f$ concentrated around singular height $a$ such that for any other homotopically non-trivial $\theta : \scA \mono f$ we have $\theta_a \mono \theta$.
\end{enumerate}
If only the second condition holds we say $f$ is \textit{non-recursively} generic.
\end{defn}

Morally this definition says that at most one interesting thing can happen at each singular height (and recursively so). In that sense, all ``events" are in ``generic position".

\section{``Bottom-up" characterisation as inductively constructed set} \label{sec:bottom_up_gcomp}

We now give a second characterisation of generic composites: for this, we will think about the operations that can be used to build new cubes from old cubes.

\subsection{Stacking of $\SI$-families}

We start with the basic operation underlying the rest of the section: namely, stacking two singular intervals on top of one another.

\begin{constr}[Stacking $\SI$-families]  \label{constr:SI_family_stacking} Let $\scA, \scB : X \to \SI$ be $\SI$-families over a poset $X$. We construct $\scA \stack \scB$ called the \textit{stacking of $\scB$ on $\scA$}. For this we note that the \textit{ordered sum functor} $( - \uplus - ) : \Pos \times \Pos$ restricts (with respect to the embedding of \autoref{rmk:SI_is_cat_of_tot_ord}) to a functor
\begin{equation}
( - \uplus - ) : \SI \times \SI \to \SI
\end{equation}
which can be illustrated as follows
\begin{restoretext}
\begingroup\sbox0{\includegraphics{test/page1.png}}\includegraphics[clip,trim=0 {.0\ht0} 0 {.0\ht0} ,width=\textwidth]{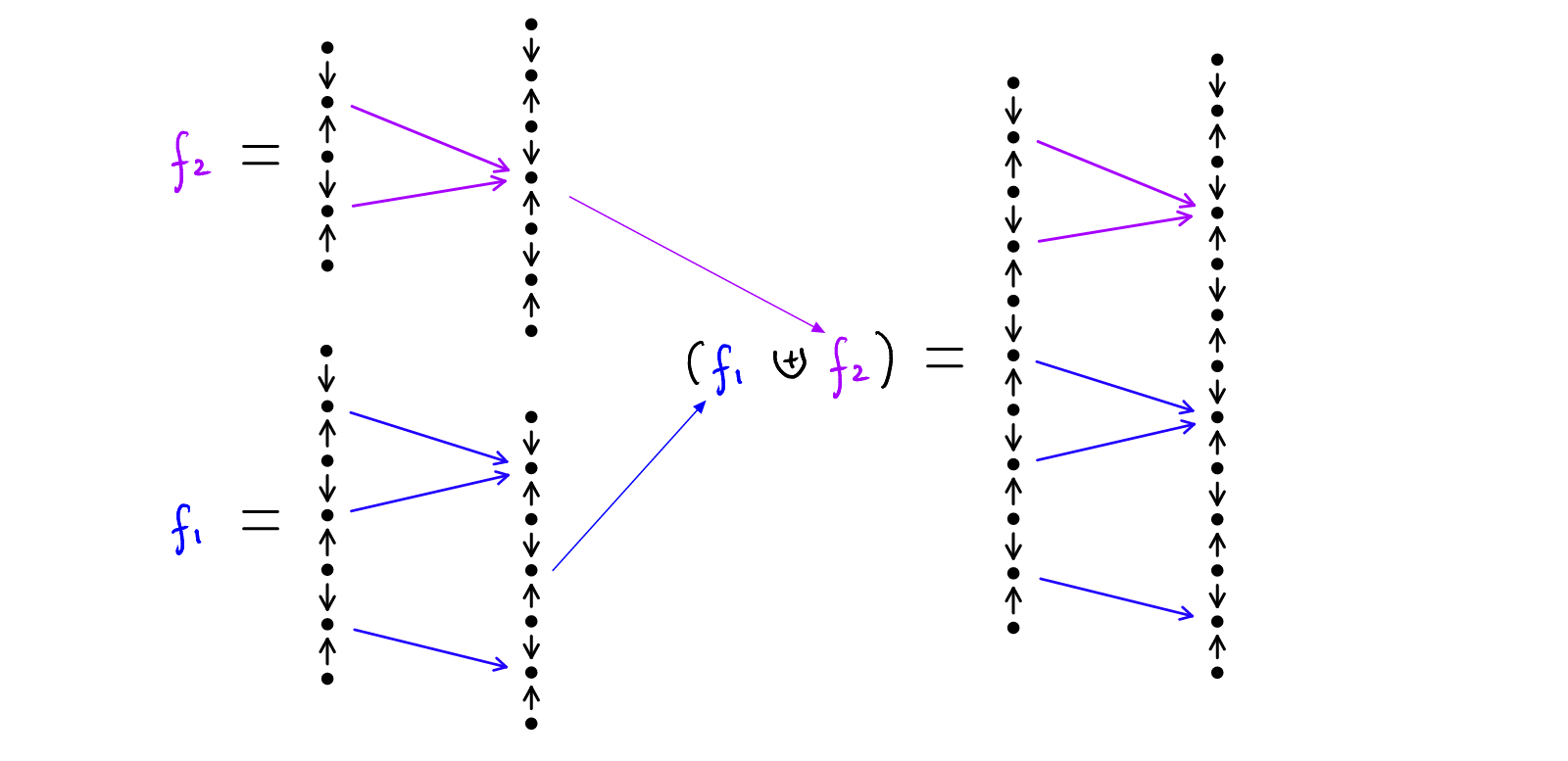}
\endgroup\end{restoretext}
Explicitly, the definition goes as follows.
\begin{itemize}
\item On objects $I,J \in \SI$, we define $I \uplus J$ to be the singular interval of height 
\begin{equation}
\iH_{(\scA \stack \scB)(x)} = \iH_{\scA(x)} + \iH_{\scB(x)}
\end{equation}
(recall from \autoref{defn:singular_intervals} that a singular interval is uniquely determined by its height.)
\item On morphism $(f : I_1 \to I_2), (g : J_1 \to J_2) \in \mor(\SI)$, we define, if $ a\in [0,\iH_{I_1}]$,
\begin{equation} \label{eq:lower_stacking_morphism}
(f \uplus g)(a) =f(a)
\end{equation}
and if $a \in [\iH_{\scA(x)} ,\iH_{\scA(x)} + \iH_{\scB(x)}]$
\begin{equation} \label{eq:upper_stacking_morphism}
(f \uplus g)(a) = \iH_{I_2} + g(a - \iH_{I_1})
\end{equation}
\end{itemize}
We then construct $(\scA \stack \scB) : X \to \SI$, the \textit{stacking of $\scA$ with $\scB$}, as the composite
\begin{equation} \label{eq:most_basic_stack}
X \xto {\Delta} X \times X \xto {\scA \times \scB} \SI \times \SI \xto {\uplus} \SI
\end{equation}
where $\Delta$ is the diagonal functor mapping $x \to (x,x)$ ($x \in X$).
\end{constr}

\begin{eg}[Stacking \SI-families] As an example, consider the families $\scA$ and $\scB$ defined by their \SI-bundles as
\begin{restoretext}
\begingroup\sbox0{\includegraphics{test/page1.png}}\includegraphics[clip,trim=0 {.15\ht0} 0 {.1\ht0} ,width=\textwidth]{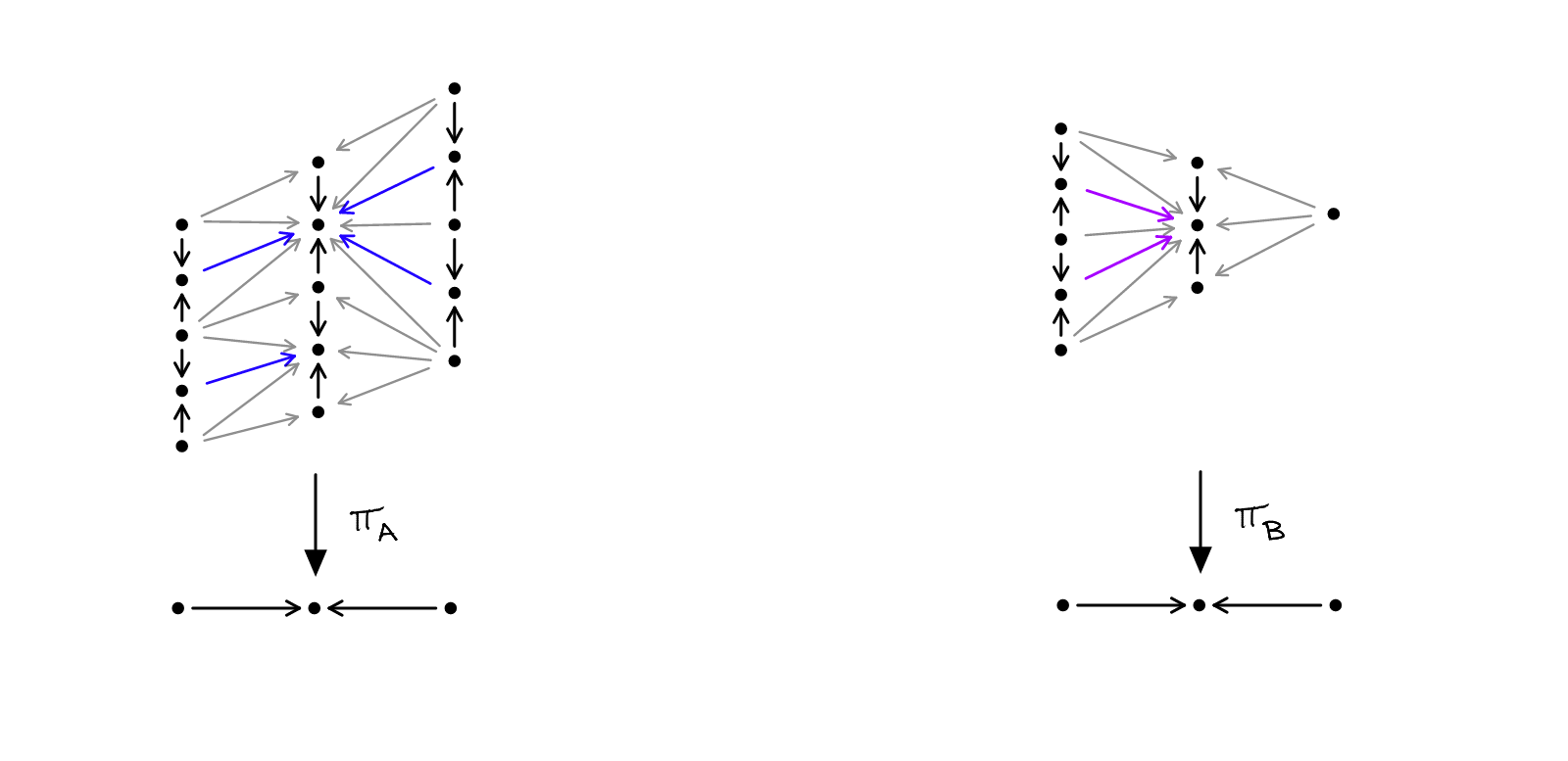}
\endgroup\end{restoretext}
Then the stacking of $\scA$ with $\scB$ is given by the \SI-family
\begin{restoretext}
\begingroup\sbox0{\includegraphics{test/page1.png}}\includegraphics[clip,trim=0 {.1\ht0} 0 {.05\ht0} ,width=\textwidth]{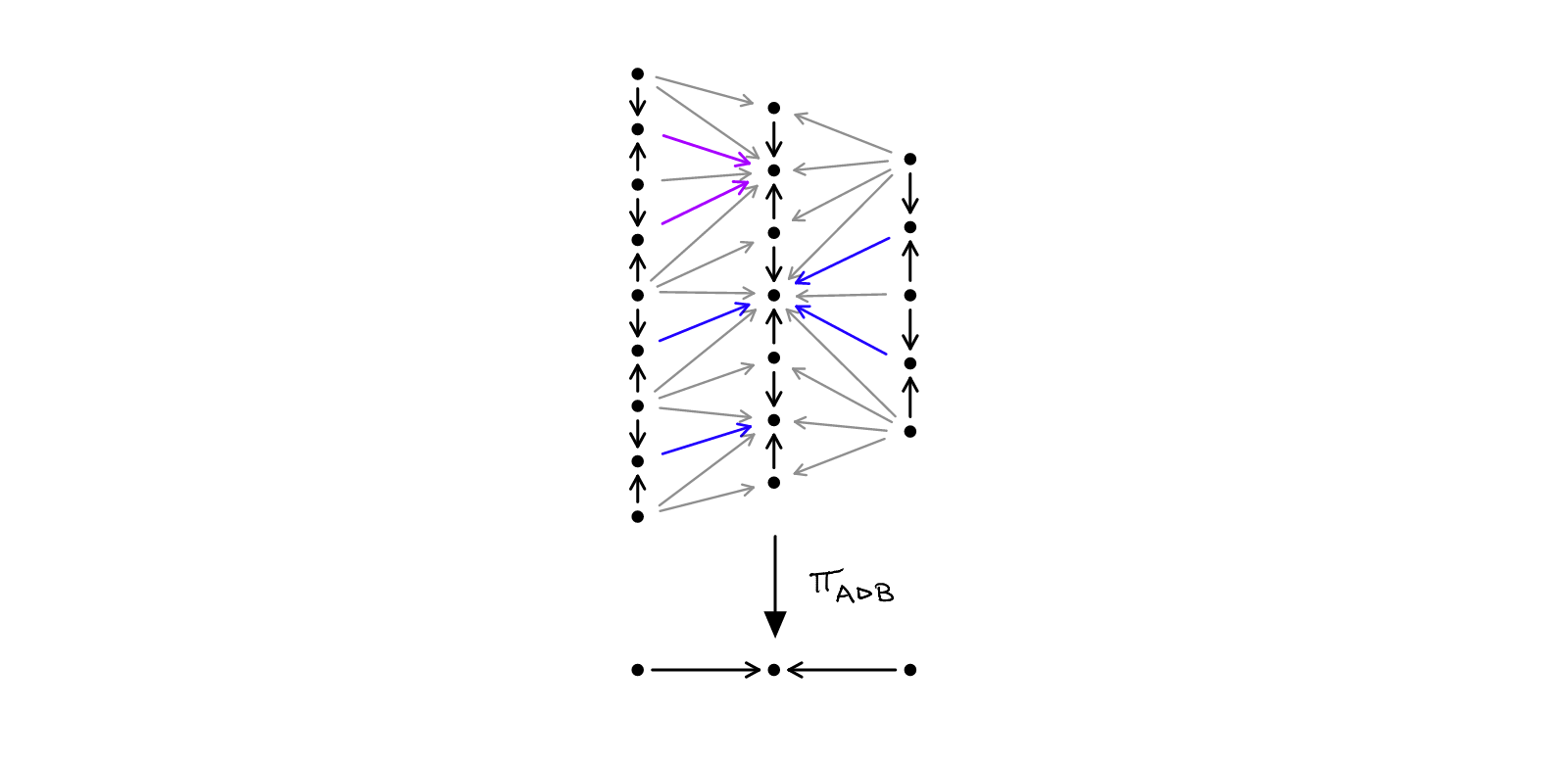}
\endgroup\end{restoretext}
\end{eg}

\begin{constr}[Isomorphism of total space of stacking and glueing of total spaces] \label{constr:stacking_gluing_iso} Let $\scA, \scB : X\to \SI$. We construct a gluing of total spaces $\sG(\scA) \cup_X \sG(\scB)$ and a canonical isomorphism of posets
\begin{equation} \label{eq:source_target_pushout2}
\sI_{  \scA,  \scB} : \sG(  \scA) \cup_X \sG(  \scB) \iso \sG(  \scA \stack   \scB) 
\end{equation}
This is \stfwd{}. The explicit construction proceeds as follows. We first define
\begin{align}
\stacklow : \sG(\scA) \into \sG(\scA \stack \scB) \\
\stackup : \sG(\scB) \into \sG(\scA \stack \scB)
\end{align}
(note that $ \scA, \scB$ are implicit in this notation). Namely, for $(x,a) \in \sG(  \scA)$ we set
\begin{equation}
\stacklow (x,a) = (x,a)
\end{equation}
and for $(x,a) \in \sG(  \scB)$ we set
\begin{equation}
\stackup (x,a) = (x,a + \iH_{  \scA(x)})
\end{equation}
We claim both maps are full inclusions. We show this for $\stackup$ (the proof for $\stacklow$ is similar). We claim
\begin{equation}
(x,b_1) \to (x,b_2) \in \sG(\scB) \iff (x,b_1 + \iH_{\scA(x)}) \to (y,b_2 + \iH_{\scA(y)}) \in \sG(\scA \stack \scB)
\end{equation}
If $b_1 \in \sing(\scB(x))$ this follows by \eqref{eq:defn_order_realisation_1} and definition of $\scA \stack \scB$. If $b_1 \in \reg(\scB(x))$ then by \eqref{eq:defn_order_realisation_2} we need to show
\begin{align}
&\wwidehat \scB (b_1 - 1) \leq b_2 \leq \wwidehat \scB (b_1 + 1) \\
 \iff & \wwidehat {\scA \stack \scB} (b_1 + \iH_{\scA(x)} - 1) \leq b_2 + \iH_{\scA(y)} \leq \wwidehat {\scA \stack \scB} (b_1 + \iH_{\scA(x)} + 1)
\end{align}
This again follows by definition of $\scA \stack \scB$, and in particular, since $\scA(x \to y)(\iH_{\scA(x)}) < \iH_{\scA(y)}$.

Note that the definitions entail
\begin{equation}
G := \stacklow  \mtgt_\scA(x) = \stackup \msrc_\scA(x) 
\end{equation}
Thus, employing \autoref{constr:pushouts_in_labelled_posets} in $\Poss {\sG(  \scA \stack   \scB)}$), we find a pushout
\begin{equation}
\xymatrix{ (X,G) \ar[r]^{\mtgt_{\scA}} \ar[d]_{\msrc_{\scB}} & (\sG(\scA), \stacklow) \ar[d]^{} 
\\ (\sG(\scB),\stackup) \ar[r]_-{} & (\sG(\scA) \cup_X \sG(\scB), \sI_{  \scA,  \scB}) \pushoutfar }
\end{equation}
However the definitions of our full inclusions $\stacklow$ and $\stackup$ further satisfy  that $\im(G) = \im(\stacklow) \cap \im(\stackup)$ as well as $\sG(\scA \stack \scB) = \im(\stacklow) \cup \im(\stackup)$. This identifies $\sG(\scA \stack \scB)$ as a pushout itself and thus $\sI_{  \scA,  \scB}$ is an isomorphism.
\end{constr}

The next remark shows that the definitions of $\stacklow$ and $\stackup$ can in fact be expressed using tools developed in \autoref{ch:emb}.

\begin{rmk}[Endpoints inclusions for stacking] \label{rmk:canonical_total_space_inclusion} With the assumptions of the previous constructions, recall the poset inclusions $\stacklow$ and $\stackup$. We illustrate their definitions by the following example
\begin{restoretext}
\begingroup\sbox0{\includegraphics{test/page1.png}}\includegraphics[clip,trim=0 {.1\ht0} 0 {.1\ht0} ,width=\textwidth]{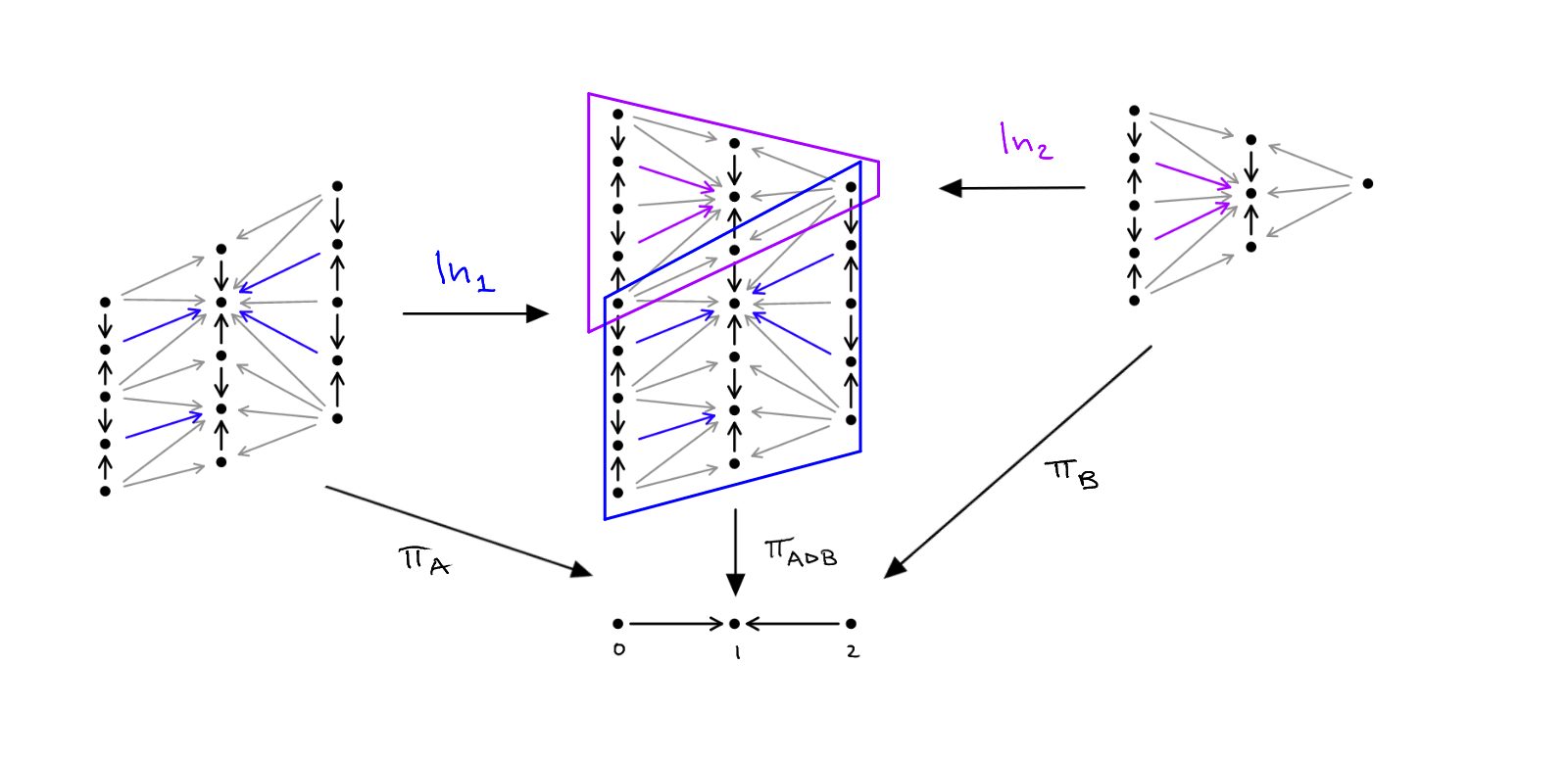}
\endgroup\end{restoretext}
Comparing \autoref{defn:stacking_labelled_families} and \autoref{constr:subfamilies_from_endpoints}, we see that these definitions actually coincide with the following previous constructions
\begin{equation}
\sJ^{\und{\scA \stack \scB}}\restsec{[\stacklow \msrc_{\und\scA}, \stacklow \mtgt_{\und\scA}]} = \stacklow
\end{equation}
and
\begin{equation}
\sJ^{\und{\scA \stack \scB}}\restsec{[\stackup \msrc_{\und\scB}, \stackup \mtgt_{\und\scB}]} = \stackup
\end{equation}
For instance for the above example, the following depicts the open sections $\stacklow \msrc_{\und\scA}, \stacklow \mtgt_{\und\scA}$ by \cblue{} and \cgreen{} circles respectively, and the open sections $\stackup \msrc_{\und\scB}, \stackup \mtgt_{\und\scB}$ by \cgreen{} and \cred{} circles respectively
\begin{restoretext}
\begingroup\sbox0{\includegraphics{test/page1.png}}\includegraphics[clip,trim=0 {.1\ht0} 0 {.1\ht0} ,width=\textwidth]{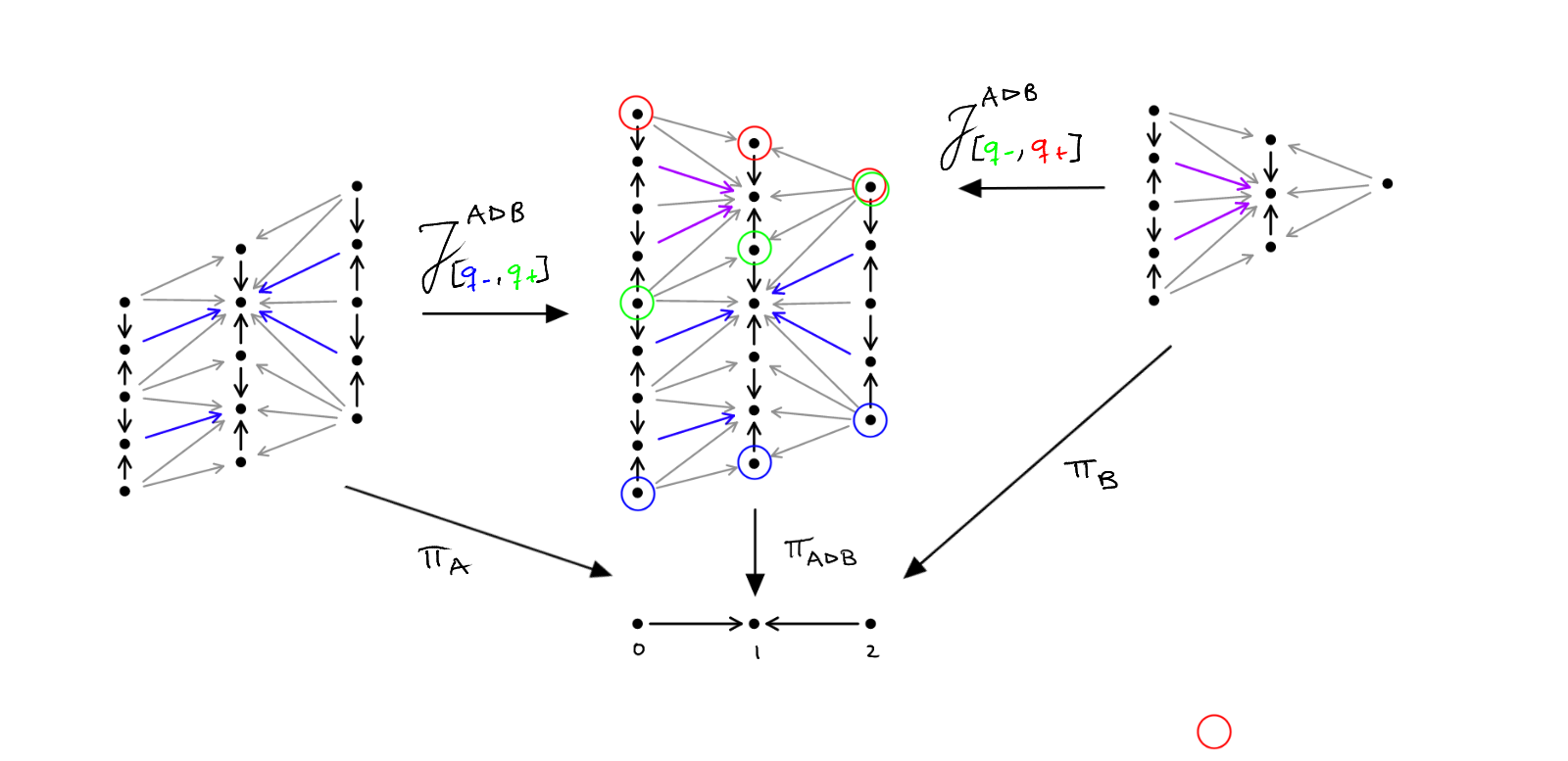}
\endgroup\end{restoretext}
\end{rmk}

\subsection{Stacking of labelled $\SI$-families}

Next, we want to generalise the stacking construction to $\SIvertone  \cC$ families. However, care needs to be taken with labelling functors in this case. For instance, consider $\SIvertone  \cC$ families $\scA$ and $\scB$ defined by
\begin{restoretext}
\begingroup\sbox0{\includegraphics{test/page1.png}}\includegraphics[clip,trim=0 {.15\ht0} 0 {.1\ht0} ,width=\textwidth]{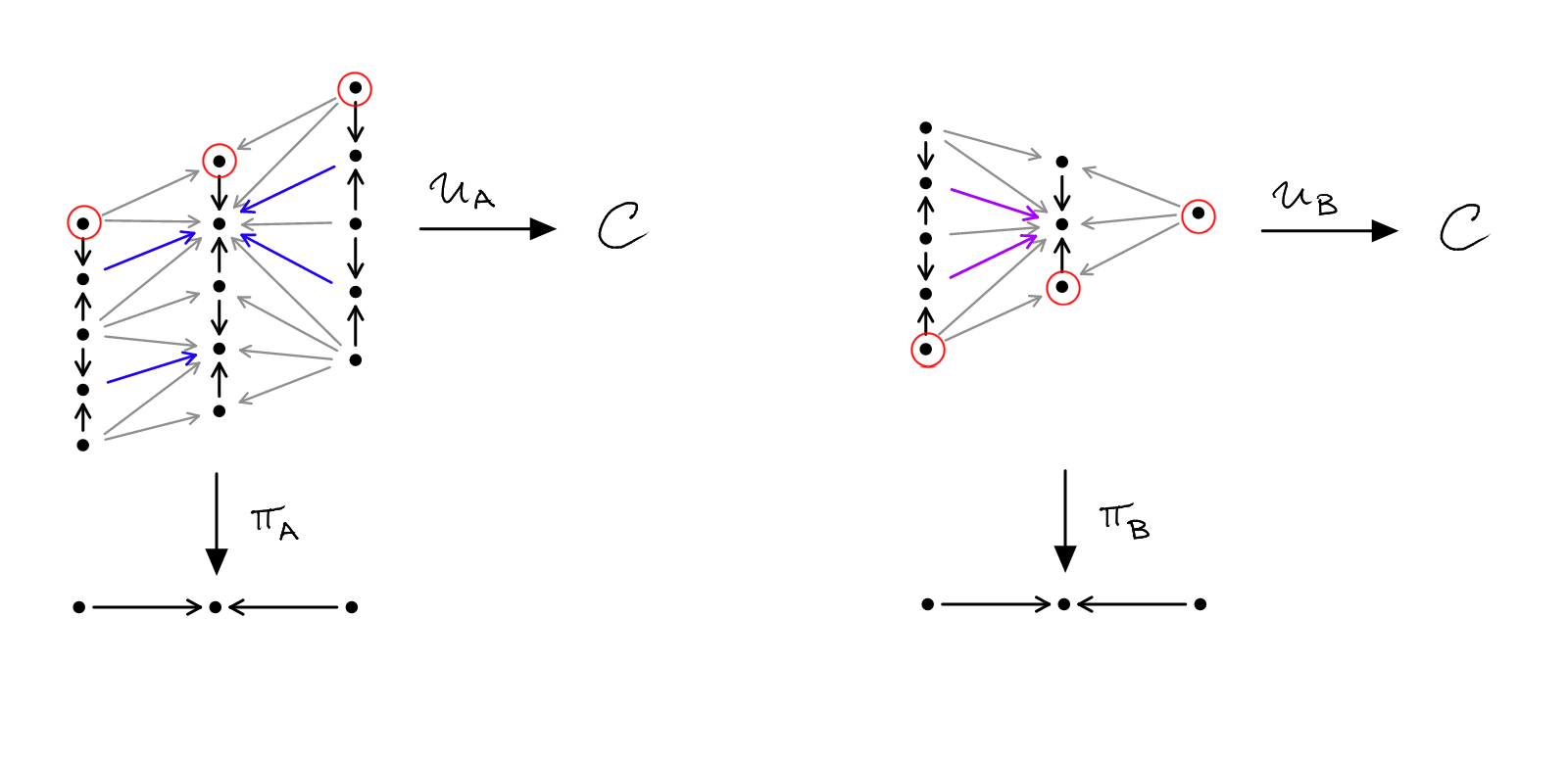}
\endgroup\end{restoretext}
We would like to define a stacking $\scA \stack \scB$ which inherits its labelling from $\scA$ and $\scB$. For this to be possible, the labelling $\sU_\scA$ on the target of $\scA$ (whose objects are marked by \cred{} circles in $\tsG 1(\scA)$ above) needs to coincide with the labelling $\sU_\scB$ on the source of $\scB$ (whose objects are also marked by \cred{} circles in $\tsG 1(\scB)$). These parts of $\sG(\scA)$ and $\sG(\scB)$ will be identified in $\sG(\scA \stack \scB)$ (they both correspond to the full subposet of $\sG(\scA \stack \scB)$ highlighted by \cred{} circles) below
\begin{restoretext}
\begingroup\sbox0{\includegraphics{test/page1.png}}\includegraphics[clip,trim=0 {.1\ht0} 0 {.09\ht0} ,width=\textwidth]{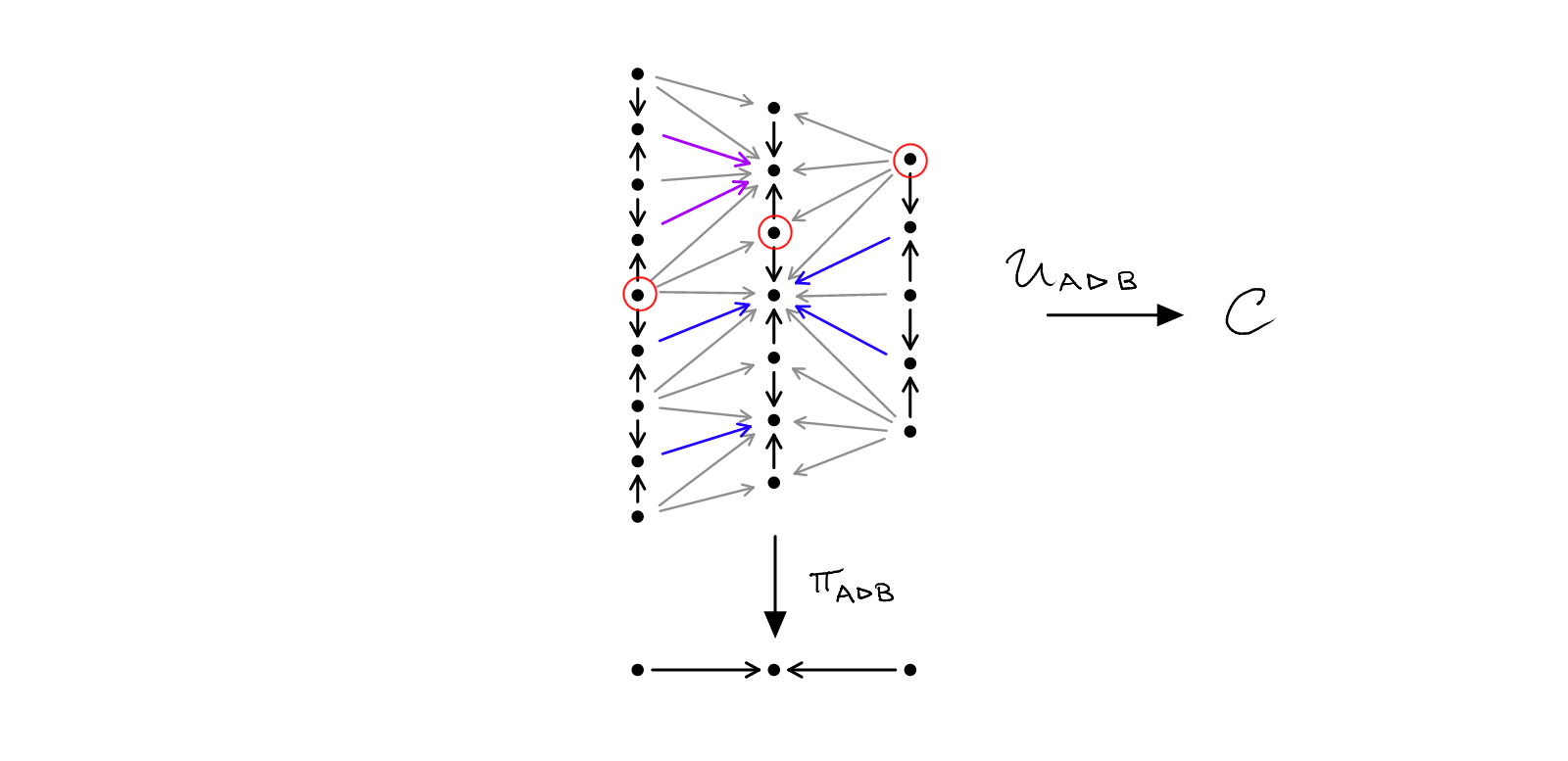}
\endgroup\end{restoretext}
Thus the assumptions on $\sU_\scA$ and $\sU_\scB$ allow us to naturally provide a labelling $\sU_{\scA \stack \scB}$.

This observation is formalised in the following construction. To define the labelling $\sU_{\scA \stack \scB}$ we use the gluing of labellings defined in \autoref{constr:pushouts_in_labelled_posets}.

\begin{constr}[Stacking $\cC$-labelled $\SI$-families]\label{defn:stacking_labelled_families} Let $\scA,\scB : X \to \SIvertone \cC$ be $\cC$-labelled $\SI$-families over a poset $X$, such that
\begin{equation} \label{eq:source_target_pushout_cond}
U := \sU_\scA \mtgt_{\und \scA} = \sU_\scB \msrc_{\und \scA} 
\end{equation}
We define the \textit{stacking of $\scB$ on $\scA$} which will be denoted by
\begin{equation}
\scA \stack \scB : X \to \SIvertone \cC
\end{equation}
Employing \autoref{constr:pushouts_in_labelled_posets} we find a pushout (in $\Poss \cC$)
\begin{equation} \label{eq:source_target_pushout}
\xymatrix{ (X,U) \ar[r]^{\mtgt_{\scA}} \ar[d]_{\msrc_{\scB}} & (\sG(\scA),\sU_\scA) \ar[d]^{ } 
\\ (\sG(\scB),\sU_\scB) \ar[r]_-{ } & (\sG(\scA) \cup_X \sG(\scB), \sU_\scA \cup_X \sU_\scB) \pushoutfar }
\end{equation}
Note that $\mtgt_{\scA}$ and $\msrc_{\scB}$ satisfy the required conditions for this construction by assumption \eqref{eq:source_target_pushout_cond} and \autoref{rmk:src_and_tgt_have_lifts}. 

We now construct the \textit{stacking of $\scB$ on $\scA$} as the $\cC$-labelled $\SI$-family denoted by $\scA \stack \scB$ and defined by
\begin{equation} \label{eq:stack_basic_def}
\scA \stack \scB := \sR_{(\und \scA \stack \und \scB), (\sU_\scA \cup_X \sU_\scB) \sI\inv_{\und \scA, \und \scB}}
\end{equation}
In other words, $\scA \stack \scB$ is the $\cC$-labelled $\SI$-family obtained from the $\SI$-family $\und \scA \stack \und \scB : X \to \SI$, together with the labelling
\begin{equation}
\sG(\und \scA \stack \und \scB) \xto {\sI\inv_{\und \scA, \und \scB}} \sG(\und \scA) \cup_{X} \sG(\und \scB) \xto {\sU_\scA \cup_X \sU_\scB} \cC
\end{equation}
which uses the definitions in \eqref{eq:source_target_pushout2} and \eqref{eq:source_target_pushout}.
\end{constr}

As an immediate consequence of the construction, we remark that from \eqref{eq:stack_basic_def} it follows that
\begin{align} \label{eq:subbund_of_stack}
\sU_{\scA \stack \scB} \stacklow &= \sU_\scA \\
\sU_{\scA \stack \scB} \stackup &= \sU_\scB
\end{align}

The following remark states the interaction of sources and targets with stacking.

\begin{rmk}[Source and target for stacking] \label{rmk:src_tgt_for_stacking} By a fibrewise comparison (and using \eqref{eq:stack_basic_def}) we have
\begin{align}
\tsU 1_{\scA \stack \scB} \msrc_{\und{\scA \stack\scB}} &= (\sU_\scA \cup_X \sU_\scB) \sI\inv_{\und \scA, \und \scB} \msrc_{\und{\scA} \stack\und{\scB}} \\
&= \sU_\scA \msrc_{\und\scA}
\end{align}
and similarly we find
\begin{equation}
\tsU 1_{\scA \stack \scB} \mtgt_{\und{\scA \stack\scB}} = \sU_\scB \mtgt_{\und\scB} 
\end{equation}
\end{rmk}

Finally, we discuss how base change can be applied to stackings.

\begin{rmk}[Base change for stacking] \label{rmk:basechange_stacking} Assume $\scA$, $\scB$ as in the previous construction, and assume $H : Y \to X$. Note that by precomposing \eqref{eq:most_basic_stack} with $H$ we obtain
\begin{equation}
(\und\scA \stack \und\scB) H = (\und\scA H)\stack (\und\scB H)
\end{equation}
Using the previous constructions, a fibrewise comparison then yields
\begin{align}
&\big(\sG((\und\scA \stack \und\scB) H) \xto {\sG(H)} \sG(\scA \stack \scB) \xto {\sI\inv_{\und\scA,\und\scB}} \sG(\und \scA) \cup_{X} \sG(\und \scB) \xto {\sU_\scA \cup_X \sU_\scB} \cC\big) \\
= ~&\big(\sG(\und\scA H\stack \und\scB H)) \xto {\sI\inv_{\und\scA H,\und\scB H}} \sG(\und \scA H) \cup_{Y} \sG(\und \scB H) \xto {\sU_\scA \sG(H) \cup_X \sU_\scB\sG(H)} \cC\big)
\end{align}
Using these two equations together with \autoref{claim:grothendieck_span_construction_basechange} we find that
\begin{align}
(\scA \stack \scB) H &= \sR_{(\und \scA \stack \und \scB), (\sU_\scA \cup_X \sU_\scB) \sI\inv_{\und \scA, \und \scB}} H \\
&= \sR_{(\und \scA H \stack \und \scB H), (\sU_{\scA H} \cup_Y \sU_{\scB H})\sI\inv_{\und \scA H, \und \scB H}} \\
&= (\scA H\stack \scB H)
\end{align}
Thus a base changed stacking of families is the stacking of the base changed families. This is illustrated in the following examples
\begin{restoretext}
\begingroup\sbox0{\includegraphics{test/page1.png}}\includegraphics[clip,trim=0 {.0\ht0} 0 {.0\ht0} ,width=\textwidth]{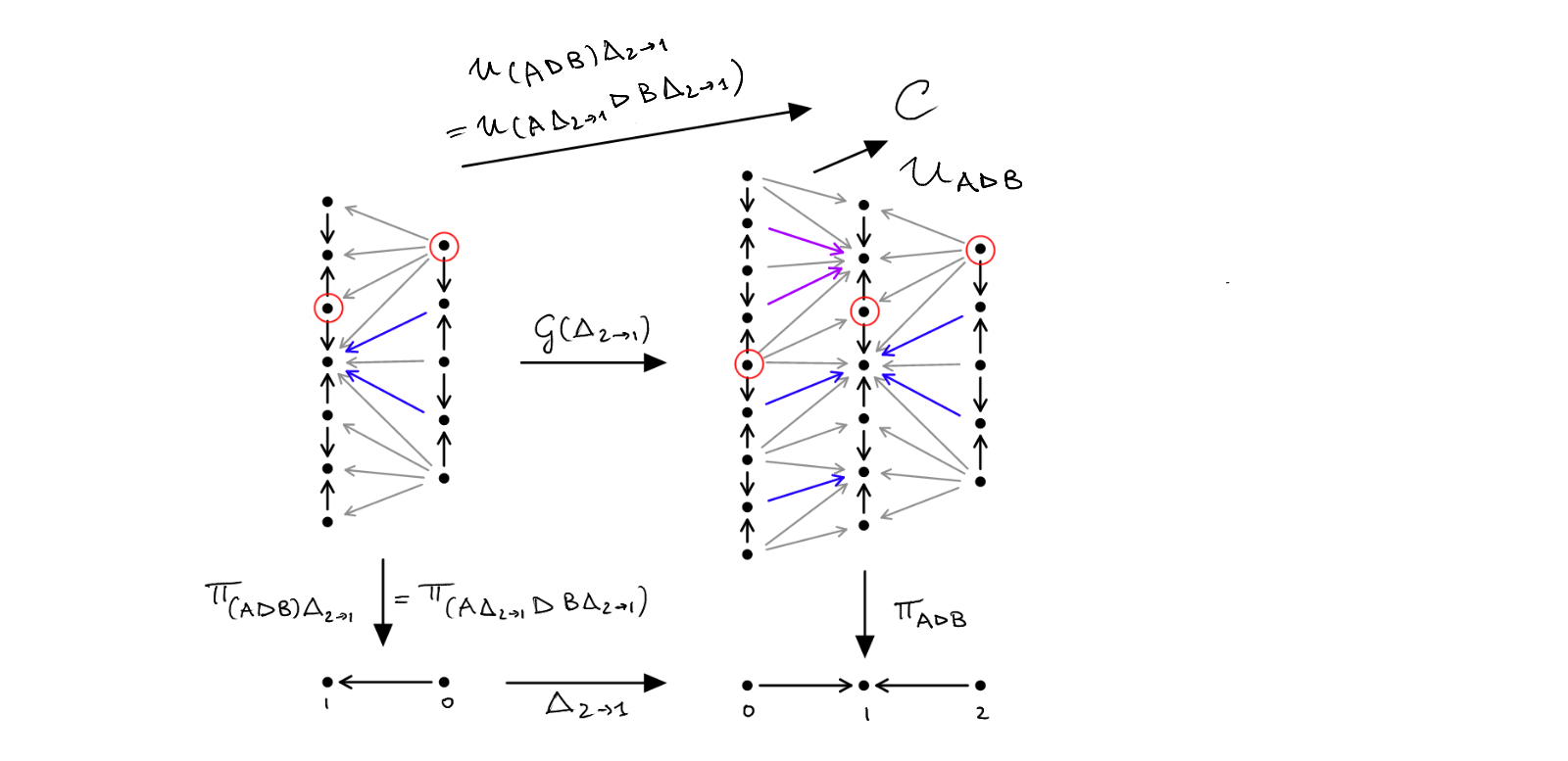}
\endgroup\end{restoretext}
In the above the gluing boundary of both families is indicated by \cred{} circles.
\end{rmk}

\subsection{Cubical sources and targets}

The following definition intuitively defines the two ``sides" of the cubes facing in the $k$th direction.

\begin{defn}[$k$-level cubical source and target] \label{constr:source_target} Let $\scA : X \to \SIvert n \cC$. For $1 \leq k \leq n$, define the $k$-level cubical source 
\begin{equation}
\csrc^k(\scA) := {\tsU {k}_\scA} \msrc_{\tusU {k-1}_\scA} : \tsG {k-1}(\scA) \to \SIvert {n-k} \cC
\end{equation}
and the $k$th cubical target
\begin{equation}
\ctgt^k(\scA) := {\tsU {k}_\scA} \mtgt_{\tusU {k-1}_\scA} : \tsG {k-1}(\scA) \to \SIvert {n-k} \cC
\end{equation} 
\noindent We usually write $\csrc^1$ and $\ctgt^1$ as $\csrc$ and $\ctgt$ respectively.
\end{defn}

\begin{eg}[Cubical sources and targets] \label{eg:sources_and_targets} \hfill
\begin{enumerate}
\item Recall the $\SIvert 2 \cC$-cube $\scA_a$ from \autoref{eg:subfamilies}. We compute $\csrc(\scA_a)$ to be the $\SIvert 1 \cC$-family
\begin{restoretext}
\begingroup\sbox0{\includegraphics{test/page1.png}}\includegraphics[clip,trim=0 {.15\ht0} 0 {.15\ht0} ,width=\textwidth]{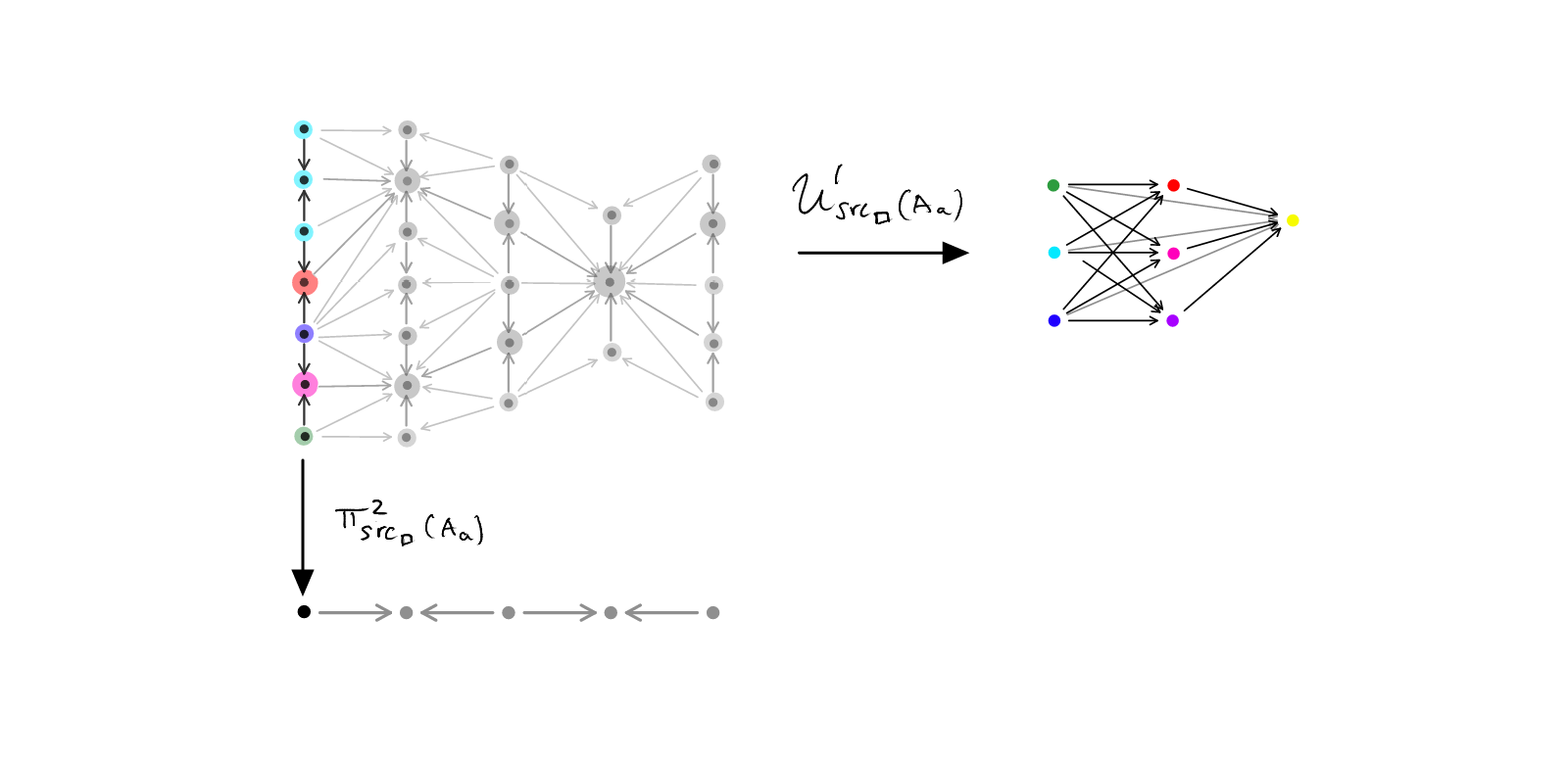}
\endgroup\end{restoretext}
and $\ctgt(\scA_a)$ to be the $\SIvert 1 \cC$-family
\begin{restoretext}
\begingroup\sbox0{\includegraphics{test/page1.png}}\includegraphics[clip,trim=0 {.2\ht0} 0 {.1\ht0} ,width=\textwidth]{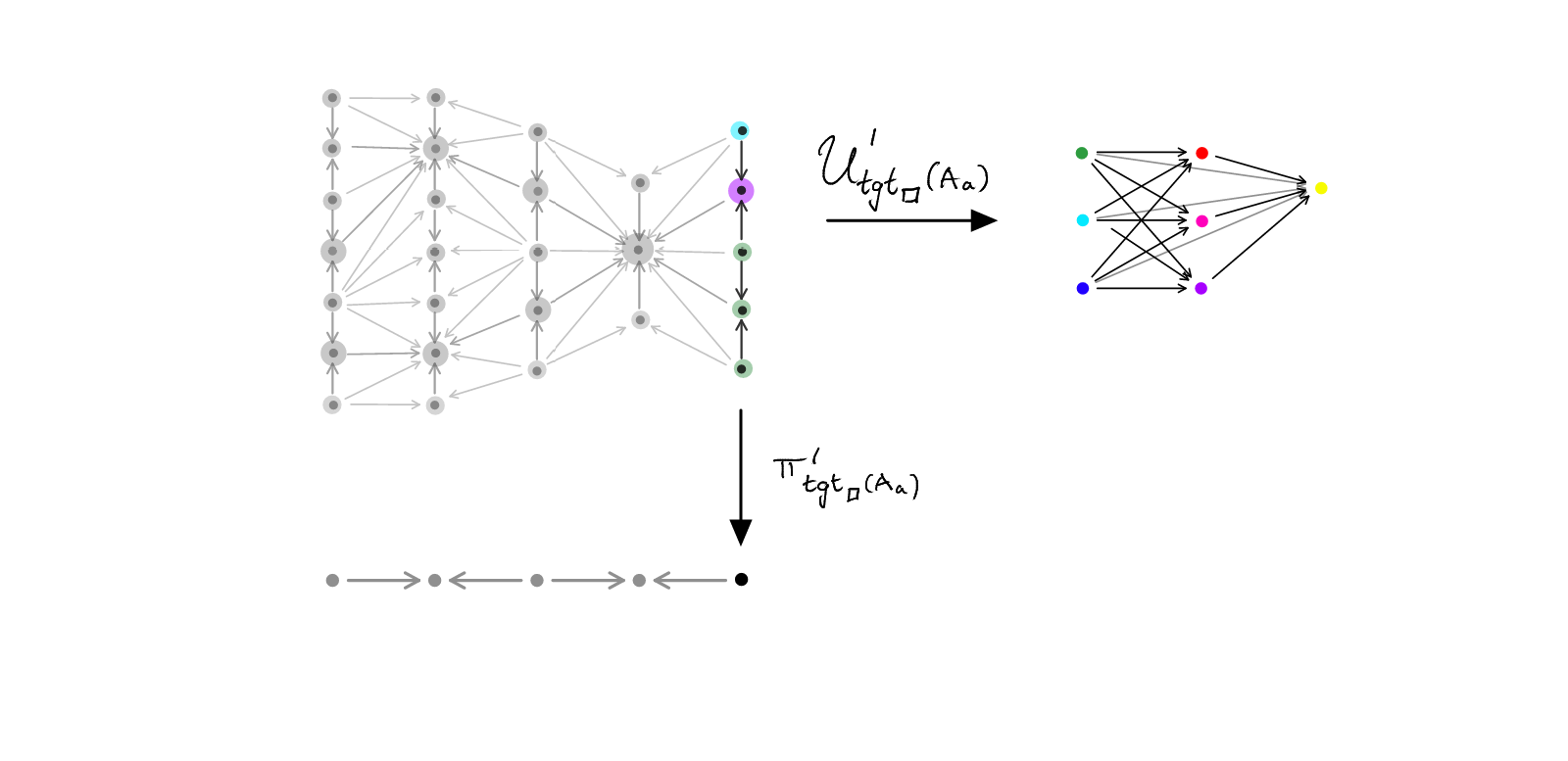}
\endgroup\end{restoretext}
We further compute $\csrc^2(\scA_a)$ to be the $\SIvert 0 \cC$-family
\begin{restoretext}
\begingroup\sbox0{\includegraphics{test/page1.png}}\includegraphics[clip,trim=0 {.3\ht0} 0 {.2\ht0} ,width=\textwidth]{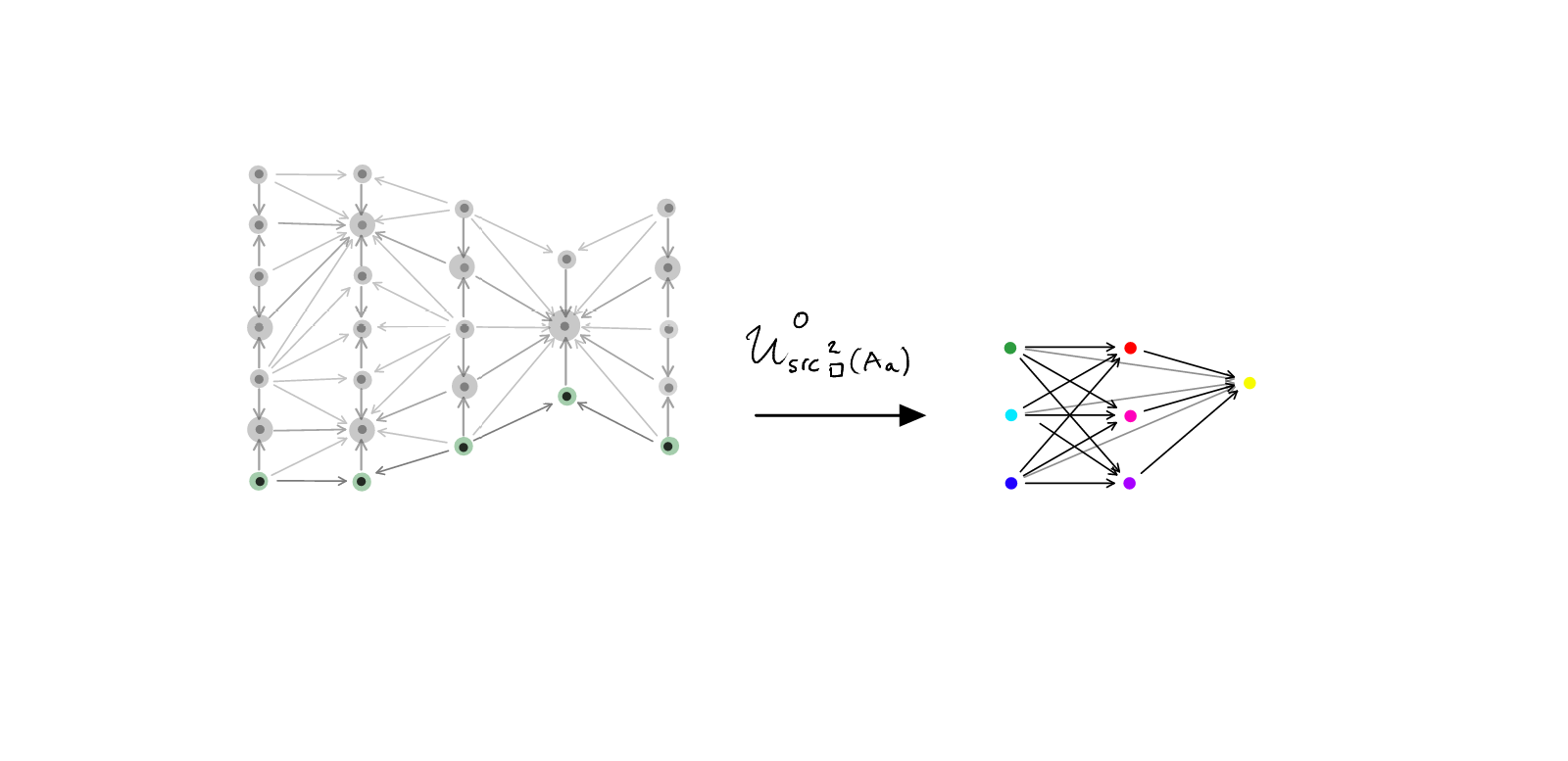}
\endgroup\end{restoretext}
and $\ctgt^2(scA_a)$ is the $\SIvert 0 \cC$-family
\begin{restoretext}
\begingroup\sbox0{\includegraphics{test/page1.png}}\includegraphics[clip,trim=0 {.2\ht0} 0 {.3\ht0} ,width=\textwidth]{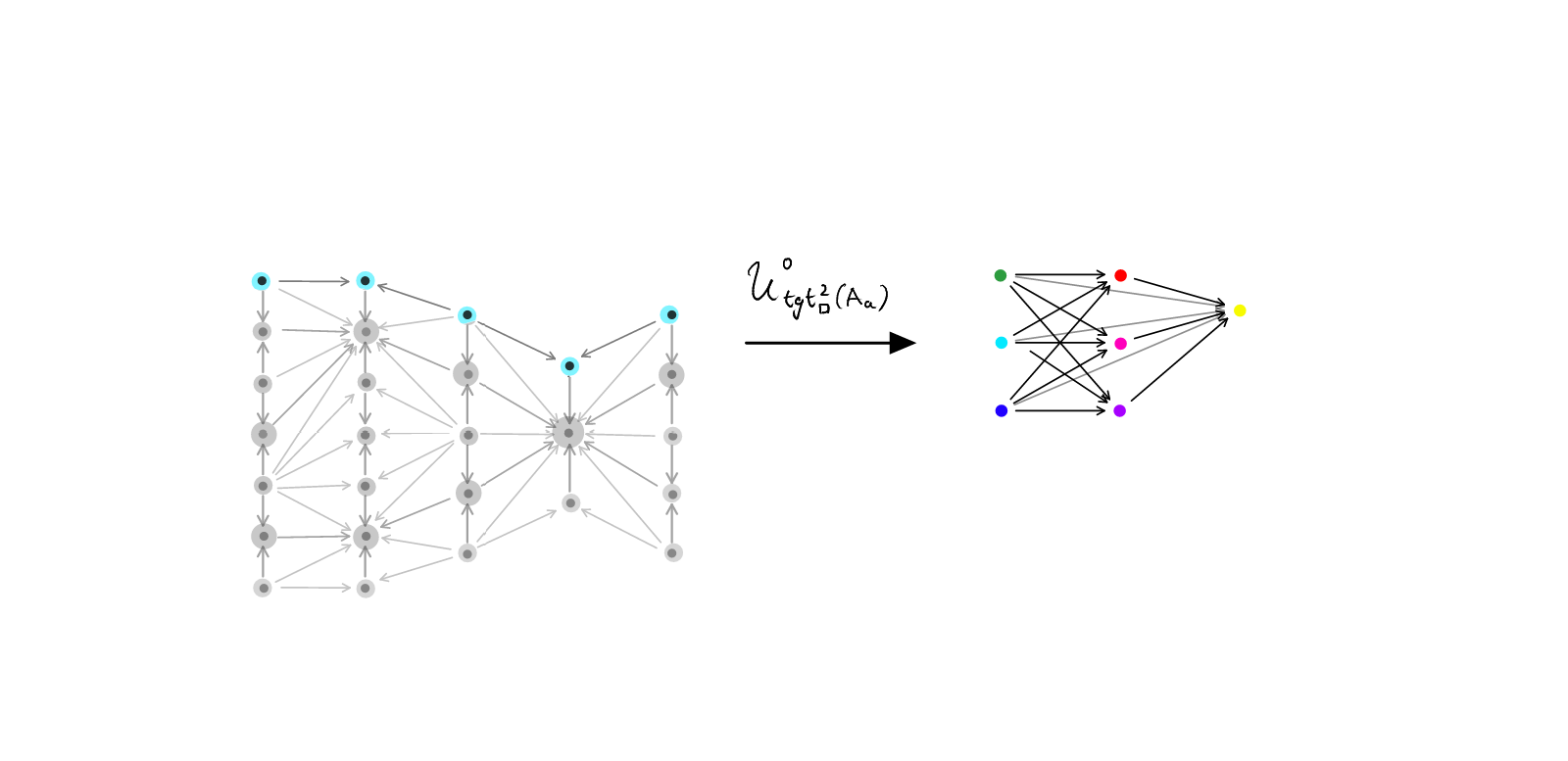}
\endgroup\end{restoretext}

\item As a second example consider the $\SIvert 3 \cC$-family $\scA_c$ defined by
\begin{restoretext}
\begin{noverticalspace}
\begingroup\sbox0{\includegraphics{test/page1.png}}\includegraphics[clip,trim=0 {.0\ht0} 0 {.0\ht0} ,width=\textwidth]{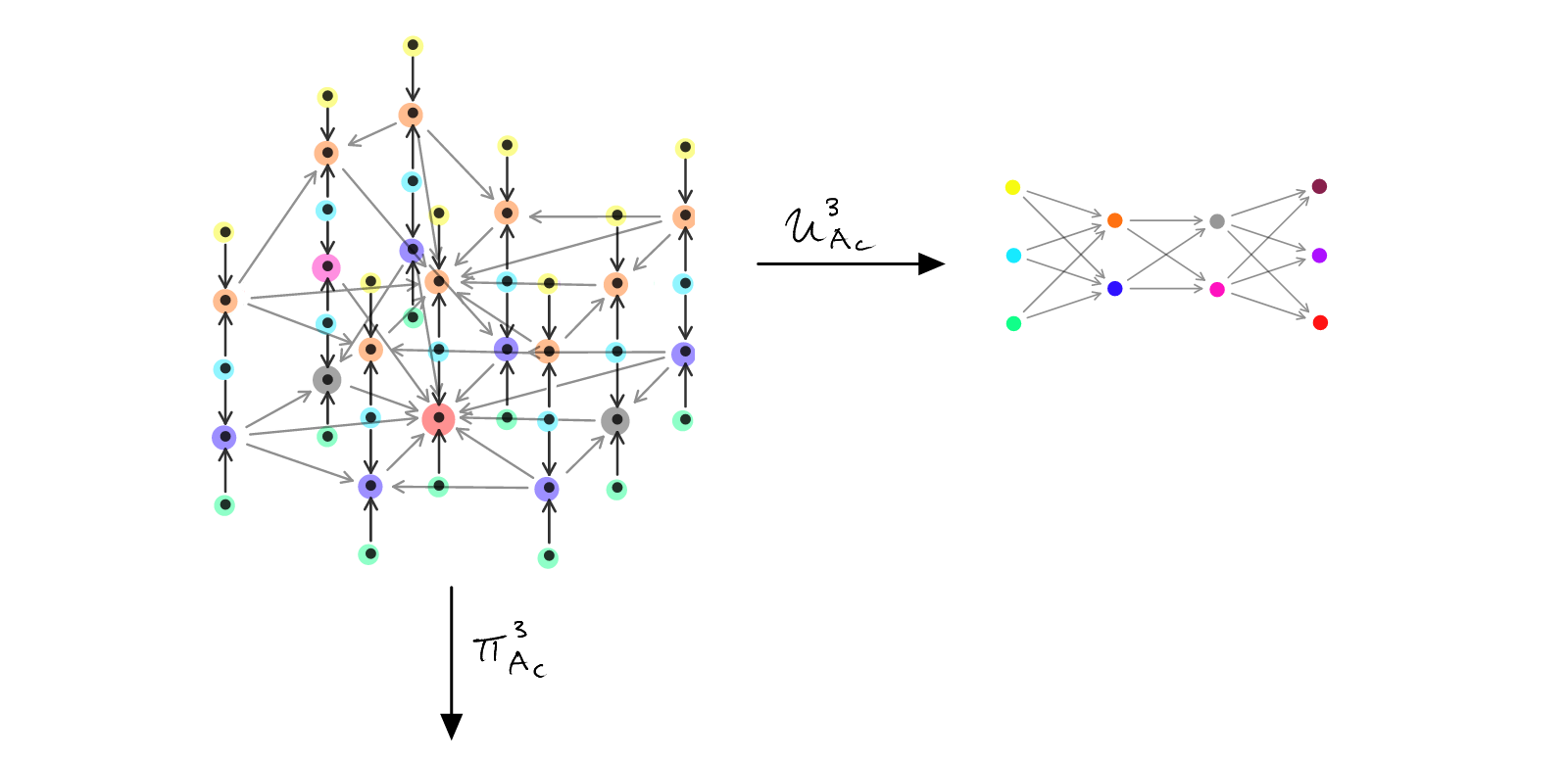}
\endgroup \\*
\begingroup\sbox0{\includegraphics{test/page1.png}}\includegraphics[clip,trim=0 {.09\ht0} 0 {.0\ht0} ,width=\textwidth]{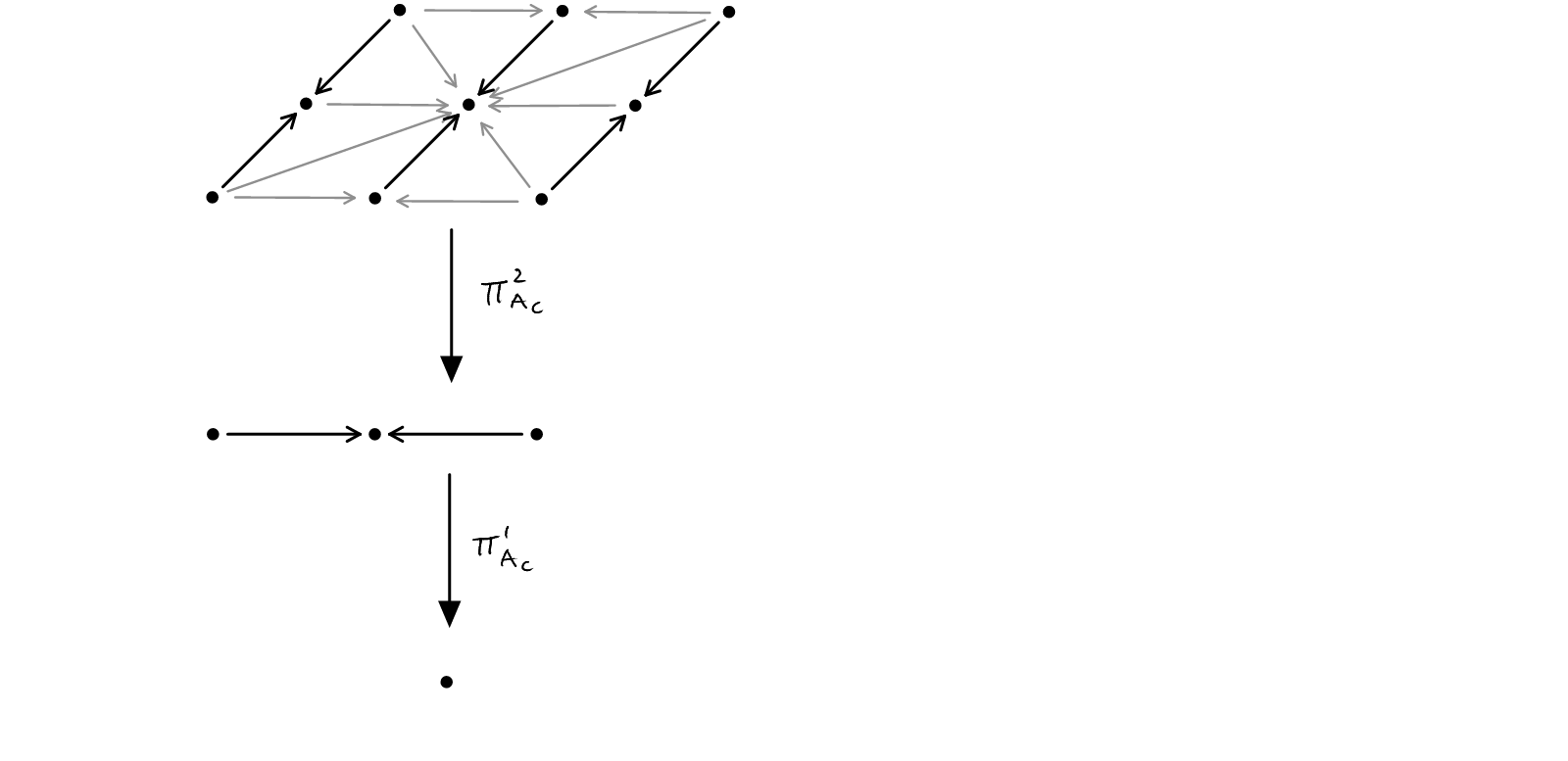}
\endgroup
\end{noverticalspace}
\end{restoretext}
We compute $\ctgt^2(\scA_c)$ to be the $\SIvert 1 \cC$-family with data
\begin{restoretext}
\begin{noverticalspace}
\begingroup\sbox0{\includegraphics{test/page1.png}}\includegraphics[clip,trim=0 {.0\ht0} 0 {.0\ht0} ,width=\textwidth]{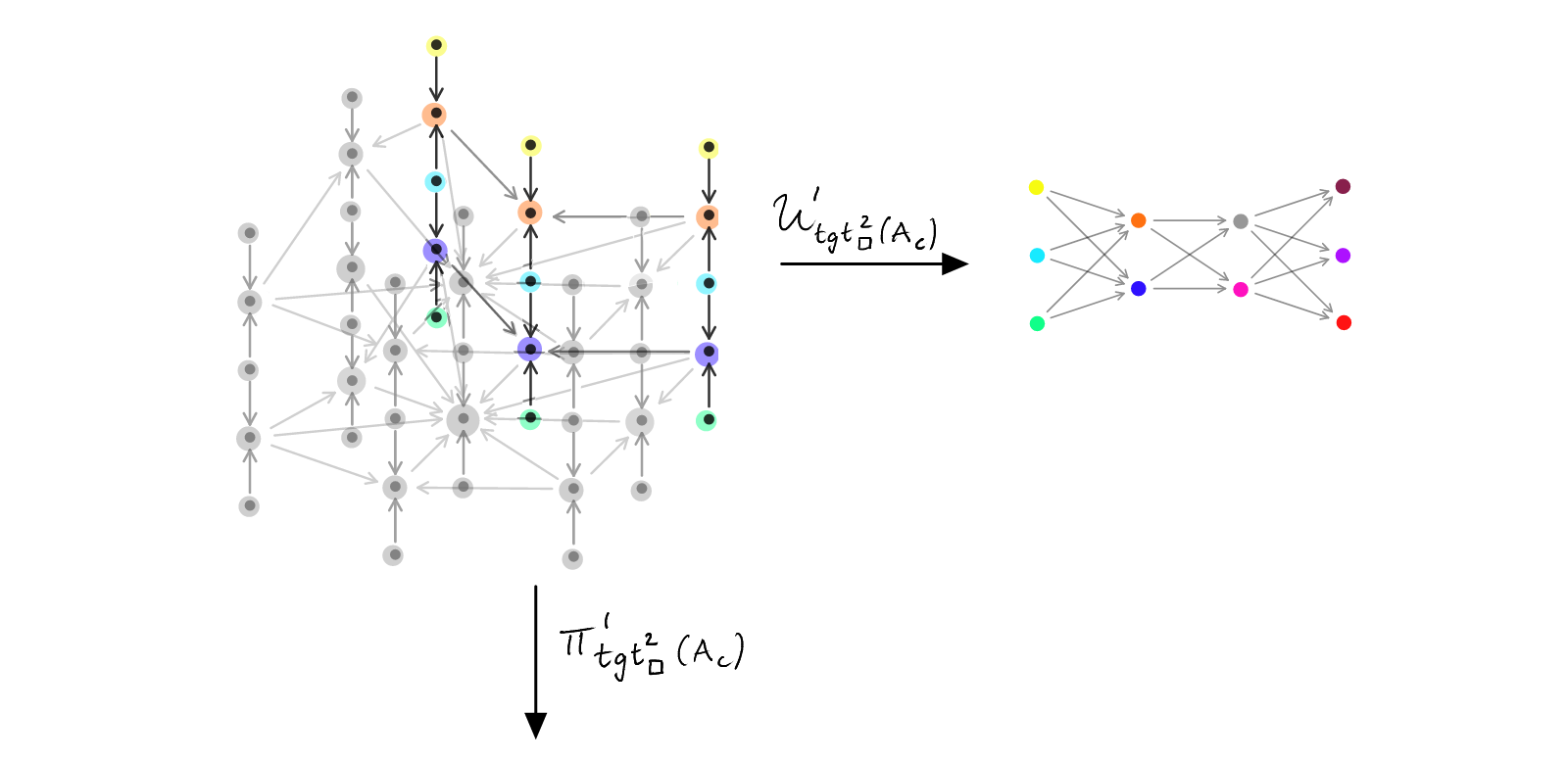}
\endgroup \\*
\begingroup\sbox0{\includegraphics{test/page1.png}}\includegraphics[clip,trim=0 {.7\ht0} 0 {.0\ht0} ,width=\textwidth]{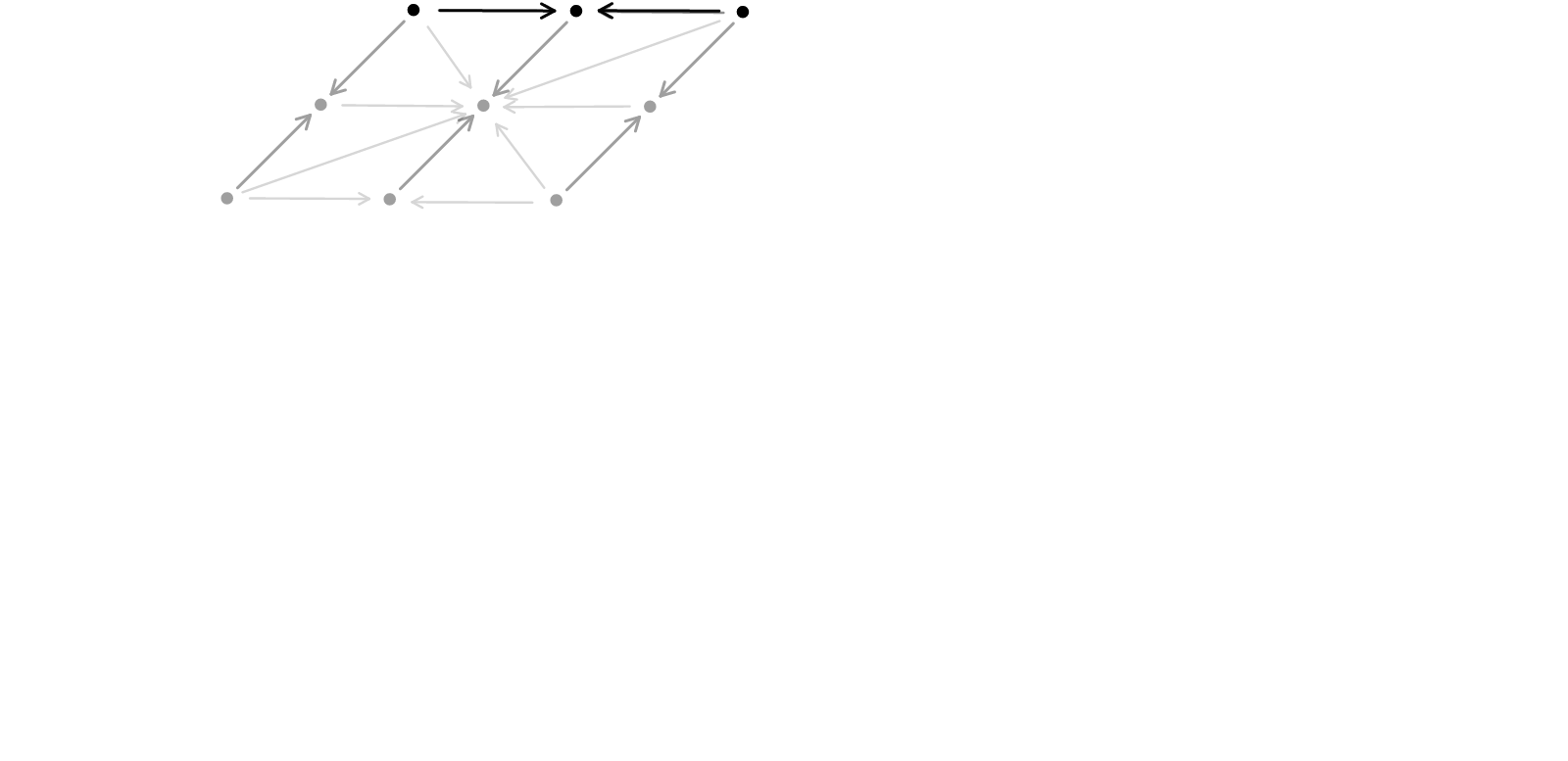}
\endgroup
\end{noverticalspace}
\end{restoretext}
and $\ctgt(\scA_c)$ is the $\SIvert 2 \cC$-family given by
\begin{restoretext}
\begin{noverticalspace}
\begingroup\sbox0{\includegraphics{test/page1.png}}\includegraphics[clip,trim=0 {.0\ht0} 0 {.0\ht0} ,width=\textwidth]{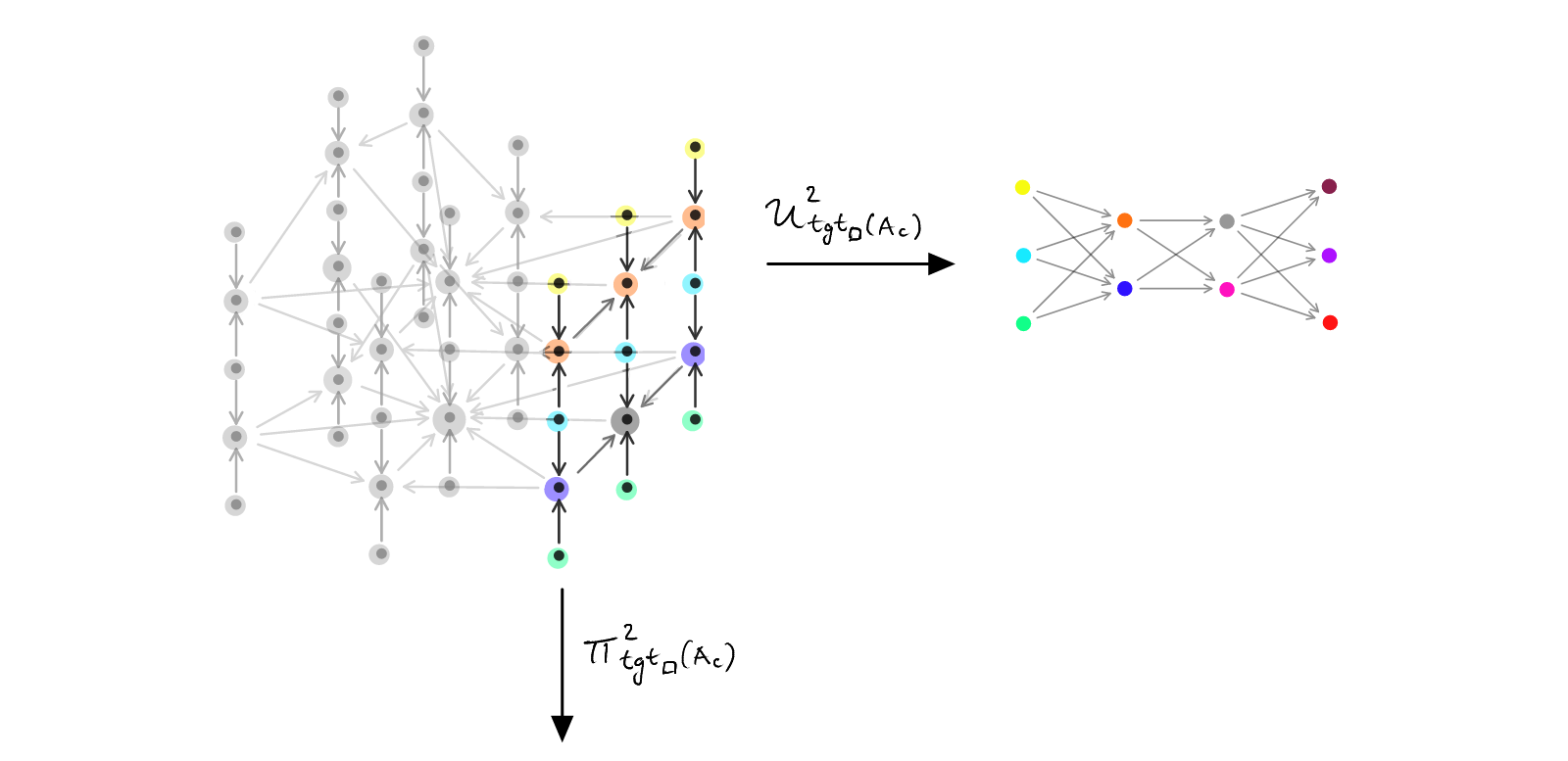}
\endgroup \\*
\begingroup\sbox0{\includegraphics{test/page1.png}}\includegraphics[clip,trim=0 {.4\ht0} 0 {.0\ht0} ,width=\textwidth]{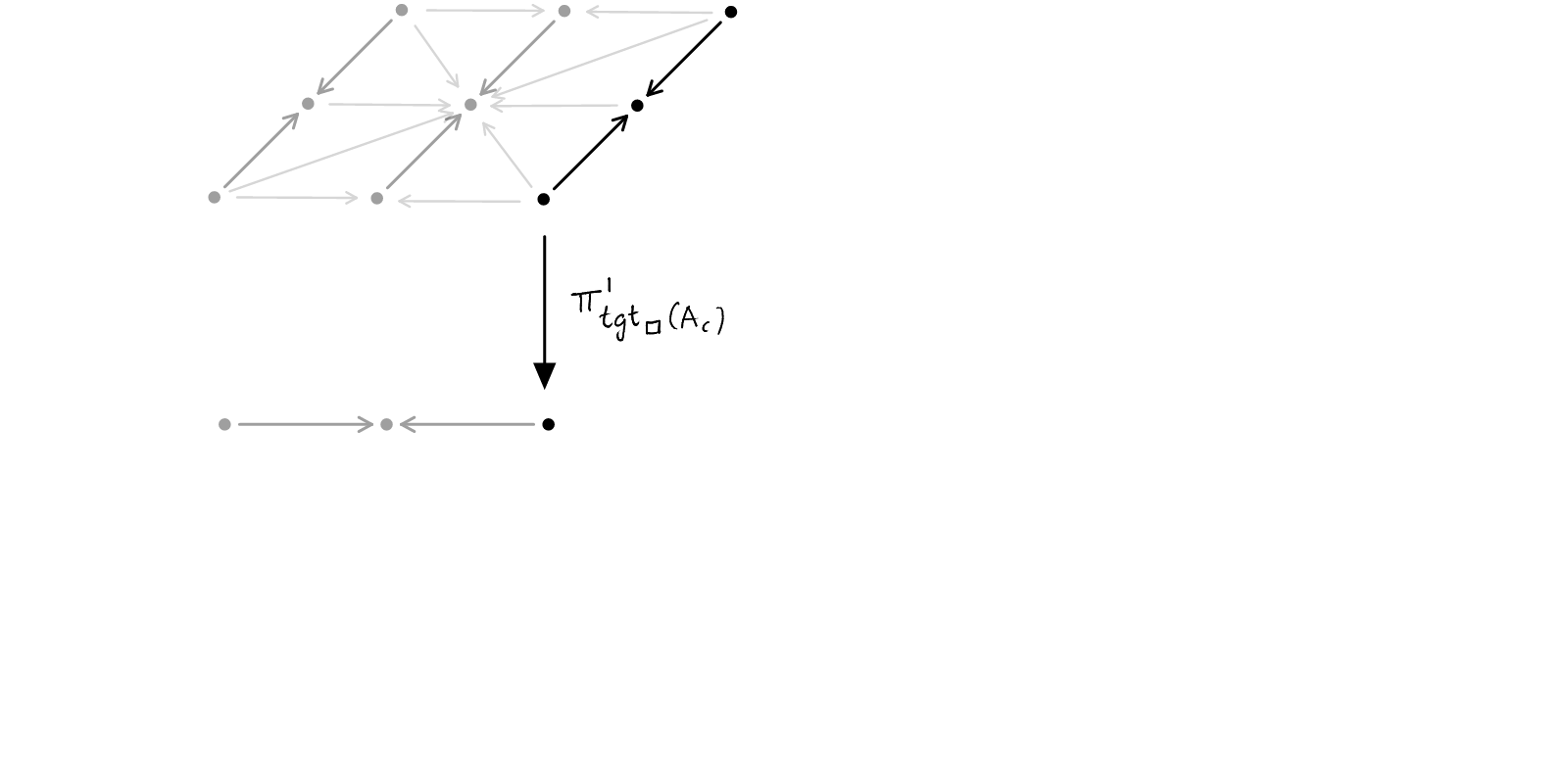}
\endgroup
\end{noverticalspace}
\end{restoretext}

\end{enumerate}

\end{eg}

\begin{rmk}[Restating the condition for stacking] Let $\scA, \scB : X \to \SIvertone \cC$. Note that $\scA \stack \scB$ as given in \autoref{defn:stacking_labelled_families} is defined if and only if $\ctgt(\scA) = \csrc(\scB)$.
\end{rmk}

\begin{rmk}[Source and target inherit connectedness]\label{rmk:connectedness_of_src_tgt}  By \autoref{rmk:connectedness} note that the domains of $\csrc^k(\scA)$ and $\ctgt^k(\scA)$ are connected if $X$ is.
\end{rmk}

\begin{rmk}[Source and target inherit normalisation] \label{lem:src_tgt_normalisation} Let $A : \bnum{1} \to \SIvert n \cC$ be globular and normalised (up to level $n$). As a corollary to \autoref{thm:normalisation_on_dc_restrictions} we find that, for $1 \leq k \leq n$, $\csrc^k(\scA)$ and $\ctgt^k(\scA)$ are normalised (up to level $(n-k+1)$).
\end{rmk}

\begin{rmk}[Source and target commute with relabelling] \label{rmk:src_tgt_label_transfer} Let $\scA : X \to \SIvert n \cC$ and $F : \cC \to \cD$. Then we compute
\begin{align}
\csrc^k (\SIvert n F \scA) &= \SIvert {n-k} F \tsU k_{\scA} \msrc_{\tusU {k-1}_\scA} \\
&= \SIvert {n-k} F \csrc (\scA)
\end{align}
and similarly
\begin{align}
\ctgt^k (\SIvert n F \scA) = \SIvert {n-k} F \ctgt^k (\scA)
\end{align}
Thus source and target commute with relabelling.
\end{rmk}

\subsection{$k$-Level stacking}

We now generalise stacking of interval families to cube families.

\begin{constr}[$k$-level stacking] \label{constr:k_level_stacking} For $1 \leq k \leq n$ we define partial binary operations $\glue {k}$ on the set of $\cC$-labelled $\SI^n$-families as follows: let $\scA, \scB : X \to \SIvert n \cC$. Then, if $\sT^{k-1}_\scA = \sT^{k-1}_\scB$ and also $\ctgt^k(\scA) = \csrc^k(\scB)$, we set
\begin{equation} \label{eq:k_level_stacking}
\scA \glue {k} \scB = \tsR {k-1}_{\sT^{k-1}_\scA,\tsU {k-1}_\scA \stack \tsU {k-1}_\scB}
\end{equation}
which is called the \textit{$k$-level stacking of $\scA$ and $\scB$}.
\end{constr}

\begin{eg}[$k$-level stacking] 
\begin{enumerate}
\item Consider the $\SIvert 3 \cC$-family $\scA_d$ defined by the data
\begin{restoretext}
\begin{noverticalspace}
\begingroup\sbox0{\includegraphics{test/page1.png}}\includegraphics[clip,trim=0 {.0\ht0} 0 {.1\ht0} ,width=\textwidth]{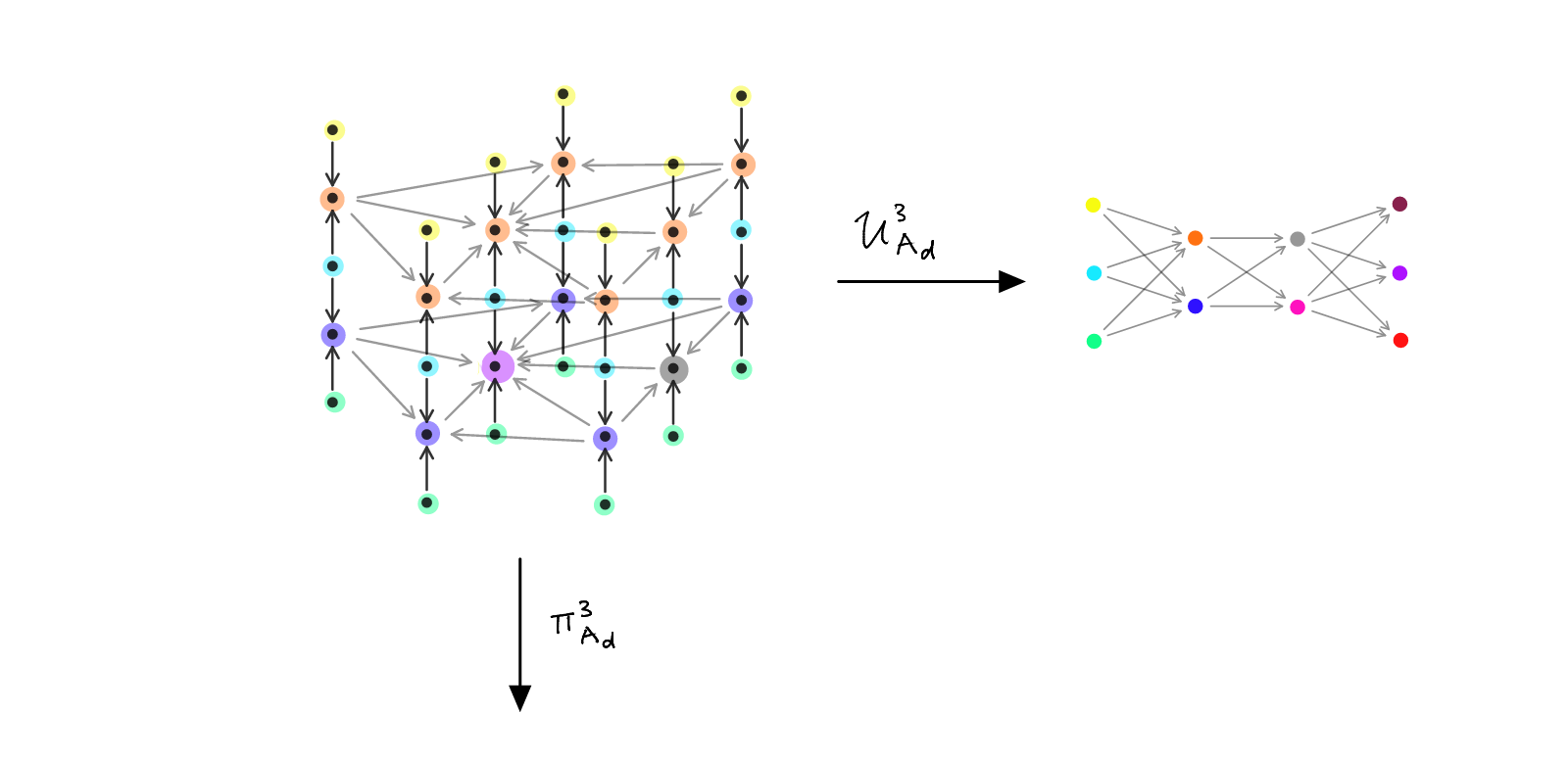}
\endgroup \\*
\begingroup\sbox0{\includegraphics{test/page1.png}}\includegraphics[clip,trim=0 {.1\ht0} 0 {.0\ht0} ,width=\textwidth]{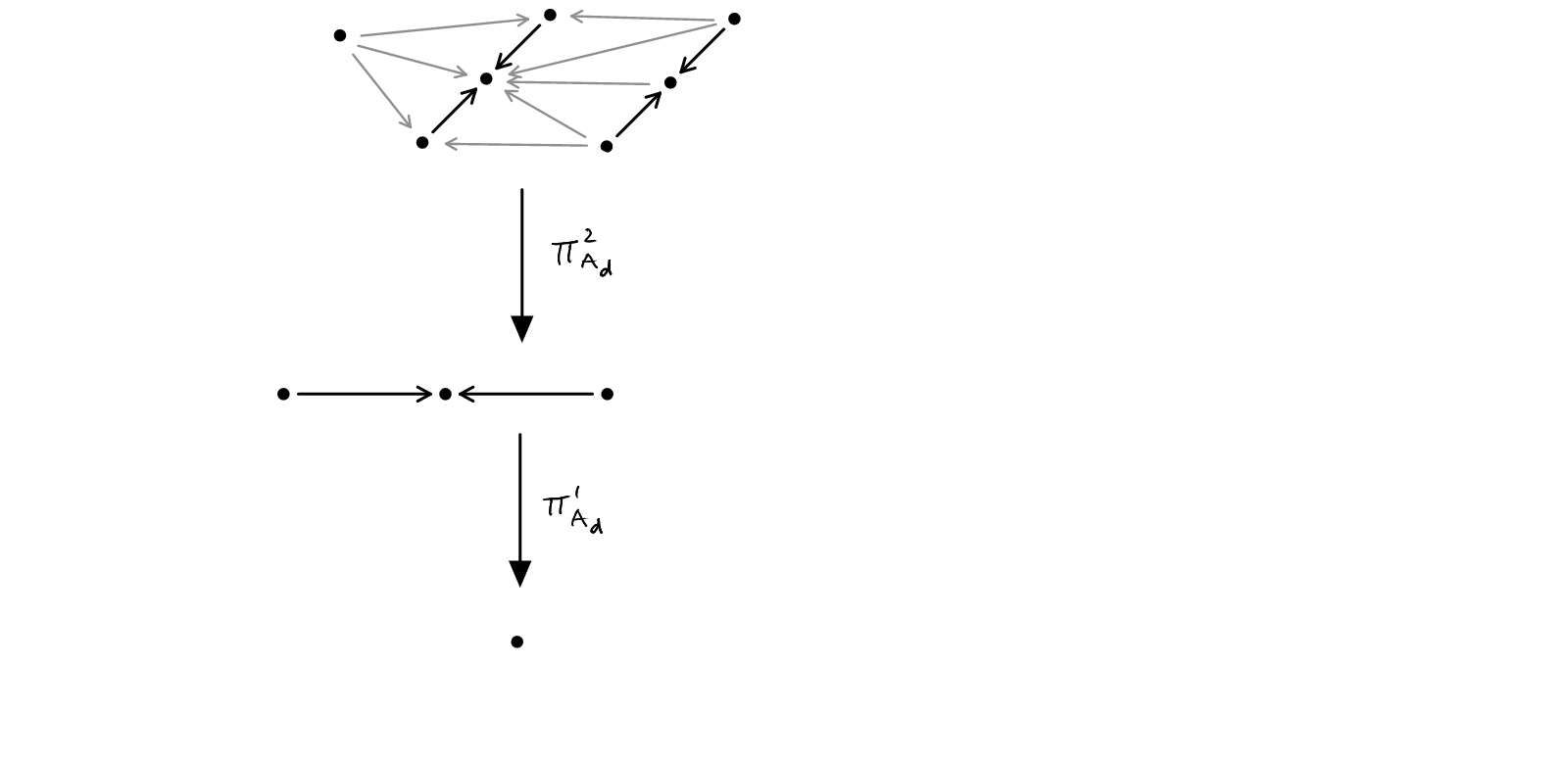}
\endgroup
\end{noverticalspace}
\end{restoretext}
Recall $\scA_c$ and $\csrc^2(\scA_c)$ from \autoref{eg:sources_and_targets}. Observe that $\sT^1_{\scA_c} = \sT^1_{\scA_d}$ (since $\tpi 1_{\scA_c} = \tpi 1_{\scA_d}$) and further that $\csrc^2(\scA_d)$ is the same family as the previously computed family $\ctgt^2(\scA_c)$. Thus $\scA_c \glue 2 \scA_d$ exists. It is the family $\scC_b$ with data
\begin{restoretext}
\begin{noverticalspace}
\begingroup\sbox0{\includegraphics{test/page1.png}}\includegraphics[clip,trim=0 {.0\ht0} 0 {.0\ht0} ,width=\textwidth]{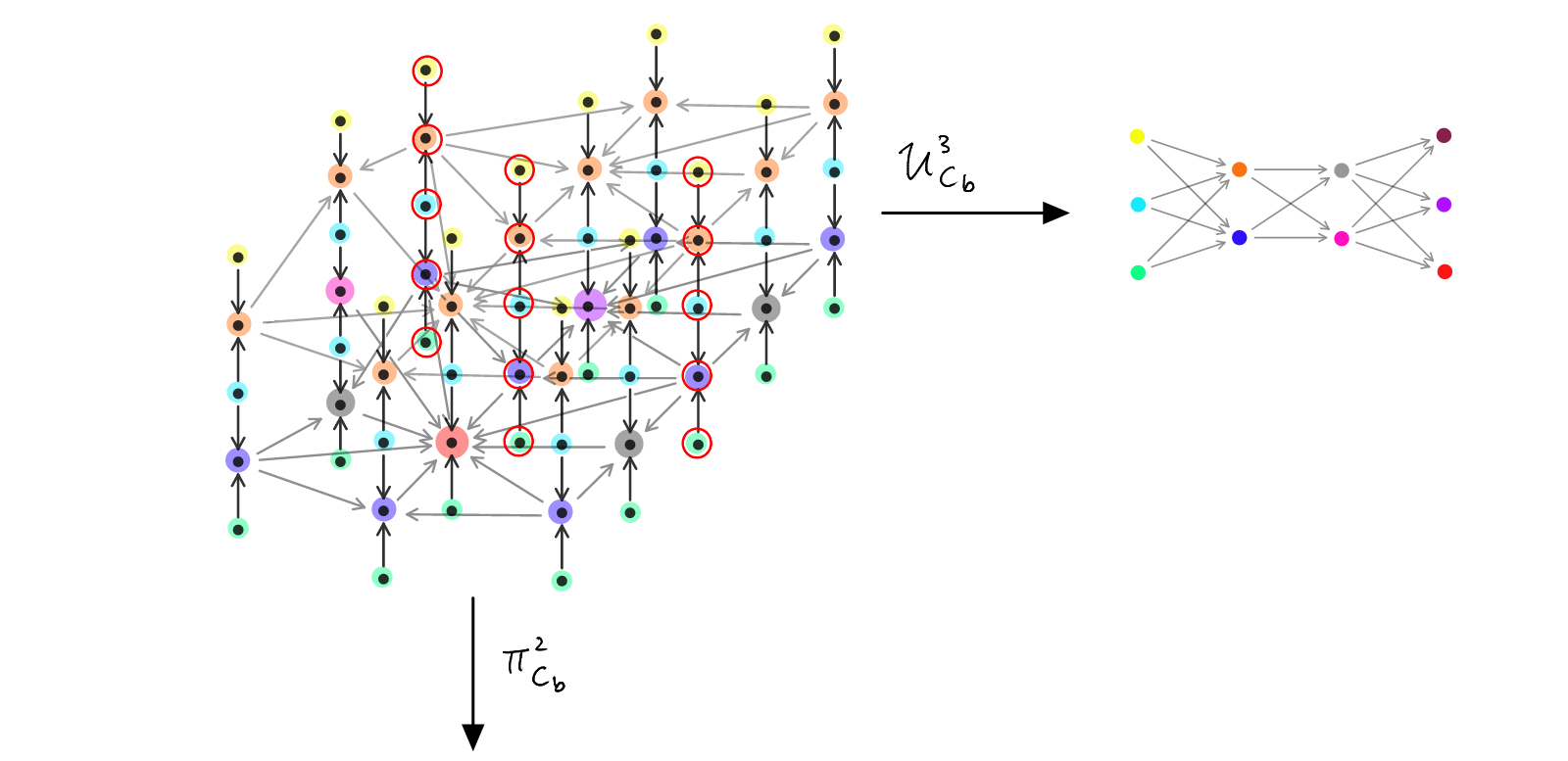}
\endgroup \\*
\begingroup\sbox0{\includegraphics{test/page1.png}}\includegraphics[clip,trim=0 {.0\ht0} 0 {.0\ht0} ,width=\textwidth]{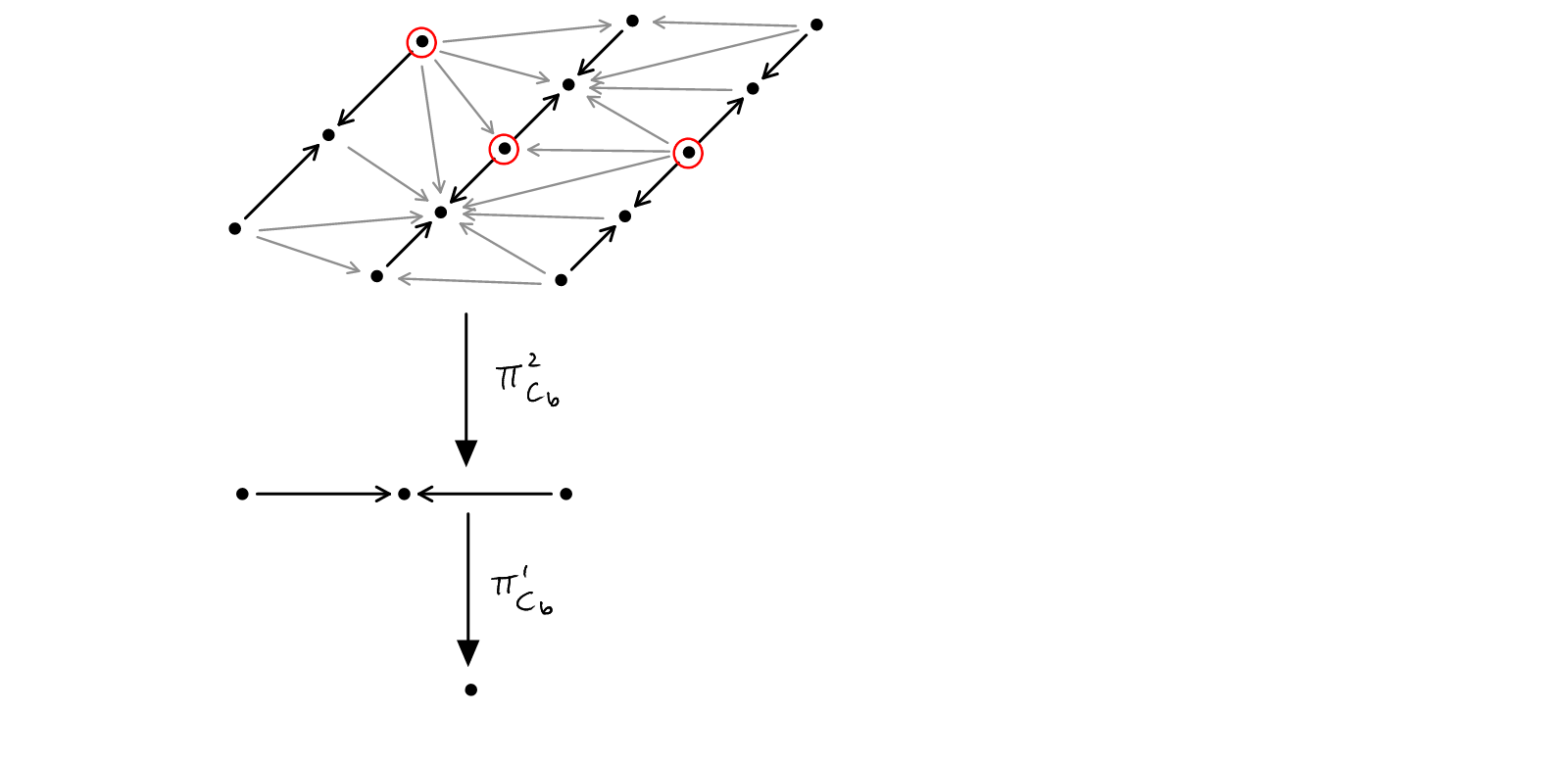}
\endgroup
\end{noverticalspace}
\end{restoretext}
where we highlighted the ``glueing" boundary (which equals both $\csrc^2(\scA_d)$ and $\ctgt^2(\scA_c)$) by \cred{} circles. Note that the stacking $\tsU 1_{\scA_c} \stack \tsU 2_{\scA_d}$ does happen at level $2$ in this case, and that $\scA_c \glue 2 \scA_d$ is then obtained from $\tsU 1_{\scA_c} \stack \tsU 2_{\scA_d}$ by adding the remaining tower of \SI-bundles below level $1$, that is
\begin{equation}
\scA_c \glue 2 \scA_d = \tsR 1_{\sT^1_{\scA_c}, \tsU 1_{\scA_c} \stack \tsU 1_{\scA_d}}
\end{equation}

\item As a second example consider the $\SIvert 3 \cC$-family $\scA_e$ defined by 
\begin{restoretext}
\begin{noverticalspace}
\begingroup\sbox0{\includegraphics{test/page1.png}}\includegraphics[clip,trim=0 {.0\ht0} 0 {.0\ht0} ,width=\textwidth]{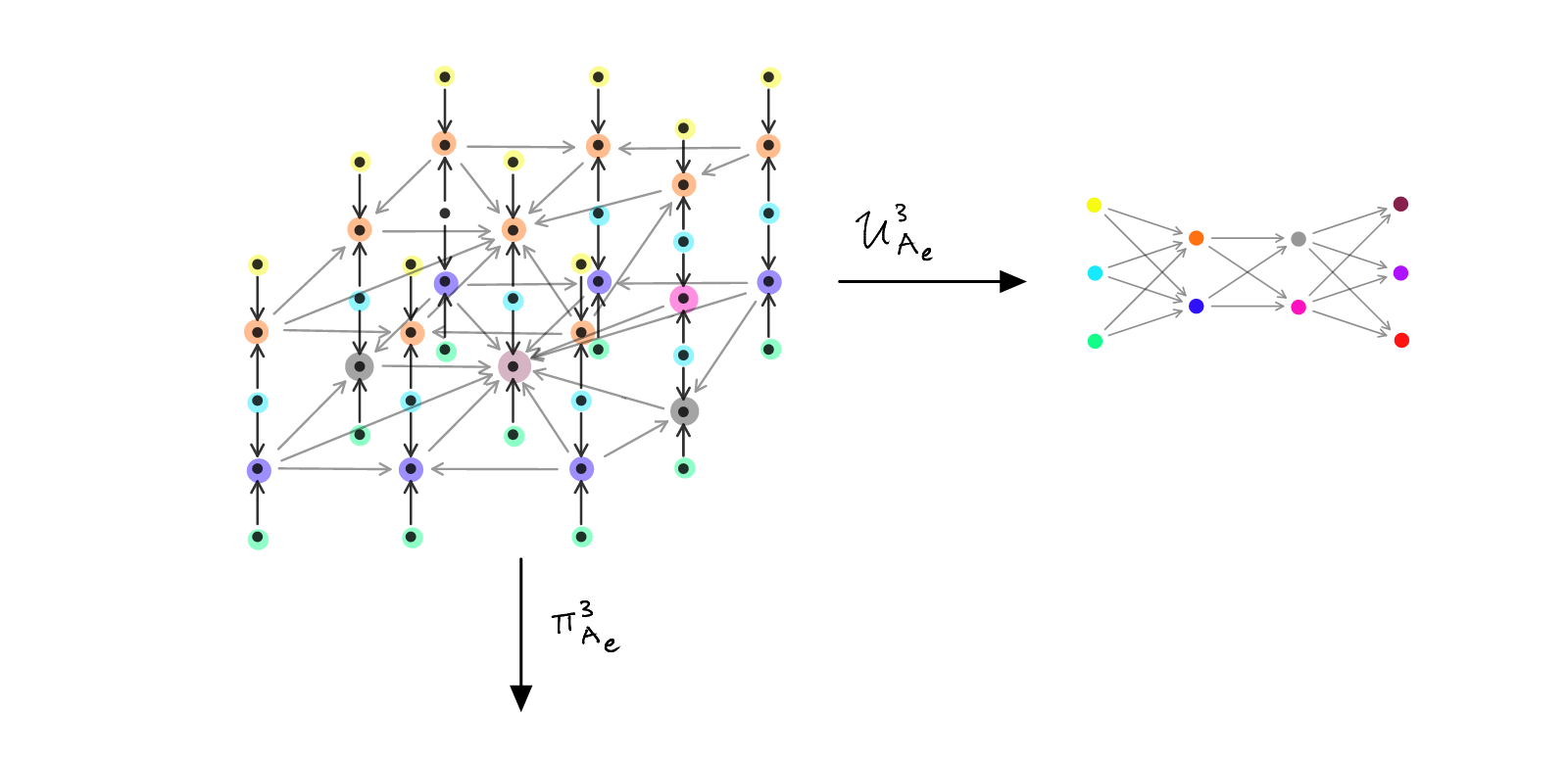}
\endgroup \\*
\begingroup\sbox0{\includegraphics{test/page1.png}}\includegraphics[clip,trim=0 {.0\ht0} 0 {.0\ht0} ,width=\textwidth]{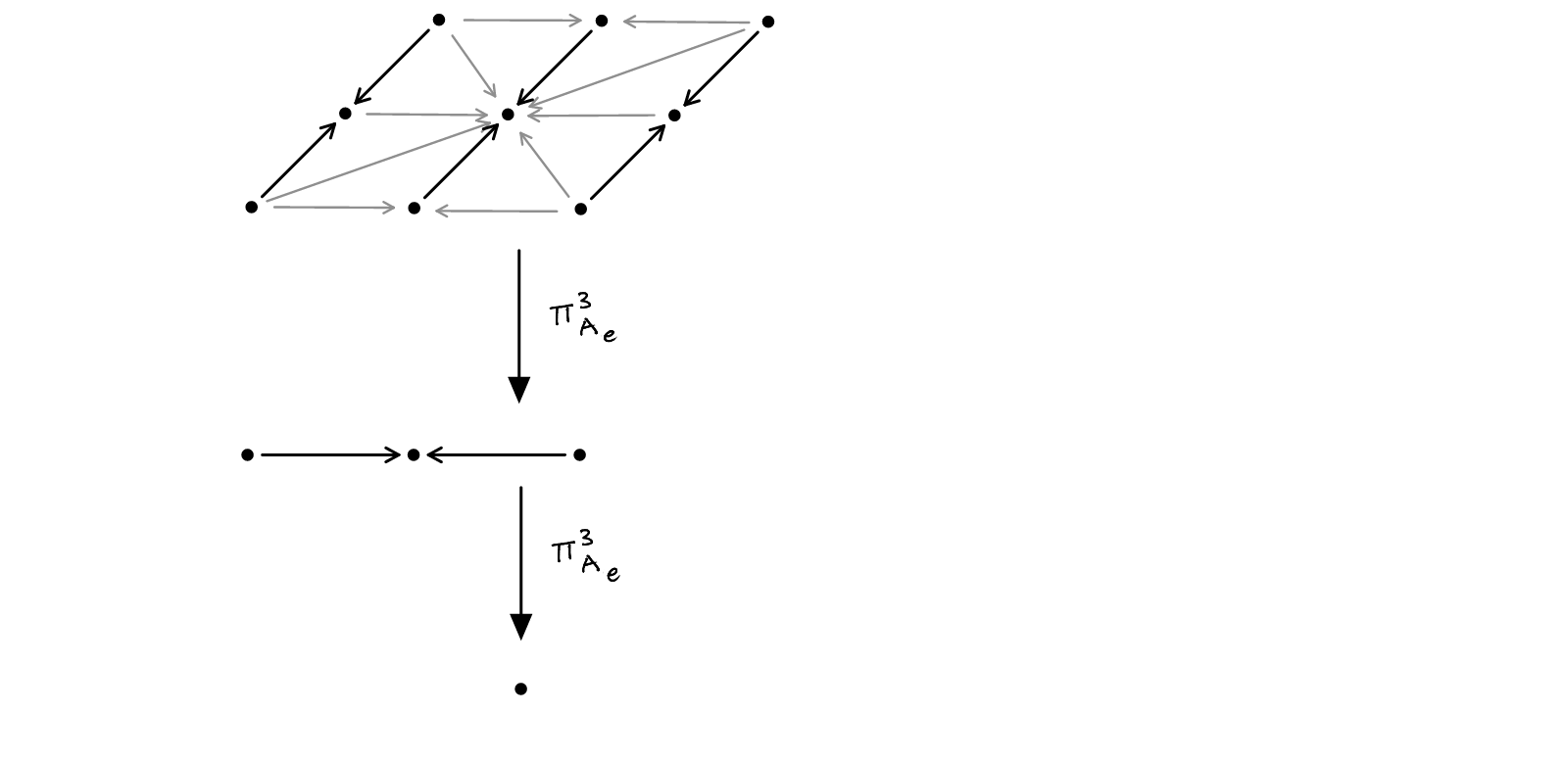}
\endgroup
\end{noverticalspace}
\end{restoretext}
We compute $\ctgt(\scA_e)$ to be the $\SIvert 2 \cC$-family
\begin{restoretext}
\begin{noverticalspace}
\begingroup\sbox0{\includegraphics{test/page1.png}}\includegraphics[clip,trim=0 {.0\ht0} 0 {.0\ht0} ,width=\textwidth]{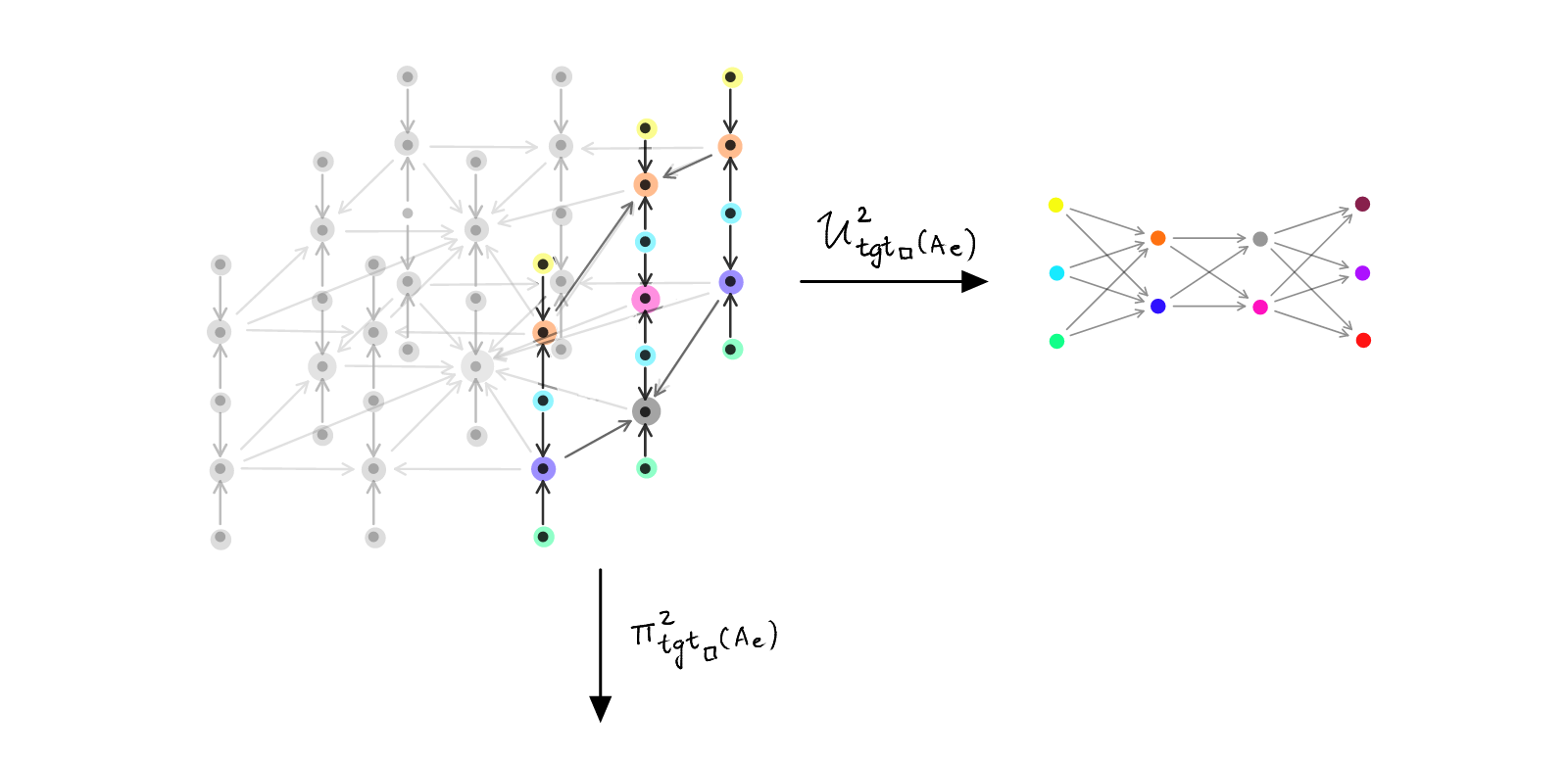}
\endgroup \\*
\begingroup\sbox0{\includegraphics{test/page1.png}}\includegraphics[clip,trim=0 {.35\ht0} 0 {.0\ht0} ,width=\textwidth]{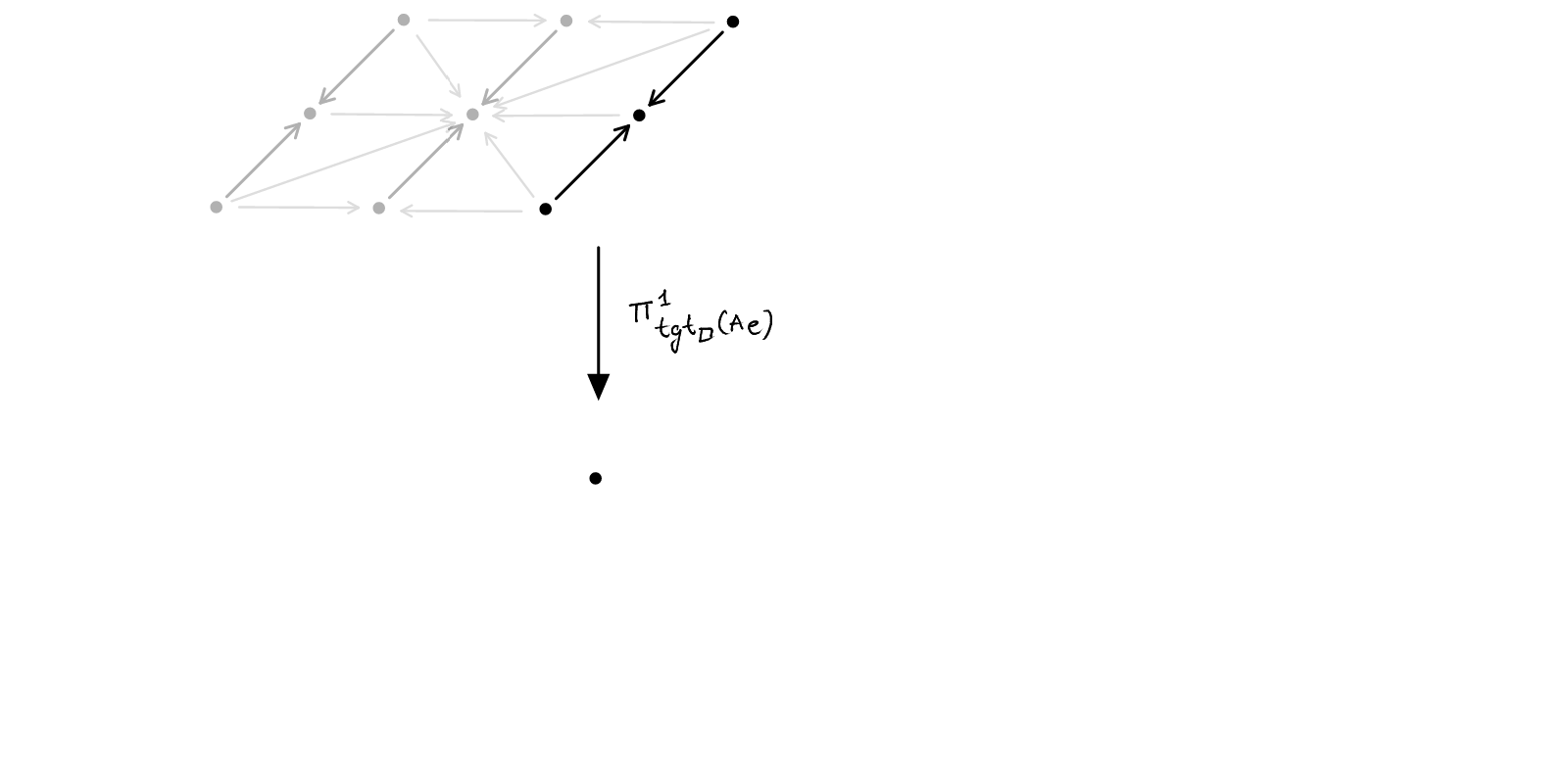}
\endgroup
\end{noverticalspace}
\end{restoretext}
This can be seen to coincide with $\csrc(\scC_b)$ from \autoref{eg:subfamilies}. Thus $\scA_e \glue 1 \scC_b$ exists, and can be computed to give $\scA_b$ from \autoref{eg:subfamilies}.
\end{enumerate}
\end{eg}

Note that in each of the above examples of $k$-level stackings $\scA \glue k \scB$, clearly $\scA$ and $\scB$ are naturally subfamilies of $\scA \glue k \scB$. More generally, we have the following.

\begin{rmk}[Canonical subfamilies of stacked cubes] \label{rmk:can_subbund_of_gluing} With the same assumptions as in the previous construction, we show that both $\scA$ and $\scB$ are canonical subfamilies of $\scA \glue k \scB$. Namely, setting $\scC := \tusU {k-1}_\scA$ and $\scD := \tusU {k-1}_\scB$ we find embeddings as follows
\begin{align}
\stackinc k_1 := \sJ^{\scA \glue k \scB,(k-1)}\restsec{[\stacklow \msrc_{\scC},\stacklow \mtgt_{\scC}]} &: \scA \subset \scA \glue k \scB \\
\stackinc k_2 := \sJ^{\scA \glue k \scB,(k-1)}\restsec{[\stackup \msrc_{\scD},\stackup \mtgt_{\scD}]} &: \scB \subset \scA \glue k \scB
\end{align}
The claim that these \textit{stacking embeddings} obtained from $(k-1)$-level endpoints indeed have domain $\scA$ resp. $\scB$ follows from \autoref{rmk:canonical_total_space_inclusion} and  \eqref{eq:subbund_of_stack} after unpacking the definitions in \autoref{constr:k_level_stacking} and \autoref{constr:subfamilies_from_endpoints}. Note that the stacking $\scA \glue k \scB$ is implicit in the notation for stacking embeddings (the index $i \in \Set{1,2}$ reflects whether the embedding refers to the first or second component of the stacking),
\end{rmk}

\subsection{Stacking preserves normalisation}

An important property of stacking cube families is the following.

\begin{lem}[$k$-Level stacking preserves normalisation] \label{lem:stacking_preserves_nomalisation} Assume $\scA, \scB : X \to \SIvert n \cC$ are normalised, satisfy $\sT^{k-1}_\scA = \sT^{k-1}_\scB$ and also $\ctgt^k(\scA) = \csrc^k(\scB)$. Then $\scA \glue {k} \scB$ is normalised.
\proof The proof is \stfwd{}. Note that the embeddings defined in \autoref{rmk:can_subbund_of_gluing} satisfy for $0 \leq l \leq n$ that
\begin{equation}
(\stackinc k_1)^l (\tsG l(\scA)) \cup (\stackinc k_2)^l (\tsG l(\scB)) = \tsG l(\scA \glue k \scB)
\end{equation}
If $\scA\glue k \scB$ is not normalised then we can choose a non-identity $\lambda : (\scA\glue k \scB) \kcoll l \scC$ for some $0 < l \leq n$ and $\scC : X \to \SIvert n k$. Since $\lambda$ is a non-identity injection $\lambda : \tsU {l-1}_\scC \into \tsU {l-1}_{\scA \glue k \scB}$ we find $x \in \tsG {l-1}(\scC) = \tsG {l-1}(\scA \glue k \scB)$ and $y \in \singcont(\tusU {l-1}_{\scA \glue k \scB}(x))$ such that $y \notin \im(\lambda_x)$. By \autoref{thm:collapse_maps_vs_injections}, this is equivalent to
\begin{equation}
\sS^\lambda (x,y) \in \regcont(\tsG l(\scC))
\end{equation}
But by our first observation
\begin{equation}
(x,y) \in \tsG l(\scA \glue k \scB) \quad \imp \quad (x,y) = (\stackinc k_1)^l (x_1,y_1) \text{~or~} (x,y) = (\stackinc k_2)^l (x_2,y_2)
\end{equation}
for some $x_i, y_i$. We assume the former case (the argument in the latter case is similar). By \autoref{constr:collapse_on_subfamilies} we find $(\stackinc k_1)\pbstar  \lambda : \scA \kcoll l \scD$. But now the pullback \eqref{eq:subfamily_and_collapse} defining $\sS^{(\stackinc k_1)\pbstar \lambda}$ implies that
\begin{equation}
\sS^{(\stackinc k_1)\pbstar  \lambda} (x_1,y_1) \in \regcont(\tsG k(\scD))
\end{equation}
which in turn by \autoref{thm:collapse_maps_vs_injections} implies that $(\stackinc k_1)\pbstar  \lambda$ is a non-identity collapse (as it has a non-identity component at $x_1$). This contradicts $\scA$ being in normal form. \qed
\end{lem}

\subsection{Whiskering} \label{ssec:glob_comp}

Having defined $k$-level stacking of cubes, we are now in the position to define a more specialised operation which only applies to normalised globular cubes. This models the usual idea of ``whiskering" in strict $n$-categories (cf. \cite{leinster-operads}) in our cubical setting. 

To begin, recall the definition of globular sources and targets from the previous chapter. We remark the following notational confluence.

\begin{rmk}[Globular source and target] \label{constr:globular_src_tgt} Let $\scA : \bnum{1} \to \SIvert n \cC$ be globular, $n > 0$. In this case, the \textit{globular source} $\gsrc(\scA)$ and \textit{target} $\gtgt(\scA)$ satisfy
\begin{align}
\gsrc(\scA) &=  \csrc(\scA)\\
\gtgt(\scA) &=  \ctgt(\scA)
\end{align}
\end{rmk}

The following example serves to demonstrate an important fact.

\begin{eg}[``Globularity" of globular families] Recall the $\SIvert 2 \cC$-family $\scA_a$ and its sources and targets from \autoref{eg:sources_and_targets}. This is a globular family as previously remarked. In this case, we make the following observation. Taking the globular source and target of the globular target of $\scA_a$ we end up with two $\SIvert 0 \cC$-families as follows
\begin{restoretext}
\begingroup\sbox0{\includegraphics{test/page1.png}}\includegraphics[clip,trim=0 {.2\ht0} 0 {.01\ht0} ,width=\textwidth]{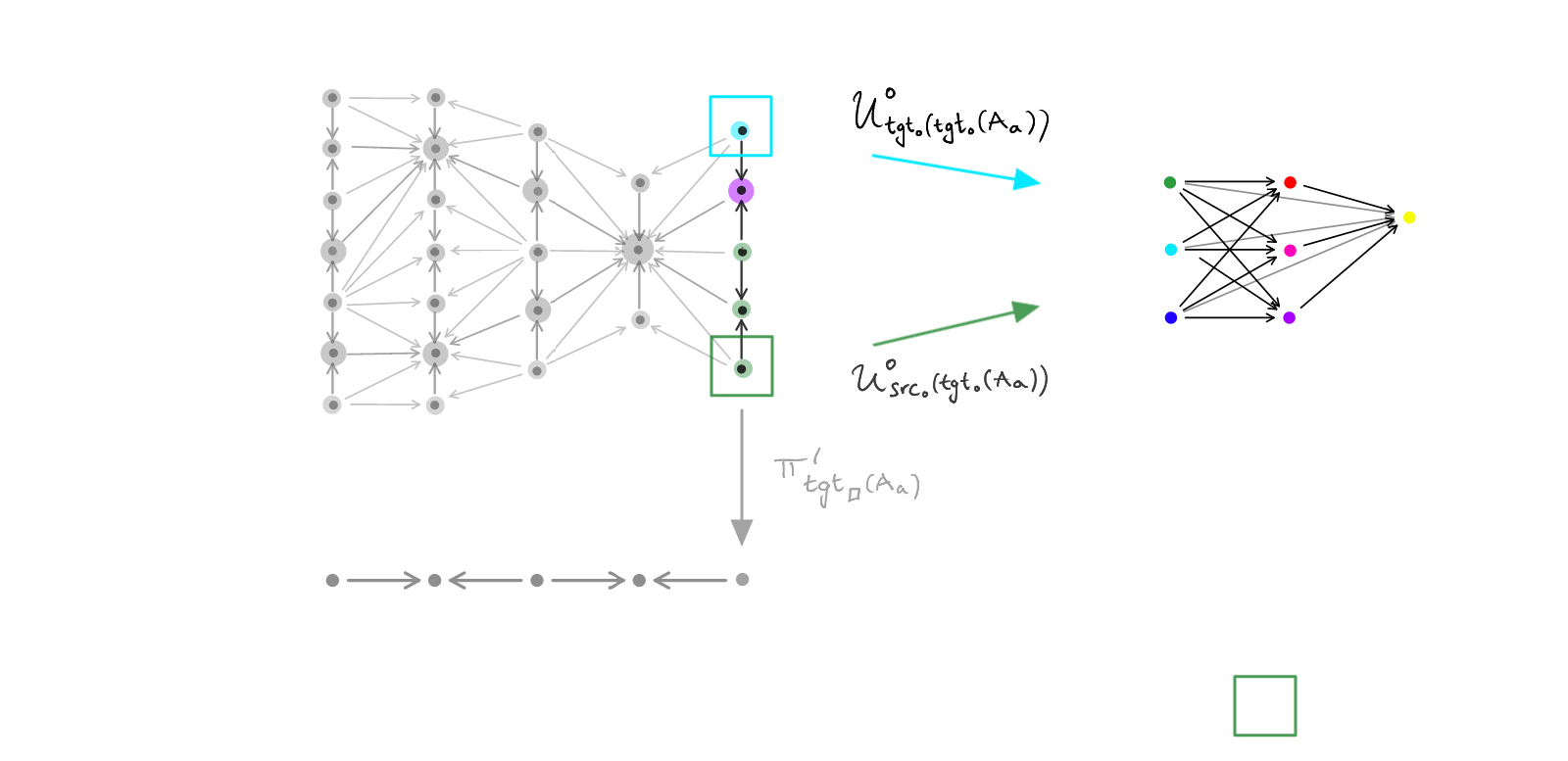}
\endgroup\end{restoretext}
Similarly, taking the globular source and target of the globular source of $\scA_a$ we end up with the $\SIvert 0 \cC$-families 
\begin{restoretext}
\begingroup\sbox0{\includegraphics{test/page1.png}}\includegraphics[clip,trim=0 {.1\ht0} 0 {.1\ht0} ,width=\textwidth]{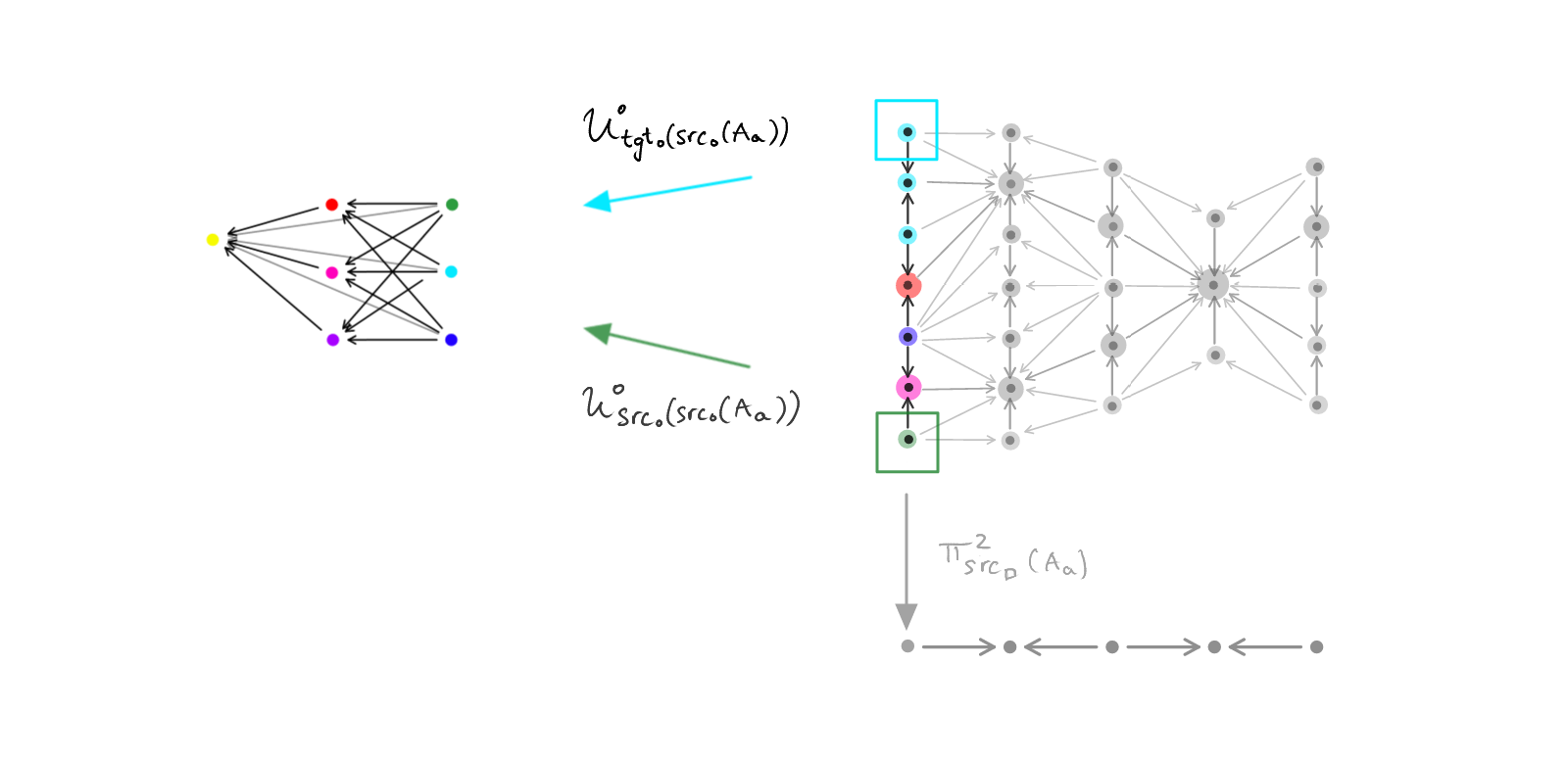}
\endgroup\end{restoretext}
Thus the source of the source equals the source of the target, and the target of the source equals the target of the target. In the first case, this is because globularity requires the entire $\csrc^2(\scA_a)$ to be (locally) trivial, as illustrated below
\begin{restoretext}
\begingroup\sbox0{\includegraphics{test/page1.png}}\includegraphics[clip,trim=0 {.3\ht0} 0 {.2\ht0} ,width=\textwidth]{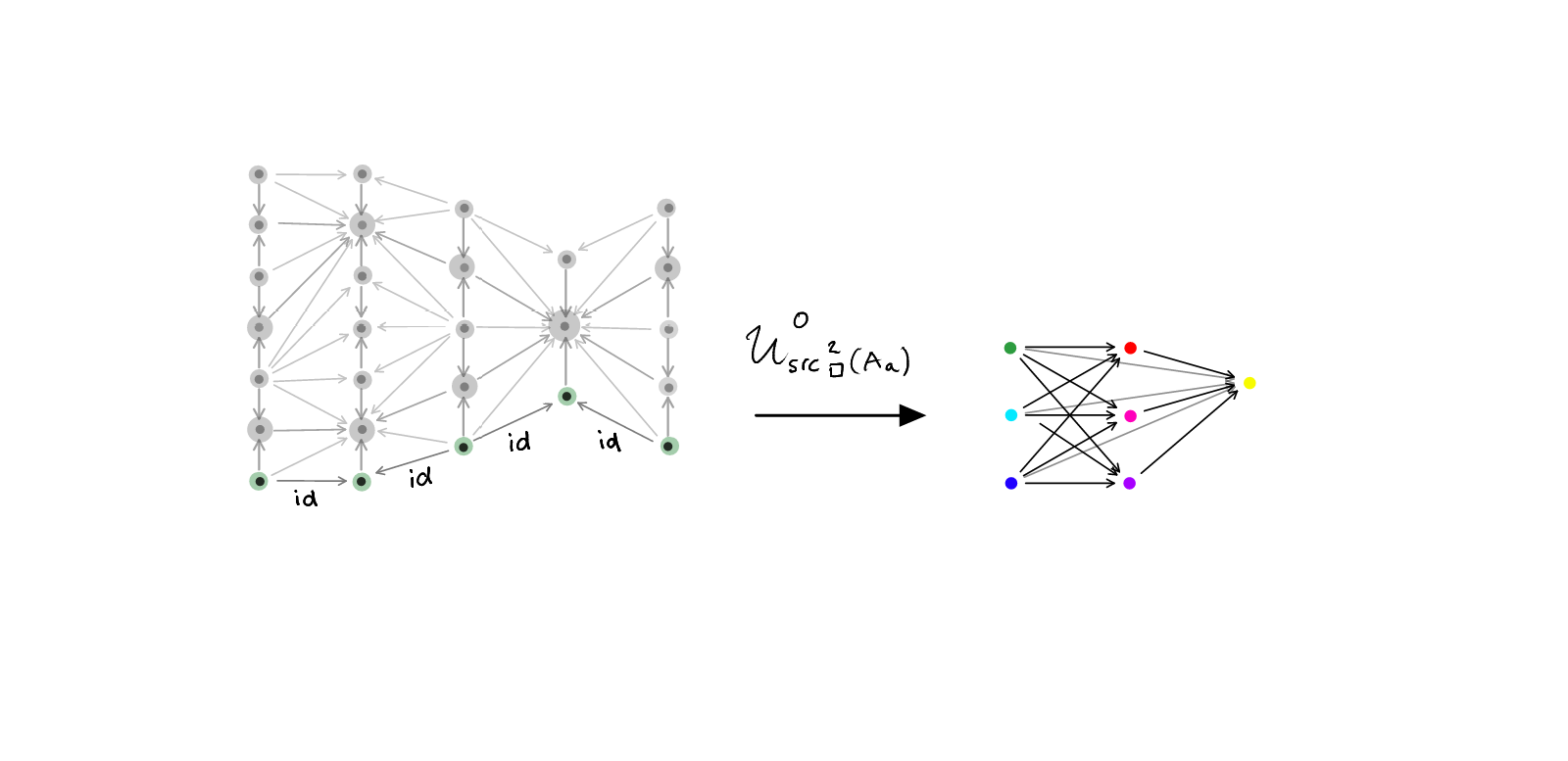}
\endgroup\end{restoretext}
Here we marked morphisms with $\id$ that have to be labelled by the identity due to globularity. Similarly, in the second case we find
\begin{restoretext}
\begingroup\sbox0{\includegraphics{test/page1.png}}\includegraphics[clip,trim=0 {.2\ht0} 0 {.3\ht0} ,width=\textwidth]{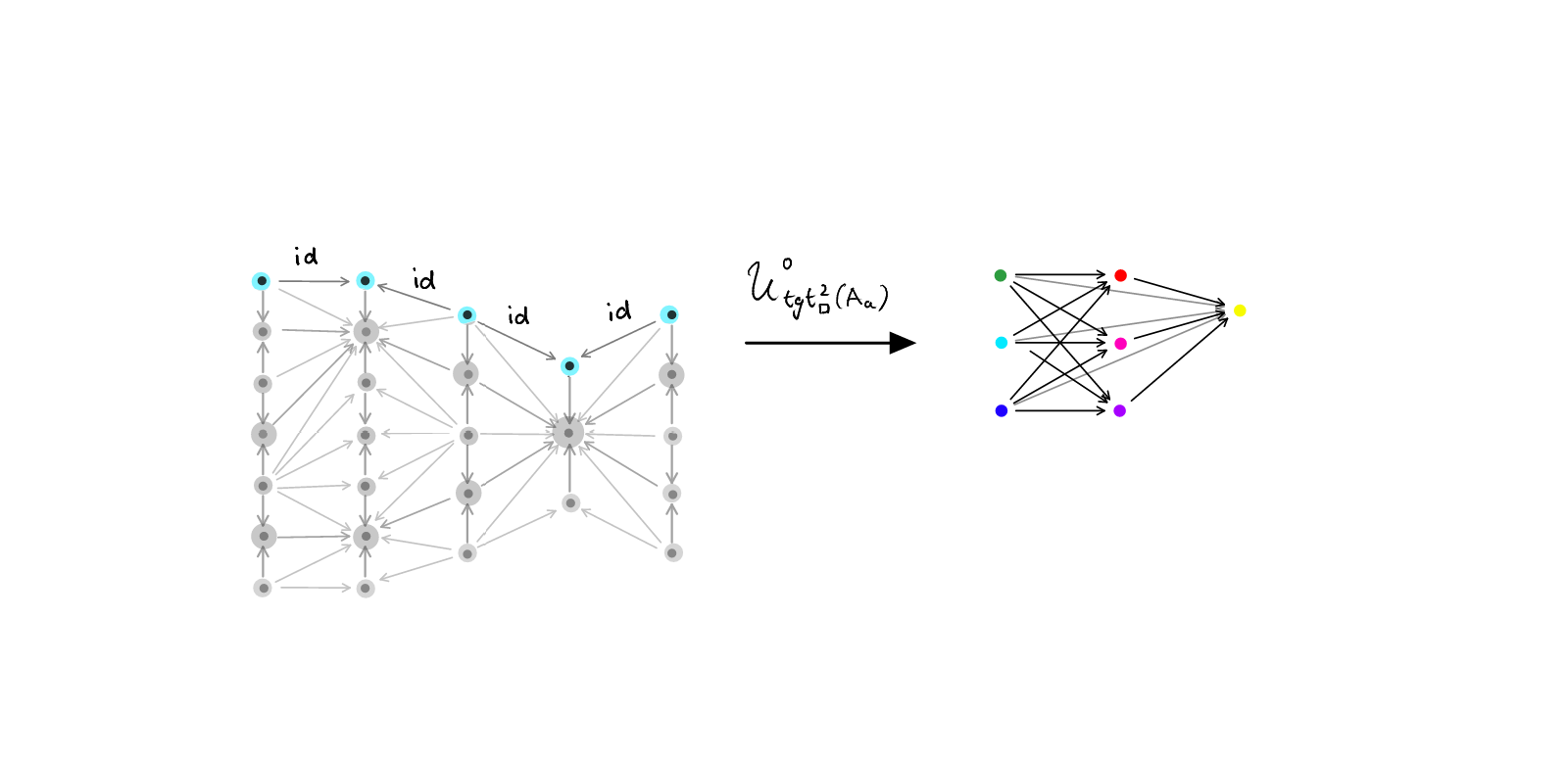}
\endgroup\end{restoretext}
\end{eg}

This observation is generalised and formalised by the following.

\begin{lem}[Globularity conditions] \label{lem:glob_cond} Let $\scA : \bnum{1} \to \SIvert n \cC$ be globular and $n > 1$. Then 
\begin{equation}
\NF{\gsrc(\gsrc(\scA))}^{n-2} = \NF{\gsrc(\gtgt(\scA))}^{n-2}
\end{equation}
and 
\begin{equation}
\NF{\gtgt(\gsrc(\scA))}^{n-2} = \NF{\gtgt(\gtgt(\scA))}^{n-2}
\end{equation}
\proof  The proof is \stfwd{}. We prove the first statement (the second follows similarly). Note by \autoref{constr:globular_src_tgt}, $\gsrc\gsrc(\scA), \gsrc\gtgt(\scA) : \bnum{1} \to \SIvert {n-2} \cC$ are given by
\begin{align}
\gsrc\gsrc(\scA) &= \tsU 2_\scA \sG(\Delta_{(0,0)}) \Delta_{((0,0),0)} \\
&= \tsU 2_\scA \Delta_{((0,0),0)}
\end{align}
and
\begin{align}
\gsrc\gtgt(\scA)(0) &= \tsU 2_\scA \sG(\Delta_{(0,\iH_{\scA}(x))}) \Delta_{((0,\iH_{\scA}(x)),0)} \\
&= \tsU 2_\scA \Delta_{((0,\iH_{\scA}(x)),0)}
\end{align}
Next, observe that both $((0,0),0)$ and $((0,\iH_{\scA}(x)),0)$ lie in (cf. \autoref{constr:source_and_target_inclusion})
\begin{equation}
Y := \im(\msrc_{\tsU 1_\scA}) \subset \tsG 2 (\scA)
\end{equation}
Using \autoref{rmk:src_tgt_are_regular_content} and the \autoref{defn:globular_families} of globularity we infer that for any $y_1 \to y_2 \in \mor(Y)$, both $\rest {\tsU 2_{\scA}} {y_1}$ and $\rest {\tsU 2_{\scA}} {y_2}$ have the same normal form. Since by \autoref{rmk:connectedness_of_src_tgt} $Y$ is connected we infer that all $\rest {\tsU 2_{\scA}} y$ (for $y \in Y$) have the same normal form. In particular $\gsrc\gsrc(\scA), \gsrc\gtgt(\scA)$ have the same normal form as claimed. \qed
\end{lem}

\begin{cor}[Globularity on the nose]\label{cor:globularity} If $A : \bnum{1} \to \SIvert n \cC$ is globular and normalised then
\begin{equation}
\gsrc(\gsrc(\scA)) = \gsrc(\gtgt(\scA))
\end{equation}
and 
\begin{equation}
\gtgt(\gsrc(\scA)) = \gtgt(\gtgt(\scA))
\end{equation} 
\proof This follows from \autoref{lem:src_tgt_normalisation} and \autoref{lem:glob_cond}.
\qed
\end{cor}

\begin{defn}[$k$-level globular source and target] \label{defn:klvl_source_target}
Let $A : \bnum{1} \to \SIvert n \cC$. Then we inductively define $\gsrc^1(\scA) = \gsrc(\scA)$ and $\gtgt^1(\scA) = \gtgt(\scA)$ and set 
\begin{equation}
\gsrc^{k+1}(\scA) = \gsrc(\gsrc^k(\scA))
\end{equation}
as well as
\begin{equation}
\gtgt^{k+1}(\scA) = \gtgt(\gtgt^k(\scA))
\end{equation}
\end{defn}

We now introduce a special form of composition that applies in the case of globular families.

\begin{constr}[Whiskering composition] \label{constr:glob_comp} Let $1 \leq k \leq n$, and let $\scA : \bnum{1} \to \SIvert n \cC$ be globular and normalised up to level $n$ and $\scB : \bnum{1} \to \SIvert {n-k+1} \cC$ be globular normalised up to level $(n-k+1)$. Under the condition that
\begin{equation} \label{eq:postwhiskering_condition}
\gtgt^k(\scA) = \gsrc(\scB)
\end{equation}
we define the \textit{$k$-level (post-)whiskering of $\scA$ with $\scB$}, denoted by $\scA \whisker {k} {n} \scB$, as follows: let $\bang : \tsG {k-1}(\scA) \to \bnum{1}$ denote the unique functor to the terminal category $\bnum{1}$. Since $\scA$ is globular normalised, \autoref{lem:local_trivial_to_global_trivial} implies that $\ctgt^k(\scA)$ is globally trivial, and thus in particular
\begin{equation}
\ctgt^k(\scA)  = \gtgt^k(\scA) \bang
\end{equation}
Using this together with \eqref{eq:postwhiskering_condition}, 
\begin{equation}
\ctgt^k(\scA) = \gsrc(\scB) \bang
\end{equation}
we can then compute
\begin{align}
\ctgt^k(\scA)  &= \csrc(\scB) \bang \\
&= \tsU 1_\scB \msrc_{\und\scB} \bang \\
&= \tsU 1_\scB \sG(\bang) \msrc_{\und\scB\bang} \\
&= \tsU 1_{\scB\bang}  \msrc_{\und\scB\bang} \\
&= \csrc(\scB \bang)
\end{align}
Here, we first unpacked the definition of the globular source, then \autoref{constr:source_target}, before using \autoref{rmk:basechange_for_src_and_tgt} in the third step. The fourth step follows from \autoref{claim:grothendieck_span_construction_basechange}, and in the last step we used \autoref{constr:source_target} again. Thus, by \autoref{constr:k_level_stacking}, we can now define 
\begin{align} \label{eq:post_whiskering}
\scA \whisker {k} {n} \scB &:= \tsR {k-1}_{\sT^{k-1}_\scA, \tsU {k-1}_\scA \stack \scB \bang}
\end{align}

Similarly, under the condition that
\begin{equation} \label{eq:prewhiskering_condition}
\gsrc^k(\scA) = \gtgt(\scB)
\end{equation}
we define the \textit{$k$-level (pre-)whiskering of $\scA$ with $\scB$}, denoted by $\scB \whisker {k} {n} \scA$, as follows:
\begin{align} \label{eq:pre_whiskering}
\scB \whisker {k} {n} \scA &:= \tsR {k-1}_{\sT^{k-1}_\scA, \scB \bang \stack \tsU {k-1}_\scA} 
\end{align}
\end{constr}

We illustrate the above construction step by step for the following example.

\begin{eg}[Whiskering] Consider the $\SIvert 1 \cC$-family $\scB_c$ defined by
\begin{restoretext}
\begingroup\sbox0{\includegraphics{test/page1.png}}\includegraphics[clip,trim=0 {.3\ht0} 0 {.15\ht0} ,width=\textwidth]{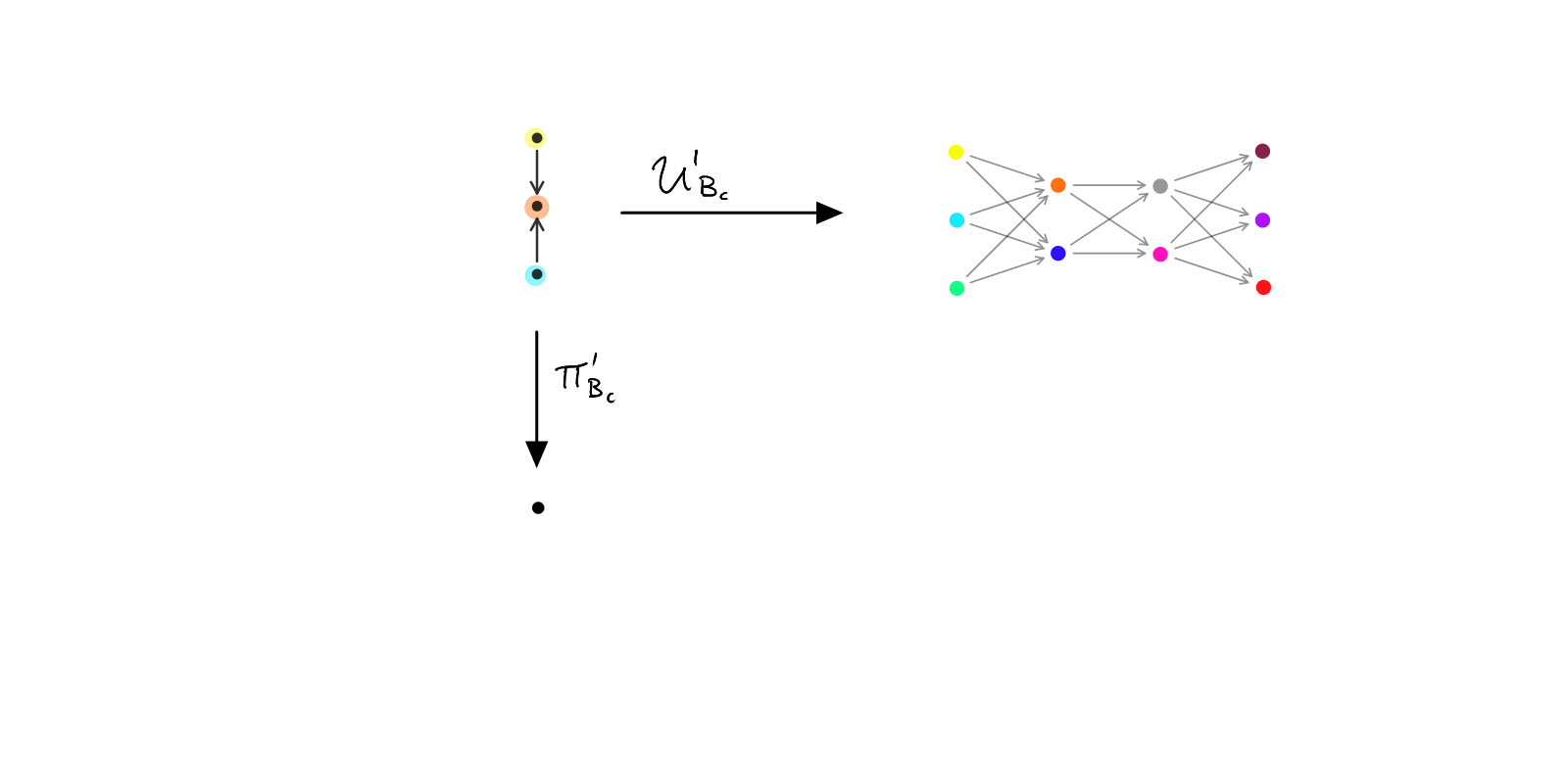}
\endgroup\end{restoretext}
and the $\SIvert 3 \cC$-family $\scA_f$ defined by
\begin{restoretext}
\begin{noverticalspace}
\begingroup\sbox0{\includegraphics{test/page1.png}}\includegraphics[clip,trim=0 {.0\ht0} 0 {.25\ht0} ,width=\textwidth]{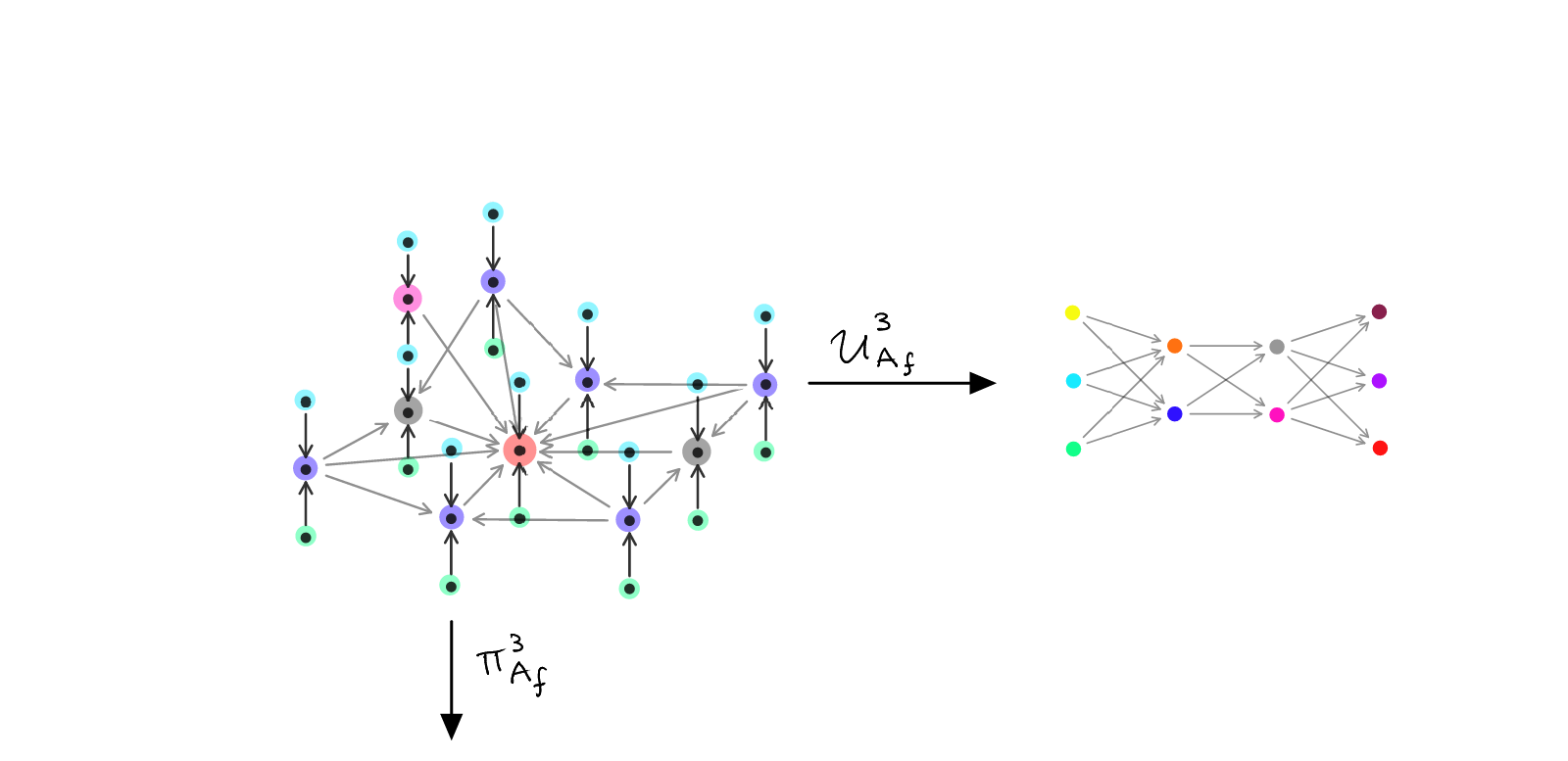}
\endgroup \\*
\begingroup\sbox0{\includegraphics{test/page1.png}}\includegraphics[clip,trim=0 {.2\ht0} 0 {.0\ht0} ,width=\textwidth]{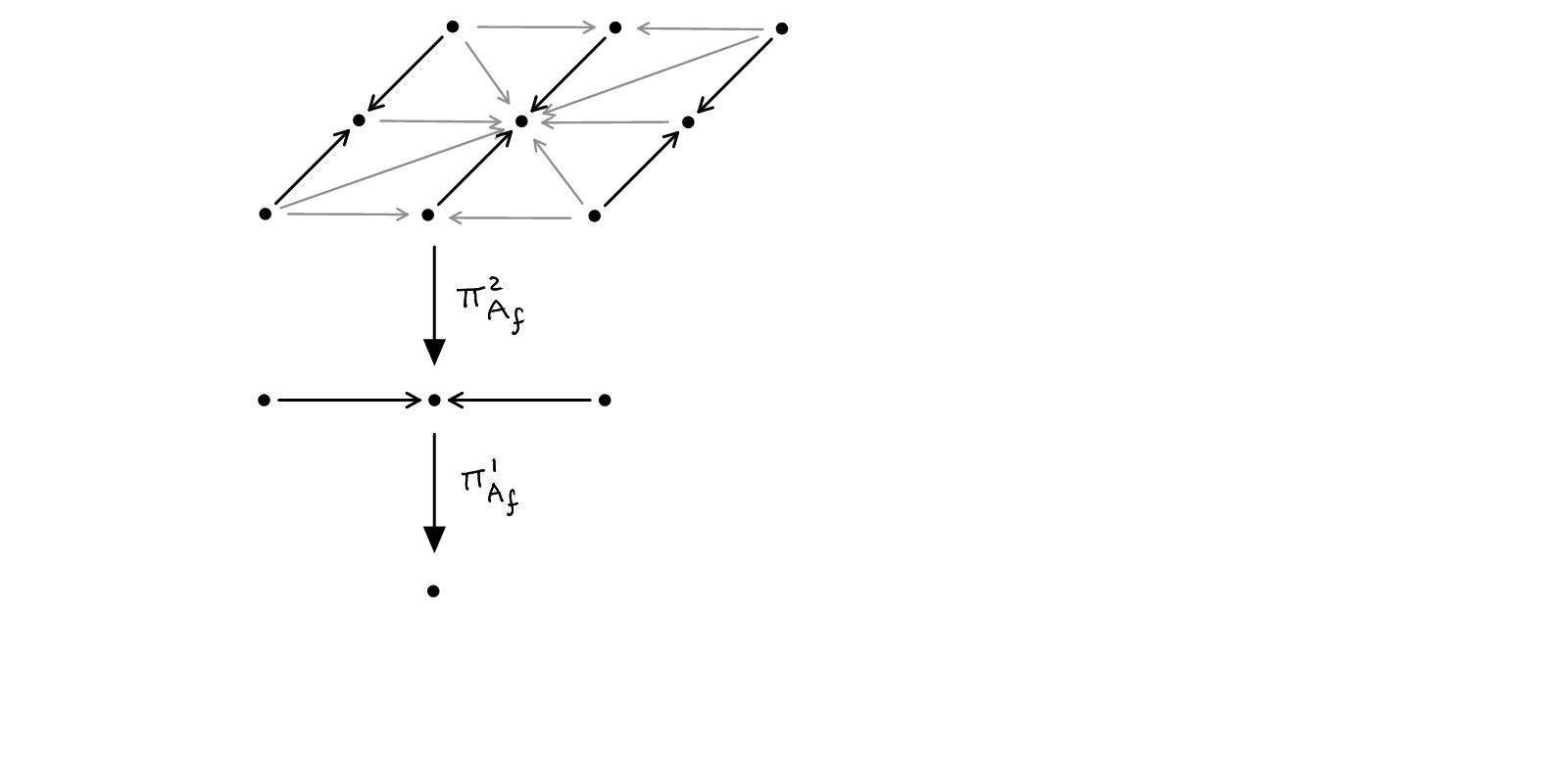}
\endgroup
\end{noverticalspace}
\end{restoretext}
Note that
\begin{equation} \label{eq:whiskering_cond_eg}
\gtgt^3(\scA_f) = \gsrc(\scB)
\end{equation}
namely, both equal the $\SIvert 0 \cC$-family
\begin{restoretext}
\begingroup\sbox0{\includegraphics{test/page1.png}}\includegraphics[clip,trim=0 {.45\ht0} 0 {.25\ht0} ,width=\textwidth]{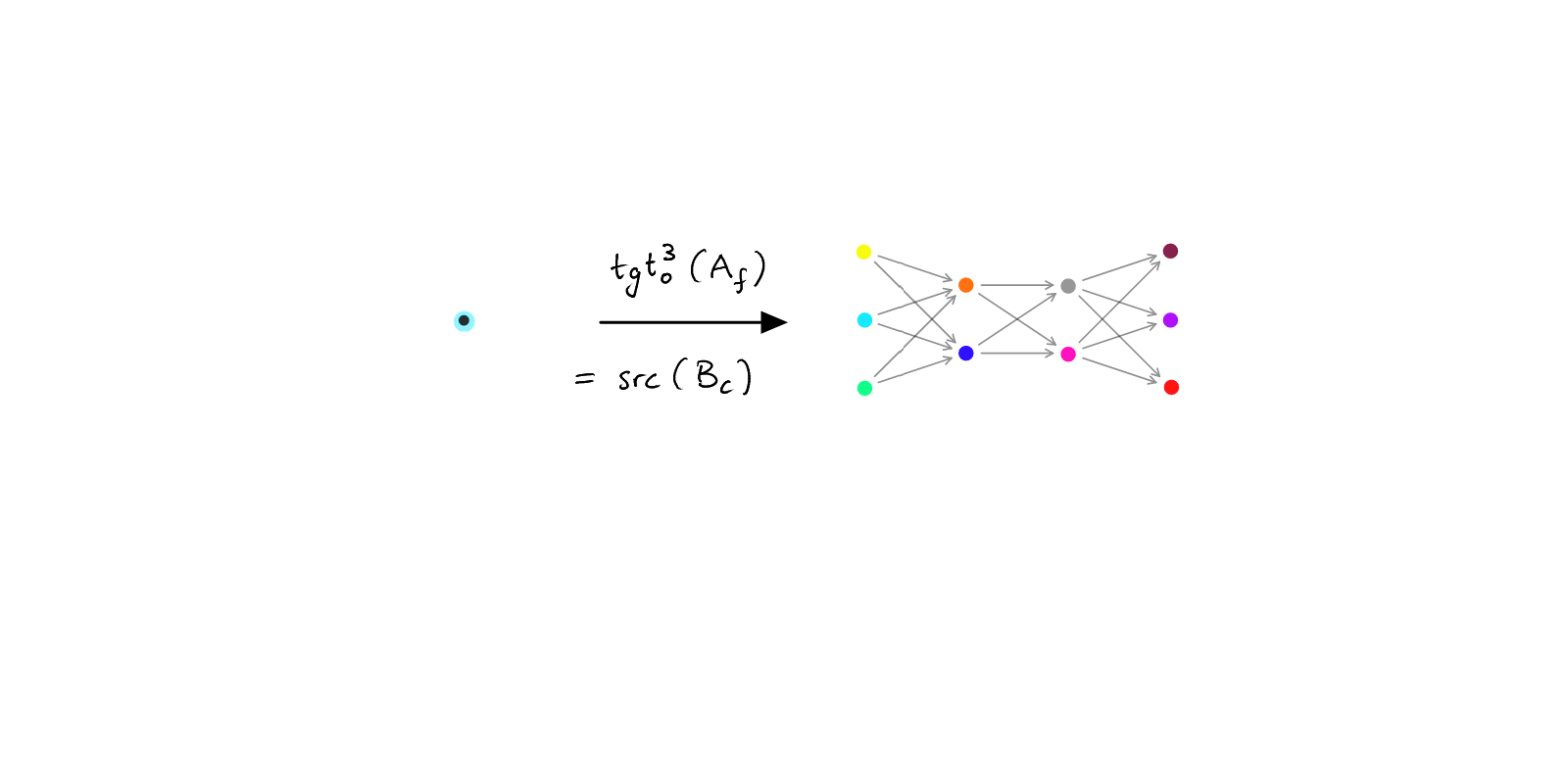}
\endgroup\end{restoretext}
Thus $\scA_f \whisker 3 3 \scB_c$ is well-defined. The functor $\bang : \tsG 2(\scA_f) \to \bnum{1}$ is the terminal functor
\begin{restoretext}
\begingroup\sbox0{\includegraphics{test/page1.png}}\includegraphics[clip,trim=0 {.2\ht0} 0 {.35\ht0} ,width=\textwidth]{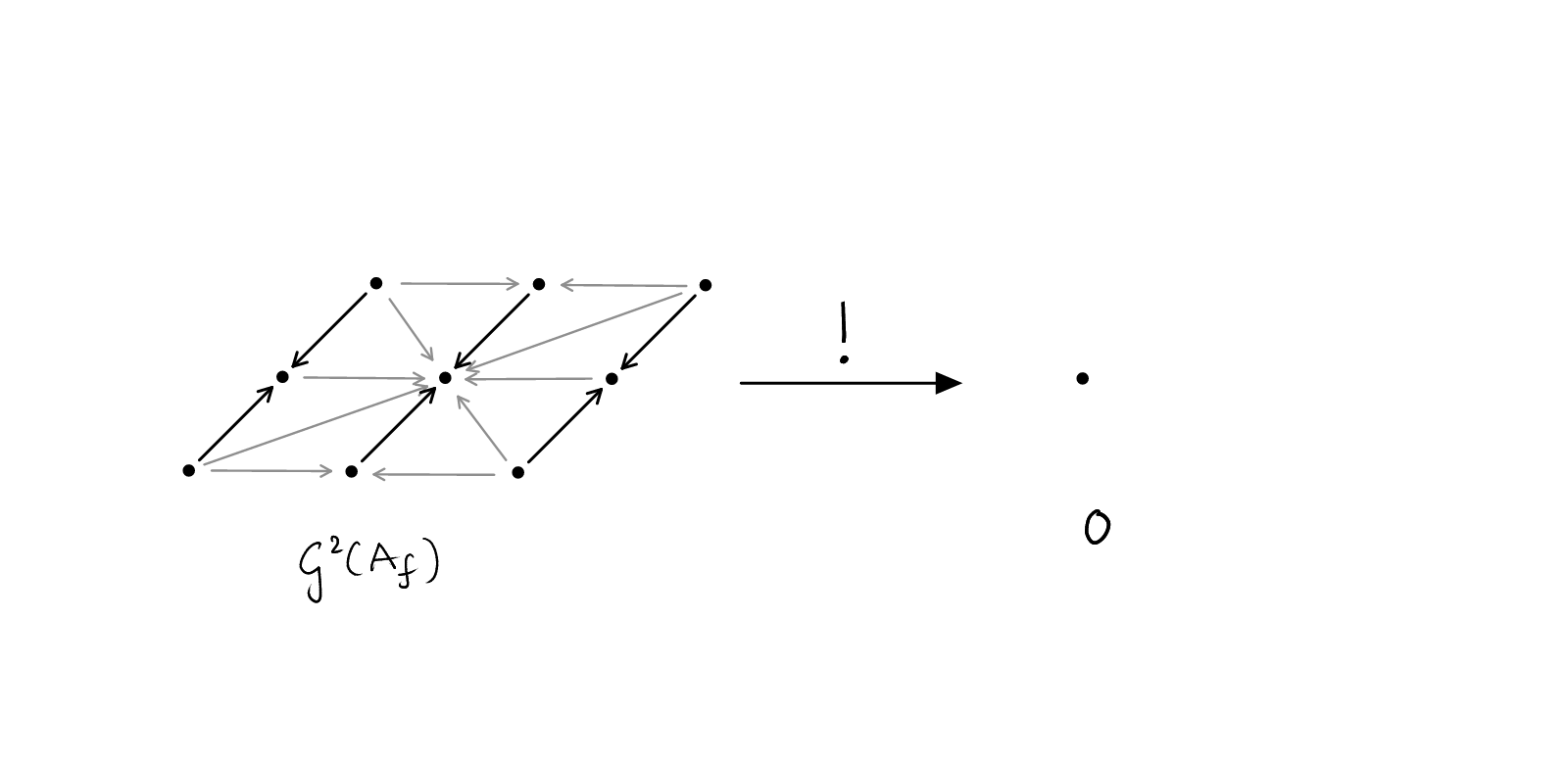}
\endgroup\end{restoretext}
And pulling back $\scB_c$ along this functor gives a $\SIvert 1 \cC$-family $\scB_c \bang$ with data
\begin{restoretext}
\begingroup\sbox0{\includegraphics{test/page1.png}}\includegraphics[clip,trim=0 {.1\ht0} 0 {.0\ht0} ,width=\textwidth]{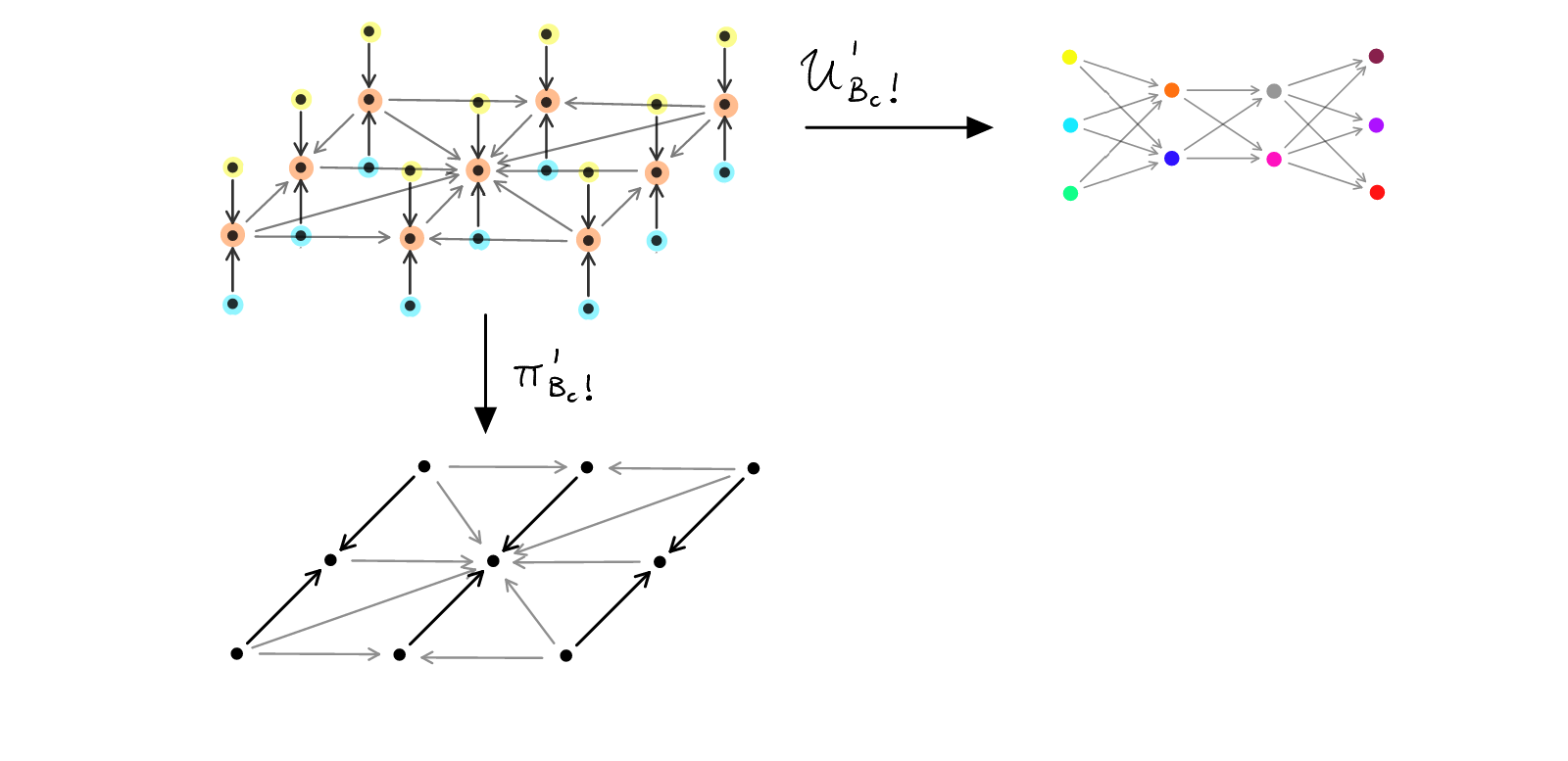}
\endgroup\end{restoretext}
Comparing to the family $\tsU 2_{\scA_f}$
\begin{restoretext}
\begingroup\sbox0{\includegraphics{test/page1.png}}\includegraphics[clip,trim=0 {.0\ht0} 0 {.0\ht0} ,width=\textwidth]{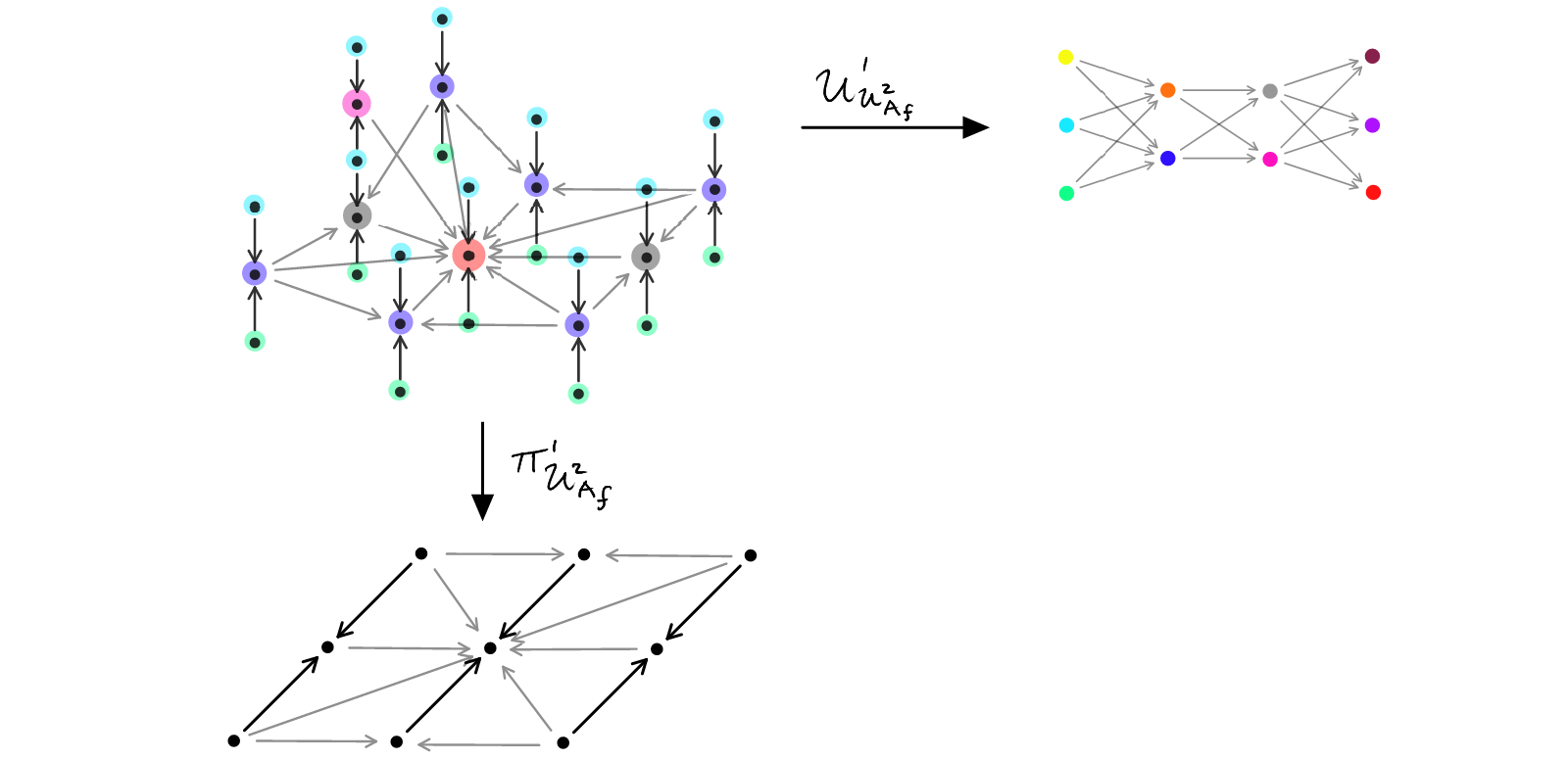}
\endgroup\end{restoretext}
we see that 
\begin{equation}
\ctgt(\tsU 2_{\scA_f}) = \csrc(\scB_c \bang)
\end{equation}
(as shown in the construction, this follows from globularity and assumption \eqref{eq:whiskering_cond_eg} in general). Namely, both equal the $\SIvert 0 \cC$-family
\begin{restoretext}
\begingroup\sbox0{\includegraphics{test/page1.png}}\includegraphics[clip,trim=0 {.35\ht0} 0 {.3\ht0} ,width=\textwidth]{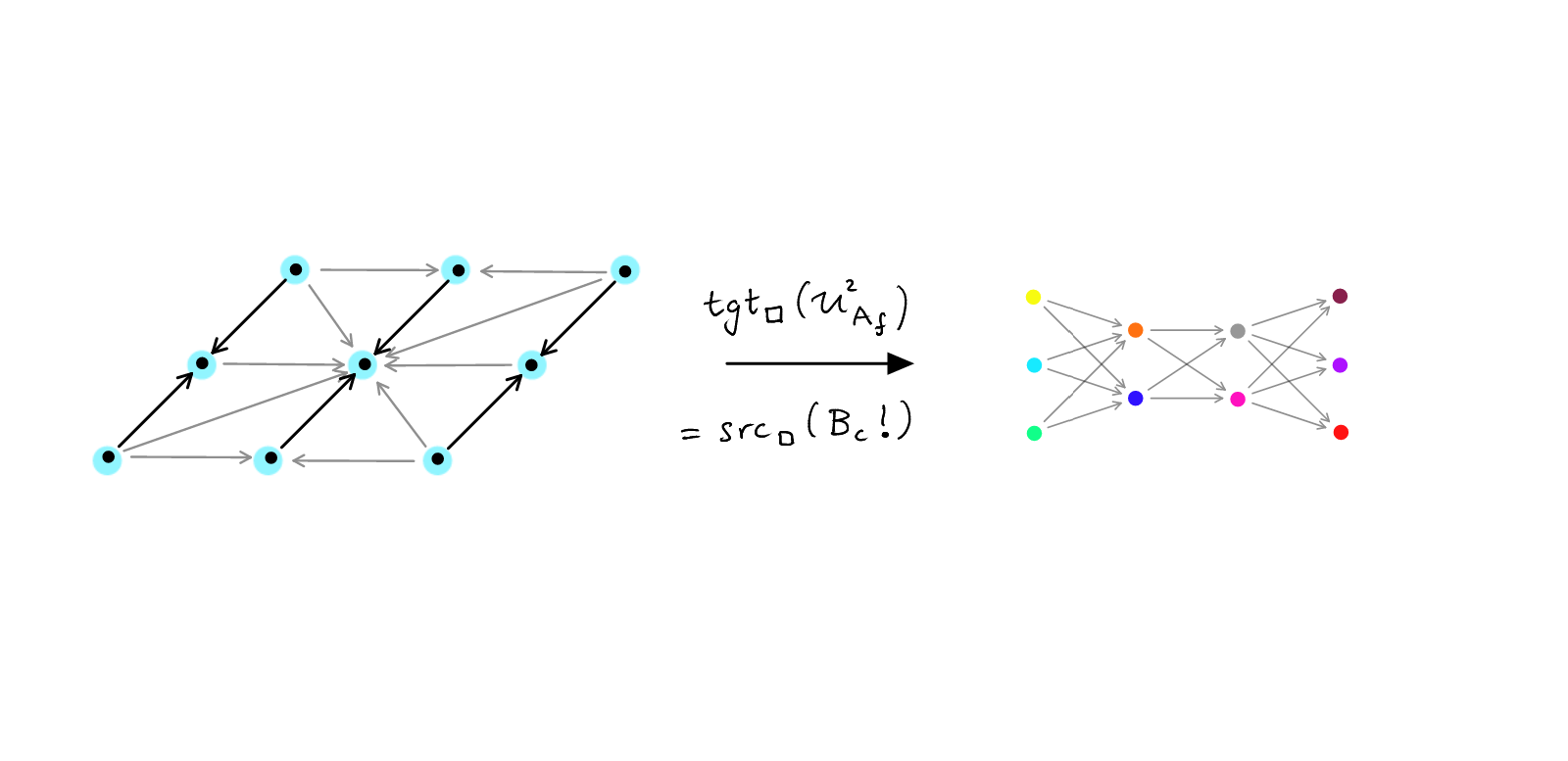}
\endgroup\end{restoretext}
Thus the $1$-level stacking $\tsU 2_{\scA_f} \glue 1 \scB_c \bang$ exists. $\scA_f \whisker 3 3 \scB_c$ is now obtained from $\tsU 2_{\scA_f} \glue 1 \scB_c \bang$ by adding the remaining tower of bundles below level $2$, that is
\begin{equation}
\scA_f \whisker 3 3 \scB_c = \tsR 2_{\sT^2_{\scA_f}, \tsU 2_{\scA_f} \glue 1 \scB_c \bang}
\end{equation}
As a result we see that $\scA_f \whisker 3 3 \scB_c$ equals the family $\scA_c$ from \autoref{eg:sources_and_targets}.
\end{eg}

\subsection{Properties of whiskering composition}

The next three lemmas are technical proofs of properties of whiskering.

\begin{lem}[Whiskering preserves normalisation] \label{lem:gen_comp_pres_norm}  $\scA \whisker {k} {n} \scB$ as defined in \eqref{eq:post_whiskering} is in normal form (up to level $n$). Similarly, $\scB \whisker {k} {n} \scA$ as defined in \eqref{eq:pre_whiskering} is in normal form (up to level $n$).
\proof The proof is \stfwd{}. We prove the first statement (the second follows similarly). Set
\begin{equation} \label{eq:gen_comp_pres_norm}
\widetilde\scB := \tsR {k-1}_{\sT^{k-1}_\scA, \scB \bang}
\end{equation}
Note that by \autoref{constr:k_level_stacking} and \autoref{constr:glob_comp} we have
\begin{align}
\scA \whisker {k} {n} \scB &:= \tsR {k-1}_{\sT^{k-1}_\scA, \tsU {k-1}_\scA \stack \scB \bang} \\
&= \scA \glue {k} \widetilde\scB
\end{align}
Assume by contradiction that $\lambda : \scA \whisker {k} {n} \scB \kcoll l \scC$.  Following the argument of \autoref{lem:stacking_preserves_nomalisation} we find a singular height $y$ over some $x \in \sG^{l-1}(\scA \whisker k n \scB)$ such that
\begin{equation}
\sS^\lambda (x,y) \in \regcont(\tsG l(\scC))
\end{equation}
As before we have
\begin{equation}
(x,y) \in \tsG l(\scA \glue k \widetilde\scB) \quad \imp \quad (x,y) = (\stackinc k_1)^l (x_2,y_2) \text{~or~} (x,y) = (\stackinc k_2)^l  (x_2,y_2)
\end{equation}
for some $x_i, y_i$. If the former holds, the contradiction can be derived as in \autoref{lem:stacking_preserves_nomalisation}. Thus assume
\begin{equation}
(x,y) \notin (\stackinc k_1)^l (\tsG l(\scA)) \text{~and~} (x,y) \in (\stackinc k_2)^l (\tsG l(\widetilde\scB))
\end{equation}
Since $\stackinc k_1= \id : \tsG j(\scA) \to \tsG j(\scA \whisker {k} {n} \scB)$ for $j < k$ (cf. \autoref{constr:subfamilies_from_endpoints}) we must have $l \geq k$, and thus for some $\scD$ we have
\begin{equation}
(\stackinc k_2)\pbstar  \lambda : \tsU {k-1}_{\widetilde \scB}  \kcoll {l - (k-1)} \scD
\end{equation}
From \eqref{eq:gen_comp_pres_norm} we deduce
\begin{equation}
\tsU {k-1}_{\widetilde \scB} = \scB\bang
\end{equation}
and thus (denoting the projection of $x_2$ to $\sG^{k-1}(\widetilde \scB)$ by $x^{k-1}_2$)
\begin{equation}
\scB = \tsU {k-1}_{\widetilde \scB} \Delta_{x^{k-1}_2}
\end{equation} 
We infer a collapse
\begin{equation}
((\stackinc k_2)\pbstar \lambda) \tsG {l-k}(\Delta_{x^{k-1}_2}) : \tsU {k-1}_{\widetilde \scB} \Delta_{x^{k-1}} \kcoll {l - (k-1)} \scD \Delta_{x^{k-1}_2}
\end{equation}
But by \autoref{defn:regions} and construction of $x^{k-1}_2$, we have $\widetilde x \in \tsG {l-k}(\scB)$ such that $\tsG {l-k}(\Delta_{x^{k-1}_2})(\widetilde x) = x_2$. Thus $((\stackinc k_2)\pbstar \lambda) \tsG {l-k}(\Delta_{x^{k-1}_2})$ has a non-identity component at $\widetilde x$, since the component $((\stackinc k_2)\pbstar \lambda)_{x_2}$ is not the identity by choice of $x_2$ and $x$ (as was argued for in \autoref{lem:stacking_preserves_nomalisation}). This contradicts $\scB$ being in normal form. \qed
\end{lem}

\begin{lem}[Whiskering preserves globularity] \label{lem:whisk_pres_glob} $\scA \whisker {k} {n} \scB$ as defined in \eqref{eq:post_whiskering} is globular. Similarly, $\scB \whisker {k} {n} \scA$ as defined in \eqref{eq:pre_whiskering} is globular. 
\proof The proof is \stfwd{}. We prove the first statement (the second follows similarly). Set
\begin{equation} \label{eq:gen_comp_pres_norm2}
\widetilde\scB := \tsR {k-1}_{\sT^{k-1}_\scA, \scB \bang}
\end{equation}
We need to show $\scA \glue {k} \widetilde \scB$ is $l$-level globular, that is, setting $\scC := \tsU l_{\scA \glue {k} \widetilde \scB}$, we need to show  $\rest {\scC} {\regcont(\tsG l(\scA \glue {k} \widetilde \scB))}$ is locally trivial, for all $0 < l \leq n$. For this in turn we need to show for $(x \to y) \in \regcont(\tsG l(\scA \glue {k} \widetilde \scB))$ that $\rest {\scC} {(x \to y)}$ normalises to the identity.

Using the embeddings defined in \autoref{rmk:can_subbund_of_gluing} we recall
\begin{equation}
(\stackinc k_1)^l (\tsG l(\scA)) \cup (\stackinc k_2)^l (\tsG l(\widetilde\scB)) = \tsG l(\scA \glue k \widetilde\scB)
\end{equation}
Since both $\scA$ and $\widetilde\scB$ are \gls{downwardclosed} embeddings of $\scA \glue k \widetilde\scB$ (cf. \autoref{cor:downward_closed_subbund}) we infer that either $(x \to y) \in (\stackinc k_1)^l(\tsG l(\scA))$ or $(x \to y) \in (\stackinc k_1)^l(\tsG l(\widetilde\scB))$. 

If $l < k$, we know (as shown in the proof of \autoref{lem:gen_comp_pres_norm}) that $(x \to y) \in \mor(\tsG l(\scA))$. Set $i := k-1-l\geq 0$. We compute
\begin{align}
\scC \Delta_{x \to y} &= \tsR {i}_{\sT^{i}_\scC, \tsU {k-1}_\scA \stack \scB\bang} \Delta_{x \to y} \\
&=  \tsR {i}_{\sT^{i}_\scC \Delta_{x \to y}, (\tsU {k-1}_\scA \stack \scB\bang)\tsG {i}(\Delta_{x \to y})} \\
&= \tsR {i}_{\sT^{i}_\scC \Delta_{x \to y}, (\tsU {k-1}_\scA \tsG {i}(\Delta_{x \to y})) \stack (\scB\bang\tsG {i}(\Delta_{x \to y}))} \\
&= (\tsU l_\scA \Delta_{x \to y}) \whisker {i+1} n \scB
\end{align}
Here, in the first step we unpacked the definition of $\scC$, and in the second step we used \autoref{cor:basechange_compact}. The third step follows from \autoref{rmk:basechange_stacking} and the last step follows from definition of whiskering (and $\bang \tsG {i}(\Delta_{x \to y}) =~\bang$ as well as $\tsU {k-1}_\scA \tsG {i}(\Delta_{x \to y}) = \tsU {i}_{\tsU l_\scA \Delta_{x \to y}}$ by \autoref{claim:grothendieck_span_construction_basechange}). 

Since $\scA$ is normalised and globular we know $\tsU l_\scA \Delta_{x \to y} = \id$ and thus $\tsU l_\scA \Delta_{x \to y}$ factors through the terminal category, that is it can be written as $\tsU l_\scA \Delta_{x \to y} = \scD \bang$ for $\bang : \bnum{2} \to \bnum{1}$ and some $\scD : \bnum{1} \to \SIvert {n-l} \cC$. We deduce
\begin{align}
\scC \Delta_{x \to y} &= \scD \bang \whisker {i+1} n \scB \\
&= \tsR {i}_{\sT^{i}_{\scD \bang}, \tsU {i}_{\scD \bang} \stack \scB\bang} \\
&= \tsR {i}_{\sT^{i}_\scD \bang, \tsU {i}_\scD \tsG i(\bang) \stack \scB\bang \tsG i(\bang)} \\
&= \tsR {i}_{\sT^{i}_\scD, \tsU {i}_\scD \stack \scB\bang} \bang
\end{align}
Here, after unwinding definitions we used \autoref{claim:grothendieck_span_construction_basechange} in the third step and \autoref{rmk:basechange_stacking} together with \autoref{cor:basechange_compact} in the forth. Thus $\scC \Delta_{x \to y} $ factors through $\bnum{1}$, that is $\scC \Delta_{x \to y} = \id$, and normalises to the identity as required, since by \autoref{lem:gen_comp_pres_norm} it is already in normal form.

Next, we look at the case that $l \geq k$. Again we know that either $(x \to y) \in (\stackinc k_1)^l(\tsG l(\scA))$ or $(x \to y) \in (\stackinc k_2)^l(\tsG l(\widetilde\scB))$. In the former case, in which $(x \to y) = (\stackinc k_1)^l (\widetilde x \to \widetilde y)$ for some $\widetilde x, \widetilde y \in \tsG l(\scA)$, we deduce from the definition of $\stackinc k_1$ that
\begin{equation}
\tsU l_{\scA}  \Delta_{\widetilde x \to \widetilde y} = \tsU l_{\scA \glue k \widetilde\scB} (\stackinc k_1)^l \Delta_{\widetilde x \to \widetilde y}
\end{equation}
thus globularity follows from globularity of $\scA$ in this case. It remains to consider the case in which $(x \to y) = (\stackinc k_2)^l (\widetilde x \to \widetilde y)$ for some $\widetilde x, \widetilde y \in \tsG l(\widetilde\scB)$. As before we have
\begin{align}
\tsU l_{\widetilde \scB}  \Delta_{\widetilde x \to \widetilde y} = \tsU l_{\scA \glue k \widetilde\scB}(\stackinc k_2)^l \Delta_{\widetilde x \to \widetilde y}
\end{align}
We compute (cf. \eqref{eq:gen_comp_pres_norm2})
\begin{align}
\tsU l_{\widetilde \scB}  \Delta_{\widetilde x \to \widetilde y} &=  \tsU l_{\tsR {k-1}_{\sT^{k-1}_\scA, \scB \bang}} \Delta_{\widetilde x \to \widetilde y} \\
&=  \tsU {l-k+1}_{\scB \bang} \Delta_{\widetilde x \to \widetilde y} \\
&=  \tsU {l-k+1}_{\scB} \tsG {l-k+1}(\bang) \Delta_{\widetilde x \to \widetilde y} 
\end{align}
But since $\tsG {l-k+1}(\bang)$ is (fibrewise) open, the claim that $\rest \scC {x \to y}$ normalises to the identity follows from $\scB$ being globular.
\qed
\end{lem}

\begin{lem}[Source and target of whiskering composition] \label{lem:gen_comp_src_tgt} For $\scA \whisker {k} {n} \scB$ as defined in \eqref{eq:post_whiskering} we have if $k = 1$
\begin{align} \label{eq:post_whisk_src_tgt_k_eq_n}
\gsrc(\scA \whisker {1} {n} \scB) &= \gsrc(\scA) \\
\gtgt(\scA \whisker {1} {n} \scB) &= \gtgt(\scB)
\end{align}
and  if $k > 1$
\begin{align} \label{eq:post_whisk_src_tgt_k_leq_n}
\gsrc(\scA \whisker {k} {n} \scB) &= \gsrc(\scA) \whisker {k-1} {n-1} \scB \\
\gtgt(\scA \whisker {k} {n} \scB) &= \gtgt(\scA) \whisker {k-1} {n-1} \scB 
\end{align}
Similarly, for $\scB \whisker {k} {n} \scA$ as defined in \eqref{eq:pre_whiskering} we have if $k = 1$
\begin{align} 
\gsrc(\scB \whisker {1} {n} \scA) &= \gsrc(\scB) \\
\gtgt(\scB \whisker {1} {n} \scA) &= \gtgt(\scA)
\end{align}
and  if $k > 1$
\begin{align} 
\gsrc(\scB \whisker {k} {n} \scA) &= \scB \whisker {k-1} {n-1} \gsrc(\scA) \\
\gtgt(\scB \whisker {k} {n} \scA) &= \scB\whisker {k-1} {n-1} \gtgt(\scA) 
\end{align}
\proof The proof is \stfwd{}. We show the first statement of \eqref{eq:post_whisk_src_tgt_k_eq_n} and the first statement of \eqref{eq:post_whisk_src_tgt_k_leq_n} (all other statements follow similarly). For the former we compute 
\begin{align}
\gsrc(\scA \whisker {1} {n} \scB) &= \gsrc(\scA \stack \scB) \\
&= \tsU 1_{\scA \stack \scB} \msrc_{\und\scA \stack \und \scB} \\
&= \tsU 1_{\scA} \msrc_{\und\scA}\\
&= \gsrc(\scA)
\end{align}
where we used \autoref{rmk:src_tgt_for_stacking} in the second step.

For the latter statement we compute
\begin{align}
\gsrc(\scA \whisker {k} {n} \scB) &= \tsU 1_{\scA \whisker {k} {n} \scB} \msrc_{\und{\scA \whisker {k} {n} \scB}} \\
&= \tsR {k-2}_{\sT^{k-2}_{\tsU 1_\scA}, \tsU {k-2}_{\tsU 1_\scA} \stack \scB\bang} \msrc_{\und{\scA}} \\
&= \tsR {k-2}_{\sT^{k-2}_{\scA} \msrc_{\und{\scA}}, (\tsU {k-2}_{\tsU 1_\scA} \stack \scB\bang) \tsG {k-2}(\msrc_{\und{\scA}})} \\
&= \tsR {k-2}_{\sT^{k-2}_{\gsrc(\scA)}, (\tsU {k-2}_{\gsrc(\scA)} \stack \scB\bang)} \\
&= \gsrc(\scA) \whisker {k-1} {n-1} \scB
\end{align}
where, in the first step we unpacked the definition of $\gsrc$, in the second step we unpacked the definition of $\whisker k n$ (for $k \geq 2$), in the third step we used \autoref{cor:basechange_compact}, in the fourth step we used \autoref{constr:unpacking_collapse} and \autoref{rmk:basechange_stacking} (and that $\bang \tsG {k-2}(\msrc_{\und{\scA}}) = ~\bang$), and in the final step we repacked the definition of $\gsrc$. \qed
\end{lem}

\subsection{Definition of $\GComp$}

We now define the endofunctor $\GComp : \globset \to \globset$. We chose a slightly different but equivalent formulation to that in \autoref{ssec:sum_gcomp} adding a lot of detail along the way.

\begin{constr}[Globular sets of generic composites] \label{constr:bottom_up_gcomp} Let $\sC$ be a presented associative $\infty$-category. We inductively (in dimension $k$) define the globular set $\GComp(\sC)$ of \textit{generic composites} as follows. We first set
\begin{equation}
\GComp(\sC)_0 = \Comp(\sC)_0
\end{equation}
Now assume we have defined $\GComp(\sC)_{k-1}$. We define $\GComp(\sC)_k$ inductively as follows
\begin{itemize}
\item For each $g \in \sC_k$, if $\gsrc\abss{g}, \gtgt\abss{g} \in \GComp(\sC)_{k-1}$ then $\abss{g} \in \GComp(\sC)_k$
\item If $f_i \in \GComp(\sC)_{k-1}$ ($0 \leq i \leq m$) and if $H : f_0 \xiso{} f_1 \xiso{} ... \xiso{} f_m$ is non-recursively generic (cf. \autoref{defn:top_down_gen_comp}), then $H \in \GComp(\sC)_k$
\item Finally, if $f_0 \in \GComp(\sC)_k$, $f_1 \in \GComp(\sC)_j$ for $j \leq k$, and if $f_0 \whisker j k f_1$ (resp. $f_1 \whisker j k f_0$) exist then $(f_0 \whisker j k f_1) \in \GComp(\sC)_k$ (resp. $(f_1 \whisker j k f_0) \in \GComp(\sC)_k$)
\end{itemize}
This completes the definition of $\GComp(\sC)$. We now show that it is a globular subset of $\Comp(\sC)$. We first show it is a globular set, for which we need to show closure under source and target maps. Let $f \in \GComp(\sC)$. We distinguish three cases
\begin{itemize}
\item If $f = \abss{g}$, $g \in \sC_k$  we must have $\gsrc(f), \gtgt(f) \in \GComp(\sC)_{k-1}$ by our construction above.
\item If $f = (H : f_0 \xiso{} f_1 \xiso{} ... \xiso{} f_m)$ is a homotopy then source $f_0$ and target $f_m$ are in $\GComp(\sC)_{k-1}$ by assumption above.
\item If $f = (f_0 \whisker j k f_1)$ then using \autoref{lem:gen_comp_src_tgt} we find that for $k = 1$ we have $\gsrc(f) = \gsrc(f_0), \gtgt(f) =\gtgt(f_1) \in \GComp(\sC)_{k-1}$ by induction. Otherwise if $k > 1$ (again using \autoref{lem:gen_comp_src_tgt}), then 
\begin{equation}
\gsrc(f) = \gsrc(f_0) \whisker {j-1} {k-1} \abss{f_1} \in \GComp(\sC)_{k-1}
\end{equation}
by induction and since $\GComp(\sC)_{k-1}$ is closed under whiskering. Similarly for $\gtgt(f)$.
\end{itemize}
\noindent This shows $\GComp(\sC)$ is a globular set

Next, we prove it is a globular subset of $\Comp(\sC)$. Let $f \in \GComp(\sC)$. We distinguish three cases.
\begin{itemize}
\item If $f = \abss{g}$, $g \in \sC_k$ then $\abss{g} \in \Comp(\sC)_k$ by definition (cf. \autoref{defn:pres_ANC} and \autoref{defn:PANC_mor}).
\item If $f = (H : f_0 \xiso{} f_1 \xiso{} ... \xiso{} f_m)$ is a homotopy then $H \in \Comp(\sC)_k$ follows from \autoref{defn:homotopies}.
\item If $f = (f_0 \whisker j k f_1)$, note that $f$ is globular and normalised by \autoref{lem:whisk_pres_glob} and \autoref{lem:gen_comp_pres_norm}. For well-typedness, set $\scA = f_0$ and $\scB = f_1$. Using notation from \autoref{lem:gen_comp_pres_norm} note that $f = \scA \glue k \widetilde\scB$. Using \autoref{rmk:can_subbund_of_gluing}, for $p \in \sG^k(f)$ the minimal neighbour $f \slash p$ by virtue of minimality either lands fully in $\scA$ or fully in $\widetilde \scB$. Since $\scA$ is well-typed, assume the latter. In this case, we can find a collapse sequence $\vec\mu$ (of identities $\vvec\mu^i = \id$ for $i > k$ and $\vvec\mu^i : \const_{\singint 0} \to \tusU {i}_{f}$ for $i \leq k$) which collapses $\widetilde \scB$ to $\Id^{k-1}(\scB)$ and thus well-typedness follows from well-typedness of $\scB$ in this case.
\end{itemize}
This proves that $\GComp(\sC) \subset \Comp(\sC)$.

The association $\sC \mapsto \GComp(\sC)$ canonically lifts to a functor
\begin{equation}
\GComp : \pCat_\infty \to [\lG\op,\SetCat]
\end{equation}
Indeed, assuming $\alpha : \sC \to \sD$ to be a map of presentations, we construct
\begin{align}
\GComp(\alpha)_k :~ \GComp(\sC)(k) &\to \GComp(\sD)(k) \\
f &\mapsto \SIvert k {\alpha} f
\end{align}
The fact that this is well-defined can again be proven inductively (the fact that it is globular follows since taking sources and targets commutes with relabelling). Assume we proved well-definedness for $l < k$ ($k = 0$ is trivial) and let $f \in \GComp(\sC)_k$. We distinguish three cases
\begin{itemize}
\item If $f = \abss{g}$, $g \in \sC_k$ then we must have $\GComp(\alpha)_k(f) \in \GComp(\sD)_k$ by the definition of maps of presentations and arguing inductively.

\item If f $f = (H : f_0 \xiso{} f_1 \xiso{} ... \xiso{} f_m)$ is a homotopy then by induction we find for each $i$
\begin{equation}
f'_i = \SIvert k {\alpha} f_i \in \GComp(\sD)_{k-1}
\end{equation}
Also note that the non-triviality property of subcubes (cf. \autoref{notn:non_identity_subcubes}) is preserved by maps of presentations since homotopies are locally constant. Consequently, we have 
\begin{equation}
(\SIvert k {\alpha} f : f'_0 \xiso{} f'_1 \xiso{} ... \xiso{} f'_k) \in \GComp(\sD)
\end{equation}
as required.

\item If $f = (f_0 \whisker j k f_1)$ then $\GComp(\alpha)_k(f) \in \GComp(\sD)(k-1)$ since
\begin{align}
\GComp(\alpha)_k(f) &= \SIvert k {\alpha} ( f_0 \whisker k n f_1 )\\
&=  (\SIvert k {\alpha} f_0) \whisker k n (\SIvert k {\alpha} f_1)
\end{align}
and thus the statement follows inductively.
\end{itemize}

Note that for a globular set $S$ we have $\kC(S) \in \pgCat$. Indeed, by \autoref{constr:panc_from_glob}, the fact that sources and targets commute with relabelling and \autoref{constr:terminal_n_globe} we find that $g \in S_m$
\begin{align}
\gsrc \abss{g} &= \abss{s_m g} \\
\gtgt \abss{g} &= \abss{t_m g}
\end{align}
which are in turn elements of $\GComp(\kC(S))(m-1)$ as required. By functoriality of $\kC$, the association $S \mapsto \GComp(\kC(S))$ then gives rise to an endofunctor
\begin{equation}
\GComp : [\lG\op,\SetCat] \to [\lG\op,\SetCat]
\end{equation}
Note that we are abusing notation slightly, writing $\GComp$ in instead of $\GComp \circ \kC$. There further is a canonical natural transformation (referred to as the \textit{unit} of $\GComp$)
\begin{equation}
\abss{-} : \id \to \GComp
\end{equation}
such that
\begin{equation}
\abss{-}_S : S \to \GComp(S)
\end{equation}
maps $g \in S_m$ to $\abss{g} \in \GComp(S) \subset \SIvert m {\GGamma{}{\kC(S)}}$. Importantly, $\GComp$ does not have further monad structure (see next remark).
\end{constr}

\begin{rmk}[Failure of monad structure] \label{rmk:monad_fail} Importantly, $\GComp$ does not have further monad structure. This is $\GComp(S)$ contains both whiskering composite and homotopies, and this can cause non-unique choices when grafting the former into the latter: consider for instance the binary interchange in $\GComp(\GComp(S))$ (that is, the easiest elementary homotopy, see next definition), and assume each $1$-stratum is in fact a $1$-level  whiskering composite in $\GComp(S)$ of generators in $S$. Then there is no unique choice for a ``grafted" homotopy in $\GComp(S)$.
\end{rmk}

\subsection{Elementary generic morphisms}

\begin{defn}[Elementary generic morphism] \label{defn:elementarity} Let $\sC \in \pCat_\infty$. We say $f \in \GComp(\sC)_n$ is \textit{elementary} if $f = f_1 \whisker k n f_2$ implies that either one of $f_1$, $f_2$ is an identity (that is, it is of the form $f_i = \Id_{\abss{g}}$ for $g \in \GComp(\sC)$). Otherwise, $f$ is said to be \textit{non-elementary}.
\end{defn}

\subsection{Examples of generic composites}

We will exemplify \autoref{constr:bottom_up_gcomp} in the case of $\kC(S)$ from \autoref{eg:panc_from_gset}. Abusing notation, we will write $z$ in place of $\abss{z}$ for $z \in \kC(S)_k = S_k$ in this section. 

To begin with, note that 
\begin{equation}
\gtgt^2 x = \gsrc h
\end{equation}
and thus we can form the generic composite $x \whisker 2 2 h$ which is given by
\begin{restoretext}
\begingroup\sbox0{\includegraphics{test/page1.png}}\includegraphics[clip,trim=0 {.0\ht0} 0 {.0\ht0} ,width=\textwidth]{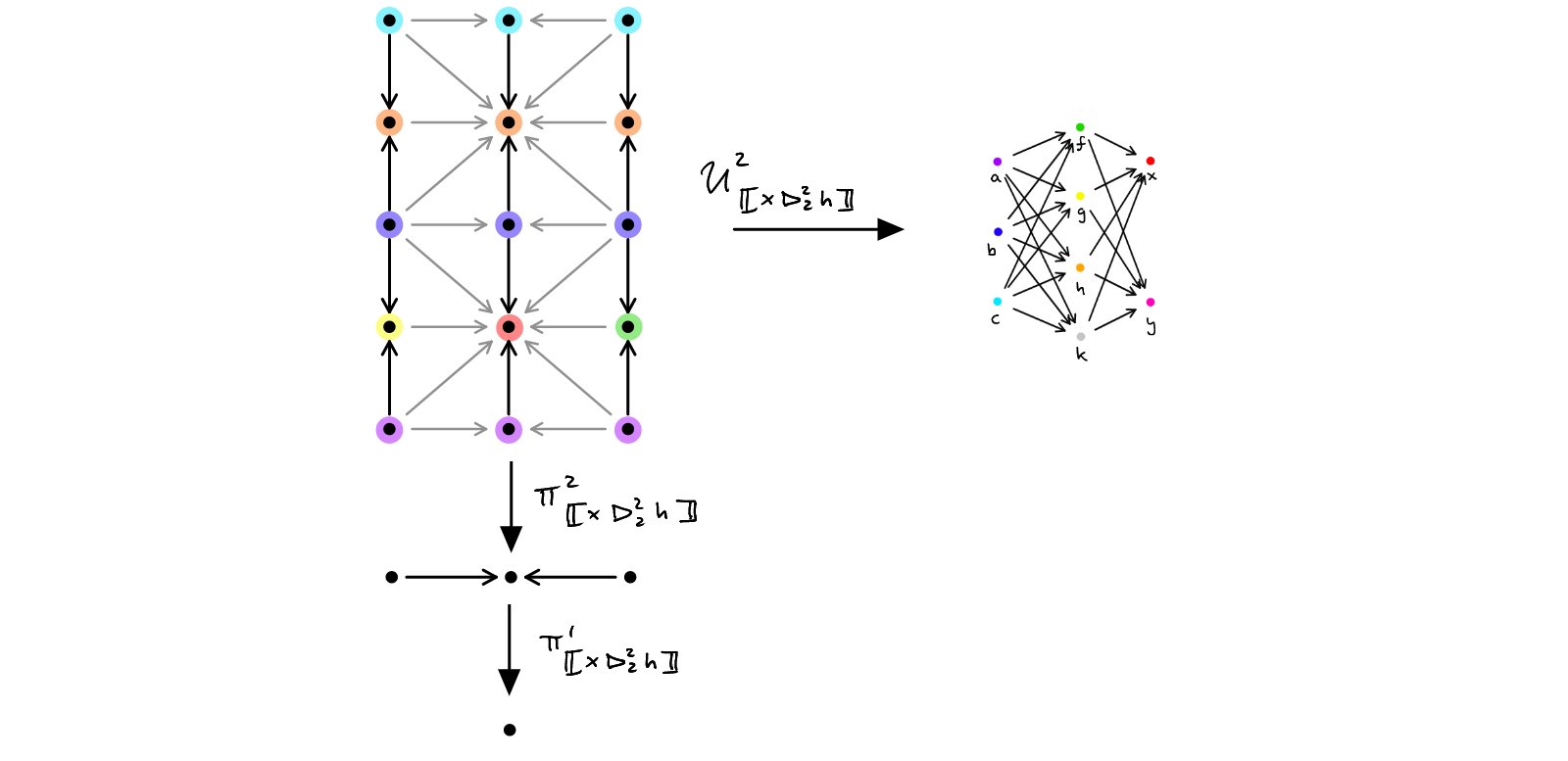}
\endgroup\end{restoretext}
Similarly, since
\begin{equation}
\gsrc^2 y = \gtgt g
\end{equation}
we can form the generic $g \whisker 2 2 y$ given by
\begin{restoretext}
\begingroup\sbox0{\includegraphics{test/page1.png}}\includegraphics[clip,trim=0 {.0\ht0} 0 {.0\ht0} ,width=\textwidth]{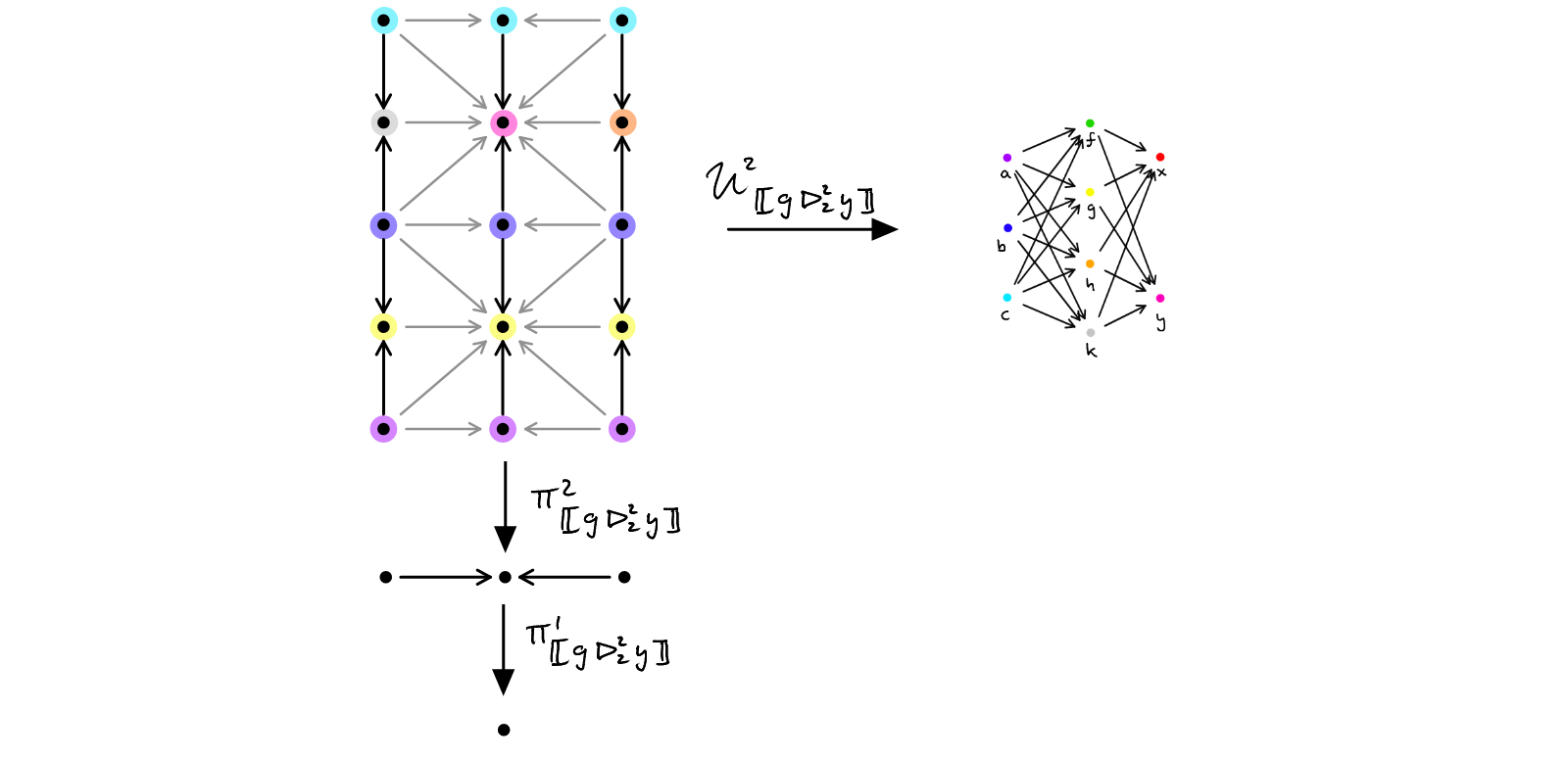}
\endgroup\end{restoretext}
From the previous two generic composites, we infer that $(g \whisker 2 2 y) \whisker 1 2 (x \whisker 2 2 h)$ is a valid generic composite given by
\begin{restoretext}
\begingroup\sbox0{\includegraphics{test/page1.png}}\includegraphics[clip,trim=0 {.0\ht0} 0 {.0\ht0} ,width=\textwidth]{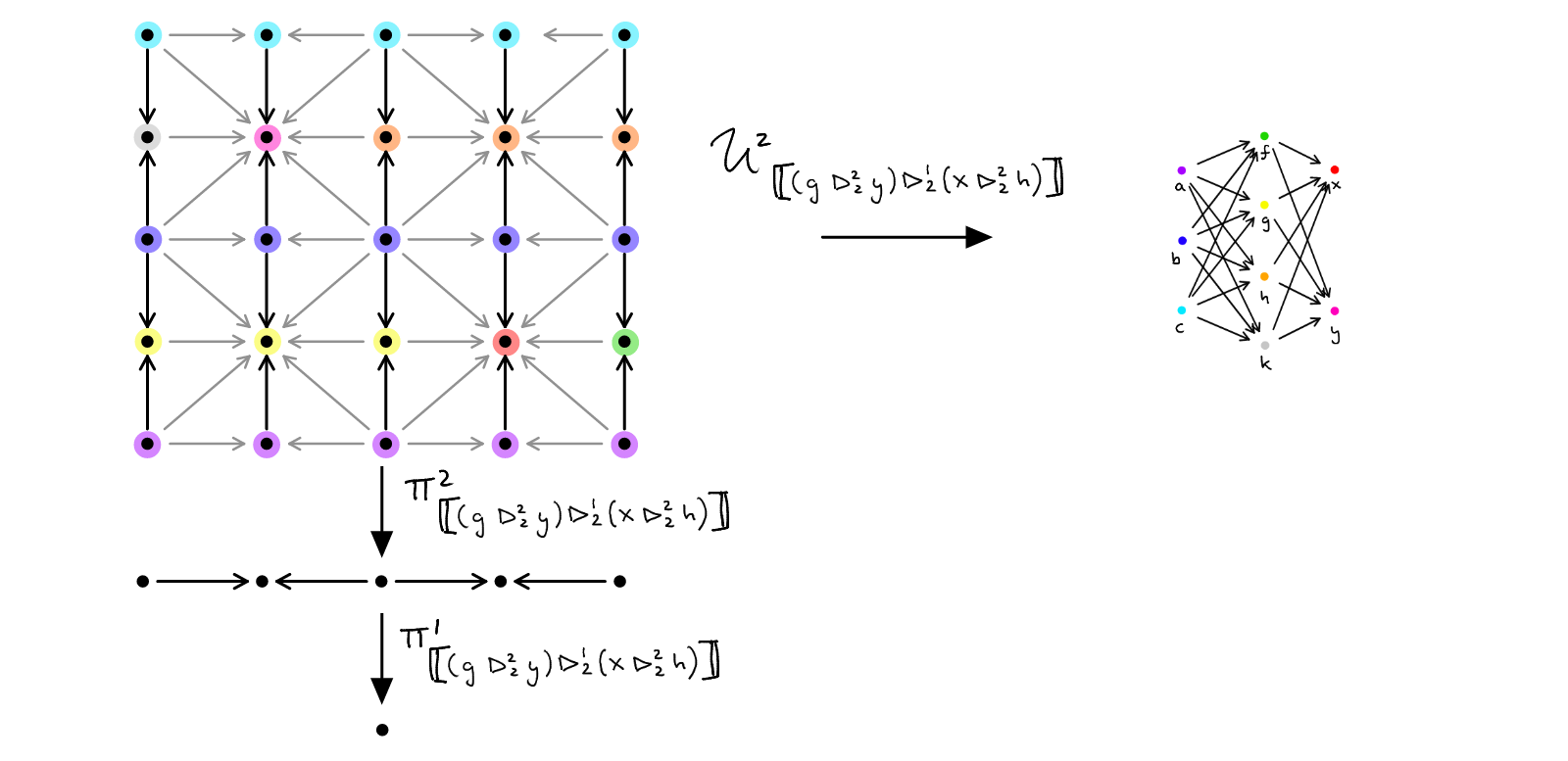}
\endgroup\end{restoretext}
Similarly, $(x \whisker 2 2 h) \whisker 1 2  (g \whisker 2 2 y)$ is a generic composite of $\kC(S)$. 

Now things get slightly more interesting in dimension $3$. Our first generic (in fact, elementary) homotopy is the binary interchange $\interchanger_{x,y}$ (abbreviated by $D$ below)
\begin{restoretext}
\begin{noverticalspace}
\begingroup\sbox0{\includegraphics{test/page1.png}}\includegraphics[clip,trim=0 {.0\ht0} 0 {.0\ht0} ,width=.9\textwidth]{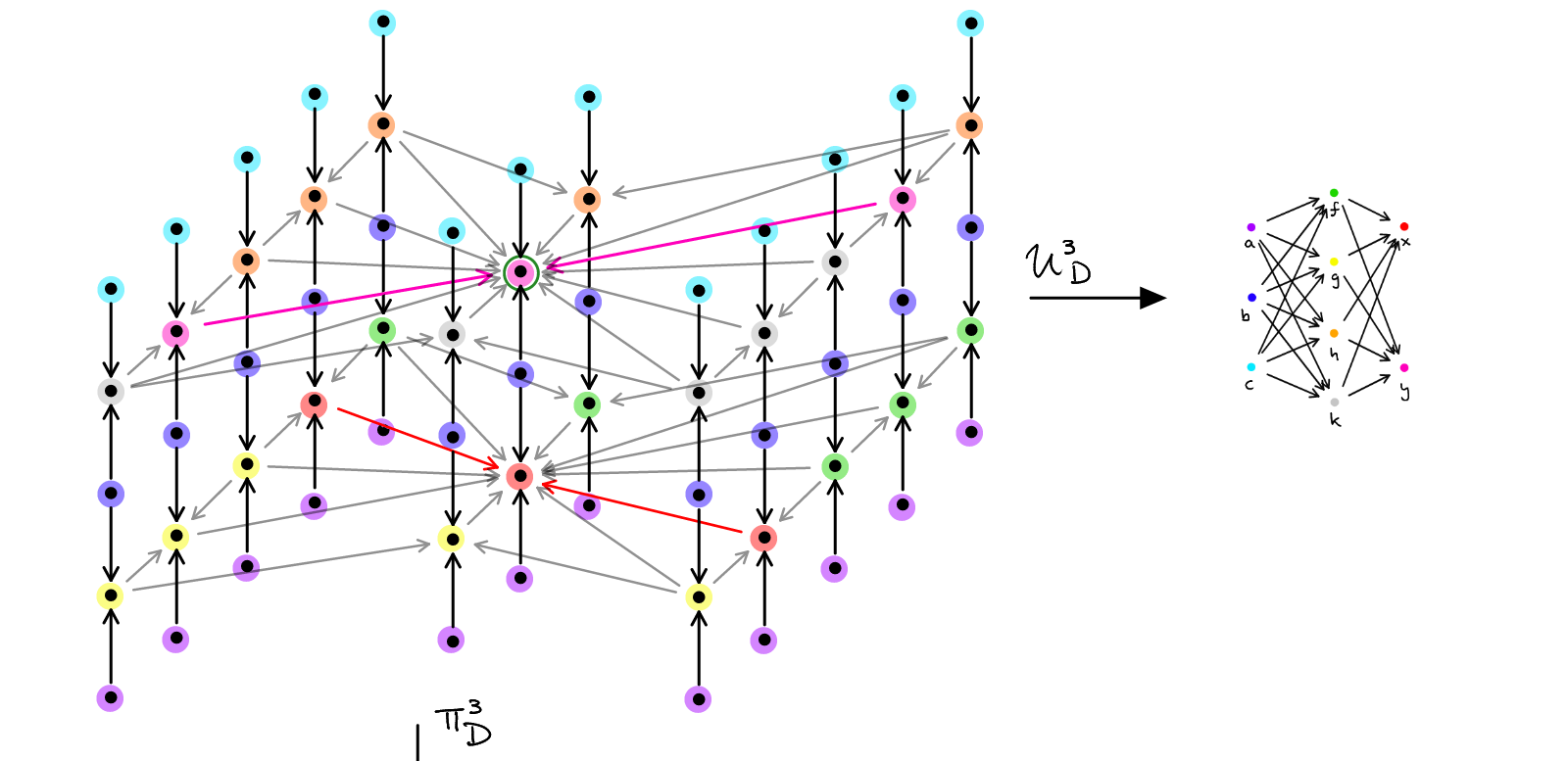}
\endgroup \\*
\begingroup\sbox0{\includegraphics{test/page1.png}}\includegraphics[clip,trim=0 {.0\ht0} 0 {.0\ht0} ,width=.9\textwidth]{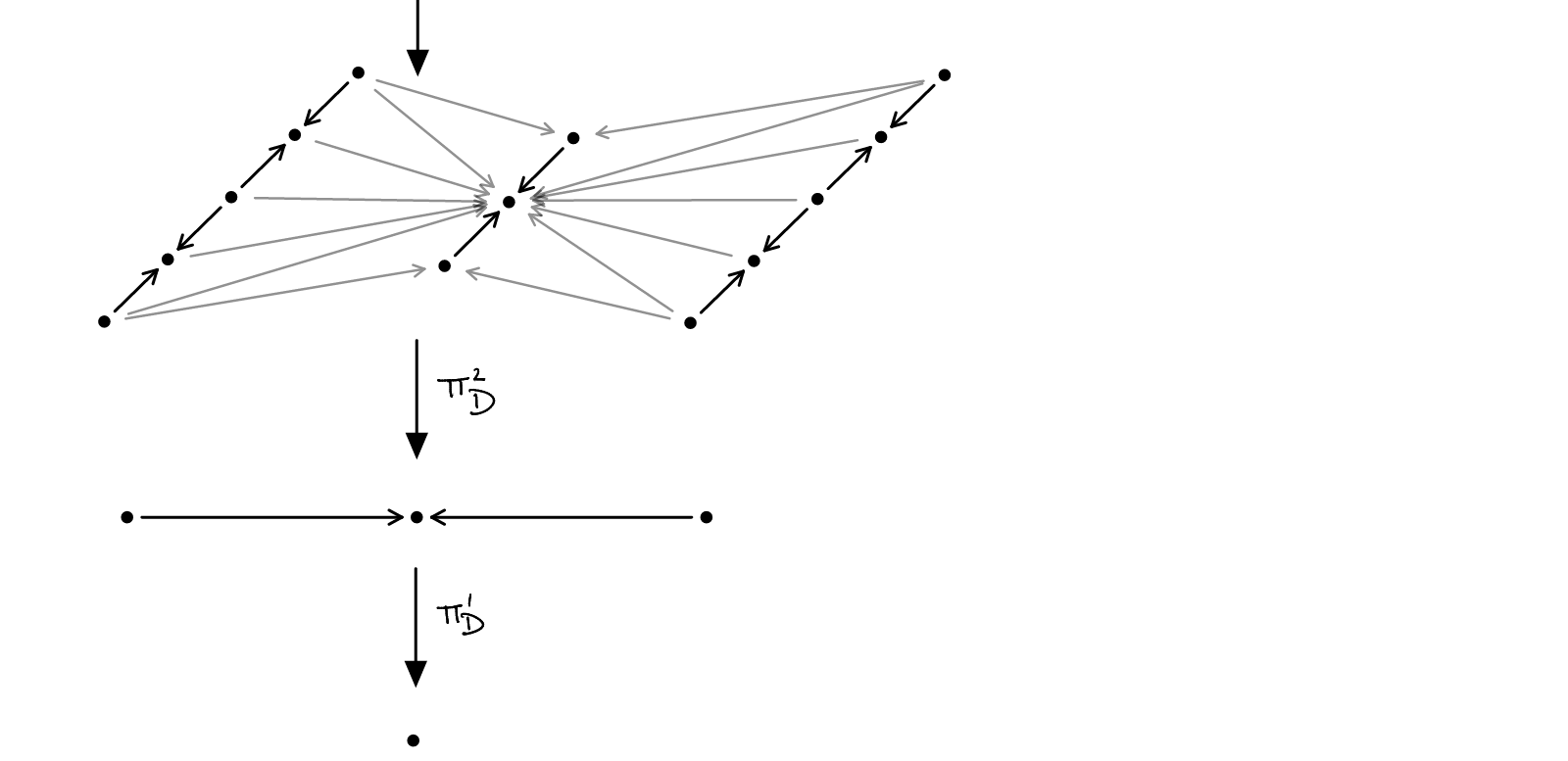}
\endgroup
\end{noverticalspace}
\end{restoretext}
note that its source and target are generic composites, namely
\begin{align}
\gsrc (\interchanger_{x,y}) =  {(g \whisker 2 2 y) \whisker 1 2 (x \whisker 2 2 h)}\\
\gtgt (\interchanger_{x,y}) =  {(x \whisker 2 2 k) \whisker 1 2 (f \whisker 2 2 y)}
\end{align}
and further that it indeed satisfies the condition of being both a homotopy and generic. For the latter note that any proper subcubes, such as
\begin{restoretext}
\begin{noverticalspace}
\begingroup\sbox0{\includegraphics{test/page1.png}}\includegraphics[clip,trim=0 {.0\ht0} 0 {.15\ht0} ,width=.9\textwidth]{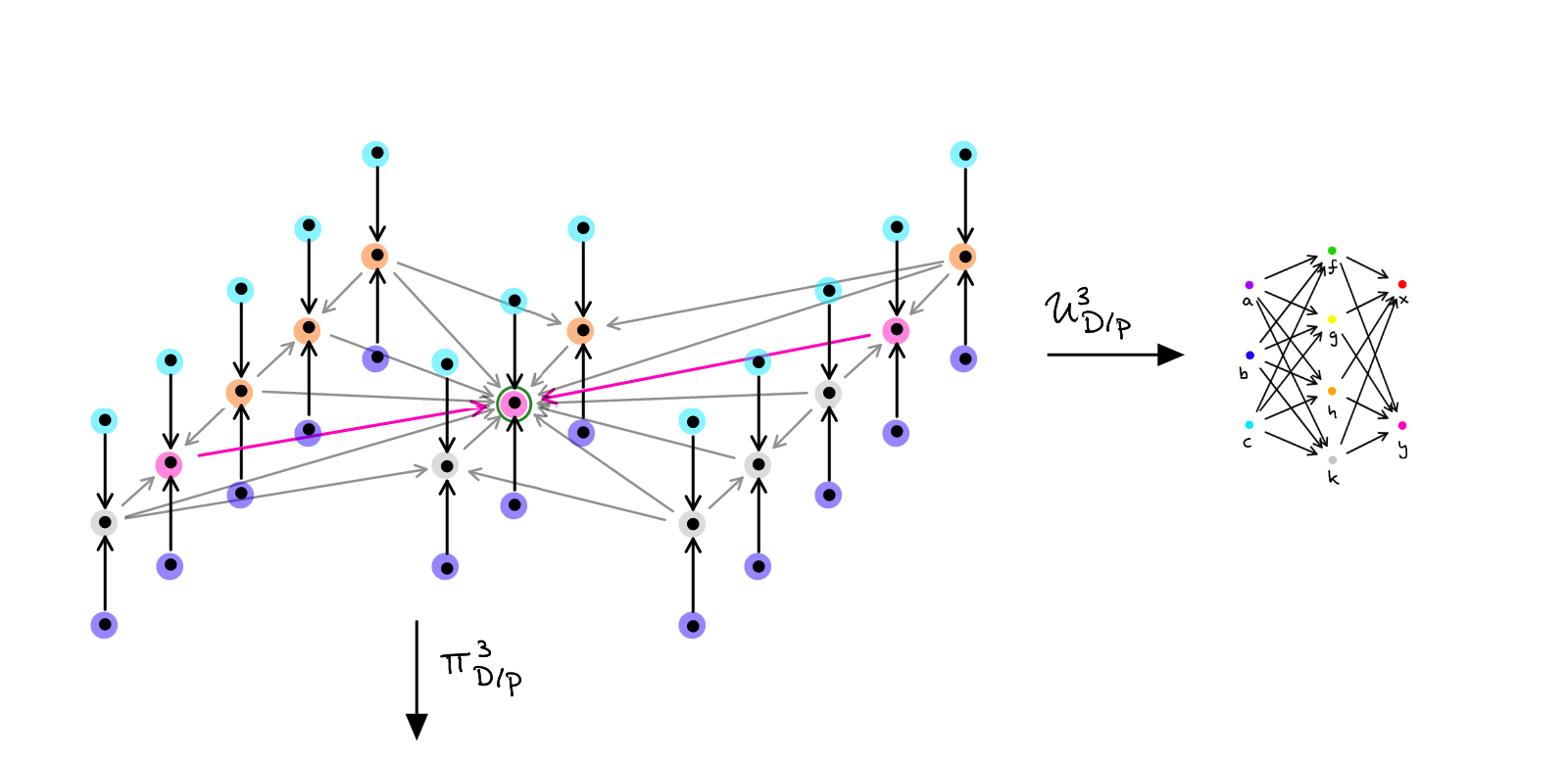}
\endgroup \\*
\begingroup\sbox0{\includegraphics{test/page1.png}}\includegraphics[clip,trim=0 {.0\ht0} 0 {.0\ht0} ,width=.9\textwidth]{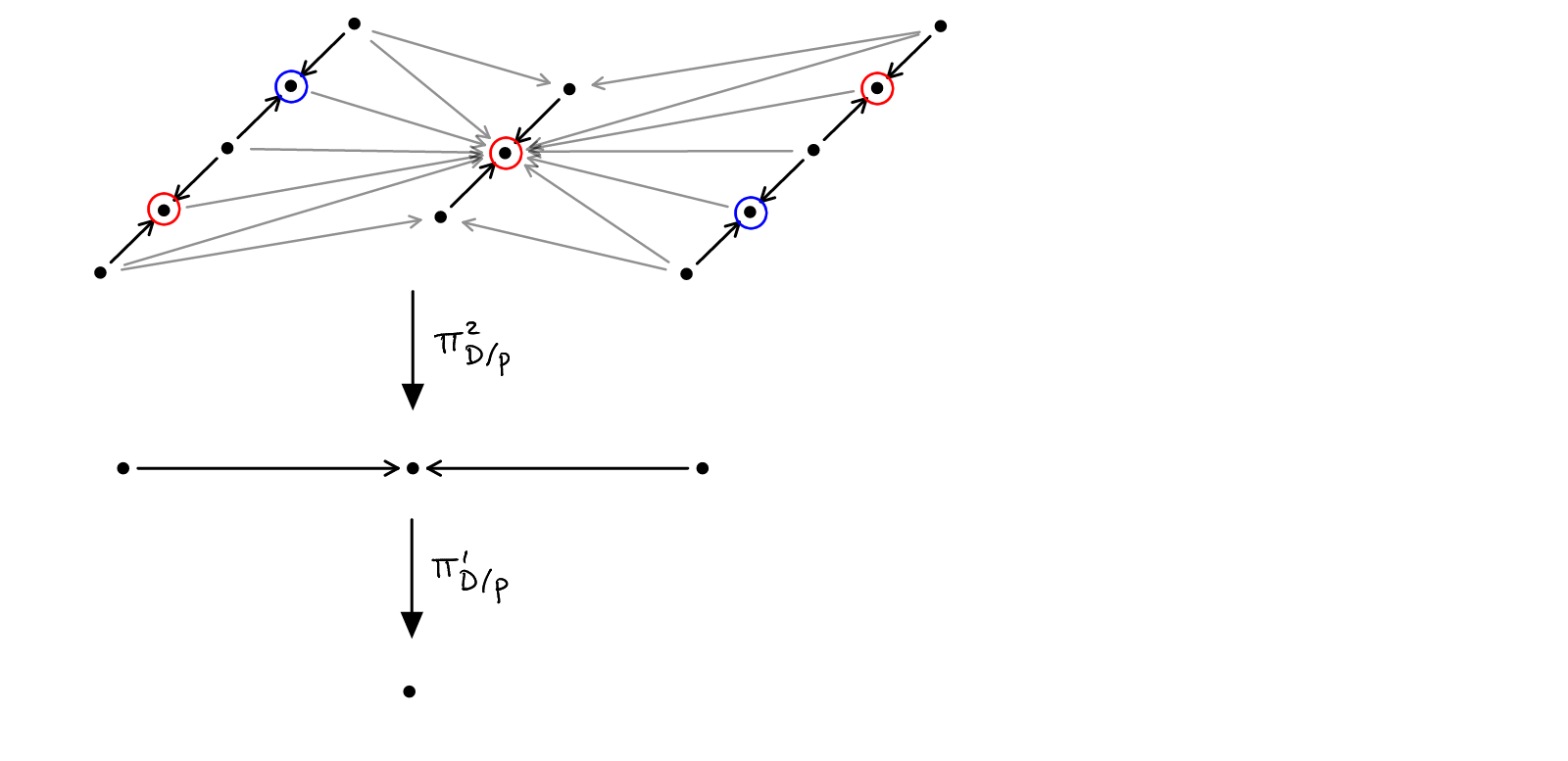}
\endgroup
\end{noverticalspace}
\end{restoretext}
is homotopically trivial (in fact, by the identity homotopy), as defined in  \autoref{notn:non_identity_subcubes}. The only homotopically non-trivial subcube is the entire cube. This shows genericity. A similar argument shows local triviality: namely, each minimal neighbourhood (such as the one above around the point circled in \cgreen{}) normalises to the identity, and this verifies that $\interchanger_{x,y}$ is indeed a homotopy. 

\begin{notn}[Theory of interchange] \label{notn:interchanger} Our chosen $S$ is in fact the ``theory of interchange", and every globular map $\alpha : S \to S'$ exhibits two interchangeable elements in $S'$ (as the images of $x$ and $y$). The resulting generic composite in $\kC(S')$ will be denoted by $\interchanger_{\alpha(x), \alpha(y)}$ (obtained by relabelling $\interchanger_{x,y}$ using $\alpha$ accordingly).
\end{notn}

The following are morphisms but not generic composites. First consider
\begin{restoretext}
\begingroup\sbox0{\includegraphics{test/page1.png}}\includegraphics[clip,trim=0 {.0\ht0} 0 {.0\ht0} ,width=\textwidth]{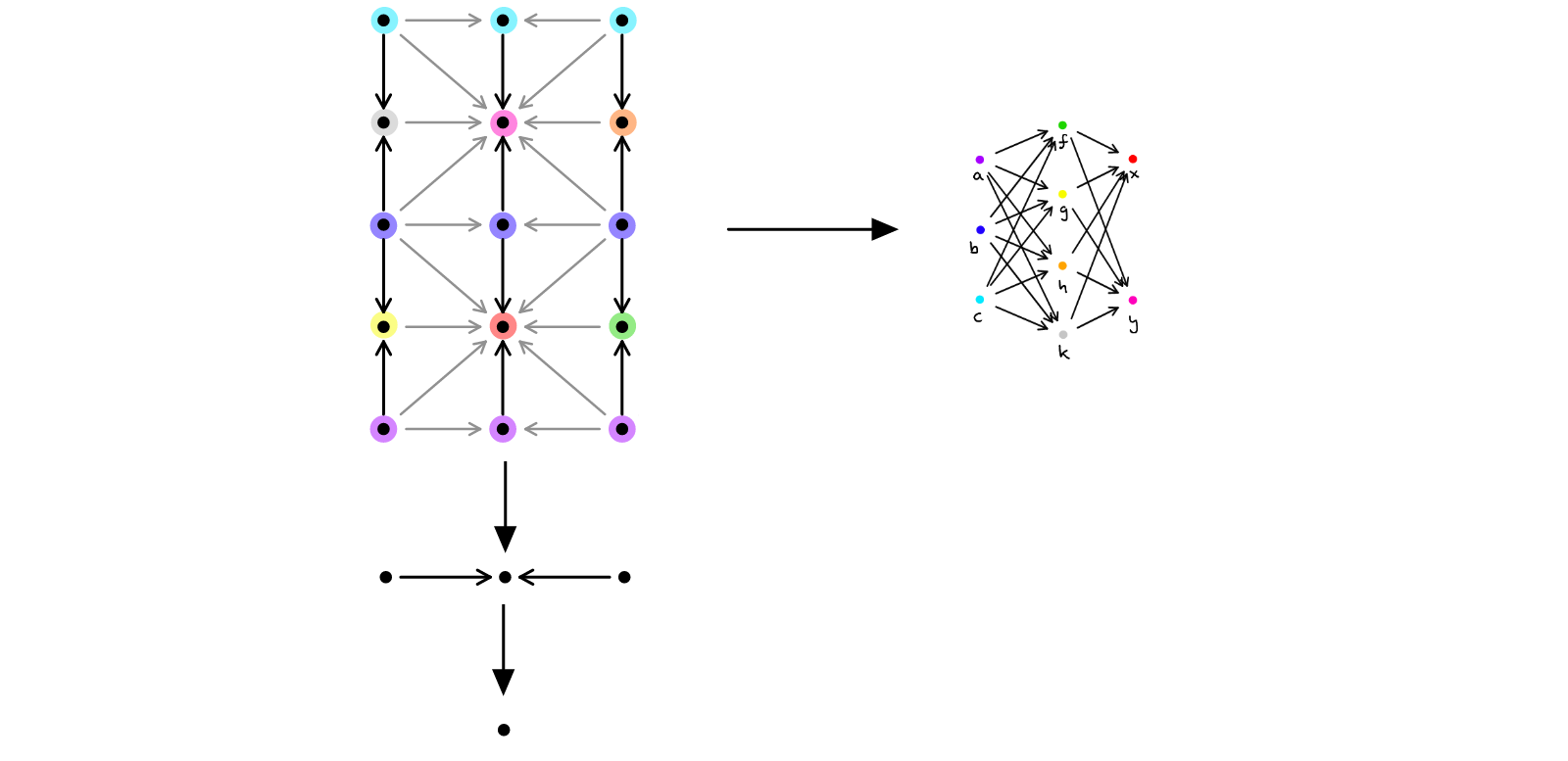}
\endgroup\end{restoretext}
which cannot be obtained by a whiskering operation. 

Next, the ``ternary interchange"
\begin{restoretext}
\begingroup\sbox0{\includegraphics{test/page1.png}}\includegraphics[clip,trim=0 {.2\ht0} 0 {.2\ht0} ,width=\textwidth]{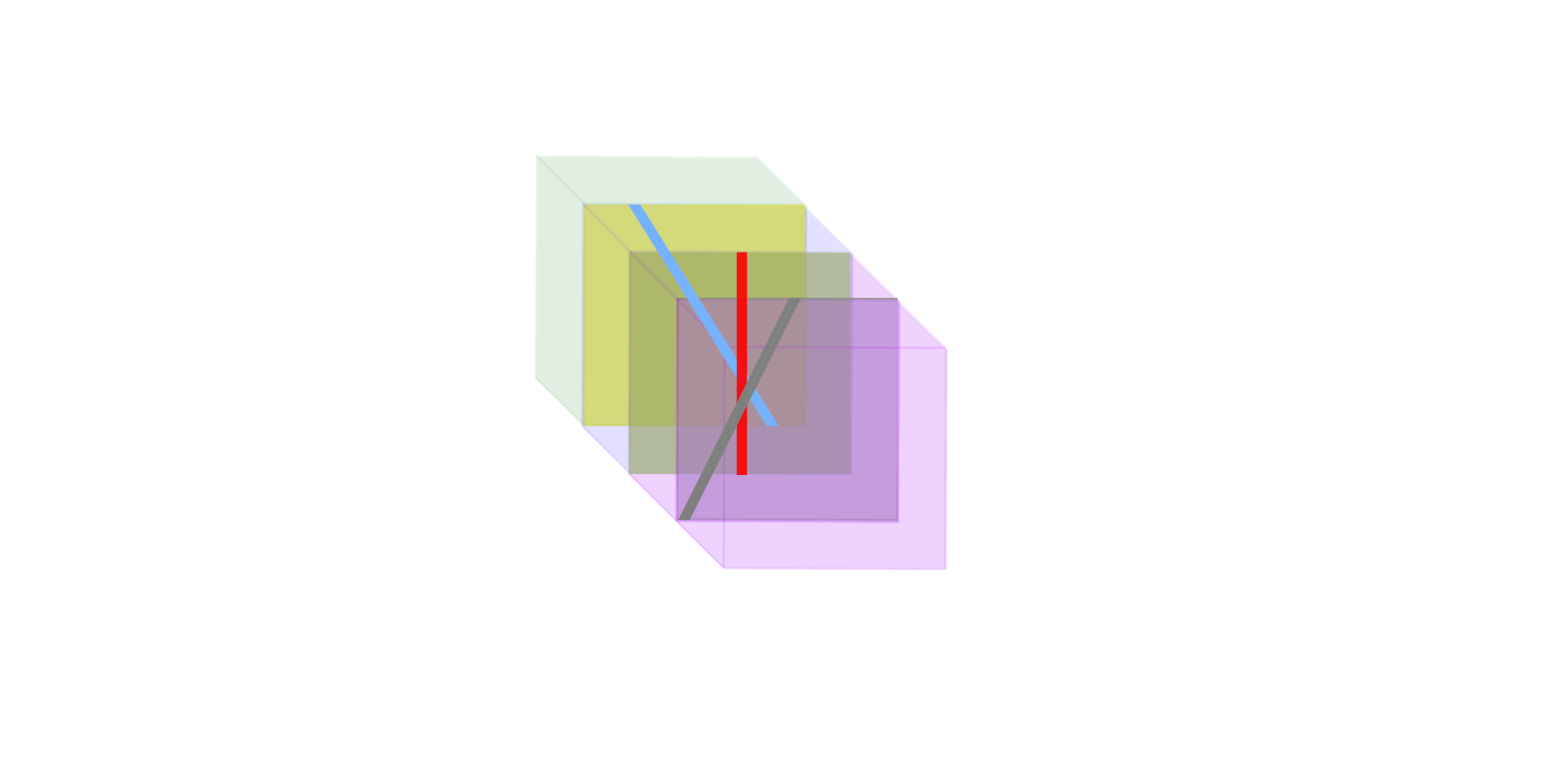}
\endgroup\end{restoretext}
is a homotopy but it is not generic (as defined in \autoref{defn:top_down_gen_comp}) since it contains the binary interchange as proper subcubes in two different ways, and these two subcubes (both not homotopic to the identity) do not factor through each other. Note however that the ternary interchange is equivalent (i.e. homotopic) to a generic composite only containing binary interchanges, namely by the 4-homotopy
\begin{restoretext}
\begingroup\sbox0{\includegraphics{test/page1.png}}\includegraphics[clip,trim=0 {.1\ht0} 0 {.2\ht0} ,width=\textwidth]{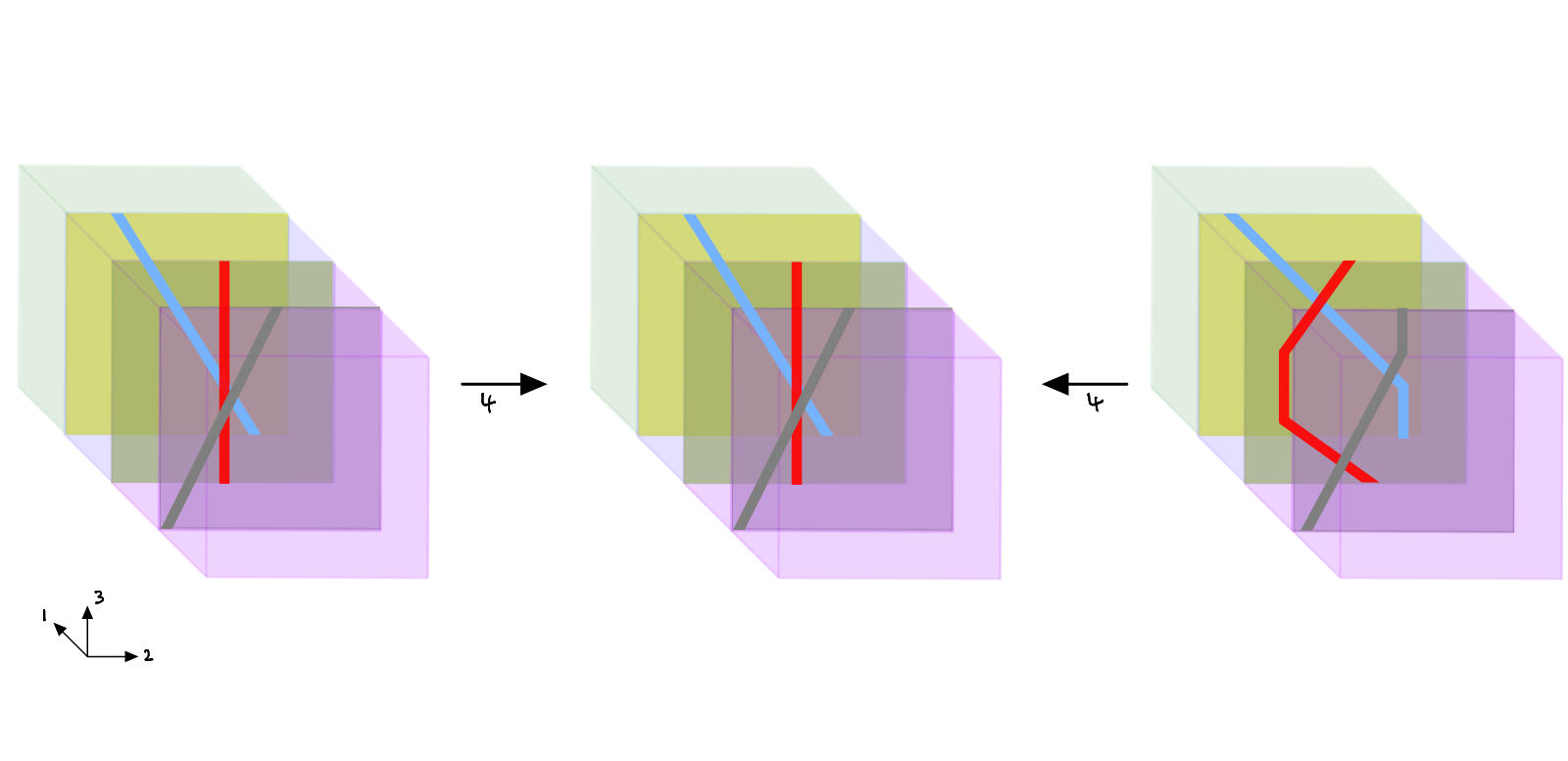}
\endgroup\end{restoretext}
To be fully explicit at this stage we also indicate the combinatorial structure of this 4-homotopy in the following tower of posets
\begin{restoretext}
\begin{noverticalspace}
\begingroup\sbox0{\includegraphics{test/page1.png}}\includegraphics[clip,trim=0 {.0\ht0} 0 {.0\ht0} ,width=\textwidth]{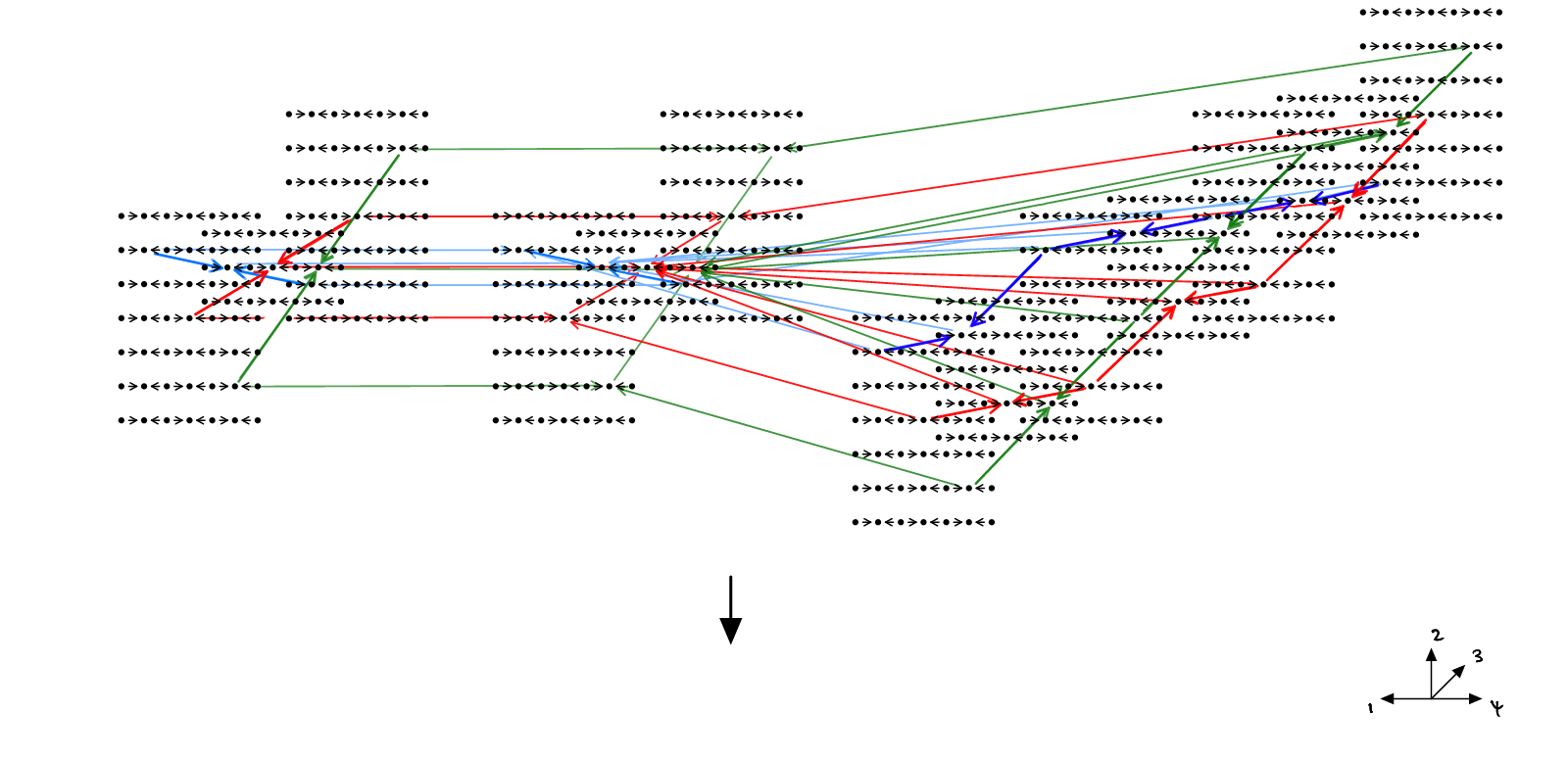}
\endgroup \\*
\begingroup\sbox0{\includegraphics{test/page1.png}}\includegraphics[clip,trim=0 {.15\ht0} 0 {.0\ht0} ,width=\textwidth]{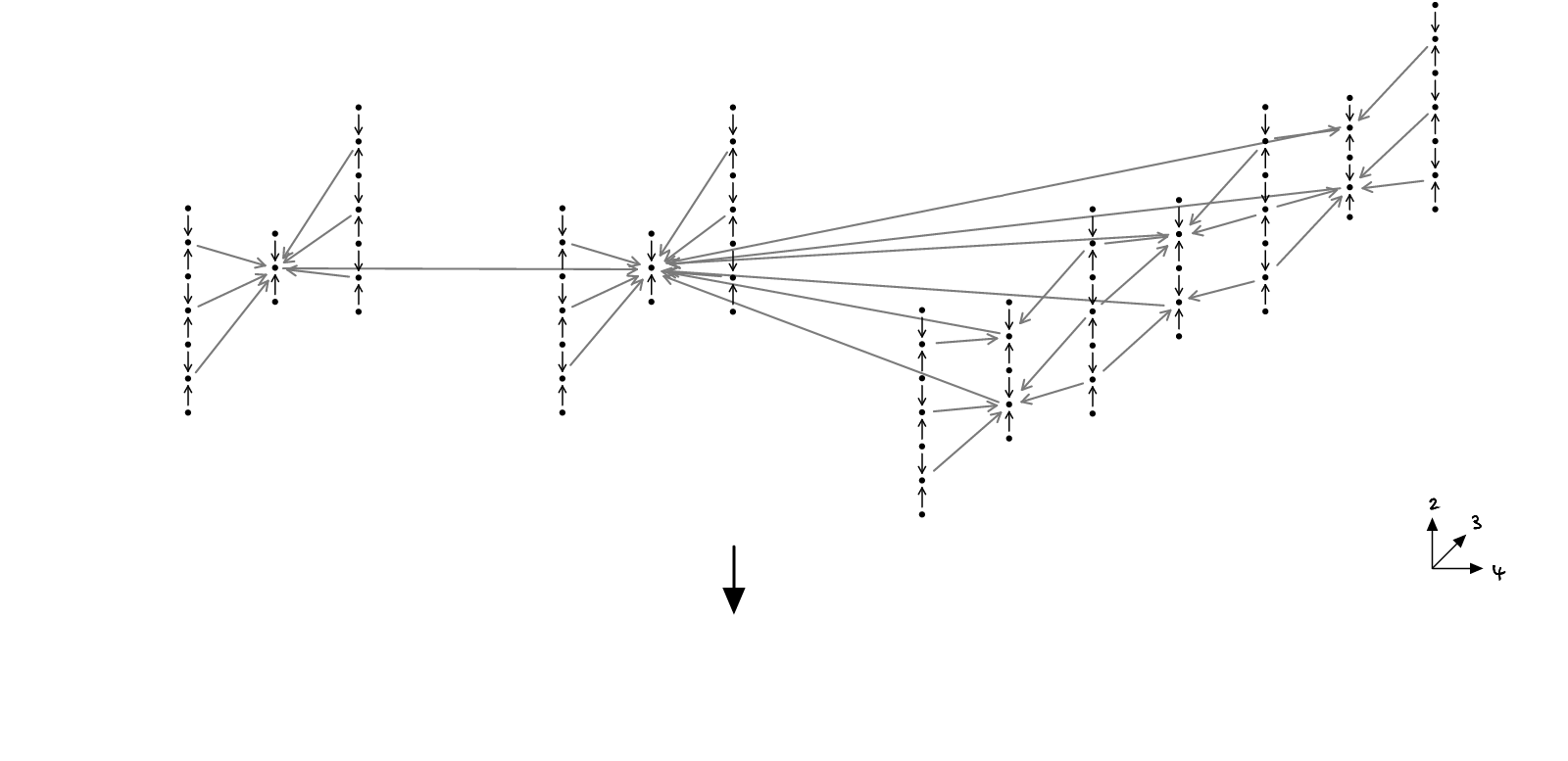}
\endgroup \\*
\begingroup\sbox0{\includegraphics{test/page1.png}}\includegraphics[clip,trim=0 {.0\ht0} 0 {.0\ht0} ,width=\textwidth]{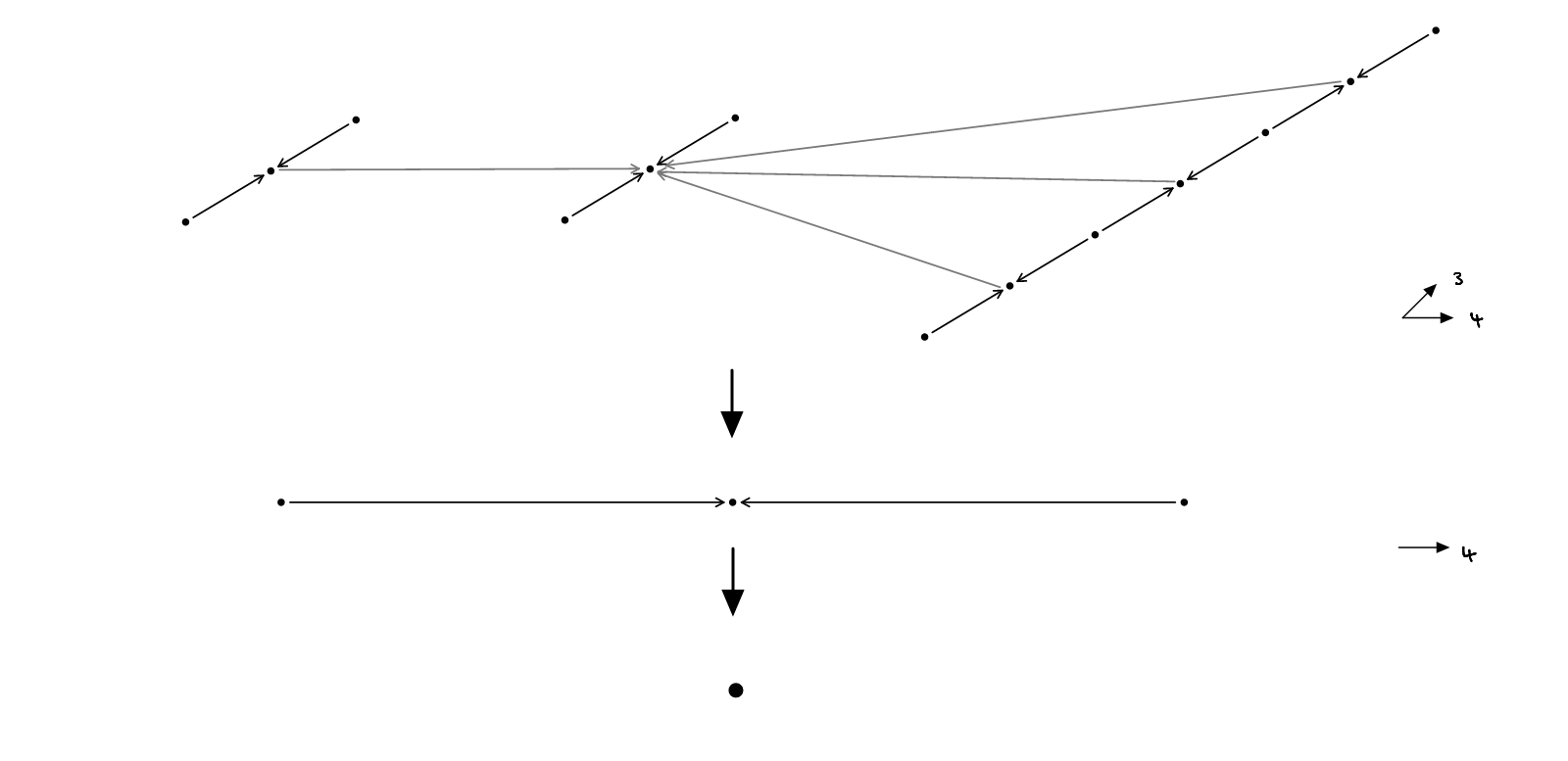}
\endgroup
\end{noverticalspace}
\end{restoretext}

Similarly, the homotopy
\begin{restoretext}
\begingroup\sbox0{\includegraphics{test/page1.png}}\includegraphics[clip,trim=0 {.27\ht0} 0 {.18\ht0} ,width=\textwidth]{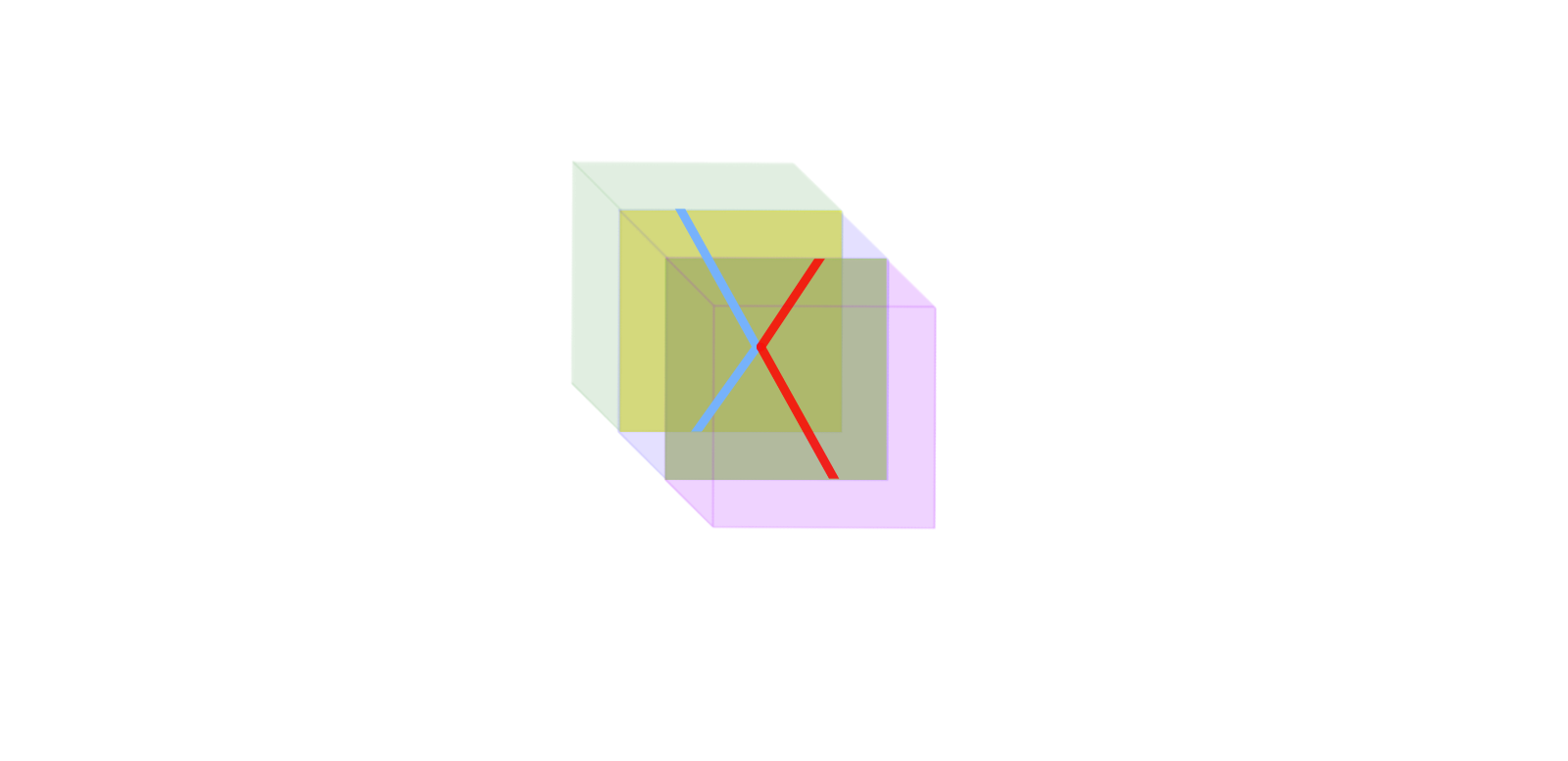}
\endgroup\end{restoretext}
is not generic either: it doesn't contain any subcube that is not homotopic to the identity. As before it is, however, equivalent to a generic composite (namely, the identity homotopy), as the following 4-homotopy shows
\begin{restoretext}
\begingroup\sbox0{\includegraphics{test/page1.png}}\includegraphics[clip,trim=0 {.1\ht0} 0 {.2\ht0} ,width=\textwidth]{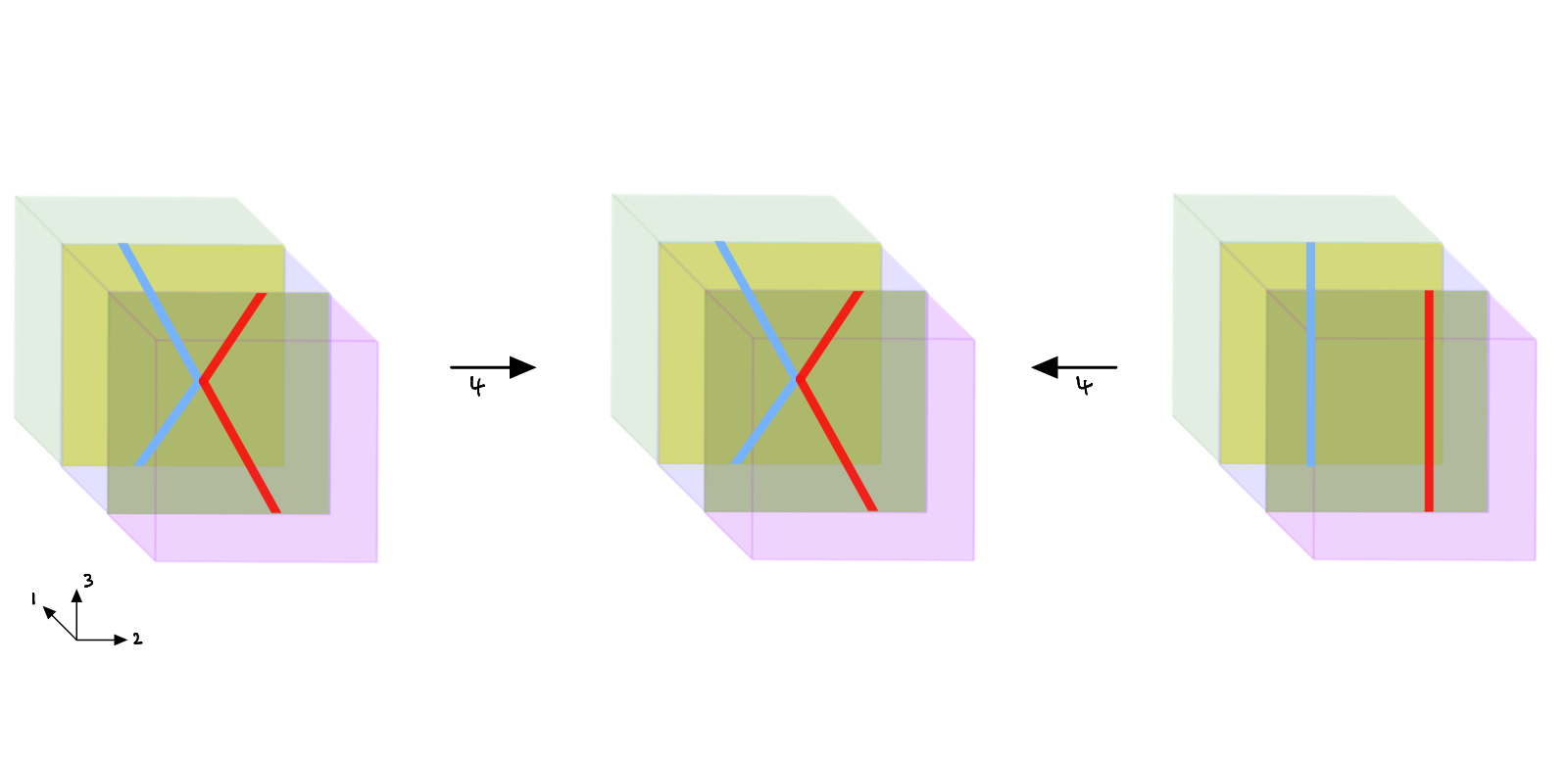}
\endgroup\end{restoretext}

To summarise, this brief discussion of examples was supposed to illustrate the following point: generic composites are a ``dense" subclass of all composites, in the sense that any composite is homotopic to a generic one. Morally, this is because a ``perturbation of a diagram" is a homotopy itself: thus, given a non-generic composite, we can perturb it until all manifolds are in generic position which constructs the required homotopy.

\subsection{Claim of equivalence of characterisations}

We make the following claim without giving its full proof, as this will not be needed for our further discussion.

\begin{claim} Let $\sC \in \pCat_\infty$. Then (cf. \autoref{defn:top_down_gen_comp}, \autoref{constr:bottom_up_gcomp})
\begin{equation}
\GComp(\sC) = \GComps(\sC)
\end{equation}

\proof[Proof sketch] For
\begin{equation}
f \in \GComps(\sC)_n \imp f \in \GComp(\sC)_n
\end{equation}
we argue inductively on the dimension of $f$. First assume $f$ has a single singular height (that is, $\sG^1(f) = \singint 1$). Then $f_1$ in $\theta_1 : f_1 \to f$ (cf. \autoref{defn:top_down_gen_comp}) can be shown to be either the type of an $n$-generator or an elementary homotopy. Further, the remainder of $f$ can be $k$-level whiskered onto $f_1$, proceeding step by step for $k = n, n-1, ... ,3,2$,  to yield $f$. Using the inductive assumption, this shows that $f$ can be written as an element in $\GComp(\sC)$. For general $f$, we can proceed in this way for each singular height each time obtaining an element in $\GComp(\sC)$, and then construct $f$ by $1$-level whiskering of those elements, showing $f$ itself is also an element of $\GComp(\sC)$.

Conversely, for 
\begin{equation}
f \in \GComp(\sC)_n \imp f \in \GComps(\sC)_n
\end{equation}
we argue inductively on dimension and the structure of $f$. If $f$ is elementary than the statement follows trivially. If $f$ is non-elementary, then one can show that whiskering preserves the conditions of \autoref{defn:top_down_gen_comp}. \qed
\end{claim}

\chapter{Foundation-dependent higher categories} \label{ch:associative}

In this chapter we will give a definition of associative $n$-categories and discuss certain variations of it. While the definition we give allows to deduce many higher categorical identities in a uniform and elegant fashion, we do not attempt to provide rigorous evidence for its correctness: the advance brought by the framework developed in this thesis lies more in providing the language to formulate such  ``associative" definitions of higher categories rather than in the definitions themselves.

Our definition will be a straight-forward step based on the work we did in the previous chapters on generic composites and coherent invertibility. Namely, an associative $n$-category will be an algebra $(\sC, \kM_\sC)$ to the generic composites endofunctors $\GComp : \globset \to \globset$ satisfying equations. These equations are encoded as follows. Each $\GComp$-algebra has a (coherent) resolution: this is the presented $\infty$-category obtained by replacing every equation of the algebra multiplication by a coherently invertible generator. Now an associative $n$-category $\sC$ is a $\GComp$-algebra such that its composition operation $\kM_\sC$ extends to generic composites $\GComp(\kiC{}(\kM_\sC))$ of its resolution in a way that is compatible with the inductive structure of $\GComp$, which was explored in-depth in the previous chapter.

In fact, we will find that the idea of resolutions can be iterated by dropping the structural compatibility condition, and replacing it with ``weak" composition operation data on the resolution itself. This iterative idea will allows us to give a ``spectrum" of definitions ranging from ``fully associative" to ``fully weak".

After defining associative $n$-categories (and their ``weaker" counterparts) in \autoref{sec:ANC_def}, we will then use \autoref{sec:ANC_egs} to discuss examples in dimension $n=0,1,2$ and $3$, and we will illustrate general observations about associative $n$-categories by means of these examples.

\section{Definitions} \label{sec:ANC_def}

\subsection{$\mathrm{GComp}$-algebras and their resolution $\mathfrak{R}$}

We will define $\GComp$-algebras as well as their resolution (a similar construction holds for $\Comp$). Resolutions are associated presented associative $\infty$-categories obtained by turning all algebra equations into coherently invertible generators. First note

\begin{notn}[Truncated globular sets] \label{notn:trunc_glob_sets} $n\in \lN\cup\Set{\infty}$. An $n$-truncated globular $S$ is a globular set with elements concentrated in degrees $m \leq n$, that is, $S_m = \emptyset$ for $m > n$. Every globular set $S$ has an $n$-truncation $S_{\truncleq n}$ defined by setting $(S_{\truncleq n})_m$ to be $S_m$ if $m \leq n$ and $\emptyset$ otherwise. Globular maps $s_m, t_m$ are inherited from those of $S$ for $m < n$ and are the empty maps otherwise. 
\end{notn}

Generalising \autoref{ssec:sum_alg_and_inf} we introduce the following definition of an $n$-truncated algebras.

\begin{defn}[Algebras of $\GComp$] A (unital) \textit{$n$-truncated $\GComp$-algebra} $(S,\iM)$ is a tuple of an $n$-truncated $S \in \globset$ and $\iM : \GComp(S)_{\truncleq n} \to S$ such that
\begin{equation}
\xymatrix{ S \ar[d]_{\abss{-}_S} \ar[dr]^{\id} \\
\GComp(S)_{\truncleq n} \ar[r]_-{\iM} & S }
\end{equation}
\end{defn}

The (coherent) resolution of $\iM$ is then constructed by replacing all algebraic equations $\iM(f) = g$ by coherently invertible generators (apart from those forced by unitality). Concretely, we have the following.

\begin{constr}[Resolutions of $\GComp$-algebras] \label{constr:resolutions} Let $(S,\iM)$ be a $\GComp$-algebra. We construct the resolution $\kiC{}(\iM)$ as the colimit of a sequence
\begin{equation}
\kiC {(0)} (\iM)\into \kiC {(1)} (\iM)\into \kiC {(2)}(\iM) \into ...
\end{equation}
whose components are inductively defined as follows: for $k = 0$, we set
\begin{equation}
\kiC {(0)} (\iM) = \kC(S)
\end{equation}
Now inductively assume we have constructed $\kiC {(k)}  (\iM)$ ($k < n$) such that for each $h\in \Comp(S)_{k-1} \setminus S_{k-1}$ we have generators
\begin{equation}
\infeq {h} {k} \in \kiC {(k)}  (\iM)_k
\end{equation}
As a convention, for $h \in S_{k-1}$ we will understand the symbol $\infeq h k$ as $\Id_h$. Inductively, $\infeq h k$, $h\in \Comp(S)_{k-1}$, satisfy the following.

\begin{itemize}
\item \textit{Inductive property of the construction}: For $0 < i < (k-1)$ set
\begin{align}
s_i &= \abss{(\infeq {\gsrc^i(h)} {k-i})\inv} \\
t_i &= \abss{\infeq {\gtgt^i(h)} {k-i}}
\end{align}
Then $\abss{\infeq h {k}}$ satisfies
\begin{equation}
\gsrc \abss{\infeq h {k}} = (s_1 \whisker 1 {k-1} ... (s_{k-2} \whisker {k-2} {k-1} h \whisker {k-2} {k-1} t_{k-2}) ... \whisker 1 {k-1} t_1)
\end{equation}
and
\begin{equation}
\gtgt \abss{\infeq h {k}} = \abss{\iM(h)}
\end{equation}
$\abss{(\infeq h {k})\inv}$ satisfies the same, with source and target exchanged.
\end{itemize}

\noindent Now we define $\kiC {(k+1)}  (\iM)$. Let $f \in \GComp(S)_k \setminus S_k$. For $0 < i < k$ set
\begin{align}
s_i &= \abss{(\infeq {\gsrc^i(f)} {k+1-i})\inv} \\
t_i &= \abss{\infeq {\gtgt^i(f)} {k+1-i}}
\end{align}
We claim there are $x^f_i \in \GComp(\kiC {(k)}  (\iM))_k$ defined inductively as follows: For $i = k,k-1, ..., 1$ we first set $x^k_{f} = f$ and (as shown below) we can then define
\begin{equation}
x^f_i := (s_i \whisker i k x^{i+1}_{f} \whisker i k t_i)
\end{equation}
For $i = k-1$ this is well-defined by the above definitions of $s_i$, $t_i$. Inductively assume $x^f_i$ to be of the form above, and argue as follows: For $j < i$, using \autoref{lem:gen_comp_src_tgt} we compute
\begin{equation}
\gsrc^{j}(x^f_i) =  s_i \whisker {i-j} {k-j}  \gsrc^{j}(x^f_{i+1}) \whisker {i-j} {k-j} t_i
\end{equation}
which inductively implies
\begin{equation}
\gsrc^{i-1} (x^f_i) = (s_i \whisker 1 {k+1-i} ... (s_{k-1} \whisker {k-i} {k+1-i} \gsrc^{i-1} f \whisker {k-i} {k+1-i} t_{k-1}) ... \whisker 1 {k+1-i} t_i)
\end{equation}
Using the inductive assumption for $\infeq {\gsrc^{i-1} f} {k+2-i}$ this implies
\begin{equation}
\gsrc^{i-1} = \gtgt(s_{i-1})
\end{equation}
Similarly, we find
\begin{equation}
\gtgt^{i-1}(x^f_i) = \gsrc(s_{i-1})
\end{equation}
This means
\begin{equation}
x^f_{i-1} := (s_{i-1} \whisker {i-1} k x^{i}_{f} \whisker {i-1} k t_{i-1})
\end{equation}
is well-defined, completing the inductive construction of $x^f_i$. We define
\begin{equation}
x^f = x^f_1
\end{equation}
Note that by its inductive definition, using \autoref{lem:gen_comp_src_tgt}, and by using inductive assumption for $s_i$, $t_i$, we find
\begin{align}
\gsrc^i(x^f) &= \gsrc(s_i) \\
&= \abss{\iM(\gsrc^i(f))}
\end{align}
and similarly
\begin{align}
\gtgt^i(x^f) &= \gtgt(t_i) \\
&= \abss{\iM(\gtgt^i(f))}
\end{align}
In particular, using globularity of $M_S$ and \autoref{constr:panc_from_glob}, $x^f$ has the same sources and targets as $y^f \in \GComp(\kiC {(k)}  (\iM))_k$ defined by
\begin{equation}
 y^f := \abss{\iM(f)}
\end{equation} 
We can thus define
\begin{equation}
\kiC {(k+1)}  (\iM)  = \kiC {(k)}  (\iM) \igadd {x^f, y^f} \Set{\infeq f {k+1} ~|~ f \in \GComp(S)_k}
\end{equation}
which completes the inductive construction.

Finally, $\kiC {}  (\iM)$ is defined as $\kiC {(n)}  (\iM)$, or if $n = \infty$, as the colimit of 
\begin{equation}
\kiC {(0)}  (\iM) \into \kiC {(1)}  (\iM)\into \kiC {(2)}  (\iM) \into \dots
\end{equation}
\end{constr}

\begin{defn}[Compositional depth] \label{defn:comp_depth} Note that there is a canonical inclusion of presentations $\infinc {\iM} : \kC(S) \into \kiC{} (\iM)$. Elements of $\GComp(\kiC{}(\iM))$ \textit{not} in the image of $\GComp(\infinc {\iM}) : \GComp(\kC(S))_{\truncleq n} \into \GComp(\kiC{}(\iM))_{\truncleq n}$ are called composites of \textit{compositional depth $1$}. Otherwise, composites are said to be of \textit{compositional depth 0}.
\end{defn}

Morally the resolution of $\iM$ is a presented higher category whose composites represent the ``free interactions of algebra laws with elements of the algebra". Composites of compositional depth $1$ are those that involve algebra laws (and their coherences). Depth $0$ composite are nothing but regular composites build from the elements of the algebra.

\begin{rmk}[Algebras and resolutions for $\Comp$] The preceding definition and construction have analogues obtained when replacing $\GComp$ by $\Comp$.
\end{rmk}

\subsection{$m$-Iterated resolutions}

We briefly mention the following generalisation of the construction in the preceding section: Namely, we will define $m$-iterated resolutions of algebras, which involves choosing new multiplication data on a given resolution, when is then used to built the next resolution.

\begin{constr}[Sketch of $m$-iterated resolutions] \label{constr:resolutions_iterated} An $m$-iterated resolution of a ($n$-truncated) $\GComp$-algebra $\iM$ is a commuting diagram of the form
\begin{equation}
\xymatrix@C=2cm{ \kC(S) \ar[d]_{\abss{-}_S} \ar[dr]^{\id} & \\
\GComp(\kC(S))_{\truncleq n} \ar[r]_{\iM^0_S} \ar[d] & S \\
\GComp(\kiC{}(\iM^0_S))_{\truncleq n} \ar[ur]_{\iM^1_S} \ar[d] & \\
\GComp(\kiC{}(\iM^1_S))_{\truncleq n} \ar[uur]_{\iM^2_S} \ar[d] & \\
\cdots \ar[d] & \\
\GComp(\kiC{}(\iM^{m-2}_S))_{\truncleq n} \ar[uuuur]_{\iM^{m-1}_S} \ar[d]  & \\
\GComp(\kiC{}(\iM^{m-1}_S))_{\truncleq n} 
}
\end{equation} 
such that $\iM^0_S = \iM$ and in each step one defines $\kiC{}(\iM^{k}_S)$ to contain $\kiC{}(\iM^{k-1}_S)$ together with invertible generators $\infeq f l$ witnessing all equations $\iM^k_S(f) = g$ for $f \in \GComp(\kiC{}(\iM^{k}_S))_l \setminus  \GComp(\kiC{}(\iM^{k-1}_S)_l$ (such an $f$ is said to be of \textit{compositional depth $(k+1)$}). The types of $\infeq f {l+1}$ are analogous to those in \autoref{constr:resolutions}. One then first \textit{chooses} a map $\iM^{k+1}_S : \Comp(\kiC{}(\iM^{k+1}_S)) \to S$ before constructing $\Comp(\kiC{}(\iM^{k+2}_S))$ in a similar fashion.
\end{constr}

\subsection{Associative higher categories}

We are now in the position to give a definition of associative $n$-categories (which was already given in the case $n = \infty$ in \autoref{ssec:sum_anc}). 

\begin{defn}[Associative $n$-categories] \label{defn:ANC} Let $n \in \lN \cup \infty$. An \textit{associative $n$-category} $\sC$ is an $n$-truncated $\GComp$-algebra denoted by\footnote{Note that we are abusing notation very slightly by re-using the symbol $\sC$ for both the associative $n$-category and the underlying $\GComp$-algebra.} $(\sC, \kM_\sC)$ such that there is an extension $\extkM_\sC$ of $\kM_\sC$ to the resolution of $\sC$, that is,
\begin{equation}
\xymatrix{ \GComp(\sC)_{\truncleq n} \ar[d] \ar[r]^-{\kM_\sC} & \sC \\
 \GComp(\kiC{}(\iM))_{\truncleq n} \ar@{-->}[ur]_-{\extkM_\sC}& }
\end{equation}
uniquely determined subject to the following conditions
\begin{enumerate}
\item \textit{\complaw{} law}: For $f = f_1 \whisker k m f_2 \in  \GComp(\kiC{}(\iM))_{\truncleq n}$ we require\footnote{Note that non-elementarity of $f$ is not required for this to hold}
\begin{equation} \label{eq:composition_law}
\extkM_\sC(f_1 \whisker k m f_2) = \kM_\sC (\abss{\extkM_\sC(f_1)} \whisker k m \abss{\extkM_\sC(f_2)})
\end{equation}

\item \textit{\cohlaw{} law}: For elementary $f \in \GComp(\kiC{}(\kM_\sC))_k$, $k \leq n + 1$, if either $k = n +1$ or $f$ is of compositional depth $1$, then
\begin{equation} \label{eq:coh_law_st}
\extkM_\sC(\gsrc(f))  = g = \extkM_\sC(\gtgt(f)) 
\end{equation}
and (if $k \leq n$)
\begin{equation}
\extkM_\sC(f) = \kM_\sC(\Id_{\abss{g}})
\end{equation}
\end{enumerate}
\end{defn}

Note that generic composites play a crucial role in this definition: Firstly, since every generic composite is either elementary or a whiskering composite, it is possible to determine $\extkM_\sC$ completely by forcing compatibility with whiskering (via the \complaw{} law) and determining its values on elementary composites (via the \cohlaw{} law).

Secondly, there are then two types of elementary composites, namely, those of depth $0$ (in $\GComp(\sC)$) and those of depth $1$ (in $\GComp(\kiC{}(\iM)) \setminus \GComp(\sC)$). The former are either generators of $\sC$ or elementary homotopies in $\sC$. Depth $0$ homotopies represent categorical \textit{coherences}. In general these are \textit{non-strict}, in the sense that $\kM_\sC$ does not necessarily evaluate them to equal identities $\kM_\sC(\Id)$. In contrast, depth $1$ elementary composites are composites involving witnesses of compositions: they are either witnesses of compositions themselves or homotopies (that is, coherences) containing them. Now, by the \cohlaw{} law we require depth $1$ coherences to translate to \textit{strict} equations of their source and target (namely by \eqref{eq:coh_law_st}). In other words, these elementary composites are parametrising exactly the equations that we want to hold strictly in an associative $n$-category $\sC$ (in addition to the equations enforced on non-elementary composites via \eqref{eq:composition_law}).

In the next section we will see how to alternatively chose depth $1$ coherences, and more generally depth $k$ coherences (which are coherences involving witnesses of depth-$(k-1)$ stuff), in a \textit{non-strict} fashion. We also refer the reader to \autoref{ch:weak} for more thoughts on ``depth-$k$" structures, including how weak associators and weak pentagonators arise as such.

We first point out the following technical remark about the preceding definition.

\begin{rmk}[The right-hand side of the \complaw{} law is defined] To see that the whiskering on the right hand side of \eqref{eq:composition_law} exists, note that by existence of $f_1 \whisker k m f_2$ we can distinguish the following two cases
\begin{enumerate}
\item In the first case we have $\gtgt^k (f_1) = \gsrc (f_2)$, in which case we compute
\begin{restoretext}
\begin{alignat}{3}
\gtgt^k \abss{ {\extkM_\sC (f_1)}} &= \abss{ {t_{n,k} \big(\extkM_\sC (f_1)\big)}} \hspace{0.5cm} && \text{\autoref{constr:panc_from_glob}, \autoref{notn:glob_set_src_tgt}}\\
&=  \abss{\extkM_\sC ( \gtgt^k (f_1))} && \text{$\extkM_\sC$ is globular} \\
&=  \abss{\extkM_\sC ( \gsrc (f_2))} && \text{case assumption} \\
&=  \abss{ s_n \big(\extkM_\sC (  (f_2) )\big)} &&  \\
&=  \gsrc \abss{\extkM_\sC (  (f_2))}  &&
\end{alignat}
\end{restoretext}
Thus, the condition for whiskering is satisfied.

\item In the second case we have $\gtgt (f_1) = \gsrc^k (f_2)$, in which case we can argue analogously.
\end{enumerate}
\end{rmk}

We will later discuss the preceding definition in low dimensions. We can summarise our discussion as follows. 
\begin{itemize}
\item An associative $0$-category is a set
\item An associative $1$-category is an unbiased category
\item An associative $2$-category is an unbiased strict $2$-category
\item An associative $3$-category is an unbiased $\mathbf{Gray}$-category
\end{itemize}
Before discussing examples, we first describe alternatives to our chosen definition.

\subsection{A spectrum of ``fully associative" to ``fully weak" definitions}

As remarked in the previous section, coherences in higher categories appear at different ``compositional depths". In the previous definition we chose all coherences to be strict starting at depth $1$. However, the following definition sketch shows that we could also require additional structure on our higher category corresponding to weak instances of those coherences.

\begin{defn}[Higher categories, weak up to depth $m$] Let $m,n \in \lN \cup \Set{\infty}$. A \textit{``weak up to depth $m$" $n$-category} $\sC$  consists of a globular set (which we denote again by $\sC$), together with maps $\kM^i_\sC$ ($0 < i < m$) of globular sets with the property that if $n < \infty$ there exists $\extkM_\sC$ in (cf. \autoref{constr:resolutions_iterated})
\begin{equation}
\xymatrix@C=2cm{ \kC(S) \ar[d]_{\abss{-}_\sC} \ar[dr]^{\id} & \\
\GComp(\kC(S))_{\truncleq n} \ar[r]^{~\kM^0_\sC = \kM_\sC} \ar[d] & S \\
\GComp(\kiC{}(\kM^0_\sC))_{\truncleq n} \ar[ur]_{\kM^1_\sC} \ar[d] & \\
\GComp(\kiC{}(\kM^1_\sC))_{\truncleq n} \ar[uur]_{\kM^2_\sC} \ar[d] & \\
\cdots \ar[d] & \\
\GComp(\kiC{}(\kM^{m-2}_\sC))_{\truncleq n} \ar[uuuur]|{\kM^{m-1}_\sC} \ar[d]  & \\
\GComp(\kiC{}(\kM^{m-1}_\sC))_{\truncleq n} \ar@{-->}[uuuuur]_-{\extkM_\sC}
}
\end{equation} 
which is fully determined subject to the \complaw{} law and the \cohlaw{} law from \autoref{defn:ANC} (the latter now applying to elementary $f$ of depth $m$).
\end{defn}

Unlike in associative $n$-categories, $\kM^k_\sC$ are now choices of data: Indeed, previously compositional and coherence associativity laws determined such choices, making the higher category associative. Now $\kM^k_\sC$ pick out coherence cells in $\sC$, which makes $\sC$ \textit{weaker}, in the sense of requiring more coherence data to be chosen. In fact, it is possible from this to retrieve usual coherence cells, such as associators and pentagonators from the above iterated resolutions, which is illustrated in \autoref{ch:weak}.

If $m = \infty$, we say $\sC$ is a \textit{fully} weak $n$-category. Note that in this case, passing to the colimit of the vertical chain of presentation inclusions, we obtain a globular set $\Comp(\kiC{}(\kM^\infty_\sC))$ and a map
\begin{equation}
\kM^\infty_\sC : \Comp(\kiC{}(\kM^\infty_\sC)) \to \sC
\end{equation}
Morally, $\Comp(\kiC{}(\kM^\infty_\sC))$ contains all composites generated by  the generators of $\sC$, the witness of their composition (depth $1$), and the iterated depth $k$ witnesses of the compositions for any $k$. This can be thought of ``bracketings of arbitrary depth" (such as $((...(...)...(...)...))$) now being accounted for in the composition operation $\kM^\infty_\sC$.

\subsection{Relation to presented associative $n$-categories}  \label{ssec:pres_conj}

We very briefly touch on the relation between presented associative $n$-categories and associative $n$-categories. The notion of presented associative $n$-category is a red herring: it is not immediately clear that a presented associative $n$-category is an associative $n$-category. This observation is related to \autoref{rmk:monad_fail}.

\begin{conj}[\Free{} categories are categories] Given a \free{} associative $n$-category $\sC$ then $\Comp(\sC)^{\quotg}_{\truncleq n}$ naturally admits structure of an associative $n$-categories (up to certain arbitrary choices for $\kM_\sC$ on homotopies).
\end{conj}

\begin{rmk}[Weakifications] A different but less direct approach of possibly deriving an associative $n$-category from a presented associative $n$-category, is to work with the globular subset of terminal-globe-like generator types of the weakification (based on techniques mentioned in \autoref{ch:weak}) of a presented associative $n$-categories.
\end{rmk}

\section{Examples} \label{sec:ANC_egs}

In this section we describe examples of \autoref{defn:ANC} in low dimensions. In each dimension ($> 0$), we will proceed as follows to define an associative $n$-category $\sC$
\begin{enumerate}
\item \textit{Morphism data}: We describe the underlying globular set of $\sC$
\item \textit{Composition data}: We define the value of $\kM_\sC$ on $\whisker k n$ operations of two (higher) morphisms in $\sC$ and on each elementary homotopy. By the \unitlaw{} law, as well as the \complaw{} and \cohlaw{} laws this then fully determines $\kM_\sC$
\item \textit{Properties of composition}: We discuss implications of the \complaw{} and \cohlaw{} laws, and how $\kM_\sC$ satisfies them
\end{enumerate}

\begin{notn}[Omitting $\abss{-}$ for types] To simplify notation, for $a \in \sC$, $\sC \in \globset$, we will write $a$ in place of $\abss{a}$ for the type of $a$ in $\kC(\sC)$.
\end{notn}

\subsection{Dimension 0}

An associative $0$-category $\sC$ is a $0$-truncated $\GComp$ satisfying the \complaw{} and \cohlaw{} laws. Since $\GComp(\sC)_0 \iso \sC$ canonically,  $\kM_\sC$ is fully determined by the  \unitlaw{} law. In other words an associative $0$-category is nothing but a set. Note also that $\extkM_\sC = \kM_\sC$ in \autoref{defn:ANC} in this case.

\subsection{Dimension 1}

An associative $1$-category $\sC$ is a $1$-truncated $\GComp$-algebra satisfying the \complaw{} and \cohlaw{} laws. As an example consider the following

\begin{enumerate}
\item \textit{Morphism data}: Let $\sC$ be the globular set with
\begin{align}
\sC_0 &= \Set{a,b} \\
\sC_1 &= \Set{k,g,h,\id_a,\id_b}
\end{align}
such that
\begin{align}
s(k) = t(k) = s(g) = s(h) = s(\id_a) = t(\id_a) &= a \\
 t(g) = t(h) = s(\id_b) = t(\id_b) &= b
\end{align}

\item \textit{Composition data}: We define $\kM_\sC$ on elementary homotopies by setting
\begin{align}
\kM_\sC(\Id_{a}) &= \id_a \\
\kM_\sC(\Id_{b}) &= \id_b \\
\end{align}
(note that identities are the only elementary homotopies in this dimension). We further define $\kM_\sC$ on whiskerings of elementary composites as follows
\begin{align}
\kM_\sC(k \whisker 1 1 k) &= k\\
\kM_\sC(k \whisker 1 1 g) &= h\\
\kM_\sC(k \whisker 1 1 h) &= h\\
\kM_\sC(k \whisker 1 1 \id_a) &= k\\
\kM_\sC(\id_a \whisker 1 1 g) &= g\\
\kM_\sC(\id_a \whisker 1 1 h) &= h\\
\kM_\sC(\id_a \whisker 1 1 k) &= k\\
\kM_\sC(g \whisker 1 1 \id_b) &= g\\
\kM_\sC(h \whisker 1 1 \id_b) &= h
\end{align}
As we will now see, only the first three of those equations need to be specified. The rest follows since $\id_x$ as been identified with the homotopy $\Id_{x}$ by $\kM_\sC$ ($x \in \Set{a,b}$). In subsequent examples we will thus not specify the whiskering composites with morphisms $\kM_\sC(\Id_x)$ (which will be consistently named $\id_x$).

\item First note that
\begin{align}
(x \whisker k n y) \whisker k n z &= (x \whisker k n y) \whisker k n z \\
&=: x \whisker k n y \whisker k n z 
\end{align}
From the \complaw{} it then follows that e.g.
\begin{align}
\kM_\sC(k \whisker 1 1 k \whisker 1 1 g) &= \kM_\sC(\kM_\sC(k \whisker 1 1 k) \whisker 1 1 \kM_\sC(g))
= \kM_\sC(\kM_\sC(k) \whisker 1 1 \kM_\sC(k \whisker 1 1 g))
\end{align}
and thus
\begin{equation}
\kM_\sC(k \whisker 1 1 g) = \kM(k \whisker 1 1 h)
\end{equation}
which is indeed satisfied by our choice for $\kM_\sC$. More generally, this argument shows that composition is associative.

Composition is also necessarily unital, since for instance
\begin{align}
\kM_\sC(g) &= \kM_\sC(\Id_a \whisker 1 1 g) \\
&= \kM_\sC(\kM_\sC(\Id_{a}) \whisker 1 1 \kM_\sC(g))
\end{align}
\end{enumerate}

In summary, $\kM_\sC$ gives $\sC$ the structure of an unbiased $1$-category. Here the predicate ``unbiased" means, that $\kM_\sC : \GComp(\sC) \to \sC$ not only specifies binary composites of morphisms but more generally $n$-ary composites (cf. \cite{leinster-operads}).

\subsection{Dimension 2}

An associative $2$-category $\sC$ is a $2$-truncated $\GComp$-algebra satisfying the \complaw{} and \cohlaw{} laws. As an example consider the following

\begin{enumerate}
\item \textit{Morphism data}: Let $\sC$ be the globular set consisting of
\begin{align}
\sC_0 &= \Set{a,b} \\
\sC_1 &= \Set{g, \id_a,\id_b} \\
\sC_2 &= \Set{\alpha, \beta, \id_g, \id_{\id_a},\id_{\id_b}}
\end{align}
with 
\begin{align}
s(g) &= a \\
t(g) &= b \\
s(\alpha) = t(\alpha) &= \id_a \\
s(\beta) = t(\beta) &= g
\end{align}
(here we omitted specifying sources and targets of morphisms with the name $\id_x$ which are clearly supposed to both be $x$).

\item \textit{Composition data}: Next, we specify $\kM_\sC$ by identifying identity homotopies $\Id_x$ with $\id_x$ as before (there are no other interesting homotopies below dimension 3), and then setting
\begin{align}
\kM_\sC(\alpha \whisker 2 2 g) &= \beta \\
\kM_\sC(\id_a \whisker 2 2 \beta) &= \beta \\
\kM_\sC(\alpha \whisker 1 2 \alpha) &=  \alpha \\
\kM_\sC(\beta \whisker 2 2 \beta) &= \beta
\end{align}

\item \textit{Properties of composition}: We have already seen associativity and unitality playing out in dimension 1 and these arguments remain valid for 2-categories. 

A first novelty in dimension 2, is that we are now in a context in which we can take identities of (composites involving) identities. These behave as expected. For instance, by the \cohlaw{} law we find
\begin{align}
\kM_\sC(\Id_{\id_a \whisker 1 1 g}) &= \kM_\sC(\Id_{\kM_\sC(\id_a \whisker 1 1 g)}) \\
&= \kM_\sC(\Id_{g}) \\
&= \id_g
\end{align}
Here, the first equation follows from identifying source and target of $\Id_{\left(\infeq {\id_a \whisker 1 1 g} 2\right)}$, the second from unitality as discussed for 1-categories, and the third by the (implicit) definition of $\kM_\sC$. Similarly, one finds
\begin{equation}
\kM_\sC(\Id_{\Id_a}) = \kM_\sC(\Id_{\id_a}) = \id_{\id_a}
\end{equation}
as expected. As a final example, we compute
\begin{align}
g &= \kM_\sC(g)\\
 &= \kM_\sC (\Id_{\Id_a} \whisker 1 2 g) \\
 &= \kM_\sC (\id_{\id_a} \whisker 1 2 g)
\end{align}
and thus $\id_{\id_a}$ behaves as a strict unit for $\whisker 1 2$ as expected.

A second new feature in dimension 2, is that since the coherence law also forces equations derived from elementary 3-dimensional homotopies,  \textit{interchangers} start to play a role. For instance, the interchanger $\interchanger_{\alpha,\beta}$, with manifold diagram
\begin{restoretext}
\begingroup\sbox0{\includegraphics{ANCimg3/empty.png}}\includegraphics[clip,trim={.0\ht0} {.3\ht0} {.0\ht0} {.25\ht0} ,width=\textwidth]{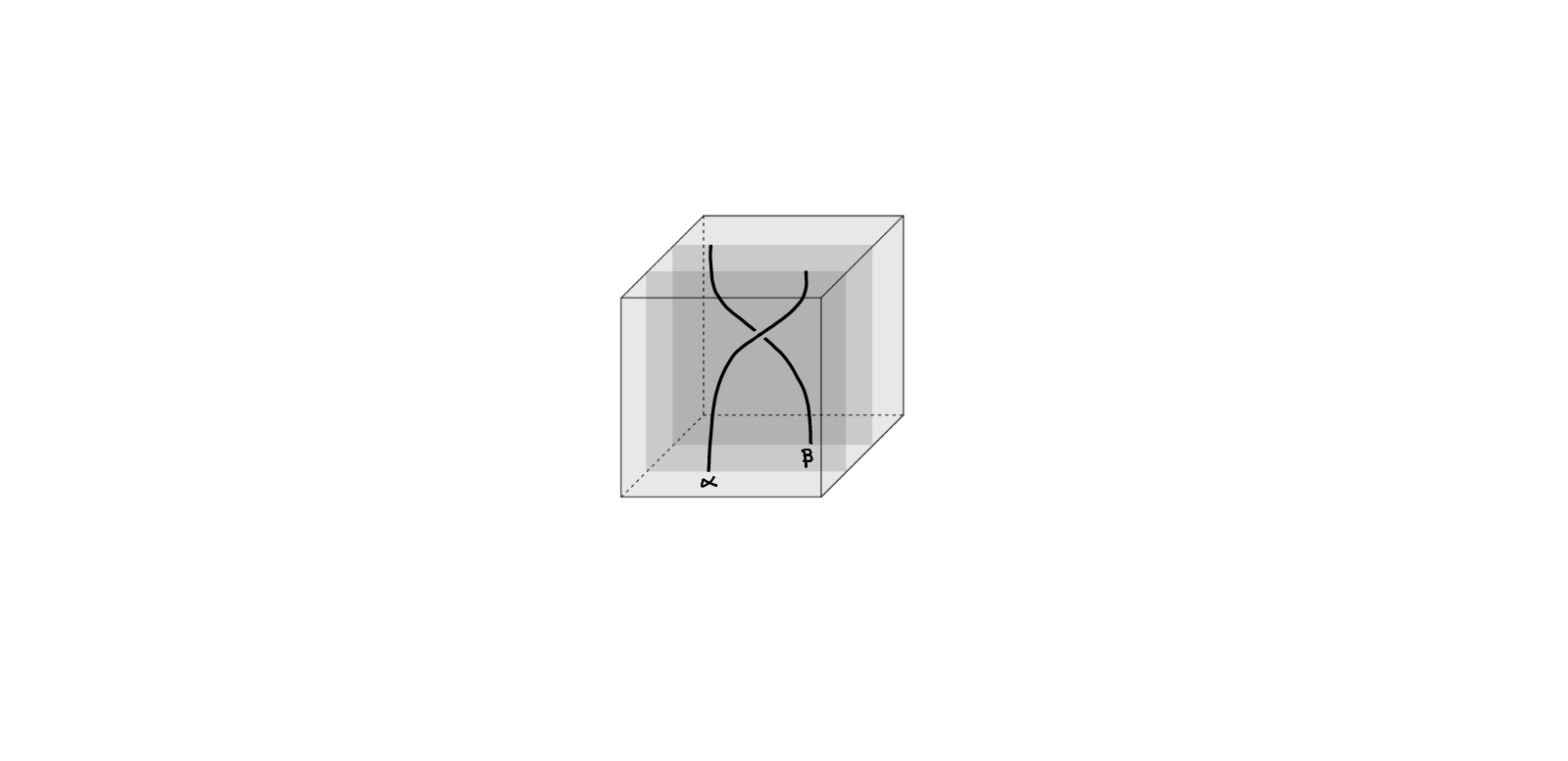}
\endgroup\end{restoretext}
forces the equation (by the \cohlaw{} law)
\begin{equation}
\kM_\sC((\alpha \whisker 2 2 g) \whisker 1 2 (\id_a \whisker 2 2 \beta)) = \kM_\sC((\id_a \whisker 2 2 \beta) \whisker 1 2 (\alpha \whisker 2 2 g))
\end{equation}
Both side can be evaluated using the \complaw{} law, and one finds that they agree as required. In general the interchange law thus holds strictly in associative 2-categories. However, it will not hold strictly in associative $n$-categories for $n > 2$.

We remark that the interchange law should not be confused with the exchange law which in fact cannot be formulated in full generality using only whiskering compositions, but to the degree that it can, it does hold strictly in all associative $n$-categories. A special instance of the exchange law for instance reads
\begin{equation}
\kM_sC( (\gamma \whisker 2 2 f) \whisker 1 2 (\gamma' \whisker 2 2 f)) = \kM_\sC((\gamma \whisker 1 2 \gamma') \whisker 2 2 f)
\end{equation}
where $\gamma, \gamma'$ are composable 2-morphisms and $f$ is a composable 1-morphism.
\end{enumerate}

In summary, associative $2$-categories are unbiased strict $2$-categories: that is, they are 2-categories with strictly associative and unital whiskering composition of arbitrary arity and strict interchange law.

\subsection{Dimension 3}

An associative $3$-category $\sC$ is a $3$-truncated $\GComp$-algebra satisfying the \complaw{} and \cohlaw{} laws. As an example consider the following

\begin{enumerate}
\item \textit{Morphism data}: Let $\sC$ be the globular set consisting of
\begin{align}
\sC_0 &= \Set{a,b,c} \\
\sC_1 &= \Set{f,g,h}\\
\sC_2 &= \Set{\alpha,\beta,\gamma,\delta}\\
\sC_3 &= \Set{\mathrm{ir}_{\alpha,\beta}, \mathrm{ir}_{\alpha,\beta}\inv}
\end{align}
together with according identity elements $\id_x \in \sC_k$ for each $x \in \sC_{k-1}$ (which have been omitted from the above sets for readability). Sources and targets of non-identity elements are given by
\begin{align}
s(f) = s(h) &= a\\
s(g) = t(f) &= b\\
t(g) = t(h) &= c\\
s(\alpha) = t(\alpha) &= f\\
s(\beta) = t(\beta) &= g\\
s(\gamma) = t(\gamma) = s(\delta) = t(\delta) &= h\\
s(\mathrm{ir}_{\alpha,\beta}) = t(\mathrm{ir}_{\alpha,\beta}\inv) &= \gamma \\
s(\mathrm{ir}_{\alpha,\beta}\inv) = t(\mathrm{ir}_{\alpha,\beta}) &= \delta
\end{align}

\item \textit{Composition data}: We next define $\kM_\sC$ by setting on binary whiskering composites
\begin{align}
\kM_\sC(f \whisker 1 1 g) &= h \\
\kM_\sC(\alpha \whisker 2 2 g) &= \gamma \\
\kM_\sC(f \whisker 2 2 \beta) &= \delta \\
\kM_\sC(\alpha \whisker 1 2 \alpha) &= \alpha \\
\kM_\sC(\beta \whisker 1 2 \beta) &= \beta \\
\kM_\sC(\delta \whisker 1 2 \gamma) &=  \delta\\
\kM_\sC(\gamma \whisker 1 2 \delta) &=  \gamma\\
\kM_\sC(\gamma \whisker 1 2 \gamma) &= \gamma \\
\kM_\sC(\delta \whisker 1 2 \delta) &=  \delta\\
\kM_\sC(\mathrm{ir}_{\alpha,\beta} \whisker 2 3 \gamma) = \kM_\sC(\mathrm{ir}_{\alpha,\beta} \whisker 2 3 \delta) &= \mathrm{ir}_{\alpha,\beta}\\
\kM_\sC(\mathrm{ir}_{\alpha,\beta}\inv \whisker 2 3 \gamma) = \kM_\sC(\mathrm{ir}_{\alpha,\beta}\inv \whisker 2 3 \delta) &=  \mathrm{ir}_{\alpha,\beta}\inv\\
\kM_\sC(\gamma \whisker 2 3 \mathrm{ir}_{\alpha,\beta}) = \kM_\sC(\gamma \whisker 2 3 \mathrm{ir}_{\alpha,\beta}\inv) &=  \id_{\gamma}\\
\kM_\sC(\delta \whisker 2 3 \mathrm{ir}_{\alpha,\beta}) = \kM_\sC(\delta \whisker 2 3 \mathrm{ir}_{\alpha,\beta}\inv) &=  \id_{\delta}\\
\kM_\sC(\mathrm{ir}_{\alpha,\beta} \whisker 1 3 \mathrm{ir}_{\alpha,\beta}\inv) &=  \id_{\gamma}\\
\kM_\sC(\mathrm{ir}_{\alpha,\beta}\inv \whisker 1 3 \mathrm{ir}_{\alpha,\beta}) &= \id_{\delta} \\
\end{align}
We still need to define $\kM_\sC$ on elementary homotopies. For identities, this is done as before (that is, we set $\kM_\sC(\Id_x) = \id_x$). However, we now have elementary homotopies which are not identity homotopies as well, namely, in the form of interchangers. For those, we make the following definitions:
\begin{align} 
\kM_\sC(\interchanger_{\alpha,\beta}) = \mathrm{ir}_{\alpha,\beta} \\
\kM_\sC(\interchanger_{\alpha,\beta}\inv) = \mathrm{ir}_{\alpha,\beta}\inv
\end{align}
Using the \complaw{} law, one can easily verify that this assignment is globular, that is, for instance
\begin{align}
\kM_\sC(\gsrc(\interchanger_{\alpha,\beta})) = s(\mathrm{ir}_{\alpha,\beta}) \\
\kM_\sC(\gtgt(\interchanger_{\alpha,\beta})) = t(\mathrm{ir}_{\alpha,\beta}) \\
\end{align}
The reader might wonder whether we also need to define $\kM_\sC$ on interchangers such as
\begin{equation}
\interchanger_{\alpha,\id_g}
\end{equation}
 As we will now see, unitality extends strictly to homotopies as well and we are forced to set
\begin{equation}
\kM_\sC(\interchanger_{\alpha,\id_g}) = \id_{\kM_\sC(\alpha \whisker 2 2 g)}
\end{equation}
in such a situation.

\item \textit{Properties of composition}: Associativity and unitality hold as they did in previous dimensions. The genuinely new phenomenon is now the laws that are being forced by higher homotopies. Consider for instance the following element of $\GComp(\kiC{}(\kM_\sC))_4$ (represented by its manifold diagram):
\begin{restoretext}
\begingroup\sbox0{\includegraphics{ANCimg3/empty.png}}\includegraphics[clip,trim={.0\ht0} {.25\ht0} {.0\ht0} {.25\ht0} ,width=\textwidth]{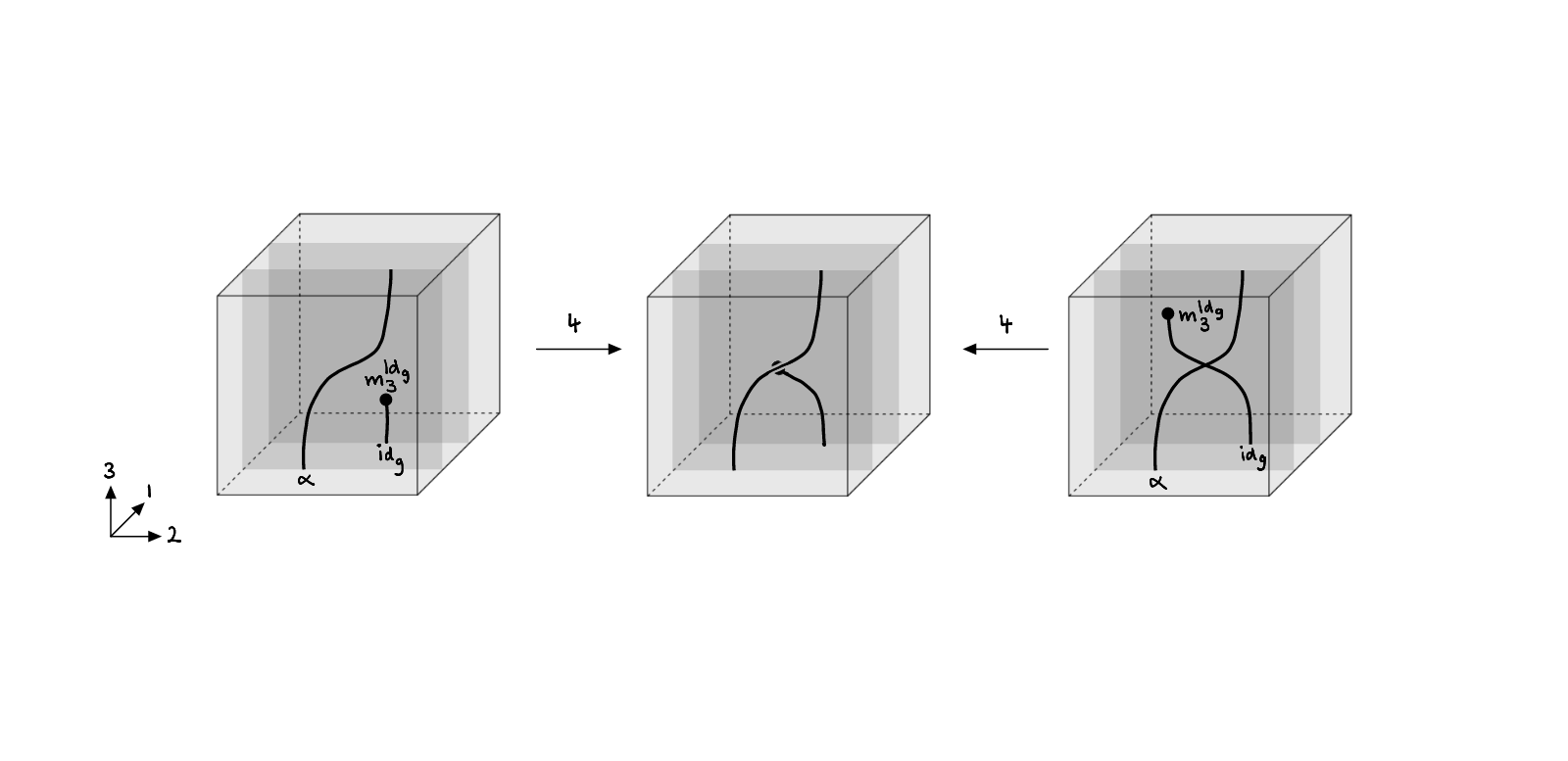}
\endgroup\end{restoretext}
Here, recall $\infeq {\Id_g} 3$ is the witness of $\kM_\sC(\Id_g) = \id_g$ and has type
\begin{restoretext}
\begingroup\sbox0{\includegraphics{ANCimg3/empty.png}}\includegraphics[clip,trim={.0\ht0} {.3\ht0} {.0\ht0} {.25\ht0} ,width=\textwidth]{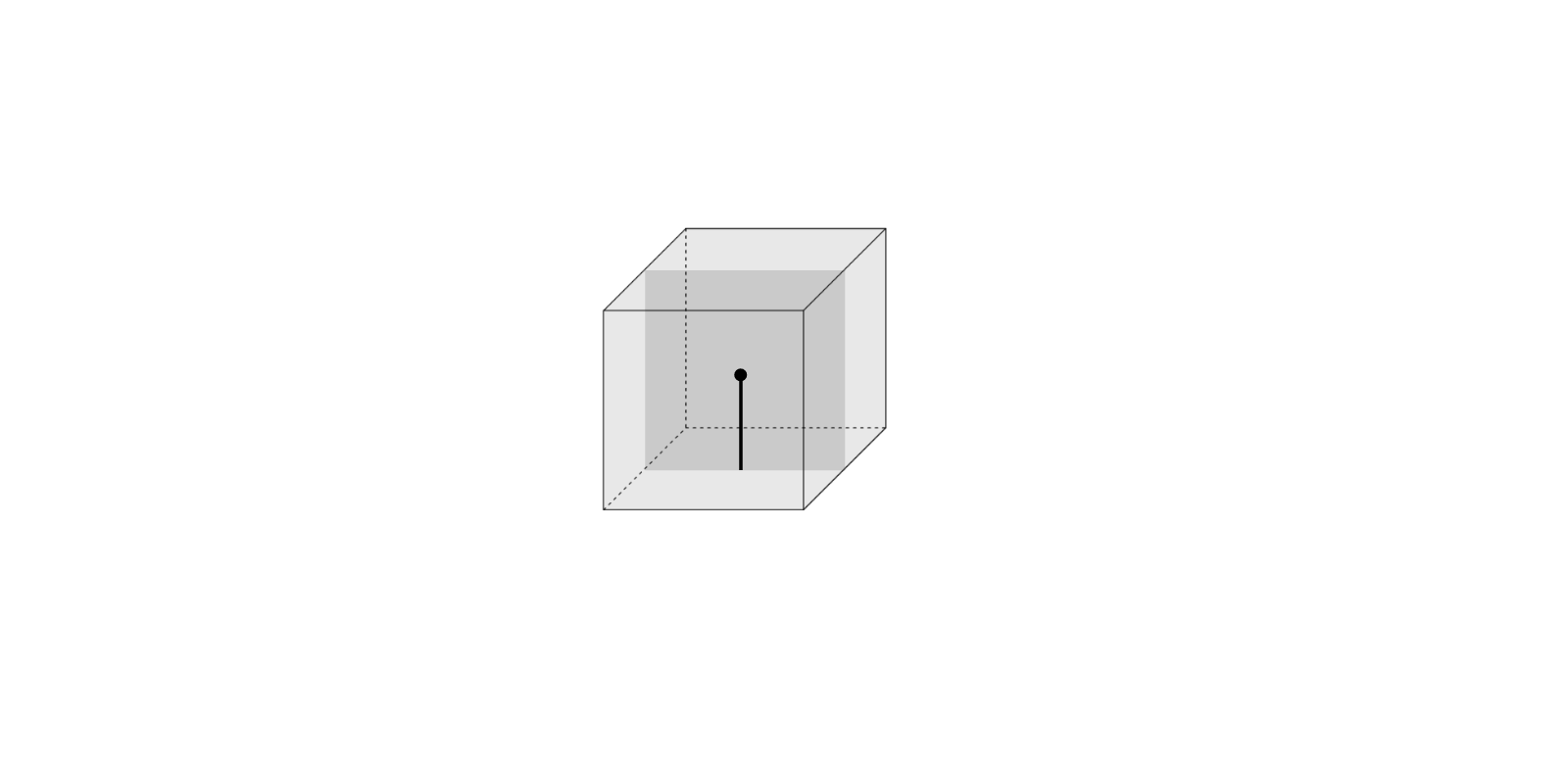}
\endgroup\end{restoretext}
As a consequence of the above $4$-dimensional diagram, the \cohlaw{} law forces an equation between its source and target, namely
\begin{restoretext}
\begingroup\sbox0{\includegraphics{ANCimg3/empty.png}}\includegraphics[clip,trim={.0\ht0} {.3\ht0} {.0\ht0} {.25\ht0} ,width=\textwidth]{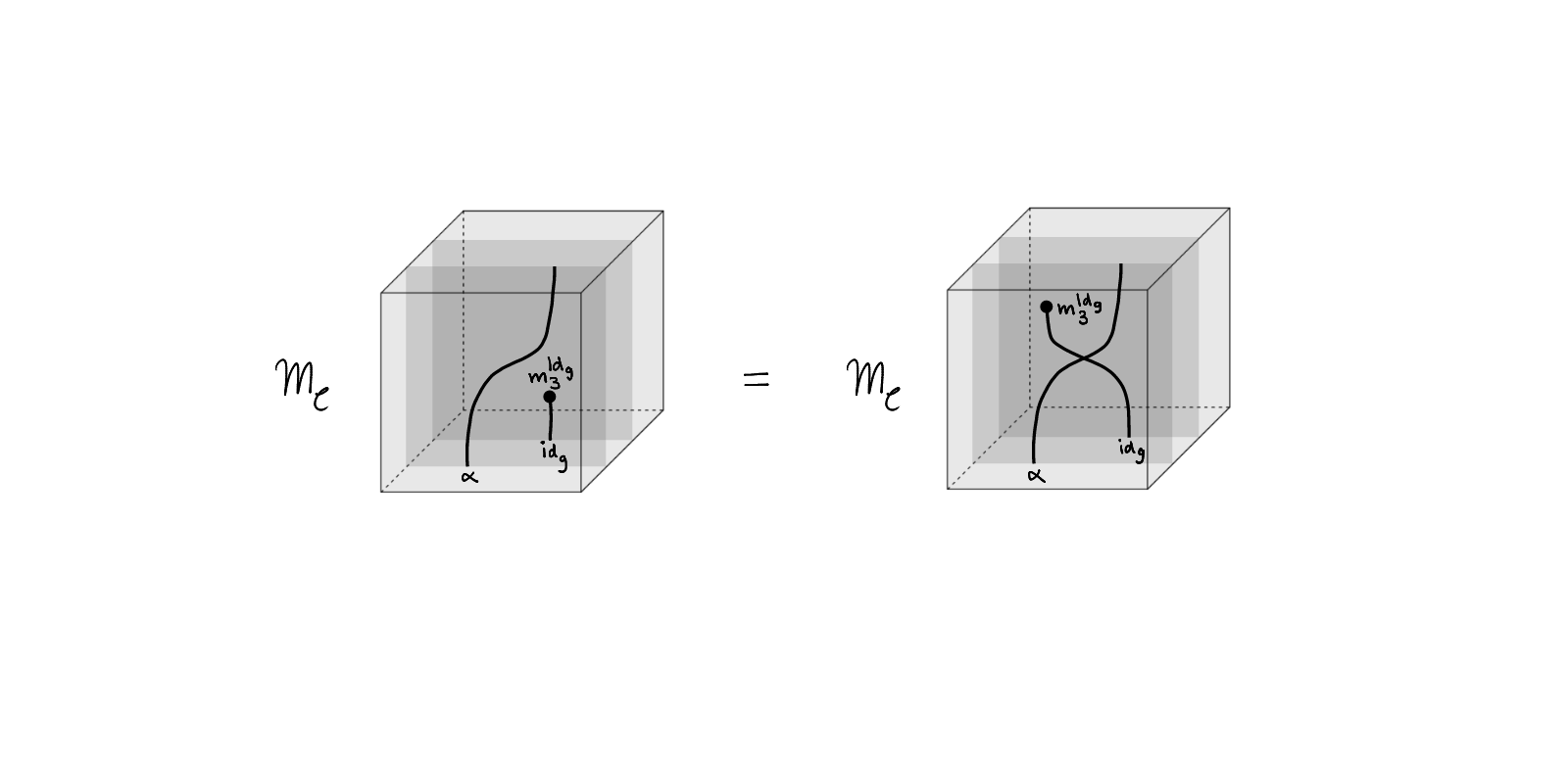}
\endgroup\end{restoretext}
or, written as a symbolic equation,
\begin{equation}
\extkM_\sC ((\alpha \whisker 2 2 g) \whisker 2 3 (f \whisker 3 3 \infeq {\Id_g} 3)) 
= \extkM_\sC(\interchanger_{\alpha,\id_g} \whisker 3 3 ((f \whisker 3 3 \infeq {\Id_g} 3) \whisker 2 3 (\alpha \whisker 2 2 g)))
\end{equation}
Now, both sides of the equation can be further simplified using the \complaw{} and \cohlaw{} laws, yielding
\begin{equation}
\kM(\Id_{\alpha \whisker 1 2 g}) = \kM(\interchanger_{\alpha,\id_g})
\end{equation}
and the left-hand side evaluates to $\id_{\gamma}$ as previously claimed.

There are other elementary homotopies in $\GComp(\kiC{}(\kM_\sC)$ which similarly encode equation that need to be satisfied by $\kM_sC$, for instance
\begin{restoretext}
\begingroup\sbox0{\includegraphics{ANCimg3/empty.png}}\includegraphics[clip,trim={.0\ht0} {.3\ht0} {.0\ht0} {.25\ht0} ,width=\textwidth]{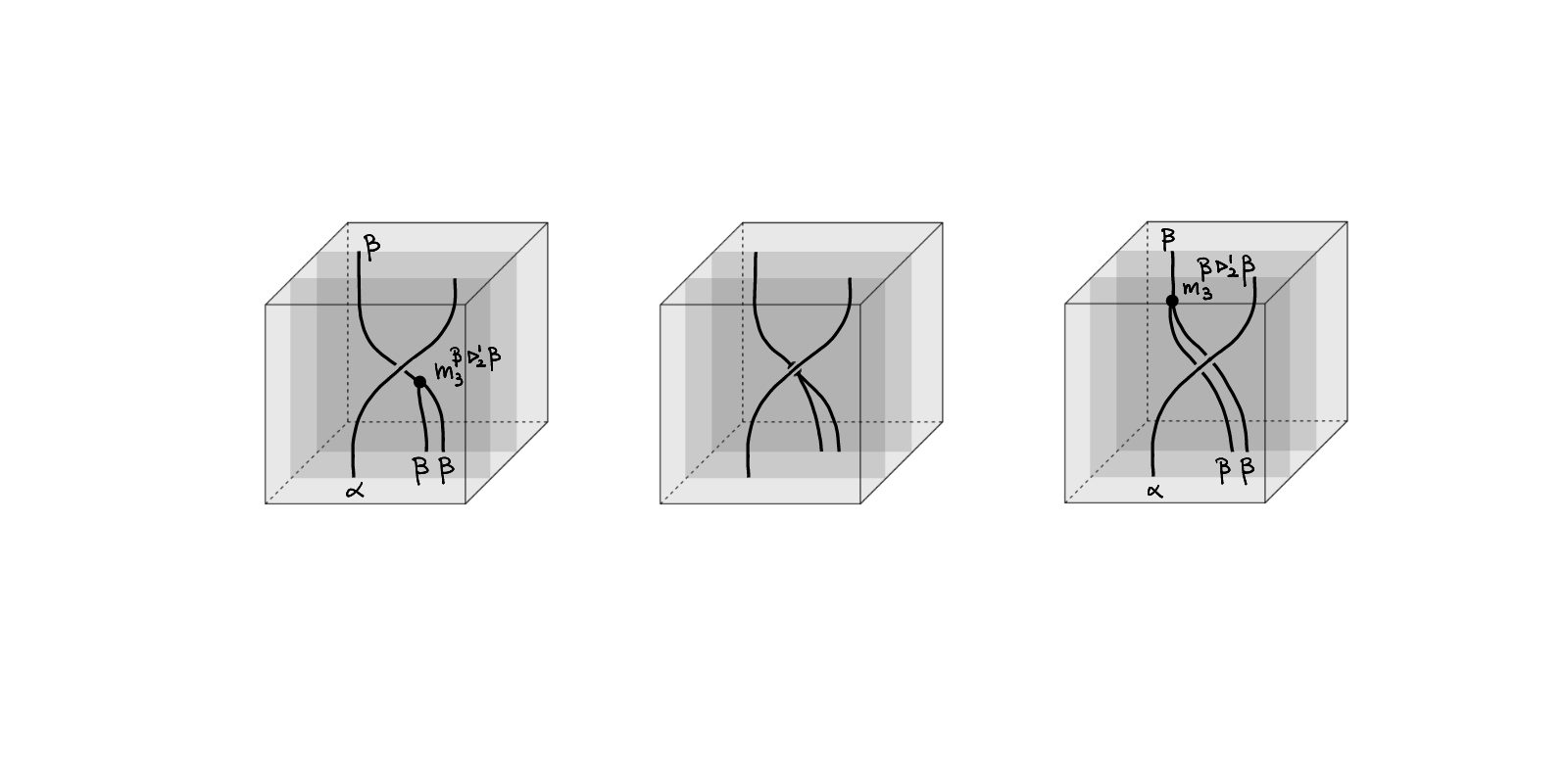}
\endgroup\end{restoretext}
in $\GComp(\kiC{}(\kM_\sC)$ means that we must have
\begin{restoretext}
\begingroup\sbox0{\includegraphics{ANCimg3/empty.png}}\includegraphics[clip,trim={.0\ht0} {.3\ht0} {.0\ht0} {.25\ht0} ,width=\textwidth]{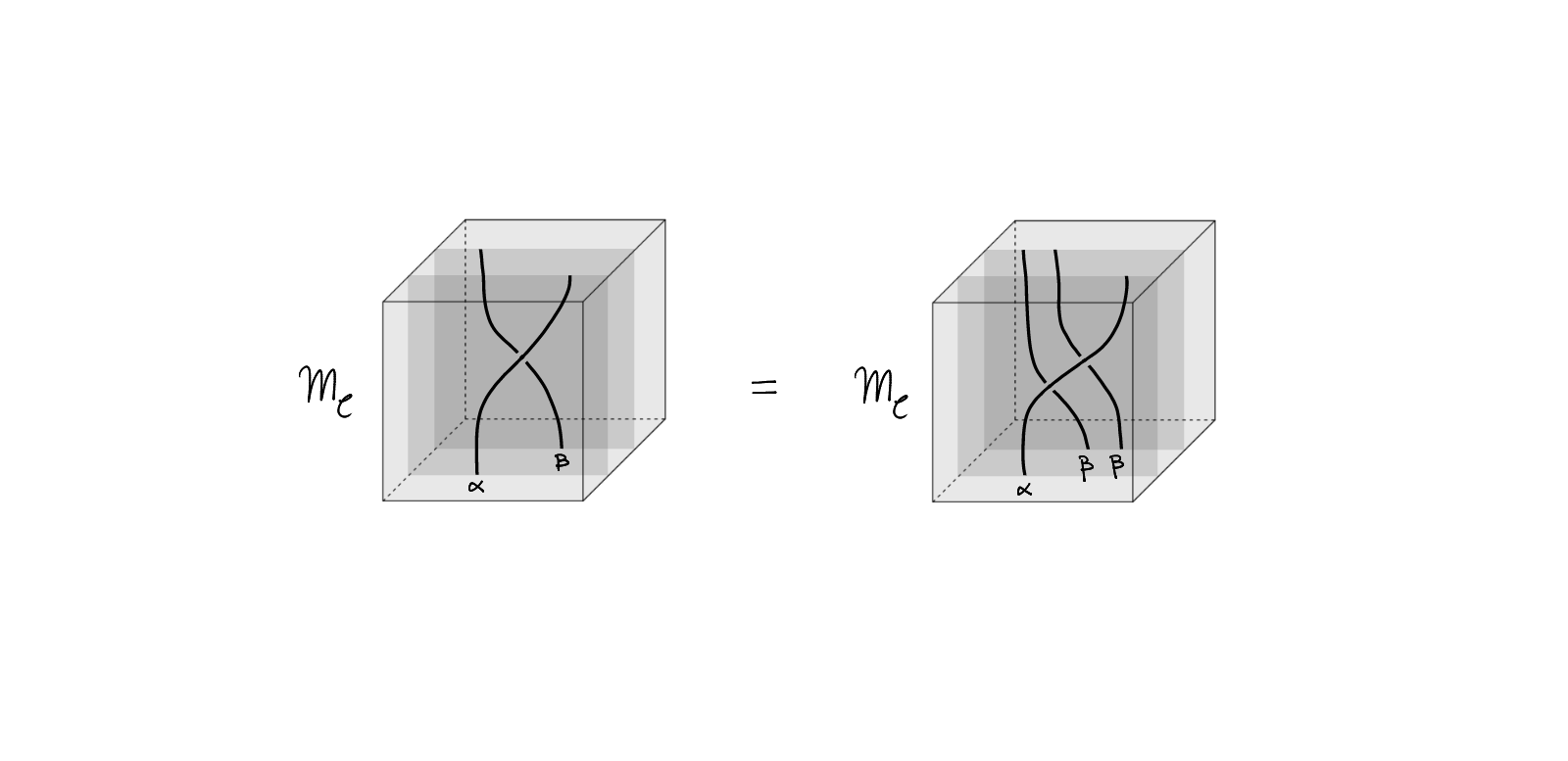}
\endgroup\end{restoretext}
which (together with the other equations) is left to be checked by the reader. As a consequence, we see that composition is not only strictly associative, but also strictly interacts with weak coherences---And in dimension 3 as for the present examples, there is only one weak coherence: the interchanger.
\end{enumerate}

The last observation also holds for another prominent model of (non-strict) 3-categories, called \textit{Gray categories}. We thus conclude that associative $3$-categories can be thought of as unbiased Gray categories. 

\begin{rmk}[Comparison to weak definitions] In contrast, we remark that for ``weak up to depth $m$" associative $3$-categories  (for $m > 0$) the above two coherences would have had to be chosen which would lead to additional weak coherence data.
\end{rmk}

\subsection{Dimension n}

More generally, associative $n$-categories are an attempt to capture the notion of ``Gray $n$-categories".

A conjectural example of an associative $n$-category is given by choosing objects to be strict $(n-1)$-categories, morphisms to be functors, and $k$-morphisms to be pseudo-natural $k$-transformations. A similar example should be obtainable when starting with associative $(n-1)$-categories as objects. Note that this example hinges upon a good definition of pseudo-natural $k$-transformation which we hope to give in future work. 

\chapter*{Conclusion and future work}

We began this thesis with our definition of flag-foliation-compatible stratifications of Euclidean space, resulting in the category $\Cubeo n \cC(\bnum 1)$ of $\cC$-labelled $n$-cubes and open maps between them. These geometric gadgets are powerful for expressing many interesting mathematical ideas, in particular the idea of manifold diagrams---which, as discussed, has been long sought after as a local model for higher categories. By using the central observation that $\cC$-labelled $n$-cubes have recursive structure, we came up with a combinatorial structure for them, phrased in a category $\Buno n \cC$. We provided a (functional, not yet functorial) geometric realisation procedure $\norm{-} : \Buno n \cC \to \Cubeo n \cC(\bnum 1)$ showing that the objects of $\Buno n \cC$ describe a large (or at the very least, interesting) subclass of objects in $\Cubeo n \cC(\bnum 1)$. We gave an extensive and elementary analysis of the objects and morphisms in $\Buno n \cC$, unravelling their exciting combinatorial properties, in particular focusing on epimorphisms (referred to as collapses) and monomorphism (referred to as embeddings). One result that certainly stands out is the normal form theorem, which gives a way to decide ``geometric equality" of combinatorial cubes by comparing their normal forms. The combinatorial theory of $\cC$-labelled $n$-cubes was then used to define several (fully) algebraic models of higher categories. 

Many allusions towards the future potential of the framework have been made: it seems possible that labelled $n$-cubes are an elegant (and minimal!) combinatorial setting for capturing notions such as tangles and extended cobordism, which in turn are both connected to many other parts of (unstable and stable) Algebraic Topology, including, at the most basic level, CW-complexes (via a generalised Thom-Pontryagin construction, sketched in \autoref{ch:geom}). There are numerous other potential directions which so far haven't been mentioned at all: for instance, the framework presented here has a natural extension to ``tensor $n$-categories" which in turn allow us to study objects that should live in the $n$-category $n$\textbf{Vect} of so-called $n$-vector spaces over some ring $k$. 

Future research will need to show how the ``combinatorial elegance" of labelled $n$-cube can be put to use for the needs of Algebraic Topology and Higher Category Theory.

\appendix
\setcounter{section}{0}
\renewcommand{\thesection}{A.{\arabic{section}}}
\renewcommand{\theHsection}{A.{\arabic{section}}}

\chapter{Connection of algebraic and geometric models} \label{ch:geom}

This appendix sketches the connections of presented associative $\infty$-groupoids, coherent invertibility and CW-complexes, using a generalised version of the Thom-Pontryagin construction (which will only be sketched, but with a reasonable level of detail).

\section{A theory of invertibility yielding manifolds} \label{ssec:com_cob}

Our starting point is \autoref{sec:group_TI}. From now on, we will take $\TI^n$ (and $\TI$) to mean $\TI^{n,\infty +1}$ (and $\TI^{\infty,\infty +1}$).

We leisurely discuss morphisms of $\TI$ (the section can be skipped by readers who feel they know what these should be). The elements of $\Comp(\TI)_k$ are called  \textit{combinatorial framed $k$-cobordism}. Note that at a first glance, a better name would be ``combinatorial $(k+1)$-framed $k$-tangles", however ultimately the morphisms of $\TI$ (via \autoref{constr:adjoin_inv_gen}) will be used in \textit{arbitrary co-dimensions} which makes the name ``cobordisms" appropriate as well. We will comment on the notion of framing only later on in this chapter.

For now, we give examples of types and morphisms of $\TI$ in low dimensions
\begin{enumerate}
\item $\TI_1$, as we saw, has exactly two elements, a $1$-morphism from $-$ to $+$ and its inverse. General morphisms are alternating sequences of these morphisms.

\item In dimension $2$, we saw $\TI_2$ has exactly four elements: two cups and two caps. General $2$-morphisms in $\Comp(\TI)_2$ are obtained as globular, normalised well-typed $2$-cubes labelled by the above generators.

Note that we cannot have generating $2$-morphisms of any different form. For instance the following are impossible types of generators
\begin{restoretext}
\begingroup\sbox0{\includegraphics{ANCimg/page1.png}}\includegraphics[clip,trim=0 {.2\ht0} 0 {.35\ht0} ,width=\textwidth]{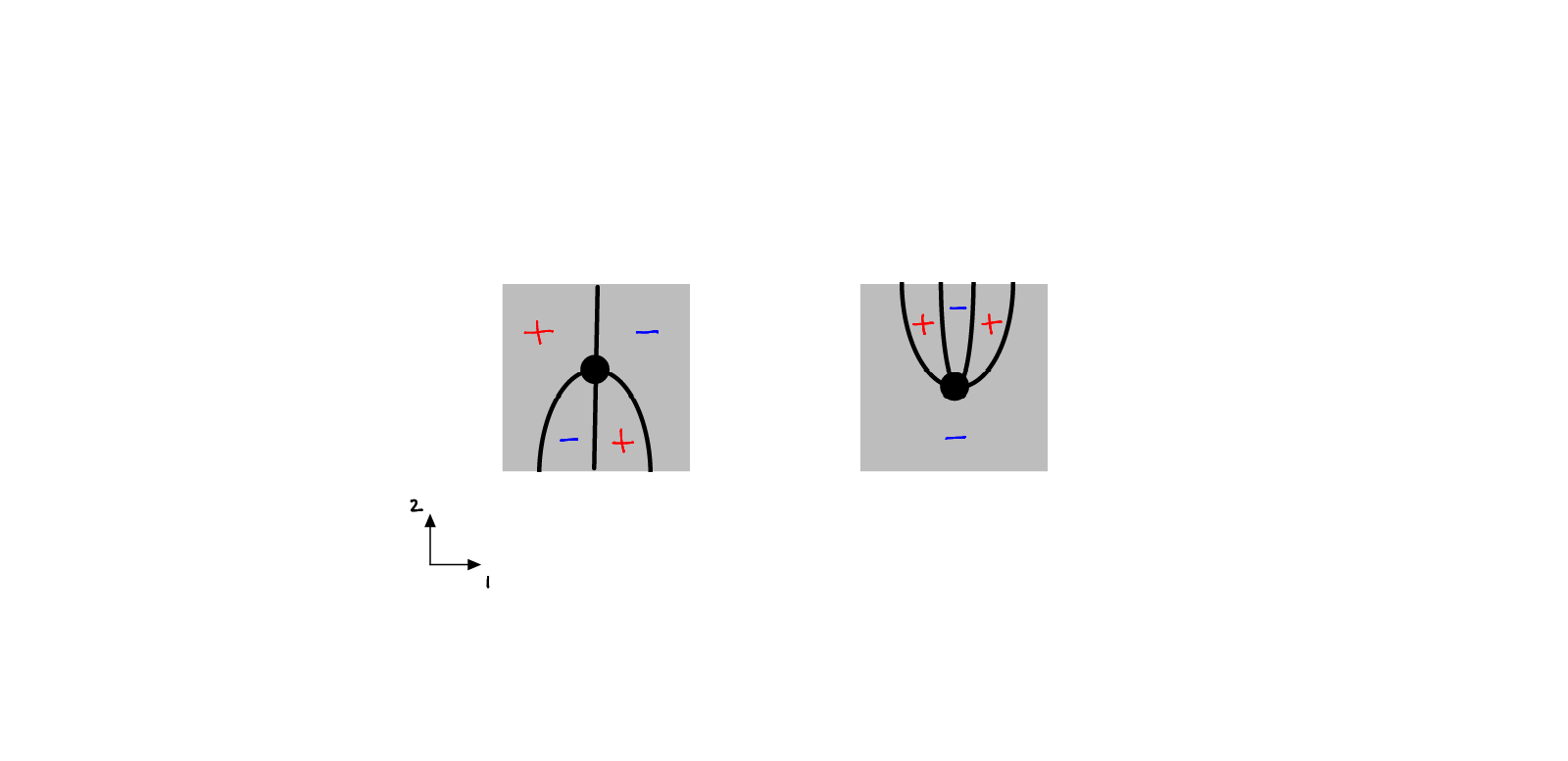}
\endgroup\end{restoretext}
In each case the $P^1_{\infty +1}$-condition from \autoref{constr:TI} is not satisfies. In general the $P^1_{\infty +1}$-condition will guarantee, that our pictures (more precisely, the complements of regions labelled by $-$ and $+$) will look like manifolds.

\item In dimension $n\geq 2$, $\TI_n$ has infinitely many elements for $n \geq 3$. For instance, for $n = 3$, the easiest elements of $\TI_3$ have types of the forms such as
\begin{restoretext}
\begingroup\sbox0{\includegraphics{ANCimg/page1.png}}\includegraphics[clip,trim=0 {.05\ht0} 0 {.1\ht0} ,width=\textwidth]{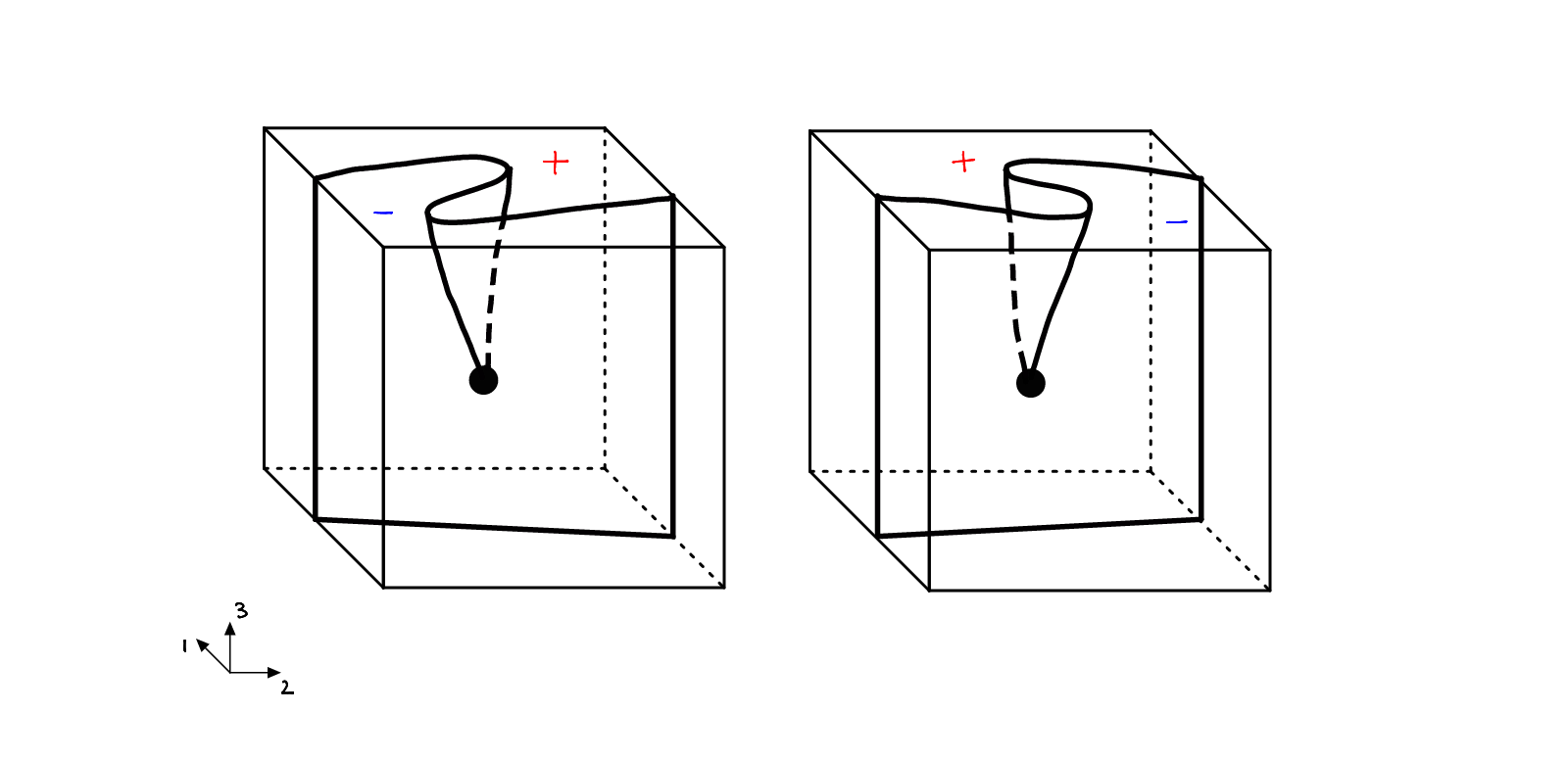}
\endgroup\end{restoretext}
(called snake singularities) and
\begin{restoretext}
\begingroup\sbox0{\includegraphics{ANCimg/page1.png}}\includegraphics[clip,trim=0 {.05\ht0} 0 {.1\ht0} ,width=\textwidth]{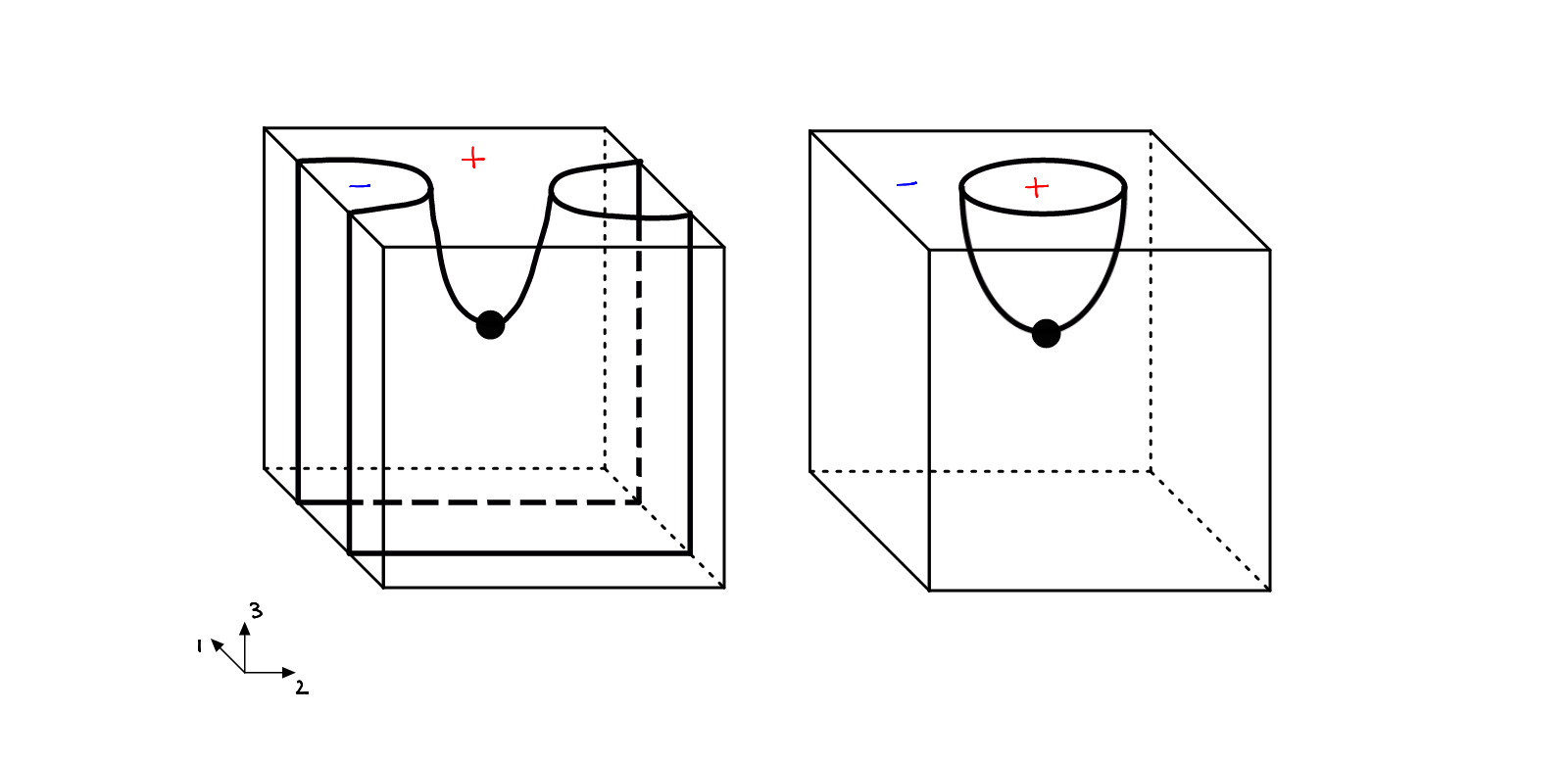}
\endgroup\end{restoretext}
(called saddle/crotch and circle death/birth singularities). Note that in each case the $3$-cube is split into two ball-homeomorphic connected components labelled by $-$ and $+$ respectively, as required in the definition of $\TI$.

The shared feature across the above singularities is that they are \textit{binary}. That is, each of them relates two generators of co-dimension 1 in their input (or their output). However, $\TI$ contains more complicated singularities than these, and based on this it would be be reasonable to call $\TI$ the \textit{unbiased} theory of invertibility. For instance, it also includes (non-binary) generating $3$-morphisms of the following types
\begin{restoretext}
\begingroup\sbox0{\includegraphics{ANCimg/page1.png}}\includegraphics[clip,trim=0 {.05\ht0} 0 {.1\ht0} ,width=\textwidth]{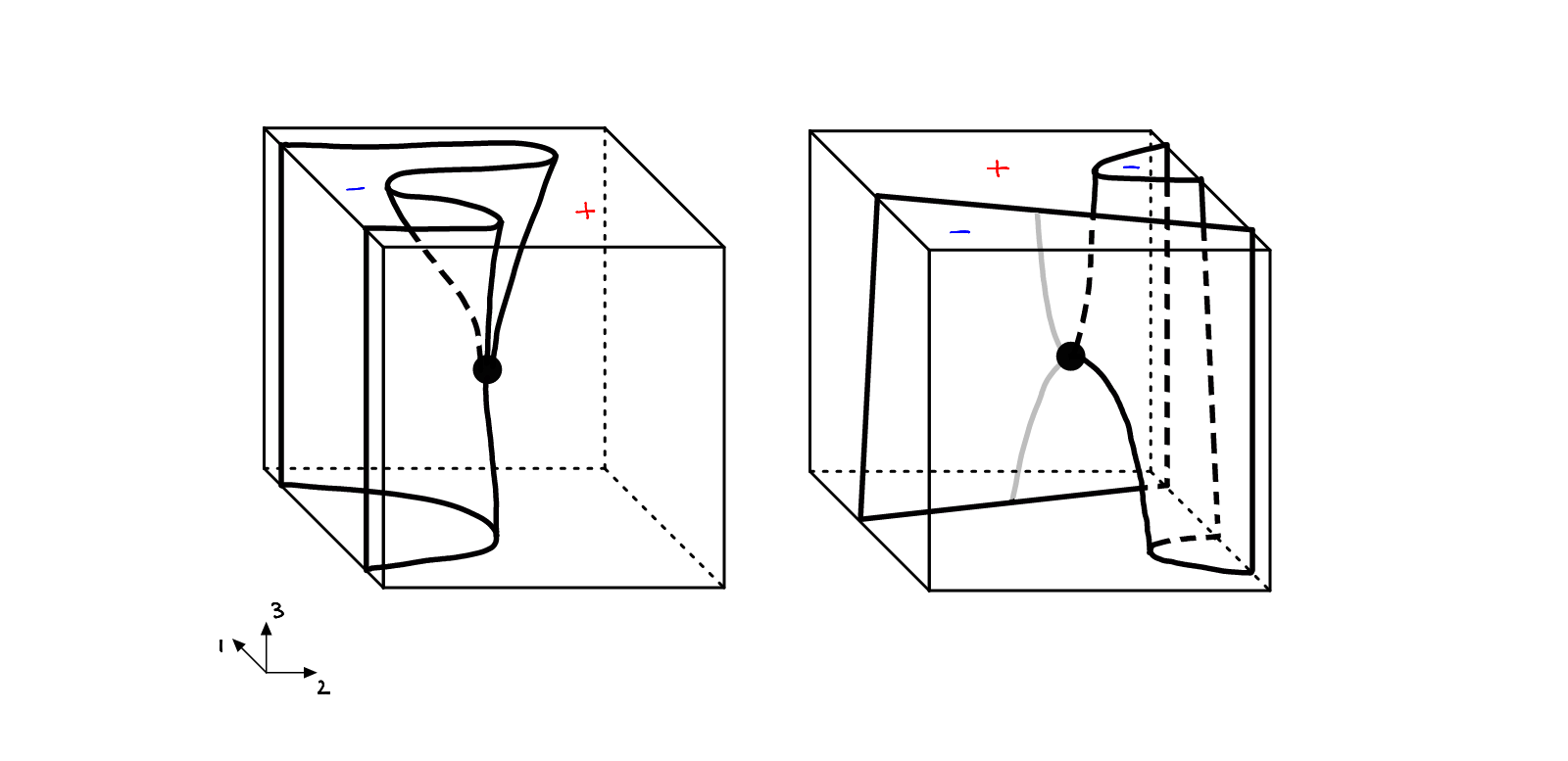}
\endgroup\end{restoretext}
The type on the left has a total of three cups and caps in its output, while the type on the right has a single cup in both its input and its output (also note that lines in \cgray{} are only meant as visual aids).
\end{enumerate}

\begin{rmk}[Finite presentations of the theory of invertibility] \label{rmk:finite_cob_rep} Importantly, the above two types can be represented by (that is, they are homotopic to) morphisms composed of the snake, saddle and crotch singularities. In fact there is finite collection of $3$-singularities that can represent any other singularity in dimension $3$. A (related) discussion of this can be found in \cite{schommer2011classification}. It is expected to hold true for general dimension $n$, and work towards a general finite presentation $\TI^{n,\mathrm{fin}}$ (in fact using only \textit{binary} generators) is in progress. The claim would then be, that the chain of inclusion
\begin{equation}
\TI^{n,\mathrm{fin}} \into \TI^{n,\infty+1} \into \TI^{n,\infty} \into .... \into \TI^{n,0} \into \TI^{n,-1}
\end{equation}
induces ``equivalences of categories" in each degree, and thus, we can represent each morphisms in $\TI^{n,k}$ equivalently as a morphism in $\TI^{n,\mathrm{fin}}$. For instance, the morphism 
\begin{restoretext}
\begingroup\sbox0{\includegraphics{ANCimg3/empty.png}}\includegraphics[clip,trim={.0\ht0} {.32\ht0} {.0\ht0} {.3\ht0} ,width=\textwidth]{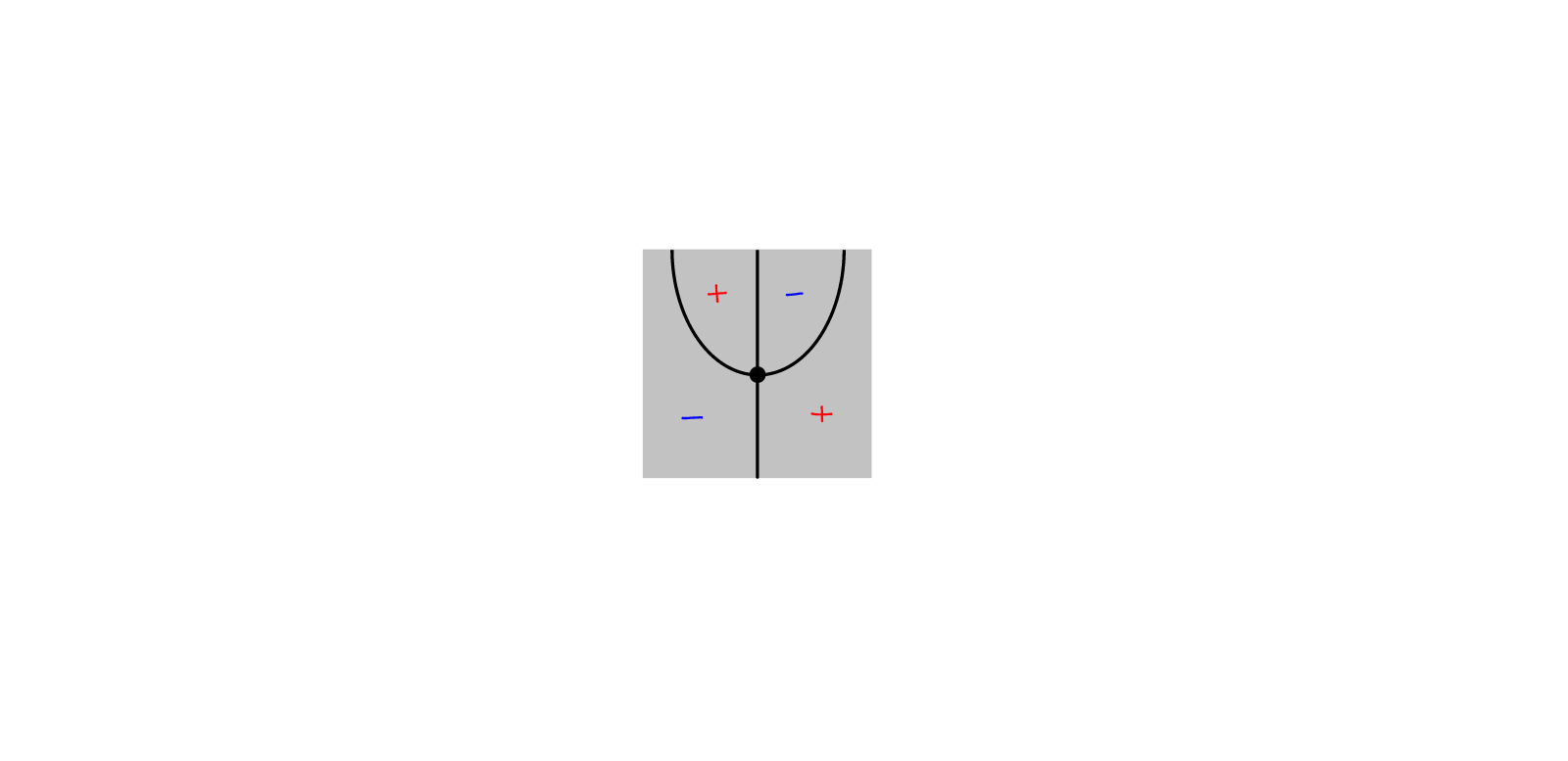}
\endgroup\end{restoretext}
in $\TI^{n,-1}$ can be represented by a morphism
\begin{restoretext}
\begingroup\sbox0{\includegraphics{ANCimg3/empty.png}}\includegraphics[clip,trim={.0\ht0} {.32\ht0} {.0\ht0} {.3\ht0} ,width=\textwidth]{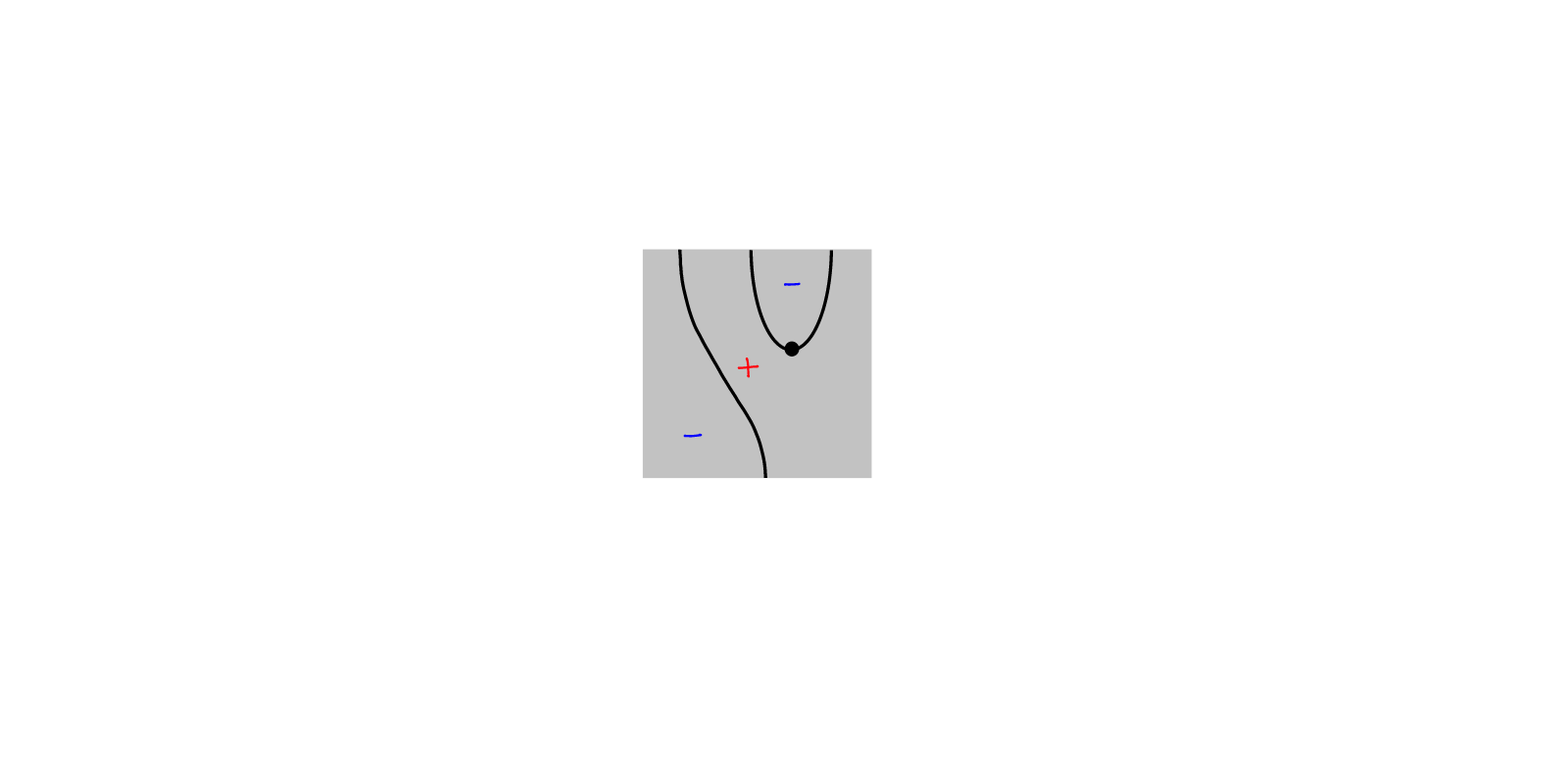}
\endgroup\end{restoretext}
in $\TI^{n,\mathrm{fin}}$ if the cap singularity were to be part of $\TI^{n,\mathrm{fin}}$ (which it is of course!).
\end{rmk}

\begin{rmk}[Composites] Adding to the pictures of generators above, we draw pictures of composites. As remarked above, in general only ``manifold-like" generators will be allowed in $\TI$. For instance, the following types for generating $3$-morphisms are not allowed
\begin{restoretext}
\begingroup\sbox0{\includegraphics{ANCimg/page1.png}}\includegraphics[clip,trim=0 {.2\ht0} 0 {.15\ht0} ,width=\textwidth]{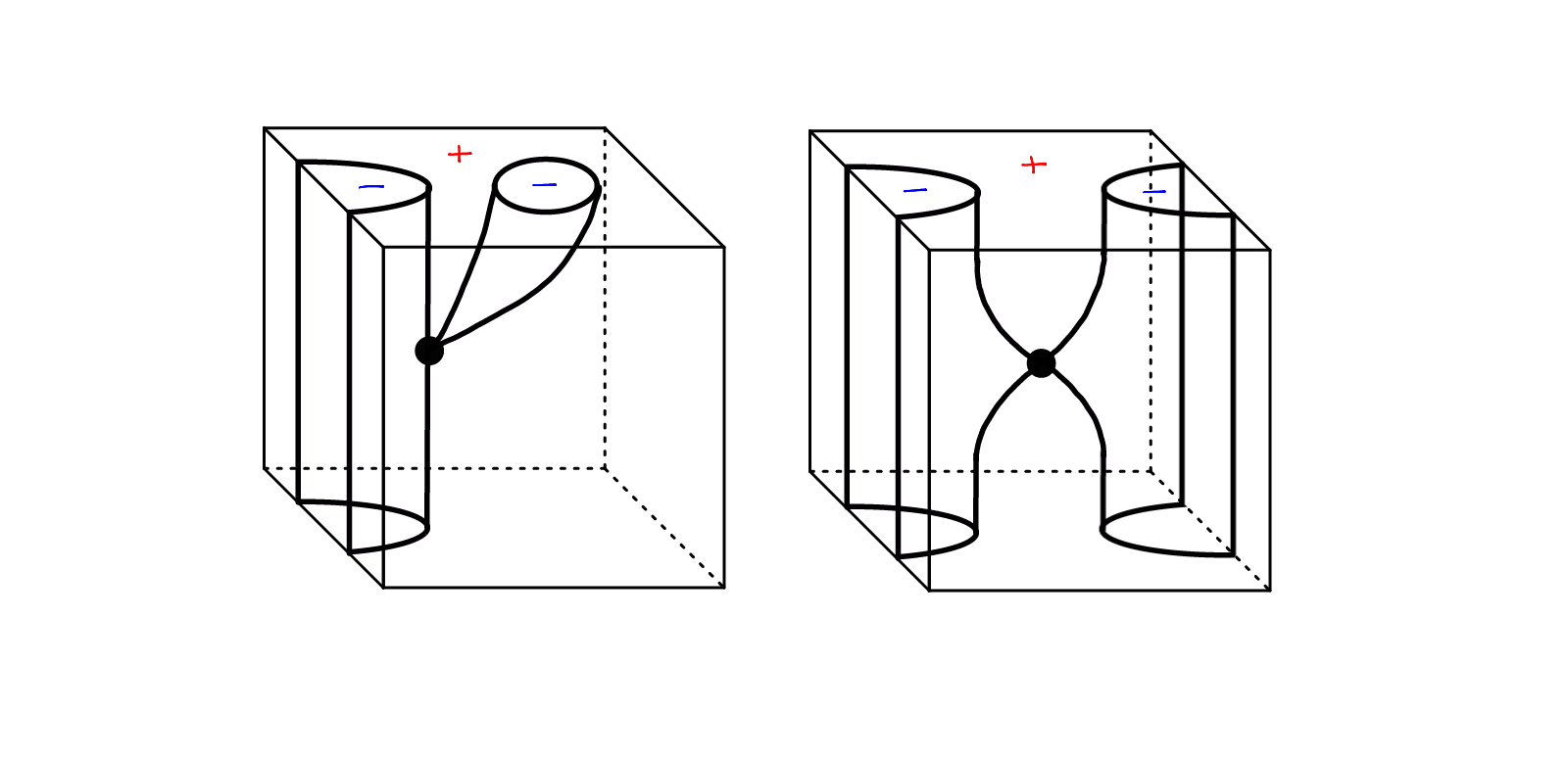}
\endgroup\end{restoretext}
In both case the condition from \autoref{constr:TI} is not satisfied. 

Morphisms $\Comp(\TI)_n$ are globular normalised well-typed $n$-cubes labelled by the generators of $\TI$. For instance, $3$-morphisms can be of the form
\begin{restoretext}
\begingroup\sbox0{\includegraphics{ANCimg/page1.png}}\includegraphics[clip,trim=0 {.0\ht0} 0 {.15\ht0} ,width=\textwidth]{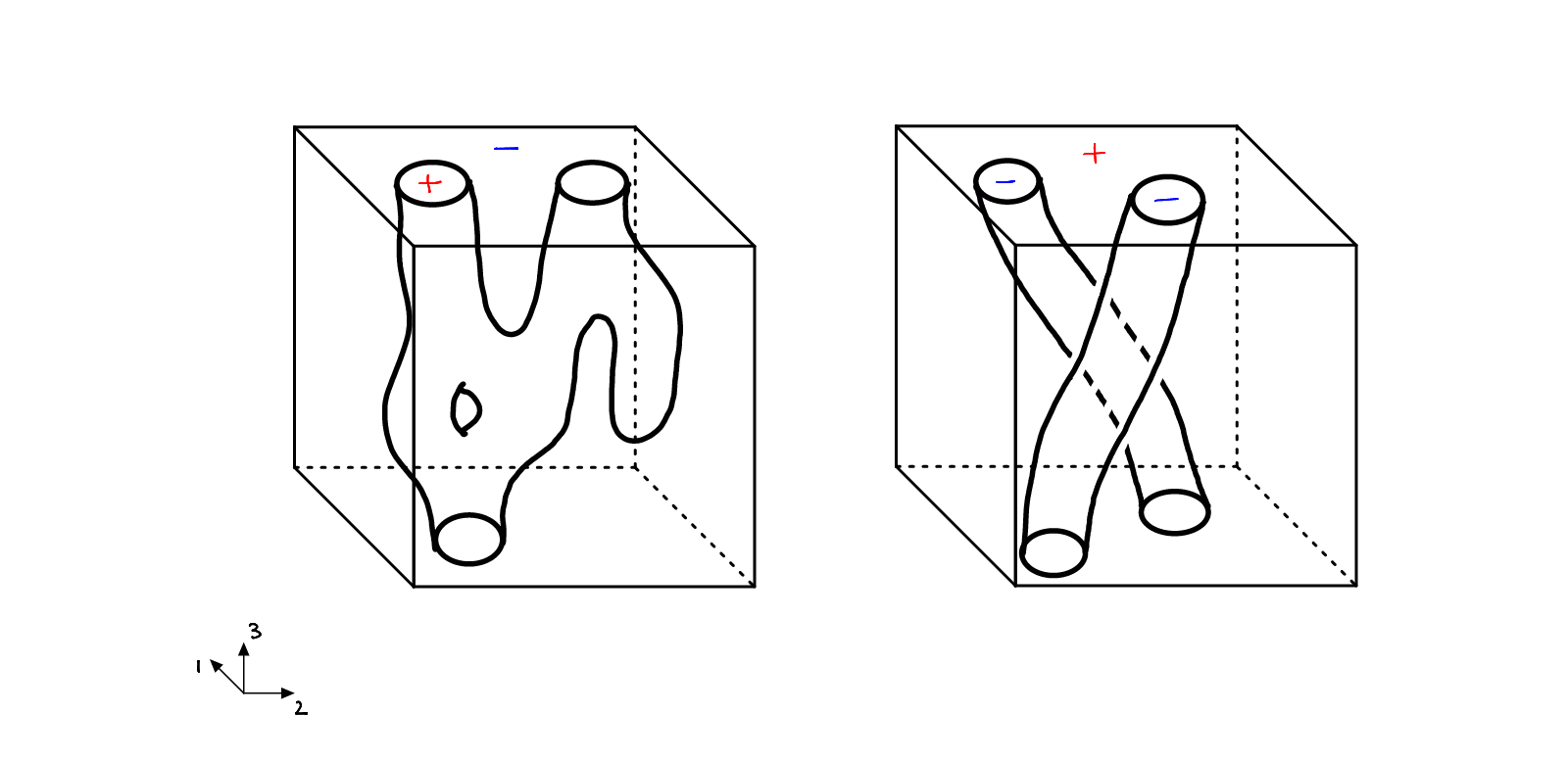}
\endgroup\end{restoretext}
As before we depicted the manifolds obtained as the \textit{union} of all submanifolds labelled by some $\ic_{S \equiv T}$ (thus the complement of the union of the regions labelled by $+$ and $-$). 
\end{rmk}

\begin{rmk}[Connections to Morse Theory] We also remark, that the combinatorial structure of manifold diagrams gives a fine-grained classification of these manifolds. This classification has the flavour of Morse theory, or Cerf Theory (cf. \cite{gay2012reconstructing},  \cite{cerf1970stratification}), which studies singularities of (Morse) functions on manifolds. To illustrate this, as a final example consider the following generating  $4$-morphism in $\TI_4$, whose source and target are (combinatorially) distinct surfaces
\begin{restoretext}
\begingroup\sbox0{\includegraphics{ANCimg/page1.png}}\includegraphics[clip,trim=0 {.0\ht0} 0 {.15\ht0} ,width=\textwidth]{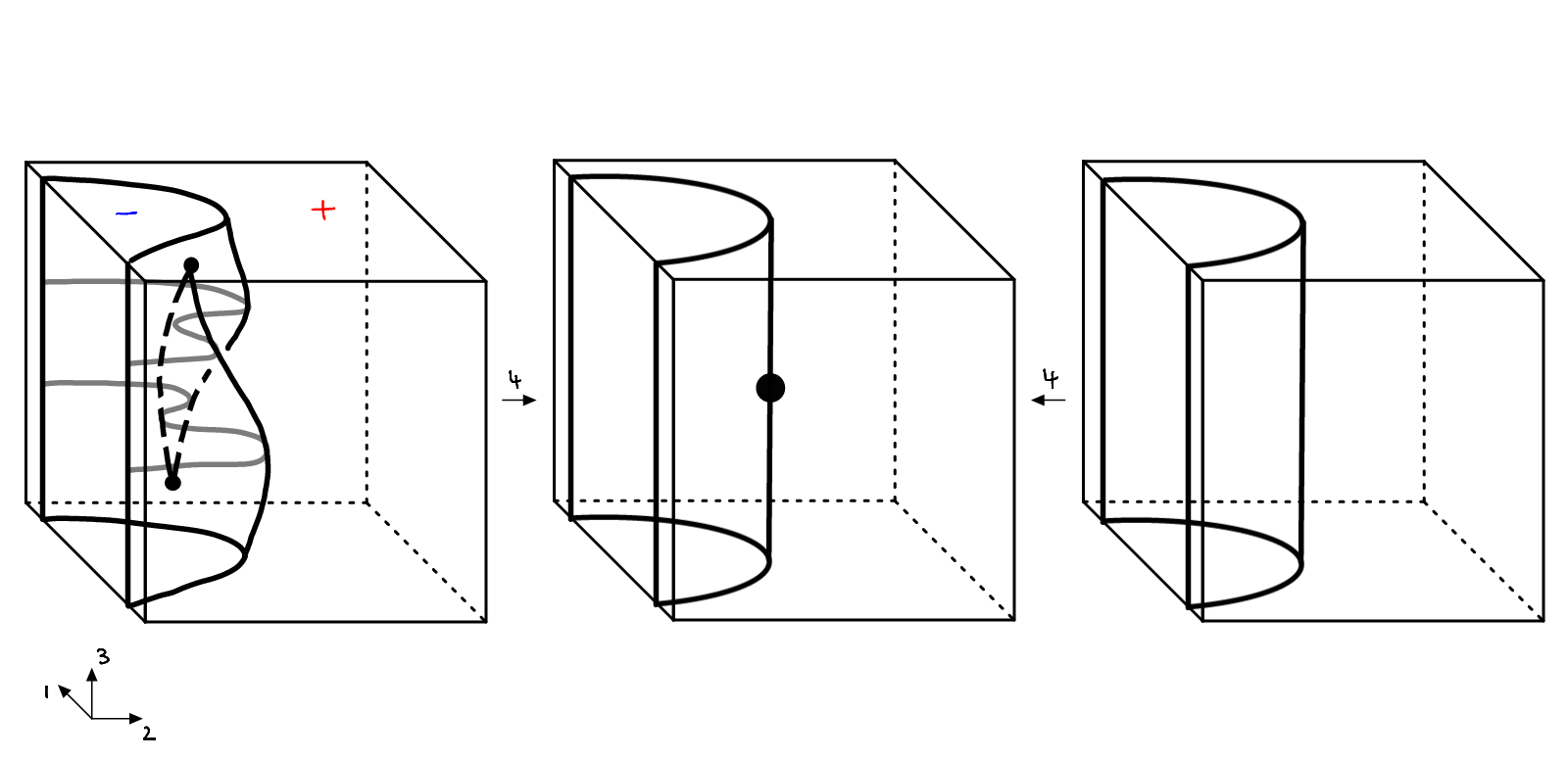}
\endgroup\end{restoretext}
The $3$-morphism on the left contains two snake singularities, while the $3$-morphism on the right is the identity on a cap. This $4$-generator is in fact known as the \textit{swallowtail} singularity \cite{vicary2018globular}. The surfaces only differ by the types of singularities of their ``height function" in direction $2$, which are marked in \cblack{} in the $3$-cube.
\end{rmk}

\section{Connection to CW complexes} \label{sec:CWcomplexes}

We now turn our attention towards classical conceptions of space, and investigate their relation to our algebraic model of presented associative $\infty$-groupoids (and $\infty$-categories). As it turns out, this relation will be phrasable in elementary terms, and provide insights about classical analogues to the algebraic constructions in this thesis. 

To pre-face this section, we re-iterate our warning that, as we now enter the realm of (set-theoretical) geometry and topology our arguments will often not be fully rigorous.

As a concrete model of spaces we choose CW-complexes, and we sketch a translation of CW-complexes into associative $n$-groupoids. The central tool in this translation will be a generalisation of the so-called Thom-Pontryagin construction, which (in both its classical and generalised form) we discuss first. 

\subsection{Generalised Thom-Pontryagin construction}

The goal of this section is to sketch the generalised Thom-Pontryagin construction, which will play a crucial role in translating CW-complexes into groupoids. We recall the following.

\begin{defn}[Framed manifolds and cobordism] Let $Y$ be a $(n+k)$-manifold. A closed $k$-submanifold $M \subset Y$ together with a trivialisation $\eps_M$ of the normal bundle $N(M)$, that is, a bundle isomorphism
\begin{equation}
\eps_M : M \times \lR^n \iso N(M)
\end{equation}
is called \textit{$(n+k)$-framed (in $Y$)}. If $Y = \lR^{n+k}$ then this in particular implies that the stabilisation $T(M) \oplus \lR^n$ of the tangent bundle can be trivialised. 

A cobordism $C : M \equiv N$ between $(n+k)$-framed $k$-manifolds $M, N$ (in $Y$) is a $(n+k+1)$-framed $(k+1)$-manifold $C \subset [0,1] \times Y$ (in $\lR \times Y$) such that firstly
\begin{equation}
\partial C = \Big(\Set{0} \times M~\Big) \sqcup \Big(\Set{1} \times N~\Big)
\end{equation}
where $\sqcup$ denotes the disjoint union, and secondly, the framing of $M$ and $N$ can be recovered from the framing of $C$ by restriction to $\Set{0} \times Y \iso Y$ respectively $\Set{1} \times Y \iso Y$.

Cobordism induces an equivalence relation on the set of $(n+k)$-framed $k$-manifolds in $Y$. Quotienting by this relation gives the set $\Omega^{\mathrm{fr}}_k(Y)$.
\end{defn}

We now briefly recall the basic idea of the Thom-Pontryagin construction (see \cite{pontrjagin2007smooth}). To differentiate this construction from subsequent ideas, we will also refer to it as \textit{based} Thom-Pontryagin construction.

\begin{constr}[Based Thom-Pontryagin construction] Let $S^m$ be the $m$-sphere. The Thom-Pontryagin construction establishes a bijection
\begin{equation}
\pi_{n+k}(S^n) \iso \Omega^{\mathrm{fr}}_k(\lR^{n+k})
\end{equation}
We will describe the mapping from left to right. Let $x \in \pi_{n+k}(S^n)$ be represented by $f : S^{n+k} \to S^n$. Note $S^{m} \iso \lR^m\cup \Set{\infty}$. Using an appropriate homotopy we can assume the following properties of $f$
\begin{enumerate}
\item \textit{Smoothness}: We assume $f$ to be smooth with respect to the usual differential structure on $S^m$.
\item \textit{Transversality}: We assume that $f$ is transversal at $p$, for a fixed choice of $p \in \lR^n \subset S^n$
\item \textit{Basepoint preservation}: We assume $f\inv(\infty) = \infty$.
\end{enumerate}
This implies $f\inv(p) \subset \lR^{n+k}$ is a $k$-dimensional submanifold of $\lR^{n+k}$. The tangent space $T_p$ of $\lR^n$ at $p$ induces a trivialisation of the normal bundle $f\inv(p)$. Thus
\begin{equation}
f\pbstar  p := f\inv(p)
\end{equation}
is a $(n+k)$-framed $k$-manifold. This is a representative of  the image of $x$ under the above mapping $\pi_{n+k}(S^n) \to \Omega^{\mathrm{fr}}_k(\lR^{n+k})$.
\end{constr}

\begin{eg}[Hopf map] \label{eg:hopf_map} We will exemplify the Thom-Pontryagin construction in the case of the Hopf map $h :  S^3 \to S^2$. The map is illustrated below
\begin{restoretext}
\begingroup\sbox0{\includegraphics{ANCimg/page1.png}}\includegraphics[clip,trim=0 {.0\ht0} 0 {.05\ht0} ,width=\textwidth]{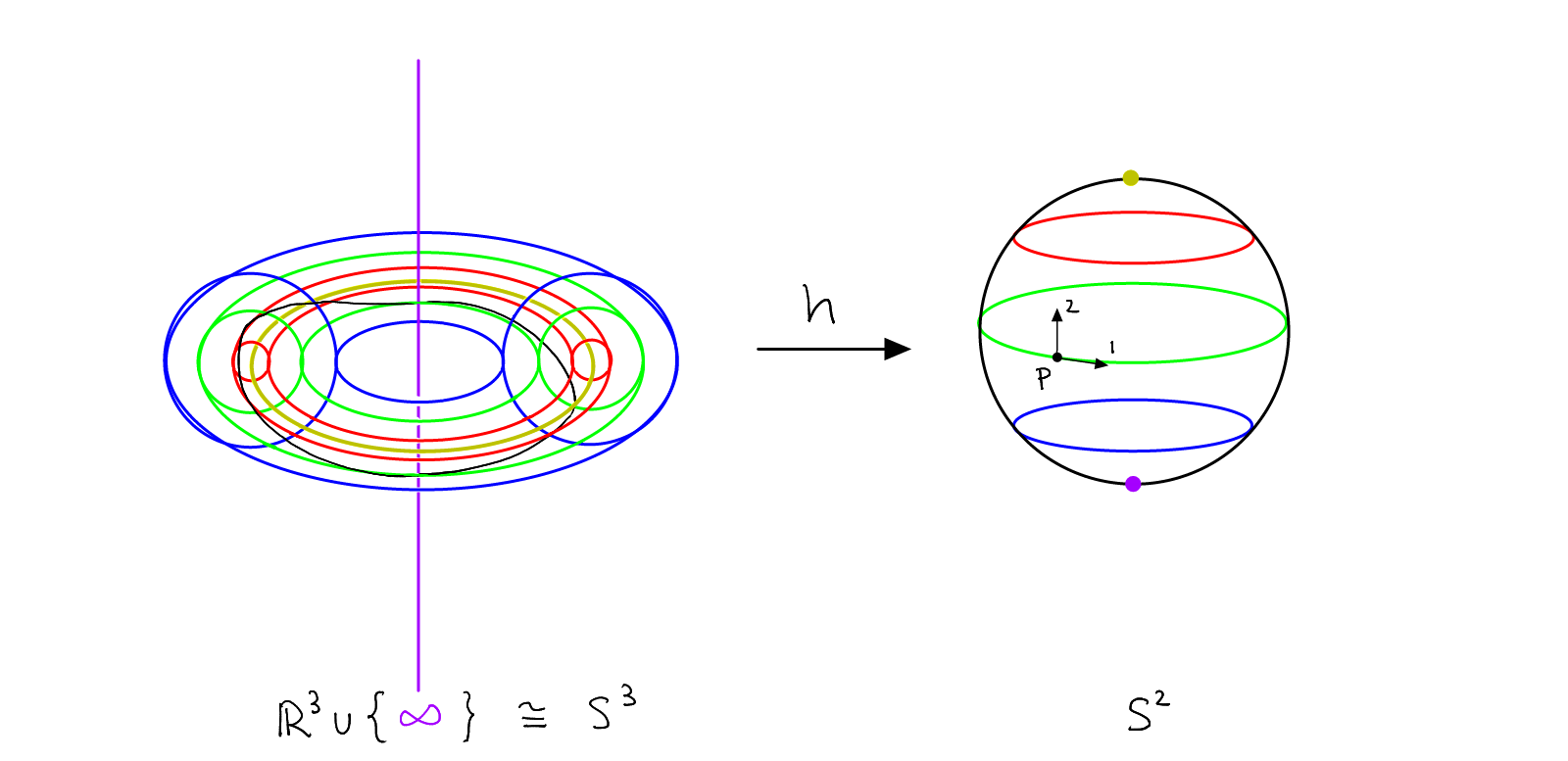}
\endgroup\end{restoretext}
On the left we identity the domain of $h$, $S^3$, with $\lR^3 \cup \Set{\infty}$. For any $q \in S^2$, the preimages $h\inv(q)$ are circles. For latitudal circles, preimages are tori. Preimages are color-coded in the above: the preimage of the \cpurple{} point (the ``south pole") is the circle running through $\infty$ indicated by the \cpurple{} line. The preimage of the \cyellow{} point (the ``north pole") is the \cyellow{} circle. The preimage of the \cred{} (resp. \cdarkgreen{} and \cblue{}) circle is the \cred{} (resp. \cdarkgreen{} and \cblue{}) torus, which we sketched by four touching circles of the same color. The preimage of the \cblack{} point $p$ (lying on the \cdarkgreen{} circle) is the thin \cblack{} circle running in the \cdarkgreen{} torus. Identifying $S^2 \iso \lR^2 \cup \Set{\infty}$ the Hopf map restricts to a smooth map $\lR^3 \to \lR^2$. Choosing tangent vector $1$ and $2$ at $p$ as indicated by arrows above, pulls back to give a trivialisation of the normal vector bundle of $h\inv(p)$ as follows
\begin{restoretext}
\begingroup\sbox0{\includegraphics{ANCimg/page1.png}}\includegraphics[clip,trim=0 {.15\ht0} 0 {.3\ht0} ,width=\textwidth]{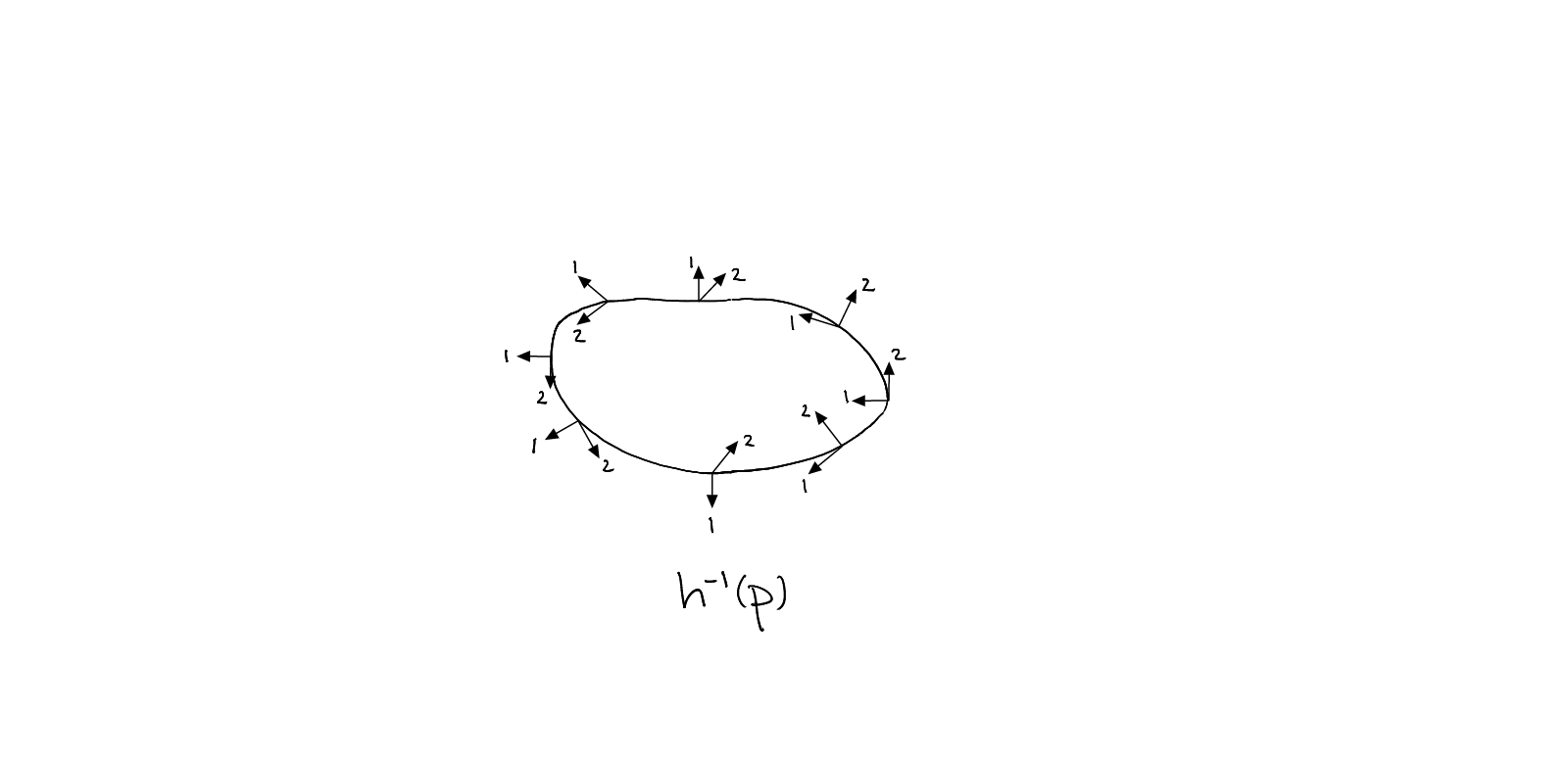}
\endgroup\end{restoretext}
Note that the normal vector frame ``rotates" by $2\pi$ when running around the circle $h\inv(p)$. A quick and informal way to see this is the following: note that vector $1$ was chosen tangential to the \cdarkgreen{} circle. Thus the corresponding normal vector is always tangential to the \cdarkgreen{} torus. Vector $2$ is orthogonal to the \cdarkgreen{} circle and ``points" towards the \cred{} circle. Thus the corresponding normal vector always is normal to the \cdarkgreen{} torus, and points from the \cdarkgreen{} torus to the \cred{} torus.

In conclusion, the Thom-Pontryagin construction associates to the Hopf map the circle in $\lR^3$ with a normal vector frame rotating by $2\pi$ (when compared to the canonical normal framing of $S^1 \subset \lR^2 \subset \lR^2 \oplus \lR^1$).
\end{eg}

While $\pi_{n+k}(S^n)$ considers based maps, we also claim the following unbased version of the Thom-Pontryagin construction.

\begin{constr}[Unbased Thom-Pontryagin construction] The unbased Thom-Pontryagin construction gives a bijection
\begin{equation}
[S^{n+k},S^n] \iso \Omega^{\mathrm{fr}}_k(S^{n+k})
\end{equation}
Here, the left-hand side denotes homotopy classes of maps from $S^{n+k}$ to $S^n$. The construction of the map from left to right is analogous to the based Thom-Pontryagin construction, but we can drop any mention of the basepoints $\infty \in S^m \iso \lR^m \cup \Set{\infty}$. Namely, given $x \in [S^{n+k},S^n]$ we chose $f \in x$ which is smooth and transversal at $p \in S^n$. Then $f\inv(p)$ represents the image of $x$ in $\Omega^{\mathrm{fr}}_k(S^{n+k})$ (and obtains its framing from $T_p S^n$ as before).
\end{constr}

The unbased Thom-Pontryagin construction can be generalised from homotopy classes of maps $f : S^{n+k}\to S^n$ of spheres to homotopy classes of maps $f : S^{n+k} \to X^{(n)}$ from a sphere to the ($n$-skeleton) of a CW-complex. For a mathematical rigorous version of the construction (in the piecewise linear case) we refer the reader to \cite{buonchristiano1976geometric}, Proposition VII.4.1. We will only give a sketch of the construction here. We recall the following.

\begin{defn}[CW complexes] \label{defn:CW_complex} Let $D^n \subset \lR^n$ be the $n$-disk, $i^\partial: \partial D^n = S^{n-1} \into D^n$ its boundary (the $(n-1)$-sphere) and $i^\circ: (D^n)^\circ \into D^n$ its interior (as a subspace of $\lR^n$). Let $X$ be a space, and $f : S^{n-1} \to X$ a continuous map. By an attachment $X \cup_f D^n$ of $D^n$ on $X$ along $f$ we mean the pushout of spaces
\begin{equation}
\xymatrix{ S^{n-1} \ar[r]^{f} \ar[d]_{i^\partial} & X \ar[d] \\
D^n \ar[r] & X \cup_f D^n \pushoutfar} 
\end{equation}
Note that the right leg of the pushout gives an inclusion $X \subset X \cup_f D^n$. Let $I_k$ be a set of labels and let $\ata_g : S^{k -1} \to X$, $g \in I_k$ be attaching maps of $k$-disks. Then $X \cup_{\ata_g} \Set{g \in I_k}$ is the space obtained as the colimit of the transfinite sequence of gluings
\begin{equation}
X \into (X \cup_{\ata_{g_1}} D^k) \into ((X \cup_{\ata_{g_1}} D^k)\cup_{\ata_{g_2}} D^k) \into \dots
\end{equation}
This can be seen to be independent of any ordering $g_1, g_2, g_2, \dots$ chosen on $I_k$.

A CW complex is a space inductively constructed as follows. $X^{(0)}$ is a discrete space (that is, a union of points $x \in \sN^X_0$). $X^{(k)}$ is obtained from $X^{(k-1)}$ by gluings of $k$-disks. That is
\begin{equation}
X^{(k)} = X^{(k-1)} \cup_{\ata_g} \Set{g \in \sN^X_k}
\end{equation}
for some choice of \textit{labels} $\sN^X_k$ and \textit{attaching maps} $\ata_g : S^{k-1} \to X^{(k-1)}$, $g \in \sN^X_k$. $X$ is the space obtained as the colimit of the sequence
\begin{equation}
X^{(0)} \into X^{(1)} \into X^{(2)} \into ...
\end{equation}
$X^{(n)}$ is called the $n$-skeleton of $X$. We also note that the pushout leg $D^k \to  X^{(k-1)} \cup_{\ata_g} D^k$ composes with $X^{({k-1})} \cup_{\ata_g} D^k \into X$ to give a map $\cell_g : D^k \into X$, called the \textit{cell in $X$ labelled by $g$}. We set
\begin{equation}
\cpartial_g := \cell_g i^\partial : S^{k-1} \to X
\end{equation}
called the boundary of $\cell_g$. Note $\cpartial_g$ factors through $X^{({k-1})}$ and this factorisation equals $\ata_g$. We further set
\begin{equation}
\ccirc_g := \cell_g i^\circ : (D^k)^\circ \to X
\end{equation}
called the interior of $\cell_g$. Note that $\ccirc_g$ has image in $X^{(k)}$. Further, note that $\ccirc_g$ is a homeomorphism with its image. In this sense, we say that the subspace $\im(\ccirc_g)\subset X$ inherits smooth structure from $D^k$.
\end{defn}

\begin{egs}[CW complexes] We give three examples of CW-complexes. The first will be the ``half-cut circle" $X$ which has the following skeleta $X^{(0)} \subset X^{(1)} \subset X^{(2)}$.
\begin{restoretext}
\begingroup\sbox0{\includegraphics{ANCimg/page1.png}}\includegraphics[clip,trim=0 {.25\ht0} 0 {.25\ht0} ,width=\textwidth]{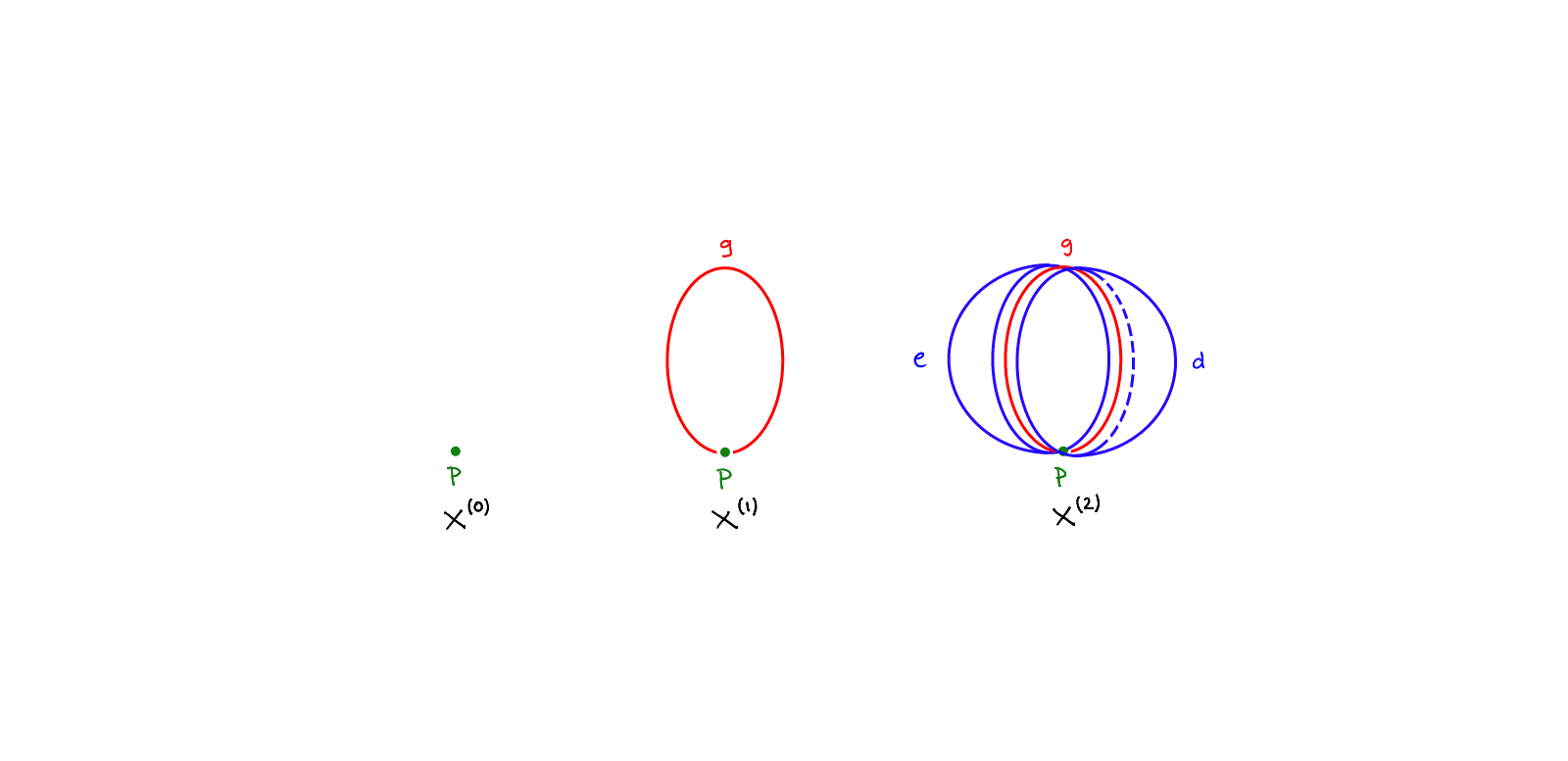}
\endgroup\end{restoretext}
Thus, $X^{(0)}$ is a point $p$ to which we glue a $1$-disk $g$ on its endpoint to obtain $X^{(1)}$. $X^{(1)}$ is thus a circle, to which we glue two disks $e$ and $d$ to obtain $X^{(2)}$.

The second example is the ``two-point torus" $X$, whose skeleta are illustrated in the following picture
\begin{restoretext}
\begingroup\sbox0{\includegraphics{ANCimg/page1.png}}\includegraphics[clip,trim=0 {.15\ht0} 0 {.2\ht0} ,width=\textwidth]{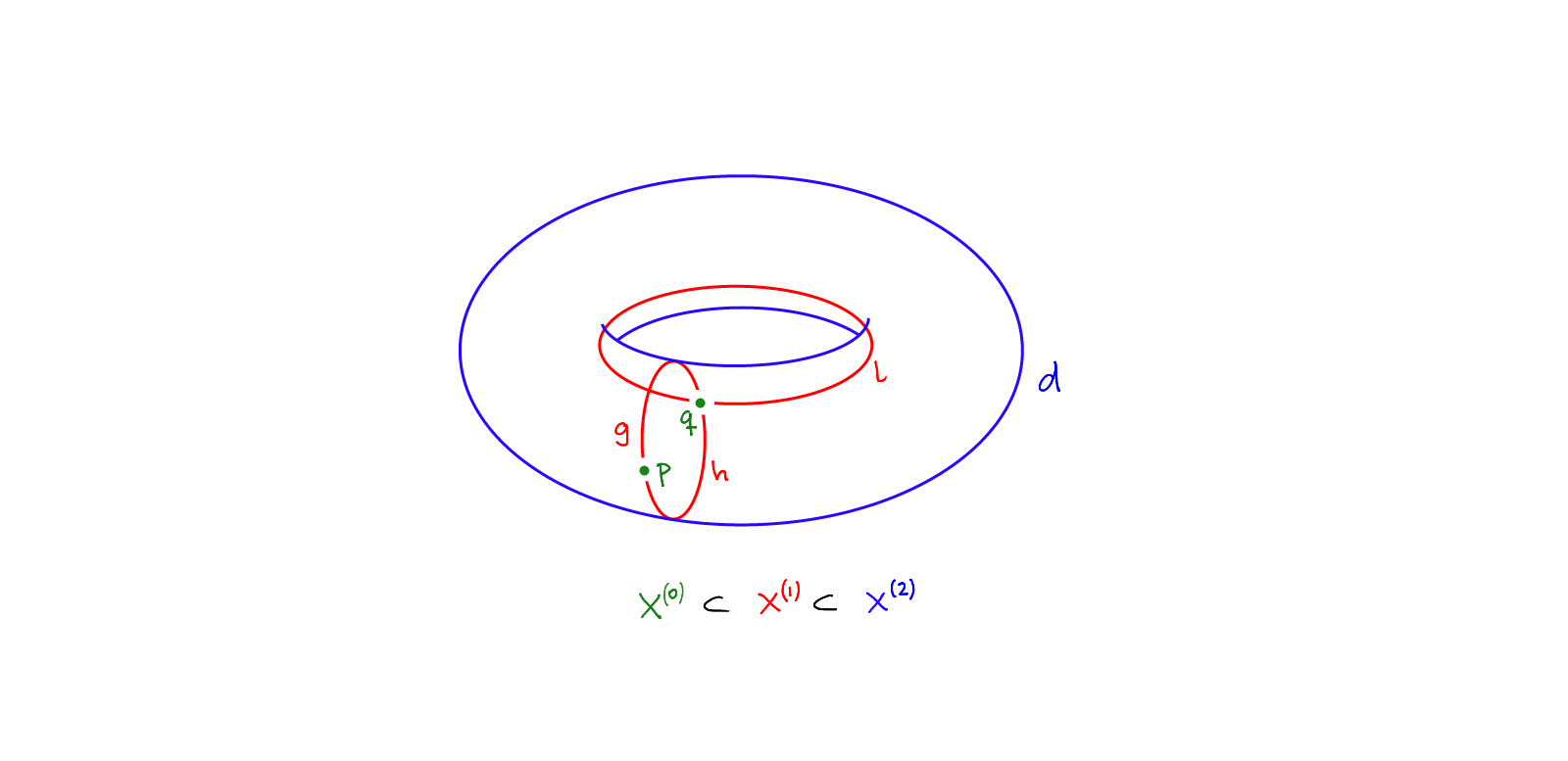}
\endgroup\end{restoretext}
Thus $X^{(0)}$ consists of two points $p$ and $q$, to which we glue three $1$-disks $g$, $h$ and $l$ to obtain $X^{(1)}$. By an appropriate gluing of a $2$-disk $d$ this yields the torus.

Our final example is called the complex projective plane. The skeleta of this space $X$ are (symbolically) illustrated as follows
\begin{restoretext}
\begingroup\sbox0{\includegraphics{ANCimg/page1.png}}\includegraphics[clip,trim=0 {.05\ht0} 0 {.0\ht0} ,width=\textwidth]{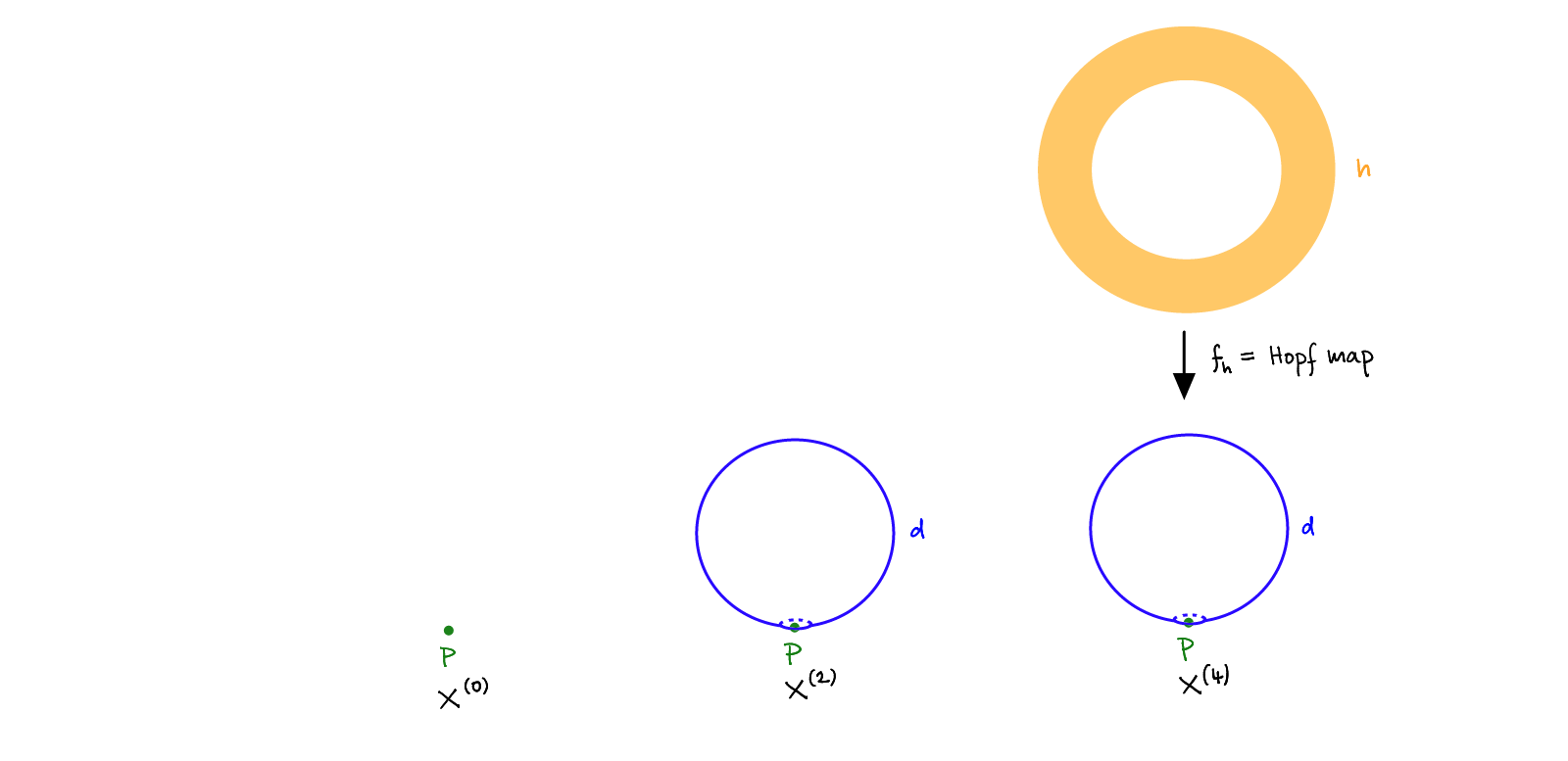}
\endgroup\end{restoretext}
$X^{(0)}$ is a single point $p$, to which we glue a $2$-disk $d$ to obtain $X^{(2)}$. $X^{(2)}$ is a $2$-sphere, and $X^{(4)}$ is then obtained by attaching a $4$-disk $h$ with attaching map $\ata_h$ being the Hopf map.
\end{egs}

Next, we sketch a generalisation of the notion of framed manifold in $X$ to that of framed stratification of $X$. The definition should be compared with that of manifold diagram in \autoref{ssec:po_mfld_diag}.

\begin{defn}(Framed $n$-stratifications) \label{defn:framed_strat} Let $Y$ be a $n$-manifold and let $\sN$ be a namescope of dimension $\infty$. A framed $n$-stratification $\bM$ of $Y$ in the namescope $\sN$ consists of sets $\bM_0$, $\bM_1$ ... $\bM_n$ of disjoint $n$-framed submanifolds of $Y$ covering $Y$,  which are identified with subsets of the name sets by injections $\mathsf{name}^k_\bM : \bM_k \subset \sN_k$. They satisfy the following
\begin{enumerate}
\item \textit{Dimension}: $f \in \bM_k$ is a $(n-k)$-dimensional open manifold
\item \textit{Boundary}: For $f \in \bM_k$ we have
\begin{equation}
\partial f = \bigcup_{h \subset \overline{f}} h
\end{equation}
\end{enumerate}
Two framed $n$-stratifications of $Y$, $\bM$, $\bN$ in the namescope $\sN$ are cobordant, written $\bM \equiv \bN$, if there is a framed $n$-stratification $\bC$ of $[0,1] \times Y$ in the namescope $\sN$, such that each framed manifold $f$ in $\sC_i$ yields a cobordism for of the corresponding elements in $g \in \bM_i$ and $h \in \bN_i$. In other words, for any $f$ we have $g$, $h$ with
\begin{equation}
\mathsf{name}^k_\bC (f) = \mathsf{name}^k_\bM (g) = \mathsf{name}^k_\bN (h)
\end{equation}
and $f$ is a cobordism $g \equiv h$. 

Quotienting the set of framed $n$-stratifications of $Y$ by this relation gives the set $\Psi^{\mathrm{fr}}_{\sN}(Y)$.
\end{defn}

We now turn to the generalised unbased Thom-Pontryagin construction. This goes hand in hand with a notion of framed dual stratification $\kD(X)$ for a given CW-complex $X$. The following construction gives a sketch of both notions in a mutually inductive fashion.

\begin{constr}[Framed dual stratifications and generalised unbased Pontraygin] \label{constr:dual_complex} The \textit{dual stratification} $\kD(X)$ of a CW-complex $X$ is a collection of disjoint subsets $g\dualdag \subset X$, for each label $g \in \sN^X_l$ of the complex, covering $X$. $g\dualdag$ is called the \textit{dual stratum of $g$}. The construction of the dual stratification $\kD(X)$ is doubly inductive: it is ordered by dimension $k$ of skeleta and (transfinitely) ordered by attaching cells 
\begin{equation}
\ig_1 \leq \ig_2 \leq ... \leq \ig_\alpha \leq \ig_{\alpha+1} \leq ...
\end{equation}
which we assume to be ordered in a chosen way at each dimension $k$. 
Let $X^{(m)}_{k,\alpha}$ denote the $m$-skeleton of the complex that has all cells attached up to $k$-cells of order number less or equal $\alpha$ (note if $m < k$ then $X^{(m)}_{k,\alpha} = X^{(m)}$). We inductively define strata
\begin{equation}
g\dualdag_{k,\alpha} = g\dualdag \cap X^{(k)}_{k,\alpha}
\end{equation}
of a stratification $\kD(X^{(k)}_{k,\alpha})$ of $X^{(k)}_{k,\alpha}$. Passing to the colimit over $\alpha$ we obtain strata
\begin{equation}
g\dualdag_{k} = g\dualdag \cap X^{(k)}
\end{equation}
of a stratification $\kD(X^{(k)})$ of $X^{(k)}$. $\kD(X)$ and its strata $g\dualdag$ will then be obtained as a colimit over $k$.

We claim the following inductive properties in the $(k,\alpha)$th step of the construction
\begin{enumerate}
\item \textit{Strata retracts}: Let $m \leq k$. Define
\begin{equation}
X^{k,\alpha}_{(m)} = \bigcup_{g \in I^X_{l \leq m}} g\dualdag_{k,\alpha}
\end{equation}
We inductively claim $X^{k,\alpha}_{(m)} \subset X^{(k)}_{k,\alpha}$ is open and that we have a retract
\begin{equation}
\retr^{k,\alpha}_m : X^{k,\alpha}_{(m)} \to X^{(m)}_{k,\alpha}
\end{equation}
such that $\retr^{k,\alpha}_m(g\dualdag_{k,\alpha}) = g\dualdag_{k,\alpha} \cap X^{(m)}_{k,\alpha}$ for $g \in I^X_{l \leq m}$.

\item \textit{Strata pullbacks}: Let $m \leq k$ and $h \in \sN^X_m$ such that $h \leq g_\alpha$. There is a framed $m$-stratification $\bM(h)$ of $D^m$ in the namescope $\sN^X$, consisting of manifolds denoted by
\begin{equation}
\bM(h)_l = \Set{h\pbstar g\dualdag ~|~ g \in \sN^X_l, g\dualdag_k \cap \im(\cell_h) \neq \emptyset }
\end{equation}
We chose (cf. \autoref{defn:framed_strat}) 
\begin{equation}
\mathsf{name}^l(h\pbstar g\dualdag) = g
\end{equation}
$h\pbstar  g\dualdag$ is called the pullback stratum of $g\dualdag$ to $h$. We inductively claim $h\pbstar h\dualdag \in (D^m)^\circ$.

\end{enumerate}

The base case $k = 0$ is trivial: $X^{(0)}_{0,\alpha}$ consists of  0-cells $\ig_\beta$, $\beta \leq \alpha$, which we define to be their own dual strata $(\ig_\beta)\dualdag_{0,\alpha} := \ig_\beta$. For strata pullbacks we set $(\ig_\beta)\pbstar(\ig_\beta)\dualdag := D^0$.

Now let $k \geq 1$. To complete the inductive construction of $\kD(X^{(k)}_{k,\alpha + 1})$ and $\retr^{k,\alpha + 1}_m$, assume $\kD(X^{(k)}_{k,\alpha})$ and $\retr^{k,\alpha}_m$ to be defined.  We first give a generalised Thom-Pontryagin construction for a closed $(k-1)$-manifold $M$, which will provide us with a mapping 
\begin{equation}
\kP_M : [M, X^{({k-1})}] \to \Psi^{\mathrm{fr}}_{\sN^X}(M)
\end{equation}
Let $f \in [M,X^{({k-1})}]$ be a homotopy class. We assume the existence\footnote{That such a choice is possible, is a conjecture which we don't attempt to prove. We only point out the analogy to the step in the classical based Thom-Pontryagin construction. Therefore, this choice might more appropriately be read as a restriction of our discussion to CW-complexes whose attaching maps can be chosen to have the stated properties of ``smoothness" and ``transversality"} of a same-named representative $f$ of this class, which, in analogy with the classical Thom-Pontryagin construction, satisfies two conditions of ``smoothness" and ``transversality" defined below. We write $f_{(m)}$ for the restriction of $f$ to the inverse image of $X_{(m)}$ under $f$, and for $g \in I^X_m$, noting $\im(\ccirc_g)$ is open in $X^{(m)}$, we write $f_{(g)}$ for the restriction of $\retr^{k-1}_m f_{(m)}$ to the inverse image of $\im(\ccirc_g)$ under $\retr^{k-1}_m f_{(m)}$. For each $g \in I^X_{m \leq k-1}$, we then require the following
\begin{enumerate}
\item \textit{Smoothness}:  We require $f_{(g)}$ to be smooth.
\item \textit{Transversality}: We require $f_{(g)}$ to be transversal at $g\pbstar g\dualdag_m$ (note $\im(\ccirc_g) \iso (D^m)^\circ$ canonically).
\end{enumerate}
As a consequence of these assumptions we find
\begin{equation}
f\inv_{(g)} (g\pbstar g\dualdag_{m}) \subset M
\end{equation}
is a $(k-1)$-framed $(k-1-m)$-manifold in $M$. Note that this manifold equals $f\inv (g\dualdag_{k,\alpha})$ because of our inductive assumptions. As another consequence of the inductive assumption, we find that the collection of $f\inv_{(g)} (g\pbstar g\dualdag_{m})$, $g \in \sN^X_l$, $l < k$, forms a framed $(k-1)$-stratification $\kP_M(f)$ of $M$ in the namescope $\sN^X$. This provides the mapping for the generalised unbased Thom-Pontryagin construction as stated above. As an aside, we claim (without proof) that analogous to the classical case this mapping is injective.

We now chose $M = S^{k-1}$ (in which case we denote $\kP_M$ by just $\kP$), and let $f$ equal the attaching map $\ata_{h} \in [S^{k-1},X^{({k-1})}]$ of the $(\alpha + 1)$th $k$-cell $h := \ig_{\alpha + 1}$. Using the generalised unbased Thom-Pontryagin construction we obtain
\begin{equation}
\kP(\ata_{h} ) \in \Psi^{\mathrm{fr}}_{\sN^X}(S^{k-1})
\end{equation}
which is a framed stratification of $S^{k-1}$.

Next consider the quotient space
\begin{equation}
D^k \iso \frac{S^{k-1} \times [0,1]}{S^{k-1} \times \Set{0}}
\end{equation}
In other words, $D^k$ is the cone of $S^{k-1}$ (with vertex point $p = S^{k-1} \times \Set{0}$). We define $h\dualdag_{k,\alpha+1} := \cell_h(p)$ and $h\pbstar  h\dualdag = p$. Thus $h\dualdag_{k,\alpha+1}$ is the point given by the image of the vertex point. A $k$-framing of the $0$-dimensional manifold $h\pbstar  h\dualdag$ in $D^k \subset \lR^k$ can be chosen arbitrarily (up to orientation there are only two choices, which categorically correspond to choosing a cell or its inverse). We further define for each $g \in \sN^X_m$, $m < k$,
\begin{equation}
h\pbstar g\dualdag := (\ata_{h})\inv_{(g)} (g\pbstar g\dualdag) \times [0,1) \subset D^k
\end{equation}
with framing inherited from the framing of $g\pbstar g\dualdag$ in $\bM(g)$. Note that the collection of $h\pbstar g\dualdag$ defined using the coning approach above together with the vertex point $h\pbstar h\dualdag$ forms a framed $k$-stratification $\bM(h)$ of $D^k$ in the namescope $\sN^X$.

We define a framed stratification $\kD^k(X^{({k})}_{k,\alpha + 1})$ of $X^{({k})}_{k,\alpha + 1} = X^{({k})}_{k,\alpha} \cup_{\ata_h} D^k$, by defining its strata as (for $g \in I^X_l, l \leq k, g < h$)
\begin{equation}
g\dualdag_{k,\alpha + 1} = \cell_{h}(h\pbstar g\dualdag) \cup g\dualdag_{k, \alpha}
\end{equation}
For $m \leq k$, we extend $\retr^{k,\alpha}_m$ to a map
\begin{equation}
\retr^{k,\alpha + 1}_m : X^{k,\alpha + 1}_{(m)} \to X^{(m)}_{k,\alpha+1}
\end{equation}
which is the identity if $m = k$, and otherwise ``retracts" the punctured cell $h$ before using the inductively defined retract. Explicitly, it first maps
\begin{equation}
x \in \im(\ccirc_h) \setminus { h\dualdag_{k,\alpha + 1} } \iso (D^k)^\circ \setminus \Set{p}  \iso S^{k-1} \times (0,1) \ni (q,r) \mapsto \ata_h (q)
\end{equation}
before applying $\retr^{k,\alpha}_m$. This completes the inductive construction of $\kD(X^{(k)}_{k,\alpha})$. By passing to the colimit over $\alpha$ it completes the construction of $\kD(X^{(k)})$. Finally, by passing to the colimit over $k$ it completes the construction of $\kD(X)$.
\end{constr}

We now give examples of the previous construction, discussing the inductive argument step-by-step.

\begin{egs}[Dual complexes] We illustrate the construction of dual complexes in the case of our three previous examples of CW complexes. The inductive steps of construction of the dual complex of the half-cut sphere can be illustrated as follows (we denote $g\dualdag_k$ by $g\dualdag$ for simplicity, as the index $k$ of $g\dualdag$ can be inferred from the $k$-skeleton it is contained in)
\begin{restoretext}
\begingroup\sbox0{\includegraphics{ANCimg/page1.png}}\includegraphics[clip,trim={.25\ht0}  {.25\ht0} {.25\ht0}  {.25\ht0} ,width=\textwidth]{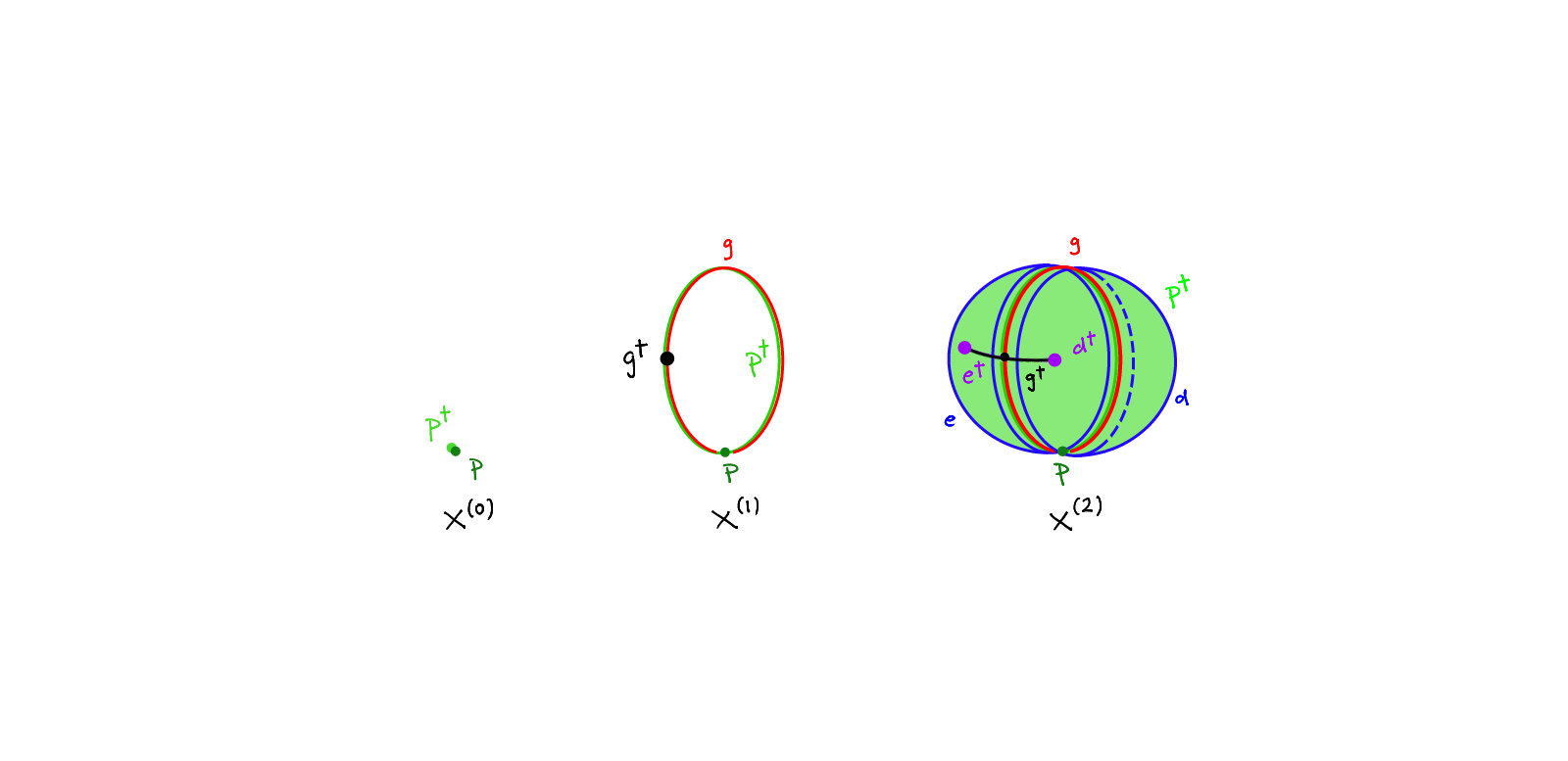}
\endgroup\end{restoretext}
In the first step, the dual $p\dualdag_0$ of the single $0$-cell $p \in X^{(0)}$ is $p$ itself (this is the inductive base case of the previous construction). In the next step, we find $\ata \inv _g (p\dualdag) = S^0$ and by the given ``coning" construction we have
\begin{align}
g\pbstar  p\dualdag &= S^0 \times [0,1) \\
&= D^1 \setminus \Set{g\pbstar  g\dualdag } 
\end{align}
where $g\pbstar  g\dualdag \in D^1$ is a point in the interior of $D^1$ corresponding to the vertex point $S^0 \times \Set{1}$ under the identification
\begin{equation}
D^1 \iso \frac{S^0 \times [0,1]}{S^0 \times \Set{1}}
\end{equation}
Consequently, $g\dualdag_1$ is the point $\cell_g(g\pbstar g\dualdag)$ on the circle $X^{(1)}$ and 
\begin{equation}
p\dualdag_1 = \cell_g(g\pbstar p\dualdag)\cup p\dualdag_{0}
\end{equation}
is the rest of the circle $X^{(1)}$. Now in the last inductive step (building $X^{(2)}$), the procedure turns $g\dualdag_2$ into an ``interval" (marked by a \cblack{} line above), whose endpoints are $e\dualdag_2$ and $d\dualdag_2$ respectively, and $p\dualdag_2$ covers the rest of the sphere. Note that $g\dualdag_2 = g\dualdag$ (and similarly for $p$, $e$, $d$) since $X$ has no cells in dimension higher than $2$. 

It is instructive to visualise the pullback stratifications $\bM(g)$, $\bM(e)$ and $\bM(d)$ of dual strata as follows
\begin{restoretext}
\begingroup\sbox0{\includegraphics{ANCimg/page1.png}}\includegraphics[clip,trim={.2\ht0}  {.3\ht0} {.2\ht0}  {.25\ht0} ,width=\textwidth]{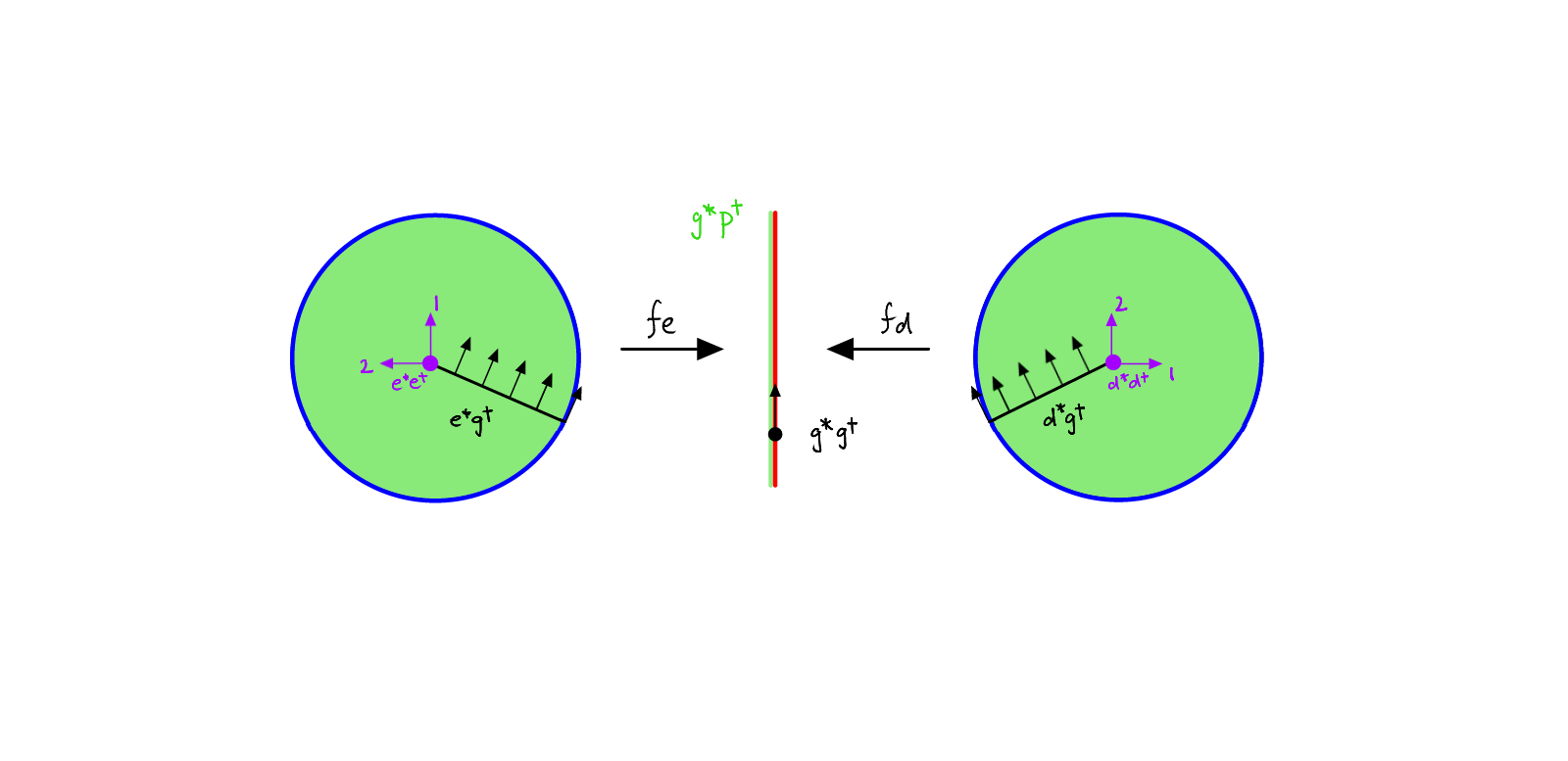}
\endgroup\end{restoretext}
On the left we depicted pullbacks of dual strata into $e$, on the right into $d$ and in the middle into $g$ (note that we over-layed the cell $g$, marked in \cred{}, and its dual strata, marked in \cgreen{} and \cblack{}). Arrows in bold from left to middle and right to middle indicate attaching maps $\ata_e$ and $\ata_d$. Restricted to the inverse image of $\im(\ccirc_g) \iso (D^1)^\circ$ these are the maps $(\ata_e)_{(g)}$ respectively $(\ata_d)_{(g)}$ which are chosen smooth and transversal at $g\pbstar  g\dualdag$. Their respective inverse images of $g\pbstar g\dualdag$ are then the $1$-framed $0$-manifold given by $e\pbstar g\dualdag$ and $d\pbstar g\dualdag$ restricted to $S^1$ (which is the circle marked in \cblue{} above). In the next step, these are extended to $D^2$ by ``coning" $S^1$, obtaining $2$-framed $1$-manifolds $e\pbstar g\dualdag$ and $d\pbstar g\dualdag$ as shown above. We finally find
\begin{equation}
g\dualdag_2 = g\dualdag_1 \cup \cell_e(e\pbstar g\dualdag) \cup \cell_d(d\pbstar g\dualdag)
\end{equation}
Similarly, one determines $p\dualdag_2$, $e\dualdag_2$ and $d\dualdag_2$. Note that, the framings of $e\pbstar e\dualdag$ and $d\pbstar d\dualdag$ can be chosen arbitrarily at this stage, and are indicated by \cpurple{} arrows above. This completes the construction of our first example.

As our second example, we consider the previously defined two-point torus. Its dual strata (namely $p\dualdag$, $q\dualdag$, $g\dualdag$, $h\dualdag$, $l\dualdag$ and $d\dualdag$) are illustrated in the following graphic
\begin{restoretext}
\begingroup\sbox0{\includegraphics{ANCimg/page1.png}}\includegraphics[clip,trim=0 {.2\ht0} 0 {.2\ht0} ,width=\textwidth]{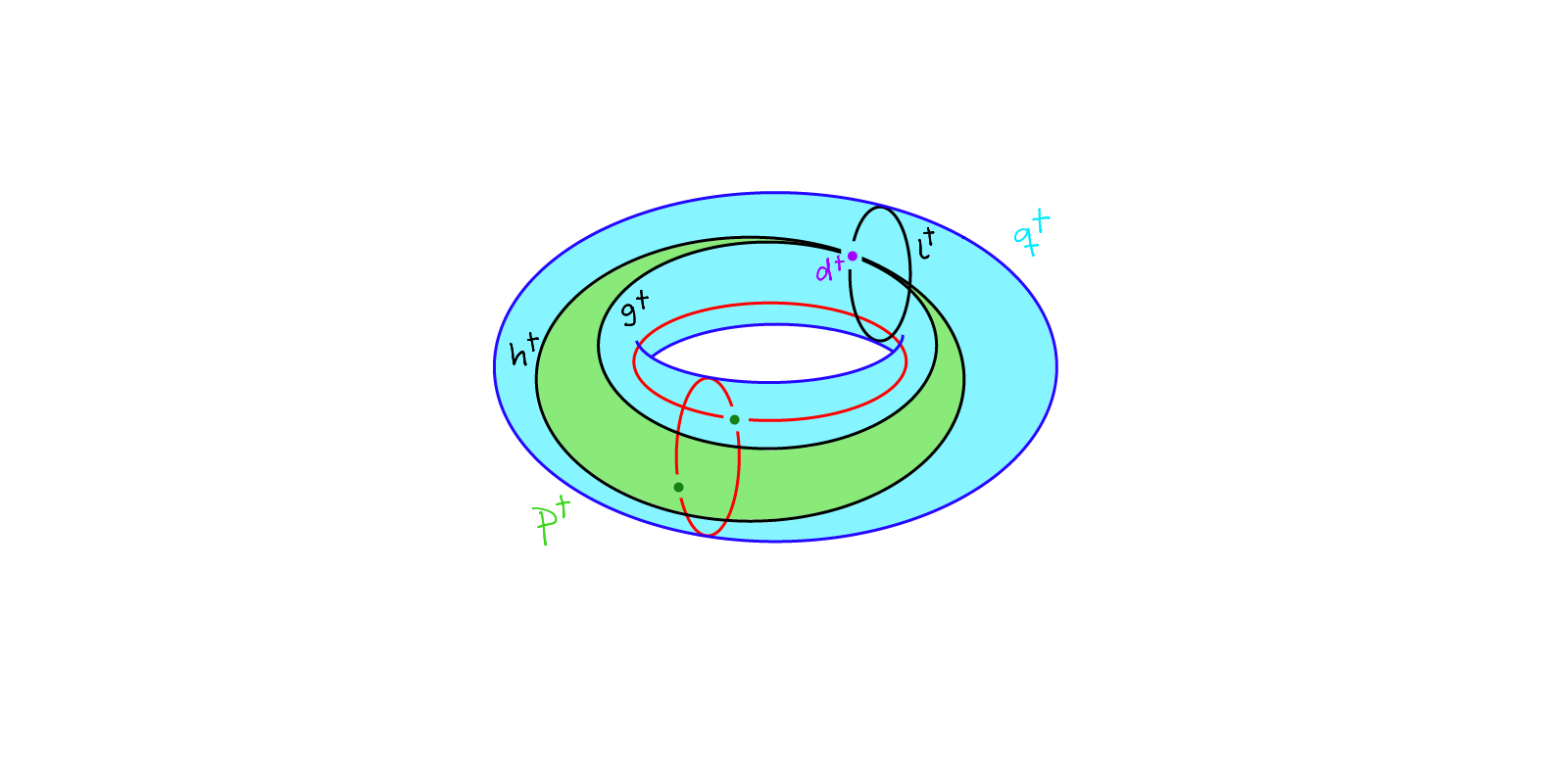}
\endgroup\end{restoretext}
Again, it is instructive to depict the pullback of dual strata together with their framing
\begin{restoretext}
\begingroup\sbox0{\includegraphics{ANCimg/page1.png}}\includegraphics[clip,trim=0 {.0\ht0} 0 {.0\ht0} ,width=\textwidth]{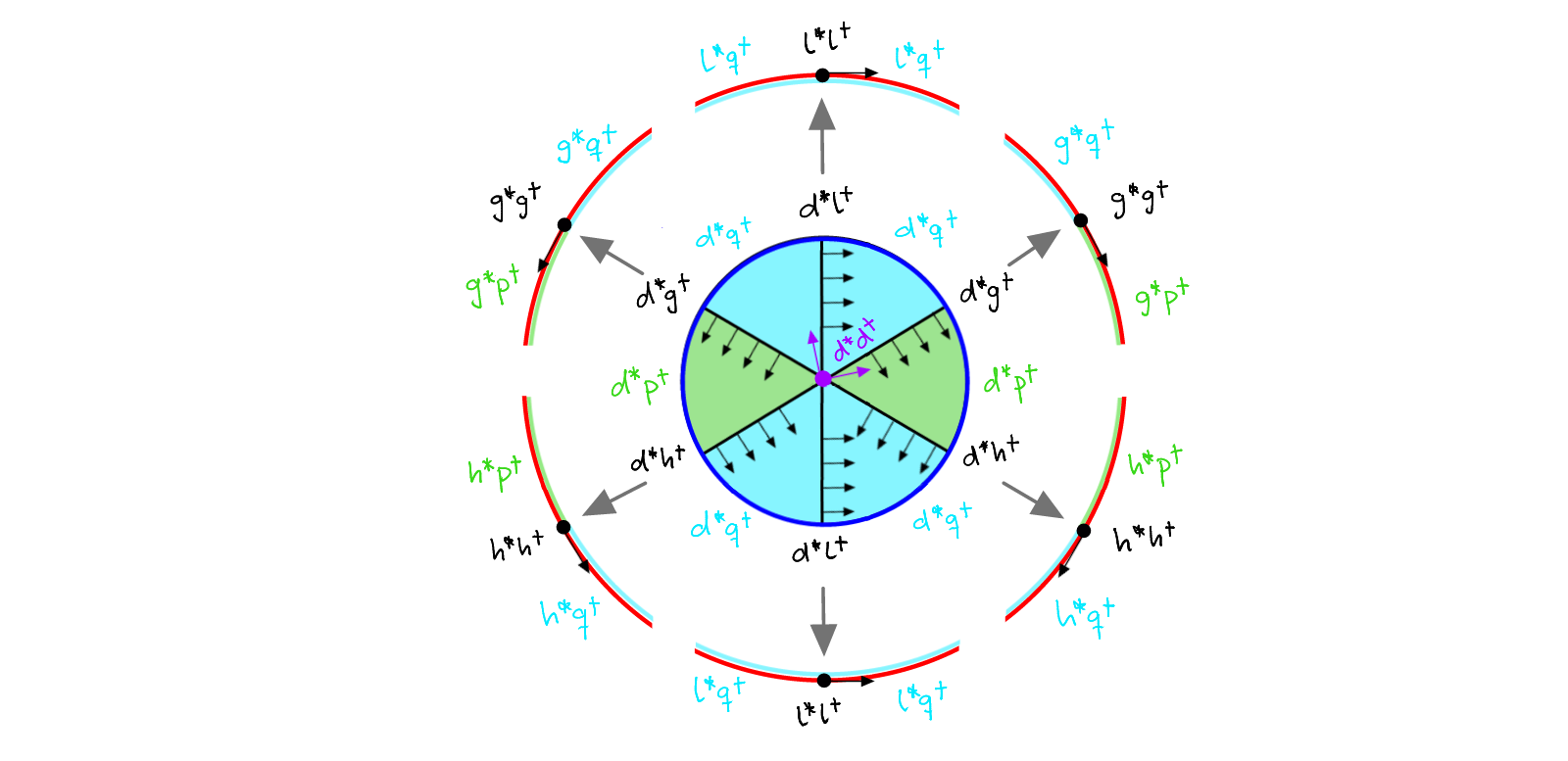}
\endgroup\end{restoretext}
In the middle, we pictured the pullback stratification $\bM(d)$ of $D^2$, and on the outside pullback stratifications $\bM(g)$, $\bM(h)$ and $\bM(l)$ of $D^1$ (note that all of latter appear twice). Bold  arrows indicate the attachment map $\ata _d$ of $d$ (each of $g$, $h$ and $l$ are covered twice by this map). Once more, we note that this attachment map can be chosen to be ``smooth" and ``transversal" as required in the construction.

Finally, we consider the complex projective plane $X$. Its (inductively constructed) dual complex can be illustrated as follows
\begin{restoretext}
\begingroup\sbox0{\includegraphics{ANCimg/page1.png}}\includegraphics[clip,trim=0 {.0\ht0} 0 {.0\ht0} ,width=\textwidth]{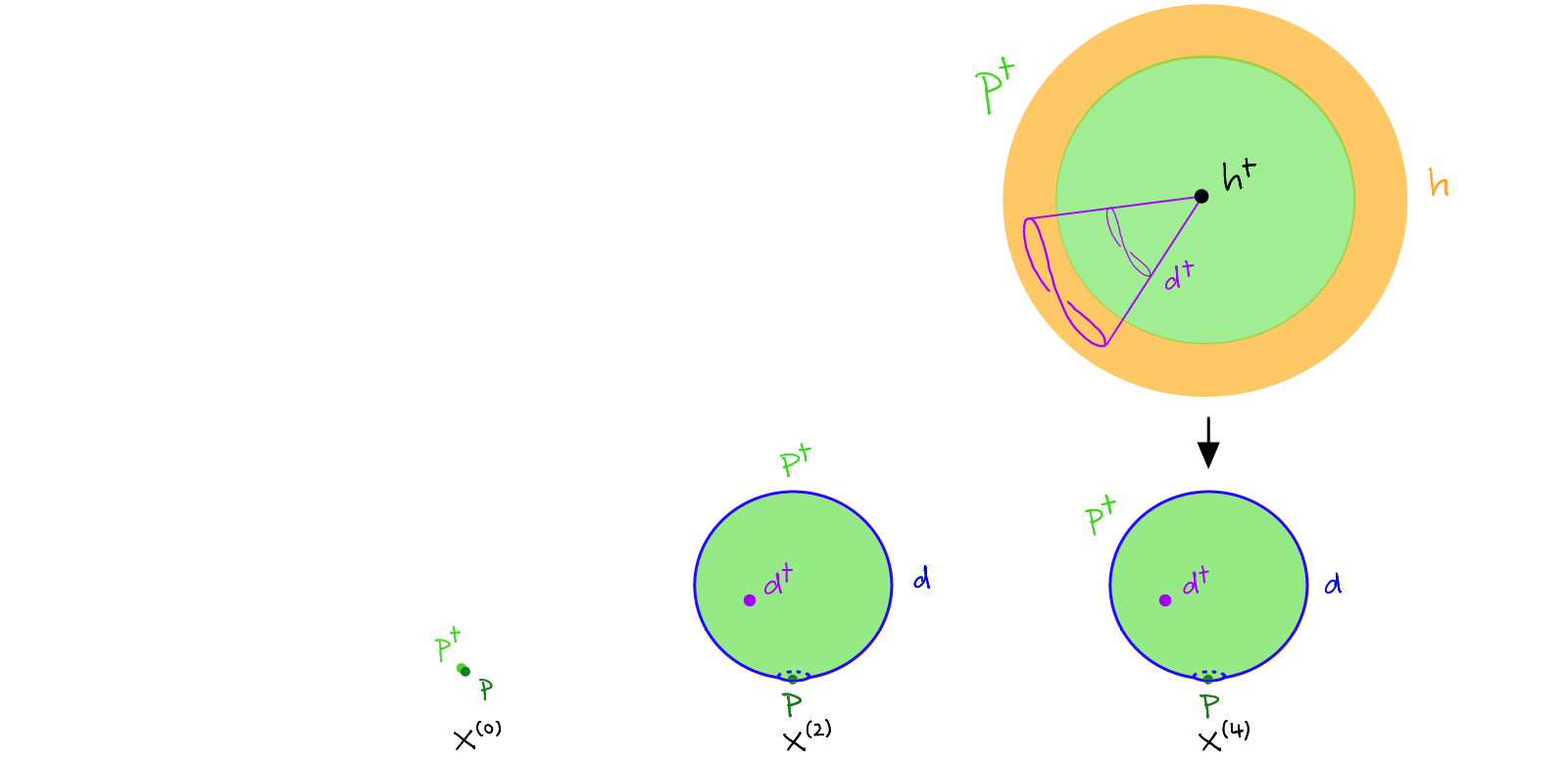}
\endgroup\end{restoretext}
The pullback stratification $\bM(d)$ of $D^2$ is just
\begin{restoretext}
\begingroup\sbox0{\includegraphics{ANCimg/page1.png}}\includegraphics[clip,trim=0 {.25\ht0} 0 {.2\ht0} ,width=\textwidth]{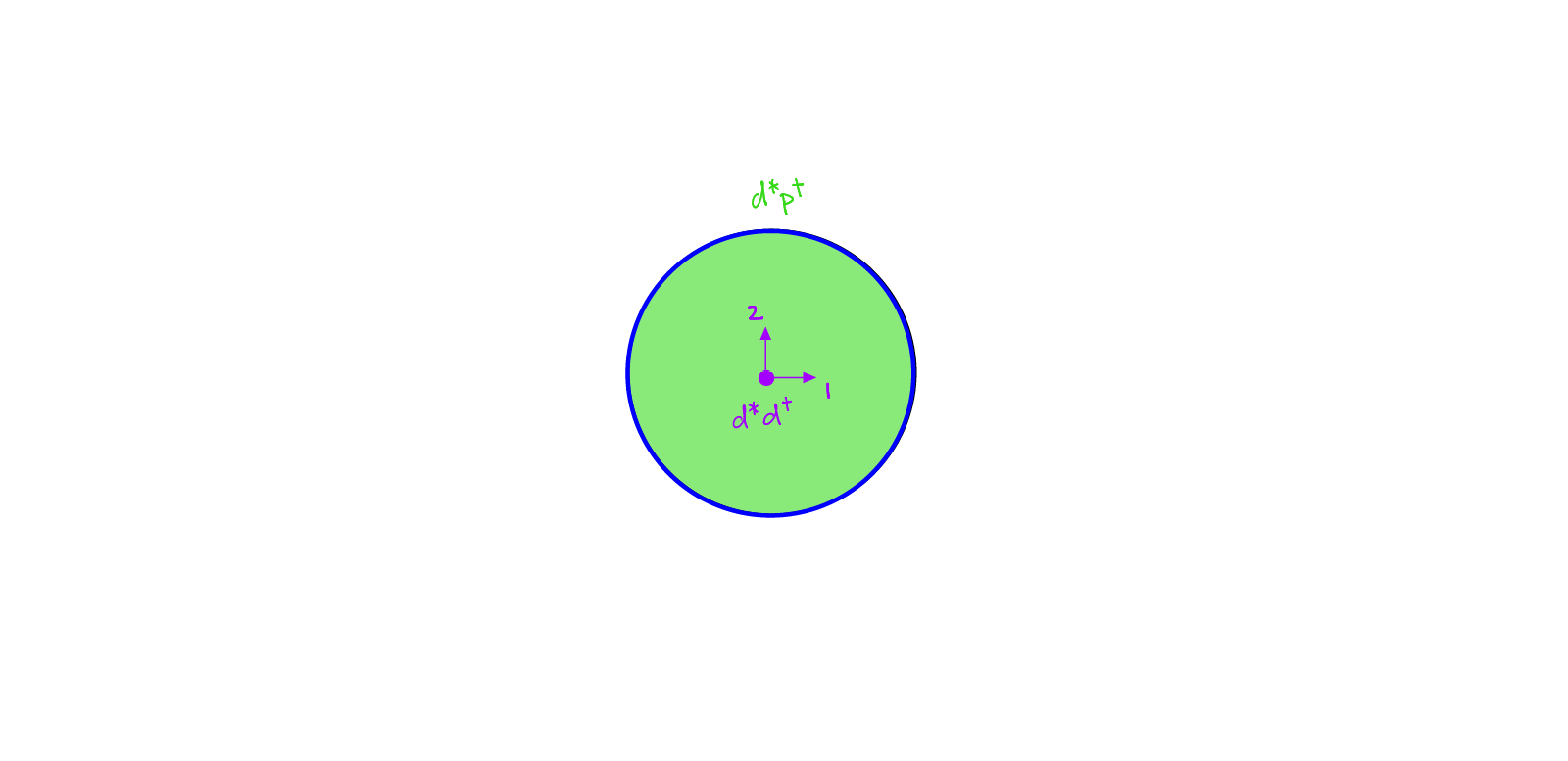}
\endgroup\end{restoretext}
Recalling our discussion of the Hopf map in \autoref{eg:hopf_map}, the pullback stratification $\bM(d)$ of $D^4$ is illustrated as follows
\begin{restoretext}
\begingroup\sbox0{\includegraphics{ANCimg/page1.png}}\includegraphics[clip,trim=0 {.0\ht0} 0 {.0\ht0} ,width=\textwidth]{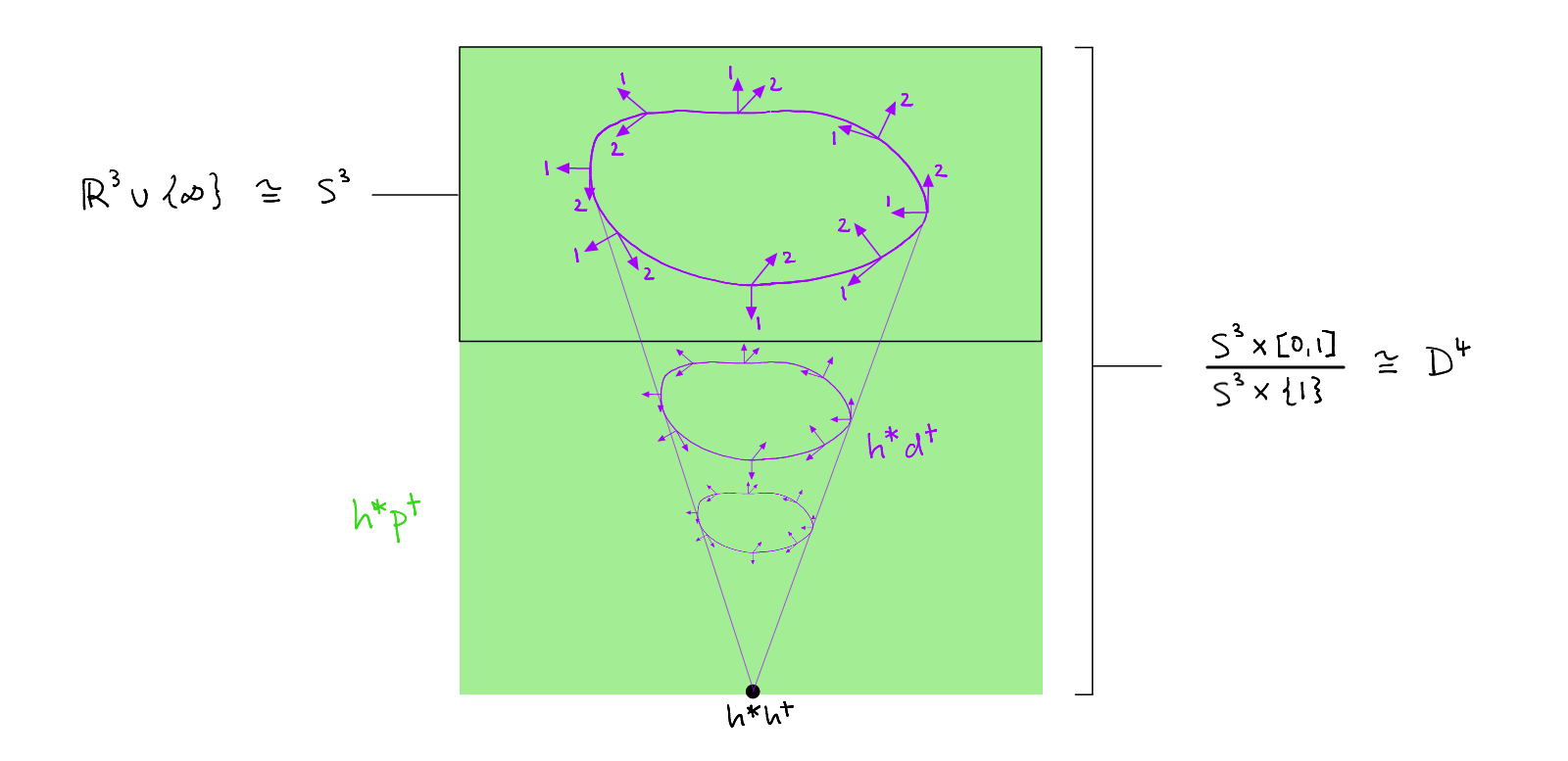}
\endgroup\end{restoretext}
This is an representation of submanifolds of $D^4$, with the top part (marked by a \cblack{} rectangle) representing $S^3$. $\ata\inv_h (d\dualdag) \subset S^3$ is (as we have found in \autoref{eg:hopf_map}) a $3$-framed $1$-manifold with frame rotating by $2\pi$. Consequently $h\pbstar d\dualdag$ is a ``coning" of this framed manifold (marked in \cpurple{}) but excludes its vertex point point $h\pbstar h\dualdag$ (marked in \cblack{}). $h\pbstar p\dualdag$ fills the rest of the space in $D^4$. Based on our discussion of the Hopf map it should be clear the attachment map can be chosen to be ``smooth" and ``transversal" as required in the preceding construction.
\end{egs}

\subsection{Framing for morphisms of groupoids}

In this section we will heuristically explain that morphisms of presented associative $n$-groupoids have framings on the manifolds in the manifold diagrams corresponding to their morphisms. While our explanation will not be mathematically precise, the underlying idea has turned out to be ``precise enough" to work out many interesting computations in the setting of presented associative $n$-groupoids.

We sketch this explanation in three steps: firstly, we will explain framings of $\Comp(\TI)$, then framings of coherent invertibility data of an invertible generator, and finally framings of general morphism in groupoids. It is important to emphasize once more that we are working with manifold diagrams of groupoids: this will allow us to regard a generator and its coherent invertibility data as a \textit{single manifold} (in contrast to the distinction of individual labelled regions as manifolds, see \autoref{ssec:po_mfld_diag}).

\begin{enumerate}
\item \textit{Framings for morphisms in \TI}:
We first explain how framings can be found for manifold diagrams corresponding to morphisms in $\TI$. Namely, combinatorial cobordisms (morphisms in $\TI$) can be given a one-dimensional normal frame on a manifold which is the \textit{union of all manifolds labelled by some $\ic_{S \equiv T}$}, which excludes manifolds labelled by $-$ and $+$. We assume this manifold can be represented smoothly\footnote{We leave a discussion of the existence of these ``smooth representations" of manifold diagrams corresponding to morphism in $\TI$ to future work, and content ourselves with an empirical argument that such representations seem to be naturally attainable in many low-dimensional examples with no evident obstructions in higher dimensions.} (that is, we assume such a representative in its equivalence class, cf. \autoref{ssec:po_mfld_diag}). Then, the framing of the resulting smooth manifold is given by having the single normal vector point from the regions labelled by ``$-$" to the regions labelled by ``$+$" everywhere. 

We illustrate this with several examples. First consider the case of $1$-generators, where we find the framing
\begin{restoretext}
\begingroup\sbox0{\includegraphics{ANCimg/page1.png}}\includegraphics[clip,trim=0 {.35\ht0} 0 {.25\ht0} ,width=\textwidth]{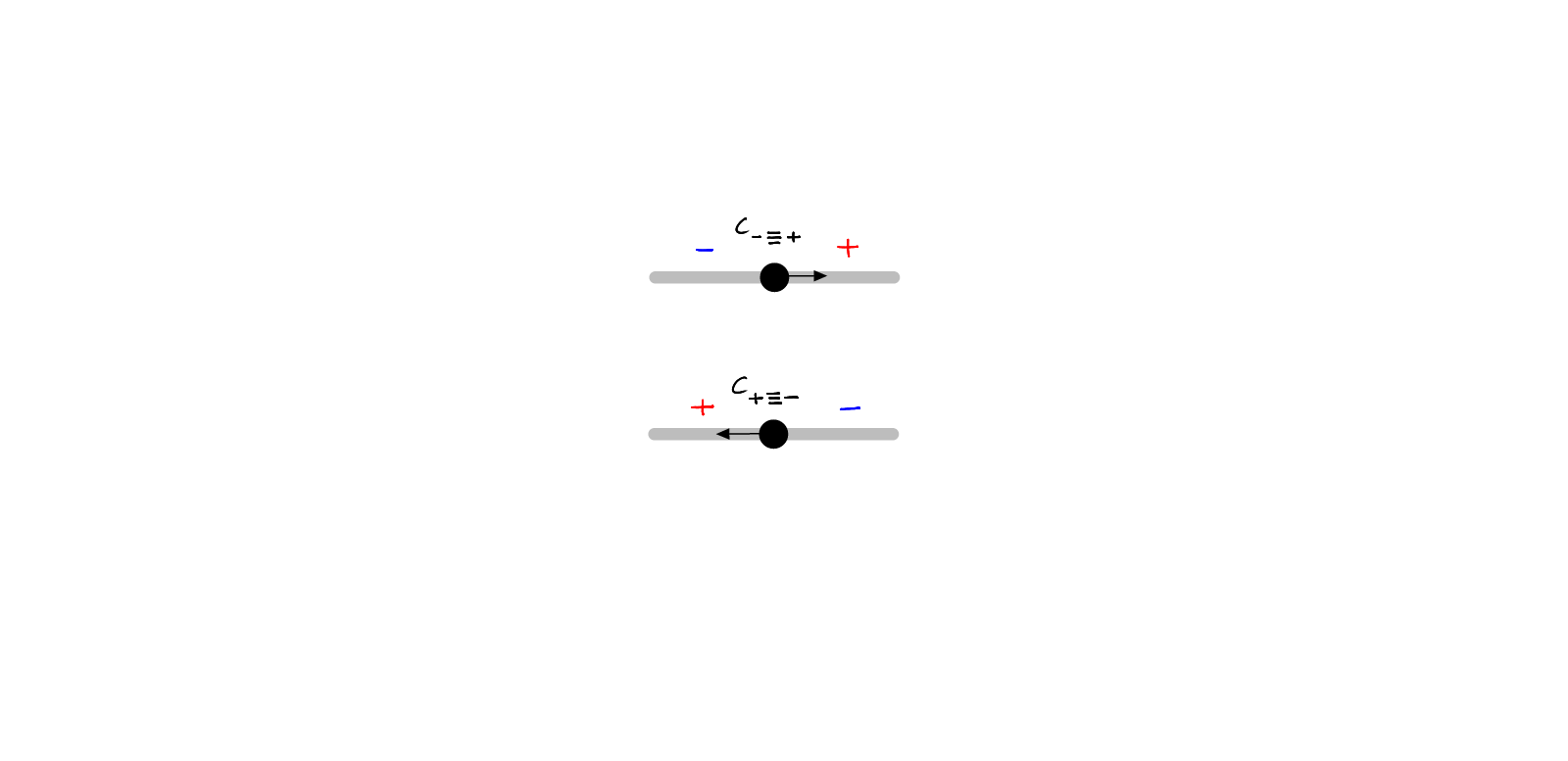}
\endgroup\end{restoretext}
In this case we thus have a 1-framed 0-manifold given by the manifold corresponding to $\ic_{-\equiv +}$ (respectively $\ic_{+\equiv -}$). The normal vector points from ``$-$" to ``$+$" in both cases.

Similarly, $2$-generators (that is, cup and cap singularities) obtain a framing given by
\begin{restoretext}
\begingroup\sbox0{\includegraphics{ANCimg/page1.png}}\includegraphics[clip,trim=0 {.0\ht0} 0 {.1\ht0} ,width=\textwidth]{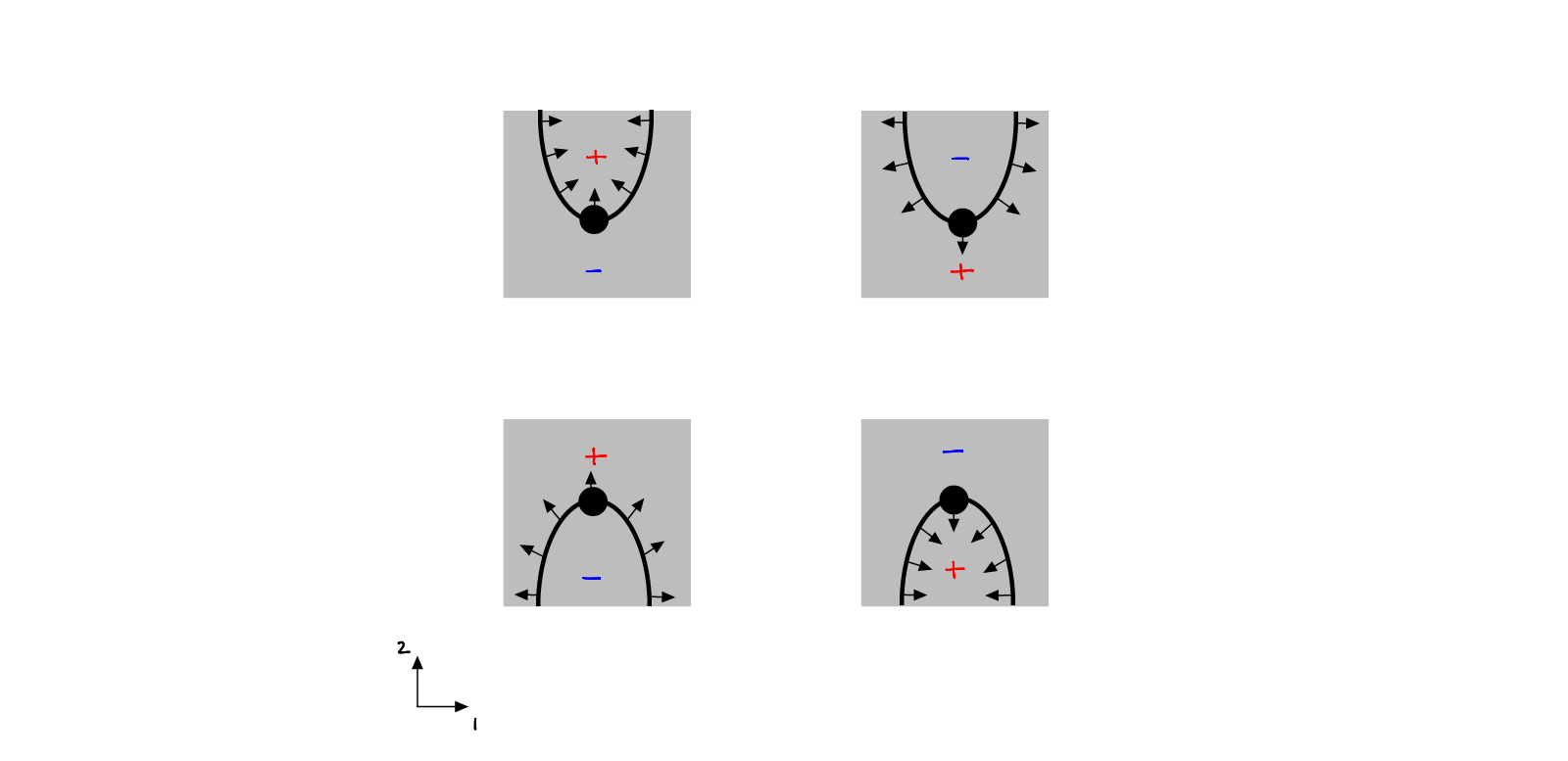}
\endgroup\end{restoretext}
Here, we obtained $1$-framed $1$-manifolds, and the manifolds are given as a union of submanifolds corresponding to generators in $\TI_1$ and $\TI_2$. For instance, the upper left manifold is the union of the three manifold labelled by $\ic_{-\equiv +}$, $\ic_{+\equiv -}$ and $\ic_{\abss{\ic_{-\equiv +}} \whisker 1 1 \abss{\ic_{+ \equiv -}} \equiv \Id_{-}}$. The normal vectors again point from ``$-$" to ``$+$".

Another example, the two following $3$-generators (which were previously called saddle, and death of circle singularities), can be framed by
\begin{restoretext}
\begingroup\sbox0{\includegraphics{ANCimg/page1.png}}\includegraphics[clip,trim=0 {.05\ht0} 0 {.1\ht0} ,width=\textwidth]{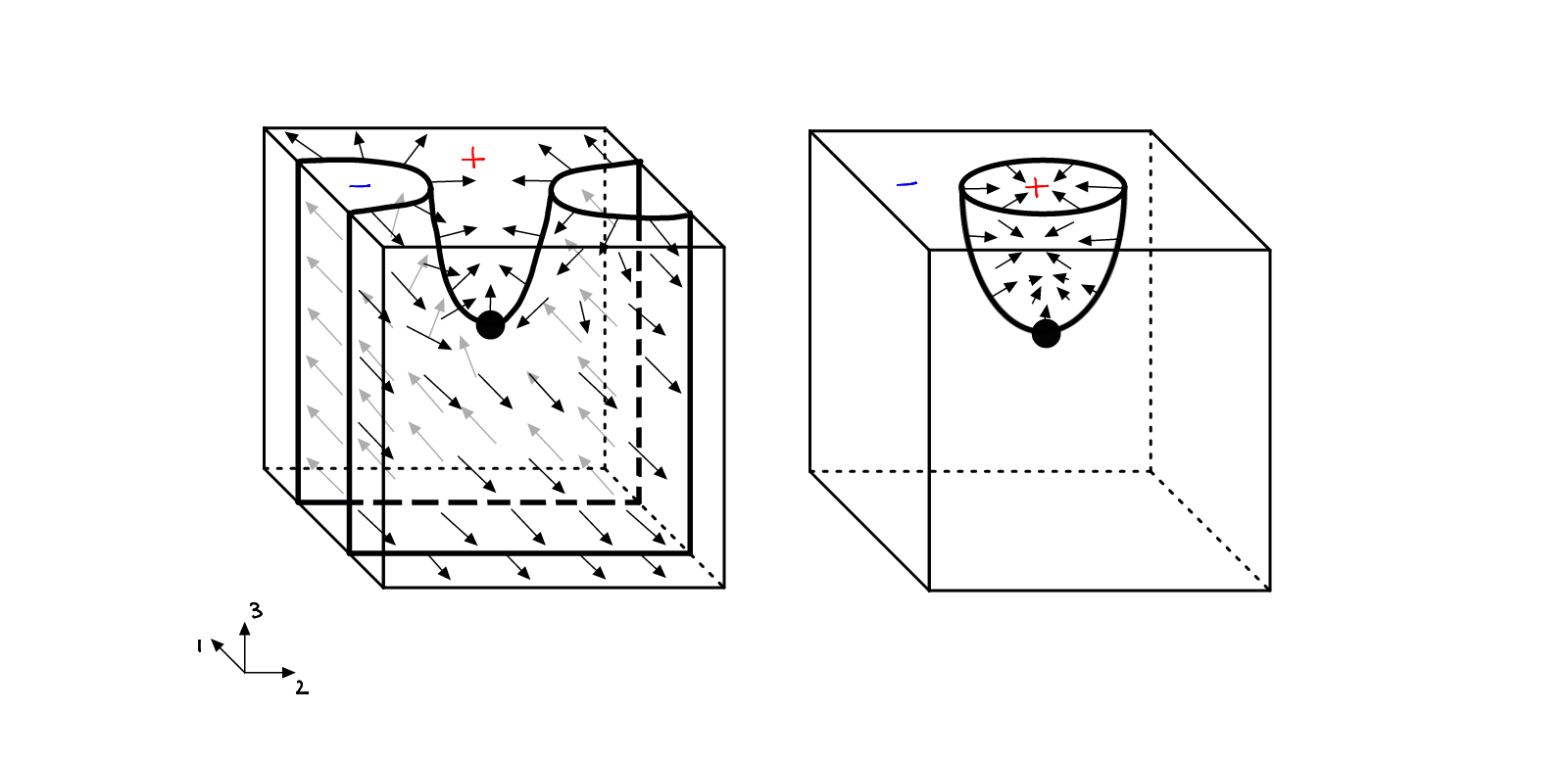}
\endgroup\end{restoretext}

More generally, we claim the idea applies in all dimensions and for all morphisms in $\TI$. The reader can verify that other previously mentioned (and smoothly represented) combinatorial cobordisms, including morphism that are not types of generators as a above, obtain a $1$-framing in the same way.

\item \textit{Framing for types of general invertible generators}:
Next we consider how framings can be obtained on $n$-manifolds containing general invertible element $\ig \in \sC_{m+1}$ with source $x$ and target $y$, together with the singularities $\icg\ig_{S \equiv T}$ belonging to its coherent invertibility data. We start with manifold diagrams of \textit{types} of $(m+k)$-generators $\abss{\icg\ig_{S \equiv T}}$. Within such a diagram, the manifold consisting of the union of all manifolds labelled by some $\icg\ig_{U \equiv V}$ can be given a $(m+1)$-dimensional normal frame which consists of $m$ constant normal vectors in the $m$ directions of the source and target cubes of $\ig$ and a single normal vector obtained analogous to the previous item for $\ic_{S \equiv T}$ after projecting along the first $m$ normal vectors. Indeed, by \eqref{eq:inv_prop2}  this projection recovers the cube of the $k$-generator $\ic_{S \equiv T}$ in $\TI$ and thus the previous item applies.

We illustrate this in the case of \autoref{eg:adjoining_inv_gen} where $m = 1$. We find that the normal framings of the types of $\icg\ig_{- \equiv +}$ and $\icg\ig_{+ \equiv -}$ are given by
\begin{restoretext}
\begingroup\sbox0{\includegraphics{ANCimg/page1.png}}\includegraphics[clip,trim=0 {.2\ht0} 0 {.35\ht0} ,width=\textwidth]{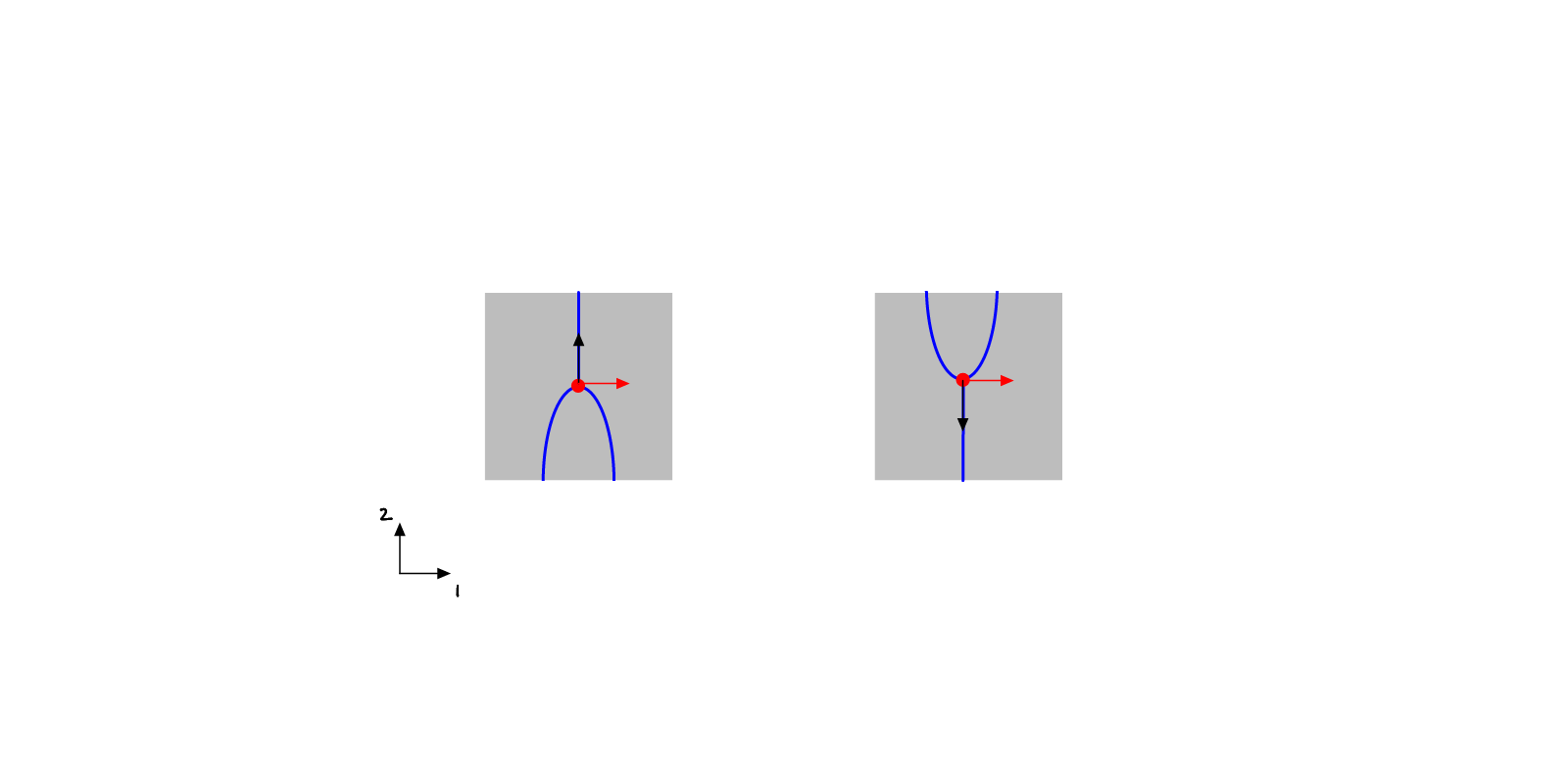}
\endgroup\end{restoretext}
Here, the \cblack{} arrow is the framing of the corresponding generator in $\TI$, and the \cred{} arrow is the normal vector into source/target directions (in this case, this is the direction of coordinate 1 in the coordinate system given above). This should be compared to our previous pictures of $1$-dimensional manifold diagrams for $\ic_{- \equiv +}$ and $\ic_{+ \equiv -}$.

Similarly, cups and caps for $\ig$ obtain framings as follows
\begin{restoretext}
\begingroup\sbox0{\includegraphics{ANCimg/page1.png}}\includegraphics[clip,trim=0 {.0\ht0} 0 {.1\ht0} ,width=\textwidth]{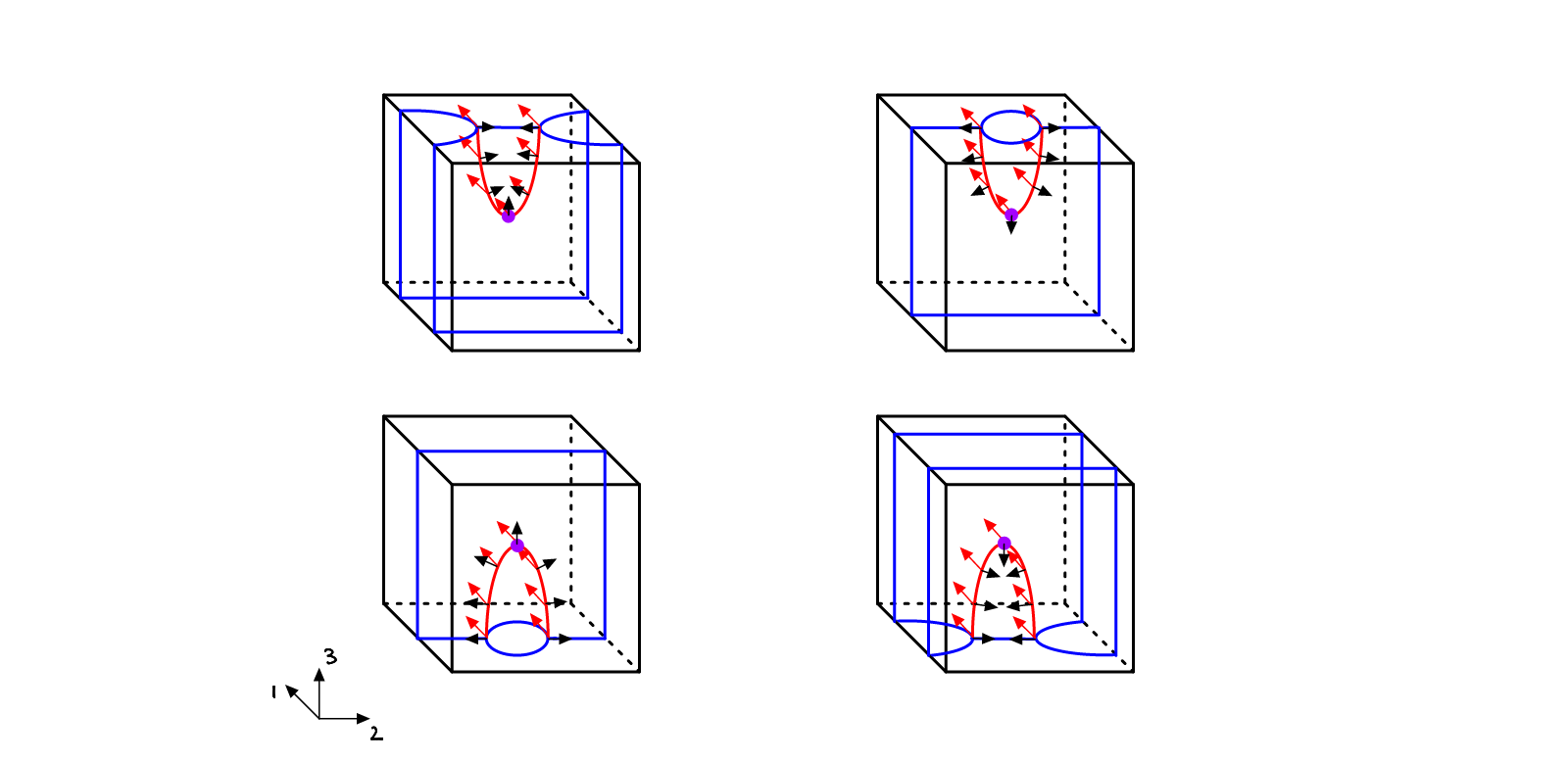}
\endgroup\end{restoretext}
where again \cblack{} arrows are inherited from $\TI$ whereas \cred{} arrows (constantly) point in the directions of source/target (again, this is coordinate direction 1). This should be compared to our pictures of $2$-dimensional manifold diagrams for cups and caps in the previous item.

The procedure extends trivially to identities on types, and thus, by well-typedness, to morphism containing the invertible generator $\ig$. However, as we will now see, for general morphisms this does not always lead to an ``exact" normal framing.

\item \textit{Framing for morphisms of general invertible generators}:
While in the above case of types and their identities the procedure leads to a normal framing, for general morphisms it only does so \textit{approximately}. This is because manifolds corresponding to an invertible generator $\ig$ (and its singularities $\icg\ig_{S\equiv T}$) can not always be confined to live in a plane for which the added constant vectors are normal to---the above picture falsely suggest this could be the case. However, it is still admissible to form normal vectors in the  above way, as we can approximate a normal framing arbitrarily close by appropriately ``scaling" any behaviour outside the plane. 

We illustrate this as follows. Recall our example $\sC$ of a \free{} associative $2$-category from \autoref{ch:presented}. This has an invertible generator $d$ with framing
\begin{restoretext}
\begingroup\sbox0{\includegraphics{ANCimg/page1.png}}\includegraphics[clip,trim=0 {.2\ht0} 0 {.35\ht0} ,width=\textwidth]{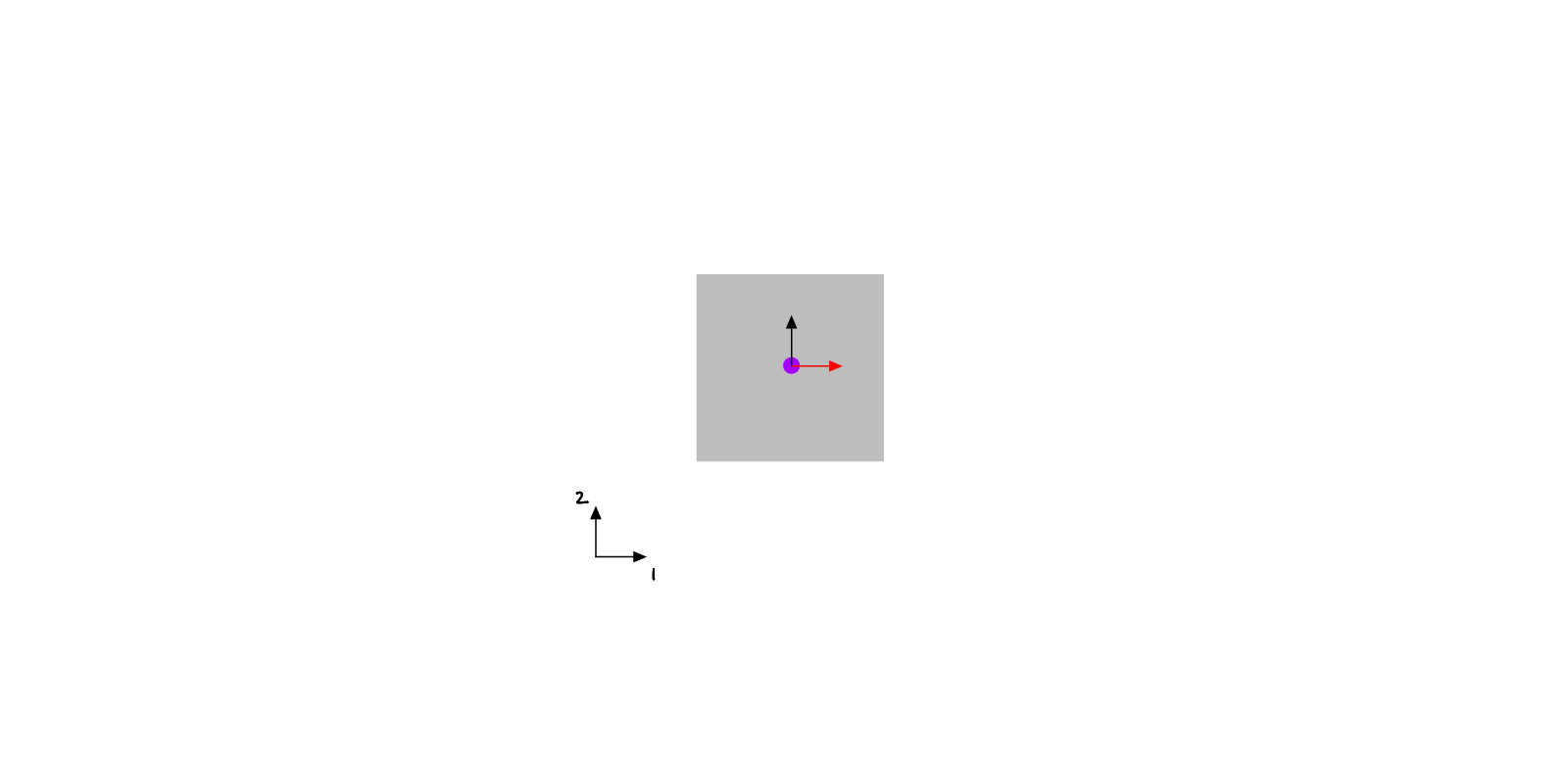}
\endgroup\end{restoretext}
This generator and its cap and cup singularities allow us to form the following morphism on the left below
\begin{restoretext}
\begingroup\sbox0{\includegraphics{ANCimg/page1.png}}\includegraphics[clip,trim=0 {.0\ht0} 0 {.0\ht0} ,width=\textwidth]{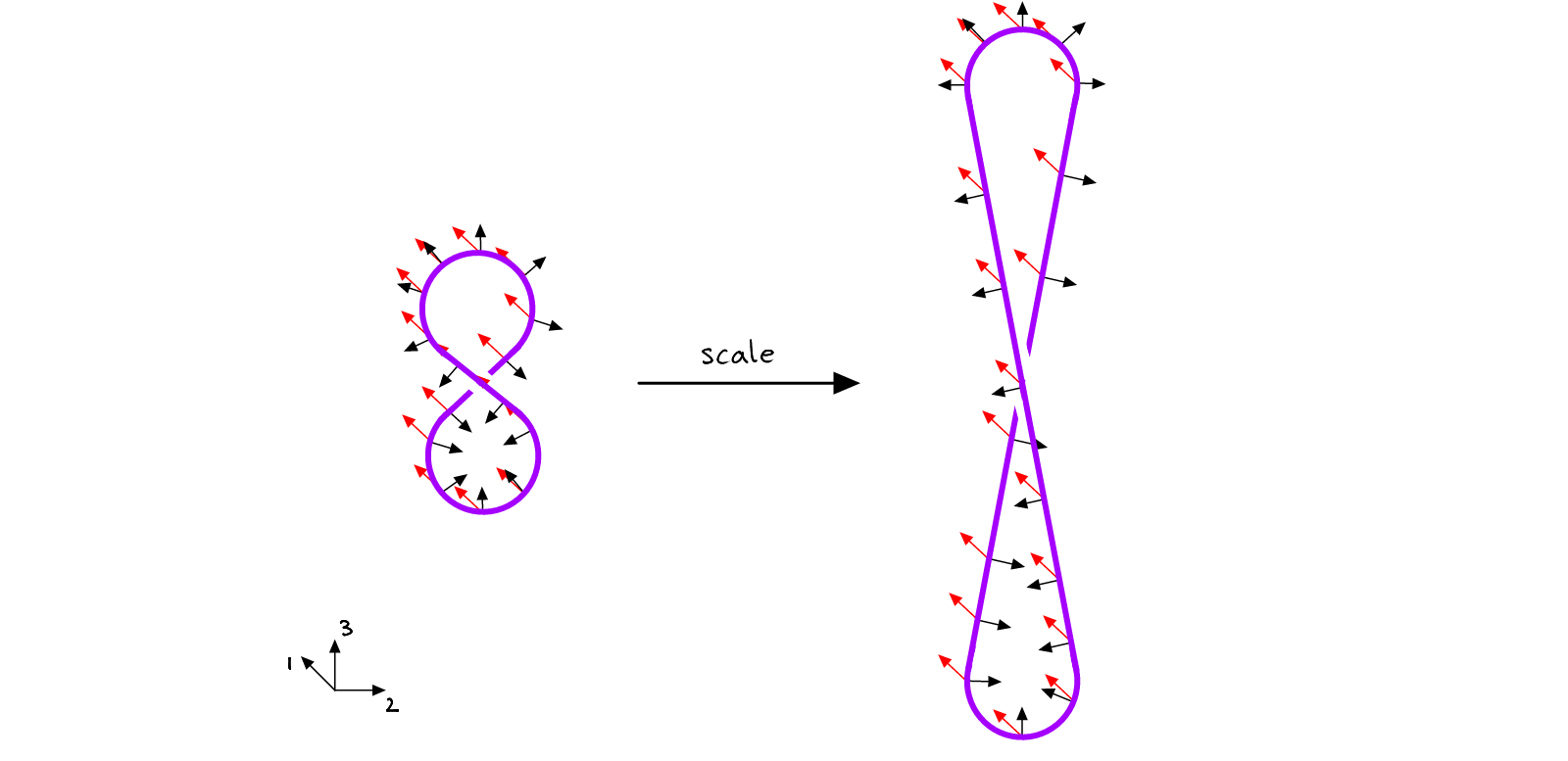}
\endgroup\end{restoretext}
In the above, we omit visually distinguishing the generator $d$ from its derived singularities, such as $d\inv$, its cups and caps etc (they are all part of the \cpurple{} line). This highlights again that $d$ and the singularities derived from its coherent invertibility form a single manifold, for which we want to find a framing. As before \cblack{} arrows above can be derived from normal vectors in $\TI$ after projecting out additional dimensions (in this case, this means projecting along coordinate 1). Red arrows are constant vectors into the directions of this projection (that is, of coordinate 1). Note that, due to its shape, namely due to the presence of the interchanger, the above morphisms cannot possibly be put into a single plane normal to coordinate 1. As a consequence the \cred{} arrows aren't truly normal. However, they are approximately so: by scaling the above morphism appropriately in the  directions normal to coordinate $1$, in parts where it deviates from lying in a plane, we can make the \cred{} arrows arbitrary close to being normal. This is depicted on the right in the above picture. 
\end{enumerate}

This motivates this section's claim that $n$-manifolds associated to an invertible $(m+1)$-generator $\ig$ (and its associated singularities $\icg\ig_{S \equiv T}$)  have a deformation equivalent representative that admits $(n+m+1)$-framing. Since every $k$-cell $g \in \sX^{k}$ in a higher groupoid $\sX^k$ is invertible,  we consequently make the central conjecture that any manifold diagram of a $k$-morphisms in $\sX$ has a deformation equivalent representative that admits framing on all its submanifolds (where again we treat generators together with their singularities as a single manifold as before). We call such a manifold diagrams \textit{framed}. In particular, up to identifying the $k$-cube and the $k$-disk, any framed manifold diagram is in fact a framed stratification in $\Psi^{\mathrm{fr}}_{\sX}(D^k)$ (cf. \autoref{defn:framed_strat}). 

\subsection{Translating CW-complexes to $\infty$-groupoids}

The translation from presentations of CW-complexes $X$ to \free{} associative $n$-groupoids $\sX$ is now straight-forward. We summarise the steps that it involves and then proceed to give examples.

The construction of $\sX$ is inductive in $k$. Assume $\sX^m$, $m < k$ has been defined. 
\begin{enumerate}
\item For each $h \in \sN^X_k$ find the framed $k$-stratification $\bM(h)$ of $D^k$ via the previous construction.

\item Choose a globular foliation of the $k$-disk (cf. \autoref{ssec:sum_CW}), such that the pullback of strata (which might result in singularities, cf. \autoref{ssec:sum_CW}) yields a ``conical" framed manifold diagram $\abss{h\dualdag}$, that is, the framed stratification of the $k$-cube is the cone of the framed stratification of its boundary. As we will see below, the resulting $\abss{h\dualdag}$ is not unique.

\item Add the generator $h\dualdag$ to $\sX^k$. Choose its type to be the labelled singular $k$-cube whose corresponding framed manifold diagram is equivalent to $\abss{h\dualdag}$.
\end{enumerate}

We exemplify the procedure for the three CW-complexes that we previously considered

\begin{egs}[Building a manifold diagram from the generalised Thom-Pontryagin construction] 
First consider $X$ to be the half-cut sphere. Then according to the above procedure $\sX^0$ contains a single generating $0$-morphisms $p\dualdag$
\begin{restoretext}
\begingroup\sbox0{\includegraphics{ANCimg/page1.png}}\includegraphics[clip,trim=0 {.6\ht0} 0 {.25\ht0} ,width=\textwidth]{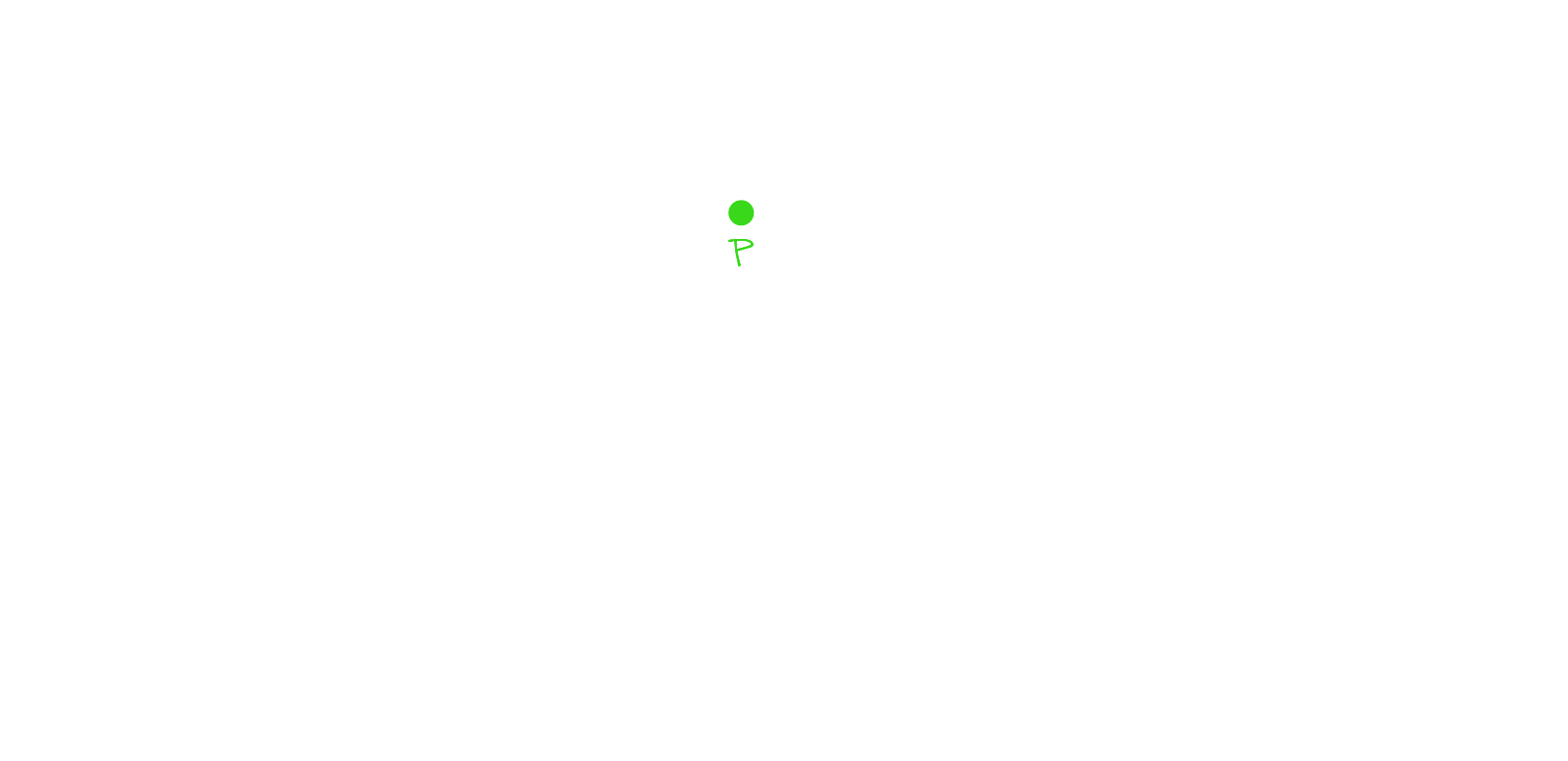}
\endgroup\end{restoretext}
Next, $g \in \sN^X_1$ translates to an invertible $1$-generator $g\dualdag \in \sX^1$ with type
\begin{restoretext}
\begingroup\sbox0{\includegraphics{ANCimg/page1.png}}\includegraphics[clip,trim=0 {.35\ht0} 0 {.4\ht0} ,width=\textwidth]{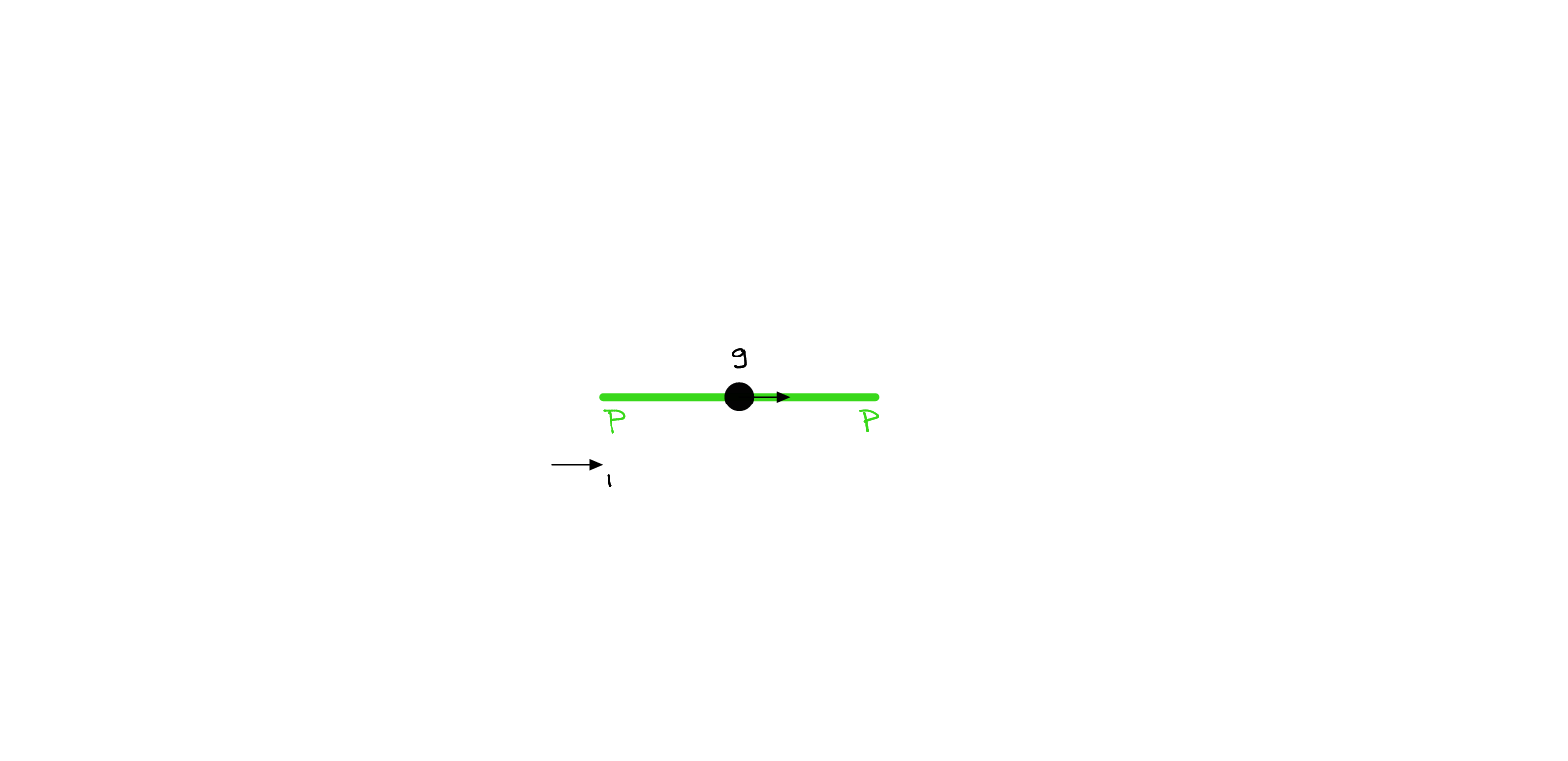}
\endgroup\end{restoretext}
Note how the fact that $\ata _g$ maps the endpoints of $D^1$ to $p$ translates into $g\dualdag \in \sX^1$ having source and target $p\dualdag \in \sX^0$. Finally, the two cells $d,e \in \sN^X_2$ translate to invertible 2-generators $d\dualdag,e\dualdag \in \sX^2$ with types
\begin{restoretext}
\begingroup\sbox0{\includegraphics{ANCimg/page1.png}}\includegraphics[clip,trim=0 {.2\ht0} 0 {.25\ht0} ,width=\textwidth]{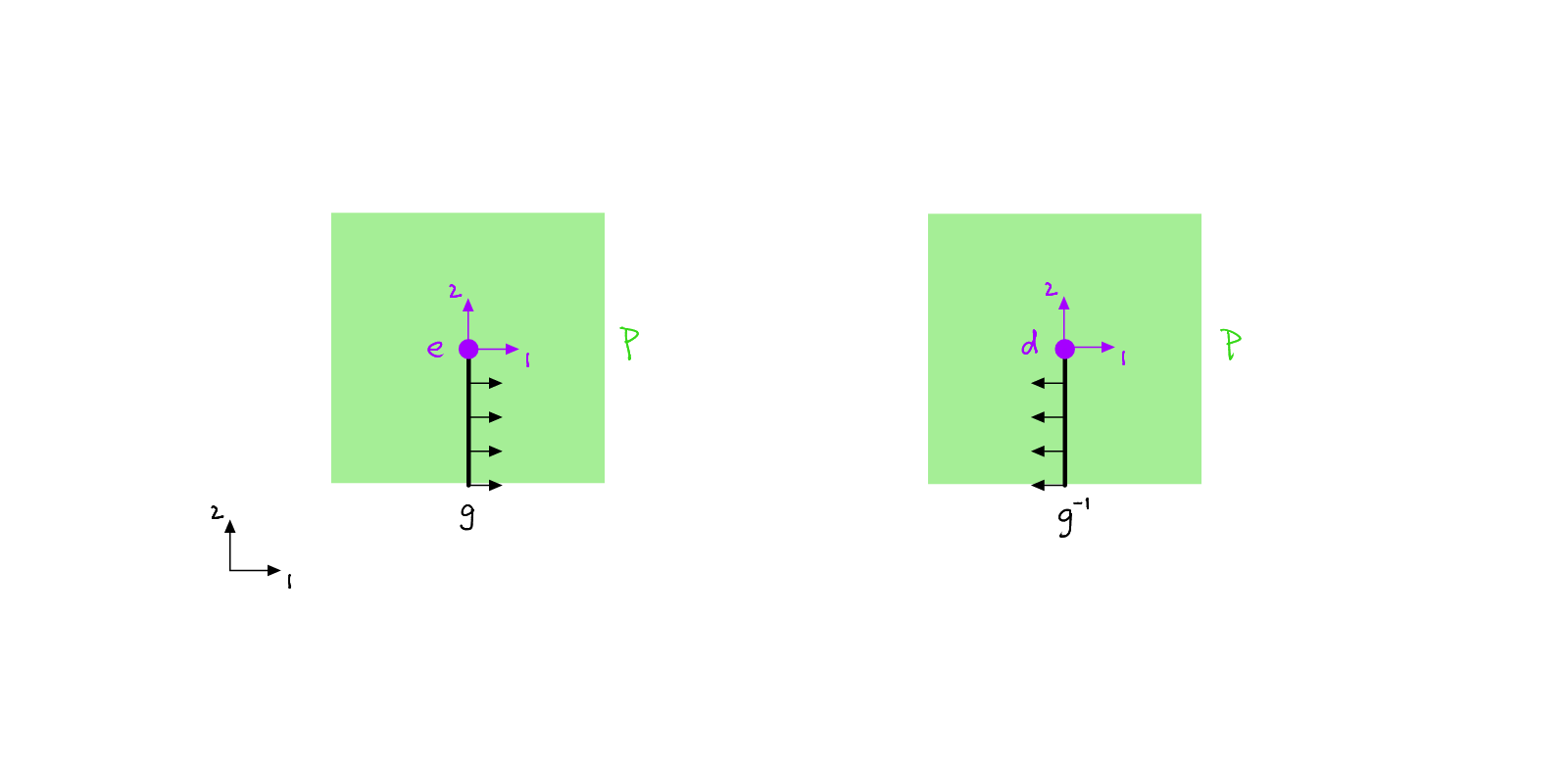}
\endgroup\end{restoretext}
Note, how in both cases the framed manifold diagrams $\abss{e\dualdag}$ (and $\abss{d\dualdag}$) equal to the pullback stratification $\bM(e)$ (and $\bM(d)$) up to an identification of the $2$-cube with the $2$-sphere (by a globular foliation, cf. \autoref{ssec:sum_CW}). A different choice of this identification however can lead to a different type $\abss{e\dualdag}$ as depicted on the left below
\begin{restoretext}
\begingroup\sbox0{\includegraphics{ANCimg/page1.png}}\includegraphics[clip,trim=0 {.2\ht0} 0 {.25\ht0} ,width=\textwidth]{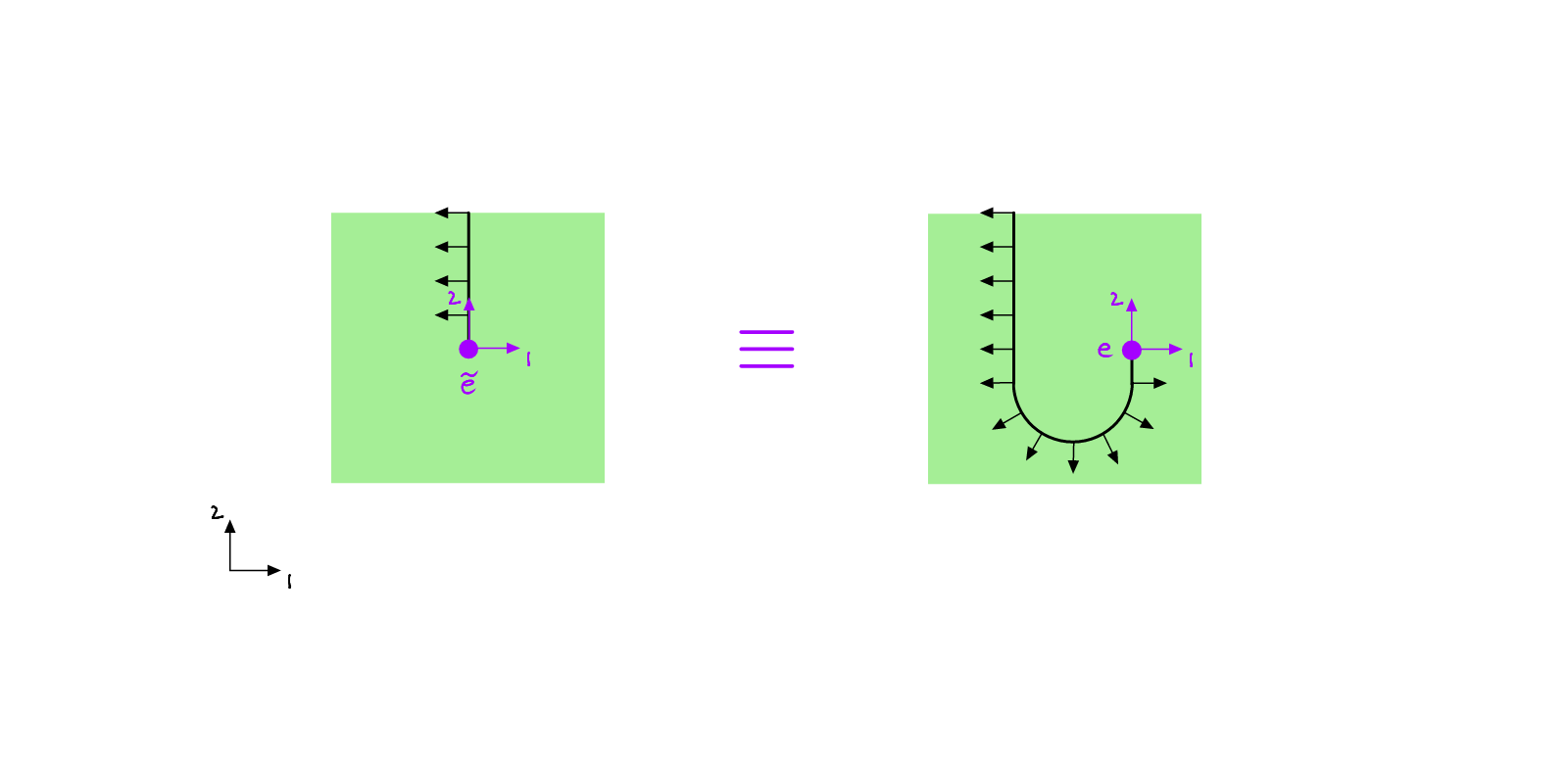}
\endgroup\end{restoretext}
Roughly speaking, this choice would result in an ``equivalent" groupoid because the behaviour of this new type $\abss{e\dualdag}$ can be mimicked by the old generator together with the singularities derived from invertibility of $g\dualdag$. This is illustrated on the right above. 

As a second example, we consider the two-point torus $X$. The two $0$-cells $p,q\in X$ translate in to $0$-generators $p\dualdag,q\dualdag \in \sX^0$
\begin{restoretext}
\begingroup\sbox0{\includegraphics{ANCimg/page1.png}}\includegraphics[clip,trim=0 {.4\ht0} 0 {.2\ht0} ,width=\textwidth]{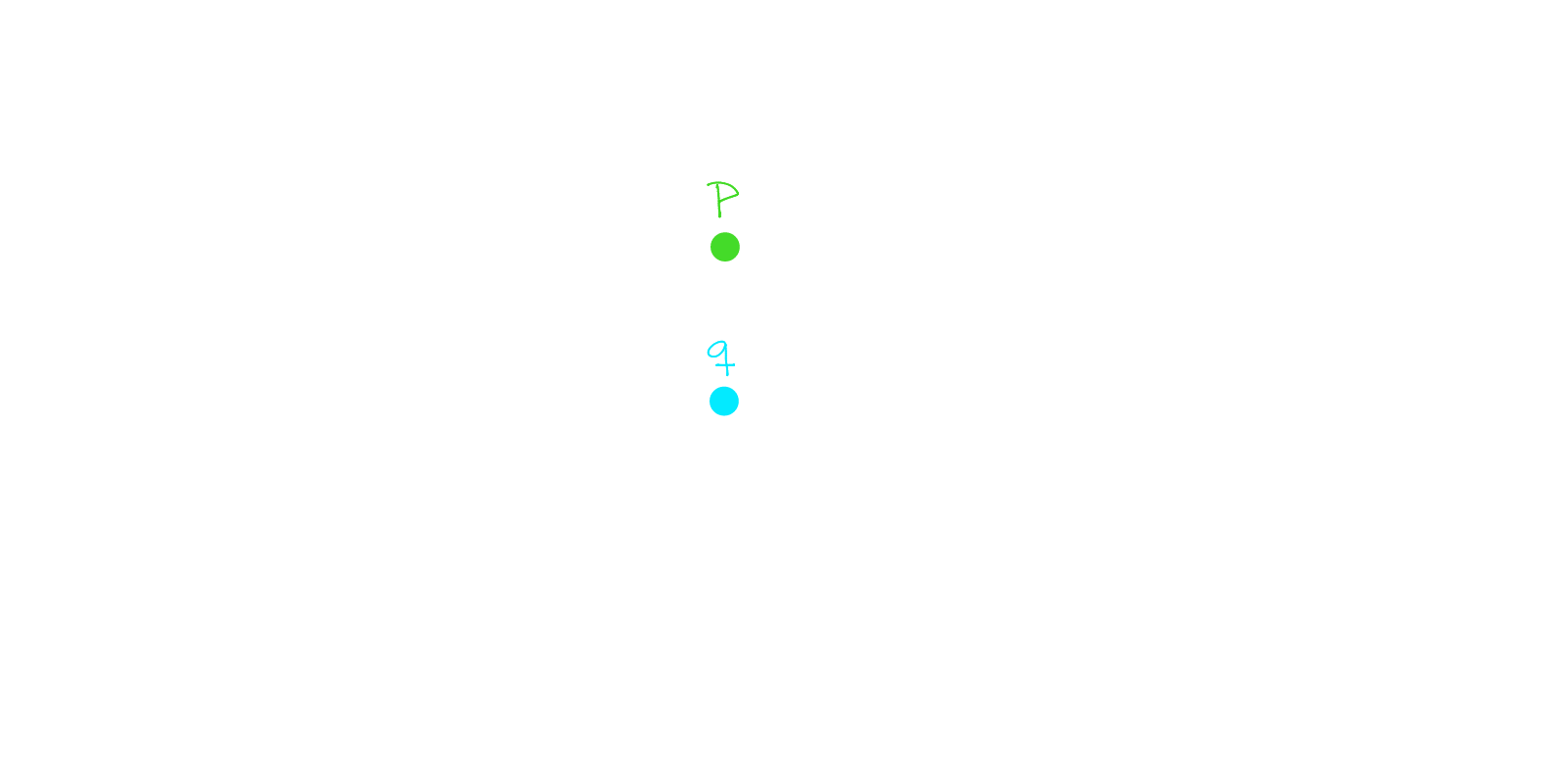}
\endgroup\end{restoretext}
The three $1$-cells $g,h,l \in \sN^X_1$i translate into invertible $1$-generators $g\dualdag,h\dualdag,l\dualdag \in \sX^1$ with types
\begin{restoretext}
\begingroup\sbox0{\includegraphics{ANCimg/page1.png}}\includegraphics[clip,trim=0 {.05\ht0} 0 {.35\ht0} ,width=\textwidth]{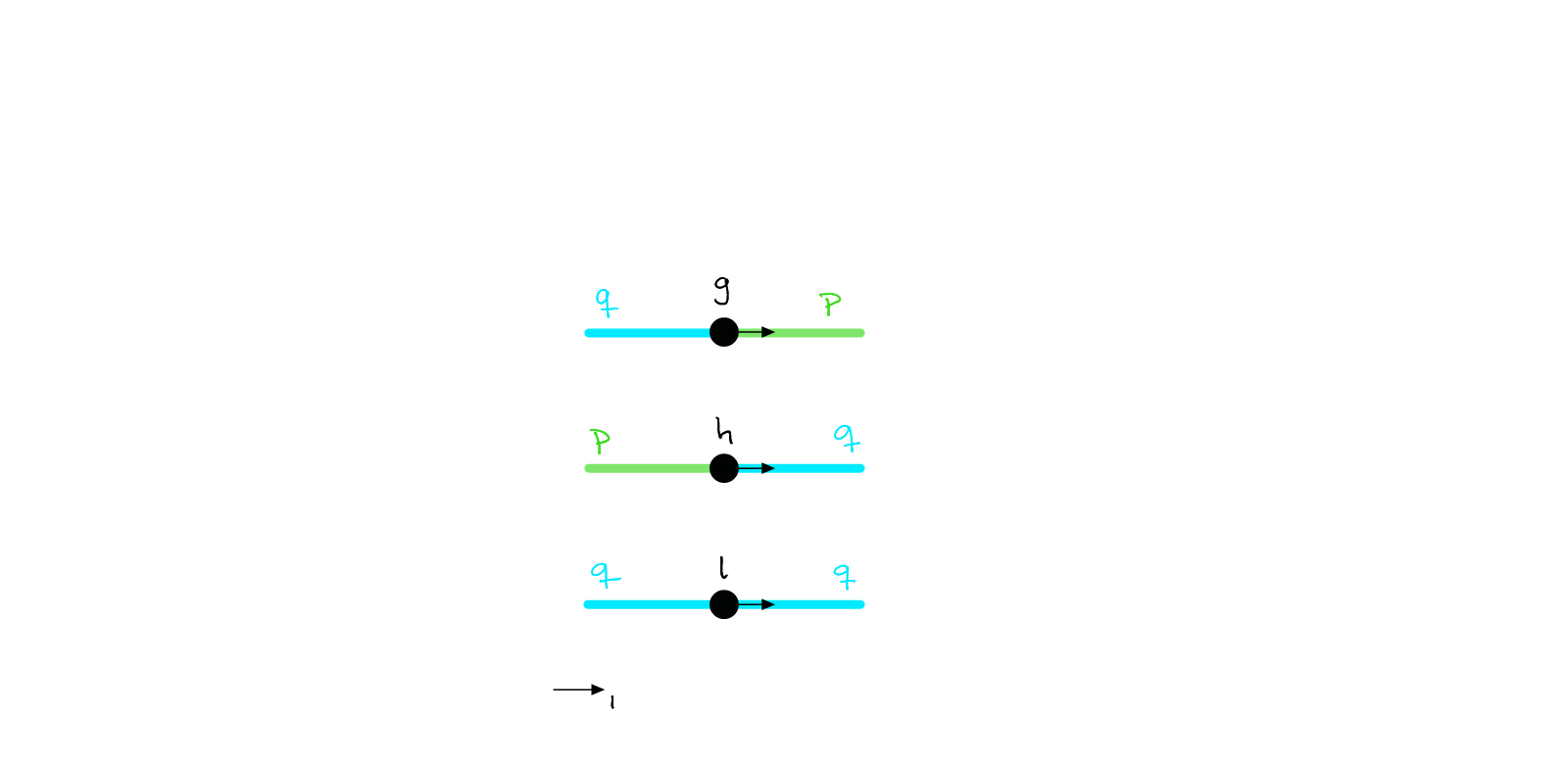}
\endgroup\end{restoretext}
Finally, the single $2$-cell $d \in \sN^X_2$ translates to an invertible $2$-generator $d\dualdag \in \sX^2$ which could be given one of the following types
\begin{restoretext}
\begingroup\sbox0{\includegraphics{ANCimg/page1.png}}\includegraphics[clip,trim=0 {.2\ht0} 0 {.2\ht0} ,width=\textwidth]{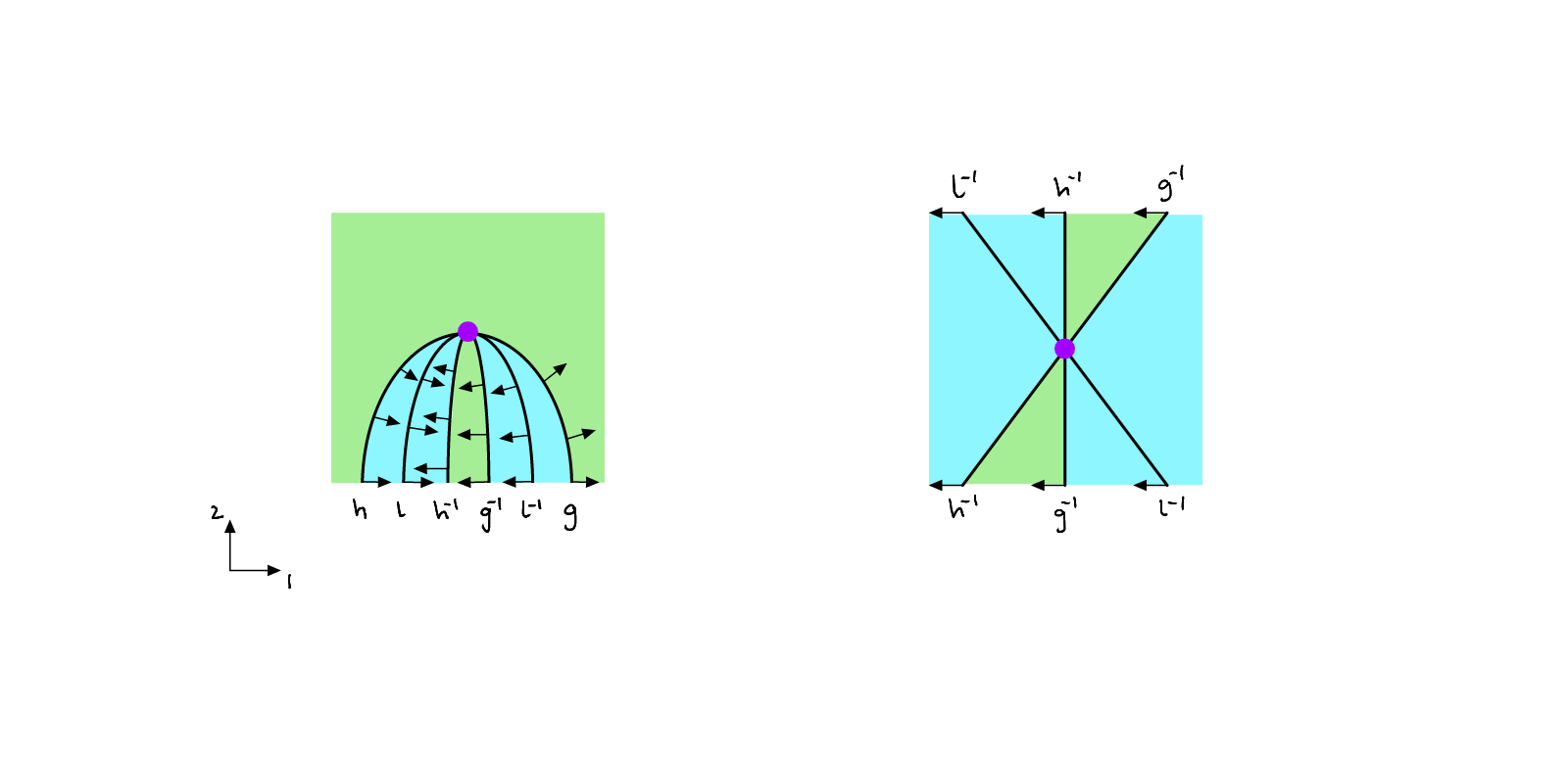}
\endgroup\end{restoretext}
Both versions can again be seen to be ``equivalent": For instance the right can be obtained from the left by bending the rightmost and two leftmost wires upwards.

Finally, we consider the complex projective plane $X$, which has a single $0$-generator
\begin{restoretext}
\begingroup\sbox0{\includegraphics{ANCimg/page1.png}}\includegraphics[clip,trim=0 {.5\ht0} 0 {.3\ht0} ,width=\textwidth]{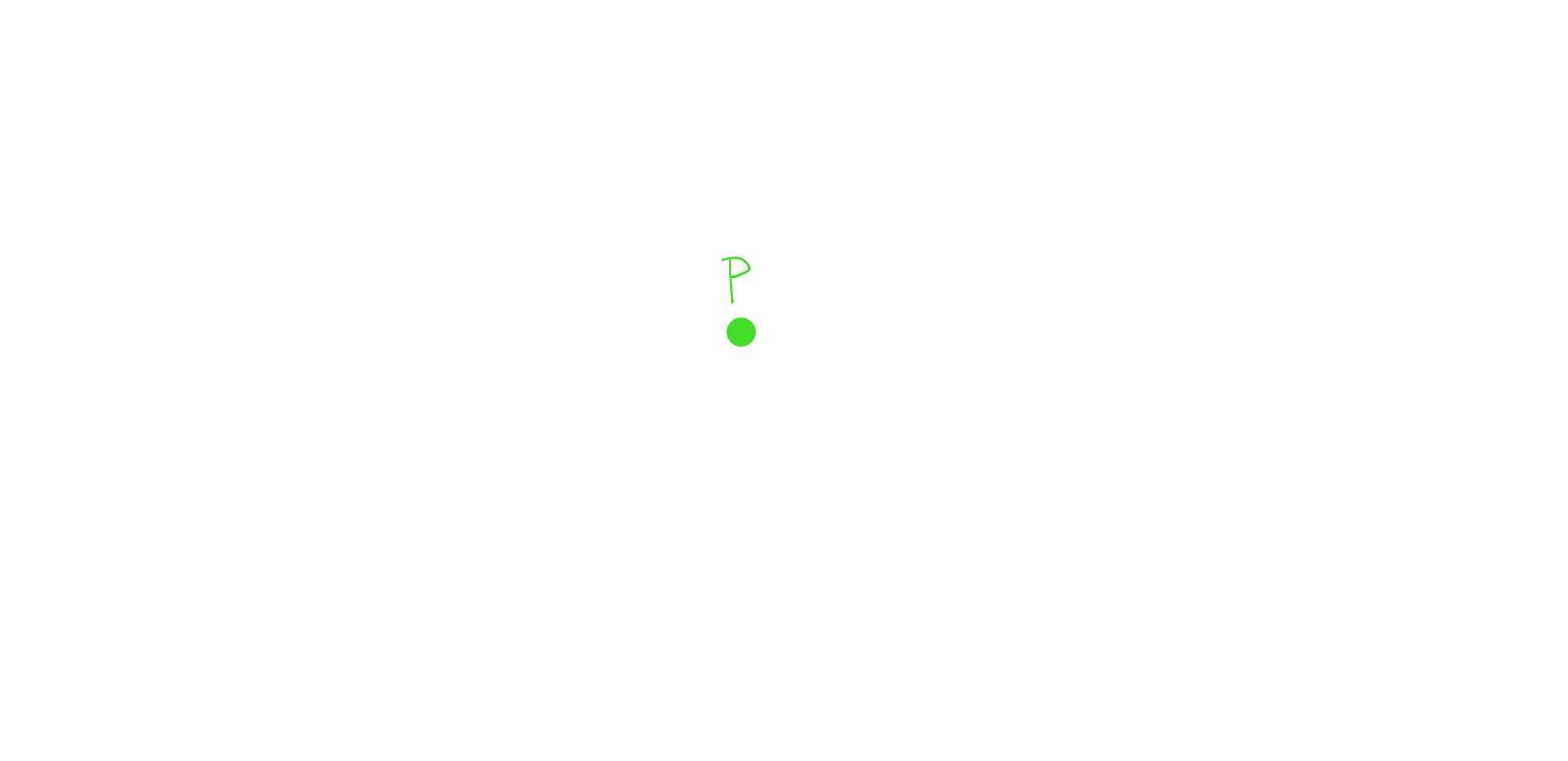}
\endgroup\end{restoretext}
in $\sX^0$ corresponding to $p \in \sN^X_0$. Next, $d\dualdag \in \sN^X_2$ corresponds to $d \in \sX^2$ with type
\begin{restoretext}
\begingroup\sbox0{\includegraphics{ANCimg/page1.png}}\includegraphics[clip,trim=0 {.2\ht0} 0 {.25\ht0} ,width=\textwidth]{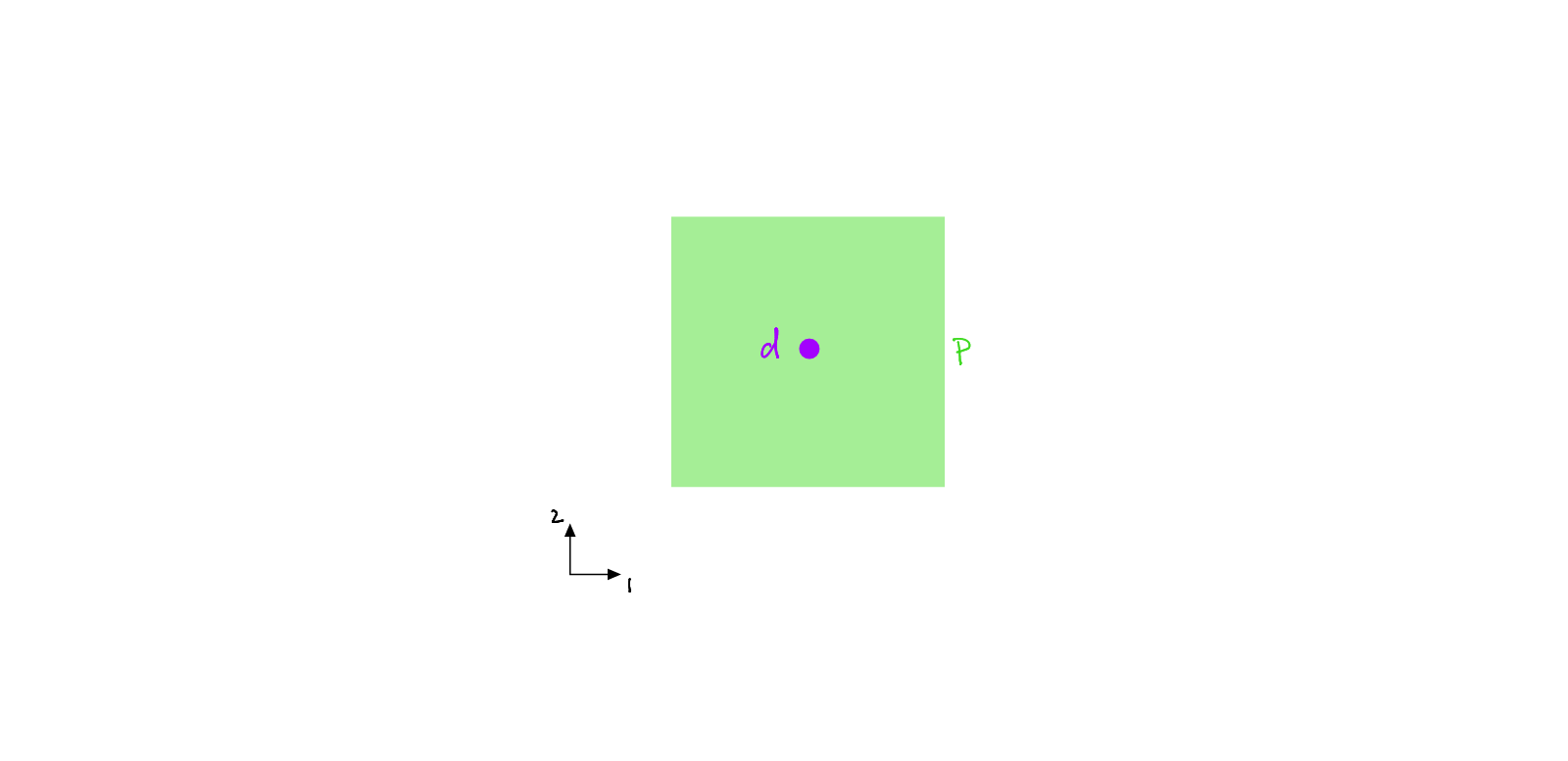}
\endgroup\end{restoretext}
Last but not least, $h \in \sN^X_4$ corresponds to $h\dualdag \in \sX^4$ with type $\abss{h\dualdag}$ given by
\begin{restoretext}
\begingroup\sbox0{\includegraphics{ANCimg/page1.png}}\includegraphics[clip,trim=0 {.05\ht0} 0 {.15\ht0} ,width=\textwidth]{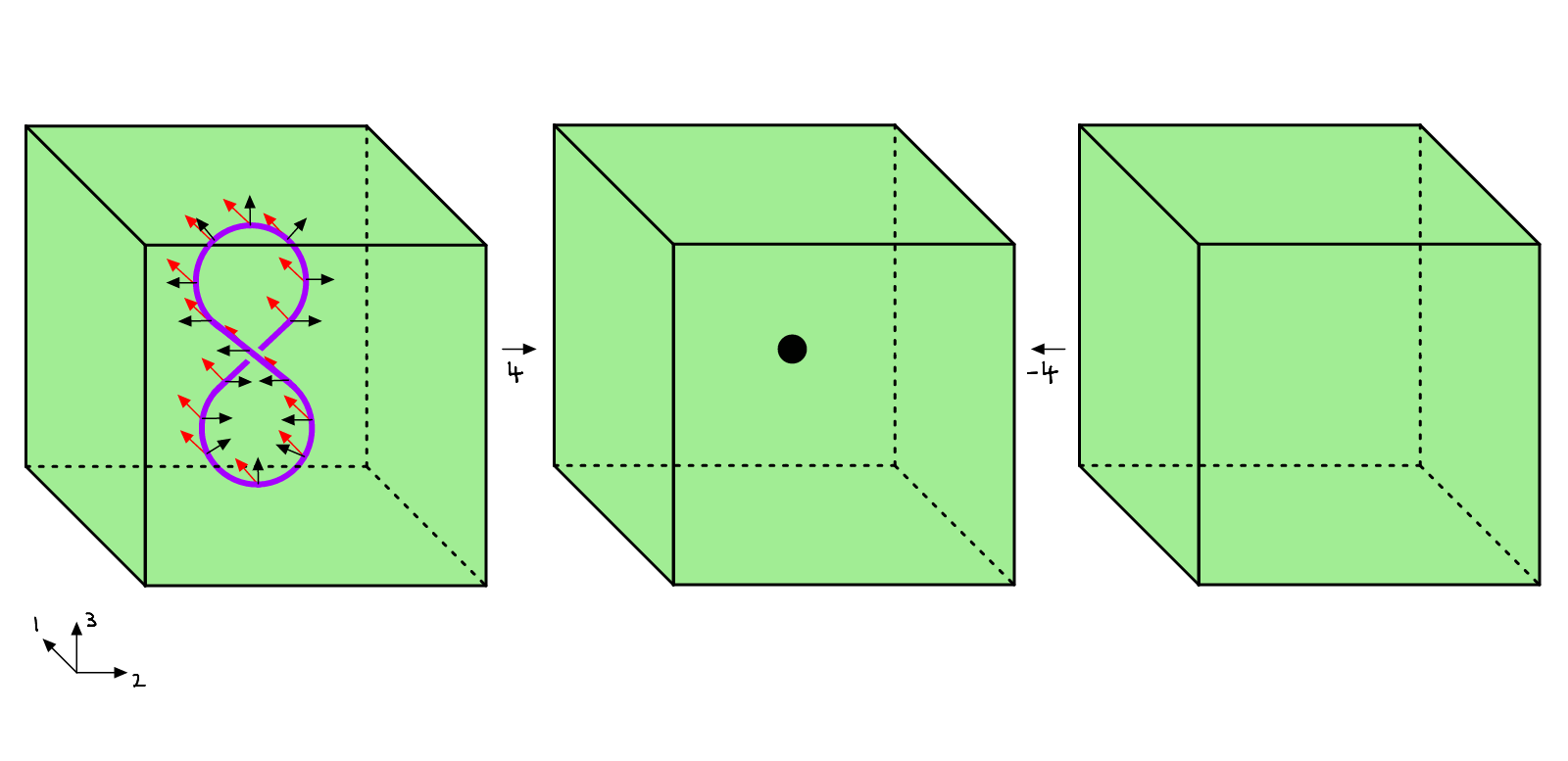}
\endgroup\end{restoretext}
Note that the source of $\abss{h\dualdag}$ contains a circle manifold whose normal framing can be obtained by previously discussed arguments. This framing rotates by $2\pi$ when traversing the circle, and equals (that is, is framed cobordant) to the output of the Thom-Pontryagin construction for the Hopf map (cf.  \autoref{eg:hopf_map}). Up to identification of the $4$-cube and the $4$-disk by an appropriate globular foliation, we see that the above $4$-cube equals $\bM(h)$ as previously constructed.
\end{egs}

\subsection{Connections to the homotopy hypothesis} \label{ssec:homotopy_hyp}

Based on the previous informal discussion, a formal procedure for translating CW-complexes into \free{} associative higher groupoids can be obtained by filling in the bottom and sides of the following square:
\begin{equation} 
\xymatrix{ \text{$\infty$-Groupoids} \quad \ar[d]_{\text{types}} \ar@{<->}[r] & \quad \text{CW-complexes} \ar[d]^{\text{cells}} \\
\text{Manifold diagrams}\quad \ar@{<->}[r] & \quad \text{Framed stratifications} 
}
\end{equation}
Here, the right-hand side arrow corresponds to the ``generalised unbased Thom-Pontryagin construction" discussed in the previous section. The left-hand side arrow notes that a (\free{} associative) $\infty$-groupoid is determined by its generating types. The lower arrow was loosely illustrated by examples in the previous section. We manifest the upper arrow in the following conjecture.
\begin{conj}[Combinatorial representation of stratifications] \label{conj:cobordism_hyp} Let $X$ be a CW-complex. There is a (non-unique) $\infty$-groupoid $\sX$ such that $\sX_k = I^X_k$, and for each $h \in I^X_k$, $\bM(h)$ is equivalent to $\abss{h}$. Conversely, if $\sX$ arises from $X$ in this way, then $X$ as a CW-complex can be uniquely recovered from $\sX$ by a mapping from groupoids to CW-complexes based on (the inverse) of the generalised Thom-Pontryagin construction.
\end{conj}
\noindent Natural tools for attempting to prove such a conjecture (which involves in particular a formalisation of the lower arrow in the above square), are Morse theory and Cerf theory. 

This is a fine-grained version of the homotopy hypothesis, as it gives a (many-to-one) correspondence at the level of presentations of groupoids and spaces. Note that the ``classical homotopy hypothesis" \cite{baez2007homotopy} asks for a correspondence of spaces and higher groupoids up to equivalence. We have neither given a definition of functors nor of equivalences of higher categories at this stage. We expect that the correspondence at the level of presentations will be a good indicator of how to formulate such notions correctly---that is, in a way that the ``classical homotopy hypothesis" holds. The reader might have noticed that the line of thought outlined in this chapter also has close ties to another important hypothesis: the so-called generalised tangle hypothesis \cite{baez1995higher}. Detailing the connections to both hypotheses is left to future work.

\setcounter{section}{0}
\renewcommand{\thesection}{B.{\arabic{section}}}
\renewcommand{\theHsection}{B.{\arabic{section}}}

\chapter{\Free{} associative $n$-fold categories} \label{ch:nfold}

We define \free{} associative $n$-fold categories. The definition of \free{} associative $n$-fold categories is a proper generalisation of the definition of \free{} associative $n$-categories which we have studied in-depth in \autoref{ch:presented}. Instead of having sets of generators $\sC_i$ indexed by $i \in (\bnum{n+1})$ we will have sets of generators $\sC_\alpha$ indexed by $\alpha \in \bnum{2}^{n+1}$. The index $\alpha$ should be understood according to the \stratatype{}  terminology which was defined (topologically) in \autoref{ssec:po_mfld_diag} and (combinatorially) in \autoref{ssec:regions} (note that the two definitions differ by an inversion of order of components).

\section{Definition}

\subsection{Presentations}

\begin{defn} A \textit{\free{} associative $n$-fold category} $\sC$, where $n \in \Set{-1}\cup \lN \cup \Set{\infty}$, is an ``$n$-cube" of sets
\begin{equation}
\Set{ \sC_\alpha ~|~ \alpha = (\alpha_0, \alpha_1, ... ,\alpha_n, \alpha_{n+1}) \in \bnum{2}^{n+1} }
\end{equation}
(where $\sC_\alpha$ is called the set of generating $\alpha$-morphisms if $\alpha_{n+1} = 0$, and set of generating $\alpha$-equalities if $\alpha_{n+1} = 1$) together with the following structure: let
\begin{equation}
\Gamma_\sC : \bnum{2}^{n+1} \to \PRel
\end{equation}
map $\alpha \mapsto \Discr\sC_\alpha$ and $(\alpha \to \beta)$ maps to $R^{\mathrm{full}}_{\Discr\sC_\alpha, \Discr\sC_\beta}$ (cf. \autoref{rmk:full_pro_rel}). Then, for each $g \in \sC_\alpha$, $\sC$ comes with data of a normalised cube.
\begin{equation}
\abss{g} : \bnum{1} \to \SIvert n {\GGamma{}\sC}
\end{equation}
which satisfies
\begin{itemize}
\item \textit{Minimality}: There is $\ip_g \in \tsG n(\abss{g})$ with $\tsU n_{\abss{g}}(\ip_g) = g$ and
\begin{equation}
\abss{g} = \abss{g} \sslash \ip_g
\end{equation}
\item \textit{Dimensionality}: The region $\ip_g$ is of type $\alpha$, that is 
\begin{equation}
\ctyp^n_{\abss{g}} (\ip_g) = \alpha
\end{equation}
\item \textit{Well-typedness}: For all $p \in \tsG n(\abss{g})$ with $\tsU n_{\abss{g}}(p) = f$,
\begin{equation}
\NF{~\abss{g} \sslash p ~}^n = \abss{f}
\end{equation}
\end{itemize}
The first two conditions can be summarised by saying that the type is ``conical" as before.

Note that for $n < \infty$ and $\alpha = (\beta, 1)$ (that is, $g$ is a generating $\alpha$-equality), then we could also add proof-irrelevance and symmetry conditions similar to those in \autoref{defn:pres_ANC}, which would guarantee for instance the existence of a unique ``inverse" $g\inv$. We leave stating these conditions to the diligent reader, as the work needed to formulate them mostly doesn't provide valuable new insights.
\end{defn}

\subsection{Morphisms}

The previous definition only defines presentations of (presented associative) $n$-fold categories. The next definition answers the question what type of morphisms can be naturally build from such presentations.

\begin{defn}[Morphisms in an $n$-fold category] Given a \free{} associative $n$-fold category $\sC$, the \textit{set $\Comp(\sC)_\alpha$ of $\alpha$-morphisms}, $\alpha \in \bnum 2^n$, is defined to contain normalised $n$-cubes $\scD : \mathbf{1} : \SIvert n \scC$ satisfying the following

\begin{enumerate}
\item \textit{Correctness of dimension}: If $\alpha_i = 0$ for $1 \leq i \leq n$ then
\begin{equation}
\tusU i_\scD = \const_{\singint 0}
\end{equation}
\item \textit{Well-typedness}: If $p \in \tsG k_\scD$ with $\tsU n_{\abss{g}}(p) = f$,
\begin{equation}
\NF{~\abss{g} \sslash p ~}^n = \abss{f}
\end{equation}
\end{enumerate}

\end{defn}

\section{The graphical calculus of double and triple categories}

Morphism can be represented by a graphical calculus. For this the discussion in \autoref{ssec:coloring} fully applies, that is, $\alpha$-morphisms of $n$-fold cubes can be represented just as any other colored $n$-cube.

The resulting graphical calculus of \free{} associative $n$-fold categories is similar but more powerful than the graphical calculus of \free{} associative $n$-categories. We give the following examples. In dimension 2, that is, in the case of ``double categories" we are now allowed to have corners in strings
\begin{restoretext}
\begingroup\sbox0{\includegraphics{ANCimg/page1.png}}\includegraphics[clip,trim=0 {.2\ht0} 0 {.2\ht0} ,width=\textwidth]{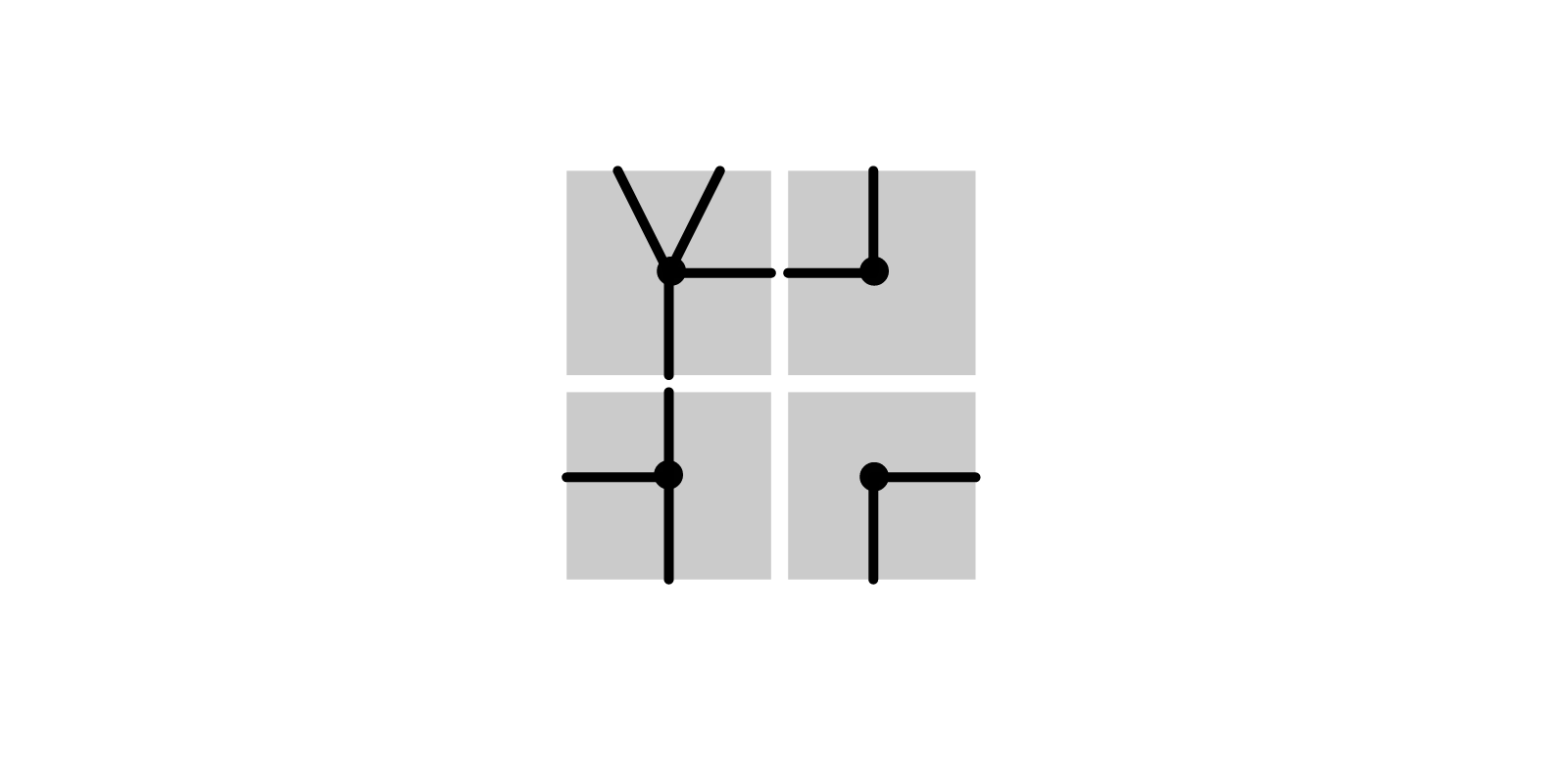}
\endgroup\end{restoretext}
In three dimensions we can have corners in surface as illustrated in the following
\begin{restoretext}
\begingroup\sbox0{\includegraphics{ANCimg/page1.png}}\includegraphics[clip,trim=0 {.1\ht0} 0 {.1\ht0} ,width=\textwidth]{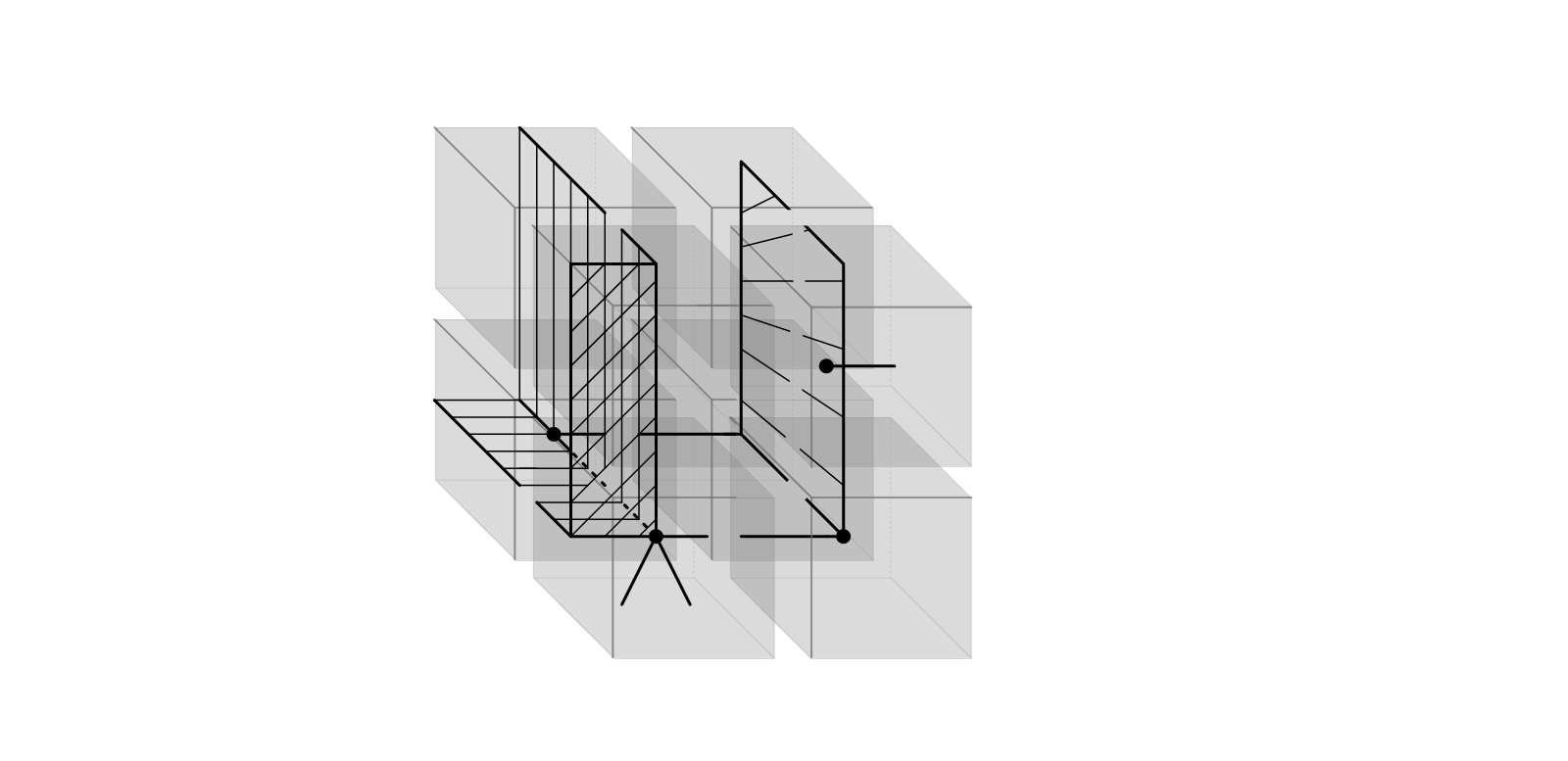}
\endgroup\end{restoretext}
The $2$-dimensional case is a formalisation of a graphical calculus that has already been used in the context of proarrow-equipments in double categories \cite{myers2016string}. The $3$-dimensional case, and the $n$-dimensional case are novel graphical calculi.

\setcounter{section}{0}
\renewcommand{\thesection}{C.{\arabic{section}}}
\renewcommand{\theHsection}{C.{\arabic{section}}}

\chapter{\Free{} fully weak $n$-categories} \label{ch:weak}

In this appendix we return to the question ``where did all the other coherence data go?". We discuss how coherence data (in the usual sense of a fully weak $n$-category) can possibly be recovered by combining homotopies with an explicit (invertible) composition operation. The idea is similar to the notion of resolution from \autoref{constr:resolutions}.

\section{Definition}

\subsection{Resolutions of associative presentations}

\begin{constr}[Resolutions of presentations] Let $\sC \in \pCat_\infty$. Recall the inclusion $\abss{-} :  \sC \in \Comp(\sC)$ of sets. Define
\begin{equation}
\kiC{0}(\sC) = \sC \gadd{\gsrc(f), \gtgt(f)} \Set{ \bracket(f) ~|~f \in \Comp(\sC) \setminus \sC }
\end{equation}
$\kiC{0}(\sC)$ contains each composite of $\sC$ as a generator, which (for non-trivial composites $f\in \Comp(\sC) \setminus \sC$) is called the \textit{bracketed version} (or the \textit{bracketing}) of $f$. We further define
\begin{equation}
\kiC{}(\sC) = \kiC{0}(\sC) \igadd{f, \bracket(f)} \Set{ \witness{f} ~|~f \in \Comp(\sC) \setminus \sC }
\end{equation}
which adds invertible \textit{bracketing witnesses} between each composite $f$ and its bracketed version $\bracket(f)$. $\kiC{}(\sC)$ is called the \textit{resolution} of $\sC$. Note that there is a canonical inclusion of presentations
\begin{equation}
\sC \into \kiC{}(\sC)
\end{equation}
\end{constr}

\subsection{Weak presentations}

\begin{defn}[\Free{} weak $n$-categories] A \textit{\free{} weak $n$-category} $\sC$ is the colimit (cf. \autoref{rmk:colimit_of_presentations}) of a sequence of associative $n$-categories
\begin{equation} \label{eq:weak_cat_colimit}
\sC^0 \xinto {\bunit^1} \sC^1 \xinto {\bunit^1} \sC^2 \xinto {\bunit^2} \dots
\end{equation}
satisfying
\begin{equation}
\sC^{k} = \kiC{}(\sC^{k-1}) \gadd{s(g_k), t(g_k)} \Set{g_k \in G_k}
\end{equation}
Here, $\sC^k$ is called the \textit{depth $k$ bracketing of $\sC$}, and $G_k$ is called the set of \textit{depth $k$ generators}.
\end{defn}

\begin{rmk}[Interpretation of bracketing and bracketing witnesses] Using usual higher category theory lingo, bracketings should be understood as \textit{candidate compositions} (cf. \cite{groth2015short}). For instance, if $f, g$ are (composable) generating $1$-morphisms in $\sC^{i-1}$ and recalling that $f \whisker 1 1 g$ denotes the $1$-cube in which $f$ and $g$ are composed, then $\bracket(f \whisker 1 1 g)$ (a generator of $\sC^i$) should be understood as the $1$-morphism which is the canonical composition candidate for composing $f$ and $g$. $\witness{\scD}$ is in turn the witness of equivalence between the composition candidate $\bracket(f \whisker 1 1 g)$ and the composite $f \whisker 1 1 g$ of $f$ and $g$. 
\end{rmk}

\subsection{Morphisms}

Let $\sC$ be a \free{} weak $n$-category. Every $\scD \in \Comp(\sC)_k$ lives in $\Comp(\sC^k)_k$ for some  large enough $k$ (up to relabelling by $\bunit^\infty_k : \sC^k \into \sC$, cf. \autoref{rmk:colimit_of_presentations}), but every such $\scD$ is then isomorphic to a generator $\bracket(\scD) \in \sC^{k+1}$. Thus, our ``old" notion of morphism contains redundant information as it records morphisms ``twice" in that way. A more appropriate notion of morphisms and composition is given by the following definition.

\begin{defn}[Morphisms in weak categories] Let $\sC$ be a \free{} weak $n$-category. A \textit{$k$-morphism $f$ in $\sC$} is an element $f \in \sC_k$. Two $k$-morphisms $f$ and $g$ are \textit{equivalent}, if there is a map of presentations
\begin{equation}
(\sC \igadd{f,g}  \ie) \to \sC
\end{equation}
which restricts to the identity on $\sC$.

In place of our previous terminology for morphisms, we introduce the following: $\scD \in \Comp(\sC)_k$ is called a \textit{compositional $k$-shape} in $\sC$. A composition candidate $f$ for this shape is a $k$-morphisms that is equivalent to $\bracket(\scD) \in \sC_k$.
\end{defn}

\subsection{Recovering coherence data}

Coherences data now arises naturally by a combination of bracketing witnesses and homotopies. We will exemplify this in the case of (non-strict) identities, associators, interchangers and pentagonators. 

Let $\sC^0$ contain two composable generating $1$-morphisms $f$ and $g$ (we don't need to specify their sources and targets for the following discussion). Then the depth $1$ bracketing $\sC^1$ contains the $1$-generator $\bracket(\abss{f} \whisker 1 1 \abss{g})$ which we will abbreviate as $(fg)$. It also contains the generating $2$-morphism $\witness{\abss{f} \whisker 1 1 \abss{g}}$ which then has the type
\begin{restoretext}
\begingroup\sbox0{\includegraphics{ANCimg/page1.png}}\includegraphics[clip,trim=0 {.2\ht0} 0 {.2\ht0} ,width=\textwidth]{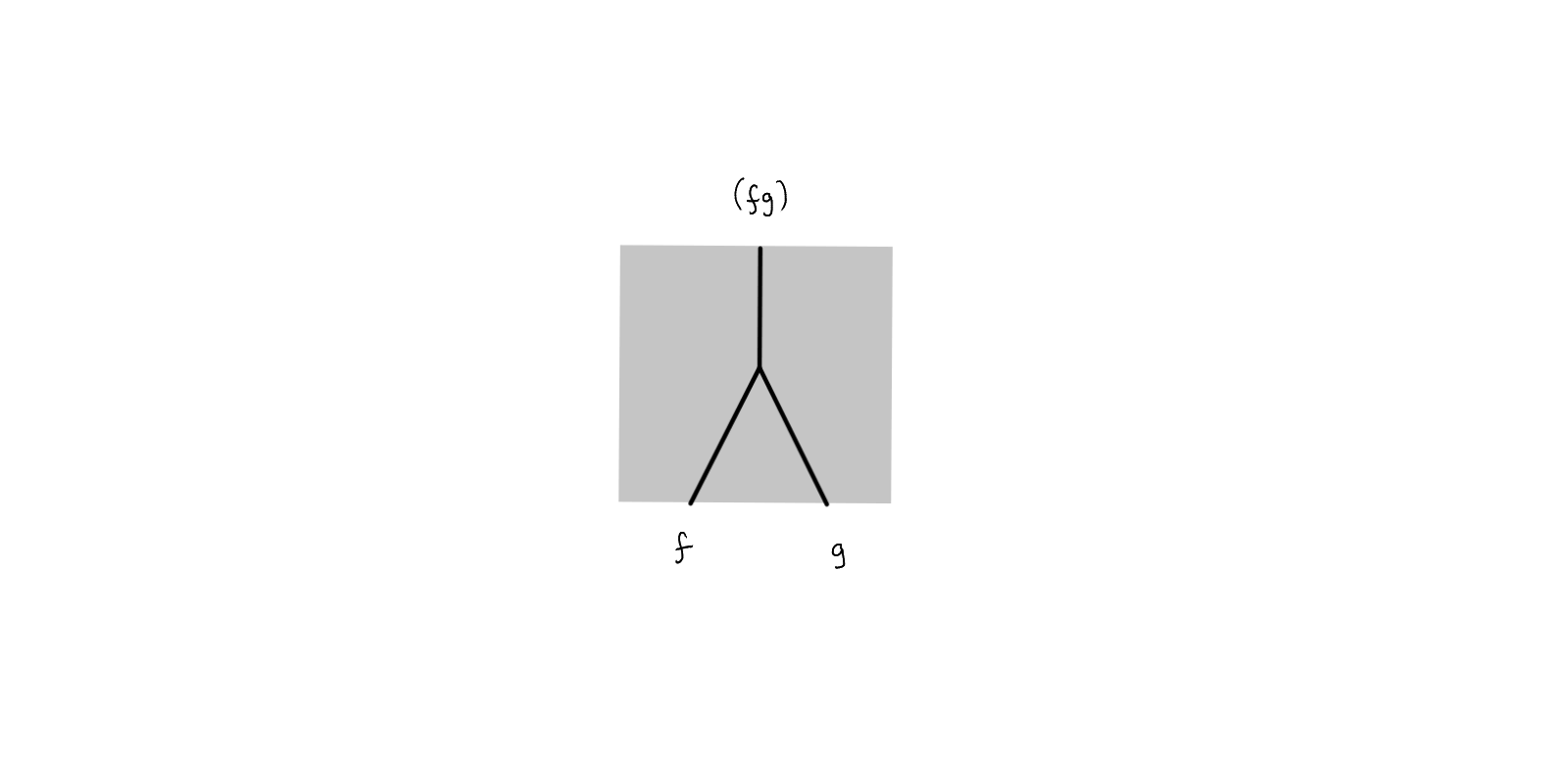}
\endgroup\end{restoretext}
Note that here we depict the central singularity of the bracketing witness simply by a bifurcation of wires (as opposed to putting an additional point on that bifurcation as we often did before). The cup singularities for the invertible generator $\witness{\abss{f} \whisker 1 1 \abss{g}}$ are then of the form
\begin{restoretext}
\begingroup\sbox0{\includegraphics{ANCimg/page1.png}}\includegraphics[clip,trim=0 {.25\ht0} 0 {.25\ht0} ,width=\textwidth]{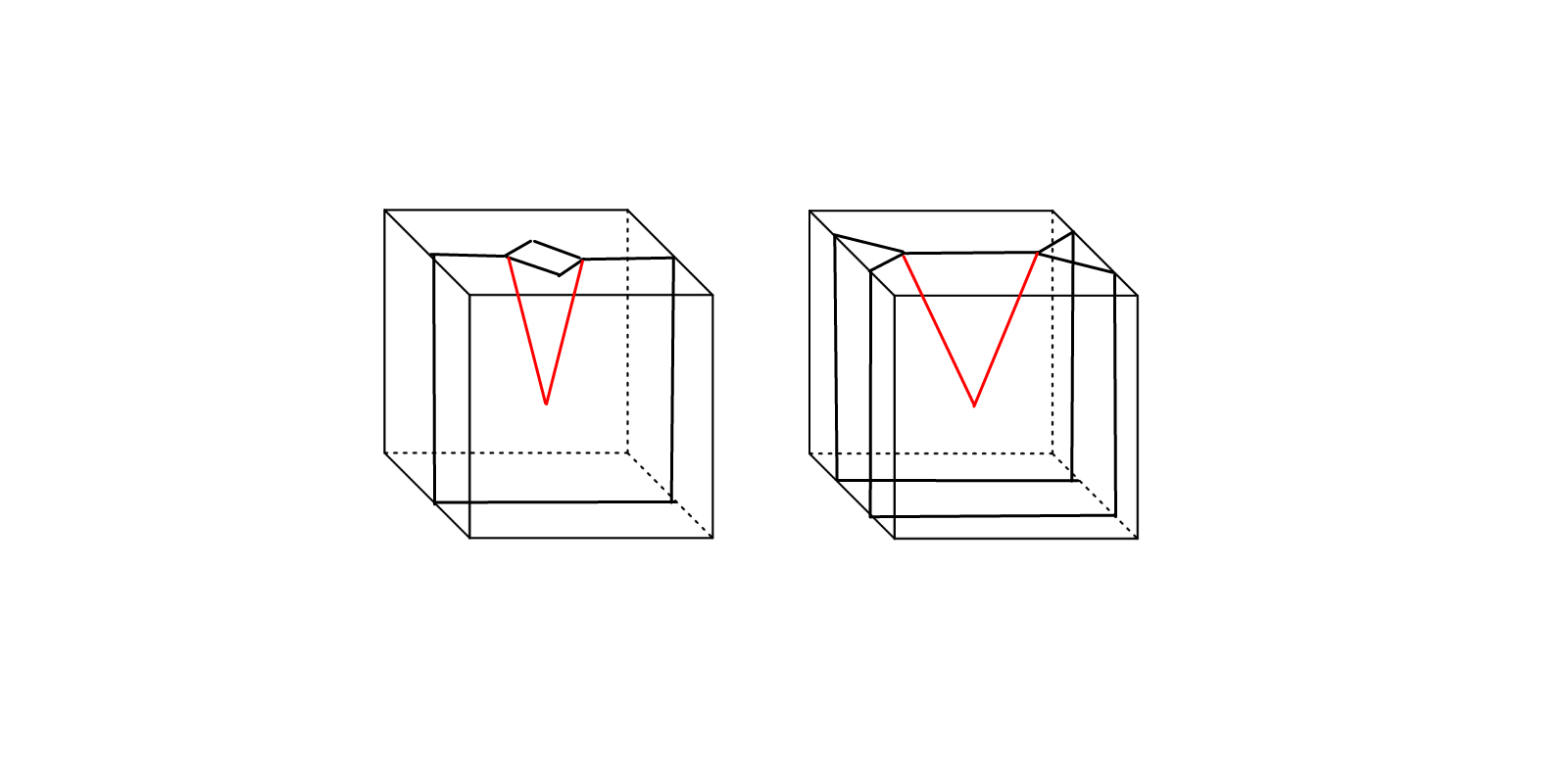}
\endgroup\end{restoretext}
and the cap singularities can be obtained by flipping these pictures vertically. 

Now assume $\sC^0$ has three composable $1$-generators $f,g,h$. As before the depth $i$ bracketing contains bracketing $(fg)$ and $(gh)$, but the depth $2$ bracketing now also contains $((fg)h)$ and $(f(gh))$, which abbreviate the $1$-generators
\begin{equation}
\bracket(\bracket(\abss{f} \whisker 1 1 \abss{g})\whisker 1 1 \bunit^1(h))
\end{equation}
and 
\begin{equation}
\bracket(\bunit^1(f) \whisker 1 1 \bracket(g \whisker 1 1 h))
\end{equation}
respectively. These generators are isomorphic by the following compositional $2$-shape in $\sC^2$
\begin{restoretext}
\begingroup\sbox0{\includegraphics{ANCimg/page1.png}}\includegraphics[clip,trim=0 {.25\ht0} 0 {.15\ht0} ,width=\textwidth]{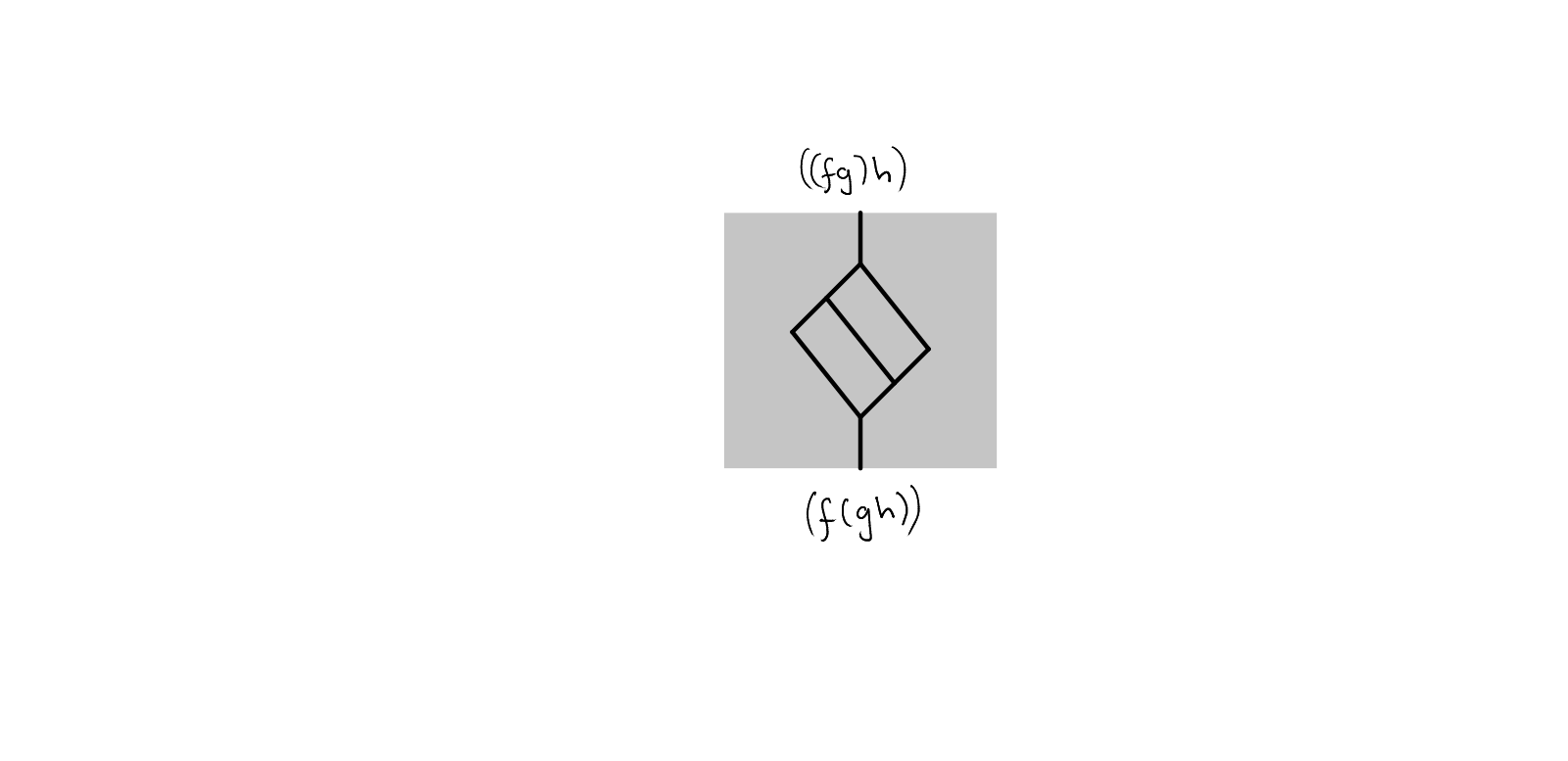}
\endgroup\end{restoretext}
Bracketing this morphism gives a $2$-generator in $\sC^3$ called the \textit{associator} (for $f,g,h$). 

This plays the role of the usual weak associator in higher category theory. In particular we will now claim that it satisfies the pentagon identity, and it does so ``weakly" giving rise to a cell called the \textit{pentagonator}.

Assume composable $1$-generators $f,g,h,k$ in $\sC^0$. Then, using bracketing of depth $3$ we find the following $3$-morphism in $\sC^3$
\begin{restoretext}
\begingroup\sbox0{\includegraphics{ANCimg/page1.png}}\includegraphics[clip,trim=0 {.05\ht0} 0 {.0\ht0} ,width=\textwidth]{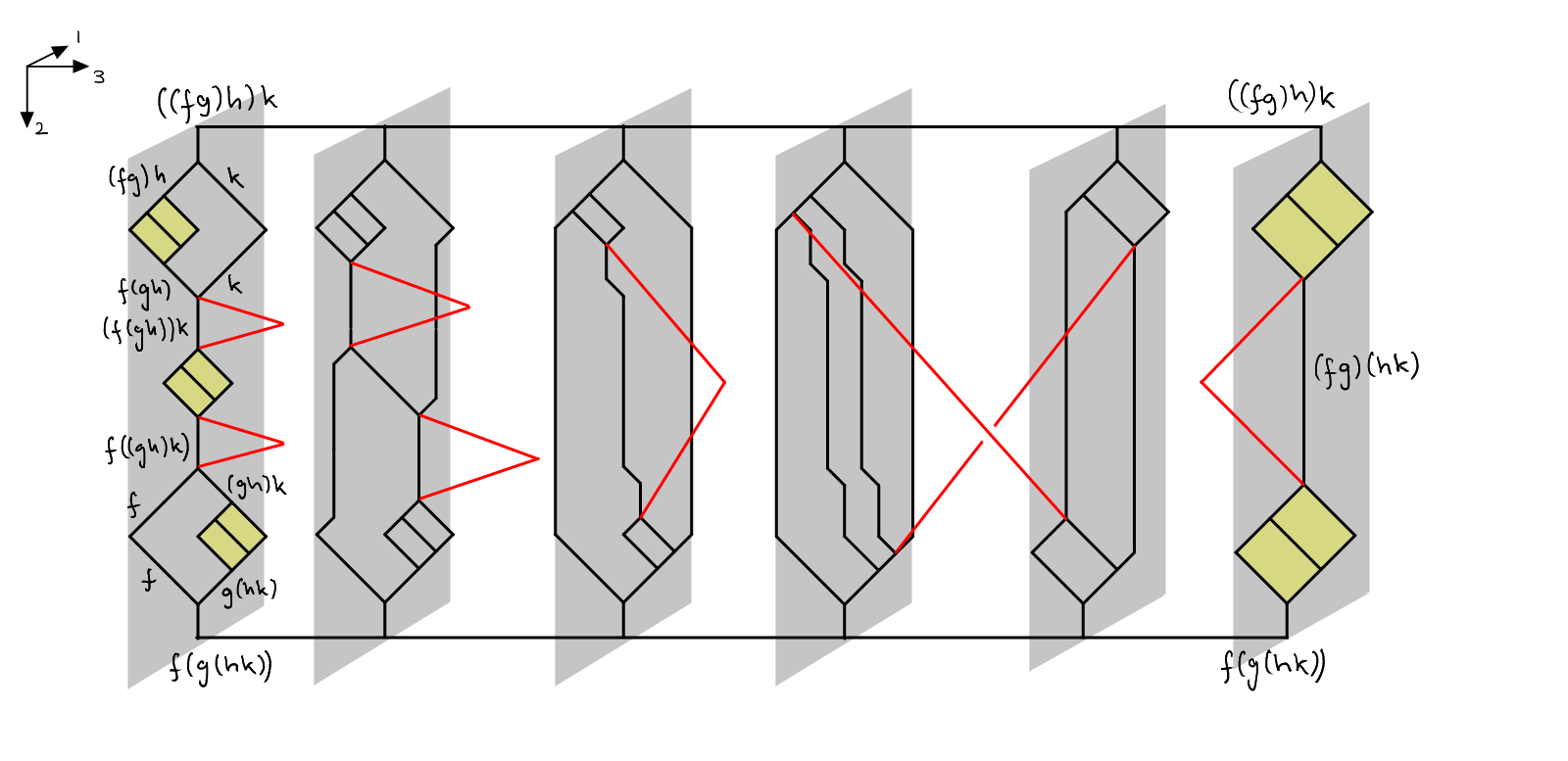}
\endgroup\end{restoretext}
The source of this $3$-morphisms is a composition of three (unbracketed) associators which transform the generator $(((fg)h)k)$ into the generator $(f(g(hk)))$ (note that outer brackets have been omitted in the picture, and that associators are highlighted by \cyellow{} coloring). The target is a composition of two (unbracketed) associators which transform the generator $(((fg)h)k)$ into the generator $(f(g(hk)))$. The lines in \cred{} highlight cup and cap singularities, as well as an interchanger, which are part of the above compositional $3$-shape. The pentagonator is the bracketing of the above morphism (after adding bracketing witnesses for associators on the left and on the right), and thus lives in $\sC^4$ as a generator.

\begin{rmk}[Coherences and homotopies] It is a crucial observation that in the above depiction of the unbracketed pentagonator one of the steps involves and interchange. This shows that in general, coherence data of higher categories requires the study of homotopies. 
\end{rmk}

Finally, we will also briefly discuss weakness of identities in \free{} weak $n$-categories. Weak identities are generators obtained by bracketing identity morphisms: that is, assume $A$ is a $0$-generator of $\sC^0$ then $\sC^1$ will contain the ``weak identity" $\bracket(\Id_A) \in \sC^1_1$ (which usually will be denoted by $\id_A$) together with a bracketing witness which has type
\begin{restoretext}
\begingroup\sbox0{\includegraphics{ANCimg/page1.png}}\includegraphics[clip,trim=0 {.25\ht0} 0 {.25\ht0} ,width=\textwidth]{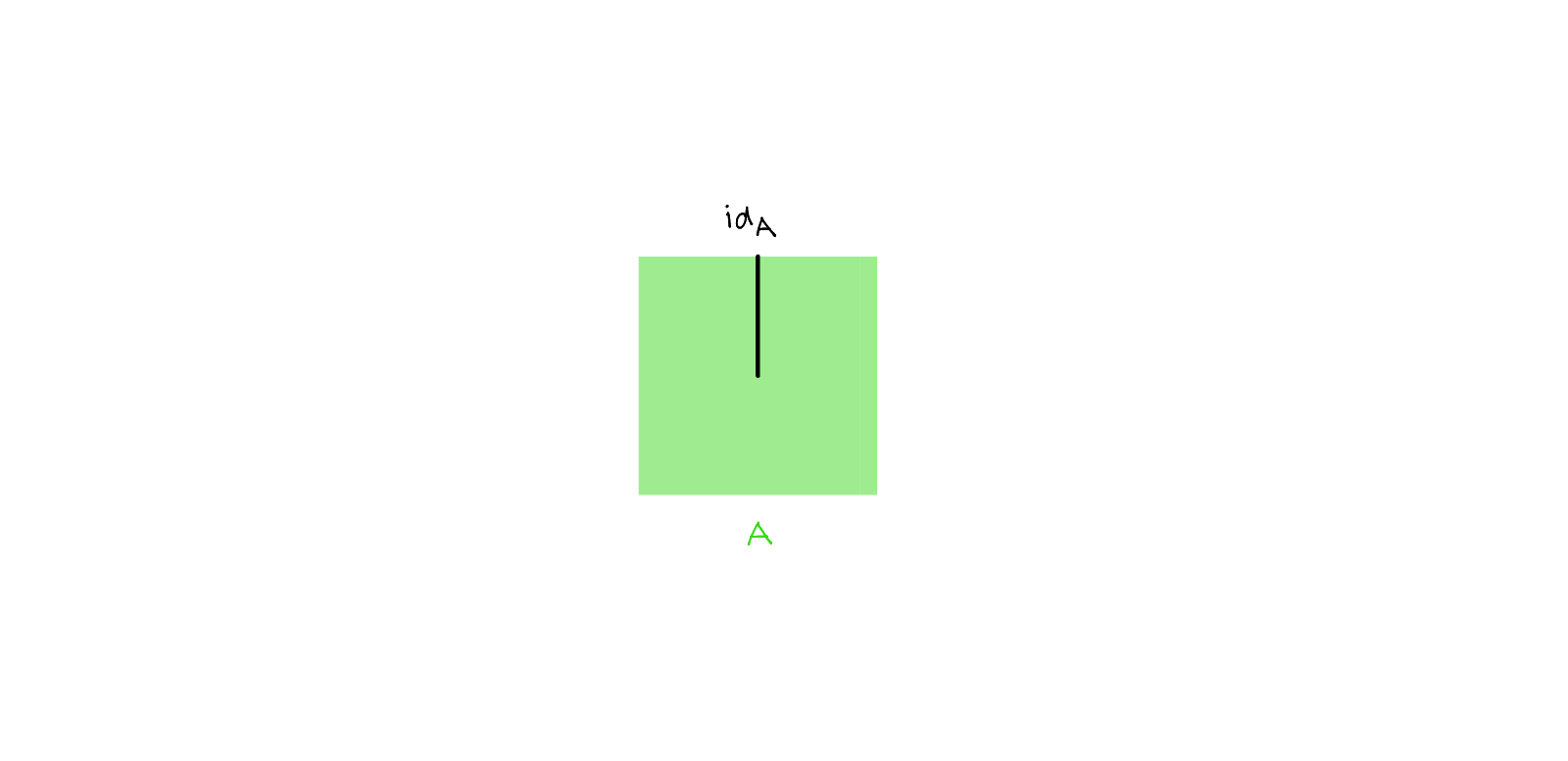}
\endgroup\end{restoretext}
This, together with witnesses of  $\whisker 1 1$ composition which were discussed above, allows us to form the following composition $2$-shape in $\sC^2$
\begin{restoretext}
\begingroup\sbox0{\includegraphics{ANCimg/page1.png}}\includegraphics[clip,trim=0 {.25\ht0} 0 {.2\ht0} ,width=\textwidth]{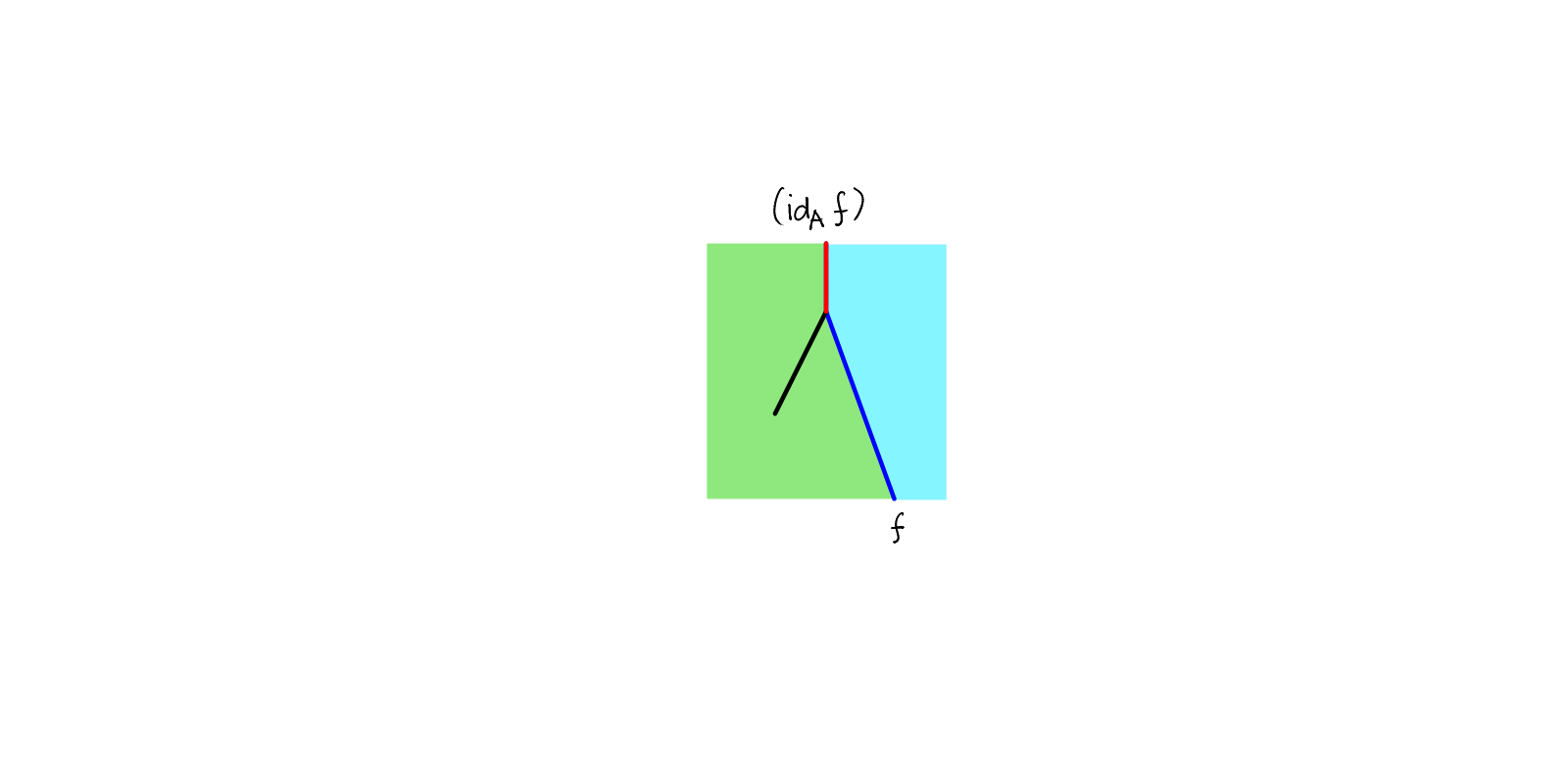}
\endgroup\end{restoretext}
This composition $2$-shape (composed of two invertible generators, and thus invertible itself) shows equivalence of $\id_A f$ with $f$, and in particular establishes that $\id_A$ acts as a weak identity as we expected.

\section{Equivalence of associative and fully weak approach} \label{sec:weak_assoc_equiv}

\subsection{Weakification and forgetful functor} \label{ssec:weakification}

Note that presented weak $\infty$-categories form a category $\pwCat$ whose morphisms are level-wise maps of presentations of their depth $k$ bracketings, commuting with the inclusions of depth $k$ into depth $(k+1)$ bracketings.

\begin{defn}[Weakification functor] There is a functor, called the \textit{weakification} functor, 
\begin{equation}
\Omega : \pCat_\infty \to \pwCat
\end{equation}
define by mapping $\sC \in \pCat$ to the chain
\begin{equation}
\kiC{(0)}(\sC) \into \kiC{(1)}(\sC) \into \kiC{(2)}(\sC) \into ...
\end{equation}
where $\kiC{(0)}(\sC) = \sC$ and $\kiC{(k+1)}(\sC) = \kiC{}(\kiC{(k)}(\sC))$. This extends naturally to maps of presentations.
\end{defn}

We further define

\begin{defn}[Forgetful functor] There is a \textit{forgetful} functor
\begin{equation}
U : \pwCat \to \pCat_\infty
\end{equation}
defined to map $\sC$ to $\sC$ (but forgetting the ``depth $k$" structure of a presented weak $n$-category).
\end{defn}

\subsection{Discussion}

Based on the above discussion, by the forgetful functor every \free{} weak category is a \free{} associative category. Conversely, from every associative category one can generate a \free{} weak category by weakification, by freely adjoining bracketings and bracketing witnesses and passing to the colimit. These processes are not inverse on the nose, and we won't show they are inverse ``up to equivalence". However, we remark that adjoining bracketings and bracketing witnesses, whether we start in the weak or associative case, leaves us with morally the same category since we are only adding more ``composition candidates".

\glsaddall
\printglossary[type=termsandsymbols,style=termdescsymb]  %
\printglossary[type=variables,style=vardesc]

\addcontentsline{toc}{chapter}{Bibliography}
\bibliography{refs}        %
\bibliographystyle{plain}  %

\end{document}